%% file: twistor.tex
\documentclass[a4,10pt]{article}
\usepackage{amsmath,amscd,amssymb}

\input{ntn}
\input{classdef}

\begin{document}

\title{Asymptotic behaviour of 
tame harmonic bundles and \\
an application to pure twistor $D$-modules
}
\author{Takuro Mochizuki}
\date{}

\maketitle

\begin{abstract}
\input{abstract}
\end{abstract}

\section{Introduction}

\input{103.2}

\input{103}

\input{103.1}

\subsection{The outline of the paper}

\input{105}
\input{105.1}

\input{105.2}

\input{105.3}

\subsection{Some remarks}
\input{d9}

\subsection{Acknowledgement}
\input{a89}

\part{Preliminary}
\label{part;b11.17.100}

\section{Preliminary}
\label{section;a11.9.50}

\subsection{Notation}
\label{subsection;a11.9.51}

\input{1}

\input{a87}

\input{b21}

\subsection{Prolongation by an increasing order}
\label{subsection;a11.9.52}

\input{b20}

\subsection{A preliminary for $\mu_c$-equivariant bundle}
\label{subsection;a11.9.53}
\input{a50.2}

\subsection{Some very elementary preliminary for convexity}
\label{subsection;a11.9.54}

\input{b15}

\subsection{An elementary remark on some distributions}
\label{subsection;a11.9.55}

\input{a34.4}

\subsection{Some preliminary from elementary linear algebra}
\label{subsection;a11.9.56}

\input{b22}

\input{c11.1}

\input{d4}

\input{c12.1}

\subsection{Preliminary from complex differential Geometry}

\label{subsection;a11.9.57}

\input{15.1}

\input{a49.1}

\input{a49.2}

\subsection{Preliminary from functional analysis}
\label{subsection;a11.9.58}
\input{15.2}

\input{a36.2}

\subsection{An estimate of the norms of Higgs field and the conjugate}
\label{subsection;a11.9.59}

\input{b25.1}

\subsection{Convergency of the sequence of harmonic bundles}
\label{subsection;a11.9.60}

\input{a36.1}

\subsection{Higgs field and twisted map}
\label{subsection;04.1.27.60}

\input{d1}

\section{Preliminary for mixed twistor structure}
\label{section;a11.9.61}

\subsection{$\proj^1$-holomorphic vector bundle over $X\times\proj^1$}
\label{subsection;a11.11.50}

\input{b30}

\input{b30.1}

\subsection{Equivariant $\proj^1$-holomorphic bundle over $X\times\proj^1$}
\label{subsection;a11.11.51}
\input{b6}

\input{b6.1}

\input{b6.2}

\subsection{Tate objects and $\nbigo(p,q)$}
\label{subsection;a11.11.52}
\input{b5}

\input{a87.1}

\input{a70}

\subsection{Equivalence of some categories}
\label{subsection;a11.11.53}

\input{b6.3}

\input{b7}

\subsection{Variation of $\proj^1$-holomorphic bundles}
\label{subsection;a11.11.54}

\input{b30.2}

\input{b30.3}

\input{b30.4}

\input{b6.4}

\subsection{The twistor nilpotent orbit}
\label{subsection;a11.11.55}

\input{a71}

\input{b10.3}

\input{b6.5}

\subsection{Split polarized mixed twistor structure and the nilpotent orbit}
\label{subsection;a11.11.56}

\input{a71.1}

\input{a48.2}

\input{a48.3}

\input{b6.6}

\subsection{The induced tuple on the divisor}
\label{subsection;a11.11.57}

\input{a48.4}

\input{a48.5}

\input{b8.1}

\input{b9}

\subsection{Translation of some results due to Kashiwara, Kawai and Saito}
\label{subsection;a11.11.58}

\input{a33}

\input{b10}

\input{b10.1}

\input{a33.1}

\input{b10.4}

\subsection{$\nbigr$-triple and twistor structure}

\label{subsection;a11.11.59}
\input{b11}

\input{a81}

\input{a81.1}

\section{Preliminary for filtrations}
\label{section;d11.12.1}

\subsection{Filtrations and decompositions on a vector space}

\label{subsection;d11.12.2}
\input{a50}

\subsection{Filtrations and decompositions on a vector bundle}
\label{subsection;d11.12.3}

\input{a50.1}

\subsection{Compatibility of the filtrations and nilpotent maps}
\label{subsection;d11.12.4}

\input{a51}
\input{a51.1}

\input{a51.2}

\subsection{Extension of splittings}

\label{subsection;d11.12.5}
\input{a50.3}
\input{a51.3}

\input{a51.4}

\subsection{Compatibility of the filtrations and nilpotent maps
 on the divisors}
\label{subsection;d11.12.6}

\input{a52}

\section{Some lemmas for generically splitted case}
\label{section;d11.12.10}

\subsection{Filtrations}
\label{subsection;d11.12.11}
\input{23.1}

\input{a72}

\subsection{Compatibility of morphisms and filtrations}
\label{subsection;d11.12.12}

\input{a72.1}

\input{b12.1}

\section{Model bundle}
\label{section;a12.10.1}

\subsection{Basic example I}

\label{subsection;a11.17.55}

\input{a53.3}

\subsection{Basic example II}

\label{subsection;a11.17.56}

\input{a53.2}

\input{a88}


\part{Prolongation of deformed holomorphic bundles}
\label{part;a11.17.101}

\section{Harmonic bundles on $\Delta^{\ast}$}
\label{section;b11.14.1}

\subsection{Simpson's Main estimate}

\label{subsection;b11.9.2}

\input{6}

\subsection{The KMS-structure of tame harmonic bundles
on a punctured disc}

\label{subsection;a11.14.2}

\input{2}

\input{2.1}

\subsection{Basic comparison due to Simpson}

\label{subsection;10.11.35}

\input{3}

\subsection{Multi-valued flat sections}
\label{subsection;a11.14.3}

\input{4}

\input{4.1}

\input{4.2}

\subsection{Family of multi-valued sections}

\label{subsection;a11.14.4}

\input{7}

\subsection{Asymptotic orthogonality}

\label{subsection;a11.14.5}

\input{5}

\subsection{Maximum principle for the distance of the harmonic metrics}
\label{subsection;04.1.27.61}

\input{d2}

\section{Harmonic bundles on $\Delta^{\ast\,l}\times\Delta^{n-l}$}

\label{section;b11.16.1}

\subsection{Preliminary}
\label{subsection;b11.16.2}
\input{10}
\input{10.3}

\input{d3}

\subsection{Simpson's Main estimate in higher dimensional case}

\label{subsection;b11.9.3}

\input{10.1}

\input{10.2}

\subsection{Prolongation in the case that $\lambda$ is generic}
\label{subsubsection;b11.16.5}
\input{11}

\input{a73}

\subsection{Extension of holomorphic sections on a hyperplane}

\label{subsubsection;b11.16.6}

\input{12}

\input{12.1}

\input{a73.1}

\subsection{Preliminary prolongation of $\nbigelambda$ (Special case)}

\label{subsection;b11.16.7}

\input{12.2}

\input{12.3}

\subsection{Prolongation of $\nbigelambda$
 and the compatibility of the parabolic filtrations}

\label{subsection;b11.16.8}

\input{14}

\input{14.1}

\subsection{Prolongation of
$\nbige_{|\Delta(\lambda_0,\epsilon_0)\times (X-D)}$}

\label{subsection;b11.11.30}

\input{15}

\input{15.3}

\input{15.4}

\subsection{KMS-structure of $\prolongg{\vecb}{\nbige}$}
\label{subsection;b11.16.10}

\input{a55}

\input{15.8}

\input{15.5}

\input{15.6}

\input{15.7}

\subsection{The induced vector bundle}

\label{subsection;b11.16.11}

\input{a57}

\section{The $KMS$-structure on
the spaces of the multi-valued flat sections}
\label{section;b11.17.1}

\subsection{The filtration $\lefttop{i}\nbigf$}

\label{subsection;b11.17.2}

\input{16}

\input{16.1}

\input{16.2}

\subsection{The compatibility of the filtrations
$\lefttop{i}\nbigf$ $(i=1,\ldots,l)$}

\label{subsection;b11.17.3}

\input{a55.1}

\input{16.3}

\input{16.4}

\input{a57.1}

\subsection{The induced objects}

\label{subsection;b11.17.4}

\input{16.6}

\input{a57.2}

\section{The filtrations and the decompositions on $\nbige$ and
 $\prolongg{\vecc}{\nbige}$}

\label{section;a11.17.20}

\subsection{The filtrations and the decompositions on $\nbigelambda$
and $\prolongg{\vecc}{\nbigelambda}$}

\label{subsection;a11.17.21}

\input{16.5}

\subsection{The decomposition $\EEzero$ and the filtration
  $\Fzero$ on $\nbige$ for $\lambda_0\neq 0$}

\label{subsection;a11.17.22}

\input{17}

\input{17.1}

\subsection{The morphisms between $\lefttop{\lbar}\nbigg_{\vecu}$ and
  $\lefttop{\lbar}\nbigg_{\vecu}\nbigh$}
\label{subsection;a11.17.23}

\input{17.2}

\part{Limiting mixed twistor theorem and some consequence}
\label{part;a11.17.102}

\section{The induced vector bundle}

\label{section;a11.17.50}

\subsection{The variation of pure twistor structures}

\label{subsection;a11.17.51}

\input{21}

\subsection{The induced objects of the conjugate and the pairing}

\label{subsection;a11.17.52}

\input{a53}

\subsection{The induced vector bundles over $\proj^1$}

\label{subsection;a11.17.53}

\input{a53.1}

\input{20.5}

\input{a57.3}

\subsection{$\Gr^W_hS^{\can}_{\vecu}(E)$ and $\Gr^W_h S_{\vecu}(E,P)$}

\label{subsection;a11.17.54}

\input{a58.2}

\input{a57.4}

\section{Limiting mixed twistor theorem}

\label{section;a11.17.60}

\subsection{Limiting mixed twistor theorem in the case of curves}

\label{subsection;c11.14.5}

\input{a88.1}

\input{a88.2}

\input{a88.3}

\input{a58}

\input{a58.1}

\subsection{Limiting mixed twistor theorem in the higher dimensional case}

\label{subsection;a11.17.61}

\input{21.4}

\input{a59.2}
\input{a59.1}

\input{a59.3}

\input{21.5}

\subsection{Some consequences}

\label{subsection;a11.17.62}

\input{22}

\input{22.1}

\input{a33.2}

\input{a84.1}

\section{Norm estimate}

\label{section;a11.17.70}

\subsection{Preliminary}

\label{subsection;a11.17.71}

\input{25}

\input{25.1}

\input{25.2}

\subsection{Preliminary norm estimate}

\label{subsection;a11.17.72}

\input{26}

\subsection{Norm estimate for holomorphic sections}

\label{subsection;a11.17.73}

\input{26.1}

\subsection{Norm estimate for flat sections}

\label{subsection;a11.17.74}

\input{a60}

\part{An application to the theory of pure twistor $D$-modules}

\label{part;a11.17.103}

\section{Nearby cycle functor for $\nbigr$-module}

\label{section;b11.18.1}

\subsection{The KMS structure of $\nbigr$-module}

\label{subsection;b11.18.2}

\input{a74}

\input{30}

\input{a34.1}

\subsection{Specialization of the pairing of Sabbah}

\label{subsection;b11.18.3}

\input{a74.1}

\input{30.1}

\input{a34.2}

\input{a35}

\section{Prolongation of $\nbigr$-module $\nbige$}

\label{section;b11.24.1}

\input{a74.2}

\subsection{Naive prolongment $\naiveprolong{\nbige}$ and the filtrations}

\label{subsection;b11.24.2}

\input{a75}

\input{a75.2}

\input{a75.5}

\subsection{Prolongment $\gbige$}
\label{subsection;a11.23.40}

\input{a75.1}

\input{a75.6}

\subsection{Comparison of $\lefttop{I}T^{(\lambda_0)}(\vecc,\vecd)$
 and $\lefttop{I}\tilde{T}^{(\lambda_0)}(\vecc,\vecd)$}

\label{subsection;b11.24.5}

\input{31.1}

\input{a75.3}

\input{a77.2}

\subsection{Relation of the filtrations of $\gbige$}

\label{subsection;b11.24.6}

\input{a76.1}

\input{a76}

\input{a76.2}
\input{31.5}

\subsection{The characterization of $\gbige$}
\label{subsection;b12.6.10}

\input{a79}

\section{The filtrations of $\gbige[\deldel_t]$}

\label{section;b11.24.10}

\input{a78}

\subsection{The filtration $\Uzero$}

\label{subsection;b11.24.11}

\input{32}

\subsection{Preliminary reductions and decompositions}

\label{subsection;b11.24.12}

\input{32.1}

\input{32.2}

\subsection{Primitive decomposition}
\label{subsection;b11.24.13}
\input{32.3}

\subsection{Strict $S$-decomposability along $t=0$}
\label{subsection;b12.6.15}

\input{b32.1}
\input{b32.2}

\subsection{The decomposition}
\label{subsection;b12.6.30}

\input{32.4}

\input{32.5}

\input{a33.3}

\input{a84}

\section{The weight filtration on $\psi_{t,u}\gbige$
and the induced $\nbigr$-triple}

\label{section;b12.2.0}

\input{a83.1}

\subsection{The weight filtration on $\lefttop{I}\nbigl$}
\label{subsection;b12.2.1}

\input{a80.1}

\input{a80.2}

\subsection{The filtration $\FFzero$}
\label{subsection;b12.2.2}

\input{a83.2}

\input{a80.3}

\input{a80.4}

\input{a80.5}

\subsection{Strict $S$-decomposability}

\label{subsection;b12.2.3}

\input{a80.6}

\input{a80.7}

\input{a81.2}

\section{The sesqui-linear pairings}
\label{section;b12.6.100}

\input{a85.1}

\subsection{The sesqui-linear pairing on $\gbige$}

\label{subsection;b11.24.7}

\input{a34}

\subsection{The sesqui-linear pairing on
$\lefttop{\nbar}\overline{\nbigg}_{\vecu}(\nbige)$}
\label{subsection;b12.6.200}

\input{a85}

\input{a83}

\input{a84.2}

\subsection{The induced sesqui-linear pairings}

\label{subsection;b12.6.210}

\input{a85.2}

\subsection{A comparison of the induced objects}

\label{subsection;b11.24.8}

\input{a85.3}

\input{a85.4}

\subsection{The specialization of the pairings}

\label{subsection;b12.2.4}

\input{a34.6}

\section{A prolongment as a regular pure twistor $D$-module}

\label{section;b11.24.30}

\subsection{Correspondence}
\label{subsection;a11.23.50}
\input{a77}

\subsection{The tameness of the corresponding harmonic bundle}
\label{subsection;04.1.27.35}
\input{d5}

\subsection{The existence of
 the prolongment as regular pure twistor $D$-module}
\label{subsection;a11.23.51}

\input{a35.1}

\input{a77.1}

\subsection{The uniqueness of the prolongment as 
   regular pure twistor $D$-module}

\label{subsection;a11.23.52}

\input{a35.2}

\subsection{The pure imaginary case}
\label{subsection;04.2.5.1}

\input{d7}

\subsection{The conjectures of Kashiwara and Sabbah}
\label{subsection;04.2.5.2}

\input{d8}

\part{Appendix}

\label{part;a11.17.104}

\section{Pure twistor $D$-modules and polarization}

\subsection{Pure twistor $D$-modules and polarization}

\input{a62}

\subsection{Pure imaginary pure twistor $D$-modules}

\label{subsection;04.2.5.100}
\input{d6}

\subsection{The pure twistor $D$-modules on a smooth curve (correspondence)}

\input{a62.1}

\section{The decomposition theorem for pure twistor $D$-modules}

\input{a62.2}

\subsection{Preliminary}
\input{a61}

\subsection{Quasi isomorphisms of complexes of sheaves}
\input{a61.1}

\subsection{Globalization of isomorphisms}
\input{a61.2}

\subsection{Family of isomorphisms}

\input{a86}

\input{a86.1}

\input{harmref.tex}
\noindent
{\it Address\\
Department of Mathematics,
Osaka City University,
Sugimoto, Sumiyoshi-ku,
Osaka 558-8585, Japan.\\
takuro@sci.osaka-cu.ac.jp\\
}

\end{document}

%% file: ntn.tex
\newcommand{\nbiga}{\mathcal{A}}
\newcommand{\nbigb}{\mathcal{B}}
\newcommand{\nbigc}{\mathcal{C}}
\newcommand{\nbigd}{\mathcal{D}}
\newcommand{\nbige}{\mathcal{E}}
\newcommand{\nbigf}{\mathcal{F}}
\newcommand{\nbigg}{\mathcal{G}}
\newcommand{\nbigh}{\mathcal{H}}
\newcommand{\nbigi}{\mathcal{I}}

\newcommand{\nbigk}{\mathcal{K}}
\newcommand{\nbigl}{\mathcal{L}}
\newcommand{\nbigm}{\mathcal{M}}
\newcommand{\nbign}{\mathcal{N}}
\newcommand{\nbigo}{\mathcal{O}}

\newcommand{\nbigq}{\mathcal{Q}}
\newcommand{\nbigr}{\mathcal{R}}
\newcommand{\nbigs}{\mathcal{S}}
\newcommand{\nbigt}{\mathcal{T}}
\newcommand{\nbigu}{\mathcal{U}}
\newcommand{\nbigv}{\mathcal{V}}
\newcommand{\nbigw}{\mathcal{W}}
\newcommand{\nbigx}{\mathcal{X}}
\newcommand{\nbigy}{\mathcal{Y}}

\newcommand{\proj}{\mathbb{P}}
\newcommand{\seisuu}{\mathbb{Z}}
\newcommand{\rnum}{{\boldsymbol Q}}

\newcommand{\cnum}{{\boldsymbol C}}
\newcommand{\real}{{\boldsymbol R}}
\newcommand{\hyperh}{\mathbb{H}}

\newcommand{\Tate}{\mathbb{T}}

\newcommand{\DD}{\mathbb{D}}
\newcommand{\EE}{\mathbb{E}}
\newcommand{\FF}{\mathbb{F}}
\newcommand{\GG}{\mathbb{G}}

\newcommand{\gbiga}{\mathfrak A}
\newcommand{\gbigb}{\mathfrak B}
\newcommand{\gbigc}{\mathfrak C}
\newcommand{\gbigd}{\mathfrak D}
\newcommand{\gbige}{\mathfrak E}

\newcommand{\gbigh}{\mathfrak H}

\newcommand{\gbigm}{\mathfrak M}

\newcommand{\gbigs}{\mathfrak S}
\newcommand{\gbigt}{\mathfrak T}

\newcommand{\gminib}{\mathfrak b}

\newcommand{\gminie}{\mathfrak e}

\newcommand{\gminik}{\mathfrak k}

\newcommand{\gminim}{\mathfrak m}

\newcommand{\gminip}{\mathfrak p}
\newcommand{\gminiq}{\mathfrak q}

\newcommand{\gminis}{\mathfrak s}
\newcommand{\gminit}{\mathfrak t}

\newcommand{\veceta}{{\boldsymbol \eta}}
\newcommand{\vecrho}{{\boldsymbol \rho}}
\newcommand{\vece}{{\boldsymbol e}}

\newcommand{\vecv}{{\boldsymbol v}}
\newcommand{\vecu}{{\boldsymbol u}}
\newcommand{\vecw}{{\boldsymbol w}}
\newcommand{\vecgamma}{{\boldsymbol \gamma}}
\newcommand{\vecl}{{\boldsymbol l}}

\newcommand{\vecalpha}{{\boldsymbol \alpha}}
\newcommand{\veca}{{\boldsymbol a}}
\newcommand{\vecb}{{\boldsymbol b}}
\newcommand{\vecbeta}{{\boldsymbol \beta}}
\newcommand{\vecdelta}{{\boldsymbol \delta}}
\newcommand{\vecs}{{\boldsymbol s}}
\newcommand{\vect}{{\boldsymbol t}}
\newcommand{\vecc}{{\boldsymbol c}}
\newcommand{\vecd}{{\boldsymbol d}}
\newcommand{\vech}{{\boldsymbol h}}
\newcommand{\veck}{{\boldsymbol k}}
\newcommand{\vecm}{{\boldsymbol m}}
\newcommand{\vecM}{{\boldsymbol M}}
\newcommand{\vecN}{{\boldsymbol N}}
\newcommand{\vecI}{{\boldsymbol I}}
\newcommand{\vecomega}{{\boldsymbol \omega}}

\newcommand{\vecf}{{\boldsymbol f}}

\newcommand{\vecF}{{\boldsymbol F}}
\newcommand{\vecE}{{\boldsymbol E}}
\newcommand{\vecW}{{\boldsymbol W}}
\newcommand{\vecn}{{\boldsymbol n}}
\newcommand{\vecp}{{\boldsymbol p}}
\newcommand{\veczeta}{{\boldsymbol \zeta}}
\newcommand{\vecC}{{\boldsymbol C}}
\newcommand{\vecU}{{\boldsymbol U}}

\newcommand{\llarr}{\longleftarrow}

\newcommand{\lrarr}{\longrightarrow}

\newcommand{\pf}{{\bf Proof}\hspace{.1in}}
\newcommand{\qed}{\mbox{\rule{1.2mm}{3mm}}}

\def\End{\mathop{\rm End}\nolimits}

\def\Cok{\mathop{\rm Cok}\nolimits}
\def\cok{\mathop{\rm Cok}\nolimits}
\def\Image{\mathop{\rm Im}\nolimits}
\def\Realpart{\mathop{\rm Re}\nolimits}

\def\Re{\mathop{\rm Re}\nolimits}

\def\Gr{\mathop{\rm Gr}\nolimits}

\def\rank{\mathop{\rm rank}\nolimits}

\def\Ker{\mathop{\rm Ker}\nolimits}
\def\ker{\mathop{\rm Ker}\nolimits}
\def\modulo{\mathop{\rm modulo}\nolimits}

\def\Gr{\mathop{\rm Gr}\nolimits}
\def\Sym{\mathop{\rm Sym}\nolimits}

\def\ad{\mathop{\rm ad}\nolimits}
\def\Res{\mathop{\rm Res}\nolimits}
\def\Rep{\mathop{\rm Rep}\nolimits}
\def\ord{\mathop{\rm ord}\nolimits}

\def\Ric{\mathop{\rm Ric}\nolimits}
\def\tr{\mathop{\rm tr}\nolimits}
\def\vol{\mathop{\rm dvol}\nolimits}
\def\dvol{\mathop{\rm dvol}\nolimits}
\def\Diff{\mathop{\rm Diff}\nolimits}
\def\vecdeg{\mathop{\rm \bf deg}\nolimits}
\def\vecord{\mathop{\rm \bf ord}\nolimits}
\def\can{\mathop{\rm can}\nolimits}
\def\var{\mathop{\rm var}\nolimits}
\def\id{\mathop{\rm id}\nolimits}
\def\Pat{\mathop{\rm Pat}\nolimits}
\def\Prim{\mathop{\rm Prim}\nolimits}
\def\Red{\mathop{\rm Red}\nolimits}
\def\gcd{\mathop{\rm g.c.d.}\nolimits}

\def\filt{\mathop{\rm Filt}\nolimits}
\def\equi{\mathop{\rm Equi}\nolimits}
\def\Supp{\mathop{\rm Supp}\nolimits}

\def\VPTgen{\mathop{\rm VPT_{gen}}\nolimits}
\def\MPT{\mathop{\rm MPT}\nolimits}
\def\RHD{\mathop{\rm RHD}\nolimits}

\newcommand{\del}{\partial}
\newcommand{\delbar}{\overline{\del}}

\newcommand{\Deltabar}{\overline{\Delta}}
\newcommand{\Deltabarast}{\Deltabar^{\ast}}

\newcommand{\nbar}{\underline{n}}

\newcommand{\jbar}{\underline{j}}
\newcommand{\mbar}{\underline{m}}
\newcommand{\kbar}{\underline{k}}

\newcommand{\mbariti}{\underline{m_1}}
\newcommand{\mminusitibar}{\underline{m-1}}
\newcommand{\lbariti}{\underline{l_1}}

\newcommand{\ibar}{\underline{i}}
\newcommand{\iitibar}{\underline{i+1}}
\newcommand{\nminusitibar}{\underline{n-1}}
\newcommand{\nminusnibar}{\underline{n-2}}

\newcommand{\lbar}{\underline{l}}
\newcommand{\pbar}{\underline{p}}
\newcommand{\itibar}{\underline{1}}
\newcommand{\nibar}{\underline{2}}

\newcommand{\alphasitabar}{\underline{\alpha}}

\newcommand{\shikaku}{\sharp}
\newcommand{\sankaku}{\triangle}
\newcommand{\twoprime}{\prime\prime}

\newcommand{\harmonicbundle}{(E,\delbar_E,\theta,h)}

\newcommand{\barz}{\bar{z}}
\newcommand{\zbar}{\barz}
\newcommand{\zetabar}{\bar{\zeta}}
\newcommand{\baralpha}{\bar{\alpha}}
\newcommand{\alphabar}{\baralpha}
\newcommand{\barlambda}{\bar{\lambda}}
\newcommand{\lambdabar}{\barlambda}

\newcommand{\tbar}{\bar{t}}
\newcommand{\xbar}{\bar{x}}
\newcommand{\sbar}{\bar{s}}

\newcommand{\modeldeform}{{\mathcal Mod}}

\newcommand{\Poin}{{\bf p}}
\newcommand{\poin}{\Poin}
\newcommand{\prolong}[1]{{}^{\diamond}{#1}}

\newcommand{\prolongg}[2]{{}_{#1}{#2}}

\newcommand{\resddlambda}{\Res(\DD^{\lambda})}

\newcommand{\DDlambda}{\DD^{\lambda}}

\newcommand{\laplacian}{\Delta''}

\newcommand{\doublelangle}{\langle\langle}
\newcommand{\doublerangle}{\rangle\rangle}

\newcommand{\nbigelambda}{\nbige^{\lambda}}
\newcommand{\nbigelambdazero}{\nbige^{\lambda_0}}

\newcommand{\nbigelor}{\nbige^{\lor}}

\newcommand{\blowup}{\widetilde}

\newcommand{\nbigxlambda}{\nbigx^{\lambda}}
\newcommand{\nbigxlambdazero}{\nbigx^{\lambda_0}}
\newcommand{\nbigdlambda}{\nbigd^{\lambda}}

\newcommand{\KMS}{{\mathcal{KMS}}}
\newcommand{\KMSoverline}{\overline{\mathcal{KMS}}}
\newcommand{\KMSE}[1]{\KMS(\nbige^{#1})}
\newcommand{\KMSEoverline}[1]{\overline{\KMS}(\nbige^{#1})}

\newcommand{\KMSEprolongg}[2]{\KMS(\prolongg{#1}{\nbige^{#2}})}

\newcommand{\Par}{{\mathcal Par}}
\newcommand{\Sp}{{\mathcal Sp}}

\newcommand{\ParE}[1]{{\mathcal Par}(\nbige^{#1})}
\newcommand{\SpE}[1]{{\mathcal Sp}(\nbige^{#1})}

\newcommand{\ParEprolongg}[2]{{\mathcal Par}(\prolongg{#1}{\nbige}^{#2})}
\newcommand{\SpEprolongg}[2]{{\mathcal Sp}(\prolongg{#1}{\nbige}^{#2})}

\newcommand{\kmsmap}{\gminik}
\newcommand{\paramap}{\gminip}
\newcommand{\eigenmap}{\gminie}
\newcommand{\multiplicity}{\gminim}

\newcommand{\twistmap}{\gminit}

\newcommand{\lefttop}[1]{{}^{#1}}
\newcommand{\leftbottom}[1]{{}_{#1}}

\def\Equi{\mathop{\rm Equi}\nolimits}
\def\Bifilt{\mathop{\rm Bifilt}\nolimits}

\newcommand{\FEzero}{(F^{(\lambda_0)},\EE^{(\lambda_0)})}

\newcommand{\Fzero}{F^{(\lambda_0)}}
\newcommand{\EEzero}{\EE^{(\lambda_0)}}
\newcommand{\FFzero}{\FF^{(\lambda_0)}}
\newcommand{\GGzero}{\GG^{(\lambda_0)}}
\newcommand{\Vzero}{V^{(\lambda_0)}}

\newcommand{\psizero}{\psi^{(\lambda_0)}}
\newcommand{\tildepsizero}{\widetilde{\psi}^{(\lambda_0)}}
\newcommand{\tildepsi}{\widetilde{\psi}}
\newcommand{\psitilde}{\tildepsi}

\newcommand{\nbigfzero}{\nbigf^{(\lambda_0)}}

\newcommand{\Ulambdazero}{U^{(\lambda_0)}}
\newcommand{\Uzero}{U^{(\lambda_0)}}

\newcommand{\Psizero}{\Psi^{(\lambda_0)}}

\newcommand{\deldel}{\eth}
\newcommand{\deldelbar}{\overline{\deldel}}

\newcommand{\lamda}{\lambda}
\newcommand{\lambdazero}{\lambda_{0}}
\newcommand{\lamdazero}{\lamda_0}

\newcommand{\AAA}{{\boldsymbol A}}
\newcommand{\distribution}{\gbigd\gminib}

\newcommand{\naiveprolong}[1]{\lefttop{\square}{#1}}
\newcommand{\naiveprolongg}[2]{{}^{\square}_{#1}{#2}}

\newcommand{\rmoduleprolong}[1]{\mathfrak #1}

\newcommand{\vecnbign}{{\mathcal {\boldsymbol {\tilde{N}}}}}
\newcommand{\nbigvecn}{{\mathcal{\boldsymbol N}}}

\newcommand{\supp}{\gminis}

\newcommand{\kakkolambdazero}{^{(\lambda_0)}}

\newcommand{\tildenbigm}{\widetilde{\nbigm}}
\newcommand{\tildenbigt}{\widetilde{\nbigt}}
\newcommand{\tildegbige}{\widetilde{\gbige}}
\newcommand{\tildenbige}{\widetilde{\nbige}}

\newcommand{\dualhenomap}{\clubsuit}
\newcommand{\rtriplecat}{\nbigr-{\rm Triples}}
\newcommand{\closedopen}[2]{[#1,#2[}

\newcommand{\PH}{{\mathcal P}{\mathcal H}}

%% file: classdef.tex
\newcommand{\bdmath}{\begin{displaymath}}
\newcommand{\edmath}{\end{displaymath}}

\newcommand{\beqn}{\begin{equation}}
\newcommand{\eeqn}{\end{equation}}

\newcommand{\beqnarray}{\begin{eqnarray}}
\newcommand{\eeqnarray}{\end{eqnarray}}

\newcommand{\bitemize}{\begin{itemize}}
\newcommand{\eitemize}{\end{itemize}}

\newcommand{\benumerate}{\begin{enumerate}}
\newcommand{\eenumerate}{\end{enumerate}}

\newcommand{\bdescriprion}{\begin{description}}
\newcommand{\edescriprion}{\end{description}}

\newtheorem{thm}{Theorem}[section]
\newtheorem{cor}{Corollary}[section]
\newtheorem{conj}{Conjecture}[section]
\newtheorem{rem}{Remark}[section]
\newtheorem{lem}{Lemma}[section]
\newtheorem{prop}{Proposition}[section]
\newtheorem{df}{Definition}[section]

\newtheorem{condition}{Condition}[section]
\newtheorem{assumption}{Assumption}[section]
\newtheorem{principle}{Principle}[section]

\def\longuparrow{\vcenter{%
   \lineskip0pt\lineskiplimit0pt\baselineskip0pt\ialign{%
   \hfil{##}\hfil\crcr\hbox to 0pt{\hss$\uparrow$\hss}\cr%
   \hbox to 0pt{\hss\vrule width.4pt depth 0pt height 1em\hss}\cr}}}
\def\longdownarrow{\vcenter{%
   \lineskip0pt\lineskiplimit0pt\baselineskip0pt\ialign{\hfil{##}\hfil%
   \crcr\hbox to 0pt{\hss\vrule width.4pt depth 0pt height 1em\hss}\cr%
   \hbox to 0pt{\hss$\downarrow$\hss}\cr}}}

%% file: abstract.tex
We study the asymptotic behaviour of tame harmonic bundles.
First of all, we prove a local freeness of the prolongation
by an increasing order.
Then we obtain the polarized mixed twistor structure.
As one of the applications, we obtain the norm estimate of holomorphic
or flat sections by weight filtrations of the monodromies.

As other application,
we establish the correspondence
of semisimple regular holonomic $D$-modules
and polarizable pure imaginary pure twistor $D$-modules
through a tame pure imaginary harmonic bundles,
which is a conjecture of Sabbah.
Then the regular holonomic version of Kashiwara's conjecture
follows from the results of Sabbah and us.

\vspace{.1in}
\noindent
Keywords: 
Higgs fields,
harmonic bundle,
variation of Hodge structure,
mixed twistor structure,
$D$-module.

\noindent
MSC: 14C30, 53C43, 32S40.

%% file: 103.2.tex
\subsection{Simpson's Meta-Theorem}

The guiding principle of our study is
the following, which we call Simpson's Meta-Theorem:

\begin{principle}
The theory of Hodge structure should be generalized 
to the theory of twistor structures.
\hfill\qed
\end{principle}

In {\rm\cite{s3}},
Simpson stated the above principle as follows:
\begin{quote}
\noindent
{\bf Meta-Theorem}
{\em
If the words ``mixed Hodge structure''
(resp. ``variation of mixed Hodge structure'')
are replaced by the words ``mixed twistor structure''
(resp. ``variation of mixed twistor structure'')
in the hypotheses and conclusions of any theorem
in Hodge theory,
then one obtains a true statement.
The proof of the new statement will be analogous to
the proof of the old statement.}
\end{quote}
We regard it as a kind of principle.
As for the study of variation of pure twistor structures
(harmonic bundle),
it may occur that
the proof of new statement is not analogous to
the proof of the old statement,
in our current understanding.

%% file: 103.tex

\subsection{The purposes in this paper}
We have two main purposes in this paper.

\begin{enumerate}
\item \label{number;11.6.1}
 In the previous paper \cite{mochi},
 we discussed the behaviour of tame harmonic bundle
 imposed the nilpotentness and the trivial parabolic structure
 conditions.
 We would like to remove the assumption.
 We also improve the argument.
 In particular,
 we use the reduction to Hodge theory more efficiently.
\item \label{number;11.6.2}
 We would like to apply the study on the behaviour of tame harmonic
 bundle to the theory of pure twistor $D$-module,
 introduced by  Sabbah.
\end{enumerate}

\subsection{On the purpose \ref{number;11.6.1}}

Our principle in the study of tame harmonic bundle
is as follows,
which is a `corollary' of Simpson's Meta-Theorem:
\begin{principle}
 The asymptotic behaviour of tame harmonic bundle
 should be similar to the asymptotic behaviour of
 variation of polarized Hodge structures.
\hfill\qed
\end{principle}
Although our goal is to show the theorems
known for complex variation of polarized Hodge structures,
we do not follow closely \cite{sch}, \cite{cks1} and \cite{k}.
Instead
we follow the more differential geometric method pioneered by Simpson.
We refer the following two difficulty
to apply the classical method in the Hodge theory directly.

\begin{description}
\item[(a)]
 The nilpotent orbit theorem for harmonic bundle
 is not known.
\item[(b)]
 In the case of harmonic bundles,
 we have non-trivial eigenvalues of the residues
 and non-trivial parabolic structures.
\end{description}

\subsubsection{The difficulty (a)}

In the study of
complex variation of polarized Hodge structures,
which will be abbreviated as CVHS in the following for simplicity,
the nilpotent orbit theorem due to Schmid
is quite important, and it is a starting point of
the later studies of Cattani, Kaplan, Kashiwara, Kawai and Schmid.
However, we do not know even the formulation of
nilpotent orbit theorem for harmonic bundles.
Since the harmonic bundle can be regarded
as a pluri-harmonic map from a complex manifold
to a symmetric space,
it would be possible that we would obtain
a generalization of the nilpotent orbit theorem
for harmonic bundle, after the study would be made progressed.
(See Remark \ref{rem;b12.5.1}, for example.)
But anyway we do not have the nilpotent orbit theorem
at the moment, and thus we will find another starting point.

\begin{rem}
Now we have understood the asymptotic behaviour of
tame harmonic bundle pretty well,
and hence the author does not think
that a generalization of nilpotent orbit theorem
is necessary as a starting point,
although it would be interesting.
\hfill\qed
\end{rem}

\subsubsection{The difficulty (b)}
\label{subsubsection;b12.5.5}

We have the non-trivial eigenvalues of the residue of
the Higgs field,
and non-trivial parabolic structures.
As an example,
we have the following simple example.
Let $\Delta^{\ast}$ denote the punctured disc
$\bigl\{z\in\cnum\,\big|\,0<|z|<1\bigr\}$.

\vspace{.1in}
\noindent
{\bf Example}
Let us consider the holomorphic bundle
$E:=\nbigo_{\Delta^{\ast}}\cdot e$ of rank $1$
over $\Delta^{\ast}$.
We have the Higgs field $\theta:=\alpha\cdot dz/z$
 $(\alpha\in\cnum)$
and the metric $h$ determined by $h(e,e):=|z|^{-2a}$
 $(a\in\real)$.
Then it is easy to check that
the tuple $(E,\theta,h)$ is a harmonic bundle.

\vspace{.1in}
In the case of variation of Hodge,
the corresponding Higgs field is always nilpotent.
Hence if $\alpha\neq 0$, the example cannot be
the variation of Hodge structures.
In the case $\alpha=0$, the example is Hodge.
However if $a$ is not rational,
then the monodromy of the corresponding local system is not
quasi unipotent. Recall that it is often assumed that
the local monodromy is quasi unipotent
in the study of variation of Hodge structures.
In this sense,
the example is far from (usual) Hodge in the case
$(a,\alpha)\not\in\rnum\times\{0\}$.

\subsubsection{A starting point in the paper \cite{mochi}}

In our paper \cite{mochi},
we discussed the problem under the assumption
that the difficulty $(b)$ does not occur,
namely, the assumption of the nilpotentness
and the trivial parabolic structure.
We recall what was our starting point
in \cite{mochi},
instead of the nilpotent orbit theorem.

We put $X:=\Delta^n$, $D_i:=\{z_i=0\}$,
$D:=\bigcup_{i=1}^n D_i$.
We put $\nbigx:=X\times\cnum_{\lambda}$
and $\nbigd:=D\times\cnum_{\lambda}$.
Let $\harmonicbundle$ be a tame nilpotent harmonic bundle
with trivial parabolic structure.
We have the deformed holomorphic bundle
$\nbige$ and the $\lambda$-connection $\DD$ on $\nbigx-\nbigd$.
We prolong the sheaf $\nbige$ on $\nbigx-\nbigd$
to the sheaf $\prolong{\nbige}$ over $\nbigx$
by imposing the condition on increasing order.
(See the subsection \ref{subsection;a11.9.52}
for $\prolong{\nbige}$ and $\prolongg{\vecb}{\nbige}$.)
Then we proved the following.

\begin{prop}[Theorem 4.1 and Proposition 4.9 in \cite{mochi}]
\label{prop;11.6.1}
Under the assumption of the nilpotentness
and the trivial parabolic structures,
the $\nbigo_{\nbigx}$-module $\prolong{\nbige}$ is locally free,
and $\DD$ is a regular $\lambda$-connection,
in the sense
$\DD f\in\prolong{\nbige}\otimes\Omega_X(\log D)$
for any section
$f\in\prolong{\nbige}$.
\hfill\qed
\end{prop}

Then we obtain the holomorphic vector bundle
$V_0:=\prolong{\nbige}_{|\{O\}\times\cnum_{\lambda}}$
on the complex plane $\cnum_{\lambda}$.
The residues $\Res_i(\DD)$ induce the nilpotent endomorphism $\nbign_i$.

On the other hand,
we have the harmonic bundle $(E,\del_E,\theta^{\dagger},h)$
on the conjugate complex manifold $X^{\dagger}-D^{\dagger}$.
We put $\nbigx^{\dagger}:=X^{\dagger}\times\cnum_{\mu}$
and $\nbigd^{\dagger}:=D^{\dagger}\times\cnum_{\mu}$.
We obtain the deformed holomorphic bundle
$\nbige^{\dagger}$ and the $\mu$-connection $\DD$,
and then the prolongment $\prolong{\nbige^{\dagger}}$.
Thus we obtain the holomorphic bundle
$V_{\infty}:=\prolong{\nbige^{\dagger}}_{|\cnum_{\mu}}$
and the nilpotent endomorphism $\nbign^{\dagger}_i$
on the complex plane $\cnum_{\mu}$.

We glue $\cnum_{\lambda}$ and $\cnum_{\mu}$
by the relation $\lambda=\mu^{-1}$,
and thus we obtain $\proj^1$.
By taking a point $P\in X-D$,
we obtain the gluing of
$(V_{0},\nbign_i)$ and 
$(V_{\infty},-\nbign_i^{\dagger})$.
Thus we obtain the holomorphic vector bundle
$S(E,P)$ and the nilpotent maps
$\nbign^{\sankaku}_i:
 S(E,P)\lrarr S(E,P)\otimes\nbigo_{\proj^1}(2)$
over $\proj^1$.
The nilpotent map
$\nbign^{\sankaku}(\nbar)=\sum\nbign_i^{\sankaku}$
induces the weight filtration $W$ on $S(E,P)$.

\begin{prop} [Theorem 7.2 \cite{mochi}]\label{prop;11.6.2}
The filtered vector bundle
$(S(E,P),W)$ is a mixed twistor structure.
\hfill\qed
\end{prop}

Propositions \ref{prop;11.6.1} and \ref{prop;11.6.2} are
the starting points of our study in the paper \cite{mochi}.
Then we obtain the constantness of the filtrations,
the compatibility of the nilpotent maps,
the norm estimate, the limiting CVHS
and the purity theorem
by using some geometric argument.

\subsubsection{When the difficulty (b) occurs}
\label{subsubsection;a12.10.5}

When the difficulty (b) occurs,
we cannot use the argument in \cite{mochi} straitforwardly.
Let us see what happens in the example
in the subsubsection \ref{subsubsection;b12.5.5}.

In the example, $\del_E$  and $\theta^{\dagger}$ are as follows:
\[
 \del_E e=e\cdot (-a)\frac{dz}{z},
\quad
 \theta^{\dagger}
=\overline{\alpha}\cdot\frac{d\overline{z}}{\overline{z}}.
\]
Then we have the frame $f$ of $\nbige$ given as follows:
\[
 f:=
 \exp\bigl(
  -\overline{\alpha}\cdot\lambda\cdot \log|z|^2
 \bigr)\cdot e.
\]
The $\lambda$-connection is as follows:
\[
 \DD f=f\cdot
 \bigl(\alpha-a\cdot\lambda-\overline{\alpha}\cdot\lambda^2\bigr)\cdot
 \frac{dz}{z}.
\]
In particular,
$\Res(\DD)=\alpha-a\cdot \lambda-\overline{\alpha}\cdot\lambda^2$.
The norm of $u$ with respect to $h$ is as follows:
\[
 |f|_h=
 |z|^{-a-2\Realpart(\overline{\alpha}\cdot\lambda)}.
\]
Then we obtain
$-\ord(|f_{|\nbigx^{\lambda}}|_h)
=a+2\Realpart(\overline{\alpha}\cdot\lambda)$,
which depends on $\lambda$,
in the case $\alpha\neq 0$.
It means that 
the sheaves $\prolong{\nbige}$ or $\prolongg{\vecb}{\nbige}$
for any $\vecb\in\real^n$ are not locally free.
Namely Proposition \ref{prop;11.6.1} does not hold
in general.
(See Remark \ref{rem;b12.5.6}
on the explanation from the view point of the curvature.)

\subsubsection{How we can modify?}

First we discuss the prolongment for fixed $\lambda$,
namely we consider the prolongment of
$\nbigelambda$ to the sheaf
$\prolongg{\vecb}{\nbigelambda}$ on $\nbigx^{\lambda}$.
In this case, we can show the local freeness
by the essentially same argument as the proof of 
Proposition \ref{prop;11.6.2}.
\begin{prop}
$\prolongg{\vecb}{\nbigelambda}$ is locally free.
\hfill\qed
\end{prop}

Then we have following two structures
of $\prolongg{\vecb}{\nbigelambda}_{|\nbigd_i^{\lambda}}$
$(i=1,\ldots,n)$.
\begin{itemize}
\item
The parabolic filtration $\lefttop{i}F$,
which is a filtration of
$\prolongg{\vecb}{\nbigelambda}_{|\nbigd_i^{\lambda}}$
in the category of vector bundles.
\item
The generalized eigen decomposition $\lefttop{i}\EE$
of $\prolongg{\vecb}{\nbigelambda}_{|\nbigd_i^{\lambda}}$
with respect to the action of the residue
$\Res_i(\DDlambda)$.
\end{itemize}
They are called the KMS-structure
(Kashiwara-Malgrange-Sabbah-Simpson).

The parabolic structure $F$ is determined by 
the increasing order, namely it is given as follows:
\[
  \lefttop{i}F_{c}(\prolongg{\vecb}{\nbigelambda}_{|\nbigd_i^{\lambda}})
:=\Image\Bigl(
 \prolongg{\vecb+(c-b_i)\vecdelta_i}{\nbigelambda}
\lrarr
 \prolongg{\vecb}{\nbigelambda}
 \Bigr).
\]
The filtration $\lefttop{i}F$
and the decomposition $\lefttop{i}\EE$
are compatible,
in the sense
$\lefttop{i}F_a
=\bigoplus_{\alpha\in\cnum}
 \lefttop{i}\EE_{\alpha}\cap \lefttop{i}F_a$.
Then we obtain the following data:
\[
 \KMS(\prolongg{\vecb}{\nbige^{\lambda}},i):=
 \bigl\{
 (a,\alpha)\in\real\times\cnum\,|\,
 \lefttop{i}\Gr^F_a\lefttop{i}\EE(\prolongg{\vecb}{\nbige}_{|D_i},\alpha)
 \neq 0
 \bigr\}\subset\real\times\cnum.
\]
We put
$\KMS(\nbige^{\lambda},i)=\bigcup_{\vecb}
 \KMS(\prolongg{\vecb}{\nbigelambda},i)
\subset\real\times\cnum$.
The elements of $\KMS(\nbigelambda,i)$
is called the KMS-spectrum at $\lambda$.
The number
$\dim\lefttop{i}\Gr^F_a
 \lefttop{i}\EE(\prolongg{\vecb}{\nbige}_{|D_i},\alpha)$
is called the multiplicity
of $(a,\alpha)\in\KMS(\prolongg{\vecb},i)$.

We have the $\seisuu$-action on $\KMS(\nbigelambda,i)$ as follows:
\begin{lem}
Let $(a,\alpha)$ be an element of $\real\times\cnum$.
Then 
$(a,\alpha)\in\KMS(\nbigelambda,i)$
if and only if
$(a+1,\alpha-\lambda)\in\KMS(\nbigelambda,i)$.
The multiplicities of
$(a,\alpha)$ and $(a+1,\alpha-\lambda)$ are same.
\hfill\qed
\end{lem}

We have the bijection
$\kmsmap(\lambda):
 \real\times\cnum\lrarr\real\times\cnum$.
For $u=(a,\alpha)\in\real\times\cnum$,
we put as follows:
\[
 \kmsmap(\lambda,u):=\bigl(
 \paramap(\lambda,u),\eigenmap(\lambda,u)
 \bigr),
\quad\quad
\left\{
\begin{array}{l}
 \paramap(\lambda,u):=a+2\Re\bigl(\lambda\cdot\overline{\alpha}\bigr),\\
\mbox{{}}\\
 \eigenmap(\lambda,u):=\alpha-a\cdot\lambda-\overline{\alpha}\cdot\lambda^2.
\end{array}
\right.
\]
We note that $\eigenmap(\lambda,u)$ is the eigenvalue of the residue
$\Res(\DD)$ in the example.
We also note that $-\paramap(\lambda,u)$ is the increasing order
of $f$ in the example 
(the subsubsections \ref{subsubsection;b12.5.5}
and \ref{subsubsection;a12.10.5}).

The following proposition is essentially due to Simpson.
\begin{prop}
The map $\kmsmap(\lambda,u)$
induces the bijection
$\KMS(\nbige^0,i)\lrarr
 \KMS(\nbigelambda,i)$.
The multiplicities are preserved.
\hfill\qed
\end{prop}
Note that the map preserves the $\seisuu$-action.

Let $I$ be a subset of $\nbar=\{1,\ldots,n\}$.
Then we have the 
filtrations $\lefttop{i}F$ $(i\in I)$
and
the decompositions $\lefttop{i}\EE$ $(i\in I)$
of
$\prolongg{\vecb}{\nbigelambda}_{|\nbigd_I^{\lambda}}$:
It can be shown that they are compatible.
Then we obtain the following subset
of $\real^I\times\cnum^I=(\real\times\cnum)^I$:
\[
 \KMS(\prolongg{\vecb}{\nbigelambda},I)
:=\bigl\{
 (\veca,\vecalpha)\in\real^I\times\cnum^I\,\big|\,
 \lefttop{I}\Gr^F_{\veca}
 \lefttop{I}\EE\bigl(\prolongg{\vecb}{\nbige}^{\lambda}_{|\nbigd_I^{\lambda}},
 \vecalpha \bigr)\neq 0
 \bigr\}.
\]
We put
$\KMS(\nbigelambda,I):=
 \bigcup_{\vecb}\KMS(\prolongg{\vecb}{\nbigelambda},I)$.
The element
$\vecu=(\veca,\vecalpha)\in
 \KMS(\prolongg{\vecb}{\nbigelambda},I)$
is called $KMS$-spectrum,
and the number
$\dim\lefttop{I}\Gr^F_{\veca}
 \lefttop{I}\EE(\prolongg{\vecb}{\nbige}^{\lambda}_{|\nbigd_I^{\lambda}},
 \vecalpha)$
is called the multiplicity of $\vecu$.

Similarly to the case of $I=\{i\}$,
we have the $\seisuu^I$-action
on $\KMS(\nbigelambda,I)$,
preserving the multiplicities.

\begin{prop}
We have the bijection
$\kmsmap(\lambda):\KMS(\nbige^0,I)\lrarr\KMS(\nbigelambda,I)$,
preserving the $\seisuu$-action and the multiplicities.
\hfill\qed
\end{prop}

Then we put as follows, for any element
$\vecu\in\KMS(\nbige^0,\nbar)$:
\[
 \lefttop{\nbar}\nbigg_{\vecu}^{\lambda}
:=\lefttop{\nbar}\Gr^{F}_{\paramap(\lambda,\vecu)}
 \lefttop{\nbar}\EE(\prolongg{\vecb}{\nbige},\eigenmap(\lambda,\vecu)).
\]

The residue $\Res_i(\DD)$ induces the endomorphism
of $\lefttop{\nbar}\nbigg_{\vecu}^{\lambda}$.
The unique eigenvalue of the endomorphism
is $\eigenmap(\lambda,u_i)$ by our construction.
The nilpotent part is denoted by $\nbign_i^{\lambda}$.

\begin{rem}
The tame harmonic bundle $\harmonicbundle$
is nilpotent and with the trivial parabolic structure,
if and only if
the set $\KMS(\nbige^0,\nbar)$ is same as
$\seisuu^n\times\{0\}$.
Due to the $\seisuu$-action,
we have only to consider the spectrum
$0\in\KMS(\nbigelambda,\nbar)$,
and we have
$\lefttop{\nbar}\nbigg_{0}^{\lambda}
=\prolong{\nbige}_{|(O,\lambda)}$.
\hfill\qed
\end{rem}

Then we obtain the family
$\bigl\{
 \lefttop{\nbar}\nbigg_{\vecu}^{\lambda}\,\big|\,
 \lambda\in\cnum \bigr\}$
of vector spaces.
We would like to give the structure of a holomorphic vector bundle
over the complex plane $\cnum_{\lambda}$.
In the case where $\harmonicbundle$ is nilpotent and
with trivial parabolic structure,
we have the holomorphic bundle
$\prolong{\nbige}$ over $\nbigx$,
and thus we obtain the holomorphic bundle
$\prolong{\nbige}_{|\cnum_{\lambda}}$,
which gives the structure of the holomorphic vector bundle
of the family
$\bigl\{\lefttop{\nbar}\nbigg_{0}^{\lambda}\,\big|\,
 \lambda\in\cnum \bigr\}$.
As we have already said,
the sheaves
$\prolong{\nbige}$ or $\prolongg{\vecb}{\nbige}$
on $\nbigx$ are not locally free in general,
we cannot apply the method directly.

\vspace{.1in}

Let us pick any point $\lambda_0\in\cnum_{\lambda}$.
Let us pick an element $\vecb\in\real^n$ such that
$b_i\not\in\KMS(\nbige^{\lambda_0},i)$
for $i\in\nbar$.
Let us take a sufficiently small positive number $\epsilon_0$.
We put
$\nbigx(\lambda_0,\epsilon_0)
 :=\Delta(\lambda_0,\epsilon_0)\times X$
and $\nbigx^{\lambda}:=\{\lambda\}\times X$.
\begin{prop}\mbox{{}}
\begin{itemize}
\item
Then $\prolongg{\vecb}{\nbige}$ is a locally free sheaf
on $\nbigx(\lambda_0,\epsilon_0)$.
\item
For any point $\lambda\in\Delta(\lambda_0,\epsilon_0)$,
we have the canonical isomorphism
$\prolongg{\vecb}{\nbige}_{|\nbigxlambda}
\simeq
 \prolongg{\vecb}{\nbigelambda}$.
\hfill\qed
\end{itemize}
\end{prop}

Let $\nbigd_i(\lambda_0,\epsilon_0)$
denote $D_i\times\Delta(\lambda_0,\epsilon_0)$,
and $\nbigd_i^{\lambda}$ denote $D_i\times\{\lambda\}$.
We have the filtration $\lefttop{i}\Fzero$
and the decomposition $\lefttop{i}\EE$
of the vector bundle
$\prolongg{\vecb}{\nbige}_{|\nbigd_i(\lambda_0,\epsilon_0)}$:
\begin{itemize}
\item
 The restriction of $\lefttop{i}\Fzero$
 to $\nbigx^{\lambda_0}$ is same as 
 the parabolic filtration $\lefttop{i}F$
 of $\prolongg{\vecb}{\nbige^{\lambda_0}}_{|\nbigd_i^{\lambda_0}}$.
\item
 The restriction $\lefttop{i}\EEzero$
 to $\nbigx^{\lambda_0}$ is same as 
 the generalized eigen decomposition $\lefttop{i}\EE$
 of $\prolongg{\vecb}{\nbige^{\lambda_0}}_{|\nbigd_i^{\lambda_0}}$.
\item
 $\lefttop{i}\Fzero$ and $\lefttop{i}\EEzero$ are compatible.
\end{itemize}

Let us explain the restriction of
$\Fzero$ and $\EEzero$
to $\nbigd_i^{\lambda}$ for a point $\lambda$
near $\lambda_0$.
The relation of the filtration $\Fzero$
and the parabolic filtration $F$ on 
$\prolongg{\vecb}{\nbige}_{|\nbigd_i^{\lambda}}$
is as follows: 
\[
 \lefttop{i}\Fzero_b\bigl(
 \prolongg{\vecb}{\nbige}_{|\nbigd_i(\lambda_0,\epsilon_0)}
 \bigr)_{|\nbigd_i^{\lambda}}
=\lefttop{i}F_{b+\epsilon}\bigl(
 \prolongg{\vecb}{\nbigelambda}_{|\nbigd_i^{\lambda}}
 \bigr).
\]
Here $\epsilon$ denotes a small number
depending on $b$ and $\lambda$.
We also have the following:
\[
 \lefttop{i}\EEzero\bigl(
 \prolongg{\vecb}{\nbige}_{|\nbigd_i(\lambda_0,\epsilon_0)}
 ,\alpha \bigr)_{|\nbigd_i^{\lambda}}
=\bigoplus_{|\beta-\alpha|<\eta}
 \lefttop{i}\EE\bigl(
 \prolongg{\vecb}{\nbigelambda}_{|\nbigd_i^{\lambda}},\beta
 \bigr).
\]

For any subset $I\subset\nbar$,
we obtain the filtrations $\lefttop{i}\Fzero$ $(i\in I)$
and $\lefttop{i}\EEzero$ $(i\in I)$
of the vector bundle
$\prolongg{\vecb}{\nbige}_{|\nbigd_I(\lambda_0,\epsilon_0)}$,
and they satisfies the relations as above.

Then we put as follows, for any element
$\vecu\in\KMS(\nbige^0,\nbar)$:
\[
 \lefttop{\nbar}\nbigg_{\vecu}^{(\lambda_0)}
=\lefttop{\nbar}\Gr^{\Fzero}_{\paramap(\lambda_0,\vecu)}
 \lefttop{\nbar}\EEzero
 \bigl(
 \prolongg{\vecb}{\nbige}_{|\nbigd_{\nbar}(\lambda_0,\epsilon_0)},
 \eigenmap(\lambda_0,\vecu)
 \bigr).
\]
Then we have
$\lefttop{\nbar}\nbigg^{(\lambda_0)}_{\vecu\,|\,\lambda}
=\lefttop{\nbar}\nbigg^{\lambda}_{\vecu}$
for any point $\lambda\in\Delta(\lambda_0,\epsilon_0)$.
When the intersection $S:=\Delta(\lambda_0,\epsilon_0)\cap
  \Delta(\lambda_1,\epsilon_1)$ is not empty,
we have the canonical isomorphism
$\lefttop{\nbar}\nbigg^{(\lambda_0)}_{\vecu\,|\,S}
\simeq
 \lefttop{\nbar}\nbigg^{(\lambda_1)}_{\vecu\,|\,S}$.
Thus we obtain the global vector bundle
$\lefttop{\nbar}\nbigg_{\vecu}$
on $\cnum_{\lambda}$
such that
$\lefttop{\nbar}\nbigg_{\vecu\,|\,\Delta(\lambda_0,\epsilon_0)}
\simeq
 \lefttop{\nbar}\nbigg_{\vecu}^{(\lambda_0)}$
and $\lefttop{\nbar}\nbigg_{\vecu\,|\,\lambda}\simeq
 \lefttop{\nbar}\nbigg^{\lambda}_{\vecu}$.
We also have the nilpotent endomorphism
$\nbign_i$ $(i\in\nbar)$.

By the same construction
for the tame harmonic bundle
$(E,\del_E,\theta^{\dagger},h)$ on $X^{\dagger}-D^{\dagger}$,
we obtain the holomorphic bundle
$\lefttop{\nbar}\nbigg^{\dagger}_{\vecu}$
for any element $\vecu\in\KMS(\nbige^{\dagger\,0},\nbar)$,
with the nilpotent map
$\nbign_i^{\dagger}$.

We have the morphism
$\real\times\cnum\lrarr\real\times\cnum$
given by $(a,\alpha)\longmapsto (-a,\overline{\alpha})$.
It induces the bijection
$\KMS(\nbige^0,\nbar)\lrarr\KMS(\nbige^{\dagger\,0},\nbar)$.
We denote the correspondence by
$\vecu\longmapsto \vecu^{\dagger}$.
Then we have the gluing
of $\lefttop{\nbar}\nbigg_{\vecu}$
and $\lefttop{\nbar}\nbigg^{\dagger}_{\vecu}$,
and thus we obtain the vector bundle
$S_{\vecu}(E,P)$.
We also obtain the nilpotent morphism
$\nbign^{\sankaku}_i:S_{\vecu}(E,P)
\lrarr S_{\vecu}(E,P)\otimes\nbigo(2)$.
As in the previous paper,
we put $\nbign^{\sankaku}(\nbar)=\sum_{i=1}^n\nbign^{\sankaku}_i$,
which induce the weight filtration $W$ on $S_{\vecu}(E,P)$.
Then it can be shown that
$(S_{\vecu}(E,P),W)$ is a mixed twistor structure,
which is called the limiting mixed twistor structure.

\begin{rem}
We have another gluing, and the resulted vector bundle
is denoted by $S^{\can}_{\vecu}(E)$.
In fact,
it is more close to
the traditional construction of the limiting mixed Hodge
than $S_{\vecu}(E,P)$.
\hfill\qed
\end{rem}

\subsubsection{Polarization of the limiting mixed twistor structures}

We discuss the naturally induced polarization $S$
on the limiting mixed twistor structure,
following Sabbah \cite{sabbah},
who considered the polarized limiting mixed twistor structures
for tame harmonic bundle on a quasi projective curve.
The following theorem is one of the main goals
in the study of Part \ref{part;a11.17.101}--\ref{part;a11.17.102}.

\begin{thm}[Theorem \ref{thm;10.26.110}]\label{thm;b12.5.10}
The tuples $\bigl(S_{\vecu}^{\can}(E),W,\vecN^{\sankaku},S\bigr)$ and
$(S_u(E,P),W,\vecN^{\sankaku},S)$ $(P\in X-D)$ are
polarized mixed twistor structure.
\hfill\qed
\end{thm}

We use Theorem \ref{thm;b12.5.10}
by taking the associated graded objects.
We have the associated graded vector bundle
$V^{(0)}:=S^{\can}(E)$,
on which we have the naturally induced
filtration $W^{(0)}$,
the nilpotent maps $\vecN^{(0)}$
and the pairing $S^{(0)}$.
Then the tuple
$\bigl(V^{(0)},W^{(0)},\vecN^{(0)},S^{(0)}\bigr)$
is again a polarized mixed twistor structure.
We can take an appropriate torus action
on the tuple,
and thus it is a polarized Hodge structure.
Moreover, it can be shown that
$\bigl(V^{(0)},W^{(0)},\vecN^{(0)},S^{(0)}\bigr)$
is a nilpotent orbit.
Since the nilpotent orbit
was studied very closely
in the theory of variation of Hodge structures,
we can say that we understand 
the tuple $\bigl(V^{(0)},W^{(0)},\vecN^{(0)},S^{(0)}\bigr)$
very well,
due to the classical results on variation of Hodge structures.
Much information on the tuple
$\bigl(S^{\can},W,\vecN,S\bigr)$ can be obtained
from
$\bigl(V^{(0)},W^{(0)},\vecN^{(0)},S^{(0)}\bigr)$.
For example,
we can obtain the compatibility of the nilpotent maps
and vanishing cycle theorem.
We can also apply the lemmas due to Kashiwara and Saito
on the nilpotent orbit to $\bigl(S_{\vecu}^{\can},W,\vecN,S\bigr)$.
In this sense, the study of tame harmonic bundle
is reduced to the study of variation of Hodge structures.

In the previous paper \cite{mochi}, we often used the argument
to take a `limit' of a sequence of tame harmonic bundles.
For example, 
we consider the morphisms
$\lefttop{\nbar}\psi_m:X-D\lrarr X-D$
given by
$(z_1,\ldots,z_n)\longmapsto (z_1^m,\ldots,z_n^m)$,
and we consider the sequence
$\bigl\{ \psi_m^{\ast}\harmonicbundle\bigr\}$
of harmonic bundles.
Under the assumption that
$\harmonicbundle$ is nilpotent and with trivial parabolic structure,
we obtained the complex variation of Hodge structure
as the `limit'.
We also considered the morphisms
$\lefttop{\itibar}\psi_m:X-D\lrarr X-D$
given by $(z_1,\ldots,z_n)\longmapsto (z_1^m,z_2,\ldots,z_n)$,
and we used a limit of the sequence
$\bigl\{\lefttop{\itibar}\psi_m^{\ast}\harmonicbundle\bigr\}$
to derive a sequential compatibility of the residues.

The argument to take a limit does not work if the residues of the Higgs field
is not nilpotent.
In the case of the example
in the subsubsection \ref{subsubsection;b12.5.5},
we have
$\lefttop{\itibar}\psi_m^{\ast}\theta=m\cdot \alpha\cdot dz/z$,
and thus it is not easy to guess what is `limit'
for $m\to\infty$.

In a sense, taking the associated graded tuple
$\bigl(V^{(0)},W^{(0)},\vecN^{(0)},S^{(0)}\bigr)$
corresponds to taking a `limit'.
Let us consider the case
that $\harmonicbundle$ is nilpotent and with trivial parabolic
structure.
As is already mentioned,
we obtain the limiting CVHS
$(E^{(\infty)},\delbar_{E^{(\infty)}},
 \theta^{(\infty)},h^{(\infty)})$ as a limit.
Then it is easy to see that
the limiting mixed twistor structure
$(E^{(\infty)},\delbar_{E^{(\infty)}},\theta^{(\infty)},h^{(\infty)})$
is naturally isomorphic to
the associated graded mixed twistor structure
$\bigl(V^{(0)},W^{(0)}\bigr)$.

Most of our argument to take a limit in the previous paper
\cite{mochi} can be replaced with the argument
to consider the associated graded tuple.
The only exception is the proof of the constantness of
the filtration, for which
we use some elementary calculus instead of 
taking a limit.

\begin{rem}
As is mentioned above,
we do not use the argument to take a limit
in this paper.
However, 
it seems significant to observe that CVHS appears
as the limit.
\hfill\qed
\end{rem}

%% file: 103.1.tex

\subsection{On the purpose \ref{number;11.6.2}}

\subsubsection{Pure twistor $D$-module and Sabbah's program}
\label{subsubsection;b12.6.300}

The author thinks that
we have already understood the asymptotic behaviour of
tame harmonic bundle pretty well.
In Part IV, we would like to apply the study in Part II--III
to the theory of pure twistor $D$-modules of Sabbah.

Following Simpson's Meta-Theorem,
it is interesting and natural to ask whether we can construct
the theory of ``twistor module'', which should be a generalization
of the theory of Hodge module of Saito.
In this direction, Sabbah has already done a remarkable work
(\cite{sabbah2}).
He gave a definition of regular pure twistor $D$-modules
and proved a decomposition theorem.
The motivation of Sabbah is to attack
a conjecture of Kashiwara,
so that we recall a part of the conjecture.

\begin{conj}[The regular holonomic version of Kashiwara's conjecture]
\label{conjecture;b12.5.25}\mbox{{}}
\begin{description}
\item[(Push-forward)]
Let $X$ and $Y$ be a quasi projective manifold over 
the complex number field $\cnum$.
Let $f:X\lrarr Y$ be a proper morphism.
Let $\nbigf$ be a semisimple regular holonomic $D$-module on $X$.
\begin{itemize}
\item
Then the push-forward $Rf_{+}\nbigf$ is isomorphic to
the direct sum $\bigoplus R^if_{+}\nbigf$
in the derived category of cohomologically holonomic complexes on $Y$.
\item
The hard Lefschetz theorem for $\bigoplus R^if_{+}\nbigf$ holds.
\item
Each $R^if_{+}\nbigf$ is semisimple.
\end{itemize}
\item[(Vanishing cycles)]
Let $X$ be a quasi projective manifold, 
and $\nbigf$ be a semisimple regular holonomic $D$-module on $X$.
Let $f$ be a holomorphic function on $X$.
We take the nearby cycle functor and the vanishing cycle functor
along $f$. Then the associated graded object with respect to the
monodromy weight filtration is semisimple.
\hfill\qed
\end{description}
\end{conj}

\begin{rem}\label{rem;b12.5.40}
The conjecture of Kashiwara is stronger than the statement above.
In fact, he conjectured that the statement is true
for semisimple holonomic $D$-modules
which are not necessarily regular.
See {\rm\cite{k5}} for more precise.
\hfill\qed
\end{rem}

Sabbah's program to attack the conjecture is as follows:
\begin{description}
\item [Step 1.]
 To establish the correspondence of
 tame harmonic bundle and semisimple local system.
\item[Step 2.]
 To give a definition of regular pure twistor $D$-module
 and to prove the decomposition theorem for 
 regular pure twistor $D$-module.
\item[Step 3.]
 To establish the correspondence of
 tame harmonic bundles and regular pure twistor $D$-modules.
\end{description}

As for Step 1,
there is known the classical result of Corlette
who proved that semisimple local system on a projective manifold
corresponds to a harmonic bundle.
The result was generalized by Jost-Zuo,
who proved there exists a pluri-harmonic metric
on a semisimple local system on a quasi projective manifold,
in other words, there exists the structure of harmonic bundle
on any semisimple local system on a quasi projective manifold.
It is refined in \cite{mochi3},
and we know that 
a flat bundle $(E,\nabla)$
on a quasi projective manifold is semisimple
if and only if there exists a tame pure imaginary pluri-harmonic metric
on $(E,\nabla)$.
Hence we can say that the Step 1 is established.

As is already remarked,
Sabbah established the step 2 in \cite{sabbah2}.
Sabbah also proved that
a harmonic bundle, without singularity,
gives a regular pure twistor $D$-module
for the step 3.
As a corollary of the results due to Corlette and himself,
Sabbah obtained the decomposition theorem and the hard Lefschetz theorem
for semisimple local system.

\begin{rem}\mbox{{}}
\begin{itemize}
\item
After the author submitted the first version of this paper
to math arXiv,
Sabbah kindly informed on the revision {\rm\cite{sabbah2}}
of his paper {\rm\cite{sabbah}},
although
we mainly refer {\rm\cite{sabbah}} in this paper.
\item
The definition of pure twistor $D$-module in {\rm\cite{sabbah2}} is 
more general than that given in {\rm\cite{sabbah}},
and it is most appropriate for attacking
Kashiwara's conjecture.
But it can also be said that
it is slightly narrow from the view point of Simpson's Meta Theorem.
See Appendix.
\item
Although
Sabbah discusses the pure twistor $D$-modules
which are not necessarily regular,
we consider only
regular pure twistor $D$-modules
in this paper,
even if we omit to distinguish ``regular''.
\hfill\qed
\end{itemize}
\end{rem}

\subsubsection{The goal of the part IV, the conjectures of
 Sabbah and Kashiwara}

Let $X$ be a complex manifold.
Let $Z$ be an irreducible closed subset of $X$.
A tame harmonic bundle generically defined on $Z$
is defined to be a tame harmonic bundle
defined over a smooth Zariski open subset of $Z$.
The purpose of the part IV is to establish the step 3,
namely
we will prove that a tame harmonic bundle generically defined 
on $Z$ gives the pure twistor $D$-module of weight $0$.
More precisely,
we will prove the following correspondence:
\begin{thm}[Theorem \ref{thm;b12.5.20}, Theorem \ref{thm;04.2.4.1}]
 \label{thm;b12.5.30}
We have the bijective correspondences:
\[
 \VPTgen(Z,w)\simeq \MPT(Z,w),
\quad\quad
 \VPTgen^{pi}(Z,w)\simeq \MPT^{pi}(Z,w).
\]
(See the subsubsection {\rm\ref{subsubsection;b12.5.21}}
for the notation.)
\hfill\qed
\end{thm}

Recall that Saito established the following correspondence.
\begin{prop}[Saito, \cite{saito2}]\label{prop;04.2.23.1}
Variations of pure polarized Hodge structures of weight $w$
which are defined over Zariski open subsets of $Z$
correspond to
pure polarized Hodge modules of weight $w$ whose strict supports are $Z$.
\hfill\qed
\end{prop}

Theorem \ref{thm;b12.5.30} is a natural generalization
of Proposition \ref{prop;04.2.23.1}.
Once we have established the resemblance
of the asymptotic behaviours of tame harmonic bundles
and CVHS,
we can use a part of Saito's idea
to prove Theorem \ref{thm;b12.5.30}.
In fact,
it can be said that
the crucial ideas can be found in
the subsections 3.19--3.21 of \cite{saito2}.
However, it is not so clear for the author how to modify
the arguments in 3.1--3.5 and the most part of 3.b of \cite{saito2}
for our purpose.
Hence we will go along the other route.

\vspace{.1in}

Sabbah conjectured that
every semisimple regular holonomic $D$-module on
a complex projective manifold
underlies a pure imaginary pure twistor $D$-module.
(See the section 4.2.c in\cite{sabbah2} or
 the subsection \ref{subsection;04.2.5.2}.)
Once Theorem \ref{thm;b12.5.30} and Step 1 above are established, 
it is easy to show that his conjecture is true.
Namely we obtain the following theorem.
(See the subsubsection \ref{subsubsection;b12.5.21}
and the subsection \ref{subsection;04.2.5.2} for the notation.)

\begin{thm}[The conjecture of Sabbah, Theorem \ref{thm;04.2.5.20}]
The map $\Xi_{Dol}:\MPT^{pi}(Z,0)\lrarr \RHD^{ss}(Z)$
is surjective.
\hfill\qed
\end{thm}

As a result, we obtain the regular holonomic version
of Kashiwara's conjecture,
combining the results of Sabbah (\cite{sabbah2}) and us.

\begin{cor}
Conjecture {\rm\ref{conjecture;b12.5.25}} is true.
\hfill\qed
\end{cor}

\begin{rem}
K. Vilonen informed the author of the work of D. Gaitsgory 
{\rm\cite{gaitsgory}},
who proved de Jong's conjecture.
Since V. Drinfeld proved that de Jong's conjecture
implies the regular holonomic version of Kashiwara's conjecture
{\rm \cite{drinfeld}},
Conjecture {\rm\ref{conjecture;b12.5.25}} has been established
also by their works.
\hfill\qed
\end{rem}

%% file: 105.tex

\subsubsection{Part \ref{part;b11.17.100}, Section \ref{section;a11.9.50}}

The Part \ref{part;b11.17.100} is a preparation
for the subsequent parts.
The author expects that the readers can skip Part \ref{part;b11.17.100}
until they need it.

In the section \ref{section;a11.9.50},
we prepare some notation and lemmas from several areas.
In the subsection \ref{subsection;a11.9.51},
we prepare the notation of some sets and the functions.
The maps $\kappa_c$ and $\nu_c$
in the subsubsection \ref{subsubsection;b11.11.10}
are used to describe the descent of sections for ramified covering.
The maps $\kmsmap$, $\eigenmap$, and $\paramap$
in the subsubsection \ref{subsubsection;b11.11.11}
are used for the control of the $KMS$-structure of tame harmonic
bundles.

In the subsection \ref{subsection;a11.9.52},
we recall the prolongment of a holomorphic vector bundle
with a hermitian metric over $X-D$ to an $\nbigo_X$-sheaf on $X$.
Here $X$ denotes a complex manifold,
and $D$ denotes a normal crossing divisor.
We see when the prolonged sheaf is locally free
(Lemma \ref{lem;9.5.10}).

In the subsection \ref{subsection;a11.9.53},
we prepare something on the $\mu_{c}$-equivariant
holomorphic bundles.
It will be useful when we consider the descent
of holomorphic bundle.

In the subsection \ref{subsection;a11.9.54},
we give a lemma for pluri-subharmonic function
and convexity.
It will be used in the proof of preliminary constantness
of the filtrations in the subsubsection \ref{subsubsection;11.9.1}.
Although the argument is elementary,
it is one of the key steps.
Hence we give some detail.

In the subsection \ref{subsection;a11.9.55},
we consider the distributions given by a polynomials 
of logarithmic functions.
The result will be used when we calculate
the specialization of the sesqui-linear pairing
of the $\nbigr$-triples obtained from tame harmonic bundle
(the subsubsection \ref{subsubsection;b11.9.1}).

In the subsection \ref{subsection;a11.9.56},
we mainly give some lemmas on a metrics on a finite dimensional vector space.
They are used in the subsection \ref{subsection;b11.9.2}
and the subsection \ref{subsection;b11.9.3}.
The notation $\EE_{\epsilon}$ $(\epsilon>0)$
is introduced.
We also recall the relation between the norm of isomorphism
and the distance of the hermitian metrics.
The result will be used in the subsection \ref{subsection;04.1.27.61}.

In the subsection \ref{subsection;a11.9.57},
we recall two kinds of lemmas for the acceptable bundles.
One is the vanishing of the higher cohomology groups
(Lemma \ref{lem;9.10.72} and corollary \ref{cor;10.11.56}).
The other is corollary \ref{cor;11.28.15},
which controls the estimate of the increasing order of
holomorphic sections.

In the subsection \ref{subsection;a11.9.58},
we give a lemma for a complex of Hilbert space bundles
over the disc (Lemma \ref{lem;9.10.75}).
It will be used to obtain a locally free prolongment
of the deformed holomorphic bundle of tame harmonic bundle
(the subsection \ref{subsection;b11.11.30}).
Although the procedure is standard,
we have to care the infinite dimensionality,
and we do not know an appropriate reference.
Thus we give some detail.
We also recall a standard lemma for embeddings of Sobolev spaces,
which will be used in the subsection \ref{subsection;a11.9.59}.

In the subsection \ref{subsection;a11.9.59},
we give an estimate of Higgs field of harmonic bundle.
In the subsection \ref{subsection;a11.9.60},
we give an improvement of the convergency of
a sequence of harmonic bundles,
which was given in our previous paper \cite{mochi}.
Since we do not use an argument to take a limit,
the reader can skip the subsection \ref{subsection;a11.9.60}.

In the subsection \ref{subsection;04.1.27.60},
we recall the relation of Higgs field and a twisted map
associated with the flat $\lambda$-connection with a hermitian metric.
It will be used in the subsubsection \ref{subsubsection;04.1.27.65}.

\subsubsection{Section \ref{section;a11.9.61}}

Following Simpson, we give some detail
on the relation of Hodge structure and twistor structure.
In the subsection \ref{subsection;a11.11.50},
we introduce some terminology and the notation.

In the subsection \ref{subsection;a11.11.51},
we recall the equivalence of the category
of equivariant holomorphic vector bundle over $\proj^1$
and the category of bi-filtered vector space.
The equivalence is compatible with real structures.

In the subsection \ref{subsection;a11.11.52},
we give the concrete description of
the Tate objects, the objects $\nbigo(p,q)$ and $\nbigo(n)$
in the category of twistor structures.
In particular, we give the isomorphism
$\iota_{(p,q)}:\sigma^{\ast}\nbigo(p,q)\lrarr\nbigo(q,p)$.
The isomorphism is fixed in the sequel.
We also compare $\Tate(n)$ with the Tate objects in the Hodge theory.

In the subsection \ref{subsection;a11.11.53},
we recall the equivalence of the category of
polarized pure Hodge structures and the category of
equivariant polarized pure twistor structures.
We also introduce polarized mixed twistor structure,
and we see the equivalence of the category
of polarized mixed Hodge structures
and the category of equivariant polarized mixed twistor structures.
Note that our choice of the signature of the nilpotent maps
is different from that in the standard Hodge theory.

In the subsection \ref{subsection;a11.11.54},
we recall the variation of twistor structures,
or more generally, the variation of $\proj^1$-holomorphic bundles.
An example given in the subsubsection \ref{subsubsection;10.5.15}
is important for our understanding of the harmonic bundles.
We also see the equivalence of the category of 
variation of Hodge and the category of the variation of
equivariant pure twistor structures,
which is compatible with some additional structures.

In the subsection \ref{subsection;a11.11.55},
we introduce the twistor nilpotent orbit.
We see that it is a generalization of nilpotent orbit
in the Hodge theory (Proposition \ref{prop;b11.11.1}).

In the subsection \ref{subsection;a11.11.56},
we see that a split polarized mixed twistor structure
gives a nilpotent orbit in the Hodge theory
(Corollary \ref{cor;a11.10.25}).
Since we can pick an appropriate torus action
on the split polarized mixed twistor structure,
we can regard it as a split polarized mixed Hodge structure.
Thus the result may be known in the Hodge theory, probably.
However it is very important for our application,
and hence we give some detail.
In particular, we can always obtain the nilpotent orbit
in the Hodge theory from a polarized mixed twistor structure,
by taking the associated graded object.
This is one of the key steps to reduce our study of
tame harmonic bundle to the classical study of Hodge structures.

In the subsection \ref{subsection;a11.11.57},
we see that the polarized mixed twistor structure
$(V,W,\vecN,S)$ induces the polarized structure
on the primitive part $P_h\Gr^{W(N)}_h(V)$
(Proposition \ref{prop;10.6.5}).
By using the result,
we see that the tuple of nilpotent maps $\vecN$ is
strongly sequentially compatible
(Lemma \ref{lem;10.26.111}).

In the subsection \ref{subsection;a11.11.58},
we translate some results known for the Hodge structure,
due to Kashiwara, Kawai and Saito,
to the results for the twistor structure
(Proposition \ref{prop;9.23.6},
 Corollary \ref{cor;9.23.7},
Proposition \ref{prop;a11.11.4}
and Lemma \ref{lem;9.26.20}).
They are crucial for our study to relate
tame harmonic bundles and pure twistor $D$-modules.

In the subsection \ref{subsection;a11.11.59},
we give the concrete correspondence of the twistor structure
in the sense of Simpson and those in the sense of Sabbah.

\subsubsection{Section \ref{section;d11.12.1}}

We give some definition of the compatibility of
decompositions, filtrations and nilpotent maps.
Although we refer the definition of 
`sequential compatibility' and `strongly sequential compatibility'
from our previous paper \cite{mochi},
we do not use the lemmas in \cite{mochi} essentially.

In the subsection \ref{subsection;d11.12.2},
we give definitions of the compatibility
of filtrations and decompositions on a vector space.
In the subsection \ref{subsection;d11.12.3},
we give definitions of the compatibility
of filtrations and decompositions on a vector bundle.
In the subsection \ref{subsection;d11.12.4},
we give definitions of the compatibility
of filtrations, decompositions and nilpotent maps.

In the subsection \ref{subsection;d11.12.5},
we give some lemmas for extending a splitting
given on a divisor.
In the subsection \ref{subsection;d11.12.6},
we give definitions of compatibility of
decompositions, filtrations and nilpotent maps
given on divisors.
By using the result in the subsection \ref{subsection;d11.12.5},
we see the existence of splitting.

\subsubsection{Section \ref{section;d11.12.10}}

We consider a compatible tuple of filtrations
on a discrete valuation ring $R$,
such that the splitting is given on the generic point.
We assume that some nice property holds on a generic point $K$.
we also assume the nice property holds
on the associated graded vector bundle on $R$.
Under such assumptions, we see that
the nice property holds on $R$.

In the subsection \ref{subsection;d11.12.11},
we discuss the sequential compatibility of
the nilpotent maps.
In the subsection \ref{subsection;d11.12.12},
we discuss the strictness of the morphism
with the filtrations.

The results in the section \ref{section;d11.12.10}
will be very useful,
when we combine them with the limiting mixed twistor theorem.
Briefly and imprecisely speaking,
we can derive some information
for the associated graded bundle of the parabolic filtrations.
Then we can obtain the information
of the original bundle
by using the results in the section \ref{section;d11.12.10}.

\subsubsection{Section \ref{section;a12.10.1}}
We give easy and basic examples of harmonic bundles
on a punctured disc, which we call model bundles.
They are fundamental for the study of the asymptotic behaviour.
In a sense, the study of 
the asymptotic behaviour of general tame harmonic bundles
can be reduced to these basic examples.
The author apologizes that we use the notation
introduced in the section \ref{section;a11.17.50}.
We also refer the subsection 3.2 in \cite{mochi}
for model bundles,
although some constants in \cite{mochi}
are different from those in this paper.

%% file: 105.1.tex

\subsubsection{Part \ref{part;a11.17.101}, Section \ref{section;b11.14.1}}

In Part \ref{part;a11.17.101},
we discuss the prolongment of the deformed holomorphic bundle
of tame harmonic bundle.

In the section \ref{section;b11.14.1},
we recall the result of Simpson on the study
of tame harmonic bundles over the punctured disc,
with minor generalization.
They play the fundamental role in the study of
tame harmonic bundles on a higher dimensional complex manifold.

In the subsection \ref{subsection;b11.9.2},
we give some detail on Simpson's Main estimate,
that is,
the estimate on the norm of Higgs field around the singularity.
Since we would like to use the result in the higher dimensional case,
we clarify the dependence of the constants.
We also see the asymptotic orthogonality 
of the generalized eigen decomposition.

In the subsection \ref{subsection;a11.14.2},
we recall the results on the prolongment
of the deformed holomorphic bundle $\nbigelambda$
for a fixed $\lambda$.
We introduce the KMS-structure of the prolongment
$\prolongg{\vecb}{\nbigelambda}$, and we see the functoriality
of the structure.
In particular, we give some detail on the functoriality
for pull backs via the ramified covering.
It will be useful for the study in the higher dimensional case.

In the subsection \ref{subsection;10.11.35},
we recall the basic comparison due to Simpson.
As a result, the KMS-structure at $\lambda$
can be controlled by the KMS-structure at $0$,
and we see that the weight filtration at $\lambda$
is equivalent to that at $0$.
We also obtain a rough relation of the frame at $\lambda$
and the frame at $0$,
which will be used to show the asymptotic orthogonality.

In the subsection \ref{subsection;a11.14.3},
we give some detail on the space $H(\nbigelambda)$
of the multi-valued flat sections
of the deformed holomorphic bundle $\nbigelambda$.
We introduce the KMS-structure of $H(\nbigelambda)$,
and we compare it with the KMS-structure of
$\prolongg{\vecb}{\nbigelambda}$.
In the subsubsection \ref{subsubsection;c11.14.1},
we introduce the notion of `generic' with respect
to the KMS-structure.
In the subsubsection \ref{subsubsection;9.11.25},
we see that the prolongment by an increasing order
is equivalent to the quasi canonical prolongment,
if $\lambda$ is generic.
In this sense,
the prolongment for generic $\lambda$ is canonically given,
even if we forget the metric $h$.

In the subsection \ref{subsection;a11.14.4},
we consider the family of the spaces
$\bigl\{H(\nbigelambda)\,\big|\,\lambda\in\cnum_{\lambda}^{\ast}\bigr\}$.
We introduce the decomposition $\EEzero$ and $\Fzero$
which are defined on a neighbourhood of
$\lambda_0\in\cnum^{\ast}_{\lambda}$.

In the subsection \ref{subsection;a11.14.5},
we see the asymptotic orthogonality of
the generalized eigen decomposition, the parabolic filtration,
and the weight filtration.
They are used in the proof of the limiting mixed twistor theorem
in the case of curves (the subsection \ref{subsection;c11.14.5}).
The asymptotic orthogonality of the generalized eigen decomposition
is also used for the local prolongment of $\nbige$
(the subsection \ref{subsection;b11.11.30}).

In the subsection \ref{subsection;04.1.27.61},
we give a maximum principle
for the distance of the harmonic metrics on a punctured disc.
The result will be used to show a characterization of tameness
in the subsubsection \ref{subsubsection;04.1.27.65}.

\subsubsection{Section \ref{section;b11.16.1}}

We give some detail on the prolongment of the deformed holomorphic
bundles of the tame harmonic bundle
on $\Delta^{\ast}\times\Delta^{n-l}$.
The section is one of the hearts of this paper.

In the subsection \ref{subsection;b11.16.2},
we give the remark on the constantness of the KMS-spectrum.
We also see that the tame harmonic bundle of rank one
is very easy to understand.
We use the facts without mention.
In the subsubsection \ref{subsubsection;04.1.27.65},
we give a simple characterization of tameness
(Corollary \ref{cor;04.1.27.51}),
which is useful when we check the tameness of a harmonic bundle.
(See the subsubsection \ref{subsubsection;04.1.27.70}).

In the subsection \ref{subsection;b11.9.3},
we give the estimate of the Higgs bundle around the singularity
in the higher dimensional case.
Since we see the dependence of the constants closely
in the subsection \ref{subsection;b11.9.2},
the argument for the generalization to the higher dimensional case
is elementary.
As a consequence, we see that the tame harmonic bundle is acceptable.
Thus we can apply the result in the subsection
\ref{subsection;a11.9.57}.

In the subsubsection \ref{subsubsection;b11.16.5},
we give the prolongment of $\nbigelambda$
in the case that $\lambda$ is generic.
In this case, the situation is very easy.
We see that the quasi canonical prolongment
gives the prolongment by an increasing order in that case.
Note that the direction of the argument is reverse
to those in the one dimensional case.
The results are used in the next subsections.

In the subsection \ref{subsubsection;b11.16.6},
we see the extension property of sections of $\nbigelambda$
defined over a hyperplane,
by using the result in the subsection \ref{subsection;a11.9.57}.
The argument is essentially given in our previous paper \cite{mochi}.
However it is one of the most technical parts
for the prolongment,
and hence we give some detail.
In the subsubsection \ref{subsubsection;d11.14.31},
we give estimates of Higgs fields by using the results
in the subsection \ref{subsection;a11.9.59}.
For a holomorphic section on a hyperplane,
we construct a cocycle in the subsubsection \ref{subsubsection;a11.15.60}.
By using the estimate in the subsubsection
\ref{subsubsection;d11.14.31},
we give an estimate of the cocycle.
The extension property is stated 
in the subsubsection \ref{subsubsection;10.11.80},
and it is proved in the subsubsection \ref{subsubsection;d11.14.40}.
We use the result in the subsection \ref{subsection;a11.9.57}.
In the subsubsection \ref{subsubsection;c11.16.20},
we also state the extension property in the codimension one.
Since this is the easier case, we give only an indication of the proof.

In the subsections \ref{subsection;b11.16.7}
and \ref{subsection;b11.16.8},
we show that the prolongment $\prolongg{\vecb}{\nbigelambda}$
is locally free.
As a preliminary, we show the claim under the assumption as in 
Lemma \ref{lem;9.9.95} in the subsection \ref{subsection;b11.16.7},
by using the result in the subsection \ref{subsubsection;b11.16.6}.
Then we show the claim without the assumption
in the subsection \ref{subsection;b11.16.8}.
We also see that the parabolic structures of the divisors
give the compatible tuple of the filtrations.
For that purpose,
we consider the pull back of $\nbigelambda$
via the ramified covering $\psi_{\vecc}$
for an appropriate $\vecc\in\seisuu_{>0}^l$
in the subsubsection \ref{subsubsection;c11.16.1}.
Due to the result in the subsection \ref{subsection;b11.16.7},
the prolongment of the pull back is locally free.
Moreover we have the action of the finite abelian group,
which induces the decompositions on the divisors.
Due to the result in the subsection \ref{subsection;a11.14.2},
we see that the decompositions give the splittings of
the parabolic filtrations.
In the subsubsection \ref{subsubsection;c11.16.2},
we take an equivariant frame $\vecv$ of $\prolong{\psi_{\vecc}\nbigelambda}$
which is compatible with the decompositions on the divisors.
Then we take the descent of $\vecv$.
We will see that the descent gives the frame of
the prolongment of $\nbigelambda$
by using the result in the subsubsection \ref{subsection;a11.14.2},
and we will obtain the local freeness of the prolongment.

In the subsection \ref{subsection;b11.11.30},
we see that the prolongment $\prolongg{\vecb}{\nbige}$
is locally free on $\nbigx(\lambda_0,\epsilon_0)$.
First we show the extension property of holomorphic sections.
For that purpose, we use the asymptotic orthogonality
(the subsection \ref{subsection;b11.11.30})
and the trivialization given by the argument in 
the subsection \ref{subsection;a11.9.58}.
Once the extension property is shown,
it is easy to show the local freeness.

In the subsection \ref{subsection;b11.16.10},
we see some structures induced on the divisors.
In particular, we obtain the filtrations
$\lefttop{i}\Fzero$ (the subsubsection \ref{subsubsection;10.18.20})
and the decompositions $\lefttop{i}\EEzero$
(the subsubsection \ref{subsubsection;10.18.20})
of the vector bundle
$\prolongg{\vecb}{\nbige}_{|\nbigd_i(\lambda_0,\epsilon_0)}$.
The tuples of the filtrations and the decompositions are compatible.
In particular, we obtain the induced vector bundle
$\lefttop{\lbar}\nbigg^{(\lambda_0)}_{\vecu}$
on $\nbigd_{\lbar}(\lambda_0,\epsilon_0)$.

In the subsection \ref{subsection;b11.16.11},
we see that 
$\bigl\{
\lefttop{\lbar}\nbigg^{(\lambda_0)}_{\vecu}\,\big|\,
 \lambda_0\in\cnum_{\lambda}
\bigr\}$ gives the holomorphic bundle over
$\nbigd_{\lbar}$,
and we give some detail on the vector bundle.

\subsubsection{Section \ref{section;b11.17.1}}

We give some detail of the $KMS$-structure
on the space of the multi-valued flat sections
of $\nbigelambda$ on $\Delta^{\ast\,l}\times\Delta^{n-l}$.

In the subsection \ref{subsection;b11.17.2},
we see some easy properties of
the filtrations $\lefttop{i}\nbigf$.

In the subsection \ref{subsection;b11.17.3},
we show the compatibility of the tuple of the filtrations
$\bigl(\lefttop{i}\nbigf\,\big|\,i\in \lbar\bigr)$.
The argument to deal the filtrations is complicated a little,
as usual.
However it is elementary.

Then we obtain the induced object
$\lefttop{\lbar}\nbigg^{(\lambda_0)}_{\vecu}(\nbigh)$.
In the subsection \ref{subsection;b11.17.4},
we see that the family
$\bigl\{
 \lefttop{\lbar}\nbigg^{(\lambda_0)}_{\vecu}(\nbigh)\,\big|\,
 \lambda_0\in\cnum_{\lambda}
 \bigr\}$ gives a vector bundle over $\cnum_{\lambda}^{\ast}$.
We also see some additional structures,
the nilpotent maps and the pairing.

\subsubsection{Section \ref{section;a11.17.20}}

In the subsection \ref{subsection;a11.17.21},
we see the compatibility of the naturally defined tuple
of the filtrations and the decompositions
on $\prolongg{\vecc}{\nbigelambda}$.

In the subsection \ref{subsection;a11.17.22},
we obtain the filtrations and the decompositions
of the prolongment $\prolongg{\vecc}{\nbige}$.
As a result, we obtain the induced object
$\lefttop{\lbar}\nbigg_{\vecu}(\nbige)$ on $\nbigx^{\shikaku}$.

Then we obtain the isomorphisms
$\Phi^{\can}_{\vecu}:
 \lefttop{\lbar}\nbigg_{\vecu}(\nbigh)\simeq
 \lefttop{\lbar}\nbigg_{\vecu\,|\,\cnum_{\lambda}^{\ast}}$
and
$\Phi_{\vecu,P,O}:
 \lefttop{\lbar}\nbigg_{\vecu}(\nbige)_{|\cnum_{\lambda}^{\ast}\times\{P\}}
\simeq
 \lefttop{\lbar}\nbigg_{\vecu\,|\,\cnum^{\ast}_{\lambda}}$,
which is described in the subsection \ref{subsection;a11.17.23}.

%% file: 105.2.tex

\subsubsection{Part \ref{part;a11.17.102}, Section \ref{section;a11.17.50}}

In Part \ref{part;a11.17.102},
we prove a limiting mixed twistor theorem.
As an application, we obtain the norm estimate
for holomorphic sections and multi-valued flat sections.

In the section \ref{section;a11.17.50},
we give some detail on the construction of
the vector bundle over $\proj^1$
with the nilpotent maps and the pairing,
from a tame harmonic bundle.

In the subsection \ref{subsection;a11.17.51},
we recall the variation of polarized pure twistor structure
induced by the harmonic bundles.
In particular, we recall about the conjugate 
of harmonic bundles.
The formalism was given by Simpson in \cite{s3}.

In the subsection \ref{subsection;a11.17.52},
we see the KMS-structure and the induced objects of the conjugate.
Briefly speaking, the conjugate is isomorphic to
the dual.
However we remark that the signature of the nilpotent maps
are reversed.

In the subsection \ref{subsection;a11.17.53},
we see the construction of the vector bundles
$S_{\vecu}^{\can}(E)$ and $S_{\vecu}(E,P)$.
We also see the induced nilpotent maps
and the pairings on them.
In particular,
we obtain the weight filtration $W$ on
$S_{\vecu}^{\can}(E)$ and $S_{\vecu}(E,P)$
in the case $\dim(X)=1$.

In the subsection \ref{subsection;a11.17.54},
we give some detail on the associated graded
vector bundles
$\Gr^WS_{\vecu}^{\can}(E)$ and $\Gr^WS_{\vecu}(E,P)$
in the case $\dim(X)=1$,
which are simple.


\subsubsection{Section \ref{section;a11.17.60}}

We prove the limiting mixed twistor theorem,
which will be very important in the study in Part \ref{part;a11.17.103}.
It is also useful to control the conjugacy classes of
the nilpotent parts of the residues.

In the subsection \ref{subsection;c11.14.5},
we prove the limiting mixed twistor theorem
in the case $\dim(X)=1$.
Although the proof is essentially same as
those in \cite{mochi} and \cite{sabbah},
it is rather complicated to state the argument precisely.

In the subsection \ref{subsection;a11.17.61},
we prove the limiting mixed twistor theorem
for higher dimensional case.
The different part from that in \cite{mochi}
is the proof of the constantness of the filtrations
(Lemma \ref{lem;9.16.2}).
In our previous paper \cite{mochi}, we used the argument to take a `limit'.
Instead we use the result in the subsection \ref{subsection;a11.9.54}.

When we take the associated graded objects
of $\bigl(S^{\can}(E),W,\vecN,S\bigr)$,
we obtain the nilpotent orbit
due to the limiting mixed twistor theorem.
In the subsection \ref{subsection;a11.17.62},
we derive some consequences
by using the results for polarized mixed twistor structures
(the subsection \ref{subsection;a11.11.57})
and the results in the section \ref{section;d11.12.10}.
We show the strongly sequential compatibility
and some decomposition.

\subsubsection{Section \ref{section;a11.17.70}}

As one of the application of the limiting mixed twistor theorem,
we give a norm estimate of the holomorphic sections
or flat sections of the deformed holomorphic bundles.
The arguments are essentially same as those in our previous paper
\cite{mochi}.
We have only to care the parabolic structure.

In the subsection \ref{subsection;a11.17.71},
we give some remark on functoriality for pull backs
of deformed holomorphic bundles.
In the subsection \ref{subsection;a11.17.72},
we show a preliminary norm estimate for holomorphic sections.
In the subsection \ref{subsection;a11.17.73},
we derive a norm estimate for holomorphic sections.
In the subsection \ref{subsection;a11.17.74},
we reduce the norm estimate for flat sections
to the norm estimate for holomorphic sections.
Contrast to the holomorphic case,
we do not discuss the norm estimate of the family of flat sections.
It seems that we need some additional argument.

%% file: 105.3.tex

\subsubsection{Part \ref{part;a11.17.103},
 Section \ref{section;b11.18.1}}

In Part \ref{part;a11.17.103},
we apply the results in Part
\ref{part;a11.17.101}--\ref{part;a11.17.102}
to the theory of pure twistor $D$-modules.

As a preparation,
in the section \ref{section;b11.18.1},
we recall the specialization of $\nbigr$-modules
and the sesqui-linear pairings introduced by Sabbah
with minor generalization.
The author recommends the reader to read the very readable paper
\cite{sabbah}.
We use some results in \cite{sabbah} without mention.

In the subsection \ref{subsection;b11.18.2},
we recall the specialization of $\nbigr$-module.
Although Sabbah considered only the local unitary case,
we give some more general definitions.
However, the basic theory of $\nbigr$-modules are 
established sufficiently generally in \cite{sabbah}.
Even if we need a generalization,
we need at most minor modification.

In the subsection \ref{subsection;b11.18.3},
we recall the specialization of the sesqui-linear pairings,
with some minor generalization.
Again, we need at most minor modification.

\begin{rem}
Sabbah kindly informed to the author on the revised version
{\rm\cite{sabbah2}} of his paper,
in which the generalization treated in this section
is already discussed essentially.
We keep this section for our reference in the later sections.
\hfill\qed
\end{rem}

\subsubsection{Section \ref{section;b11.24.1}}

We give the prolongments of the deformed holomorphic bundle
$\nbige$ over $\nbigx-\nbigd$ to the $\nbigr$-modules
on $\nbigx$.

First we give a naive prolongment
$\naiveprolong{\nbige}$ in the subsection \ref{subsection;b11.24.2}.
Although it is not coherent,
$\naiveprolong{\nbige}$ has many nice properties,
and it is easy to understand $\naiveprolong{\nbige}$ algebraically.
We will use $\naiveprolong{\nbige}$ as the ambient sheaf.
We introduce the sheaf $\lefttop{I}\tilde{T}^{(\lambda_0)}(\vecc,\vecd)$,
and we prepare a lemma in the subsubsection
\ref{subsubsection;c11.24.1},
which will be used in the subsection \ref{subsection;b11.24.5}.

In the subsection \ref{subsection;a11.23.40},
we give the prolongment $\gbige$,
which we really need.
To understand $\gbige$ more closely,
we introduce the filtrations $\lefttop{I}\Vzero$ $(I\subset\lbar)$,
and we obtain the sheaves  $\lefttop{I}T^{(\lambda_0)}(\vecc,\vecd)$.
In the subsection \ref{subsection;b11.24.5},
we show $\lefttop{I}T^{(\lambda_0)}(\vecc,\vecd)$
and $\lefttop{I}\tilde{T}^{(\lambda_0)}(\vecc,\vecd)$
are naturally isomorphic
(Lemma \ref{lem;9.18.11}).
Although it looks a little long and complicated,
the arguments are elementary.

In the subsection \ref{subsection;b11.24.6},
we show a kind of compatibility of the filtrations
$\lefttop{i}\Vzero(\gbige)$ $(i\in \lbar)$.
First goal is to show
$\lefttop{\lbar}\nbigv_S\bigl(\gbige\bigr)=
 \lefttop{\lbar}\nbigv_S\bigl(\naiveprolong{\nbige}\bigr)
\cap \gbige$
(Proposition \ref{prop;b11.21.5}).
Once we prove Proposition \ref{prop;b11.21.5},
we can easily translate some nice properties
of $\lefttop{i}\Vzero\bigl(\naiveprolong{\nbige}\bigr)$
to $\lefttop{i}\Vzero\bigl(\gbige\bigr)$.
In particular, we obtain the strictly $S$-decomposability
of $\gbige$ along $z_i=0$ for $i\in \lbar$.
We also obtain the primitive decomposition of sections
of $\gbige$.

In the subsection \ref{subsection;b12.6.10},
we give a characterization of $\gbige$
as the prolongment of the deformed holomorphic bundle $\nbige$
obtained from harmonic bundle.
It will be useful
when we consider the specialization of
$\gbige$ along $\prod z_i^{m_i}$.

\subsubsection{Section \ref{section;b11.24.10}}

We give some detail
on the push-forward $\gbige[\deldel_t]$
for the graph of the holomorphic functions
$\prod_{i=1} z_i^{m_i}$.
Following Saito \cite{saito2},
we introduce the filtration $\Uzero$
in the subsection \ref{subsection;b11.24.11}.
Our purpose is to show that
$\Uzero$  gives the $V$-filtration along $t=0$ at $\lambda_0$
whose associated graded module is strict,
namely,
$\gbige[\deldel_t]$ is strictly specializable along $t_0$.
We will also show that the naturally induced
filtrations $\lefttop{i}\Vzero$
on $\tildepsi_{t,u}(\gbige[\deldel_t])$
gives the $V$-filtration
along $z_i=0$ at $\lambda_0$,
whose associated graded module is strict,
namely $\tildepsi_{t,u}(\gbige[\deldel_t])$
is strictly specializable along $z_i=0$ $(i\in\lbar)$.

In the subsection \ref{subsection;b11.24.12},
we give some algorithms to describe sections of
$\gbige[\deldel_t]$ in a normal (but not unique) way.

In the subsection \ref{subsection;b11.24.13},
we obtain the primitive decompositions
of sections of $\Gr^{\Uzero}_b\gbige[\deldel_t]$.
The first goal is Proposition \ref{prop;9.18.70}.
Once Proposition \ref{prop;9.18.70} is established,
the rest are rather formal.

We obtain the strict $S$-decomposability of
$\gbige[\deldel_t]$
in the subsection \ref{subsection;b12.6.15}.
Even if it looks complicated,
it is a rather formal consequence of
the results in the subsection \ref{subsection;b11.24.13}.

In the subsection \ref{subsection;b12.6.30},
we would like to see the form of
$\tildepsi_{t,u}(\gbige[\deldel_t])$
briefly.
Since it is not easy to see it directly,
we see
$\lefttop{I}\Gr^{\Vzero}_{\vecc}
 \lefttop{J}\Vzero_{\vecd}(\tildepsi_{t,u}\gbige[\deldel_t])$,
instead.
Our goal is to relate them with
the construction in the subsubsection
\ref{subsubsection;9.13.50}.
Compare the formulas (\ref{eq;b12.6.50})
and (\ref{eq;b12.6.80}).

\subsubsection{Section \ref{section;b12.2.0}}

We consider the nilpotent map
$N=t\deldel_t+\eigenmap(\lambda,u)$
on $\tildepsi_{t,u}\gbige[\deldel_t]$.
Then we obtain the weight filtration $W(N)$.
We would like to see
$\Gr^{W(N)}\tildepsi_{t,u}\gbige[\deldel_t]$.
It is not easy to see it directly,
we introduce the filtration $\FFzero$
on $\tildepsi_{t,u}\gbige[\deldel_t]$,
such that $\Gr^{\FFzero}_m$ is a direct sum 
of the locally free sheaves $\lefttop{I}\nbigl$
on $\nbigd_I$ $(|I|=n-m)$.

In the subsection \ref{subsection;b12.2.1},
we see the relation of the weight filtration $W(N)$
and the induced filtrations $\lefttop{\lbar}\Vzero$
on $\lefttop{I}\nbigl$.
We also see the decomposition of
$P_h\Gr^{W(N)}_h$.
The results easily follow from
the limiting mixed twistor theorem
and the results in the section \ref{section;d11.12.10}.

In the subsection \ref{subsection;b12.2.2},
we introduce the filtration $\FF$.
We see that the exact sequences
associated to $\FF$ is strict
with respect to the weight filtration $W(N)$.
A key observation is given
in the proof of Lemma \ref{lem;a12.2.130}.
As a result,
we can understand the filtrations
$\lefttop{I}\Vzero$ on $\Gr^{W(N)}(\psi_{t,u}\gbige[\deldel_t])$
sufficiently well,
and we obtain the strictly specializability
of $\Gr^{W(N)}(\psi_{t,u}\gbige[\deldel_t])$
along $z_i=0$ $(i\in \nbar)$.

In the subsection \ref{subsection;b12.2.3},
we obtain the strict $S$-decomposability
of $\Gr^{W(N)}(\psi_{t,u}\gbige[\deldel_t])$
along $z_i=0$ $(i\in \nbar)$,
by using the results in the previous subsubsections
and the lemmas of Kashiwara and Saito
prepared in the subsection \ref{subsection;a11.11.58}.
As a result, we obtain the decomposition
by the supports,
as in Proposition \ref{prop;b12.6.150}.
We see some properties of the components.

\subsubsection{Section \ref{section;b12.6.100}}

In the subsection \ref{subsection;b11.24.7},
we give the sesqui-linear pairing $\gbigc$
of $\gbige$,
which is the unique prolongation of
the pairing of $\nbige$.

In the subsection \ref{subsection;b12.6.200},
we introduce $\lefttop{\nbar}\overline{\nbigg}_{\vecu}$
for $\vecu\in\KMSoverline(\nbige^0,\nbar)$,
which is a family of flat bundles on
$(X-D)\times\cnum^{\ast}$.
We see that the sesqui-linear pairing $C_0$ of $\nbige$
induces the pairing on $\lefttop{\nbar}\overline{\nbigg}_{\vecu}$.
By using it, we obtain
the sesqui-linear pairing on
$\lefttop{\nbar}\nbigg_{\vecu'}$ for $\vecu'\in\KMS(\nbige^0,\nbar)$
in the subsection \ref{subsection;b12.6.210}.
We see that it is essentially obtained
by the formalism in the subsection \ref{subsection;b11.18.3}.

In the subsection \ref{subsection;b11.24.8},
we compare two $\nbigr$-triples in $0$-dimension.
One is obtained as the specialization of the $\nbigr$-triple.
The other is obtained from
the vector bundle constructed in the section 
\ref{section;a11.17.50}.

In the subsection \ref{subsection;b12.2.4},
we show that
the component whose support is $\{0\}$ is
a polarized pure twistor structure
(Corollary \ref{cor;a11.30.50}),
by using the lemma of Kashiwara in the subsubsection
\ref{subsubsection;9.13.50}.
Once we know Corollary \ref{cor;a11.30.50},
we immediately know that 
the smooth part of the components
of $P_h\Gr^{W(N)}\psi_{t,u}\gbige[\deldel_t]$
together with the induced sesqui-linear pairing
is a variation of pure twistor structures.
By using the characterization of
the prolongment,
the each components of the decomposition
of $\tildepsi_{t,u}(\gbige[\deldel_t],\gbige[\deldel_t],\gbigc)$
is isomorphic to
a $\nbigr$-triple obtained from a tame harmonic bundle.
As a result,
we arrive at the stage
we can use an induction
to show the existence of the prolongment
as the pure twistor $D$-module.

\subsubsection{Section \ref{section;b11.24.30}}

We see the correspondence of
tame harmonic bundles and pure twistor $D$-modules
are given.
The statement is given in the subsection \ref{subsection;a11.23.50}.

In the subsection \ref{subsection;04.1.27.35},
we show that every regular pure twistor $D$-module gives
a tame harmonic bundle which is defined on a Zariski open subset
of the strict support.

The existence of the prolongment
as pure twistor $D$-module is shown in \ref{subsection;a11.23.51}.
This is a formal consequence
of the result in the subsection \ref{subsection;b12.2.4}.

The uniqueness is shown in the subsection \ref{subsection;a11.23.52}.
We need some consideration by using
the uniqueness of the intermediate extension
and the Riemann-Hilbert correspondence.

In the subsection \ref{subsection;04.2.5.1},
we see that
the correspondence of Theorem \ref{thm;b12.5.20}
preserves the pure imaginary property.
By the correspondence,
we can show that Sabbah's conjecture is true,
i.e.,
semisimple regular holonomic $D$-modules
correspond to pure imaginary pure twistor $D$-modules.
As a result, we obtain the regular holonomic version
of Kashiwara's conjecture
by combining the results of Sabbah and us.

\subsubsection{Part \ref{part;a11.17.104}, Appendix}

As is noted in the subsubsection \ref{subsubsection;b12.6.300},
the definition of pure twistor $D$-module given by Sabbah
(see \cite{sabbah} and \cite{sabbah2})
is slightly narrow from the view point of Simpson's Meta Theorem,
although his definition is natural to attack
the conjecture of Kashiwara. 
Hence we give some minor complement for pure twistor $D$-modules
as an appendix. 
In particular, we give some detail on the decomposition theorem
of the pure twistor $D$-module on a smooth projective curve.

%% file: d9.tex

\subsubsection{On the contents and the length}

Although the paper is long,
most of the contents are rather standard.
It is the intension of the author
to organize the facts for our purpose
and to give a proof when he does not know an appropriate reference.
He does not pretend that most of the contents are original.
He just hopes that everything contained in this paper
is clear for everyone.

\subsubsection{The difference of the version 1 and the version 2}

Many typos and arguments are improved.
We only refer the main difference.
\begin{itemize}
\item
Because the paper \cite{mochi3} has been written,
the discussion on  pure imaginary pure twistor $D$-module
is added.
We give the correspondence
of semisimple regular holonomic $D$-modules and
of pure imaginary pure twistor $D$-modules
through tame pure imaginary harmonic bundles.
Namely, the conjectures of Sabbah and Kashiwara are established.
The related part is mainly
the subsections \ref{subsection;04.2.5.1}--\ref{subsection;04.2.5.2},
and the subsection \ref{subsection;04.2.5.100}.

\item
The explanation for
the tameness of harmonic bundle
obtained from a pure twistor $D$-modules
is added.
The related part is mainly
the subsubsections
\ref{subsubsection;04.2.16.2}--\ref{subsubsection;04.2.16.1},
the subsection \ref{subsection;04.1.27.60},
the subsection \ref{subsection;04.1.27.61},
the subsubsection \ref{subsubsection;04.1.27.65},
and the subsection \ref{subsection;04.1.27.35}.
\end{itemize}

%% file: a89.tex

The author is grateful to the colleagues in Osaka City University.
He specially thanks Mikiya Masuda for his encouragement and supports.

The author thanks Yoshifumi Tsuchimoto and Akira Ishii
for their constant encouragements.

The author thanks Tomohide Terasoma who informed various things
to the author.
He is also grateful to Terasoma's family for their kindness.

The author thanks Masaki Kashiwara for some useful information.

The author thanks Mark Andrea de Cataldo,
who attracted the author's attention
to Hodge modules and a decomposition theorem of perverse sheaves.
The author also thanks his great tolerance for the author's lack of
communication ability.

The author thanks Claude Sabbah who kindly informed
the revised version of his paper.

The author thanks Kari Vilonen who kindly informed 
the work of D. Gaitsgory \cite{gaitsgory}.

The author prepared the paper during his stay at the Institute
for Advanced Study. The author is sincerely grateful
to their excellent hospitality.
He also acknowledges National Scientific Foundation
for a grant DMS 9729992,
although any opinions, findings and conclusions or recommendations
expressed in this material are those of the author
and do not necessarily reflect the views of the National Science Foundation.

The author thanks the financial supports by Japan Society for
the Promotion of Science and the Sumitomo Foundation.

%% file: 1.tex

\subsubsection{Sets}

We will use the following notation:

\begin{tabular}{llll}
$\seisuu$: & the set of the integers, &
$\seisuu_{>0}$: &the set of the positive integers,\\
$\rnum$: &the set of the rational numbers,&
$\rnum_{>0}$: &the set of the positive rational numbers,\\
$\real$: & the set of the real numbers, &
$\real_{>0}$: &the set of the positive real numbers,\\
$\cnum$: &the set of the complex numbers, &
$\nbar$: & the set $\{1,2,\ldots,n\}$,\\
$M(r)$: & the set of $r\times r$-matrices,&
$\nbigh_r$: & the set of $r\times r$-hermitian matrices,\\
$\gbigs_l$: & the $l$-th symmetric group,\\

\end{tabular}

We denote the set of positive hermitian metric
of $V$ by $\PH(V)$.
We often identify it with the set of the positive hermitian matrices
by taking an appropriate base of $V$.

\subsubsection{A disc, a punctured disc and some products}

For any positive number $C>0$ and $z_0\in\cnum$,
the open disc $\bigl\{z\in\cnum\,\big|\,|z-z_0|<C\bigr\}$
is denoted by $\Delta(z_0,C)$,
and the punctured disc $\Delta(z_0,C)-\{z_0\}$
is denoted by $\Delta^{\ast}(z_0,C)$.
When $z_0=0$, $\Delta(0,C)$ and $\Delta^{\ast}(0,C)$
are often denoted by $\Delta(C)$ and $\Delta^{\ast}(C)$.
Moreover, if $C=1$, $\Delta(1)$ and $\Delta^{\ast}(1)$
are often denoted by $\Delta$ and $\Delta^{\ast}$.
If we emphasize the variable,
we describe as $\Delta_z$, $\Delta_i$.
For example,
$\Delta_z\times \Delta_w=\{(z,w)\in\Delta\times \Delta\}$,
and $\Delta_1\times \Delta_2=\{(z_1,z_2)\in\Delta\times \Delta\}$.
We often use the notation $\cnum_{\lambda}$ and $\cnum_{\mu}$
to denote the complex planes
$\bigl\{\lambda\in\cnum\bigr\}$
and $\bigl\{\mu\in\cnum\bigr\}$.

Unfortunately, the notation $\Delta$ is also used to denote the
Laplacian. The author hopes that there will be no confusion.

We put $\Deltabar(C):=\bigl\{z\,\big|\,|z|\leq C\bigr\}$
and $\Deltabarast(C):=\bigl\{z\,\big|\,0<|z|\leq C\bigr\}$.
In the case $C=1$,
we use $\Deltabar$ and $\Deltabarast$
instead of $\Deltabar(1)$ and $\Deltabarast(1)$.

For a complex manifold $X$, a point $\lambda_0\in\cnum_{\lambda}$
and a positive number $\epsilon_0$,
we often consider the product
$X\times\Delta(\lambda_0,\epsilon_0)$
and $X\times\Delta^{\ast}(\lambda_0,\epsilon_0)$.
For simplicity, we denote them by
$\nbigx(\lambda_0,\epsilon_0)$ and $\nbigx^{\ast}(\lambda_0,\epsilon_0)$
respectively.
We also use the notation $\nbigx$ and $\nbigx^{\shikaku}$
to denote the $X\times\cnum_{\lambda}$ and $X\times
\cnum_{\lambda}^{\ast}$.

For a complex manifold $X$, we have the conjugate complex manifold,
which is denoted by $X^{\dagger}$.


\subsubsection{Projections}

Let $I$ be a finite set and $J$ be a subset of $I$.
In general, $q_J:X^I\lrarr X^J$ denotes the naturally induced projection
taking the $i$-th components for $i\in J$.
Similarly, $\pi_J:X^I\lrarr X^{I-J}$
denotes the projection omitting the $j$-th component for $j\in J$.
However, we will often use $\pi$ to denote some other projections.
If $J$ consists of the unique element $j$,
We often use the notation $q_j$ and $\pi_j$
instead of $q_{\{j\}}$ and $\pi_{\{j\}}$.

%% file: a87.tex

\subsubsection{The order on $\real^n$}

We have the natural order on $\real$.
Let $n$ be a positive integer.
We often use the order on $\real^n$
given as follows:
For elements $\veca,\vecb\in\real^n$,
we say $\veca\leq \vecb$ 
if and only if
$q_i(\veca)\leq q_i(\vecb)$ for any $i\in\nbar$.
When we consider such order,
we say $\veca<\vecb$ if and only if
$q_i(\veca)< q_i(\vecb)$ for any $i\in\nbar$,
and we say $\veca\lneq\vecb$ if and only if
$\veca\leq\vecb$ and $\veca\neq\vecb$.

%% file: b21.tex

\subsubsection{$\kappa_c$ and $\nu_c$ }
\label{subsubsection;b11.11.10}

Let $c$ be a real number.
The maps $\kappa_c:\real\lrarr ]c-1,c]$
and $\nu_c:\real\lrarr\seisuu$ are defined by the following condition:
\begin{quote}
For a real number $x\in \real$,
the equality $\kappa_c(x)+\nu_c(x)=x$ holds.
\end{quote}
In the case $c=0$,
we use the notation $\kappa$ and $\nu$
instead of $\kappa_c$ and $\nu_c$.

Let $S$ be a finite subset of $]-1,0]$
and $b$ be a positive number.
We obtain the map
$\phi_b:S\lrarr \seisuu$ given by
$\phi_b(x)=\nu(b\cdot x)$.
\begin{df} \label{df;9.5.11}
A real number $b$ is sufficiently large with respect to $S$,
if the map $\phi_b$ is injective.
\hfill\qed
\end{df}

\subsubsection{$\paramap$, $\eigenmap$ and $\kmsmap$}
\label{subsubsection;b11.11.11}

For any element $u=(a,\alpha)\in \real\times\cnum$,
we put as follows:
\[
 \begin{array}{l}
 \paramap(\lambda,u):=
 a+2\cdot\Realpart(\lambda\cdot\bar{\alpha}),\\
 \mbox{{}}\\
 \eigenmap(\lambda,u):=
 \alpha-a\cdot\lambda-\bar{\alpha}\cdot\lambda^2.
 \end{array}
\]
Then we obtain the following morphism:
\[
 \kmsmap(\lambda)=\bigl(\paramap(\lambda),\eigenmap(\lambda)\bigr):
 \real\times\cnum\lrarr\real\times\cnum.
\]
The following lemma is checked by a direct calculation.
\begin{lem} \label{lem;10.11.25}
$\kmsmap(\lambda)$ is bijective.
\end{lem}
\pf
Let us consider the equation
$\paramap\bigl(\lambda,(a,\alpha)\bigr)=A$
and $\eigenmap\bigl(\lambda,u\bigr)=B$
for $(B,A)\in\cnum\times\real$.
Then we have the unique solution:
\[
 \alpha=\frac{\lambda\cdot A+B}{|\lambda|^2+1},
\quad
 a=\frac{(-|\lambda|^2+1)A-2\Realpart(\lambda\overline{B})}
 {|\lambda|^2+1}.
\]
It may be useful to use the relation
$\lambda\cdot\paramap(\lambda,(a,\alpha))+\eigenmap(\lambda,(a,\alpha))
=(|\lambda|^2+1)\cdot\alpha$.
\hfill\qed

\begin{rem}
We will use Lemma {\rm\ref{lem;10.11.25}}
to control the KMS-structure of the prolongment
of deformed holomorphic bundles.
\hfill\qed
\end{rem}

For any element $u=(\alpha,a)\in\cnum\times\real$,
we put $u^{\dagger}=(\overline{\alpha},-a)\in\cnum\times\real$.

\begin{lem} \label{lem;10.16.6}
Let $u$  be an element of $\cnum\times\real$.
We have the following formula:
\begin{equation} \label{eq;10.18.25}
-\paramap(\lambda,u)
=\paramap\bigl(-\overline{\lambda},u^{\dagger}\bigr),
\quad\quad
 \overline{\eigenmap(\lambda,u)}
=\eigenmap\bigl(-\overline{\lambda},u^{\dagger}\bigr).
\end{equation}
\end{lem}
\pf
It can be checked by a direct calculation.
\hfill\qed

\vspace{.1in}
For any element $u=(a,\alpha)\in\real\times\cnum$,
we put as follows:
\begin{equation}\label{eq;a11.15.1}
\begin{array}{l}
 \paramap^f(\lambda,u):=
  \Realpart\bigl(\lambda\cdot\overline{\alpha}+\lambda^{-1}\cdot\alpha\bigr)
 =\paramap(\lambda,u)
 +\Realpart\bigl(\lambda^{-1}\cdot\eigenmap(\lambda,u)\bigr),
 \\
 \mbox{{}}\\
 \eigenmap^f(\lambda,u):=
 \exp\bigl(
 -2\pi\sqrt{-1}\cdot(\lambda^{-1}\cdot\alpha-a-\lambda\cdot\overline{\alpha})
 \bigr)
=\exp\bigl(
 -2\pi\sqrt{-1}\lambda^{-1}\cdot\eigenmap(\lambda,u)
 \bigr).
\end{array}
\end{equation}

\begin{lem}\label{lem;a12.1.10}
Let $u$ be an element of $\real\times\cnum$.
We have the following formula:
\[
 \eigenmap^f(\lambda,u)=\eigenmap^f\bigl(\lambda^{-1},u^{\dagger}\bigr)^{-1},
\quad
 \paramap^f(\lambda,u)=\paramap^f\bigl(\lambda^{-1},u^{\dagger}\bigr).
\]
\end{lem}
\pf
We have the following equalities:
\[
 \eigenmap^{f}(\lambda,u)
=\exp\bigl(-2\pi\sqrt{-1}(\lambda^{-1}\alpha-a-\lambda\bar{\alpha})\bigr)
=\eigenmap^f\bigl(\lambda^{-1},u^{\dagger}\bigr)^{-1}.
\]
We also have the following equalities:
\[
 \paramap^f(\lambda,u)
=\Realpart\bigl(\lambda\bar{\alpha}+\lambda^{-1}\alpha\bigr)
=\paramap^f\bigl(\lambda^{-1},u^{\dagger}\bigr).
\]
Thus we are done.
\hfill\qed

\vspace{.1in}

\noindent
{\bf Notation}
To denote the element $(1,0)\in\real\times\cnum$,
we often use the notation $\vecdelta_{0}$.
We often identify $(\real\times\cnum)^I$
and $\real^I\times\cnum^I$.

\vspace{.1in}

\label{subsubsection;10.16.5}
For any element $\vecu=(\vecalpha,\vecb)\in\cnum^l\times\real^l$,
we put $\vecu^{\dagger}:=(\overline{\vecalpha},-\vecb)$.

%% file: b20.tex

\subsubsection{Notation}

Let $X$ be a complex manifold
and $D=\bigcup_{i=1}^N D_i$ be a normal crossing divisor.
Let $E$ be a holomorphic vector bundle with
a hermitian metric $h$ over $X-D$.

Let $U$ be an open subset of $X$, which is admissible
with respect to $D$,
i.e., we have a coordinate $(z_1,\ldots,z_n)$ satisfying the following:
\[
   D\cap U=\bigcup_{k=1}^l(D_{i_k}\cap U),
 \quad D_{i_k}=\{z_k=0\}.
\]

For any section $f\in\Gamma\bigl(U\cap(X-D),E\bigr)$,
let $|f|_h$ denote the norm function of $f$ with respect to 
the metric $h$.
We describe $|f|_h=O\bigl(\prod_{i=1}^l|z_k|^{-b_k}\bigr)$,
if there exists a positive number $C$
such that $|f|_h\leq C\cdot \prod_{k=1}^l|z_k|^{-b_k}$.

Recall that `$-\ord(f)\leq \vecb$' means the following:
\[
 |f|_h=O\Bigl(\prod_{k=1}^l|z_k|^{-b_k-\epsilon}\Bigr)
\mbox{ for any positive number }\epsilon.
\]

For any $\vecb\in\real^N$,
the sheaf $\prolongg{\vecb}E$ is defined as follows:
\[
 \Gamma(U,\prolongg{\vecb}{E}):=
 \big\{f\in \Gamma\bigl(U\cap(X-D),E\bigr)\,\big|\,
 -\ord(f)\leq \vecb
 \big\}.
\]
The sheaf $\prolongg{\vecb}E$ is called the prolongment of $E$
by an increasing order $\vecb$.
In particular, we use the notation $\prolong{E}$
in the case $\vecb=(0,\ldots,0)$.

\subsubsection{Adaptedness and adaptedness up to log order}

Let $X$ be a $C^{\infty}$-manifold,
and $E$ be a $C^{\infty}$-vector bundle with a hermitian metric $h$.
Let $\vecv=(v_1,\ldots,v_r)$ be a $C^{\infty}$-frame of $E$.
We obtain the $\nbigh(r)$-valued function $H(h,\vecv)$,
whose $(i,j)$-component is given by $h(v_i,v_j)$.
Recall that the frame $\vecv$ is called adapted,
if $H(h,\vecv)$ and $H(h,\vecv)^{-1}$ are bounded.

Let $E=\bigoplus_i E_i$ be a $C^{\infty}$-decomposition of $E$.
The hermitian metric $h$ of $E$ induces the metric $h_i$ on $E_i$.
Then we obtain the metric $\bigoplus_i h_i$ of $E$.
\begin{df} \label{df;10.11.10}
The decomposition $E=\bigoplus_i E_i$ is quasi adapted
with respect to $h$,
if $h$ and $\bigoplus_i h_i$ are mutually bounded.
\hfill\qed
\end{df}

Let us consider the case $X=\Delta^{\ast\,l}\times \Delta^{n-l}$.
We have the coordinate $(z_1,\ldots,z_n)$.
Let $E$, $h$ and $\vecv$ be as above.
\begin{df}
A frame $\vecv$ is called adapted up to log order,
if the following inequalities hold over $X$,
for some positive numbers $C_i$ $(i=1,2)$ and $M$:
\[
0<
 C_1\cdot \Bigl(
 -\sum_{i=1}^l \log|z_i|
 \Bigr)^{-M}
\leq H(h,\vecv)
\leq C_2\cdot\Bigl(
 -\sum_{i=1}^l \log|z_i|
 \Bigr)^M.
\]
\hfill\qed
\end{df}

\subsubsection{Lemmas for local freeness and the parabolic filtration}

We put $X:=\Delta^n$
and $D:=\bigcup_{i=1}^l D_i$,
where $D_i=\bigl\{(z_1,\ldots,z_n)\in X\,\big|\,z_i=0\bigr\}$.
Let $E$ be a holomorphic bundle over $X-D$
and $h$ be a hermitian metric of $E$.

\begin{lem} \label{lem;9.5.10}
We assume that
we have a holomorphic frame
$\vecv=\bigl(v_j\,\big|\,j=1,\ldots,\rank(E)\bigr)$
of $E$ satisfying the following conditions:
\begin{itemize}
\item
 There exists $b_i(v_j)\in ]-1,0]$
 for $1\leq j\leq \rank(E)$ and for $1\leq i\leq l$.
\item
 The $C^{\infty}$-frame
 $\vecv'=\bigl(v_j'\,\big|\,j=1,\ldots,\rank(E)\bigr)$,
given as follows, is adapted up to log order.
 \[
  v_j':=v_j\cdot\prod_{i=1}^l |z_i|^{b_i(v_j)}.
 \]
\end{itemize}
Then the following holds:
\begin{enumerate}
\item
The $\nbigo_X$-sheaf $\prolong{E}$ is locally free.
\item
Each $v_i$ is a section of $\prolong{E}$,
and $\vecv$ gives a frame of $\prolong{E}$.
\end{enumerate}
\end{lem}
\pf
It is easy to see that $v_j$ are sections of $\prolong{E}$.
Let $f$ be a section of $\prolong{E}$ over 
an open subset of $U\subset X$.
We have the following development on $U$:
\[
 f=\sum f_j\cdot v_j
  =\sum f_j\cdot \prod_{i=1}^l|z_i|^{-b_i(v_j)}\cdot v_j'.
\]
Here $f_j$ are holomorphic functions on $U\cap (X-D)$.
Since $\vecv'$ is adapted up to log order,
we have the following estimate:
\[
 \Big|f_j\cdot \prod_{i=1}^l|z_i|^{-b_i(v_j)}\Big|
 =O\Bigl(
 \prod_{i=1}^{l}|z_i|^{-\epsilon}
 \Bigr),
\quad
\mbox{for any positive number }\epsilon.
\]
In fact, the left hand side can be dominated
by a polynomial of $-\log |z_i|$ $(i=1,\ldots,l)$.
Hence we obtain the following estimate, for any $\epsilon>0$,
\[
 |f_j|=O\Bigl(
 \prod_{i=1}^l |z_i|^{b_i(v_j)-\epsilon}\Bigr).
\]
Note that $b_i(v_j)-\epsilon>-1$ for any sufficiently small 
positive number $\epsilon$.
Hence the functions $f_j$ are holomorphic on $U$.
Thus $\vecv$ is a frame of $\prolong{E}$,
and $\prolong{E}$ is locally free.
\hfill\qed

\vspace{.1in}

Let $\vecdelta_i$ denote the element
$(\overbrace{0,\ldots,0}^{i-1},1,0,\ldots,0)$.
For any real number $b\leq 0$,
we have the natural morphism
$\prolongg{b\vecdelta_i}{E}\lrarr \prolong{E}$.
Thus we have the parabolic filtration $\lefttop{i}F$ 
on $D_i$ given as follows:
\[
 \lefttop{i}F_b:=\Image\bigl(
 \prolongg{b\vecdelta_i}{E}_{|D_i}\lrarr E_{|D_i}
 \bigr).
\]

\begin{lem} \label{lem;10.10.31}
We impose the same assumption in Lemma {\rm \ref{lem;9.5.10}}.
Then the following holds.
\begin{itemize}
\item
For each $i$, $\lefttop{i}F$ is a filtration
in the category of vector bundles on $D_i$,
namely the associated graded sheaf $\lefttop{i}\Gr^F$
is locally free on $D_i$.
\item
The tuple of the filtrations
$\bigl(\lefttop{i}F\,\big|\,i=1,\ldots,l\bigr)$ is compatible,
in the sense of Definition {\rm\ref{df;10.11.90}}.
\item
 We have
 $\lefttop{i}\deg^F(v_j)=b_i(v_j)$.
\end{itemize}
\end{lem}
\pf
Let $f$ be a section of $\prolongg{b\cdot\vecdelta_i}{E}$,
we have the description
$f=\sum f_j\cdot v_j$ for some holomorphic functions $f_j$.
By an argument similar to the proof of Lemma \ref{lem;9.5.10},
we obtain the vanishing
$ f_{j\,|\,D_i}=0$ in the case $b_i(v_j)>b$.

Thus $\lefttop{i}F_b$ is the vector subbundle
of $\prolong{E}_{|D_i}$ generated by
$v_i$ such that $b_i(v_j)\leq b$.
It implies all the claims.
\hfill\qed

%% file: a50.2.tex

\subsubsection{The action of the group $\mu_c$}

For any positive integer $c$,
we put $\mu_c:=\bigl\{z\in\cnum\,\big|\,z^c=1\bigr\}$.
We pick a generator $\omega(c)$ of $\mu_c$.
For any element $\vecc=(c_1,\ldots,c_n)\in\seisuu_{>0}^m$,
we put $\mu_{\vecc}:=\prod_{i=1}^m\mu_{c_i}$.
We denote the element
$\bigl(\overbrace{1,\ldots,1}^{i-1},\omega(c_i),1,\ldots,1\bigr)$
by $\omega(c_i)$ for simplicity.
We have the natural inclusion
$\mu_{c_i}\lrarr\mu_{\vecc}$.
The image is also denoted by $\mu_{c_i}$

We put $X=\Delta^n$, $D_i:=\{z_i=0\}$ and $D=\bigcup_{i=1}^n D_i$
for some $l\leq n$.
We have the natural $G_m^n$-action on $X$,
given by the componentwise multiplication.
Let $\vecc$ be an element of $\seisuu_{>0}^m$.
If we take a homomorphism $\rho:\mu_{\vecc}\lrarr G_m^n$,
we obtain the $\mu_{\vecc}$-action $\rho$ on $X$.
In the following,
we consider only such $\mu_{\vecc}$-actions on $X$.

\subsubsection{An equivariant section and an equivariant lift}

Let $\rho$ be a $\mu_{\vecc}$-action on $X$.
Let $E$ be a $\rho$-equivariant holomorphic vector bundle on $X$.
\begin{df}
A section $f$ of $E$ is called $\rho$-equivariant,
if there exists a homomorphism
$\chi:\mu_{\vecc}\lrarr\cnum^{\ast}$
such that
$g^{\ast}(f)=\chi(g)\cdot f$ for any element $g\in\mu_{\vecc}$.
\hfill\qed
\end{df}

Let $\Gamma(X,E)$ denote the space of holomorphic sections of $E$
over $X$.
We have the natural $\mu_{\vecc}$-action on $\Gamma(X,E)$.
Since $\mu_{\vecc}$ is a finite group,
we have the canonical decomposition:
\[
  \Gamma(X,E)\simeq\bigoplus_{\eta\in Rep(\mu_{\vecc})}
 Hom\bigl(V_{\eta},\Gamma(X,E)\bigr)\otimes V_{\eta}.
\]
Here $\Rep(\mu_{\vecc})$ denote the set of the equivalence classes
of the irreducible representations of $\mu_{\vecc}$,
and $V_{\eta}$ $(\eta\in\Rep(\mu_{\vecc}))$ denotes
an irreducible representation corresponding to $\eta$.
Then a section $f\in\Gamma(X,E)$ is equivariant if and only if
$f$ is contained in one of the components
in the canonical decomposition.

Let $I$ be a subset of $\nbar$.
We put $D_I:=\bigcap_{i\in I}D_i$.
\begin{lem} \label{lem;10.10.1}
Let $f_0$ be a holomorphic equivariant section of $E_{|D_I}$,
i.e.,
there exists a homomorphism $\chi:\mu_{\vecc}\lrarr\cnum^{\ast}$
such that $g^{\ast}(f_0)=\chi(g)\cdot f_0$ for any $g\in\mu_{\vecc}$.
Then there exists a holomorphic equivariant section $f$
of $E$ on $X$ satisfying
$f_{|D_I}=f_0$ and $g^{\ast}(f)=\chi(g)\cdot f$.
\end{lem}
\pf
We have the equivariant surjection
$\Gamma(X,E)\lrarr \Gamma(D_I,E_{|D_I})$.
Since $\mu_{\vecc}$ is finite,
we have the canonical decompositions
of $\Gamma(X,E)$ and $\Gamma(D_I,E_{|D_I})$.
Then we obtain the surjections
of the components of the canonical decompositions.
Hence we are done.
\hfill\qed

\vspace{.1in}
Let $\vecu=(u_i)$ be an equivariant base of $E_{|O}$,
i.e.,
there exist $\chi_i:G\lrarr\cnum^{\ast}$
such that $g^{\ast}(u_i)=\chi_i(g)\cdot u_i$
for each $i$.
\begin{cor} \label{cor;10.7.10}
There exists an equivariant frame $\vecv$ of $E$
on a neighbourhood of $O$,
such that $g^{\ast}(v_i)=\chi_i(g)\cdot v_i$
and $v_{i\,|\,O}=u_i$.
\hfill\qed
\end{cor}

\begin{df}
A frame as in Corollary {\rm\ref{cor;10.7.10}}
is called an equivariant frame.
\hfill\qed
\end{df}

\begin{lem} \label{lem;10.7.15}
Let $E_j$ $(j=1,2)$ be $\mu_{\vecc}$-vector bundles over $X$,
and $\pi:E_1\lrarr E_2$ be equivariant surjection.
Let $N$ be a $\mu_{\vecc}$-subbundle of $E_2$
and $M$ be a $\mu_{\vecc}$-subbundle of $E_{1\,|\,D}$.
Assume that the restriction of $\pi_{|D}$ to $M$
gives an isomorphism of $M$ and $N_{|D}$.

Then there exists a $\mu_{\vecc}$-subbundle $\widetilde{M}$ of $E_1$
defined around $O$
satisfying the following:
\begin{itemize}
\item
The restriction of $\pi_{|D}$ to $\widetilde{M}$
gives an isomorphism of $\widetilde{M}$ and $N$.
\end{itemize}
Such $\widetilde{M}$ is called an equivariant lift of $N$
extending $M$.
\end{lem}
\pf
We have the following exact sequence:
\begin{equation}\label{eq;8.18.1}
 \Gamma(E_1)
\lrarr \Gamma(E_2)\oplus\Gamma(E_{1\,|\,D})
\lrarr \Gamma(E_{2\,|\,D})\lrarr 0.
\end{equation}
Since $\mu_{\vecc}$ is finite,
$\mu_{\vecc}$-vector spaces in the complex (\ref{eq;8.18.1})
have the canonical decomposition.
On each component of the canonical decomposition,
the complex is exact.

Let take an equivariant frame $\vecv_N$ of $N$
on a neighbourhood of $O$.
By using an isomorphism $\pi_{|D}:M\lrarr N_{|D}$,
we take an equivariant frame $\vecv_M$ of $M$
around $O$
such that $\pi_{|D}(\vecv_M)=\vecv_{N|D}$.
By using the equivariant exact sequence (\ref{eq;8.18.1}),
we can take an equivariant sections $v_{\widetilde{M},i}$ of $\widetilde{M}$
around $O$ such that
$\pi(v_{\widetilde{M},i})=v_{N,i}$
and $v_{\widetilde{M},i\,|\,D}=v_{M,i}$.
Then $\vecv_{\widetilde{M}}:=(v_{\widetilde{M},i})$ gives
a $\mu_{\vecc}$-subbundle $\widetilde{M}$ of $E_1$,
which has desired properties.
\hfill\qed

\subsubsection{Compatible frames}

We put $X=\Delta^n$, $D_i=\{z_i=0\}$ and $D=\bigcup_{i=1}^n D_i$.
Let $V$ be a $\mu_{\vecc}$-vector bundle on $X$.
Assume that we are given $\mu_{\vecc}$-subbundles
$H_{I}\subset V_{|D_I}$ for any subset $I\subset\nbar$,
satisfying
$H_{I'\,|\, D_I}\supset H_I$ $(I'\subset I)$.

\begin{lem} \label{lem;10.10.2}
Let $v$ be a $\mu_{\vecc}$-equivariant section of $H_{\lbar}$.
Then there exists a $\mu_{\vecc}$-equivariant section $\tilde{v}$
of $V$ on a neighbourhood of $D_{\nbar}$,
satisfying
$\tilde{v}_{|D_I}\in H_{I}$ and
$\tilde{v}_{|D_{\lbar}}=v$.
\end{lem}
\pf
We construct $\tilde{v}_{|D_I}$ descending inductively on $|I|$.
Assume that we have already took $\tilde{v}_{|D_I'}$
for any $I'\supsetneq I$.
We put $\del D_I:=\bigcup_{I'\supsetneq} D_{I'}$.
Then $\tilde{v}_{|\del D_I}$ is an equivariant section of $H_{I\,|\,\del D_I}$.
By an argument  similar to the proof of Lemma \ref{lem;10.10.1},
we can extend it to an equivariant section
$\tilde{v}_{|D_I}$ of $H_I$.
Thus the inductive construction can proceed.
\hfill\qed

\vspace{.1in}

Let $S$ be a set. For any subset $I\subset\nbar$,
we have the set $S^I:=\{f:I\lrarr S\}$.
For any pair $I\subset I'$, we have the naturally defined projection
$q^I:S^{I'}\lrarr S^I$.
Let us consider decompositions
$V=\bigoplus_{u\in S^I}\lefttop{I}U_s$.
We denote the tuple $\bigl(\lefttop{I}U_s\,\big|\,s\in S^I\bigr)$
by $\lefttop{I}U$.
We assume the following, for any subset $I\subset\nbar$
and for any element $u\in S^I$:
\[
  \lefttop{I}U_{s\,|\,D_{I'}}
=\bigoplus_{q_I(u')=u}\lefttop{I'}U_{u'}.
\]

\begin{lem}\label{lem;c11.12.2}
Let $\vecv$ be an equivariant frame of $V_{|D_{\lbar}}$
compatible with $\lefttop{\nbar}U$.
Then we have an equivariant frame $\tilde{\vecv}$ of $V$
on a neighbourhood of $D_{\nbar}$,
satisfying that $\tilde{\vecv}_{|D_{\nbar}}=\vecv$
and that $\tilde{\vecv}_{|D_I}$ is compatible with
$\lefttop{I}U$.
\end{lem}
\pf
From any element $u\in S^{\lbar}$,
we obtain the elements $q_I(u)\in S^I$.
Then we obtain the subbundles
$\lefttop{I}U_{q_I(u)}\subset V_{|D_I}$.
Hence we can take a tuple of sections $\tilde{\vecv}$
satisfying that $\tilde{\vecv}_{|D_{\nbar}}=\vecv$
and that $\tilde{\vecv}_{|D_I}$ is compatible with
$\lefttop{I}U$,
by using Lemma \ref{lem;10.10.2}.
On a neighbourhood of $D_{\nbar}$,
$\tilde{\vecv}$ gives a frame.
\hfill\qed

%% file: b15.tex

\subsubsection{Preliminary}

\label{subsubsection;10.26.95}

Let $T^n$ be an $n$-dimensional torus,
i.e.,
$T^{n}=\bigl\{(z_1,\ldots,z_n)\in\cnum^n\,\big|\,|z_i|=1\bigr\}$.
We use the coordinate $z_i=\exp\bigl(\sqrt{-1}\theta_i\bigr)$.
Let $\vecalpha=(\alpha_1,\ldots,\alpha_n)$ be an
element of $\real^n$
such that $\alpha_1,\ldots,\alpha_n$ are linearly independent
over $\rnum$.
Let us consider the morphism $\psi_{\vecalpha}:\real\lrarr T^n$,
given  as follows:
\[
 \psi_{\vecalpha}(x)=
 \bigl(
 \exp\bigl(\sqrt{-1}\alpha_i\cdot x\bigr)
 \,\big|\, i=1,\ldots,n \bigr).
\]

Let $f$ be an $\real$-valued function on $T^n$.
Then we have the Fourier decomposition of $f$:
\[
 f=\sum_{\vecm\in \seisuu^n} a_{\vecm}\cdot
 \exp\Bigl(
 \sqrt{-1}\cdot\sum_{i=1}^n m_i\cdot\theta_i
 \Bigr).
\]
Here $m_i$ denotes the $i$-th component of $\vecm$,
and $a_{\vecm}$ are complex numbers.
Since $f$ is $\real$-valued, we have the relation
$\overline{a}_{\vecm}=a_{-\vecm}$.
Thus we have the following equality:
\[
 f=\sum_{\vecm\in\seisuu^n}
 \frac{1}{2}
 \Bigl(
 a_{\vecv}\cdot\exp\Bigl(
 \sqrt{-1}\sum m_i\cdot\theta_i
 \Bigr)
+\overline{a}_{\vecm}\cdot
 \exp\Bigl(
 -\sqrt{-1}\sum m_i\cdot\theta_i
 \Bigr)
 \Bigr).
\]
In the following, we put
$\vecm\cdot \vecalpha:=\sum_{i=1}^n m_i\cdot \alpha_i$.
Then we have the following:
\[
 \psi_{\vecalpha}^{-1}(f)(x)
=\sum_{\vecm\in\seisuu^n}
 \frac{1}{2}
 \Bigl(
  a_{\vecm}\cdot\exp\bigl(
 \sqrt{-1}\vecm\cdot\vecalpha\cdot x
 \bigr)
+\overline{a}_{\vecm}
 \exp\bigl(
 -\sqrt{-1}\vecm\cdot\vecalpha\cdot x
 \bigr)
 \Bigr).
\]
Let $b_{\vecm}$ and $c_{\vecm}$ denote the real part and
the imaginary part of $a_{\vecm}$ respectively.
Then we obtain the following:
\[
 \psi_{\vecalpha}^{-1}(f)(x)
=\sum_{\vecm\in\seisuu^n}
 b_{\vecm}\cdot \cos\bigl(
 \vecm\cdot\vecalpha\cdot x
 \bigr)
-\sum_{\vecm\in\seisuu_n}
 c_{\vecm}\cdot\sin\bigl(
 \vecm\cdot\vecalpha\cdot x
 \bigr).
\]

\begin{lem}
We have the finiteness
$\sum_{\vecm}|a_{\vecm}|<\infty$.
\end{lem}
\pf
We put $\seisuu^{\ast}:=\seisuu-\{0\}$.
We have the natural inclusion
$\seisuu^{\ast\,I}\lrarr \seisuu^n$.
Let $T_I$ be the sub-torus of $T^n$
determined by the condition
$z_j=1$ for any $j\in\nbar-I$.
Then the restriction $f_{|T_I}$ is $C^{\infty}$.
Hence we have the following:
\[
 \sum _{\vecm\in \seisuu^{\ast\,I}}
 \prod_{i\in I}|m_i|^2\cdot |a_{\vecm}|^2<\infty.
\]
Thus we obtain the following:
\[
 \sum_{\vecm\in\seisuu^{\ast\,I}}|a_{\vecm}|
<
 \Bigl(
 \sum_{\vecm\in\seisuu^{\ast\,I}}
 \prod_{i\in I}|m_i|^{-2}
 \Bigr)
\cdot
 \Bigl(
 \sum_{\vecm\in\seisuu^{\ast\,I}}
 \prod_{i\in I}|m_i|^{2}
\cdot
 |a_{\vecm}|^2
 \Bigr)
<\infty.
\]
On the other hand, we have the following:
\[
 \sum_{\vecm\in\seisuu^n}|a_{\vecm}|
=\sum_{I\subset\nbar}
 \sum_{\vecm\in\seisuu^{\ast\,I}}
 |a_{\vecm}|<\infty.
\]
Thus we are done.
\hfill\qed

\begin{cor} \label{cor;10.2.6}
We have $\sum_{\vecm\in\seisuu^n}|b_{\vecm}|<\infty$
and $\sum_{\vecm\in\seisuu^n}|c_{\vecm}|<\infty$.
\hfill\qed
\end{cor}

\begin{lem}
We have the following equalities:
\[
 \int^N_{-N} 1\cdot\frac{dx}{2N}=1,
\quad\quad
 \int^N_{-N}\cos(\vecm\cdot\vecalpha\cdot x)\cdot
 \frac{dx}{2N}
=\frac{\sin (\vecm\cdot\vecalpha\cdot N)}{N\cdot(\vecm\cdot\vecalpha)},
\quad\quad
 \int^N_{-N}
 \sin (\vecm\cdot\vecalpha\cdot x)\cdot\frac{dx}{2N}=0.
\]
\end{lem}
\pf
It can be checked by direct calculations.
\hfill\qed

\vspace{.1in}

We put as follows:
\[
 \Xi_N(f):=
 \int^{N}_{-N}
 \psi_{\vecalpha}^{-1}(f)\cdot \frac{dx}{2N}.
\]

\begin{lem}
We have the following equality:
\[
 \Xi_N(f)=
 b_{0}+
 \sum_{\vecm\neq 0}
 \frac{\sin(\vecm\cdot\vecalpha\cdot N)}{N\cdot(\vecm\cdot\vecalpha)}
\cdot b_{\vecm}.
\]
\end{lem}
\pf
We have the following:
\[
 \psi_{\vecalpha}^{-1}(f)=
 \lim_{M\to\infty}
 \Bigl(
 \sum_{|\vecm|\leq M}
 b_{\vecm}\cdot\cos(\vecm\cdot\vecalpha\cdot x)
-
 \sum_{|\vecm|\leq M}
 c_{\vecm}\cdot\sin(\vecm\cdot\vecalpha\cdot x)
 \Bigr).
\]
Due to Corollary \ref{cor;10.2.6},
we can change the order of the integral and the summation.
Thus we obtain the following:
\[
 \Xi_N(f)=
 \sum_{\vecm\in\seisuu^n}
 \Bigl(
 \int^N_{-N}
 \cos(\vecm\cdot\vecalpha\cdot x)\frac{dx}{2N}\cdot b_{\vecm}
-\int^N_{-N}
 \sin(\vecm\cdot\vecalpha\cdot x)\frac{dx}{2N}\cdot c_{\vecm}
 \Bigr).
\]
Thus we obtain the result.
\hfill\qed

\begin{lem}
When $N\to \infty$,
the sequence of the numbers $\{\Xi_N(f)\}$ is convergent.
We put $\Xi(f):=\lim_{N\to\infty}\Xi_N(f)$.
Then we have $\Xi(f)=b_0$.
\end{lem}
\pf
We have only to show the following:
\[
 \lim_{N\to\infty}
 \Bigl(
 \sum_{\vecm\neq 0}
 \frac{\sin(\vecm\cdot\vecalpha\cdot N)}{N\cdot\vecm\cdot\vecalpha}
 \cdot b_{\vecm}
 \Bigr)=0.
\]
We have the following inequality:
\[
 \left|
 \frac{\sin(\vecm\cdot\vecalpha\cdot N)}{N\cdot\vecm\cdot\vecalpha}
\cdot b_{\vecm}
 \right|
\leq |b_{\vecm}|.
\]
We also have $\sum_{\vecm\neq 0}|b_{\vecm}|<\infty$.
Then we obtain the following:
\[
 \lim_{N\to\infty}\Bigl(
 \sum_{\vecm\neq 0}
 \frac{\sin(\vecm\cdot\vecalpha\cdot N)}{N\cdot\vecm\cdot\vecalpha}
\cdot b_{\vecm}
 \Bigr)
=
 \sum_{\vecm\neq 0}\Bigl(
 \lim_{N\to\infty}
 \frac{\sin(\vecm\cdot\vecalpha\cdot N)}{N\cdot\vecm\cdot\vecalpha}
\cdot b_{\vecm}\Bigr)
=0.
\]
Thus we are done.
\hfill\qed

\vspace{.1in}

The morphism $\psi_{\vecalpha}$
induces the morphism $\real\times\real_{>0}\lrarr T^n\times\real_{>0}$,
which we denote also by $\psi_{\vecalpha}$.
Let $f$ be an $\real$-valued $C^{\infty}$-function
on $T^n\times\real_{>0}$.
Then we have the Fourier decomposition as before:
\[
 f=\sum_{\vecm} a_{\vecm}(y)
 \cdot\exp\Bigl(
 \sqrt{-1}\sum_{i=1}^nm_i\cdot\theta_i
 \Bigr).
\]
Pick $y\in\real_{>0}$.
Then we put as follows:
\[
 \Xi_N(f)(y):=
 \int_{-N}^N(\psi_{\vecalpha}^{-1}f)(x,y)\frac{dx}{2N}.
\]
We have the limit $\Xi(f)(y):=\lim_{N\to\infty}\Xi_N(f)(y)$,
and thus we obtain the functions
$\Xi(f),\Xi_N(f):\real_{>0}\lrarr \real$.
If we decompose $a_{\vecm}(y)$ into the
real part $b_{\vecm}(y)$ and the imaginary part $c_{\vecm}(y)$,
then we have
$\Xi(f)(y)=b_0(y)$.
In particular,
$\Xi(f)$ is $C^{\infty}$
and we have the following:
\[
 \Bigl(
 \frac{\del}{\del y}
 \Bigr)^l
 \Xi(f)
=\frac{d^l b_0(y)}{dy^l}
=\Xi\Bigl(
 \frac{\del^l f}{\del y^l}
 \Bigr).
\]

\subsubsection{Convexity}

We continue to use the setting in the subsubsection 
\ref{subsubsection;10.26.95}.
We put $F:=\psi_{\vecalpha}^{-1}(f)$.
\begin{lem} \label{lem;10.26.91}
Assume $F$ is subharmonic,
i.e., the following inequality holds:
\[
 -\frac{\del^2 F}{\del x^2}
 -\frac{\del^2 F}{\del y^2}\leq 0.
\]
Then we have the following convexity of $\Xi(f)$:
\[
 -\frac{\del^2}{\del y^2}\Xi(f)\leq 0.
\]
\end{lem}
\pf
We have the following equality:
\[
 -\frac{\del^2}{\del y^2}\Xi(f)
=\Xi\Bigl(
 -\frac{\del^2f}{\del y^2}
 \Bigr)
=\lim_{N\to\infty}
 \int_{-N}^N
 \Bigl(
 -\frac{\del^2 F}{\del y^2}
 \Bigr)
 \frac{dx}{2N}.
\]
We have the following inequality due to the subharmonicity
of $F$:
\[
 \int_{-N}^N
 \Bigl(
 -\frac{\del^2 F}{\del y^2}
 \Bigr)\frac{dx}{2N}
\leq
 \int_{-N}^N
 \frac{\del^2F}{\del x^2}
 \frac{dx}{2N}
=\frac{1}{2N}
\Bigl(
 \frac{\del F}{\del x}(N,y)
-\frac{\del F}{\del x}(-N,y)
\Bigr).
\]
Let $V$ denote the vector field on $T^n$,
given as follows:
\[
 V=\sum_{i=1}^n\alpha_i\frac{\del}{\del \theta_i}.
\]
Then we have the following equality:
\[
 \frac{\del F}{\del x}=
\psi_{\vecalpha}^{-1}(Vf).
\]
If we fix $y$,
then $Vf$ is bounded function on a compact set $T^n$,
and thus $\del F/\del x$ is a bounded function on $\real$.
Hence we have the following:
\[
 \lim_{N\to\infty}
 \frac{1}{2N}
 \Bigl(
 \frac{\del F}{\del x}(N,y)
-\frac{\del F}{\del x}(-N,y)
 \Bigr)=0.
\]
Thus we obtain the result.
\hfill\qed

\subsubsection{An elementary boundedness of a convex function}

Lemma \ref{lem;10.26.91} will be used in the proof of 
preliminary constantness of the filtration
in the subsubsection \ref{subsubsection;11.9.1},
together with the following lemma.

\begin{lem} \label{lem;10.26.92}
Let $f:\real_{\geq 1}\lrarr\real$ be a $C^{\infty}$-function.
Assume the following:
\begin{enumerate}
\item \label{8.27.1}
 ${\displaystyle-\frac{d^2 f}{d y^2}\leq 0}$.
\item \label{8.27.2}
 There exist positive numbers $C_1$ and $C_2$
such that
 $f\leq C_1+C_2\cdot \log y$.
\end{enumerate}
Then there exists a positive number $C_3$
such that $f\leq C_3$.
\end{lem}
\pf
Due to the condition \ref{8.27.1},
the function $f$ is convex below.
On the other hand,
the right hand side in the condition \ref{8.27.2}
is convex above.
They imply that $f$ is dominated by a constant.
\hfill\qed

%% file: a34.4.tex

\subsubsection{Some integrals}

In this subsubsection,
$\phi$ denotes a $C^{\infty}$-function on $\real$
whose support is compact.
We denote the differential $\frac{d\phi}{dr}$ by $\phi'$.
We put as follows:
\[
 L_n:=\frac{(\log r^2)^n}{n!}.
\]

\begin{lem}
We have the following equality:
\begin{equation} \label{eq;9.25.30}
 \int_{0}^{\infty} \phi \cdot L_n \cdot r^{2s-1}\cdot dr
=\frac{1}{2}
\sum_{i=0}^n
 (-s)^{-i-1}\int_0^{\infty}\phi'\cdot L_{n-i}\cdot r^{2s}dr.
\end{equation}
\end{lem}
\pf
We have the following equality:
\[
 L_n\cdot r^{2s-1}
=\frac{1}{2}
\frac{d}{dr}
\Bigl(
 \sum_{i=0}^n (-1)^i\cdot s^{-i-1}
 L_{n-i}\cdot r^{2s}
\Bigr).
\]
Then (\ref{eq;9.25.30}) immediately follows.
\hfill\qed

\vspace{.1in}

Since $\phi'$ is $C^{\infty}$ around $r=0$,
$ \int_{0}^{\infty}\phi'\cdot L_k r^{2s}\cdot dr $
gives a entire function for the variable $s$ around $s=0$.
\begin{lem}
We have the following Taylor development at $s=0$:
\begin{equation} \label{eq;9.25.31}
 \int_{0}^{\infty}\phi'\cdot L_k \cdot r^{2s}\cdot dr
=\sum_{l=0}^{\infty}
 s^l\left(
 \begin{array}{c} k+l\\l\end{array}
 \right)
 \cdot\int_{0}^{\infty}\phi'\cdot L_{k+l}\cdot dr.
\end{equation}
\end{lem}
\pf
We have the following:
\[
 \frac{1}{l!}\left(
 \frac{d}{ds}\right)^l
\int_{0}^{\infty}\phi'\cdot L_k\cdot r^{2s}\cdot dr
=\frac{1}{l!}
\int_{0}^{\infty}
\phi'\cdot\frac{(\log r^2)^{k+l}}{k!}r^{2s}dr
=\left(
 \begin{array}{c}k+l\\l\end{array}
\right)
\cdot \int_{0}^{\infty}\phi'\cdot L_{k+l}\cdot r^{2s}dr.
\]
Then (\ref{eq;9.25.31}) immediately follows.
\hfill\qed

\begin{lem}
We have the following equality:
\begin{equation} \label{eq;9.25.32}
 \int_0^{\infty} \phi\cdot L_n\cdot r^{2s-1}\cdot dr
=\frac{(-1)^n}{2}\cdot s^{-n-1}\cdot\phi(0)
+ \sum_{l=n+1}^{\infty}
 X_{n,l}\cdot s^{l-n-1}\int_{0}^{\infty}\phi'\cdot L_l\cdot dr.
\end{equation}
Here we put as follows:
\begin{equation} \label{eq;11.9.2}
 X_{n,l}:=
 \frac{(-1)^{n-1}}{2}
 \sum_{h=0}^n (-1)^h
 \left(\begin{array}{c}l\\ h\end{array} \right)
\in\real.
\end{equation}
\end{lem}
\pf
We have the following:
\begin{equation}
 \int \phi\cdot L_n\cdot r^{2s-1}\cdot dr
=\frac{1}{2}
\sum_{i=0}^n(-s)^{-i-1}
 \sum_{l=0}^{\infty}s^l\cdot
 \left(\begin{array}{c}n-i+l\\ l\end{array}
 \right)
\int_{0}^{\infty}\phi'\cdot L_{n-i+l}\cdot dr.
\end{equation}
By putting $h=n-i$ and $m=n-i+l$,
the right hand side can be rewritten as follows:
\begin{multline}
 \frac{1}{2}
\sum_{m=0}^{\infty}
\sum_{\substack{h+l=m\\ h\leq n}}
(-1)^{-n+h-1}\cdot s^{-n+m-1}\cdot
\left(\begin{array}{c}m \\ n \end{array}\right)
\cdot \int_{0}^{\infty}\phi'\cdot L_m\cdot dr\\
=\frac{1}{2}\sum_{m=0}^{\infty}
\Bigl(
\sum_{\substack{h+l=m\\ h\leq n}}
(-1)^h\left(\begin{array}{c}m\\ n\end{array}
 \right)
\Bigr)\cdot(-1)^{n-1}s^{-n+m-1}
\int_{0}^{\infty}\phi'\cdot L_m\cdot dr.
\end{multline}
The term $m=0$ is as follows:
\[
 \frac{1}{2}\cdot (-1)^{n-1}s^{-n-1}
\cdot\int_0^{\infty}\phi'\cdot dr=
 \frac{(-1)^n}{2}\cdot s^{-n-1}\cdot \phi(0).
\]
It is easy to see that
the terms $1\leq m\leq n$ vanish.
Then we obtain (\ref{eq;9.25.32}).
\hfill\qed

\begin{cor}
We have the following development:
\[
 \int\phi\cdot L_n\cdot r^{2ms-1}dr
=
\frac{(-1)^n}{2}(ms)^{-n-1}\cdot\phi(0)
+\frac{1}{2}\sum_{l=n+1}^{\infty}
 X_{n,l}\cdot (m\cdot s)^{l-n-1}\int_{0}^{\infty}\phi'\cdot L_k \cdot dr.
\]
Here $X_{n,l}$ is given in {\rm(\ref{eq;11.9.2})}.
\hfill\qed
\end{cor}

\begin{cor}\label{cor;b11.23.65}
Let $\phi$ be a test function the complex plane $\cnum$.
We have the following formula:
\[
 \Res_{s=0}\int \phi\cdot L_n(|z|^2)\cdot |z|^{2s-2}\cdot
 \frac{\sqrt{-1}}{2\pi}dz\wedge d\bar{z}
=\left\{
 \begin{array}{ll}
 1 & (n=0),\\
 \mbox{{}}\\
 0 & (n\neq 0)\\
 \end{array}
 \right.
\]
\hfill\qed
\end{cor}

\subsubsection{Some distributions}

In this subsubsection $\phi$ denotes
a test function on $\cnum^n$.
Let us consider the function $\Phi$ on $\cnum^{\ast\,n}$
of the following form:
\[
 \Phi=
\sum_{k=1}^N
 \sum_{\vecn\in S}s^k\cdot a_{\vecn,k}
\cdot\prod_{i=1}^nL_{n_i}(|z_i|)^{n_i}.
\]
Here $S$ denotes a finite subset of
$\seisuu^n$,
$a_{\vecn,k}$ denote complex numbers,
and $n_i$ $(i=1,\ldots,n)$ denote the $i$-th component of $\vecn$.
We have the distribution $\hat{\Phi}$
defined as follows:
\[
 \hat{\Phi}(\phi):=
\Res_{s=0}\int_{\cnum^n}
 \Phi\cdot\phi\cdot
 \prod_{i=1}^n|z_i|^{2m_i\cdot s-2}\cdot 
 \frac{\sqrt{-1}}{2\pi}dz_i\wedge d\bar{z}_i.
\]

\begin{lem}\label{lem;a11.30.20}
Assume that the support of the distribution $\hat{\Phi}$
is contained in $\{O\}$.
Then we have the following formula:
\[
 \hat{\Phi}(\phi)
=\sum_{(\vecn,k)\in S_1} a_{\vecn,k}
 \cdot(-1)^{\sum n_i}\cdot\phi(0)\cdot\prod m_i^{-n_i-1}.
\]
Here $S_1$ denotes the set of the elements
 $(\vecn,k)\in S\times\seisuu$
satisfying $\sum n_i=k-n+1$.
\end{lem}
\pf
We have only to consider the test functions
of the form $\phi(z_1,\ldots,z_n)=\prod_{i=1}^n\phi_i(z_i)$.
Note the following equality:
\[
 \int \Phi\cdot\phi\cdot\prod_{i=1}^n |z_i|^{2m_is-2}
 \cdot\frac{\sqrt{-1}}{2\pi}dz_i\wedge d\bar{z}_i
=\sum_{i=1}^N \sum_{\vecn\in S}
 s^k\cdot a_{\vecn,k}\cdot 2^n
 \cdot\prod_{i=1}^n \int_{0}^{\infty}\phi_i\cdot L_{n_i}\cdot
 r^{2m_is-1}_idr_i.
\]
We have the following equalities:
\begin{multline} \label{eq;10.2.7}
 \prod_{i=1}^n
 \int_{0}^{\infty}\phi_i\cdot L_{n_i}\cdot r_i^{2m_is-1}dr_i\\
=\sum_{I\sqcup J=\nbar}
 \prod_{i\in I}m_i^{-n_i-1}
\cdot s^{-|\vecn|-n}
\cdot \prod_{i\in I}\Bigl(
 \frac{(-1)^{n_i}}{2}\phi_i(0)\Bigr)
\cdot\prod_{j\in J}\Bigl(
 \sum_{l_j=n_j+1}^{\infty}
 X_{n_j,l_j}\cdot (m_j\cdot s)^{l_j}
\cdot \int_{0}^{\infty}\phi'\cdot L_{l_j}\cdot dr
 \Bigr)\\
=
\sum_{I\sqcup J=\nbar}
 \prod_{i\in I}m_i^{-n_i-1}
 \cdot \prod_{j\in J}m_j^{l_j}
 \cdot s^{-|\vecn|-n+|\vecl|}
\cdot\prod_{i\in I}\Bigl(\frac{(-1)^{n_i}}{2}\phi_i(0) \Bigr)
\cdot\sum_{\vecl\in\seisuu^J}
 \prod_{j\in J}X_{n_j,l_j}\cdot
 \int \phi_j'\cdot L_{l_j}\cdot dr_j
\end{multline}
Here we put $|\vecn|=\sum_{i=1}^n n_i$ and $|\vecl|=\sum_{j\in J}l_j$.
We put as follows:
\[
 A(\vecm,\vecn,I,J,\vecl):=
 \prod_{i\in I}m_i^{-n_i-1}
 \cdot \prod_{j\in J}m_j^{l_j}
 \cdot\prod_{j\in J}X_{n_j,l_j}\cdot\prod_{i\in I}\frac{(-1)^{n_i}}{2}
\in\real.
\]
In the case $I=\nbar$,
we have the following equality:
\[
 A(\vecm,\vecn,\nbar,\emptyset,0)
=\prod_{i=1}^nm_i^{-n_i-1}\cdot(-1)^{|\vecn|}\cdot 2^{-n}.
\]

Then the right hand side can be rewritten as follows:
\[
 \sum_{I\sqcup J=\nbar}
 \sum_{\vecl\in\seisuu^J}
 A(\vecm,\vecn,I,J,\vecl)
\cdot s^{-|\vecn|-n+|\vecl|}
\cdot\prod_{i\in I}\phi_i(0)\cdot
\prod_{j\in J}\int \phi_j'\cdot L_{l_j}dr_j.
\]

We have the following:
\[
 \Res_{s=0}s^k\prod_{i=1}^n \int\phi_i\cdot L_{n_i}\cdot
 r_i^{2m_is-1}dr_i
=\sum_{I\sqcup J=\nbar}
 \sum_{\vecl\in S(\vecn,n,k)}
 A(\vecm,\vecn,I,J,\lbar)
\cdot\prod_{i\in I}\phi_i(0)\cdot
\prod_{j\in J}\int\phi'_j\cdot L_{l_j}\cdot dr_j.
\]
Here we put
$S(\vecn,n,k):=
 \bigl\{\vecl\in\seisuu^J_{\geq\,0}\,\big|\,-|\vecn|-n+|\vecl|+k=-1\bigr\}$.
Then we obtain the following:
\begin{multline}
 \hat{\Phi}(\phi)=
\sum_{\vecn,k}a_{\vecn,k}
\sum_{I\sqcup J=\nbar}\sum_{\vecl\in S(\vecn,n,k)}
A(\vecm,\vecn,I,J,\vecl)
\cdot \prod_{i\in I}\phi_i(0)\cdot
\prod_{j\in J}\int \phi_j'\cdot L_{l_j}\cdot dr_j\\
=\sum_{I\sqcup J=\nbar}\sum_{\vecl\in \seisuu^J}
 Y(I,J,\vecl)\prod_{i\in I}\phi_i(0)\cdot
 \prod_{j\in J}\int\phi_j'\cdot L_{l_j}\cdot dr_j
\end{multline}
Here we put as follows:
\[
 Y(I,J,\vecl)=\sum_{\vecn,k}a_{\vecn,k}\cdot
 A(\vecm,\vecn,I,J,\vecl)\in\real.
\]
Note the tuple
$\bigl\{\prod_{j\in J}L_{l_j}\,\big|\,\vecl\in\seisuu_{>0}^J\bigr\}$
of $C^{\infty}$-functions on $\cnum^{\ast\,J}$
is linearly independent over $\real$.
Since we have assumed that
the support of $\hat{\Phi}$ is contained in $\{O\}$,
we obtain the vanishings
of the constants $Y(I,J,\vecl)$ $(J\neq\emptyset)$.
Then we obtain the following:
\[
 \hat{\Phi}(\phi)
=\sum_{(\vecn,k)\in S_1}
 a_{\vecn,k}\cdot 2^n\cdot 2^{-n}\cdot (-1)^{|\vecn|}
 \cdot\prod_{i=1}^n m_i^{-n_i-1}\phi_i(0)
=\sum_{(\vecn,k)\in S_1}
 a_{\vecn,k}\cdot (-1)^{|\vecn|}\cdot
 \phi(0)\cdot \prod_{i=1}^nm_i^{-n_i-1}.
\]
Then we obtain the result.
\hfill\qed

%% file: b22.tex

\subsubsection{The generalized eigen decomposition}

Let $V$ be a finite dimensional vector space over $\cnum$,
and $f$ be an endomorphism of $V$.
We often denote the set of eigenvalues of $f$
by $\Sp(f)$.
For any element $\alpha\in\Sp(f)$,
we denote the generalized eigenspace corresponding to $\alpha$
by $\EE(f,\alpha)$.
We often denote it by $\EE(\alpha)$ or $\EE(V,\alpha)$,
if there are no confusion.
Note the compatibility of the filtration
and the generalized eigen decomposition.
\begin{lem}
Let $V$ be a finite dimensional vector space,
and $f$ be an endomorphism of $V$.
Let $F$ be a filtration of $V$
such that $f$ preserves $F$.
Then the generalized eigen decomposition of $f$ and $F$ is compatible.
\end{lem}
\pf
If a subspace $W$ of $V$ is preserved by $f$,
then we have $W=\bigoplus_{\alpha\in\Sp(f)} W\cap \EE(f,\alpha)$.
Then the lemma immediately follows.
\hfill\qed

\subsubsection{A lemma for the boundedness of the hermitian metrics}
\label{subsubsection;a12.2.5}

Let $V$ be a finite dimensional vector space over $\cnum$
and $h$ be a hermitian metric of $V$.
Let $S\subset \cnum$ be a finite subset.
Let $\eta$ be a positive number
such that $B_a(\eta)\cap B_b(\eta)=\emptyset$
for $a\neq b\in S$.
Then we put as follows:
\begin{equation}\label{eq;a12.9.1}
 \nbigs(S,\eta,C):=
 \bigl\{
 f\in End(V)\,\big|\,
 |f|_h\leq C,\,\,
 \Sp(f)\subset \bigcup_{a\in S}B_a(\eta)
 \bigr\}.
\end{equation}

\begin{lem}
The subset $\nbigs(S,\eta,C)\subset End(V)$ is compact.
\end{lem}
\pf
Due to the condition $|f|_h\leq C$,
the set $\nbigs(S,\eta,C)$ is bounded.
We also have the closedness of the defining condition
of $\nbigs(S,\eta,C)$.
\hfill\qed

\vspace{.1in}

For any $f\in\nbigs(S,\eta,C)$,
we have the decomposition of $V$:
\[
 V=\bigoplus_{a\in S}\EE_{\eta}(f,a),
\quad
 \EE_{\eta}(f,a):=
 \bigoplus_{\substack{\alpha\in\Sp(f),\\ |\alpha-a|<\eta}}
 \EE(f,\alpha).
\]
Then we obtain the hermitian metric
$h(f)$ given as follows:
\[
 h(f):=\bigoplus_{a\in S}h_{|\EE_{\eta}(f,a) }.
\]

\begin{lem} \label{lem;8.16.1}
The set
$\bigl\{h(f)\,\big|\,f\in\nbigs(S,\eta,m,C)\bigr\}$
is compact.
\end{lem}
\pf
It immediately follows from the compactness of $\nbigs(S,\eta,C)$
\hfill\qed

\vspace{.1in}

\subsubsection{A lemma for $\epsilon$-orthogonality}

Let $V$ be a finite dimensional vector space over $\cnum$.
Let $S$ be a finite set,
and $V=\bigoplus_{a\in S}V_a$ be a decomposition.
Let $h$ be a hermitian metric of $V$.
\begin{df}
Let $\epsilon$ be a positive number such that $\epsilon\leq 1$.
The decomposition $V=\bigoplus_{a\in S} V_a$
is called $\epsilon$-orthogonal,
if the inequalities
$|h(u,v)|\leq \epsilon\cdot |u|_h\cdot |v|_h$
hold for any elements $u,v\in V$.
\hfill\qed
\end{df}

Assume the following:
\begin{itemize}
\item
There exist positive constants $C_1$ and $C_2$
such that the following holds for any 
element $v=\sum_{a\in S} v_a$ of $V$:
\[
 C_1\cdot\sum |v_b|_h\leq |v|_h\leq C_2\cdot\sum |v_b|_h.
\]
\item
We assume that the decomposition $V=\bigoplus_{a\in S}V_a$
is $\epsilon$-orthogonal with respect to $h$.
\end{itemize}

Let $g$ be an element of $\bigoplus_{a\in S}End(V_a)$,
and $g^{\dagger}$ be the adjoint of $g$ with respect to $h$.
We have the decomposition:
\[
 g^{\dagger}=\sum (g^{\dagger})_{a\,b},
\quad\quad
 (g^{\dagger})_{a\,b}\in Hom(V_a,V_b).
\]

\begin{lem} \label{lem;9.7.120}
There exists a positive constant $C$
such that the following holds:
\begin{itemize}
\item $C$ is independent of $\epsilon$.
\item 
The inequalities
$ \bigl|(g^{\dagger})_{a\,b}\bigr|_h
\leq
 C\cdot\epsilon\cdot|g|_h$
hold for any elements $g\in \bigoplus_{a\in S}End(V_a)$
and for any $a\neq b\in S$.
\end{itemize}
\end{lem}
\pf
Let $v$ be an element of $V_a$.
We put $w_b:=g^{\dagger}_{a\,b}(v)$.
Let $u$ be an element of $V_c$ for $c\neq a$.
Then we have the following equality:
\[
 \bigl(v,g(u)\bigr)_h
=\bigl(g^{\dagger}(v),u\bigr)_h
=\sum (w_b,u)_h.
\]
Hence we have the following:
\[
 (w_c,u)_h=\bigl(v,g(u)\bigr)_h-\sum_{b\neq c}(w_b,u)_h.
\]
Then we obtain the following due to the $\epsilon$-orthogonality:
\begin{multline}
 |(w_c,u)_h|
\leq
 \epsilon\cdot |v|_h\cdot|g|_h\cdot|u|_h
+\sum_{b\neq c}\epsilon\cdot|w_b|_h\cdot|u|_h
\leq
 \epsilon\cdot |v|_h\cdot|g|_h\cdot|u|_h
+\epsilon\cdot C_1^{-1}|g^{\dagger}(v)|_h\cdot |u|_h \\
\leq
 \epsilon
 \cdot
 \bigl(
 1+C_1^{-1}
 \bigr)
|v|_h\cdot|g|_h\cdot|u|_h.
\end{multline}
In particular,
we consider the case $u=w_c$:
\[
 |w_c|_h\leq
 \epsilon\cdot \bigl(
 1+ C_1^{-1}\bigr)
 \cdot|v|_h\cdot |g|_h
\]
Then we obtain the following inequality for any $v\in V_a$:
\[
 |g^{\dagger}_{a\,c}(v)|_h
\leq \epsilon\cdot \bigl(
 1+C_1^{-1}\bigr) 
 \cdot|v|_h\cdot|g|_h.
\]
It implies the claim.
\hfill\qed


%% file: c11.1.tex

\subsubsection{Norm of endomorphisms}
\label{subsubsection;04.2.16.2}

Let $V$ be an $r$-dimensional vector space,
and $h$ be a hermitian metric of $V$.
Let $\vecv$ be a base of $V$.
We put $H:=H(h,\vecv)$.
Let $f$ be an endomorphism of $V$,
and then we have the matrix $A$
such that $f\cdot\vecv=\vecv\cdot A$.
Let $f^{\dagger}$ denote the adjoint of $f$
with respect to the metric $h$.
Then we have
$f^{\dagger}\vecv=
 \vecv\cdot \overline{H}^{-1}\cdot\lefttop{t}\overline{A}\cdot\overline{H}$.

We have the norm $|f|_h$ of $f$ with respect to the metric $h$.
\begin{lem}\label{lem;04.1.6.1}
The following holds:
\[
 |f|_h^2
=\tr\bigl(A\cdot\overline{H}^{-1}\cdot
  \lefttop{t}\overline{A}\cdot \overline{H}\bigr).
\]
\end{lem}
\pf
We have $|f|_h^2=\tr(f\cdot f^{\dagger})$.
Then the claim immediately follows.
\hfill\qed

%% file: d4.tex

\subsubsection{The $GL(r)$-invariant metric}
\label{subsubsection;04.1.27.5}

Let $\PH(r)$ denote the set of
$(r\times r)$-positive definite hermitian matrices.
We have the standard left action $\kappa$ of $GL(r)$ on $\PH(r)$:
\[
 GL(r)\times \PH(r)\lrarr \PH(r),
\quad
 (g,H)\longmapsto \kappa(g,H)=g\cdot H\cdot \lefttop{t}\bar{g}.
\]
For any point $H\in\PH(r)$,
the tangent space $T_H\PH(r)$ is naturally identified with
the vector space $\nbigh(r)$.
Let $I_r$ denote the identity matrix.
We have the positive definite metric of $T_{I_r}\PH(r)$
given by $(A,B)_{I_r}=\tr(A\cdot B)=\tr(A\cdot \lefttop{t}\bar{B})$.
It is easy to see the metric is invariant
with respect to the $U(r)$-action on $T_{I_r}\PH(r)$.

Let $H$ be any point of $\PH(r)$
and let $g$ be an element of $GL(r)$
such that $H=g\cdot \lefttop{t}\bar{g}$.
Then the metric of $T_{H}\PH(r)$ is given as follows:
\begin{equation}\label{eq;04.1.4.2}
(A,B)_{H}=\bigl(\kappa(g^{-1})_{\ast}A,\kappa(g^{-1})_{\ast}B\bigr)_{I_r}
=\tr\Bigl(
 g^{-1}A\lefttop{t}\bar{g}^{-1}\cdot g^{-1}B\lefttop{t}\bar{g}^{-1}
 \Bigr)
=\tr\bigl(H^{-1}AH^{-1}B\bigr).
\end{equation}
Since $(\cdot,\cdot)_{I_r}$ is $U(r)$-invariant,
the metric $(\cdot,\cdot)_{H}$ on $T_H\PH(r)$ is well defined.
Thus we have the $GL(r)$-invariant Riemannian metric
of $\PH(r)$.

It is well known that $\PH(r)$ with the metric above is
a symmetric space with non-positive curvature.
We denote the induced distance by $d_{\PH(r)}$.
We often use the simple notation $d$ to denote $d_{\PH(r)}$,
if there are no confusion.

Let $X$ be a manifold,
and $\Psi:X\lrarr\PH(r)$ be a differentiable map.
Let $P$ be a point of $X$,
and $v$ be an element of the tangent space $T_PX$.
\begin{lem}\label{lem;04.1.13.1}
We have the following formula:
\[
 \bigl|d\Psi(v)\bigr|^2_{T_{\Psi(P)}\PH(r)}
=\tr\Bigl(
 \Psi(P)^{-1}\cdot d\Psi(v)\cdot \Psi(P)^{-1}\cdot d\Psi(v)
 \Bigr)
=\tr\Bigl(
 \overline{\Psi(P)}^{-1}\cdot
 \overline{d\Psi(v)}
 \cdot 
 \overline{\Psi(P)}^{-1}\cdot
 \overline{d\Psi(v)}
 \Bigr).
\]
\end{lem}
\pf
It follows from (\ref{eq;04.1.4.2}).
\hfill\qed

%% file: c12.1.tex

\subsubsection{Comparison of the norm and the distance}
\label{subsubsection;04.2.16.1}

For any elements $H_1$ and $H_2$ of $\PH(r)$,
we have an element $g\in GL(r)$
such that $\kappa(g,H_1)=I_r$
and that $\kappa(g,H_2)$ is the diagonal matrix.
The set of the eigenvalues of $\kappa(g,H_2)$ is independent
of a choice of $g$.
Let $e^{\alpha_1},\ldots,e^{\alpha_r}$ be 
the eigenvalues of $\kappa(g,H_2)$.
We put as follows:
\[
 \delta(H_1,H_2):=
\left(
 \sum_{i=1}^r \Bigl(
 \frac{e^{\alpha_i}-e^{-\alpha_i}}{2}
 \Bigr)^2
\right)^{1/2}.
\]
On the other hand,
we have the distance
$d_{\PH(r)}(H_1,H_2)=\bigl(\sum \alpha_i^2\bigr)^{1/2}$.

For any real number $R$, we put as follows:
\[
 C(R):=\frac{e^{R}-e^{-R}}{2R}.
\]
If $0\leq x\leq R$,
we have $x\leq C(R)\cdot x$.
We also note that $C(R)\to 1$ when $R\to 0$.

\begin{lem}\mbox{{}}\label{lem;a12.28.12}
\begin{itemize}
\item
 We have the inequality:
\[
 d_{\PH(r)}(H_1,H_2)
\leq
 \delta(H_1,H_2).
\]
\item
 If $d(H_1,H_2)\leq R$,
 we have the inequality:
\[
 \delta(H_1,H_2)\leq
 C(R)\cdot d_{\PH(r)}(H_1,H_2).
\]
\end{itemize}
\end{lem}
\pf
It can be checked elementarily.
\hfill\qed

\vspace{.1in}
We reformulate Lemma \ref{lem;a12.28.12} as follows:
Let $V$ be an $r$-dimensional vector space,
and let $h_1$ and $h_2$ be hermitian metrics of $V$.
The identity map induces the map
$\Phi:(V,h_1)\lrarr (V,h_2)$.
We have the norms $|\Phi|$ and $\big|\Phi^{-1}\big|$.

\begin{lem}\label{lem;04.1.27.6}\mbox{{}}
\begin{itemize}
\item
 The following inequality holds:
 \[
  d_{\PH(r)}(h_1,h_2)^2
\leq
 \frac{
 \big|\Phi\big|^2
+\big|\Phi^{-1}\big|^{2}
-2r}{4}.
 \]
\item
 If $d_{\PH(r)}(h_1,h_2)\leq R$,
 the following inequality holds:
\[
 \frac{
  \big|\Phi\big|^2+\big|\Phi^{-1}\big|^2-2r}{4}
\leq
 C(R)^2\cdot
 d_{\PH(r)}(h_1,h_2)^2.
\]
\end{itemize}
\end{lem}
\pf
If we take an appropriate base of $V$,
$h_1$ and $h_2$ are represented by 
the identity matrix $I_r$
and the diagonal matrices whose diagonal entries are
$e^{\alpha_1},\ldots,e^{\alpha_r}$.
It is easy to check that
$|\Phi|^2=\sum_{i=1}^r e^{2\alpha_i}$
and $|\Phi^{-1}|^2=\sum_{i=1}^r e^{-2\alpha_i}$.
Then it immediately follows
$4^{-1}\bigl(
 |\Phi|^2+|\Phi^{-1}|^2-2r
 \bigr)=\delta(H_1,H_2)^2$.
Thus the claims follow from Lemma \ref{lem;a12.28.12}.
\hfill\qed

%% file: 15.1.tex

\subsubsection{Some results of Andreotti-Vesentini}

We recall some results of Andreotti-Vesentini in \cite{av}.
Let $(Y,g)$ be a complete Kahler manifold, not necessarily
compact.
We denote the natural volume form of $Y$ by $\vol$.
Let $(E,\delbar_E,h)$ be a hermitian holomorphic bundle
over $Y$.
The hermitian metric $h$ and the Kahler metric $g$
induces the fiberwise hermitian metric
of $E\otimes\Omega^{p,q}_Y$,
which we denote by $(\cdot,\cdot)_{h,g}$.
The space of $C^{\infty}$ $(p,q)$-forms with compact support
is denoted by $A_c^{p,q}(E)$.
For any $\eta_1,\eta_2\in A_c^{p,q}(E)$,
we put as follows:
\[
 \langle
 \eta_1,\eta_2
 \rangle_h=
 \int (\eta_1,\eta_2)_{h,g}\cdot\vol,
\quad
 ||\eta||^2_h=\langle \eta,\eta\rangle_h.
\]
The completion of $A_c^{p,q}(E)$
with respect to the norm $||\cdot||_h$
is denoted by $A_h^{p,q}(E)$.

We have the operator
$\delbar_E:A_c^{p,q}(E)\lrarr A_c^{p,q+1}(E)$,
and the formal adjoint
$\delbar_E^{\ast}:A_c^{p,q}(E)\lrarr A_c^{p,q-1}(E)$.
We use the notation
$\laplacian=\delbar_E^{\ast}\delbar_E+\delbar_E\delbar_E^{\ast}$.
We have the maximal closed extensions
$\delbar_E:A_h^{p,q}(E)\lrarr A_h^{p,q+1}(E)$
and 
$\delbar_E^{\ast}:A_h^{p,q}(E)\lrarr A_h^{p,q-1}(E)$.
We denote the domains of $\delbar_E$
and $\delbar_E^{\ast}$ by $Dom(\delbar_E)$
and $Dom(\delbar_E^{\ast})$ respectively.

\begin{prop}[Proposition 5 of \cite{av}] \label{prop;11.28.6}
In $W^{p,q}:=Dom(\delbar_E)\cap Dom(\delbar_E^{\ast})$,
the space $A_c^{p,q}(E)$ is dense
with respect to the the graph norm:
$||\eta||_h^2+||\delbar_E\eta||^2_h+||\delbar_E^{\ast}\eta||^2_h$.
(See also {\rm\cite{cg}}).
\hfill\qed
\end{prop}

\begin{prop}[Theorem 21 of \cite{av}] \label{prop;12.11.1}
Assume that 
there exists a positive number $c>0$ satisfying the following:
\begin{quote}
Then, for any $\eta\in W^{p,q}$,
we have $||\delbar_E\eta||_h^2+||\delbar_E^{\ast}\eta||_h^2
 \geq c\cdot ||\eta||_{h}^2$.
\end{quote}
For any $C^{\infty}$-element $\eta\in A_h^{p,q}(E)$
such that $\delbar_E(\eta)=0$,
we have a $C^{\infty}$-solution $\rho\in A_h^{p,q-1}(E)$
satisfying the equation $\delbar_E(\rho)=\eta$.
\hfill\qed
\end{prop}

\subsubsection{Kodaira identity}

For a Kahler manifold $Y$,
we have the operator $\Lambda:\Omega^{p,q}\lrarr \Omega^{p-1,q-1}$
(see  62 page of \cite{koba}).
For a section $f$ of $End(E)\otimes \Omega^{p_0,q_0}_Y$,
we have the natural morphism $A_c^{p,q}(E)\lrarr A_c^{p+p_0,q+q_0}(E)$,
defined by $\eta\longmapsto f\wedge \eta$.
We denote the morphism by $e(f)$.

Let $E$ be a holomorphic vector bundle with a hermitian metric $h$
over $Y$.
We have the metric connection of $E$
induced by the holomorphic structure $\delbar_E$
and the hermitian metric $h$.
We denote the curvature by $R(h)$.

We have the Levi-Civita connection of the tangent bundle
of $Y$.
It induces the connection of $E\otimes\Omega^{0,q}$:
\[
\nabla:
 A_c^{0,0}(E\otimes\Omega^{0,q})\lrarr
 A_c^{0,1}(E\otimes\Omega^{0,q})\oplus
 A_c^{1,0}(E\otimes \Omega^{0,q}).
\]
We denote the $(0,1)$-part of $\nabla$ by $\nabla''$
to distinguish with $\delbar_E:A_c^{0,q}(E)\lrarr A_c^{0,q+1}$.
The $(1,0)$-part of $\nabla$ is same as
$\del$ of $E\otimes\Omega^{0,q}$.
We denote the curvature of $\nabla$
by $R(\nabla)$.

We have the Laplacian $\Delta''$ on $A^{0,q}_c(E)$
and the equalities of the operators:
\[
 \Delta''=\delbar_E\delbar_E^{\ast}+\delbar_E^{\ast}\delbar_E
=\del_E\del_E^{\ast}+\del_E^{\ast}\del_E
+\sqrt{-1}\bigl[e(R(h)),\Lambda\bigr].
\]
In particular,
we have the equality
$\Delta''
 =\del_E^{\ast}\del_E-\sqrt{-1}\Lambda\circ e(R(h))$
on $A_c^{0,q}(E)$.

On the other hand,
we have the following Laplacian
on $A_c^{0,0}(E\otimes\Omega^{0,q})$:
\[
 \Delta''=
 \nabla''\nabla^{\prime\prime\ast}
+\nabla^{\prime\prime\ast}\nabla''
=\del_E^{\ast}\del_E-\sqrt{-1}\Lambda\circ e(R(\nabla)).
\]

For an element $\eta$ of $A_c^{0,0}(E\otimes\Omega^{0,q})$,
note the following equality:
\[
 \Lambda\circ e(R(\nabla))(\eta)
=\Lambda(R(\nabla))\cdot \eta.
\]
Since we have $R(\nabla)=R(h)+R(\Omega^{0,q})$,
we have the following:
\[
 \Lambda\circ e(R(\nabla))(\eta)
=\Lambda(R(h))\cdot\eta
+\Lambda(R(\Omega^{0,q}))\cdot\eta.
\]
Then we have the following identity:
\begin{multline}
\big\langle
 \delbar_E\eta,\,\delbar_E\eta
\big\rangle_h
=\big\langle
  \del_E\eta,\,\del_E\eta
 \big\rangle_h
-\sqrt{-1}\cdot\big\langle
 \Lambda(R(h)\cdot\eta),\,\eta
 \big\rangle_h \\
=\big\langle
  \nabla''\eta,\,\nabla''\eta
 \big\rangle_h
+\sqrt{-1}\cdot\big\langle
 \Lambda(\Omega^{0,q})\cdot\eta,\,\eta
 \big\rangle_h
-\sqrt{-1}\cdot
 \big\langle
  \Lambda\bigl(R(h)\cdot \eta\bigr)
 -\Lambda R(h)\cdot \eta,\,
 \eta
 \big\rangle_h
\end{multline}

Let $\chi$ be a $\real$-valued $C^{\infty}$-function.
If we put $\tilde{h}:=h\cdot e^{-\chi}$,
we obtain the following equality:
\begin{multline} \label{eq;8.18.16}
 \bigl\langle
 \delbar_E\eta,\,\,\delbar_E\eta
 \bigr\rangle_{\tilde{h}}
+\bigl\langle
 \del_E\eta,\,\,\del_E\eta
 \bigr\rangle_{\tilde{h}} \\
=\bigl\langle
 \nabla''\eta,\,\,
 \nabla''\eta
 \bigr\rangle_{\tilde{h}}
+\sqrt{-1}\cdot \bigl\langle
 \Lambda(\Omega^{0,q})\eta,\,\,
 \eta
 \bigr\rangle_{\tilde{h}} 
-\sqrt{-1}\cdot \bigl\langle
 \Lambda(R(h)\eta)-\Lambda (R(h))\cdot\eta,\,\,
 \eta
 \bigr\rangle_{\tilde{h}} \\
+\sqrt{-1}\bigl\langle
 \Lambda(\delbar\del\chi\cdot\eta)-\Lambda(\delbar\del\chi)\cdot\eta,\,\,
 \eta
 \bigr\rangle_{\tilde{h}}.
\end{multline}

\subsubsection{Kodaira identity for sections of $A^{0,1}_c(E)$}

We denote the Ricci curvature of the Kahler metric $g$
by $\Ric(g)$.
We can naturally
regard $\Ric(g)$ as a section of $End(E)\otimes \Omega^{1,1}$,
by the natural diagonal inclusion $\cnum\lrarr End(E)$.

Let $f$ be a section of $End(E)\otimes\Omega^{1,1}_Y$,
and $\eta$ be an element of $A_c^{0,1}(E)$.
Then we put as follows: \label{page;133}
\begin{equation} \label{eq;10.11.60}
 \doublelangle
 f,\eta
 \doublerangle_h
:=
 -\sqrt{-1}
\bigl(
\xi,\eta
\bigr)_h,
\quad
 \xi:=
 \Bigl(
 \Lambda\circ e(f)
 -e\bigl(\Lambda(f)\bigr)
 \Bigr)
 (\eta)
=
 \Lambda
 \bigl(
 f\cdot\eta
 \bigr)
-\Lambda(f)\cdot\eta
\end{equation}

We recall the following special case.
\begin{prop}[Kodaira \cite{ko1}, Cornalba-Griffiths \cite{cg}]
Let $\eta$ be an element of $A_c^{0,1}(E)$.
We have the following equality:
\[
 ||\delbar_E(\eta)||_h^2
+||\delbar_E^{\ast}(\eta)||_h^2
=
 ||\nabla''\eta||^2+
\int\!\!
\big\langle\big\langle
R(h)+Ric(g),
\eta
\big\rangle\big\rangle_h\!
\vol.
\]
\end{prop}
\pf
See Proposition 4.5 in \cite{mochi}.
\hfill\qed

\begin{cor}\label{cor;12.11.2}
Let $\eta$ be an element of
$Dom(\delbar_E)\cap Dom(\delbar_E^{\ast})$
in $A^{0,1}_h(E)$.
Then we have the following inequality:
\[
 ||\delbar_E(\eta)||^2_h
+||\delbar_E^{\ast}(\eta)||_h^2
\geq
\int\!\!
 \doublelangle
 R(h)+\Ric(g),\eta
 \doublerangle_h
\vol.
\]
\hfill\qed
\end{cor}

\subsubsection{Acceptable bundle}

\label{subsubsection;d11.14.35}

We put $X=\Delta^n$ and $D=\bigcup_{i=1}^l D_i$.
We have the Poincar\'{e} metric $g_{\poin}$ on $X-D$:
\[
 g_{\poin}:=
 \sum_{j=1}^l q_j^{\ast}g_0
+\sum_{j=l+1}^n q_j^{\ast}g_1.
\]
Here we put as follows:
\[
 g_1=\frac{2\cdot dz\cdot d\bar{z}}{(1-|z|^2)^2},
\quad\quad
 g_0=\frac{2\cdot dz\cdot d\bar{z}}{|z|^2(-\log|z|^2)^2}.
\]
We have the corresponding Kahler form:
\[
 \omega_{\poin}:=
 \sum_{j=1}^lq_j^{\ast}\omega_0
+\sum_{j=l+1}^{n}q_j^{\ast}\omega_1.
\]
Here we put as follows:
\[
 \omega_1:=
 \frac{\sqrt{-1}dz\wedge d\bar{z}}{(1-|z|^2)^2},
\quad\quad
 \omega_0:=
 \frac{\sqrt{-1}dz\wedge d\bar{z}}{|z|^2(-\log|z|^2)^2}.
\]

Let $P$ be a point of $X$
and $(\nbigu,\varphi)$ be an admissible coordinate around $P$.
By the isomorphism
$\varphi:\nbigu-D\simeq \Delta^{\ast\,l}\times\Delta^{n-l}$,
we take the Poincar\'e metric $g_{\Poin}$ on $\nbigu-D$.
The metric $h$ of $E$ and the metric $g_{\Poin}$ on $T(\nbigu-D)$
induce the metric $(\cdot,\cdot)_{h,g_{\Poin}}$
of $End(E)\otimes \Omega^{p,q}$
over $\nbigu-D$.
Recall the following definition.
\begin{df} \label{df;11.28.5}
We say that $(E,\delbar_E,h)$ is acceptable at $P$,
if the following holds:
\begin{itemize}
\item
 Let $(\nbigu,\varphi)$ be an admissible coordinate around $P$.
  The norms of the curvature $R(h)$ with respect to
 the metric $(\cdot,\cdot)_{h,g_{\poin}}$
 is bounded over $\nbigu-D$.
\end{itemize}
When $(E,\delbar_E,h)$ is acceptable at any point $P$,
then we say that it is acceptable.
\hfill\qed
\end{df}

Let $(E,\delbar_E,h)$ be an acceptable bundle over $X-D$.
We apply (\ref{eq;8.18.16}) for
$\chi=\tau(\veca,N)$,
where the function $\tau(\veca,N)$ is as follows:
\begin{equation}\label{eq;a11.9.7}
 \tau(\veca,N):=
 \sum_{i=1}^la_i\log|z_i|^2
 +N\cdot
 \Bigl(
 \sum_{i=1}^l\log(-\log|z_i|^2)
+\sum_{i=l+1}^n\log(1-|z_i|^2)
 \Bigr).
\end{equation}
We use the notation
$|\cdot|_{\veca,N}$,
$||\cdot||_{\veca,N}$,
$(\cdot,\cdot)_{\veca,N}$ and
$\doublelangle\cdot,\cdot\doublerangle_{\veca,N}$
instead of
$|\cdot|_{h_{\veca,N}}$,
$||\cdot||_{h_{\veca,N}}$,
$(\cdot,\cdot)_{h_{\veca,N}}$ and
$\doublelangle\cdot,\cdot\doublerangle_{h_{\veca,N}}$
for simplicity.
We also use the notation $A_{\veca,N}^{p,q}(E)$
instead of $A_{h_{\veca,N}}^{p,q}(E)$.

\label{subsubsection;10.12.1}

\begin{lem}
We have the following equality for any section $\eta\in A^{0,q}_c(E)$:
\[
 \sqrt{-1}
 \Bigl(
 \bigl(
 \Lambda\delbar\del\tau(\veca,N)\eta
 \bigr)
-\Lambda( \delbar\del\tau(\veca,N)\cdot\eta)
 \Bigr)
 =
-N\cdot q\cdot \eta.
\]
\end{lem}
\pf
We have $\delbar\del\tau(\veca,N)=-N\cdot\sqrt{-1}\cdot\omega_{\poin}$.
Then we obtain the following:
\[
 {\rm L.H.S.}
 =\sqrt{-1}\cdot \bigl(-N\sqrt{-1}\bigr)\cdot
 \Bigl(
 \Lambda(\omega_{\poin}\eta)
-(\Lambda\omega_{\poin})\cdot\eta
 \Bigr)
=N\cdot
\bigl(-(n-q)\cdot\eta-n\cdot\eta\bigr)
=-q\cdot N\cdot\eta={\rm R.H.S.}.
\]
Hence we are done.
\hfill\qed

\subsubsection{The vanishing of the cohomology of acceptable bundle
 (When $N$ is sufficiently negative)}

We have the following inequality:
\begin{equation}\label{eq;11.9.3}
 \bigl\langle
 \delbar_E\eta,\,\,\delbar_E\eta
 \bigr\rangle_{\veca,N}
+\bigl\langle
 \del_E\eta,\,\,\del_E\eta
 \bigr\rangle_{\veca,N}
\geq
 \sqrt{-1}\cdot
 \bigl\langle
  \Lambda(\Omega^{0,q})\cdot\eta,\,\,\eta
 \bigr\rangle_{\veca,N}
-\sqrt{-1}\cdot
 \Bigl\langle
 \Lambda\bigl(R(h)\cdot\eta\bigr)
-\Lambda R(h)\cdot\eta,\,\,
\eta
 \Bigr\rangle_{\veca,N}
-N\cdot q\cdot||\eta||_{\veca,N}^2.
\end{equation}

\begin{lem} \label{lem;11.9.4}
When $(E,\delbar_E,h)$ is an acceptable bundle,
there exists a positive constant $C>0$
satisfying the following for any $q=1,\ldots,n$
and for any $\eta\in A_c^{0,q}(E)$:
\[
 \Bigl|
 \Bigl\langle
 \Lambda(\Omega^{0,q})\cdot\eta-\Lambda\bigl(R(h)\cdot\eta\bigr)
 +\Lambda(R(h))\cdot\eta,\,\,
\eta
 \Bigr\rangle_{h}
 \Bigr|
\leq
 C\cdot|\eta|_{h}^2
\]
\end{lem}
\pf
Recall that $R(h)$ is dominated by $\omega_{\poin}$.
Then the claim follows from a direct calculation
of the curvature of Poincar\'{e} metric $g_{\poin}$.
\hfill\qed

\vspace{.1in}

If we take a sufficiently negative integer $N$
such that $N<-C-1$ for the constant $C$ as in Lemma \ref{lem;11.9.4},
then we obtain the following inequalities for any $q\geq 1$
and for any $\eta\in A_{c}^{0,q}(E)$,
due to the inequality (\ref{eq;11.9.3}):
\begin{equation}\label{eq;11.9.5}
 \bigl\langle
  \delbar_E\eta,\,\,\delbar_E\eta
 \bigr\rangle_{\veca,N}
+\bigl\langle
 \del_E\eta,\,\,\del_E\eta
 \bigr\rangle_{\veca,N}
\geq ||\eta||_{\veca,N}^2.
\end{equation}

\begin{lem} \label{lem;9.10.72}
Let $C$ be a positive constant as in Lemma {\rm\ref{lem;11.9.4}}.
If $N<-C-1$,
we have the vanishings of
any higher cohomology group
$H^i\bigl(A^{0,\cdot}_{\veca,N}\bigl(
 \nbigelambda\bigr),\delbar_{\nbigelambda}\bigr)$ $(i>0)$.
\end{lem}
\pf
It follows from Proposition \ref{prop;11.28.6}
and Proposition \ref{prop;12.11.1}
and (\ref{eq;11.9.5}).
\hfill\qed

\subsubsection{A lemma for an increasing order
 (When $N$ is sufficiently positive)}

Let $(E,\delbar_E,h)$ be an acceptable hermitian holomorphic bundle
over $X-D$.
The following lemma will be used in many times.
It says, we obtain the increasing order of
a holomorphic section over $\Delta^{\ast\,l}\times \Delta^{n-l}$
from the increasing order of the restriction
to the curves, which are transversal with the smooth part of
the singularity $D$.

Let $\pi_j$ denote the projection of $X$ onto $D_j$.
Let $Y(C)$ denote the set
$\bigl\{(z_1,\ldots,z_n)\,\big|\,
 0<|z_i|<C,\,(i=1,\ldots,l),\,\,|z_i|<C\,\,(i=l+1,\ldots,n)\bigr\}$.
\begin{cor} \label{cor;11.28.15}
Let $F$ be a holomorphic section of $\nbige^{\lambda}$.
Let $a_j$ and $k_j$ be real numbers $(j=1,\ldots,l)$.

For any point $p\in\Delta^{\ast}(C)^{l-1}\times \Delta(C)^{n-l}$
and any $1\leq j\leq l$,
we assume that we are given numbers
$C_1(p,j)$, $C_2(p,j)$, $a(p,j)$ and $k(p,j)$
satisfying the following:
\begin{enumerate}
\item $C_1(p,j)$ and $C_2(p,j)$ are positive numbers.
\item
 $a(p,j)$ and $k(p,j)$ are real numbers satisfying
 $a(p,j)\leq a_j$ and $k(p,j)\leq k_j$.
\item
The following inequality holds on $\pi_j^{-1}(p)$:
\[
 0<C_1(p,j)\leq \bigl|F_{|\pi_j^{-1}(p)}\bigr|_h\cdot |z_j|^{-a(p,j)}\cdot
 \bigl(-\log|z_j|\bigr)^{-k(p,j)}
 \leq C_2(p,j).
\]
\item $C_1(p,j)$, $C_2(p,j)$, $a(p,j)$ and $k(p,j)$ may depend on $p$ and $j$.
\end{enumerate}
Then there exists a positive constant $C_3$ and a large number $M$,
satisfying the following:
\begin{itemize}
\item
 The inequality
 $|s|_h\leq C_3\cdot \prod_{j=1}^l|z_j|^{a_j}(-\log |z_j|)^M$
 holds over $Y(C)$.
\item
The constant $C_3$ depends only on the values of $|s|_h$
on the following compact set:
\[
 \bigl\{(z_1,\ldots,z_n)\,\big|\,|z_j|=C,\,\,(j=1,\ldots,l),\,\,
 |z_j|\leq C\,\,(j=l+1,\ldots,n)\bigr\}.
\]
\end{itemize}
\end{cor}
\pf
See Corollary 4.12 of \cite{mochi}.
\hfill\qed

%% file: a49.1.tex

\subsubsection{The metrics and the curvatures of $\nbigo(i)$ on $\proj^1$}

\label{subsubsection;10.11.70}

In the next few subsubsecions,
we recall the subsubsection 4.7.3 in \cite{mochi}.
Let $\proj^1$
denote the one dimensional projective space over $\cnum$.
We use the homogeneous coordinate $[t_0:t_1]$.
The points $[0:1]$ and $[1:0]$ are denoted by $0$ and $\infty$
respectively.
We use the coordinates $t=t_0/t_1$ and $s=t_1/t_0$.
We have the line bundle $\nbigo(i)$ over $\proj^1$.
The coordinates of $\nbigo(i)$ is given as follows:
$(t,\zeta_1)$ over $\proj^1-\{\infty\}$,
and $(s,\zeta_2)$ over $\proj^1-\{0\}$.
The relations are given by
$s=t^{-1}$ and $t^{-i}\cdot \zeta_1=\zeta_2$.

Recall that we have the smooth metric $h_i$ of $\nbigo(i)$.
Let $\xi=(t,\zeta_1)=(s,\zeta_2)$ be an element of $\nbigo(i)$.
\[
 h_i(\xi,\xi):=
 |\zeta_1|^2\cdot \bigl(1+|t|^2\bigr)^{-i}
=|\zeta_2|^2\cdot \bigl( 1+|s|^2\bigr)^{-i}.
\]
For any real numbers $a$ and $b$,
we have the possibly singular metrics $h_{i,(a,b)}$
of $\nbigo(i)$:
Let $\xi=(t,\zeta_1)=(s,\zeta_2)$ be an element of $\nbigo(i)$.
\[
 h_{i,(a,b)}(\xi,\xi):=
 h_i(\xi,\xi)\cdot \bigl(1+|t|^{-2}\bigr)^a\cdot \bigl(1+|t|^2\bigr)^b
=h_i(\xi,\xi)\cdot \bigl(1+|s|^{2}\bigr)^a\cdot \bigl(1+|s|^{-2}\bigr)^b.
\]
Around $|t|=0$,
the order of $h_{i,(a,b)}$ is equivalent to $|t|^{-2a}$.
Around $|s|=0$,
the order of $h_{i,(a,b)}$ is equivalent to $|s|^{-2b}$.
The curvature $R(h_{i,a,b})$ is as follows:
\begin{equation}\label{eq;a11.9.6}
 R(h_{i,a,b})=
 (-a-b+i)\cdot
 \frac{dt\cdot d\tbar}{\bigl(1+|t|\bigr)^2}.
\end{equation}

Let take a point $P\in \proj^1$.
Then we obtain a morphism $\nbigo(i)\lrarr \nbigo(i+1)$
of coherent sheaves.
\begin{lem}
The morphism is bounded with respect to the metrics
$h_{i,(a,b)}$ and $h_{i+1,(a,b)}$.
\hfill\qed
\end{lem}

\subsubsection{Some open subset of the line bundle $\nbigo(-1)$
 with the complete Kahler metric}
\label{subsubsection;10.11.50}

We are mainly interested in the case $i=-1$.
We regard $\nbigo(-1)$ as a complex manifold.
The open submanifold
$Y$ is defined to be
$\bigl\{\xi\in\nbigo(-1)\,\big|\,h_{-1,(0,0)}(\xi,\xi)<1\bigr\}$.
We denote the naturally defined projection of $Y$
onto $\proj^1$ by $\pi$.
We denote the image of the $0$-section $\proj^1\lrarr Y$
by $\proj^1$.
Then we have the normal crossing divisor
$D'=\proj^1\cup \pi^{-1}(0)\cup \pi^{-1}(\infty)$
of $Y$.
The manifold $Y-D'$ is same as 
$\bigl\{(t,x)\in \cnum^{\ast\,2}\,\big|\,
 |x|^2\bigl(1+|t|^2\bigr)<1\bigr\}$.

We have the complete Kahler metric $g:=g_1+g_2+g_3$ of $Y-D'$
given as follows:
As a contribution of the $0$-section $\proj^1$,
we put $\tau_1=-\log\Bigl[\bigl(1+|t|^2\bigr)\cdot |x|^2\Bigr]$,
and as follows:
\[
 g_1:=\frac{1}{\tau_1^2}
 \left(
 \frac{\tbar\cdot dt}{1+|t|^2}
+\frac{dx}{x}
 \right)
\cdot
 \left(
 \frac{t\cdot d\tbar}{1+|t|^2}
+\frac{d\xbar}{\xbar}
 \right)
+\frac{1}{\tau_1}
 \frac{dt\cdot d\tbar}{\bigl(1+|t|^2\bigr)^2}.
\]
\begin{lem}
Note that $g_1$ gives the complete Kahler metric of 
$Y-\proj^1$.
It is equivalent to the Poincar\`{e} metric
around the divisor $\proj^1$.
\hfill\qed
\end{lem}

As a contribution of $\pi^{-1}(\infty)$,
we put $\tau_2=\log\bigl(1+|t|^2\bigr)$, and as follows:
\[
 g_2=\frac{1}{\tau_2}
 \Bigl(
 -1+\frac{|t|^2}{\tau_2}
 \Bigr)
 \cdot
 \frac{dt\cdot d\tbar}{\bigl(1+|t|^2\bigr)^2}.
\]

\begin{lem}\label{lem;b11.10.1}
We have $-1+|t|^2\cdot \tau_2^{-1}>0$.
Around $|t|=\infty$, or equivalently, around $|s|=0$,
the $g_2$ is similar to $(-|s|\log|s|)^{-2}ds\cdot d\sbar$.
Around $|t|=0$,
we have $g_2=\bigl(2^{-1}+o(|t|^2)\bigr)\cdot dt\cdot d\tbar$.
\end{lem}
\pf
It can be checked by a direct calculation.
\hfill\qed

\vspace{.1in}

As the contribution of the divisor $\pi^{-1}(0)$,
we put $\tau_3:=\log(1+|t|^2)-\log|t|^2=\log(1+|s|^2)$,
where we use $s=t^{-1}$.
And we put as follows:
\[
 g_3=\frac{1}{\tau_3}
 \cdot
 \Bigl(
 -1+\frac{|s|^2}{\tau_3}
 \Bigr)
 \frac{ds\cdot d\sbar}{\bigl(1+|s|^2\bigr)^2}.
\]
By the symmetry, the behaviour of $g_3$ is similar to $g_2$.
(See Lemma \ref{lem;b11.10.1})

The following lemma can be checked directly.
\begin{lem}
The metric $g$ gives the complete Kahler metric of 
the complex manifold $Y-D'$.
Around the divisors $\proj^1$,
$\pi^{-1}(0)$ and $\pi^{-1}(\infty)$,
the behaviours of the metric $g$
are equivalent to the Poincar\'{e} metric.
\hfill\qed
\end{lem}

We note the following formulas:
\[
 \begin{array}{l}
 {\displaystyle
 \delbar\del\log \tau_1=
 \frac{1}{\tau_1^2}
 \Bigl(
  \frac{\tbar\cdot dt}{1+|t|^2}+\frac{dx}{x}
 \Bigr)
 \wedge
 \Bigl(
  \frac{t\cdot d\tbar}{1+|t|^2}+\frac{d\xbar}{\xbar}
 \Bigr)
 +\frac{1}{\tau_1}
 \frac{dt\wedge d\tbar}{\bigl(1+|t|^2\bigr)^2}=:\omega_1,}\\
\mbox{{}}\\
{\displaystyle
 \delbar\del\log \tau_2=
 \frac{1}{\tau_2}
 \Bigl(
 -1+ \frac{|t|^2}{\tau_2}
 \Bigr)
\cdot
 \frac{dt\wedge d\tbar}{\bigl(1+|t|^2\bigr)^2}=:\omega_2,}\quad
{\displaystyle
 \delbar\del\log \tau_3=
 \frac{1}{\tau_3}
 \Bigl(
 -1+ \frac{|s|^2}{\tau_3}
 \Bigr)
\cdot
 \frac{ds\wedge d\sbar}{\bigl(1+|s|^2\bigr)^2}=:\omega_3.}
 \end{array}
\]
We put $\omega=\omega_1+\omega_2+\omega_3$.
We put as follows:
\[
 H_0=
 \frac{1}{\tau_1}
+\frac{1}{\tau_2}
 \Bigl(
   -1+\frac{|t|^2}{\tau_2}
 \Bigr)
+\frac{1}{\tau_3}
 \Bigl(
  -1+\frac{|s|^2}{\tau_3}
 \Bigr)>0.
\]
Then we have the following:
\[
 \omega^2=\det(g)\cdot dt\wedge d\tbar \wedge dx\wedge d\xbar
= \Bigl(
 \frac{1}{\tau_1^2\cdot |x|^2\cdot \bigl(1+|t|^2\bigr)^2}
 \times H_0
 \Bigr)
 \cdot dt\wedge d\tbar \wedge dx\wedge d\xbar.
\]
We put as follows:
\[
 H_1:=\frac{H_0}{\bigl(1+|t|^2\bigr)\cdot \bigl(1+|s|^2\bigr)}
\]
Recall that we have $\Ric(g)=\delbar\del(\det(g))$.
\begin{lem}\mbox{{}}
\begin{itemize}
\item
Let $C$ be a number such that $0<C<1$.
On the domain
$\bigl\{(t,x)\in\cnum^{\ast\,2}\,\big|\,
 |x|^2\cdot\bigl(1+|t|^2\bigr)\leq C\bigr\}$,
we have the following similarity of the behaviour:
\[
\begin{array}{ll}
 H_1\sim (\log|t|)^{-2}, & (|t|\to\infty, \mbox{\rm or, } |t|\to 0),\\
\mbox{{}}\\
 H_1\sim (-\log |x|)^{-1}, & (|x|\to 0).
\end{array}
\]
\item
We have the equality:
$Ric(g)-\delbar\del\log(H_1)=-\delbar\del\log\tau_1^2$.
\hfill\qed
\end{itemize}
\end{lem}

%% file: a49.2.tex

\subsubsection{The inequality and the vanishing for the acceptable bundle}

\label{subsubsection;10.11.71}

We put $\Delta_z^2=\{(z_1,z_2)\,|\,|z_i|<1\}$
and $D_i'=\{z_i=0\}\subset \Delta_z^2$.
Let $\varphi:\blowup{\Delta_z^2}\lrarr \Delta_z^2$
denote the blow up of $\Delta_z^2$ at the origin $O=(0,0)$.
We have the exceptional divisor $\varphi^{-1}(O)$,
the proper transforms $\blowup{D'}_i$ of $D'_i$.

We put $\blowup{X}=\blowup{\Delta_z^2}\times\Delta_w^{n-2}$.
Then we have the composite $\psi$ of the natural morphisms:
\[
\begin{CD}
 \blowup{X}@>{\varphi\times id}>>
 \Delta^2_z\times\Delta^{n-2}_w
 @>>>\Delta_z^{n}.
\end{CD}
\]
Here the latter morphism is the natural isomorphism
given by $w_{i}=z_{i+2}$.
We put $\blowup{D}:=\psi^{-1}(D)$, which is same as the following:
\[
 \bigl(
 \varphi^{-1}(0,0)\cup \blowup{D_1'}\cup\blowup{D_2'}
 \bigr)\times \Delta_w^{n-2}
\cup
 \Delta^2_z\times
 \Bigl(
 \bigcup_{i=1}^{n-2}\{w_i=0\}
 \Bigr).
\]
The restriction of $\psi$ to $\blowup{X}-\blowup{D}$
gives an isomorphism $\blowup{X}-\blowup{D}\simeq X-D$.

We can take a holomorphic embedding $\iota$
of $Y$, given in the subsubsection \ref{subsubsection;10.11.50},
to $\widetilde{\Delta}^2$
satisfying the following:
\begin{itemize}
\item The image of the $0$-section $\proj^1$ is the exceptional divisor
 $\phi^{-1}(O)$.
\item
We have $\iota^{-1}(D_1')=\pi^{-1}(\infty)$
and $\iota^{-1}(D_2')=\pi^{-1}(0)$.
\end{itemize}

We put $\overline{X}:=Y\times\Delta_w^{n-2}$.
Then we have the naturally induced morphism
$\overline{X}\lrarr \blowup{\Delta_z^2}\times\Delta_w^{n-2}$,
which we also denote by $\iota$.
We put as follows:
\[
 \overline{D}:=\iota^{-1}(\blowup{D}),\quad
 \overline{X^{(1)}}:=\pi^{-1}(P)\times
 \Delta_w^{n-2},\quad
 \overline{D^{(1)}}:=\overline{X^{(1)}}\cap \overline{D}.
\]
Note that $\iota(\overline{X})$ gives a neighborhood
of $\blowup{X_0}$ in $\blowup{X}$.
The composite $\psi\circ\iota$ is denoted by $\psi_1$.

Let $(E,\delbar_E)$ be a holomorphic bundle
with a hermitian metric $h$ over $X-D$.
We assume that $(E,\delbar_E,h)$ is acceptable
in the sense of Definition \ref{df;11.28.5}.
We denote the curvature of $\psi_1^{\ast}(E,\delbar_E,h)$
by $\psi_1^{\ast}R(h)$.
It is dominated by
$\psi_1^{\ast}\delbar\del\log\tau(\veca,N)$ for sufficiently
negative number $N$.
(See (\ref{eq;a11.9.7}) for the definition of $\tau(\veca,N)$).

Let $\epsilon_i$ $(i=1,2)$ be positive numbers
such that $\epsilon_1+\epsilon_2<1$.
\begin{lem} \label{lem;10.11.55}
We can pick positive numbers $\epsilon$
and negative numbers
$a$ and $b$ satisfying the following:
\[
 -a-b=1,
\quad
 0<-a<1-\epsilon_1-\epsilon,
\quad
 0<-b<1-\epsilon_2-\epsilon.
\]
\end{lem}
\pf
It can be checked elementarily.
\hfill\qed

\vspace{.1in}

Let $\epsilon$, $a$ and $b$ be as in Lemma \ref{lem;10.11.55}.
The metric $\tilde{h}$ of
$\psi_1^{\ast}(E)(-\overline{X}^{(1)}):=
\psi_1^{\ast}(E)\otimes\pi^{\ast}\nbigo_{\proj}(-1)$
over the complex manifold $Y-D'$
is given as follows:
\begin{equation} \label{eq;10.11.75}
 \tilde{h}_{N,\epsilon,a,b}:=
 \psi_1^{\ast}h_{\veca,N}\cdot
 h_{-1,a,b}\cdot
 H_1^{-1}\cdot
  \tau_1^{2+\epsilon}
 \bigl(
 \tau_2\cdot\tau_3
 \bigr)^{\epsilon}.
\end{equation}
For simplicity, we use the notation $\tilde{h}$
instead of $\tilde{h}_{N,\epsilon,a,b}$.

\begin{lem} \label{lem;5.22.30}
Let $\epsilon$, $a$ and $b$ be as in Lemma {\rm\ref{lem;10.11.55}}.
When $N$ is sufficiently smaller than $0$,
then the following inequality holds for any
$\eta\in A^{0,1}_c(\psi_1^{\ast}E(-\overline{X}^{(1)}))$:
\[
 \doublelangle
 R(\tilde{h})+\Ric(g),\eta
 \doublerangle_{\tilde{h}}
\geq \epsilon||\eta||_{\tilde{h}}^2.
\]
(See the formula {\rm(\ref{eq;10.11.60})}
for the definition of
$\doublelangle\cdot,\cdot\doublerangle_{\tilde{h}}$.)
\end{lem}
\pf
We have the following equality:
\begin{multline}
 R(\tilde{h})+\Ric(g)
=R(\psi_1^{\ast}h_{\veca,N})+R(h_{-1,a,b})
-\delbar\del\log H_1
+(2+\epsilon)\cdot\delbar\del\log \tau_1
+\epsilon \cdot\delbar\del(\log \tau_2+\log\tau_3)
+\Ric(g)\\
= R(\psi_1^{\ast}h_{\veca,N})
+\epsilon(\omega_1+\omega_2+\omega_3).
\end{multline}
Here we have used 
$R(h_{-1,a,b})=0$ due to (\ref{eq;a11.9.6})
and our choice of $a$ and $b$.
By taking sufficiently negative $N$,
we can assume the following inequality
for any $\eta\in A_c^{0,1}(E)$ on $X-D$:
\[
 \doublelangle
 R(h_{\veca,N}),\eta
 \doublerangle_{\veca,N}\geq\,0.
\]
Then, by a fiberwise linear algebraic argument,
it is easy to see that the following inequality holds
for any $\eta\in A_c^{0,1}(\psi_1^{\ast}(E))$:
\[
 \doublelangle
 \psi^{\ast}R(h_{\veca,N}),\eta
 \doublerangle_{\tilde{h}}\geq\,0.
\]
On the other hand, we obtain
$\big\langle\big\langle
 \omega_1+\omega_2+\omega_3,\eta
 \big\rangle\big\rangle_{\tilde{h}}
\geq \epsilon\cdot||\eta||_{\tilde{h}}$,
which can be checked directly from definition.
Thus we are done.
\hfill\qed

\begin{cor} \label{cor;10.11.56}
Let $\epsilon$, $a$ and $b$ be as in Lemma {\rm\ref{lem;10.11.55}}.
If $N$ is sufficiently negative,
then the first cohomology
$H^1\bigl(
 A^{0,\cdot}_{\tilde{h}}\bigl(\psi_1^{-1}E(-\overline{X}^{(1)})
 \bigr)\bigr)$
vanishes.
\end{cor}
\pf
It immediately follows from Lemma \ref{lem;5.22.30},
Proposition \ref{prop;11.28.6}
and Proposition \ref{prop;12.11.1}.
\hfill\qed

\vspace{.1in}
We remark the following.
\begin{lem}
The contribution of
$h_{-1,a,b}\cdot H_1^{-1}\cdot\tau_1^{2+\epsilon}
 \cdot \bigl(\tau_2\cdot\tau_3\bigr)^{\epsilon}$
to the metric $\tilde{h}$
is equivalent to the following, around $|t|=0$:
\[
 |t|^{-2a}\cdot\bigl(-\log |t|\bigr)^{2}\cdot |t|^{2\epsilon}
 \bigl(-\log|t|\bigr)^{2\epsilon}
=|t|^{-2(a-\epsilon)}\cdot \bigl(-\log|t|\bigr)^{2+2\epsilon}.
\]
We have a similar estimate around $|s|=0$.
The contribution is equivalent to $(-\log |x|)^{3+\epsilon}$
around $x=0$.
\end{lem}
\pf
We have the following equivalences:
\[
 h_{-1,a,b}\sim |t|^{-2a},
\quad
 H_1^{-1}\sim (-\log|t|)^{2},
\quad
 \tau_1^{2+\epsilon}\sim \bigl(-\log |x|^2\bigr)^{2+\epsilon}
\quad
 \tau_2^{\epsilon}\sim |t|^{2\epsilon}
\quad
 \tau_3^{\epsilon}\sim (-\log |t|)^{2\epsilon}.
\]
Thus we are done.
\hfill\qed

%% file: 15.2.tex

\subsubsection{On a family of complexes of Hilbert spaces}
\label{subsubsection;9.10.51}

\begin{lem} \label{lem;9.10.50}
Let $H_i$ $(i=1,2)$ be Hilbert spaces.
Let $F(\lambda):H_1\lrarr H_2$ be bounded morphisms
depending on $\lambda\in\Delta(\lambda_0,\epsilon)$.
Assume the following:
\begin{enumerate}
\item \label{8.18.21}
 $F(\lambda)$ is bounded for any $\lambda\in\Delta(\lambda_0,\epsilon)$.
\item \label{8.18.22}
 There exists a positive constant $C$
 such that
 $||F(\lambda)-F(\lambda')||\leq C\cdot|\lambda-\lambda'|$
 for any $\lambda,\lambda'\in\Delta(\lambda_0,\epsilon)$.
\item \label{8.18.23}
 $F$ is holomorphic with respect to $\lambda$
 in the following sense:
For any $v\in H_1$, $F(\lambda)(v)$ gives a holomorphic
function from $\Delta(\lambda_0,\epsilon)\lrarr H_2$.
(See {\rm\cite{ma}} for Hilbert space valued holomorphic functions,
for example.)
 \item \label{8.18.24}
 $F(\lambda_0)$ is surjective.
\end{enumerate}
Then there exists $\eta>0$ such that the following holds:
\begin{itemize}
\item
 $F(\lambda)$ is surjective for any $\lambda\in\Delta(\lambda_0,\eta)$.
\item
 There exists a family of bounded
 morphisms $\Psi(\lambda):H_1\lrarr H_2$
 for $\lambda\in\Delta(\lambda_0,\eta)$,
 and the following holds:
\begin{enumerate}
 \item
 $\Psi(\lambda)$ is homeomorphic
 for any $\lambda\in\Delta(\lambda_0,\eta)$.
 \item
 There exists a positive constant $C'$
 such that the following holds:
\[
 \max\Bigl\{
 \bigl|\bigl|\Psi(\lambda)-\Psi(\lambda')\bigr|\bigr|,\quad
 \bigl|\bigl|\Psi(\lambda)^{-1}-\Psi(\lambda')^{-1}\bigr|\bigr|
\Bigr\}
  \leq 
 |\lambda-\lambda'|\cdot C'.
\]
\item
 $\Psi(\lambda)$ is holomorphic with respect to $\lambda$.
\item
 The following diagramm is commutative:
\[
 \begin{CD}
 H_1 @>{F(\lambda_0)}>> H_2\\
 @V{\Psi}VV @V{id}VV \\
 H_1 @>{F(\lambda)}>> H_2.
 \end{CD}
\]
\end{enumerate}
\end{itemize}
\end{lem}
\pf
We put $C_1:=\ker F(\lambda_0)$ and $C_2=\ker F(\lambda_0)^{\bot}$.
We put $\varphi:=F(\lambda_0)_{|C_2}:C_2\lrarr H_2$ is bijective
and bounded.
Hence it is homeomorphic,
i.e.,
$\varphi^{-1}$ is also bounded.

We have the bounded morphisms
$a(\lambda):C_1\lrarr H_2$
and $b(\lambda):C_2\lrarr H_2$
defined by
$a(\lambda)=F(\lambda)_{|C_1}$
and $b(\lambda)=F(\lambda)_{|C_2}-\varphi$.
It is easy to check the following:
\begin{itemize}
\item
 $a$ and $b$ are holomorphic with respect to $\lambda$.
\item
 There exists positive constants $C_2$ such that
 $|a(\lambda)-a(\lambda')|\leq C_2\cdot|\lambda-\lambda'|$
 and
 $|b(\lambda)-b(\lambda')|\leq C_2\cdot|\lambda-\lambda'|$.
\item
 $a(\lambda_0)=b(\lambda_0)=0$.
\end{itemize}
Since $\varphi$ is homeomorphic,
there exists $\eta>0$
such that $\varphi+b(\lambda)$ is homeomorphic
for any $\lambda\in\Delta(\lambda_0,\eta)$.
In particular, $F(\lambda)$ is surjective.

The morphism
$\Psi:C_1\oplus C_2\lrarr C_1\oplus C_2$ can be given as follows:
\[
 \left(
 \begin{array}{cc}
 1 & (\varphi+b(\lambda))^{-1}\circ a(\lambda)\\
 0 & (\varphi+b(\lambda))^{-1}\circ \varphi(\lambda)
 \end{array}
 \right).
\]
Then it is easy to check that $\Psi$ has the desired properties.
\hfill\qed

\begin{cor} \label{cor;9.10.70}
Let $F$ and $H_i$ $(i=1,2)$ be as in Lemma {\rm\ref{lem;9.10.50}}.
There exists a positive constant $\eta>0$
and the linear morphism $G(\lambda):\ker F(\lambda_0)\lrarr H_1$
depending on $\lambda\in\Delta(\lambda_0,\eta)$
with the following properties:
\begin{itemize}
\item
It satisfies the conditions {\rm\ref{8.18.21}}, {\rm\ref{8.18.22}}
and {\rm\ref{8.18.23}}
in Lemma {\rm\ref{lem;9.10.50}}.
\item
For any $\lambda\in\Delta(\lambda_0,\eta)$,
$G(\lambda)$ gives a homeomorphism
$\ker F(\lambda_0)\simeq \ker F(\lambda)$.
\end{itemize}
Namely it gives the trivialization of the family
of Hilbert spaces
$\bigl\{\ker(F(\lambda))\,\bigr|\,\lambda\in \Delta(\lambda_0,\eta)\bigr\}$.
\end{cor}
\pf
We have only to restrict $\Psi$ to $\Ker(F(\lambda_0))$.
\hfill\qed

\vspace{.1in}

\noindent
{\bf Reduction procedure}
We have the standard reduction procedure of the family of the complexes
of Hilbert spaces, as is explained in the following.
Let us consider the following family of complexes of Hilbert spaces
$H_i$ $(i=0,1,2)$
depending on $\lambda\in\Delta(\lambda_0,\epsilon)$:
\begin{equation}
\begin{CD}
 H_0
 @>{F_0(\lambda)}>>
 H_1
 @>{F_1(\lambda)}>>
 H_2.
\end{CD}
\end{equation}
Assume the family satisfies the following conditions:
\begin{itemize}
\item
$F_0$ and $F_1$ satisfy the conditions
{\rm\ref{8.18.21}}, {\rm\ref{8.18.22}} and {\rm\ref{8.18.23}}
in Lemma {\rm\ref{lem;9.10.50}}.
\item
$F_1$ satisfies the condition {\rm\ref{8.18.24}}
in Lemma {\rm\ref{lem;9.10.50}}.
\item
The complex at $\lambda_0$ is exact,
i.e.,
$\ker(F_1(\lambda_0))=\Image F_0(\lambda_0)$.
\end{itemize}

Then we can take the family of the morphisms
$\Psi(\lambda):H_1\lrarr H_1$
depending on $\lambda\in\Delta(\lambda_0,\eta)$
as in Lemma \ref{lem;9.10.50}.
We put $F_0'(\lambda):=\Psi^{-1}\circ F_0(\lambda)$.

\begin{lem}
The image of $F_0'(\lambda)$ is contained in
$\Ker F_1(\lambda_0)$.
\end{lem}
\pf
We have $F_1(\lambda_0)\circ F_0'(\lambda)=
F_1(\lambda)\circ \Psi\circ \Psi^{-1}\circ F_0(\lambda)
=F_1(\lambda)\circ F_0(\lambda)=0$.
\hfill\qed

\vspace{.1in}

Hence we obtain the family of morphisms
$F_0''(\lambda):H_0\lrarr \ker(F(\lambda_0))$
depending on $\lambda\in\Delta(\lambda_0,\eta)$.
By our construction,
$F_0''$ satisfies the conditions
\ref{8.18.21}, \ref{8.18.22} and \ref{8.18.23}
in Lemma \ref{lem;9.10.50}.
Since the complex is exact at $\lambda_0$,
the morphism $F_0''$ satisfies the condition \ref{8.18.24}
in Lemma \ref{lem;9.10.50}, too.

\begin{lem} \label{lem;9.10.75}
Let $H_i$ $(i=0,\ldots,n)$ be Hilbert spaces.
For simplicity, we put $H_{n+1}=0$.
Let $d_i(\lambda):H_i\lrarr H_{i+1}$ be family of linear morphisms
satisfying the conditions 
{\rm\ref{8.18.21}}, {\rm\ref{8.18.22}} and {\rm\ref{8.18.23}}
in Lemma {\rm\ref{lem;9.10.50}}.
Assume the following:
\begin{itemize} \label{lem;9.10.71}
\item The higher cohomology groups vanishes at $\lambda_0$.
Namely, we have $\ker (d_{i+1}(\lambda_0))=\Image(d_{i}(\lambda_0))$
for $i=0,\ldots,n-1$.
\end{itemize}
Then there exists a positive constant $\eta$
and the family of linear morphisms
$G(\lambda):\ker d_0(\lambda_0)\lrarr H^0$
$(\lambda\in\Delta(\lambda_0,\eta))$
satisfying the following:
\begin{itemize}
\item The higher cohomology groups vanishes at 
any $\lambda\in\Delta(\lambda_0,\eta)$.
\item
$G(\lambda)$
satisfies the conditions 
{\rm\ref{8.18.21}}, {\rm\ref{8.18.22}} and {\rm\ref{8.18.23}}
in Lemma {\rm\ref{lem;9.10.50}},
and it gives the trivialization of
the family
$\bigl\{\ker d_0(\lambda)\,\big|\,\lambda\in\Delta(\lambda_0,\eta)\bigr\}$,
namely,
$G(\lambda)$ gives the homeomorphism
of $\ker d_0(\lambda_0)$ and $\ker d_0(\lambda)$
for $\lambda\in\Delta(\lambda_0,\eta)$.
\end{itemize}
\end{lem}
\pf
We have only to use the above reduction procedure successively
and Corollary \ref{cor;9.10.70}.
\hfill\qed

%% file: a36.2.tex

\subsubsection{Sobolev spaces}

We recall the following theorem
(See Theorem 9.1 in \cite{p}).

\begin{prop} \label{prop;9.9.14}
Let $p$ and $q$ be real numbers such that $p,q\geq 1$.
Let $k$ and $l$ be real numbers satisfying the following inequality:
$k-\frac{m}{p}\geq l-\frac{m}{q}$.
Then we have the natural inclusion
$L^p_k(\real^m)\subset L_l^q(\real^m)$,
and it is continuous.
If the inequality $k>l$ and $k-\frac{m}{p}>l-\frac{m}{q}$ hold,
then the inclusion is compact.
\end{prop}

\begin{lem} \label{lem;9.9.16}
Let $p$ be a real number such that $p>2d$.
We put as follows for $i=1,\ldots,d$:
\[
 q_i:=
 \left\{
 \begin{array}{ll}
 {\displaystyle \frac{2d}{d-i}
 } &
 {\displaystyle (i<d)
 } \\
 \mbox{{}}\\ 
 p & (i=d)\\
 \end{array}
 \right.
\]
Then we have the continuous inclusion
$L_2^{q_{i-1}}(\real^{2d})\subset L_1^{q_i}(\real^{2d})$.
\end{lem}
\pf
We have the following relation from our choice of $q_i$:
\[
 2-\frac{2n}{q_{i-1}}
 \geq 1-\frac{2n}{q_i}.
\]
Then the lemma immediately follows from \ref{prop;9.9.14}.
\hfill\qed

%% file: b25.1.tex

\subsubsection{Preliminary}

Let $X_0$ be an open subset of $\cnum_z^{m}$.
Assume that we are given the holomorphic one forms
$\eta_i\in \Omega_{X_0}^{1,0}$ $(i=1,\ldots,m)$, which gives are frame.
For example, we consider $dz_1,\ldots,dz_m$
or $z_1^{-1}\cdot dz_1,\ldots,z_m^{-1}\cdot dz_m$.
We put $X=\Delta^d_{\zeta}\times X_0$.

Let $\harmonicbundle$ be a harmonic bundle over $X$ of rank $r$.
Then we have the deformed holomorphic bundle $\nbige^{\lambda}$
with the holomorphic structure $\delbar_{\nbigelambda}$
and the $\lambda$-connection $\DD^{\lambda}$.
(See the subsection 3.1 in our previous paper \cite{mochi}).

Let $\vecv$ be a $C^{\infty}$-frame of $E=\nbigelambda$.
The $(0,1)$-form $K\in C^{\infty}(X,\Omega_X^{0,1}\otimes M(r))$
is determined by the following relation:
\[
 \delbar_{\nbigelambda}\vecv
=\vecv\cdot K.
\]
The $(1,0)$-form $A\in C^{\infty}(X,\Omega^{1,0}\otimes M(r))$
is determined by the following relation:
\[
 \DD^{\lambda}\vecv=\vecv\cdot (A+K).
\]

We also have $\Theta\in C^{\infty}(X,M(r)\otimes\Omega^{1,0}_X)$
and $\Theta^{\dagger}\in C^{\infty}(X,M(r)\otimes\Omega^{0,1}_X)$
given by the following relation:
\[
 \theta\vecv=\vecv\cdot \Theta,
\quad\quad
 \theta^{\dagger}\vecv=\vecv\cdot\Theta^{\dagger}.
\]
We decompose as follows:
\[
 \Theta=\sum \Theta_{\zeta_i}\cdot d\zeta_i+
   \sum \Theta_{\eta_j} \cdot\eta_j,
\quad\quad 
 \Theta^{\dagger}=\sum \Theta^{\dagger}_{\zeta_i}\cdot d\overline{\zeta}_i
+\sum \Theta^{\dagger}_{\eta_j}\cdot \overline{\eta}_j.
\]

\begin{lem} \label{lem;9.9.35}
We have the following relation:
\[
 \delbar \Theta+\bigl[K-\lambda\Theta^{\dagger},\Theta\bigr]=0.
\]
\end{lem}
\pf
It follows from the relations
$\delbar_E=\delbar_{\nbigelambda}-\lambda\cdot\theta^{\dagger}$
and 
$\delbar_E\theta=0$.
\hfill\qed

\begin{lem}\label{lem;9.9.36}
We have the following relation:
\[
 \del \Theta^{\dagger}
+\lambda^{-1}\bigl[
 A-\Theta,\Theta^{\dagger} \bigr]=0.
\]
\end{lem}
\pf
It follows from the relation
$\del_E=\lambda^{-1}\cdot
\big(\DD^{\lambda}-\delbar_{\nbigelambda}-\theta\big)$
and $\del_E\theta^{\dagger}=0$.
\hfill\qed

\vspace{.1in}

Let $\varphi$ be a $C^{\infty}$-function.
\begin{lem}
We have the following formula:
\begin{equation} \label{eq;9.9.38}
 \begin{array}{l}
 {\displaystyle
 \delbar\bigl(
 \varphi^{2b}\cdot\Theta
 \bigr)
=-\bigl[
 \varphi^b K-\lambda\!\cdot\!\varphi^b\Theta^{\dagger},
 \,\,
 \varphi^b\Theta
 \bigr]
-2\cdot\delbar\varphi^b\cdot\varphi^b\Theta,
 }\\
 \mbox{{}}\\
 {\displaystyle
 \del(\varphi^{2b}\Theta^{\dagger})
 =-\lambda^{-1}\cdot
 \bigl[
 \varphi^{b} A-\varphi^b\Theta,\,\,
 \varphi^b\Theta^{\dagger}
 \bigr]
-2\cdot\del \varphi^b\cdot\varphi^b\Theta^{\dagger}
 }
\end{array}
\end{equation}
\end{lem}
\pf
It immediately follows from Lemma \ref{lem;9.9.35}
and Lemma \ref{lem;9.9.36}.
\hfill\qed

\begin{lem}
We have the following relation for 
$\zeta=\zeta_j$ and $a=\zeta_j,z_k$:
\begin{equation} \label{eq;9.9.39}
\begin{array}{l}
{\displaystyle
 \delbar_{\zeta}\bigl(
 \varphi^{2b}\Theta_a
 \bigr)
=-\bigl[ \varphi^bK_{\zeta}-\lambda\cdot \varphi^b\Theta^{\dagger}_{\zeta},
 \,\,
 \varphi^b\Theta_a
 \bigr]
-2\cdot\delbar_{\zeta}\varphi^b\cdot\varphi^b\Theta_a,}\\
\mbox{{}}\\
 {\displaystyle
 \del_{\zeta}\bigl(
 \varphi^{2b}\Theta^{\dagger}_a
 \bigr)
=-\lambda^{-1}\cdot\bigl[
 \varphi^bA_{\zeta}-\varphi^b\Theta_{\zeta},\,\,
 \varphi^b\Theta^{\dagger}_a
 \bigr]
-2\cdot\delbar_{\zeta}\varphi^b\cdot\varphi^b\Theta^{\dagger}_a.
 }
\end{array}
\end{equation}
\end{lem}
\pf
It immediately follows from (\ref{eq;9.9.38})
\hfill\qed

\begin{lem}
We have the following formulas:
\begin{equation}\label{eq;9.9.40}
\begin{array}{ll}
{\displaystyle
 \del_{\zeta}\delbar_{\zeta}
 \bigl(\varphi^{2b}\Theta
 \bigr)
=  }&
 {\displaystyle
 -\bigl[
 \del_{\zeta}\bigl(\varphi^b K_{\zeta}\bigr)
 -\lambda\cdot\del_{\zeta}\varphi^b\Theta^{\dagger}_{\zeta},\,\,
 \varphi^b\Theta_a
 \bigr]
-\bigl[
 \varphi^bK_{\zeta}-\lambda\cdot\varphi^b\Theta^{\dagger}_{\zeta},\,\,
 \del_{\zeta}\varphi^b\Theta_a
 \bigr]
-2\cdot\del_{\zeta}\delbar_{\zeta}\varphi^{b}\cdot
 \varphi^b\Theta_a } \\
 &
 {\displaystyle
-2\cdot\delbar_{\zeta}\varphi^b\cdot
 \del_{\zeta}\big(\varphi^b\Theta_a\big),
 }\\
\mbox{{}}\\
 {\displaystyle
 \delbar_{\zeta}\del_{\zeta}
 \bigl(\varphi^{2b}\Theta^{\dagger}
 \bigr)
= }&
{\displaystyle
 -\lambda^{-1}\cdot\bigl[
 \delbar_{\zeta}(\varphi^b A_{\zeta})
 -\delbar_{\zeta}(\varphi^b\Theta_{\zeta}),\,\,
 \varphi^b\Theta^{\dagger}
 \bigr]
-\lambda^{-1}\cdot\bigl[
 \varphi^b A_{\zeta}-\varphi^b\Theta_{\zeta},\,\,
 \delbar_{\zeta}(\varphi^b\Theta^{\dagger})
 \bigr] } \\
 &{\displaystyle
-2\cdot\delbar_{\zeta}\del_{\zeta}\varphi^b\cdot\varphi^b\Theta^{\dagger}
-2\cdot\del_{\zeta}\varphi^b\cdot\delbar_{\zeta}(\varphi^b\Theta^{\dagger}).
 }
\end{array}
\end{equation}
\end{lem}
\pf
It immediately follows from (\ref{eq;9.9.39}).
\hfill\qed

\begin{lem} \label{lem;9.9.11}
Let $f$ be a compact support $C^{\infty}$-function on $X$.
Then we have the following equality:
\[
 \int (\delbar f,\delbar f)
=\int(\del f,\del f).
\]
\end{lem}
\pf
It follows from the following:
\[
 \int(\delbar f,\delbar f)
=-\int(f,\Delta(f))
=\int(\del f,\del f).
\]
Here we have used the Kahler identity
$\delbar^{\ast}\delbar=\Delta=\del^{\ast}\del$.
\hfill\qed

\subsubsection{An estimate}
\label{subsubsection;10.2.16}

Let $P$ be a point of $X_0$.
A $L_k^p$-function spaces on $\Delta^d_{\zeta}\times\{P\}$
is denoted by $L_{k,P}^p$.
If $k=0$, we use the notation $L_P^p$.
For a $C^{\infty}$-function $F$ on $X$,
we obtain the restriction
$F_{|P}:=F_{|\Delta^d_{\zeta}\times \{P\}}$,
and the norm $\big|\big|F_{|P}\big|\big|_{L_{k}^p}$.
For simplicity, we denote it
by $\big|\big|F\big|\big|_{L^p_{k\,P}}$.
We put as follows:
\begin{equation} \label{eq;10.2.9}
 Q(P):=
 ||\Theta||_{L^{\infty}_{P}}
+||\Theta^{\dagger}||_{L^{\infty}_P}
+||K||_{L^{\infty}_P}
+||A||_{L^{\infty}_P}
+1.
\end{equation}

Let us pick a $C^{\infty}$-function $\varphi$
on $\Delta_{\zeta}^d$
satisfying the following:
\[
 0\leq \varphi(P_1)\leq 1,
\quad\quad
 \varphi(\zeta)
 =\left\{
 \begin{array}{ll}
 1 & \bigl(P_1\in \Delta_{\zeta}^{d}(1/3)\bigr) \\
 \mbox{{}}\\
 0 & \bigl(P_1\not\in \Delta_{\zeta}^{d}(2/3)\bigr) \\
 \end{array}
 \right.
\]

\begin{lem} \label{lem;9.9.42}
For any $b\in \seisuu_{>0}$,
we have the following:
\[
\max\Big\{
 ||\varphi^b\Theta||_{L^{\infty}_P},\,\,
 ||\varphi^b\Theta^{\dagger}||_{L^{\infty}_P}
\Big\}
\leq Q(P).
\]
\end{lem}
\pf
It immediately follows from our choice of $Q(P)$.
\hfill\qed

\begin{lem} \label{lem;9.9.43}
For any $b\in 2\cdot \seisuu_{>0}$,
there exists a positive constant $C_1(b)$ such that
the following holds:
\[
 \max\left\{
 ||\delbar_{\zeta}(\varphi^b\Theta)||_{L_{P}^2},\,\,
 ||\del_{\zeta}(\varphi^b\Theta^{\dagger})||_{L_{P}^2}
 \right\}
\leq C_1(b)\cdot Q(P)^2.
\]
\end{lem}
\pf
It follows from (\ref{eq;9.9.39}) and Lemma \ref{lem;9.9.42}.
\hfill\qed

\begin{lem} \label{lem;9.9.45}
For any $b\in 2\cdot \seisuu_{>0}$,
there exists a positive constant $C_1(b)$ such that
the following holds:
\[
 \max\left\{
 ||\del_{\zeta}(\varphi^b\Theta)||_{L_{P}^2},\,\,
 ||\delbar_{\zeta}(\varphi^b\Theta^{\dagger})||_{L_{P}^2}
 \right\}
\leq C_1(b)\cdot Q(P)^2.
\]
\end{lem}
\pf
It follows from Lemma \ref{lem;9.9.43} and Lemma \ref{lem;9.9.11}.
\hfill\qed

\begin{lem}
For any $b\in 2^2\cdot \seisuu_{>0}$,
there exists a positive constant $C_2(b)>0$ such that the following
 holds:
\[
 \max\Bigl\{
 ||\varphi^b\Theta||_{L^2_{2,P}},\,\,
 ||\varphi^b\Theta^{\dagger}||_{L^2_{2,P}}
 \Bigr\}
\leq C_3(b)\cdot Q(P)^3.
\]
\end{lem}
\pf
From (\ref{eq;9.9.40}), 
 Lemma \ref{lem;9.9.45} and Lemma \ref{lem;9.9.42},
we obtain the following for any $b\in 2^2\cdot\seisuu_{>0}$
\[
 \max\Bigl\{
 || \del_{\zeta}\delbar_{\zeta}
 \bigl(\varphi^b\Theta\bigr)_a||_{L^2},\,\,
 ||\del_{\zeta}\delbar_{\zeta}
 \bigl(\varphi^b\Theta^{\dagger}_a\bigr)||_{L^2}
 \Bigr\}
\leq 
C'(b)\cdot Q(P)^3.
\]
Thus we are done.
\hfill\qed

\begin{lem} \label{lem;9.9.48}
Let $b$ be an integer contained in $2^{2+i}\cdot\seisuu_{>0}$.
Let $q_i$ $(i=1,\ldots,d)$ be the numbers given in 
Lemma {\rm\ref{lem;9.9.16}}.
Then there exist positive constants $C_{2+i}(b)$
satisfying the following:
\begin{equation}\label{eq;9.9.46}
 \max\Bigl\{
 || \varphi^{b}\cdot\Theta_a||_{L^{q_i}_{2,P}},\,\,
 || \varphi^{b}\cdot\Theta_a^{\dagger}||_{L^{q_i}_{2,P}}
 \Bigr\}
\leq C_{2+i}(b)\cdot Q(P)^{3+i}.
\end{equation}
\end{lem}
\pf
We use an induction on $i$.
We have already known the claim for $i=0$ holds.
We assume the claim for $i$ holds,
and we will show that the claim for $i+1$ holds.

Due to our hypothesis of our induction,
we have the following inequality:
\[
 \max\Bigl\{
 || \varphi^{b}\cdot\Theta_a||_{L^{q_i}_{2,P}},\,\,
 || \varphi^{b}\cdot\Theta_a^{\dagger}||_{L^{q_i}_{2,P}}
 \Bigr\}
\leq C_{2+i}(b)\cdot Q(P)^{3+i}.
\]
Recall that we have the continuous inclusion
$L_{2}^{q_{i}}\subset L_1^{q_{i+1}}$
due to Lemma \ref{lem;9.9.16}.
Thus we obtain the following:
\[
 \max\Bigl\{
 || \varphi^{b}\cdot\Theta_a||_{L^{q_{i+1}}_{1,P}},\,\,
 || \varphi^{b}\cdot\Theta_a^{\dagger}||_{L^{q_{i+1}}_{1,P}}
 \Bigr\}
\leq C'(b)\cdot Q(P)^{3+i}.
\]
In particular, we obtain the following:
\[
 \max\Bigl\{
 || \del_{\zeta}\varphi^{b}\cdot\Theta_a||_{L^{q_{i+1}}_{P}},\,\,
 || \delbar_{\zeta}\varphi^{b}\cdot\Theta_a^{\dagger}||_{L^{q_{i+1}}_{P}}
 \Bigr\}
\leq C''(b)\cdot Q(P)^{3+i}.
\]
Then we obtain the following
due to Lemma \ref{lem;9.9.42} and (\ref{eq;9.9.40}):
\[
  \max\Bigl\{
 || \delbar\del_{\zeta}\varphi^{2b}\cdot\Theta_a||_{L^{q_{i+1}}_{P}},\,\,
 || \del\delbar_{\zeta}\varphi^{2b}\cdot\Theta_a^{\dagger}||_{L^{q_{i+1}}_{P}}
 \Bigr\}
\leq C^{(3)}(b)\cdot Q(P)^{3+i+1}.
\]
Thus we obtain the inequality (\ref{eq;9.9.46}) for $i+1$,
and the induction can proceed.
\hfill\qed

\begin{cor}
Let $p$ be a real number such that $p>2d=2\dim_{\cnum}X$.
Let $b$ be an integer contained in $2^{2+d}\cdot\seisuu_{>0}$.
Then there exists a positive constant $C_{2+d}$
satisfying the following:
\[
  \max\Bigl\{
 || \varphi^{b}\cdot\Theta_a||_{L^{p}_{2,P}},\,\,
 || \varphi^{b}\cdot\Theta_a^{\dagger}||_{L^{p}_{2,P}}
 \Bigr\}
\leq C_{2+d}(b)\cdot Q(P)^{3+d}.
\]
\end{cor}
\pf
It immediately follows from Lemma \ref{lem;9.9.48}.
\hfill\qed

\begin{lem} \label{lem;9.9.50}
For any $b\in 2^{d+k}\cdot\seisuu_{>0}$,
there exists a positive constant $C_{d+k}(b)$
satisfying the following:
\[
 \max\Bigl\{
 ||\varphi^b\Theta||_{L^p_{k,P}},
 ||\varphi^b\Theta^{\dagger}||_{L^p_{k,P}}
 \Bigr\}
\leq C_{k+d}(b)\cdot Q(P)^{d+k+1}.
\]
\end{lem}
\pf
We can show it by a standard boot strapping argument.
\hfill\qed

\vspace{.1in}

We put $K_0:=\overline{\Delta}_{\zeta}(1/3)^{d}$.
\begin{cor} \label{cor;9.9.61}
Let $k$ be an integer, and $p$ be any sufficiently large number.
There exist positive constants $M$ and $C$
such that
the $L_k^p$-norms and $C^k$-norms of $\Theta_{|K_0\times \{P\}}$ 
and $\Theta^{\dagger}_{|K_0\times\{P\}}$ are dominated
by $C\cdot Q(P)^M$.
\end{cor}
\pf
Let us pick a real number $k_1$
such that we have the continuous inclusion
$L_k^p(\real^d)\subset C^k(\real^d)$.
Then we can pick $b_1$, $C_1$ and $M_1$ such that the following holds,
due to Lemma \ref{lem;9.9.50}:
\[
 \max\Bigl\{
 || \varphi^{b_1}\cdot\Theta||_{L^p_{k_1,P}},
 || \varphi^{b_1}\cdot\Theta^{\dagger}||_{L^p_{k_1,P}}
 \Bigr\}
\leq C_1\cdot Q(P)^{M_1}.
\]
Then we obtain the estimate for $C^k$-norms
of $\varphi^{b_1}\Theta$ and $\varphi^{b_1}\Theta^{\dagger}$.
Since $\varphi$ is identically $1$
on a neighbourhood of $K_0$ by our choice of $\varphi$.
Thus we obtain the result.
\hfill\qed

\subsubsection{Estimate for the differential of $\Theta^{\dagger}_{\zeta}$}

Let us consider the case
$\eta_i=dz_i/z_i$ $(i=1,\ldots,m)$.
We have the vector fields $V_i:=z_i\del/\del z_i$,
and $\del f=\sum V_i(f)\wedge \eta_i$
and $\delbar f=\sum \bar{V}_i(f)\wedge \overline{\eta}_i$.

\begin{lem} \label{lem;9.9.63}
We have the following relation:
\[
 \delbar \Theta^{\dagger}
+\bigl[K-\lambda\cdot\Theta^{\dagger},\,\,
 \Theta^{\dagger}\bigr]
=0.
\]
\end{lem}
\pf
It follows from the relation
$\delbar_E\theta^{\dagger}=0$
and $\delbar_E=\delbar_{\nbigelambda}-\lambda\theta^{\dagger}$.
\hfill\qed

\vspace{.1in}

We obtain the following relation for any $\zeta=\zeta_j$
and $i=1,\ldots,m$:
\[
 \overline{V}_i\Theta^{\dagger}_{\zeta}
-\delbar_{\zeta}\Theta^{\dagger}_{\eta_i}
+\bigl[
 K_{\eta_i}-\lambda\Theta^{\dagger}_{\eta_i},\,\,
 \Theta^{\dagger}_{\zeta}
 \bigr]
-\bigl[
 K_{\zeta}-\lambda\Theta^{\dagger}_{\zeta},\,\,
 \Theta^{\dagger}_{\eta_i}
 \bigr]
=0.
\]
Hence we obtain the following:
\begin{lem}\label{lem;d11.14.30}
There exist positive constants
$C$ and $M$ such that
the functions
$\overline{V}_i\Theta^{\dagger}_{\zeta_j}$
$(i=1,\ldots,m,\,\,j=1,\ldots,d)$
are dominated by $C\cdot Q(P)^M$
on $K_0\times X_0$.
\hfill\qed
\end{lem}

%% file: a36.1.tex

Although
we do not use the argument to take a `limit'
of a sequence of harmonic bundles,
contrary to the previous paper \cite{mochi},
the author thinks that such convergency seems significant
for the study of harmonic bundles.
We can improve
the argument for the convergency in \cite{mochi}
by using the estimate in the subsubsection \ref{subsubsection;10.2.16}.
We explain it in this subsection.
In this subsection, we assume that $\lambda\neq 0$.

\subsubsection{Convergency of 
the sequences of the Higgs fields $\Theta^{(n)}$ and
the adjoint maps $\Theta^{\dagger\,(n)}$}

Let $X$ be $\Delta^{\ast\,l}\times \Delta^{d-l}$,
and $(E^{(n)},\theta^{(n)},h^{(n)})$ be a harmonic bundles on $X$
such that $\rank(E^{(n)})=r$.
Recall that
we have the deformed holomorphic bundles
$(\nbige^{\lambda\,(n)},d^{\twoprime\,\lambda\,(n)},
 \DD^{\lambda\,(n)},h^{(n)})$
on $\nbigx^{\lambda}=\{\lambda\}\times X\subset \nbigx$.
In this subsection, the metric and the measure of $X$
are $\sum_{i=1}^ndz_i\cdot d\zbar_i $
and $\prod_{i=1}^n |dz_i\cdot d\zbar_i|$.

Assume that we are given holomorphic frames
$\vecw^{(n)}=(w_1^{(n)},\ldots,w_r^{(n)})$
of $\nbige^{\lambda\,(n)}$.
Then we have the elements
$\Theta^{(n)}\in C^{\infty}\bigl(X,M(r)\otimes\Omega^{1,0}\bigr)$
determined by
$\theta^{(n)}\vecw^{(n)}=\vecw^{(n)}\cdot\Theta^{(n)}$.
We also have the elements $\Theta^{\dagger\,(n)}\in
 C^{\infty}(X,M(r)\otimes\Omega^{0,1}_X)$ determined by
$\theta^{(n)\dagger}\vecw^{(n)}=
 \vecw^{(n)}\Theta^{\dagger\,(n)}$.

We assume the following condition.
\begin{condition}\mbox{{}}\label{condition;10.2.8}
\begin{enumerate}
\item
We have the holomorphic $\lambda$-connection forms
$A^{(n)}\in \Gamma(X,M(r)\otimes\Omega_X^{1,0})$
determined by $\DD^{\lambda\,(n)}\vecw^{(n)}=\vecw^{(n)}\cdot A^{(n)}$.
Then the sequence $\bigl\{ A^{(n)} \bigr\}$ converges to
$A^{(\infty)}\in \Gamma(X,M(r)\otimes \Omega_X^{1,0})$
on any compact subset $Y\subset X$.
\item
On any compact subset $Y\subset X$,
$\bigl\{
     \Theta^{(n)}\bigr\}$ 
and $\bigl\{\Theta^{\dagger\,(n)}\bigr\}$
are bounded independently of $n$.
\hfill\qed
\end{enumerate}
\end{condition}

Let $Y$ be a compact subset of $X$.
We put as follows:
\[
 Q^{(n)}_Y:=
 ||\Theta^{(n)}||_{L^{\infty}(Y)}
+||\Theta^{\dagger\,(n)}||_{L^{\infty}(Y)}
+||A^{(n)}||_{L^{\infty}(Y)}.
\]
Then we have a positive constant $C_Y$ independent of $n$
such that $Q^{(n)}_Y\leq C_Y$ for any integer $n$
due to Condition \ref{condition;10.2.8}.

\begin{lem} \label{lem;10.2.9}
Let $k$ be a positive number.
There exist positive constants
$C$ and $M$ such that
$L_k^p$-norms and $C^k$-norms of
$\Theta_{|Y}^{(n)}$ and $\Theta^{(n)\dagger}_{|Y}$
are dominated by $C\cdot C_Y^M$.
\end{lem}
\pf
Since we have the boundedness of $Q^{(n)}_Y$,
we can use Corollary \ref{cor;9.9.61}.
Note that we use the holomorphic frames $\vecw^{(n)}$,
the term $||K||_{L^{\infty}}$ is trivial.
\hfill\qed

\begin{cor} \label{cor;10.2.10}
There exists a subsequence $\{n_i\}$
such that
$\bigl\{\Theta_{|Y}^{(n_i)}\bigr\}$ and
$\bigl\{
 \Theta^{(n_i)\dagger}_{|Y}\bigr\}$
are convergent in the $C^{\infty}$-sense.
\end{cor}
\pf
We have only to use Lemma \ref{lem;10.2.9}
and an easy diagonal argument.
\hfill\qed

\begin{prop} \label{prop;10.2.11}
There exists a subsequence $\{n_i\}$
such that
$\bigl\{\Theta_{|Y}^{(n_i)}\bigr\}$ and
$\bigl\{
 \Theta^{(n_i)\dagger}_{|Y}\bigr\}$
are convergent on any compact subset $Y\subset X$
in the $C^{\infty}$-sense.
\end{prop}
\pf
We have only to use Corollary \ref{cor;10.2.10}
and an easy diagonal argument.
\hfill\qed

\subsubsection{The convergency of $\bigl\{H^{(n)}\bigr\}$}
\label{subsubsection;9.9.23}

Besides Condition \ref{condition;10.2.8},
we consider the following additional condition:
\begin{condition}\label{condition;10.2.15}\mbox{{}}
\begin{itemize}
\item
We put $H^{(n)}:=H(h^{(n)},\vecw^{(n)})\in
 C^{\infty}\bigl(X,\nbigh(r)\bigr)$.
On any compact subset $Y\subset X$,
$\bigl\{H^{(n)}\bigr\}$ and $\bigl\{H^{(n)\,-1}\bigr\}$ are bounded
independently of $n$.
Namely we have a constant $C_Y$ depending on $Y$
such that $\big|H^{(n)}\big|<C_Y$ and $\big|H^{(n)\,-1}\big|<C_Y$.
\hfill\qed
\end{itemize}
\end{condition}

Due to Proposition \ref{prop;10.2.11},
we may assume that $\bigl\{\Theta^{(n)}\bigr\}$
and $\bigl\{\Theta^{\dagger\,(n)}\bigr\}$ are convergent
on any compact subset $Y$ in the $C^{\infty}$-sense,
by picking a subsequence.

We have the unitary connection 
$\nabla^{(n)}:=
 \delbar_{\nbige^{\lambda\,(n)}}+\del_{\nbige^{\lambda\,(n)}}$
of $\nbige^{\lambda\,(n)}$.
Let $B^{(n)}$ denote the connection form
of $\nabla^{(n)}$,
i.e.,
$\nabla^{(n)}\vecw^{(n)}=\vecw^{(n)}\cdot B^{(n)}$.
Then we have the following by a standard theory:
\begin{equation} \label{eq;9.9.25}
B^{(n)}=H^{(n)\,-1}\del H^{(n)}. 
\end{equation}
On the other hand,
we have $\nabla^{(n)}=\DD^{\lambda\,(n)}-2\theta^{(n)}$.
Hence we obtain the following relation:
\begin{equation}\label{eq;9.9.26}
 B^{(n)}=A^{(n)}-\lambda^{-1}\bigl(1+|\lambda|^2\bigr)\cdot \Theta^{(n)}.
\end{equation}

Let $Y$ be a compact subset of $X$,
and we pick the function $\varphi$ as in the previous subsubsection.
Let us pick a compact subset $Y\subset X$,
and let pick compact subsets $Y_1$ and $Y_2$ of $X$
such that 
$Y$ is contained in the interior of $Y_1$,
and that $Y_1$ is contained in the interior of $Y_2$.
We can pick an element $\varphi\in C^{\infty}(X,\real)$,
satisfying the following:
\[
 0\leq \varphi(x)\leq 1,\quad
 \varphi(x)=
 \left\{
 \begin{array}{ll}
 1 & (x\in Y)\\
 0 & (x\not\in Y_2)
 \end{array}
 \right.
\]

\begin{lem}
We have the following formula:
\begin{equation} \label{eq;a9.9.26}
 \del\bigl(
 \varphi^{2b}H^{(n)}
 \bigr)
=\varphi^b\cdot H^{(n)}
\cdot\bigl(\varphi^bA^{(n)}
-\varphi^b\lambda^{-1}\bigl(1+|\lambda|^2\bigr)\Theta^{(n)}
 \bigr)
+2\del \varphi^b\cdot\varphi^b H^{(n)}.
\end{equation}
\end{lem}
\pf
We obtain the following relation due to
(\ref{eq;9.9.25}) and (\ref{eq;9.9.26}):
\[
  \del H^{(n)}
=H^{(n)}\cdot 
 \bigl(
 A^{(n)}-2\cdot \lambda^{-1}\bigl(1+|\lambda|^1\bigr)\Theta^{(n)}
 \bigr).
\]
Then (\ref{eq;a9.9.26}) follows immediately.
\hfill\qed

\begin{lem} \label{lem;9.9.27}
Let $\{n_i\}$ be a subsequence of $\{n\}$.
\begin{itemize}
\item
Assume 
$\bigl\{\varphi^b\cdot H^{(n_i)} \bigr\}$ is bounded in $L_k^p$,
and $\bigl\{ \del\bigl(\varphi^b\cdot H^{(n_i)} \bigr)
 \bigr\}$ is bounded in $L_k^p$.
Then $\bigl\{\varphi^b\cdot H^{(n_i)} \bigr\}$ is bounded in $L_k^p$
is bounded in $L_{k+1}^p$.
\item
Assume 
$\bigl\{\varphi^b\cdot H^{(n_i)} \bigr\}$ is bounded in $L_k^p$,
and $\bigl\{ \del\bigl(\varphi^b\cdot H^{(n_i)} \bigr)
 \bigr\}$ is convergent in $L_k^p$.
Then $\bigl\{\varphi^b\cdot H^{(n_i)} \bigr\}$ is bounded in $L_k^p$
is convergent in $L_{k+1}^p$.
\end{itemize}
\end{lem}
\pf
Note that we have
$\overline{\varphi\cdot H_{i\,j}}
=\varphi\cdot H_{j\,i}$.
Hence we have the following:
\[
 \overline{\delbar \bigl(
\varphi\cdot H_{i\,j}\bigr)}
=\del \bigl( \varphi\cdot H_{j\,i}\bigr).
\]
Thus the estimate for $\del H$ implies the estimate
for $\delbar H$.
Thus we obtain the result.
\hfill\qed

\begin{lem} \label{lem;9.9.28}
There exists a subsequence $\{n_i\}$
such that
$\{\varphi^{2b}\cdot H^{(n_i)}\}$ are convergent in $L^p_0$
for any $b\in 2\cdot \seisuu$.
\end{lem}
\pf
Due to (\ref{eq;a9.9.26}) and Lemma \ref{lem;9.9.27},
we have the boundedness of the family
$\{\varphi^{2b}\cdot H^{(n_i)}\}$ are bounded in $L^p_1$.
Then we obtain the result due to the compactness
of the inclusion $L^p_1\subset L^p_0$.
\hfill\qed

\begin{lem} \label{lem;10.2.12}
Let $\{n_i\}$ be the subsequence as in Lemma {\rm \ref{lem;9.9.28}}.
Then $\bigl\{H^{(n_i)}\bigr\}$ are convergent on the compact subset $Y$
in the $L_l^p$-sense for any $l$.
\end{lem}
\pf
Let $k$ be a positive integer.
For any sufficiently large integer $b$,
the family $\{\varphi^b\cdot H^{(n_i)}\}$
are convergent in $L_k^p$-sense,
by using (\ref{eq;a9.9.26}) and the bootstrapping argument.
Since $\varphi=1$ on a neighbourhood of $K$,
we obtain the result.
\hfill\qed

\begin{prop}
Assume the conditions in Condition {\rm\ref{condition;10.2.8}}
and Condition {\rm\ref{condition;10.2.15}} are satisfied.
There exists a subsequence $\{n_i\}$ such that
$\bigl\{H^{(n_i)}\bigr\}$,
$\bigl\{\Theta^{(n_i)}\bigr\}$
and $\bigl\{\Theta^{\dagger\,(n_i)}\bigr\}$
are convergent on any compact subset $Y$
in the $C^{\infty}$-sense.
\end{prop}
\pf
We have already seen the convergency of
$\bigl\{\Theta^{(n)}\bigr\}$
and $\bigl\{\Theta^{\dagger\,(n)}\bigr\}$
for some subsequence in Proposition \ref{prop;10.2.11}.
We have seen the convergency of $\{H^{(n)}\}$
on any compact subset for some subsequence in Lemma \ref{lem;10.2.12}.
Then we have only to use the diagonal argument.
\hfill\qed

%% file: d1.tex

\subsubsection{A description of Higgs field}
\label{subsubsection;04.2.5.150}

Let $X$ be a complex manifold,
and let $\harmonicbundle$ be a harmonic bundle over $X$.
Let $(\nbigelambda,\DDlambda)$ be the corresponding
$\lambda$-connection, and $\DD^{\lambda\,f}$ denote the associated
flat connection.
We have the decomposition $\DD^{\lambda\,f}=d'+d''$ into
the $(1,0)$-part and the $(0,1)$-part.
Let $\delta'$ denote the $(1,0)$-operator
such that $d''+\delta'$ is a unitary connection.
Let $\delta''$ denote the $(0,1)$-operator
such that $d'+\delta''$ is a unitary connection.
Then we have the following formula:
\[
 d'=\del+\lambda^{-1}\cdot\theta,
\quad
 d''=\delbar+\lambda\cdot\theta^{\dagger},
\quad
 \delta'=\del-\overline{\lambda}\cdot\theta,
\quad
 \delta''=\delbar-\overline{\lambda}^{-1}\cdot\theta^{\dagger}.
\]
It can be rewritten as follows:
\[
\begin{array}{ll}
 \del=\lambda\cdot\bigl(1+|\lambda|^2\bigr)^{-1}
\cdot\bigl(\overline{\lambda}\cdot d'+\lambda^{-1}\cdot\delta'\bigr)
 &
 \delbar=\overline{\lambda}\cdot\bigl(1+|\lambda|^2\bigr)^{-1}
 \cdot\bigl(\overline{\lambda}^{-1}\cdot d''+\lambda\cdot\delta''\bigr),\\
 \theta=\lambda\cdot\bigl(1+|\lambda|^2\bigr)^{-1}
 \cdot\bigl(d'-\delta'\bigr),
 &
 \theta^{\dagger}=\overline{\lambda}\cdot\bigl(1+|\lambda|^2\bigr)^{-1}
 \cdot\bigl(d''-\delta''\bigr).
\end{array}
\]

Let $\vecv$ be a flat frame with respect to
the flat connection $\DD^{\lambda\,f}$.
We put $H=H(h,\vecv)$.
We have the $C^{\infty}$-section $\Theta$
and $\Theta^{\dagger}$
of $M(r)\otimes\Omega^{1,0}$ and $M(r)\otimes\Omega^{0,1}$
respectively,
determined by the relations
$\theta\cdot\vecv=\vecv\cdot\Theta$ and
$\theta^{\dagger}\cdot\vecv=\vecv\cdot\Theta^{\dagger}$
respectively.
Then we have the following relation:
\[
 \delbar H_{i\,j}=\delbar h(v_i,v_j)
=h(v_i,\delta'v_j)=-h\bigl(v_i,(d'-\delta')v_j\bigr)
=-h\bigl(v_i,\lambda^{-1}\cdot(1+|\lambda|^2)\theta v_j\bigr)
=-\overline{\lambda}^{-1}\cdot \bigl(1+|\lambda|^2\bigr)
 \cdot H_{i\,j}\cdot\overline{\Theta_{k\,j}}.
\]
Hence we obtain the relation
$\Theta=
 -\lambda\cdot\bigl(1+|\lambda|^2\bigr)^{-1}\cdot
 \overline{H}\cdot\del \overline{H}$.
Similarly, we obtain the relation
$\Theta^{\dagger}=
 -\overline{\lambda}\cdot\bigl(1+|\lambda|^1\bigr)^{-1}
\cdot \overline{H}^{-1}\cdot \delbar\overline{H}$.

Hence we have the following relation:
\[
  \overline{\lambda}\cdot\Theta
+\lambda\cdot\Theta^{\dagger}
=|\lambda|^2\cdot \bigl(1+|\lambda|^2\bigr)^{-1}
\cdot\overline{H}^{-1}\cdot d\overline{H}.
\]

Let $v$ be a real vector field on $X$.
Hence we obtain the following equality, 
by using Lemma \ref{lem;04.1.6.1}:
\begin{multline} \label{eq;04.1.27.1}
 \bigl|
\overline{\lambda}\cdot\theta(v)
+\lambda\cdot\theta^{\dagger}(v)
 \bigr|^2_h
=\bigl|\lambda\bigr|^2\cdot\bigl(1+|\lambda|^2\bigr)^{-1}
\cdot
\tr\Bigl( \overline{H}^{-1}\cdot d\overline{H}(v)
 \cdot \overline{H}^{-1}\cdot d\overline{H}(v)
 \cdot\overline{H}^{-1}\cdot\overline{H}
 \Bigr) \\
=\bigl|\lambda\bigr|^2\cdot\bigl(1+|\lambda|^2\bigr)^{-1}
\cdot
 \tr\Bigl(
 H^{-1}\cdot d H(v)\cdot H^{-1}\cdot dH(v)
 \Bigr).
\end{multline}

\subsubsection{Higgs field as the differential of the twisted map}
\label{subsubsection;04.2.5.151}

Let $\pi:\tilde{X}\lrarr X$ be a universal covering.
Recall that $\PH(r)$ denote the sets of
the $(r\times r)$-positive definite hermitian matrices.
By taking a flat trivialization of the flat bundle
$\pi^{-1}\bigl(\nbigelambda,\DD^{\lambda\,f}\bigr)$,
we obtain the equivariant map $\Psi_h:\tilde{X}\lrarr\PH(r)$,
which is essentially independent of a choice of a flat trivialization.

We have the $\pi_1(X)$-action on $\tilde{X}$.
The monodromy induces endomorphism $\rho:\pi_1(X)\lrarr GL(r)$,
which induces the $\pi_1(X)$-action on $\PH(r)$.
The map $\Psi_h$ is equivariant with respect to
the actions of the fundamental group $\pi_1(X)$.

The map $\Psi_h$ is called the twisted map associated with
$(\nbigelambda,\DD^{\lambda\,f},h)$.
We often regard it as a map
$\Psi_h:X\lrarr \PH(r)/\pi_1(X)$,
although $\PH(r)\big/\pi_1(X)$ is not a good topological space
in general.

Then the equality (\ref{eq;04.1.27.1}) can be reformulated
as follows, by using Lemma \ref{lem;04.1.13.1}:
\begin{equation}\label{eq;04.1.27.20}
 \bigl|
 \overline{\lambda}\cdot\theta(v)
+\lambda\cdot\theta^{\dagger}(v)
 \bigr|^2
=|\lambda|^4\cdot \bigl(1+|\lambda|^2\bigr)^{-1}
\cdot
 \bigl|d\Psi_h(v)\bigr|^2.
\end{equation}

%% file: b30.tex

\subsubsection{$Y$-holomorphic structure}

Let $X$ be a $C^{\infty}$-manifold
and $Y$ be a complex manifold.
Let $V$ be a $C^{\infty}$-vector bundle over $X\times Y$.
Let $p_Y$ denote the projection of $X\times Y$ onto $Y$.

\begin{df}
A $Y$-holomorphic structure of $V$ is defined to be
a differential operator
$ d''_{V,Y}:C^{\infty}(V)
\lrarr C^{\infty}(V\otimes p_Y^{\ast}\Omega^{0,1}_{Y})$
satisfying the following:
\begin{itemize}
\item
 $(d''_{V,Y})^2=0$.
\item
 For any function $f\in C^{\infty}\bigl(X\times Y\bigr)$
and any section $s\in C^{\infty}\bigl(X\times Y,V\bigr)$,
the following holds:
\[
 d''_{V,Y}(f\cdot s)
= \delbar_Y(f)\cdot s +f\cdot d''_{V,Y}s.
\]
\end{itemize}
The pair $(V,d''_{V,Y})$ or $V$ is called
a $Y$-holomorphic vector bundle.
We often denote $d''$ instead of $d''_{V,Y}$.
\hfill\qed
\end{df}

\begin{df}
Let $(V,d'')$ be a $Y$-holomorphic vector bundle.
A $Y$-holomorphic section of $V$ is a $C^{\infty}$-section $s$
of $V$ such that $d''(s)=0$.
\hfill\qed
\end{df}

Let $(V,d''_{V,Y})$ be a $Y$-holomorphic vector bundle.
Then the $Y$-holomorphic structure
on the dual $V^{\lor}$ is naturally defined.
Let $(V^{(i)},d''_{V^{(i)},Y})$ $(i=1,2)$ be $Y$-holomorphic vector bundles
over $X\times Y$.
Then the tensor product and the direct sum
of $(V^{(i)},d''_{V^{(i)},Y})$ $(i=1,2)$ are naturally defined.

\begin{df}
Let $(V^{(i)},d''_{V^{(i)},Y})$ $(i=1,2)$ be $Y$-holomorphic vector bundles
over $X\times Y$.
A morphism of $(V^{(1)},d''_{V^{(1)},Y})$ to
$(V^{(2)},d''_{V^{(2)},Y})$ is defined to be
a $Y$-holomorphic section of
$Hom(V^{(1)},V^{(2)})$.
\hfill\qed
\end{df}

\subsubsection{Some description of
$\proj^1$-holomorphic vector bundle over $X\times\proj^1$}

We are particularly interested in the $\proj^1$-holomorphic vector
bundle over $X\times \proj^1$.
It is often useful to consider the following object.
In the following, we regard $\proj^1$ as the gluing
of $\cnum_{\lambda}$ and $\cnum_{\mu}$
by the relation $\lambda=\mu^{-1}$.

\begin{df}
A patched object on $X$ is defined to be
a tuple $(V_{\lambda},V_{\mu},\varphi)$ as follows:
\begin{itemize}
\item
$V_{a}$ is a $\cnum_{a}$-holomorphic vector bundle over
$X\times \cnum_{a}$ 
for $a=\lambda,\mu$ .
\item
We have the $\cnum_{\lambda}^{\ast}$-holomorphic vector bundles
$V_{\lambda\,|\,X\times\cnum_{\lambda}^{\ast}}$
and $V_{\mu\,|\,X\times\cnum_{\mu}^{\ast}}$
over $X\times\cnum_{\lambda}^{\ast}$.
Then $\varphi$ is a $\cnum_{\lambda}^{\ast}$-holomorphic
isomorphism of 
$V_{\lambda\,|\,X\times\cnum_{\lambda}^{\ast}}$
to $V_{\mu\,|\,X\times\cnum_{\mu}^{\ast}}$.
\hfill\qed
\end{itemize}
\end{df}

\begin{df}
Let $(V^{(i)}_{\lambda},V^{(i)}_{\mu},\varphi^{(i)})$ $(i=1,2)$
be patched objects over $X$.
A morphism $f$ of
$(V^{(1)}_{\lambda},V^{(1)}_{\mu},\varphi^{(1)})$
to
$(V^{(2)}_{\lambda},V^{(2)}_{\mu},\varphi^{(2)})$
is defined to be a tuple
$(f_{\lambda},f_{\mu})$, where $f_a$ $(a=\lambda,\mu)$
is a $\cnum_{a}$-holomorphic morphism
$V^{(1)}_a\lrarr V^{(2)}_{a}$ 
which is compatible with $\varphi^{(1)}$ and $\varphi^{(2)}$.
\hfill\qed
\end{df}

A direct sum, a tensor product and a dual are naturally defined.

\vspace{.1in}
\noindent 
{\bf Equivalence}\\
Let $(V_{\lambda},V_{\mu},\varphi)$ be a patched object over $X$.
By gluing them, we obtain
the $\proj^1$-holomorphic vector bundle over $X\times\proj^1$.
On the other hand,
let $V$ be a $\proj^1$-holomorphic vector bundle
over $X\times\proj^1$,
we put $V_{a}:=V_{|\cnum_{a}}$ for $a=\lambda,\mu$.
We have the naturally defined isomorphism
$\id:
 V_{\lambda\,|\,\cnum_{\lambda}^{\ast}}\lrarr
 V_{\mu\,|\,\cnum_{\mu}^{\ast}}$.
Thus we obtain the patched object.
It is easy to check the equivalence
of the category of the $\proj^1$-holomorphic vector bundle
over $X\times\proj^1$
and the category of the patched object over $X$.
The equivalence preserves a tensor product, a direct sum
and a dual.

\vspace{.1in}
\noindent
{\bf Another description}\\
We can also consider another kind of patched objects
$(V_0,V_1,V_{\infty};\alpha_0,\alpha_{\infty})$ as follows:
For simplicity,
we put $Y_0:=\cnum_{\lambda}$, $Y_{\infty}:=\cnum_{\mu}$
and $Y_1:=\cnum_{\lambda}^{\ast}$.
\begin{itemize}
\item
 $V_a$ be a $Y_{a}$-holomorphic vector bundle
 over $X\times Y_a$.
\item
 We have the $Y_1$-holomorphic
 vector bundles $V_{a\,|\,Y_1}$ $(a=0,1,\infty)$.
 Then $\alpha$ is a $Y_1$-holomorphic morphism
 $V_{a\,|\,Y_1}\lrarr V_1$ $(a=0,\infty)$.
\end{itemize}
We can naturally define the morphism of such patched objects.
We also have a tensor product, a direct sum and a dual of such objects.
As before we have the naturally defined equivalence
of the category of such patched objects and the category
of the $\proj^1$-holomorphic vector bundle over $X\times\proj^1$.

\vspace{.1in}

Due to the equivalences explained above,
we will often use the descriptions
$(V_{\lambda},V_{\mu})$ or
$(V_0,V_{\infty},V_1;\alpha,\beta)$
to denote a $\proj^1$-holomorphic vector bundle
over $X\times\proj^1$,
in the following.

%% file: b30.1.tex

\subsubsection{The involution $\sigma$ and the induced bundles}

Let $\sigma$ denote the morphism $\proj^1\lrarr\proj^1$
given by $[z_0:z_{\infty}]\longmapsto [\bar{z}_{\infty}:-\bar{z}_0]$.
We also use the notation $\sigma$ to denote the following induced
morphisms:
\begin{itemize}
\item
$\cnum_{\lambda}\lrarr\cnum_{\mu}$ given by
$\sigma^{\ast}\mu=-\overline{\lambda}$.
\item
$\cnum_{\mu}\lrarr\cnum_{\lambda}$ given by
$\sigma^{\ast}\lambda=-\overline{\mu}$.
\item
$\cnum_{\lambda}^{\ast}\lrarr\cnum_{\lambda}^{\ast}$ given by
$\lambda\longmapsto -\overline{\lambda}^{-1}$.
\end{itemize}

We take the anti-linear isomorphism
$\varphi:C^{\infty}(X\times\cnum_{\mu})\lrarr
 C^{\infty}(X\times\cnum_{\lambda})$
by $\varphi(f)=\sigma^{\ast}(\overline{f})$. \label{eq;9.19.2}
It is induced by the conjugate map
$\varphi_0:
 \sigma^{\ast}\cnum_{X\times\cnum_{\mu}}\lrarr
 \cnum_{X\times\cnum_{\lambda}}$
given by $a\longmapsto \overline{a}$,
where $\cnum_Y$ denotes the trivial line bundle over
a $C^{\infty}$-manifold $Y$.

Since the morphism $\sigma$ is anti-holomorphic,
we have the naturally induced isomorphism
$\sigma^{\ast}\Omega^{0,1}_{\cnum_{\mu}}
\lrarr \Omega^{1,0}_{\cnum_{\lambda}}$
given by
$\sigma^{\ast}(a\cdot d\overline{\mu})\longmapsto 
 -\sigma^{\ast}(a)\cdot d\lambda$.
On the other hand, 
we have the conjugate
$\Omega^{1,0}_{\cnum_{\lambda}}\lrarr 
 \Omega^{0,1}_{\cnum_{\lambda}}$
given by $a\cdot d\lambda\longmapsto \overline{a}\cdot d\overline{\lambda}$,
which is anti-linear.
As the composite,
we obtain the linear morphism
$\varphi_0:\sigma(\Omega^{0,1}_{\cnum_{\mu}})
  \lrarr \Omega^{0,1}_{\cnum_{\lambda}}$.
They induces the morphism
$\varphi:
C^{\infty}\bigl(
 X\times\cnum_{\mu}, \Omega^{0,1}_{\cnum_{\mu}}
\bigr)
\lrarr
 C^{\infty}\bigl(
 X\times\cnum_{\lambda}, \Omega^{0,1}_{\cnum_{\lambda}}
\bigr)$
given by $f\cdot d\overline{\mu}\longmapsto
 -\sigma^{\ast}(\overline{f})\cdot d\overline{\lambda}$.

Let $f$ be an element of
$C^{\infty}(X\times\cnum_{\mu})$.
Then we obtain the elements
$\varphi(f)\in C^{\infty}(X\times\cnum_{\lambda})$
and
$\delbar_{\lambda}
 \varphi(f)
 \in
 C^{\infty}
 \bigl(
  X\times \cnum_{\lambda},\Omega^{0,1}_{\cnum_{\lambda}}
  \bigr)$.
On the other hand,
we have the elements
$\delbar_{\mu}f\in C^{\infty}\bigl(
 X\times\cnum_{\mu},
 \Omega^{0,1}_{\cnum_{\mu}}
 \bigr)$.
We also have the element
$\varphi\bigl(\delbar_{\mu }f
 \bigr)
\in C^{\infty}\bigl(
  X\times \cnum_{\lambda},\Omega^{0,1}_{\cnum_{\lambda}}
  \bigr)$.

\begin{lem}
We have the equality
$\delbar_{\lambda}\varphi(f)
=\varphi\big(\delbar_{\mu} f \big)$.
Namely, we have 
$\delbar_{\lambda}\circ\varphi=\varphi\circ\delbar_{\mu}$.
\end{lem}
\pf
We have the following:
\begin{equation}\label{eq;9.12.1}
 \varphi\Bigl(
 \frac{\del f}{\del \overline{\mu}}d\overline{\mu}
 \Bigr)
=-\sigma^{\ast}
 \overline{\Bigl(\frac{\del f}{\del \overline{\mu}}\Bigr)}
\cdot
 d\overline{\lambda}
=-\sigma^{\ast}
 \left(
 \frac{\del \overline{f}}{\del \mu}
 \right) d\overline{\lambda}.
\end{equation}
We have the following equalities:
\begin{equation}\label{eq;9.12.2}
 \sigma^{\ast}\left(
 \frac{\del \overline{f}}{\del \mu}
 \right)(z,\lambda)
=\frac{\del \overline{f}}{\del \mu} (z,-\overline{\lambda})
=-\frac{\del \overline{f}(z,-\overline{\lambda})}{\del \overline{\lambda}}
=-\frac{\del}{\del \overline{\lambda}}
 \bigl(
 \varphi(f)(z,\lambda)\bigr).
\end{equation}
From (\ref{eq;9.12.1}) and (\ref{eq;9.12.2}),
we obtain the result.
\hfill\qed

\vspace{.1in}

Let $V$ be a $C^{\infty}$-vector bundle over $\cnum_{\mu}\times X$.
Then we obtain the pull back $\sigma^{\ast}V$
over $\cnum_{\lambda}\times X$.
We change the $\cnum$-vector bundle structure of $\sigma^{\ast}V$
as follows:
For a complex number $a\in \cnum$
and an element $\sigma^{\ast}(v)\in \sigma^{\ast}V_{|(\lambda,P)}$,
$a\cdot \sigma^{\ast}(v)$ is defined to be
$\sigma^{\ast}(\overline{a}\cdot v)$.
We denote the resulted $C^{\infty}$-vector bundle
by $\sigma(V)$.
The anti-linear morphisms $\varphi_0$ given above can be 
regarded as the linear morphisms
$\sigma\cnum_{X\times\cnum_{\mu}}\lrarr
 \cnum_{X\times\cnum_{\lambda}}$
or
$\sigma\Omega^{0,1}_{\cnum_{\mu}}
\lrarr \Omega^{0,1}_{\cnum_{\lambda}}$.

Then the morphism $\varphi_0$ induces
the following isomorphism
\[
 \varphi_0:
 \sigma\bigl(
 V\otimes\Omega^{0,1}_{\cnum_{\mu}}
 \bigr)
\lrarr
 \sigma\bigl(V\bigr)
\otimes
 \Omega^{0,1}_{\cnum_{\lambda}},
\quad
 \varphi_0\bigl(
 \sigma(v\otimes d\overline{\mu})
 \bigr)
=-\sigma(v)\otimes d\overline{\lambda}.
\]
It induces the morphism $\varphi_0:
 C^{\infty}\Bigl(X\times\cnum_{\lambda},
  \sigma\bigl(
 V\otimes\Omega^{0,1}_{\cnum_{\mu}}
 \bigr)
 \Bigr)
\lrarr
 C^{\infty}\Bigl(X\times\cnum_{\lambda},
 \sigma\bigl(V\bigr)
\otimes
 \Omega^{0,1}_{\cnum_{\lambda}}
\Bigr)$.

Let $(V,\delbar_{\mu})$ be a $\cnum_{\mu}$-holomorphic vector bundle
over $X\times\cnum_{\mu}$.
Then the differential operator $\delbar_{\lambda}$
on $\sigma(V)$ is defined as follows:
\[
 \delbar_{\lambda}(\sigma(v))
=\varphi_0\bigl(
 \sigma(\delbar_{\mu}v)
 \bigr).
\]
\begin{lem}
The operator $\delbar_{\lambda}$
gives the $\cnum_{\lambda}$-holomorphic structure
of $\sigma(V)$.
\end{lem}
\pf
We have the following:
\begin{equation}
 \delbar_{\lambda}\bigl(
f\cdot \sigma(v)
 \bigr)
=\delbar_{\lambda}\bigl(
 \sigma(\varphi(f)\cdot v)
 \bigr)
=\varphi_0\bigl(
 \sigma\bigl(
 \delbar_{\mu}(\varphi(f)\cdot v)
 \bigr)
 \bigr)
=\varphi_0\bigl(
 \sigma\bigl(
 \delbar_{\mu}\varphi(f)\cdot v
+\varphi(f)\cdot\delbar_{\mu}v
 \bigr)
 \bigr)
\end{equation}
We have the following:
\begin{equation}
 \varphi_0\sigma
 \bigl(
 \delbar_{\mu}\varphi(f)\cdot v
 \bigr)
=\varphi_0\sigma \bigl(
 \varphi(\delbar_{\lambda}f)
 \bigr)
\cdot\sigma(v)
=\delbar_{\lambda}f\cdot\sigma(v).
\end{equation}
We also have the following:
\[
 \varphi_0\bigl(
 \sigma\bigl(
 \varphi(f)\cdot\delbar_{\mu}v
 \bigr)\bigr)
=f\cdot \varphi_0\bigl(
 \sigma\bigl(\delbar_{\mu}v
 \bigr)
 \bigr).
\]
It is easy to check that $\delbar_{\lambda}^2=0$.
Thus we are done.
\hfill\qed

\vspace{.1in}
Namely we obtain the induced $\cnum_{\lambda}$-holomorphic
vector bundle $\sigma(V)$ from a $\cnum_{\mu}$-holomorphic
vector bundle $V$.
Similarly, a $Y$-holomorphic vector bundle $V$
induces the $\sigma(Y)$-holomorphic vector bundle
$\sigma(V)$ for $Y=\cnum_{\lambda}$, $\cnum_{\mu}$,
$\cnum_{\lambda}^{\ast}$ and $\proj^1$.

%% file: b6.tex

\subsubsection{Torus action}

Let $\rho_0$ denote the $G_m$-action on $\proj^1$
given as follows:
\[
 \rho_0\bigl(t,[z_0:z_{\infty}]\bigr)=[t\cdot z_0:z_{\infty}].
\]
It induces the $G_m$-action on $X\times\proj^1$,
which we also denote by $\rho_0$.
We have the open subsets
$X\times\cnum_{\lambda}$,
$X\times\cnum_{\mu}$
and $X\times\cnum_{\lambda}^{\ast}$
of $X\times\proj^1$.
They are stable with respect to $\rho_0$.

Let $V$ be a $C^{\infty}$-vector bundle over $X\times Y$
for $Y=\cnum_{\lambda}$, $\cnum_{\mu}$, $\cnum_{\lambda}^{\ast}$
or $\proj^1$

When we consider the $G_m$-action on $V$,
we consider only the action which is a lift of $\rho_0$.
Assume that $V$ is holomorphic along the $\proj^1$-direction.
Let $\rho$ be a $G_m$-action on $V$.
It is called holomorphic if
$t^{\ast}\delbar_{\lambda}=\delbar_{\lambda}$
for any $t\in G_m$.

%% file: b6.1.tex

\subsubsection{Rees bundle (one filtration)}

Let $p_0$, $p_{\infty}$ and $p_{1}$
denote the projection of
$X\times \cnum_{\lambda}$,
$X\times\cnum_{\mu}$
and $X\times\cnum_{\lambda}^{\ast}$
onto $X$.
Let $H$ be a $C^{\infty}$-bundle over $X$,
and $F$ be a decreasing $C^{\infty}$-filtration of $H$.
Let us pick a point of $P$
and a small neighbourhood $U$ of $P$.
We may assume that we have a frame $\vecv=(v_i)$
of $H_{|U}$ which is compatible with $F$.
We put $b_i=\deg_F(v_i)=\max\{h\,|\,v_i\in F_h\}$.
We put $\tilde{v}_i:=\lambda^{-b_i}\cdot p_1^{\ast}(v_i)$.
Then we obtain the frame $\tilde{\vecv}:=(\tilde{v}_i)$
of $p_1^{\ast}H_{|U}$
over $\cnum_{\lambda}^{\ast}\times U$.
The frame $\tilde{\vecv}$ gives a prolongation of $p_1^{\ast}H_{|U}$
to a vector bundle $\xi(H_{|U};F)$ over $U\times\cnum_{\lambda}$.

\begin{lem}
The $C^{\infty}$-bundle $\xi(H_{|U},F)$
is independent of a choice of compatible frame $\vecv$.
\end{lem}
\pf
Assume that $\vecu$ is other compatible frame.
We have the following relation:
\[
 u_i=
 \sum b_{j\,i}\cdot v_j.
\]
Here we have $b_{j\,i}=0$ in the case $\deg_F(u_i)>\deg_F(v_j)$.
Then we obtain the following relation.
\[
 \tilde{u}_i
=\lambda^{-b(u_i)}\cdot u_i
=\sum b_{j\,i}\cdot \lambda^{-b(u_i)+b(v_j)}\cdot\lambda^{-b(v_j)}\cdot v_j
=\sum b_{j\,i}\cdot\lambda^{-b(u_i)+b(v_j)}\cdot\tilde{v}_j.
\]
It implies the well definedness 
of $\xi(H_{|U};F)$.
\hfill\qed

\begin{cor}
We obtain the global $C^{\infty}$-vector bundle
$\xi(H;F)$, which is $\cnum_{\lambda}$-holomorphic.
\hfill\qed
\end{cor}

It is characterized locally as follows:
We may assume that $U=X$.
We have the ring $C^{\infty}(X)[\lambda]$.
We have the $C^{\infty}(X)[\lambda]$-module
$C^{\infty}(X,H)[\lambda,\lambda^{-1}]$.
We have the submodule of $C^{\infty}(X,H)[\lambda,\lambda^{-1}]$
given as follows:
\[
 \sum_{p\in\seisuu}
 \lambda^{-p}\cdot C^{\infty}(X,F^p)[\lambda].
\]
By taking $\tilde{\vecv}$ as above,
we can show that it is locally free,
and $\xi(H_{|U},F)$ is the corresponding vector bundle.
The restriction $\xi(H;F)_{|\nbigx^0}$ is naturally isomorphic
to the associated graded vector bundle
$\Gr_F(H)=\bigoplus_p \Gr_F^p(H)$.

The $\rho_0$ can be naturally lifted to
the action on $p_1^{\ast}(H)$.
It is easy to check that the action can be prolonged
to the action on $\xi(H;F)$.
Since $\nbigx^0:=X\times\{0\}$ is fixed by the torus action,
we obtain the weight decomposition
of $\xi(H;F)_{|\nbigx^0}$.
It is given by the decomposition $\Gr_F(H)=\bigoplus_p \Gr_F^p(H)$,
i.e.,
the weight on $\Gr_F^p(H)$ is $p$.

\vspace{.1in}

Let $\filt(X)$ be the category of filtered
$C^{\infty}$-vector bundle over $X$.
For two filtered bundles $(H_i,F_i)$ $(i=1,2)$,
a morphism $F:(H_1,F_1)\lrarr (H_2,F_2)$ is defined to
be a morphism $H_1\lrarr H_2$ preserving filtrations.
Let $\equi(X\times\cnum_{\lambda})$
be the $C^{\infty}$-vector bundle with $G_m$-action
over $X\times\cnum_{\lambda}$.
For equivariant bundles $V_i$ $(i=1,2)$,
a morphism $f:V_1\lrarr V_2$ is defined to be
an equivariant $\cnum_{\lambda}$-holomorphic section
of the equivariant bundle $Hom(V_1,V_2)$.

Let $f:(H_1,F)\lrarr (H_2,F)$ be a morphism.
Let $v$ be a section of $F^p(H_1)$.
Then $f(v)$ is contained in $F^{p}(H_2)$.
Hence $f(\lambda^{-p}\cdot v)$
gives a section of $\xi(H_2,F)$.
Thus we obtain the section $\xi(f)$ of $Hom(V_1,V_2)$.
It is easy to check that $\xi(f)$ is equivariant with respect to
the torus action.
Therefore we obtain the equivariant morphism
$\xi(f):V_1\lrarr V_2$.

Let $(H_i,F)$ $(i=1,2)$ be a filtered vector bundle.
Recall that the filtration of the tensor product $H_1\otimes H_2$
is defined as follows:
\[
 F^p(H_1\otimes H_2):=
 \sum_{r+s\geq p}F^r(H_1)\otimes F^s(H_2).
\]
Let $(H,F)$ be a filtered vector bundle.
Recall that the dual of $(H,F)$ is defined as follows:
\[
 F^p(H^{\lor})=\ker\bigl(H^{\lor}\lrarr (F^{-p+1})^{\lor}\bigr).
\]

\begin{lem} \label{lem;10.5.1}
$\xi$ gives the fully faithful functor from $\filt(X)$ to $\equi(X)$.
It preserves direct sums, tensor products and duals.
\end{lem}
\pf
It is clear that $\xi$ is a functor
and that $\xi$ preserves direct sums.

We pick a splittings $H_1=\bigoplus_p U_{1}^p$
and $H_2=\bigoplus_p U_2^p$ of the filtrations.
Then the decomposition
$H_1\otimes H_2=\bigoplus_{p}\bigoplus_{r+s=p}U_1^r\otimes U_2^s$
gives the splitting of the filtration of $H_1\otimes H_2$.
Then it is easy to see that
the tensor product 
$\xi(H_1,F)\otimes \xi(H_2,F)$
is isomorphic to
$\xi\bigl((H_1,F)\otimes (H_2,F)\bigr)$.

We have the perfect pairing
$H\otimes H^{\lor}\lrarr \cnum$.
Due to the definition of $F(H^{\lor})$,
the composite of the following morphisms is trivial
for $i\geq -p+1$:
\[
 F^i(H)\otimes F^{p}(H^{\lor})
\lrarr
 F^{-p+1}(H)\otimes F^p(H^{\lor})\lrarr \cnum.
\]
For any $v\in F^i(H)$ and $v^{\lor}\in F^p(H^{\lor})$,
the pairing 
$\bigl\langle\lambda^{-p}\cdot v^{\lor},\lambda^{-i}\cdot v
 \bigr\rangle$ is contained in $\cnum[\lambda]$.
It implies we have the equivariant morphism
$\xi(H^{\lor},f)\lrarr \xi(H,f)^{\lor}$.
By using the perfectness of $\Gr_F(H)\otimes\Gr_F(H^{\lor})\lrarr\cnum$,
we obtain that $\xi(H^{\lor},F)$ is isomorphic to $\xi(H,F)^{\lor}$.

Let $\phi:\xi(H_1,F)\lrarr \xi(H_2,F)$ be an equivariant morphism.
Then we obtain the morphism $f=\phi_{|\{1\}\times X}:H_1\lrarr H_2$.
Since $\phi$ is equivariant,
we have $\phi(\lambda^{-p}\cdot v)=\lambda^{-p}f(v)$.
Let $v$ be a section of $F^p(H_1)$.
Since $\phi(\lambda^{-p}\cdot v)=\lambda^{-p}\cdot f(v)$ is a section of
$\xi(H_2,F)$,
$f(v)$ is contained in $F^p(H_2)$.
Thus we are done.
\hfill\qed

\begin{rem}
In fact,
the functor $\xi$ gives the equivalence of two categories.
\hfill\qed
\end{rem}

\subsubsection{Rees bundle (bi-filtration)}

Let $H$ be a $C^{\infty}$-vector bundle over $X$.
Let $F$ and $G$ be filtrations of $H$
in the $C^{\infty}$-category.
As we have already seen,
we obtain the equivariant bundle
$\xi(H;F)$ over $X\times\cnum_{\lambda}$.
By a similar way,
we obtain the equivariant bundle
$\xi(H;G)$ over $X\times\cnum_{\mu}$.
Note that they are isomorphic to $p_1^{\ast}H$
on $X\times\cnum_{\lambda}^{\ast}$.
Hence we can glue them,
and we obtain the equivariant vector bundle
$\xi(H,F,G)$ over $X\times\proj^1$.
Or, we can say that we obtain the patched object
$\bigl(
 \xi(H,F),\xi(H,G),p_1^{\ast}H;
\id,\id
\bigr)$.

Let $\Bifilt(X)$ be the category of bi-filtered $C^{\infty}$-vector
bundle over $X$.
Let $(H_i,F,G)$ $(i=1,2)$ be bi-filtered vector bundles over $X$.
A morphism $f:(H_1,F,G)\lrarr (H_2,F,G)$ is defined to be
a morphism $f:H_1\lrarr H_2$ preserving the filtrations $F$ and $G$.

Let $\Equi(X\times\proj^1)$
be the category of equivariant $\proj^1$-holomorphic
vector bundle over $X\times\proj^1$.
Let $V_i$ $(i=1,2)$ be equivariant vector bundles
over $X\times\proj^1$.
A morphism $f:V_1\lrarr V_2$ is defined to be
an equivariant section of the equivariant bundle
$Hom(V_1,V_2)$.

Let $f:(H_1,F,G)\lrarr (H_2,F,G)$ be a morphism
of bi-filtered vector bundles.
Then we obtain the morphism $\xi(f):\xi(H_1,F,G)\lrarr \xi(H_2,F,G)$.

\begin{prop} \label{prop;9.12.15}
The functor $\xi$ gives an equivalence of two categories.
It preserves direct sums, tensor products and duals.
It is functorial for a $C^{\infty}$-morphism $Y\lrarr X$.
\end{prop}
\pf
We only show the equivalence of the categories,
because the rests are easy.
It follows easily 
from Lemma \ref{lem;10.5.1}
that $\xi$ is fully faithful.
Thus we have only to show that $\xi$ is essentially surjective.

Let $V$ be an equivariant vector bundle over $X\times\proj^1$.
We will construct the bi-filtered vector bundle
$(H,F,G)$.
We put $H=V_{|\nbigx^1}$, where we put $\nbigx^1:=X\times\{1\}$.
We will construct the filtrations $F$ and $G$.
We have only to construct them locally on $X$.

Let us pick a point $P$ of $X$.
The following lemmas are standard.
\begin{lem} \label{lem;9.12.10}
There exists a number $n_0$,
such that
the following holds for any $n\geq n_0$:
\begin{enumerate}
\item \label{9.12.11}
The following morphism is surjective:
\[
 H^0\bigl(V\otimes\nbigo(0,n)_{|\{P\}\times\proj^1}\bigr)
\lrarr
 V\otimes \nbigo(0,n)_{|\{P\}\times\{0\}}=
V_{|\{P\}\times\{0\}}.
\]
\item \label{9.12.12}
$H^1\bigl(V\otimes\nbigo(0,n)_{|\{P\}\times\proj^1}\bigr)=0$.
\hfill\qed
\end{enumerate}
\end{lem}

\begin{lem}
We have a neighbourhood $U$ of $P$ in $X$
satisfying the following:
\begin{itemize}
\item
The properties {\rm\ref{9.12.11}} and {\rm\ref{9.12.12}}
in Lemma {\rm\ref{lem;9.12.10}}
for any point $Q\in U$.
\item
$\bigl\{
 H^0\bigl(V\otimes\nbigo_{\proj^1}(0,n)_{|Q\times\proj^1}\bigr)\,\big|\,
 Q\in U
 \bigr\}$
forms a $C^{\infty}$-vector bundle $\nbige_n$ on $U$.
\end{itemize}
\end{lem}
\pf
It can be shown by arguments
similar to those in the subsubsection
\ref{subsubsection;9.10.51}.
\hfill\qed

\vspace{.1in}

Let us return to the proof of Proposition \ref{prop;9.12.15}.
Then we have the equivariant surjective morphism
$\pi:\nbige_n\lrarr V_{|\nbigx^0}$
over $X$.
If $n<n'$, then we have the commutativity
$\nbige_n\subset\nbige_{n'}\lrarr V_{|\nbigx^0}$.

Let $\vecv=(v_i)$ be a $C^{\infty}$-frame of $V_{|\nbigx^0}$,
which is compatible with the weight decomposition
$V_{|\nbigx^0}=\bigoplus U_h$.
We denote the weight of $v_i$ by $w_i$.
Let $\tilde{v}_i$ be a $C^{\infty}$-section
of the weight $w_i$-part of $\nbige_n$
such that $\pi(\tilde{v}_i)=v_i$.
Note that $\tilde{v}_i$ naturally give 
sections of $V\otimes\nbigo(0,n)$,
and they are $\cnum_{\lambda}$-holomorphic.
Then there exists a neighbourhood $U_1$ of $O$ in $\cnum_{\lambda}$,
and a neighbourhood $U_2$ of $P$ in $X$,
such that $\tilde{\vecv}=(\tilde{v}_i)$ gives
a frame of $V\otimes\nbigo(0,n)_{|U_1\times U_2}$.

Let $\nbige_{n,h}$ denote the weight $h$-part of $\nbige_n$.
The decreasing filtration $F$ is defined as follows:
\[
 F^h:=\Image\Bigl(
 \bigoplus_{k\geq h}\nbige_{n,k}\lrarr H=V_{|\{1\}\times X}
 \Bigr).
\]
\begin{lem} \label{lem;a11.9.10}
We have $\Gr_F^h\simeq U_h$.
\end{lem}
\pf
Let $f$ be a section of the weight $h$-part of $\nbige_n$.
Then there exist $C^{\infty}$-functions
$a_i(x,\lambda)$ which are $\cnum_{\lambda}$-holomorphic,
such that the following holds:
\[
 f=\sum a_i(x,\lambda)\cdot\tilde{v}_i.
\]
Due to the condition on the weight,
the functions $a_i(x,\lambda)$ are described
as $\tilde{a}_i(x)\cdot\lambda^{w_i-h}$
for $C^{\infty}$-functions $\tilde{a}_i(x)$.
Note that $\tilde{a}_i(x)=0$
if $w_i-h<0$.
Then Lemma \ref{lem;a11.9.10} follows.
\hfill\qed

\begin{lem}
The construction of the subbundle $F^h$ of $H$
is independent of a choice of $n$.
\end{lem}
\pf
Recall that we have $\nbige_n\lrarr\nbige_{n'}\lrarr
V_{\nbigx^{\lambda}}$ for any $n\leq n'$ and for any $\lambda$.
The lemma follows easily.
\hfill\qed

\vspace{.1in}

Thus we obtain the well defined filtration $F$ of $H$.
We denote the image of $\tilde{v}_i$ by $\bar{v}_i$.
We have the equivariant morphism
$p_1^{\ast}H\lrarr V_{|X\times\cnum_{\lambda}^{\ast}}$,
given by
$\lambda^{-w_i}\cdot \bar{v}_i\longmapsto
\tilde{v}_{i\,|\,\cnum_{\lambda}^{\ast}\times X}$.
It is prolonged to the equivariant isomorphism
$\xi(H,F)\lrarr V_{|\cnum_{\lambda}\times X}$
over $\cnum_{\lambda}\times X$.

Similarly,
we obtain the filtration $G$
by considering $V\otimes\nbigo(n,0)$
and $V_{|\lambda=\infty}$,
and we have the natural equivariant isomorphism
$\xi(H,G)\lrarr V_{|X\times\cnum_{\mu}}$.
Then we obtain the isomorphism
$\xi(H,F,G)\lrarr V$.
Thus the proof of Proposition \ref{prop;9.12.15}
is accomplished.
\hfill\qed

\begin{lem}
An equivariant subbundle $W$ of $V=\xi(H,F,G)$
corresponds to
a subbundle $H'$ of $H$ with strict filtrations
$F'=F\cap H'$ and $G'=G\cap H'$.
\end{lem}
\pf
Let $H'$ be a subbundle of $H$.
We put $F'=F\cap H'$ and $G'=G\cap H'$.
Then it is easy to check that
$\xi(H',F',G')$ gives a subbundle of $V$.

Let $W$ be a subbundle of $V$.
We put $H'=W_{|\lambda=1}$,
and then we obtain $\xi(H',F',G')$ as above.
Due to the equivariance,
the restrictions of
$W$ and $\xi(H',F',G')$ to $X\times\cnum_{\lambda}^{\ast}$
are same.
It implies $W$ and $\xi(H',F',G')$ are same.
Thus we are done.
\hfill\qed

%% file: b6.2.tex

\subsubsection{Real structure of equivariant bundles}

Let $(V,\rho)$ be an equivariant vector bundle over $X\times\proj^1$.
We have the $G_m$-action $\sigma(\rho)$
on $\sigma(V)$ defined as follows:
\[
 \sigma(\rho)(t)\cdot\bigl(\sigma(v)\bigr):=
 \sigma\bigl(\rho(\overline{t})^{-1}\cdot v\bigr).
\]
On the other hand,
we have the bi-filtered vector space $(H,F,G)$
and the conjugate $(H^{\dagger},G^{\dagger},F^{\dagger})$.

\begin{lem}
If $(V,\rho)=\xi(H,F,G)$,
then we have $\bigl(\sigma(V),\sigma(\rho)\bigr)
=\xi\bigl(H^{\dagger},G^{\dagger},F^{\dagger}\bigr)$.
\end{lem}
\pf
It directly follows from the definition.
\hfill\qed

\begin{df}
A real structure on an equivariant bundle
$(V,\rho)$ is defined to be an equivariant isomorphism
$\iota_V:\bigl(\sigma(V),\sigma(\rho)\bigr)\simeq (V,\rho)$,
such that
$\iota_V\circ\sigma(\iota_V)=\id_V$.
\hfill\qed
\end{df}

\begin{rem}
In this paper,
we do not consider the real structure of non-equivariant
$\proj^1$-holomorphic vector bundle.
See the section {\rm 2} in {\rm\cite{s3}}.
\hfill\qed
\end{rem}

The real structure of $(H,F,G)$ is 
defined to be the isomorphism
$\iota_H:(H,F,G)\lrarr (H^{\dagger},G^{\dagger},F^{\dagger})$.
In other words,
$\iota_H$ is anti-isomorphism $H\lrarr H$
and we have $G=F^{\dagger}$.

\begin{lem}
In the case $(V,\rho)=\xi(H,F,G)$,
a real structure $\iota_V$ of $(V,\rho)$
and a real structure $\iota_H$ of $\xi(H,F,G)$
corresponds
by the relation $\iota_V=\xi(\iota_H)$.
\hfill\qed
\end{lem}

\begin{df}
Let $(W,\rho_W)$ be a vector subbundle of $(V,\rho)$.
If we have $\iota_V\bigl(\sigma(W)\bigr)\subset W$,
the subbundle $W$ is called defined over $\real$.
\hfill\qed
\end{df}

\begin{lem}
An equivariant subbundle $(W,\rho_W)$ is defined over $\real$
if and only if
the  subbundle $H'=W_{|\lambda=1}$ of $H$
is defined over $\real$.
\end{lem}
\pf
We have the action of $\iota_V$
on $C^{\infty}(X,H)[\lambda,\lambda^{-1}]$,
it is equivalent with $\iota_H\otimes \iota_0$.
Here $\iota_0$ is given by $\iota_0(\lambda)=-\lambda^{-1}$.
Let us consider the following property:
\begin{quote}
 $C^{\infty}(X,H')[\lambda,\lambda^{-1}]$
is preserved by $\iota_V=\iota_H\otimes\iota_0$.
\end{quote}
Then it is easy to see that three properties are equivalent.
\hfill\qed

\begin{df}
Let $(V_i,\rho_i,\iota_{i})$ $(i=1,2)$ be equivariant vector bundles
defined over $\real$.
Let $f:(V_1,\rho_1)\lrarr (V_2,\rho_2)$ be an equivariant morphism.
It is called defined over $\real$ if
$\iota_{2}\circ\sigma(f)=f\circ\iota_{1}$.
\hfill\qed
\end{df}

\begin{lem}
Let $(H_i,F_i,G_i)$ $(i=1,2)$ be bi-filtered vector bundles
with real structures $\iota_{H_i}$.
We put $(V_i,\rho_i)=\xi(H_i,F_i,G_i)$.
A morphism $f:(V_1,\rho_1)\lrarr (V_2,\rho_2)$ is defined over $\real$
if and only if
$f=\xi(N)$ and 
$N:(H_1,F_1,G_1)\lrarr (H_2,F_2,G_2)$ is defined over $\real$.
\end{lem}
\pf
From $f$, we obtain the morphism
$f':C^{\infty}(X,H_1)[\lambda,\lambda^{-1}]
\lrarr
 C^{\infty}(X,H_2)[\lambda,\lambda^{-1}]$.
Then the both properties are equivalent
to the property that $f'$ is compatible with $\real$-structure.
\hfill\qed

%% file: b5.tex
\subsubsection{Tate object}

The following patched object is called the Tate object:
\[
 \Tate(n)=
\Bigl(
 \nbigo_{\cnum_{\lambda}}\!\cdot t_0^{(n)}\!,\,\,
 \nbigo_{\cnum_{\mu}}\!\cdot t_{\infty}^{(n)}\!,\,\,
 \nbigo_{\cnum^{\ast}_{\lambda}}\!\cdot t_{1}^{(n)}\!,\,\,
 \alpha_0^{(n)}\!,\,\,\alpha_{\infty}^{(n)}
\Bigr).
\]
Here the morphisms $\alpha_0^{(n)}$ and $\alpha_{\infty}^{(n)}$
are given as follows:
\[
 \alpha_0^{(n)}(t_0^{(n)})=(\sqrt{-1}\lambda)^n\cdot t_1^{(n)},
\quad
 \alpha_{\infty}^{(n)}(t_{\infty}^{(n)})
=\bigl(-\sqrt{-1}\mu\bigr)^n\cdot t_1^{(n)}.
\]
The vector bundle corresponding to the patched object above
can be regarded
as the gluing of $\nbigo_{\cnum_{\lambda}}\cdot t_0^{(n)}$
and $\nbigo_{\cnum_{\mu}}\cdot t_{\infty}^{(n)}$
by the relation
$(\sqrt{-1}\lambda)^{-n}\cdot t_0^{(n)}
=(-\sqrt{-1}\mu)^{-n}\cdot t_{\infty}^{(n)}$,
i.e.,
$(\sqrt{-1}\lambda)^{-2n}\cdot t_0^{(n)}=t_{\infty}^{(n)}$.
Hence it is isomorphic to $\nbigo_{\proj^1}(-2n)$.
We denote the corresponding vector bundle
also by $\Tate(n)$ or $\nbigo(-2n)$ for simplicity.

By the projection $\pi:X\times\proj^1\lrarr\proj^1$,
we obtain the patched object $\pi^{\ast}\Tate(n)$
over $X$,
which we denote by $\Tate(n)_X$ or $\Tate(n)$ for simplicity.

We have the torus action $\rho_{\Tate(n)}$ on $\Tate(n)$
given as follows:
\[
 (t,t^{(n)}_0)\longmapsto t^n\cdot t_0^{(n)},
\quad
 (t,t^{(n)}_{\infty})\longmapsto t^{-n}\cdot t_{\infty}^{(n)},
\quad
 (t,t_1^{(n)})\longmapsto t_1^{(n)}.
\]
We have the isomorphism
$\iota_{\Tate(n)}:\sigma^{\ast}\Tate(n)\lrarr \Tate(n)$
given as follows:
\[
 \sigma^{\ast}(t_0^{(n)})\longmapsto (-1)^n\cdot t_{\infty}^{(n)},
\quad
 \sigma^{\ast}(t_{\infty}^{(n)})\longmapsto (-1)^n\cdot t_{0}^{(n)},
\quad
 \sigma^{\ast}(t_1^{(n)})\longmapsto t_1^{(n)}.
\]
Note that it is well defined as it can be checked as follows:
\[
\begin{array}{l}
 \sigma^{\ast}t_{\infty}^{(n)}
=\sigma\bigl(
 (\sqrt{-1}\lambda)^{-2n}\cdot t_0^{(n)}
 \bigr)
=(\sqrt{-1}\mu)^{-2n}\cdot\sigma^{\ast}t_{0}^{(n)}
\longmapsto
 (\sqrt{-1}\mu)^{-2n}\cdot (-1)^n\cdot t_{\infty}^{(n)}
=(-1)^n\cdot t_0^{(n)},\\
\mbox{{}}\\
 \sigma^{\ast}t_1^{(n)}
=\sigma^{\ast}\bigl(
 (\sqrt{-1}\lambda)^{-n}\cdot t_0^{(n)}
 \bigr)
=(\sqrt{-1}\mu)^{-n}\cdot\sigma^{\ast}t_0^{(n)}
\longmapsto
 (\sqrt{-1}\mu)^{-n}\cdot (-1)^n\cdot t_{\infty}^{(n)}
=t_1^{(n)}.
\end{array}
\]
\begin{lem}
We have
$\iota_{\Tate(n)}\circ\sigma^{\ast}\bigl(\iota_{\Tate(n)}\bigr)=
 \id_{\Tate(n)}$.
Namely, the morphism $\iota_{\Tate(n)}$ gives the real structure
of $\Tate(n)$.
\end{lem}
\pf
It can be checked by a direct calculation.
\hfill\qed

\vspace{.1in}

We have
the corresponding real structure on $\Tate(n)_{|\lambda=1}$.
In the following,
the real base $t^{(n)}_{1\,|\,\lambda=1}$
is fixed.

\subsubsection{$\nbigo(p,q)$ and $\nbigo(n)$}

\label{subsubsection;a11.11.15}

We have the following patched object:
\[
 \nbigo_{\proj^1}(p,q)=
 \Bigl(
 \nbigo_{\cnum_{\lambda}}\!\cdot\! f_0^{(p,q)}\!,\,\,\,
 \nbigo_{\cnum_{\mu}}\!\cdot\! f^{(p,q)}_{\infty}\!,\,\,\,
 \nbigo_{\cnum^{\ast}_{\lambda}}\!\cdot\! f_1^{(p,q)}\!,\,\,\,
 \alpha_0^{(p,q)}\!,\,\,\alpha_{\infty}^{(p,q)}
 \Bigr).
\]
Here the morphism $\alpha_0^{(p,q)}$ and $\alpha_{\infty}^{(p,q)}$
are given as follows:
\[
 \alpha_0^{(p,q)}\bigl(f_0^{(p,q)}\bigr)=
\bigl(\sqrt{-1}\lambda\bigr)^{-p}\cdot f_1^{(p,q)},
\quad\quad
 \alpha_{\infty}^{(p,q)}(f_{\infty}^{(p,q)})
=\bigl(-\sqrt{-1}\mu\bigr)^{-q}\cdot f_1^{(p,q)}.
\]
Since $\nbigo(p,q)$ is the vector bundle over $\proj^1$
obtained as the gluing of
$\nbigo_{\cnum_{\lambda}}\cdot f^{(p,q)}_{0}$
and $\nbigo_{\cnum_{\mu}}\cdot f^{(p,q)}_{\infty}$
by the relation
$\bigl(\sqrt{-1}\lambda\bigr)^{p+q}\cdot f_0^{p,q}=f_{\infty}^{p,q}$,
it is isomorphic to $\nbigo_{\proj^1}(p+q)$.
We also note that we have the canonical isomorphism
$\phi_{(p,q),(p',q')}:\nbigo(p,q)\simeq \nbigo(p',q')$
in the case $p+q=p'+q'$,
which is given as follows:
\begin{equation}\label{eq;a12.5.1}
 f_0^{(p,q)}\longmapsto f_0^{(p',q')},
\quad\quad
 f_{\infty}^{(p,q)}\longmapsto f_{\infty}^{(p',q')}.
\end{equation}
In this sense, we may also use $f_0^{(n)}$ and $f_{\infty}^{(n)}$
instead of $f_0^{(p,q)}$ and $f_{\infty}^{(p,q)}$,
when we forget the torus action and we have $p+q=n$.

We denote $\nbigv(\nbigo(p,q))$ by $\nbigo(p,q)$
or $\nbigo(p+q)$ for simplicity.
Let $\pi_X:X\times \proj^1\lrarr\proj^1$ denote the projection.
The induced patched objects
$\pi_X^{\ast}\nbigo(p,q)$ are denoted by $\nbigo(p,q)_X$
or simply by $\nbigo(p,q)$.

We have the torus action $\rho_{(p,q)}$ on
$\nbigo(p,q)$ given as follows:
\[
 (t,f_0^{(p,q)})\longmapsto t^p\cdot f_0^{(p,q)},
\quad
 (t,f_{\infty}^{(p,q)})\longmapsto t^{-q}\cdot f_{\infty}^{(p,q)},
\quad
 (t,f_1^{(p,q)})\longmapsto f_1^{(p,q)}.
\]
The isomorphism $\phi_{(p,q),(p',q')}$ is not compatible
with the torus actions,
if $(p,q)\neq (p',q')$.

We have the isomorphism
$\iota_{(p,q)}:\sigma^{\ast}\nbigo(p,q)\simeq \nbigo(q,p)$
given as follows:
\[
\begin{array}{l}
 \sigma^{\ast}(f_0^{(p,q)})\longmapsto
\sqrt{-1}^{p+q} f_{\infty}^{(q,p)},\\
 \mbox{{}}\\
 \sigma^{\ast}(f_{\infty}^{(p,q)})\longmapsto
 (-\sqrt{-1})^{p+q}f_0^{(q,p)},\\
\mbox{{}}\\
 \sigma^{\ast}(f_1^{(p,q)})\longmapsto 
\sqrt{-1}^{q-p}f_1^{(q,p)}.
\end{array}
\]

\begin{lem}
The morphism $\iota_{(p,q)}$ is well defined.
\end{lem}
\pf
We have only to check that
the second and the third correspondences are induced
by the first correspondence as is checked.
We have the following equalities:
\begin{multline}
 \sigma^{\ast}(f_{\infty}^{(p,q)})
=\sigma^{\ast}\bigl(
 (\sqrt{-1}\lambda)^{p+q}\cdot f_0^{(p,q)}
 \bigr)
=(\sqrt{-1}\mu)^{p+q}\cdot\sigma^{\ast}(f_{0}^{(p,q)}) \\
\longmapsto
 (-1)^{p+q} (\sqrt{-1}\lambda)^{-p-q}\cdot 
 (-\sqrt{-1})^{p+q}\cdot f_{\infty}^{(q,p)}
=\sqrt{-1}^{p+q}\cdot f_0^{(q,p)}
\end{multline}
We also have the following:
\[
  \sigma^{\ast}(f_1^{(p,q)})
=\sigma^{\ast}\bigl(
 (\sqrt{-1}\lambda)^p\cdot f_0^{(p,q)}
 \bigr)
=(\sqrt{-1}\mu)^{p}\cdot\sigma^{\ast}f_{0}^{(p,q)}
\longmapsto
 (-1)^p\cdot(\sqrt{-1}\lambda)^{-p}\cdot
 \sqrt{-1}^{p+q}\cdot
 f_{\infty}^{(q,p)}
=\sqrt{q-p} \cdot f_1^{(q,p)}.
\]
Thus the morphism is well defined.
\hfill\qed

\begin{lem}
The morphism $\iota_{(p,q)}$ is compatible with
$\phi_{(p,q),(p',q')}$,
i.e.,
\[
 \iota_{(p',q')}\circ \sigma^{\ast}\bigl(\phi_{(p,q),(p',q')}\bigr)
=\phi_{(p,q),(p',q')}\circ \iota_{(p',q')}.
\]
\end{lem}
\pf
It can be checked by a direct calculation.
\hfill\qed

\begin{lem}
We have the natural isomorphism
\[
 \nbigo(p_1,q_1)\otimes \nbigo(p_2,q_2)
\longmapsto
 \nbigo(p_1+p_2,q_1+q_2).
\]
It is given by
$f_a^{(p_1,q_1)}\otimes f_a^{(p_2,q_2)}\longmapsto
 f_a^{(p_1+p_2,q_1+q_2)}$
for $a=0,1,\infty$.

It is compatible with the morphisms $\phi_{(p,q),(p',q')}$,
$\rho_{(p,q)}$ and $\iota_{(p,q)}$.
\hfill\qed
\end{lem}

\begin{lem}
We have the isomorphism
$\nbigo(-n,-n)\simeq \Tate(n)$
given by
$f_a^{(-n,-n)}\longmapsto t_a^{(n)}$ for $a=0,1,\infty$.
\hfill\qed
\end{lem}

\subsubsection{The description as the Rees bundle}

We have the following description as the Rees bundle.
We put $\cnum(p,q)=\cnum\cdot e^{(p,q)}$.
The decreasing filtrations
$\leftbottom{a}F_{(p,q)}$ $(a=1,2)$ are defined as follows:
\[
 \leftbottom{1}F_{(p,q)}^i=
 \left\{
 \begin{array}{ll}
 0 & (i>p)\\
 \cnum(p,q) & (i\leq p)\\
 \end{array}
 \right.
\quad\quad
  \leftbottom{2}F_{(p,q)}^i=
 \left\{
 \begin{array}{ll}
 0 & (i>q)\\
 \cnum(p,q) & (i\leq q)\\
 \end{array}
 \right.
\]
We have the Rees bundle
$\xi\bigl(\cnum(p,q),\leftbottom{1}F_{(p,q)},\leftbottom{2}F_{(p,q)}
    \bigr)$.
The correspondence
 $e^{(p,q)}\longmapsto f_1^{(p,q)}$
induces
the equivariant isomorphism
$\xi\bigl(\cnum(p,q),\leftbottom{1}F_{(p,q)},\leftbottom{2}F_{(p,q)}\bigr)
\lrarr
 \nbigo_{\proj^1}(p,q)$

We have the isomorphism
$\cnum(p,q)^{\dagger}\lrarr \cnum(q,p)$
given by
$a\cdot\overline{e^{p,q}}
 \longmapsto
\sqrt{-1}^{q-p}\cdot a\cdot e^{q,p}$.
It induces the isomorphism
$\iota_{(p,q)}:\sigma^{\ast}\nbigo(p,q)\lrarr\nbigo(q,p)$.

We have the pairing
$\langle\cdot,\cdot\rangle_{(p,q)}:
 \cnum(p,q)\otimes\cnum(p,q)^{\dagger}\lrarr\cnum(n,n)$,
which is given as follows:
\[
 e^{(p,q)}\otimes \overline{e^{(p,q)}}\longmapsto
\sqrt{-1}^{q-p}\cdot e^{(n,n)}.
\]
Since the real base of $\cnum(n,n)$ is given
by $e^{(n,n)}$,
the pairing $\sqrt{-1}^{p-q}\cdot\langle\cdot,\cdot\rangle_{(p,q)}$
is positive definite.
Note that the pairing $\langle\cdot,\cdot\rangle$ corresponds
to $S_{(p,q)}$ below.

%% file: a87.1.tex

\subsubsection{Polarization of twistor structure}

In the case $p+q=n$,
the canonical pairing
$S_{(p,q)}:\nbigo(p,q)\otimes\sigma^{\ast}\nbigo(p,q)
\lrarr \Tate(-n)$ is defined to be
the composite of the morphisms
$\nbigo(p,q)\otimes\sigma^{\ast}\nbigo(p,q)
\lrarr
 \nbigo(p,q)\otimes\nbigo(q,p)
\lrarr \nbigo(p+q,p+q)=\Tate(-n)$.
More precise correspondence is as follows:
\begin{equation}\label{eq;a12.5.2}
\begin{array}{l}
 f_1^{(p,q)}\otimes \sigma^{\ast}(f_1^{(p,q)})
\longmapsto \sqrt{-1}^{q-p}\cdot f_1^{(p,q)}\otimes f_1^{(q,p)}
\longmapsto \sqrt{-1}^{q-p}\cdot t_1^{(-n)},\\
 \mbox{{}}\\
 f_0^{(p,q)}\otimes \sigma^{\ast}(f_{\infty}^{(p,q)})
\longmapsto(-\sqrt{-1})^{n}\cdot f_0^{(p,q)}\otimes f_0^{q,p}
\longmapsto (-\sqrt{-1})^n\cdot t_0^{(-n)},\\
 \mbox{{}}\\
f_{\infty}^{(p,q)}\otimes\sigma^{\ast}(f_0^{(p,q)})
\longmapsto \sqrt{-1}{}^{n}\cdot f_{\infty}^{(p,q)}\otimes f_{\infty}^{(q,p)}
\longmapsto \sqrt{-1}{}^n \cdot t_{\infty}^{(-n)}.
\end{array}
\end{equation}

\begin{lem}
Let consider the case $p+q=p'+q'=n$.
Under the isomorphism $\nbigo(p,q)\simeq\nbigo(p',q')$
given by {\rm (\ref{eq;a12.5.1})},
we have $S_{(p,q)}=S_{(p',q')}$.
\end{lem}
\pf
It immediately follows from (\ref{eq;a12.5.2}).
\hfill\qed

\vspace{.1in}

When we forget the torus action,
we use the notation $S_{(n)}$ instead of $S_{(p,q)}$.

\begin{lem}
The pairing $S_{(n)}$ is $(-1)^n$-symmetric.
\end{lem}
\pf
We have
$S_{(n)}\bigl(f_0^{(n)}\otimes \sigma(f_{\infty}^{(n)})\bigr)
=(\sqrt{-1})^{-n}\cdot t_0^{(-n)}$
and
$S_{(n)}\bigl(f_{\infty}^{(n)}\otimes\sigma(f_0^{(n)})\bigr)
=(\sqrt{-1})^nt_{\infty}^{(-n)}$.
Hence we have the following equalities:
\[
 \sigma\Bigl(
 S_{(n)}\bigl(f_{\infty}^{(n)}\otimes\sigma(f_0^{(n)})\bigr)
 \Bigr)
=(\sqrt{-1})^{-n}\cdot (-1)^n\cdot t_0^{(-n)}
=(-1)^n\cdot
S_{(n)}\bigl(f_0^{(n)}\otimes \sigma(f_{\infty}^{(n)})\bigr).
\]
Thus we are done.
\hfill\qed

\vspace{.1in}

Let $V$ be a $\proj^1$-holomorphic bundle,
and let $S:V\otimes\sigma(V)\lrarr \Tate(-n)$
be a pairing.
Then the pairing $S_{(i)}$ induces the following pairing:
\[
 S:\bigl(V\otimes\nbigo(i)\bigr)
 \otimes
 \sigma\bigl(V\otimes\nbigo(i)\bigr)
\lrarr \Tate(-n+i).
\]

\begin{df}[Simpson]
Let $V$ be a pure twistor of weight $n$,
and let $S:V\otimes\sigma(V)\lrarr\Tate(-n)$ is a pairing.
\begin{itemize}
\item
In the case $n=0$,
the pairing $S$ is called a polarization
if the induced hermitian pairing
on $H^0(\proj^1,V)$ is positive definite.
\item
For any $n$,
the pairing $S$ is called polarization
if the induced pairing
$S:\bigl(V\otimes\nbigo(-n)\bigr)\otimes
\sigma\bigl(V\otimes\nbigo(-n)\bigr)
\lrarr\Tate(0)$ is polarization
of the pure twistor structure $V\otimes\nbigo(-n)$ of weight $0$.
\hfill\qed
\end{itemize}
\end{df}

\subsubsection{The polarization of dual and the conjugate}

The polarization
$S:\nbigo(n)\otimes\sigma\bigl(\nbigo(n)\bigr)\lrarr\Tate(-n)$
and the isomorphism $\sigma(\Tate(-n))\simeq\Tate(-n)$
induce the pairing
$\sigma(S):\sigma\bigl(\nbigo(n)\bigr)\otimes\nbigo(n)\lrarr\Tate(-n)$.
\begin{lem}
The pairing $\sigma(S)$ is a polarization.
\end{lem}
\pf
The pairing $\sigma(S)$ is the composite of  the following
correspondences:
\[
 \begin{array}{l}
 \sigma(f_0^{(n)})\otimes f_{\infty}^{(n)}
 \longmapsto
 \sigma \bigl((\sqrt{-1})^n\cdot t_0^{(-n)}\bigr)
 \longmapsto (\sqrt{-1})^n\cdot (-1)^n\cdot t_{\infty}^{(-n)}
 =(\sqrt{-1})^{-n}\cdot t_{\infty}^{(-n)}, \\
 \mbox{{}}\\
 \sigma(f_{\infty}^{(n)})\otimes f_0^{(n)}
 \longmapsto
 \sigma\bigl((\sqrt{-1})^n\cdot t_{\infty}^{(-n)}\bigr)
 \longmapsto (\sqrt{-1})^{-n}\cdot (-1)^n\cdot t_0^{(-n)}
 =(\sqrt{-1})^n\cdot t_0^{(-n)}.
 \end{array}
\]
Hence the pairing
$\sigma(S)(-n) $ is given by the composite of the following
correspondences:
\[
 \begin{array}{l}
 \bigl(
 \sigma(f_0^{(n)})\otimes f_{\infty}^{(-n)}\bigr)
 \otimes
 \bigl(f_{\infty}\otimes\sigma(f_0^{(-n)})\bigr)
 \longmapsto
 (\sqrt{-1})^{-n}\cdot t_{\infty}^{(-n)}
 \cdot (\sqrt{-1})^{-n}\cdot t_{\infty}^{(n)}
 =(-1)^n\cdot t_{\infty}^{(0)},\\
 \mbox{{}}\\
 \bigl(
 \sigma(f_{\infty}^{(n)})\otimes f_0^{(-n)}
 \bigr)
\otimes
 \bigl(
 f_0^{(n)}\otimes \sigma(f_{\infty}^{(-n)})
 \bigr)
\longmapsto
 (\sqrt{-1})^n\cdot t_0^{(-n)}\cdot (\sqrt{-1})^n\cdot t_0^{(n)}
=(-1)^n\cdot t_{\infty}^{(0)}.
 \end{array}
\]
A global section $s$ of $\sigma\bigl(\nbigo(n)\bigr)\otimes\nbigo(-n)$
is given as follows:
\[
 s=\sigma\bigl(f_0^{(n)}\bigr)\otimes f_{\infty}^{(-n)}
  =\sigma\bigl((\sqrt{-1}\lambda)^{-n}\cdot f_{\infty}^{(n)}\bigr)
  =(-1)^n\cdot \sigma(f_{\infty}^{(n)})\otimes f_0^{(-n)}.
\]
Hence we obtain $\sigma(S)(-n)(s,\sigma(s))=1$,
which means $\sigma(S)$ is a polarization.
\hfill\qed

\vspace{.1in}

As for the dual,
the induced pairing
$S^{\lor}:\nbigo(-n)\otimes\sigma\bigl(\nbigo(-n)\bigr)\lrarr\Tate(n)$
is clearly the polarization.

\begin{lem}
Let $V$ be a pure twistor of weight $n$,
and $S:V\otimes\sigma\bigl(V\bigr)\lrarr\Tate(-n)$ is a polarization.
Then the following induced pairings are also the polarization:
\[
 \sigma(S):\sigma(V)\otimes V\lrarr\Tate(-n),
\quad\quad
 S^{\lor}:V^{\lor}\otimes\sigma(V^{\lor})\lrarr\Tate(n).
\]
\end{lem}
\pf
Since any polarized pure twistor of weight $n$ is
isomorphic to a direct sum of
$\bigl(\nbigo(n),S\bigr)$,
the lemma can be reduced to the case
$\bigl(\nbigo(n),S\bigr)$.
It has been already checked above.
\hfill\qed

%% file: a70.tex

\subsubsection{Some remarks}

Let $(H,F,G)$ be a bi-filtered vector bundle
and we put $(V,\rho)=\xi(H,F,G)$.
A pairing $S:(V,\rho)\otimes \sigma(V,\rho)\lrarr \Tate(-n)$
corresponds to the pairing
$\langle\cdot,\cdot\rangle:
 (H,F,G)\otimes(H^{\dagger},G^{\dagger},F^{\dagger})
\lrarr \cnum(n,n)$.
The correspondence can be regarded as follows:
$S$ induces the morphism
$C^{\infty}(X,H)[\lambda,\lambda^{-1}]\otimes
 C^{\infty}(X,H^{\dagger})[\lambda,\lambda^{-1}]
\lrarr \cnum(n,n)[\lambda,\lambda^{-1}]$,
and $H\otimes H^{\dagger}\lrarr\cnum(n,n)$.

\begin{lem}
$S$ is defined over $\real$
if and only if $\langle\cdot,\cdot\rangle$ is defined over $\real$.
\hfill\qed
\end{lem}

\begin{lem}
A morphism $V\lrarr V\otimes \Tate(-1)$
corresponds to
a morphism $(H,F,G)\lrarr (H,F,G)\otimes \cnum(1,1)$.
\end{lem}
\pf
It follows from $\Tate(-1)=\xi(\cnum(1,1))$
and the equivalence of the categories
(Proposition \ref{prop;9.12.15}).
\hfill\qed

\vspace{.1in}

Let $f:(V,\rho)\lrarr (V,\rho)\otimes \Tate(-1)$,
which corresponds to
$f_{|1}:=f_{|\lambda=1}:(H,F,G)\lrarr (H,F,G)\otimes \cnum(1,1)$.
\begin{lem} \label{lem;10.5.7}
$S(f\otimes\id)+S\bigl(\id\otimes \sigma(f)\bigr)=0$
if and only if
$\langle f_{|1}\otimes \id\rangle+\langle \id\otimes f_{|1}\rangle=0$.
\hfill\qed
\end{lem}

Recall that we have the isomorphism
$\cnum(1,1)\otimes\cnum(-1,-1)\lrarr\cnum(0,0)$,
given by the following:
\[
 \bigl(
 t^{(-1)}_{1\,|\,\lambda=1}
 \bigr)
\otimes
 \bigl(
 t^{(1)}_{1\,|\,\lambda=1}
 \bigr)
\longmapsto
 t^{(0)}_{1\,|\,\lambda=1}.
\]
The multiplication
$t^{(1)}_{1\,|\,\lambda=1}$
gives the isomorphism
of the vector spaces
$\cnum(1,1)$ and $\cnum(0,0)$,
which preserves the real bases.
Under the isomorphism,
we can identify
the morphisms
$\bigl(
 f\cdot t^{(1)}
 \bigr)_{1\,|\,\lambda=1}$
and $f_{|\lambda=1}$.
We put $N:=-(f\cdot t^{(1)})_{|\lambda=1}$.
\begin{lem} \label{lem;10.5.5}
Under the identification,
we have the following identities:
\[
 \bigl\langle f_{|1}\otimes \id\bigr\rangle
+\bigl\langle \id\otimes \sigma(f)_{|1}\bigr\rangle
=-\Bigl(
 \bigl\langle N\otimes\id\bigr\rangle
+\bigl\langle\id\otimes \sigma(N)\bigr\rangle
 \Bigr).
\]
\[
 (-1)^h\cdot \bigl\langle
 \id\otimes\sigma\bigl(f^h\bigr)_{|\lambda=1}
 \bigr\rangle
=\bigl\langle
 \id\otimes\sigma(N^h)
 \bigr\rangle.
\]
\end{lem}
\pf
We have the following:
\[
   \bigl\langle f_{|1}\otimes \id\bigr\rangle
+\bigl\langle \id\otimes \sigma(f)_{|1}\bigr\rangle
=
 -\bigl\langle\big(f\otimes t_1^{(1)}\big)_{|1}
 \otimes \id\bigr\rangle
-\bigl\langle \id\otimes
 \sigma(f\otimes t_1^{(1)})_{|1}\bigr\rangle \\
=-\bigl(
\bigl\langle N
 \otimes \id\bigr\rangle
+\bigl\langle \id\otimes
 \sigma(N)
 \bigr\rangle
\bigr).
\]
We also have the following:
\[
  (-1)^h\cdot\bigl\langle
 \id\otimes\sigma\bigl(f\bigr)_{|\lambda=1}
 \bigr\rangle
=\bigl\langle
 \id\otimes\sigma\bigl(
 \bigl(-f\cdot t^{(1)}_1\bigr)^h
 \bigr)_{|\lambda=1}
 \bigr\rangle \\
=\bigl\langle
 \id\otimes\sigma(N^h)
 \bigr\rangle.
\]
Thus we are done.
\hfill\qed

%% file: b6.3.tex
\subsubsection{complex Hodge structure and the equivariant twistor structure}

We put $(V,\rho):=\xi(H,F,G)$
for a bi-filtered vector bundle $(H,F,G)$
bundle over a $C^{\infty}$-manifold $X$.

\begin{lem}
$V$ is pure twistor bundle of weight $n$
if and only if $(H,F,G)$ is complex Hodge structure bundle,
i.e.,
$F$ and $G$ are $n$-opposed.
\end{lem}
\pf
If $F$ and $G$ are $n$-opposed,
then we have the decomposition $H=\bigoplus_{p+q=n}H^{p,q}$,
where $H^{p,q}=F^p\cap G^q$.
Then $\xi(H,F,G)\simeq \bigoplus_{p+q=n}H^{p,q}\otimes\nbigo(p,q)$.
Hence it is pure twistor of weight $n$.

We can assume $X$ is a point.
By considering $V\otimes\nbigo(-n,0)$,
we can reduce the problem to the case $n=0$.
We have the weight decomposition
\[
 H^0(\proj^1,V)=\bigoplus_h U_h.
\]
Here $U_h$ denotes the weight $h$-space.
Then $U_h$ gives the subbundle
$U_h\otimes\nbigo_{\proj^1}\subset V$,
which is isomorphic to a direct sum of $\nbigo(h,-h)$,
and $(V,\rho)=\bigoplus_h U_h\otimes\nbigo(h,-h)$.
Then it can be checked that the corresponding filtrations
$F$ and $G$ are $n$-opposed.
\hfill\qed

\begin{cor}
The functor $\xi$ gives the equivalence
of the following two categories:
\begin{itemize}
\item
 The category of equivariant pure twistor bundle of weight $n$.
\item
 The category of complex pure Hodge structure bundle of weight $n$.
\end{itemize}
The functor $\xi$ gives the equivalence
of the following two categories:
\begin{itemize}
\item
 The category of equivariant pure twistor bundle of weight $n$
defined over $\real$.
\item
 The category of real pure Hodge structure bundle of weight $n$.
\hfill\qed
\end{itemize}
\end{cor}

\begin{cor}
 The functor $\xi$ gives the equivalence of the following two
categories:
\begin{itemize}
\item
 The category of equivariant mixed twistor bundle.
\item
 The category of complex mixed Hodge structure bundle.
\end{itemize}
 The functor $\xi$ gives the equivalence of the following two
categories:
\begin{itemize}
\item
 The category of equivariant mixed twistor bundle
 defined over $\real$.
\item
 The category of real mixed Hodge structure bundle.
\hfill\qed
\end{itemize}
\end{cor}

Let $S:(V,\rho)\otimes \sigma(V,\rho)\lrarr \Tate(-n)$
be an equivariant pairing.
Let
$\langle\cdot,\cdot\rangle:
 (H,F,G)\otimes (H^{\dagger},G^{\dagger},F^{\dagger})
 \lrarr \cnum(n,n)$
be the corresponding pairing.

\begin{lem}
The pairing $S$ is a polarization of twistor structure
if and only if
the induced pairing $\langle\cdot,\cdot\rangle$ is a polarization
of Hodge structure.
\end{lem}
\pf
We may assume that $X$ is a point.
We have the following:
\[
 \sigma^{\ast}S\bigl(a\cdot u\otimes \sigma^{\ast}(b\cdot v)\bigr)
=\sigma^{\ast}\Bigl(
 a\cdot\overline{\sigma^{\ast}(b)}\cdot\langle u,v\rangle
 \Bigr)
=\overline{\sigma^{\ast}(a)}\cdot b\cdot
 \overline{\langle u,v\rangle}.
\]
We also have the following equality:
\[
 S\bigl(b\cdot v\otimes \sigma^{\ast}(a\cdot u)\bigr)=
 \overline{\sigma^{\ast}(a)}\cdot b\cdot\langle v,u\rangle.
\]
Hence $S$ is $(-1)^n$-symmetric
if and only if $\langle\cdot,\cdot\rangle$
is $(-1)^n$-hermitian symmetric.

To compare the positive definiteness conditions,
note the following:
\[
 H^0\bigl(\proj^1,V\otimes\nbigo(-n)\bigr)
=\bigoplus_{p+q=n}
 \bigl\{\lambda^{-p+n}\cdot u\,\big|\,
 u\in H^{p,q}
 \bigr\}.
\]
Then we have the following:
\begin{multline}
 S(-n)\Bigl(\lambda^{-p+n}\cdot u\otimes f_1^{(-n,0)}\otimes
 \sigma\bigl(\lambda^{-p+n}\cdot u\otimes f_1^{(-n,0)}\bigr)
 \Bigr) \\
=\langle u,u\rangle
\cdot\lambda^{-p+n}\cdot
 (-\lambda)^{p-n}\cdot
 f_1^{(-n,0)}\otimes \bigl(\sqrt{-1}^n\cdot f_1^{(0,-n)}\bigr)
=\langle u,u\rangle\cdot t_1^{(n)}\cdot
\sqrt{-1}^{p-q}.
\end{multline}
Here we have used the equality
$(-1)^{p-n}\cdot \sqrt{-1}^n=\sqrt{-1}^{p-q}$.
Thus $H^0(S(-n))$ is positive definite
if and only if
$\sqrt{-1}^{p-q}\cdot\langle\cdot,\cdot\rangle$
is positive definite.
Note that the real base of $\cnum(n,n)$ is fixed as $t_1^{(-n)}$.
\hfill\qed

\begin{cor} \label{cor;10.5.6}
The functor $\xi$ gives an equivalence of the following categories:
\begin{itemize}
\item
 The category of the equivariant polarized pure twistor bundle
 of weight $n$.
\item
 The category of polarized complex pure Hodge structure bundle of weight $n$.
\end{itemize}
The functor $\xi$ gives an equivalence of the following categories:
\begin{itemize}
\item
 The category of the equivariant polarized pure twistor bundle
 of weight $n$ defined over $\real$.
\item
 The category of polarized real pure Hodge structure bundle of weight $n$.
\hfill\qed
\end{itemize}
\end{cor}

%% file: b7.tex

\subsubsection{Polarized mixed twistor structures}

\begin{df}
A tuple $(V,W,N,S)$ as
follows is called a polarized mixed twistor structure
of weight $n$ in one variable:
\begin{enumerate}
\item
 $(V,W)$ is a mixed twistor structure.
\item
 The morphism $f:V\lrarr V\otimes \Tate(-1)$ is
 a morphism of mixed twistors, and it is nilpotent.
 The weight filtration of $N$ is denoted by $W(f)$.
\item
 The pairing
 $S:V\otimes \sigma(V)\lrarr \Tate(-n)$ is a morphism
 of mixed twistor satisfying the following:
\[
 S(f\otimes \id)+S(\id\otimes\sigma(f))=0.
\]
Note that we obtain the induced morphism:
\[
 S:\Gr^W_{h+n}\otimes \Gr^W_{-h+n}\lrarr \Tate(-n).
\]
\item
 We have $W_h=W(f)_{h-n}$ for any $h$.
\item
 The induced pairing
 $S(f^h\otimes\id)=
 (-1)^h\cdot S\bigl(\id\otimes \sigma(f^h)\bigr)$
 gives the polarization  of  the primitive part
 $P\Gr^W_{h+n}=\ker\bigl(N^{h+1}:\Gr^W_{h+n}\lrarr \Gr^{W}_{n-h-2}\bigr)$.
\end{enumerate}
If a tuple $(V,W,f,S)$ satisfies the first four conditions,
then it is called a pseudo-polarized mixed twistor structure
of weight $n$ in one variable.
\hfill\qed
\end{df}

Assume $(V,W,f,S)$ is a pseudo-polarized mixed twistor structure
of weight $n$ in one variable.
We put as follows:
\[
 V^{(0)}:=\Gr^W(V),
\quad
 W_h^{(0)}=\bigoplus_{i\leq h}\Gr^W_i(V).
\]
We obtain the induced morphism
$f^{(0)}:\Gr^{W}_i\lrarr \Gr^{W}_{i-2}\otimes \Tate(-1)$,
and $f^{(0)}:V^{(0)}\lrarr V^{(0)}$.
We also obtain the induced morphism
$S^{(0)}:V^{(0)}\otimes\sigma(V^{(0)})
\lrarr \Tate(-n)$.
Then it is easy to check
that $\bigl(V^{(0)},W^{(0)},f^{(0)},S^{(0)}\bigr)$ is also
a pseudo-polarized mixed twistor structure
of weight $n$ in one variable.

The following lemma is clear.
\begin{lem}\label{lem;9.13.5}
Let $(V,W,f,S)$ be a pseudo polarized mixed twistor structure
of weight $n$ in one variable.
Then it is a polarized mixed twistor structure of weight $n$
in one variable,
if and only if
the induced tuple $(V^{(0)},W^{(0)},f^{(0)},S^{(0)})$
is a polarized mixed twistor
structure of weight $n$ in one variable.
\hfill\qed
\end{lem}

\begin{df}
A tuple $(V,W,\vecf,S)$ is called polarized mixed twistor structure
of weight $n$ in $l$ variable,
if
\begin{itemize}
\item
 $(V,W)$ is a mixed twistor structure.
\item
 $\vecf$ is a tuple of nilpotent morphisms
 $f_i:V\lrarr V\otimes\Tate(-1)$ $(i=1,\ldots,l)$
 of mixed twistor structures.
\item
 Let $S:V\otimes\sigma(V)\lrarr \Tate(-n)$
 be a morphism of mixed twistor structures
 satisfying the following:
\[
 S(f_i\otimes\id)+S\bigl(\id\otimes\sigma(f_i)\bigr)=0.
\]
\item
 For any element $\veca\in\cnum^l$,
 we put $f(\veca):=\sum a_i \cdot f_i$.
 Then 
 $\bigl(V,W,f(\veca),S\bigr)$ is a polarized mixed twistor structure
 in one variable of weight $n$
 for any element $\veca\in\real^l_{>0}$.
\end{itemize}
For simplicity,
polarized mixed twistor structure of weight $n$ in $l$-variable
is abbreviated to Pol-MTS of $(n,l)$-type.

If a tuple $(V,W,\vecf,S)$ satisfies the first three conditions,
then it is called pseudo polarized mixed twistor structure
of weight $n$ in $l$-variable.
Similarly, pseudo polarized mixed twistor structure
of weight $n$ in $l$-variable
is abbreviated to
$\Psi$-Pol-MTS of $(n,l)$-type.
\hfill\qed
\end{df}

Let $(V,W,\vecf,S)$ is a $\Psi$-Pol-MTS of $(n,l)$-type.
Then we put as follows:
\[
 V^{(0)}:=\Gr^W(V),
\quad
 W^{(0)}_h:=\bigoplus_{i\leq h}\Gr^W_i(V).
\]
We have the induced morphism
$f^{(0)}_j:\Gr^W_i(V)\lrarr \Gr^{W}_{i-2}(V)\otimes\Tate(-1)$
and $f^{(0)}_j:V^{(0)}\lrarr V^{(0)}\otimes\Tate(-1)$.
We also have the induced morphism:
\[
\begin{array}{l}
 S^{(0)}:V^{(0)}\otimes \sigma(V^{(0)})\lrarr \Tate(-n),\\
 \mbox{{}}\\
 S^{(0)}:\Gr^W_{i+n}(V)\otimes \sigma\bigl(\Gr^{W}_{-i+n}(V)\bigr)
 \lrarr \Tate(-n).
\end{array}
\]
Then it is easy to check
that $\bigl(V^{(0)},W^{(0)},\vecf^{(0)},S^{(0)}\bigr)$
is also a $\Psi$-Pol-MTS of $(n,l)$-type.

\begin{lem}\label{lem;a12.1.35}
Let $(V,W,\vecf,S)$ is a $\Psi$-Pol-MTS of $(n,l)$-type.
It is a Pol-MTS of $(n,l)$-type,
if and only if
the induced tuple
$(V^{(0)},W^{(0)},\vecf^{(0)},S^{(0)})$ is a Pol-MTS of $(n,l)$.
\end{lem}
\pf
It follows from Lemma \ref{lem;9.13.5}.
\hfill\qed

\begin{lem}
The functor $\xi$ gives the following categories:
\begin{itemize}
\item
 The category of polarized equivariant mixed twistor structures.
\item
 The category of polarized mixed Hodge structures 
 defined over $\cnum$.
\end{itemize}
The functor $\xi$ gives the following categories:
\begin{itemize}
\item
 The category of equivariant mixed twistor structures
 defined over $\real$.
\item
 The category of polarized mixed Hodge structures
 defined over $\real$.
\end{itemize}
\end{lem}
\pf
It follows from Lemma \ref{lem;10.5.7},
Lemma \ref{lem;10.5.5}
and Corollary \ref{cor;10.5.6}.
\hfill\qed

%% file: b30.2.tex

\subsubsection{Definition and some functorial properties}

We put as follows:
\[
 \xi\Omega^1_X:=
 \Omega_X^{1,0}\otimes\nbigo_{\proj^1}(1,0)
 \oplus
 \Omega_X^{0,1}\otimes\nbigo_{\proj^1}(0,1),
\quad\quad
 \xi\Omega^h_X:=
 \bigwedge^h\Bigl(
 \xi\Omega^1_X
 \Bigr).
\]
We have the differential operator $\DD^{\sankaku}$
defined as follows:
\[
 \DD^{\sankaku}_X:
 C^{\infty}(X)
\lrarr
 C^{\infty}(X,\xi\Omega^1_X),
\quad
 g\longmapsto
 \del_X(g)\otimes f_{\infty}^{(1,0)}
+\delbar_X(g)\otimes \bigl(
 \sqrt{-1}\cdot f_0^{(0,1)}\bigr).
\]
When we forget the torus action,
we can use the notation $f^{(1)}_x$ $(x=0,1,\infty)$
instead of $f^{(1,0)}_x$ or $f^{(0,1)}_x$.
On the open subsets $X\times\cnum_{\lambda}$,
$X\times\cnum_{\mu}$ and $X\times\cnum_{\lambda}^{\ast}$,
it can be regarded as follows:
On $X\times\cnum_{\lambda}$,
we take the base $f^{(1)}_0$ of $\nbigo(1)_{|\cnum_{\lambda}}$.
Then $\DD_X^{\sankaku}$ induces the operator
$(\delbar_X+\lambda\cdot\del_X)\otimes (\sqrt{-1}\cdot f_0^{(1)})$.
On $X\times\cnum_{\mu}$,
we take the base $f^{(1)}_{\infty}$ of $\nbigo(1)_{|\cnum_{\mu}}$.
Then $\DD_X^{\sankaku}$ induces the operator
$(\del_X+\mu\delbar_X)\otimes f_{\infty}^{(1)}$.
On $X\times\cnum_{\lambda}^{\ast}$,
we take the base of $f^{(1,0)}_{\infty}$ of $\nbigo(1,0)$
and $\sqrt{-1}\cdot f^{(0,1)}_{0}$ of $\nbigo(0,1)$,
and then $\DD^{\sankaku}$ induces the operator
$d_X=\delbar_X+\del_X$.

The following lemma can be checked easily.
\begin{lem}\mbox{{}}
\begin{enumerate}
\item
The Leibniz rule holds in the following sense:
\[
 \DD_X^{\sankaku}(f\cdot g)
=f\cdot\DD_X^{\sankaku}(g)
+\DD_X^{\sankaku}(f)\cdot g.
\]
\item
We have the induced operator
$\DD^{\sankaku}$ on $C^{\infty}(X\times\proj^1,\xi\Omega^{\cdot}_X)$,
and we have the flatness $(\DD_X^{\sankaku}+\delbar_{\proj^1})^2=0$.
\hfill\qed
\end{enumerate}
\end{lem}

Let $V$ be a $\proj^1$-holomorphic vector bundle over $X\times\proj^1$.
\begin{df}
A variation of
$\proj^1$-holomorphic vector bundle over $X\times\proj^1$
is defined to be a differential operator
$\DD^{\sankaku}_V:
 C^{\infty}(X\times\proj^1,V)
 \lrarr
 C^{\infty}\bigl(X\times \proj^1,V\otimes \xi\Omega_X^1\bigr)$
satisfying the following conditions:
\[
 \bigl(
 \DD_V^{\sankaku}+d''_V
 \bigr)^2=0,
\quad\quad
 \DD_V^{\sankaku}(f\cdot v)
=f\cdot \DD_V^{\sankaku}(v)
+(\DD_X^{\sankaku}f)\cdot v.
\]
\hfill\qed
\end{df}

A tensor product, a direct sum and a dual
for variations of $\proj^1$-holomorphic vector bundles
are naturally defined.
Let $(V^{(i)},\DD^{\sankaku}_{V^{(i)}})$ $(i=1,2)$ be 
variation of $\proj^1$-holomorphic vector bundles.
A morphism of $(V^{(1)},\DD^{\sankaku}_{V^{(1)}})$
to $(V^{(2)},\DD^{\sankaku}_{V^{(2)}})$
is defined to be a $\proj^1$-holomorphic and $\DD^{\sankaku}_V$-flat
section of $Hom(V^{(1)},V^{(2)})$.

\subsubsection{Some description of a variation of 
 $\proj^1$-holomorphic bundles}

Let $V$ be a $\proj^1$-holomorphic
vector bundle with variation $\DD^{\sankaku}_V$.
On $X\times\cnum_{\lambda}$, we take the base $f_0^{(1)}$
of $\nbigo(1)$, and then $\DD^{\sankaku}_V$
induces the $\lambda$-connection
$\DD_V:
 C^{\infty}(X\times\cnum_{\lambda},V_0)\lrarr
 C^{\infty}(X\times\cnum_{\lambda},V_0\otimes\Omega_X^{1})$,
and we have the flatness $(\DD_{V_0}+d'')^2=0$.
On $X\times\cnum_{\mu}$, we take the base $f_{\infty}^{(1)}$
of $\nbigo(1)$, and then $\DD^{\sankaku}_{V_{\infty}}$ induces
the $\mu$-connection
$\DD_{V_{\infty}}^{\dagger}:
 C^{\infty}\bigl(X^{\dagger}\times\cnum_{\mu},V_{\infty}\bigr)
 \lrarr
 C^{\infty}\bigl(X^{\dagger}\times\cnum_{\mu},
 V_{\infty}\otimes\Omega^{1}_{X^{\dagger}}\bigr)$,
and we have the flatness $(\DD_{V_{\infty}}^{\dagger}+d'')^2=0$.
On $X\times\cnum_{\lambda}^{\ast}$,
we take the bases $f_{\infty}^{(1,0)}$ of $\nbigo(1,0)$
and $f^{(0,1)}_0$ of $\nbigo(0,1)$,
and then $\DD^{\sankaku}$ induces
the family of holomorphic connections
$\DD^{f}_{V_1}:
 C^{\infty}\bigl(X\times\cnum_{\lambda}^{\ast},
 V_{1}\bigr)
 \lrarr
 C^{\infty}\bigl(X\times\cnum_{\lambda}^{\ast},
 V_{1}\otimes\Omega^{1}_{X}\bigr) $,
and we have the flatness
$(\DD^f_{V_1}+d'')^2=0$.

\label{subsubsection;10.18.2}

On the other hand, we can consider a patched object
$(V_0,V_{\infty},V_1;\alpha_0,\alpha_{\infty})$:
\begin{itemize}
\item
 $V_0$ is a $\cnum_{\lambda}$-holomorphic vector bundle
 over $X\times\cnum_{\lambda}$,
which is equipped with the $\lambda$-connection $\DD_{V_0}$
such that $(\DD_{V_0}+d'')^2=0$.
\item
 $V_{\infty}$ is a $\cnum_{\mu}$-holomorphic vector bundle
 over $X\times\cnum_{\mu}$,
 which is equipped with the $\mu$-connection $\DD_{V_{\infty}}$
 such that $(\DD_{V_{\infty}}+d'')^2=0$.
\item
$V_1$ is a $\cnum_{\lambda}^{\ast}$-holomorphic vector bundle
 over $X\times \cnum_{\lambda}^{\ast}$,
which is equipped with the holomorphic family of the flat connections
$\DD^f_{V_1}$.
\item
 We have the induced families of flat connections
 $(V_{0\,|\,X\times\cnum_{\lambda}^{\ast}},\DD^f_{V_0})$,
 $(V_{\infty\,|\,X\times\cnum_{\lambda}^{\ast}},
 \DD^{\dagger\,f}_{V_{\infty}})$.
 Then $\alpha_a$ $(a=0,\infty)$ are isomorphisms
 $V_{a\,|\, X\times\cnum_{\lambda}^{\ast}}
 \lrarr
 V_{1\,|\,X\times \cnum_{\lambda}^{\ast}}$,
 which are compatible with the $\cnum_{\lambda}^{\ast}$-holomorphic
 structure and the family of flat connections.
\end{itemize}

Once we are given such a patched object
$(V_0,V_{\infty},V_1;\alpha_0,\alpha_{\infty})$,
then we have the $\proj^1$-holomorphic vector bundle
$V$ by gluing.
\begin{lem}\label{lem;a12.1.1}
We have the well defined differential operator
$\DD^{\sankaku}_V$ given as follows:
\[
 \DD^{\sankaku}_V(f)=
\left\{
\begin{array}{ll}
 \DD_{V_0}(f)\otimes \bigl(\sqrt{-1}\cdot f_{0}^{(1)}\bigr),
 & (\mbox{on }X\times\cnum_{\lambda})\\
 \mbox{{}}\\
 \DD_{V_{\infty}}^{\dagger}(f)\otimes f_{\infty}^{(1)},
 & (\mbox{on }X\times\cnum_{\mu}).
\end{array}
\right.
\]
The operator $\DD^{\sankaku}_V$ gives a variation of
$\proj^1$-holomorphic vector bundles.
\end{lem}
\pf
The well definedness follows from the compatibility
of $\alpha_a$ $(a=0,\infty)$ with the flat connections.
It is easy to see that $\DD^{\sankaku}_{V}$ gives a variation.
\hfill\qed

\begin{cor}
The category of the patched objects above
and the category of the variations of $\proj^1$-holomorphic vector
 bundles are equivalent,
by the correspondence given above.
\hfill\qed
\end{cor}

We can also consider patched objects
$(V_0,V_{\infty};\psi)$:
\begin{itemize}
\item
$V_0$ is a $\cnum_{\lambda}$-holomorphic vector bundle
equipped with $\lambda$-connection $D_{V_0}$.
\item
$V_{\infty}$ is a $\cnum_{\mu}$-holomorphic vector bundle
equipped with $\mu$-connection $D^{\dagger}_{V_{\infty}}$.
\item
$\psi$ is an isomorphism
of the induced family of the flat bundles
$V_{0\,|\,\cnum_{\lambda}^{\ast}\times X}$
and $V_{1\,|\,\cnum_{\lambda}^{\ast}\times X}$.
\end{itemize}
As before, the category of such patched objects
and the category of variations of $\proj^1$-holomorphic vector bundles
are naturally equivalent.

\vspace{.1in}

In the following,
we often use the descriptions
$\bigl(V_0,V_{\infty},V_1;\alpha_0,\alpha_{\infty}\bigr)$
or $(V_0,V_{\infty},\psi)$
to denote the variation of $\proj^1$-holomorphic vector bundles.

%% file: b30.3.tex

\subsubsection{An example of variation of
the $\proj^1$-holomorphic vector bundle}

\label{subsubsection;10.5.15}

Let $V$ be a vector bundle over $\proj^1$,
and $f_i:V\lrarr V\otimes \Tate(-1)$ $(i=1,\ldots,n)$
be nilpotent morphisms
such that $f_i$ and $f_j$ are commutative.

We put $X:=\cnum^n$, $D_i:=\{z_i=0\}$,
and $D=\bigcup_{i=1}^lD_i$.
We will construct the 
$\proj^1$-holomorphic vector bundle $\nbigv$
and the variation $\DD^{\sankaku}_{\nbigv}$
over $X-D$.

We put $V_0:=V_{|\cnum_{\lambda}}$
and $V_{\infty}:=V_{|\cnum_{\mu}}$.
We put as follows:
\[
 \nbigv_0:=
 V_{0}\otimes \nbigo_{\nbigx-\nbigd},
\quad\quad
 \nbigv_{\infty}:=
 V_{\infty}\otimes\nbigo_{\nbigx^{\dagger}-\nbigd^{\dagger}}.
\]

From the morphism $f_i\in Hom(V,V\otimes\Tate(-1))$,
we obtain
the morphism $f_i\in Hom(V_0,V_0\otimes \Tate(-1))$
on $\cnum_{\lambda}$,
and then we have the endomorphism
$f_i\otimes t_0^{(1)}\in \End(V_{0})$.
Then the $\lambda$-connection $\DD_{\nbigv_0}$
is given.
Namely,
for any $v\in\Gamma(\cnum_{\lambda},V_0)$
and $g\in C^{\infty}(\nbigx-\nbigd)$,
we put as follows:
\begin{equation} \label{eq;10.5.10}
 \DD_{\nbigv_0}(g\cdot v):=
\sum_{i=1}^lg\cdot f_i(v)\cdot t_0^{(1)}\cdot\frac{dz_i}{z_i}
+\bigl(\lambda\cdot \del_X(g)+\delbar_X(g)\bigr)\cdot v.
\end{equation}

We also have the endomorphism
$f_i\otimes t_{\infty}^{(1)}
\in \End(V_{\infty})$
on $\cnum_{\mu}$.
The $\mu$-connection $\DD^{\dagger}_{\nbigv_{\infty}}$
is given.
Namely,
for any $v^{\dagger}\in\Gamma(\cnum_{\mu},V_{\infty})$
and $g\in C^{\infty}(\nbigx^{\dagger}-\nbigd^{\dagger})$,
we put as follows:
\begin{equation}\label{eq;10.5.11}
 \DD^{\dagger}_{\nbigv_{\infty}} (g\cdot v^{\dagger})
:=\sum_{i=1}^l
 g\cdot f_i(v^{\dagger})\cdot t_{\infty}^{(1)}\cdot
 \frac{d\bar{z}_i}{\bar{z}_i}
+\bigl(\mu\cdot\delbar_X(g)+\del_X(g)\bigr)\cdot v^{\dagger}.
\end{equation}

We will give the isomorphism
$\Psi:
 \nbigv_{0|\nbigx^{\shikaku}-\nbigd^{\shikaku}}
\lrarr
 \nbigv_{\infty|\nbigx^{\dagger\,\shikaku}-\nbigd^{\dagger\,\shikaku}}$.
Let $\lambda=\mu^{-1}$ be an element of $\cnum_{\lambda}^{\ast}$
and $P$ be a point of $X-D$.
Let $v$ be an element of $V_{|\lambda}$.
It gives the elements
$v_0\in V_{0\,|\,\lambda}$ and
$v_{\infty}\in V_{\infty\,|\,\mu}$.
Then we naturally obtain the following elements:
\[
 v_{0\,|\,P}\in \nbigv_{0|(\lambda,P)},
\quad\quad
 v^{\dagger}_{\infty\,|\,P}\in \nbigv_{\infty|(\mu,P)}.
\]
Then $\Psi$ is defined as follows:
\begin{equation}\label{eq;10.5.16}
 \Psi\Bigl(
 \exp\Bigl(
 -\sum_{i=1}^n \lambda^{-1}\cdot\log|z_i(P)|^2
\cdot f_i\otimes t_0^{(1)}
 \Bigr)
\cdot v_{0|(\lambda,P)}
 \Bigr)
=v_{\infty|(\mu,P)},
\end{equation}
or equivalently,
\[
 \Psi(v_{0\,|\,P})=
 \exp\Bigl(
 -\sum_{i=1}^n\mu^{-1}\cdot
 \log|z_i(P)|^2
\cdot f_i\otimes t_{\infty}^{(1)}
 \Bigr)\cdot v_{\infty\,|\,P}.
\]
Here we have used the relation
$\lambda^{-1}\cdot t_0^{(1)}=-\mu^{-1}\cdot t_{\infty}^{(1)}$.

Corresponding to $v\in V_{|\lambda}$,
we have the flat sections
of $\nbigv_{0\,|\,\nbigx^{\lambda}-\nbigd^{\lambda}}$
and $\nbigv_{\infty\,|\,\nbigx^{\dagger\,\mu}-\nbigd^{\dagger\,\mu}}$:
\[
 \begin{array}{l}
{\displaystyle
 \tilde{s}_0:=
 \exp\Bigl(
 -\sqrt{-1}\sum_{i=1}^n \log z_i\cdot f_i\otimes t^{(1)}_1
 \Bigr)\cdot v_0,}
 \\
 \mbox{{}}\\
 {\displaystyle
 \tilde{s}_{\infty}:=
 \exp\Bigl(
 \sqrt{-1}\sum_{i=1}^n\log \overline{z}_i\cdot f_i\otimes t^{(1)}_1
 \Bigr)\cdot v_{\infty}.}
\end{array}
\]
\begin{lem}\label{lem;10.17.6}
We have the relation $\tilde{s}_0=\tilde{s}_{\infty}$.
The gluing $\Psi$ is characterized by this property.
\end{lem}
\pf
It can be easily checked by a direct calculation.
\hfill\qed

\begin{lem}
$\Psi$ is compatible with $\DD_{\nbigv_0}$
and $\DD^{\dagger}_{\nbigv_{\infty}}$.
\end{lem}
\pf
We put as follows:
\[
 \tilde{v}:=
 \exp\Bigl(
 -\sum_{i=1}^n\lambda^{-1}\cdot\log|z_i|^2\cdot
 f_i\otimes t_0^{(1)}
 \Bigr)v_0.
\]
Then we have the following:
\begin{multline}
 \DD^f_{\nbigv_0}(\tilde{v})
=\sum \lambda^{-1}\cdot f_i(\tilde{v})
\otimes t_0^{(1)}
\frac{dz_i}{z_i}
-\sum \lambda^{-1}
 \left(
 \frac{dz_i}{z_i}+\frac{d\bar{z}_i}{\bar{z}_i}
 \right)\cdot f_i(\tilde{v})\otimes t_0^{(1)}\\
=-\sum \lambda^{-1}f_i(\tilde{v})
  \otimes t_0^{(1)}\frac{d\bar{z}_i}{\bar{z}_i}
=\mu^{-1}\sum f_i(\tilde{v})\otimes t_{\infty}^{(1)}
 \frac{d\bar{z}_i}{\bar{z}_i}.
\end{multline}
Hence it is compatible with the definition
of
$\DD^{\dagger\,f}_{\nbigv_{\infty}}$
given in (\ref{eq;10.5.11}).
\hfill\qed

\vspace{.1in}
We put $N_i:=-f_i\otimes t_1^{(1)}$.
\begin{lem} \label{lem;10.17.1}
The endomorphism $\exp\bigl(2\pi\cdot N_i\bigr)$
is the monodromy of the loop $\{0\leq t\leq 1\}\lrarr \cnum^{\ast\,n}$
given by the following:
\[
 t
\longmapsto
 \bigl(z_1,\ldots, z_{i-1},\exp(-2\pi\sqrt{-1}t)\cdot z_i,
 z_{i+1},\ldots,z_{n}\bigr).
\]
\end{lem}
\pf
It can be checked by a direct calculation.
\hfill\qed

\begin{rem} \label{rem;10.17.5}
When we consider the monodromy in this paper,
we usually use the loop with the inverse direction:
\[
 t
\longmapsto
 \bigl(z_1,\ldots,\cdot z_{i-1},\exp(2\pi\sqrt{-1}t)\cdot z_i,
 z_{i+1},\ldots,z_{n}\bigr).
\]
However, Cattani-Kaplan-Schmid use the loop
given in Lemma {\rm\ref{lem;10.17.1}}.
The author apologize the inconvenience,
and he hopes that there are no confusion.
\hfill\qed
\end{rem}

%% file: b30.4.tex

\subsubsection{The involution and the induced variation}

\label{subsubsection;10.18.3}

We have the isomorphisms
$\iota_{1,0}:\sigma\nbigo(1,0)\lrarr \nbigo(0,1)$
and $\iota_{0,1}:\sigma\nbigo(0,1)\lrarr\nbigo(1,0)$.
We take the isomorphisms
$\sigma\Omega^{1,0}_X\lrarr \Omega^{0,1}_X$
and
$\sigma\Omega^{0,1}_X\lrarr\Omega^{1,0}_X$
given by the following:
\[
 dz_i\longmapsto -d\bar{z}_i,
\quad\quad
 d\bar{z}_i\longmapsto dz_i.
\]
Then we obtain the morphisms
$ \sigma^{\ast}\bigl(
 \Omega_X^{1,0}\otimes\nbigo(1,0)
 \bigr)\simeq 
 \Omega_X^{0,1}\otimes\nbigo(0,1)$
and
$\sigma^{\ast}\bigl(
 \Omega_X^{0,1}\otimes\nbigo(0,1)
 \bigr)
\simeq
 \Omega_X^{1,0}\otimes\nbigo(1,0)$.
We denote them by $\varphi_0$.

\begin{lem} \label{lem;10.18.6}
We have the following equalities:
\[
\begin{array}{ll}
\varphi_0\bigl(
 \sigma^{\ast}\bigl(dz_i\otimes f_{\infty}^{(1,0)}\bigr)
\bigr)
=d\zbar_i\otimes\sqrt{-1}f_0^{(0,1)},
 &
 \varphi_0\bigl(
 \sigma^{\ast}\bigl(d\zbar_i\otimes f_{\infty}^{(1,0)}\bigr)
\bigr)
=-dz_i\otimes\sqrt{-1}f_0^{(0,1)},\\
 \mbox{{}}\\
\varphi_0\bigl(
 \sigma^{\ast}\bigl(d\zbar_i\otimes \sqrt{-1}\cdot f_0^{(0,1)}\bigr)
 \bigr)
=dz_i\otimes f^{(1,0)}_{\infty},
 &
 \varphi_0\bigl(
 \sigma^{\ast}\bigl(dz_i\otimes \sqrt {-1}\cdot f_0^{(0,1)}\bigr)
 \bigr)
=-d\zbar_i\otimes f^{(1,0)}_{\infty}.
\end{array}
\]
\end{lem}
\pf
It can be checked by a direct calculation.
\hfill\qed

\vspace{.1in}

The morphism $\varphi_0$ induces the following morphisms:
\[
 \begin{array}{ll}
C^{\infty}\bigl(X\times\proj^1,
 \Omega_X^{0,1}\otimes\nbigo(1,0)
 \bigr)\simeq 
C^{\infty}\bigl(X\times\proj^1,
 \Omega_X^{0,1}\otimes\nbigo(0,1)
 \bigr), \\
 \mbox{{}}\\
C^{\infty}\bigl(X\times\proj^1,
 \Omega_X^{0,1}\otimes\nbigo(0,1)
 \bigr)
\simeq
C^{\infty}\bigl(X\times\proj^1,
 \Omega_X^{1,0}\otimes\nbigo(1,0)
\bigl).
\end{array}
\]
We denote them by $\varphi$.

\begin{lem} \label{lem;10.17.2}
We have the following:
\[
\begin{array}{l}
 \varphi\bigl(
 g\cdot dz\otimes f^{(1,0)}_{\infty}
 \bigr)
 =\sigma^{\ast}(\bar{g})\cdot d\bar{z}\otimes \sqrt{-1}\cdot f_0^{(0,1)},
\quad\quad
 \varphi\bigl(
 g\cdot d\bar{z}\otimes \sqrt{-1}\cdot f_0^{(0,1)}\bigr)
 =\sigma^{\ast}(\bar{g})\cdot dz\otimes f_{\infty}^{(1,0)}.
\end{array}
\]
\end{lem}
\pf
It can be checked by a direct calculation.
\hfill\qed

\vspace{.1in}
Recall that we put $\varphi(f):=\sigma^{\ast}(\overline{f})$
for a function $f$.

\begin{lem}
We have
$\varphi\circ\DD^{\sankaku}_X=
 \DD^{\sankaku}\circ\varphi$.
\end{lem}
\pf
We have the following:
\[
 \varphi\bigl(
 \DD^{\sankaku}_Xf
 \bigr)
=\varphi(\del_Xf\otimes f_{\infty}^{(1)}
 +\delbar_X f\otimes \sqrt{-1}f_0^{(1)}).
\]
We have the following, by using Lemma \ref{lem;10.17.2}:
\[
 \varphi\Bigl(
 \frac{\del f}{\del z_i}dz_i\otimes f^{(1)}_{\infty}
 \Bigr)
=\varphi\Bigl(
 \frac{\del f}{\del z_i}\Bigr)
 \cdot d\bar{z}_i\otimes \sqrt{-1}f^{(1)}_0
=\frac{\del \varphi(f)}{\del \bar{z}_i}\cdot
 d\bar{z}_i\otimes \sqrt{-1}f_0^{(1)}.
\]
Similarly we have the following:
\[
 \varphi\Bigl(
 \frac{\del f}{\del \bar{z}_i}\cdot d\bar{z}_i\otimes \sqrt{-1}\cdot f_0^{(1)}
 \Bigr)
=\frac{\del\varphi(f)}{\del z_i}\cdot dz_i\otimes f^{(1)}_{\infty}.
\]
It implies the commutativity
$\varphi\circ\DD_X^{\sankaku}=\DD_X^{\sankaku}\circ\varphi$.
\hfill\qed

\vspace{.1in}
Let $V$ be a $\proj^1$-holomorphic vector bundle
over $X\times \proj^1$.
Let $\DD^{\sankaku}_V$ be a variation of $\proj^1$-holomorphic
vector bundles on $V$.
We have the $\proj^1$-holomorphic bundle $\sigma(V)$.
Then we have the operator $\DD^{\sankaku}_{\sigma(V)}$
on $\sigma(V)$ defined as follows:
\[
 \DD^{\sankaku}_{\sigma(V)}\bigl(\sigma(v)\bigr)
=\varphi_0\circ\sigma\bigl(
 \DD^{\sankaku}_Vv
 \bigr).
\]
\begin{lem}
$\DD^{\sankaku}_{\sigma(V)}$ is a variation of $\proj^1$-holomorphic
vector bundles.
\end{lem}
\pf
We have the following:
\begin{multline}
 \DD^{\sankaku}_{\sigma(V)}\bigl(
 f\cdot\sigma(v)
 \bigr)
=\DD^{\sankaku}_{\sigma(V)}
 \bigl(\sigma\bigl(\varphi(f)\cdot v\bigr)\bigr)
=\varphi_0\circ\sigma\bigl(
 \DD^{\sankaku}_V(\varphi(f)\cdot v)
 \bigr) \\
=\varphi_0\circ\sigma\Bigl(
 \DD^{\sankaku}_X(\varphi(f))\cdot v
+\varphi(f)\cdot\DD^{\sankaku}_V(v)
 \Bigr)
=\DD_X^{\sankaku}(f)\cdot\sigma(v)
+f\cdot \DD^{\sankaku}_{\sigma(V)}(\sigma(v)).
\end{multline}
Thus we are done.
\hfill\qed

%% file: b6.4.tex

\subsubsection{Equivariant variation of $\proj^1$-holomorphic vector bundles}

The following lemma is easy to see.
\begin{lem}
$\DD_X^{\sankaku}$ is equivariant with respect to
the torus action.
We have
$\DD_X^{\sankaku}(\sigma^{\ast}g)=\sigma^{\ast}\bigl(
 \DD_X^{\sankaku}g\bigr)$.
\hfill\qed
\end{lem}

Let $(V,\rho)$ be an equivariant $\proj^1$-holomorphic vector bundle
over $X\times\proj^1$.
Let $\DD^{\sankaku}_V$ be the variation of
$\proj^1$-holomorphic vector bundles.
We put $H:=V_{|\{1\}\times X}$.
Then we have the flat connection
$\DD_0:=\DD^{\sankaku}_{V\,|\,\{\lambda\}\times X}$.
Since $V$ is equivariant, we have the two filtrations
$F$ and $G$ on $H$,
such that $\xi(H,F,G)\simeq (V,\rho)$.

\begin{lem}
 $\DD^{\sankaku}_V$ is equivariant,
if and only if $\DD_0$ satisfies the Griffiths transversality
in the following sense:
\[
 \begin{array}{ll}
 \DD_0^{(0,1)}\Bigl(
 C^{\infty}(F^p)\Bigr)
 \subset
 C^{\infty}(F^p\otimes\Omega_X^{0,1}), &
 \DD_0^{(1,0)}\Bigl(
 C^{\infty}(F^p)
 \Bigr)\subset
 C^{\infty}(F^{p-1}\otimes\Omega_X^{1,0}),\\
\mbox{{}}\\
 \DD_0^{(0,1)}\Bigl(
 C^{\infty}(G^p)\Bigr)\subset
 C^{\infty}(G^{p-1}\otimes\Omega_X^{0,1}), &
 \DD_0^{(1,0)}\Bigl(
 C^{\infty}(G^p)\Bigr)\subset
 C^{\infty}(G^p\otimes\Omega_X^{1,0}).
 \end{array}
\]
\end{lem}
\pf
Let $s$ be a section of $F^p$.
Then $\lambda^{-p}\cdot s$ is a section of $V$
on $X\times\cnum_{\lambda}$.
We have the following:
\[
 \DD_V^{\sankaku}(\lambda^{-p}\cdot s)
=\lambda^{-p}\cdot
 \bigl(\DD^{(1,0)}s\otimes f_{\infty}^{(1,0)}+
 \DD^{(0,1)}s\otimes f_0^{(0,1)} \bigr)
=\lambda^{-p+1}\cdot\DD_0^{(1,0)}(s)\otimes f_0^{(1,0)}
+\lambda^{-p}\cdot\DD_0^{(0,1)}(s)\otimes f_0^{(0,1)}.
\]
Thus $\DD_V^{\sankaku}(\lambda^{-p}\cdot s)$ is a section
of $V\otimes\xi\Omega^1_X$ on $X\times\cnum_{\lambda}$
if and only if
the following is satisfied:
\[
 \DD_0^{(1,0)}(s)\in F^{p-1}\otimes\Omega_X^{1,0},
\quad\quad
 \DD_0^{(0,1)}(s)\in F^{p}\otimes\Omega_X^{0,1}.
\]
Similar things for $G$ hold.
Thus we are done.
\hfill\qed

\vspace{.1in}

Let $(W,\rho)$ be an equivariant subbundle of $(V,\rho)$.
We have the corresponding vector bundle
$H_W$ and $H_V$.
\begin{lem}
We have
$\DD_V^{\sankaku}\bigl(
 C^{\infty}(X,W)\bigr)\subset
 C^{\infty}(X,W\otimes\xi\Omega_X^1)$
if and only if
$H_W$ is a flat vector subbundle of $H$
with respect to the flat connection $\DD_0$.
\end{lem}
\pf
$\DD_V$ induces the morphism
$C^{\infty}(X,H)[\lambda,\lambda^{-1}]
\lrarr C^{\infty}(X,H\otimes\Omega_X^1)[\lambda,\lambda^{-1}]$,
which is same as $\DD_0\otimes \id$.
The both claims are equivalent to the following:
\[
 \DD_V\Bigl(
 C^{\infty}(X,H_W)[\lambda,\lambda^{-1}]
 \Bigr)
\subset
 C^{\infty}(X,H_W\otimes\Omega^1_X)[\lambda,\lambda^{-1}].
\]
Then the lemma follows.
\hfill\qed

\vspace{.1in}
Let $f$ be an equivariant morphism
$(V_1,\rho_1)\lrarr (V_2,\rho_2)$
corresponding to
$f_{|1}:(H_1,F_1,G_1)\lrarr (H_2,F_2,G_2)$.
\begin{lem}
$\DD_{V_2}^{\sankaku}\circ f=f\circ\DD_{V_1}^{\sankaku}$
if and only if
$\DD_{0\,H_2}\circ f_{|1}=f_{|1}\circ \DD_{0\,H_1}$.
\end{lem}
\pf
$f$ induces the morphism
$\tilde{f}:
 C^{\infty}(X,H_1)[\lambda,\lambda^{-1}]
\lrarr
 C^{\infty}(X,H_2)[\lambda,\lambda^{-1}]$.
Then the both claims are equivalent to
the compatibility of $\tilde{f}$,
$\DD_{V_2}^{\sankaku}$ and $\DD_{V_1}^{\sankaku}$.
\hfill\qed

\begin{cor}\mbox{{}}
\begin{itemize}
\item
$\iota_V$ is flat with respect to $\DD_{V}^{\sankaku}$
if and only if
$\iota_H$ is flat with respect to $\DD_{0}$.
\item
 $S$ is flat if and only if $\langle\cdot,\cdot\rangle$ is flat.
\end{itemize}
\end{cor}
\pf
It follows from the previous lemma.
\hfill\qed

\begin{cor} \label{cor;10.5.20}
$\xi$ gives the equivalence of the following categories:
\begin{itemize}
\item
 The category of variation of pure twistor structures.
\item
 The category of variation of complex pure Hodge structure.
\end{itemize}
It is compatible with real structures and polarizations.
\hfill\qed
\end{cor}

\begin{cor}
$\xi$ gives the equivalence of the following categories:
\begin{itemize}
\item
 The category of variation of mixed twistor structures.
\item
 The category of variation of complex mixed Hodge structure.
\end{itemize}
It is compatible with the real structures and the polarizations.
\hfill\qed
\end{cor}

%% file: a71.tex

\subsubsection{Pairing}

We put $X=\cnum^n$,
$D_i:=\bigl\{z_i=0\bigr\}$
and $D=\bigcup_{i=1}^n D_i$.
Let $V$ be a holomorphic vector bundle over $\proj^1$
and $\vecf$ be a tuple of nilpotent maps
$f_i:V\lrarr V\otimes\Tate(-1)$ $(i=1,\ldots,n)$.
From $(V,\vecf)$, we obtain the $\proj^1$-holomorphic
vector bundle $\nbigv$ 
and the variation $\DD^{\sankaku}_{\nbigv}$
over $\proj^1\times (X-D)$
(the subsubsection \ref{subsubsection;10.5.15}).

Let $S:V\otimes\sigma(V)\lrarr \Tate(0)$ be a pairing
such that
$S(f_i\otimes\id)+S\bigl(\id\otimes\sigma(f_i)\bigr)=0$.
Then the pairing
$\tilde{S}:\nbigv_0\otimes\sigma\nbigv_{\infty}
\lrarr \nbigo_{\nbigx-\nbigd}$
is given.
Namely,
for any $u\in V_0$, $v\in V_{\infty}$,
$a\in C^{\infty}(\nbigx-\nbigd)$
and $b\in C^{\infty}(\nbigx^{\dagger}-\nbigd^{\dagger})$,
we put as follows:
\[
 \tilde{S}\bigl(a\cdot u\otimes \sigma(b\cdot v)\bigr)
=a\cdot\overline{\sigma^{\ast}(b)}\cdot
 S\bigl(u\otimes\sigma(v)\bigr).
\]
Then we obtain the morphism
$\tilde{S}:\nbigv(V,\vecf)_0\otimes\sigma\nbigv(V,\vecf)_{\infty}
 \lrarr \nbigo_{\nbigx}$.
\begin{lem} \label{lem;10.17.9}
$\tilde{S}$ is a morphism of $\lambda$-connections.
\end{lem}
\pf
Let $u$ and $v$ be sections of
$V_0$ and $\sigma(V_{\infty})$ respectively.
We have $\DD_X\tilde{S}\bigl(u,\sigma^{\ast}(v)\bigr)=0$.
On the other hand,
we have the following equality on $\nbigx-\nbigd$:
\[
{\displaystyle
 \tilde{S}\bigl(\DD^{\sankaku} u\otimes\sigma(v)\bigr)
=\sum_i S\bigl(f_i(u)\otimes t_0^{(1)}\otimes v\bigr)
 \cdot\frac{dz_i}{z_i}\otimes \sqrt{-1}\cdot f_0^{(1,0)}}\\
\]
We also have the following:
\begin{multline}
 \tilde{S}\bigl(u\otimes \DD^{\sankaku}\sigma (v)\bigr)
=\sum_i\tilde{S}\Bigl(u\otimes
  \varphi_0\sigma\Bigl(f_i(v)\otimes
  t_{\infty}^{(1)}\otimes f_{\infty}^{(0,1)}
 \cdot\frac{d\overline{z}_i}{\overline{z}_i}
 \Bigr)\Bigr) \\
=\sum_i\tilde{S}\Bigl(
 u\otimes\varphi_0\sigma\bigl(f_i(v)\bigr)
 \Bigr)\otimes (-t_{0}^{(1)})
 \otimes(-\sqrt{-1})
\cdot f_0^{(1,0)}\cdot\frac{dz_i}{z_i}
=\sum_i
 S\Bigl(u\otimes \sigma\bigl(f_i(v)
 \bigr)\Bigr)
\cdot t_0^{(1)}
 \cdot\frac{dz_i}{z_i} \otimes \sqrt{-1}\cdot f_0^{(1,0)}.
 \end{multline}
Thus we obtain the equality
$\tilde{S}(\DD_{\nbigv_0}\otimes\id)
+\tilde{S}(\id\otimes\DD_{\sigma(\nbigv_{\infty})})=
 \DD_{X}\circ\tilde{S}$.
\hfill\qed

\vspace{.1in}

On the plane $\nbigx^{\dagger}-\nbigd^{\dagger}$,
we have the pairing
$\tilde{S}:\nbigv_{\infty}\otimes \sigma \nbigv_0
 \lrarr\nbigo_{\nbigx^{\dagger}-\nbigd^{\dagger}}$.
\begin{lem}
$\tilde{S}$ is a morphism of $\mu$-connections.
\end{lem}
\pf
Note the following equality on $\nbigx^{\dagger}-\nbigd^{\dagger}$:
\[
\tilde{S}\bigl(\DD^{\sankaku}_{\nbigv_{\infty}}
 u\otimes \sigma (v)\bigr)=
\sum_i
S\bigl(f_i(u)\otimes t_{\infty}^{(1)}
 \otimes \sigma (v)
 \bigr)\otimes
 f_{\infty}^{(0,1)}\cdot
 \frac{d\overline{z}_i}{\overline{z}_i}
\]
We also have the following:
\begin{multline}
 \tilde{S}\bigl(u\otimes \DD^{\dagger}_{\sigma(\nbigv_0)} \sigma(v)\bigr)
=
\sum_i
S\Bigl(u\otimes \varphi_0\sigma \Bigl(
 f_i(v)\otimes t_0^{(1)}\cdot
 \sqrt{-1}f_0^{(1,0)}\cdot\frac{dz_i}{z_i}
 \Bigr)
 \Bigr) \\
=
\sum_i
 S\bigl(u\otimes \sigma\bigl(f_i(v)\bigr)\bigr)
\otimes \bigl(-t_{\infty}^{(1)}\bigr)
 \cdot f_{\infty}^{(0,1)}
 \cdot(-\sqrt{-1})\cdot (\sqrt{-1}f_{\infty}^{(0,1)})
\cdot
\Bigl(-\frac{d\overline{z}_i}{\overline{z}_i}\Bigr)
=\sum_i
  S\bigl(u\otimes \sigma\bigl(f_i(v)\bigr)\bigr)
\otimes t_{\infty}^{(1)}\cdot f_{\infty}^{(0,1)}
 \cdot \frac{d\zbar_i}{\zbar_i}.
\end{multline}
Thus we obtain the equality
$\tilde{S}\bigl(\DD^{\dagger}_{\nbigv_{\infty}}(u)\otimes v\bigr)
+\tilde{S}\bigl(u\otimes\DD_{\sigma(\nbigv_0)}\sigma(v)\bigr)=0$
for any $u,v\in \Gamma(\cnum_{\mu},V_{\infty})$.
\hfill\qed

\begin{lem}
We obtain the pairing
$S_{\nbigv}:\nbigv\otimes\sigma\nbigv\lrarr  \Tate(0)$.
\end{lem}
\pf
We have only to check the pairings on the planes
$\nbigx-\nbigd$ and $\nbigx^{\dagger}-\nbigd^{\dagger}$
are preserved by the gluing morphism.
Note we have the following:
\[
 S\Bigl(
 \bigl(-\lambda^{-1}\cdot f_i\otimes t_0^{(1)}\bigr)
\otimes
 \id
 \Bigr)
+S\Bigl(
 \id\otimes\sigma\bigl(-\mu^{-1}\cdot f_i\otimes t_{\infty}^{(1)}
 \bigr)
 \Bigr)=0.
\]
Then we obtain the following compatibility.
\[
 S(u\otimes\sigma^{\ast}v)
=S\Bigl(
\exp \Bigl(
 -\sum\lambda^{-1}\cdot\log|z_i|^2 f_i\otimes t_0^{(1)}
 \Bigr)u
\otimes
\sigma^{\ast}
\exp
 \Bigl(
 -\sum\mu^{-1}\cdot\log|z_i|^2 f_i\otimes t_{\infty}^{(1)}
 \Bigr)v
 \Bigr).
\]
Thus we are done.
\hfill\qed

\subsubsection{Definition of twistor nilpotent orbit}

\begin{df}
 $(V,\vecf,S)$ is called a twistor nilpotent orbit,
if there exists a positive constant $C>0$ such that
$\bigl(\nbigv,\DD^{\sankaku}_{\nbigv},S_{\nbigv}\bigr)$
is a variation of polarized pure twistor structure
over $\Delta(C)^{\ast\,n}$.
\hfill\qed
\end{df}

\begin{lem}
The resulted harmonic bundle over $\Delta(C)^{\ast\,n}$
is tame. 
The eigenvalues of the residues of Higgs field are trivial.
The parabolic structure is trivial.
\end{lem}
\pf
The first two claims are clear from our construction of the variation.
By our construction, it is clear
that the eigenvalue is trivial.
By seeing the eigenvalues of $\lambda$-connections
for any $\lambda$,
we obtain the triviality of the parabolic structures.
\hfill\qed


\begin{lem}
The tuple
$\bigl(S^{\can}(Pat(V,\vecf,S)),\Res_i,S\bigr)$
is isomorphic to the original $(V,\vecf,S)$.
(See the subsubsection {\rm\ref{subsubsection;a12.5.5}}.)
\end{lem}
\pf
We have only to note that
the prolongment of $M(V,\vecf)$ is 
$V_0\otimes\nbigo_{\nbigx}$.
On $X\times\cnum^{\ast}_{\lambda}$,
it is clear.
Then by using the Hartogs Theorem,
we obtain the coincidence
$\prolong{M(V,\nbigf)}=V_0\otimes\nbigo_{\nbigx}$.
Similarly,
we obtain $\prolong{A(V,\nbigf)}=
 V_{\infty}\otimes\nbigo_{\nbigx^{\dagger}}$.
Then we obtain the isomorphisms:
\[
 S^{\can}_{|\cnum_{\lambda}}\simeq V_0,
\quad\quad
 S^{\can}_{|\cnum_{\mu}}\simeq V_{\infty}.
\]

Let us compare the gluing.
Let $\lambda$ be a point of $\cnum_{\lambda}^{\ast}$.
Let $v$ be an element of $V_{|\lambda}$.
Let us consider the multi-valued flat section $\tilde{v}_1$
of $M(V,\vecf)_{|\lambda}$
and the multi-valued flat section $\tilde{v}_2$
of $A(V,\vecf)_{|\lambda^{-1}}$ are given as follows:
\[
 \tilde{v}_1:=
\exp\Bigl(-\sqrt{-1}\sum_{i=1}^n\log z_i\cdot f_i\otimes t_1^{(1)}\bigr)v,
\quad
 \tilde{v}_2:=
\exp\Bigl(\sqrt{-1}\sum_{i=1}^n
 \log \overline{z}_i\cdot f_i\otimes t_1^{(1)}\Bigr)v.
\]
Then $\tilde{v}_1$ gives a flat section of
$V_{0\,|\,\lambda}$,
and $\tilde{v}_2$ gives a flat section of
$V_{\infty\,|\,\lambda^{-1}}$.
Then the gluing of $S^{\can}$ is
obtained by the following relations:
\[
 \tilde{v}_1=\tilde{v}_2,
\quad
 \tilde{v}_1\longmapsto v\in V_{0\,|\,\lambda},
\quad
 \tilde{v}_2\longmapsto v\in V_{\infty\,|\,\lambda^{-1}}.
\]
Here recall $\Phi^{\can}$ is obtained
by taking the degree $0$-part
of the polynomials $\sum v_J \cdot(\log z)^J$,
which gives the second correspondence.
The third correspondence can be obtained similarly.
Hence the gluing of $S^{\can}$ is same as
the gluing of $V$.

The comparison of $f_i$ and $\Res_i$
and the comparison of the pairings
should be checked only on $\nbigx$,
and it is easy.
\hfill\qed


%% file: b10.3.tex
\subsubsection{A lemma for the restriction of twistor nilpotent orbit}

We put $X=\cnum^l$ and $D=\bigcup_i D_i$.
Let $\pi:X\lrarr D_{\mbar}$ denote the projection
$(z_1,\ldots,z_l)\longmapsto (z_{m+1},\ldots,z_l)$.
For any point $Q\in D^{\circ}_{\mbar}$,
$\pi^{-1}(Q)$ is naturally isomorphic to $\Delta^m$.

From $(V,\vecf)$,
we obtain the variation $(\nbigv,\DD)$
of $\proj^1$-holomorphic bundles over $X-D$.
Let us consider the restriction
$\nbigv_{|\pi^{-1}(Q)}$.

We put $\tilde{Q}=(\overbrace{1,\ldots}^m,Q)\in X-D$.
We put $V':=\nbigv_{|\{\tilde{Q}\}\times\proj^1}$.
The vector bundle $V'$ is a twist of $V$ by the following
endomorphism of $V_{|\cnum_{\lambda}^{\ast}}$:
\[
 \exp\Bigl(
 -\sum_{i=m+1}^l f_i\cdot \log|z_i(Q)|^2\cdot t_1^{(1)}
 \Bigr).
\]
The tuple $\vecf'$ of the morphisms
$f_i':V'\lrarr V'\otimes\Tate(-1)$ $(i=1,\ldots,m)$
are naturally defined.
Then it is easy to see that we have the isomorphism
\[
 \nbigv(V',\vecf')
\simeq
 \nbigv(V,\vecf)_{|\pi^{-1}(Q)}.
\]
Thus we obtain the following lemma.
\begin{lem} \label{lem;9.23.3}
Let $\harmonicbundle$ be 
a harmonic bundle over $\Delta^{\ast\,l}$
corresponding a twistor nilpotent orbit.
Then the restriction
$\harmonicbundle_{|\pi^{-1}(Q)}$ is also a harmonic bundle
corresponding to a twistor nilpotent orbit.
\hfill\qed
\end{lem}

%% file: b6.5.tex

\subsubsection{Twist of Rees bundles}

We would like to see the relation of twistor nilpotent orbit
and the nilpotent orbit in the Hodge theory.

Let $H$ be a vector space with a decreasing filtration $F=(F^p)$.
For an endomorphism $g$ of $H$,
we put $g\cdot F:=(g\cdot F^p)$,
which is called the twist of $F$ by $g$.
Then we obtain the left $Aut(H)$-action on
the set of filtrations of $H$.

We have the natural isomorphism
$i_F:\xi(H,F)_{|X\times \cnum_{\lambda}^{\ast}}\simeq p_1^{\ast}H$.
For an element $g\in Aut(H)$,
we have the natural isomorphism
$g:(H,F)\lrarr (H,g\!\cdot\! F)$.
Then it induces the isomorphism
$\phi_g:\xi(H,F)\lrarr \xi(H,g\!\cdot\! F)$.

Let $Aut_{eq}(p_1^{\ast}H)$
be equivariant automorphisms of $p_1^{\ast}H$.
Clearly we have the natural isomorphism
$Aut_{eq}(p_1^{\ast}H)\simeq Aut(H)$.
We do not distinguish them.

Let $g$ be an element of $Aut(H)$.
Then we have the twist of $i_F$,
i.e.,
$g\circ i_F:\xi(H,F)_{|X\times\cnum_{\lambda}^{\ast}}
\simeq p_1^{\ast}(H)$.
Then we have the following commutative diagramm:
\[
 \begin{CD}
 \xi(H,F)_{|X\times\cnum_{\lambda}^{\ast}}
 @>{g\circ i_F}>> p_1^{\ast}H\\
 @V{\phi_g}VV @V{\id}VV \\
 \xi(H,g\!\cdot\! F)_{|X\times\cnum_{\lambda}^{\ast}}@>{i_F}>> p_1^{\ast}H.
 \end{CD}
\]

Let $(H,F,G)$ be bi-filtered vector bundle.
The Rees bundle $\xi(H,F,G)$ is obtained by the following
gluing:
\[
\begin{CD}
 \xi(H,F)_{|X\times\cnum_{\lambda}^{\ast}}
 @>{i_F}>>
 p_1^{\ast}H
 @<{i_G}<<
 \xi(H,G)_{|X\times\cnum_{\mu}^{\ast}}
\end{CD}
\]
Let $g_i$ $(i=1,2)$ be element of $Aut(H)$.
Then the vector bundle $\xi(H,F,G,g_1,g_2)$,
is obtained as the twisting $\xi(H,F,G)$ by $g_i$:
\[
 \begin{CD}
 \xi(H,F)_{|X\times\cnum_{\lambda}^{\ast}}
 @>{g_1\circ i_F}>>
 p_1^{\ast}H
 @<{g_2\circ i_G}<<
 \xi(H,G)_{|X\times\cnum_{\mu}^{\ast}}
\end{CD}
\]

\begin{lem} \label{lem;10.5.18}
$\xi(H,F,G,g_1,g_2)$ is naturally isomorphic to
$\xi(H,g_1\!\cdot\! F,g_2\!\cdot\! G)$.
\end{lem}
\pf
We have the following commutative diagram:
\[
 \begin{CD}
 \xi(H,F)_{|X\times\cnum_{\lambda}^{\ast}}
 @>{g_1\circ i_F}>>
 p_1^{\ast}H
 @<{g_2\circ i_G}<<
 \xi(H,G)_{|X\times\cnum_{\mu}^{\ast}}\\
 @V{\phi_{g_1}}VV @V{\id}VV @V{\phi_{g_2}}VV \\
 \xi(H,g_1\cdot F)_{|X\times\cnum_{\lambda}^{\ast}}
 @>{i_F}>>
 p_1^{\ast}H
 @<{i_G}<<
 \xi(H,g_2\cdot G)_{|X\times\cnum_{\mu}^{\ast}}
 \end{CD}
\]
It gives the isomorphism desired.
\hfill\qed

\subsubsection{The induced variation from an equivariant nilpotent tuple}

We put $X=\cnum^n$, $D_i=\{z_i=0\}$ and $D=\bigcup_{i=1}^n D_i$.
We put $\tilde{X}=\cnum^n$,
and then we have the universal covering
$\pi:\tilde{X}\lrarr X-D$,
given by the correspondence
$\zeta_i\longmapsto \exp\bigl(\sqrt{-1}\zeta_i\bigr)$.

From $(V,\vecf)$, we obtain the $\proj^1$-holomorphic
vector bundle $\nbigv$ 
and the variation $\DD^{\sankaku}_{\nbigv}$
over $\proj^1\times (X-D)$
(the subsubsection \ref{subsubsection;10.5.15}).

\begin{lem}
If $(V,\vecf)$ is equivariant,
then the vector bundle $\nbigv$ and the variation
$\DD^{\sankaku}_{\nbigv}$ are naturally equivariant.
\end{lem}
\pf
Since $-\sqrt {-1}f_i\otimes t_{1}^{(1)}$ is equivariant,
the gluing (\ref{eq;10.5.16}) is equivariant. Thus we have the natural
torus action on $\nbigv$.
The following sections, appearing in (\ref{eq;10.5.10})
and (\ref{eq;10.5.11}),
are invariant with respect to the torus action:
\[
 f_i\otimes t_0^{(1)}\otimes \frac{dz_i}{z_i}\cdot f_0^{(1,0)},
\quad
 f_i\otimes t_{\infty}^{(1)}\otimes\frac{d\overline{z}_i}{\overline{z}_i}
 \cdot f_{\infty}^{(0,1)}.
\]
Hence $\DD^{\sankaku}_{\nbigv}$ is also equivariant.
\hfill\qed

\vspace{.1in}

We put $\nbigh:=\nbigv_{|\{1\}\times (X-D)}$.
Then we obtain the two filtrations
$\nbigf$ and $\nbigg$
such that
$\xi(\nbigh,\nbigf,\nbigg)\simeq
 \nbigv$.
We also have the flat connection $\DD_0$ on $\nbigh$.

On $P_0:=\overbrace{(1,\ldots,1)}^n\in X-D$,
the fiber $\nbigv_{|\proj^1\times\{P_0\}}$ is naturally identified with $V$,
and $\nbigh_{|(1,P_0)}\simeq H$.

Let us consider $\pi^{\ast}\nbigh$.
We have the flat connection $\pi^{\ast}\DD_0$
and the isomorphism $\pi^{\ast}\nbigh_{|O}\simeq H$.
They induce the isomorphism
$\pi^{\ast}\nbigh\simeq p^{\ast}H$,
where $p$ denotes the natural morphism
$\cnum^n\lrarr \{1\}\subset \proj^1$.
We have the two filtrations
$\pi^{\ast}\nbigf$ and $p^{\ast}F$.
We also have
$\pi^{\ast}\nbigg$ and $p^{\ast}G$.

We put $N_i:=-(f_i\otimes t_1^{(1)})_{|\lambda=1}$.
Recall Lemma \ref{lem;10.17.1}
and Remark \ref{rem;10.17.5}.

\begin{lem} \label{lem;10.17.8}
The following equalities of filtrations hold:
\[
 \pi^{\ast}\nbigf
=\exp\Bigl(
 \sum \zeta_i N_i
 \Bigr)\cdot p^{\ast}F,
\quad\quad
  \pi^{\ast}\nbigg
=\exp\Bigl(
 \sum \overline{\zeta}_i N_i
 \Bigr)\cdot p^{\ast}G.
\]
\end{lem}
\pf
We put as follows:
\[
 g_0:=\exp\bigl(\sum_i \zeta_i\cdot N_i\bigr),
\quad
 g_{\infty}:=\exp\bigl(
 \sum_i \zetabar_i\cdot N_i
 \bigr).
\]
Let $v$ be an element of $H\simeq \nbigv_{|(O,\lambda)}$.
We put $\tilde{v}_0:=g_0\cdot v$.
Then it is the flat section of
$\pi^{\ast}\nbigv_{0\,|\,\nbigx^{\lambda}}$
such that $\tilde{v}_{0|(O,\lambda)}=v$.

We put $\tilde{v}_{\infty}:=g_{\infty}\cdot v$ .
Then it is the flat section of
$\pi^{\ast}\nbigv_{\infty\,|\,\nbigx^{\dagger\,\lambda}}$
such that $\tilde{v}_{2\,|\,(O,\lambda)}=v$.

The construction of the vector bundle $\nbigv$
is given by the relation
$\tilde{v}_{0}=\tilde{v}_{\infty}$,
due to Lemma \ref{lem;10.17.6}.
Hence the vector bundle $\nbigv$ is given by the following gluing:
\[
 \begin{CD}
 \xi\bigl(p^{\ast}H,p^{\ast}F\bigr)_{\tilde{X}\times\cnum_{\lambda}^{\ast}}
 @>g_0^{-1}\circ i_F>>
 p_1^{\ast}H @<{g_{\infty}^{-1}\circ i_G}<<
 \xi\bigl(p^{\ast}H,p^{\ast}G\bigr)_{\tilde{X}\times\cnum_{\mu}^{\ast}}.
 \end{CD}
\]
Here $p_1$ denote the canonical morphism of
$\tilde{X}\times\cnum_{\lambda}^{\ast}$ to a point.
Then we obtain the following equality,
due to Lemma \ref{lem;10.5.18}
\[
 \xi\bigl(\pi^{\ast}\nbigh,\pi^{\ast}\nbigf,\pi^{\ast}\nbigf\bigr)
\simeq
 \xi\Bigl(p^{\ast}H,\,
 g_0
 \!\cdot\! 
  p^{\ast}F,\,
 g_{\infty}
 \cdot\! p^{\ast}G
 \Bigr).
\]
Thus we are done.
\hfill\qed

\subsubsection{Reword}

Let $(H,F,G)$ be a bi-filtered vector space,
and $N_i:(H,F,G)\lrarr (H,F,G)\otimes \cnum(1,1)$
be a morphism.
We put $V=\xi(H,F,G)$,
and then we have $f_i:V\lrarr V\otimes\Tate(-1)$.

We have the trivial local system $p^{\ast}H$
on $\cnum^n$.
Then we obtain the $C^{\infty}$-bundle
$\nbigh^{(1)}$
with the natural flat connection $\DD^{(1)}$.
We have the $\seisuu^n$-action on $\cnum^n$
by $\veczeta\longmapsto \veczeta+\vecn\cdot 2\pi\sqrt{-1}$.
We lift it to the action on $\nbigh^{(1)}$
as follows:
\[
 (\veczeta,v)\longmapsto
 \Bigl(\veczeta+2\pi\cdot\vecn,\,\,
 \prod \exp\bigl(2\pi\cdot n_i\cdot N_i\bigr)\cdot v
 \Bigr).
\]
Here we regard $N_i$ as the endomorphism of
$H$ by using the isomorphism of $\cnum(1,1)\simeq \cnum(0,0)$
given by $t^{(-1)}_{1\,|\,\lambda=1}
 \longmapsto t^{(0)}_{1\,|\,\lambda=1}$.
We have two filtrations on $\nbigh^{(1)}$:
\[
 \exp\Bigl(\sum_i\zeta_i \cdot N_i\Bigr)\cdot p^{\ast}F,
\quad\quad
 \exp\Bigl(\sum \overline{\zeta}_i \cdot N_i\Bigr)\cdot p^{\ast}G.
\]
We obtain the $C^{\infty}$-bundle $\nbigh$
with the flat connection
on $\cnum^{\ast\,n}$,
and the filtrations $\nbigf$ and $\nbigg$.
They satisfy the Griffiths transversality,
and $\xi(\nbigh,\nbigf,\nbigg)\simeq Pat(V,\vecf)$.

Let $(H,F,\overline{F})$ be a bi-filtered vector space
defined over $\real$.
We put $V:=\xi(H,F,\overline{F})$.
Let $N_i:(H,F,\overline{F})\lrarr (H,F,\overline{F})\otimes \cnum(1,1)$
be morphisms.
The morphisms $N_i$ induce the morphisms
from $V\lrarr V\otimes\Tate(-1)$,
which we also denote by $N_i$.
We put $f_i:=-N_i$.
Then we obtain the tuple of endomorphisms $\vecf=(f_i)$.
\begin{cor} \label{cor;10.5.21}
We have the isomorphism
$\xi(\nbigh,\nbigf,\overline{\nbigf})
\simeq \nbigv(V,\vecf)$.
\end{cor}
\pf
This is a reformulation of 
Lemma \ref{lem;10.17.8}.
\hfill\qed

\begin{prop}\label{prop;b11.11.1}
The functor $\xi$ gives an equivalence of the following categories:
\begin{itemize}
\item
 The equivariant twistor nilpotent orbit defined over $\real$.
\item
 The nilpotent orbit in the category of Hodge theory
 in the sense of Schmid
 (Definition {\rm1.14} in {\rm\cite{cks1}}, for example).
\end{itemize}
\end{prop}
\pf
It follows from the various equivalences
(Corollary \ref{cor;10.5.20} and Corollary \ref{cor;10.5.21}).
\hfill\qed

%% file: a71.1.tex

\subsubsection{Definition}

\begin{df}
Let $(V,W,\vecf,S)$ be a Pol-MTS of type $(n,l)$.
Assume that the grading $V=\bigoplus V_h$ is given,
such that the following holds:
\begin{itemize}
\item $V_h$ is pure twistor of weight $h$.
\item $W_h=\bigoplus_{i\leq h}V_i$.
\item $f_j$ preserves the grading.
\item The restriction of $S$ to $V_i\otimes V_j$ is $0$
 unless $i+j=n$.
\end{itemize}
In that case, $(V,W,\vecf,S)$ is called
a split Pol-MTS of type $(n,l)$.

If $(V,W,\vecf,S)$ is a $\Psi$-Pol-MTS of type $(n,l)$,
and if the grading satisfying the above conditions is given,
then $(V,W,\vecf,S)$ is called a split $\Psi$-Pol-MTS of type $(n,l)$
\hfill\qed
\end{df}

%% file: a48.2.tex

\subsubsection{Preliminary on the split Pol-MTS in one variable of rank 2}

We put $V^{[2]}:=\nbigo(1,0)\oplus\nbigo(0,-1)$.
We have the filtration given by
$W_{-1}=\nbigo(0,-1)\subset W_1=V^{[2]}$.
Then $(V^{[2]},W)$ is a mixed twistor structure.

Let $F^{[2]}:V^{[2]}\lrarr V^{[2]}\otimes\Tate(-1)$ be
the morphism given by
$f_x^{(1,0)}\longmapsto f_x^{(0,-1)}\otimes t_x^{(-1)}$
and $f_x^{(0,-1)}\longmapsto 0$
for $x=0,1,\infty$.
We put $N=F^{[2]}\otimes t_1^{(1)}$.

We put $V^{[2]}_0:=V^{[2]}_{|\cnum_{\lambda}}$
and $V^{[2]}_{\infty}:=V^{[2]}_{|\cnum_{\mu}}$.
We have the frames $\bigl(f_0^{(1,0)},f_0^{(0,-1)}\bigr)$
and $\bigl(f_{\infty}^{(1,0)},f_{\infty}^{(0,-1)}\bigr)$
of $V^{[2]}_0$ and $V^{[2]}_{\infty}$ respectively.
We have the frame
$\bigl(
 f_1^{(1,0)},f_{1}^{(0,-1)}
 \bigr)$ of
$V^{[2]}_{0\,|\,\cnum_{\lambda}^{\ast}}
=V^{[2]}_{\infty\,|\,\cnum_{\mu}^{\ast}}$.
We use the notation
$f_1^{(1,0)\,\dagger}$ and $f_1^{(0,-1)\,\dagger}$,
when we consider them
as the frame of $V^{[2]}_{\infty\,|\,\cnum_{\mu^{\ast}}}$.

We twist the gluing as follows:
Let $\lambda$ be a point of $\cnum_{\lambda}^{\ast}$.
Let $v$ be an element of $V^{[2]}_{|\lambda}$.
It induces the element of $v\in V^{[2]}_{0\,|\,\lambda}$
and  $v^{\dagger}\in V^{[2]}_{\infty\,|\,\mu}$, where we put
$\mu=\lambda^{-1}$.
The original gluing is given by $v=v^{\dagger}$.
The twisted gluing is given by the following relation:
\[
 \exp\bigl(
 \sqrt{-1}y\cdot N
 \bigr)v=v^{\dagger}.
\]
Note that $y\cdot N$ gives the following correspondence:
\[
 f_1^{(1,0)}\longmapsto y\cdot f_1^{(0,-1)},
\quad\quad
 f_1^{(0,-1)}\longmapsto 0.
\]
Thus the gluing is given by the following:
\[
 \bigl(f_1^{(1,0)},f_1^{(0,-1)}\bigr)
\cdot
 \left(
 \begin{array}{cc}
 1          & 0 \\
 \sqrt{-1}y & 1
 \end{array}
 \right)
=\bigl(
 f_1^{(1,0)\,\dagger},
 f_1^{(0,-1)\,\dagger}
 \bigr).
\]
Then we obtain the following relation:
\[
 \bigl(
 f_0^{(1,0)},f_0^{(0,-1)}
 \bigr)
\cdot
 \left(
 \begin{array}{cc}
 \sqrt{-1}\lambda & 0 \\
 \sqrt{-1} y & -\sqrt{-1}\mu
 \end{array}
 \right)
=\bigl(
 f_{\infty}^{(1,0)},f_{\infty}^{(0,-1)}
 \bigr).
\]
Let $\tilde{V}^{[2]}_y$ denote the resulted vector bundle.

\begin{lem}
Assume that $y>0$.
The vector bundle $\tilde{V}^{[2]}_y$  is a pure twistor of weight $0$
of rank $2$.
The tuple of global sections $(\tilde{s}_1,\tilde{s}_2)$,
which are given as follows,
is a base of the space of the global sections:
\begin{equation}\label{eq;10.5.25}
 \begin{array}{l}
 \tilde{s}_1:=\sqrt{-1}\lambda\cdot f_0^{(1,0)}
+\sqrt{-1}y\cdot f_0^{(0,-1)}
=f_{\infty}^{(1,0)},\\
 \mbox{{}}\\
 \tilde{s}_2:=f_0^{(1,0)}
=-\sqrt{-1}\mu\cdot f_{\infty}^{(1,0)}
 -\sqrt{-1}y\cdot f_{\infty}^{(0,-1)}.
 \end{array}
\end{equation}
\end{lem}
\pf
It can be checked by direct calculations.
\hfill\qed

\vspace{.1in}

The pairing $\eta_1:\nbigo(1,0)\otimes\sigma\bigl(\nbigo(0,-1)\bigr)
\lrarr \Tate(0)$ is given as the composite of the following
naturally defined morphisms:
\[
 \nbigo(1,0)\otimes\sigma\bigl(\nbigo(0,-1)\bigr)
\lrarr
 \nbigo(1,0)\otimes\nbigo(-1,0)
\lrarr
 \Tate(0).
\]
The pairing $\eta_2:\nbigo(0,-1)\otimes\sigma\bigl(\nbigo(1,0)\bigr)
\lrarr\Tate(0)$ is given similarly.

Let us consider the pairing
$S^{[2]}:V^{[2]}\otimes V^{[2]}\lrarr \Tate(0)$
given as follows:
The restriction of $S^{[2]}$ to
$\nbigo(1,0)\otimes\sigma\bigl(\nbigo(1,0)\bigr)
 \oplus \nbigo(0,-1)\otimes\sigma\bigl(\nbigo(0,-1)\bigr)$
is defined to be trivial.
The restriction of $S^{[2]}$ to 
$\nbigo(1,0)\otimes\sigma\bigl(\nbigo(0,-1)\bigr)$
is defined to be $-\eta_1$.
The restriction of $S^{[2]}$
to $\nbigo(0,-1)\otimes\sigma\bigl(\nbigo(1,0)\bigr)$
is defined to be $\eta_2$.
\begin{lem}
We have
$S^{[2]}\bigl(f_0^{(1,0)}\otimes
  \sigma\bigl(f_{\infty}^{(0,-1)}\bigr)\bigr)
=-\sqrt{-1}t_0^{(0)}$
and
$S^{[2]}\bigl(
 f_0^{(0,-1)}\otimes\sigma\bigl(
 f_{\infty}^{(1,0)}
 \bigr)
 \bigr)=-\sqrt{-1}t_0^{(0)}$.
\end{lem}
\pf
It can be checked by a direct calculation.
\hfill\qed

\begin{lem} \mbox{{}}\label{lem;10.5.30}
\begin{itemize}
\item
 The pairing $S^{[2]}$ is symmetric.
\item
 We have the relation
 $S^{[2]}\bigl(\id\otimes\sigma(F^{[2]})\bigr)
 +S^{[2]}\bigl(F^{[2]}\otimes\id\bigr)$.
\item
The induced pairing $-S^{[2]}(\id\otimes\sigma(F^{[2]}))$
on $\nbigo(1,0)$ gives a polarization of weight $1$.
\end{itemize}
\end{lem}
\pf
By a direct calculation,
we have 
$ S^{[2]}\bigl(
 f_0^{(1,0)}\otimes \sigma(f_{\infty}^{(0,-1)})
 \bigr)
=-\bigl(
 \sqrt{-1}\cdot t_0^{(0)}\bigr)$
and
$S^{[2]}\bigl(
 f_{\infty}^{(0,-1)}
\otimes\sigma\bigl(
 f_{0}^{(1,0)}
 \bigr)
 \bigr)=\sqrt{-1}t_{\infty}^{(0)}$.
Thus we obtain
$\sigma\bigl(
 S^{[2]}\bigl(f_0^{(1,0)}\otimes \sigma(f_{\infty}^{(0,-1)})
 \bigr)
 \bigr)
=S^{[2]}\bigl(
 f_{\infty}^{(0,-1)}\otimes\sigma\bigl(
 f_0^{(1,0)}
 \bigr)
 \bigr)$.
Hence the pairing $S^{[2]}$ is symmetric.

Let us show the second claim.
The morphism $S^{[2]}\bigl(\id\otimes \sigma(F^{[2]})\bigr)$
gives the composite of the following correspondence:
\[
 f_0^{(1,0)}\otimes \sigma\bigl(f_{\infty}^{(1,0)}\bigr)
\longmapsto
 f_0^{(1,0)}\otimes
 \sigma\bigl(f_{\infty}^{(0,-1)}\otimes t_{\infty}^{(-1)}\bigr)
\longmapsto
 -\sqrt{-1}\cdot
 \bigl( -t_0^{(-1)} \bigr)=\sqrt{-1}\cdot t_0^{(-1)}.
\]
On the other hand, 
$S^{[2]}\bigl(F^{[2]}\otimes \id\bigr)$
gives the composite of the following correspondence:
\[
 f_0^{(1,0)}\otimes \sigma\bigl(f_{\infty}^{(1,0)}
 \bigr)
\longmapsto
 f_0^{(0,-1)}\otimes t_0^{(-1)}\otimes
 \sigma\bigl(f_{\infty}^{(1,0)}\bigr)
\longmapsto
 -\sqrt{-1}\cdot t_0^{(-1)}.
\]
Thus we obtain the second claim.

Let us show the third claim.
A base of the space $H^0\bigl(\nbigo(1,0)\otimes\nbigo(-1,0)\bigr)$
is given by
$s=f_0^{(1,0)}\otimes f_0^{(-1,0)}
=f_{\infty}^{(1,0)}\otimes f_{\infty}^{(-1,0)}$.
We have the following:
\[
 S^{[2]}\bigl(F^{[2]}(s),\sigma(s)\bigr)=
 S^{[2]}\bigl(
 f_{\infty}^{(0,-1)}\otimes t_0^{(-1)}\otimes f_{\infty}^{(-1,0)},
 \sigma\bigl(
 f_0^{(1,0)}\otimes f_0^{(-1,0)}
 \bigr)
 \bigr)
=1.
\]
Thus we are done.
\hfill\qed

\begin{lem}
The tuple $(V^{[2]},W,F^{[2]},S^{[2]})$
is a split polarized mixed twistor structure
of weight $0$.
\end{lem}
\pf
It immediately follows from Lemma \ref{lem;10.5.30}.
\hfill\qed

\vspace{.1in}

On the other hand,
the pairing $S^{[2]}$ induces
the pairing $\tilde{S}^{[2]}$ on $\tilde{V}^{[2]}_y$,
for we have the relation:
\[
 S^{[2]}\bigl(\id\otimes\sigma(F^{[2]}\otimes t_1^{(1)})\bigr)
+S^{[2]}\bigl(F^{[2]}\otimes t_1^{(1)}\otimes\id\bigr)=0.
\]
\begin{lem}
The pairing $\tilde{S}^{[2]}$ is a polarization of $\tilde{V}^{[2]}_y$
of weight $0$.
\end{lem}
\pf
We have only to show the positivity
$\tilde{S}^{[2]}(\tilde{s}_i,\tilde{s}_i)>0$
for the sections $\tilde{s}_i$ $(i=1,2)$
given in (\ref{eq;10.5.25}).
As for $\tilde{s}_1$, we have the following:
\[
\tilde{S}^{[2]}
\bigl(
 \tilde{s}_1,\tilde{s}_1
\bigr)
=
 S^{[2]}\bigl(
 \sqrt{-1}\lambda\cdot f_0^{(1,0)}
+\sqrt{-1}y\cdot f_0^{(0,-1)},\,
 \sigma\bigl(
 f_{\infty}^{(1,0)}
 \bigr)
 \bigr)
=\sqrt{-1}y\cdot S^{(0)}\bigl(
 f_0^{(0,-1)},
 \sigma\bigl(
 f_{\infty}^{(1,0)}\bigr)
 \bigr)
=y\cdot\sqrt{-1}\cdot\sqrt{-1}^{-1}=y.
\]
As for $\tilde{s}_2$, we have the following:
\begin{multline}
  \tilde{S}^{[2]}(\tilde{s}_2,\tilde{s}_2)
=S^{[2]}\bigl(
 -\sqrt{-1}\mu\cdot f_{\infty}^{(1,0)}
 -\sqrt{-1}y\cdot f_{\infty}^{(0,-1)},
 \sigma\bigl(
 f_0^{(1,0)} 
 \bigr)
 \bigr)
=S^{[2]}\bigl(
 -\sqrt{-1}y\cdot f_{\infty}^{(0,-1)}, 
 \sigma\bigl(
 f_0^{(1,0)}
 \bigr)
 \bigr) \\
=-\sqrt{-1}y\cdot S^{[2]}\bigl(
 f_{\infty}^{(0,-1)},
 \sigma\bigl(
 f_{0}^{(1,0)}
 \bigr)
 \bigr)=-\sqrt{-1}y\cdot \sqrt{-1}=y.
\end{multline}
Hence we have $\tilde{S}^{[2]}(\tilde{s}_i,\tilde{s}_i)>0$
for $i=1,2$.
Thus we are done.
\hfill\qed

\subsubsection{Preliminary on the split Pol-MTS of rank h}

Let $h>1$ be an integer.
We put $V^{(1)}:=V^{[2]\,\otimes\,h-1}$.
We have the naturally defined pairing
$S^{(1)}:V^{(1)}\otimes\sigma(V^{(1)})\lrarr \Tate(0)$,
given as follows:
\[
 S^{(1)}\Bigl(
 \bigotimes_{i=1}^{h-1} f_x^{(p_i,q_i)},\,\,
 \bigotimes_{i=1}^{h-1} \sigma\bigl(f_x^{(\tilde{p}_i,\tilde{q}_i)}\bigr)
 \Bigr)
=\prod_{i=1}^{h-1}
 S^{[2]}\bigl(
 f_x^{(p_i,q_i)},
 \sigma(f_x^{(\tilde{p}_i,\tilde{q}_i)})
 \bigr).
\]
Here $(p_i,q_i)$ and $(\tilde{p}_i,\tilde{q}_i)$
denote $(1,0)$ or $(0,-1)$.
We also have the morphism
$F^{(1)}:V^{(1)}\lrarr V^{(1)}\otimes\Tate(-1)$,
which is induced from $F^{[2]}$ by the Leibniz rule.

We have the natural $\gbigs_{h-1}$-action on $V^{(1)}$.
It preserves $S^{(1)}$ and $F^{(1)}$.
Then we obtain the invariant part $V^{[h]}=\Sym^{h-1}(V^{(0)})$
and the induced pairing $S^{[h]}$ and the induced
morphism $F^{[h]}$.
We also have the induced filtration $W$ of $V^{[h]}$.
The following lemma is clear.
\begin{lem}
We have the natural grading:
\begin{equation} \label{eq;10.17.10}
 V^{[h]}\simeq
 \bigoplus_{\substack{p+q=h-1,\\p,q\geq 0}}
 \nbigo(1,0)^{\otimes\,p}\otimes\nbigo(0,-1)^{\otimes\,q}.
\end{equation}
The filtration of the left hand side in {\rm (\ref{eq;10.17.10})}
is isomorphic to the following filtration of the right hand side:
\[
 W_a=\bigoplus_{p-q\leq a}
 \nbigo(1,0)^{\otimes\,p}\otimes\nbigo(0,-1)^{\otimes\,q}.
\]
In particular, the vector bundle $V^{[h]}$ with the filtration $W$
is a mixed twistor structure.
\hfill\qed
\end{lem}

\label{subsubsection;10.17.11}
\begin{lem}
The tuple $\bigl(V^{[h]},W,F^{[h]},S^{[h]}\bigr)$ is
a split polarized mixed twistor structure of weight $0$.
\end{lem}
\pf
The condition
$S^{(1)}\bigl(\id\otimes\sigma\bigl(F^{(1)}\bigr)\bigr)
+S^{(1)}\bigl(F^{(1)}\otimes\id\bigr)=0$
can be checked easily.
Then 
$S^{[h]}\bigl(\id\otimes\sigma\bigl(F^{[h]}\bigr)\bigr)
+S^{[h]}\bigl(F^{[h]}\otimes\id\bigr)=0$
immediately follows.

It is easy to see that there exists a positive number $B$
such that  the following holds:
\[
 \big(F^{[h]}\big)^h(f_x^{(1,0)\,\otimes\,h-1})
=B\cdot f_x^{(0,-1)\,\otimes\,h-1}.
\]
A base of the space of the global sections
of $\Gr^{W}_{h-1}\otimes\nbigo(-h+1,0)$ is given by the following:
\[
 s=
 f_0^{(1,0)\,\otimes\,h-1}
 \otimes f^{(-h+1,0)}_0
=f_{\infty}^{(1,0)\otimes\,h-1}
 \otimes f^{(-h+1,0)}_{\infty}.
\]
We have the following:
\[
  S^{[h]}\bigl(
 F^{[h]\,h-1}(s),\sigma(s)
 \bigr)
=
B\cdot
 S^{[h]}\Bigl(
 f_{\infty}^{(0,-1)\,\otimes\,h-1}
\otimes
 f_{\infty}^{(-h+1,0)},\,\,
\sigma\bigl(
 f_0^{(1,0)\,\otimes\,h-1}
\otimes f_0^{(-h+1,0)}
 \bigr)
\Bigr)
=B.
\]
Thus $S^{[h]}\bigl(F^{[h]\,h-1}\otimes\id\bigr)$ is a polarization
on $\Gr^W_{h-1}=P\Gr^W_{h-1}$.
\hfill\qed

\vspace{.1in}
On the other hand,
we have the induced pairing $\tilde{S}^{[h]}$
on $\tilde{V}^{[h]}$.
\begin{lem} \label{lem;10.5.31}
$\tilde{S}^{[h]}$
is a polarization of $\tilde{V}^{[h]}$.
In particular,
$\bigl(\tilde{V}^{[h]},W,F^{[h]},\tilde{S}^{[h]}\bigr)$
is a polarized twistor structure of weight $0$.
\end{lem}
\pf
It is easy to see that $\tilde{S}^{(1)}$
is a polarization of $\tilde{V}^{(1)}$.
Then the lemma immediately follows.
\hfill\qed

\subsubsection{Classification of split Pol-MTS in one variable}

Let us consider the vector bundle
$V=\bigoplus_{p+q=h,0\leq p,q\leq h}\nbigo(p,-q)$
over $\proj^1$.
The filtration $W$ is given as follows:
\[
 W_a:=\bigoplus _{p-q\leq a} \nbigo(p,-q).
\]
Then $(V,W)$ is a mixed twistor structure.

\begin{lem} \label{lem;10.5.32}
Let $(V,W)$ be as above.
Let $F:V\lrarr V\otimes\Tate(-1)$ be a morphism of
mixed twistor structure preserving the grading.
Let $S:V\otimes\sigma(V)\lrarr\Tate(0)$ be a morphism
of mixed twistor structures.
Assume that $(V,W,F,S)$ be a split polarized mixed twistor structure
of weight $0$.
Then $(V,W,F,S)$ are determined uniquely up to isomorphisms.
\end{lem}
\pf
Up to isomorphisms,
we may assume that
$F:\nbigo(p,-q)\lrarr\nbigo(p-1,-q-1)\otimes\Tate(-1)$
is given by
$f_x^{(p,-q)}\longmapsto f_x^{(p-1,-q-1)}\otimes t^{(-1)}_x$
for $x=0,1,\infty$.
Since we have the relation
$S\bigl(\id\otimes\sigma(F)\bigr)+S\bigl(F\otimes\id\bigr)=0$,
we obtain the following:
\[
 S\bigl(f_{x}^{(p-1,-q-1)},\,\,\sigma(f_x^{(q+1,-p+1)})\bigr)
+S\bigl(f_x^{(p,-q)},\,\,\sigma(f_x^{(q,-p)})\bigr)=0.
\]
Hence $S$ is determined by the number
$C=S\bigl(f_{\infty}^{(0,-h)}\otimes\sigma(f_{0}^{(h,0)})\bigr)$.

We have the following:
\[
 S\Bigl(
 F^h\bigl(
 f^{(h,0)}_{\infty}\otimes f^{(-h,0)}\bigr),\,\,
 \sigma\bigl(
 f^{(h,0)}_{0}\otimes f^{(-h,0)}
 \bigr)
 \Bigr)
=C.
\]
Hence we obtain $C>0$.
We may assume that $C=1$ up to isomorphisms again.
Thus we are done.
\hfill\qed

\begin{cor} \label{cor;10.6.1}
Let $V$, $W$, $F$ and $S$ as above.
If $(V,W,F,S)$ is a split polarized mixed twistor structure,
then it is isomorphic to a polarized mixed twistor structure
of the form $(V^{[h]},W,F^{[h]},S^{[h]})$,
given in the subsubsection {\rm\ref{subsubsection;10.17.11}}.
\end{cor}
\pf
It immediately follows from Lemma \ref{lem;10.5.31}
and Lemma \ref{lem;10.5.32}.
\hfill\qed

\begin{cor} \label{cor;a11.10.1}
Let $(V,W,F,S)$ be a split polarized mixed twistor structure
of weight $n$.
Then it is isomorphic to a split Pol-MTS of weight $n$
of the following form:
\[
 \bigoplus_{i} \bigl(V^{[h_i]},W,F^{[h_i]},S^{[h_i]}\bigr)\otimes\nbigo(n)
\]
\end{cor}
\pf
We have only to consider the case $n=0$.
By taking the primitive decomposition,
we can reduce the problem to the case
$V=\bigoplus_{p+q=h,0\leq p,q\leq h}\nbigo(p,-q)$.
Then the lemma immediately follows from Corollary \ref{cor;10.6.1}.

\mbox{{}}\hfill\qed

\begin{lem} \label{lem;10.6.2}
Let $(V,W,F,S)$ be a split polarized mixed twistor structure
of weight $n$.
Then the twisted vector bundle $\tilde{V}$ is 
a pure twistor of weight $n$.
\end{lem}
\pf
It immediately follows from Corollary \ref{cor;a11.10.1}
and Lemma \ref{lem;10.5.31}.
\hfill\qed


%% file: a48.3.tex

\subsubsection{Pol-MTS and nilpotent orbit in one variable}

\begin{lem} \label{lem;10.6.15}
Let $(V,W,f,S)$ be a Pol-MTS of type $(n,1)$.
Then it is a nilpotent orbit of type $(n,1)$.
\end{lem}
\pf
We put $V^{(0)}:=\Gr^{W}(V)$.
Then we have the induced tuple
$(V^{(0)},W^{(0)},f^{(0)},S^{(0)})$.
We put as follows:
\[
 V_0:=V_{|\cnum_{\lambda}},\quad
 V_{\infty}:=V_{|\cnum_{\mu}},\quad
 V^{(0)}_{0}:=V^{(0)}_{|\cnum_{\lambda}},\quad
 V^{(0)}_{\infty}:=V^{(0)}_{|\cnum_{\mu}}.
\]
Let us take frames
$\vecu^{(0)}_0=(u^{(0)}_{0\,i})$ and
$\vecu^{(0)}_{\infty}=(u^{(0)}_{\infty\,i})$
of $V^{(0)}_0$ and $V^{(0)}_{\infty}$ respectively,
which are compatible with the natural grading.
We take frames $\vecu_0=(u_{0\,i})$ and
$\vecu_{\infty}=(u_{\infty\,j})$ of $V_0$ and $V_{\infty}$ respectively,
such that they induce $\vecu^{(0)}_0$ and $\vecu^{(0)}_{\infty}$
respectively.
We put as follows:
\[
 K(i):=\deg^W(u_{0\,i})=\deg^{W}(u^{(0)}_{0\,i}),
\quad\quad
 L(i):=\deg^{W}(u_{\infty\,i})=\deg^W(u^{(0)}_{\infty\,i}).
\]
We have the relations:
\[
 u^{(0)}_{\infty\,i}
=\sum B^{(0)}_{j\,i}\cdot u^{(0)}_{0\,j},
\quad
 u_{\infty\,i}
=\sum B_{j\,i}\cdot u_{0\,j}.
\]
\begin{lem}\mbox{{}} \label{lem;10.6.12}
\begin{itemize}
\item
We have $B_{j\,i}^{(0)}=0$
unless $L(i)=K(j)$.
\item
We have $B_{j\,i}=0$
unless $L(i)\geq K(j)$.
\item
We have $B_{j\,i}=B_{j\,i}^{(0)}$
in the case $L(i)=K(j)$.
\end{itemize}
\end{lem}
\pf
It immediately follows from our choice of the frames.
\hfill\qed

\vspace{.1in}
We have the following relations:
\[
 f^{(0)}\otimes t_1^{(-1)}(u_{0\,i}^{(0)})
=\sum A_{j\,i}^{(0)}\cdot u_{0\,j}^{(0)},
\quad
 f\otimes t_1^{(-1)}(u_{0\,i})
=\sum A_{j\,i}\cdot u_{0\,j}.
\]
\begin{lem}\mbox{{}} \label{lem;10.6.11}
\begin{itemize}
\item
 We have $A_{i\,j}^{(0)}=0$
unless $K(j)=K(i)-2$.
\item
 We have $A_{i\,j}=0$ unless $K(j)\leq K(i)-2$.
\item
 In the case $K(j)=K(i)-2$,
 we have $A_{i\,j}^{(0)}=A_{i\,j}$.
\end{itemize}
\end{lem}
\pf
It immediately follows from our choice of the frames.
\hfill\qed

\vspace{.1in}

Let us pick a point $P\in\Delta^{\ast}$,
and we put $y:=-\log|z(P)|^2>0$.
The restrictions of $\Pat(V^{(0)},f^{(0)})$ 
and $\Pat(V,f)$ to
$\proj^1\times\{P\}$ are given by the following gluings:
\[
  \vecu^{(0)}_{\infty}
=\vecu^{(0)}_0\cdot B^{(0)}\cdot \exp\bigl(\sqrt{-1}y\cdot A^{(0)}\bigr),
\quad
 \vecu_{\infty}
=\vecu_0\cdot B\cdot \exp\bigl(\sqrt{-1}y\cdot A\bigr).
\]
Let us consider the frames
$\vecu_{\infty}(y)$, $\vecu_0(y)$,
$\vecu_{\infty}^{(0)}(y)$ and $\vecu_0^{(0)}(y)$ given as follows:
\[
\begin{array}{ll}
 u_{\infty\,i}(y):=y^{-L(i)/2}\cdot u_{\infty\,i},
 &
 u_{0\,i}(y):=y^{-K(i)/2}\cdot u_{0\,i},\\
 \mbox{{}}\\
 u^{(0)}_{\infty\,i}(y):=y^{-L(i)/2}\cdot u^{(0)}_{\infty\,i},
 &
 u^{(0)}_{0\,i}(y):=y^{-K(i)/2}\cdot u^{(0)}_{0\,i}.
\end{array}
\]
Then it is easy to see that we have the following relation:
\begin{equation} \label{eq;10.6.10}
 \vecu_{\infty}^{(0)}(y)=\vecu_0^{(0)}(y)\cdot B^{(0)}\cdot 
 \exp\bigl(\sqrt{-1}A^{(0)}\bigr).
\end{equation}
Since $\Pat(V^{(0)},f^{(0)})$ is a variation of pure twistor structures,
the vector bundle whose gluing is given by (\ref{eq;10.6.10})
is pure twistor of weight $0$.

On the other hand, we have the following relation:
\[
 \vecu_{\infty}(y)=\vecu_0(y)\cdot
 B(y)\cdot \exp\bigl(\sqrt{-1}A(y)\bigr).
\]
Here $B(y)_{i\,j}$ and $A(y)_{i\,j}$ are given as follows:
\[
 B(y)_{i\,j}:=y^{(K(j)-L(i))/2}\cdot B_{i\,j},
\quad
 A(y)_{i\,j}:=y^{(K(j)+2-L(i))/2}\cdot A_{i\,j}.
\]
\begin{lem}
We have the following:
\[
 \lim_{y\to\infty} B(y)\cdot \exp\bigl(\sqrt{-1}A(y)\bigr)
=B^{(0)}\cdot \exp\bigl(\sqrt{-1}A^{(0)}\bigr).
\]
\end{lem}
\pf
It immediately follows from Lemma \ref{lem;10.6.12}
and Lemma \ref{lem;10.6.11}.
\hfill\qed

\vspace{.1in}
In particular, there exists a positive constant $\epsilon$
such that
the restriction of $\Pat(V,f)$
to $\Delta^{\ast}(\epsilon)$ is a variation of pure twistor structures.

Let us consider the pairing $\tilde{S}$ 
on $\Pat(V,f)$.
Note that it gives the non-degenerate hermitian pairing
of $H^0\bigl(\Pat(V,f)_{|\proj^1\times\{P\}}\bigr)$
for any point $P\in\Delta^{\ast}(\epsilon)$.

We have the following:
\[
 S\bigl(u_{0\,i}(y),\sigma\bigl(u_{\infty,j}(y)\bigr)\bigr)
=y^{-(L(i)+K(j))/2}\cdot S\bigl(u_{0\,i},\sigma(u_{\infty,j})\bigr).
\]
Hence the limit of $S$ is same as $S^{(0)}$
when we take the limit $y\to\infty$.
Hence we obtain the positive definiteness of $\tilde{S}$.
Thus the proof of \ref{lem;10.6.15} is accomplished.
\hfill\qed

\begin{lem} \label{lem;10.6.16}
Let $(V,W,f,S)$ be a split nilpotent orbit of weight $n$.
Then $(V,W,f,S)$ is a split Pol-MTS of type $(n,1)$.
\end{lem}
\pf
By taking the primitive decomposition,
we may assume
$V=\bigoplus_{p+q=n,0\leq p,q}\nbigo(p,-q)$,
and $W_a=\bigoplus_{p-q\leq a}\nbigo(p,-q)$.
By taking an appropriate isomorphisms,
we may assume that $f$ is naturally defined morphism.
Recall that $S$ is determined up to constant multiplication.
(See the proof of Lemma \ref{lem;10.5.32}.)
Then we may derive that
$(V,W,f,S)$ is isomorphic to 
a mixed twistor of the form $(V^{[h]},W^{[h]},f^{[h]},S^{[h]})$,
given in the subsubsection \ref{subsubsection;10.17.11}.
Thus $(V,W,f,S)$ is a Pol-MTS.
\hfill\qed

\begin{lem}
Let $(V,W,f,S)$ be a nilpotent orbit of type $(n,1)$.
Then $(V,W,f,S)$ is a Pol-MTS of type $(n,1)$.
\end{lem}
\pf
We put $V^{(0)}:=\Gr^{W}(V)$,
and then we have the induced tuple
$(V^{(0)},W^{(0)},f^{(0)},S^{(0)})$.
We use the notation in the proof of Lemma \ref{lem;10.6.15}.
We may assume that
the restriction of $\Pat(V,f)$ to $\Delta^{\ast}(C)$
is a variation of pure twistor
for some $0<C<1$.
Let pick a point $P\in\Delta^{\ast}(C)$.
Let us consider the frames
$\vecu_{\infty}(n)$, $\vecu_0(n)$,
$\vecu_{\infty}^{(0)}(n)$ and $\vecu_0^{(0)}(n)$ given as follows:
\[
\begin{array}{ll}
 u_{\infty\,i}(n):=n^{-L(i)/2}\cdot u_{\infty\,i},
 &
 u_{0\,i}(n):=n^{-K(i)/2}\cdot u_{0\,i},\\
 \mbox{{}}\\
 u^{(0)}_{\infty\,i}(n):=n^{-L(i)/2}\cdot u^{(0)}_{\infty\,i},
 &
 u^{(0)}_{0\,i}(n):=n^{-K(i)/2}\cdot u^{(0)}_{0\,i}.
\end{array}
\]
By an argument similar to the proof of Lemma \ref{lem;10.6.15},
the limit vector bundle is naturally isomorphic to
$\Pat(V^{(0)},f^{(0)})$,
when we take the limit $n\to\infty$.
Note that $\Pat(V^{(0)},f^{(0)})$ on $\Delta^{\ast}$
is always a variation of pure twistors.

The pairings $\tilde{S}(n)$ and $\tilde{S}^{(0)}$
induce the perfect hermitian product
of $H^0(\Pat(V,f)_{|P})$ and $H^0(\Pat(V^{(0)},f^{(0)})_{|P})$.
Since the limit of $\tilde{S}(n)$ is $\tilde{S}^{(0)}$,
and since $\tilde{S}(n)$ is a positive definite,
we obtain the positive definiteness of $\tilde{S}^{(0)}$.
It means that
$(V^{(0)},W^{(0)},f^{(0)},S^{(0)})$
is a nilpotent orbit.
Due to Lemma \ref{lem;10.6.16},
$(V^{(0)},W^{(0)},f^{(0)},S^{(0)})$
is a Pol-MTS.
It implies that $(V,W,f,S)$ is a Pol-MTS.
\hfill\qed

\vspace{.1in}

In all, we obtain the following:
\begin{prop}
Let $(V,W,f,S)$ be a $\Psi$-Pol-MTS of type $(n,1)$.
It is a nilpotent orbit of type $(n,1)$
if and only if
it is a Pol-MTS of type $(n,1)$.
\hfill\qed
\end{prop}

%% file: b6.6.tex
\subsubsection{The twistor nilpotent orbit of split type}

\begin{prop}\label{prop;9.13.1}
Let $(V,W,\vecf,S)$ be a split $\Psi$-Pol-MTS of type $(n,l)$.
In this case,
the induced variation of $\proj^1$-holomorphic bundles
$\bigl(\nbigv,\DD^{\sankaku}_{\nbigv}\bigr)$
a variation of pure twistors of weight $0$
over $\Delta^{\ast\,l}$.

Moreover,
let $S$ be a pairing of $V$ above such that
$\bigl(V,S,\sum a_i f_i\bigr)$ is a split Pol-MTS.
Then the induced pairing $S_{\nbigv}$ of $\nbigv$
gives a polarization of $\bigl(\nbigv,\DD^{\sankaku}_{\nbigv}\bigr)$.
\end{prop}
\pf
We have only to show the following:
\begin{itemize}
\item
 Let $P$ be a point of $\Delta^{\ast\,l}$.
 Then  $\nbigv_{|\{P\}\times\proj^1}$ is
 isomorphic to a trivial bundle.
\item
 Let $S$ be as above. Then $S_{\nbigv\,|\,\{P\}\times\proj^1}$
 gives a polarization of the pure twistor
 $\nbigv_{|\{P\}\times\proj^1}$.
\end{itemize}
Recall that $\nbigv_{|\{P\}\times\proj^1}$
is obtained as the twisting of the gluing of $V$ by
$\exp\Bigl(
 \sqrt{-1}\sum_{i=1}^l -\log |z_i(P)|^2\cdot f_i\otimes t_1^{(1)}
  \Bigr)$.
Note that we have $-\log|z_i(P)|^2>0$.
Then the proposition immediately follows from
Lemma \ref{lem;10.6.2}.
\hfill\qed

\begin{cor}\label{cor;b11.11.2}
Let $(V,W,\vecf,S)$ be a split Pol-MTS of type $(n,l)$.
Then it is a twistor nilpotent orbit.
\hfill\qed
\end{cor}

\subsubsection{A split Pol-MTS and a nilpotent orbit in Hodge theory}

Let $(V,W,\vecf,S)$ be a split Pol-MTS of type $(n,l)$.
Then we can pick a $G_m$-action on $(V,\vecf,S)$.
For example, it is given as follows:
For an integer $h$,
the integers $p(h)$ and $q(h)$ are defined as follows:
If we have $h=2m$ for some integer $m$,
then we put $p(h)=q(h):=m$.
If we have $h=2m+1$ for some integer $m$,
then we put $p(h):=m+1$ and $q(h):=m$.
Then the $G_m$-action on $V$
can be given by the isomorphism
$V\simeq \bigoplus_h V_h\otimes\nbigo\bigl(p(h),q(h)\bigr)$.
We denote the action by $\rho_1$.

\begin{lem}
$\vecf$ and $S$ are equivariant
with respect to $\rho_1$.
If we have a real structure $\iota$ of $V$
preserving the grading,
then it is also equivariant.
\end{lem}
\pf
The equivariance of $\vecf$ and $\iota$ is clear.
To see that $S$ is equivariant,
we have only to check the equivariance of the following morphism:
\[
 \nbigo\bigl(p(h_1),q(h_1)\bigr)\otimes
 \nbigo\bigl(p(h_2),q(h_2)\bigr)\lrarr
 \Tate(-n)
\]
Here we put $2n=h_1+h_2$.
But it is clear, for we have $p(h_1)+q(h_2)=p(h_2)+q(h_1)=n$.
\hfill\qed

\begin{cor}\label{cor;9.13.3}
Let $(V,\vecf,S)$ be a split Pol-MTS.
Let $\iota$ be a real structure of $(V,\vecf,S)$.
It gives a nilpotent orbit in the Hodge theory,
when we  take an appropriate $G_m$-action of $(V,\vecf,S)$.
\hfill\qed
\end{cor}

\begin{cor}\label{cor;a11.10.25}
Let $(V,\vecf,S)$ be a split Pol-MTS.
Then the tuple $(V,\vecf,S)\oplus \sigma(V,\vecf,S)$
gives a nilpotent orbit in the Hodge theory,
when we take an appropriate $G_m$-action.
\end{cor}
\pf
Since we have the canonical real structure
on $(V,\vecf,S)\oplus \sigma(V,\vecf,S)$,
it immediately follows from Corollary \ref{cor;9.13.3}.
\hfill\qed

%% file: a48.4.tex

\subsubsection{The nilpotent orbit on the divisors}

Let $(V,W,\vecf,S)$ be a nilpotent orbit of type $(n,l)$.
We put $V^{(1)}_h:=P\Gr^{W(f_1)}_h(V)$.
Then we have the induced filtration
$W^{(1)}$,
the induced morphisms
$f_i^{(1)}:V^{(1)}_h\lrarr V^{(1)}_h\otimes\Tate(-1)$ $(i=2,\ldots,l)$.
The pairing $S\bigl(f_1^h\otimes\id\bigr)$
induces the pairing $S^{(1)}_h$.
Thus we obtain the induced tuple
$\bigl(V^{(1)}_h,W^{(1)},\vecf^{(1)},S^{(1)}_h\bigr)$,
which is a $\Psi$-Pol-MTS of type $(n+h,l-1)$.

\begin{lem} \label{lem;10.17.20}
The tuple $\bigl(V^{(1)}_h,W^{(1)},\vecf^{(1)},S^{(1)}_h\bigr)$
is a nilpotent orbit of type $(n+h,l-1)$.
\end{lem}
\pf
We may assume that $n=0$.
We have only to show
that $\Pat(V^{(1)}_h,W^{(1)},\vecf^{(1)},S^{(1)}_h)$
is a variation of polarized pure twistor of weight $h$,
on $\Delta^{\ast\,n-1}$.
We identify $\Delta^{\ast\,n-1}$
and $D_1^{\circ}$, naturally.

Let us take a point $P\in D_1^{\circ}$.
Let $\tilde{P}$ denote the point of $\cnum^l$
such that $\{\tilde{P}\}=q_1^{-1}(1)\cap \pi_1^{-1}(P)$.
We put $V(P):=\Pat(V,\vecf)_{|\tilde{P}}$.
We have the induced filtration $W(P)$,
the induced morphism $f_1$,
and the induced pairing $S$
on $V(P)$.
Due to Lemma \ref{lem;9.23.3},
the tuple
$(V(P),W,f_1,S)$ is a nilpotent orbit.
In particular,
$(V(P),W,f_1,S)$ is a polarized MTS.
It implies
that
$\Pat(V^{(1)}_h,\vecf^{(1)})_{|P}$ is pure twistor of weight $h$,
and $S^{(1)}_h$ is a polarization
of $\Pat(V^{(1)}_h,\vecf^{(1)})_{|P}$.
Thus we are done.
\hfill\qed

%% file: a48.5.tex

\subsubsection{The Pol-MTS on the divisor}

Let $(V,W,\vecN,S)$ be a Pol-MTS of type $(n,l)$.
We put
$V^{(1)}_h:=P\Gr^{W(N_1)}_h(V)$,
and we obtain the mixed twistor structure
$W^{(1)}$ on $V^{(1)}$,
and the tuple of induced nilpotent morphisms
$N^{(1)}_2,\ldots,N^{(1)}_l$.
We also obtain the pairing
$ S^{(1)}_h:V_h^{(1)}\otimes \sigma(V_h^{(1)})
\lrarr \Tate(-n-h)$,
by putting $S^{(1)}_h:=S(N_1^h\otimes \id)$.
Then we obtain the tuple
$\bigl(V^{(1)}_h,W^{(1)},\vecN^{(1)},S^{(1)}_h\bigr)$.

\begin{lem} \label{lem;10.17.21}
If the tuple $(V,W,\vecN,S)$ is a split Pol-MTS of type $(n,l)$,
the induced tuple $(V^{(1)}_h,W^{(1)},\vecN^{(1)},S^{(1)})$ above
is a split Pol-MTS of type $(n+h,l-1)$.
\end{lem}
\pf
It is easy to check that
$\bigl(V^{(1)}_h,W^{(1)},\vecN^{(1)},S^{(1)}\bigr)$
is a split $\Psi$-Pol-MTS of type $(n+h,l-1)$.
Since $(V,W,\vecN,S)$ is a nilpotent orbit,
the tuple is $(V^{(1)}_h,W^{(1)},\vecN^{(1)},S^{(1)})$
is also nilpotent orbit, due to Lemma \ref{lem;10.17.20}.
In particular, it is a Pol-MTS of type $(n+h,l-1)$.
Thus we are done.
\hfill\qed

\vspace{.1in}

Let $(V,W,\vecN,S)$ be a $\Psi$-Pol-MTS of type $(n,l)$.
Let $(V^{(1)}_b,W^{(1)},\vecN^{(1)},S^{(1)})$ be
the induced $\Psi$-Pol-MTS of type $(n+b,l-1)$.
Then we obtain the following induced
$\Psi$-Pol-MTS of type $(n+b,l-1)$:
\[
 \Bigl(
 (V^{(1)}_b)^{(0)},(W^{(1)})^{(0)},(\vecN^{(1)})^{(0)},(S^{(1)})^{(0)}
 \Bigr).
\]
On the other hand,
we have the induced
$\Psi$-Pol-MTS $(V^{(0)},W^{(0)},\vecN^{(0)},S^{(0)})$
of type $(n,l)$ obtained from
$(V,W,\vecN,S)$.
Then we obtain the following induced $\Psi$-Pol-MTS
of type $(n+b,l-1)$ as is given above:
\[
 \Bigl(
 (V^{(0)})^{(1)}_b,
  (W^{(0)})^{(1)},
 (\vecN^{(0)})^{(1)},
 (S^{(0)})^{(1)} 
 \Bigr)
\]

\begin{lem} \label{lem;9.13.26}
We have the natural isomorphism:
\[
 \Bigl(
(V^{(1)}_b)^{(0)},(W^{(1)})^{(0)},(\vecN^{(1)})^{(0)},(S^{(1)})^{(0)}
 \Bigr)
\simeq
 \Bigl(
 (V^{(0)})^{(1)}_b,
  (W^{(0)})^{(1)},
 (\vecN^{(0)})^{(1)},
 (S^{(0)})^{(1)} 
 \Bigr).
\]
\end{lem}
\pf
We have the two filtrations on $\Gr^W_h$,
i.e.,
$W(N_1^{(0)})$ and $W(N_1)^{(0)}$.
We use the following lemma.
\begin{lem} \label{lem;9.13.25}
We have $W(N_1^{(0)})=W(N_1)^{(0)}$.
\end{lem}
\pf
We have only to check that
the filtration $W(N_1)^{(0)}$
satisfies the axioms of weight filtrations
for $N_1^{(0)}$.
It is easy to see
$N_1^{(0)}W(N_1)^{(0)}_h\subset
 W(N_1)^{(0)}_{h-2}\otimes\Tate(-1)$.
We have
$\Gr^{W(N_1)^{(0)}}_b\bigl(
 \Gr^W_h\bigr)
\simeq
 \Gr^{W}_h\bigl(
 \Gr^{W(N_1)}_b\bigr)$.
Since $N_1$ is a morphism of mixed twistor structures,
the following morphism is isomorphic:
\[
 N_1^b:\Gr^{W}_h\bigl(
 \Gr^{W(N_1)}_b\bigr)
\lrarr
 \Gr^{W}_{h-2b} \bigl(
 \Gr^{W(N_1)}_{-b}
 \bigr)\otimes\Tate(-b).
\]
Hence 
the following morphism is isomorphic:
\[
 (N_1^{(0)})^b:
\bigoplus_h
 \Gr^{W(N_1)^{(0)}}_b\bigl(
 \Gr^W_h\bigr)
\simeq
\bigoplus_h
 \Gr^{W(N_1)^{(0)}}_{-b}
 \bigl(
 \Gr^W_h\otimes\Tate(-b)\bigr).
\]
Thus we obtain Lemma \ref{lem;9.13.25}.
\hfill\qed

\vspace{.1in}
Thus we have the canonical isomorphism
$\Gr^{W(N_1)}_b\bigl(\Gr^W\bigr)\simeq
 \Gr^W\Gr^{W(N_1)}_b$.
We have the following commutative diagramm:
\[
 \begin{CD}
 \Gr^{W(N_1)}_b\bigl(\Gr^W\bigr)
 @>{\simeq}>>
 \Gr^W\Gr^{W(N_1)}_b\\
 @V{N_1^{b+1}}VV @V{N_1^{b+1}}VV \\
 Gr^{W(N_1^{(0)})}_{-b-2}\Gr^W\otimes\Tate(-h)
 @>{\simeq}>>
 \Gr^W\Gr^{W(N_1)}_{-b-2}\otimes\Tate(-b).
 \end{CD}
\]
The kernel of the left vertical arrow
is the primitive part $P\Gr^{W(N_1)}_b\Gr^W$
by definition.
On the other hand,
it is easy to see that the kernel of the right vertical arrow
is naturally isomorphic to
$\Gr^W P\Gr_b^{W(N_1)}$,
by using the primitive decomposition of $\Gr_b^{W(N_1)}$.
Thus we obtain the canonical isomorphism
$\bigl(V^{(0)}\bigr)_b^{(1)}\simeq
 \bigl(V^{(1)}_b\bigr)^{(0)}$.
Once we obtain the isomorphism of vector bundles,
then it is easy to see the coincidence
of the filtration, the nilpotent maps and the pairings.
Thus we obtain Lemma \ref{lem;9.13.26}.
\hfill\qed

\begin{prop}\label{prop;10.6.5}
Let $(V,W,\vecN,S)$ be a Pol-MTS of type $(n,l)$.
Then $(V_h^{(1)},W^{(1)},\vecN^{(1)},S^{(1)})$ is a Pol-MTS 
of type $(n+h,l-1)$.
\end{prop}
\pf
We have only to show that the tuple
$\bigl(
(V_h^{(1)})^{(0)},
(W^{(1)})^{(0)},
(\vecN^{(1)})^{(0)},
(S^{(1)})^{(0)}\bigr)$ is a Pol-MTS of type $(n+h,l-1)$.
Since the tuple $(V^{(0)},W^{(0)},\vecN^{(0)},S^{(0)})$
is a split Pol-MTS of type $(n,l)$,
the induced tuple
$((V_h^{(0)})^{(1)},
(W^{(0)})^{(1)},
(\vecN^{(0)})^{(1)},
(S^{(0)})^{(1)})$ is a Pol-MTS of type $(n+h,l-1)$
due to Corollary Lemma \ref{lem;10.17.21}.
Then we obtain the result from Lemma \ref{lem;9.13.26}.
\hfill\qed

%% file: b8.1.tex

\subsubsection{Strongly sequentially compatibility}

\begin{lem}\label{lem;10.26.111}
Let $(V,W,\vecN,S)$ be a Pol-MTS of type $(n,l)$.
Then $\vecN$ is strongly sequentially compatible.
\end{lem}
\pf
We use an induction on $l$.
Due to the hypothesis of the induction and Proposition
\ref{prop;10.6.5},
the tuple $\vecN^{(1)}=(N_2^{(1)},\ldots,N_l^{(1)})$
on $V^{(1)}=\bigoplus V_h^{(1)}$
is sequentially compatible.

Since $(V^{(0)},W^{(0)},\vecN^{(0)},S^{(0)})$ is a nilpotent orbit,
$\vecN^{(0)}$ is strongly sequentially compatible,
which was shown in our previous paper \cite{mochi}.
Hence the conjugacy classes of
$N^{(0)}(\veca)$ $(\veca\in\real_{>0}^I)$ 
are constant if we fix a subset $I\subset \lbar$.
Since the conjugacy classes of $N(\veca)$
and $N^{(0)}(\veca)$ are same,
we obtain the constantness of the conjugacy classes
of $N(\veca)$ for $\veca\in\real_{>0}^I$.

\begin{lem}\label{lem;a11.10.20}
We have
$W_{h+b}\bigl(N(\ibar)\bigr)^{(1)}\cap \Gr^{W(N_1)}_h
=W_b\bigl(N(\ibar)^{(1)}\bigr)\cap \Gr^{W(N_1)}_h$.
\end{lem}
\pf
Let $W^{(1)}$ denote the induced mixed twistor structure
on $\Gr^{W(N_1)}_h$.
Since $W^{(1)}$ gives the mixed twistor structure,
we have the following:
\begin{equation}\label{eq;a11.10.14}
\begin{array}{l}
 \Gr^{W^{(1)}}\bigl(
 W_{h+b}\bigl(N(\ibar)^{(1)}\bigr)
\cap
 \Gr^{W(N_1)}_h
 \bigr)
=W_{h+b}\bigl(
 N(\ibar)^{(1)\,(0)}
 \bigr)\cap \Gr^{W^{(1)}}\Gr^{W(N_1)}_h,\\
\mbox{{}}\\
 \Gr^{W^{(1)}}\bigl(
 W_b\bigl(N(\ibar)^{(1)}\bigr)\cap \Gr_h^{W(N_1)}
 \bigr)
=W_b\bigl(N(\ibar)^{(1)}\bigr)^{(0)}
\cap \Gr^{W^{(1)}}\Gr^{W(N_1)}_h.
\end{array}
\end{equation}
Since $\vecN^{(0)}$ is strongly sequentially compatible,
we have the following:
\begin{equation}\label{eq;a11.10.15}
 W_{h+b}\bigl(N(\ibar)^{(0)}\bigr)^{(1)}
\cap \Gr^{W(N_1)^{(0)}}_h\Gr^W
=W_b\bigl(N(\ibar)^{(0)\,(1)}\bigr)
\cap \Gr^{W(N_1)^{(0)}}_h\Gr^W.
\end{equation}

\begin{lem}\label{lem;a11.10.21}
Under the isomorphism $\Gr^W\Gr^{W(N_1)}_b\simeq \Gr^{W(N_1)^{(0)}}_h\Gr^W$,
we have the isomorphisms
\begin{equation}\label{eq;a11.10.16}
W_{h+b}(N(\ibar)^{(0)})^{(1)}\simeq
 W_{h+b}(N(\ibar))^{(1)\,(0)},
\quad\quad
W_b\bigl(N(\ibar)^{(0)\,(1)}\bigr)
\simeq
 W_b\bigl(N(\ibar)^{(1)}\bigr)^{(0)}.
\end{equation}
\end{lem}
\pf
Since the filtration $W$ gives mixed twistor structure,
we have
$W_{h+b}\bigl(N(\ibar)^{(0)}\bigr)^{(1)}
=W_{h+b}\bigl(N(\ibar)\bigr)^{(0)\,(1)}$,
due to Lemma \ref{lem;a12.5.10}.
Then the first isomorphism
is equivalent to the isomorphism
$\Gr^WGr^{W(N_1)}_h\bigl(W_{h+b}(N(\ibar))\bigr)
\simeq
 \Gr^{W(N_1)}_h\Gr^W\bigl(W_{h+b}\bigl(N(\ibar)\bigr)\bigr)$,
which always holds.

Since the filtration $W^{(1)}$ gives mixed twistor structure,
we have
$W_b\bigl(N(\ibar)^{(1)}\bigr)^{(0)}
=W_b\bigl(N(\ibar)^{(1)\,(0)}\bigr)$,
due to Lemma \ref{lem;a12.5.10}.
Thus the second isomorphism is equivalent
to $W_b\bigl(N(\ibar)^{(0)\,(1)}\bigr)\simeq
 W_b\bigl(N(\ibar)^{(1)\,(0)}\bigr)$.
It follows from the equality
$N(\ibar)^{(0)\,(1)}=N(\ibar)^{(1)(0)}$
under the isomorphism
$Gr^{W^{(1)}}\Gr^{W(N_1)}_h\simeq 
 Gr^{W(N_1)^{(0)}}_h\Gr^W$.
Thus we obtain Lemma \ref{lem;a11.10.21}.
\hfill\qed

\vspace{.1in}
From (\ref{eq;a11.10.14}), (\ref{eq;a11.10.15}) and
(\ref{eq;a11.10.16}),
we obtain the following:
\begin{equation}\label{eq;a11.10.18}
 \Gr^{W^{(1)}}\bigl(
 W_{h+b}\bigl(N(\ibar)^{(1)}\bigr)
\cap
 \Gr^{W(N_1)}_h
 \bigr)
=\Gr^{W^{(1)}}\bigl(
 W_b\bigl(N(\ibar)^{(1)}\bigr)\cap \Gr_h^{W(N_1)}
 \bigr)
\end{equation}
Then Lemma \ref{lem;a11.10.20} follows 
from (\ref{eq;a11.10.18}) and Lemma \ref{lem;a11.10.22} below.
\hfill\qed

\vspace{.1in}
\begin{lem} \label{lem;a11.10.23}
We have
$W_{h+b}(N(\ibar))^{(1)}\cap V^{(1)}_h=
 W_{b}(N(\ibar)^{(1)})\cap V^{(1)}_h$.
\end{lem}
\pf
Since $N(\ibar)^{(1)}$ is compatible with the primitive decomposition
of $\Gr^{W(N_1)}_b$,
the filtration $W(N(\ibar)^{(1)})$ is compatible 
with the primitive decomposition of $\Gr^{W(N_1)}_b$.
Then Lemma \ref{lem;a11.10.23} follows from
Lemma \ref{lem;a11.10.20}.
\hfill\qed

\vspace{.1in}

We would like to show the surjectivity of the morphism
\begin{equation}\label{eq;9.13.30}
 \bigcap_{i=1}^l W_{h_i}(\ibar)
\lrarr \bigcap_{i=2}^l W_{h_i}^{(1)}(\ibar)\cap \Gr_{h_i}^{W(N_1)}.
\end{equation}
The morphism (\ref{eq;9.13.30})
induces the morphism on
the associated graded vector bundles
for $W=W(\lbar)$, and the induced morphism is surjective.
Thus the morphism (\ref{eq;9.13.30}) itself is surjective.
Similarly, we can show the surjectivity of the following morphism:
\[
 \Ker N_1^{h_1}\cap
\bigcap_{i=1}^l W_{h_i}(\ibar)
\lrarr
 \bigcap_{i=2}^l W_{h_i}^{(1)}(\ibar)
 \cap P_{h_1}\Gr^{W(\itibar)}_{h_1}.
\]
Thus we obtain the strongly sequential compatibility,
i.e.,
the proof of Lemma \ref{lem;10.26.111} is accomplished.
\hfill\qed

%% file: b9.tex

\subsubsection{Lemmas for mixed twistor structure}

\begin{lem} \label{lem;a12.5.10}
Let $(V,W)$ be a mixed twistor structure.
Let $f:V\lrarr V\otimes\Tate(-1)$ be a nilpotent morphism of mixed twistor
structure.
Let $W(f)$ be the weight filtration of $f$.
We have the induced mixed twistor structure $V^{(0)}:=\Gr^{W}(V)$
and the induced morphism
$f^{(0)}:V^{(0)}\lrarr V^{(0)}\otimes\Tate(-1)$.
We also have the induced filtration $W(f)^{(0)}$.
Then we have $W(f^{(0)})=W(f)^{(0)}$.
\end{lem}
\pf
In the case $f^{h+1}=0$,
we have $W_{-h}(f)=\Image(f^{h})$
and $W_h(f)=\Ker(f^{h})$.
Since $f$ is strict with respect to the filtration $W$,
we obtain
$W_{-h}(f)^{0}=W_{-h}(f^{(0)})$
and $W_h(f)^{(0)}=W_h(f^{(0)})$.
Due to the recursive construction of the weight filtration
(see \cite{d3}),
we obtain $W_k(f)^{(0)}=W_k(f^{(0)})$ for any $k$.
\hfill\qed

\begin{lem}\label{lem;a11.10.22}
Let $(V,W)$ be a MTS.
Let $V_i\subset V$ $(i=1,2)$ be sub MTS.
Assume that $\Gr^W(V_1)=\Gr^W(V_2)$ in $\Gr^W(V)$.
Then we have $V_1=V_2$.
\end{lem}
\pf
We have only to show $W_h\cap V_1=W_h\cap V_2$ for any $h$.

1. We put $h_0:=\min\{h\,|\,\dim W_h\cap V_1\neq 0\}$.
Then $W_{h_0}\cap V_i\simeq \Gr^W_{h_0}(V_i)\subset \Gr_{h_0}^W(V)$.
We put $U:=\Gr^W_{h_0}(V_1)=\Gr^W_{h_0}(V_2)$.
Then $W_{h_0}\cap V_i$ give the isomorphism
$\phi_i:U\lrarr W_{h_0}\cap V_i\subset V$.
Then $(\phi_1-\phi_2)(U)\subset W_{h_0-1}(V)$.
Since $U$ is pure twistor of weight $h_0$,
we obtain $\phi_1-\phi_2=0$,
i.e.,
$W_{h_0}\cap V_1=W_{h_0}\cap V_2$.

2. We assume that the claim holds for $h-1$,
and then we will show the claim for $h$.
We put $K:=W_{h-1}\cap V_1=W_{h-1}\cap V_2$.
It is sub MTS of $V$.
Thus $V/K$ is also a MTS,
and $V_i/K\subset V/K$ $(i=1,2)$ satisfies the condition.
Thus the problem is reduced to the previous case.
\hfill\qed

%% file: a33.tex

\subsubsection{Vanishing cycle theorem due to Kashiwara-Kawai}

Let $(V,W,\vecN,S)$ be a Pol-MTS of type $(n,l)$.
We put $V^{(1)}:=\Image(N_1)$.
Since $N_1:V\lrarr V\otimes\Tate(-1)$ is a morphism of mixed twistors,
$V^{(1)}$ is a sub mixed twistor of $V\otimes\Tate(-1)$.
The filtrations $W^{(1)}$
and the tuple of the nilpotent maps $\vecN^{(1)}$
on $V^{(1)}$ are naturally induced.
Since $V^{(1)}$ is a subbundle of $V\otimes\Tate(-1)$,
we have the naturally induced pairing
$S':V\otimes\sigma(V^{(1)})\lrarr \Tate(-n-1)$.
It is easy to see that $S'$ vanishes
on $\Ker(N_1)\otimes\sigma(V^{(1)})$.
Hence we obtain the induced pairing
$S^{(1)}:V^{(1)}\otimes\sigma(V^{(1)})\lrarr \Tate(-n-1)$.

\begin{lem} \label{lem;9.23.1}
Let $(V,W,\vecN,S)$ be a split Pol-MTS of type $(n,l)$.
Then the induced tuple
$(V^{(1)},W^{(1)},\vecN^{(1)},S^{(1)})$
is a split Pol-MTS of type $(n+1,l)$.
\end{lem}
\pf
Recall that we can pick a torus action
on $(V,W,\vecN,S)$,
and then it is regarded as a nilpotent orbit
in the Hodge theory
(See Corollary \ref{cor;a11.10.25}).
Thus the lemma is a consequence of the vanishing cycle theorem
due to Kashiwara-Kawai
(Theorem 2.15 in \cite{k3}).
\hfill\qed

\vspace{.1in}
For a polarized mixed twistor structure $(V,W,\vecN,S)$,
we put $V^{(0)}:=\Gr^{W}(V)$.
Then we have the induced split Pol-MTS 
$(V^{(0)},W^{(0)},\vecN^{(0)},S^{(0)})$
of type $(n,l)$.
By applying the construction above,
we obtain $(V^{(0)\,(1)},W^{(0)\,(1)},\vecN^{(0)\,(1)},S^{(0)\,(1)})$.

On the other hand,
we obtain the tuple $(V^{(1)},W^{(1)},\vecN^{(1)},S^{(1)})$.
It is easy to check that the tuple
is a $\Psi$-Pol MTS of type $(n,l)$.
Hence we obtain the induced tuple
$(V^{(1)\,(0)},W^{(1)\,(0)},\vecN^{(1)\,(0)},S^{(1)\,(0)})$.

\begin{lem} \label{lem;9.23.2}
The tuples
$(V^{(0)\,(1)},W^{(0)\,(1)},\vecN^{(0)\,(1)},S^{(0)\,(1)})$
and
$(V^{(1)\,(0)},W^{(1)\,(0)},\vecN^{(1)\,(0)},S^{(1)\,(0)})$
are naturally isomorphic.
\end{lem}
\pf
It is clear from our constructions.
\hfill\qed

\vspace{.1in}

\begin{prop} \label{prop;9.23.6}
Let $(V,W,\vecN,S)$ be a Pol-MTS of type $(n,l)$.
We put $V^{(1)}:=\Image(N_1)$.
Then the naturally induced tuple
$(V^{(1)},W^{(1)},\vecN^{(1)},S^{(1)})$
is a Pol-MTS of type $(n,l)$.
\end{prop}
\pf
It follows from Lemma \ref{lem;9.23.1} and Lemma \ref{lem;9.23.2}.
\hfill\qed

%% file: b10.tex

\subsubsection{Kashiwara's Lemma}
\label{subsubsection;9.13.50}

Let $V$ be a vector bundle over $\proj^1$
and $\vecN$ be a tuple of
nilpotent maps $N_i:V\lrarr V\otimes \Tate(-1)$ $(i=1,\ldots,l)$.
Let $N$ be a variable,
and we put as follows:
\[
 V[N]:=
 \bigoplus_{i=0}^{\infty}
 V\otimes \Tate(i)\cdot N^i.
\]
We have the natural mixed twistor structure $W$
on the infinite dimensional vector bundle $V[N]$.
We have the natural morphism of mixed twistor structures:
\[
 N:V[N]\lrarr V[N]\otimes \Tate(-1),
\quad
 u\cdot N^i\longmapsto u\cdot N^{i+1}.
\]
For a subset $I\subset\lbar$,
we have the morphism:
\[
 \prod_{i\in I}(N-N_i):
 V[N]\otimes\Tate(n)\lrarr V[N].
\]
Then we put as follows:
\begin{equation}\label{eq;b12.6.80}
 \nbigv_{I}(V,\vecN):=V[N]\Big/\prod_{i\in I}(N-N_i)
\simeq
 \bigoplus_{i=0}^{|I|-1}V\otimes \Tate(i)\cdot N^i.
\end{equation}
In the case $I=\underline{l_1}$, we use the notation
$\nbigv_{l_1}(V,\vecN)$ instead of $\nbigv_{\underline{l_1}}(V,\vecN)$.
We often omit to denote $\vecN$.

Let $(V,W)$ be a mixed twistor structure
and $\vecN$ be a tuple of morphisms $N_i:V\lrarr V\otimes\Tate(-1)$
of mixed twistors.
We have the induced filtration $W$
and the induced morphisms $N,N_i$ on $(\nbigv_{I},W)$.
The following lemma is clear from the construction.
\begin{lem}
The morphisms $N$ and $N_i$ are morphisms of 
mixed twistor structures.
\hfill\qed
\end{lem}

\begin{cor}
 The conjugacy classes of $aN+\sum b_i N_i$ are independent of
$\lambda\in\proj^1$.
\hfill\qed
\end{cor}

We put $\vecnbign:=(N,N_1,\ldots,N_l)$.
Thus we obtain the mixed twistor structure
$(\nbigv_I,W)$ and the tuple of nilpotent maps $\vecnbign$.

Let $(V,W,\vecN,S)$ be a Pol-MTS$(n,l)$.
We have the induced object
$(\nbigv_{I},\nbigw,\vecnbign)$.
We have the induced pairing
$V(i)\otimes\sigma\bigl(V(j)\bigr)\lrarr \Tate(-n+i+j)$,
which we also denote by $S$.
We have the induced morphism
$S':V[N,N^{-1}]\otimes\sigma (V[N,N^{-1}])
\lrarr \Tate(-n)[N,N^{-1}]$, defined as follows:
\[
 S\bigl(u\cdot N^i,\sigma(v\cdot N^j)\bigr)
:=(-1)^j\cdot S\bigl(u,\sigma(v)\bigr)\cdot N^{i+j}.
\]
Since $N_i$ are nilpotent,
we have the morphism:
\[
  \prod_{i\in I}(N-N_i)^{-1}:
 V[N]\lrarr V[N]\otimes\Tate(|I|).
\]
Then we have the following element 
of $\Tate(-n+|I|-k)[N,N^{-1}]$
for any $f_i\in V\otimes\Tate(i)$ and $g_j\in V\otimes\Tate(j)$,
and for any $k\in\seisuu$:
\[
 S'\Bigl(
 \prod_{i\in I}(N-N_i)^{-1}f_i\otimes N^{i+k},
 \,\,
 \sigma\bigl(g_j\cdot N_j\bigr)
 \Bigr)
=\sum_{\vecn\in\seisuu^I}
 (-1)^j\cdot N^{-|I|+i+j+k-|\vecn|}
 S\Bigl(\prod_{i\in I}N_a^{n_a}f_i,\sigma(g_j)\Bigr)
\in \Tate(-n+|I|-k)[N].
\]

By taking the residue,
i.e.,
the $N^{-1}\cdot \Tate(-n+|I|-k-1)$-component,
we obtain the morphism:
$ \nbigv_{I}(S):
 V[N,N^{-1}]\otimes \sigma(V[N,N^{-1}])
\lrarr \Tate(-n+|I|-k-1)$:
\begin{equation}\label{eq;a11.10.30}
 \nbigv_{I}(S)\bigl(f(N),\sigma\bigl(g(N)\bigr)\bigr)
=\Res_{N=0}S'\Bigl(
 \tilde{N}^{-1}_I\cdot f(N),\sigma\bigl(g(N)\bigr)
 \Bigr).
\end{equation}
Here we put $\tilde{N}_I:=\prod_{i\in I}(N-N_i)$.

\begin{lem}
We have the following formula for
$f=\sum f_i\cdot N^i$ and $\sum g_j\cdot N^j$
and for any $k\in\seisuu$:
\[
 \nbigv_{I}(S)\bigl(
 f(N)\cdot N^k,\,\,\sigma\bigl(g(N)\bigr)
 \bigr)
=\sum_{i,j}
 \sum_{-|I|-|\vecn|+i+j+k=-1}
 (-1)^j\cdot
 S\Bigl(
 \prod_{a\in I}N_a^{n_a}\cdot f_i,\,\sigma(g_j)
 \Bigr)\in\Tate(-n+|I|-1-k).
\]
\end{lem}
\pf
It can be checked by a direct calculation.
\hfill\qed

\vspace{.1in}

Then we obtain
$\nbigv_{I}(S):\nbigv_{I}(V)\otimes\sigma(\nbigv_{I}(V))
\lrarr \Tate(-n+|I|-1)$.
Thus we obtain the tuple
$\bigl(\nbigv_{I}(V),\nbigw,\nbigvecn,\nbigv_{I}(V)\bigr)$.

The signature of the construction in {\rm\cite{saito2}}
looks slightly different from that in the above construction
at a sight.
Let us see that they are same, in fact.
Recall the construction in \cite{saito2}.
To distinguish, we use the variable $s$.
We also denote the given nilpotent maps of $V$
by $s_i$.
Note the relation $s_i=-N_i$.
(See Remark {\rm\ref{rem;10.17.5}}, for example.)

We obtain the vector bundle
$\nbigv_I(V)^{(1)}:=V[s]/\prod_{i\in I}(s-s_i)$.
The extended pairing $S^{(1)}:V[s,s^{-1}]\otimes\sigma(V)[s,s^{-1}]
 \lrarr \Tate(-n)[s,s^{-1}]$
given in \cite{saito2}
is as follows:
\[
 S^{(1)}(u\cdot s^i,v\cdot s^j)=(-1)^i\cdot S(u,v)\cdot s^{i+j}.
\]
Then we obtain the pairing $\tilde{S}^{(1)}$,
given as follows:
\[
 \tilde{S}^{(1)}(s^i\cdot u,s^j\cdot v):=
 \Res_{s=0}S^{(1)}\bigl(u\!\cdot\! s^i,
 \,\,\tilde{s}^{-1}_I\!\cdot\! v\!\cdot\! s^j\bigr).
\]
Here $\tilde{s}_I$ denote $\prod_{i\in I}(s-s_i)$.

Since we have the relation $s_i=-N_i$,
the correspondence $N\longmapsto -s$ induces the
isomorphism $V[N]\simeq V[s]$,
and $\nbigv_I(V)\simeq \nbigv_I(V)^{(1)}$.
\begin{lem}\label{lem;a11.10.42}
Under the isomorphism,
we have $\tilde{S}(N^iu,N^jv)=\tilde{S}^{(1)}(s^iu,s^jv)$.
\end{lem}
\pf
We put
$\nbigs:=\bigl\{\vecl\in\seisuu^I\,\big|\,
 -|I|-|\vecl|+i+j=-1
 \bigr\}$.
By a direct calculation,
we have the following:
\begin{equation}\label{eq;a11.10.40}
 \tilde{S}(N^iu,N^jv)=
 \sum_{\vecl\in\nbigs}
 (-1)^j\cdot S\Bigl(\prod N_k^{l_k}u,v\Bigr).
\end{equation}
We also have the following:
\begin{equation}\label{eq;a11.10.41}
 \tilde{S}^{(1)}(s^iu,s^jv)=
 \sum_{\vecl\in\nbigs}
 (-1)^i\cdot S\Bigl(u,\prod s_k^{l_k}v\Bigr).
\end{equation}
Let us substitute $s=-N$ and $s_i=-N_i$ in
(\ref{eq;a11.10.41}),
then we obtain the following:
\[
 (-1)^{i+j}\tilde{S}^{(1)}\bigl(
 N^iu,N^jv
 \bigr)
=\sum_{\vecl\in\nbigs}
 (-1)^i\cdot S\Bigl(u,\prod (-N_k)^{l_k}v\Bigr)
=\sum_{\vecl\in\nbigs}
 (-1)^i\cdot S\Bigl(\prod N_k^{l_k}u,v \Bigr).
\]
Since it is same as (\ref{eq;a11.10.40}),
Lemma \ref{lem;a11.10.42} is shown.
\hfill\qed

\begin{lem}
Let $(V,W,\vecN,S)$ be a split Pol-MTS of type $(n,ln)$.
Then the induced tuple
$\bigl(\nbigv_{I}(V),\nbigw,\nbigvecn,\nbigv_{I}(V)\bigr)$
is a split Pol-MTS of type $(n-|I|+1,l+1)$.
\end{lem}
\pf
Clearly we have a splitting.
We put $(V',W',\vecN',S')=(V,W,\vecN,S)\oplus\sigma(V,W,\vecN,S)$,
which gives an $\real$-nilpotent orbit.
Then
$\bigl(\nbigv_{l_1}(V'),\nbigw,\nbigvecn,\nbigv_{l_1}(V')\bigr)$
is an $\real$-nilpotent orbit.
Since we have the following:
\[
 \bigl(\nbigv_{l_1}(V'),\nbigw,\nbigvecn,\nbigv_{l_1}(S)\bigr)
=\bigl(\nbigv_{l_1}(V),\nbigw,\nbigvecn,\nbigv_{l_1}(S)
  \bigr)
\oplus
 \sigma\bigl(\nbigv_{l_1}(V),\nbigw,\nbigvecn,\nbigv_{l_1}(S)\bigr).
\]
Thus it follows from the original Kashiwara's lemma
(Proposition 3.19 \cite{saito2} and Appendix \cite{saito2}).
Or, it is not difficult to apply Kashiwara's argument
in Appendix in \cite{saito2} to the nilpotent twistor orbits.
\hfill\qed

\begin{cor} \label{cor;9.23.7}
Let $(V,W,\vecN,S)$ be a Pol-MTS of type $(n,l)$.
Then $\bigl(\nbigv_{I}(V),\nbigw,\nbigvecn,\nbigv_{I}(S)\bigr)$
is a Pol-MTS of type $(n-|I|+1,l+1)$.
\end{cor}
\pf
It is easy to check that
the tuple
$\bigl(\nbigv_{I}(V),\nbigw,\nbigvecn,\nbigv_{I}(S)\bigr)$
satisfies the first three conditions.
Then the tuple
$\bigl(\nbigv_I(V)^{(0)}_{I},\nbigw^{(0)},\nbigvecn^{(0)},
  \nbigv_I(S)^{(0)}
 \bigr)$
is naturally isomorphic to
the tuple obtained from 
$\bigl(\nbigv_I(V^{(0)}),\nbigw,\nbigvecn^{(0)},\nbigv_I(S^{(0)})
 \bigr)$.
Thus we are done.
\hfill\qed

%% file: b10.1.tex

\subsubsection{Preliminary lemma}
\label{subsubsection;a11.10.50}

In the subsubsections
\ref{subsubsection;a11.10.50}--\ref{subsubsection;a11.10.51},
we recall the argument given by Saito in the page 302 in \cite{saito2}.
See also the original argument.

Let $(V,W,\vecN,S)$ be an equivariant $\real$-nilpotent orbit 
of of $(n,l)$.
Let $(\hat{V},\hat{W},\hat{\vecN},\hat{S})$ be
an equivariant $\real$-nilpotent orbit of type $(n-1,l)$.
Let $f:V\lrarr \hat{V}\otimes\Tate(-1)$
and $g:\hat{V}\lrarr V$
be morphisms of MTS,
such that
\begin{itemize}
\item
 $g\circ f=N_2: V\lrarr V\otimes \Tate(-1)$.
\item
 $f\circ g=\hat{N}_2:\hat{V}\lrarr\hat{V}\otimes\Tate(-1)$.
\item
 $\hat{S}(f\otimes\id)=S(\id\otimes g)$.
\end{itemize}

Then we have the induced morphisms
\[
 \begin{array}{l}
 f^{(1)}:
 P_h\Gr^{W(N_1)}_h(V)
 \lrarr
 P_h\Gr_h^{W(\hat{N}_1)}(\hat{V})\otimes\Tate(-1),\\
 \mbox{{}}\\
 g^{(1)}:
 P_h\Gr^{W(\hat{N}_1)}_h(\hat{V})
 \lrarr
 P_h\Gr^{W(N_1)}_h(V).
 \end{array}
\]

\begin{lem} \label{lem;9.13.35}
In the case $l=2$,
we have the following decomposition:
\[
 P_h\Gr^{W(\hat{N}_1)}(\hat{V})
=
 \Ker g^{(1)}
\oplus
\Bigl(
\Image(f^{(1)})\otimes\Tate(1)
\Bigr)
\]
\end{lem}
\pf
It is easy to see that
$f^{(1)}$ and $g^{(1)}$ are morphisms of mixed twistor structures,
and that
we have $f^{(1)}\circ g^{(1)}=\hat{N}_2^{(1)}$
and $g^{(1)}\circ f^{(1)}=N_2^{(1)}$.
We also have
$\hat{S}^{(1)}(f^{(1)}\otimes 1)=\hat{S}(1\otimes g^{(1)})$.

We have the induced morphisms
\[
 \begin{array}{l}
 f^{(2)}:
 \Gr^W_a P_h\Gr^{W(N_1)}_{h}(V)
 \lrarr
 \Gr^{W}_{a-2}\bigl(P_h\Gr_h^{W(\hat{N}_1)}(\hat{V})\bigr)
  \otimes\Tate(-1),\\
 \mbox{{}}\\
 g^{(2)}:
 \Gr^{W}_aP_h\Gr^{W(\hat{N}_1)}_h(\hat{V})
 \lrarr
 \Gr^W_aP_h\Gr^{W(N_1)}_h(V).
 \end{array}
\]
In all, we have the following:
\[
 \begin{array}{l}
 f^{(2)}:
 \Gr^W P_h\Gr^{W(N_1)}_{h}(V)
 \lrarr
 \Gr^{W}(P_h\Gr_h^{W(\hat{N}_1)}(\hat{V}))\otimes\Tate(-1),\\
 \mbox{{}}\\
 g^{(2)}:
 \Gr^{W} P_h\Gr^{W(\hat{N}_1)}_h(\hat{V})
 \lrarr
 \Gr^W P_h\Gr^{W(N_1)}_h(V).
 \end{array}
\]
We have only to show
$\Image(f^{(2)})\otimes\Tate(1)\oplus
 \ker g^{(2)}=
 \Gr^WP_h\Gr_h^{W(\hat{N}_1)}(\hat{V})$.

Note the following equalities:
\[
\begin{array}{l}
 \Gr^W_a P_h \Gr^{W(N_1)}_h(V)
=\Gr^{W(N_2^{(1)})}_{a-n-h}
 P_h\Gr^{W(N_1)}_h(V),\\
 \mbox{{}}\\
 \Gr^W_a P_h\Gr^{W(\hat{N}_1)}_h(\hat{V})
=\Gr^{W(\hat{N}_2^{(1)})}_{a-n+1-h}
 P_h\Gr^{W(\hat{N}_1)}_h(V).
\end{array}
\]
Thus $f^{(2)}$ and $g^{(2)}$ can be regarded as the morphisms:
\[
 \begin{array}{l}
 f^{(2)}:
 \Gr^{W(N_2^{(1)})}_a P_h\Gr_h^{W(N_1)}(V)
\lrarr
 \Gr^{W(\hat{N}_2)}_{a-1}(P_h\Gr^{W(\hat{N}_2)}_h(\hat{V}) )\otimes\Tate(-1)
 \\
 \mbox{{}}\\
 g^{(2)}:
 \Gr^{W(\hat{N}_2^{(1)})}_aP_h\Gr^{W(N_1)}_h(\hat{V})
 \lrarr
 \Gr^{W(N_2)}_{a-1}(P_h\Gr_h^{W(N_1)}(V)).
 \end{array}
\]
We put $n+h+1=n'$ and as follows:
\[
 \begin{array}{ll}
 H^a:=\Gr_a^{W(N_2^{(1)})}P_h\Gr_h^{W(N_1)}(V)
 & \mbox{ weight }a+n'-1\\
 \mbox{{}}\\
 \hat{H}^a:=\Gr^{W(N_2^{(1)})}P_h\Gr_h^{W(\hat{N}_1)}(\hat{V})
 \otimes\Tate(-1) & \mbox{ weight }a+n'
 \end{array}
\]
Then we obtain the following morphisms:
\[
\begin{array}{l}
 f^{(2)}:H^{a}\lrarr \hat{H}^{a-1},\\
\mbox{{}}\\
 g^{(2)}:\hat{H}^{a}\lrarr H^{a-1}\otimes \Tate(-1),
\end{array}
\]
We also have the pairings:
\[
\begin{array}{l}
 S^{(2)}:H^a\otimes\sigma H^{-a}\lrarr \Tate(-n'+1),\\
 \mbox{{}}\\
 \hat{S}^{(2)}:\hat{H}^a\otimes\sigma\hat{H}^{-a}
\lrarr \Tate(-n').
\end{array}
\]
We have $\hat{S}^{(2)}(f^{(2)}\otimes\id)=S^{(2)}(1\otimes g^{(2)})$,
and 
we have the polarizations
$S^{(2)}(N^{(2)\,a}\otimes\id)$
and $\hat{S}^{(2)}(\hat{N}^{(2)\,a}\otimes\id)$.
Then the problem is reduced to
Lemma 5.2.15 of \cite{saito1}.
\hfill\qed

\begin{prop} \label{prop;9.23.5}
We have the following decomposition
for any $l$:
\[
 P_h\Gr^{W(\hat{N}_1)}(\hat{V})
=
 \Ker g^{(1)}
\oplus
\Bigl(
\Image(f^{(1)})\otimes\Tate(1)
\Bigr)
\]
\end{prop}
\pf
We put $X=\Delta^l$, $D_i=\{z_i=0\}$ and $D=\bigcup_{i=1}^l D_i$.
Let $(E,\delbar_E,h,\theta)$
and $(\hat{E},\delbar_{\hat{E}},\hat{h},\hat{\theta})$
be the harmonic bundles over $X-D$,
corresponding to
$(V,W,\vecN,S)$ and $(\hat{V},\hat{W},\hat{\vecN},\hat{S})$
respectively.
Since we discuss in the category of mixed twistor structures,
such decomposition holds at $\lambda=1$
if and only if
such decomposition holds over $\proj^1$.
Let $(\nbige,\DD)$ and $(\hat{\nbige},\hat{\DD})$
be the associated deformed holomorphic bundles,
and $\prolong{\nbige}$ and $\prolong{\hat{\nbige}}$
be the canonical prolongment.

The morphisms $f$ and $g$
induce the morphisms $F:\prolong{\nbige}\lrarr \prolong{\hat{\nbige}}$
and $G:\prolong{\hat{\nbige}}\lrarr\prolong{\nbige}$,
which are flat with respect to the $\lambda$-connections
$\DD$ and $\hat{\DD}$.

Let us pick a point
$Q\in D_1\cap D_2-\bigcup_{j>2}\bigl(D_i\cap D_1\cap D_2\bigr)$.
By using the normalizing frame,
we obtain the isomorphism:
\[
 \bigl(\prolong{\nbige}_{|\,(1,O)},
 \Res_{1}(\DD),\Res_{2}(\DD)\bigr)
\simeq
 \bigl(\prolong{\hat{\nbige}}_{|\,(1,Q)},
 \Res_{1}(\DD),\Res_{2}(\DD)\bigr)
\]
We have a similar isomorphism for
$\hat{\nbige}$.
Hence we have only to show the following decomposition:
\begin{equation} \label{eq;9.23.4}
 \Gr^{W(N)}(\prolong{\hat{\nbige}}_{|(1,Q)})
=\Image(F_{|(1,Q)})\oplus \Ker(G_{|(1,Q)}).
\end{equation}

Let $\pi:X\lrarr D_1\cap D_{2}$ denote the projection
onto the first two components.
We obtain the surface $\pi^{-1}(Q)\simeq \Delta^2\subset X$.
The restrictions of
$(E_i,\delbar_{E_i},h_i,\theta_i)$
to $\pi^{-1}(Q)$ are nilpotent orbits
(Lemma \ref{lem;9.23.3}, for example).
Thus we obtain the decomposition (\ref{eq;9.23.4})
from Lemma \ref{lem;9.13.35}.
Thus the proof of Proposition \ref{prop;9.23.5}
is accomplished.
\hfill\qed

%% file: a33.1.tex

\subsubsection{The preliminary decomposition}

\label{subsubsection;9.23.11}

Let us consider a pair $\vecI=(J,K)$ 
of subsets $J$ and $K$ of $\lbar$,
such that $J\cap K=\emptyset$ and $J\neq\emptyset$.
For such pair, we put $|\vecI|:=|J|+|K|$.
For two pairs $\vecI_i=(J_i,K_i)$ $(i=1,2)$,
we mean 
$J_1\subset J_2$ and $K_1\subset K_2$
by $\vecI_1\subset \vecI_2$.
In the case $\vecI_1\subset\vecI_2$,
we put $\vecI_1-\vecI_2:=(J_2-J_1)\cup (K_2-K_1)$.

Let $(V,W,\vecN,S)$ be a Pol-MTS of type $(n,l)$.
We put $N_K:=\prod_{k\in K}N_k$.
We have the naturally induced filtration $W$
the tuple of the nilpotent maps $\vecN$,
and the pairing $S$ on $\Image(N_K)$.
Due to Proposition \ref{prop;9.23.6},
the tuple
$(\Image(N_K),W,\vecN,S)$ is a Pol-MTS of type $(n+|K|,l)$.
By the construction given in the subsubsection
\ref{subsubsection;9.13.50},
we obtain the Pol-MTS $\nbigv_J(\Image(N_K))$
of type $(n+|K|-|J|+1,l+1)$
(See Corollary \ref{cor;9.23.7}).

Let us pick an element $i$ of $K$,
and we put $K':=K-\{i\}$.
The inclusion $\Image(N_{K})\subset \Image(N_{K'})\otimes\Tate(-1)$
is denoted by $\var_i$.
The morphism $N_i$ induces the morphism
$\Image(N_{K'})\lrarr \Image(N_K)$,
which is denoted by $\can_i$.
The morphism $\var_i$ and $\can_i$
induce the morphisms
$\var_i:
 \nbigv_J(\Image(N_{K}))\lrarr
 \nbigv_J(\Image(N_{K'}))\otimes\Tate(-1)$
and
$\can_i:
 \nbigv_J(\Image(N_{K'}))\lrarr \nbigv_J(\Image(N_{K}))$.

Let us pick an element $i\in J$,
and we put $J':=J-\{i\}$.
The identity of $\Image(N_K)[N]$
induces the morphism
$\var_i:\nbigv_J(\Image(N_K))\lrarr\nbigv_{J'}(\Image(N_{K}))$.
The morphism $(N-N_i):\Image(N_K)[N]\lrarr\Image(N_K)[N]\otimes\Tate(-1)$
induces the morphism
$\can_i:\nbigv_{J'}(\Image(N_{K}))\otimes\Tate(1)
\lrarr\nbigv_J(\Image(N_K))$.

We introduce the following notation.
For a pair $\vecI=(J,K)$ as above,
we put
$\nbigu_{\vecI}
 :=\nbigv_J(\Image(N_K))\otimes\Tate\bigl(|K|\bigr)$.
Then we obtain the morphisms
$\can_i:\nbigu_{\vecI'}\lrarr \nbigu_{\vecI}$
and $\var_i:\nbigu_{\vecI}\lrarr\nbigu_{\vecI'}$
for $\vecI'\subset\vecI$ such that $\vecI-\vecI'=\{i\}$.
They induce the morphisms
$\can_i:
 P_h\Gr^{W(N)}_h\nbigu_{\vecI'}\lrarr P_h\Gr^{W(N)}_h\nbigu_{\vecI}$
and
$\var_i:
 P_h\Gr^{W(N)}_h\nbigu_{\vecI}\lrarr P_h\Gr^{W(N)}_h\nbigu_{\vecI'}$.


\begin{prop} \label{prop;9.23.10}
Let $\vecI=(J,K)$ be a pair of subsets
$J$ and $K$ of $\lbar$
such that $J\cap K=\emptyset$ and $J\neq\emptyset$.
Let $\vecI'$ be a pair of subsets
such that $\vecI'\subset\vecI$ and
$\vecI-\vecI'=\{i\}$.
Then we have the following decomposition:
\[
 P_h\Gr^{W(N)}_h\bigl(
 \nbigu_{\vecI}
 \bigr)
=\Ker(\var_i)\oplus\Image\bigl(\can_i\otimes\Tate(1)\bigr).
\]
\end{prop}
\pf
It follows from Proposition \ref{prop;9.23.5}.
\hfill\qed

%% file: b10.4.tex

\subsubsection{The decomposition}

\label{subsubsection;a11.10.51}

For pairs $\vecI'\subset \vecI$,
we have the naturally defined morphisms
$\can(\vecI',\vecI):\nbigu_{\vecI'}\lrarr\nbigu_{\vecI}$
and
$\var(\vecI',\vecI):\nbigu_{\vecI}\lrarr\nbigu_{\vecI'}$:
\[
 \can(\vecI',\vecI):=
 \prod_{i\in\vecI-\vecI'}\can_i,
\quad\quad
 \var(\vecI',\vecI):=
 \prod_{i\in \vecI-\vecI'}\var_i.
\]

Let $\vecI_i=(J_i,K_i)$ $(i=0,1,2,3)$ be pairs of subsets of $\lbar$.
We assume 
$\vecI_1\cap\vecI_2:=(J_1\cap J_2,K_1\cap K_2)=\vecI_0$,
and $\vecI_1\cup\vecI_2:=(J_1\cup J_2,K_1\cup K_2)=\vecI_3$.
\begin{lem} \label{lem;a11.11.1}
The following morphisms are commutative:
\[
 \begin{CD}
\nbigu_{\vecI_0} @>{\can(\vecI_0,\vecI_1)}>> \nbigu_{\vecI_1} \\
 @V{\can(\vecI_0,\vecI_2)}VV @V{\can(\vecI_1,\vecI_3)}VV \\
 \nbigu_{\vecI_2} @>{\can(\vecI_2,\vecI_3)}>> \nbigu_{\vecI_3},
\end{CD}
\quad\quad\quad\quad\quad
 \begin{CD}
\nbigu_{\vecI_0} @>{\var(\vecI_0,\vecI_1)}>> \nbigu_{\vecI_1} \\
 @V{\can(\vecI_0,\vecI_2)}VV @V{\can(\vecI_1,\vecI_3)}VV \\
 \nbigu_{\vecI_2} @>{\var(\vecI_2,\vecI_3)}>> \nbigu_{\vecI_3},
\end{CD}
\quad\quad\quad\quad\quad
 \begin{CD}
\nbigu_{\vecI_0} @>{\var(\vecI_0,\vecI_1)}>> \nbigu_{\vecI_1} \\
 @V{\var(\vecI_0,\vecI_2)}VV @V{\var(\vecI_1,\vecI_3)}VV \\
 \nbigu_{\vecI_2} @>{\var(\vecI_2,\vecI_3)}>> \nbigu_{\vecI_3}
\end{CD}
\]
\end{lem}
\pf
It can be directly checked from the definition.
\hfill\qed

\vspace{.1in}

For $\vecI'\subset \vecI$,
we obtain the morphisms
$\can(\vecI',\vecI):
 P_h\Gr^{W(N)}_h(\nbigu_{\vecI'})
 \lrarr P_h\Gr^{W(N)}_h(\nbigu_{\vecI})$
and
$\var(\vecI',\vecI):
 P_h\Gr^{W(N)}_h(\nbigu_{\vecI})
 \lrarr P_h\Gr^{W(N)}_h(\nbigu_{\vecI'})$.

\begin{prop} \label{prop;a11.11.4}
For $\vecI'\subset \vecI$,
we have the following decomposition:
\begin{equation}\label{eq;a11.11.3}
  P_h\Gr^{W(N)}_h\bigl(
 \nbigu_{\vecI}
 \bigr)
=\Ker\bigl(\var(\vecI',\vecI)\bigr)
 \oplus
 \Image\bigl(\can(\vecI',\vecI)\otimes\Tate(|\vecI|-|\vecI'|)\bigr).
\end{equation}
\end{prop}
\pf
We use an induction on $|\vecI|-|\vecI'|$.
The case $|\vecI|-|\vecI'|=1$ is shown in Proposition \ref{prop;9.23.10}.
We assume that the claim holds
in the case $|\vecI|-|\vecI'|<a$,
and we will prove that the claim also holds
in the case $|\vecI|-|\vecI'|=a$.

We take $\vecI_i$ $(i=0,1,2,3)$ as in Lemma \ref{lem;a11.11.1}
satisfying the following:
\[
  \vecI_0=\vecI',
\quad\quad
 \vecI_3=\vecI,
\quad\quad
 \vecI_0\subsetneq \vecI_i\subsetneq \vecI_3,\,\,(i=1,2).
\]
Due to the hypothesis of the induction,
we have the following decompositions:
\begin{equation}\label{eq;a11.11.2}
 P_h\Gr^{W(N)}_h\bigl(\nbigu_{\vecI_3}\bigr)
=\Image\bigl(\can(\vecI_2,\vecI_3)\bigr)
 \oplus
 \Ker\bigl(\var(\vecI_2,\vecI_3)\bigr),
\quad
 P_h\Gr^{W(N)}_h\bigl(\nbigu_{\vecI_1}\bigr)
=\Image\bigl(\can(\vecI_0,\vecI_1)\bigr)
 \oplus
 \Ker\bigl(\var(\vecI_0,\vecI_1)\bigr).
\end{equation}
Here we have omitted to denote the Tate twist, for simplicity.

Due to the commutativity in Lemma \ref{lem;a11.11.1},
the morphism $\can(\vecI_1,\vecI_3)$ and $\var(\vecI_1,\vecI_3)$
are compatible with the decompositions in (\ref{eq;a11.11.2}).
Then it is easy to derive the decomposition
(\ref{eq;a11.11.3}) for $\vecI=\vecI_3$ and $\vecI'=\vecI_0$.
\hfill\qed

\vspace{.1in}

We put as follows:
\[
 \nbigc_{\vecI}:=
 \bigcap_{\vecI'\subsetneq\vecI}
 \Ker\bigl(\var(\vecI',\vecI)\bigr)
\subset
 P_h\Gr^{W(N)}_h\bigl(\nbigu_{\vecI}\bigr).
\]

\[
 \nbigc_{\vecI,\vecI'}:=
 \can(\vecI',\vecI)\bigl(
 \nbigc_{\vecI'}
 \bigr)
 \subset P_h\Gr^{W(N)}_h\bigl(\nbigc_{\vecI}\bigr).
\]

\begin{lem}
We have the following:
\[
 \nbigc_{\vecI}
=\bigcap_{\substack{\vecI'\subset\vecI,\\|\vecI|-|\vecI'|=1}}
  \Ker\bigl(\var(\vecI',\vecI)\bigr)
\subset
 P_h\Gr^{W(N)}_h\bigl(\nbigu_{\vecI}\bigr).
\]
\end{lem}
\pf
It can be shown by an argument similar to the proof of
Proposition \ref{prop;a11.11.4}.
\hfill\qed

\begin{lem}
Let us consider pairs $\vecI_i\subset \vecI$ $(i=1,2)$.
Assume $\vecI_1\not\subset \vecI_2$.
Then we have $\var(\vecI_2,\vecI)\bigl(\nbigc_{\vecI_1,\vecI}\bigr)=0$.
\end{lem}
\pf
It immediately follows from the commutativity
(Lemma \ref{lem;a11.11.1}) as follows:
\[
 \var(\vecI_2,\vecI)\bigl(\nbigc_{\vecI_1,\vecI}\bigr)
=\var(\vecI_2,\vecI)\circ \can(\vecI_1,\vecI)\bigl(\nbigc_{\vecI_1}\bigr)
=\can(\vecI_1\cap \vecI_2,\vecI_2)\circ
 \var(\vecI_1\cap \vecI_2,\vecI_1)
 \bigl(\nbigc_{\vecI_1}\bigr)
 =0.
\]
\hfill\qed

\begin{lem} \label{lem;9.13.42}
Let us consider sub-pairs $\vecI_1\subset \vecI_2\subset\vecI$.
Then the restriction of the morphism
$\var(\vecI_2,\vecI)_{|\nbigc_{\vecI,\vecI_1}}$
is injective.
\end{lem}
\pf
It follows from
$\nbigc_{\vecI,\vecI_1}\subset \Image \bigl(\can(\vecI_2,\vecI)\bigr)$
and the decomposition in Proposition \ref{prop;a11.11.4}.
\hfill\qed

\begin{lem}\label{lem;a11.11.5}
The naturally defined morphism
$\phi:\bigoplus_{\vecI'\subset \vecI}
 \nbigc_{\vecI,\vecI'}\lrarr P_h\Gr^{W(N)}_h\nbigu_{\vecI}$ is onto.
\end{lem}
\pf
It follows from the surjectivity
of $\sum \Image\bigl(\can(\vecI',\vecI)\bigr)
 \oplus\nbigc_{\vecI}\lrarr P_h\Gr^{W(N)}_h(\nbigu_{\vecI})$.
\hfill\qed

\begin{prop} \label{prop;9.13.41}
We have the following decomposition:
\[
 P_h\Gr^{W(N)}_h(\nbigu_{\vecI})
\simeq
 \bigoplus_{\vecI'\subset \vecI}\nbigc_{\vecI,\vecI'}.
\]
\end{prop}
\pf
We have only to show the injectivity of $\phi$
in Lemma \ref{lem;a11.11.5}.
Assume that $\sum v_{\vecI'}=0$,
where $v_{\vecI}'$ are elements of $\nbigc_{\vecI,\vecI'}$.
We put $\nbigs=\{\vecI'\,|\,v_{\vecI'}\neq 0\}$.
Assume that $\nbigs\neq \emptyset$,
and we will derive a contradiction.
Let $\vecI_0'$ be the minimal element.
Then we obtain the following:
\[
 0=\var(\vecI_0',\vecI)\Bigl(\sum v_{\vecI'}\Bigr)
=\var(\vecI_0',\vecI)(v_{\vecI_0'})
\]
It implies $v_{\vecI_0'}=0$ due to Lemma \ref{lem;9.13.42}.
But it contradicts our choice of $\vecI_0'$.
Thus we obtain Proposition \ref{prop;9.13.41}.
\hfill\qed

\begin{lem} \label{lem;9.26.20}
For a pair $\vecI=(J,K)$,
$\nbigc_{\vecI}$ is pure twistor of weight $n-|J|+1+h$.
The pairing of $P_h\Gr^{W(N)}_h(\nbigu_{\vecI})$
induce the polarization on $\nbigc_{\vecI}$.
\end{lem}
\pf
$\nbigc_I$ is a direct summand of $P_h\Gr^{W(N)}_h(\nbigu_{\vecI})$.
On $\nbigc_I$, the nilpotent maps $\nbign_i$ $(i\in \vecI)$ are $0$,
by definition of $\nbigc_I$.
Thus the weight filtration $\nbigw(I)$ induced by
$\sum_{i\in I}\nbign_i$ is trivial on $\nbigc_I$.
Recall that the filtration
$\nbigw(I)$ is same as the induced mixed twistor structure
of $P_h\Gr^{W(N)}_h\nbigu_{\vecI}$,
up to shift of the degree.
Thus we obtain the first claim.
The second claim can be shown similarly.
\hfill\qed

%% file: b11.tex

\subsubsection{Perfect strict $\nbigr$-triple and vector bundle
over $\proj^1$}
\label{subsubsection;b11.23.1}

\begin{df}
$\nbigt=(\nbigm',\nbigm'',C)$ be $\nbigr$-triple in $0$-dimension
(See {\rm\cite{sabbah}}).
It is called perfect and strict if the following holds:
\begin{itemize}
\item
 The $\nbigo_{\cnum_{\lambda}}$-modules
 $\nbigm'$ and $\nbigm''$ are locally free.
\item
 The sesqui-linear pairing
 $C:\nbigm'(\AAA)\otimes \overline{\nbigm''({\AAA})}\lrarr\nbigo({\AAA})$
 is perfect.
\hfill\qed
\end{itemize}
\end{df}

Let $V$ be a vector bundle over $\proj^1$.
We put
$V_0:=V_{|\cnum_{\lambda}}$
and $V_{\infty}:=V_{|\cnum_{\mu}}$.
We have the natural perfect pairing
$C_V:V^{\lor}_0(\AAA)\otimes \sigma V_{\infty}(\AAA)
 \lrarr\nbigo(\AAA)$.
We put
$\Theta(V)=\bigl(V_0^{\lor},\sigma(V_{\infty}),C_V\bigr)$,
which is a strict perfect $\nbigr$-triple.

On the other hand,
if $\nbigt=(\nbigm',\nbigm'',C)$ is a strict perfect $\nbigr$-triple,
then $C$ induces the isomorphism
$\nbigm^{\prime\,\lor}_{|\AAA}\simeq
 \sigma^{-1}(\nbigm'')_{|\AAA}$.
Thus we obtain the vector bundle $V(\nbigt)$.

Let $V^{(i)}$ $(i=1,2)$ be a locally free $\nbigo_{\proj^1}$-modules,
and let $f:V^{(1)}\lrarr V^{(2)}$ be a morphism
of $\nbigo$-sheaves.
Then we obtain the morphism
$\bigl(f_0^{\lor},\sigma(f_{\infty})\bigr):
 \Theta(V^{(1)})\lrarr\Theta(V^{(2)})$.

\begin{lem}
By the correspondence above,
we obtain the equivalence of the following categories:
\[
 \Theta:
 (\mbox{vector bundle over }\proj^1)
\lrarr
 (\mbox{strict perfect }\nbigr\mbox{-triple})
\]
\hfill\qed
\end{lem}

For perfect strict $\nbigr$-triples
$\nbigt_i=(\nbigm_i',\nbigm_i'',C_i)$ $(i=1,2)$,
the tensor product
$\nbigt_1\otimes\nbigt_2$
is given by
$(\nbigm_1'\otimes\nbigm_2',\nbigm_1''\otimes\nbigm_2'',C_1\otimes C_2)$.
The direct sum is also naturally defined.

Let $\nbigt=(\nbigm',\nbigm'',C)$
be a perfect strict $\nbigr$-triple.
The pairing $C$ induces the perfect pairing
$C':\nbigm^{\prime\lor}({\AAA})
\otimes \nbigm^{\prime\prime\lor}({\AAA})
\lrarr\nbigo({\AAA})$.
Then we obtain 
the dual
$\nbigt^{\lor}
=(\nbigm^{\prime\,\lor},\nbigm^{\prime\prime\,\lor},C')$.

\begin{lem}
$\Theta$ preserves tensor products, duals and direct sums.
\end{lem}

We have Sabbah's hermitian adjoint
$\nbigt^{\ast}=(\nbigm'',\nbigm',C^{\ast})$
in the category of $\nbigr$-triples
(\cite{sabbah}).
The hermitian adjoint in the category of vector bundles over $\proj^1$
is given by $V^{\ast}:=\sigma(V)^{\lor}$.

\subsubsection{Tate objects}

The $\nbigr$-triple corresponding to
$\nbigo_{\proj^1}(p,q)$ is as follows:
\[
 \Theta(\nbigo_{\proj^1}(p,q))
=\Bigl(\nbigo_{\cnum_{\lambda}}\!\cdot\! f_0^{(-p,-q)},\,\,\,
 \nbigo_{\cnum_{\lambda}}\!\cdot\! \sigma(f_{\infty}^{(p,q)}),\,\,\, C
 \Bigr).
\]
Here the pairing $C$ is given as follows:
\[
 C\bigl(f_0^{(-p,-q)},\sigma(f_{\infty}^{(p,q)})\bigr)
=\bigl(\sqrt{-1}\lambda\bigr)^{p+q}.
\]
For simplicity, we denote it by $\bigl(\sqrt{-1}\lambda\bigr)^{p+q}$.
In particular,
if we forget the torus action,
$\Theta(\nbigo(n))$ is given as follows:
\[
 \Theta(\nbigo_{\proj^1}(n))
=\Bigl(\nbigo_{\cnum_{\lambda}}\!\cdot\! f^{(-n)}_0,\,\,\,\,
  \nbigo_{\cnum_{\lambda}}\!\cdot\! \sigma(f^{(n)}_{\infty}),\,\,\,\,
 \bigl(\sqrt{-1}\lambda\bigr)^n
 \Bigr).
\]
The $\nbigr$-triple
$\Theta(\sigma(\nbigo_{\proj}(n)))$
is given by
$ \Bigl(
 \nbigo_{\cnum_{\lambda}}\!\cdot\! \sigma(f_{\infty}^{(n)}),\,\,\,\,
 \nbigo_{\cnum_{\lambda}}\!\cdot\! f_0^{(-n)},\,\,\,\,
 (-1)^n\cdot\bigl(\sqrt{-1}\lambda\bigr)^n
 \Bigr)$.
The latter is same as the hermitian adjoint
$\Theta\bigl(\nbigo_{\proj^1}(-n)\bigr)^{\ast}$
in the category of $\nbigr$-triples.
Hence the isomorphism $\sigma^{\ast}\nbigo(-n)\lrarr \nbigo(-n)$
(the subsubsection \ref{subsubsection;a11.11.15})
induces the isomorphism
$\Theta\bigl(\nbigo_{\proj^1}(n)\bigr)^{\ast}
\lrarr
 \Theta\bigl(\nbigo_{\proj^1}(-n)\bigr)$.

\begin{lem}
The induced isomorphism
$\Theta\bigl(\nbigo_{\proj^1}(n)\bigr)^{\ast}
\lrarr
 \Theta\bigl(\nbigo_{\proj^1}(-n)\bigr)$
is given by the pair of the maps
$(\varphi_1,\varphi_2)$:
\[
\begin{array}{ll}
 \nbigo_{\cnum_{\lambda}}\cdot \sigma\bigl(f_{\infty}^{(n)}\bigr)
 \stackrel{\varphi_1}{\llarr}
 \nbigo_{\cnum_{\lambda}}\cdot f_0^{(n)},
 & \varphi_1(f_0^{(n)})=\sqrt{-1}{}^n\sigma(f_{\infty}^{(n)}),\\
 \mbox{{}}\\
 \nbigo_{\cnum_{\lambda}}\cdot f_0^{(-n)}
 \stackrel{\varphi_2}{\lrarr}
 \nbigo_{\cnum_{\lambda}}\cdot\sigma\bigl(f_{\infty}^{(-n)}\bigr),
 & \varphi_2\bigl(f_0^{(-n)}\bigr)
    =\sqrt{-1}{}^n\sigma\bigl(f_{\infty}^{(-n)}\bigr).
\end{array}
\]
\end{lem}
\pf
The morphism $\varphi_1$ is the dual
of the isomorphism
$\iota_{-n}:\sigma\bigl(\nbigo(-n)_{\infty}\bigr)\lrarr 
 \nbigo(-n)_0$.
The morphism $\varphi_2$ is given 
by $\sigma(\iota_{-n})$.
Then the claim can be checked by a direct calculation.
\hfill\qed

\vspace{.1in}

For any half integers $k\in\frac{1}{2}\seisuu$,
we put 
 $\Tate^{s}(k):=
 \Bigl(\nbigo_{\cnum_{\lambda}}\!\cdot\! e_0^{(2k)},\,\,\,\,
 \nbigo_{\cnum_{\lambda}}\!\cdot\! e^{(2k)}_{\infty},\,\,\,\,
 (\sqrt{-1}\lambda)^{-2k}
 \Bigr)$,
which is the Tate object in the category of $\nbigr$-triples,
namely,
the Tate twist in the category of $\nbigr$-triple 
is given by the tensor product with $\Tate^s(k)$.
Recall that the canonical isomorphism
$\Tate^s(k)\lrarr \Tate^{s}(-k)^{\ast}$ is given by 
the pair of maps $\bigl((-1)^{2k},1\bigr)$
(see the section 1.6.a in \cite{sabbah}).

We fix the complex number $a$ such that $a^2=\sqrt{-1}$.
Then we take the isomorphism
$\Phi^{n}=(\Phi_1^{(n)},\Phi_2^{(n)}):
 \Theta\bigl(\nbigo_{\proj^1}(n)\bigr)
\lrarr \Tate^s(-n/2)$ given as follows:
\[
 \begin{array}{ll}
 \nbigo_{\cnum_{\lambda}}\cdot  f_0^{(-n)}
 \stackrel{\Phi^{(n)}_1}{\llarr}
 \nbigo_{\cnum_{\lambda}}\cdot e_0^{(-n)}, &
 \Phi_1(e_0^{(-n)})=a^n\cdot f_0^{(-n)}, \\
 \mbox{{}}\\
 \nbigo_{\cnum_{\lambda}}\cdot \sigma(f_{\infty}^{(n)})
 \stackrel{\Phi^{(n)}_2}{\lrarr}
 \nbigo_{\cnum_{\lambda}}\cdot e_{\infty}^{(-n)},
 & \Phi_2(f_{\infty}^{(n)})=a^{-n}\cdot e_{\infty}^{(-n)}.
 \end{array}
\]
In particular,
the isomorphism $\Theta(\Tate(n))\lrarr \Tate^{S}(n)$
is given by the pair $\bigl(\Phi^{(-2n)}_1,\Phi^{(-2n)}_2\bigr)$:
\[
 \Phi^{(-2n)}_1(e_0^{(2n)})=(\sqrt{-1})^n \cdot t_0^{(-n)},
\quad\quad
 \Phi^{(-2n)}_2\bigl(\sigma(t_{\infty}^{(n)})\bigr)
 =(\sqrt{-1})^{-n}\cdot e_{\infty}^{(2n)}.
\]

\begin{lem} \label{lem;9.13.45}
The following diagramm is commutative:
\begin{equation} \label{eq;a11.30.1}
 \begin{CD}
 \Theta\bigl(\nbigo_{\proj^1}(n)\bigr)^{\ast}
 @>>>
 \Theta\bigl(\nbigo_{\proj^1}(-n)\bigr) \\
 @AAA @VVV \\
 \Tate^s(-n/2)^{\ast} @>>> \Tate^s(n/2).
 \end{CD}
\end{equation}
\end{lem}
\pf
We have the following diagramms:
\begin{equation}\label{eq;a11.30.2}
 \begin{CD}
 \nbigo\cdot \sigma(f_{\infty}^{(n)})
 @<{(\sqrt{-1})^{n}}<< \nbigo\cdot f_{0}^{(n)}\\
 @VV{a^{-n}}V @A{a^{-n}}AA\\
 \nbigo\cdot e_{\infty}^{(-n)}
 @<{(-1)^n}<<
 \nbigo\cdot e_0^{(n)},
 \end{CD}
\quad\quad\quad\quad
 \begin{CD}
 \nbigo\cdot f_0^{(-n)} @>{(\sqrt{-1})^n}>>
 \nbigo\cdot\sigma\bigl(f_{\infty}^{(-n)}\bigr) \\
 @AA{a^n}A @VV{a^n}V \\
 \nbigo\cdot e_0^{(-n)} @>{1}>>
 \nbigo\cdot e_{\infty}^{(n)}.
 \end{CD}
\end{equation}
Since we have $a^{-2n}\cdot(\sqrt{-1})^n=(-1)^n$
and $a^{2n}\cdot (\sqrt{-1})^n=1$ due to our choice of $a$,
we obtain the commutativity of the diagramms in (\ref{eq;a11.30.2})
and thus the commutativity of the diagramm (\ref{eq;a11.30.1}).
\hfill\qed

\begin{lem}
The functor $\Theta$ essentially commutes with Tate twists
of the both categories.
The functor $\Theta$ preserves the compatibility of the Tate twist
and the adjunction.
\end{lem}
\pf
It follows from Lemma \ref{lem;9.13.45}.
\hfill\qed

\vspace{.1in}
Recall that we have the canonical isomorphism
$\Theta\bigl(\sigma(V^{\lor})\bigr)\simeq \Theta(V)^{\ast}$.
\begin{cor}
The pairing $S:V\otimes\sigma^{\ast}(V)\lrarr \Tate(-n)$
is a polarization of pure twistor structure of weight $n$
if and only if
the induced morphism
$S:\Theta(V)\lrarr\Theta(V)^{\ast}\otimes\Tate^S(-n)$
is a polarization of pure twistor $\nbigr$-triple of weight $n$
in $0$-dimension.
\end{cor}
\pf
In the case $n=0$, it is clear.
Since $\Theta$ is compatible with the Tate twist,
we can reduce the problem to the case $n=0$, as follows:

\begin{tabular}{ll}
 & $S:V\otimes\sigma^{\ast}(V)\lrarr\Tate(-n)$ is a polarization\\ 
 \mbox{{}}\\
$\Longleftrightarrow$ &
$S(-n):\bigl(
 V\otimes\nbigo_{\proj^1}(-n)\bigr)
 \otimes
 \sigma^{\ast}\bigl(V\otimes\nbigo(-n)\bigr)
 \lrarr\Tate(0)$ is a polarization\\
 \mbox{{{}}}\\
$\Longleftrightarrow$ &
$\Theta\bigl(V\otimes\nbigo_{\proj^1}(-n)\bigr)
\lrarr\Theta\bigl(V\otimes\nbigo_{\proj^1}(-n)\bigr)^{\ast}$
 is a polarization.\\
 \mbox{{{}}}\\
$\Longleftrightarrow$ &
$\Theta(V)\otimes\Tate(n/2)
\lrarr
 \bigl(\Theta(V)\otimes \Tate(n/2)\bigr)^{\ast}$
 is a polarization.\\
 \mbox{{}}\\
$\Longleftrightarrow$ &
$\Theta(V)\lrarr\Theta(V)^{\ast}\otimes\Tate^S(-n)$ is a polarization.
\end{tabular}

\mbox{{}}
\hfill\qed

\vspace{.1in}

\begin{df}
Let $\nbigt$ be a strict perfect $\nbigr$-triple in $0$-dimension,
$W$ be a filtration of $\nbigt$,
and $\vecN$ be a tuple of nilpotent morphisms
$N_i:\nbigt\lrarr\nbigt\otimes\Tate^S(-1)$ $(i=1,\ldots,l)$.
Let $S$ be a hermitian duality of $\nbigt$ of weight $n$.
The tuple
$(\nbigt,W,\vecN,S)$ is called a Pol-MTS of type $(n,l)$
in the sense of Sabbah,
if it is isomorphic to
$\Theta(V,W,\vecN,S)$ for some Pol-MTS $(V,W,\vecN,S)$
of type $(n,l)$.
\hfill\qed
\end{df}

%% file: a81.tex

\subsubsection{The induced polarized pure twistor}
\label{subsubsection;a11.30.3}

Let $V$ be a pure twistor of weight $n$
and $S:V\otimes\sigma(V)\lrarr \Tate(-n)$ be a polarization,
which induces the pairing $\sigma(S):\sigma(V)\otimes V\lrarr \Tate(-n)$.
Then we obtain the isomorphism
$\rho_S:\sigma(V)\lrarr V^{\lor}\otimes\Tate(-n)$,
and we have the induced polarization
$(\rho_S',\rho_S''):
\Theta\bigl(\sigma(V)\bigr)
 \lrarr\Theta\bigl(\sigma(V)\bigr)^{\ast}\otimes\Tate^S(-n)$.
We can naturally regarded $\rho_S'$ and $\rho_S''$
as morphisms $V_0\lrarr \sigma\bigl(V_{\infty}^{\lor}\bigr)$.
Let $C':V_0(\AAA)\otimes\overline{V_0(\AAA)}\lrarr \nbigo(\AAA)$
be the sesqui-linear pairing given as follows:
\[
 C'(u,\bar{v}):=
 C_{\sigma(V)}\bigl(
 \rho_S''(u),\bar{v}\bigr).
\]
Then we obtain the $\nbigr$-triple
$(V_0,V_0,C')$,
and recall that
the morphism
$\bigl((-1)^n,1\bigr):(V_0,V_0,C')\lrarr 
 (V_0,V_0,C')\otimes\Tate^{S}(-n)$ gives the polarization.

Let us calculate the pairing $C'$.
First we consider the case $V=\nbigo(n)$.
We denote the pairing of $V$ by $S$.
The pairing $\sigma(S)$ gives the correspondence
$\sigma(f_0^{(n)})\otimes f_{\infty}^{(n)}
\longmapsto (\sqrt{-1})^{-n}\cdot t_{\infty}^{(-n)}$.
Thus we have 
$\rho_S(\sigma(f_0^{(n)}))=
 (\sqrt{-1})^{-n}\cdot f_{\infty}^{(-n)}
 \otimes t_{\infty}^{(-n)}$.
Hence we have the following:
\begin{multline}
 \rho_S''(f_0^{(n)})=
 \sigma\bigl(
 \rho_S(\sigma(f_0^{(n)}))
 \bigr)
=(\sqrt{-1})^n\cdot \sigma(f_{\infty}^{(-n)})
 \otimes \sigma(t_{\infty}^{(-n)}) \\
\longmapsto
 (\sqrt{-1})^n\cdot \sigma(f_{\infty}^{(-n)})
 \cdot(\sqrt{-1})^n\cdot e_{\infty}^{(-2n)}
 =(-1)^n\cdot \sigma(f_{\infty}^{(-n)})\cdot e_{\infty}^{(-2n)}.
\end{multline}
Then the pairing
$C':\nbigo(n)_0(\AAA)\otimes\nbigo(n)_0(\AAA)\lrarr\nbigo(\AAA)$
is given as follows:
\begin{multline}
 C'\bigl(f_0^{(n)}\otimes \overline{f_0^{(n)}}\bigr)
=C_{\sigma(\nbigo(n))}
 \bigl(
 \sigma(f_{\infty}^{(-n)})\cdot(-1)^n,
 \overline{f_0^{(n)}}
 \bigr)
=(-1)^n\cdot \bigl\langle
 \sigma(f_{\infty}^{(-n)}),
 \sigma(f_{0}^{(n)})
 \bigr\rangle
=(-1)^n\cdot\sigma^{\ast}
 \bigl\langle f_{\infty}^{(-n)},f_0^{(n)}
 \bigr\rangle \\
=(-1)^n\cdot
 \sigma^{\ast}\bigl(
 (\sqrt{-1}\lambda)^{-n}
 \bigr)
=(\sqrt{-1}\lambda)^n.
\end{multline}
On the other hand,
we have the following:
\[
 S\bigl(f_0^{(n)},\sigma\bigl(f_0^{(n)}\bigr)\bigr)
=S\bigl(
 f_0^{(n)},
 \sigma\bigl((\sqrt{-1}\lambda)^{-n}\cdot f_{\infty}^{(n)}\bigr)
 \bigr)
=(\sqrt{-1}\lambda)^n\cdot(-1)^n\cdot
 (\sqrt{-1})^{-n}\cdot t_0^{(-n)}
=(\sqrt{-1}\lambda)^{n}\cdot (\sqrt{-1})^n\cdot t_0^{(-n)}.
\]
Hence we have the following formula:
\[
 C'\bigl(f_0^{(n)},\overline{f_0^{(n)}}\bigr)
\cdot t_0^{(-n)}
=(\sqrt{-1})^{-n}\cdot
 S\bigl(f_0^{(n)},\sigma^{\ast}(f_0^{(n)})\bigr).
\]

\begin{lem}\label{lem;a11.30.10}
Let $(V,S)$ be a polarized pure twistor of weight $n$.
Then the $\nbigr$-triple
$(V_0,V_0,C')$ is a pure twistor of weight $n$,
where the sesqui-linear pairing $C'$ is given as follows:
\[
 C'(u,\overline{v})\cdot t_0^{(-n)}
=(\sqrt{-1})^{-n}\cdot S(u\otimes\sigma^{\ast}(v)).
\]
\end{lem}
\pf
Since any polarized pure twistor of weight $n$
is isomorphic to a direct sum of $(\nbigo(n),S)$,
we can reduce the lemma to the case
$V=\nbigo(n)$.
It has been already shown above.
\hfill\qed

%% file: a81.1.tex

\subsubsection{A lemma}

Let $(V,W,\vecN,S)$ be a Pol-MTS of type $(0,l)$.
Recall that we obtain the
Pol-MTS $(\nbigv_{\lbar}(V),W,\vecnbign,\nbigv_I(S))$
of type $(1-l,l+1)$,
and we obtain the Pol-MTS 
$(P_k\Gr^{W(N)}_k\nbigv_{\lbar}(V),W,\vecN,S_k)$
of type $(1-l+k,l)$.
We have the subbundle
$\nbigc_{\lbar}\subset P_k\Gr^{W(N)}_k \nbigv_{\lbar}(V)$,
which is pure twistor of weight $1-l+k$.
Moreover $S_k$ induces the polarization of
$\nbigc_{\lbar}$.
(See the subsubsection \ref{subsubsection;9.13.50}
and the subsubsection \ref{subsubsection;a11.10.51}).
Hence we obtain the pure twistor
$\nbigc_{\lbar}\otimes\nbigo(l-1)$
and the polarization of weight $k$.

We put $\big(\nbigc_{\lbar}\otimes\nbigo(l-1)\bigr)_{0}:=
 \bigl(\nbigc_{\lbar}\otimes\nbigo(l-1)\bigr)_{|\cnum_{\lambda}}$.
Due to the result in the subsubsection \ref{subsubsection;a11.30.3},
we obtain the pure twistor
$\bigl(
 (\nbigc_{\lbar}\otimes\nbigo(l-1))_{0},
 (\nbigc_{\lbar}\otimes\nbigo(l-1))_{0},
 C'
 \bigr)$.
Let us calculate the sesqui-linear pairing $C'$.
The endomorphism $\bar{N}_a\in\End(V_0,V_0)$ is determined by
$N_a=\bar{N}_a\cdot t_0^{(-1)}$.

We have the following:
\[
 \nbigv_{l}(V)_0
=\bigoplus_{i=0}^{\infty}V_0\otimes t_0^{(i)}\cdot N^i
 \Big/ \prod_{a=1}^l(N-N_a).
\]
For any $u_i,v_j\in V(\AAA)$ and for any $k\in\seisuu$,
we have the following:
\begin{multline}
 S'\bigl(u_i\cdot t_0^{(i)}\cdot N^{i+k},
 \sigma\bigl(v_j\cdot t_0^{(j)}\cdot N^j\bigr)
 \bigr)
=(-1)^j\cdot
 S\bigl(u_i,\sigma(v_j)\bigr)
 \cdot (\sqrt{-1}\lambda)^{-2j}\cdot (-1)^j\cdot t_0^{(i+j+k)} \\
=S\bigl(u_i,\sigma(v_j)\bigr)
 \cdot (\sqrt{-1}\lambda)^{-2j}\cdot t_0^{(i+j+k)}.
\end{multline}
Here we have used
$\sigma(t_0^{(j)})=
 \sigma\bigl((\sqrt{-1}\lambda)^{2j}\cdot t_{\infty}^{(j)}\bigr)
=(\sqrt{-1}\lambda)^{-2j}\cdot(-1)^j\cdot t_0^{(j)}$.
Then we have the following:
\begin{multline}
 S'\bigl(
 \tilde{N}^{-1}_{\lbar}
 \bigl(u_i\cdot N^{i+k}\cdot t_0^{(i)}\bigr),
 \sigma\bigl(v_j\cdot t_0^{(j)}\cdot N^j\bigr)
 \bigr)
=\sum_{\vecn}(-1)^j
 \cdot N^{-l-|\vecn|+i+j+k}
 S\Bigl(\prod_{a=1}^lN_a^{n_a}u_i\otimes t_0^{(i)},
 \sigma\bigl(v_j\cdot t_0^{(j)}\bigr)
 \Bigr)\\
=\sum_{\vecn}(\sqrt{-1}\lambda)^{-2j}
 \cdot
 N^{-l-|\vecn|+i+j+k}
 \cdot
 S\Bigl(\prod_{a=1}^l\bar{N}_a^{n_a}\cdot u_i,\sigma(v_j)\Bigr)
 \cdot t_0^{(-|\vecn|+i+j)}.
\end{multline}
Hence we obtain the following:
\[
 \nbigv_{\lbar}(S)\bigl(u_i\cdot t_0^{(i)}\cdot N^{i+k},
 \sigma\bigl(v_j\cdot t_0^{(j)}\cdot N^j\bigr)
 \bigr)
=\sum_{\vecn\in S(k)}
 (\sqrt{-1}\lambda)^{-2j}\cdot
 S\Bigl(
 \prod_{a=1}^n \bar{N}_a^{n_a}\cdot u_i,\sigma(v_j)
 \Bigr)\cdot t_0^{(l-1-k)}.
\]
Here we put 
$S(k):=\bigl\{\vecn\in\seisuu^l,\,\big|\,
 -|\vecn|-l+i+j+k=-1\bigr\}$.

\begin{lem}
When we twist by $\nbigo(l-1)$,
we have the following:
\[
 \nbigv_{\lbar}(S)\Bigl(
 u_i\cdot t_0^{(i)}\cdot N^{i+k}\cdot f_0^{(l-1)},\,\,
 \sigma\bigl(v_j\cdot t_0^{(j)}\cdot N^j\cdot f_0^{(l-1)}\bigr)
 \Bigr)
=\sum_{\vecn\in S(k)}
 (\sqrt{-1}\lambda)^{l-1-2j}
 \cdot (\sqrt{-1})^{l-1}
 \cdot S\Bigl(\prod_{a=1}^l N_a^{n_a}\cdot u_i,\sigma(v_j)\Bigr)
 \cdot t_0^{(-k)}.
\]
\end{lem}
\pf
We have used the following correspondence
for the canonical polarization of $\nbigo(n-1)$:
\begin{multline}
 f_0^{(l-1)}\otimes\sigma\bigl(f_0^{(l-1)}\bigr)
 =f_0^{(l-1)}\otimes
  \sigma\bigl((\sqrt{-1}\lambda)^{-l+1}\cdot f_{\infty}^{(l-1)}\bigr)
 \longmapsto \\
 (-1)^{-l+1}\cdot
 (\sqrt{-1}\lambda)^{l-1}\cdot (\sqrt{-1})^{-l+1}
=(\sqrt{-1})^{l-1}\cdot (\sqrt{-1}\lambda)^{l-1}.
\end{multline}
\hfill\qed

Then we obtain the following formula:
\begin{multline}
 \nbigv_{\lbar}S\Bigl(
 u_i\cdot t_0^{(i)}\cdot (-\sqrt{-1}N)^{i}\cdot N^k\cdot f_0^{(l-1)},\,
 \sigma\bigl(
 v_j\cdot t_0^{(j)}\cdot (\sqrt{-1}N)^{j}\cdot f_0^{(l-1)}
 \bigr) 
 \Bigr) \\
=\sum_{\vecn\in S(k)}
 (\sqrt{-1}\lambda)^{l-1-2j}
 \cdot {\sqrt{-1}}^{k}\cdot
 S\Bigl(
 \prod_{a=1}^l (-\sqrt{-1}\bar{N}_a)\cdot u_i,
 \sigma(v_j)
 \Bigr)\cdot t_0^{(-k)}.
\end{multline}

Hence we obtain the following.
\begin{lem}\label{lem;a11.30.41}
We put $Y:=\nbigc_{\lbar}\otimes\nbigo(l-1)$.
Then the $\nbigr$-triple
$(Y_0,Y_0,C')$ is pure twistor of weight $k$.
Here $C'$ is given by the following:
\[
 C'\bigl(u_i\cdot (-\sqrt{-1}N)^i,
 \overline{v\cdot (\sqrt{-1}N)^j}\bigr)
=\sum_{\vecn\in S(k)}
 (\sqrt{-1}\lambda)^{l-1-2j}
 \cdot S\Bigl(
 \prod_{a=1}^l(-\sqrt{-1}\bar{N}_a)^{n_a}\cdot u_i,
 \sigma(v_j)
 \Bigr)
\]
\end{lem}
\pf
We have only to apply Lemma \ref{lem;a11.30.10}.
\hfill\qed

%% file: a50.tex

\subsubsection{Filtration}

Let $k$ be a field.
Let $V$ be a finite dimensional vector space over $k$.
\begin{df}
An increasing filtration $F$ of $V$ indexed by $\real$ is defined to be
a family of subspaces
$\bigl\{F_{\eta}\subset V\,\big|\,\eta\in\real\bigr\}$
satisfying the conditions
$F_{\eta}\subset F_{\eta'}$ $(\eta\leq \eta')$
and $F_{\eta}=V$ for any sufficiently large $\eta$.
In this paper,
we mainly use the increasing filtration.
So `filtration' will mean `increasing filtration'
if we do not mention.
\hfill\qed
\end{df}

When we consider a tuple of filtrations,
we often use the notation $\vecF=(\lefttop{i}F\,|\,i\in I)$.
We also use the notation $\lefttop{i}\Gr^F$
to denote the associated graded vector space
$\Gr^{\lefttop{i}F}$.

Let $\vecF:=\bigl(\lefttop{i}F\,\big|\,i\in I\bigr)$
be a tuple of filtrations indexed by $\real$.
For any subset $J\subset I$ and for any $\veceta\in\real^J$,
we put as follows:
\[
 \lefttop{J}F_{\veceta}=\bigcap_{j\in J}\lefttop{j}F_{\eta_j}.
\]
\begin{df} \label{df;b11.12.1}
$\vecF:=\bigl(\lefttop{i}F\,\big|\,i\in I\bigr)$
be a tuple of filtrations indexed by $\real$.
It is called compatible,
if we have a decomposition of
 $V=\bigoplus_{\veceta\in\real^I}U_{\veceta}$
satisfying the following:
\begin{equation} \label{eq;9.10.1}
 \lefttop{I}F_{\vecrho}
=\bigoplus_{\veceta\in\real^l,\veceta\leq \vecrho}
 U_{\veceta}.
\end{equation}
If $\vecF$ is compatible,
a decomposition of $V$ satisfying {\rm(\ref{eq;9.10.1})}
is called a splitting of $\vecF$.
\hfill\qed
\end{df}

Recall that 
$\vecF=(\lefttop{i}F\,|\,i\in I)$ is called sequentially compatible
if the following is satisfied:
\begin{itemize}
\item
 We have the induced filtrations $\lefttop{i}F^{(1)}$
 $(i=2,\ldots,l)$ on $\lefttop{1}\Gr^{F}$.
 Then $\lefttop{i}F^{(1)}$ are sequentially compatible.
\item
 The following map is surjective:
 \[
 \bigcap_{i=1}^l \lefttop{i}F_{h_i}
\lrarr
 \lefttop{1}\Gr^{F}_{h_1}\cap
 \bigcap_{i=2}^l
 \lefttop{i}F^{(1)}_{h_i}.
 \]
\end{itemize}

\begin{lem}
Compatibility and sequential compatibility
are equivalent.
\end{lem}
\pf
We saw that sequential compatibility implies
compatibility in our previous paper \cite{mochi}.
It is easy to see that compatibility implies
sequential compatibility.
\hfill\qed

\vspace{.1in}

\noindent
{\bf Notation}
Let $\vecF=(\lefttop{i}F\,|\,i\in I)$ be a tuple of filtrations.
For any subset $J\subset I$
and $\veceta\in\real^J$,
we put as follows:
\[
 {\displaystyle
 \lefttop{J}\Gr^F_{\veceta}(V):=
 \frac{\lefttop{J}F_{\veceta}}{
 \sum_{\veceta'\lneq \veceta}
 \lefttop{I}F_{\veceta'}  }
 }.
\]

\begin{lem}\label{lem;10.12.31}
Let $\vecF=(\lefttop{i}F\,\big|\,i\in I)$ be a tuple of filtrations
on $V$.
Then we have the following inequality:
\[
 \sum_{\veceta} \dim \lefttop{I}\Gr^F_{\veceta}\geq \dim V.
\]
If the equality holds, then
the filtrations $\lefttop{i}F$ $(i=1,\ldots,l)$ are compatible.
\end{lem}
\pf
We have the surjection
$\pi_{\veceta}:\lefttop{I}F_{\veceta}\lrarr \lefttop{I}\Gr^F_{\veceta}$.
We pick subspace $U_{\veceta}\subset\lefttop{I}F_{\veceta}$
such that the restriction of $\pi_{\veceta}$ to $U_{\veceta}$
gives the isomorphism of $U_{\veceta}$ and $\lefttop{I}\Gr^F_{\veceta}$.
Then we have the naturally defined surjection
$f:\bigoplus_{\veceta}U_{\veceta}\lrarr V$.
It implies the inequality.
If the equality holds, then $f$ is isomorphic,
and thus the decomposition $V=\bigoplus_{\veceta}U_{\veceta}$
gives a splitting of the filtrations.
\hfill\qed

\subsubsection{Decomposition}

Let $k$ be a field.
Let $V$ be a finite dimensional vector space over $k$.
\begin{df}
A decomposition $\vecE$ of $V$ indexed by a set $S$
is defined to be a family of subspaces
$\bigl\{\EE_{\alpha}\subset V\,\big|\,\alpha\in S\bigr\}$
such that $V=\bigoplus_{\alpha\in S}\EE_{\alpha}$.
\hfill\qed
\end{df}

When we consider a tuple of decompositions,
we often use the notation
$\vecE=\bigl(\lefttop{i}\EE\,\big|\,i\in I\bigr)$.
For any subset $J\subset I$
and for any $\vecalpha\in\cnum^J$,
we put as follows:
\[
 \lefttop{J}\EE_{\vecalpha}:=
 \bigcap_{j\in J}\lefttop{j}\EE_{\vecalpha_j}.
\]
\begin{df}
Let $\vecE=\bigl(\lefttop{i}\EE\,\big|\,i\in I\bigr)$
be a tuple of decompositions indexed by a set $S$.
It is called compatible,
if the following holds:
\[
 V=\bigoplus_{\vecalpha\in S^I}\lefttop{I}\EE_{\vecalpha}.
\]
\hfill\qed
\end{df}

\subsubsection{Filtrations and decompositions}
\label{subsubsection;a11.12.1}

Let $V$ be a vector space over $k$.
\begin{df} \label{df;b11.12.2}
Let $\vecE=\bigl(\lefttop{i}\EE\,\big|\,i\in I\bigr)$
be a tuple of decompositions of $V$.
Let $\vecF=\bigl(\lefttop{j}F\,\big|\,j\in J\bigr)$
be a tuple of filtrations of $F$.
Then the tuple $(\vecE,\vecF)$ is called compatible,
if the following holds:
\begin{itemize}
\item
 The tuple $\vecE$ of decompositions is compatible.
\item
 The tuple $\vecF$ of the filtrations is compatible,
 in the sense of Definition {\rm\ref{df;b11.12.1}}.
\item
 We have
$\lefttop{j}F=
 \bigoplus \bigl(
 \lefttop{i}F_b\cap \lefttop{I}\EE_{\vecalpha}\bigr)$.
\hfill\qed
\end{itemize}
\end{df}

Let $J_1$ be a subset of $J$
and $\veceta$ be an element of $\real^{J_1}$.
The filtrations $\lefttop{i}F$ $(i\not\in J_1)$ induce the filtration
on $\lefttop{J}\Gr_{\veceta}^F$.
We denote the induced filtration by $\lefttop{i}F^{J_1}$.
Similarly, we use the notation $\lefttop{J_2}F^{J_1}$
and $\lefttop{I}\EE^{J_1}$.
We use the notation
$\vecF^{J_1}$ and
$\vecE^{J_1}$
to denote the induced tuples
$\bigl(\lefttop{j}F^{J_1}\,\big|\,j\in J-J_1\bigr)$
and 
$\bigl(\lefttop{i}\EE^{J_1}\,\big|\,i\in I\bigr)$.

%% file: a50.1.tex
\subsubsection{Filtrations on a vector bundle}

Let $X$ be a complex manifold or scheme
equipped with an action of a finite group $G$.
Let $V$ be a $G$-equivariant vector bundle.

\begin{df}\mbox{{}}
\begin{itemize}
\item
A $G$-filtration $F$ of $V$ indexed by $\real$
in the category of $G$-equivariant
vector bundles is defined to be a family of
$G$-equivariant subbundles
$\bigl\{F_{\eta}\subset V\,\big|\,\eta\in\real\bigr\}$
satisfying the conditions
$F_{\eta}\subset F_{\eta'}$ $(\eta\leq \eta')$
and $F_{\eta}=V$ for any sufficiently large $\eta$.
\item
A $G$-decomposition $\EE$ of $V$ indexed by $\cnum$
in the category of $G$-equivariant vector bundles is defined to be
a family of $G$-equivariant subbundles
$\bigl\{\EE_{\alpha}\,\big|\,\alpha\in\cnum\bigr\}$
such that $V=\bigoplus_{\alpha\in\cnum}\EE_{\alpha}$.
\hfill\qed
\end{itemize}
\end{df}
We often omit to denote ``$G$-'',
if there are any confusion.

Let $\vecF=\bigl(\lefttop{i}F\,\big|\,i\in I\bigr)$ be a tuple of
$G$-filtrations of $V$.
Let $P$ be a point of $X$.
Then we obtain the tuple of filtrations
$\bigl(\lefttop{i}F_{|P}\big|\,i\in I\bigr)$
of the vector space $V_{|P}$.
We denote the tuple by $\vecF_{|P}$.
Similarly, we use the notation $\vecE_{|P}$
to denote the restriction of $\vecE$ to the fibers $V_{|P}$.

\begin{df} \label{df;b11.12.3}
Let $\vecE=\bigl\{\lefttop{i}\EE\,\big|\,i\in I\bigr\}$
be a tuple of $G$-decompositions of $V$.
Let $\vecF=\bigl\{\lefttop{j}F\,\big|\,j\in J\bigr\}$
be a tuple of $G$-filtrations of $V$.
The tuple $(\vecE,\vecF)$ is called compatible,
if the following holds:
\begin{itemize}
\item
 For any point $P\in X$,
 the tuple
 $(\vecE_{|P},\vecF_{|P})$  is compatible
 in the sense of Definition {\rm\ref{df;b11.12.2}}.
\item
 Let $I_1$ and $J_1$ be subsets of $I$ and $J$ respectively.
 Let $\vecalpha$ and $\veceta$ be elements of
 $\real^{J_1}$ and $\cnum^{I_1}$.
 Then 
 $\Bigl\{\lefttop{J_1}F_{\veceta\,|\,P}\cap
 \lefttop{I_1}\EE_{\vecalpha\,|\,P}\,\Big|\,P\in X
 \Bigr\}$ forms a vector bundle over $X$.
\hfill\qed
\end{itemize}
\end{df}

The second condition can be reworded as follows:
\begin{lem}
We put as follows:
\[
 \bar{d}(\vecalpha,\veceta,P):=
 \dim\Bigl(
 \lefttop{J_1}F_{\veceta\,|\,P}\cap
 \lefttop{I_1}\EE_{\vecalpha\,|\,P}
 \Bigr).
\]
Then
 $\Bigl\{\lefttop{J_1}F_{\veceta\,|\,P}\cap
 \lefttop{I_1}\EE_{\vecalpha\,|\,P}\,\Big|\,P\in X
 \Bigr\}$ forms a vector bundle over $X$,
if and only if the numbers $\bar{d}(\vecalpha,\veceta,P)$
 are independent of a choice of $P\in X$.
\hfill\qed
\end{lem}

\begin{rem}
We put as follows:
\[
 d(\vecalpha,\veceta,P):=
 \dim\Bigl(
 \lefttop{J}\Gr^{F}_{\veceta\,|\,P}\cap\lefttop{I}\EE^J_{\vecalpha\,|\,P}
 \Bigr).
\]
If the tuple $(\vecE_{|P},\vecF_{|P})$ is compatible,
the family of the numbers
$\bigl\{d(\vecalpha,\veceta,P)\,\big|\,
 \vecalpha\in\cnum^I,\,\veceta\in\real^J
 \bigr\}$
can be reconstructed from
the family of the numbers
$\bigl\{\bar{d}(\vecalpha,\veceta,P)\,\big|\,
 \vecalpha\in\cnum^I,\,\veceta\in\real^J
 \bigr\}$.
Hence if $(\vecE,\vecF)$ is compatible,
we obtain the vector bundle
$\lefttop{J}\Gr^F_{\veceta}\cap\lefttop{I}\EE^{J}_{\vecalpha}$
on $X$.

Such an argument will be used in many times without mention.
\hfill\qed
\end{rem}

Let $\vecE=\bigl(\lefttop{i}\EE\,\big|\,i\in I\bigr)$ be a
tuple of $G$-decompositions of $V$,
and let $\vecF=\bigl(\lefttop{j}F\,\big|\,j\in J\bigr)$ be a tuple
of $G$-filtrations of $V$.
For any subset $I_1\subset I$,
we denote the tuple 
$\bigl(\lefttop{i}\EE\,\big|\,i\in I_1\bigr)$
by $q_{I_1}(\vecE)$.
We also use the notation
$q_{J_1}(\vecF)$ for any subset $J_1\subset J$.

\begin{lem}
Let $\vecE$, $\vecF$ be as above.
Assume that $(\vecE,\vecF)$ is compatible.
Then $(q_{I_1}(\vecE),q_{J_1}(\vecF))$ is compatible
for any subsets $I_1\subset I$ and $J_1\subset J$.
\hfill\qed
\end{lem}

\begin{df} \label{df;b11.12.4}
Let $\vecE=\bigl(\lefttop{i}\EE\,\big|\,i\in I\bigr)$ be a
tuple of $G$-decompositions of $V$,
and let $\vecF=\bigl(\lefttop{j}F\,\big|\,j\in J\bigr)$ be a tuple
of $G$-filtrations of $V$.
Assume that $(\vecE,\vecF)$ is compatible.
A decomposition
$V=
 \bigoplus_{
 (\veca,\vecalpha)\in\real^J\times\cnum^I}
 U_{(\veca,\vecalpha)}$ is called a splitting of
$(\vecE,\vecF)$, if the following holds
for any $(\veca,\vecalpha)\in\real^J\times\cnum^I$:
\[
 \lefttop{J}F_{\veca}\cap\lefttop{I}\EE_{\vecalpha}
 =\bigoplus_{\vecb\leq\veca}
 U_{(\vecb,\vecalpha)}.
\]
\hfill\qed
\end{df}

\subsubsection{Compatible tuple $(\vecE,\vecF,\vecW)$}
\label{subsubsection;a12.9.2}

Let $V$ be a $G$-equivariant vector bundle over $X$.
Let $\vecE=\bigl(\lefttop{s}\EE\,\big|\,s\in S\bigr)$
be a compatible tuple of $G$-decompositions of $V$.
Let $\vecF=\bigl(\lefttop{i}F\,\big|\,i\in \lbar\bigr)$
be a compatible tuple of $G$-filtrations of $V$.
Let $m$ be an integer such that $1\leq m\leq l$.
For any element $(\vecalpha,\veceta)\in\cnum^S\times\real^m$,
we have the $G$-vector bundle
$\lefttop{\mbar}\Gr^F_{\veceta}\cap
 \lefttop{S}\EE^{\mbar}_{\vecalpha}$.

Let $l_1$ be an integer such that $1\leq l_1\leq l$.
Let us consider a tuple
$\vecW=\bigl(W(\mbar)\,\big|\,m\in\lbariti\bigr)$
such that
each $W(\mbar)$
is a $G$-filtration 
of $\lefttop{\mbar}\Gr^F_{\veceta}$
indexed by $\seisuu$
which is compatible
with the tuple $(\vecE^{\mbar},\vecF^{\mbar})$.
(See the subsubsection \ref{subsubsection;a11.12.1}
 for $\vecE^{\mbar}$ and $\vecF^{\mbar}$).
In that case,
the filtrations $W(\ibar)$ $(i\leq m)$
induce the filtrations on $\lefttop{\mbar}\Gr^F_{\veceta}$
for any $\veceta\in\real^m$.
We denote the induced filtrations by
$W^{\mbar}(\ibar)$.
We denote the tuple
$\bigl(W^{\mbar}(\ibar)\,\big|\,i\leq m\bigr)$
by $\vecW^{\mbar}$.

\begin{df}\label{df;10.7.16}
Let $V$, $\vecE$, $\vecF$, $l_1$ and $\vecW$ be as above.
The tuple
$(\vecE,\vecF,\vecW)$ is called compatible
if $(\vecE^{\mbar},\vecF^{\mbar},\vecW^{\mbar})$
is compatible for any $m\leq l$,
in the sense of Definition {\rm \ref{df;b11.12.3}}.
\hfill\qed
\end{df}

Let $J$ be a subset of $\lbar$.
Let $m(J)$ denote the number determined by the conditions
$\underline{m(J)}\subset J\cap \lbariti$ and
 $m(J)+1\not\in J\cap \lbariti$.
We put $q_J(\vecW):=\bigl(W(\ibar)\,\big|\,i=1,\ldots,m(J)\bigr)$.
Then the following lemma is easy to see.
\begin{lem}
Let $S'$ and $J$ be subsets of $S$ and $I$.
Then we obtain the tuple
$\bigl(q_{S'}(\vecE),q_{J}(\vecF),q_J(\vecW)\bigr)$.
It is compatible.
\hfill\qed
\end{lem}

\subsubsection{Splitting of $(\vecE,\vecF,\vecW)$}

Let $l'$ be any integer such that $1\leq l'\leq l$,
and we put $l'_1:=\min\{l',l_1\}$.
We have the projection
$\real^{l}\lrarr\real^{l'}$, taking the first $l'$-components.
Then we obtain the projection
$\pi_1:
 \cnum^S\times\real^l\times\seisuu^{l_1}
 \lrarr
 \cnum^S\times\real^{l'}$.
We also have the projection
$\pi_3:\cnum^S\times\real^l\times\seisuu^{l_1'}
\lrarr\cnum^S\times\real^{l'}$.

On the other hand,
we have the projection $\seisuu^{l_1}\lrarr\seisuu^{l_1'}$,
taking the first $l_1'$-components.
Then we obtain the projection
$\pi_2:\cnum^S\times\real^l\times\seisuu^{l_1}
\lrarr\cnum^S\times\real^l\times\seisuu^{l_1'}$.
Note we have $\pi_1=\pi_3\circ\pi_2$.

Let us consider a $G$-equivariant decomposition of $V$:
\begin{equation} \label{eq;10.7.1}
 V=\bigoplus_{\vecu\in\cnum^S\times\real^l\times\seisuu^{l_1}}
 U_{\vecu}.
\end{equation}
For elements $\vecv\in\cnum^S\times\real^{l'}$
and $\vecu_1\in\cnum^S\times\real^{l}\times\seisuu^{l_1'}$,
we put as follows:
\[
 \lefttop{l'}C_{\vecv}:=
 \bigoplus_{\pi_1(\vecu)=\vecv}
 U_{\vecu},
\quad
 \lefttop{l_1'}B_{\vecu_1}:=
 \bigoplus_{\pi_2(\vecu)=\vecu_1}
 U_{\vecu}.
\]
Then we have the following:
\[
 \lefttop{l'}C_{\vecv}
=\bigoplus_{\pi_3(\vecu_1)=\vecv}\lefttop{l'_1}B_{\vecu_1}.
\]

\begin{df} \label{df;a11.12.2}
Assume that $(\vecE,\vecF,\vecW)$ is compatible.
The decomposition {\rm(\ref{eq;10.7.1})} is called a splitting
of $(\vecE,\vecF,\vecW)$, if the following holds:
\begin{enumerate}
\item
 The decomposition
 $V=\bigoplus_{\vecv\in\cnum^S\times\real^{l'}}
 \lefttop{l'}C_{\vecv}$
 is a splitting of
 the tuple  $(\vecE,q_{\lbar'}(\vecF))$ for any $l'$,
 in the sense of Definition {\rm \ref{df;b11.12.4}}.
\item
 For any element $\vecv=(\vecalpha,\veca)\in\cnum^S\times\real^{l'}$,
 we obtain the following decomposition
 via the isomorphism
 $\lefttop{\lbar'}\Gr^F_{\veca}
 \cap
 \lefttop{S}\EE^{\lbar'}_{\vecalpha}
 \simeq \lefttop{l'}C_{\vecv}$ given by the previous condition:
\begin{equation} \label{eq;10.7.2}
 \lefttop{\lbar'}\Gr^F_{\veca}
\cap
 \lefttop{S}\EE^{\lbar'}_{\vecalpha}
\simeq
 \bigoplus_{\pi_3(\vecu_1)=\vecv}
 \lefttop{l'_1}B_{\vecu_1}.
\end{equation}
 Then the decomposition {\rm (\ref{eq;10.7.2})}
 is a splitting of the tuple
 $(\vecF^{\lbar'},\vecW^{\lbar'})$
 in the sense of Definition {\rm\ref{df;b11.12.4}}.
\end{enumerate}
We also say that
$\bigl(U_{\vecu}\,|\,
 \vecu\in\cnum^S\times\real^l\times \seisuu^{l_1}\bigr)$
is a splitting.
\hfill\qed
\end{df}

Let $(\vecE,\vecF,\vecW)$ be a compatible tuple,
and let $\bigl(U_{\vecu}\,\big|\,
 \vecu\in\cnum^S\times\real^l\times \seisuu^{l_1}\bigr)$
be a splitting of $(\vecE,\vecF,\vecW)$.
Let $S'$ and $J$ be subsets of $S$ and $I$ respectively.
We have the naturally defined projection:
\[
 q_{S',J}:\cnum^S\times\real^l\times\seisuu^{l_1}
\lrarr\cnum^{S'}\times\real^{J}\times\seisuu^{m(J)}.
\]
Then we put as follows, for any element
$\vecu_1\in\cnum^{S'}\times\real^{J}\times\seisuu^{m(J)}$:
\[
 U'_{\vecu_1}:=
 \bigoplus_{q_{S',J}(\vecu)=\vecu_1} U_{\vecu}.
\]
The following lemma is easy to see.
\begin{lem}
The decomposition
$\bigl( U'_{\vecu_1}\,\big|\,
 \vecu_1\in \cnum^{S'}\times\real^{J}\times\seisuu^{m(J)}
 \bigr)$
gives a splitting of the compatible tuple
$\bigl(q_{S}(\vecE),q_J(\vecF),q_J(\vecW)\bigr)$.
\hfill\qed
\end{lem}

%% file: a51.tex

\subsubsection{A lemma}

Let $V$ be a vector space over $k$.
\begin{lem}\label{lem;9.16.25}
Let $F$ be a filtration on $V$,
and $N$ be a nilpotent map on $V$ preserving the filtration $F$.
On the associated graded vector space $\Gr^F(V)$,
we have the induced filtration $W^F(N)$
and the induced nilpotent map $N^F$.
Assume $W(N^F)=W^F(N)$.

Then the induced filtration $F^{(1)}$
and the primitive decomposition of $\Gr^{W(N)}$
are compatible.
\end{lem}
\pf
Let us consider the induced isomorphism:
\begin{equation} \label{eq;9.16.10}
N^h:\Gr^{W(N)}_h\lrarr \Gr^{W(N)}_{-h}.
\end{equation}
It preserves the filtration $F^{(1)}$.
First of all,
we would like to show that the morphism (\ref{eq;9.16.10})
is strict with respect to the filtration $F^{(1)}$.

Since the isomorphism (\ref{eq;9.16.10}) preserves the filtration,
the equality
$\dim F^{(1)}_a\cap \Gr^{W(N)}_h=
 \dim F^{(1)}_a\cap\Gr^{W(N)}_{-h}$
implies the strictness of (\ref{eq;9.16.10}).
So we have only to show the equality.
For that purpose,
we have only to show the following equality 
for any $a$ and any $h$:
\begin{equation}\label{eq;9.16.11}
\dim \Gr^{F^{(1)}}_a \bigl(
 \Gr^{W(N)}_h\bigr)
=\dim \Gr^{F^{(1)}}_a\bigl(
 \Gr^{W(N)}_{-h}\bigr).
\end{equation}

Note we have the following equalities for any $a$ and $h$:
\begin{equation}\label{eq;9.16.12}
 \dim \Gr_a^{F^{(1)}}\bigl(
 \Gr^{W(N)}_h\bigr)
=\dim \Gr_h^{W^F(N)}\bigl(
 \Gr^{F}_a\bigr)
=\dim \Gr^{W(N^F)}_h\bigl(
 \Gr_a^F\bigr).
\end{equation}
By definition of the weight filtrations,
we have the following equality:
\begin{equation}\label{eq;9.16.13}
 \dim\Gr^{W(N^F)}_h\bigl(
 \Gr_a^{F}\bigr)
=\dim\Gr^{W(N^F)}_{-h}\bigl(
 \Gr_a^{F}\bigr).
\end{equation}
From (\ref{eq;9.16.12}) and (\ref{eq;9.16.13}),
we obtain the equality (\ref{eq;9.16.11}).
Thus we obtain the strictness of the morphism of (\ref{eq;9.16.10}).

Recall that the primitive part is given as follows:
\[
 P_{h}\Gr_h^{W(N)}=\ker(N^{h+1}:\Gr^{W(N)}_h\lrarr \Gr^{W(N)}_{-h-2}),
\]
We put
$P_h\Gr^{W(N)}_{h-2a}=N^aP_h\Gr^{W(N)}_h$
for any $0\leq a\leq 2h$,
and then we have the decomposition.
\[
 \Gr^{W(N)}_h
=\bigoplus_{a\geq 0}
 P_{|h|+2a}\Gr^{W(N)}_{h}.
\]
Let $x$ be an element of $F_b^{(1)}\cap \Gr^{W(N)}_h$.
We have the primitive decomposition $x=\sum x_a$,
where $x_a\in P_{|h|+2a}\Gr^{W(N)}_{h}$.
We would like to show each element $x_a$ is contained in $F_b^{(1)}$.

Since the isomorphism (\ref{eq;9.16.10}) is strict,
we have only to consider the case $h\geq 0$.
We assume that $\bigl\{a\,\big|\,x_a\not\in F_b^{(1)}\bigr\}$,
and we will derive the contradiction.
We put $a_0:=\max\bigl\{a\,\big|\,x_a\not\in F_b^{(1)}\bigr\}$.
We may assume that
$x=\sum_{a\leq a_0} x_a$.
Then we have the following:
\[
 N^{h+a_0}x=N^{h+a_0}x_{a_0}
\in F^{(1)}_b.
\]
Due to the strictness of (\ref{eq;9.16.10}),
there exists the element $y\in F_b^{(1)}\cap \Gr^{W(N)}_{h+2a_0}$
such that
$N^{h+2a_0}y=N^{h+a_0}x_{a_0}$.
Due to the property of the primitive decomposition,
we have the equality $N^{a_0}y=x_{a_0}\in F_b^{(1)}$,
which contradicts $x_{a_0}\not\in F_b^{(1)}$.
Thus we are done.
\hfill\qed

\subsubsection{Sequential compatibility}

\begin{df} \label{df;9.16.20}
Let $V$ be a finite dimensional vector space.
Let $\nbign_1,\ldots,\nbign_m$ be commuting tuple
of nilpotent maps of $V$.
Let $\lefttop{1}F,\ldots , \lefttop{l}F$
be filtrations on $V$.
$(\nbign_1,\ldots,\nbign_m;\lefttop{1}F,\ldots,\lefttop{l}F)$
is called sequentially compatible,
if the following conditions hold:
\begin{enumerate}
\item \label{9.16.21}
 $\nbign_j$ preserves the filtration $\lefttop{i}F$.
\item \label{9.16.22}
 We put $\nbign(\jbar):=\sum_{i\leq j}\nbign_i$.
 We denote the weight filtration of $\nbign(\jbar)$
 by $W(\jbar)$.
Then the filtrations
 $(W(\itibar),W(\nibar),\ldots,W(\mbar),\lefttop{1}F,\ldots,\lefttop{l}F)$
are compatible,
 in the sense of Definition {\rm \ref{df;b11.12.1}}.
\item \label{9.16.23}
 On the associated graded vector spaces
 $\lefttop{I}\Gr^F_{\veca}$ $(\veca\in\real^I)$,
 we have
 the induced filtrations $\lefttop{I}W(\jbar)$
 and the induced nilpotent maps $\lefttop{I}\nbign(\jbar)$.
 Then we have $W(\lefttop{I}\nbign(\jbar))=\lefttop{I}W(\jbar)$.
\item \label{9.16.24}
 On the associated graded vector spaces
 $\lefttop{I}\Gr^F_{\veca}$ $(\veca\in\real^I)$,
 the tuple of nilpotent maps
 $\lefttop{I}\nbign_1,\ldots,\lefttop{I}\nbign_m$
 are sequentially compatible
 (see Definition {\rm 2.7} in {\rm\cite{mochi}}).
\hfill\qed
\end{enumerate}
\end{df}

\begin{lem} \label{lem;9.16.26}
Assume a tuple
$\nbigs=(\nbign_1,\ldots,\nbign_m,\lefttop{1}F,\ldots,\lefttop{l}F)$
is sequentially compatible.
\begin{description}
\item[(A)]
The induced tuple
 $\nbigs(I):=\bigl(\lefttop{I}\nbign_1,\ldots,\lefttop{I}\nbign_m,
  \lefttop{i}F\,\,(i\in I^c)\bigr)$ on $\lefttop{I}\Gr^F$
 is sequentially compatible.
\item[(B)]
 The induced tuple
 $\nbigs^{(1)}:=
 \bigl(\nbign^{(1)}_2,\ldots,\nbign^{(1)}_m,
 \lefttop{i}F\,\,\,(i\in \lbar)\bigr)$ on $\Gr^{W(\itibar)}$
 is sequentially compatible.
 Moreover, it is compatible with the primitive decomposition
 of $\Gr^{W(\itibar)}$.
\end{description}
\end{lem}
\pf
First we see the claim $(A)$.
The conditions \ref{9.16.21}, \ref{9.16.23} and \ref{9.16.24}
in Definition \ref{df;9.16.20} are clear.
The condition \ref{9.16.22} and \ref{9.16.23} for
$\nbigs$ implies the condition \ref{9.16.22} for $\nbigs(I)$.

Let us show the claim $(B)$.
The condition \ref{9.16.21} is clear.
Note that we have
$W(\nbign^{(1)}(\jbar))_h=
 W^{(1)}(\jbar)_{h+a}$ on $\Gr^{W(\itibar)}_a$
due to the condition \ref{9.16.24} for $\nbigs$.
Then the condition \ref{9.16.22} for $\nbigs^{(1)}$
follows from the same condition for $\nbigs$.
We have the following on
$\Gr^{W^F(\itibar)}_a\lefttop{I}\Gr^F
\simeq
\lefttop{I}\Gr^{F^{(1)}}\Gr^{W(\itibar)}$:
\[
 \lefttop{I}W(\jbar)^{(1)}_{h+a}
=W\bigl(\lefttop{I}\nbign^{(1)}(\jbar)\bigr)_h.
\]
Thus the conditions \ref{9.16.23} and \ref{9.16.24} for $\nbigs^{(1)}$
follows from the conditions \ref{9.16.22}, \ref{9.16.23}
and \ref{9.16.24}
for $\nbigs$.

The compatibility with the primitive decomposition follows
from Lemma \ref{lem;9.16.25}.
\hfill\qed

\begin{cor}
Assume that a tuple
$(\nbign_1,\ldots,\nbign_m,\lefttop{1}F,\ldots,\lefttop{l}F)$
is sequentially compatible.
Then the induced tuple
 $(\nbign^{(1)}_2,\ldots,\nbign^{(1)}_m,
 \lefttop{1}F^{(1)},\ldots,\lefttop{l}F^{(1)})$
on the primitive part $P\Gr^{W(\itibar)}$
is sequentially compatible.
\end{cor}
\pf
It follows from the claim $(B)$ in Lemma \ref{lem;9.16.26}.
\hfill\qed

\begin{df} \label{df;9.16.30}
Let $X$ be a complex manifold or scheme.
Let $V$ be a vector bundle over $X$.
Let $\nbign_1,\ldots,\nbign_m$ be nilpotent maps on $V$.
Let $\lefttop{1}F,\ldots,\lefttop{l}F$ be filtrations of $V$
in the category of vector bundles.

A tuple $(\nbign_1,\ldots,\nbign_m,\lefttop{1}F,\ldots,\lefttop{l}F)$
is called sequentially compatible,
if the following conditions hold:
\begin{enumerate}
\item \label{8.28.1}
 For any point $P\in X$,
 the tuple
 $(\nbign_{1\,|\,P},\ldots,\nbign_{m\,|\,P},
 \lefttop{1}F_{|P},\ldots,\lefttop{l}F_{|P} )$
is sequentially compatible,
in the sense of Definition {\rm\ref{df;9.16.20}}.
\item \label{8.28.2}
 The family
 $\Bigl\{
 \bigcap_{j=1}^m W(\jbar)_{h_j\,|\,P}\cap \bigcap_{i=1}^lF_{a_i\,|\,P}
 \,\Big|\, P\in X
 \Bigr\}$
 forms a vector bundle over $X$ for each $(\vech,\veca)$.
\hfill\qed
\end{enumerate}
\end{df}

%% file: a51.1.tex

\subsubsection{Compatibility of the nilpotent maps defined on the
   associated graded bundles}

Let $V$ be a vector bundle over $X$.
Let $\vecE=\big(\lefttop{s}\EE\,\big|\,s\in S\big)$ be a compatible
tuple of the decompositions.
Let $\vecF=\bigl(\lefttop{i}F\,\big|\,i\in \lbar\bigr)$ be a compatible
tuple of the filtrations.
Let us consider a tuple of nilpotent maps
$\vecN=(N_i)$,
where $N_i$ are defined on $\lefttop{\ibar}\Gr^{F}(V)$.
In that case,
we have the induced nilpotent maps
$N^{\mbar}_i$ on $\lefttop{\mbar}\Gr^F$ for any $i\leq m$.
We denote the tuple $(N_1^{\mbar},\ldots,N^{\mbar}_m)$
by $\vecN^{\mbar}$.

\begin{df} \label{df;10.9.1}
The tuple $(\vecE,\vecF,\vecN)$ is called sequentially compatible,
if the following holds:
\begin{itemize}
\item $N_i$ preserves the decomposition $\vecE^{\ibar}$.
\item The tuple $(\vecN^{\mbar},\vecF^{\mbar})$ 
 on $\lefttop{\mbar}\Gr^F$ is
 sequentially compatible for any $m\leq l$,
 in the sense of Definition {\rm\ref{df;9.16.30}}.
 (See the subsubsection {\rm\ref{subsubsection;a11.12.1}}
 for $\vecF^{\mbar}$.)
\hfill\qed
\end{itemize}
\end{df}

When we are given a compatible tuple $(\vecE,\vecF,\vecN)$,
we obtain the weight filtration
$W(\mbar)$ on $\lefttop{\mbar}\Gr^F$,
and thus the tuple $\vecW=\bigl(W(\mbar)\,\big|\,m\in\lbar\bigr)$.
Hence we obtain the induced tuple
$(\vecE,\vecF,\vecW)$
as in the subsubsection \ref{subsubsection;a12.9.2}.

\begin{lem}
If $(\vecE,\vecF,\vecN)$ is compatible
in the sense of Definition {\rm\ref{df;10.9.1}},
then the induced tuple $(\vecE,\vecF,\vecW)$ is compatible
in the sense of Definition {\rm\ref{df;10.7.16}}.
\end{lem}
\pf
It immediately follows from the definition.
\hfill\qed

\begin{df} \label{df;10.9.3}
A splitting of a sequentially compatible tuple $(\vecE,\vecF,\vecN)$
is defined to be a splitting of the induced tuple
$(\vecE,\vecF,\vecW)$,
in the sense of Definition {\rm\ref{df;a11.12.2}}.
\hfill\qed
\end{df}

%% file: a51.2.tex

\subsubsection{Strongly sequential compatibility of nilpotent maps}

Let $V$ be a vector bundle on $X$.
Let $\vecE=\bigl(\lefttop{s}\EE\,\big|\,s\in S\bigr)$ be a compatible
tuple of $G$-decompositions.
Let $\vecN=\bigl(N_i\,|\,i=1,\ldots,l\bigr)$ be a tuple of nilpotent maps
of $V$.

\begin{df} \label{df;b11.12.10}
The tuple $(\vecE,\vecN)$ is called strongly sequentially compatible,
if the following holds:
\begin{itemize}
\item
 Each nilpotent map $N_i$ preserves the decompositions $\lefttop{s}E$.
\item
 The tuple of nilpotent maps $\vecN$ is strongly sequentially compatible
 (see Definition {\rm 2.9} in {\rm\cite{mochi}}).
\hfill\qed
\end{itemize}
\end{df}

\begin{df}\label{df;10.9.5}
Let $(\vecE,\vecN)$ be a strongly sequentially compatible tuple.
A splitting of $(\vecE,\vecN)$ is defined to be
decompositions
\[
 V=\bigoplus_{\substack{
 \vecalpha\in\cnum^S,\\
 \vech\in\seisuu^m
 }}
 U_{\vecalpha,\vech},
\quad\quad
 U_{\vecalpha,\vech}
=\bigoplus_{a\geq 0} P_{|q_1(\vech)|+2a}U_{\vecalpha,\vech}.
\]
They are assumed to satisfy the following conditions:
\begin{itemize}
\item
 $\nbign(\itibar)\bigl(
 P_h U_{\vecalpha,\vech}\bigr)=
 P_{h}U_{\vecalpha,\vech-2\vecdelta}$.
\item
 $\nbign(\itibar)^{h+1}\bigl(
 P_hU_{\vecalpha,\vech}\bigr)=0$.
\item
 $\nbign(\itibar)^{h}:
 P_hU_{\vecalpha,\vech}\lrarr
 P_hU_{\vecalpha,\vech-2h\vecdelta}
 $ is isomorphic
if $q_1(\vech)=h\geq 0$.
\end{itemize}
Here we put $\vecdelta:=(1,\ldots,1)\in\seisuu^m$.
\hfill\qed
\end{df}

For any subset $S'\subset S$,
the tuple $\bigl(\lefttop{s}\EE\,\big|\,s\in S'\bigr)$
is denoted by $q_{S'}(\vecE)$.
For any integer $m'$ such that $1\leq m'\leq m$,
the tuple $(N_i\,|\,i=1,\ldots,m')$ is denoted by
$q_{\underline{m}'}(\vecN)$.

The following lemma is easy to see.
\begin{lem}\mbox{{}}
\begin{itemize}
\item
For any subset $S'\subset S$ and for any integer $m'$ such that
$1\leq m'\leq m$,
the induced tuple
$\bigl(q_{S'}(\vecE),q_{\mbar'}(\vecN)\bigr)$ 
is compatible.
\item
Assume $(U_{\vecalpha,\vech},P_hU_{\vecalpha,\vech}\,|\,
 \vecalpha\in\cnum^S,\,\vech\in\seisuu^m)$
is a splitting of $(\vecE,\vecN)$.
For any elements $\vecbeta\in\cnum^{S'}$ and $\veck\in\seisuu^{m'}$,
we put as follows:
\[
 U'_{\vecbeta,\veck}
:=
 \bigoplus_{\substack{
 q_{S'}(\vecalpha)=\vecbeta\\
 q_{\mbar'}(\vech)=\veck
 }}
 U'_{\vecalpha,\vech},
\quad\quad
 P_{h}U_{\vecbeta,\veck}
:=\bigoplus_{\substack{
 q_{S'}(\vecalpha)=\vecbeta,\\
 q_{\mbar'}(\vech)=\veck
}}
 P_hU_{\vecalpha,\vech}.
\]
Then 
$(U'_{\vecbeta,\veck},P_{h}U'_{\vecbeta,\veck})$ is a splitting.
\hfill\qed
\end{itemize}
\end{lem}

\begin{rem}
We can consider the strongly sequential compatibility
of the tuple $\bigl(\vecE,\vecF,\vecW\bigr)$.
However it is rather complicated, and it is useless for our purpose.
So we do not omit it.
\hfill\qed
\end{rem}

%% file: a50.3.tex

\subsubsection{A lift of a splitting of a compatible filtrations}

We put $X:=\Delta^n$, $D_i:=\{z_i=0\}$ and
$D:=\bigcup_{i=1}^pD_i$ for some $p\leq n$.
Let $\vecc$ be an element of $\seisuu^m_{>0}$,
and $\rho$ be a $\mu_{\vecc}$-action on $X$.
Let $V$ be a $\mu_{\vecc}$-bundle,
and let
$\vecF:=\bigl(\lefttop{i}F\,\big|\,i\in I \bigr)$
be a compatible tuple of equivariant filtrations of $V$.

\begin{lem} \label{lem;10.7.21}
Assume that we have an equivariant splitting 
$\big(U^D_{\veceta}\,\big|\,\veceta\in\real^I\bigr)$
of $\vecF_{|D}$,
i.e., we have a decomposition of $V_{|D}$ as follows:
\[
 V_{|D}=
 \bigoplus_{\veceta\in\real^l}
 U^D_{\veceta},
\quad\quad
 \lefttop{I}F_{\vecrho|D}
=\bigoplus_{\veceta\leq\vecrho}U_{\veceta}.
\]
Then there exists an equivariant splitting
$\bigl(U_{\veceta}\,\big|\,\veceta\in\real^I
 \bigr)$ of
$\vecF$ on a neighbourhood of $O$
such that $U_{\veceta|D}=U^D_{\veceta}$.
\end{lem}
\pf
We have the equivariant surjection
$\lefttop{I}F_{\veceta}\lrarr
 \lefttop{I}\Gr^F_{\veceta}$
over $X$.
On the divisor $D$, 
we have the subbundle
$U_{\veceta}^D\subset \lefttop{I}F_{\veceta\,|\,D}$
is given.
Then we may extend it to 
the subbundle $\tilde{U}_{\veceta}\subset \lefttop{I}F_{\veceta}$
on a neighbourhood of $O$,
by using Lemma \ref{lem;10.7.15}.
\hfill\qed

\subsubsection{A lift of equivariant splitting}

We put $X=\Delta^n$, $D_i:=\{z_i=0\}$
and $D=\bigcup_{i=1}^p D_i$.
We take a $\mu_{\vecc}$-action on $X$.
Let $V$ be a $\mu_{\vecc}$-equivariant bundle over $X$.
Let $\vecF=\big(\lefttop{i}F\,\big|\,i\in\lbar\big)$ be
a compatible tuple of equivariant filtrations of a vector bundle $V$.

We have the vector bundle $\lefttop{l}\Gr^{F}$,
and we have the tuple of the induced filtrations
$\vecF^{(1)}_a:=
 \bigl(\lefttop{l}\Gr^F_a\cap\lefttop{i}F^{(1)}\,
 \big|\,i=1,\ldots,l-1\bigr)$
on $\lefttop{l}\Gr^F_a$ for any $a\in\real$.

We have the natural isomorphism
$\real^l\simeq \real^{l-1}\times\real$.
We use the notation $(\veceta,a)$ to denote an element of
$\real^l$, where $\veceta$ and $a$ denote elements
of $\real^{l-1}$ and $\real$ respectively.

Let 
$\bigl(U^D_{(\veceta,a)}\,\big|\,(\veceta,a)\in\real^l \bigr)$
be a splitting of the restriction $\vecF_{|D}$.
For any $a\in\real$,
let $\bigl(U^{(1)}_{(\veceta,a)}\,\big|\,\veceta\in\real^{l-1} \bigr)$
be a splitting of
the filtration $\vecF^{(1)}_a$ of $\lefttop{l}\Gr^F_a$.

We have the naturally defined surjection:
\[
 \pi_a:\lefttop{\lbar}F_{(\veceta,a)}
 \lrarr
 \lefttop{l}\Gr^F_a\cap\,
 \lefttop{\underline{l-1}}F^{(1)}_{\veceta}.
\]
We assume 
$\pi_{a|D}\bigl(U^{D}_{(\veceta,a)}\bigr)
 =U^{(1)}_{(\veceta,a)}$.
\begin{lem} \label{lem;10.7.30}
We have a splitting 
$\bigl(U_{(\veceta,a)}\,
 \big|\,\veceta\in\real^{l-1} \bigr)$
of $\vecF$ satisfying the following:
\[
 U_{(\veceta,a)\,|\,D}=U^D_{(\veceta,a)},
\quad\quad
 \pi_a\bigl(
 U_{(\veceta,a)}
 \bigr)
=U^{(1)}_{(\veceta,a)}.
\]
Such a splitting
$\bigl(U_{(\veceta,a)}\,
 \big|\,\veceta\in\real^{l-1} \bigr)$
is called a lift of
$\bigl(U^{(1)}_{(\veceta,a)}\,\big|\,
 (\veceta,a)\in\real^l \bigr)$
extending
$\bigl(U^D_{(\veceta,a)}\,
 \big|\,\veceta\in\real^{l-1} \bigr)$.
\end{lem}
\pf
We have only to take an equivariant lift of
$U^{(1)}_{(\veceta,a)}$ extending
$U^D_{(\veceta,a)}$
by applying Lemma \ref{lem;10.7.15}.
\hfill\qed

\subsubsection{Extension of a splitting of compatible tuple
 $(\vecE,\vecF,\vecW)$}

We put $X=\Delta^n$, $D_i:=\{z_i=0\}$
and $D=\bigcup_{i=1}^p D_i$ for some $p\leq n$.
We take a $\mu_{\vecc}$-action on $X$.
Let $V$ be an equivariant holomorphic bundle on $X$.
Let $(\vecE,\vecF,\vecW)$ be a compatible tuple
as in Definition \ref{df;10.7.16}.
Assume that we are given a splitting
$\bigl(
 U^D_{\vecu}\,\big|\,
 \vecu\in\cnum^S\times\real^l\times\seisuu^{l_1}
 \bigr)$
of the restriction $(\vecE,\vecF,\vecW)_{|D}$
(see Definition \ref{df;a11.12.2}).
We would like to extend it to
a splitting of $(\vecE,\vecF,\vecW)$ on a neighbourhood of $O$.

Let $l'$ be integer such that $1\leq l'\leq l_1$.
We have the projection $\real^{l}\lrarr\real^{l'}$,
taking the first $l'$-components.
Then we obtain the projections
$\pi_{1,l'}:
  \cnum^S\times\real^l\times\seisuu^{l_1}\lrarr
  \cnum^S\times\real^{l'}$
and
$\pi_{3,l'}:\cnum^S\times\real^{l}\times\seisuu^{l'}
  \lrarr\cnum^S\times\real^{l'}$.
We have the projection
$\seisuu^{l_1}\lrarr\seisuu^{l'}$,
taking the first $l'$-components.
Then we obtain the projection
$\pi_{2,l'}:\cnum^S\times\real^{l}\times\seisuu^{l_1}\lrarr
 \cnum^S\times\real^l\times\seisuu^{l'}$.

For any element $\vecu\in\cnum^S\times\real^l\times\seisuu^{l_1}$,
we put $\pi_{1,l'}(\vecu)=(\vecalpha,\veca)$,
and we put as follows:
\[
 \lefttop{\lbar'}U^D_{\vecu}:=
 \pi_{\veca}\bigl(
 U^D_{\vecu}
 \bigr)
\subset
\bigl(
 \lefttop{\lbar'}\Gr^F_{\veca}\cap
 \lefttop{S}\EE_{\vecalpha}
\bigr)_{|D}.
\]
For an element $\vecv\in\cnum^S\times\real^l\times\seisuu^{l'}$,
we put as follows:
\[
 \lefttop{\lbar'}B^D_{\vecv}
:=\bigoplus_{\pi_{2,l'}(\vecu)=\vecv}
 \lefttop{\lbar'}U^D_{\vecu}
\]
Then we have the decomposition
\begin{equation}\label{eq;a12.9.3}
 \bigl(
 \lefttop{\lbar'}\Gr^F_{\veca}
 \cap \lefttop{S}\EE_{\vecalpha}
 \bigr)_{|D}
=\bigoplus_{\pi_{3,l'}(\vecv)=(\veca,\vecalpha)}
 \lefttop{\lbar'}B_{\vecv}^D.
\end{equation}
Then
the decomposition (\ref{eq;a12.9.3})
gives a splitting of
$(\vecF^{\lbar'},\vecW^{\lbar'})_{|D}$,
by definition.

\begin{lem} \label{lem;a12.7.1}
Let $l'$ be any integer such that $0\leq l'\leq l_1$,
and let $\vecu$ be any element of
$\cnum^S\times\real^l\times\seisuu^{l_1}$.
We can take subbundles
$\lefttop{\lbar'}U'_{\vecu}$
of $\lefttop{\lbar'}\Gr^F_{\pi_{1,l'}(\vecu)}$
satisfying the following:
\begin{enumerate}
\item
We have the following decomposition,
for any element $(\vecalpha,\veca)\in\cnum^S\times\real^{\,l'}$:
\[
 \lefttop{S}\EE_{\vecalpha}\lefttop{\lbar'}\Gr^F_{\veca}
=\bigoplus_{\pi_{1,l'}(\vecu)=(\vecalpha,\veca)}
 \lefttop{\lbar'}U'_{\vecu}.
\]
\item
 For any element $\vecv\in\cnum^S\times\real^{l'}\times\seisuu^{l'}$,
 we put as follows:
\begin{equation}\label{eq;a11.12.3}
 \lefttop{\lbar'}B'_{\vecv}:=
 \bigoplus_{\pi_{2,l'}(\vecu)=\vecv}
 \lefttop{\lbar'}U'_{\vecu}.
\end{equation}
Due to the first condition,
we have the decomposition
$\lefttop{S}\EE_{\vecalpha}\lefttop{\lbar'}\Gr^F_{\veca}=
 \bigoplus_{\pi_{3,l'}(\vecv)=(\vecalpha,\veca)}
  \lefttop{\lbar'}B'_{\vecv}$
for any $\veca\in\real^{l'}$.
Then the decomposition gives a splitting
of $(\vecF^{\lbar'},\vecW^{\lbar'})$.
\item
 We have $\lefttop{\lbar'}U_{\vecu\,|\,D}=\lefttop{\lbar'}U^D_{\vecu}$.
 In particular,
 we have $\lefttop{\lbar'}B_{\vecv\,|\,D}=\lefttop{\lbar'}B^D_{\vecv}$.
\end{enumerate}
\end{lem}
\pf
We can restrict our attention to each component
$\lefttop{S}\EE_{\vecalpha}$
 $(\vecalpha\in\cnum^S)$.
Thus we may assume that $|S|=1$.
In the following, we omit to denote $\lefttop{S}\EE_{\vecalpha}$.
We use a descending induction on $l'$.

\noindent
1. First we construct such decomposition in the case $l'=l_1$.
We have the following decomposition,
for any $\veca\in\real^{l_1}$:
\begin{equation} \label{eq;10.7.20}
 \lefttop{\lbariti}\Gr^F_{\veca\,|\,D}
=\bigoplus_{\pi_{1,l_1}(\vecu)=\veca}
 \lefttop{\underline{l_1}}
 U^D_{\vecu}.
\end{equation}
Then it is a splitting of the compatible tuple
$(\vecF^{\lbariti},\vecW^{\lbariti})_{|D}$.
By using Lemma \ref{lem;10.7.21},
we can obtain a splitting
$\bigl(
 \lefttop{\lbariti}U'_{\vecu}\,\big|\,
 \vecu\in \real^l\times\seisuu^{l_1},\,\,\pi_{1,l_1}(\vecu)=\veca
 \bigr)$
extending the splitting
$\bigl(
 \lefttop{\lbariti}U^{D}_{\vecu}\,\big|\,
 \vecu\in\real^l\times\seisuu^{l_1}\,\,\pi_{1,l_1}(\vecu)=\veca
 \bigr)$.
Thus the claim has been shown in the case $l'=l_1$.

\vspace{.1in}
\noindent
2. We assume that we have already obtained
the vector subbundles
$\lefttop{\underline{l'+1}}U'_{\vecu}$,
and we will construct $\lefttop{\lbar'}U'_{\vecu}$.

We have the projection $\seisuu^{l'+1}\lrarr \seisuu^{l'}$
by taking the first $l'$-components.
We denote the induced projection 
$\real^{l}\times\seisuu^{l'+1}\lrarr\real^{l}\times\seisuu^{l'}$
by $\pi_{4,l'}$.
Similarly the projection $\real^{l}\lrarr\real^{l'+1}$ induces
the morphism $\pi_{5,l'}:\real^l\times\seisuu^{l'}\lrarr\real^{l'+1}$.

We take $\lefttop{\underline{l'+1}}B'_{\vecv}$
$(\vecv\in\real^l\times \seisuu^{l'+1})$
as in (\ref{eq;a11.12.3}).
Then we put as follows,
for any element $\vecv\in\real^l\times\seisuu^{l'}$:
\[
 \lefttop{\underline{l'+1}}B''_{\vecv}
=\bigoplus_{\pi_{4,l'}(\vecv_1)=\vecv}
 \lefttop{\underline{l'+1}}B'_{\vecv_1}
\subset
 \lefttop{l'+1}\Gr^F_{\pi_{5,l'}(\vecv)}.
\]
Then we have the decomposition,
for any $(\veca,b)\in\real^{l'+1}=\real^{l'}\times\real$:
\[
 \lefttop{\underline{l'+1}}\Gr^{F}_{(\veca,b)}
=\bigoplus_{\pi_{5,l'}(\vecv)=(\veca,b)}
 \lefttop{\underline{l'+1}}B''_{\vecv}.
\]
Recall that
$\bigl(
 \lefttop{\underline{l'+1}}B'_{\vecv}\,\big|\,
 \pi_{2,l'+1}(\vecv)=(\veca,b)
 \bigr)$
gives the splitting of the following tuple of the filtrations
on $\lefttop{\underline{l'+1}}\Gr^F_{(\veca,b)}$:
\[
 \bigl(
 \vecF^{\underline{l'+1}},\vecW^{\underline{l'+1}}
 \bigr)
=\bigl(
 \lefttop{l'+2}F^{\underline{l'+1}},\ldots,
 \lefttop{l}F^{\underline{l'+1}}\,;\,
  W^{\underline{l'+1}}(\itibar),\ldots,
  W^{\underline{l'+1}}(\underline{l'+1})
 \bigr).
\]
Hence the decomposition
$\bigl(
 \lefttop{\underline{l'+1}}B''_{\vecv}\,\big|\,
 \pi_{5,l'+1}(\vecv)=(\veca,b)
 \bigr)$
gives the splitting of the following compatible tuple of the
filtrations
on $\lefttop{\underline{l'+1}}\Gr^F_{(\veca,b)}$:
\[
 \bigl(
 \lefttop{l'+2}F^{\underline{l'+1}},
 \ldots,\lefttop{l}F^{\underline{l'+1}}\,;\,
 W^{\underline{l'+1}}(\itibar),
\ldots,
 W^{\underline{l'+1}}(\underline{l'})
 \bigr).
\]

We have the compatible tuple of filtrations
$(\vecF^{\lbar'},\vecW^{\lbar'})$
on $\lefttop{\lbar'}\Gr^F_{\veca}$.
We have the isomorphism
$\lefttop{l'+1}\Gr_b\lefttop{\lbar'}\Gr^F_{\veca}
 \simeq
 \lefttop{\underline{l'+1}}\Gr^F_{(\veca,b)}$.

Then we can take a lift
$\bigl(
 \lefttop{\underline{l'}}B'_{\vecv}\,\big|\,
 \vecv\in\real^{l}\times\real^{l'_1}
 \bigr)$
of
$\bigl(
 \lefttop{\underline{l'+1}}B''_{\vecv}\,\big|\,
 \vecv\in\real^{l}\times\real^{l'_1}
 \bigr)$
extending 
$\bigl(
 \lefttop{\underline{l'}}B^D_{\vecv}\,\big|\,
 \vecv\in\real^{l}\times\real^{l'_1}
 \bigr)$,
by applying Lemma \ref{lem;10.7.30}.

Then we have the naturally defined isomorphism
$\lefttop{\underline{l'}}B'_{\vecv}
\lrarr
 \lefttop{\underline{l'+1}}B''_{\vecv}$
for any $\vecv\in\real^l\times\seisuu^{l'}$.
For any $\vecu\in\real^l\times\seisuu^{l_1}$
such that $\pi_{2,l'}(\vecu)=\vecv$,
we can lift $\lefttop{\lbar'}U'_{\vecu}$
of $\lefttop{\underline{l'+1}}U'_{\vecu}$
extending $\lefttop{\lbar'}U_{\vecu}^D$,
by applying Lemma \ref{lem;10.7.15}.
Thus the induction can proceed.
\hfill\qed

\begin{cor} \label{cor;10.9.2}
Assume that we are given a splitting
$\bigl(
 U^D_{\vecu}\,\big|\,
 \vecu\in\cnum^S\times\real^l\times\seisuu^{l_1}
 \bigr)$
of the restriction $(\vecE,\vecF,\vecW)_{|D}$.
Then we have a splitting
$\bigl(
 U_{\vecu}\,\big|\,
 \vecu\in\cnum^S\times\real^l\times\seisuu^{l_1}
 \bigr)$
of $(\vecE,\vecF,\vecW)$
defined around the origin $O$,
such that
$U_{\vecu\,|\,D}=U^D_{\vecu}$.
\end{cor}
\pf
We have only to put
$U_{\vecu}:=\lefttop{\underline{0}}U_{\vecu}$.
\hfill\qed

%% file: a51.3.tex

\subsubsection{Extension of splitting of 
sequentially compatible tuple $(\vecE,\vecF,\vecN)$}

We put $X=\Delta^n$, $D_i:=\{z_i=0\}$
and $D=\bigcup_{i=1}^p D_i$ for some $p\leq n$.
We take a $\mu_{\vecc}$-action on $X$.
Let $V$ be an equivariant holomorphic bundle on $X$.
Let $(\vecE,\vecF,\vecN)$ be a sequentially compatible tuple
as in Definition \ref{df;10.9.1}.

\begin{lem}
Assume that we are given a splitting
$\bigl(
 U^D_{\vecu}\,\big|\,
 \vecu\in\cnum^S\times\real^l\times\seisuu^{l_1}
 \bigr)$
of the restriction $(\vecE,\vecF,\vecN)_{|D}$.
Then we have a splitting
$\bigl(
 U_{\vecu}\,\big|\,
 \vecu\in\cnum^S\times\real^l\times\seisuu^{l_1}
 \bigr)$
of $(\vecE,\vecF,\vecN)$,
such that
$U_{\vecu\,|\,D}=U^D_{\vecu}$.
\end{lem}
\pf
It immediately follows from
Definition \ref{df;10.9.3} and Corollary \ref{cor;10.9.2}.
\hfill\qed

%% file: a51.4.tex

\subsubsection{Extension of splitting of 
strongly sequentially compatible tuple $(\vecE,\vecN)$}

We put $X=\Delta^n$, $D_i:=\{z_i=0\}$
and $D=\bigcup_{i=1}^p D_i$ for some $p\leq n$.
We take a $\mu_{\vecc}$-action on $X$.
Let $V$ be an equivariant holomorphic bundle on $X$.
Let $(\vecE,\vecN)$ be a strongly sequentially compatible tuple
as in Definition \ref{df;10.9.5}.

\begin{lem}
Assume that we are given a splitting
$\bigl(
 U^D_{\vecu},\,P_hU^D_{\vecu}\,\big|\,
 \vecu\in\cnum^S\times\seisuu^{l_1}, h\in\seisuu_{\geq\,0}
 \bigr)$
of the restriction $(\vecE,\vecN)_{|D}$.
Then we have a splitting
$\bigl(
 U_{\vecu},\,P_hU_{\vecu}\,\big|\,
 \vecu\in\cnum^S\times\seisuu^{l_1}, h\in\seisuu_{\geq\,0}
 \bigr)$
of $(\vecE,\vecN)$
such that
$U_{\vecu\,|\,D}=U^D_{\vecu}$ and
$P_hU_{\vecu\,|\,D}=P_hU_{\vecu}^D$.
\end{lem}
\pf
We have only to consider each component $\lefttop{S}\EE_{\vecalpha}$.
Thus we may assume $|S|=1$.
In the following, we omit to denote $\vecalpha$.

We put $V^{(1)}:=P_hGr_h^{W(N_1)}(V)$.
We have the surjection
$\pi_h:\ker\bigl(N_1^{h+1}\bigr)\lrarr P_h\Gr^{W(N_1)}_h(V)$.
The image of $P_hU^D_{\vech}$ via the restriction of the morphism
$\pi_{h\,|\,D}$ is denoted by $P_hU^{\prime\,D}_{\vech}$.

We obtain the induced sequentially compatible tuple
$\vecN^{(1)}:=(N_2^{(1)},\ldots,N_{l_1}^{(1)})$
on $V^{(1)}$,
and we obtain the splitting
$\bigl(P_hU^{\prime\,D}_{\vech}\,\big|\,
 \vech\in\seisuu^{l_1},q_1(\vech)=h\bigr)$.
We can extend it to 
a splitting
$\bigl(P_hU'_{\vech}\,\big|\,
 \vech\in\seisuu^{l_1},q_1(\vech)=h
 \bigr)$
by using Lemma \ref{lem;10.7.21}.
We can take an equivariant lift $P_hU_{\vech}$
of $P_hU'_{\vech}$ extending $P_hU^D_{\vech}$.
We have only to put as follows,
for any $a\in\seisuu_{\geq\,0}$ and for
any $\vech\in\seisuu^{l_1}$:
\[
 P_{|q_1(\vech)|+2a}U_{\vech}
=\left\{
 \begin{array}{ll}
 N_1^{a}\bigl(
 P_{q_1(\vech)+2a  }U_{\vech+2a\vecdelta}
 \bigr),
 & (q_1(\vech)\geq 0),\\
 \mbox{{}}\\
 N_1^{a+|q_1(\vech)|}\bigl(
 P_{q_1(\vech)+2a  }U_{\vech+2a\vecdelta}
 \bigr),
 & (q_1(\vech)< 0).
 \end{array}
 \right.
\]
Then $\bigl(P_{h}U_{\vech}\,\
 \big|\,h\in\seisuu_{\geq\,0},\,
 \vech\in\seisuu^{l_1} \bigr)$
gives a desired decomposition.
\hfill\qed

%% file: a52.tex

\subsubsection{Compatibility of filtrations and decompositions}

We put $X=\Delta^n$ and $D=\bigcup_{i=1}^l D_i$.
We take a $\mu_{\vecc}$-action on $X$.
Let $V$ be an equivariant vector bundle over $X$.
Let $\vecE=\big(\lefttop{i}\EE,\,\big|\,i=1,\ldots,l\bigr)$
be a tuple of
equivariant decompositions $\lefttop{i}\EE$ of $V_{|D_i}$.
Let $\vecF=\bigl(\lefttop{i}F\,\big|\,i=1,\ldots,l\bigr)$
be a tuple of equivariant filtrations of $V_{|D_i}$.

Let $I$ be a subset of $\lbar$.
Then we obtain the tuple
$\lefttop{I}\vecE:=\bigl(\lefttop{i}\EE\,\big|\,i\in I \bigr)$
of the equivariant decompositions of $V_{|D_I}$
and the tuple
$\lefttop{I}\vecF:=\bigl(\lefttop{i}F\,\big|\,i\in I
 \bigr)$
of the equivariant filtrations of $V_{|D_I}$.

\begin{df} \label{df;10.11.90}
The tuple $(\vecE,\vecF)$ is called compatible
if $(\lefttop{I}\vecE,\lefttop{I}\vecF)$ are compatible
for any $I\subset\lbar$,
in the sense of Definition {\rm\ref{df;b11.12.3}}.
\hfill\qed
\end{df}

Let $\vecE=\bigl(\lefttop{i}\EE\,\big|\,i\in\lbar\bigr)$
and $\vecF=\bigl(\lefttop{i}F\,\big|\,i\in\lbar\bigr)$
be as above.
Assume that $(\vecE,\vecF)$ is compatible.
Let us consider splittings
$\lefttop{I}\vecU
=\bigl(\lefttop{I}U_{\vecu}\,\big|\,\vecu\in\cnum^I\times\real^I\bigr)$
for any subset $I\subset\lbar$
of $(\lefttop{I}\vecE,\lefttop{I}\vecF)$,
(Definition \ref{df;b11.12.3}).

For any subset $I\subset I'$,
let $q_I$ denote the projection
$\cnum^{I'}\times\real^{I'}\lrarr\cnum^{I}\times\real^I$.

\begin{df}\label{df;b11.12.15}
A tuple of splittings
$\bigl(\lefttop{I}\vecU\,\big|\,I\subset\lbar\bigr)$
is called a splitting of the tuple
$\bigl(\vecE,\vecF\bigr)$,
if $\lefttop{I}U_{\vecu\,|\,D_{I'}}
=\bigoplus_{q_I(\vecu_1)=\vecu}
 \lefttop{I'}U_{\vecu_1}$ hold
for any $I\subset I'$.
\hfill\qed
\end{df}

\subsubsection{Compatibility of $(\vecE,\vecF,\vecW)$}

We put $X=\Delta^n$, $D_i=\{z_i=0\}$ and $D=\bigcup_{i=1}^l D_i$
for some $l\leq n$.
Let $\vecE=\bigl(\lefttop{i}\EE\,\big|\,i\in\lbar\bigr)$
be a tuple of decompositions of $V_{|D_i}$ $(i\in\lbar)$.
Let $\vecF=\bigl(\lefttop{i}F\,\big|\,i\in\lbar\bigr)$
be a tuple of filtrations of $V_{|D_i}$ $(i\in\lbar)$.
We assume that $(\vecE,\vecF)$ is compatible.

Let us consider a tuple
$\vecW=\big(W(\ibar)\,\big|\, i\in \lbariti \big)$,
where $W(\ibar)$ are filtrations of
the vector bundle $\lefttop{\ibar}\Gr^F(V)$ on $D_{\ibar}$.
We have the induced filtrations
$W^{J}(\ibar)$ of $\lefttop{J}\Gr^F(V)$
on $D_{J}$ for any subset $J$ such that $\ibar\subset J$.

For any subset $J\subset\lbar$,
we have the number $m(J)$ determined by the condition
$\underline{m(J)}\subset J$ and $m(J)+1\not\in J$.
We denote the tuple
$\bigl( W^J(\ibar)\,\big|\,i\in\underline{m(J)} \bigr)$
by $\lefttop{J}\vecW$.
Then we obtain the tuple
$(\lefttop{J}\vecE,\lefttop{J}\vecF,\lefttop{J}\vecW)$.

\begin{df} \label{df;b11.12.16}
The tuple $(\vecE,\vecF,\vecW)$ is called compatible,
if the induced tuples 
$(\lefttop{J}\vecE,\lefttop{J}\vecF,\lefttop{J}\vecW)$
on $D_J$
are compatible for any subset $J\subset\lbar$,
in the sense of Definition {\rm\ref{df;10.7.16}}.
\hfill\qed
\end{df}

Let us consider splittings
$\lefttop{J}\vecU=\bigl(
 \lefttop{J}U_{\vecu}\,\big|\,
 \vecu\in \cnum^J\times\real^J\times
 \seisuu^{m(J)} \bigr)$
of the tuple
$\bigl(
 \lefttop{J}\vecE,\lefttop{J}\vecF,\lefttop{J}\vecW
\bigr)$,
in the sense of Definition \ref{df;a11.12.2}.
For any pair of subsets $I\subset I'$ of $\lbar$,
let $q_I$ denote the projection
$\cnum^{I'}\times\real^{I'}\times\seisuu^{m(I')}
\lrarr
 \cnum^{I}\times\real^I\times\seisuu^{m(I)}$.
\begin{df}\label{df;b11.12.17}
A tuple of splittings
$\vecU=\bigl(\lefttop{J}\vecU\,\big|\,J\subset\lbar\bigr)$
is called a splitting of $(\vecE,\vecF,\vecW)$,
if
$\lefttop{I}U_{\vecu\,|\,D_{I'}}=
 \bigoplus_{q_I(\vecu')=\vecu}\lefttop{I'}U_{\vecu'}$
hold for any subset $I\subset I'$
and for any $\vecu\in \cnum^I\times\real^I\times\seisuu^{m(I)}$.
\hfill\qed
\end{df}

\begin{prop} \label{prop;10.9.10}
Let $(\vecE,\vecF,\vecW)$ be as above,
and we assume that it is compatible.
Then there exists a splitting of
$(\vecE,\vecF,\vecW)$,
in the sense of Definition {\rm\ref{df;b11.12.17}}.
\end{prop}
\pf
On $D_{\lbar}$,
we have only to take a splitting 
$\lefttop{\lbar}\vecU$ of the compatible tuple
$(\lefttop{\lbar}\vecE,\lefttop{\lbar}F,\lefttop{\lbar}W)$.

We construct the splittings $\lefttop{I}\vecU$
descending inductively on $|I|$.
We assume that we have the splittings $\lefttop{I'}\vecU$
for any $I'\supsetneq I$.
Then we will construct the splitting $\lefttop{I}\vecU$.

We put $\del D_I:=\bigcup_{I'\supsetneq I}D_{I'}$.
From the given splittings $\lefttop{I'}\vecU$ $(I'\supsetneq I)$,
we obtain the splitting 
$\lefttop{I}\vecU^{\del D_I}$
of the restriction
$(\lefttop{I}\vecE,\lefttop{I}\vecF,\lefttop{I}\vecW)_{|\del D_I}$.
Then we can extend $\lefttop{I}\vecU^{\del D_I}$
to $\lefttop{I}\vecU$,
due to Lemma \ref{lem;a12.7.1}.
Thus the inductive construction can proceed,
and hence we are done.
\hfill\qed

\begin{cor}
Let $(\vecE,\vecF)$ be a compatible tuple
as in Definition {\rm\ref{df;10.11.90}}.
Then we have a splitting of $(\vecE,\vecF)$
in the sense of Definition {\rm\ref{df;b11.12.15}}.
\end{cor}
\pf
We have only to consider Proposition \ref{prop;10.9.10}
in the case where the filtrations $W(\mbar)$ are trivial.
\hfill\qed

\subsubsection{Sequential compatibility}

\label{subsubsection;10.26.123}

We put $X=\Delta^n$, $D_i=\{z_i=0\}$ and $D=\bigcup_{i=1}^l D_i$
for some $l\leq n$.
Let $\vecE=\bigl(\lefttop{i}\EE\,\big|\,i\in\lbar\bigr)$
be a tuple of decompositions of $V_{|D_i}$ $(i\in\lbar)$.
Let $\vecF=\bigl(\lefttop{i}F\,\big|\,i\in\lbar\bigr)$
be a tuple of filtrations of $V_{|D_i}$ $(i\in\lbar)$.
We assume that $(\vecE,\vecF)$ is compatible.

Let us consider a tuple
$\vecN=\big(N_i\,\big|\, i\in \lbariti \big)$,
where $N_i$ are nilpotent maps of
the vector bundle $\lefttop{\ibar}\Gr^F(V)$ on $D_{\ibar}$.
We have the induced filtrations
$N^{J}_i$ of $\lefttop{J}\Gr^F(V)$
on $D_{J}$ for any subset $J$ such that $\ibar\subset J$.

For any subset $J\subset\lbar$,
we have the number $m(J)$ determined by the condition
$\underline{m(J)}\subset J$ and $m(J)+1\not\in J$.
We denote the tuple
$\bigl(N_i\,\big|\,i\in\underline{m(J)} \bigr)$
by $\lefttop{J}\vecN$.
Then we obtain the tuple
$(\lefttop{J}\vecE,\lefttop{J}\vecF,\lefttop{J}\vecN)$.

\begin{df} \label{df;b11.12.20}
The tuple $(\vecE,\vecF,\vecN)$ is called sequentially compatible,
if the induced tuples
$\bigl(\lefttop{J}\vecE,\lefttop{J}\vecF,\lefttop{J}\vecN\bigr)$
are sequentially compatible for any subset $J\subset\lbar$,
in the sense of Definition {\rm\ref{df;10.9.1}}.
\hfill\qed
\end{df}

The nilpotent endomorphism $N(\ibar)$ induces
the filtration $W(\ibar)$ of
the vector bundle $\lefttop{\ibar}\Gr^F(V)$
over $D_{\ibar}$.
Thus we obtain the tuple
$\vecW=\bigl(W(\ibar)\,\big|\,i\in\lbariti\bigr)$.
\begin{lem} \label{lem;c11.12.1}
The tuple $(\vecE,\vecF,\vecW)$ is compatible,
in the sense of Definition {\rm\ref{df;b11.12.16}}.
\end{lem}
\pf
It immediately follows from the definition
of compatibility.
\hfill\qed

\vspace{.1in}

Let us consider splittings
$\lefttop{J}\vecU=\bigl(
 \lefttop{J}U_{\vecu}\,\big|\,
 \vecu\in \cnum^J\times\real^J\times
 \seisuu^{m(J)} \bigr)$
of the tuple
$\bigl(
 \lefttop{J}\vecE,\lefttop{J}\vecF,\lefttop{J}\vecN
\bigr)$,
in the sense of Definition \ref{df;10.9.3}.
For any pair of subsets $I\subset I'$ of $\lbar$,
let $q_I$ denote the projection
$\cnum^{I'}\times\real^{I'}\times\seisuu^{m(I')}
\lrarr
 \cnum^{I}\times\real^I\times\seisuu^{m(I)}$.

\begin{df}\label{df;b11.12.21}
A tuple of splittings
$\vecU=\bigl(\lefttop{J}\vecU\,\big|\,J\subset\lbar\bigr)$
is called a splitting of $(\vecE,\vecF,\vecN)$,
if
$\lefttop{I}U_{\vecu\,|\,D_{I'}}=
 \bigoplus_{q_I(\vecu')=\vecu}\lefttop{I'}U_{\vecu'}$
hold for any subset $I\subset I'$
and for any $\vecu\in \cnum^I\times\real^I\times\seisuu^{m(I)}$.
\hfill\qed
\end{df}

\begin{prop} \label{prop;10.9.11}
Let $(\vecE,\vecF,\vecN)$ be sequentially compatible
in the sense of Definition {\rm\ref{df;b11.12.20}}.
and we assume that it is compatible.
Then there exists a splitting of
$(\vecE,\vecF,\vecN)$,
in the sense of Definition {\rm\ref{df;b11.12.21}}.
\end{prop}
\pf
It immediately follows from
Lemma \ref{lem;c11.12.1} and Proposition \ref{prop;10.9.10}.
\hfill\qed

\begin{cor} \label{cor;10.9.12}
Assume that $(\vecE,\vecF,\vecN)$
is sequentially compatible.
Then there exists a frame $\vecv$
compatible with the decompositions
$\lefttop{i}\EE$, the filtrations $\lefttop{i}F$ on $D_i$,
and the filtrations $W(\mbar)$ on $D_{\mbar}$.
\end{cor}
\pf
We have a splitting of $(\vecE,\lefttop{i}F,\nbign_i)$,
due to Proposition \ref{prop;10.9.11}.
Then we have only to take a frame compatible with
splittings (Lemma \ref{lem;c11.12.2}).
\hfill\qed

\vspace{.1in}
For any frame $\vecv$ compatible with
$(\vecE,\vecF,\vecN)$,
we have the decomposition
$\vecv=\bigcup \vecv_{\veca,\vecalpha,\vech}$:
\[
 \vecv_{\veca,\vecalpha,\vech}
=\Bigl(
 v_i\,\Big|\,
 \vecdeg^{\EE,F}(v_i)=(\vecalpha,\veca),\,
 \deg^{W(\mbar)}(v_i)=h_m
 \Bigr)
\]

Let us consider the special case $N_i=0$ $(i=1,\ldots,l)$.
Then we obtain the following corollary.
\begin{cor} \label{cor;10.11.95}
Let $(\vecE,\vecF)$ be a compatible tuple,
in the sense of Definition {\rm\ref{df;10.11.90}}.
Then there exists an equivariant frame $\vecv$ compatible with
the decompositions $\lefttop{i}\EE$ and the filtrations
$\lefttop{i}F$.
\hfill\qed
\end{cor}

\subsubsection{Strongly sequential compatibility}

We put $X=\Delta^n$, $D_i=\{z_i=0\}$ and $D=\bigcup_{i=1}^l D_i$
for some $l\leq n$.
Let $\vecE=\bigl(\lefttop{i}\EE\,\big|\,i\in\lbar\bigr)$
be a tuple of decompositions of $V_{|D_i}$ $(i\in\lbar)$.
Let $\vecN=\bigl(N_i\,\big|\,i\in\lbar\bigr)$
be a tuple of nilpotent maps of $V_{|D_i}$ $(i\in\lbar)$.

Let $I$ be a subset of $\lbar$.
Then we obtain the tuple
$\lefttop{I}\vecE$ of the decompositions of $V_{|D_I}$.
We also obtain the tuple
$\lefttop{I}\vecN:=\bigl(
 N_{i\,|\,D_I}\,\big|\,
 i\in I
 \bigr)$ of nilpotent maps of $V_{|D_I}$.

\begin{df}\label{df;b11.12.25}
The tuple $(\vecE,\vecN)$ is called strongly sequentially compatible
if $\bigl(\lefttop{I}\vecE,\lefttop{I}\vecN\bigr)$
are strongly sequentially compatible
for any subset $I\subset \lbar$,
in the sense of Definition {\rm\ref{df;b11.12.10}}
\hfill\qed
\end{df}

Let $\vecE$ and $\vecN$ be as above.
We assume that they are strongly sequentially compatible.
Let us consider splittings
$\lefttop{I}\vecU$
of $(\lefttop{I}\vecE,\lefttop{I}\vecN)$,
in the sense of Definition \ref{df;10.9.5}.
For any subset $I\subset I'$,
let $q_I$ denote the projection
$\cnum^{I'}\times\seisuu^{I'}\lrarr\cnum^{I}\times\seisuu^I$.

\begin{df}
A tuple of splittings
$\bigl(\lefttop{I}\vecU\,\big|\,I\subset\lbar\bigr)$
is called a splitting of the tuple
$(\vecE,\vecN)$,
if 
$\lefttop{I}U_{\vecu\,|\,D_{I'}}
=\bigoplus_{q_I(\vecu')=\vecu}\lefttop{I'}U_{\vecu'}$
hold for any $I\subset I'$.
\hfill\qed
\end{df}

\begin{prop}
Let $(\vecE,\vecN)$ be as above,
and we assume that it is compatible.
Then there exists a splitting of $(\vecE,\vecN)$.
\end{prop}
\pf
It can be shown by an argument
similar to the proof of Proposition \ref{prop;10.9.10}.
\hfill\qed

\vspace{.1in}
Let $(\vecE,\vecN)$ be strongly sequentially compatible tuple.
From the nilpotent maps $N(\mbar)=\sum_{i=1}^m N_i$
of $V_{|D_{\mbar}}$,
we have the weight filtrations $W(\mbar)$ of $V_{|D_{\mbar}}$.
We denote the tuple by $\vecW$.

\begin{cor}
We can take a frame $\vecv$ compatible with
$(\vecE,\vecW)$.
\end{cor}
\pf
Similar to Corollary \ref{cor;10.9.12}.
\hfill\qed

\vspace{.1in}

We have the decomposition
$\vecv=\bigcup_{(\veca,\vech)\in \cnum^l\times\seisuu^l}
 \vecv_{(\veca,\vech)}$.

\begin{cor} \label{cor;10.27.26}
We can take a compatible frame $\vecv$ 
satisfying the following:
\begin{itemize}
\item
 We have the decomposition
 $\vecv_{(\veca,\vech)}
 =\bigcup_{a\geq 0}
 P\vecv_{(\veca,\vech),h(a)}$,
 where we put $h(a)=|q_1(\vech)|+2a$.
\item
 $P\vecv_{(\veca,\vech),h}$ consists of
 sections $v_{(\veca,\vech),h},i$
 $\bigl(i=1,\ldots,d(\veca,\vech,h)\bigr)$,
 and the following holds:
\[
 N(\itibar)\bigl(
 v_{(\veca,\vech),h,i}
 \bigr)
=\left\{
\begin{array}{ll}
 v_{(\veca,\vech-2\vecdelta),h,i}, & (-h+2\leq q_1(\vech)\leq h),\\
 \mbox{{}}\\
 0 & (\mbox{\rm otherwise}).
\end{array}
 \right.
\]
\end{itemize}
Such $\vecv$ is called strongly compatible
with $(\vecE,\vecN)$.
\hfill\qed
\end{cor}

%% file: 23.1.tex

\subsubsection{One nilpotent map}
\label{subsubsection;9.20.10}

Let $R$ be a discrete valuation ring,
$K$ be the quotient field,
and $k$ be the residue field.
Let $V$ be a free $R$-module of finite rank.
Assume that we are given the following data, in this subsubsection.
\begin{condition}\label{condition;9.16.35}
\mbox{{}}
\begin{enumerate}
\item
Let $\lefttop{i}F$ $(i=1,\ldots,l)$ be compatible filtrations
of $V$ in the category of $R$-free modules.
We denote the tuple of the filtrations $(\lefttop{i}F\,|\,i\in I)$
by $\vecF$.
\item
We have the decomposition
$V_{|K}=\bigoplus_a U_a$, which gives a splitting of
$\vecF_{|K}$.
\item
Let $N$ be a nilpotent endomorphism of $V$
preserving $\vecF$.
\item
The restriction
$N_{|K}$ preserves $U_{\vecb}$ for any $\vecb$.
\item \label{9.16.36}
The endomorphism
$N$ induces the nilpotent endomorphism $\lefttop{\lbar}N^{F}_{\veca}$
of $\lefttop{\lbar}\Gr^F_{\veca}(V)$.
Then the conjugacy classes of
$\lefttop{\lbar}N^F_{\veca\,|\,x}$ $(x=k,K)$ are same.
\hfill\qed
\end{enumerate}
\end{condition}

\begin{lem}\mbox{{}} \label{lem;c11.12.20}
Let $m$ be an integer such that $1\leq m\leq l$.
\begin{itemize}
\item
The conjugacy classes
of $\lefttop{\mbar^c}N^F$ on
$\lefttop{\mbar^c}Gr^F$ are constant,
i.e.,
the conjugacy classes of $\lefttop{\mbar^c}N^F_{|x}$
$(x=k,K)$ are same.
Here we put $\mbar^c=\lbar-\mbar$.
\item
We put
$ H_{m,\veca,h}:=
 W_h(\lefttop{\mbar^c}N^F)
\cap
 \lefttop{\mbar}F_{\rho_m(\veca)}
 \cap
 \lefttop{\mbar^c}Gr^F_{\eta_m(\veca)}$.
Here we put $\rho_{m}(\veca)=(a_1,\ldots,a_m)$
and $\eta_m(\veca)=(a_{m+1},\ldots,a_l)$.
Then $H_{m,\veca,h}$ forms a vector subbundle
of $\lefttop{\mbar^c}\Gr^F_{\eta_m(\veca)}$,
and the image of
$\phi_{m,\veca,h}:
 H_{m,\veca,h}
\lrarr \lefttop{\mminusitibar^c}\Gr^F_{\eta_{m-1}(\veca)}$
is same as $H_{m-1,\veca,h}$.
\end{itemize}
\end{lem}
\pf
We have the sequence of degeneration
$\lefttop{\mbar^c}N^F_{|K}
 \Longrightarrow
 \lefttop{\mbar^c}N^F_{|k}
 \Longrightarrow
 \lefttop{\lbar}N^F_{|k}$.
Since the conjugacy classes of the first one and the last one
are same,
we obtain the first claim.
In particular,
$W(\lefttop{\mbar^c}N^F)\cap \lefttop{\mbar^c}\Gr^F_{\veca}$
gives a filtration of $\lefttop{\mbar^c}\Gr^F_{\veca}$
in the category of the vector bundles.

As a preparation of the proof of the second claim,
we put as follows for $x=k,K$:
\[
 H_{m\,\veca,h\,|\,x}
:=
 W_h(\lefttop{\mbar^c}N^F_{\eta_m(\veca)})_{|x}
\cap
 \lefttop{\mbar}F_{\rho_m(\veca)\,|\,x}
\cap
 \lefttop{\mbar^c}\Gr^F_{\eta_m(\veca)\,|\,x}
\]
\begin{rem}
The author apologizes that the notation is not so appropriate.
\hfill\qed
\end{rem}

To see that $H_{m,\veca,h}$ forms a vector subbundle,
we have only to show the following equality:
\begin{equation} \label{eq;8.28.10}
 \dim H_{m,\veca,h\,|\,k}
=\dim H_{m,\veca,h\,|\,K}.
\end{equation}

Let us show the second claim in Lemma \ref{lem;c11.12.20}
by an induction on $(l,m)$.
Let $P(l,m)$ denote the second claim for $(l,m)$.
Note $0\leq m\leq l$.

In the case $(l,m)=(0,0)$,
the claim $P(0,0)$ is clear.
In the case $(l,m)=(l,0)$,
the claim $P(l,0)$ follows from
the condition \ref{9.16.36} in Condition \ref{condition;9.16.35}.
Thus we have only to show
that $P(l,m-1)+P(l-1,m-1)$ implies $P(l,m)$.

See the following naturally defined morphism:
\[
\varphi_{m,\veca,h}:
 H_{m,\veca,h}\lrarr
 \lefttop{m}\Gr^F_{a_m}
 \lefttop{\mbar^c}\Gr^F_{\eta_{m}(\veca)}
=\lefttop{\mminusitibar^c}\Gr^F_{\eta_{m-1}(\veca)}.
\]
On a generic point $K$,
we have the splitting of the filtration $F$ compatible with $N$.
Hence we have
$\Image(\varphi_{m,\veca,h\,|\,K})=H_{m-1,\veca,h\,|\,K}$.
In particular,
we have $\dim\Image(\varphi_{m,\veca,h\,|\,K})
=\dim H_{m-1\,\veca\,h\,|\,K}$.
Since $H_{m-1,\veca,h}$ is a subbundle
of $\lefttop{\mminusitibar^c}\Gr^F_{\eta_{m-1}(\veca)}$
due to the hypothesis $P(l,m-1)$ of the induction.
We also obtain the following:
\begin{equation}\label{eq;8.28.5}
 \Image\bigl(\varphi_{m,\veca,h}\bigr)\subset H_{m-1,\veca,h}.
\end{equation}
We have the morphisms
$\varphi_{m,\veca,h\,|\,k}:
 H_{m,\veca,h\,|\,k}\lrarr
 \lefttop{\mminusitibar^c}\Gr^F_{\eta_m(\veca)\,|\,k} $.
Due to (\ref{eq;8.28.5}),
we have the implication
$\Image(\varphi_{m,\veca,h\,|\,k})\subset
 H_{m-1,\veca,h\,|\,k}$.
Thus we obtain the following inequality:
\begin{equation}
 \dim \Image(\varphi_{m,\veca,h\,|\,k})
\leq
 \dim H_{m-1,\veca,h\,|\,k}.
\end{equation}

Due to the hypothesis $P(l,m-1)$ of the induction,
we have $\dim H_{m-1,\veca,h\,|\,K}= \dim H_{m-1\,\veca,h\,|\,k}$.
Hence we obtain the following:
\begin{equation} \label{eq;c11.12.21}
 \dim \Image(\varphi_{m,\veca,h\,|\,k})
\leq
 \dim \Image(\varphi_{m,\veca,h\,|\,K}).
\end{equation}

We have the following for some $\epsilon>0$:
\[
 \ker \varphi_{m,\veca,h\,|\,x}
=H_{m,\veca,h\,|\,x}\cap
 \lefttop{m}F_{<a_m\,|\,x}
=H_{m,\veca-\epsilon\vecdelta_m,h\,|\,x}.
\]

For any elements $\vecb\in\real^{m-1}$
and $\vecc\in \real^{\mbar^c}$,
we put as follows:
\[
 I_{\vecb,\vecc,h}:=
 \lefttop{\mminusitibar}F_{\vecb}
 \cap
 \lefttop{\mbar^c}\Gr^F_{\vecc}
 \cap
 W_h(\lefttop{\mbar^c}N^F).
\]
Note that the tuple
$\bigl(N,\lefttop{1}F,\ldots,\lefttop{m-1}F,
 \lefttop{m+1}F,\ldots\lefttop{l}F\bigr)$
also satisfies Condition \ref{condition;9.16.35}.
Due to the hypothesis $P(l-1,m-1)$ of the induction,
$I_{\vecb,\vecc,h}$ is a vector subbundle,
and we have the following:
\begin{equation}\label{eq;c11.12.23}
 I_{\vecb,\vecc,h\,|\,x}=
 \lefttop{\mminusitibar}F_{\vecb\,|\,x}
\cap \lefttop{\mbar^c}\Gr^F_{\vecc\,|\,x}
\cap W_h(\lefttop{\mbar^c}N^F)_{|x},
\,\,(x=k,K)
\quad\quad
 \dim(I_{\vecb,\vecc,h\,|\,k})=
 \dim(I_{\vecb,\vecc,h\,|\,K}).
\end{equation}

We have the induced filtration $\lefttop{m}F(I_{\vecb,\vecc,h\,|\,x})$
on $I_{\vecb,\vecc,h\,|\,x}$,
and we have the following:
\begin{equation}
 \dim \Image\varphi_{m,\veca,\vech\,|\,x}
=\dim \Gr^{\lefttop{m}F}_{a_m}
  \bigl(I_{\rho_{m-1}(\veca),\eta_m(\veca),h\,|\,x}\bigr).
\end{equation}
Hence we have the following equality for $x=k,K$:
\begin{equation}\label{eq;c11.12.24}
 \sum_{\substack{
 \rho_{m-1}(\veca)=\vecb,\\
 \eta_{m}(\veca)=\vecc
 }}\dim\Image(\varphi_{m,\veca,h\,|\,x})
=\dim(I_{\vecb,\vecc\,|\,x}).
\end{equation}
Thus we obtain the following equality,
from (\ref{eq;c11.12.21}), (\ref{eq;c11.12.23})
and (\ref{eq;c11.12.24}):
\begin{equation}\label{eq;8.28.11}
 \dim \Image(\varphi_{m\,\veca\,\vech\,|\,k})
=\dim \Image(\varphi_{m\,\veca\,\vech\,|\,K}).
\end{equation}
We have the following:
\begin{equation}\label{eq;c11.12.25}
 \dim H_{m,\veca,h\,|\,x}
=\sum_{\substack{
 \pi_m(\vecb)=\veca,\\
 q_m(\vecb)\leq q_m(\veca)
 }} \dim\Image(\varphi_{m,\vecb,h\,|\,x}).
\end{equation}
Here $\pi_m$ denote the projection $\real^l\lrarr\real^{l-1}$,
omitting the $m$-th component.

Thus we obtain (\ref{eq;8.28.10})
from (\ref{eq;8.28.11}) and (\ref{eq;c11.12.25}).
Due to (\ref{eq;8.28.11}) and
$\dim\varphi_{m,\veca,h\,|\,K}
=\dim H_{m-1,\veca,h\,|\,K}$,
we obtain
$\dim\Image(\varphi_{m,\veca,h\,|\,k})=\dim H_{m-1,\veca,h\,|\,k}$,
which implies
$\Image(\varphi_{m,\veca,h})=H_{m-1,\veca,h}$.
Thus the induction can proceed.
\hfill\qed

%% file: a72.tex

\subsubsection{A tuple of nilpotent maps}

Let $V$ be a free $R$-module of finite rank.
Assume that we are given the following data.

\begin{condition}\label{condition;9.16.40}
\begin{enumerate}
\item \label{c11.12.40}
 Let $\lefttop{i}F$ $(i=1,\ldots,l)$ be compatible filtrations,
 and $N_j$ $(j=1,\ldots,\alpha)$ be a commuting tuple of nilpotent maps.
 We denote the tuple of filtrations by $\vecF$,
and we put $N(\ibar)=\sum_{j\leq i}N_j$.
\item \label{c11.12.41}
 We have the splitting $V_{|K}=\bigoplus U_{\veca}$
 of $\vecF$.
\item \label{c11.12.42}
 $N_{j\,|\,K}$ preserves $U_{\veca}$,
 and $N_j$ preserves the filtration $\vecF$.
\item \label{9.16.41}
 We have the induced morphisms
 \,$\lefttop{\lbar}N_j^F$ $(j=1,\ldots,\alpha)$
on $\lefttop{\lbar}\Gr^F$.
Then $(\lefttop{\lbar}N_1^F,\ldots,\lefttop{\lbar}N^F_{\alphasitabar})$
is sequentially compatible.
\hfill\qed
\end{enumerate}
\end{condition}

We put as follows:
\[
 J_{\vech,\veca,m}:=
\bigcap_{j=1}^{\alpha}W(\jbar)_{h_j}
\cap
\lefttop{\mbar}F_{\rho_m(\veca)}
\cap
\lefttop{\mbar^c}\Gr_{\eta_m(\veca)}.
\]
Let us consider the following morphism:
\[
 \psi_{\vech,\veca,m}:
 J_{\vech,\veca,m}
\lrarr
 \lefttop{\mminusitibar^c}\Gr_{\eta_{m-1}(\veca)}.
\]
\begin{lem} \label{lem;c11.12.30}
$J_{\vech,\veca,m}$ is a subbundle
of $\lefttop{\mbar^c}\Gr_{\eta_m(\veca)}$,
and we have
$\Image(\psi_{\vech,\veca,m})
=J_{\vech,\veca,m-1}$.
\end{lem}
\pf
The argument is essentially same as the proof of the previous one.
We put as follows:
\[
 J_{\vech,\veca,m\,|\,x}:=
\bigcap_{j=1}^{\alpha}W(\jbar)_{h_j\,|\,x}
\cap
\lefttop{\mbar}F_{\rho_m(\veca)\,|\,x}
\cap
\lefttop{\mbar^c}\Gr_{\eta_m(\veca)\,|\,x}.
\]
Then we have only to show the following equalities:
\begin{description}
\item[(A)]
$\dim J_{\vech,\veca,m\,|\,k}
=\dim J_{\vech,\veca,m\,|\,K}$.
\item[(B)]
$\dim\Image\psi_{\vech,\veca,m\,|\,k}
=\dim J_{\vech,\veca,m-1\,|\,k}$.
\end{description}

We denote the claims for $(l,m)$ by $P(l,m)$.
We show $P(l,m)$ by an induction.
Note that $0\leq m\leq l$.
The claim $P(0,0)$ is trivial,
and the claim $P(l,0)$ follows from the condition \ref{9.16.41}
in Condition \ref{condition;9.16.40}.
Thus we have only to show that
$P(l,m-1)+ P(l-1,m-1)$ implies $P(l,m)$.

Let us consider the following morphisms:
\[
 \begin{array}{l}
 \psi_1:
 W(\jbar)_{h_j}\cap \lefttop{\mbar^c}\Gr_{\eta_m(\veca)}
\lrarr
 \lefttop{\mminusitibar^c}\Gr_{\eta_{m-1}(\veca)}.
 \end{array}
\]

\begin{lem}
We have the following:
\[
 \begin{array}{l}
 \Image(\psi_1)=
 W(\jbar)_{h_j}\cap
 \lefttop{\mminusitibar^c}\Gr_{\eta_{m-1}(\veca)}.
 \end{array}
\]
\end{lem}
\pf
It immediately follows from Lemma \ref{lem;c11.12.20}.
\hfill\qed

\vspace{.1in}

Hence we obtain
$\Image(\psi_{\vech,\veca,m\,|\,x})
\subset  J_{\vech,\veca,m-1\,|\,x}$.
Since we are given the splitting 
of the filtrations $\vecF$ compatible with the nilpotent maps 
on the generic point $K$,
we have the following:
\[
 \Image(\psi_{\vech,\veca,m\,|\,K})=J_{\vech,\veca,m-1\,|\,K}.
\]
Due to the hypothesis $P(l,m-1)$ of the induction,
we have the following equality:
\begin{equation} \label{eq;8.28.15}
\dim J_{\vech,\veca,m-1\,|\,k}
=\dim J_{\vech,\veca,m-1\,|\,K}.
\end{equation}
Hence we obtain the following inequality:
\begin{equation}
 \dim \Image\psi_{\vech,\veca,m\,|\,k}
\leq
 \dim \Image\psi_{\vech,\veca,m\,|\,K}.
\end{equation}

On the other hand,
we have
$\ker(\psi_{\vech,\veca,m\,|\,x})=
 J_{\vech,\veca-\epsilon\vecdelta_m,m\,|\,x}$
for some small positive number $\epsilon$.

For any elements $\vecb\in\real^{m-1}$ and $\vecc\in\real^{\mbar^c}$,
we put as follows:
\[
 I'_{\vech,\vecb,\vecc}:=
 \bigcap_{j=1}^{\alpha}W(\jbar)_{h_j}
\cap
 \lefttop{\mminusitibar}F_{\rho_{m-1}(\veca)}
\cap
 \lefttop{\mbar^c}\Gr_{\eta_m(\veca)}.
\]
Due to the hypothesis $P(l-1,m-1)$ of the induction,
$I'_{\vech,\vecb,\vecc}$ is a vector subbundle.
We also have the following
\begin{equation}\label{eq;8.28.16}
 I'_{\vech,\vecb,\vecc\,|\,x}:=
 \bigcap_{j=1}^{\alpha}W(\jbar)_{h_j\,|\,x}
\cap
 \lefttop{\mminusitibar}F_{\rho_{m-1}(\veca)\,|\,x}
\cap
 \lefttop{\mbar^c}\Gr_{\eta_m(\veca)\,|\,x},
\quad
 \dim(I'_{\vech,\vecb,\vecc\,|\,k})
=\dim(I'_{\vech,\vecb,\vecc\,|\,K}).
\end{equation}

We have the induced filtration
$\lefttop{m}F\bigl(I'_{\vech,\vecb,\vecc\,|\,x}\bigr)$,
and we have the following equality:
\[
 \dim \Image(\psi_{\vech,\veca,m\,|\,x})
=\dim\Gr^{\lefttop{m}F}_{a_m}
 \bigl(I'_{\vech,\rho_{m-1}(\veca),\eta_m(\veca)\,|\,x}\bigr).
\]
Thus we have the following:
\begin{equation} \label{eq;8.28.17}
 \sum_{\substack{
 \rho_{m-1}(\veca)=\vecb,\\
 \eta_m(\veca)=\vecc
 }}\dim \Image(\psi_{\vech,\veca,m\,|\,x})
=\dim (I'_{\vech,\vecb,\vecc\,|\,x}).
\end{equation}
Due to (\ref{eq;8.28.15}), (\ref{eq;8.28.16})
and (\ref{eq;8.28.17}),
we obtain the equality
$\dim \Image\psi_{\vech,\veca,m\,|\,k}
=\dim \Image\psi_{\vech,\veca,m\,|\,K}$,
which implies $(B)$
in $P(l,m)$.

We have the following:
\[
 \dim J_{\vech,\veca,m\,|\,x}
=\sum_{\substack{
 \pi_m(\vecb)=\pi_m(\veca),\\
 q_m(\vecb)\leq q_m(\veca)
 }}
 \dim \Image(\psi_{\vech,\vecb,m\,|\,x}).
\]
Here $\pi_m$ denotes the projection
$\real^l\lrarr\real^{l-1}$,
forgetting the $m$-th component.
Then we obtain the claim $(A)$ in $P(l,m)$.
\hfill\qed

\begin{cor} \label{cor;c11.12.50}
The tuple of the filtrations
 $\bigl(
W(\itibar),\ldots,W(\alphasitabar),\lefttop{1}F,\ldots,\lefttop{l}F
\bigr)$
are compatible,
in the sense of Definition {\rm\ref{df;b11.12.1}}.
\end{cor}
\pf
We have the following morphisms:
\[
\begin{CD}
 \bigcap_{j=1}^{\alpha}W(\jbar)_{h_j}
\cap
 \lefttop{\lbar}F_{\veca}
 @>{\psi_{\vech,\veca}}>>
 \bigcap_{j=1}^{\alpha}W(\jbar)_{h_j}
\cap
 \lefttop{\lbar}\Gr^F_{\veca}
 @>{\phi_{\vech,\veca}}>>
 \Gr^W_{\vech}\lefttop{\lbar}\Gr^F_{\veca}.
\end{CD}
\]
The morphism $\psi_{\vech,\veca}$ is surjective
due to Lemma \ref{lem;c11.12.20},
and the morphism $\phi_{\vech,\veca}$
is surjective due to the condition \ref{9.16.41}
in Condition \ref{condition;9.16.40}.

Let us pick subbundles
$C_{\vech,\veca}\subset
 \bigcap_{j=1}^{\alpha}W(\jbar)_{h_j}
\cap \lefttop{\lbar}F_{\veca}$
such that
the restriction of $\phi_{\vech,\veca}\circ\psi_{\vech,\veca}$
to $C_{\vech,\veca}$ is isomorphic.

Since the filtrations
$\lefttop{\lbar}W(\itibar),\ldots\lefttop{\lbar}W(\alphasitabar)$
are compatible,
we obtain the following:
\[
 \lefttop{\lbar}\Gr^F_{\veca}
=\bigoplus_{\vech}\psi_{\vech,\veca}
 (C_{\vech,\veca}).
\]
Since $\lefttop{i}F$ $(i=1,\ldots,l)$ are compatible,
we obtain the following:
\[
 V=\bigoplus_{\veca}\bigoplus_{\vech} C_{\vech,\veca}.
\]
We have only to show the following:
\begin{equation}\label{eq;c11.12.35}
\bigcap_{j=1}^{\alpha}W(\jbar)_{h_j}
\cap
\lefttop{\lbar}F_{\veca}
=\bigoplus_{(\veck,\vecb)\leq (\vech,\veca)}
 C_{\vech,\veca}.
\end{equation}
Since we have the splitting on the generic point $K$,
it is easy to see that the restriction of (\ref{eq;c11.12.35})
holds over the generic point $K$.
We have already known that the both sides of (\ref{eq;c11.12.35})
are subbundles of $V$ (Lemma \ref{lem;c11.12.20}).
Then we can conclude that (\ref{eq;c11.12.35})
holds on $R$.
\hfill\qed

\begin{cor}
The conjugacy classes of $\nbign(\jbar)$ are constant over $R$,
and $\bigcap_{j=1}^{\alpha}W(\jbar)_{h_j}\cap
 \lefttop{\lbar}F_{\veca}$ are vector bundles
for any $\veca$ and any $\vech$.
\hfill\qed
\end{cor}

\begin{prop} \label{prop;10.26.120}
 $(\nbign_1,\ldots,\nbign_{\alpha},\lefttop{1}F,\ldots,\lefttop{l}F)$
are sequentially compatible,
in the sense of Definition {\rm\ref{df;9.16.30}}.
\end{prop}
\pf
 We have only to show that
 $(\nbign_{1\,|\,k},\ldots,\nbign_{\alpha\,|\,k},
 \lefttop{1}F_{|k},\ldots,\lefttop{l}F_{|k})$
 are sequentially compatible
 in the sense of Definition \ref{df;9.16.20}.

 The condition \ref{9.16.21} in Definition \ref{df;9.16.20}
 follows from 
 the condition \ref{c11.12.42} in Condition \ref{condition;9.16.40}.
 The condition \ref{9.16.22} in Definition \ref{df;9.16.20}
 follows from Corollary \ref{cor;c11.12.50}.
 The condition \ref{9.16.23} follows Lemma \ref{lem;c11.12.20}.

 Let us see the condition \ref{9.16.24}.
 We use an induction on $\alpha$.
 In the case $\alpha=1$, there remain nothing to prove.
 We assume that the claim holds in the case $\alpha-1$,
 and we will prove the claim also holds in the case $\alpha$.

 On the vector bundle $\Gr^{W(N_1)}(V)$,
 we have the induced filtrations
 $\lefttop{1}F^{(1)},\ldots \lefttop{l}F^{(1)}$
 and the nilpotent maps $\nbign^{(1)}_2,\ldots,\nbign^{(1)}_{\alpha}$.
 Due to the hypothesis of the induction,
 the tuple
 $\bigl(\nbign^{(1)}_2,\ldots,\nbign^{(2)}_{\alpha},
 \lefttop{1}F^{(1)},\ldots,\lefttop{l}F^{(1)}
 \bigr)$ is sequentially compatible.
 
 We put $\nbign^{(1)}(\ibar)=\sum_{j\leq i}\nbign^{(1)}_j$,
 which is same as the induced morphism by $\nbign(\ibar)$
 on $\Gr^{W(N_1)}$.
 Let $W^{(1)}(\ibar)$ denote the induced filtration
 on $\Gr^{W(N_1)}$  by $W(\ibar)$.
We have only to show the following:
\begin{equation} \label{eq;c11.12.55}
 W^{(1)}(\ibar)_{h+a}\cap
 \Gr^{W(N_1)}_a\lefttop{I}\Gr^F
=W(\nbign^{(1)}(\ibar))_h\cap \Gr^{W(N_1)}_a\lefttop{I}\Gr^F.
\end{equation}
 Since we have the splitting on the generic point $K$,
 it is easy to see that (\ref{eq;c11.12.55}) 
 holds when it is restricted to the generic point $K$.
We have already known that the both sides are vector subbundles
of $\Gr^{W(N_1)}_a\lefttop{I}\Gr^F$,
we can conclude that (\ref{eq;c11.12.55}) holds on $R$.
\hfill\qed

%% file: a72.1.tex

\subsubsection{A compatible tuple of filtrations and a morphism}

\label{subsubsection;c11.12.80}

Let $R$ be a discrete valuation ring,  $K$ be the quotient field,
and $k$ be the residue field.
Let $V^{(a)}$ be free $R$-modules $(a=1,2)$,
and let $\vecF:=\bigl(
 \lefttop{i}F(V^{(a)})\,\big|\, i\in I\bigr)$
be a compatible tuple of filtrations of $V^{(a)}$ in the category of
free $R$-modules.
Let $f:V^{(1)}\lrarr V^{(2)}$ be the morphism
preserving the filtrations.
We have the induced morphism
$\lefttop{I}\Gr^F_{\veceta}(f):
 \lefttop{I}\Gr^F_{\veceta}\bigl(V^{(1)}\bigr)
 \lrarr \lefttop{I}\Gr^F_{\veceta}\bigl(V^{(2)}\bigr)$
for any element $\veceta\in\real^I$.

Let $S$ be a finite subset of $\real^I$.
For simplicity, we use the following notation:
\[
 \lefttop{I}F_{S}(V^{(a)}):=
 \sum_{\veceta\in S}\lefttop{I}F_{\veceta}(V^{(a)}),
\quad
 \lefttop{I}\Gr^F_S(V^{(a)}):=
 \bigoplus_{\veceta\in S}\lefttop{I}\Gr_{\veceta}(V^{(a)}).
\]
We have the naturally defined projection
$\pi_{S}:\lefttop{I}F_S\bigl(V^{(a)}\bigr)
\lrarr \lefttop{I}\Gr^{F}_S\bigl(V^{(a)}\bigr)$.

\begin{prop} \label{prop;c11.12.60}
Assume the following:
\begin{itemize}
\item
We have a splitting
 $V^{(a)}_{|K}=\bigoplus_{\veceta\in\real^l}U_{\veceta}^{(a)}$
of the tuple of filtrations $\vecF_{|K}$
satisfying
$f_K(U_{\veceta}^{(a)})\subset U^{(2)}_{\veceta}$.
\item
 The image
 $\Image\lefttop{I}\Gr^F_{\veceta}(f)$
 is the vector subbundle
 of $\lefttop{I}\Gr^F_{\veceta}(V^{(2)})$,
 for any $\veceta\in\real^I$.
\end{itemize}
Then the following claims hold.
\begin{enumerate}
\item
 For any finite subset $S\subset\real^I$,
 the image $f\bigl(\lefttop{I}F_{S}(V^{(1)})\bigr)$
 is a vector subbundle of $V^{(2)}$.
\item
 We have the following:
\[
 f\bigl( \lefttop{I}F_{S}(V^{(1)})
 \bigr)
=\Image(f)\cap\lefttop{I}F_S(V^{(2)}).
\]
\end{enumerate}
\end{prop}
\pf
For any finite subset $S$ of $\real^I$,
we put as follows:
\[
 L(S):=\max\Bigl\{
 \sum_{i\in I} q_i(\veca)\,\Big|\,
 \veca\in S
 \Bigr\}.
\]
The following claim is denoted by $P(r)$:
\begin{description}
\item[$P(r)$]
 The claim of Proposition \ref{prop;c11.12.60} holds
 in the case $L(S)\leq r$.
\end{description}

\begin{lem}\mbox{{}}\label{lem;c11.12.61}
\begin{itemize}
\item
 The claim $P(r)$ holds for any sufficiently negative $r$.
\item
 If the claim $P(r)$ holds for some $r\in\real$,
 then there exists a positive number 
 such that $P(r')$ holds for any $r'$ such that $r\leq r'\leq r+\epsilon$.
\end{itemize}
\end{lem}
\pf
It follows from the finiteness of the set
$\bigl\{b\in\real\,\big|\, 
 \exists i,\,\exists a,
 \lefttop{i}\Gr^F_b\bigl(V^{(a)}\bigr)\neq 0,
 \bigr\}$.
\hfill\qed

\vspace{.1in}

Due to Lemma \ref{lem;c11.12.61},
we have only to show 
that $P(r)$ holds under the assumption
$P(r')$ holds for any real numbers $r'<r$,
which we will show in the following.

\begin{lem}
Let $S$ be a finite subset of $\real^I$.
Let us consider the projection
$\pi_{S}:
 \lefttop{I}F_{S}\bigl(V^{(a)}\bigr)
 \lrarr
 \lefttop{I}\Gr^{F}_{S}\bigl(V^{(a)}\bigr)$.
The kernel of $\pi_S$ is described as the form
$\lefttop{I}F_{S'}\bigl(V^{(a)}\bigr)$
for some finite subset $S'\subset\real^I$
such that $L(S')<L(S)$.
\end{lem}
\pf
We have the following:
\[
 \ker \pi_S
=\sum_{\veca\in S}
 \sum_{\veceta\lneq\veca}
 \lefttop{I}F_{\veceta}\bigl(V^{(a)}\bigr).
\]
Then it is clear from the compatibility of 
the tuple $\vecF$.
\hfill\qed

\begin{lem}\label{lem;c11.12.65}
Let us consider the morphism
$\pi'_S:f\bigl(\lefttop{I}F_{S}(V^{(1)})\bigr)
\lrarr \lefttop{I}\Gr^F_S(V^{(2)})$,
induced by the projection $\pi_{S}$.
Then we have the following:
\[
 \Image \pi_S'
=\Image \lefttop{I}\Gr^F_{S}(f)
 \subset \lefttop{I}\Gr^F_S(V^{(2)}).
\]
\end{lem}
\pf
We have the following commutative diagramm:
\[
 \begin{CD} 
 \lefttop{I}F_S(V^{(1)})
 @>>>
 \lefttop{I}F_S(V^{(2)}) \\
 @V{\pi_{S}}VV @V{\pi_S}VV\\
 \lefttop{I}\Gr^F_S(V^{(1)})
 @>{\lefttop{I}\Gr^F_S(f)}>>
 \lefttop{I}\Gr^F_S(V^{(2)}).
 \end{CD}
\]
Then Lemma \ref{lem;c11.12.65} immediately follows.
\hfill\qed

\vspace{.1in}

We always have
$f\bigl(\lefttop{\lbar}F_{S}(V^{(1)})\bigr)
\subset
 \Image(f)\cap \lefttop{I}F_{S}(V^{(2)})$.
Hence we always have the following:
\begin{equation}\label{eq;c11.12.66}
 \pi_{S}\Bigl(
 f\bigl(\lefttop{\lbar}F_{S}(V^{(1)})\bigr)
 \Bigr)
\subset
 \pi_S\Bigl(
 \Image(f)\cap \lefttop{I}F_{S}(V^{(2)})
 \Bigr).
\end{equation}
\begin{lem}
In {\rm(\ref{eq;c11.12.66})},
the equality holds.
\end{lem}
\pf
Since we have the splitting on the generic point $K$,
it is easy to see that the equality holds,
when we restrict (\ref{eq;c11.12.66}) to the generic point.
Since the left hand side is a vector subbundle
in $\lefttop{I}\Gr^F_S(V^{(2)})$ due to Lemma \ref{lem;c11.12.65},
we obtain the equality on $R$.
\hfill\qed

\vspace{.1in}
Let us pick an appropriate finite subset $S'\subset\real^I$
such that $L(S')<L(S)$
and $\ker\pi_S=\lefttop{I}F_{S'}\bigl(V^{(a)}\bigr)$ $(a=1,2)$.
We have the following:
\begin{equation}\label{eq;c11.12.70}
 \Image(f)\cap \lefttop{I}F_{S}(V^{(2)})\cap\ker\pi_S
=\Image(f)\cap \lefttop{I}F_{S'}\bigl(V^{(2)}\bigr).
\end{equation}
We have the following implication:
\begin{equation} \label{eq;c11.12.71}
 f\bigl(
 \lefttop{I}F_{S'}\bigl(V^{(1)}\bigr)
 \bigr)
\subset
 f\bigl(
 \lefttop{I}F_{S}\bigl(V^{(1)}\bigr)
 \bigr)
\cap \ker\pi_S
\subset
 \Image(f)\cap
 \lefttop{I}F_S\bigl(V^{(1)}\bigr)
 \cap\ker\pi_S.
\end{equation}
Due to the assumption $P(r')$ $(r'<r)$,
we have 
$f\bigl(
 \lefttop{I}F_{S'}\bigl(V^{(1)}\bigr)
 \bigr)=
\Image(f)\cap \lefttop{I}F_{S'}\bigl(V^{(2)}\bigr)$.
Then we obtain the following equality
from (\ref{eq;c11.12.70}) and (\ref{eq;c11.12.71}):
\begin{equation}\label{eq;c11.12.72}
 f\bigl(
 \lefttop{I}F_{S}\bigl(V^{(1)}\bigr)
 \bigr)
\cap \ker\pi_S
=f\bigl(
  \lefttop{I}F_{S'}\bigl(V^{(1)}\bigr)
 \bigr)
=\Image(f)\cap
 \lefttop{I}F_{S'}\bigl(V^{(2)}\bigr)
=\Image(f)\cap \ker\pi_S.
\end{equation}
Hence we obtain the equality
$f\bigl(\lefttop{I}F_S\bigl(V^{(1)}\bigr)\bigr)
=\Image(f)\cap\lefttop{I}F_S\bigl(V^{(2)}\bigr)$.

Due to the assumption $P(r')$ $(r'<r)$,
$ f\bigl(
 \lefttop{I}F_{S}\bigl(V^{(1)}\bigr)
 \bigr)
\cap \ker\pi_S$ is a vector subbundle of $V^{(2)}$.
It follows that
$ f\bigl(
 \lefttop{I}F_{S}\bigl(V^{(1)}\bigr)
 \bigr)$ is a vector subbundle of $V^{(2)}$.
Thus the proof of Proposition \ref{prop;c11.12.60}
is accomplished.
\hfill\qed

\begin{cor}
Under the assumption of Proposition {\rm\ref{prop;c11.12.60}},
we have
$\lefttop{I}\Gr^F_{\veceta}\bigl(\Image(f)\bigr)
=\Image\lefttop{I}\Gr^F_{\veceta}(f)$,
for any element $\veceta\in\real^I$.
\hfill\qed
\end{cor}

%% file: b12.1.tex

\subsubsection{Decomposition}

Let $R$, $K$, $k$, $V^{(a)}$ and 
$\vecF=\bigl(\lefttop{i}F\,\big|\,i\in I\bigr)$
be as in the subsubsection \ref{subsubsection;c11.12.80}.
Let $f:V^{(1)}\lrarr V^{(2)}$ and $g:V^{(2)}\lrarr V^{(1)}$
be morphisms preserving the filtrations.

\begin{lem} \label{lem;9.23.15}
Assume the following:
\begin{itemize}
\item
We have splittings
$V^{(a)}_{|K}=\bigoplus U_{\veceta}^{(a)}$ 
of the filtrations $\vecF\bigl(V^{(a)}\bigr)_{|K}$
on the generic point $K$
\item
We have
$f_{|K}\bigl(U^{(1)}_{\veceta\,|\,K}\bigr)
\subset U^{(2)}_{\veceta\,|\,K}$,
and
$ g_{|K}\bigl(U^{(2)}_{\veceta\,|\,K}\bigr)
\subset
 U^{(1)}_{\veceta\,|\,K}$
on the generic point $K$.
\item
We have
$\lefttop{I}\Gr_{\veceta}\bigl(V^{(2)}\bigr)
=\Image \bigl(\lefttop{I}\Gr_{\veceta}(f)\bigr)
\oplus
 \Ker\bigl(\lefttop{I}\Gr_{\veceta}(g)\bigr)$.
In particular,
$\Image\bigl(\lefttop{I}\Gr_{\veceta}(f)\bigr)$,
$\Ker\bigl(\lefttop{I}\Gr_{\veceta}(g)\bigr)$
and $\Image(\lefttop{I}\Gr_{\veceta}(g))$
are the vector subbundles.
\end{itemize}

Then we have the decomposition
$V^{(2)}=\Image(f)\oplus\Ker(g)$.
\end{lem}
\pf
$\Image(f)$ and $\Ker(g)$ are vector subbundles of $V^{(2)}$,
due to Proposition \ref{prop;c11.12.60}.
The tuple of filtrations $\vecF\bigl(V^{(2)}\bigr)$
induce the tuple of filtrations
$\vecF\bigl(\Image(f)\bigr)$
and $\vecF\bigl(\ker(g)\bigr)$ on $R$
The decomposition
$\bigl(U^{(a)}_{\veceta}\,\big|\,\veceta\in\real^I\bigr)$
of $V^{(2)}_{|K}$ induces the decomposition
of $\Image(f)$ and $\Ker(g)$ on the generic point.
Let us consider the naturally defined morphism
$\Phi:\Image(f)\oplus\Ker(g)\lrarr V^{(2)}$.
Then it is easy to see that the assumption of
Proposition \ref{prop;c11.12.60} is satisfied.
Thus the image $\Image\Phi$ is a vector subbundle of $V^{(2)}$.
Since $\Image(\Phi_{|K})$ and $V^{(2)}_{|K}$ are same,
we obtain $\Image(\Phi)=V^{(2)}$ on $R$.
\hfill\qed

%% file: a53.3.tex

\subsubsection{Preliminary}
\label{subsubsection;b12.9.1}

We put $X:=\Delta$ and $D=\{O\}$.
We put $E:=\nbigo\cdot e$.
We consider the hermitian metric $h$ given by
$h(e,e)=|z|^{-2a}$.
We consider the Higgs field $\theta$ given by
$\theta(e)=e\cdot \alpha\cdot dz/z$.
Then we have the following:
\[
 \delbar_Ee=0,\quad
 \del_E e=e\cdot\Bigl(-a\cdot\frac{dz}{z}\Bigr),\quad
 \theta^{\dagger}(e)=e\cdot \overline{\alpha}\cdot\frac{d\zbar}{\zbar}.
\]

We put $v=\exp\bigl(-\bar{\alpha}\cdot\lambda\cdot \log|z|^2\bigr)\cdot e$.
Then we obtain the following:
\[
 \bigl(
 \delbar_E+\lambda\theta^{\dagger}
 \bigr)v=0,
\quad
 \DD v=\bigl(\lambda\cdot\del_E+\theta\bigr)\cdot v
=\bigl(-\alphabar\cdot\lambda^2-a\cdot\lambda+\alpha\bigr)v\cdot
 \frac{dz}{z}.
\]
Then $v$ gives a holomorphic frame
of $\nbige$ over $\nbigx-\nbigd$.
We put as follows:
\[
 s:=\exp\bigl(
 \bigl(\overline{\alpha}\cdot\lambda+a-\alpha\cdot\lambda^{-1}\bigr)
 \cdot\log z
 \bigr)\cdot v
=\exp\bigl(
 -\alpha\cdot\lambda^{-1}\log z
 +a\cdot\log z-\overline{\alpha}\cdot\lambda\cdot\log \overline{z}
 \bigr)\cdot e.
\]
Then $s$ is a frame of the holomorphic bundle
$\nbigh(E)$ over $\cnum^{\ast}$.

We put $e^{\dagger}:=|z|^{2a}\cdot e$.
Then we have $\del_E e^{\dagger}=0$.
Namely $e^{\dagger}$ gives the frame of
$(E,\del_E)$.
We also have $h(e^{\dagger},e^{\dagger})=|z|^{2a}$,
and $\delbar_Ee^{\dagger}=e^{\dagger}\cdot a\cdot d\zbar/\zbar$.

We put
$v^{\dagger}
:=\exp\bigl(-\alpha\cdot\mu\cdot\log|z|^2\bigr)\cdot e^{\dagger}$.
Then we obtain the following:
\[
 \bigl(
 \del_E+\mu\cdot\theta
 \bigr)\cdot v^{\dagger}=0,
\quad
 \DD^{\dagger}v^{\dagger}
=\bigl(\mu\cdot\delbar_E+\theta^{\dagger}\bigr)\cdot v^{\dagger}
=\bigl(-\alpha\cdot\mu^2+a\cdot\mu+\overline{\alpha}\bigr)
 \cdot v^{\dagger}\cdot\frac{d\zbar}{\zbar}.
\]
Namely $v^{\dagger}$ is a frame of $\nbige^{\dagger}$ over
$\nbigx^{\dagger}-\nbigd^{\dagger}$.
We put as follows:
\[
 s^{\dagger}=
 \exp\bigl(
 \bigl(\mu\cdot\alpha-a-\overline{\alpha}\cdot\mu^{-1}\bigr)
 \log \zbar
 \bigr)\cdot v^{\dagger}.
\]
Then $s^{\dagger}$ is a frame of $\nbigh^{\dagger}(E)$.

\begin{lem}
We have $s=s^{\dagger}$.
\end{lem}
\pf
It can be shown by a direct calculation, as follows:
\[
 s^{\dagger}=
 \exp\bigl(
 -\overline{\alpha}\cdot\mu^{-1}\cdot\log \zbar
 -a\cdot \log \zbar
 -\alpha\cdot\mu\cdot\log z
 \bigr)\cdot e^{\dagger}\\
=\exp\bigl(
 -\overline{\alpha}\cdot\mu^{-1}\cdot\log \zbar
 +a\cdot\log z
 -\alpha\cdot\mu\cdot\log z
 \bigr)e=s.
\]
Thus we are done.
\hfill\qed

\vspace{.1in}

We put $w:=v_{|\cnum_{\lambda}\times\{O\}}$
and $w^{\dagger}:=v^{\dagger}_{|\cnum_{\mu}\times\{O\}}$.
The induced objects $\nbigg_u$ and $\nbigg^{\dagger}_{u^{\dagger}}$
are generated by $w$ and $w^{\dagger}$ respectively.

\subsubsection{The gluings of $S^{\can}_u(E)$ and $S_u(E,P)$}

By definition,
we have $\Phi^{\can}_u(s)=w_{|\cnum^{\ast}_{\lambda}}$ and
$\Phi^{\dagger\,\can}_{u^{\dagger}}(s^{\dagger})
 =w^{\dagger}_{|\cnum_{\mu}^{\ast}}$.
We also have $s=s^{\dagger}$.
Then the gluing of $S^{\can}_u(E)$ is given by the relation
$w_{|\cnum_{\lambda}^{\ast}}
=w^{\dagger}_{\cnum_{\mu}^{\ast}}$.
In particular,
$S^{\can}_u(E)$ is isomorphic to $\nbigo_{\proj^1}$.

Let us see the gluing of $S(E,P)$.
We denote $v_{|\cnum_{\lambda}^{\ast}\times\{P\}}$
by $v_{|P}$ for simplicity of notation.
We will also use the similar convention.
We have the following:
\[
 v_{|P}=\exp\bigl(
 -\overline{\alpha}\cdot\lambda\cdot
 \log|z(P)|^2
 \bigr)\cdot e_{|P}.
\]
Then we have the following:
\begin{multline}
 v^{\dagger}_{|P}
=\exp\bigl(
 -\alpha\cdot\mu\cdot\log|z(P)|^2
 \bigr)\cdot e^{\dagger}_{|P}
=\exp\bigl(
 -\alpha\cdot\mu\cdot\log|z(P)|^2
 +a\cdot \log|z(P)|^2
 \bigr)\cdot e_{|P} \\
=\exp\bigl(
 \bigl(
 -\alpha\cdot\mu+a+\overline{\alpha}\cdot\lambda
 \bigr)\cdot \log|z(P)|^2
 \bigr)\cdot v_{|P}.
\end{multline}

We put as follows:
\[
 \tilde{w}:=
 \exp\bigl(
 \overline{\alpha}\cdot\lambda\cdot
 \log|z(P)|^2
 \bigr)\cdot w,
\quad
 \tilde{w}^{\dagger}:=
 \exp\bigl(
 \alpha\cdot\mu\cdot\log|z(P)|^2
 \bigr)
\cdot w^{\dagger}.
\]
Then $\nbigg_u$ and $\nbigg^{\dagger}_u$
are generated by $w$ and $w^{\dagger}$ respectively,
and the gluing of $S_u(E,P)$ is given by the relation
$\tilde{u}^{\dagger}=\tilde{u}\cdot |z(P)|^{2a}$.
In particular, $S_u(E,P)$ is isomorphic to
$\nbigo_{\proj^1}$.

\begin{cor}\mbox{{}}
\begin{itemize}
\item
The vector bundles $S^{\can}_u(E)$ and $S(E,P)$
are pure twistors of weight $0$.
\item
A frame of $H^0\big(\proj^1,S^{\can}_u(E)\big)$
is given by $w=w^{\dagger}$.
\item
A frame of $H^0\big(\proj^1,S_u(E,P)\big)$
is given by $\tilde{u}=|z(P)|^{-2a}\cdot\tilde{u}^{\dagger}$.
\hfill\qed
\end{itemize}
\end{cor}

\subsubsection{Pairing}

The pairing
$S:\nbige\otimes\sigma^{\ast}\nbige\lrarr\nbigo_{\nbigx-\nbigd}$
is given as follows:
\[
 S\bigl(v,\sigma^{\ast}(v^{\dagger})\bigr)
=h\bigl(
 v(\lambda,x),
 v^{\dagger}(-\bar{\lambda},x)
 \bigr)
=\exp\bigl(-\overline{\alpha}\cdot\lambda\cdot\log|z|^2
 \bigr)
\cdot
 \overline{
 \exp\bigl(
 -\alpha\cdot\bigl(-\overline{\lambda}\bigr)
 \cdot \log|z|^2
 \bigr)
 }
=1.
\]
Thus the induced pairing
$S:\nbigg_u\otimes\sigma^{\ast}\nbigg^{\dagger}_{u^{\dagger}}
\lrarr\nbigo_{\cnum_{\lambda}}$ is given by
$S\bigl(w,\sigma^{\dagger}w^{\dagger}\bigr)=1$.

For the global section $v_1=w=w^{\dagger}$ of $S^{\can}(E)$,
we have $S(v_1,\sigma^{\ast}v_1)=S(w,\sigma^{\ast}w^{\dagger})=1>0$.
In particular,
the pairing
$S:S^{\can}(E)\otimes\sigma^{\ast}S^{\can}(E)
\lrarr \Tate(0)$ induces the polarization of $S^{\can}(E)$.

For the global section
$v_2=\tilde{u}=|z(P)|^{-2a}\cdot \tilde{u}^{\dagger}$,
we have the following:
\begin{multline}
 S\bigl(v_2,\sigma^{\ast}v_2\bigr)
=S\bigl(\tilde{u},\sigma^{\ast}|z|^{-2a}\tilde{u}^{\dagger}\bigr)\\
=|z|^{-2a}\cdot\exp\bigl(
 \overline{\alpha}\cdot\lambda\cdot\log|z(P)|^2
 \bigr)
\cdot
 \overline{
 \exp\bigl(
 -\alpha\cdot\overline{\lambda}\cdot\log|z(P)|^2
 \bigr)
 }\cdot S(w,\sigma^{\ast}w^{\dagger})=|z(P)|^{-2a}>0.
\end{multline}
In particular,
the pairing $S:S_u(E,P)\otimes\sigma^{\ast}S_u(E,P)\lrarr \Tate(0)$
is the polarization.

\subsubsection{The canonical frame}

Let $\lambda_0$ be an element of $\cnum_{\lambda}$.
Let $u=(a,\alpha)$ be an element of $\real\times\cnum$.\\

\vspace{.1in}
\noindent
{\bf The case $\alpha\neq 0$}\\
Let $b$ be a real number such that
$b\not\in \paramap(\lambda_0,u)+\seisuu$.
In the case we have the integer
$\nu$ determined by the condition
$b-1<\nu+\paramap(\lambda_0,u)<b$.
Let $v$ be as in the subsubsection \ref{subsubsection;b12.9.1}.
Then the section $z^{-\nu}v$ is a frame of
$\prolongg{b}{\nbigl(u)}$,
which we call the canonical frame at $\lambda_0$.

\vspace{.1in}
\noindent
{\bf The case $\alpha=0$}.
In the case, we recall $\paramap(\lambda,u)=a$ for any $\lambda$,
and we have $v=e$.
Let $b$ be any real number.
We have the integer
$\nu$ determined by the condition $b-1<\nu+a\leq b$.
Then the frame $z^{-\nu}\cdot v$ is called the canonical frame
of $\nbigl(u)$ at $\lambda_0$.

\subsubsection{The higher dimensional case}

Let $\vecu=(u_1,\ldots,u_n)$ be an element of $(\real\times\cnum)^n$,
where $u_i=(a_i,\alpha_i)$.
We put $X=\Delta^n$ and $D=\bigcup_{i=1}D_i$,
where $D_i=\{z_i=0\}$.
Then we put $L(\vecu)=\nbigo_{X-D}\cdot e$.
We have the Higgs field $\theta$ and the metric $h$
determined as follows:
\[
 \theta\cdot e=e\cdot\sum\alpha_i\cdot\frac{dz_i}{z_i},
\quad
 h(e,e)=\prod_{i=1}^n|z_i|^{-2a_i}.
\]
Let $q_i$ denote the projection
of $X$ onto the $i$-th component.
Then we have the isomorphism
$L(\vecu)\simeq \bigoplus_{i=1}^n q_i^{\ast}L(u_i)$
compatible with the metric and the Higgs field.
Hence $(L(\vecu),\theta,h)$ is a harmonic bundle.
We have the holomorphic frame of the deformed holomorphic bundle
$\nbigl(\vecu)$ over $\nbigx-\nbigd$:
\[
 v=\exp\Bigl(
 -\sum \alpha_i\cdot \log|z_i|^2
 \Bigr)\cdot e.
\]
Let $\lambda_0$ be an element of $\cnum_{\lambda}$.
We have the canonical frame $f_i$ of $\nbigl(u_i)$
around $\lambda_0$.
Then the section
$\prod_{i=1}^n q_i^{\ast}f_i$ gives the canonical frame
of $\prolongg{\vecb}{\nbigl(\vecu)}$ around $\lambda_0$.

%% file: a53.2.tex

\subsubsection{Preliminary}
\label{subsubsection;b12.9.5}

We put $X:=\Delta$ and $D:=\{O\}$.
We put $y:=-\log|z|^2$.
We put $E:=\nbigo_{X-D}\!\cdot\! e_1\oplus\nbigo_{X-D}\!\cdot \!e_{-1}$.
The metric $h$ is given as follows:
\[
 H(h,\vece)=
 \left(
 \begin{array}{cc}
 y & 0 \\ 0 & y^{-1}
 \end{array}
 \right).
\]
The Higgs field $\theta$ is given as follows:
\[
 \theta\cdot\vece=
 \vece\cdot
 \left(
 \begin{array}{cc}
 0 & 0 \\ 1 & 0
 \end{array}
 \right)
\frac{dz}{z}.
\]
Then $Mod(2)=(E,\delbar_E,\theta,h)$ is a harmonic bundle.
(See our previous paper \cite{mochi}, for example).
We take a frame $\vecv=(v^{1,0},v^{0,1})$ of
the deformed holomorphic bundle $\nbige$
given as follows:
\begin{equation}\label{eq;b12.9.10}
 \vecv=\vece\cdot
 \left(
 \begin{array}{cc}
 1 & -\lambda\cdot y^{-1}\\
 0 & 1
 \end{array}
 \right).
\end{equation}
Then $\vecv$ gives a normalizing frame of $\prolong{\nbige}$,
i.e.,
it satisfies the following condition:
\[
 \DD\vecv=
 \vecv\cdot
 \left(
 \begin{array}{cc}
 0 & 0 \\ 1 & 0
 \end{array}
 \right)\frac{dz}{z}.
\]
We also put as follows:
\[
 \vece^{\dagger}=\vece\cdot
 \left(
 \begin{array}{cc}
 y^{-1} & 0 \\ 0 & y
 \end{array}
 \right),
\quad
 \vecv^{\dagger}=
 \vece^{\dagger}\cdot
 \left(
 \begin{array}{cc}
 1 & 0 \\ -\mu\cdot y^{-1} & 1
 \end{array}
 \right).
\]
Then $\vecv^{\dagger}$ gives a normalizing frame
of $\prolong{\nbige^{\dagger}}$:
\[
 \DD^{\dagger}\vecv^{\dagger}
=\vecv^{\dagger}\cdot
 \left(
 \begin{array}{cc}
 0 & 1 \\ 0 & 0
 \end{array}
 \right)\frac{d\bar{z}}{\bar{z}}.
\]

\subsubsection{The induced objects}

It is easy to see that
$\KMS(\nbige^0)=\bigl\{(n,0)\,\big|\,n\in\seisuu\bigr\}
 \subset\real\times\cnum$.
We have the natural isomorphisms
$S^{\can}_{(n,0)}(E)\simeq S^{\can}_{(0,0)}(E)$
and
$S_{(n,0)}(E,P)\simeq S_{(0,0)}(E,P)$.
Hence we only consider the case $\vecu=(0,0)$.
In the following, we omit to denote the subscript $\vecu$.
We also omit to denote `$\itibar$'.

We put 
$u^{1,0}:=v^{1,0}_{|O\times\cnum_{\lambda}}$
and 
$u^{0,1}=v^{0,1}_{|O\times\cnum_{\lambda}}$.
We have 
$\nbigg=
 \nbigo_{\cnum_{\lambda}}\!\cdot\! u^{1,0}
 \oplus
 \nbigo_{\cnum_{\lambda}}\!\cdot\! u^{0,1}$,
by definition.
We also have
$\nbigg(\nbige)=
 \nbige_{|\nbigx^{\shikaku}}$.
Thus we have the following:
\[
 \nbigg(\nbige)_{|P}=\nbige_{|P}
=\nbigo_{\cnum_{\lambda}}\!\cdot\! v^{1,0}_{|P}
 \oplus
 \nbigo_{\cnum_{\lambda}}\!\cdot\! v^{0,1}_{|P}
=\nbigo_{\cnum_{\lambda}}\!\cdot\! e_{1\,|\,P}
 \oplus
 \nbigo_{\cnum_{\lambda}}\!\cdot\! e_{-1\,|\,P}.
\]
Here we use
the notation ``$|P$''
instead of ``${|\{P\}\times\cnum_{\lambda}^{\ast}}$''.
In the following in this subsubsection,
we use this abbreviation of the notation for simplicity.

We have the following equality on the plane $\{P\}\times\cnum_{\lambda}$:
\[
 \vecv_{|P}=\vece_{|P}\cdot
 \left(
 \begin{array}{cc}
 1 & -\lambda\cdot y(P)^{-1}\\
 0 & 1
 \end{array}
 \right).
\]

On the conjugate side,
we have the following:
\[
 \nbigg^{\dagger}
=\nbigo_{\cnum_{\mu}}\!\cdot\! u^{\dagger\,1,0}\,
\oplus\,
 \nbigo_{\cnum_{\mu}}\!\cdot\! u^{\dagger\,0,1},
\quad\quad
 u^{\dagger\,1,0}:=v^{\dagger\,1,0}_{|O\times\cnum_{\mu}},
\quad
 u^{\dagger\,0,1}:=v^{\dagger\,0,1}_{|O\times\cnum_{\mu}}.
\]
We also have the following:
\[
 \nbigg^{\dagger}(\nbige^{\dagger})_{|P}
=\nbigo_{\cnum_{\mu}}\!\cdot\! v^{\dagger\,1,0}_{|P}
\oplus
 \nbigo_{\cnum_{\mu}}\!\cdot\! v^{\dagger\,0,1}_{|P}
=\nbigo_{\cnum_{\mu}}\!\cdot\! e^{\dagger}_{1\,|\,P}
\oplus
 \nbigo_{\cnum_{\mu}}\!\cdot\! e^{\dagger}_{-1\,|\,P}.
\]
Then we have the following:
\begin{equation} \label{eq;9.13.60}
 \vecv^{\dagger}_{|P}
=\vece^{\dagger}_{|P}\cdot
 \left(
 \begin{array}{cc}
 1 & 0 \\
 -\mu\cdot y(P)^{-1} & 1
 \end{array}
 \right),
\quad\quad
 \vece^{\dagger}_{|P}
=\vece_{|P}\left(
 \begin{array}{cc}
 y(P)^{-1} & 0 \\ 0 & y(P)
 \end{array}
\right).
\end{equation}

\subsubsection{The gluing of $S(P)$ and the nilpotent maps}

Then we have the following:
\begin{equation} \label{eq;9.13.61}
 \Phi_{(P,O)}\bigl(\vece_{|P}\bigr)
=\vecu\cdot
 \left(\begin{array}{cc}
 1 & \lambda\cdot y(P)^{-1}\\ 0 & 1
 \end{array}
 \right),
\quad\quad
 \Phi^{\dagger}_{(P,O)}\bigl(\vece^{\dagger}_{|P}\bigr)
=\vecu^{\dagger}\cdot
 \left(\begin{array}{cc}
 1 & 0 \\ \mu\cdot y(P)^{-1} & 1
 \end{array}
 \right).
\end{equation}
From (\ref{eq;9.13.61}) and the first equation in (\ref{eq;9.13.60}),
the vector bundle $S(E,P)$
is obtained by the following gluing:
\[
 \vecu=\vecu^{\dagger}
 \left(
 \begin{array}{cc}
 1 & 0 \\ \mu\cdot y(P)^{-1} & 1
 \end{array}
 \right)
\cdot
 \left(
 \begin{array}{cc}
 y(P) & 0 \\ 0 & y(P)^{-1}
 \end{array}
 \right)\cdot
 \left(
 \begin{array}{cc}
 1 & -y(P)^{-1}\cdot\lambda \\
 0 & 1
 \end{array}
 \right)
=\vecu^{\dagger}
 \left(
 \begin{array}{cc}
 y(P) & -\lambda \\
 \mu & 0
 \end{array}
 \right).
\]
The nilpotent map $\nbign^{\sankaku}$ is given as follows:
\[
 \nbign^{\sankaku}_{|\cnum_{\lambda}}:\quad
 u^{1,0}\longmapsto u^{0,1}\otimes t_0,
\quad
 u^{0,1}\longmapsto 0,
\quad
\mbox{\rm on the plane $\cnum_{\lambda}$,}
\]
\[
 \nbign^{\sankaku}_{|\cnum_{\mu}}:\quad
 u^{\dagger\,1,0}\longmapsto 0,
\quad
 u^{\dagger\,0,1} \longmapsto -u^{\dagger\,1,0}\otimes t_{\infty},
\quad
\mbox{\rm on the plane $\cnum_{\mu}$.}
\]
Hence 
the vector bundle
$\Gr^{W^{\sankaku}}_1$ is given by the gluing
$u^{1,0}=\mu\cdot u^{\dagger\,0,1}$.
The vector bundle
$\Gr^{W^{\sankaku}}_{-1}$
is given by $u^{0,1}=-\lambda\cdot u^{\dagger\,1,0}$.
Therefore $\Gr^{W^{\sankaku}}_a$ is a pure twistor of weight $a$,
and $(S(E,P),W^{\sankaku})$ is a mixed twistor.

\subsubsection{The pairing}

The induced pairing
$S:\Gr^{W^{\sankaku}}_{1}\otimes\sigma^{\ast}\Gr^{W^{\sankaku}}_{-1}
\lrarr \Tate(0)$ is given as follows:
\[
 u^{1,0}\otimes\sigma^{\ast}(u^{\dagger\,1,0})\longmapsto 1,
\quad
 u^{\dagger\,0,1}\otimes\sigma^{\ast}(u^{0,1})\longmapsto 1.
\]

\begin{rem}
Note that
we have the following equality on $\cnum_{\lambda}^{\ast}$:
\begin{multline}
 u^{1,0}\otimes \sigma^{\ast}(u^{\dagger\,1,0})
=\mu \cdot u^{\dagger\,0,1}\otimes
 \sigma^{\ast}(-\lambda\cdot u^{\dagger\,0,1})
=\mu\cdot \varphi^{-1}(-\lambda)\cdot u^{\dagger\,0,1}
\otimes\sigma^{\ast}(u^{\dagger\,0,1})   \\ 
=\mu\cdot
 \overline{\bigl(-(-\overline{\lambda})\bigr)}\cdot u^{\dagger\,0,1}
\otimes\sigma^{\ast}(u^{\dagger\,0,1})
=u^{\dagger\,0,1}\otimes\sigma^{\ast}u^{\dagger\,0,1}.
\end{multline}
Thus the pairing $S$ is well defined, of course.
\hfill\qed
\end{rem}

We have induced the pairing
$S\bigl(\nbign^{\sankaku}\otimes\id\bigr):
\Gr_1^{W^{\sankaku}}\otimes \sigma^{\ast}\Gr^{W^{\sankaku}}_{1}
\lrarr \Tate(-1)$.

\begin{lem}
It gives a polarization.
\end{lem}
\pf
Let us consider the induced pairing
$\Gr^{W^{\sankaku}}_1\otimes\nbigo(-1)\otimes
\sigma^{\ast}\bigl(\Gr^{W^{\sankaku}}_1\otimes\nbigo(-1)\bigr)
\lrarr \Tate(0)$.
A base $s$ of the one dimensional vector space
$H^0\bigl(\proj^1,\Gr_1^{W^{\sankaku}}\otimes\nbigo(-1)\bigr)$
is given by the following:
\[
s=
 -\sqrt{-1}\cdot u^{1,0}\otimes f_0^{(-1)}
=u^{\dagger\,0,1}\otimes f_{\infty}^{(-1)}.
\]
Via the pairing
$S\bigl(\nbign^{\sankaku}\otimes\id\bigr)$,
we have the composite of the following correspondences:
\begin{multline}
 \bigl(-\sqrt{-1}\cdot u^{1,0}\otimes f^{(-1)}_0\bigr)
\otimes
 \sigma^{\ast}
 \bigl(u^{\dagger\,0,1}\otimes f^{(-1)}_{\infty}\bigr)
\longmapsto
 -\sqrt{-1}\bigl(u^{0,1}\otimes t_0^{(-1)}\otimes f_0^{(-1)}\bigr)
\otimes
 \sigma^{\ast}
 \bigl(u^{\dagger\,0,1}\otimes f^{(-1)}_{\infty}\bigr)\\
\longmapsto
 -\sqrt{-1}\cdot f_0^{(1)}\otimes\sigma^{\ast}(f_{\infty}^{(-1)})
\longmapsto 1.
\end{multline}
Thus we are done.
\hfill\qed

\vspace{.1in}
In all, we obtain the following.
\begin{cor}
The tuple $\bigl(S(E,P),W,\nbign^{\sankaku},S\bigr)$
is a polarized mixed twistor structure.
\hfill\qed
\end{cor}

\subsubsection{The case $S^{\can}(E)$}

Next we consider the vector bundle $S^{\can}(E)$.
Let $\pi:\hyperh\lrarr\Delta^{\ast}$ be the universal covering
given by $\zeta\longmapsto \exp\bigl(\sqrt{-1}\zeta\bigr)=z$.
Pick a point $\tilde{P}$ such that $\pi(\tilde{P})=P$.
We put as follows:
\[
 \vecs=\vecv\cdot\exp\Bigl(
 -\lambda^{-1}\cdot\bigl(\log z-\log z(\tilde{P})\bigr)
 \left(
 \begin{array}{cc}
  0 & 0 \\ 1 & 0
 \end{array}
 \right)
 \Bigr).
\]
Then $\vecs$ is a frame of the space of multi-valued flat sections,
such that $\vecs(\tilde{P})=\vecv(P)$.

We put as follows:
\[
 \vecs^{\dagger}
=\vecv^{\dagger}\cdot
 \exp\Bigl(
 -\mu^{-1}\cdot\bigl(\log \bar{z}-\log \bar{z}(\tilde{P})\bigr)
 \left(
 \begin{array}{cc}
 0 & 1 \\ 0 & 0
 \end{array}
 \right)
 \Bigr).
\]
Then $\vecs^{\dagger}$ be a frame of the space of
multi-valued sections
such that $\vecs^{\dagger}(\tilde{P})=\vecv^{\dagger}(P)$.
Then we have the following relation:
\[
 \vecs=\vecs^{\dagger}
\left(
 \begin{array}{cc}
 y & -\lambda \\ \mu & 0
 \end{array}
\right).
\]

To see the morphism $\Phi^{\can}$,
we develop $\vecs$ as a polynomial of $\log z$,
and we take the degree $0$-part.
Therefore the morphism $\nbigh\lrarr\nbigg_{|\cnum_{\lambda}^{\ast}}$
is given as follows:
\[
 \vecs\longmapsto
 \vecu\cdot \exp\Bigl(
 \lambda^{-1}\log z(\tilde{P})
 \left(
 \begin{array}{cc}
 0 & 0 \\ 1 & 0
 \end{array}
 \right)
 \Bigr).
\]
Similarly,
we obtain the morphism
$\nbigh\lrarr\nbigg^{\dagger}_{|\cnum_{\lambda}^{\ast}}$:
\[
 \vecs^{\dagger}\longmapsto
 \vecu^{\dagger}\cdot \exp\Bigl(
 \mu^{-1}\log \bar{z}(\tilde{P})
\left( \begin{array}{cc} 0 & 1 \\ 0 & 0
\end{array}
\right)
 \Bigr).
\]
Then the gluing of $S^{\can}(E)$ is given as follows:
\begin{multline}
 \vecu
=\vecu^{\dagger}\cdot
 \exp\left(\mu^{-1}\log \bar{z}(\tilde{P})
 \left( \begin{array}{cc}
 0 & 1 \\ 0 & 0
 \end{array}
 \right)
 \right)\cdot
\left(
\begin{array}{cc}
 y & -\lambda \\ \mu & 0
\end{array}
\right)\cdot
\exp\left(
 -\lambda^{-1}\log z(\tilde{P})
 \left(
 \begin{array}{cc}
 0 & 0 \\ 1 & 0
 \end{array}
 \right) \right) \\
=\vecu^{\dagger}\cdot
 \exp\left(
 \mu^{-1}\log |z(P)|^2
 \left(
 \begin{array}{cc}
 0 & 1 \\ 0 & 0
 \end{array}
 \right)
 \right)\cdot
\left(
 \begin{array}{cc}
 y & -\lambda \\ \mu & 0
 \end{array}
\right) 
=\vecu^{\dagger}\cdot
 \left(
 \begin{array}{cc}
 1 & -\mu^{-1}y \\ 0 & 1
 \end{array}
 \right)\cdot
 \left(
 \begin{array}{cc}
 y & -\lambda \\ \mu & 0
 \end{array}
 \right)
=\vecu^{\dagger}
\left(
 \begin{array}{cc}
 0 & -\lambda \\ \mu & 0
 \end{array}
\right).
\end{multline}
Thus $S^{\can}(E)$ is naturally isomorphic
to $\Gr^{W^{\sankaku}}_1\oplus\Gr^{W^{\sankaku}}_{-1}$.

The filtration $W^{\sankaku}$ gives the mixed twistor structure,
and the pairing $S$ gives the polarization,
which has been already shown.
In all, we obtain the following.
\begin{lem}
The tuples $\bigl(S^{\can},W,\nbign^{\sankaku},S\bigr)$
are polarized mixed twistor structures.
\hfill\qed
\end{lem}

%% file: a88.tex

\subsubsection{The rank $l$ case}

By taking the $(l-1)$-th symmetric product
of the harmonic bundle $Mod(2)$
given in the subsubsection \ref{subsubsection;b12.9.5},
we obtain the harmonic bundle
$Mod(l)$, whose rank is $l$.
Let $\modeldeform(l)$ denote the deformed holomorphic bundle of
$Mod(l)$.

The frame $\vecw$ of
$\modeldeform(l)$ is called canonical,
if the $\lambda$-connection form with respect to $\vecw$
is of the form $A\cdot dz/z$ for some constant nilpotent matrix $A$.
The frame $\vecv$ of $\modeldeform(2)$ given in (\ref{eq;b12.9.10})
induces a canonical frame of $\modeldeform(l)$.

The following lemma is clear.
\begin{lem}\label{lem;b12.10.1}
Let $V$ be a vector space over $\cnum$,
and let $N$ be a nilpotent map of $V$.
We have a harmonic bundle of the form
$\bigoplus_{i}Mod(l_i)$ such that
the residue is isomorphic to $(V,N)$.
Such a harmonic bundle is denoted by $E(V,N)$,
although it is not uniquely determined.
\hfill\qed
\end{lem}

\subsubsection{General model bundles}
\label{subsubsection;b12.9.50}

In general,
the harmonic bundle of the following form is called 
a model bundle:
\[
 \bigoplus_i
 Mod(l_i)\otimes L(u_i).
\]

Let $\lambda_0$ be a point of $\cnum_{\lambda}$.
We have canonical frames of $\nbigl(u_i)$
and $\modeldeform(l_i)$
around $\lambda_0$.
The tensor products of them induce
the frame of $\bigoplus\modeldeform(l_i)\otimes\nbigl(u_i)$.
The induced frame is called a canonical frame.

%% file: 6.tex

\label{subsection;10.10.5}

\subsubsection{Statement}

Let $R>0$ be a positive number.
We put $X=\Delta(R)$, and $D=\{O\}$.
Let $\harmonicbundle$ be a tame harmonic bundle
defined over a neighbourhood of $X-D$ in $\cnum^{\ast}$.
We denote $\theta=f_0\cdot dz/z$.
Assume that we have the decomposition
 $E=\bigoplus_{a\in S_0} E_a$
satisfying the following conditions:
\begin{condition} \mbox{{}}\label{condition;9.6.60}
\begin{enumerate}
\item \label{number;10.10.8}
The endomorphism $f_0$ preserves the decomposition,
that is,
$f_0=\bigoplus_a f_{0\,a}$
for $f_{0\,a}\in End(E_a)$.
\item \label{number;9.6.61}
There exists $C_0>0$ and $\epsilon_0>0$,
such that
 $|b-a|<C_0\cdot |z|^{\epsilon_0}$
 for any eigenvalue $b$ of $f_{0\,a}(z)$.
\item \label{number;9.6.85}
We put $\xi:=\sum_{a\in\Sp(\nbige^0)}\rank(E_a)\cdot|a|^2$.
Then we have $\xi\neq 0$.
\hfill\qed
\end{enumerate}
\end{condition}

We take a total order $\leq_1$ of $S_0:=\Sp(\nbige^0)$.
Then we obtain the filtration $F$ defined as follows:
\[
 F_aE:=\bigoplus_{b\leq_1 a}E_b.
\]
Let $E_a'$ denote the orthogonal complement
of $F_{<a}(E)$ in $F_a(E)$.
We have $\dim(E_a)=\dim(E_a')$,
but $E_a\neq E_a'$ in general.
We put as follows:
\begin{equation} \label{eq;9.15.5}
 \rho:=\bigoplus_{a\in S_0}a\cdot id_{E_a}
\in \bigoplus_{a\in S_0}\End(E_a),
\quad
 \rho':=\bigoplus_{a\in S_0}a\cdot id_{E_a'}
\in \bigoplus_{a\in S_0}\End(E_a').
\end{equation}
Note we have $|\rho'|_h^2=\xi$.

The following proposition is our main purpose
in the subsection \ref{subsection;10.10.5}.
It will be proved in the subsubsections
\ref{subsubsection;10.10.6}--\ref{subsubsection;10.10.7}.

\begin{prop} \label{prop;9.6.77}\mbox{{}}
\begin{description}
\item[(I)]
 Let $R_1$ be a positive number such that $R_1<R$.
 There exists $C_1>0$ satisfying the following inequality
 on $\Delta^{\ast}(R_1)$:
\begin{equation} \label{eq;9.6.75}
 |f_0-\rho'|_h\leq C_1\cdot\Bigl(-\log \frac{|z|}{R}\Bigr)^{-1}
\end{equation}
\item[(II)]
 There exist positive constants $C_2>0$, $\epsilon_2>0$ and $R_2>0$,
 satisfying the following inequality on $\Delta^{\ast}(R_2)$:
\begin{equation}
 \big|\rho-\rho'\big|_h\leq C_2\cdot |z|^{\epsilon_2}
\end{equation}
\end{description}
Here $C_1$, $C_2$, $\epsilon_2$, $R_2$ depends only on
the constants
$R$, $C_0$, $\epsilon_0$, $R_1$, $S_0$, $\rank(E)$
and the set $\bigl\{\rank E_a\,\big|\,a\in S_0\bigr\}$.
\end{prop}

For simplicity, we introduce the following terminology
which is used only in this subsection.
\begin{df}
A constant is called good,
if it depends only on the constants
$R$, $C_0$, $\epsilon_0$, $R_1$, $S_0$, $\rank(E)$
and the set $\bigl\{\rank E_a\,\big|\,a\in S_0\bigr\}$.
\hfill\qed
\end{df}

\begin{rem}
Proposition {\rm \ref{prop;9.6.77}} was proved
by Simpson (Theorem {\rm 1} in {\rm \cite{s2}}).
In fact, we closely follow his argument.
However we would like to care the dependence of constants,
to use the proposition in higher dimensional case.
It is the reason why we give some detail.
\hfill\qed
\end{rem}

\subsubsection{Preliminary}
\label{subsubsection;10.10.6}

We put $f:=f_0\cdot z^{-1}$.
Recall we have the following fundamental inequality due to Simpson.
\begin{lem}[Simpson, the page 729 in \cite{s2}]
The following inequality holds:
\[
 \Delta\log|f|_h^2
\leq
 \frac{-\big|[f,f^{\dagger}]\big|_h^2}{|f|_h^2}.
\]
\hfill\qed
\end{lem}

Due to the condition \ref{number;10.10.8} in Condition \ref{condition;9.6.60},
we have the endomorphisms $f_{0\,a\,|\,Q}\in End(E_{a\,|\,Q})$.
We have the $\EE$-decomposition of $E_{a\,|\,Q}$:
\[
 E_{a\,|\,Q}=
 \bigoplus_{\alpha\in\Sp(f_{0\,a\,|\,Q})}
 \EE(f_{0\,a\,|\,Q},\alpha).
\]
We have the natural bijection
$ \Sp(f_{0\,|\,Q})
\simeq\bigl\{
 (a,\alpha)\,\big|\,
 a\in S_0,\,\,
 \alpha\in \Sp(f_{0\,a\,|\,Q})
 \bigr\}$.
We pick a total order $\leq_2$ on $\Sp(f_{0\,a\,|\,Q})$
on each $a$.
Then we obtain the total order $\leq_3$
on $\Sp(f_{0\,|\,Q})$,
which is given by the lexicographic order
of $\leq_1$ and $\leq_2$.

We obtain the filtration $F^{(1)}$ on $E_{|Q}$
defined as follows:
\[
 F^{(1)}_{(a,\alpha)}(E_{|Q})=
\bigoplus_{(b,\beta)\leq_3 (a,\alpha)}
 \EE(f_{0\,b\,|\,Q},\beta).
\]
Let $H_{(a,\alpha)}$
denote the orthogonal complement of
$F^{(1)}_{<(a,\alpha)}$ in $F^{(1)}_{(a,\alpha)}$.
We have the following:
\[
 E'_{a\,|\,Q}
=\bigoplus_{\alpha\in\Sp(f_{0\,a\,|\,Q})}H_{(a,\alpha)}.
\]
We put as follows:
\[
 \tilde{\rho}_Q:=
 \bigoplus_{(a,\alpha)\in \Sp(f_{0\,|\,Q})}
 \alpha\cdot id_{H_{(a,\alpha)}}.
\]

\begin{lem}\label{lem;9.6.78}
There exists a good constant $C$ and $\epsilon_0$
satisfying the following:
\[
 \big|\tilde{\rho}_Q-\rho'_{|Q}\big|_h\leq C\cdot|z(Q)|^{\epsilon_0}.
\]
\end{lem}
\pf
It follows from the condition \ref{number;9.6.61}
in Condition \ref{condition;9.6.60}.
\hfill\qed

\begin{cor}
There exists a good constant $A_2$
satisfying the following:
\[
 |\tilde{\rho}_Q|_h^2\leq A_2.
\]
\end{cor}
\pf It immediately follows from Lemma \ref{lem;9.6.78}.
\hfill\qed

\vspace{.1in}
We decompose as $f_{|Q}:=\tilde{\rho}_Q\cdot z(Q)^{-1}+g_Q$.
Then we have the following equality:
\[
 |f_{|Q}|_h^2=
 |\tilde{\rho}_Q|_h^2\cdot |z(Q)|^{-2}
+|g_Q|_h^2.
\]
\begin{lem}
There exists a good constant $A_1$ satisfying the following:
\begin{equation} \label{eq;9.6.90}
 \big|\big[f_{|Q},\,\,f^{\dagger}_{|Q}\big]\big|_h
 \geq
 A_1\cdot |g_Q|^2.
\end{equation}
\end{lem}
\pf
It can be shown by an elementary linear algebraic argument.
\hfill\qed

\subsubsection{Step 1 of the proof of (I)}

\begin{lem} \label{lem;9.6.70}
There exist good constants $A_3$ and $A_4$
satisfying the following:
\begin{quote}
 For any point $Q\in\Delta^{\ast}(R)$,
 one of the following holds:
\begin{itemize}
\item
 $|f(Q)|_h^2\leq A_3\cdot|z(Q)|^{-2}$.
\item
 $\Delta \log|f|_h^2(Q)\leq -A_4\cdot |f|_h^2(Q)$.
\end{itemize}
\end{quote}
\end{lem}
\pf
Assume that $|f_{|Q}|_h^2\geq 2\cdot A_2\cdot |z(Q)|^{-2}$.
Then we obtain
$|g_Q|_h^2\geq 2^{-1}\cdot |f_{|Q}|_h^2$.
Hence we obtain the following:
\[
  \bigl( \Delta\log|f|_h^2\bigr)(Q)
\leq
 \frac{-\big|\big[f_{|Q},\,f^{\dagger}_{|Q}\big]\big|_h^2}{|f_{|Q}|_h^2}
\leq
 \frac{-A_1^2\cdot \big|g_Q\big|_h^4}{|f_{|Q}|_h^2}
\leq
 \frac{-A_1^2}{4}\cdot |f_{|Q}|_h^2.
\]
Thus we may take $A_3=2A_2$ and $A_4=4^{-1}A_1^2$.
\hfill\qed

\vspace{.1in}

Let $\eta$ and $B$ be positive numbers.
We put as follows:
\[
 m_{\eta,B}:=
 \frac{B}{(|z|-\eta)^2\cdot(|z|-R)^2}.
\]
It is a $C^{\infty}$-function on the region
$\bigl\{z\,\big|\,\eta<|z|<R\bigr\}$.
Note that $m_{\eta,B}$ is $\infty$
on the boundary
$\bigl\{z\,\big|\,|z|=\eta\bigr\}\cup\bigl\{z\,\big|\,|z|=R\bigr\}$.

\begin{lem}
Let us take a positive number $B$ satisfying the following inequalities:
\[
 B\geq \frac{4R^2}{A_4},
\quad
 B>A_3\cdot R^2.
\]
In particular, $B$ is a good constant.
Then we obtain the following inequalities:
\[
 \Delta\log m_{\eta,B}\geq -A_4\cdot m_{\eta,B},
\quad
 m_{\eta,B}\geq A_3\cdot|z|^{-2}.
\]
\end{lem}
\pf
We have the following formula:
\begin{equation} \label{eq;9.6.65}
 \Delta\log m_{\eta,B}
=\frac{-2\eta}{|z|\cdot \bigl(|z|-\eta\bigr)^2}
+\frac{-2R}{|z|\cdot\bigl(|z|-R\bigr)^2}
=\frac{-2}{B}
 \Bigl(
 \frac{\bigl(|z|-R\bigr)^2}{|z|}\cdot\eta
 +\frac{\bigl(|z|-\eta\bigr)^2}{|z|}\cdot R
 \Bigr)
\cdot m_{\eta,B}.
\end{equation}
From $\eta<|z|<R$,
we have the following inequalities:
\[
  \frac{\eta}{|z|}<1,
\quad
 (|z|-R)^2\leq R^2,
\quad
 \frac{(|z|-\eta)^2}{|z|}
\leq |z|-\eta\leq R.
\]
Thus we have the following inequality:
\begin{equation} \label{eq;9.6.66}
 \frac{\eta}{|z|}\cdot\bigl(|z|-R\bigr)^2
+\frac{\bigl(|z|-\eta\bigr)^2}{|z|}\cdot R
 \leq 2R^2.
\end{equation}
From (\ref{eq;9.6.65}) and (\ref{eq;9.6.66}),
we obtain the following inequality
on the region $\bigl\{z\,\big|\,\eta<|z|<R\bigr\}$:
\[
 \Delta\log m_{\eta,B}\geq
-\frac{4}{B}\cdot R^2\cdot m_{\eta,B}.
\]
For any $B\geq \frac{4R^2}{A_4}$ and for any $\eta>0$,
we obtain the first inequality
on the region $\{\eta<|z|<R\}$.

We also have the following inequalities:
\[
 m_{\eta,B}
=\frac{B}{(|z|-\eta)^2\cdot(|z|-R)^2}
\geq
 \frac{A_3\cdot R^2}{(|z|-\eta)^2\cdot(|z|-R)^2}
\geq
 \frac{A_3}{(|z|-\eta)^2}
\geq \frac{A_3}{|z|^2}.
\]
Hence we obtain the second inequality.
\hfill\qed

\vspace{.1in}
We put
$S_1:=
 \bigl\{z\in\Delta^{\ast}(R)\,\big|\,
  |f(z)|_h^2>m_{\eta,B}(z) \bigr\}$.
\begin{lem}
The set $S_1$ is empty.
In other words,
we have the inequality
$|f(z)|_h^2\leq m_{\eta,B}(z)$
for any point $z\in\Delta^{\ast}(R)$ such that $|z|>\eta$.
\end{lem}
\pf
Assume that $S_1$ is not empty,
and we will derive a contradiction.
On the region $S_1$,
we have the inequality
$|f|_h^2>m_{\eta,B}\geq A_3\cdot |z|^{-2}$.
Hence we have $\Delta\log|f|_h^2\leq -A_4\cdot|f|_h^2$
due to Lemma \ref{lem;9.6.70}.
Then we obtain the following inequality on $S_1$:
\[
 \Delta\log\bigl(|f|_h^2\big/m_{\eta,B}\bigr)
\leq
 -A_4\cdot \bigl(|f|_h^2-m_{\eta,B}\bigr)<0.
\]
It implies that the function $\log\bigl(|f|_h^2/m_{\eta,B}\bigr)$
cannot take the maximal in the region $S_1$.

Let $\bar{S}_1$ denote the closure of $S_1$ in $\cnum$.
Since we have $m_{\eta}=\infty$ 
on the boundary $ \{|z|=\eta\}\cup\{|z|=R\}$,
the intersection 
of the sets $\bar{S}_1$
and $\bigl\{z\,|\,|z|=\eta\bigr\}\cup \bigl\{|z|=R\bigr\}$ is empty.
Hence we have $|f|_h^2=m_{\eta,B}$
on the boundary of $\bar{S}_1$.
Hence we obtain $|f|_h^2\leq m_{\eta,B}$
on $S_1$,
but it contradicts the definition of $S_1$.
Hence we obtain $S_1=\emptyset$.
\hfill\qed

\vspace{.1in}

When we take a limit $\eta\to 0$,
we obtain the inequality on $\Delta^{\ast}(R)$:
\[
 \bigl|f(z)\bigr|_h^2\leq
 \frac{B}{|z|^2\cdot\bigl(|z|-R\bigr)^2}.
\]

Hence, we have arrived at the following:
\begin{lem}\label{lem;a11.13.1}
For any $R_3$ such that $R_1<R_3<R$,
there exists a good constant $A_5$
such that the following inequality
holds on $\Delta^{\ast}(R_3)$:
\[
 \bigl|f(z)\bigr|_h^2\leq \frac{A_5}{|z|^2}.
\]
\hfill\qed
\end{lem}

\subsubsection{Step 2 for the proof of (I)}

We put as follows:
$k:=\log|f|_h^2-\log\bigl(\xi\cdot|z|^{-2}\bigr)$.
Then we have the following:
\[
 k(Q)=
\log\Bigl(
 \frac{|z(Q)|^2}{\xi}\cdot |f_{|Q}|^2
\Bigr)
=\log
 \Bigl(
 \frac{|z(Q)|^2}{\xi}
 \cdot
 \bigl(
  |\tilde{\rho}_Q|_h^2\cdot|z(Q)|^{-2}
+|g_Q|_h^2
 \bigr)
 \Bigr).
\]
We put $b_Q:=|\tilde{\rho}_Q|^2-\xi$.

\begin{lem} \label{lem;9.6.79}
We have a good constant $A_6$
such that $|b_Q|\leq A_6\cdot |z(Q)|^{\epsilon_0}$.
\end{lem}
\pf
This is a reformulation of Lemma \ref{lem;9.6.78}.
\hfill\qed

\begin{lem}
There exists a good constant $A_7$
satisfying the following for any point $Q\in\Delta^{\ast}(R_3)$:
\begin{equation}\label{eq;9.6.81}
 A_7\cdot\Bigl(
 \xi^{-1}\cdot b_Q+\frac{|z(Q)|^2}{\xi}\cdot|g_Q|_h^2
 \Bigr)
\leq k(Q)
\leq
 \xi^{-1}\cdot b_Q+\frac{|z(Q)|^2}{\xi}\cdot|g_Q|_h^2.
\end{equation}
\end{lem}
\pf
We have the following description:
\[
 k(Q)=\log
 \Bigl(
  1+\xi^{-1}\cdot b_Q+\frac{|z(Q)|^2}{\xi}\cdot|g_Q|_h^2
 \Bigr).
\]
Then the right inequality is obvious.
We have only to obtain the left inequality.

Recall we have obtained the following estimate on $\Delta^{\ast}(R_3)$
in Step 1:
\[
 |g_Q|_h^2\leq |f_{|Q}|_h^2
\leq A_5\cdot |z(Q)|^{-2}.
\]
Hence we have the following inequality:
\[
 0\leq \frac{|z(Q)|^2}{\xi}\cdot |g_Q|_h^2
\leq \frac{A_5}{\xi}.
\]
Then
we obtain the left inequality for some good constant,
for example,
by using the convexity of the logarithmic function.
\hfill\qed

\vspace{.1in}
We will show Lemma \ref{lem;9.6.80} later.
\begin{lem} \label{lem;9.6.80}
There exists a good constant $A_8$
satisfying the following:
\[
 k\leq A_8\cdot \left(-\log\frac{|z|}{R_3}\right)^{-2}.
\]
\end{lem}

\begin{lem}
Lemma {\rm\ref{lem;9.6.80}}
implies the claim (I) of Proposition {\rm\ref{prop;9.6.77}}.
\end{lem}
\pf
Assume that we have shown Lemma \ref{lem;9.6.80},
then we obtain the following inequality on $\Delta^{\ast}(R_1)$:
\[
 |f_0-\rho'|_h^2
=b_Q+|z(Q)|^2\cdot|g_Q|_h^2
\leq \xi\cdot A_7^{-1}\cdot k
\leq
 \xi\cdot A_7^{-1}\cdot
 A_8\cdot \Bigl(-\log\frac{|z|}{R_3}\Bigr)^{-2}.
\]
It implies the desired inequality (\ref{eq;9.6.75})
in (I).
Thus we have reduced the proof of the claim (I)
to Lemma \ref{lem;9.6.80}.
\hfill\qed

\vspace{.1in}
\subsubsection{Proof of Lemma \ref{lem;9.6.80}}
Let us prove Lemma \ref{lem;9.6.80}.
There exists a good constant $A_{10}$
satisfying the following inequality on $\Delta^{\ast}(R_3)$:
\[
 \xi^{-1}\cdot A_6\cdot|z|^{\epsilon_0}
\leq
 \frac{1}{2}A_{10}
 \left(-\log\frac{|z|}{R_3}\right)^{-2}.
\]
Here $A_6$ appeared in Lemma \ref{lem;9.6.79}.

Take a good constant $A_{11}$ satisfying the following:
\begin{equation} \label{eq;9.6.92}
 A_{11}<\frac{A_1^2\cdot\xi^2}{4\cdot A_5},
\quad
A_{11}<\frac{6}{A_{10}}.
\end{equation}
The first condition will be used in Lemma \ref{lem;9.6.91}
and the second condition will be used to obtain
the inequality (\ref{eq;9.6.93}).

\begin{lem} \label{lem;9.6.91}
One of the following holds:
\begin{itemize}
\item
 ${\displaystyle
 k(Q)<A_{10}\cdot\Bigl(\log\frac{|z(Q)|}{R_3}\Bigr)^{-2}
 }$
\item
 $(\Delta k)(Q)< -A_{11}\cdot k(Q)^2\cdot |z(Q)|^{-2}$.
\end{itemize}
\end{lem}
\pf
Assume $k(Q)\geq A_{10}\cdot\big(-\log\frac{|z(Q)|}{ R_3}\big)^{-2}$.
Then we obtain the following:
\[
 \frac{1}{2}k(Q)\geq \xi^{-1}\cdot A_6\cdot |z|^{\epsilon_0}
 \geq \xi^{-1}\cdot |b_Q|.
\]
Namely we obtain the following:
\[
 k(Q)-\xi^{-1}\cdot b_Q
\geq \frac{k(Q)}{2}.
\]
By using the right inequality in (\ref{eq;9.6.81}),
we obtain the following inequality:
\[
 \frac{|z(Q)|^2}{\xi}\cdot |g_Q|_h^2
\geq
 \frac{k(Q)}{2}.
\]
Thus we obtain the following:
\begin{equation}\label{eq;9.6.91}
 -|g_Q|_h^4\leq \frac{-\xi^2}{4\cdot|z(Q)|^4}k(Q)^2.
\end{equation}

Note the following:
\[
 \Delta k=\Delta\log|f|_h^2
\leq \frac{-\big|\big[f,f^{\dagger}\big]\big|_h^2}{|f|_h^2}.
\]
Then we obtain the following inequality
by using Lemma \ref{lem;a11.13.1},
the inequality (\ref{eq;9.6.90}),
the left inequality in (\ref{eq;9.6.92}),
and the inequality (\ref{eq;9.6.91}):
\[
 (\Delta k)(Q)\leq
-A_1^2\cdot\frac{\xi^2}{|z(Q)|^4\cdot 4}\cdot k(Q)^2\cdot
 \frac{|z(Q)|^2}{A_5}
=\frac{-A_1^2\cdot\xi^2\cdot k(Q)^2}{4\cdot A_5\cdot|z(Q)|^2}
< -A_{11}\cdot \frac{k(Q)^2}{|z(Q)|}.
\]
Thus we are done.
\hfill\qed

\vspace{.1in}

For any positive numbers $\epsilon$ and $B$,
we put as follows:
\[
 p_{B,\epsilon}=
 B\cdot \Bigl(
 -\log\frac{|z|}{R_3}\Bigr)^{-2}
+\epsilon\cdot \left(-\log\frac{|z|}{R_3}\right).
\]

\begin{lem}
We put $B=6\cdot A_{11}^{-1}$.
We have the following inequalities:
\begin{equation} \label{eq;9.6.100}
 \Delta p_{B,\epsilon}\geq -A_{11}\cdot\frac{p_{B,\epsilon}^2}{|z|^2},
\quad\quad
 p_{B,\epsilon}>A_{10}\cdot\Bigl(-\log \frac{|z|}{R_3}\Bigr)^{-2}. 
\end{equation}
\end{lem}
\pf
We have the following formula,
where we use the real coordinate $z=r\cdot e^{\sqrt{-1}\theta}$:
\[
  \Delta p_{B,\epsilon}=
-\left(
 \frac{\del}{\del r}+\frac{1}{r}
 \right)
 \left(
 \frac{2B}{(-\log r/R_3)^3\cdot r}
 \right)
=\frac{-6B}{(-\log r/R_3)^4\cdot r^2}
\geq
 \frac{-6p_{B,\epsilon}^2}{B\cdot r^2}.
\]
Here we have used the inequality:
\[
 \frac{p_{B,\epsilon}^2}{r^2}
\geq
 \frac{B^2}{(-\log r/R_3)^2\cdot r^2}.
\]
Since we put $B=6\cdot A_{11}^{-1}$,
we obtain the first inequality in (\ref{eq;9.6.100}).

Note that we have the following, due to the second inequality
in (\ref{eq;9.6.92}):
\begin{equation}\label{eq;9.6.93}
B=6\cdot A_{11}^{-1}>A_{10}.
\end{equation}
Thus we obtain the second equality
in (\ref{eq;9.6.100}).
\hfill\qed

\vspace{.1in}

We put $ S_2:=
 \bigl\{Q\,\big|\,  k(Q)> p_{\epsilon}(Q)\bigr\}$.
\begin{lem} \label{lem;10.10.10}
The set $S_2$ is empty.
In other words, we have $k(Q)\leq p_{\epsilon}(Q)$
for any point $Q\in\Delta^{\ast}$.
\end{lem}
\pf
We assume that $S_2$ is not empty,
and we will derive a contradiction.
Let $Q$ be a point of $S_2$,
and then 
we have the inequality
$k(Q)>
 p_{\epsilon}(Q)\geq A_{10}\cdot\big(-\log\frac{|z(Q)|}{R_3}\big)^{-2}$.
Then we obtain the following inequality
due to Lemma \ref{lem;9.6.91}:
\[
 \Delta k(Q)< -A_{11}\frac{|k(Q)|^2}{|z(Q)|^2}.
\]
Hence we have the following inequality:
\[
 \Delta(k-p_{\epsilon})(Q)
<
 -A_{11}\frac{k(Q)^2-p_{\epsilon}(Q)^2}{|z(Q)|^2}< 0.
\]
It means that the function $k-p_{\epsilon}$
does not have any maximal point in the region $S_2$.
Since we have $p_{\epsilon}=\infty$ 
on the boundary $\{|z|=0\}\cup \{|z|=R_3\}$,
the intersection of the sets
$S_2$ and $\{|z|=0\}\cup \{|z|=R_3\}$
is empty.
Hence we have $k(Q)=p_{\epsilon}(Q)$
on the boundary of $S_2$.
Thus
we obtain the inequality $k\leq p_{\epsilon}$ on the region $S_2$,
which contradicts the definition of $S_2$.
\hfill\qed

\vspace{.1in}
Let us return to the proof of Lemma \ref{lem;9.6.80}.
We obtain $k(Q)\leq p_{B,\epsilon}$
for any positive number $\epsilon$ due to Lemma \ref{lem;10.10.10}.
Thus we obtain the following inequality
for a good constant:
\[
 k\leq B\cdot\Bigl(-\log\frac{|z|}{R_3}\Bigr)^{-2}.
\]
Therefore the proof of Lemma \ref{lem;9.6.80} 
and the proof of the claim (I) are accomplished.
\hfill\qed

\subsubsection{The proof of (II)}
\label{subsubsection;10.10.7}

For any element $l\in\bigoplus_{a\leq_1 b}Hom(E_b',E_a')$,
the adjoint $\ad(l)$ induces the endomorphism of
$\bigoplus_{a<_1b}Hom(E_b',E_a')$,
which we denote by $F_l$.
Recall that we put $f_0=f\cdot z$.

\begin{lem}\label{lem;9.7.1}
There exist good constants $R_4$ and $A_{14}$
such that $F_{f_0}$ is invertible
on $\Delta^{\ast}(R_4)$,
and the norms of $F_{f_0}$ and $F_{f_0}^{-1}$
are  dominated by $A_{14}$.
\end{lem}
\pf
We put $\tilde{g}:=f_0-\rho'$,
which is a section of $\bigoplus_{a\leq_1 b}Hom(E_b',E_a')$.
Then we have the following inequality on $\Delta^{\ast}(R_1)$
due to the claim (I):
\begin{equation} \label{eq;9.7.3}
 |\tilde{g}|_h
\leq
 C_1\cdot\left(-\log\frac{|z|}{R}\right)^{-1}.
\end{equation}
Hence the norm of the endomorphism $F_{\tilde{g}}$ is dominated
by $A_{13}\cdot \bigl(-\log(|z|\cdot R^{-1})\bigr)^{-1}$
for a good constant $A_{13}$.
On the other hand,
the endomorphism $F_{\rho'}$ is invertible
and the norms of $F_{\rho'}$ and $F_{\rho'}^{-1}$
are dominated by a good constant $A_{12}$.
Thus we obtain Lemma \ref{lem;9.7.1}.
\hfill\qed

\vspace{.1in}

We put $q:=\rho-\rho'$,
which is an element of $\bigoplus_{a<_1b}Hom(E_b',E_a')$.

\begin{lem} \label{lem;9.7.2}
There exists a good constant $A_{15}$ such that
$|q|_h\leq A_{15}\cdot \bigl(-\log(|z|\cdot R^{-1})\bigr)^{-1}$
on the region $\Delta^{\ast}(R_4)$.
\end{lem}
\pf
We have the following equality by using $[\rho',\rho']=0$:
\[
0=[f_0,\rho]=[f_0,\rho'+q]
 =F_{f_0}(q)+[\rho'+\tilde{g},\rho']
 =F_{f_0}(q)+[\tilde{g},\rho'].
\]
Then we obtain Lemma \ref{lem;9.7.2}
by using Lemma \ref{lem;9.7.1} and (\ref{eq;9.7.3}).
\hfill\qed

\begin{lem} \label{lem;9.7.4}
There exist good constants $R_5$ and $A_{16}$
such that the following inequality holds
on the region $\Delta^{\ast}(R_5)$:
\[
 \big|[\rho,f^{\dagger}]\big|_h^2
\geq
 A_{16}\cdot|q|_h^2\cdot|z|^{-2}.
\]
\end{lem}
\pf
We have the following equality by a direct calculation:
\[
 [\rho,f^{\dagger}]=
 [\rho',f^{\dagger}]
+[q,\,\bar{z}^{-1}\!\cdot\!\bar{\rho}']
+[q,\,\bar{z}^{-1}\!\cdot\!\tilde{g}^{\dagger}].
\]
Here we put $\bar{\rho}'=\bigoplus_{a\in S_0}\bar{a}\cdot id_{E_a'}$,
and we have used the relation $\rho^{\prime\,\dagger}=\bar{\rho}'$.
Note the following:
\[
 \bigl[\rho',f^{\dagger}\bigr]
\in
 \bigoplus_{a>_1b}
 Hom(E'_{b},E'_a),
\quad\quad
 \bigl[q,\,\bar{z}^{-1}\!\cdot\!\bar{\rho}'\bigr]
\in
 \bigoplus_{a<_1b}
 Hom(E'_{b},E'_a).
\]
There exists good constants $C>0$ and $C'>0$ satisfying the following:
\[
\begin{array}{l}
  \big|[q,\bar{z}^{-1}\cdot\tilde{g}^{\dagger}]\big|_h
\leq
 C'\cdot|z|^{-1}\cdot|\ad(\tilde{g})|_h\cdot|q|_h
\leq
 C\cdot \bigl(-|z|\log\frac{|z|}{R}\bigr)^{-1}\cdot |q|_h,\\
 \mbox{{}}\\
 \big|[q,\bar{z}^{-1}\bar{\rho}']\big|_h
\geq C\cdot|z|^{-1}|q|_h.
\end{array}
\]
In all, we obtain Lemma \ref{lem;9.7.4}.
\hfill\qed

\vspace{.1in}

Recall the following inequality
due to Simpson (in the page 731 of \cite{s2}):
\begin{equation}\label{eq;9.7.5}
 \Delta\log|\rho|_h^2
\leq
 \frac{-\big|[\rho,f^{\dagger}]\big|_h^2}{|\rho|_h^2}.
\end{equation}

\begin{lem}
There exist good constants
$R_6$ and $A_{17}$ such that
the following inequality holds on $\Delta^{\ast}(R_{6})$:
\[
 \Delta\log|\rho|_h^2
\leq
 -A_{17}\cdot|q|_h^2\cdot |z|^{-2}.
\]
\end{lem}
\pf
It follows from (\ref{eq;9.7.5}) and Lemma \ref{lem;9.7.4}.
\hfill\qed

\vspace{.1in}

Since we have $|\rho|_h^2=\xi +|q|_h^2$,
we obtain the following inequality on $\Delta^{\ast}(R_6)$:
\begin{equation} \label{eq;9.7.10}
 \Delta\log(1+\xi^{-1}\cdot|q|_h^2)
\leq
 -A_{17}\cdot|q|_h^2\cdot |z|^{-2}.
\end{equation}

\begin{lem}
There exists a good constant $A_{18}$
such that the following holds:
\begin{equation} \label{eq;9.7.6}
 A_{18}\cdot\xi^{-1} \cdot|q|_h^2
\leq
 \log(1+\xi^{-1}\cdot|q|_h^2)
\leq
 \xi^{-1}\cdot|q|_h^2.
\end{equation}
\end{lem}
\pf
The right in the inequality (\ref{eq;9.7.6}) is clear.
Since we have $|q|_h\leq A_{15}\cdot\bigl(-\log(|z|\cdot R^{-1})\bigr)^{-1}$
on $\Delta^{\ast}(R_6)$
due to Lemma \ref{lem;9.7.2},
we have the following:
\begin{equation} \label{eq;9.7.9}
|q|_h^2\cdot\xi^{-1}
\leq \xi^{-1}\cdot A_{15}^2\cdot\Bigl(-\log\frac{|z|}{R}\Bigr)^{-2}.
\end{equation}
In particular,
we have a good constant $C>0$ such that
$0\leq |q|_h^2\cdot\xi^{-1}\leq C$ on $\Delta^{\ast}(R_6)$.
Thus there exists a good constant for the left inequality
in (\ref{eq;9.7.6}).
\hfill\qed

\vspace{.1in}

We put $k:=\log(1+|z|^2\cdot\xi^{-1}\cdot|q|_h^2)$.
\begin{lem}
There exists a good constant $A_{19}$
such that the following holds:
\begin{equation} \label{eq;9.7.16}
 \Delta k\leq -A_{19}\cdot k\cdot|z|^{-2}.
\end{equation}
\end{lem}
\pf
It follows from (\ref{eq;9.7.10})
and the left inequality in (\ref{eq;9.7.6}).
\hfill\qed

\begin{lem}
We have the following inequality on $\Delta^{\ast}(R_6)$:
\[
  k\leq |z|^2\cdot\xi^{-1}\cdot|q|_h^2\leq
 A_{15}\cdot\xi^{-1}\cdot\bigl(-\log(|z|\cdot R^{-1})\bigr)^{-2}.
\]
\end{lem}
\pf
It immediately follows from the right inequality in (\ref{eq;9.7.6})
and (\ref{eq;9.7.9}).
\hfill\qed

\begin{cor}
There exists a good constant $A_{20}$
such that $k\leq A_{20}$.
\hfill\qed
\end{cor}

For positive numbers $B$, $\epsilon$ and $u$,
we put as follows:
\[
 p_{B,\epsilon,u}
:=B\cdot\Bigl(|z|^u+\epsilon\cdot\Big(-\log\frac{|z|}{R}\Big)\Bigr).
\]
It is easy to check
$\Delta p_{B,\epsilon,u}=-u^2\cdot|z|^{u-2}\cdot B$.
\begin{lem} \label{lem;9.7.15}
\mbox{{}}
\begin{itemize}
\item
Take $u>0$ satisfying $u^2<A_{19}$.
Then we have the following inequality:
\begin{equation} \label{eq;9.7.17}
 \Delta p_{B,\epsilon,u}
>-A_{19}\cdot |z|^{-2}\cdot |z|^u\cdot B
>-A_{19}\cdot|z|^{-2}\cdot p_{B,\epsilon,u}.
\end{equation}
\item
Fix $u>0$.
Take $B>0$ as $B\cdot R_6^u>A_{20}$.
Then we obtain the following:
\begin{equation} \label{eq;9.7.19}
 p_{B,\epsilon,u}(R_6)=
 B\cdot \Big(R_6^u+\epsilon\cdot \Big(-\log\frac{|z|}{R_6}\Big)\Big)
>A_{20}>k(R_6).
\end{equation}
\hfill\qed
\end{itemize}
\end{lem}

Let us fix good constants $u$ and $B$ as in Lemma 
\ref{lem;9.7.15}.
We use the notation $p_{\epsilon}$ instead of $p_{B,\epsilon,u}$.
Then we obtain the following inequality
from (\ref{eq;9.7.16}) and (\ref{eq;9.7.17}):
\begin{equation}\label{eq;9.7.18}
 \Delta(k-p_{\epsilon})
< -\frac{A_{19}}{|z|^2}\cdot(k-p_{\epsilon}).
\end{equation}

We put
$ S_{3}:=
 \bigl\{Q\in\Delta^{\ast}(R_0)\,\big|\,
 k(Q)-p_{\epsilon}> 0
 \bigr\}$.
\begin{lem}
The set $S_3$ is empty.
In other words,
we have $k(Q)\leq p_{\epsilon}(Q)$ for any point
$Q\in\Delta^{\ast}$.
\end{lem}
\pf
Assume that $S_3$ is not empty,
and we will derive a contradiction.
The function $k-p_{\epsilon}$ has no maximal point
in the region $S_3$ due to the inequality (\ref{eq;9.7.18}).
Since we have $p_{\epsilon}(R_6)>k(R_6)$ due to 
the inequality (\ref{eq;9.7.19}),
and since we have $p_{\epsilon}(0)=\infty$ by definition,
the intersection of the sets $S_3$
and $\{|z|=0\}\cup \{|z|=R_6\}$ is empty.
Hence we have $k=p_{\epsilon}$
on the boundary of the closure $\bar{S}_3$.
Hence we obtain $k\leq p_{\epsilon}$ on $S_3$,
which is a contradiction.
Hence the set $S_3$ is empty.
\hfill\qed

\vspace{.1in}

When we take limit $\epsilon\to 0$,
we obtain the inequality
$k\leq B\cdot r^u$
on $\Delta^{\ast}(R_6)$.
Then there exists a good constant $A_{21}$
such that the following inequality holds on $\Delta^{\ast}(R_6)$:
\[
 |q|_h\leq A_{21}\cdot |z|^u.
\]
Thus the proof of the claim (II) is accomplished.
\hfill\qed

\subsubsection{Some consequences and the asymptotic orthogonality}

\begin{cor} \label{cor;9.7.110}
 For any $R_1<R$,
 there exists a good constant $C$ such that the following holds
 on $\Delta^{\ast}(R_1)$:
\[
 |f_0-\rho|_h\leq C\cdot\Big(-\log\frac{|z|}{R}\Big)^{-1}.
\]
\end{cor}
\pf
It follows from the estimate of $q=\rho-\rho'$
and (I) in Proposition \ref{prop;9.6.77}.
\hfill\qed

\vspace{.1in}

We have the decomposition
$f_0=\sum_{a\geq_1 b}\tilde{f}_{0,a,b}$,
where $\tilde{f}_{0,a,b}\in Hom(E_a',E_b')$.
\begin{cor} \label{cor;9.7.20}
 There exist good constants $C>0$ and $\epsilon>0$
 such that the following holds:
  \begin{itemize}
 \item
   $|\tilde{f}_{0\,a\,b}|_h\leq C\cdot|z|^{\epsilon}$
   in the case $a\neq b$.
 \item
   ${\displaystyle
     |\tilde{f}_{0\,a\,a}-a\cdot id_{E_a'}|_h
     \leq C\cdot\Bigl(-\log\frac{|z|}{R}\Bigr)^{-1}  }$.
  \end{itemize}
\end{cor}
\pf
The second inequality immediately follows 
from (I) in Proposition \ref{prop;9.6.77}.
Due to the commutativity $[f_0,\rho]=0$,
we have $[f_0,\rho']+[f_0,q]=0$.
We have the following:
\[
 [f_0,\rho']=\sum_{a>_1b}(b-a)\cdot\tilde{f}_{0\,a\,b}.
\]
On the other hand, we have the estimate
\[
 \big|[f_0,q]\big|_h\leq C\cdot |f_0|_h\cdot|q|_h
 \leq C'|z|^{\epsilon_2}.
\]
Hence we obtain the estimates for $\tilde{f}_{0\,a\,b}$.
\hfill\qed

\vspace{.1in}

We have the decomposition
$f_0^{\dagger}=\sum_{a,b}\tilde{f}^{\dagger}_{0\,a\,b}$,
where $\tilde{f}^{\dagger}_{0\,a\,b}\in Hom(E_a',E_b')$.
\begin{cor}\label{cor;9.7.21}
 There exists good constants $C$ and $\epsilon$
 such that the following holds:
 \begin{itemize}
 \item
   $|\tilde{f}^{\dagger}_{0\,a\,b}|_h\leq C\cdot|z|^{\epsilon}$
   in the case $a\neq b$.
 \item
   $|\tilde{f}^{\dagger}_{0\,a\,a}-\bar{a}\cdot id_{E_a'}|_h
     \leq C\cdot(-\log|z|)^{-1}$.
 \end{itemize}
\end{cor}
\pf
It immediately follows from Corollary \ref{cor;9.7.20}.
\hfill\qed

\vspace{.1in}

\label{eq;9.15.7}
We put $\bar{\rho}:=\bigoplus_{a\in S_0} \bar{a}\cdot id_{E_a}$.
We also put
$\bar{\rho}':=\bigoplus_{a\in S_0}\bar{a}\cdot id_{E_a'}$.
Note that $\bar{\rho}'$ is adjoint of $\rho'$,
and that
$\bar{\rho}-\bar{\rho}'\in \bigoplus_{a<_1b}Hom(E_{b}',E_{a}')$.

\begin{lem} \label{lem;9.7.25}
There exists good constants $C'$, $R'$ and $\epsilon'$
such that the following holds on $\Delta^{\ast}(R')$:
\[
 \bigl|
 \bar{\rho}-\bar{\rho}'
 \bigr|_h
\leq C'\cdot |z|^{\epsilon'}.
\]
\end{lem}
\pf
The argument is similar to the proof of Corollary \ref{cor;9.7.20}.
We have the following formula:
\[
0=[\bar{\rho},\rho]
 =[\bar{\rho},\rho']+[\bar{\rho},q]
 =\bigl[\bar{\rho}-\bar{\rho}',\bar{\rho'}
  \bigr]
 +[\bar{\rho},q].
\]
There exist good positive constants $R''$, $C''$, $C'''$ and $\epsilon''$
satisfying the following inequalities on $\Delta^{\ast}(R'')$:
\[
 \big|[\bar{\rho},q]\big|_h\leq C''\cdot|z|^{\epsilon''},
\quad\quad
 \bigl|
 \bigl[
 \bar{\rho}-\bar{\rho}',
 \rho'
 \bigr]
 \bigr|_h
\geq C''\cdot
 \bigl|\bar{\rho}-\bar{\rho}' \bigr|_h.
\]
Thus we are done.
\hfill\qed

\begin{cor}
There exist good positive constants $C'$ and $R'$
such that the following holds on $\Delta^{\ast}(R')$:
\[
 \big|f^{\dagger}_0-\bar{\rho}\big|_h
 \leq
 C\cdot|z|^{-1}\cdot\Bigl(-\log\frac{|z|}{R}\Bigr)^{-1}.
\]
\end{cor}
\pf
We have
 $f_0^{\dagger}-\bar{\rho}=(f_0^{\dagger}-\bar{\rho}')
 +(\bar{\rho}'-\bar{\rho})
=(f_0^{\dagger}-\rho^{\prime\,\dagger})+(\bar{\rho}'-\bar{\rho})$.
Then we obtain the result
from (I) in Proposition \ref{prop;9.6.77}
and Lemma \ref{lem;9.7.25}.
\hfill\qed

\vspace{.1in}

In general,
let us consider an element $g\in\bigoplus_a End(E_a)$.
Then we have the decomposition:
\[
 g=\sum_{a\geq_1 b}g_{a\,b},
\quad\quad
 g_{a\,b}\in Hom(E_a',E_b').
\]
\begin{lem} \label{lem;9.7.26}
We have the estimate
$|g_{a\,b}|_h\leq C\cdot|z|^{\epsilon}\cdot|g|_h$
for $a\neq b$.
\end{lem}
\pf
We have $0=[g,\rho]=[g,\rho']+[g,q]$.
We have an estimate
$\big|[g,q]\big|_h\leq C'\cdot |g|_h\cdot |z|^{\epsilon'}$
on $\Delta^{\ast}(R')$
for some good positive constants $C'$, $\epsilon'$ and $R'$.
On the other hand,
we have
\[
 [g,\rho']
=\sum_{a>_1b}(b-a)\cdot g_{a,b}.
\]
Thus we are done.
\hfill\qed

\vspace{.1in}
For the endomorphism $g$ above,
we also obtain the adjoint
$g^{\dagger}\in End(E)$,
and we have the decomposition:
\[
 g^{\dagger}=
\sum_{a\leq_1 b}(g^{\dagger})_{a\,b},
\quad\quad
 (g^{\dagger})_{a\,b}\in Hom(E_a',E_b').
\]
\begin{lem}
We have
 $|(g^{\dagger})_{a\,b}|_h\leq C\cdot|z|^{\epsilon}\cdot|g|_h$
 if $a\neq b$.
\end{lem}
\pf
It immediately follows from Lemma \ref{lem;9.7.26}.
\hfill\qed

\vspace{.1in}

We have the following asymptotic orthogonality.
\begin{prop} \label{prop;10.11.30}
There exist good constants
$C_3>0$, $\epsilon_3>0$ and $R_{10}>0$.
Let $a_1$ and $a_2$ be elements of $S_0$
 such that $a_1\neq a_2$.
 Then $E_{a_1}$ and $E_{a_2}$ are $|z|^{\epsilon_3}$-asymptotically
 orthogonal.
 More precisely,
 let $v_i$ be $C^{\infty}$-sections of $E_i$.
 Then it holds
 $\big|(v_1,v_2)_h\big|\leq
 C_3\cdot |z|^{\epsilon_3}\cdot|v_1|_h\cdot|v_2|_h$
 on $\Delta^{\ast}(R_{10})$.
\end{prop}
\pf
Let $v$ be a $C^{\infty}$-section of $E_a$.
We have the following decomposition:
\[
 v=\sum_{b\leq_1 a}v_b,\quad\quad
 v_b\in C^{\infty}\bigl(X-D,E_b'\bigr).
\]
We have the equalities
$\rho(v)=a\cdot v=\sum_{b\leq a}a\cdot v_{b}$.
On the other hand,
we have the following equalities:
\[
 \rho(v)=\rho'(v)+q(v)
=\sum_{b\leq_1 a}b\cdot v_b+q(v).
\]
Hence we obtain the following:
\[
 \sum_{b<_1a}(a-b)\cdot v_b=q(v).
\]
Therefore there exists a good constant $A_{22}$
such that the following holds for any $v$:
\[
 |v_b|_h\leq A_{22}\cdot|z|^u\cdot |v|_h.
\]
Let $v$ be a $C^{\infty}$-section of $E_a$
and $w$ be a $C^{\infty}$-section of $E_c$
such that $c<a$.
Then we have the following:
\[
 |(v,w)_h|
=\Bigl|
  \sum_{b\leq_1 a}(v_b,w)_h
 \Bigr|
=\Bigl|
  \sum_{b\leq_1 c}(v_b,w)_h
 \Bigr|
\leq
 \sum_{b\leq_1 c}|v_b|_h\cdot|w|_h
\leq
 A_{22}\cdot|z|^u\cdot \sum _{b\leq_1 c}|v|_h\!\cdot\!|w|_h.
\]
Hence there exists a good constant $A_{23}>0$
such that the following holds:
\[
 |(v,w)_h|
\leq
 A_{23}\cdot|z|^u\cdot|v|\cdot |w|.
\]
Thus we are done.
\hfill\qed

%% file: 2.tex

\subsubsection{Prolongment of $\nbigelambda$}
We put $X=\Delta$ and $D=\{O\}$.
Let $\harmonicbundle$ be a tame harmonic bundle
over $X-D$.
We have the deformed holomorphic bundle
$(\nbigelambda,\DDlambda)$ over $\nbigxlambda-\nbigdlambda$
with the metric $h$.
Let us recall the result on the prolongment of
$\nbigelambda$.
(See the section 10 of \cite{s1} and the section 3 of \cite{s2}.
See also the subsubsection 4.3.1--4.3.3 in \cite{mochi}.)

\begin{enumerate}
\item
 For any real number $b\in\real$,
 the $\nbigo_X$-module
 $\prolongg{b}{\nbigelambda}$ is locally free.
\item
 For any real numbers $a<b$,
 we have the canonical inclusion
 $\prolongg{a}{\nbigelambda}\lrarr\prolongg{b}{\nbigelambda}$
 of $\nbigo_X$-modules.
 Then we obtain the parabolic filtration of
 $\prolongg{b}{\nbigelambda}_{|O}$.
 Namely 
 we put 
  $F_a(\prolongg{b}{\nbigelambda}):=
 \Image(
  \prolongg{a}{\nbigelambda}_{|O}\lrarr
  \prolongg{b}{\nbigelambda}_{|O}
 )$. 
Then we have the following inclusions
for any $b-1\leq a\leq b$:
\[
0=F_{b-1}(\prolongg{b}{\nbigelambda})
\subset
 F_a(\prolongg{b}{\nbigelambda})
 \subset
 F_b(\prolongg{b}{\nbigelambda})
 =\prolongg{b}{\nbigelambda}_{|O}.
\]
\item
 We put as follows:
 \[
  F_{<a}(\nbigelambda)=\sum_{c<a}F_c(\nbigelambda)
=\bigcup_{c<a}F_c(\nbigelambda),
\quad
 \Gr_a^F(\prolongg{b}{\nbigelambda})
 :=
 \frac{F_a(\prolongg{b}{\nbigelambda})}
 {F_{<a}(\prolongg{b}{\nbigelambda})}.
 \]
If $a\leq b<a+1$,
we have the canonical isomorphism
$\Gr_a^F(\prolongg{a}{\nbigelambda})\lrarr
 \Gr_a^F(\prolongg{b}{\nbigelambda})$.
Hence we omit to denote $b$ in this case.

On the contrary,
if $b<a$ or $b\geq  a+1$,
then we have $\Gr^F_a(\prolongg{b}{\nbigelambda})=0$
by definition.

\item
Let $\vecv=(v_i)$ be a holomorphic frame of
$\prolongg{b}{\nbigelambda}$ over $X$
compatible with parabolic filtration on $D$.
We put $b_i:=\deg^F(v_i)$.
We put $v_i':=|z|^{b_i}\cdot v_i$,
and then we obtain the $C^{\infty}$-frame $\vecv':=(v_i')$
of $\nbigelambda$ over $\nbigxlambda-\nbigdlambda$.
Then $\vecv'$ is adapted up to log order

\item
$\DDlambda$ is logarithmic in the following sense:
if $f$ is a holomorphic section of $\prolongg{b}{\nbigelambda}$,
then $\DDlambda f$ is a holomorphic section of
$\prolongg{b+1}{\nbigelambda}\otimes\Omega^{1,0}
=\prolongg{b}{\nbigelambda}\otimes\Omega^{1,0}(\log O)$.
In particular,
we obtain the residue
$\resddlambda\in
 \End(\prolongg{b}{\nbigelambda}_{|O})$,
which preserves the parabolic filtration $F$.

\end{enumerate}

We have the $\EE$-decomposition of $\prolongg{b}{\nbigelambda}_{|O}$
for $\Res(\DD)$:
\[
 \prolongg{b}{\nbigelambda}_{|O}
=\bigoplus_{\alpha\in\cnum}
 \EE\bigl(\prolongg{b}{\nbigelambda}_{|O},\alpha\bigr).
\]

\begin{lem}
The decomposition $\EE$ and the parabolic filtration $F$
is compatible.
\end{lem}
\pf
Since $F_a(\prolongg{b}{\nbigelambda})$ is stable
under the action of $\resddlambda$,
we have only to apply Lemma \ref{lem;9.5.10}.
\hfill\qed

\vspace{.1in}
For any $u=(a,\alpha)\in\real\times\cnum$
and for any $b$ such that $a\leq b<a+1$,
we put as follows:
\[
 \Gr^{F,\EE}_u(\nbige^{\lambda}):=
 \EE(\Gr^F_{a}(\prolongg{b}{\nbigelambda}),\alpha)
=\Gr^F_a\EE\bigl(\prolongg{b}{\nbigelambda}_{|O},\alpha\bigr).
\]
It is independent of a choice of $b$.

 For any $u\in\real\times\cnum$,
 we have the induced morphism
 $\Gr_{u}^{F,\EE}(\resddlambda)$
 on $Gr^{F,\EE}_u(\nbigelambda)$.
The nilpotent part is denoted by $\nbign_u$.
Then we obtain the weight filtration $W$ of $\nbign_u$.

\subsubsection{KMS-spectrum}

\begin{df}
For a harmonic bundle $\harmonicbundle$,
the set $\KMSE{\lambda}{}$ is defined as follows:
\[
 \KMSE{\lambda}{}:=
 \bigl\{ u\in \real\times\cnum\,\big|\,
 \dim \Gr^{F,\EE}_u(\nbige)\neq 0
 \bigr\}.
\]
It is called the KMS-spectrum set of $\nbige$ at $(\lambda,O)$.
For any $u\in \KMSE{\lambda}$,
the number $\multiplicity(\lambda,u)$ is defined as follows:
\[
 \multiplicity(\lambda,u):=\dim \Gr^{F,\EE}_u(\nbige^{\lambda}).
\]
It is called the multiplicity of the KMS-spectrum $u$.
\hfill\qed
\end{df}

The natural morphisms $\KMSE{\lambda}\lrarr \real$
and $\KMSE{\lambda}\lrarr \cnum$
are denoted by $\pi^{\paramap}$ and $\pi^{\eigenmap}$.
We put $\ParE{\lambda}:=\Image(\pi^{\paramap})$
and $\SpE{\lambda}:=\Image(\pi^{\eigenmap})$.

\begin{prop} \label{prop;10.10.30}
We have the isomorphism
$\Gr^{F,\EE}_u(\nbigelambda)\simeq
  \Gr^{F,\EE}_{u+(1,-\lambda)}(\nbigelambda)$.
\end{prop}
\pf
Let $\vecv=(v_i)$ be a holomorphic frame
of $\prolongg{a}{\nbigelambda}$
compatible with $\EE$ and $F$.
We put $b_i:=\deg^F(v_i)$.
We put $\tilde{v}_i:=z^{-1}\cdot v_i$
and $\widetilde{\vecv}:=(\tilde{v}_i)$.

\begin{lem}
$\widetilde{\vecv}$ gives a holomorphic frame
of $\prolongg{a+1}{\nbigelambda}$.
\end{lem}
\pf
We put $\tilde{v}_i':=\tilde{v}_i\cdot |z|^{b_i-1}$,
and $\tilde{\vecv}'=(\tilde{v}_i)$.
Then it is $C^{\infty}$-frame of $\nbigelambda$
over $\nbigxlambda-\nbigdlambda$,
and it is adapted up to log order.
Then $\widetilde{\vecv}$ gives a frame of $\prolongg{a}{\nbigelambda}$
compatible with the parabolic filtration
due to Lemma \ref{lem;9.5.10} and Lemma \ref{lem;10.10.31}.
\hfill\qed

\vspace{.1in}

Let us return to the proof of Proposition \ref{prop;10.10.30}.
Let $\nbiga$ be the $\lambda$-connection form of $\DD^{\lambda}$
with respect to the frame $\vecv$,
i.e. $\DDlambda\vecv=\vecv\cdot\nbiga$ holds.
Then we have the following:
\[
 \DD\widetilde{\vecv}
=\DD\bigl(z^{-1}\cdot\vecv\bigr)
=\widetilde{\vecv}\cdot\Bigl(
 \nbiga-\lambda\frac{dz}{z}
 \Bigr).
\]
We obtain the isomorphism
$\prolongg{a}{\nbigelambda}_{|O}\lrarr
 \prolongg{a+1}{\nbigelambda}_{|O}$
defined by the correspondence
$v_i(0)\longmapsto \tilde{v}_i(0)$.
Then it induces the isomorphism:
\[
 \Gr^{\EE,F}_u(\nbigelambda)
\lrarr
 \Gr^{\EE,F}_{u+(1,-\lambda)}(\nbigelambda).
\]
Thus we are done.
\hfill\qed

\begin{cor} \label{cor;9.5.1}
We have the equality
$ \multiplicity(\lambda,u)
=\multiplicity\big(\lambda,
        u+(1,-\lambda)\big)$.
\hfill\qed
\end{cor}

\vspace{.1in}

We have the free $\seisuu$-action on
$\KMSE{\lambda}$:
\[
 \seisuu\times\KMS(\nbigelambda)\lrarr
 \KMS(\nbigelambda),\quad
 (n,u)\longmapsto
 u+n\cdot(1,-\lambda).
\]
It preserves the multiplicities.

\begin{df}
We put $\KMSEoverline{\lambda}:= \KMSE{\lambda}/\seisuu$.
Note that
the multiplicity of any element $u\in\KMSEoverline{\lambda}$
is naturally defined,
due to Corollary {\rm \ref{cor;9.5.1}}.
\hfill\qed
\end{df}

\begin{df}
We put as follows:
\[
 \KMSEprolongg{b}{\lambda}
:=
 \bigl\{u\in \KMSE{\lambda}\,\big|\,b-1<\pi^{\paramap}(u)\leq b\bigr\}
=\bigl\{u\in\real\times\cnum\,\big|\,
  \Gr^{\EE,F}_u(\prolongg{b}{\nbigelambda})\neq 0\bigr\}.
\]
\hfill\qed
\end{df}

We have the natural morphism
$\pi:\KMSEprolongg{b}{\lambda}\lrarr \KMSEoverline{\lambda}$.
The following lemma is clear.
\begin{lem}
The morphism $\pi$ is bijective.
\hfill\qed
\end{lem}

The restriction of $\pi^{\paramap}$
to $\KMSEprolongg{b}{\lambda}$
gives the morphisms
$\KMSEprolongg{b}{\lambda}\lrarr \real$.
The image $\pi^{\paramap}\bigl(\KMSEprolongg{b}{\lambda}\bigr)$
is denoted by $\ParEprolongg{b}{\lambda}$.
The following lemma is clear.
\begin{lem}
We have
$ \ParEprolongg{b}{\lambda}=
 \bigl\{a\in\real\,\big|\,
 \Gr_a^F(\prolongg{b}{\nbigelambda})\neq 0\bigr\}$.
\hfill\qed
\end{lem}

The restriction of the morphism $\pi^{\eigenmap}$
to $\KMSEprolongg{b}{\lambda}$ gives the morphism
$\pi^{\eigenmap}:\KMSEprolongg{b}{\lambda}\lrarr \cnum$.
The image
$\pi^{\eigenmap}\bigl(\KMSEprolongg{b}{\lambda}\bigr)$ is denoted by
$\SpEprolongg{b}{\lambda}$.
The following lemma is clear.
\begin{lem}
We have 
$ \SpEprolongg{b}{\lambda}=
 \bigl\{\alpha\in\cnum\,\big|\,
 \EE\bigl(\prolongg{b}{\nbigelambda}_{|O},\alpha\bigr)\neq 0\bigr\}$.
\hfill\qed
\end{lem}

\subsubsection{The functoriality of the KMS structure for pull back}

Let $c$ be a positive integer,
and $\psi_c:X\lrarr X$ be the morphism
given by $z\longmapsto z^c$.

Let $f$ be a holomorphic section of $\prolong{\nbigelambda}$
over $X$.
We put $b=-\ord(f)$. Assume $-1<b\leq 0$,
i.e., $f_{|O}\neq 0$ in $\prolong{\nbigelambda}_{|O}$.

We have the holomorphic section $\psi_c^{-1}(f)$
of $\psi_c^{-1}(\nbigelambda)$ over $\nbigxlambda-\nbigdlambda$.
Then we have $\ord(\psi_c^{-1}(f))=-c\cdot b$.
We put as follows:
\[
 \tilde{f}:=
 z^{\nu(c\cdot b)}\cdot \psi_c^{-1}(f).
\]
Then we have the following:
\[
 -\ord(\tilde{f})
=-\nu(c\cdot b)-\ord(\psi_c^{-1}(f))
=-\nu(c\cdot b)+c\cdot b=\kappa(c\cdot b).
\]
Hence $\tilde{f}$ is a section of 
$\prolongg{\kappa(cb)}{\psi_c^{-1}\nbigelambda}$.
In particular, it gives a section of
$\prolong{\psi_c^{-1}\nbigelambda}$.

Let $\vecv$ be the holomorphic frame of
$\prolong{\nbigelambda}$ compatible with $F$.
We put $b_i:=\deg^F(v_i)$.
We put $v_i':=|z|^{b_i}\cdot v_i$
and $\vecv':=(v_i')$.
Recall that $\vecv'$ is a $C^{\infty}$-frame 
of $\prolong{\nbigelambda}$ over $\nbigxlambda-\nbigdlambda$,
which is adapted up to log order.
For each $v_i$, take the holomorphic section
$\tilde{v}_i$ of $\prolong{\psi_c^{-1}\nbigelambda}$
as above.
Then we obtain the tuple of sections $\widetilde{\vecv}=(\tilde{v}_i)$
of $\prolong{\psi_c^{-1}\nbigelambda}$.

\begin{lem}
$\widetilde{\vecv}$ is a holomorphic frame
of $\prolong{\psi_c^{-1}\nbigelambda}$.
\end{lem}
\pf
We put as follows:
\[
  \tilde{v}_i':=|z|^{\kappa(cb_i)}\tilde{v}_i=
 C_i(z)\cdot \psi_c^{-1}(v_i)'.
\]
Here we have $|C_i(z)|=1$.
We obtain the tuple of
the $C^{\infty}$-sections $\tilde{\vecv}'=(\tilde{v}_i')$
of $\psi_c^{-1}\nbigelambda$.
It is a $C^{\infty}$-frame of $\psi_c^{-1}\nbigelambda$
over $\nbigxlambda-\nbigdlambda$,
and it is adapted up to log order.
Hence we have only to apply Lemma \ref{lem;9.5.10}.
\hfill\qed

\begin{cor}
We have the surjective morphism
$\psi_c^{\ast}:
\KMS(\prolong{\nbigelambda})
\lrarr
\KMS(\prolong{\psi_c^{-1}\nbigelambda})$
given as follows:
\[
 (b,\beta)
\longmapsto
\bigl(\kappa(c\cdot b),\,\,c\cdot\beta+\nu(c\cdot b)\cdot\lambda\bigr)
=c\cdot (b,\beta)+\nu(c\cdot b)\cdot(-1,\lambda).
\]
We have isomorphisms:
\[
 \Gr^{\EE,F}_u\bigl(\prolong{\psi_c^{-1}\nbigelambda}\bigr)
\simeq
 \bigoplus_{\psi_c^{\ast}(u)=u'}
 \Gr^{\EE,F}_{u'}\bigl(\prolong{\nbigelambda}\bigr),
\quad\quad
 \Gr_a^F(\prolong{\psi_c^{-1}\nbigelambda})
\simeq
 \bigoplus_{\kappa(cb)=a}
 \Gr_b(\prolong{\nbigelambda}).
\]
The isomorphism are given by $\vecv$ and $\widetilde{\vecv}$.
\end{cor}
\pf
We have only to note that $\tilde{\vecv}$ is compatible with
$\EE$ and $F$.
The compatibility of $\tilde{\vecv}$ and $F$
follows from the fact that $\tilde{\vecv}'$ is adapted 
up to log order
(Lemma \ref{lem;10.10.31}).

Let $\nbiga$ be the $\lambda$-connection form of $\DD$
with respect to the frame $\vecv$,
namely we have 
$\DD\vecv=\vecv\cdot\nbiga$.
Then we obtain
$\tilde{\DD}^{\lambda}\tilde{\vecv}
=\tilde{\vecv}\cdot\tilde{\nbiga}$,
where $\tilde{A}=\psi_c^{-1}\nbiga+B\cdot dz/z$,
and $B$ denotes the diagonal matrix such that
$B_{j\,j}=\nu(c\cdot b_j)$.
Thus we obtain the compatibility of $\tilde{\vecv}$ and 
the decomposition $\EE$.
\hfill\qed

\begin{cor}
The following holds:
\[
 \Par\bigl(\psi_c^{-1}(\nbige)\bigr)=
 \bigcup_{a\in\Par(\prolong{\nbigelambda})}
 \big(c\cdot a+\seisuu\big).
\]
\hfill\qed
\end{cor}

%% file: 2.1.tex

\subsubsection{The action}
\label{subsubsection;9.10.7}

Assume that $c$ is a positive integer
which is sufficiently large with respect to
$\Par(\prolong{\nbige^{\lambda}})$.
(See Definition \ref{df;9.5.11}).
We have the action of $\mu_{c}$ on $X$,
given by $\omega^{\ast}z=\omega\cdot z$.
It can be naturally lifted to the action on $\psi_c^{-1}\nbigelambda$.
Since $v_i$ are invariant under the action of $\mu_c$,
we have the following:
\[
 \omega^{\ast}\tilde{v}_i=\omega^{\nu(cb_i)}\cdot \tilde{v}_i.
\]

We have the weight decomposition:
\[
 \prolong{\psi_c^{-1}\nbige}_{|O}
=\bigoplus_{h}U_h.
\]
Here $\omega^{\ast}=\omega^h$ on $U_h$
for $-c+1\leq h\leq 0$.
The following lemma is clear.
\begin{lem}
$U_h=\bigl\langle \tilde{v}_i\,\big|\, \nu(c\cdot b_i)=h\bigr\rangle$.
\hfill\qed
\end{lem}

We put $S:=\bigl\{h\,\big|\,-c+1\leq h\leq 0,\,\,U_h\neq 0\bigr\}$.
Then we have
$S=\bigl\{\nu(c\cdot b)\,\big|\, b\in \Par(\prolong{\nbigelambda})\bigr\}$.
Since $c$ is sufficiently large with respect
to $\Par(\prolong{\nbigelambda})$,
any element
$b\in \Par(\prolong{\nbigelambda})$ is uniquely determined
by the number $\nu(c\cdot b)\in S$.
Thus we have the map
$\varphi:S\lrarr \Par(\prolong{\psi_c^{-1}\nbigelambda})$
given by the following correspondence:
\[
 \nu(c\cdot b)\longmapsto \kappa(c\cdot b).
\]

Let us consider the filtration $F'$ given as follows:
\[
 F'_{\tilde{b}}
:=\bigoplus_{\substack{h\in S, \\
     \varphi(h)\leq \tilde{b}}}U_h.
\]
\begin{lem} \label{lem;9.10.11}
We have
$F'_{b}(\prolong{\psi_c^{-1}\nbigelambda})
=F_b(\prolong{\psi_c^{-1}\nbigelambda})$.
\end{lem}
\pf
It follows from the following equalities:
\[
 F_{{b}}(\prolong{\psi_c^{-1}}\nbigelambda)
=\big\langle 
 \tilde{v}_i\,\big|\,\kappa(c\cdot b_i)\leq b
 \big\rangle
=F_b'(\prolong{\psi_c^{-1}\nbigelambda}).
\]
Thus we are done.
\hfill\qed

\begin{cor}
The decomposition
$\prolong{\psi_c^{-1}\nbigelambda}_{|O}=
 \bigoplus_h U_h$ gives a splitting of the parabolic filtration
in the following sense:
\[
 F_b=\bigoplus_{\substack{h\in S\\ \varphi(h)\leq b}}U_h.
\]
In particular,
$U_h$ is naturally
isomorphic to $\Gr_b(\prolong{\psi_c^{-1}\nbigelambda}_{|O})$
$(\varphi(h)=b)$.
\hfill\qed
\end{cor}

\subsubsection{Descent of the frame}

On the other hand,
we can descend the equivariant frame.
Let $f$ be a holomorphic section
of $\prolong{\psi_c^{-1}\nbigelambda}$,
such that $\omega^{\ast}(f)=\omega^h\cdot f$
for some integer $h$ such that $-c< h\leq 0$.
We put $f_1:=z^{-h}\cdot f$,
and then we have $\omega^{\ast}(f_1)=f_1$.
Hence $f_1$ induces the holomorphic section
$\tilde{f}$ of $\nbigelambda$ over $X-D$.
We have the following:
\[
 -\ord(\tilde{f})
 =-c^{-1}\cdot\ord(f_1)
 =c^{-1}\cdot\bigl(h-\ord(f)\bigr)\leq 0.
\]
Hence $\tilde{f}$ gives a holomorphic section
of $\prolong{\nbigelambda}$.

Let $\vecv=(v_i)$ be a holomorphic frame
of $\prolong{\psi_c^{-1}\nbigelambda}$
satisfying the following:
\begin{itemize}
\item
It is equivariant in the sense
$\omega^{\ast}\tilde{v}_i=\omega^{h_i}\cdot\tilde{v}_i$
for $-c<h_i\leq 0$.
\item
It is compatible with the parabolic filtration $F$.
\end{itemize}

We put $b_i:=\deg^F(v_i)$,
and then we have $-1<b_i\leq 0$.

\begin{lem}
 We have $-1< c^{-1}\cdot(h+b_i)\leq 0$.
\end{lem}
\pf
Since we have $-c+1\leq h\leq 0$ and $-1< b_i\leq 0$,
we obtain $-c\leq h+b_i<0$.
\hfill\qed

\vspace{.1in}

Let us take the section $\tilde{v}_i$ of $\prolong{\nbigelambda}$
for each $v_i$ as above.
Then we obtain the tuple of sections
$\tilde{\vecv}=(\tilde{v}_i)$.

\begin{lem} \label{lem;9.10.15}
$\tilde{\vecv}$ is a holomorphic frame of
$\prolong{\nbigelambda}$,
compatible with the parabolic filtration.
\end{lem}
\pf
We put $\tilde{b}_i:=c^{-1}\cdot(b_i+h-c)$.
We put as follows:
\[
 \tilde{v}_i':=
 |z|^{\tilde{b}_i}\cdot\tilde{v}_i,
\quad
 \tilde{\vecv}=(\tilde{v}_i).
\]
Then it can be checked that $\tilde{\vecv}$ 
is adapted up to log order.
Thus we have only to apply Lemma \ref{lem;9.5.10}
and Lemma \ref{lem;10.10.31}.
\hfill\qed

\subsubsection{Functoriality for tensor product}

Let $(E^{(a)},\delbar_{E^{(a)}},\theta^{(a)},h^{(a)})$
$(a=1,2)$
be tame harmonic bundles over $X-D$.
We obtain the prolonged deformed holomorphic bundles
$\prolong{\nbige}^{(a)\,\lambda}$.
Let $\vecv^{(a)}$ be holomorphic frames of
$\prolong{\nbige}^{(a)\,\lambda}$
compatible with the parabolic filtration $F$.
We put
$b^{(a)}_i:=\deg^F(v^{(a)}_i)$,
and $v^{(a)\,\prime}_i:=|z|^{b^{(a)}_i}v_i^{(a)}$.
The tuple of sections
$\vecv^{(a)\,\prime}=(v^{(a)\,\prime}_i)$
gives a $C^{\infty}$-frame of $\nbige^{\lambda\,(a)}$
over $X-D$,
which is adapted up to log order.

Then we obtain the $C^{\infty}$-frame
$\vecv^{(1)\,\prime}\otimes\vecv^{(2)\,\prime}$
of $\nbige^{(1)\,\lambda}\otimes\nbige^{(2)\,\lambda}$,
given as follows:
\[
 \vecv^{(1)\,\prime}\otimes\vecv^{(2)\,\prime}
=\Bigl(
 v^{(1)\,\prime}_i
 \otimes
 v^{(2)\,\prime}_j\,\Big|\,
 1\leq i\leq \rank E^{(1)},\,\,
 1\leq j\leq\rank E^{(2)}
 \Bigr).
\]
It is adapted up to log order.
Hence we put
$w_{i\,j}:=
 z^{-\epsilon(i,j)}\cdot
 v_i^{(1)}\otimes v_j^{(2)}
 $,
where $\epsilon(i,j)$ are given as follows:
\[
 \epsilon(i,j):=
 \left\{
 \begin{array}{ll}
 1 & (b_i^{(1)}+b_j^{(2)}\leq -1)\\
 0 & (\mbox{otherwise, i.e., } -1< b_i^{(1)}+b_j^{(2)}\leq 0).
 \end{array}
 \right.
\]
Then we obtain the tuple of holomorphic sections
$\vecw=
 \bigl(w_{i\,j}\,
    \big|\,1\leq i\leq \rank E^{(1)},\,\,
           1\leq j\leq \rank E^{(2)}
 \bigr)$,
and it gives the holomorphic frame of
$\prolong{\bigl(
 \nbige^{(1)\,\lambda}\otimes\nbige^{(2)\,\lambda}\bigr)}$,
compatible with the filtration.

\begin{cor}
We have the surjective morphism
$\psi:\KMS(\prolong{\nbige^{(1)\,\lambda}})
\times
 \KMS(\prolong{\nbige^{(2)\,\lambda}})
\lrarr
 \KMS\bigl(\prolong{\bigl(
 \nbige^{(1)\,\lambda}\otimes\nbige^{(2)\,\lambda}\bigr)}\bigr)$.
For elements $u_i=(b_i,\beta_i)$ for $i=1,2$,
the element $\psi(u_1,u_2)$ is given as follows:
\[
 \psi(u_1,u_2)
=
 \bigl(\kappa(b_1+b_2),\,\,\beta_1+\beta_2-\nu(b_1+b_2)\cdot\lambda
 \bigr).
\]
We have the equality of the multiplicities:
\[
 \multiplicity(\lambda,u)
=\sum_{\psi(u_1,u_2)=u}
 \multiplicity(\lambda,u_1)\cdot
 \multiplicity(\lambda,u_2)
\]
\hfill\qed
\end{cor}

\begin{cor}
We have the isomorphisms:
\[
 \Gr^F_b(
 \prolong{\bigl(
 \nbige^{(1)\,\lambda}\otimes
 \nbige^{(2)\,\lambda}\bigr)}
 )
\simeq
 \bigoplus_{\kappa(b_1+b_2)=b}
 \Gr^{F}_{b_1}(\prolong{\nbige^{(1)\,\lambda}})
 \otimes
 \Gr^{F}_{b_2}(\prolong{\nbige^{(2)\,\lambda}}).
\]

\[
  \Gr^{F,\EE}_u(
 \prolong{\bigl(
 \nbige^{(1)\,\lambda}\otimes
 \nbige^{(2)\,\lambda}\bigr)}
 )
\simeq
 \bigoplus_{\psi(u_1,u_2)=u}
 \Gr^{F,\EE}_{u_1}(\prolong{\nbige^{(1)\,\lambda}})
 \otimes
 \Gr^{F,\EE}_{u_2}(\prolong{\nbige^{(2)\,\lambda}}).
\]
\hfill\qed
\end{cor}

\begin{cor}
We have the isomorphism:
\[
 \Gr_b^F(\prolong{\Sym^h\nbigelambda})
\simeq
 \bigoplus_{(\vecb,\vecm)\in\nbigs(b,h)}
 \bigotimes_i
 \Sym^{m_i}(\Gr^F_{b_i}\prolong{\nbigelambda}).
\]
Here we put as follows:
\[
 \nbigs(b,h):=
 \Big\{(\vecb,\vecm)\,\Big|\,
 \sum m_i=h,\quad
 \kappa\Bigl(\sum m_ib_i\Bigr)=b
 \Big\}
\]
In all,
we have an isomorphism
$\Gr^F(\prolong{\Sym^{\cdot}\nbigelambda})
\simeq
 \Sym^{\cdot}(\Gr^F\prolong{\nbigelambda})
 $.

We also have an isomorphism:
\[
 \Gr_b^F\Bigl(\prolong{\bigwedge^h\nbigelambda}\Bigr)
\simeq
 \bigoplus_{(\vecb,\vecm)\in\nbigs(b,h)}
 \bigotimes_i
 \bigwedge^{m_i}(\Gr^F_{b_i}\prolong{\nbigelambda}).
\]
We have
$\Gr^F(\prolong{\bigwedge^{\cdot}\nbigelambda})
\simeq
 \bigwedge^{\cdot}(\Gr^F\prolong{\nbigelambda})
 $.
\hfill\qed
\end{cor}

\subsubsection{Functoriality for dual}

Let $\vecv=(v_i)$ be a holomorphic frame of
$\prolong{\nbigelambda}$ compatible with $F$.
We put $b_i:=\deg^F(v_i)$
$v_i':=|z|^{b_i}\cdot v_i$,
and $\vecv'=(v_i')$.
Then we obtain the dual frame
$\vecv^{\prime\,\lor}$
of $\nbige^{\lambda\,\lor}$ over $X-D$.
It is adapted up to log order.

Let $\vecv^{\lor}=(v_i^{\lor})$ be the dual frame of $\vecv$
over $X-D$.
Then it gives a holomorphic frame
of $\prolongg{1-\epsilon}{\nbige^{\lambda\,\lor}}$
for some $\epsilon>0$.

Then we put $w_i:=z^{\epsilon(i)}\cdot v_i^{\lor}$
and $\vecw=(w_i)$,
where $\epsilon(i)$ is given as follows:
\[
 \epsilon(i):=
 \left\{
 \begin{array}{ll}
 1 & (b_i\neq 0)\\
 0 & (b_i=0).
 \end{array}
 \right.
\]
Then $\vecw$ is a holomorphic frame
of $\prolong{\nbige^{\lor\,\lambda}}$
compatible with $F$.

\begin{cor}
We have the bijection
$\psi:\KMS(\prolong{\nbigelambda})
\lrarr
 \KMS(\prolong{\nbige^{\lambda\,\lor}})$.
For any $u=(b,\beta)$,
$\psi(u)$ is given by $\bigl(\kappa(-b),-\beta-\nu(-\beta)\bigr)$.
We also have the isomorphism
$ \Gr^{F\,\EE}_u(\prolong{\nbigelambda})
\simeq
 \Gr^{F\,\EE}_{\psi(u)}(\prolong{\nbige^{\lambda\,\lor}})$.
\hfill\qed
\end{cor}

%% file: 3.tex

\subsubsection{The statement}

We put $X=\Delta$ and $D=\{O\}$.
Let $\harmonicbundle$ be a tame harmonic bundle
over $X-D$.
As is already seen,
we obtain the vector space
$V_{u}:=\Gr^{\EE,F}_{u}(\prolong{E}_{|O})$
and $N_u$.
We have the model bundle $E(V,N)$ as in Lemma \ref{lem;b12.10.1}:
\[
 (E_0,\delbar_{E_0},h_0,\theta_0)=
 \bigoplus_{u\in\KMS(\prolong{E})}
 E(V_u,N_u)\otimes L(u).
\]
We can pick an isomorphism
$\Phi:\prolong{E}_0\lrarr\prolong{E}$
satisfying the following:
\begin{itemize}
\item
 $\Phi$ preserves the parabolic filtrations.
\item
 The induced morphism
 $\Gr^F_a(\Phi)\in
 Hom\bigl(\Gr^F_a(\prolong{E_0}),\Gr^F_a(\prolong{E})\bigr)$
 is compatible with
 the morphisms
 $\Gr_a^F(\Res(\theta_0))$ and $\Gr^F_a(\Res(\theta))$.
\end{itemize}

\begin{prop}[Simpson]
$\Phi$ and $\Phi^{-1}$ are bounded.
\end{prop}
\pf
See the subsubsection 4.3.3 in the previous paper \cite{mochi},
for example.
\hfill\qed

\vspace{.1in}

Since we have $\nbige^{\lambda}=E$ and $\nbige^{\lambda}_0=E_0$
as $C^{\infty}$-bundles,
we obtain the $C^{\infty}$-isomorphism
$\Phi:\nbige^{\lambda}_0\lrarr \nbigelambda$
on $\nbigxlambda-\nbigdlambda$.
Let us take holomorphic frames
$\vecv_0$ and $\vecv$ of
$\prolong{\nbigelambda_0}$
and $\prolong{\nbigelambda}$,
which are compatible with generalized eigen decompositions $\EE$,
parabolic filtrations $F$ and
the weight filtrations $W$.
We put as follows:
\[
\begin{array}{lll}
 \deg^F(v_j)=b_j,
 &
 \deg^{\EE}(v_j)=\beta_j,
 &
 {\displaystyle\frac{\deg^{W}(v_j)}{2}}=k_j, \\
\mbox{{}}\\
 \deg^{F}(v_{0\,i})=b_{0\,i}
 &
 \deg^{\EE}(v_{0\,i})=\beta_{0\,i},
 &
 {\displaystyle\frac{\deg^{W}(v_{0\,i})}{2} }=k_{0\,i}.
\end{array}
\]
We also put as follows:
\[
\begin{array}{ll}
 v_j':=v_j\cdot |z|^{b_j}\cdot(-\log|z|)^{-k_j},
 &
 \vecv'=(v_j'),\\
\mbox{{}}\\
 v'_{0\,i}:=v_{0\,i}\cdot |z|^{b_{0\,i}}\cdot(-\log|z|)^{-k_{0\,i}},
 &
 \vecv'_0=(v_{0\,i}').
\end{array}
\]
Then we obtain the $C^{\infty}$-functions
$I$ and $I'$ of
$X-D$ to $M(r)$ defined as follows:
\[
 \Phi(v_{0\,i})=\sum I_{j,i}\cdot v_j,
\quad
 \Phi(v_{0\,i}')=\sum I_{j,i}'\cdot v_j'.
\]
The following lemma can be checked by a direct calculation.
\begin{lem}
We have the equality
$ I_{j,i}'=I_{j,i}\cdot|z|^{b_{0,i}-b_j}
 \cdot\bigl(-\log|z|\bigr)^{-k_{0,i}+k_j}$.
\hfill\qed
\end{lem}

By the isomorphism $\Phi$,
we identify $E$ and $E_0$.
Let $\theta^{\dagger}$ be the conjugate of $\theta$
with respect to $h$,
and $\theta^{\dagger}_0$
be the conjugate of $\theta_0$
with respect to $h_0$.
Recall the following lemma.
\begin{lem}[Simpson, Lemma 7.3, Lemma 7.7 \cite{s2}]
\label{lem;9.15.2}
We have the following inequalities:
\[
 |\theta^{\dagger}-\theta^{\dagger}_0|_h\leq 
C\cdot |z|^{-1}(-\log|z|)^{-1},
\quad\quad
 \int
 |\theta^{\dagger}-\theta_0^{\dagger}|_h
 \cdot |z|\cdot\bigl(-\log |z|\bigr)
\frac{|dz\cdot d\bar{z}|}{|z|^2\cdot\bigl(-\log |z|\bigr)}
 <\infty. 
\]
Here $C$ denotes some positive constant.
\hfill\qed
\end{lem}

We also recall an outline of the proof of the following lemma.
\begin{prop}[Simpson]\mbox{{}} \label{prop;9.5.15}
\begin{enumerate}
\item \label{8.7.1}
 $I'$ and $I^{\prime\,-1}$ are bounded.
\item \label{8.7.6}
We have the inequality
 $|I_{j\,i}'|\leq C\cdot (-\log|z|)^{-1}$
in the case $(b_j,\beta_j)\neq (b_{0\,i},\beta_{0\,i})$.
\item \label{8.7.7}
We have the inequality
 $||I_{j\,i}'||_W<\infty$
in the case $k_j\neq k_{0\,i}$.
(See the page {\rm 764} of {\rm\cite{s2}} or the subsubsection {\rm 4.3.4}
 of {\rm \cite{mochi}} for the norm $||\cdot||_W$.)
\end{enumerate}
\end{prop}

In the case $\lambda=1$, the proposition is given by Simpson
(the section 7 in \cite{s2}).
His argument clearly works in general case.
Hence we only indicate an outline.
See loc.cit. for more detail.

\subsubsection{Outline of the proof of Proposition \ref{prop;9.5.15}}

The claim \ref{8.7.1} immediately follows from the
boundedness of $\Phi$ and $\Phi^{-1}$,
and the adaptedness of $\vecv'$ and $\vecv'_0$.
Note that
$\delbar\Phi=\lambda\cdot(\theta_0^{\dagger}-\theta^{\dagger})$.
We can apply the argument of Simpson,
and we obtain the claim \ref{8.7.7}.

Let us see the outline of the proof of the claim \ref{8.7.6}.
The holomorphic sections $\psi,\psi^{(1)}$
of $End(\nbige^{\lambda})\otimes\Omega^{1,0}$
are given as follows:
\begin{equation} \label{eq;9.6.30}
 \psi(v_j):=v_j\cdot (\beta_j+\lambda \cdot b_j)\cdot
 \frac{dz}{z},
\quad\quad
 \psi^{(1)}(v_j):=v_j\cdot(-\lambda\cdot b_j)\frac{dz}{z}.
\end{equation}
For any holomorphic section $f$ of $\nbigelambda$
such that $|f|_h\leq C_1\cdot|z|^{-b}\cdot(-\log|z|)^k$,
it is easy to see that
there exists a positive constant $C_2$
such that the following holds:
\[
 \bigl|\DDlambda f-\psi(f)-\psi_1^{(1)}(f)\bigr|_h
\leq
 C_2\cdot|z|^{-b-1}(-\log|z|)^{k-1}.
\]
\begin{lem}[Lemma 7.2 of \cite{s2}]\label{lem;10.10.35}
We have the following finiteness:
\[
 \int \bigl|\psi-(1+|\lambda|^2)\theta\bigr|_h^2
 \cdot(-\log|z|)^{1-\epsilon}\cdot
 |dzd\bar{z}|<\infty.
\]
\end{lem}
\pf
We have only to show the following claim:
If $|f|\leq C_3\cdot|z|^{-b}\cdot(-\log|z|)^{k}$,
then we have the finiteness:
\begin{equation}\label{eq;8.7.2}
\int \big|\psi(f)-(1+|\lambda|^2)\theta(f)\big|_h\cdot
 |z|^{2b}\cdot
 \big(-\log|z|\big)^{1-\epsilon-2k}\cdot
 |dzd\bar{z}|<\infty.
\end{equation}

We have only to show the inequality (\ref{eq;8.7.2})
in the case $b=0$.
Note the following inequality:
\begin{multline}\label{eq;8.7.3}
 \big|\psi(f)-(1+|\lambda|^2)\cdot\theta(f)\big|^2
\leq
 2\cdot\big|\DDlambda(f)-(1+|\lambda|^2)\cdot\theta(f)\big|^2
+C\cdot |z|^{-2}\big(-\log|z|\big)^{2(k-1)}\\
=2\cdot|\lambda|^2\cdot
 |\del_{\nbigelambda}(f)|^2
+C\cdot |z|^{-2}\cdot\bigl(-\log|z|\bigr)^{2(k-1)}.
\end{multline}

We have the Weitzenbeck formula:
\begin{equation} \label{eq;a12.6.1}
 \Delta|f|_h^2=
-|\del_{\nbigelambda}f|^2+(1+|\lambda|^2)\cdot
\bigl(|\theta f|_h^2-|\theta^{\dagger}f|_h^2\bigr).
\end{equation}
Hence we obtain the following, for some positive constant $C$:
\begin{equation}\label{eq;8.7.4}
\Delta\bigl(  |f|_h^2\cdot(-\log|z|)^{1-\epsilon-2k}\bigr)
\leq
-\frac{1}{2}
 \bigl|\del_{\nbigelambda}f\bigr|_h^2\cdot(-\log|z|)^{1-\epsilon-2k}
+C\cdot |f|_h^2\cdot |z|^{-2}\cdot\bigl(-\log|z|\bigr)^{-1-\epsilon-2k}.
\end{equation}

By using (\ref{eq;8.7.4})
and the equivalence of the norms $|f|\sim (-\log|z|)^k$,
we obtain the following:
\[
 \Delta\Bigl(
 \bigl(
 |f|_h^2\cdot(-\log|z|)^{-2k}-C
 \bigr)\cdot(-\log|z|)^{1-\epsilon}
 \Bigr)
\leq
 -\frac{1}{2}|\del_{\nbigelambda}f|_h^2
 \cdot(-\log|z|)^{1-\epsilon-2k}.
\]
We put as follows:
\[
 F:=\Bigl(
 |f|_h^2\cdot(-\log|z|)^{-2k}-C
 \Bigr)\cdot
 (-\log|z|)^{1-\epsilon}.
\]
Then we obtain the following:
\[
 \lim_{|z|\to 0}\frac{F}{-\log|z|}=0.
\]
Hence we obtain the inequality:
\begin{equation}\label{eq;8.7.5}
 \int|\del_{\nbigelambda}f|_h^2\cdot
 (-\log|z|)^{1-\epsilon-2k}\cdot
 |dzd\bar{z}|<\infty.
\end{equation}
From (\ref{eq;8.7.3}) and (\ref{eq;8.7.5}),
we obtain Lemma \ref{lem;10.10.35}.
\hfill\qed

\vspace{.1in}

Since Simpson's proof of Lemma 7.4 and Lemma 7.5 in \cite{s2}
(for the case $\lambda=1$)
can be also easily applied to the general case $\lambda\neq 1$,
we omit to give a proof of the following lemma.
\begin{lem}[Lemma 7.4 and Lemma 7.5 in \cite{s2}]
We have the following inequalities:
\begin{equation}\label{eq;9.6.32}
 |(1+|\lambda|^2)\theta-\psi|_h
\leq
 |z|^{-1}(-\log|z|)^{-1}.
\end{equation}

\begin{equation}\label{eq;a11.13.5}
 |\psi-\psi_0|_h\leq |z|^{-1}(-\log|z|)^{-1}.
\end{equation}
\hfill\qed
\end{lem}

\begin{cor}[Lemma 7.11 in \cite{s2}] \label{cor;a11.13.7}
In the case
$\beta_j+\lambda\cdot b_j\neq
 \beta_{0\,i}+\lambda\cdot b_i$,
we have the inequality $|I_{i,j}'|\leq (-\log|z|)^{-1}$.
\end{cor}
\pf
It is not difficult to derive the claim from (\ref{eq;a11.13.5}).
See the proof of Lemma 7.11 of \cite{s2}.
\hfill\qed

\vspace{.1in}

Note the relation
$\delbar\Phi=\lambda\cdot(\theta_0^{\dagger}-\theta^{\dagger})$.
Then we obtain the following:
\[
 |\delbar I_{i\,j}|
\cdot |z|^{b_{0\,i}-b_j}
\cdot (-\log|z|)^{-k_{0,i}+k_j}
\leq
 C\cdot |z|^{-1}\cdot(-\log|z|)^{-1}.
\]

In the case $b_{0,i}-b_{j}\neq 0$,
we can pick the $C^{\infty}$-function $h(\delbar I_{i\,j})$
on $\Delta^{\ast}$ satisfying the following:
\begin{itemize}
\item
 $\delbar h \bigl(\delbar I_{i\,j}\bigr)=\delbar I_{i\,j}$.
\item
 $|h (\delbar I_{i\,j})|\cdot |z|^{b_{0,i}-b_j}
 \cdot (-\log|z|)^{-k_{0,i}+k_j}
 \leq C\cdot\bigl(-\log|z|\bigr)^{-1}$.
\end{itemize}
Then we obtain the following:
\begin{itemize}
\item
 $h\bigl(\delbar I_{i\,j}\bigr)-I_{i\,j}$ is holomorphic
 on $X-D$.
\item
 $\big|h\bigl(\delbar I_{i\,j}\bigr)-I_{i\,j}\big|\cdot
 |z|^{b_{0\,i}-b_j}\cdot(-\log|z|)^{-k_{0\,i}+k_j}$
 bounded.
\end{itemize}

Hence
we obtain the following inequality
in the case $b_{0\,i}-b_j\not\in\seisuu$:
\begin{equation}\label{eq;a11.13.6}
 \big|h\bigl(\delbar I_{i\,j}\bigr)-I_{i\,j}\big|\cdot
 |z|^{b_{0\,i}-b_j}\cdot(-\log|z|)^{-k_{0\,i}+k_j}
\leq C\cdot(-\log|z|)^{-1}.
\end{equation}
Hence we have the inequality $|I_{i\,j}|\leq C\cdot (-\log|z|)^{-1}$
in the case $b_{0\,i}\neq b_j$.
Thus we obtain the claim \ref{8.7.6}
from Corollary \ref{cor;a11.13.7} and (\ref{eq;a11.13.6}),
and the outline of the proof of Proposition \ref{prop;9.5.15}
is finished.
\hfill\qed

\subsubsection{Some consequences}

Recall the bijection
$\kmsmap(\lambda):\cnum\times\real\lrarr\cnum\times\real$
is defined in the subsubsection \ref{subsubsection;b11.11.11}.

\begin{cor} \label{cor;9.11.1}
Let $u$ be an element of $\KMSE{0}$.
Then $\kmsmap(\lambda,u)\in \KMSE{\lambda}$,
and we have the equality:
\[
 \dim \Gr^W_k(\Gr^{\EE,F}_u(\nbige^0))
=\dim \Gr^W_k(\Gr^{\EE,F}_{\kmsmap(\lambda,u)}(\nbigelambda)).
\]
In particular,
we obtain the bijective morphism
$\kmsmap(\lambda):\KMS(\nbige^0)\lrarr\KMS(\nbigelambda)$
and the following equality:
\[
\dim \bigl(\Gr^{\EE,F}_u(\nbige^0)\bigr)
=\dim \bigl(\Gr^{\EE,F}_{\kmsmap(\lambda,u)}(\nbigelambda)\bigr),
\quad
 \multiplicity(0,u)=
 \multiplicity\bigl(\lambda,\kmsmap(\lambda,u)\bigr).
\]
\end{cor}
\pf
The claims for the model bundles can be checked
by direct calculations.
Then the claims for general tame harmonic bundles
follows from Proposition \ref{prop;9.5.15}.
\hfill\qed

\begin{cor}\label{cor;a11.14.3}
Let $\vece$ be a holomorphic frame of $\prolong{E}$
compatible with $\EE$, $F$ and $W$.
Let $\vecv$ be a holomorphic frame of $\prolong{\nbigelambda}$
compatible with $\EE$, $F$ and $W$.
We put as follows:
\[
 b(e_j):=\deg^F(e_j),
\quad
 k(e_j):=\frac{\deg^W(e_j)}{2},
\quad
 b(v_i):=\deg^F(v_i),
\quad
 k(v_i):=\frac{\deg^W(v_i)}{2}.
\]
We put as follows:
\[
 e_j':=e_j\cdot |z|^{b(e_j)}\cdot(-\log|z|)^{-k(e_j)},
\quad
 v_i':=v_i\cdot |z|^{b(v_i)}\cdot(-\log|z|)^{-k(v_i)}.
\]
The $C^{\infty}$-function $B:X-D\lrarr M(r)$
is given as follows:
\[
 v_i'=\sum B_{j\,i}'\cdot e_j'.
\]
Then we have the following:
\begin{itemize}
\item
 $B'$ and $B^{\prime\,-1}$ are bounded.
\item
 $|B'_{j\,i}|\leq C\cdot(-\log|z|)^{-1}$
in the case $\deg^{F,\EE}(v_i)\neq \kmsmap(\lambda,\deg^{F,\EE}(e_j))$.
\end{itemize}
\end{cor}
\pf
For the case of model bundles,
we pick $\vece_0$, $\vece_0'$, $\vecv_0$ and $\vecv_0'$
similarly,
and then we obtain $B'_{0,j\,i}$.
In this case,
we may assume the following, due to the construction of
model bundles:
\begin{quote}
$(A)$: 
$B_{0,j\,i}'=0$
if
$\deg^{F,\EE}(v_{0\,i})\neq 
 \kmsmap(\lambda,\deg^{F,\EE}(e_{0\,j}))$.
\end{quote}
From our construction of $\Phi$,
we may assume $\Phi(\vece_0)=\vece$.
Then we obtain our claims
due to the assumption $(A)$ above and
the claim \ref{8.7.6} in Proposition \ref{prop;9.5.15}.
\hfill\qed

%% file: 4.tex


\subsubsection{Order of multi-valued flat sections}

Let us fix $\lambda\in\cnum^{\ast}$.
We have the $\lambda$-connection
$(\nbigelambda,\DDlambda)$.
We have the associated flat connection $\DD^{\lambda,f}$.
Then we obtain the space of multi-valued flat sections,
which we denote by $H(\nbigelambda)$.

\begin{lem} \label{lem;8.13.1}
For any $s\in H(\nbigelambda)$
and positive number $C_0$,
there exist positive constants $C_1$ and $b$
satisfying the following:
\[
 |s|\leq C_1\cdot |z|^{-b}
\quad
 \mbox{ on }
\bigl\{z\,\big|\,|\arg z|<C_0\bigr\}.
\]
\end{lem}
\pf
It follows from tameness of our harmonic bundles.
(See Remark in the page 732 of \cite{s2}, for example.
Or it is not difficult to show directly.)
\hfill\qed

\vspace{.1in}

We have the universal covering map
$\pi:\hyperh\lrarr \Delta^{\ast}$,
given by $\zeta=x+\sqrt{-1}y\longmapsto \exp(\sqrt{-1}\zeta)=z$.
We may regard $s\in H(\nbigelambda)$ as a flat section
of $\pi^{\ast}\nbigelambda$.
We have the following equality:
\[
 \frac{\del h(s,s)}{\del x}
 =
 2\Realpart\Bigl(
 h\bigl( \nabla^{\lambda,u}_{\del_x}s,s \bigr)
 \Bigr)
=2\Realpart \Bigl(
 h\bigl(
 (\nabla^{\lambda,u}_{\del_x}-\DD^{\lambda,f}_{\del_x})s,\,s
 \bigr)\Bigr).
\]
Here $\nabla^{\lambda,u}$ denote the unitary connection
for $(\nbigelambda,h)$,
and $\del_x$ denote the vector field $\del/\del x$.
The difference $\DD^{\lambda,f}-\nabla^{\lambda,u}$
is given by
$a\cdot \theta+b\cdot \theta^{\dagger}$ $(a,b\in\cnum)$.
We have the description
$\theta=\theta^{\zeta}\cdot d\zeta$
and $\theta^{\dagger}=\theta^{\zeta\,\dagger}d\bar{\zeta}$.
Due to Simpson's main estimate (Proposition \ref{prop;9.6.77} [I]),
we have the boundedness
$|\theta^{\zeta}|_h\leq C$ and $|\theta^{\zeta\,\dagger}|_h\leq C$.
Hence we obtain the following inequality 
for some positive constant $C$:
\begin{equation} \label{eq;9.5.20}
 \left|
 \frac{\del}{\del x}\log|s|_h^2
 \right|
=
 \left|
 \frac{\del h(s,s)}{\del x}\cdot h(s,s)^{-1}
 \right|
\leq
 2|a|\cdot\big|\theta^{\zeta}\big|_h+
 2|b|\cdot\big|\theta^{\zeta\,\dagger}\big|_h
\leq C.
\end{equation}
\begin{lem} \label{lem;9.5.25}
For any positive number $C_1$, there exists 
a positive constants $C_2$
such that the following holds:
\begin{itemize}
\item For any $x_i\in\real$ $(i=1,2)$ such that $|x_i|<C_1$,
and for any $y>1$,
the following inequality holds:
\[
 \left|
 \log\big|s(x_1,y)\big|^2_h
-\log\big|s(x_2,y)\big|^2_h
 \right|
\leq C_2
\]
\end{itemize}
\end{lem}
\pf
It immediately follows from the inequality (\ref{eq;9.5.20}).
\hfill\qed

\begin{df}
Let $s$ be an element of $H(\nbigelambda)$, 
$b$ be any real number.
Then `$-\ord(s)\leq b$' means the following:
\begin{itemize}
\item
Pick any real number $x_1\in\real$.
For any positive number $\epsilon>0$,
there exists a positive constant $C$ such that
$|s(x_1,y)|\leq C\cdot e^{(b+\epsilon)\cdot y}$.
\end{itemize}
Note that such property does not depend on a choice of $x_1$,
due to Lemma {\rm\ref{lem;9.5.25}}.
\hfill\qed
\end{df}

\begin{df}
We put as follows:
\[
 \nbigf_b(H(\nbigelambda))
:=\bigl\{s\in H(\nbigelambda)\,\big|\,
 -\ord(s)\leq b\bigr\}.
\]
Thus we obtain the filtration $\nbigf$ 
on $H(\nbigelambda)$.
\hfill\qed
\end{df}

Let $M^{\lambda}$ be the monodromy on $H(\nbigelambda)$.
\begin{lem}
The filtration $\nbigf$ is preserved by $M^{\lambda}$.
In particular, $\nbigf$ is compatible with
the generalized eigen decomposition of $M^{\lambda}$.
\end{lem}
\pf
It is clear from our definition of the filtration $\nbigf$.
\hfill\qed

\subsubsection{Compatibility with the order for holomorphic sections}

Let us consider the $\EE$-decomposition
of $H(\nbigelambda)$ with respect to the monodromy $M^{\lambda}$:
\[
 H(\nbigelambda)=
 \bigoplus_{\omega\in Sp(M^{\lambda})}
 \EE\bigl(H(\nbigelambda),\omega\bigr).
\]
Let $M_{\omega}^{\lambda}$ denote the restriction of $M^{\lambda}$
to $\EE\bigl(H(\nbigelambda),\omega\bigr)$.

\label{page;9.11.15}
Pick a real number $b\in\real$.
Then there exists the unique complex number
$\alpha=\alpha(b,\omega)$ satisfying the following:
\begin{itemize}
\item
 $\exp\bigl(
 -2\pi\sqrt{-1}\cdot\alpha
 \bigr)=\omega$.
\item
 $b\leq \Realpart(\alpha)<b+1$.
\end{itemize}
Note that the number $b-\Realpart\bigl(\alpha(b,\omega)\bigr)$
is independent of $b$.

Let $M^{\lambda,u}_{\omega}$ denote the unipotent part of
$M_{\omega}^{\lambda}$,
and we put as follows:
\[
 N_{\omega}^{\lambda}:=
 \frac{-1}{2\pi\sqrt{-1}}
 \log M^{\lambda,u}_{\omega}
=\frac{-1}{2\pi\sqrt{-1}}
 \sum_{n=1}^{\infty}
 \frac{(-1)^{n-1}}{n}(M_{\omega}^{\lambda,u}-1)^n.
\]
Then we have the following:
\[
 \exp\Bigl(
 -2\pi\sqrt{-1}
 \bigl(
 \alpha(b,\omega)+N_{\omega}^{\lambda}
 \bigr)
 \Bigr)
=M^{\lambda}_{\omega}.
\]

Let $s$ be an element of $\EE\big(H(\nbige^{\lambda}),\omega\big)$.
We put as follows:
\[
 F(s,b):=
 \exp\Bigl(
  \sqrt{-1}\zeta\cdot
 \bigl(
 \alpha(b,\omega)+N_{\omega}^{\lambda}
 \bigr)
 \Bigr)\cdot s
=\exp\Bigl(
  \log z\cdot
 \bigl(
 \alpha(b,\omega)+N_{\omega}^{\lambda}
 \bigr)
 \Bigr)\cdot s.
\]
Then $F(s,b)$ induces the holomorphic section of $X-D$.

\begin{lem} \label{lem;9.6.3}
Let $s$ be an element of $\EE\bigl(H(\nbigelambda),\omega\bigr)$.
We have the following equality:
\begin{equation}\label{eq;9.6.2}
-\ord\bigl(F(s,b)\bigr)= -\ord(s)-\Realpart\bigl(\alpha(b,\omega)\bigr).
\end{equation}
\end{lem}
\pf
We take $s_1,\ldots,s_l$
as $s_h=(N_{\omega}^{\lambda})^h\cdot s$.
We have the following equality:
\begin{equation} \label{eq;9.6.1}
 F(s,b)
=z^{\alpha(b,\omega)}\cdot
 \Bigl(
 s+\sum_{i=1}^l\frac{1}{i!}s_i
 \Bigr).
\end{equation}
Hence we obtain $-\ord(F(s,b))\leq -\ord(s)-\Realpart(\alpha(b,\omega))$.

Let us consider the case
that $s\in \nbigf_b$ and $N_{\omega}s\in \nbigf_{<b_1}$.
Then $s_i$ $(i=1,\ldots,l )$ above are contained in $\nbigf_{<b_1}$,
and thus we have
$-\ord(F(s,b))=-\ord(s)-\Realpart\bigl(\alpha(b,\omega)\bigr)$.

Assume $s\in \nbigf_{b_1}$.
The number $i(s)$ is determined for $s$ by the following condition:
\[
 s_{i(s)}\not\in \nbigf_{<b_1},
\quad
 s_{i(s)+1}\in\nbigf_{<b_1}.
\]
To show the equality (\ref{eq;9.6.2}),
we use an induction on $i(s)$.
If $i(s)=0$, then we have already shown the claim.
We assume that the claim holds for any $s$ such that $i(s)<i_0$,
and we will show the claim for $s$ such that $i(s)=i_0$.
Note that $i(N_{\omega}s)=i_0-1$,
and we have the equality
$-\ord(F(N_{\omega}s))=-\ord(N_{\omega}s)-\Realpart(\alpha(b,\omega))=b_3$
due to the hypothesis of the induction.
Note the following equality:
\[
 \DD^{\lambda,f}(F(s,b))
=\alpha(b,\omega)\cdot F(s,b)\cdot\frac{dz}{z} 
+F(N_{\omega}\!\cdot\! s,b)\cdot\frac{dz}{z}.
\]
Assume $b_2:=-\ord(F(s,b))<-\ord(s)-\Realpart(\alpha(b,\omega))=:b_3$,
and we will derive the contradiction.
Note that $\DD^{\lambda,f}F(s,b)$
and $\alpha(b,\omega)\cdot F(s,b)\frac{dz}{z}$
are sections of $\prolongg{b_2+1}{\nbige}\otimes\Omega^{1,0}_X$.
On the other hand,
$F(N_{\omega}\!\cdot\! s,b)\frac{dz}{z}$
is a section of $\prolongg{b_3+1}{\nbige}\otimes\Omega^{1,0}_X$
such that it is not $0$ in $\Gr^F_{b_3+1}$.
Hence we have arrived at the contradiction.
It implies 
$-\ord(F(s,b))=-\ord(s)-\Realpart(\alpha(b,\omega))$,
and thus the induction on $i(s)$ can proceed.
Hence we are done.
\hfill\qed

\subsubsection{The compatibility of the KMS-structures}

Let us pick real numbers $c$ and $a$,
and a complex number $\omega\in\cnum$.
Recall we have the following inequality, by definition:
\[
 c-a-1<-\Realpart\bigl(\alpha(a-c,\omega)\bigr)\leq c-a.
\]
In the case $-\ord(s)\leq a$,
we obtain the following:
\[
 -\ord\bigl(F(s,a-c)\bigr)=
 -\ord(s)-\Realpart\bigl(
 \alpha(a-c,\omega)\bigr)
 \leq a-\Realpart\bigl(
 \alpha(a-c,\omega)\bigr)
 \leq c.
\]
Hence $F(s,a-c)$ gives a section of $\prolongg{c}{\nbigelambda}$.

We put
$d(a,\omega):=a-\Realpart\bigl(\alpha(a,\omega)\bigr)<0$.
Then the morphism
$\nbigf_a\EE\bigl(H(\nbigelambda),\omega \bigr)
\lrarr
 \prolongg{d(a,\omega)}{\nbigelambda}$
is given by the correspondence
$s\longmapsto F(s,a)$,
and we obtain the following induced morphism:
\[
 \varphi_{(a,\omega)}:
 \Gr^{\nbigf}_a\EE\bigl(
 H(\nbigelambda),\omega
 \bigr)
\lrarr
 \Gr^F_{d(a,\omega)}(\prolong{\nbigelambda}_{|O}).
\]
\begin{lem} \label{lem;10.11.5}
The morphism $\varphi_{(a,\omega)}$ is injective.
\end{lem}
\pf
It follows from Lemma \ref{lem;9.6.3}.
\hfill\qed

\vspace{.1in}
We have the action of $\Gr^{\nbigf}_a(N_{\omega}^{\lambda})$
on
$\Gr^{\nbigf}_a\EE(H(\nbigelambda),\omega)$
for each $a$.
Then we obtain the following endomorphism:
\[
 \alpha(a,\omega)+\Gr^{\nbigf}_a(N_{\omega}^{\lambda})
\in
 \End\Bigl(
 \Gr^{\nbigf}_{a}
  \EE\bigl(
  H(\nbige^{\lambda}),\omega
 \bigr)
 \Bigr).
\]
On the other hand,
the endomorphism $\Res(\DD^{\lambda,f})$
on $\prolong{\nbigelambda}_{|O}$
induces the following endomorphism:
\[
 \Gr^F_{d(a,\omega)}\bigl(\Res(\DD^{\lambda,f})\bigr)
\in
 \End\Bigl(
 \Gr^F_{d(a,\omega)}\bigl(\prolong{\nbigelambda}_{|O}\bigr)
 \Bigr).
\]

\begin{lem} \label{lem;a11.13.50}
We have the following equality:
\[
\varphi_{(a,\omega)}\circ
\Bigl(
 \alpha(a,\omega)+\Gr^{\nbigf}_a(N_{\omega}^{\lambda})
\Bigr)
=
 \Gr^F_{d(a,\omega)}(\Res(\DD^{\lambda,f}))
\circ
 \varphi_{(a,\omega)}.
\]
\end{lem}
\pf
We have the following morphism due to $F(\cdot,a)$:
\[
 \Bigl(
 \nbigf_a\EE(H(\nbigelambda),\omega),\,\,
 \alpha(a,\omega)+\Gr^{\nbigf}_a(N_{\omega}^{\lambda})
 \Bigr)
\lrarr
 \Bigl(
 \Gamma(X,\prolongg{d}{\nbigelambda}),\,\,\,\DD^{\lambda,f}
 \Bigr)
\]
Here we put $d:=d(a,\omega)$.

We also have the following morphisms:
\[
 \Gamma(X,\prolongg{d}{\nbigelambda})
\lrarr
 F_d(\prolong{\nbigelambda}_{|O})
\lrarr
 \Gr^F_d(\prolong{\nbigelambda}_{|O}).
\]
These morphisms are equivariant
with respect to the operators
$\DD^{\lambda,f}$,
$\Res(\DD^{\lambda,f})$ and
$\Gr^F_d(\Res(\DD^{\lambda,f}))$
respectively.
Thus we are done.
\hfill\qed

\begin{cor} \label{cor;9.11.60}
We have the following implication
in $\Gr^F_d(\prolong{\nbigelambda}_{|O})$:
\[
 \Image(\varphi_{(a,\omega)})
\subset
 \EE\Bigl(
 \Gr^F_d\bigl(\Res(\DD^{\lambda,f})\bigr),\,\,
 \alpha(a,\omega)
 \Bigr)
=\EE\bigl(
 \Gr^F_d\bigl(\Res(\DD^{\lambda})\bigr),\,\,
 \lambda\cdot\alpha(a,\omega)
 \bigr).
\]
Here we put $d:=d(a,\omega)=a-\Realpart(\alpha(a,\omega))<0$.
\end{cor}
\pf
It immediately follows from Lemma \ref{lem;a11.13.50}.
\hfill\qed

\vspace{.1in}

For any $u_1=(a,\omega)\in\real\times\cnum^{\ast}$,
we put as follows:
\[
 \twistmap(u_1)=
 \bigl(d(a,\omega),\,\,\lambda\cdot\alpha(a,\omega)\bigr).
\]
Thus we obtain the map
$\twistmap:\real\times\cnum^{\ast}\lrarr ]-1,0]\times\cnum$.

\begin{lem}
$\twistmap$ is bijective.
\end{lem}
\pf
It can be checked by a direct calculation.
\hfill\qed

\vspace{.1in}

For any $u_1=(a,\omega)\in \real\times\cnum^{\ast}$,
we put as follows:
\[
 \Gr^{\nbigf,\EE}_{u_1}\bigl(H(\nbige^{\lambda})\bigr)
:=
 \Gr^{\nbigf}_{a}\EE\bigl(H(\nbigelambda),\omega\bigr).
\]
Then we obtain the injection:
\[
 \varphi_{u_1}:
 \Gr^{\nbigf,\EE}_{u_1}\bigl(H(\nbige^{\lambda})\bigr)
\lrarr
 \Gr^{F,\EE}_{\gminit(u_1)}(\prolong{\nbigelambda})
\]
Then we obtain the following injection:
\[
 \bigoplus_{u_1\in\real\times\cnum^{\ast}}
 \varphi_{u_1}:
\bigoplus_{u_1\in\real\times\cnum^{\ast}}
 \Gr^{\nbigf,\EE}_{u_1}\bigl(H(\nbige^{\lambda})\bigr)
\lrarr
 \bigoplus_{u\in \KMS(\prolong{\nbigelambda})}
 \Gr^{F,\EE}_{u}(\prolong{\nbigelambda}).
\]

\begin{prop} \label{prop;9.6.10}
The morphism $\bigoplus_{u_1\in\real\times\cnum^{\ast}}
 \varphi_{u_1} $ is isomorphic.
Each $\varphi_{u_1}$ is isomorphic.
\end{prop}
\pf
We have already known that $\varphi_{u_1}$ are injective
(Lemma \ref{lem;10.11.5}).
We have the following equalities:
\[
 \sum_{u_1}\dim\Gr^{\nbigf,\EE}_{u_1}H(\nbigelambda)
=\rank\nbigelambda
=\sum_{u\in\KMS(\prolong{\nbigelambda})}
 \dim\Gr^{F,\EE}_u(\prolong{\nbigelambda}).
\]
The claims follow from the equalities above.
\hfill\qed

\vspace{.1in}

We have the weight filtration $W$
on $\Gr_u^{\nbigf,\EE}H(\nbigelambda)$
induced by the nilpotent map $\Gr^{\nbigf}_a(N_{\omega})$,
where we have $u=(a,\omega)$.
\begin{lem}
The morphism $\varphi_{u}$ preserves the weight filtrations $W$.
\end{lem}
\pf
It immediately follows from Lemma \ref{lem;a11.13.50}.
\hfill\qed

%% file: 4.1.tex

\subsubsection{Norm estimate for the multi-valued flat sections}


Let $\vecs=(s_i)$ be a frame of $H(\nbigelambda)$
satisfying the following conditions:
\begin{condition}\mbox{{}}
\begin{enumerate}
\item \label{8.13.2}
 It is compatible with $\EE$.
 We put $\deg^{\EE}(s_i)=\omega_i$.
\item \label{8.13.3}
 It is compatible with $\nbigf$ on
 $\EE\bigl(H(\nbigelambda),\omega\bigr)$ for any $\omega$.
 We put $\deg^{\nbigf}(s_i)=a_i$.
\item
 From the conditions {\rm\ref{8.13.2}} and {\rm\ref{8.13.3}} above,
 we obtain the induced frame $\vecs^{(1)}$
 on $\Gr^{\nbigf,\EE}H(\nbigelambda)$.
 The frame $\vecs^{(1)}$ is compatible with the weight filtration $W$.
 We put $\deg^{W}(s^{(1)}_i)=k_i$.
\end{enumerate}
If the conditions above are satisfied,
we say that $\vecs$ is compatible with
$\nbigf$, $\EE$ and $W$.
\hfill\qed
\end{condition}

The matrix $(b^{(n)}_{j\,i})$ is determined by the following condition:
\[
 (N^{\lambda})^n s_i
=\sum_j b_{j\,i}^{(n)}\cdot s_j.
\]

\begin{lem}
We have $b_{j\,i}=0$ in the following cases:
\begin{itemize}
\item
$\omega_0\neq \omega_j$.
\item
$\omega_i=\omega_j$ and $a_i<a_j$.
\item
$(\omega_i,a_i)=(\omega_j,a_j)$
and $k_i-n< k_j$.
\end{itemize}
\end{lem}
\pf
It is clear from our choice of $\vecs$.
\hfill\qed

\vspace{.1in}

We put $v_i:=F(s_i,a_i)$,
and $\vecv=(v_i)$.
Then $\vecv$ is a tuple of $\prolong{\nbigelambda}$.
We put $\alpha_i=\alpha(a_i,\omega_i)$.

\begin{lem}\label{lem;d11.14.20}
The tuple of sections
$\vecv$ is a frame of $\prolong{\nbige}$ which is compatible with
$\EE$, $F$ and $W$.
\end{lem}
\pf
It immediately follows from Proposition \ref{prop;9.6.10}
and Lemma \ref{lem;a11.13.50}.
\hfill\qed

\begin{lem} \label{lem;9.6.15}
We have the following equality:
\[
 (N^{\lambda})^nv_i
=\sum_{\substack{\omega_i=\omega_j,\\a_i\geq a_j }}
 b_{j\,i}^{(n)}\cdot z^{\alpha_i-\alpha_j}\cdot v_j.
\]
Note that if $\omega_i=\omega_j$ and $a_i\geq a_j$,
then $\alpha_i-\alpha_j$ is a non-negative integer.
\end{lem}
\pf
We have the following equalities:
\[
 (N^{\lambda})^nv_i
=F\bigl( (N^{\lambda})^ns_i,-a_i\bigr)
=\sum_{\substack{\omega_i=\omega_j,\\a_i\geq a_j }}
 b_{j\,i}^{(n)}\cdot F(s_j,-a_i)
=\sum_{\substack{\omega_i=\omega_j,\\a_i\geq a_j }}
 b_{j\,i}^{(n)}\cdot z^{\alpha_i-\alpha_j}\cdot F(s_j,-a_j).
\]
Thus we are done.
\hfill\qed

\vspace{.1in}

Let $C$ be any positive number.
We put $\alpha_i:=\alpha(a_i,\omega_i)$.
On the region $\bigl\{z\,\big|\,|\arg z|<C\bigr\}$,
we have the following equalities:
\begin{equation} \label{eq;10.11.6}
 s_i
 =z^{-\alpha_i}\cdot \exp(-\log z\cdot N_{\omega}^{\lambda})\cdot v_i
 =z^{-\alpha_i}
 \sum_{n=0}^{\infty}
 \frac{1}{n!}(-\log z)^{n}\cdot (N^{\lambda})^n\cdot v_i.
\end{equation}
It can be described as 
$s_i=z^{-\alpha_i}\sum f_j\cdot v_j$
for some multi-valued holomorphic functions $f_j$.

\begin{lem} \label{lem;10.11.7}
Let $C$ be a positive number.
We have the following,
on the region $\bigl\{|\arg z|<C\bigr\}$:
\begin{itemize}
\item
$f_i=1$.
\item
$|f_j|\leq C\cdot(-\log|z|)^{(k_i-k_j)/2}$
for some positive constant $C$,
in the case $(a_i,\omega_i)=(a_j,\omega_j)$
   and $k_i>k_j$.
\item
 $|f_j|\leq C\cdot (-\log|z|)^M$
for some positive constants $C$ and $M$,
in the case $\omega_i=\omega_j$ and $a_i>a_j$.
\item
Otherwise,
 $f_j$ vanishes identically.
\end{itemize}
\end{lem}
\pf
It immediately follows from Lemma \ref{lem;9.6.15}
and (\ref{eq;10.11.6}).
\hfill\qed

\vspace{.1in}

We put $s_i':=s_i\cdot |z|^{a_i}\cdot (-\log|z|)^{-k_i/2}$,
and $\vecs':=(s_i')$.
\begin{prop} \label{prop;9.7.55}
$\vecs'$ is adapted on the region
$\bigl\{s\,\big|\,|\arg z|<C\bigr\}$ for any positive constant $C$.
\end{prop}
\pf
It is a direct corollary of Lemma \ref{lem;10.11.7}
and the adaptedness of $\vecv'$.
\hfill\qed


\subsubsection{The decomposition and the filtration of
the flat bundle $\nbigelambda$}
\label{subsubsection;a11.15.3}

We have the generalized eigen decomposition over $X-D$
for the monodromy $M^{\lambda}$:
\begin{equation}\label{eq;8.13.10}
\nbigelambda=
 \bigoplus_{\omega\in\Sp(M^{\lambda})}\EE(\nbigelambda,\omega).
\end{equation}

\begin{cor}
The decomposition is quasi adapted (Definition {\rm\ref{df;10.11.10}}).
\end{cor}
\pf
Let $\vecv$ be the frame of $\prolong{\nbigelambda}$
obtained from $\vecs$.
Then $\vecv$ is compatible with the generalized eigen decomposition
above,
and $\vecv'$ is adapted.
Thus the decomposition is quasi adapted.
\hfill\qed

\vspace{.1in}


Let $c$ be any real number.
Let $\vecs$ be a frame of $H(\nbigelambda)$,
compatible with $\EE$, $\nbigf$ and $W$.
We put ${}_cv_i:=F(s_i,a_i-c)$,
and ${}_c\vecv=({}_cv_i)$.
We put ${}_c\alpha_i:=\alpha(a_i-c_i,\omega_i)$.

\begin{lem}
The tuple of sections ${}_c\vecv$ is a frame of $\prolongg{c}{\nbige}$
compatible with $F$, $\EE$ and $W$.
We have the following:
\begin{equation}\label{eq;9.6.20}
 \DD^{\lambda,f}({}_cv_i)
=\Bigl(
 {}_c\alpha_i\cdot{}_cv_{i}
+\sum_{\substack{
 \omega_i=\omega_j,\\ a_i\geq a_j
 }}
 b^{(1)}_{j\,i}\cdot z^{{}_c\alpha_i-{}_c\alpha_j}\cdot
 {}_cv_j
 \Bigr)
\cdot\frac{dz}{z}.
\end{equation}
\end{lem}
\pf
Since $\vecv$ is a frame of $\prolong{\nbigelambda}$
compatible with the parabolic filtration,
it is easy to check ${}_c\vecv$ is a frame of $\prolongg{c}{\nbige}$,
compatible with the parabolic filtration.
The equality (\ref{eq;9.6.20}) follows from the following:
\begin{multline}
 \DD^{\lambda,f}({}_cv_i)
=\DD^{\lambda,f}\bigl(
 z^{{}_c\alpha_i}\cdot\exp\bigl(
 \log z\cdot N^{\lambda}_{\omega_i}
 \bigr)\cdot s_i
 \bigr)
=\bigl(
 {}_c\alpha_i+N^{\lambda}_{\omega_i}
 \bigr)\cdot{}_cv_i\cdot\frac{dz}{z}\\
=\Bigl(
 {}_c\alpha_i\cdot {}_cv_i
+\sum_{\substack{\omega_i=\omega_j,\\ a_i\geq a_j}}
 b_{j\,i}^{(1)}\cdot z^{{}_c\alpha_i-{}_c\alpha_j}
 \cdot
 {}_cv_j
 \Bigr)
\frac{dz}{z}.
\end{multline}
The compatibility of ${}_c\vecv$ with $\EE$ and $W$
follows from the formula (\ref{eq;9.6.20}).
\hfill\qed

\begin{cor}
The decomposition {\rm(\ref{eq;8.13.10})} is prolonged to
the following:
\[
 \prolongg{c}{\nbigelambda}=
\bigoplus_{\omega\in\Sp(M^{\lambda})}
 \prolongg{c}{\EE(\nbigelambda,\omega)}.
\]
In particular, $\prolongg{c}{\EE(\nbigelambda,\omega)}$
is locally free.
\end{cor}
\pf
It is easy to check the claim
by using the frame ${}_c\vecv$.
\hfill\qed

\vspace{.1in}
We obtain two decomposition of $\prolongg{c}{\nbigelambda}_{|O}$:
\[
 \prolongg{c}{\nbigelambda}_{|O}
=\bigoplus_{\beta}\EE\bigl(\Res\DD^{\lambda},\beta\bigr)
=\bigoplus_{\omega\in\Sp(M^{\lambda})}
   \prolongg{c}{\EE(\nbigelambda,\omega)}_{|O}.
\]

\begin{cor}\label{cor;a11.17.20}
The following holds:
\[
 \prolongg{c}{\EE(\nbigelambda,\omega)}_{|O}
=\bigoplus_{\exp(-2\pi\sqrt{-1}\lambda^{-1}\cdot\beta)=\omega}
 \EE\bigl(\Res\DD^{\lambda},\beta\bigr).
\]
\end{cor}
\pf
It immediately follows from the formula (\ref{eq;9.6.20}).
\hfill\qed

\vspace{.1in}

The filtration $\nbigf$ on $\EE(H(\nbigelambda),\omega)$
induces the filtration $\nbigf\bigl(\EE(\nbigelambda,\omega)\bigr)$
on $\EE(\nbigelambda,\omega)$.
\begin{lem}
The filtration $\nbigf\bigl(\EE(\nbigelambda,\omega)\bigr)$
can be prolonged to the filtration
$\prolongg{c}{\nbigf\bigl(
 \EE(\nbigelambda,\omega)})\bigr)$.
We have the following equality:
\[
 \prolongg{c}{\nbigf_a\bigl(
 \EE(\nbigelambda,\omega)\bigr)}
=\big\langle
 {}_cv_i\,\big|\,
 \omega_i=\omega,
 a_i\leq a
 \big\rangle.
\]
We also have 
$\prolongg{c}{\nbigf_a (\nbigelambda)}=
 \bigoplus_{\omega}\prolongg{c}{\nbigf_a\bigl(
 \EE(\nbigelambda,\omega)\bigr)}$.
\end{lem}
\pf
It is clear from our definition.
\hfill\qed

\vspace{.1in}

Then we obtain the two filtrations on $\prolongg{c}{\nbigelambda}_{|O}$,
one is
$F(\prolongg{c}{\nbigelambda}_{|O})$
and the other is
$\prolongg{c}{\nbigf(\nbigelambda)}_{|O}$.

\begin{lem}
We have the following relation:
\[
 \prolongg{c}{\nbigf_a\bigl(
 \EE(\nbigelambda,\omega)\bigr)}_{|O}
=F_{d(a,\omega)}\bigl(
 \EE(\prolongg{c}{\nbigelambda}_{|O},
 \lambda\cdot\alpha(a-c,\omega))\bigr)
\oplus
\bigoplus_{
 \substack{
  \exp(-2\pi\sqrt{-1}\lambda^{-1}\beta)=\omega,\\
  \Re(\lambda^{-1}\beta)<\Re(\alpha(a-c,\omega))
 }
 }\!\!
 \EE(\prolongg{c}{\nbigelambda}_{|O},\beta).
\]
\end{lem}
\pf
It can be shown by using the frame ${}_c\vecv$.
\hfill\qed

\vspace{.1in}

The nilpotent morphism $N^{\lambda}_{\omega}$ on
$\EE(H(\nbigelambda),\omega)$
induces the endomorphism $N^{\lambda}_{\omega}$
on $\EE(\nbigelambda,\omega)$ over $X-D$.
It is prolonged to the endomorphism
of $\prolongg{c}{\EE(\nbigelambda,\omega)}$ over $X$.
It preserves the filtration
$\nbigf\bigl(\EE(\nbigelambda,\omega)\bigr)$.
Then we obtain the nilpotent morphism
$\Gr^{\nbigf,\EE}_u(N^{\lambda})$
on $\Gr^{\nbigf,\EE}_u(\nbigelambda)$ over $X$.

\begin{lem}
The conjugacy class of $\Gr^{\nbigf,\EE}_u(N^{\lambda})_{|P}$
is independent of $P\in X$.
Hence we obtain the weight filtration $W$
of $\Gr^{\nbigf,\EE}_u(\nbigelambda)$
in the category of the vector bundles.
\hfill\qed
\end{lem}

%% file: 4.2.tex

\subsubsection{Functoriality for tensor products}

Let $(E_i,\delbar_{E_i},h_i,\theta_i)$ $(i=1,2)$
be harmonic bundles over $X-D$.
We denote the deformed holomorphic bundle
by $\nbigelambda_i$ $(i=1,2)$.
We have the natural following isomorphism:
\[
 H(\nbigelambda_1\otimes\nbigelambda_2)
\simeq
 H(\nbigelambda_1)
   \otimes
 H(\nbigelambda_2).
\]

\begin{lem}
We have the natural isomorphisms:
\[
 \EE\bigl(H(\nbigelambda_1\otimes\nbigelambda_2),\omega\bigr)=
\bigoplus_{\omega_1\times\omega_2=\omega}
 \EE\bigl(H(\nbigelambda_1),\omega_1\bigr)
\otimes
 \EE\bigl(H(\nbigelambda_2),\omega_2\bigr).
\]
\hfill\qed
\end{lem}

We obtain the corresponding decomposition:
\begin{equation} \label{eq;8.13.11}
 \EE(\nbigelambda,\omega)=
 \bigoplus_{\omega_1\times\omega_2=\omega}
 \EE(\nbigelambda_1,\omega_1)
\otimes
 \EE(\nbigelambda_2,\omega_2).
\end{equation}

\begin{lem} \mbox{{}}
\begin{enumerate}
\item
 The decomposition {\rm (\ref{eq;8.13.11})} is quasi adapted
(Definition {\rm\ref{df;10.11.10}}).
\item
 We have the following:
 \[
  \nbigf_a\EE(H(\nbigelambda_1\otimes\nbigelambda_2),\omega)
  =
 \bigoplus_{\omega_1\times\omega_2=\omega}
 \sum_{a_1+a_2\leq a}
 \nbigf_{a_1}\EE(H(\nbigelambda_1),\omega_1)
\otimes
 \nbigf_{a_2}\EE(H(\nbigelambda_2),\omega_2).
 \]
\end{enumerate}
\end{lem}
\pf
 Let $\vecs_i$ be a frame of $H(\nbigelambda_i)$
 compatible with $\EE$, $\nbigf$ and $W$.
 Then we obtain the adapted frame $\vecs_i'$.
 By using $\vecs_1'$ and $\vecs_2'$,
 we obtain the first claim.
 The second claim follows from the first claim.
\hfill\qed

\vspace{.1in}

We have the following product of
the abelian group $\real\times\cnum^{\ast}$:
For $u_i=(a_i,\alpha_i)$,
we put $u_1\cdot u_2=(a_1+a_2,\alpha_1\times \alpha_2)$.

\begin{cor}
We have the isomorphism:
\[
 \Gr^{\nbigf,\EE}_u\bigl(H(\nbigelambda)\bigr)
\simeq
 \bigoplus_{u_1\cdot u_2=u}
 \Gr^{\nbigf,\EE}_{u_1}\bigl(
 H(\nbigelambda_1)
 \bigr)
\otimes
 \Gr^{\nbigf,\EE}_{u_2}\bigl(
 H(\nbigelambda_2)
 \bigr).
\]
\hfill\qed
\end{cor}

We have the natural isomorphisms:
\begin{equation}\label{eq;10.11.15}
 \EE\bigl(\bigwedge^a H(\nbigelambda),\omega\bigr)
\simeq
 \bigoplus_{f\in \nbigs(a,\omega)}
 \bigotimes_{\omega'\in\Sp(M^{\lambda})}
 \bigwedge^{f(\omega')}
 \EE(H(\nbigelambda),\omega'),
\end{equation}
\begin{equation} \label{eq;10.11.16}
 \EE\bigl(\Sym^aH(\nbigelambda),\omega\bigr)
 \simeq
 \bigoplus_{f\in \nbigs(a,\omega)}
 \bigotimes_{\omega'\in\Sp(M^{\lambda})}
 \Sym^{f(\omega')}
 \EE(H(\nbigelambda),\omega'),
\end{equation}
\begin{equation} \label{eq;10.11.17}
 \nbigs(a,\omega):=
 \Bigl\{
 f:\Sp(M^{\lambda})\lrarr\seisuu_{\geq\,0}
 \,\Big|\,
 \prod \omega^{\prime\,f(\omega')}=\omega,
 \,\,
 \sum f(\omega')=a
 \Bigr\}.
\end{equation}

\begin{cor}\mbox{{}}
\begin{itemize}
\item
 The decompositions {\rm(\ref{eq;10.11.15})},
 {\rm(\ref{eq;10.11.16})} and {\rm(\ref{eq;10.11.17})}
 are quasi adapted.
\item
 The parabolic filtrations on the left hand sides
 of {\rm(\ref{eq;10.11.15})},
 {\rm(\ref{eq;10.11.16})} and {\rm(\ref{eq;10.11.17})}
 are isomorphic to the induced filtrations on the right hand side.
\item
 The weight filtrations are also isomorphic.
\hfill\qed
\end{itemize}
\end{cor}

\subsubsection{Functoriality for dual}

\label{subsubsection;d11.14.21}

We have the natural isomorphism:
$H(\nbige^{\lor\,\lambda})\simeq H(\nbigelambda)^{\lor}$.

\begin{lem} \label{lem;9.11.23}
Under the isomorphism,
$\nbigf_aH(\nbige^{\lor\,\lambda})$ is same as the following:
\[
 \Bigl\{
 f\in H(\nbigelambda)^{\lor}\,\Big|\,
 f(\nbigf_{b}H(\nbigelambda))\subset
 \nbigf_{b+a}H(\nbigelambda)
 \Bigr\}.
\]
\end{lem}
\pf
Let $\vecs=(s_i)$ be a base of $H(\nbigelambda)$,
compatible with $\EE$, $\nbigf$, $W$.
Let $\vecs^{\lor}=(s_i^{\lor})$ denote the dual base.

We put $\deg^{\nbigf}(s_i)=a_i$.
We put $s_i':=s_i\cdot |z|^{a_i}$,
and $\vecs':=(s_i')$.
Then $\vecs'$ is adapted up to log order.
We put $s_i^{\lor\,\prime}:=s_i^{\lor}\cdot |z|^{-a_i}$,
and $\vecs^{\lor\,\prime}=(s_i^{\lor\,\prime})$.
Then $\vecs^{\lor\,\prime}$ is the dual base of
$\vecs'$,
and $\vecs^{\lor\,\prime}$ is adapted up to log order.
Then the claim follows easily.
\hfill\qed

\subsubsection{Functoriality for pull back}

Let $\psi_c:X\lrarr X$ given by $z\longmapsto z^c$.
We have the natural isomorphism:
\[
 \psi_c^{-1}:H(\nbigelambda)\simeq
 H(\psi_c^{-1}\nbigelambda).
\]
Let $M^{\lambda}_1$ denote the monodromy
of $\psi_c^{-1}\nbigelambda$.
We obtain the following isomorphism
for any $\omega_1\in\Sp(M^{\lambda_1})$:
\[
 \EE\bigl(H(\psi_c^{-1}\nbigelambda),\omega_1\bigr)
\simeq\bigoplus_{
   \substack{\omega\in \Sp(M^{\lambda})\\
   \omega^c=\omega_1  } }
 \EE(H(\nbigelambda),\omega).
\]

\begin{lem}\label{lem;a11.17.1}
We have the following, for any element $\omega_1\in \Sp(M^{\lambda}_1)$:
\begin{equation}\label{eq;8.14.1}
 \nbigf_{c\cdot a}\bigl(
 \EE\bigl(H(\psi_c^{-1}\nbigelambda),\omega_1\bigr)
 \bigr)
=\bigoplus_{
   \substack{\omega\in \Sp(M^{\lambda})\\
   \omega^c=\omega_1  } }
 \nbigf_a\bigl(
 \EE\bigl(H(\nbigelambda),\omega\bigr)\bigr).
\end{equation}
The weight filtrations are compatible.
\end{lem}
\pf
The compatibility for the weight filtration is clear.
Let $\vecs$ be a frame of $H(\nbigelambda)$
compatible with $\EE$, $\nbigf$ and $W$.
We put $a_i:=\deg^{\nbigf}(s_i)$,
and $s_i':=s_i\cdot |z|^{a_i}$.
Then $\vecs'=(s_i')$ is adapted up to log order.

Let consider $\psi_c^{-1}(s_i')=\psi_c^{-1}(s_i)\cdot |z|^{c\cdot a_i}$.
The frame $\psi_c^{-1}(\vecs')=\bigl(\psi_c^{-1}(s_i')\bigr)$
is adapted up to log order.
Then the equality (\ref{eq;8.14.1}) follows.
\hfill\qed

\vspace{.1in}

\subsubsection{The correspondence of KMS-spectrum}

We put as follows:
\[
 \overline{\KMS}^f(\nbigelambda):=
 \Bigl\{
  u\in \real\times\cnum^{\ast}\,
 \Big|\,
 \dim \Gr^{\nbigf,\EE}_uH(\nbigelambda)\neq 0
 \Bigr\}.
\]
The number
$\multiplicity^f(\lambda,u):=
 \dim \Gr^{\nbigf,\EE}_uH(\nbigelambda)$
is called the multiplicity.

\vspace{.1in}

The maps
$\paramap^f(\lambda):\real\times\cnum\lrarr \real$,
$\eigenmap^f(\lambda):\real\times\cnum\lrarr\cnum^{\ast}$
and $\kmsmap^f:\real\times\cnum\lrarr  \real\times\cnum^{\ast}$
are defined in the subsubsection \ref{subsubsection;b11.11.11}.

On the other hand, we put as follows, for any real number $c$:
\begin{equation}
 \nbigk(\nbige,\lambda,c)
:=\bigl\{
 u\in\KMS\bigl(\nbige^0\bigr)\,\big|\,
 c-1\leq \paramap(\lambda,u)<c
  \bigr\}.
\end{equation}
Then we have the isomorphism
$\kmsmap(\lambda):\nbigk(\nbige,\lambda,0)\lrarr
 \KMS\bigl(\prolongg{c}{\nbigelambda}\bigr)$.
Let us consider the case $c=0$.
\begin{lem}\label{lem;a11.15.2}
Let $u$ be an element of $\nbigk(\nbige,\lambda,0)$.
We have the relation
$\kmsmap(\lambda,u)=\gminit\bigl(\kmsmap^f(\lambda,u)\bigr)$.
\end{lem}
\pf
From the formula \ref{eq;a11.15.1} and the inequality
$-1<\paramap(\lambda,u)\leq 0$, we have the the following:
\[
 \eigenmap^f(\lambda,u)=
 \exp\bigl(
 -2\pi\sqrt{-1}\lambda^{-1}\cdot\eigenmap(\lambda,u),
\quad
 \paramap^f(\lambda,u)
\leq 
 \Re\bigl(\lambda^{-1}\cdot\eigenmap(\lambda,u)\bigr)
<\paramap^f(\lambda,u)+1.
\]
Thus we obtain the following, by definition of $\alpha(b,\omega)$:
\[
 \lambda^{-1}\cdot \eigenmap(\lambda,u)
=\alpha\bigl(\paramap^f(\lambda,u),\eigenmap^f(\lambda,u)\bigr).
\]
We also obtain the following:
\[
 \paramap(\lambda,u)=
 \paramap^f(\lambda,u)
-\Re\bigl(\alpha(\paramap^f(\lambda,u)),\eigenmap^f(\lambda,u)\bigr)
=d\bigl(\paramap^f(\lambda,u),\eigenmap^f(\lambda,u) \bigr).
\]
It means $\kmsmap(\lambda,u)=\gminit\bigl(\kmsmap^f(\lambda,u)\bigr)$.
\hfill\qed

\begin{lem} \label{lem;10.11.20}
The image of $\KMS(\nbige^0)$ via the morphism $\kmsmap^f(\lambda)$ 
is $\KMSoverline^f(\nbigelambda)$,
and we have the equality:
$\multiplicity^f(\lambda,\kmsmap^f(\lambda,u))=
 \multiplicity(\lambda,\kmsmap(\lambda,u))$.
\end{lem}
\pf
From Proposition \ref{prop;9.6.10}
and Lemma \ref{lem;a11.15.2},
the image of $\nbigk(\nbige,\lambda,0)$ via the morphism
$\kmsmap^f(\lambda)$ is same as $\KMSoverline^f(\nbigelambda)$,
and we have the equality of the multiplicity.
Note that we have 
the equalities $\kmsmap^f(\lambda,u+(1,0))=\kmsmap^f(\lambda,u)$
and
$\multiplicity(\lambda,\kmsmap(\lambda,u))
=\multiplicity\bigl(\lambda,\kmsmap(\lambda,u+(1,0))\bigr)$.
Thus we are done.
\hfill\qed

\begin{lem}\mbox{{}}
\begin{enumerate}
\item
We have the $\seisuu$-action on $\KMSE{0}$.
The map $\kmsmap^f(\lambda)$ induces the isomorphism
$\KMSEoverline{0}\simeq \overline{\KMS}^f(\nbigelambda)$,
which is also denoted by $\kmsmap^f(\lambda)$.
\item
We have
$\multiplicity(0,u)
=\multiplicity^f(\lambda,\kmsmap(\lambda,u))$.
\end{enumerate}
\end{lem}
\pf
It follows from Corollary \ref{cor;9.11.1}
and Lemma \ref{lem;10.11.20}.
\hfill\qed

\vspace{.1in}

Let $\vecs$ be a frame of $H(\nbigelambda)$,
which is compatible with $\EE$ and $\nbigf$.
Then we have the numbers
$a_i:=\deg^{\nbigf}(v_i)$ and $\omega_i:=\deg^{\EE}(v_i)$.
Let $c$ be any real number.
Then we have the elements $u_i\in\KMS(\nbige,\lambda,c)$
such that $\kmsmap^f(\lambda,u_i)=(a_i,\omega_i)$.
We put $v_i:=F(s_i,a_i-c)$,
and then we obtain the frame $\vecv=(v_i)\in\prolongg{c}{\nbigelambda}$
(See the subsubsection \ref{subsubsection;a11.15.3}).
\begin{lem}\label{lem;a11.15.5}
We have 
$\deg^{F,\EE}(v_i)=\kmsmap(\lambda,u_i)$.
We also have the following:
\[
 \DD\vecv=\vecv\cdot \bigl(C+N\bigr)\cdot\frac{dz}{z}.
\]
Here $C$ denotes the diagonal matrix
whose $(i,i)$-component is $\eigenmap(\lambda,u_i)$,
and $N$ denotes the nilpotent matrix.
\end{lem}
\pf
As in Lemma \ref{lem;a11.15.2},
we can show
$\alpha(a_i-c,\omega_i)
=\lambda^{-1}\cdot\eigenmap(\lambda,u_i)$
and $-\ord(v_i)=\paramap(\lambda,u_i)$.
Then the claims follow from the results in the subsubsection 
\ref{subsubsection;a11.15.3}.
\hfill\qed

\subsubsection{Genericity}
\label{subsubsection;c11.14.1}

We have the induced morphisms
$\paramap^f(\lambda):\KMSEoverline{0}\lrarr \Par^f(\nbigelambda)$
and $\eigenmap^f(\lambda):\KMSEoverline{0}\lrarr \Sp(M^{\lambda})$.

\begin{df} \label{df;a11.14.1}
$\lambda$ is called generic with respect to $\harmonicbundle$,
if the map
$\eigenmap^f(\lambda):\KMSEoverline{0}\lrarr\Sp(M^{\lambda})$
is bijective.
\hfill\qed
\end{df}

\begin{rem}
We can consider $\KMS(\prolong{\nbige}^0)$
instead of $\KMSEoverline{0}$.
\hfill\qed
\end{rem}

\begin{lem} \label{lem;10.11.21}
Let $S$ be the set of $\lambda\in\cnum$,
which are generic with respect to $\harmonicbundle$.
Then the set $\cnum^{\ast}-S$ is discrete in $\cnum^{\ast}$.
In particular, it is countable.
\end{lem}
\pf
Let pick $u=(a,\alpha)$ and $v=(b,\beta)$ be elements
of $\KMSE{0}$,
such that $(a,\alpha)\neq(b,\beta)$.
Let consider the following condition for $\lambda$:
\[
 \eigenmap^{f}(\lambda,u)=\eigenmap^f(\lambda,v).
\]
It is equivalent to the following:
\[
 \lambda^{-1}\cdot(\alpha-\beta)-(a-b)
-\lambda\cdot(\bar{\alpha}-\bar{\beta})
\in\seisuu.
\]
Let $n$ be an integer.
Let consider the following equation:
\begin{equation}\label{eq;8.14.2}
 \lambda^{-1}\cdot(\alpha-\beta)-(a-b)
-\lambda\cdot(\bar{\alpha}-\bar{\beta})=n.
\end{equation}
Let $\lambda_i(n)$ $(i=1,2)$ be the solutions of the equation
(\ref{eq;8.14.2}).
Then we have the following relation:
\begin{equation}\label{eq;8.14.4}
 |\lambda_1(n)\cdot\lambda_2(n)|=
 \left|
 -\frac{\alpha-\beta}{\bar{\alpha}-\bar{\beta}}
 \right|=1.
\end{equation}
We also have the following:
\begin{equation}\label{eq;8.14.3}
 \lambda_1(n)+\lambda_2(n)
= \frac{-(a-b+n)}{\bar{\alpha}-\bar{\beta}}.
\end{equation}
Then we obtain the following:
\[
 \lim_{|n|\to \infty}
 n^{-1}\cdot \big|\lambda_1(n)+\lambda_2(n)\big|=
 \big|\bar{\alpha}-\bar{\beta}\big|^{-1}\neq 0.
\]
Hence the set of solutions of (\ref{eq;8.14.4})
is discrete in $\cnum^{\ast}$.
Since $\KMS(\prolong{\nbige^{0}})$ is finite,
the claim follows.
\hfill\qed

\vspace{.1in}

Assume $\lambda$ is generic.
Then for any $\omega\in\Sp(M^{\lambda})$,
there exists the unique element $u_0\in \KMSEoverline{0}$
satisfying $\eigenmap^f(\lambda,u_0)=\omega$.
Note that the parabolic structure
of $\EE(H(\nbigelambda),\omega)$ is trivial in the following sense:
$\Gr^{\nbigf}_b\EE(H(\nbigelambda),\omega)\neq 0$
if and only if
$b=\paramap^f(\lambda,u_0)$.

\begin{lem} \label{lem;9.7.100}
Let $Z$ be a countable subset of $\cnum_{\lambda}^{\ast}$.
Assume the following:
\begin{quote}
 For any $\lambda\in\cnum^{\ast}-Z$,
 we know the set $\Sp(M^{\lambda})$ and the multiplicity
 of each element $\alpha\in\Sp(M^{\lambda})$.
\end{quote}
Then we know the set $\KMSE{\lambda}$ and the multiplicities
$\multiplicity(\lambda,u)$
($\lambda\in\cnum$, $u\in\KMSE{\lambda}$).
\end{lem}
\pf
It follows from Lemma \ref{lem;10.11.21}.
\hfill\qed

\subsubsection{Quasi canonical prolongment}

\label{subsubsection;9.11.25}

Let $b$ be a real number.
We have the quasi canonical prolongment
$QC_b(\nbigelambda)$ of $\nbigelambda$,
that is,
a holomorphic vector bundle over $X$
satisfying the following:
\begin{itemize}
\item
 The restriction  $QC_b(\nbigelambda)_{|X-D}$ is isomorphic
 to $\nbigelambda$.
\item
 Let $g$ be a holomorphic section of $QC_b(\nbigelambda)$.
 Then  $\DD^{\lambda,f}$ gives a holomorphic section
 of $QC_b(\nbigelambda)\otimes\Omega_X(\log D)$.
\item
 Let $\beta$ be an eigenvalue of the residue $\Res(\DD^{\lambda,f})$
 on $QC_b(\nbigelambda)_D$.
 Then the inequality $b\leq \Re(\beta)<b+1$ holds.
\end{itemize}
Recall that $QC_b(\nbigelambda)$ is uniquely determined
as the subsheaf of $j_{\ast}\nbigelambda$,
where $j$ denotes the inclusion $X-D\lrarr X$.

We have the decomposition:
\[
 QC_b(\nbigelambda)=
 \bigoplus_{\omega}QC_b\bigl(\EE(\nbigelambda,\omega)\bigr).
\]
We have the natural filtration given by the following:
\[
 QC_b\bigl(F_a\EE(\nbigelambda,\omega)\bigr).
\]

Let $\vecs$ be a base of $H(\nbigelambda)$.
We put $v_i:=F(s_i,b)$ and $\vecv=(v_i)$.
Then $\vecv$ gives a holomorphic frame of $\nbigelambda$
over $X-D$.
Let us consider the prolongment of $\nbigelambda$
by $\vecv$.
Then it satisfies the conditions above.
Hence $\vecv$ gives the holomorphic frame
of $QC_b(\nbigelambda)$.

Assume that $\lambda$ is generic.
Let us consider the following:
\[
 QC_0(\nbigelambda):=
 \bigoplus_{u\in\KMSEoverline{0}}
 QC_0\bigl(\EE\bigl(\nbigelambda,\eigenmap^f(\lambda,u)\bigr)\bigr).
\]

\begin{lem}
$QC_0\bigl(\EE\bigl(\nbigelambda,\eigenmap^f(\lambda,u)\bigr)\bigr)
=\prolongg{d}{\EE\bigl(\nbigelambda,\eigenmap^f(\lambda,u)\bigr)}$.
Here we put
$d:=\paramap^f(\lambda,u)
-\Re\bigl(\alpha\bigl(0,\eigenmap^f(\lambda,u)\bigr)\bigr)$.
\end{lem}
\pf
Let $s$ be a non-zero element of
$\EE\bigl(H(\nbigelambda),\eigenmap^f(\lambda,u)\bigr)$.
Then we have the following:
\[
 -\ord(F(s,\eigenmap^f(\lambda,u)))
=-\ord(s)-\Re(\alpha(0,\eigenmap^f(\lambda,u))).
\]
Then the claim follows,
from the uniqueness of the quasi canonical prolongment
\cite{d}.
\hfill\qed

\begin{lem}
We have the following relation:
\[
 \prolongg{c}{\nbigelambda}
=\bigoplus_{u\in \KMS(\nbigelambda)}
 QC_0\bigl(\EE(\nbigelambda,\eigenmap^f(\lambda,u))\bigr)
 \cdot z^{N(u)}.
\]
Here we put
$N(u):=
 \nu_c\Bigl(\paramap^f(\lambda,u)
 -\Re\bigl(\alpha\bigl(0,\eigenmap^f(\lambda,u)\bigr)\bigr)\Bigr)$.
(See the subsubsection {\rm\ref{subsubsection;b11.11.10}}
for $\nu_c$).
\end{lem}
\pf
It follows from $c-1<-N(u)+d\leq c$.
\hfill\qed

\begin{rem}
The lemma says that quasi canonical prolongment
is essentially same as the prolongment
by an increasing order of the norms,
in the case that $\lambda$ is generic.
\hfill\qed
\end{rem}

\begin{rem}
For any point $\lambda\in \cnum^{\ast}$, not necessarily generic,
the vector bundle $\prolongg{c}{\nbigelambda}$ is obtained from
$QC_0(\nbigelambda)$
by a sequence of elementary transformations.
\hfill\qed
\end{rem}

%% file: 7.tex

\subsubsection{The structure of holomorphic bundle}

Let $\harmonicbundle$ be a tame harmonic bundle over $X-D$,
and $\nbige$ be the deformed holomorphic bundle
with $\lambda$-connection $\DD$
over $\nbigx-\nbigd$.
Let us consider the family of the multi-valued sections
$\bigl\{H(\nbigelambda)\,\big|\,\lambda\in\cnum_{\lambda}^{\ast}\bigr\}$.

Let $\pi:\hyperh\lrarr X-D$ is the universal covering,
and let $P$ be a point of $\hyperh$
We have the holomorphic vector bundle
$\pi^{-1}\nbige_{|\cnum_{\lambda}^{\ast}\times\{P\}}$
over $\cnum_{\lambda}^{\ast}$.
Since we can pick the isomorphism
$\nbige_{|(\lambda,P)}\simeq H(\nbigelambda)$,
we obtain the structure of holomorphic vector bundle
on the family
$\{H(\nbigelambda)\,|\,\lambda\in\cnum_{\lambda}^{\ast}\}$.
The holomorphic vector bundle is denoted 
by $\nbigh(\nbige)$ or simply by $\nbigh$.
Clearly the structure does not depend
on choices of $P$ and the isomorphism.

\subsubsection{The $\EEzero$-decomposition}

Pick $\lambda_0\in\cnum_{\lambda}^{\ast}$.
We have the monodromy $M^{\lambda_0}$ on $H(\nbige^{\lambda_0})$.
We put $S_0:=\Sp(M^{\lambda_0})$.
Pick a positive number $\epsilon_1$ satisfying the following:
\[
 \epsilon_1<\min\big\{|a-b|\,\,\big|\,\,a\neq b\in S_0\big\}.
\]
Pick sufficiently small $\epsilon_0>0$
such that we have the following decomposition
on $\Delta(\lambda_0,\epsilon_0)$:
\begin{equation} \label{eq;9.7.31}
  \nbigh_{|\Delta(\lambda_0,\epsilon_0)}
=\bigoplus_{\omega\in S_0}\EEzero(\nbigh,\omega),
\quad
 \EEzero(\nbigh,\omega):=
 \EE_{\epsilon_1}\bigl(H(\nbigelambda),\omega\bigr).
\end{equation}
See (\ref{eq;a12.9.1}) for the notation $\EE_{\epsilon_1}$.
The subset $\nbigs(\omega)\subset \KMS(\prolong{\nbige^0})$ is given
as follows:
\[
 \nbigs(\omega):=
 \eigenmap^f(\lambda_0)^{-1}(\omega)
=\Bigl\{
  u\in \KMSoverline(\nbige^0)\,\,
 \Big|\,\,
 \eigenmap^f(\lambda_0,u)=\omega
 \Bigr\}.
\]

We may assume that any point
$\lambda\in\Delta^{\ast}(\lambda_0,\epsilon_0)$ are generic,
due to Lemma \ref{lem;10.11.21}.
Then we have the following decomposition on
the punctured disc
$\Delta^{\ast}(\lambda_0,\epsilon_0)$:
\begin{equation} \label{eq;9.7.30}
 \nbigh_{|\Delta^{\ast}(\lambda_0,\epsilon_0)}
=\bigoplus_{u\in \KMS(\prolong{\nbige^0})}
 \nbigh_u,
\quad
 \nbigh_{u\,|\,\lambda}:=
 \EE\bigl(H(\nbigelambda),\eigenmap^{f}(\lambda,u)\bigr).
\end{equation}
\begin{lem}
We have the following decomposition:
\[
 \EEzero(\nbigh,\omega)_{|\Delta^{\ast}(\lambda_0,\epsilon_0)}
=
 \bigoplus_{u\in\nbigs(\omega)}\nbigh_u.
\]
\end{lem}
\pf
It immediately follows from the definition
of $\EEzero$ in (\ref{eq;9.7.31})
and the decomposition (\ref{eq;9.7.30}).
\hfill\qed

\subsubsection{The filtration $\nbigf^{(\lambda_0)}$}

We remark the following.
\begin{lem} \label{lem;10.11.26}
 The map $\paramap^f(\lambda_0):\nbigs(\omega)\lrarr\real$
 is injective.
\end{lem}
\pf
It follows from Lemma \ref{lem;10.11.25}.
\hfill\qed

\vspace{.1in}

On the vector bundle
$\EEzero(\nbigh,\omega)_{|\Delta^{\ast}(\lambda_0,\epsilon_0)}$,
we have the filtration $\nbigf^{(\lambda_0)}$
defined as follows:
\[
 \nbigf_d^{(\lambda_0)}
 \EEzero(\nbigh,\omega)_{|\Delta^{\ast}(\lambda_0,\epsilon_0)}
:=
 \bigoplus_{
 \substack{
 u\in\nbigs(\omega),\\
 \paramap^f(\lambda_0,u)\leq d
 }
 }\nbigh_{u}.
\]

\begin{lem}
The above filtration $\nbigf^{(\lambda_0)}$ of
$\nbigh_{\omega\,|\,\Delta^{\ast}(\lambda_0,\epsilon_0)} $
can be prolonged to the filtration
of $\nbigh_{\omega} $ over $\Delta(\lambda_0,\epsilon_0)$.
\end{lem}
\pf
Since $\nbigf^{(\lambda_0)}$ is defined by using
the generalized eigen-decomposition of holomorphic
endomorphisms $M^{\lambda}$, the claim holds.
\hfill\qed

\vspace{.1in}
We denote the prolonged filtration also by $\nbigf^{(\lambda_0)}$.

\begin{prop} \label{prop;9.7.80}
We have the following equality:
\begin{equation} \label{eq;9.7.40}
 \Bigl(
 \nbigf^{(\lambda_0)}_d\EEzero(\nbigh,\omega)
 \Bigr)_{|\lambda_0}
=\nbigf_d\Bigl(
 \EE(H(\nbige^{\lambda_0}),\omega)
 \Bigr).
\end{equation}
\end{prop}
\pf
We use the induction as is explained in the following.
For any sufficiently small $d$,
both of the sides in (\ref{eq;9.7.40}) are $0$.
Hence the equality (\ref{eq;9.7.40}) holds trivially.
If the equality (\ref{eq;9.7.40}) holds for $d$,
then (\ref{eq;9.7.40})
holds for $d+\eta$ for any sufficiently small $\eta>0$.
Hence we have only to show the following claim:
\begin{quote}
 (C):
 Assume that (\ref{eq;9.7.40}) holds for any $d<d_0$.
 Then the equality (\ref{eq;9.7.40}) holds for $d_0$.
\end{quote}

We assume that $\paramap^f(\lambda_0,u_0)=d_0$
for $u_0\in \nbigs(\omega)$.
Note that the element $u_0$ is uniquely determined
due to Lemma \ref{lem;10.11.26}.
We put $R:=\rank \nbigf^{(\lambda_0)}_d$.
We have the natural isomorphism:
\begin{equation} \label{eq;9.7.41}
 \bigwedge^R \nbigh(\nbige)\simeq
 \nbigh\bigl(\bigwedge^R\nbige\bigr).
\end{equation}
We do not distinguish them in the following argument.
We put as follows:
\[
 u_1=\sum_{\substack{ u\in\nbigs(\omega)\\
 \paramap^f(\lambda_0,u)\leq d_0
 }}
 \multiplicity(u,0)\cdot u.
\]
We put $d_1:=\paramap^{f}(\lambda_0,u_1)$.

\begin{lem} \label{lem;9.7.52}
Let $a:\nbigs(\omega)\lrarr\seisuu_{\geq\,0}$
be a map satisfying the following conditions:
\[
 a(u)\leq \multiplicity(u,0),
\quad\quad
 \sum a(u)=R.
\]
Then we have the following inequality:
\begin{equation} \label{eq;9.7.51}
 \sum_{u\in\nbigs(\omega)}a(u)\cdot \paramap^f(\lambda_0,u)
\geq d_1.
\end{equation}
The equality in {\rm (\ref{eq;9.7.51})} holds
if and only if $\{a(u)\,|\,u\in\nbigs(\omega)\}$ satisfies the
 following:
\[
 a(u)=\left\{
 \begin{array}{ll}
 \multiplicity(u,0) &
 \bigl(\paramap^f(\lambda_0,u)\leq \paramap^f(\lambda_0,u_0)\bigr),
 \\
 \mbox{{}}\\
 0   & \bigl(\paramap^f(\lambda_0,u)> \paramap^f(\lambda_0,u_0)\bigr).
 \end{array}
 \right.
\]
\end{lem}
\pf
It immediately follows from our choice of $d_1$.
\hfill\qed

\begin{lem} \label{lem;9.7.60}
We have the following equality:
\begin{equation}\label{eq;9.7.50}
 \nbigf_{d_1}
\cap
\bigwedge^R
 \EE\bigl(
 H(\nbige^{\lambda_0}),\omega \bigr)
=\bigwedge^R
 \bigl(
 \nbigf_{d_0}
 \EE(
 H(\nbige^{\lambda_0}),{\omega})
 \bigr),
\quad\quad\quad
 \nbigf_{<d_1}
\cap
\bigwedge^R
 \EE\bigl(
 H(\nbige^{\lambda_0}),\omega \bigr)
=0.
\end{equation}
In particular,
the rank of 
$ \nbigf_{d_1}
\cap
\bigwedge^R
 \EE\bigl(
 H(\nbige^{\lambda_0}),\omega \bigr)$
is one.
\end{lem}
\pf
We put $R':=\rank \EE(H(\nbigelambdazero),\omega)$.
Let $\vecs=(s_1,\ldots,s_{R'})$ be a frame of
$\EE(H(\nbigelambdazero),\omega)$
which is compatible with $\nbigf$ and $\EE$.
We may assume that $(s_1,\ldots,s_R)$ be a frame
of $\nbigf_{d_0}\EE\big(H(\nbigelambdazero),\omega\big)$.
For any subset $I\subset\{1,\ldots,R'\}$,
we put $s_I:=\bigwedge_{i\in I}s_i$.
Due to the norm estimate for the multi-valued sections
(Proposition \ref{prop;9.7.55}),
the tuple
$\big\{s_I\,\big|\,|I|=R\big\}$ gives the frame of $\bigwedge^R
\EE(H(\nbigelambdazero),\omega)$,
which is also compatible with the filtration $\nbigf$,
and we have the inequality
$-\ord(s_I)=-\sum_{i\in I}\ord(s_i)$.
Thus we obtain $-\ord(s_I)\geq d_1$,
and the equality holds if and only if $I$ is same as the set
$\{1,\ldots,R\}$ due to Lemma \ref{lem;9.7.52}.
Then we obtain the equalities (\ref{eq;9.7.50}).
\hfill\qed

\vspace{.1in}

We have the line subbundle
$\nbigl:=\bigwedge^R
 \Bigl(
 \nbigf^{(\lambda_0)}_{d_0}
 \EEzero(
 \nbigh(\nbige),{\omega})
 \Bigr)$
of $\bigwedge^R\nbigh(\nbige)$.
Let us pick a non-trivial section $s$
of $\nbigl$
such that $s_{|\lambda_0}\neq 0$.

\begin{lem} \label{lem;9.7.75}
We have $-\ord(s_{|\lambda_0})\leq d_1$.
\end{lem}
\pf
For any $\lambda\in\Delta(\lambda_0,\epsilon_0)$,
the element $s_{|\lambda}\in H(\nbigelambda)$ is an eigenvector of
$M^{\lambda}$,
and the eigenvalue is $\eigenmap^{f}(\lambda,u_1)$.
We put as follows:
\[
 v:=
 \exp\Bigl(
 \log z\cdot \eigenmap^f(\lambda,u)
 \Bigr)\cdot s.
\]
Then it gives a holomorphic section of
$\nbige$ defined
over $\Delta(\lambda_0,\epsilon_0)\times (X-D)$.

We may assume that any $\lambda\in\Delta^{\ast}(\lambda_0,\epsilon_0)$
is generic.
Then we obtain the following equality
for any point $\lambda\in\Delta^{\ast}(\lambda_0,\epsilon_0)$:
\[
 -\ord\bigl(v_{|\{\lambda\}\times (X-D)}\bigr)
=\paramap(\lambda,u_1).
\]
We also note that $|v|$
is bounded over the compact set
$\Delta(\lambda_0,\epsilon_0)\times\{|z|=1/2\}$.
Since $(\nbigelambda,h)$ is acceptable for any $\lambda$
(Theorem 1 in \cite{s2}, or Corollary \ref{cor;10.11.80} in this paper),
there exists a positive constant $M>0$
satisfying the following,
due to Corollary \ref{cor;11.28.15}:
\begin{itemize}
\item
 For any $\epsilon>0$, there exists $C_{\epsilon}>0$
 such that the inequality
 $\big|v_{|\lambda}\big|_h
\leq
 C_{\epsilon}\cdot |z|^{-\paramap(\lambda,u_1)}
 \cdot\big(-\log|z|\big)^{M}$
 holds
 for any $\lambda\in\Delta^{\ast}(\lambda_0,\epsilon_0)$.
\end{itemize}
Then we obtain the inequality for $\lambda=\lambda_0$:
\[
 \big|v_{|\lambda_0}\big|_h
\leq
 C_{\epsilon}\cdot |z|^{-\paramap(\lambda_0,u_1)}
 \cdot\big(-\log|z|\big)^{M}
\]
Hence we obtain the inequality
$-\ord(v_{|\lambda_0})\leq \paramap(\lambda_0,u_1)$,
and thus
$-\ord(s_{|\lambda_0})\leq \paramap^f(\lambda_0,u_1)=d_1$.
\hfill\qed

\vspace{.1in}

Let us return to the proof of Proposition \ref{prop;9.7.80}.
Due to Lemma \ref{lem;9.7.60} and \ref{lem;9.7.75},
we obtain the following equality:
\[
 \bigwedge^R\Bigl(
 \nbigf^{(\lambda_0)}_{d_0}
 \EEzero\big(\nbigh(\nbigelambda),\omega\big)
 \Bigr)_{|\lambda_0}
=\bigwedge^R\Bigl(
 \nbigf_{d_0}
 \EE\big(H(\nbigelambdazero),\omega\big)
 \Bigr).
\]
It implies the equality (\ref{eq;9.7.40}) for $d_0$.
Thus the proof of Proposition \ref{prop;9.7.80} is accomplished.
\hfill\qed

\subsubsection{The filtration $\nbigf^{(\lambda_0)}$
and the decomposition $\EEzero$ on $\nbige$}

The filtrations and the decomposition for $\nbigh$ on
$\Delta(\lambda_0,\epsilon_0)$ induce those for $\nbige$
on $\Delta(\lambda_0,\epsilon_0)\times(X-D)$.
We only summarize the result.

We have the following
decomposition of the family of the $\lambda$-connections
on $\Delta(\lambda_0,\epsilon_0)\times (X-D)$:
\[
 \nbige=\bigoplus_{\omega\in S_0}
 \EE^{(\lambda_0)}(\nbige,\omega),
 \quad
 \EE^{(\lambda_0)}(\nbige,\omega)
:=\EE_{\epsilon_2}(\nbige,\omega).
\]
Moreover we have the following decomposition on
$\Delta^{\ast}(\lambda_0,\epsilon_0)\times (X-D)$:
\[
 \EE^{(\lambda_0)}(\nbige,\omega)_{|\Delta^{\ast}(\lambda_0,\epsilon_0)}
=\bigoplus_{u\in \nbigs(\omega)}
 \EE(\nbige,\eigenmap^f(\lambda,u)).
\]

On the vector bundle $\EE^{(\lambda_0)}(\nbige,\omega)$,
we have  the filtration
$\nbigf^{(\lambda_0)}$ satisfying the following conditions:
\begin{itemize}
\item
 On $\Delta^{\ast}(\lambda_0,\epsilon_0)\times (X-D)$,
 we have the following splitting:
\[
 \nbigf_d^{(\lambda_0)}\bigl( \EE^{(\lambda_0)}(\nbige,\omega)
 \bigr)_
 {|\Delta^{\ast}(\lambda_0,\epsilon_0)\times (X-D) }
=\bigoplus_{
 \substack{u\in \nbigs(\omega)\\
 \paramap^f(\lambda_0,u)\leq d
 }
 }\EE\bigl(\nbige,\eigenmap^f(\lambda,u)\bigr).
\]
\item
On $\{\lambda_0\}\times(X-D)$,
we have the following:
\[
 \nbigf_d^{(\lambda_0)}\bigl(
 \EE^{(\lambda_0)}(\nbige,\omega) \bigr)_
 {|\{\lambda_0\}\times (X-D) }
=\nbigf_d( \EE(\nbige^{\lambda_0},\omega)).
\]
\end{itemize}

%% file: 5.tex

\subsubsection{Asymptotic orthogonality
 for $\EE$-decomposition of $\nbige^0$}
\label{subsubsection;a12.1.10}

Let $\harmonicbundle$ be a tame harmonic bundle
over $\Delta^{\ast}$.
Let $\epsilon_1$ be a sufficiently small number,
and $C$ be a sufficiently small number.
Then we have the following decomposition over $\Delta^{\ast}(C)$:
\begin{equation}
 \prolong{\nbige}^0
=\bigoplus_{\alpha\in\Sp(\Res(\theta))}
 \EE_{\epsilon_1}(\prolong{\nbige}^0,\alpha).
\end{equation}
See (\ref{eq;a12.9.1}) for the notation $\EE_{\epsilon_1}$.
\begin{lem}\label{lem;a11.14.5}
There exists a positive constant $\epsilon_2$ such that
the decomposition is $|z|^{\epsilon_2}$-asymptotically orthogonal.
\end{lem}
\pf
It is shown in Proposition \ref{prop;10.11.30}.
\hfill\qed

\subsubsection{An asymptotically orthogonal decomposition
 of $\nbige^{\lambda}$}

Let $\lambda\in \cnum_{\lambda}^{\ast}$ be generic.
Let $\vecv$ be a normalizing frame of $\prolong{\nbige}$.
We have the decomposition:
\begin{equation} \label{eq;8.14.5}
 \prolong{\nbigelambda}
=\bigoplus_{-1<b\leq 0}
 \Bigl(
 \bigoplus_{
  \substack{
  u\in\KMS(\nbige^0)\\
  \kappa(\paramap(\lambda,u))=b
  }  }
 \prolong{\EE(\nbigelambda,\eigenmap^f(\lambda,u))}
 \Bigr).
\end{equation}
We would like to show
that the decomposition above is $(-\log|z|)^{-1}$-asymptotically orthogonal.

The flat connection $\tilde{\nabla}^{\lambda}$
is given by $\DD^{\lambda,f}-\lambda^{-1}\cdot\psi$,
where $\psi=\psi_0dz/z$ is defined as in the formula (\ref{eq;9.6.30}).
On the other hand,
we have the unitary connection
$\nabla^{\lambda}
=\delbar_E+\lambda\theta^{\dagger}+\del_E-\bar{\lambda}\theta$.

We put $\Phi^{\lambda}
:=\tilde{\nabla}^{\lambda}-\nabla^{\lambda}$.
Then we have the following formula:
\begin{equation} \label{eq;9.6.33}
  \Phi^{\lambda}=
 \Bigl(
 \delbar_E+\lambda\theta^{\dagger}+\del_E+\lambda^{-1}\theta
-\lambda^{-1}\psi\Bigr)
-\Bigl(
  \delbar_E+\lambda\theta^{\dagger}+\del_E-\bar{\lambda}\theta
 \Bigr) 
=\lambda^{-1}
 \Bigl(
 (1+|\lambda|^2)\cdot\theta-\psi
 \Bigr).
\end{equation}

Let $\eta$ be a positive number.
Let $M_{\eta}(\tilde{\nabla}^{\lambda})$
denote the monodromy of $\tilde{\nabla}^{\lambda}$
along the circle $|z|=\eta$.

\begin{lem} \label{lem;9.6.35}
There exists $C>0$ satisfying the following:
\begin{equation} \label{eq;9.6.31}
 \Bigl|
 \bigl(M(\tilde{\nabla}^{\lambda})u,
       M(\tilde{\nabla}^{\lambda})v
 \bigr)_h
-(u,v)_h \Bigr|
\leq
 C\cdot(-\log|z|)^{-1}\cdot|u|\cdot|v|.
\end{equation}
\end{lem}
\pf
We use the real coordinate $(x,y)$ $(0\leq x<2\pi,\,\,0<y)$
given by $z=\exp(\sqrt{-1}x-y)$.
Let $V$ denote the vector field given by ${\del}/{\del x}$.
Let $s$ be a flat section with respect to $\tilde{\nabla}^{\lambda}$
over $|z|=\eta$.
It satisfies the following:
\[
 \frac{d}{dx}h(s,s)
=2\Re \bigl(
 h(\nabla^{\lambda}_V(s),s)\bigr)
=-2\Re \bigl(
 h(\Phi^{\lambda}_V(s),s)\bigr).
\]
Hence we obtain the following:
\[
 \frac{d}{dx}\log|s|_h^2
=-2\frac{\Re \bigl( h(\Phi^{\lambda}_V(s),s)\bigr)}{|s|_h^2}.
\]
Due to the estimate (\ref{eq;9.6.32}),
we have
$|\Phi^{\lambda}_V|_h\leq C_0\cdot y_0$
for some $C_0>0$.
Thus there exists a positive constant $C_1>0$,
which is independent of $x$ and $y$,
and satisfying the following:
\[
 \left|
 \frac{d}{dx}\log|s(x,y)|_h^2
 \right|
\leq C_1\cdot y^{-1}.
\]
Hence there exists $C_2>0$ satisfying the following
for any $0\leq x\leq 2\pi$ and for any $y>0$:
\begin{equation} \label{eq;9.6.34}
 \big|s(x,y)\big|_h\leq 
 C_2\cdot\bigl|s(0,y)\bigr|_h
\end{equation}

On the other hand,
we also have the following equality:
\[
 \frac{d}{dx}(s_1,s_2)_h=
-(\Phi_V^{\lambda}s_1,s_2)_h
-(s_1,\Phi_V^{\lambda}s_2)_h.
\]
Then we obtain the following inequality
from (\ref{eq;9.6.32}), (\ref{eq;9.6.33}) and (\ref{eq;9.6.34}):
\[
 \left|
 \frac{d}{dx}(s_1,s_2)_h
 \right|
\leq
 C_3\cdot y^{-1}\cdot|s_1|_h\cdot |s_2|_h
\leq
 C_4\cdot y^{-1}\cdot|s_1(0)|_h\cdot|s_2(0)|_h.
\]
The inequality (\ref{eq;9.6.31}) immediately follows.
\hfill\qed

\vspace{.1in}

We have the $\EE$-decomposition of the flat bundles for
the monodromy $M^{\lambda}$:
\begin{equation}\label{eq;8.14.6}
 \nbigelambda
=\bigoplus_{-1<b\leq 0}
 \Bigl(
  \bigoplus_{
  \substack{
  u\in\KMS(\nbige^0)\\
  \paramap(\lambda,u)=b
 } }\EE(\nbigelambda,\eigenmap^f(\lambda,u))
 \Bigr).
\end{equation}

\begin{lem}
The decomposition {\rm (\ref{eq;8.14.6})} is
the generalized eigen decomposition
with respect to the monodromy $M_{\eta}(\tilde{\nabla}^{\lambda})$.
Namely we have the following equality:
\[
 \bigoplus_{
  \substack{
  u\in\KMS(\nbige^0)\\
  \paramap(\lambda,u)=b
 } } \EE\bigl(\nbigelambda,\eigenmap^f(\lambda,u)\bigr)
=\EE\Big(
  M_{\eta}(\tilde{\nabla}^{\lambda}),\,\,
  \exp(2\pi\sqrt{-1}b)
 \Big).
\]
\end{lem}
\pf
Since the endomorphism $\psi_0$ is flat with respect to
$\DD^{\lambda,f}$,
we have the following:
\[
 M(\tilde{\nabla}^{\lambda})
=M(\DD^{\lambda,f})\circ 
 \exp\bigl(2\pi\lambda^{-1}\sqrt{-1}\psi_0\bigr).
\]
Hence the eigenvalue of $M(\tilde{\nabla}^{\lambda})$
corresponding to $v_i$ is as follows:
\[
 \exp\bigl(-2\pi\sqrt{-1}\lambda^{-1}\cdot\beta_i\bigr)
\cdot \exp\bigl(2\pi\sqrt{-1}(\lambda^{-1}\cdot\beta_i+b_i)\bigr)
=\exp\bigl(2\pi\sqrt{-1}b_i\bigr).
\]
It implies the lemma.
\hfill\qed

\begin{lem} \label{lem;9.6.40}
There exists a positive number $C>0$
such that
the following holds for any $i$:
\[
 \Bigl|
 M_{\eta}(\tilde{\nabla}^{\lambda})v_i
-\exp(2\pi \sqrt{-1}b_i)v_i
 \Bigr|_h
\leq C \cdot |v_i|\cdot (-\log|z|)^{-1}.
\]
\end{lem}
\pf
We have the following equality:
\begin{multline}
 M_{\eta}(\tilde{\nabla}^{\lambda})v_i
-\exp(2\pi \sqrt{-1}b_i)\cdot v_i
=
 \exp\bigl(2\pi\sqrt{-1}(b_i+\lambda^{-1}\cdot\beta_i)\bigr)
\cdot
 \Bigl(
 M(\DD^{\lambda,f})\cdot v_i
-\exp(-2\pi\sqrt{-1}\lambda^{-1}\cdot\beta_i)\cdot v_i
 \Bigr) \\
=\exp(2\pi\sqrt{-1}b_i)
\cdot
\sum_{n=1}^{\infty}\frac{(-2\pi)^n}{n!}
 (N_{\omega_i}^{\lambda})^nv_i.
\end{multline}
Since we have
$|(N_{\omega_i}^{\lambda})^nv_i|
 \leq C\cdot |v_i|\cdot(-\log|z|)^{-1}$,
we obtain the result.
\hfill\qed

\begin{lem}
There exists a positive number $C$
satisfying the following:
\begin{quote}
Let $\gamma_i$ $(i=1,2)$ be elements of
$\Sp(M_{\eta}(\tilde{\nabla}^{\lambda}))$,
and $u_i$ be elements of
 $\EE(M_{\eta}(\tilde{\nabla}^{\lambda}),\gamma_i)$.
If $\gamma_1\neq \gamma_2$,
then the following holds:
\begin{equation} \label{eq;9.6.41}
 \bigl|
 (u_1,u_2)_h
 \bigr|
\leq 
 C\cdot (-\log|z|)^{-1}\cdot|u_1|_h\cdot |u_2|_h.
\end{equation}
\end{quote}
Namely the generalized eigen-decomposition of the monodromy of
$\tilde{\nabla}^{\lambda}$ is $(-\log|z|)^{-1}$-asymptotically orthogonal.
\end{lem}
\pf
Due to Lemma \ref{lem;9.6.35},
we have the following inequality:
\begin{equation} \label{eq;9.6.36}
 \Bigl|\bigl(
 M(\tilde{\nabla}^{\lambda})\cdot u_1,
 M(\tilde{\nabla}^{\lambda})\cdot u_2
 \bigr)_h
-(u_1,u_2)_h
 \Bigr|
\leq C\cdot (-\log|z|)^{-1}\cdot|u_1|_h\cdot|u_2|_h.
\end{equation}

On the other hand,
we have the following inequality due to Lemma \ref{lem;9.6.40}:
\begin{multline} \label{eq;9.6.42}
\Bigl|\bigl(
  M(\tilde{\nabla}^{\lambda})u,\,\,
 M(\tilde{\nabla}^{\lambda})v
 \bigr)_h
-\gamma_1\cdot\bar{\gamma}_2\cdot(u,v)_h
\Bigr|
=
 \Bigl|
 \bigl(
  (M(\tilde{\nabla}^{\lambda})-\gamma_1)u,\,\,
 M(\tilde{\nabla}^{\lambda})v
 \bigr)_h
+\bigl(\gamma_1\!\cdot\! u,\,\,
  (M(\tilde{\nabla}^{\lambda})-\gamma_2\bigr)v
 )
 \Bigr| \\
 \leq C'\cdot (-\log|z|)^{-1}\cdot|u|\cdot|v|.
\end{multline}
Since we have $|\gamma_i|=1$,
the condition $\gamma_1\neq\gamma_2$
implies
$\gamma_1\cdot\bar{\gamma}_2\neq 1$.
Hence we obtain the inequality (\ref{eq;9.6.41})
from (\ref{eq;9.6.36}) and (\ref{eq;9.6.42}).
\hfill\qed

\subsubsection{Asymptotic orthogonality for the parabolic filtration
 of $\nbige^0$}
\label{subsubsection;a12.1.11}

\begin{lem}\label{lem;a11.14.2}
We pick $\lambda\in\cnum_{\lambda}^{\ast}$ satisfying the following:
\begin{enumerate}
\item
$\lambda$ is generic (Definition {\rm\ref{df;a11.14.1}}).
\item
 $u\neq u'\Longrightarrow
 \paramap(\lambda,u)\neq \paramap(\lambda,u')$.
\end{enumerate}
\end{lem}
\pf
Note $|\lambda|$ is sufficiently small,
we may always assume the second condition holds.
\hfill\qed

\vspace{.1in}

Let us pick $\lambda$ as in Lemma \ref{lem;a11.14.2}.
Let $\vece$ be a frame of $\prolong{\nbige^0}$ compatible with
$\EE$, $F$ and $W$.
Let $\vecv$ be a frame of $\prolong{\nbigelambda}$ compatible with
$\EE$, $F$ and $W$.
We put $b(e_i)=\deg^{F}(e_i)$ and $b(v_i)=\deg^F(v_i)$.
We put $k(e_i)=\deg^W(e_i)$ and $k(v_i)=\deg^{W}(v_i)$.
We put $u(e_i)=\deg^{F,\EE}(e_i)$
and $u(v_i)=\deg^{F,\EE}(v_i)$.

We put $e_i':=e_i\cdot|z|^{b(e_i)}\cdot(-\log|z|)^{-k(e_i)/2}$
and $v_i':=v_i\cdot|z|^{b(v_i)}\cdot(-\log|z|)^{-k(v_i)/2}$.
We have the bounded $C^{\infty}$-function
$B:X-D\lrarr M(r)$ determined by the following:
\[
 e_i'=\sum_j B_{j\,i}'\cdot v_j'.
\]
Recall that
$ |B_{j\,i}'|\leq C\cdot (-\log|z|)^{-1}$
unless $\kmsmap\bigl(\lambda,u(e_i)\bigr)=u(v_j)$
(Corollary \ref{cor;a11.14.3}).

\begin{lem} \label{lem;9.6.50}
There exists a positive constant $C>0$ satisfying the following
condition:
\begin{quote}
In the case $u(e_i)\neq u(e_j)$,
the inequality $\big|(e_i,e_j)_h\big|
\leq C\cdot(-\log|z|)^{-1}\cdot|e_i|_h\cdot|e_j|_h$
holds for some positive constant $C$.
\end{quote}
\end{lem}
\pf
We have already obtained stronger estimate
in the case $\deg^{\EE}(e_i)\neq \deg^{\EE}(e_j)$
(Lemma \ref{lem;a11.14.5}).
We have only to prove the following estimate,
for some positive number $C_1$
in the case $b(e_i)\neq b(e_j)$:
\[
 \big|(e_i',e_j')_h\big|\leq C_1\cdot(-\log|z|)^{-1}.
\]
We have the following formula:
\begin{equation} \label{eq;8.14.8}
 (e_i',e_j')_h=
 \sum B_{k\,i}'\cdot \bar{B}_{l\,j}'\cdot
 (v_k',v_l')_h.
\end{equation}
We have the inequality
$ |B_{k\,i}'|\cdot |\bar{B}_{l\,j}'|
 <C\cdot(-\log|z|)^{-1}$,
unless the following condition holds:
\begin{equation} \label{eq;8.14.7}
\kmsmap\bigl(\lambda,u(e_i)\bigr)=u(v_k),
\quad
\kmsmap\bigl(\lambda,u(e_j)\bigr)=u(v_l).
\end{equation}

Assume that $b(e_i)\neq b(e_j)$ and the equalities (\ref{eq;8.14.7}) hold.
Then we have $b(v_k)\neq b(v_l)$ due to our choice of $\lambda$.
Thus we obtain the inequality $(v_k',v_l')_h\leq C\cdot(-\log|z|)^{-1}$
for some positive constant $C$.
Hence the right hand side in (\ref{eq;8.14.8}) is
dominated by $(-\log|z|)^{-1}$.
\hfill\qed

\vspace{.1in}

Let us take a decomposition of the vector bundle
$\EE_{\epsilon_2}(\prolong{\nbige^0},\alpha)$:
\begin{equation}\label{eq;8.14.10}
 \EE_{\epsilon_2}(\prolong{\nbige^0},\alpha)
=\bigoplus_{a\in\Par(\prolong{\nbige^0})}
 V_{(a,\alpha)}.
\end{equation}
We assume that the decomposition gives 
a splitting of the parabolic filtration $F$,
in the following sense:
\begin{equation}\label{eq;8.14.11}
 F_a(\prolong{\nbige^0}_{|O})
=\bigoplus_{\alpha}\bigoplus_{b\leq a}V_{{b,\alpha}\,|\,O}.
\end{equation}
For example,
we can pick $V_a$ as the vector subbundle of
$\EE_{\epsilon_2}(\prolong{\nbige},\alpha)$
generated by $\bigl\{e_i\,\big|\,u(e_i)=(a,\alpha)\bigr\}$.
On the other hand,
if we are given such decomposition,
then we can pick the frame $\vece$ compatible with
$\EE$, $F$ and $W$
such that $V_a$ is generated by
$\bigl\{e_i\,\big|\,u(e_i)=(a,\alpha)\bigr\}$.

\begin{prop}
If the condition {\rm (\ref{eq;8.14.11})} is satisfied,
then the decomposition {\rm (\ref{eq;8.14.10})} is
$(-\log|z|)^{-1}$-asymptotically orthogonal.
\end{prop}
\pf
It follows from Lemma \ref{lem;9.6.50}.
\hfill\qed

\subsubsection{Asymptotic orthogonality for the weight filtration}
\label{subsubsection;a12.1.12}

Let $(E_0,\delbar_{E_0},\theta_0,h_0)$
be a model bundle for $\harmonicbundle$.
We use the notation in the subsection \ref{subsection;10.11.35}.
Recall the finiteness due to Simpson (Lemma \ref{lem;9.15.2}):
\[
 \int |\theta^{\dagger}-\theta_0^{\dagger}|_h^2\cdot
 \bigl(-|z|\cdot \log|z|\bigr)\cdot
 \frac{|dz\cdot d\bar{z}|}{|z|^2(-\log|z|)}<\infty.
\]
Hence, for any $\epsilon>0$
there exists a subset $Z_{\epsilon}\subset X-D$
satisfying the following:
\begin{itemize}
\item
 The volume of $(X-D)-Z_{\epsilon}$ with respect
 to the measure $|dz\cdot d\bar{z}|\cdot |z|^{-2}(-\log|z|)^{-1}$
 is finite.
\item
 We have the estimate
 $|\theta^{\dagger}-\theta_0^{\dagger}|_h\cdot
 \bigl(-|z|\log|z|\bigr)\leq \epsilon$ on $Z_{\epsilon}$.
\end{itemize}

The endomorphisms $A_0$ and $A$ are determined by the following:
\[
 A_0\cdot \frac{dz\cdot d\bar{z}}{|z|^2(-\log|z|)^2}
=\theta_0\cdot\theta_0^{\dagger}
+\theta_0^{\dagger}\cdot\theta_0,
\quad\quad
 A\cdot \frac{dz\cdot d\bar{z}}{|z|^2(-\log|z|)^2}
=\theta\cdot\theta^{\dagger}
+\theta^{\dagger}\cdot\theta.
\]
Then $A_0$ is self dual with respect to $h_0$,
and $A$ is self dual with respect to $h$.
By a direct calculation,
it can be checked that the eigenvalues of $A_0$ are integers.

We have the decomposition
$E_0=\bigoplus_{k\in\seisuu} U_k$
satisfying the following:
\begin{itemize}
\item
 $W_k=\bigoplus_{h\leq k}U_h$.
\item
  $U_k$ is eigen space of $A_0$ corresponding
 to the integer $k$.
\end{itemize}

The following lemma is clear.
\begin{lem}\mbox{{}}
\begin{itemize}
\item
There exists $C_4>0$ such that
$|\theta-\theta_0|_h\cdot(-|z|\log|z|)\leq C_4\cdot|z|^{1/2}$.
\item
There exists $C_5>0$ the inequality
$|A_0-A|_h\leq C_5\cdot\epsilon_1$
holds on $Z_{\epsilon_1}\cap \Delta^{\ast}(\epsilon_2)$.
\hfill\qed
\end{itemize}
\end{lem}

Let us take sufficiently small positive number $\epsilon_1$
such that $5\cdot\epsilon_1<|\phi(k)-\phi(k-1)|$.
Let $v_i$ be a $C^{\infty}$-section of $U_{k_i}$
for $k_1\neq k_2$.
Since $A$ is anti-self dual with respect to $h$,
we have $(Av_1,v_2)_h+(v_1,Av_2)_h=0$.

On the other hand,
we have the following:
\[
 |A\cdot v_i-k_i\cdot v_i|_h=|A\cdot v_i-A_0\cdot v_i|_h
 \leq C_5\cdot \epsilon\cdot|v|_h.
\]
Then we obtain the following:
\begin{multline}
 \bigl|
 k_1-k_2
 \bigr|\cdot
 \Bigl|
 (v_1,v_2)_h
 \Bigr|
=
 \Bigl|
 \big(k_1\cdot v_1,v_2\big)_h-\big(v_1,k_2\cdot v_2\big)_h
 \Bigr|\\
=
 \Bigl|
 \big((k_1-A)v_1,v_2\big)_h-\big(v_1,(k_2-A)v_2\big)_h
 \Bigr|
\leq 2C_5\cdot \epsilon\cdot|v_1|_h\cdot|v_2|_h.
\end{multline}
Here we have used the self-duality of $A$
with respect to $h$.
Hence we obtain the inequality
$\big|(v_1,v_2)_h\big|\leq C_6\cdot\epsilon\cdot|v_1|_h\cdot|v_2|_h$.
In all, we obtain the following.
\begin{prop}
 For any $\epsilon>0$,
 there exists a subset $Z_{\epsilon}\subset X-D$
satisfying the following conditions:
\begin{itemize}
\item
 The measure of $(X-D)-Z_{\epsilon}$ with respect to
 $|dz\cdot d\bar{z}|\cdot |z|^{-2}(\log|z|)^{-1}$
 is finite.
\item
 Let $v,w\in \prolong{\nbige}^0$ such that
 $\deg^{F,\EE}(v)=\deg^{F,\EE}(w)$
 and $\deg^{W}(v)\neq \deg^{W}(w)$.
 On $Z_{\epsilon}$, we have the estimate
 $|(v,w)|\leq \epsilon\cdot|v|_h\cdot|w|_h$.
\hfill\qed
\end{itemize}
\end{prop}

%% file: d2.tex

Let $(E^{(i)},\delbar_{E^{(i)}},h^{(i)},\theta^{(i)})$
$(i=1,2)$
be tame harmonic bundles on $\Deltabarast$ of rank $r$.
We denote the corresponding $\lambda$-connections
by $\bigl(\nbige^{(i)\,\lambda},\DD^{\lambda}\bigr)$.
Assume the following:
\begin{itemize}
\item
The set of $KMS$-spectrum are same,
i.e.,
$\KMS(\nbige^{(1)\,0})=\KMS(\nbige^{(2)\,0})$.
\item
Let $\lambda$ be a complex number
which is generic with respect to $\KMS(\nbige^{(i)\,0})$,
and we are given a flat isomorphism
$\Phi:\bigl(\nbige^{(1)\,\lambda},\DD^{\lambda}\bigr)
\lrarr \bigl(\nbige^{(2)\,\lambda},\DD^{\lambda}\bigr)$.
\end{itemize}

\begin{lem}\label{lem;04.1.27.9}
The functions $\log|\Phi|^2$ and $\log|\Phi^{-1}|^2$ are subharmonic
on the disc $\Delta$.
\end{lem}
\pf
First, let us see the boundedness of $\Phi$ and $\Phi^{-1}$.
Since $\lambda$ is generic,
the harmonic metrics $h_i$ of
$\bigl(\nbige^{(i)\,\lambda},\DD^{\lambda}\bigr)$
are determined by the monodromy, up to boundedness.
Use the result in the subsubsection 4.3.3 in \cite{mochi},
and note that the splitting of the parabolic structure
is given by the generalized eigen-decomposition of the residues.

We have the inequality
$\Delta\log|\Phi|^2\leq 0$ and $\Delta\log|\Phi^{-1}|^2\leq 0$
on the punctured disc $\Delta^{\ast}$,
due to the Weitzenbeck formula of Simpson
(Lemma 4.1 in \cite{s2}).
Since $\log|\Phi|^2$ and $\log|\Phi^{-1}|^2$ are bounded,
the inequalities hold on the disc $\Delta$.
(See Lemma 2.2 in \cite{s2}.)
Thus we are done.
\hfill\qed

\vspace{.1in}

We identify $\bigl(\nbige^{(1)},\DD^{\lambda}\bigr)$
and $\bigl(\nbige^{(2)},\DD^{\lambda}\bigr)$
via the morphism $\Phi$.
For any point $P\in \Deltabarast$,
we can consider the distance of
$h_{1\,|\,P}$ and $h_{2\,|\,P}$ in 
the space $\PH(r)$
(the subsubsection \ref{subsubsection;04.1.27.5}).

\begin{lem}\label{lem;04.1.27.10}
Let $R$ be a real number such that
the inequality 
$d_{\PH(r)}\bigl(h_{1\,|\,Q},h_{2\,|\,Q}\bigr)\leq R$
hold for any point $Q\in \del\Deltabar$.
Then the following inequalities hold for any point $P\in \Deltabarast$:
\begin{equation}\label{eq;04.1.27.71}
 d_{\PH(r)}\bigl(h_{1\,|\,P},h_{2\,|\,P}\bigr)
\leq 
 \left(\frac{e^{R}-e^{-R}}{2R}
 \right)\cdot
 \max\Bigl\{
 d_{\PH(r)}\bigl(h_{1\,|\,Q},h_{2\,|\,Q}\bigr)
\,\Big|\,Q\in\del\Deltabar
 \Bigr\}.
\end{equation}
\end{lem}
\pf
We always have the following, due to Lemma \ref{lem;04.1.27.6}:
\begin{equation}\label{eq;04.1.27.7}
 d_{\PH(r)}(h_{1\,|\,P},h_{2\,|\,P})
\leq
 \frac{|\Phi_{|P}|^2+|\Phi^{-1}_{|P}|^2-2r}{2}.
\end{equation}
For any point $Q\in \del X$,
we have the following, due to Lemma \ref{lem;04.1.27.6} again:
\begin{equation}\label{eq;04.1.27.8}
 \frac{|\Phi_{|Q}|^2+|\Phi^{-1}_{|Q}|^2-2r}{2}
\leq
 \frac{e^R-e^{-R}}{2R}
\cdot d_{\PH(r)}\bigl(h_{1\,|\,Q},h_{2\,|\,Q}\bigr).
\end{equation}
Then (\ref{eq;04.1.27.71})
follows from the inequalities (\ref{eq;04.1.27.7}), (\ref{eq;04.1.27.8})
and Lemma \ref{lem;04.1.27.9}.
\hfill\qed

%% file: 10.tex

\subsubsection{$KMS$-Spectrum}

We put $X:=\Delta^n$, $D_i:=\{z_i=0\}$ and $D=\bigcup_{i=1}^lD_i$.
Let $\pi_i$ denote the projection $X\lrarr D_i$,
forgetting the $i$-th component.
Let $\harmonicbundle$ be a tame harmonic bundle over $X-D$.
We have the deformed holomorphic bundle
$(\nbigelambda,\DDlambda)$.
Let $i$ be an element of $\lbar$,
and $P$ be a point of
$D_i^{\circ}:=D_i-\bigcup_{j\neq i,\,\,1\leq j\leq l}D_i\cap D_j$.
We put as follows:
\[
 \KMS(\nbigelambda,i,P):=\KMS\bigl(\nbigelambda_{|\pi_i^{-1}(P)}\bigr).
\]
\begin{lem} \label{lem;9.7.101}
The set $\KMS(\nbigelambda,i,P)$ and the multiplicities
of any element $u\in \KMS(\nbigelambda,i,P)$
are independent of $P\in D_i^{\circ}$.
\end{lem}
\pf
Let $P_i$ $(i=1,2)$ be points of $D_i^{\circ}$.
For any $\lambda\in\cnum_{\lambda}^{\ast}$,
we have $\Sp(\nbigelambda,i,P_1)=\Sp(\nbigelambda,i,P_2)$,
and the multiplicities are same.
Then we obtain the result due to Lemma \ref{lem;9.7.100}.
\hfill\qed

\vspace{.1in}

In the following,
we use the notation $\KMS(\nbigelambda,i)$
instead of $\KMS(\nbigelambda,i,P)$.
Similarly,
the sets $\KMSoverline(\nbigelambda,i)$,
$\Sp(\nbigelambda,i)$
and $\Par(\nbigelambda,i)$ are obtained.

For any element $\vecb\in\real^l$,
the sets
$\KMS(\prolongg{\vecb}{\nbigelambda},i):=
 \KMS\big(\prolongg{b_i}{\bigl(\nbigelambda_{|\pi_i^{-1}(P)}\bigr)}\big)$
are also obtained,
where $b_i$ denotes the $i$-th component of $\vecb$.


Let us consider the Higgs field:
\[
 \theta
=\sum_{i=1}^l f_i\cdot\frac{dz_i}{z_i}
+\sum_{j=l+1}^n g_j\cdot dz_j.
\]
We have the characteristic polynomials
$\det(t-f_i)$ and $\det(t-g_j)$,
whose coefficients are holomorphic functions defined over $X-D$.
Since our harmonic bundle is tame,
the coefficients are prolonged to the holomorphic functions
defined over $X$,
by definition.
We denote them by the same notation.

\begin{lem}
We have $\det(t-f_i)_{|P_1}=\det(t-f_i)_{|P_2}$
if $P_a\in D_i$ $(a=1,2)$.
\end{lem}
\pf
If $P\in D_i^{\circ}$,
then we have
$\det(t-f_i)_{|P}=\det\bigl(t-f_{i\,|\,\pi_i^{-1}(P)}\bigr)_{|P}$.
Due to Lemma \ref{lem;9.7.101},
the right hand side is independent of
a choice of a point $P\in D_i^{\circ}$.
Then we obtain the result for $D_i$.
\hfill\qed


%% file: 10.3.tex

\subsubsection{Rank 1}

Let $\harmonicbundle$ be a tame harmonic bundle
of rank 1 over $X-D$.
In this case,
the set $\KMS(\nbige^0,i)$ consists of only one element $u_i$
for $i=1,\ldots,l$.
We have the model bundle $L(\vecu)$
for $\vecu=(u_1,\ldots,u_l)$.
Then $\harmonicbundle\otimes L(\vecu)$
is tame and nilpotent.
Moreover the parabolic structure is trivial.
Hence it is the restriction of the harmonic bundle
$(E_0,\delbar_{E_0},\theta_0,h_0)$ of rank 1 over $X$,
due to Corollary 4.10 in \cite{mochi}.
Namely we obtain the following.
\begin{prop}
Let 
$\harmonicbundle$ be a tame harmonic bundle
over $X-D$ of rank $1$.
Then it is isomorphic to the tensor product of
a model bundle
$L(\vecu)$
and a harmonic bundle $(E_0,\delbar_{E_0},h_0,\theta_0)$
over $X$.
\hfill\qed
\end{prop}

%% file: d3.tex

\subsubsection{A characterization of tameness}
\label{subsubsection;04.1.27.65}

\begin{lem}\label{lem;04.1.27.30}
Let $X$ be a simply connected compact region of $\cnum^n$.
Let $\pi:X\times\Deltabarast\lrarr X$ denote the projection.
Let $\harmonicbundle$ be a harmonic bundle
on $X\times\Deltabarast$.
Assume the following:
\begin{itemize}
\item
 For any point $P\in X$,
 the restriction $\harmonicbundle_{|\pi^{-1}(P)}$ is tame
 with respect to $(P,O)$.
\end{itemize}
Then the harmonic bundle $\harmonicbundle$ is tame
with respect to the divisor $X\times\{O\}$.
\end{lem}
\pf
Let $(z_1,\ldots,z_n)$ be the coordinate of $X$,
and $w$ be the coordinate of $\Deltabar$.
We have the description
$\theta=\sum_{i=1}^n g_i\cdot dz_i+f\cdot dw/w$.
We would like to show that
the coefficients of 
$\det(t-f)$ and $\det(t-g_i)$ are holomorphic on $X\times\Deltabar$.
Due to our assumption,
the coefficients of $\det(t-f)_{|\pi^{-1}(P)}$ is holomorphic 
for any $P\in X$.
Then it is easy to derive that 
the coefficients of $\det(t-f)$ is holomorphic on $X\times\Deltabar$.
So we have only to see $\det(t-g_i)$.

We remark that
$\KMS(\nbige^0_{|\pi^{-1}(P)})$ is independent of a choice
of $P\in X$, which can be shown by the same argument
as the proof of Lemma \ref{lem;9.7.101}.
Let us take a complex number $\lambda$,
which is generic with respect to
$\KMS(\nbige^0_{\pi^{-1}(P)})$.
We have the flat bundle
$\bigl(\nbigelambda,\DD^{\lambda,f}\bigr)$
on $X\times\Deltabarast$.
In particular, we have the following flat isomorphism
for any points $P_1,P_2\in X$:
\begin{equation} \label{eq;04.1.27.15}
\bigl(\nbigelambda,\DD^{\lambda,f}\bigr)_{|\pi^{-1}(P_1)}
\simeq \bigl(\nbigelambda,\DD^{\lambda,f}\bigr)_{\pi^{-1}(P_2)}.
\end{equation}

We put $h_{P}:=h_{|\pi^{-1}(P)}$.
For any point $Q\in \Deltabarast$,
we can consider the distance
$d_{\PH(r)}\bigl(h_{P_1\,|\,Q},h_{P_2\,|\,Q}\bigr)$
for any points $P_i\in X$,
via the isomorphism (\ref{eq;04.1.27.15}).

Let $P_0$ be a point of $X$,
and $v$ be a real tangent vector of $X$ at $P_0$.
It naturally gives the tangent vector $v_Q$
of $X\times\Deltabarast$ at $(P,Q)$ for any point $Q\in \Deltabarast$.
From Lemma \ref{lem;04.1.27.10}
and the formula (\ref{eq;04.1.27.20}),
we obtain the following maximum principle:
\begin{equation}\label{eq;04.1.27.21}
 \bigl|
\overline{\lambda}\cdot\theta(v_Q)
+\lambda\cdot\theta^{\dagger}(v_Q)
 \bigr|^2
\leq
 \max\Bigl\{
  \bigl|
\overline{\lambda}\cdot\theta(v_{Q'})
+\lambda\cdot\theta^{\dagger}(v_{Q'})
 \bigr|^2\,
 \Big|\, Q'\in \del\Deltabar
 \Bigr\}.
\end{equation}
Note $(2R)^{-1}\cdot (e^R-e^{-R})$ goes to $1$
when $R$ goes to $0$.

Let us use the real coordinate $z_i=x_i+\sqrt{-1}y_i$.
Let $\del_{x_i}$ and $\del_{y_i}$ denote
the vector fields $\del/\del x_i$ and $\del/\del y_i$.
Since $X\times\del\Deltabar$ is compact,
we obtain the boundedness of the following,
from (\ref{eq;04.1.27.21}):
\[
 \overline{\lambda}\cdot\theta\bigl(\del_{x_i}\bigr)
+\lambda\cdot \theta^{\dagger}\bigl(\del_{x_i}\bigr)
=\overline{\lambda}\cdot g_i+\lambda\cdot g_i^{\dagger}.
\]
\[
  \overline{\lambda}\cdot\theta\bigl(\del_{y_i}\bigr)
+\lambda\cdot \theta^{\dagger}\bigl(\del_{y_i}\bigr)
=\overline{\lambda}\cdot g_i-\lambda\cdot g_i^{\dagger}.
\]
Hence we obtain the boundedness of $g_i$.

Then we obtain the boundedness of the coefficients of $\det(t-g_i)$,
which implies that
they are holomorphic on $X\times\Deltabar$.
\hfill\qed

\begin{lem}
We put $X=\Delta^n$, $D_i:=\{z_i=0\}$ and $D=\bigcup_{i=1}^l D_i$.
Let $\pi_i$ denote the naturally defined projection $X\lrarr D_i$.
Let $\harmonicbundle$ be a harmonic bundle on $X-D$.
Assume the following:
\begin{itemize}
\item
 For any point $Q\in D_i\setminus \bigcup_{j\neq i} D_j$,
 the restriction $\harmonicbundle_{|\pi^{-1}(Q)}$ is tame.
\end{itemize}
Then the harmonic bundle is tame with respect to the divisor $D$.
\end{lem}
\pf
We put $D^{[2]}:=\bigcup_{i\neq j}D_i\cap D_j$,
which is of codimension $2$ in $X$.
We have the description:
\[
 \theta=\sum_{i=1}^l f_i\cdot \frac{dz_i}{z_i}+\sum_{i=l+1}^n g_i\cdot dz_i.
\]
Due to Lemma \ref{lem;04.1.27.30},
the coefficients of $\det(t-f_i)$ and $\det(t-g_i)$ are holomorphic
on $X-D^{[2]}$.
Then we can conclude that they are holomorphic on $X$
due to the theorem of Hartogs.
\hfill\qed

\vspace{.1in}
We have a straightforward corollary.

\begin{cor} \label{cor;04.1.27.51}
Let $X$ be a complex manifold, and $D$ be a normal crossing divisor
of $X$.
Let $\harmonicbundle$ be a harmonic bundle on $X-D$.
Assume the following:
\begin{itemize}
\item
 For any smooth curve $C$ be contained in $X$,
 which intersects with the smooth parts of $D$ transversely,
 the restriction $\harmonicbundle_{|C}$ is tame.
\end{itemize}
Then the harmonic bundle $\harmonicbundle$ is tame.
\hfill\qed
\end{cor}

%% file: 10.1.tex

\subsubsection{Preliminary decomposition}
\label{subsubsection;a11.9.17}

Pick a positive number $\epsilon_1$.
If we replace $X$ by a sufficiently small neighbourhood
of the origin $O$,
and if we replace the coordinate appropriately,
then we may assume that
there exist positive constants $\epsilon$, $\epsilon_1$ and $C$,
satisfying the following conditions:
\begin{itemize}
\item
 We have the decomposition
 $E=\bigoplus_{a\in\Sp(\nbige^0,i)}\EE_{\epsilon_1}(f_i,a)$
 over $X-D$.
\item
 Then we have the decomposition
 $f_i=\bigoplus_{a\in\Sp(\nbige^0,i)}f_{i\,a}$.
 Then we have the inequality
 $|\alpha-a|<C\cdot|z|^{\epsilon}$
 for any eigenvalue of $f_{i\,a}$.
\end{itemize}
By tensoring some model bundle $L(\vecu)$ of rank 1,
we may also assume that $\sum_{a\in\Sp(\nbige^0,i)}|a|^2\neq 0$
for any $i$.

Since $f_k$ $(k=1,\ldots,l)$ and $g_j$ $(j=l+1,\ldots,n)$
are commutative,
$\EE_{\epsilon_1}(f_i,a)$ are preserved by $f_k$ and $g_j$.
By an inductive procedure,
we may assume the following:
\begin{condition} \mbox{{}}\label{condition;9.10.25}
\begin{itemize}
\item
 The subset
 $\Sp(\theta)
  \subset\prod_{i=1}^l\Sp(\nbige^0,i)$
 is given.
\item
 The holomorphic decomposition
 $E=\bigoplus_{\veca\in\Sp(\theta)}E_{\veca}$
 is given.
 \item
 Each $E_{\veca}$ is preserved by
 $f_k$ $(k=1,\ldots,l)$ and $g_j$ $(j=l+1,\ldots,n)$.
 Hence we have the decompositions
 \[
\begin{array}{ll}
  f_i=\bigoplus_{\veca\in\Sp(\theta)}f_{i\,\veca},
 & f_{i\,\veca}\in End(E_{\veca}),\\
 \mbox{{}}\\
 g_j=\bigoplus_{\veca\in\Sp(\theta)}g_{j\,\veca},
 & g_{j\,\veca}\in End(E_{\veca}).
\end{array}  
 \]
\item
 There exist positive constants $\epsilon$ and $C$
 such that the inequality
 $|\alpha-q_i(\veca)|\leq C\cdot|z|^{\epsilon}$ holds
 for any eigenvalue $\alpha$ of $f_{i\,\veca}$.
 Here $q_i$ denotes the projection onto the $i$-th component.
\item
 We have $\sum_{a\in \Sp(\nbige^0,i)}|a|^2\neq 0$.
\hfill\qed
\end{itemize}
\end{condition}

Using Lemma \ref{lem;8.16.1} inductively,
we obtain the following.
\begin{lem} \mbox{{}}
\begin{itemize}
\item
 We put $h_0:=\bigoplus_{\veca\in \Sp(\theta)}h_{|E_{\veca}}$.
 Then $h_0$ and $h$ are mutually bounded.
\item
 Namely there exists positive constants $C_1$ and $C_2$
 such that the following inequalities hold
 for any section $v=\sum_{\veca\in\Sp(\theta)} v_{\veca}$:
\[
 C_1\sum_{\veca\in\Sp(\theta)} |v_{\veca}|_{h}
\leq
 |v|_{h}
\leq
 C_2\sum_{\veca\in\Sp(\theta)}|v_{\veca}|_{h}.
\]
\end{itemize}
\hfill\qed
\end{lem}

We put
$\lefttop{i}E_{a}=\bigoplus_{q_i(\veca)=a}E_{\veca}$.
We have two metrics on $\lefttop{i}E_{a}$:
that is,
$h_{|\lefttop{i}E_a}$
and $\bigoplus_{q_i(\veca)=a}h_{|E_{\veca}}$.
\begin{lem}\label{lem;d11.14.1}
The metrics $h_{|\lefttop{i}E_a}$
and $\bigoplus_{q_i(\veca)=a}h_{|E_{\veca}}$
are mutually bounded.
\end{lem}
\pf
It follows from Lemma \ref{lem;8.16.1}.
\hfill\qed

\subsubsection{Inequalities}
\label{subsubsection;10.11.40}

The restriction of $f_i$ to $\lefttop{j}E_{a}$
is denoted by $\lefttop{j}f_{i\,a}$.
\begin{lem}\label{lem;9.7.105}
We have the following inequalities
for a positive constant $C$
and for any $\veca$:
\begin{equation} \label{eq;9.7.107}
 |\lefttop{i}f_{i\,a}-a\cdot id_{\lefttop{i}E_{a}}|_h
\leq
 C\cdot(-\log|z_i|)^{-1}.
\end{equation}
\end{lem}
\pf
It follows from Corollary \ref{cor;9.7.20}.
\hfill\qed

\begin{lem}
We have the following inequality for some positive constant $C'$:
\begin{equation} \label{eq;9.7.106}
 |f_{i\,\veca}-q_i(\veca)\cdot id_{E_{\veca}}|_h
\leq
 C'\cdot(-\log|z_i|)^{-1}.
\end{equation}
\end{lem}
\pf
The estimate (\ref{eq;9.7.106})
follows from (\ref{eq;9.7.107}) and Lemma \ref{lem;d11.14.1}.
\hfill\qed

\vspace{.1in}

We have the decomposition $g_j=\bigoplus g_{j\,\veca}$,
where $g_{j\,\veca}\in End(E_{\veca})$.
The following lemma is easy.
\begin{lem}
We have $|g_{j\,\veca}|_h\leq C$
for some positive constant $C$.
\hfill\qed
\end{lem}

We have the adjoint maps $f_i^{\dagger}$ and the decomposition
$f^{\dagger}_i=\sum (f^{\dagger}_i)_{\veca\,\vecb}$,
where $(f_i^{\dagger})_{\veca\,\vecb}\in Hom(E_{\veca},E_{\vecb})$.
We also have the adjoint $g_j^{\dagger}$
and the decomposition
$g^{\dagger}_j=\sum (g_j^{\dagger})_{\veca\,\vecb}$,
where $(g_j^{\dagger})_{\veca\,\vecb}\in Hom(E_{\veca},E_{\vecb})$.

For any elements
$\veca$ and $\vecb$ of $\cnum^l$,
we put $\Diff(\veca,\vecb):=\{j\,|\,q_j(\veca)\neq q_j(\vecb)\}$.

\begin{lem}\mbox{{}}
There exist positive constant $C$ and $\epsilon>0$
such that the following holds:
\begin{description}
\item[(A)]
$\big|(f_i^{\dagger})_{\veca\,\veca}
 -\overline{q_i(\veca)}\cdot id_{E_{\veca}}\big|_h
\leq C\cdot(-\log|z_i|)^{-1}$.
\item[(B)]
If $q_i(\veca)=q_i(\vecb)$ and $\veca\neq\vecb$,
we have the following inequality:
\[
 \big|(f^{\dagger}_i)_{\veca\,\vecb}\big|_h
\leq
 C\cdot(-\log|z_i|)^{-1}\cdot
 \prod_{\substack{j\in\Diff(\veca,\vecb)}}
  |z_j|^{\epsilon}.
\]
\item[(C)]
 If $q_i(\veca)\neq q_i(\vecb)$,
then we have the following inequality:
\[
 \big|(f_i^{\dagger})_{\veca\,\vecb}\big|_h
\leq
 C\cdot \prod_{j\in\Diff(\veca,\vecb)}|z_j|^{\epsilon}.
\]
\item[(D)]
We have the following inequality:
\[
 \big|(g_j^{\dagger})_{\veca\,\vecb}\big|_h
\leq
 C\cdot\prod_{k\in\Diff(\veca,\vecb)}|z_k|^{\epsilon}.
\]
\end{description}
\end{lem}
\pf
\label{eq;9.15.17}
We put $\rho_i:=\sum_a a\cdot id_{\lefttop{i}E_a}$,
and $\bar{\rho}_i:=\sum_a \bar{a}\cdot id_{\lefttop{i}E_a}$.
We denote the adjoint of $\rho_i$ by $\rho_i^{\dagger}$.
Due to the result in the case of curves,
we obtain the following:
\begin{itemize}
\item
 $|\rho_i^{\dagger}-\bar{\rho}_i|_h\leq C\cdot|z_i|^{\epsilon}$,
 (Lemma \ref{lem;9.7.25}).
\item
 $|f_i-\rho_i|_h\leq C\cdot(-\log|z_i|)^{-1}$,
 (Corollary \ref{cor;9.7.110}).
\end{itemize}
In particular,
we have $|f_i^{\dagger}-\bar{\rho}_i|_h\leq C\cdot(-\log|z_i|)^{-1}$.
It implies $(A)$.

We also have the following decomposition for $j\neq i$:
\[
 (f_i-\rho_i)^{\dagger}=
\sum \lefttop{j}(f_i-\rho_i)^{\dagger}_{a,b},
\quad\quad
 \lefttop{j}(f_i-\rho_i)^{\dagger}_{a,b}
\in Hom(\lefttop{j}E_{a},\lefttop{j}E_b).
\]
We have the following inequalities
due to Lemma \ref{lem;8.16.1} and Lemma \ref{lem;9.7.120}:
\[
 \begin{array}{ll}
 \big|\lefttop{j}(f_i-\rho_i)^{\dagger}_{a\,b}\big|_h
\leq
 C\cdot|z_j|^{\epsilon}\cdot|f_i-\rho_i|_h
\leq
 C\cdot|z_j|^{\epsilon}\cdot(-\log|z_i|)^{-1},
 & (a\neq b)\\
 \mbox{{}}\\
 \big|\lefttop{j}(f_i-\rho_i)^{\dagger}_{a\,a}\big|_h
\leq C\cdot|f_i-\rho_i|_h
\leq
 C\cdot(-\log|z_i|)^{-1},
 & (a=b)\\
 \end{array}
\]
Due to Lemma \ref{lem;d11.14.1},
there exist positive constants $C$
and $\epsilon$
such that the following holds:
\begin{equation}\label{eq;d11.14.3}
 \big|(f_i-\rho_i)^{\dagger}_{\veca\,\vecb}\big|_h
\leq
 \left\{
 \begin{array}{ll}
  C\cdot|z_j|^{\epsilon}\cdot(-\log|z_i|)^{-1},
 &(q_j(\veca)\neq q_j(\vecb))\\
\mbox{{}}\\
 C\cdot(-\log|z_i|)^{-1},
 & (q_j(\veca)= q_j(\vecb)).
 \end{array}
 \right.
\end{equation}
For any subset $I\subset\lbar$
and for any positive number $\epsilon$,
there exists a positive number $\epsilon'$
such that the following inequality holds:
\begin{equation}\label{eq;d11.14.2}
\min\bigl\{|z_i|^{\epsilon}\,\big|\,i\in I\bigr\}
\leq
 \prod_{i\in I}|z_i|^{\epsilon'}.
\end{equation}
From (\ref{eq;d11.14.3}) and (\ref{eq;d11.14.2}),
we obtain the following inequality,
for some sufficiently small $\epsilon>0$:
\begin{equation} \label{eq;9.7.130}
 \big|(f_i-\rho_i)^{\dagger}_{\veca\,\vecb}\big|_h
\leq
 C\cdot(-\log|z_i|)^{-1}\cdot
 \prod_{\substack{j\neq i,\\ j\in\Diff(\veca,\vecb)}}
 |z_j|^{\epsilon}.
\end{equation}

Similarly we have the following:
\[
 |(\rho_i^{\dagger})_{\veca,\vecb}|_h
\leq
 \left\{
 \begin{array}{ll}
 C\cdot|z_j|^{\epsilon} & \mbox{ if }q_j(\veca)\neq q_j(\vecb)\\
\mbox{{}}\\
 C    & \mbox{ if }q_j(\veca)=q_j(\vecb).
 \end{array}
 \right.
\]
By definition,
we have $\bar{\rho}_{i\,\veca\,\vecb}=0$
if $\veca\neq\vecb$.
We also have the following for any $\veca$ and $\vecb$
from Lemma \ref{lem;9.7.25}:
\[
 \big|(\rho_i^{\dagger}-\bar{\rho}_i)_{\veca\,\vecb}\big|_h
\leq
 C\cdot|z_i|^{\epsilon}.
\]
In all,
we obtain the following for some positive constants $C$ and $\epsilon$:
\begin{equation} \label{eq;9.7.131}
 |(\rho_i^{\dagger}-\bar{\rho})_{\veca\,\vecb}|_h
\leq
 C\cdot|z_i|^{\epsilon}\cdot
 \prod_{\substack{j\neq i\\ j\in \Diff(\veca,\vecb)}}
 |z_j|^{\epsilon}.
\end{equation}
Hence, in the case $\veca\neq\vecb$,
we obtain the following inequality
from (\ref{eq;9.7.130}), (\ref{eq;9.7.131})
and $\bar{\rho}_{i\,\veca,\vecb}=0$ ($\veca\neq\vecb$):
\[
 |(f_i^{\dagger})_{\veca\,\vecb}|_h\leq
 C\cdot(-\log|z_j|)^{-1}\cdot
\prod_{\substack{j\neq i\\ j\in \Diff(\veca,\vecb)}}
 |z_j|^{\epsilon}.
\]
It implies the claim $(B)$.

In the case $j\neq i$, we have the following for some positive constants
$C$ and $\epsilon$:
\begin{equation}\label{eq;d11.14.4}
 \bigl|\lefttop{j}(f^{\dagger}_{i})_{a,b}\bigr|_h
\leq
 \left\{
 \begin{array}{ll}
 C\cdot|z_j|^{\epsilon}, & (a\neq b),\\
 \mbox{{}}\\
 C,   & (a=b).
 \end{array}
 \right.
\end{equation}
From (\ref{eq;d11.14.4}) and (\ref{eq;d11.14.2}),
we obtain the following estimate, for some positive constants
$C$ and $\epsilon$:
\[
 \bigl|(f^{\dagger}_{i})_{\veca\,\vecb}\bigr|_h
\leq
 \left\{
 \begin{array}{ll}
 C\cdot|z_j|^{\epsilon}, & \bigl(q_j(\veca)\neq q_j(\vecb)\bigr)\\
 \mbox{{}}\\
 C, & \bigl(q_j(\veca)=q_j(\vecb)\bigr).
 \end{array}
 \right.
\]
On the other hand,
in the case $q_i(\veca)\neq q_i(\vecb)$,
we have the following for some positive constants
$C$ and $\epsilon$:
\[
 \big|(f^{\dagger}_i)_{\veca\,\vecb}\big|_h
\leq
 C\cdot|z_i|^{\epsilon}.
\]
In all, we obtain the following inequality
for some positive constants $C$ and $\epsilon$,
if $q_i(\veca)\neq q_i(\vecb)$:
\[
 |(f_i^{\dagger})_{\veca\,\vecb}|_h
\leq
 C\cdot\prod_{j\in\Diff(\veca,\vecb)}|z_j|^{\epsilon}
\]
Hence we obtain the claim $(C)$.

The claim $(D)$ can be obtained similarly.
\hfill\qed

%% file: 10.2.tex

\subsubsection{Some consequences}

Let $\veca$ and $\vecb$ be elements of $\Sp(\theta)$.
We put as follows:
\[
 F_{(i,\veca),(j,\vecb)}:=
 f_{i\,\veca}\circ (f^{\dagger}_j)_{\vecb\,\veca}
-(f^{\dagger}_j)_{\vecb\,\veca}\circ f_{i\,\vecb}.
\]
\begin{lem}
There exist positive constants $C$ and $\epsilon$
such that the following inequalities hold:\\
The case $q_i(\veca)=q_i(\vecb)$ and $q_j(\veca)=q_j(\vecb)$:
\[
\quad
  |F_{(i,\veca),(j,\vecb)}|_h\leq
 C\cdot(-\log|z_i|)^{-1}\cdot(-\log|z_j|)^{-1}
 \cdot
 \prod_{\substack{k\in\Diff(\veca,\vecb)}}
 |z_k|^{\epsilon}.
\]
The case 
$q_i(\veca)=q_i(\vecb)$
and $q_j(\veca)\neq q_j(\vecb)$:
\[
 |F_{(i,\veca),(j,\vecb)}|_h\leq C\cdot(-\log|z_i|)^{-1}
 \cdot
 \prod_{\substack{k\in\Diff(\veca,\vecb)}}
 |z_k|^{\epsilon}.
\]
The case
$q_i(\veca)\neq q_i(\vecb)$
and $q_j(\veca)\neq q_j(\vecb)$:
\[
 |F_{(i,\veca),(j,\vecb)}|_h\leq C \cdot
 \prod_{\substack{k\in\Diff(\veca,\vecb)}}
 |z_k|^{\epsilon}.
\]
\end{lem}
\pf
Let us consider the case
$q_i(\veca)=q_i(\vecb)$
and $q_j(\veca)=q_j(\vecb)=d$.
We put $\delta_{\veca,\vecb}:=1$ $(\veca\neq \vecb)$,
or $\delta_{\veca,\vecb}:=0$ $(\veca\neq \vecb)$.
Then we have the following:
\[
 F_{(i,\veca),(j,\vecb)}=
\Bigl(f_{i\,\veca}-c\cdot\delta_{\veca\,\vecb}\cdot id_{E_{\veca\,\vecb}}
\Bigr)
\circ
\Bigl(
 (f^{\dagger}_j)_{\vecb,\veca}
-\bar{d}\cdot\delta_{\veca\,\vecb}\cdot id_{E_{\veca\,\vecb}}
\Bigr)
-
\Bigl(
 (f^{\dagger}_j)_{\vecb,\veca}
-\bar{d}\cdot\delta_{\veca\,\vecb}\cdot id_{E_{\veca\,\vecb}}
\Bigr)
\circ
 \Bigl(f_{i\,\vecb}-c\cdot\delta_{\veca\,\vecb}\cdot id_{E_{\veca\,\vecb}}
\Bigr).
\]
Therefore we obtain the following inequality:
\[
 |F_{(i,\veca),(j,\vecb)}|_h\leq
 C\cdot(-\log|z_i|)^{-1}\cdot(-\log|z_j|)^{-1}
 \cdot
 \prod_{\substack{k\in\Diff(\veca,\vecb)}}
 |z_k|^{\epsilon}.
\]
Thus we are done in the case
$q_i(\veca)=q_i(\vecb)$ and $q_j(\veca)=q_j(\vecb)=d$.
The rest cases can be shown similarly.
\hfill\qed

\begin{prop} \label{prop;9.7.140}
The two form
$\theta\wedge\theta^{\dagger}+\theta^{\dagger}\wedge\theta$
is dominated by the Kahler form of the Poincar\'{e} metric:
\[
\omega_{\poin}=
 \sum_{i=1}^l \frac{dz_i\cdot d\bar{z}_i}{|z_i|^2(-\log|z_i|)^2}
+\sum_{i=l+1}^n\frac{dz_i\cdot d\bar{z}_i}{(1-|z_i|^2)^2}.
\]
\end{prop}
\pf
The contributions of $f_i\cdot{dz_i}/{z_i}$
and $f_j^{\dagger}\cdot d\bar{z}_j/d\bar{z}_j$ to
the form $\theta\wedge\theta^{\dagger}+\theta^{\dagger}\wedge\theta$
is a sum of the following forms:
\[
 F_{(i,\veca),(j,\vecb)}
 \cdot\frac{dz_i}{z_i}\wedge \frac{d\bar{z}_j}{\bar{z}_j}.
\]
Then it is dominated by the form
$q_i^{\ast}\omega_{\Delta^{\ast}}+q_j^{\ast}\omega_{\Delta^{\ast}}$.
The rest terms can be dominated similarly.
\hfill\qed

\vspace{.1in}

We have the unitary connection of $\nbigelambda$
given by the following:
\[
 \delbar_{\nbigelambda}+\del_{\nbigelambda}
=\delbar_E+\del_E
+\lambda\cdot\theta^{\dagger}
-\bar{\lambda}\cdot\theta.
\]
The curvature $R(\delbar_{\nbigelambda}+\del_{\nbigelambda})$
is as follows:
\[
 \bigl(1+|\lambda|^2\bigr)\cdot
\bigl(\theta\wedge\theta^{\dagger}
+\theta^{\dagger}\wedge\theta\bigr).
\]

\begin{cor}\label{cor;10.11.80}
The curvature $R(\delbar_{\nbigelambda}+\del_{\nbigelambda})$
is dominated by $\omega_{\poin}$.
Namely, the hermitian holomorphic bundle
$(\nbigelambda,h)$ is acceptable.
\end{cor}
\pf
It follows from Proposition \ref{prop;9.7.140}.
\hfill\qed

\begin{rem}\label{rem;b12.5.6}
The curvature of the hermitian holomorphic bundle $(\nbige,h)$
over $\nbigx-\nbigd$
contains the terms of the following form:
\[
-d\overline{\lambda}\wedge\theta
+d\lambda\wedge\theta^{\dagger}.
\]
Hence $(\nbige,h)$ is not acceptable,
unless $\theta$ is nilpotent.
\hfill\qed
\end{rem}

%% file: 11.tex
\subsubsection{The quasi canonical prolongation} \label{subsubsection;9.8.5}

Definition \ref{df;a11.14.1} is generalized as follows.
\begin{df}
$\lambda$ is called generic with respect to $\harmonicbundle$
if the maps
$\kmsmap(\lambda):\KMS(\nbige^0,i)\lrarr \Sp(\nbigelambda,i)$
are bijective for any $i$.
\hfill\qed
\end{df}

Assume that $\lambda$ is generic.
Let $P$ be a point of $X-D$.
Then we have the endomorphisms
$M^{\lambda}_i$ on $\nbigelambda_{|P}$.
Here $M^{\lambda}_i$
is the monodromy with respect to $D_i$.
We have $\Sp(M^{\lambda}_i)=\Sp^f(\nbigelambda,i)$.
We put $\vecM^{\lambda}=(M^{\lambda}_1,\ldots,M^{\lambda}_l)$.
We have the decomposition for each $P$:
\[
 \nbigelambda_{|P}=\bigoplus_{\vecomega\in\Sp(\vecM^{\lambda})}
 \EE(\vecM^{\lambda},\vecomega).
\]
Thus we obtain the decomposition
$\nbigelambda=\bigoplus_{\vecomega\in\Sp(\vecM^{\lambda})}
 \EE(\nbigelambda,\vecomega)$.
It is compatible with the flat connection $\DD^{\lambda\,f}$.
We put $q_i(\vecomega)=\omega_i$.
We put
$\nbigk(\nbige,\lambda,0,i):=
 \bigl\{u\in\nbigk(\nbige^0,i)\,\big|\,-1<\paramap(\lambda,u)\leq 0
 \bigr\}$.
Then we have the unique element
$u_i=u_i(\vecomega)\in\nbigk(\nbige,\lambda,0,i)$
such that $\eigenmap^f(\lambda,u_i(\vecomega))=\omega_i$.
Then we obtain the number
$b_i=b_i(\vecomega):=\paramap^f(\lambda,u_i(\vecomega))$.
We put $\vecb(\vecomega):=(b_1,\ldots,b_l)$.

For $s\in\EE\bigl(H(\nbigelambda),\vecomega\bigr)$,
we put as follows
(see the page \pageref{page;9.11.15} for $\alpha(b,\omega)$):
\[
 F\bigl(s,\vecb(\vecomega)\bigr):=
\exp\Bigl(
 \sum \log z_i\cdot
 \bigl(
 \alpha(b_i,\omega_i)+N_{i\,\vecomega}
 \bigr)
 \Bigr)\cdot s.
\]
Note that we have the equality
$-\ord(s_{|\pi_i^{-1}(P)})=b_i$
for any point $P\in D_i^{\circ}$.

\begin{lem} \label{lem;9.8.1}
We have the following, for some positive constants $C$ and $M$:
\[
 |F(s,\vecb(\vecomega))|_h
\leq
 C\cdot
 \prod_{i=1}^l |z_i|^{-\paramap(\lambda,u_i)}
\cdot
 \Bigl(
 -\sum_{i=1}^l \log|z_i|
 \Bigr)^M
\]
\end{lem}
\pf
For each $P\in D_i^{\circ}$,
we have the equality
$-\ord\bigl(
 F(s,\vecb(\vecomega)\bigr)_{|\pi_i^{-1}(P)})
=\paramap(\lambda,u_i)$
(Lemma \ref{lem;a11.15.5}).
Then we obtain the result
due to Corollary \ref{cor;11.28.15}
and Corollary \ref{cor;10.11.80}.
\hfill\qed

\begin{cor}
$F(s,\vecb(\vecomega))$ are sections of $\prolong{\nbigelambda}$.
\hfill\qed
\end{cor}

Let $\vecs=(s_j)$ be a base of $H(\nbigelambda)$
compatible with $\EE$.
We put $\vecomega_j:=\vecdeg^{\EE}(s_j)$.
We put $v_j:=F(s_j,\vecb(\vecomega_j))$,
and $\vecv=(v_j)$.
It is a tuple of sections of $\prolong{\nbigelambda}$.

\begin{lem} \label{lem;d11.14.51}
If $\lambda$ is generic,
the $\nbigo_X$-module $\prolong{\nbigelambda}$ is 
coherent and locally free,
and $\vecv$ is a frame of $\prolong{\nbigelambda}$.
\end{lem}
\pf
For any point $P\in D_i^{\circ}$,
the tuple $\vecv_{|\pi_i^{-1}(P)}$ is a frame of
$\prolong{(\nbigelambda_{|\pi_i^{-1}(P)})}$,
due to Lemma \ref{lem;d11.14.20}.
We put as follows:
\[
 c_i:=
\sum_{a\in
 \Par(\prolong{\nbigelambda},i)}
\multiplicity(\lambda,a)\cdot a.
\]
We have the line bundle
$\leftbottom{\vecc}\det(\nbigelambda)$ over $X$.
We have the following on $\pi_i^{-1}(P)$:
\[
 \bigl(\leftbottom{\vecc}\det(\nbigelambda)\bigr)_{|\pi_i^{-1}(P)}
=\leftbottom{c_i}\bigl(\det(\nbigelambda\bigr)_{|\pi_i^{-1}(P)})
=\det\bigl(\prolong{\bigl(\nbigelambda_{|\pi_i^{-1}(P)} \bigr)}\bigr).
\]

We put $\Omega(\vecv)=v_1\wedge\cdots\wedge v_r$.
Then $\Omega(\vecv)_{|\pi_i^{-1}(P)}$
gives a frame of 
$\det\bigl(\prolong{\bigl(\nbigelambda_{|\pi_i^{-1}(P)} \bigr)}\bigr)$.
Hence $\Omega(\vecv)$
is a frame of $\leftbottom{\vecc}\det(\nbigelambda)$
over $X$.

Let $f$ be a holomorphic section
of $\prolong{\nbigelambda}$.
Then $f$ is described as $\sum f_j\cdot v_j$
for holomorphic functions $f_i$ on $X-D$.
Let us consider 
$f\wedge v_2\wedge\cdots\wedge v_r$.
On the curve $\pi_i^{-1}(P)$,
it gives a section of
$\leftbottom{\vecc}\det(\nbigelambda)_{|\pi_i^{-1}(P)}$.
Hence $f\wedge v_2\wedge\cdots\wedge v_r$
is a holomorphic section of
$\leftbottom{\vecc}\det(\nbigelambda)$
over $X$.
We have
$f\wedge v_2\wedge\cdots\wedge v_r=f_1\cdot\Omega(\vecv)$.
Hence $f_1$ is holomorphic over $X$.
Similarly $f_j$ is holomorphic for each $j$.
\hfill\qed

\begin{rem}
The last argument can be found in {\rm\cite{cg}}.
\hfill\qed
\end{rem}

\subsubsection{Weak norm estimate}
\label{subsubsection;a11.15.50}

In the case $s_j\in \EE\bigl(H(\nbigelambda),\vecomega\bigr)$,
we put $u_i(v_j):=u_i(\vecomega)$.
We put as follows:
\[
 v_j':=v_j\cdot\prod_{i=1}^l|z_i|^{\paramap(\lambda,u_i(v_j))}.
\]
We put $\vecv':=(v_j')$.
Then there exist positive constants $M$ and $C$
satisfying the following, due to Lemma \ref{lem;9.8.1}:
\begin{equation} \label{eq;9.8.3}
 H(h,\vecv')\leq
 C\cdot
 \Bigl(
 -\sum_{i=1}^l\log|z_i|
 \Bigr)^M.
\end{equation}

Let $\vecv^{\lor}=(v_j^{\lor})$
be the dual frame of $\vecv$ over $X-D$.
\begin{lem}\label{lem;9.8.2}
There exits positive constants $C$ and $M$
satisfying the following:
\[
 |v_j^{\lor}|_h
\leq
 C\cdot \prod_i |z_i|^{\paramap(\lambda,u_i(v_j))}
 \cdot
 \Bigl(
 -\sum_{i=1}^l\log |z_i|
 \Bigr)^M.
\]
\end{lem}
\pf
We have the estimates for the restrictions
$v^{\lor}_{j\,|\,\pi_i^{-1}(P)}$
due to the result in the subsubsection \ref{subsubsection;d11.14.21}.
Then we can derive the estimate desired,
due to Corollary \ref{cor;11.28.15}.
\hfill\qed

\vspace{.1in}

We put as follows:
\[
 v_j^{\lor\prime}:=
 v_j^{\lor}\cdot
 \prod_i |z_i|^{-\paramap(\lambda,u_i(v_j))}.
\]
We put $\vecv^{\lor\,\prime}:=(v_j^{\lor\,\prime})$.
Then $\vecv^{\lor\prime}$ is the dual frame
of $\vecv'$ over $X-D$,
and there exist positive constants $C$ and $M$
satisfying the following, due to Lemma \ref{lem;9.8.2}:
\begin{equation}\label{eq;9.8.4}
 H\bigl(h^{\lor},\vecv^{\lor\,\prime}\bigr)
\leq C\cdot
 \Bigl(
-\sum_{i=1}^l\log|z_i|
 \Bigr)^M.
\end{equation}
Hence we obtain the following.
\begin{lem} \label{lem;9.8.10}
The frame $\vecv'$ is adapted up to log order.
Namely,
there exist positive constants $C_1$, $C_2$
and $M$ satisfying the following:
\[
 C_1\cdot
 \Bigl(
 -\sum_{i=1}^l\log|z_i|
 \Bigr)^{-M}
\leq
 H(h,\vecv')
\leq
 C_2\cdot
 \Bigl(
 -\sum_{i=1}^l\log|z_i|
 \Bigr)^M.
\]
\end{lem}
\pf
The right inequality is given in (\ref{eq;9.8.3}).
The left inequality immediately follows from (\ref{eq;9.8.4}).
\hfill\qed

%% file: a73.tex

\subsubsection{Minor generalization and a restriction to a diagonal curve}

Let $\vecc$ be an element of $\real^l$.
We put
$\nbigk(\nbige,\lambda,\vecc,i):=
 \bigl\{
  u\in\KMS(\nbige^0,i)\,\big|\,
 c_i-1<\paramap(\lambda,u)\leq c_i
 \bigr\}$.
Let $\vecs$ be a frame of $H(\nbigelambda)$
as in the subsubsection \ref{subsubsection;9.8.5}.
We put
$\prolongg{\vecc}{v_j}:=F\bigl(s_j,\vecb(\vecomega_j)-\vecc\bigr)$,
and we put $\prolongg{\vecc}{\vecv}:=\bigl(\prolongg{\vecc}{v_j}\bigr)$.
\begin{lem}
If $\lambda$ is generic,
the sheaf $\prolongg{\vecc}{\nbigelambda}$ is locally free,
and $\prolongg{\vecc}{\vecv}$ gives a frame of
$\prolongg{\vecc}{\nbigelambda}$.
\end{lem}
\pf
It can be shown by an argument similar to the proof
of Lemma \ref{lem;d11.14.51}.
\hfill\qed

\vspace{.1in}
In the case $s_j\in \EE\bigl(H(\nbigelambda),\vecomega\bigr)$,
$u_i(\prolongg{\vecc}{v_j})$ denotes the unique element
of $\nbigk(\nbige,\lambda,\vecc,i)$ such that
$\eigenmap^f(\lambda,u_i(\prolongg{\vecc}{v_j}))
=\omega_i$.
We put as follows:
\[
 \prolongg{\vecc}{v_j}':=
 \prolongg{\vecc}{v_j}\cdot 
 \prod_{i=1}^l |z_i|^{\paramap(\lambda,u_i(\prolongg{\vecc}{v_j}))}.
\]
We put $\prolongg{\vecc}{\vecv}'=(\prolongg{\vecc}{v_j}')$.
\begin{lem}\label{lem;a11.15.10}
The frame $\prolongg{\vecc}{\vecv}'$ is adapted up to log order.
\end{lem}
\pf
It can be shown by an argument similar to
Lemma \ref{lem;9.8.10}.
\hfill\qed

\begin{lem} \label{lem;a11.15.11}
We have the following relation:
\[
 \DDlambda\prolongg{\vecc}{\vecv}
=\prolongg{\vecc}{\vecv}
\cdot
 \sum_j\bigl(C_j+N_j\bigr)\cdot\frac{dz_j}{z_j}.
\]
Here $C_j$ denote the diagonal matrices
whose $(i,i)$-components are
$\eigenmap\bigl(\lambda,u_j(\prolongg{\vecc}{v_i})\bigr)$,
and $N_j$ denote the nilpotent matrices.
\end{lem}
\pf
It follows from Lemma \ref{lem;a11.15.5}.
\hfill\qed

\vspace{.1in}

Let $\epsilon$ is a sufficiently small positive number,
and we put $\vecc:=(\overbrace{\epsilon,\ldots,\epsilon}^l)$.
and $\epsilon_i$ $(i=1,2)$ be positive numbers such that
$\epsilon_1+\epsilon_2<1/2$.
Assume the following:
\begin{itemize}
\item
 Let $a$ be elements of
 $\Par(\prolongg{\vecc}{\nbigelambda},i)$.
 Then $0<\epsilon-a<\epsilon_i$.
\end{itemize}

Pick an element
$(\lefttop{0}z_3,\ldots,\lefttop{0}z_{n})
 \in(\Delta^{\ast})^{n-2}$.
Let us consider the curve
$C_0:=\{(z,z,\lefttop{0}z_3,\ldots,\lefttop{0}z_n)\in X\}$.
We have the restriction of
$\harmonicbundle$ to $C_0$.

\begin{lem}\mbox{{}} \label{lem;a11.15.30}
\begin{enumerate}
\item
We have the following implication:
\[
 \KMS\bigl(
 \prolongg{2\epsilon}{\bigl(\nbige^{\lambda}_{|C_0}\bigr)}
 \bigr)
\subset
 \bigl\{
 u_1+u_2\,\big|\,
 u_i\in\
 KMS\bigl(\prolongg{\vecc}{\nbigelambda},i\bigr)
 \bigr\}.
\]
\item
We have the following equality:
\begin{equation}\label{eq;a11.15.25}
\sum_{
 b\in\Par\bigl(
 \prolongg{2\epsilon}{\bigl(\nbige^{\lambda}_{|C_0}\bigr)}
 \bigr)}
 b\cdot \multiplicity(\lambda,b)
=\sum_{i=1,2}
 \sum_{b_i\in\Par
 \bigl(\prolongg{\vecc}{\nbigelambda},i\bigr)}
 b_i\cdot\multiplicity(\lambda,b_i).
\end{equation}
\end{enumerate}
\end{lem}
\pf
Let us take the frame
$\vecw:=\prolongg{\vecc}{\vecv}$ as above.
For each $w_j$, we have
the elements
$u_j(w_i)\in\nbigk\bigl(\nbige,\lambda,\vecc,j\bigr)$.

Let us consider the restriction $\tilde{w}_j:=w_{j\,|\,C_0}$.
We put as follows:
\[
 \widetilde{w}_j':=\tilde{w}_j\cdot
 |z|^{\paramap\bigl(\lambda,u_1(w_j)+u_2(w_j)\bigr)},
\quad
 \widetilde{\vecw}':=
 \bigl(\widetilde{w}_j'\bigr).
\]
Then the $C^{\infty}$-frame $\widetilde{\vecw}'$
is adapted up to log order,
due to Lemma \ref{lem;a11.15.10}.
Thus $\widetilde{\vecw}$ gives a frame of
$\prolongg{2\epsilon}{\bigl(\nbige^{\lambda}_{|C_0}\bigr)}$.
We also have the following relation,
due to Lemma \ref{lem;a11.15.11}:
\begin{equation} \label{eq;a11.15.12}
 \DD\widetilde{\vecw}=
\widetilde{\vecw}\cdot
 \bigl(C_1+C_2+N_1+N_2\bigr)\frac{dz}{z}.
\end{equation}
Here $C_i$ and $N_j$ denote the matrices
given in Lemma \ref{lem;a11.15.11} for
$\vecw=\prolongg{\vecc}{\vecv}$.
From the adaptedness of $\widetilde{\vecw}'$ up to log order
and (\ref{eq;a11.15.12}),
we obtain the following:
\[
 \deg^{F,\EE}\bigl(\widetilde{w}_j\bigr)
=\kmsmap\bigl(\lambda,u_1(v_j)\bigr)
+\kmsmap\bigl(\lambda,u_2(v_j)\bigr).
\]
Thus we obtain the first claim.

We also obtain the following equality:
\begin{equation}\label{eq;a11.15.26}
 \sum_j
 -\ord\bigl(\widetilde{w}_j\bigr)
=\sum_j\sum_{i=1,2}
 -\lefttop{i}\ord\bigl(v_j\bigr).
\end{equation}
The left (resp. right) hand side of
(\ref{eq;a11.15.26})
are same as the left (resp. right) hand side of (\ref{eq;a11.15.25}).
Thus we obtain the second claim.
\hfill\qed

%% file: 12.tex

\subsubsection{Preliminary I, Estimates of Higgs fields}

\label{subsubsection;d11.14.31}

We put $X=\Delta_{\zeta}\times\Delta_z^{n-1}$
and $X^{(0)}=\{0\}\times\Delta_z^{n-1}$.
Let $D'_i$ denote the divisor of $\Delta_z^{n-1}$
defined by $z_i=0$.
We put $D_i=\Delta_{\zeta}\times D'_i$
and $D=\bigcup_{i=1}^l D_i$.
We put $D^{(0)}_i=\{0\}\times D'_i$
and $D^{(0)}=\bigcup_{i=1}^l D^{(0)}_i$.

Let $\lambda_0\in\cnum_{\lambda}$ be {\em generic}.
Let $\vecv$ be a frame of $\prolong{\nbige^{\lambda_0}}$
as in the subsubsection \ref{subsubsection;9.8.5}.
We put 
$ \vecb(v_i)
:=\vecdeg^F(v_i)
 \in ]-1,0]^l$.
Let $b_h(v_j)$ denotes the $h$-th component of $\vecb(v_j)$.
We put as follows:
\[
 v_i':=v_i\cdot\prod |z_h|^{b_h(v_i)},
\quad
 \vecv'=(v_i').
\]
We define the function $B:X-D\lrarr M(r)$
as follows:
\[
 B_{i\,j}:=
 \left\{
 \begin{array}{ll}
  \prod_{h=1}^l|z_h|^{b_h(v_i)}, &(i=j),\\
\mbox{{}}\\
 0, & (i\neq j).
 \end{array}
 \right.
\]
Then we have $\vecv'=\vecv\cdot B$.

We put $\delbar_{\nbigelambdazero}=\delbar_E+\lambda_0\cdot\theta^{\dagger}$,
which is the holomorphic structure of $\nbigelambdazero$.
We define $K\in C^{\infty}\bigl(X-D,\Omega^{1,0}_X\otimes M(r)\bigr)$
by $K=B^{-1}\cdot\delbar B=\delbar B\cdot B^{-1}$.
The following lemma is clear.
\begin{lem} \label{lem;9.8.11}
We have the following formula:
\[
 K_{i\,j}:=
 \left\{
 \begin{array}{ll}
 {\displaystyle
 \sum_{h=1}^l\frac{b_h(v_i)}{2}\frac{d\bar{z}_h}{\bar{z}_h},
 } & (i=j),\\
 \mbox{{}}\\
 0, & (i\neq j).
 \end{array}
 \right.
\]
We have $\delbar_{\nbigelambdazero}\vecv'=\vecv'\cdot K$.
\hfill\qed
\end{lem}

The one forms
$\Theta\in C^{\infty}\bigl(X-D,\Omega^{1,0}\otimes M(r)\bigr)$
and $\Theta^{\dagger}\in 
 C^{\infty}\bigl(X-D,\Omega^{0,1}_X\otimes M(r)\bigr)$
are given as follows:
\[
 \theta\cdot\vecv'=\vecv'\cdot\Theta,
\quad\quad
 \theta^{\dagger}\cdot\vecv'=
 \vecv'\cdot\Theta^{\dagger}.
\]
The functions $\Theta^{\zeta}$, $\Theta^k$,
$\Theta^{\dagger\,\zeta}$
and $\Theta^{\dagger\,k}$
are defined as follows:
\[
\begin{array}{l}
{\displaystyle
 \Theta=
 \Theta^{\zeta}\cdot d\zeta
+\sum_{k=1}^l\Theta^k\cdot\frac{dz_k}{z_k}
+\sum_{k=l+1}^n\Theta^k\cdot dz_k,}\\
\mbox{{}}\\
 {\displaystyle
  \Theta^{\dagger}=
 \Theta^{\dagger\zeta}\cdot d\bar{\zeta}
+\sum_{k=1}^l\Theta^{\dagger\,k}\cdot\frac{d\bar{z}_k}{\bar{z}_k}
+\sum_{k=l+1}^n\Theta^{\dagger\,k}\cdot d\bar{z}_k.}
\end{array}
\]

\begin{lem} \label{lem;9.8.15}
The functions $\Theta^{\zeta}$, $\Theta^k$, $\Theta^{\dagger\,\zeta}$
and $\Theta^{\dagger\,k}$
are dominated by polynomials of
$-\log|z_i|$ $(i=1,\ldots,l)$.
\end{lem}
\pf
For the decompositions
$\theta=\sum_{i=1}^l f_i\cdot {dz_i}/{z_i}
 +\sum_{i=l+1}^n g_i\cdot dz_i$,
the norms of $f_i$ and $g_j$ with respect to $h$
are bounded
(see the subsubsection \ref{subsubsection;10.11.40}).
We also know that $\vecv'$ is adapted up to log order
(Lemma \ref{lem;9.8.10}).
Then the result follows.
\hfill\qed

\vspace{.1in}

We have the $\lambda_0$-connection $\DD^{\lambda_0}$.
Let $A$ denote the $\lambda_0$-connection one form
of $\DD^{\lambda_0}$ with respect to $\vecv$,
i.e.,
$A$ is a holomorphic section of $M(r)\otimes\Omega^{1,0}_X(\log D)$
determined by the condition:
\[
 \DD^{\lambda_0}\vecv=\vecv\cdot A.
\]
Due to our choice of $\vecv$,
we have the decomposition
$A=\bigoplus A_{\vecomega}$
which corresponds the decomposition
$\nbige^{\lambda_0}
=\bigoplus_{\vecomega}\EE(\nbigelambdazero,\vecomega)$,
namely,
$A_{\omega}$ is a holomorphic section of
$\End\bigl(\EE(\nbigelambdazero,\vecomega)\bigr)
  \otimes\Omega^{1,0}(\log D)$.

\begin{lem}
$A$ and $B$ are commutative.
\end{lem}
\pf
It follows from the decomposition
$B=\bigoplus_{\vecomega}
   B(\vecomega)\cdot Id_{\EE(\nbigelambdazero,\vecomega)}$
for $B(\vecomega)\in C^{\infty}(X-D)$.
\hfill\qed

\vspace{.1in}

We have the following formula:
\begin{equation}\label{eq;9.8.13}
\DD^{\lambda_0}\vecv'=\DD^{\lambda_0}(\vecv\cdot B)
=\vecv\cdot A\cdot B+\vecv\cdot \lambda_0\del(B)
=\vecv'\cdot(B^{-1}AB+\lambda_0B^{-1}\del B)
=\vecv'(A+\lambda_0\overline{K}).
\end{equation}
We have the description:
\[
 A+\lambda_0\cdot \overline{K}=
 \sum_{i=1}^l g_i\cdot \frac{dz_i}{z_i}.
\]
\begin{lem}
The functions $g_i$ are bounded.
\end{lem}
\pf
Since $A$ is a holomorphic section of $M(r)\otimes\Omega^{1,0}(\log D)$,
the lemma follows from Lemma \ref{lem;9.8.11}.
\hfill\qed

\vspace{.1in}

\begin{lem}
We have the following equalities:
\begin{equation} \label{eq;9.8.12}
 \del \Theta^{\dagger}+
 \Bigl[
  \lambda_0^{-1}\cdot(A+\overline{K}-\Theta),
  \Theta^{\dagger}
 \Bigr]=0,
\quad
 \delbar\Theta+
 \Bigl[
  K-\lambda_0\Theta^{\dagger},\Theta
 \Bigr]=0.
\end{equation}
\end{lem}
\pf
Recall the equalities
$\del_E(\theta^{\dagger})=0$
and $\del_E=\lambda_0^{-1}(\DD^{\lambda_0}-\theta)$.
We obtain the first equality in (\ref{eq;9.8.12})
from (\ref{eq;9.8.13}) and the definition of $\Theta$.

As for the second equality in (\ref{eq;9.8.12}),
we have only to use $\delbar_E(\theta)=0$,
$\delbar_E=\delbar_{\nbigelambdazero}-\lambda_0\theta^{\dagger}$
and $\delbar_{\nbigelambdazero}\vecv'=\vecv'\cdot K$.
\hfill\qed

\vspace{.1in}

\begin{lem}\label{lem;9.9.62}
We regard 
$\Theta^{\dagger}(\zeta,z)$ and $\Theta(\zeta,z)$
as functions of $\zeta$.
Then the following holds:
\begin{itemize}
\item
They are $L^p_k$ for any $k\in\seisuu_{\geq 0}$
and for any $p$.
\item
The $L^p_k$-norms of $\Theta^{\dagger}(\zeta,z)$
and $\Theta(\zeta,z)$
are dominated by polynomials
of $(-\log|z_j|)$ $(j=1,\ldots,l)$.
\end{itemize}
\end{lem}
\pf
It follows from 
Lemma \ref{lem;9.8.15},
the formula (\ref{eq;9.8.12}),
the estimate for $A$ and $K$,
and Corollary \ref{cor;9.9.61}.
\hfill\qed

\begin{lem} \label{lem;9.9.71}
We have the descriptions as follows:
\[
 \del_k\Theta^{\dagger\,\zeta}
=C_k\cdot \frac{dz_k}{z_k}\cdot d\bar{\zeta}.
\]
Here, $C_k$ 
denotes an element of $C^{\infty}\bigl(X-D,M(r)\bigr)$.
Their $L^{\infty}$-norms
are dominated by polynomials of $-\log|z_j|$ $(j=1,\ldots,l)$.
\end{lem}
\pf
It follows from
the estimates for $\Theta$, $A$ and $K$,
due to Lemma \ref{lem;d11.14.30}.
\hfill\qed

\subsubsection{Preliminary II, Estimate of cocycles}
\label{subsubsection;a11.15.60}

Let $\lambda$ be any element of $\cnum_{\lambda}$,
and $\lambda_0$ be generic
as in the subsubsection \ref{subsubsection;d11.14.31}.
Let $f$ be a holomorphic section of
the sheaf $\prolong{\bigl(
 \nbigelambda_{|X^{(0)}-D^{(0)}}
 \bigr) }$ over $X$.
Since the tuple of sections $\vecv'$ gives $C^{\infty}$-frame of 
$C^{\infty}$-bundle $\nbigelambda=\nbige^{\lambda_0}$
over $X-D$,
we have the following description:
\[
 f=\sum_{i} f_i(z)\cdot v_i'(0,z),
\quad\quad
 f_i\in C^{\infty}(X^{(0)}-D^{(0)}).
\]

\begin{lem} \label{lem;10.11.41}
Assume $-\vecord(f)\leq 0$.
Then there exist positive constants $C$ and $M$
satisfying the following:
\[
 |f|_h\leq C\cdot\Bigl(-\sum_{i=1}^l\log|z_i|\Bigr)^M.
\]
\end{lem}
\pf
We obtain such estimate for the restrictions
$f_{|\pi_i^{-1}(P)}$ for any $P\in D_i^{\circ}$.
Then we obtain the result due to Corollary \ref{cor;11.28.15}.
\hfill\qed

\vspace{.1in}
Due to the adaptedness of $\vecv'$ up to log order,
we obtain the following inequalities for
some positive constants $C$ and $M$ and for any $i=1,\ldots,l$:
\begin{equation}\label{eq;9.9.60}
 |f_i|\leq C\cdot \Bigl(-\sum_{j=1}^l \log|z_j|\Bigr)^M.
\end{equation}

Note that we have the following relation:
\begin{equation} \label{eq;9.9.72}
 \delbar_{\nbigelambda}\vecv'=
\Bigl(
\delbar_{\nbigelambdazero}
+(\lambda-\lambda_0)\cdot\theta^{\dagger}
\Bigr)\vecv'=
\vecv'\cdot\Bigl(
 K+(\lambda-\lambda_0)\cdot\Theta^{\dagger}
 \Bigr).
\end{equation}
Hence we obtain the following equality on $X^{(0)}-D^{(0)}$:
\begin{equation} \label{eq;9.9.70}
 0=\delbar f(z)=
\sum_i \Bigl(
  \delbar f_i(z)+
\sum_{k,j}
  \bigl(
   \Theta^{k\,\dagger}_{i\,j}(0,z)
  +\frac{b_k(v_i)}{2}\cdot\delta_{i\,j}
  \bigr)\cdot\bar{\eta}_k\cdot f_j
 \Bigr)\cdot v_i(0,z).
\end{equation}
Here we put as follows:
\[
 \delta_{i\,j}:=
 \left\{
  \begin{array}{ll}
   1 & (i=j)\\ \mbox{{}}\\
   0 & (i\neq j)
  \end{array}
 \right.
\quad\quad
 \bar{\eta}_k=
 \left\{
 \begin{array}{ll}
  {d\bar{z}_k}/{\bar{z}_k} &(k\leq l)\\ \mbox{{}}\\
  d\bar{z}_k & (k\geq l+1).
 \end{array}
 \right.
\]

We have the line bundle $\nbigo_X(-X_0)$.
We put $\nbigelambda(-X_0):=\nbigelambda\otimes \nbigo(-X_0)$.

Let $\epsilon$ be any sufficiently small positive number.
Let $\vecdelta$ denote the element 
$(\overbrace{1,\ldots,1}^l)$.
\begin{prop} \label{prop;9.9.85}
There exists an element $\rho\in C^{\infty}(X-D,\nbigelambda)$
satisfying the following:
\begin{enumerate}
\item \label{8.17.10}
 $\rho\in A^{0,0}_{-\epsilon\cdot \vecdelta,N}(\nbigelambda)$ for 
 any real number $N$.
\item \label{8.17.11}
 $\delbar_{\nbigelambda}\rho\in
 A^{0,1}_{-\epsilon\vecdelta,N}(\nbigelambda(-X_0))$
 for any real number $N$.
\item \label{8.17.12}
 $\rho_{|X^{(0)}-D^{(0)}}=f$.
\item \label{8.17.13}
 The support of $\rho$ is contained in the region
 $\bigl\{\zeta\,\big|\,|\zeta|<1/3\bigr\}
 \times \bigl(X^{(0)}-D^{(0)}\bigr)$.
\end{enumerate}
See the subsubsection {\rm\ref{subsubsection;d11.14.35}}
for the notation.
\end{prop}
\pf
Let $\chi$ be any function over $\Delta_{\zeta}$
satisfying the following:
\[
 \chi(\zeta)=
\left\{
 \begin{array}{ll}
 1 & |\zeta|\leq 1/3 \\ \mbox{{}}\\
 0 & |\zeta|>2/3.
 \end{array}
\right.
\]
We put as follows:
\[
 f^1:=\sum_i f_i(z)\cdot v_i(\zeta,z),
\quad\quad
 f^2:=\overline{\zeta}\cdot g^2,
\quad\quad
g^2:=
 \sum_{i,j}\Theta^{\zeta\,\dagger}_{i\,j}(0,z)\cdot
 f_j(z)\cdot v_i(\zeta,z).
\]
We put $\rho:=\chi\cdot(f_1-f_2)$.
Then the claims \ref{8.17.12} and \ref{8.17.13} are clear.
Let us show the claims \ref{8.17.10}  and \ref{8.17.11}.

\begin{lem}
Let $\epsilon$ be a real number such that
$0<\epsilon<1/2$,
and $N$ be any real number.
Then we have
$\rho\in L^2_{-\epsilon\cdot\vecdelta,N}$.
\end{lem}
\pf
It follows from (\ref{eq;9.9.60}) and
Lemma \ref{lem;9.8.15}.
\hfill\qed

\vspace{.1in}

We have the following formula:
\begin{multline}
\delbar_{\nbigelambda}f^1=
 \sum_i\Bigl(
  \delbar f_i+
 \sum_j
 \Bigl(
  \sum_k
 \bigl(
 \Theta^{k\,\dagger}_{i\,j}(\zeta,z)+\frac{b_k(v_i)}{2}\delta_{i\,j}
 \bigr)\cdot\bar{\eta}_k
+\Theta^{\zeta\,\dagger}_{i\,j}(\zeta,z) \cdot d\bar{\zeta}
 \Bigr)\cdot f_j
 \Bigr)\cdot v_i(\zeta,z)\\
=\sum_{i,k,j}
 \Bigl(
  \Theta^{k\,\dagger}_{i\,j}(\zeta,z)
 -\Theta^{k\,\dagger}_{i\,j}(0,z)
 \Bigr)\cdot\bar{\eta}_k\cdot f_j\cdot v_i(\zeta,z)
+\sum_{i,j}
  \Theta^{\zeta\,\dagger}_{i\,j}(\zeta,z)\cdot d\overline{\zeta}
 \cdot f_j\cdot v_i(\zeta,z).
\end{multline}
We also have the following:
\[
 \delbar_{\nbigelambda}f^2=
\sum_{i,j}
 \Theta^{\zeta\,\dagger}_{i\,j}(0,z)
\cdot f_j(z)\cdot v_i(\zeta,z)\cdot d\bar{\zeta}
+\bar{\zeta}\cdot\delbar g^2.
\]
Hence we obtain the following:
\begin{multline}
\delbar_{\nbigelambda}\bigl(
 \chi\cdot (f_1-f_2)\bigr)=
 \delbar\chi\cdot(f_1-f_2)+\chi\cdot\delbar(f_1-f_2)
=\delbar\chi\cdot(f_1-f_2)\\
+\chi\times
\Bigl(
 \sum_{i,k,j}
 \bigl(\Theta^{k\,\dagger}_{i\,j}(\zeta,z)-\Theta^{k\,\dagger}_{i\,j}(0,z)
 \bigr)\cdot
\bar{\eta}_k\cdot f_j\cdot v_i(\zeta,z)
+\sum_{i\,j}
 \bigl(\Theta^{\zeta\,\dagger}_{i\,j}(\zeta,z)
 -\Theta^{\zeta\,\dagger}_{i\,j}(0,z)\bigr)\cdot f_j(z)
  \cdot v_i(\zeta,z)\cdot d\bar{\zeta}
+\bar{\zeta}\cdot \delbar g^2
\Bigr).
\end{multline}

\begin{lem}
$\delbar\chi\cdot(f^1-f^2)$ is $L^2$-section
of $\nbigelambda(-X_0)\otimes\Omega^{0,1}$.
\end{lem}
\pf
Note that $\delbar \chi$ vanishes on $|\zeta|<1/3$.
Then the claim follows from Lemma \ref{lem;9.8.15}
and Lemma \ref{lem;10.11.41}.
\hfill\qed

\begin{lem}
The following is $L^2$-section of
$\nbigelambda\otimes(-X_0)\otimes\Omega^{0,1}$:
\[
 \chi\times
\Bigl(
 \sum_{i,k,j}
 \bigl(\Theta^{k\,\dagger}_{i\,j}(\zeta,z)-\Theta^{k\,\dagger}_{i\,j}(0,z)
 \bigr)\cdot
\bar{\eta}_k\cdot f_j\cdot v_i(\zeta,z)
+\sum_{i\,j}
 \bigl(\Theta^{\zeta\,\dagger}_{i\,j}(\zeta,z)
 -\Theta^{\zeta\,\dagger}_{i\,j}(0,z)\bigr)\cdot f_j(z)
  \cdot v_i(\zeta,z)\cdot d\bar{\zeta}
\Bigr)
\]
\end{lem}
\pf
We put $K_0=\bigl\{\zeta\,\big|\,|\zeta|<2/3\bigr\}$.
On $K_0\times (X_0-D_0)$,
we have the following for $a=1,\ldots,n-1$, or $\zeta$:
\[
 \Bigl|
 \Theta^{a\,\dagger}_{i\,j}(\zeta,z)
 -\Theta^{a\,\dagger}_{i\,j}(0,z)
\Bigr|
\leq
 \bigl|\zeta\bigr|\times
 \bigl|\bigl|
 \Theta^{a\,\dagger}_{i\,j}(\cdot,z)
 \bigr|\bigr|_{C^1(K_0\times\{z\})}.
\]
Here $\bigl|\bigl|
 \Theta^{a\,\dagger}_{i\,j}(\cdot,z)
 \bigr|\bigr|_{C^1(K_0\times\{z\})} $
denotes the $C^1$-norm
of the restriction of $\Theta^{a\,\dagger}_{i\,j}$
to $K_0\times\{z\}$.
Then we obtain the result from Lemma \ref{lem;9.9.62}.
\hfill\qed

\begin{lem}
$\chi\cdot\bar{\zeta}\cdot \delbar g^2$ is $L^2$-section
of $\nbige(-X_0)\times \Omega^{0,1}$.
\end{lem}
\pf
We have the following:
\begin{multline}
\chi\cdot\bar{\zeta}\cdot \delbar g^2
=\chi\cdot \bar{\zeta}\cdot(\lambda-\lambda_0)\times\\
 \sum_{i\,j}
\Bigl(
\delbar \Theta_{i\,j}^{\zeta\,\dagger}(0,z)
 \cdot f_j(z)\cdot v_i'(\zeta,z)
+ \Theta_{i\,j}^{\zeta\,\dagger}(0,z)
 \delbar f_j(z) \cdot v_i(\zeta,z)
+\Theta_{i\,j}^{\zeta\,\dagger}(0,z)
 \cdot f_j(z)\cdot \delbar v_i(\zeta,z)
\Bigr).
 \end{multline}
The first summand can be dominated due to
Lemma \ref{lem;9.9.71}.
The second summand can be dominated due to (\ref{eq;9.9.70}).
The third summand can be dominated due to (\ref{eq;9.9.72}).
\hfill\qed

\vspace{.1in}

Hence the claims \ref{8.17.11} and \ref{8.17.10}
are obtained,
and the proof of Proposition \ref{prop;9.9.85}
is accomplished.
\hfill\qed

%% file: 12.1.tex

\subsubsection{Extension of holomorphic sections}

\label{subsubsection;10.11.80}

We put $X=\Delta^n$ and $D=\bigcup_{i=1}^l D_i$.
We put as follows:
\[
 X^{(1)}:=\bigl\{(z_1,\ldots,z_n)\in X\,\big|\,z_1=z_2\bigr\},
\quad
 D^{(1)}:=X^{{(1)}}\cap D,
\quad
 X_0:=\bigl\{(z_1,\ldots,z_n)\in X\,\big|\,z_1=z_2=0\bigr\}.
\]

\begin{condition} \label{condition;10.11.45}
Let $\harmonicbundle$ be a tame harmonic bundle over $X-D$.
Let $\epsilon_1$ and $\epsilon_2$ be positive numbers
such that $\epsilon_1+\epsilon_2<1$.
Assume that $\KMS\bigl(\prolong{\nbigelambda},i\bigr)$
is $\epsilon_i$-small $(i=1,2)$.
\hfill\qed
\end{condition}

We will show the following
in the subsubsection \ref{subsubsection;d11.14.40}.
\begin{prop} \label{prop;9.9.81}
Let $\harmonicbundle$ be a tame harmonic bundle
satisfying Condition {\rm\ref{condition;10.11.45}}.
Let $f$ be a holomorphic section of
$\prolong{(\nbigelambdazero_{X^{(1)}-D^{(1)}})}$
on $X^{(1)}$.
Then there exists a neighbourhood of $X_0$
in $X$,
and there exists a holomorphic section
$\tilde{f}\in \Gamma(U,\prolong{\nbigelambdazero})$,
such that $\tilde{f}_{X^{(1)}\cap U}=f_{|X^{(1)}\cap U}$.
\end{prop}

Let $\varphi:\widetilde{\Delta}^2_z\lrarr \Delta^2$
be a blow up.
We put $\widetilde{X}:=\widetilde{\Delta}^2_z\times\Delta^{n-2}_w$.
Let $\psi$ denote the composite of the following morphisms:
\[
 \widetilde{X}\stackrel{\varphi\times id}\lrarr
 \Delta_z^2\times\Delta_w^{n-2}
\lrarr \Delta^n_z.
\]
Here the second morphism is given by
$z_i=z_i$ $(i=1,2)$ and $w_i=z_{i-2}$ $(i=3,\ldots,n)$.
We put as follows:
\[
 \widetilde{D}:=
  \psi^{-1}(D)=
 \Bigl[
 \bigl(\varphi^{-1}(0,0)\cup\widetilde{D}'_1\cup \widetilde{D}'_2\bigr)
\times\Delta_w^{n-2}
\Bigr]
\cup \Bigl[
 \Delta_z^2\times
 \Bigl(
  \bigcup_{i=1}^{n-2}\bigl\{w_i=0\bigr\}
 \Bigr)\Bigr].
\]
Then the restriction of $\psi$ to $\widetilde{X}-\widetilde{D}$
gives the isomorphism of $\widetilde{X}-\widetilde{D}$
and $X-D$.

We denote the closure of
$\psi^{-1}(X^{(1)}-D^{(1)})$ by $\widetilde{X}^{(1)}$.
We put $\widetilde{D}^{(1)}:=\psi^{-1}(D)\cap \widetilde{X}^{(1)}$.
We put $\widetilde{X}_0=\psi^{-1}(X_0)$.
The restriction of $\psi$ to $\widetilde{X}^{(1)}$
gives an isomorphism of $\widetilde{X}^{(1)}$ and $X^{(1)}$.

\begin{lem}\label{lem;9.9.80}
Let $f$ be a holomorphic section of
$\prolong{\psi^{\ast}\nbigelambda}$ over $\widetilde{X}^{(1)}$.
Then there exists a neighbourhood $\widetilde{U}$ of $\widetilde{X}_0$
in $\widetilde{X}$,
and there exists a holomorphic section
$\tilde{f}\in \Gamma(\widetilde{U},\prolong{\psi^{\ast}\nbigelambda})$,
such that
$\tilde{f}_{|\tilde{X}^{(1)}\cap \tilde{U}}=
 f_{|\tilde{X}^{(1)}\cap \tilde{U}}$.
\end{lem}

It is easy to see that Lemma \ref{lem;9.9.80} implies
Proposition \ref{prop;9.9.81}.

\label{subsubsection;9.9.82}

\subsubsection{Proof}
\label{subsubsection;d11.14.40}

We essentially use the argument in the subsubsection 4.7.4
in \cite{mochi}.
To apply Corollary \ref{cor;10.11.56},
recall the setting of the subsubsections
\ref{subsubsection;10.11.70}--\ref{subsubsection;10.11.71}.
Let consider the blow up
$\widetilde{\Delta^2_z}\lrarr\Delta_z^2=\{(z_1,z_2)\}$ at $O=(0,0)$
as in the subsubsection \ref{subsubsection;9.9.82}.
We can take a holomorphic embedding $\iota$
of $Y$, given in the previous subsubsection,
to $\widetilde{\Delta}^2$
satisfying the following:
\begin{itemize}
\item The image of the $0$-section $\proj^1$ is the exceptional divisor
 $\phi^{-1}(O)$.
\item
We have $\iota^{-1}(D_1')=\pi^{-1}(\infty)$
and $\iota^{-1}(D_2')=\pi^{-1}(0)$.
\end{itemize}
We may assume that $\pi^{-1}(P)=\iota^{-1}\bigl(\blowup{C}(1,1)\bigr)$
for $P=[1:1]\in \proj^1$.

We put $\overline{X}:=Y\times\Delta_w^{n-2}$.
Then we have the naturally induced morphism
$\overline{X}\lrarr \blowup{\Delta_z^2}\times\Delta_w^{n-2}$,
which we also denote by $\iota$.
We put as follows:
\[
 \overline{D}:=\iota^{-1}(\blowup{D}),\quad
 \overline{X^{(1)}}:=\iota^{-1}(\blowup{X^{(1)}})=\pi^{-1}(P)\times
 \Delta_w^{n-2},\quad
 \overline{D^{(1)}}:=\overline{X^{(1)}}\cap \overline{D}.
\]
Note that $\iota(\overline{X})$ gives a neighborhood
of $\blowup{X_0}$ in $\blowup{X}$.
The composite $\psi\circ\iota$ is denoted by $\psi_1$.

Let $\harmonicbundle$ be a tame harmonic bundle over $X-D$
satisfying Condition \ref{condition;10.11.45}.
Recall that the holomorphic bundle $\nbigelambda$
with the hermitian metric $h$ is acceptable
(Corollary \ref{cor;10.11.80}).
We take 
the numbers $\epsilon$, $a$ and $b$ be as in Lemma \ref{lem;10.11.55},
and we take the metric $\tilde{h}$ of
$\psi_1^{-1}\nbigelambda$ as in (\ref{eq;10.11.75}).

Let $f$ be a holomorphic section of
$\prolong{(\nbige^{\lambda}_{|X_0-D_0})}$, or equivalently,
$\prolong{\psi^{\ast}\nbige^{\lambda}_{|\blowup{X}_0-\blowup{D}_0}}$.
Clearly Lemma \ref{lem;9.9.80} can be reduced to the following lemma.
\begin{lem}\label{lem;5.22.22}
There exists a holomorphic section $\tilde{f}$
of $\prolong{\psi_1^{\ast}\nbige^{\lambda} }$
over $\overline{X^{(1)}}$ such that
$\tilde{f}_{|\overline{X}^{(1)}}=f_{|\overline{X}^{(1)}}$.
\end{lem}
\pf
Take an embedding $\kappa:\Delta_{\zeta}\lrarr \proj^1-\{0,\infty\}$
such that $\kappa(0)=P$.
By using Proposition \ref{prop;9.9.85},
we can take a $C^{\infty}$-function $\rho$ whose support
is contained in $\pi^{-1}(\kappa(\Delta_{\zeta}))\times\Delta_w^{n-2}$,
and satisfying the following:
\begin{itemize}
\item $\rho$ is an element of
 $A^{0,0}_{\tilde{h}}\bigl(\psi_1^{\ast}\nbige^{\lambda}\bigr)$.
\item
 $\delbar_E(\rho)$ is an element of
 $A^{0,1}_{\tilde{h}}
\bigl(\psi_1^{\ast}\nbige^{\lambda}(-\overline{X^{(1)}})\bigr)$.
\item
 We have $\rho_{|\overline{X}^{(1)}-\overline{D}^{(1)}}=f$.
\end{itemize}
Here $\psi_1^{\ast}\nbige^0(-\overline{X^{(1)}})$
is same as $\psi_1^{\ast}\nbige^0\otimes\nbigo_{\proj^1}(-1)$.

Due to Corollary \ref{cor;10.11.56},
we can pick an element $G$ of 
$A^{0,0}_{\tilde{h}}
 \bigl(\psi_1^{\ast}\nbige^0(-\overline{X^{(1)}})\bigr)$
such that $\delbar G=\delbar\rho$.
Then we put $\tilde{f}:=\rho-G$.
Then it satisfies 
$\delbar\tilde{f}=0$,
$\tilde{f}\in A_{\tilde{h}}^{0,0}(\psi_1^{\ast}\nbige^0)$,
and
$\tilde{f}_{|\overline{X}^{(1)}-\overline{D}^{(1)}}=
 f_{|\overline{X}^{(1)}-\overline{D}^{(1)}}$.

Let us check that $\tilde{f}$ gives a section of
$\prolong{\psi_1^{\ast}\nbige^{\lambda}}$.
We regard $\tilde{f}$ as a section of $\nbigelambda$
over an open subset $\psi_1(\overline{X}-\overline{D})$
of $X-D$.
We have only to check that
$\tilde{f}$ gives a section of $\prolong{\nbige^{\lambda}}$
over $\psi_1(\overline{X})$.
To show it,
let us consider the norm of $\tilde{f}$
with respect to our original metric $h$.
We consider the restriction of  $\tilde{f}$ to
$\pi_j^{-1}(P)$, for any point $P\in D_j^{\circ}$
and for $1\leq j\leq l$.
First we consider the case $j=1$.
On the curve $\pi_1^{-1}(P)$,
the metric $\tilde{h}$ is equivalent
to $h\cdot |z_1|^{a-\epsilon}\times Q$,
where $Q$ denotes a polynomial of $\log|z_1|$.
Thus we obtain
$-\ord\bigl(\tilde{f}_{|\pi_1^{-1}(P)}\bigr)\leq -a+\epsilon<1-\epsilon_1$,
due to our choice of $\epsilon,a,b$
(Lemma \ref{lem;10.11.55}).
Since $\KMS(\prolong{\nbigelambda},1)$ is $\epsilon_1$-small
(Condition \ref{condition;10.11.45}),
the intersection of
the sets $\Par(\nbige^0,1)$ and
$\bigl\{c\,\big|\,0\leq c<1-\epsilon_1\bigr\}$ is empty.
Therefore we can conclude that
$-\ord\bigl(\tilde{f}_{|\pi_1^{-1}(P)}\bigr)\leq 0$
with respect to the original metric $h$.
Similarly we can show that
$-\ord\bigl(\tilde{f}_{|\pi_2^{-1}(P)}\bigr)\leq 0$
for any point $P\in D_2^{\circ}$.
If $j>2$, then the metrics $h$ and $\tilde{h}$
are equivalent on the curve $\pi_j^{-1}(P)$.
Thus we obtain
$-\ord\bigl(\tilde{f}_{|\pi_2^{-1}(P)}\bigr)\leq 0$
with respect to the metric $h$,
also in this case.
Thus we obtain that $\tilde{f}$ is, in fact,
a section of $\prolong{\psi_1^{\ast}\nbige^0}$,
due to Corollary \ref{cor;11.28.15}.
\hfill\qed

%% file: a73.1.tex

\subsubsection {Extension property in the codimension one case}

\label{subsubsection;c11.16.20}

We put $X:=\Delta_z^n\times\Delta_w$,
$D_i:=\{z_i=0\}$ and $D=\bigcup_{i=1}^l D_i$.
We put $X^{(2)}:=\Delta_z^n\times\{0\}$
and $D^{(2)}=D\cap X^{(2)}$.
We have the origin $(O,0)\in X^{(2)}\subset X$.
Let $\harmonicbundle$ be any tame harmonic bundle over $X-D$.
Let us consider the restriction of $\nbigelambda$
to $X^{(2)}-D^{(2)}$.

\begin{lem}\label{lem;a11.15.61}
Let $f$ be a holomorphic section of
$\prolongg{\vecb}{\bigl(\nbigelambda_{|X^{(2)}-D^{(2)}}\bigr)}$
defined on a neighbourhood of $(O,0)$ in $X^{(2)}$.
Then there exists a holomorphic section $\tilde{f}$
of $\prolongg{\vecb}{\nbigelambda}$ defined on a neighbourhood
of $(O,0)$, which satisfies
$\tilde{f}_{|X^{(2)}}=f$.
\end{lem}
\pf
The argument is essentially same as the proof of
Lemma \ref{lem;5.22.22}.
In fact,
we can show the claim more simply,
by using the results in the subsubsections
\ref{subsubsection;d11.14.31}--\ref{subsubsection;a11.15.60}
and the vanishing in Lemma \ref{lem;9.10.72}.
\hfill\qed

\subsubsection{Local freeness in codimension one}
We put $X=\Delta\times \Delta^{n}_w$,
$D=\{0\}\times \Delta^n_w$.
Let $\harmonicbundle$ be any tame harmonic bundle on $X-D$.
Let $\pi$ denote the projection of $X$ to $D$.
Let $P$ be a point of $D$.
Then we obtain the smooth curve $\pi^{-1}(P)$.

\begin{cor} \label{cor;a11.15.62}
Let $b$ be any real number.
Let $f$ be a holomorphic section of
$\prolongg{b}{\bigl(\nbige^{\lambda}_{|\pi^{-1}(P)}\bigr)}$
defined on a neighbourhood of $(0,P)$ in $\pi^{-1}(P)$.
Then there exists a holomorphic section $\tilde{f}$ of
$\prolongg{b}{\nbigl(\nbigelambda)}$ defined 
on a neighbourhood of $(0,P)$ in $X$,
such that $\tilde{f}_{|\pi^{-1}(P)}=f$.
\end{cor}
\pf
We have only to use Lemma \ref{lem;a11.15.61} inductively.
\hfill\qed

\vspace{.1in}

The following corollary will be used without mention.
\begin{cor} \label{cor;a11.15.65}
The sheaf $\prolongg{b}{\nbigelambda}$ is locally free.
The restriction 
$\prolongg{b}{\nbigelambda}\lrarr
 \prolongg{b}{\nbigelambda}_{|\pi^{-1}(P)}$ is surjective.
\end{cor}
\pf
The second claim is shown in Corollary \ref{cor;a11.15.62}.
Let us show the first claim.
We have only to prove the case $b=0$.
We have the set $\Par(\prolong{\nbigelambda},1)$.
Let $\vecv$ be a frame of
$\prolong{\nbige^{\lambda}_{|\pi_i^{-1}(P)}}$,
which is compatible with the parabolic filtration $F$.
For each $v_i$, we have the number $b_i:=\deg^{F}(v_i)$.

Then we can pick sections $\tilde{v}_i$
of $\prolongg{b_i}{\nbigelambda}$ such that
$\tilde{v}_{i\,|\,\pi_i^{-1}(P)}=v_i$,
by using Corollary \ref{cor;a11.15.62}.
Thus we obtain the tuple $\tilde{\vecv}:=(\tilde{v}_i)$
of sections of $\prolong{\nbigelambda}$.
We would like to show that $\tilde{\vecv}$ gives
a frame of $\prolong{\nbigelambda}$.

We put $\tilde{b}:=
\sum_{b\in \Par(\nbigelambda,1)}\multiplicity(\lambda,b)\cdot b$.
Then we have the natural isomorphism
$\prolongg{\tilde{b}}{\bigl(\det\nbigelambda\bigr)}_{|\pi^{-1}(P)}
\simeq \det\bigl(\prolong{\nbige^{\lambda}_{|\pi^{-1}(P)}}\bigr)$.
Let us consider the section
$\Omega(\tilde{\vecv}):=
\tilde{v}_1\wedge\cdots\wedge \tilde{v}_{\rank E}$.
Since $\vecv$ is a frame of $\prolong{\nbige^{\lambda}_{|\pi^{-1}(P)}}$,
we have $\Omega(\tilde{\vecv})_{|(0,P)}\neq 0$.
Hence $\Omega(\tilde{\vecv})$ gives a frame of
$\prolongg{\tilde{b}}{\bigl(\det\nbigelambda\bigr)}_{|\pi^{-1}(P)}$
around $(0,P)$.
Then we can conclude that $\tilde{\vecv}$
gives a frame of $\prolong{\nbigelambda}$
on a neighbourhood of $(0,P)$,
due to the last argument in the proof of Lemma \ref{lem;d11.14.51}.
\hfill\qed

%% file: 12.2.tex

\subsubsection{Preliminary}

Assume that $\Par(\prolong{\nbigelambda},i)$ are $\epsilon_i$-small
and $\sum_i \epsilon_i<1$.
Pick an element
$(\lefttop{0}z_3,\ldots,\lefttop{0}z_{n})
 \in(\Delta^{\ast})^{n-2}$.
Let us consider the curve
$C_0:=\{(z,z,\lefttop{0}z_3,\ldots,\lefttop{0}z_n)\in X\}$.
We have the restriction of
$\harmonicbundle$ to $C_0$.

\begin{lem}  \label{lem;d11.14.50}
We have the following claims
for the KMS-structure of  $\prolong{\nbige}_{|C_0}$.
\begin{itemize}
\item
 $\Par(\prolong{\nbigelambda}_{|C_0})$ is
 $(\epsilon_1+\epsilon_2)$-small.
\item
 We have the following equality:
\[
 \sum_{b\in \Par(\prolong{\nbigelambda_{|C_0}})}
 b\cdot \multiplicity(\lambda,b)
=\sum_{i=1,2}
 \sum_{b_i\in \Par(\prolong{\nbigelambda},i)}
 b_i\cdot \multiplicity(\lambda,b_i).
\]
\end{itemize}
\end{lem}
\pf
If $\epsilon>0$ is sufficiently small,
we have the following:
\[
 \Par(\leftbottom{\epsilon\cdot\vecdelta}\nbigelambda,i)
=\Par(\prolong{\nbigelambda},i),
\quad\quad
  \KMS(\leftbottom{\epsilon\cdot\vecdelta}\nbigelambda,i)
=\KMS(\prolong{\nbigelambda},i).
\] 

We use the following lemma.
\begin{lem}
If $\eta$ is sufficiently small,
the set
$\KMS(\leftbottom{\epsilon\cdot\vecdelta}\nbige^{\lambda'},i)$
depends on $\lambda'\in\Delta(\lambda,\eta)$ continuously.
\end{lem}
\pf
We put as follows:
$ \nbigk(\nbige^0,\lambda,0,i):=\bigl\{
 u\in\KMS(\nbige^0,i)\,\big|\, 
 \kmsmap(\lambda,u)\in\KMS(\prolong{\nbigelambda})
 \bigr\}$.
Then we have the following
\[
 \KMS(\prolongg{\epsilon\cdot\vecdelta}{\nbige^{\lambda'}},i)
=
 \bigl\{
 \kmsmap(\lambda',u)\,\big|\,
 u\in \nbigk(\nbige^0,\lambda,0,i)
 \bigr\}.
\]
Thus we are done.
\hfill\qed

\vspace{.1in}

Let $\lambda'\in\Delta(\lambda,\eta)$ be generic.
Due to Lemma \ref{lem;a11.15.30}, the following holds:
\[
 \KMS\bigl(\prolongg{2\epsilon}{(\nbige^{\lambda'}_{|C_0})}\bigr)
\subset
 \bigl\{u_1+u_2\,\big|\,
  u_i\in \KMS(\prolongg{\epsilon\cdot\vecdelta}{\nbige^{\lambda'}},i)
 \bigr\}
\]
Then the following holds
for any $\lambda'\in \Delta(\lambda,\eta)$:
\[
 \KMS\bigl(\prolongg{2\epsilon}{(\nbige^{\lambda'}_{|C_0})}\bigr)
 \subset
 \bigl\{u_1+u_2\,\big|\,
  u_i\in \KMS(\prolongg{\epsilon\cdot\vecdelta}{\nbige^{\lambda'}},i)
 \bigr\},
\quad
 \Par\bigl(\prolongg{2\epsilon}{(\nbige^{\lambda'}_{|C_0})}\bigr)
\subset
 \bigl\{a_1+a_2\,\big|\,
 a_i\in\Par(\prolongg{\epsilon\cdot\vecdelta}{\nbige^{\lambda'}},i)
 \bigr\}.
\]
In particular, we obtain the following:
\[
 \Par\bigl(\prolongg{2\epsilon}{(\nbige^{\lambda}_{|C_0})}\bigr)
\subset
 \bigl\{a_1+a_2\,\big|\,
 a_i\in\Par(\prolongg{\epsilon\cdot\vecdelta}{\nbige^{\lambda}},i)
 \bigr\}
=\bigl\{a_1+a_2\,\big|\,
 a_i\in\Par(\prolong{\nbige^{\lambda}},i)
 \bigr\}.
\]
Hence we obtain the first claim.
We can show the second claim similarly
by using Lemma \ref{lem;a11.15.30}.
\hfill\qed

%% file: 12.3.tex

\subsubsection{Preliminary prolongation}

\begin{lem} \label{lem;9.9.95}
Let $\epsilon_1,\ldots,\epsilon_l$ be
positive numbers such that $\sum_{i=1}^l \epsilon_i<1$.
Assume that $\Par(\prolong{\nbigelambda},i)$ is $\epsilon_i$-small.
Then the $\nbigo_X$-sheaf $\prolong{\nbigelambda}$ is locally free.
\end{lem}
\pf
We use an induction on the dimension of $X$.
As the hypothesis of the induction,
we assume the following:
\begin{quote}
The $\nbigo_X$-sheaf $\prolong{\nbigelambda}$ is locally free,
if the following holds:
\begin{itemize}
\item
 $\dim(X)\leq n-1$.
\item
$\Par(\prolong{\nbigelambda},i)$ is $\epsilon_i$-small.
\item
$\sum \epsilon_i<1$.
\end{itemize}
\end{quote}

We use the notation in the subsubsection \ref{subsubsection;10.11.80}.
Due to the hypothesis of the induction
and the first claim of Lemma \ref{lem;d11.14.50},
we have the local freeness of
$\prolong{(\nbigelambda_{|X^{(1)}-D^{(1)}})}$.
Pick a frame $\vecv=(v_i)$ of $\prolong{(\nbigelambda_{|X^{(1)}-D^{(1)}})}$
over $X^{(1)}$.
For each $v_i$,
we pick a section $\tilde{v}_i$ of $\prolong{\nbigelambda}$ over 
a neighbourhood $U$ of $X_0$
such that $\tilde{v}_{i\,|\,U\cap X^{(1)}}=v_{i\,|\,U\cap X^{(1)}}$.
We may assume that $v_i$ are defined over $X$.
Clearly Lemma \ref{lem;9.9.95}
can be reduced to the following lemma.

\begin{lem} \label{lem;10.11.81}
$\tilde{\vecv}$ gives a frame of $\prolong{\nbigelambda}$
around $X_0$.
\end{lem}
\pf
We put as follows:
\[
 \tilde{\vecb}:=(\tilde{b}_1,\ldots,\tilde{b}_l),
\quad
 \tilde{b}_i:=
 \sum_{b\in \Par(\prolong{\nbigelambda},i)}b\cdot \multiplicity(\lambda,b).
\]

The restriction $\tilde{\vecv}_{|\pi_i^{-1}(P)}$
gives a tuple of holomorphic sections of
$\prolong{(\nbigelambda_{|\pi_i^{-1}(P)})}$
for any $P\in D_i^{\circ}$.
Hence $\Omega(\tilde{\vecv})_{|\pi_i^{-1}(P)}$ is a holomorphic section
of $\det\bigl(\prolong{(\nbigelambda_{|\pi_i^{-1}(P)})} \bigr)=
 \prolongg{\tilde{b}_i}{\bigl(\det(\nbigelambda)_{|\pi_i^{-1}(P)}\bigr)}$.
It implies $\Omega(\tilde{\vecv})$ is a holomorphic section
of $\prolongg{\tilde{\vecb}}{\det(\nbigelambda)}$.

We have the natural isomorphism
$(\prolongg{\tilde{\vecb}}{\det(\nbigelambda)})_{|X^{(1)}}
\simeq
 \det\bigl(\prolong{(\nbigelambda_{|X^{(1)}})}\bigr)$.
Since $\Omega(\tilde{\vecv})_{|X^{(1)}}$
gives a frame of
$ \det\bigl(\prolong{(\nbigelambda_{|X^{(1)}})}\bigr)$,
we obtain $\Omega(\tilde{\vecv})_{|O}\neq 0$.
It is standard argument to conclude that $\vecv$ gives a frame
around the origin $O$.
(See the proof of Lemma \ref{lem;d11.14.51}).
Hence we obtain Lemma \ref{lem;10.11.81},
and thus Lemma \ref{lem;9.9.95}.
\hfill\qed

%% file: 14.tex

\subsubsection{Statements of the theorems}

In this subsection,
we will show the following theorems.

\begin{thm} \label{thm;9.10.2}
For any $\vecb\in\real^l$,
the $\nbigo_{X}$-module
$\prolongg{\vecb}{\nbigelambda}$
are coherent and locally free.
In particular, $\prolong{\nbigelambda}$ is locally free.
\end{thm}

Let us pick an element $\vecb=(b_1,\ldots,b_l)\in\real^l$.
Let $\vecdelta_i$ denote the element
$(\overbrace{0,\ldots,0}^{i-1},1,0,\ldots,0)$.
For any $\vecb'\leq \vecb$,
we have the naturally defined morphism
$\prolongg{\vecb'}{\nbigelambda}\lrarr 
 \prolongg{\vecb}{\nbigelambda}$.
For $b_i-1\leq b\leq b_i$, we put $\vecb'=\vecb+(b-b_i)\vecdelta_i$
and as follows:
\[
 \lefttop{i}F_{b}\bigl(\prolongg{\vecb}{\nbigelambda}\bigr)
:=
\Image\bigl(
 \prolongg{\vecb'}{\nbigelambda}_{|D_i}
\lrarr
 \prolongg{\vecb}{\nbigelambda}_{|D_i}
 \bigr).
\]
Then we obtain the filtration
$\lefttop{i}F(\prolongg{\vecb}{\nbigelambda}):=
 \bigl\{\lefttop{i}F_b(\prolongg{\vecb}{\nbigelambda})
 \,\big|\,b_i-1\leq b\leq b_i\bigr\}$
of $\nbigo_{D_i}$-modules.

\begin{thm} \label{thm;9.10.3}\mbox{{}}
\begin{itemize}
\item
$\lefttop{i}F(\prolongg{\vecb}{\nbigelambda})$ is a filtration
in the category of vector bundles on $D_i$.
\item
The tuple of the filtrations
$\big(\lefttop{i}F\,\big|\,i=1,\ldots,l\big)$ 
on the divisors are compatible.
(Definition {\rm\ref{df;10.11.90}}).
\item
Let $\veceta$ be an element of $\prod_{i\in I} [b_i-1,b_i]$.
We put $\tilde{\veceta}:=\vecb+\veceta-q_I(\vecb)$.
Then we have the following:
\[
 \lefttop{I}F_{\veceta}\bigl(\prolongg{\vecb}{\nbigelambda}\bigr):=
 \bigcap_{i\in I}\lefttop{i}F_{\eta_i\,|\,D_I}
 \bigl(\prolongg{\vecb}{\nbigelambda}\bigr)
=\Image(
 \prolongg{\tilde{\veceta}}{\nbigelambda}_{|D_I}
\lrarr
 \prolongg{\vecb}{\nbigelambda}_{|D_I}
 ).
\]
\end{itemize}
\end{thm}

Before entering the proof, we remark the following.
\begin{lem} \label{lem;9.10.21}
We have only to show the local freeness
of $\prolong{\nbigelambda}$ for any tame harmonic bundle
$(E,\delbar_E,h,\theta)$
to prove Theorem {\rm\ref{thm;9.10.2}} and Theorem {\rm\ref{thm;9.10.3}}.
\end{lem}
\pf
For a tame harmonic bundle $(E,\delbar_{E},h,\theta)$,
and $\vecb\in \real^l$,
we have the harmonic bundle
$(E',\delbar_{E'},h',\theta')
=(E,\delbar_E,h,\theta)\otimes L(-\vecb)$.
We denote the deformed holomorphic bundle
of $(E',\delbar_{E'},h',\theta')$ by $\nbige^{\prime\,\lambda}$.
Then we have the natural isomorphism
$\prolongg{\vecb}{\nbigelambda}\simeq
 \prolong{\nbige^{\prime\,\lambda}}$
by definition.
Thus we are done.
\hfill\qed

\vspace{.1in}
We also remark that we use Corollary \ref{cor;a11.15.65}
without mention.

\subsubsection{Step 1}

\label{subsubsection;c11.16.1}

\begin{condition} \label{condition;9.10.6}
Let us take an element $\vecc\in\seisuu_{>0}^l$
as follows for any $i$:
\begin{enumerate}
\item
 $c_i$ are sufficiently large with respect to
 $\KMS(\prolong{\nbigelambda},i)$.
\item
 There exists a number $b_i\in ]-1,0]$
 such that
 $\{-b_i+\kappa(c_i\cdot a_i)\,|\,
  a_i\in\Par(\nbigelambda,i) \}\subset ]-1,0]$
and that it is  $(2l)^{-1}$-small.
We may also assume that $b_i-(2l)^{-1}>-1$.
We put $\vecb=(b_1,\ldots,b_l)$.
\end{enumerate}
\end{condition}

We put
$(E_1,\delbar_{E_1},h_1,\theta_1):=
 \bigl(
 \psi_{\vecc}^{-1}(E,\delbar_E,h,\theta)
 \bigr)
 \otimes L(-\vecb)$.
We denote the deformed holomorphic bundle
of $(E_1,\delbar_{E_1},h_1,\theta_1)$
by $\nbigelambda_1$.
Then we have the natural isomorphism
$\prolong{\bigl(
 \psi_{\vecc}^{-1}\nbigelambda\bigr)}
=
\prolongg{\vecb}{\bigl( \psi_{\vecc}^{-1}\nbigelambda\bigr)}
\simeq
 \prolong{\nbigelambda_1}$.
We also have the following:
\begin{equation} \label{eq;9.10.5}
 \KMS(\prolong{\nbigelambda_1},i)
=\bigl\{
 -b_i+\kappa(c_i\cdot a_i)\,|\,
  a_i\in\Par(\nbigelambda,i)
 \bigr\}.
\end{equation}

\begin{lem}
$\prolong{\nbigelambda_1}$ is locally free.
\end{lem}
\pf
Due to (\ref{eq;9.10.5}) and our choice of $\vecc$
(see 2 in Condition \ref{condition;9.10.6}),
the set $\KMS(\prolong{\nbigelambda_1},i)$ are $(2l)^{-1}$-small
for any $i=1,\ldots,l$.
Note that $l\times (2l)^{-1}<1$.
Then we obtain the result by using the preliminary prolongation
(Lemma \ref{lem;9.9.95}).
\hfill\qed

\begin{cor}
The sheaf $\prolong{\bigl(\psi_{\vecc}^{\ast}\nbigelambda\bigr)}$
is a locally free $\nbigo_X$-module.
\hfill\qed
\end{cor}

\vspace{.1in}

We have the natural $\mu_{\vecc}$-action on
$\psi_{\vecc}^{\ast}\nbigelambda$,
which is prolonged to the action
on $\prolong{\bigl(\psi_{\vecc}^{\ast}\nbigelambda\bigr)}$.
In particular, we obtain the $\mu_{c_i}$-action on
$\prolong{\bigl(\psi_{\vecc}^{\ast}\nbigelambda\bigr)}_{|D_i}$.
Since the action of $\mu_{c_i}$ on $D_i$ is trivial,
we have the decomposition:
\[
 \prolong{\psi_{\vecc}^{\ast}\nbigelambda}_{|\,D_i}
=\bigoplus_{c_i-1\leq h\leq 0}V_h.
\]
Here the generator $\omega$ of $\mu_{c_i}$
acts as $\omega^h$ on $V_h$.

Let us pick a point $P$ of $D_i^{\circ}$.
We have the following morphism due to our choice of $\vecc$
(1 in Condition \ref{condition;9.10.6})
and the result in the subsubsection \ref{subsubsection;9.10.7}:
\[
 \varphi:
 \bigl\{h\,\big|\,-c_i+1\leq h\leq 0,\,\,V_h\neq 0\bigr\}
\lrarr
 \bigl\{\tilde{b}\,\big|\,-1<\tilde{b}\leq 0,\,\,
 \Gr^F_{\tilde{b}}
 \bigl(\prolong{\psi_{\vecc}^{\ast}\nbigelambda
 _{|\,\pi_i^{-1}(P)}}\bigr)    
 \neq 0\bigr\}.
\]

We consider the filtration $\lefttop{i}F'_b$
of $\prolong{\bigl(
 \psi_{\vecc}^{\ast}\nbigelambda\bigr)}_{|\,D_i}$
in the category of vector bundles on $D_i$,
given as follows:
\[
 \lefttop{i}F'_{b}:=
 \bigoplus_{\substack{\varphi(h)\leq b}}
 V_h.
\]
Due to the construction,
it is easy to see that the filtrations
$(\lefttop{i}F'\,|\,i=1,\ldots,l)$
are compatible in the sense of 
Definition \ref{df;10.11.90}.

We consider the subsheaf
$\prolongg{b\vecdelta_i}{\bigl(\psi_{\vecc}^{\ast}\nbigelambda\bigr)}'$
of $\prolong{\psi_{\vecc}^{\ast}\bigl( \nbigelambda\bigr)}$,
given as follows:
\[
 \prolongg{b\vecdelta_i}{\bigl(
 \psi_{\vecc}^{\ast}\nbigelambda\bigr)}'
=\Ker\Bigl(\pi:
 \prolong{\bigl(\psi_{\vecc}^{\ast}
 \nbigelambda\bigr)}
\lrarr
 \frac{\prolong{\bigl(
 \psi_{\vecc}^{\ast}\nbigelambda\bigr)}_{|D_i}}
  {\lefttop{i}F'_b}
 \Bigr).
\]
Here $\pi$ denotes the naturally defined morphism.

\begin{lem} \label{lem;9.10.13}
We have the following:
\[
 \prolongg{b\vecdelta_i}{\bigl(\psi_{\vecc}^{\ast}
\nbigelambda\bigr)}'
=\prolongg{b\vecdelta_i}{\bigl(\psi_{\vecc}^{\ast}
 \nbigelambda_1\bigr)},
\quad\quad
 \lefttop{i}F'_b=\lefttop{i}F_b
\]
\end{lem}
\pf
Let $f$ be a holomorphic section
of $\prolongg{b\vecdelta_i}{\bigl(\psi_{\vecc}^{\ast}\nbigelambda\bigr)}$.
It can be also regarded as a section of
$\prolong{\psi_{\vecc}^{\ast}\nbigelambda}$.
Let $P$ be a point of $D_i^{\circ}$.
We have the element $f(P)$ of
$\prolong{\psi_{\vecc}^{\ast}\nbigelambda}_{|P}
 =\prolong{
 \bigl({\psi_{\vecc}^{\ast}\nbigelambda}_{|\pi_i^{-1}(P)}\bigr)}_{|P}$.
Due to Lemma \ref{lem;9.10.11},
we obtain $f(P)\in \lefttop{i}F'_{b\,|\,P}$.

Let $\bar{f}$ denote the image of $f$
via the projection $\pi$.
Then $\bar{f}(P)=0$ for any $P\in D_i^{\circ}$.
It implies $\bar{f}=0$ on $D_i$.
Hence we obtain
$f\in
 \prolongg{b\vecdelta_i}{\bigl(\psi_{\vecc}^{\ast}\nbigelambda\bigr)}'$.

On the other hand, pick a section 
$f\in
 \prolongg{b\vecdelta_i}{\bigl(\psi_{\vecc}^{\ast}\nbigelambda\bigr)}'$.
Due to Lemma \ref{lem;9.10.11},
we obtain the following inequality for any $P\in D_i^{\circ}$:
\[
 -\ord\bigl(f_{|\pi_i^{-1}(P)}\bigr)\leq b.
\]
Then we obtain $f\in
 \prolongg{b\vecdelta_i}{\bigl(\psi_{\vecc}^{\ast}\nbigelambda\bigr)}$
due to Corollary \ref{cor;11.28.15}
and Corollary \ref{cor;10.11.80}.

In all, we obtain
$\prolongg{b\vecdelta_i}{\bigl(\psi_{\vecc}^{\ast}\nbigelambda\bigr)}
=\prolongg{b\vecdelta_i}{\bigl(\psi_{\vecc}^{\ast}\nbigelambda\bigr)}'$.
It implies that
$\prolongg{b\vecdelta_i}{\bigl(\psi_{\vecc}^{\ast}\nbigelambda\bigr)}$
is locally free,
and $\lefttop{i}F_b=\lefttop{i}F'_b$.
\hfill\qed

\begin{lem}
We have the following:
\[
 \lefttop{I}F_{\vecb}
=\Image\bigl(
 \prolongg{\vecb}{{\psi_{\vecc}^{\ast}\nbigelambda}_{|D_I}}
\lrarr
 \prolong{{\psi_{\vecc}^{\ast}\nbigelambda}_{|D_I}}
 \bigr).
\]
\end{lem}
\pf
It can be shown by an argument similar to
the proof of Lemma \ref{lem;9.10.13}.
\hfill\qed

\subsubsection{Step 2}
\label{subsubsection;c11.16.2}

Let $f$ be a section of
$\prolong{\psi_{\vecc}^{\ast}\nbigelambda}$.
Assume the following:
\begin{itemize}
\item
$f(0)\neq 0$,
and
$f$ is compatible with
the filtration $\lefttop{i}F$ $(i=1,\ldots,l)$,
i.e.,
there exists a splitting of $\lefttop{i}F$ $(i=1,\ldots,l)$
which is compatible with $f$.
\item
$f$ is equivariant,
i.e.,
$g^{\ast}(f)=\prod\omega_i^{h_i}\cdot f$
for some $-c_i+1\leq h\leq 0$.
Here $g=(\omega_1,\ldots,\omega_n)\in\mu_{\vecc}$.
\end{itemize}

We put $f_1:=\prod z_i^{-h_i}\cdot f\otimes e$.
Then it is a section of $\psi_{\vecc}^{\ast}\nbigelambda$,
and it is $\mu_{\vecc}$-invariant,
i.e.,
$g^{\ast}(f_1)=f_1$
for any $g\in\mu_{\vecc}$.
Hence there exists the unique section
$\bar{f}$ of $\nbigelambda$ on $X-D$,
such that $\phi_{\vecc}^{\ast}\bar{f}=f_{1\,|\,X-D}$.
Note the following:
\[
 -\lefttop{i}\ord(\bar{f})
=c_i^{-1}(h_i-\lefttop{i}\ord(f))\leq 0.
\]
Hence $\bar{f}$ gives the section of $\prolong{\nbigelambda}$.
If $f(O)\neq 0$, then 
$-\lefttop{i}\ord(f)> -1$,
and thus $-\lefttop{i}\ord(\bar{f})>-1$.

Let us take a frame $\vecv=(v_i)$ of
$\prolong{\phi_{\vecc}^{\ast}\nbigelambda}$
satisfying the following conditions (Corollary \ref{cor;10.11.95}):
\begin{itemize}
\item It is equivariant.
\item It is compatible with the filtrations
 $(\lefttop{1}F,\ldots,\lefttop{l}F)$.
\end{itemize}
Then we obtain a tuple 
$\bar{\vecv}=(\bar{v}_1,\ldots,\bar{v}_r)$ 
of sections of $\prolong{\nbigelambda}$
by the procedure above.

\begin{lem} \label{lem;9.10.20}
$\prolong{\nbigelambda}$ is locally free,
and $\bar{\vecv}$ gives a local frame on a neighbourhood of 
the origin $O$.
\end{lem}
\pf
Recall that
if $\dim(X)=1$ then we have already known the result
(Lemma \ref{lem;9.10.15}).
Let us consider
the element $\tilde{\vecb}=(\tilde{b}_1,\ldots,\tilde{b}_l)$
of $\real^l$, given as follows:
\[
 \tilde{b}_i:=\sum_{b\in\Par(\prolong{\nbigelambda},i)}
 b\cdot\multiplicity(\lambda,b).
\]
Let $P$ be a point of $D_i^{\circ}$.
Then $\Omega(\bar{\vecv})_{|\pi_i^{-1}(P)}$
is a frame of
$\det\bigl(\prolong{\bigl(\nbigelambda_{|\pi_i^{-1}(P)}\bigr)}\bigr)
=\prolongg{\tilde{b}_i}{\det\bigl(\nbigelambda_{|\pi_i^{-1}(P)}\bigr)}$.
Thus we obtain that $\Omega(\bar{\vecv})$ is a frame
of $\prolongg{\vecb}{\det(\nbigelambda)}$.

Let $f=\sum f_i\!\cdot\!\bar{v}_i$ be a holomorphic section of
$\prolong{\nbigelambda}$.
As usual we can show that
$f_i$ are holomorphic over $X$,
and thus $\tilde{\vecv}$ gives a frame of $\prolong{\nbigelambda}$.
(see the proof of Lemma \ref{lem;d11.14.51}).
In particular, 
the sheaf $\prolong{\nbigelambda}$ is locally free.
\hfill\qed

\vspace{.1in}

We consider the filtration 
$\lefttop{i}F'_b$
of $\prolong{\nbigelambda}_{|D_i}$
in the category of the vector bundles over $D_i$,
given as follows:
\[
 \lefttop{i}F'_{b}:=
 \big\langle
 \bar{v}_{j\,|\,D_i}
\,\big|\,-\lefttop{i}\ord(\bar{v}_j)\leq b
 \big\rangle
\]
For any $-1<b\leq 0$,
we consider the subsheaf $\prolongg{b\cdot\vecdelta_i}{\nbigelambda}'$
of $\prolong{\nbigelambda}$ given as follows:
\[
 \prolongg{b\cdot\vecdelta_i}{\bigl(\nbigelambda\bigr)}'
:=
 \ker\Bigl(\pi:
 \prolong{\nbigelambda}\lrarr
 \frac{\prolong{\nbigelambda}_{|D_i}}{\lefttop{i}F'_{b}}
 \Bigr).
\]
Here $\pi$ denote the naturally defined morphism.
Then $\prolongg{b\cdot\vecdelta_i}{\bigl(\nbigelambda\bigr)}'$
is locally free.

\begin{lem} \label{lem;9.10.16}
We have
$\prolongg{b\vecdelta_i}{\nbigelambda}=
 \prolongg{b\vecdelta_i}{\bigl(\nbigelambda\bigr)}'$
and
$\lefttop{i}F'_b=\lefttop{i}F_b$.
\end{lem}
\pf
We have already known that
the claim holds if $\dim X=1$ (Lemma \ref{lem;9.10.15}).

Let $f$ be a holomorphic section of
$\prolongg{b\vecdelta_i}{\nbigelambda}$.
We can also regard it as a section of $\prolong{\nbigelambda}$.
By applying Lemma \ref{lem;9.10.15}
to $f_{|\pi_i^{-1}(P)}\in
   \prolong{\bigl(\nbigelambda_{|\pi_i^{-1}(P)}\bigr)}$,
we obtain that $f(P)\in \lefttop{i}F'_{|P}$
for any $P\in D_i^{\circ}$.
Then it is easy to derive
that $f$ is contained in
$\prolongg{b\vecdelta_i}{\bigl(\nbigelambda\bigr)}'$.

On the other hand,
let $f$ be a holomorphic section of 
$\prolongg{b\vecdelta_i}{\bigl(\nbigelambda\bigr)}'$.
Applying Lemma \ref{lem;9.10.15}
to $f_{|\pi_i^{-1}(P)}$,
we obtain $-\ord(f_{|\pi_i^{-1}(P)})\leq b$.
Then we obtain $f\in\prolongg{b\vecdelta_i}{\nbigelambda}$
due to Corollary \ref{cor;11.28.15}.
Therefore  we obtain
$\prolongg{b\vecdelta_i}{\nbigelambda}
=\prolongg{b\vecdelta_i}{\bigl(\nbigelambda\bigr)}'$,
and thus $\lefttop{i}F_b=\lefttop{i}F'_b$.
\hfill\qed

\begin{lem} \label{lem;9.10.23}
The filtration $\lefttop{i}F$ is a filtration
in the category of the vector bundles over $D_i$.
The filtrations $(\lefttop{i}F\,|\,i=1,\ldots,l)$ are compatible.
\end{lem}
\pf
By our construction,
$\lefttop{i}F'$ is the filtration in the category of
the vector bundles over $D_i$,
and $(\lefttop{i}F'\,|\,i=1,\ldots,l)$ are compatible.
Then the lemma follows from Lemma \ref{lem;9.10.16}.
\hfill\qed

\begin{lem} \label{lem;9.10.22}
We have the following equality:
\[
 \lefttop{I}{F_{\vecb}}\bigl(\prolong{\nbigelambda}\bigr)
=\Image(
 \prolongg{\vecb}{\nbigelambda}_{|D_I}
\lrarr
 \prolong{\nbigelambda}_{|D_I}
 ).
\]
\end{lem}
\pf
It can be shown by an argument similar to
the proof of Lemma \ref{lem;9.10.16}.
\hfill\qed

\vspace{.1in}

Then Theorem \ref{thm;9.10.2} follows from Lemma \ref{lem;9.10.20}
and Lemma \ref{lem;9.10.21},
and Theorem \ref{thm;9.10.3} follows from
Lemma \ref{lem;9.10.23}, Lemma \ref{lem;9.10.22}
and Lemma \ref{lem;9.10.21}.
\hfill\qed

%% file: 14.1.tex

\subsubsection{Weak norm estimate of holomorphic sections}

Let $\vecv$ be a frame of $\prolongg{\vecb}{\nbigelambda}$
compatible with the parabolic filtrations
$(\lefttop{i}F\,|\,i=1,\ldots,l)$.
We obtain the numbers
$\lefttop{i}b(v_j):=\lefttop{i}\deg(v_j)$.
We put as follows:
\[
 v_j':=v_j\cdot\prod_{i=1}^l|z_i|^{\lefttop{i}b(v_j)},
\quad\quad
 \vecv'=(v_j').
\]
Then $\vecv'$ is a $C^{\infty}$-frame of $\nbigelambda$
over $X-D$.

\begin{prop}\label{prop;9.11.10}
$\vecv'$ is adapted up to log order.
\end{prop}
\pf
The argument is essentially same as 
the proof of Lemma \ref{lem;9.8.10}.
By our construction of $\vecv'$, the following is clear:
\[
 H(h,\vecv')\leq
 C_1\cdot\bigl(-\sum \log|z_i|\bigr)^M.
\]
Let $\vecv^{\lor}$ denote the dual frame of $\vecv$.
Then $\vecv^{\lor}$ gives a tuple of sections
of $\prolongg{-\vecb+(1-\epsilon)\vecdelta}{\nbige^{\lor\,\lambda}}$
for some $\epsilon>0$.
Let $P$ be a point of $D_i^{\circ}$.
Due to the result in the case of curves,
$\vecv^{\lor}_{|\pi_i^{-1}(P)}$ gives a frame of
$\prolongg{-b_i+(1-\epsilon)}
   {\nbige^{\lor\,\lambda}}_{|\pi^{-1}(P)}$,
which is compatible with the parabolic filtration.
We have $\lefttop{i}\deg^F(v^{\lor}_j)=-\lefttop{i}b(v_j)$
on $\pi_i^{-1}(P)$ for any point $P\in D^{\circ}_i$.
We put as follows:
\[
 \vecv^{\lor\,\prime}=(v_j^{\lor\,\prime}),
\quad\quad
 v_j^{\lor\,\prime}
:=v_j^{\lor}\cdot\prod_{i=1}^l|z_i|^{-\lefttop{i}b(v_j)}.
\]
Due to Corollary \ref{cor;11.28.15},
we obtain the following
(see the subsubsection \ref{subsubsection;a11.15.50}):
\[
 H(h^{\lor},\vecv^{\lor\,\prime})
\leq
 C_2\cdot
 \Bigl(
 -\sum_{i=1}^l\log |z_i|
 \Bigr)^M.
\]
It implies the following:
\[
 C_3\cdot \Bigl(
 -\sum_{i=1}^l\log |z_i|
 \Bigr)^{-M}
\leq H(h,\vecv').
\]
Thus we are done.
\hfill\qed

%% file: 15.tex
\subsubsection{Preliminary}

Recall that we may assume to 
have the following decomposition
(Condition \ref{condition;9.10.25}):
\[
 \nbige^0
=\bigoplus_{\veca\in\Sp(\theta)}
 E_{\veca}.
\]
Recall that there exists a positive constant $\epsilon_1$
such that $E_{\veca}$ and $E_{\vecb}$
are $\prod_{j\in\Diff(\veca,\vecb)}|z_j|^{\epsilon_1}$-asymptotically
orthogonal.

Let consider the $\lambda$-dependent
section $g(\lambda)$ of $End(E)$ over $X-D$,
given as follows:
\begin{equation}
 g(\lambda):=
 \bigoplus_{\veca}
 \exp\Bigl(
 \lambda\sum_{i=1}^l
 \bar{a}_i\cdot\log|z_i|^2
\Bigr)\cdot id_{E_{\veca}}.
\end{equation}

\begin{lem} 
We have the following equality:
\begin{equation} \label{eq;9.10.35}
 g(\lambda-\lambda_0)\cdot
 \delbar g(\lambda-\lambda_0)^{-1}
=-(\lambda-\lambda_0)\cdot
 \sum_{\veca}
 \Bigl(
 \sum_{i}\bar{a}_i\frac{d\bar{z}_i}{\bar{z}_i}
 \Bigr)\cdot id_{E_{\veca}}.
\end{equation}
\end{lem}
\pf
It can be checked by a direct calculation.
\hfill\qed

\vspace{.1in}

We have the decomposition
$\theta^{\dagger}=\phi_1+\phi_2+\phi_3$
satisfying the following:
\begin{itemize}
\item
 $\phi_1
=\sum_{\veca}
 \bigl(\sum_{i=1}^l\bar{a}_i\cdot \bar{z}_i^{-1}d\bar{z}_i\Bigr)
 \cdot id_{E_{\veca}}$.
\item
 $\phi_2=\sum_{\veca}\phi_{2\,\veca}$
where $\phi_{2\,\veca}\in End(E_{\veca})\otimes\Omega^{0,1}_X$,
and $|\phi_{2\,\veca}|_{h,\poin}$ is bounded.
Here $|\cdot|_{h,\poin}$ denotes the norm with respect to
$h$ and the Poincar\'{e} metric.
\item
 $\phi_3=\sum_{\veca\neq \vecb}\phi_{3\,\veca,\vecb}$,
 where
 $\phi_{3\,\veca\,\vecb}$ are sections of
 $Hom(E_{\veca},E_{\vecb})\otimes\Omega^{0,1}_X$.
 We have the following estimate for some positive constants
 $\epsilon_2$ and $C$:
\[
 \bigl|\phi_{3\,\veca,\vecb}\bigr|_{h,\poin}
\leq C\cdot \!\!\!\prod_{i\in\Diff(\veca,\vecb)}|z_i|^{\epsilon_2}.
\]
\end{itemize}

The following lemma is clear
\begin{lem} \label{lem;9.10.36}
$g(\lambda)$ and  $\phi_i$ ($i=1,2$) are commutative.
\hfill\qed
\end{lem}

\begin{lem}
We have the following formula:
\begin{multline} \label{eq;9.10.30}
g(\lambda-\lambda_0)\cdot
 (\delbar_E+\lambda\theta^{\dagger})\cdot
 g(\lambda-\lambda_0)^{-1}=\\
\delbar_E+\lambda_0\cdot\theta^{\dagger}
+(\lambda-\lambda_0)\cdot\bigl(\phi_2+\phi_3\bigr)
+\lambda\cdot\Bigl(
 g(\lambda-\lambda_0)\cdot \phi_3\cdot g(\lambda-\lambda_0)^{-1}
 -\phi_3
 \Bigr).
\end{multline}
\end{lem}
\pf
We have the following equality:
\begin{multline} \label{eq;9.10.31}
 g(\lambda-\lambda_0)\cdot
 (\delbar+\lambda\theta^{\dagger})\cdot
 g(\lambda-\lambda_0)^{-1}
=\delbar_E+g(\lambda-\lambda_0)\cdot\delbar g(\lambda-\lambda_0)^{-1}
+\lambda\cdot g(\lambda-\lambda_0)\cdot
  \theta^{\dagger} \cdot g(\lambda-\lambda_0)^{-1}\\
=\delbar_E+\lambda_0\cdot\theta^{\dagger}
+g(\lambda-\lambda_0)\cdot\delbar g(\lambda-\lambda_0)^{-1}
+(\lambda-\lambda_0)\cdot\theta^{\dagger}
+\lambda\cdot
 \Bigl(g(\lambda-\lambda_0)\cdot
 \theta^{\dagger}\cdot g(\lambda-\lambda_0)^{-1}
 -\theta^{\dagger}\Bigr).
\end{multline}
We obtain the following from (\ref{eq;9.10.35}):
\begin{equation} \label{eq;9.10.32}
 g(\lambda-\lambda_0)\cdot\delbar g(\lambda-\lambda_0)^{-1}
+(\lambda-\lambda_0)\cdot\theta^{\dagger}
=(\lambda-\lambda_0)\cdot\bigl(\phi_2+\phi_3\bigr).
\end{equation}
We also obtain the following from Lemma \ref{lem;9.10.36}:
\begin{equation} \label{eq;9.10.33}
 g(\lambda-\lambda_0)\cdot\theta^{\dagger}\cdot g(\lambda-\lambda_0)^{-1}
-\theta^{\dagger}
=g(\lambda-\lambda_0)\cdot \phi_3 \cdot g(\lambda-\lambda_0)^{-1}
-\phi_3.
\end{equation}
Then we obtain (\ref{eq;9.10.30})
from (\ref{eq;9.10.31}), (\ref{eq;9.10.32}) and (\ref{eq;9.10.33}).
\hfill\qed

\begin{lem} \label{lem;9.10.40}
We put
$\psi(\lambda):=
 g(\lambda-\lambda_0)\cdot \phi_3\cdot g(\lambda-\lambda_0)^{-1}
 -\phi_3$.
There exist positive constants $\eta$, $\epsilon'$ and $C$
such that the following holds
for any $\lambda,\lambda'\in\Delta(\lambda_0,\eta)$:
\[
 \big|
 \bigl(
 \psi(\lambda)-\psi(\lambda')
 \bigr)_{\veca,\vecb}
 \big|_{h,\poin}
\leq
 C\cdot |\lambda-\lambda'|\cdot
 \prod_{i\in\Diff(\veca,\vecb)}|z_i|^{\epsilon'}.
\]
Here $\bigl( A \bigr)_{\veca,\vecb} $
denotes the $(\veca,\vecb)$-component
of $A$,
and $|\cdot|_{h,\poin}$ denotes the norm
with respect to $h$ and the Poincare metric.
\end{lem}
\pf
By definition,
the $(\veca,\vecb)$-component of $\psi(\lambda)-\psi(\lambda')$
is as follows:
\[
 \phi_{3\,\veca,\vecb}
 \cdot
 \Bigl(
 \prod_{i=1}^l|z_i|^{2\lambda(a_i-b_i)}
-\prod_{i=1}^l|z_i|^{2\lambda'(a_i-b_i)}
 \Bigr).
\]
Hence the norm is dominated by the following:
\[
 |\phi_{3\,\veca\,\vecb}|_{h,\poin}
\cdot\prod_{i=1}^l |z_i|^{2\Re(\lambda'(a_i-b_i))}
\cdot
 \Bigl(\prod_{i=1}^l|z_i|^{2(\lambda-\lambda')(a_i-b_i)}-1
 \Bigr)
\leq
 C\cdot\prod_{i\in\Diff(\veca,\vecb)}
  |z_i|^{\epsilon-\eta|a_i-b_i|}(-\log|z_i|)
 \cdot |\lambda-\lambda'|.
\]
Thus we are done.
\hfill\qed

\vspace{.1in}

Let us pick a point $\lambda_0\in\cnum_{\lambda}$.
We put as follows for any $\lambda\in \cnum_{\lambda}$:
\begin{equation} \label{eq;9.10.88}
 d''(\lambda):=
 g(\lambda-\lambda_0)\cdot
 (\delbar_E+\lambda\cdot\theta^{\dagger})\cdot
 g(\lambda-\lambda_0)^{-1}.
\end{equation}
Then we have the following equality due to (\ref{eq;9.10.30}):
\[
 d''(\lambda)
=\delbar_E+\lambda\cdot \theta^{\dagger}
 +(\lambda-\lambda_0)\cdot (\phi_2+\phi_3)
 +\lambda\cdot\psi(\lambda).
\]
It gives the holomorphic structure of $C^{\infty}$-bundle $E$ over $X-D$,
and it is equivalent to
$\delbar_{\nbigelambda}=\delbar_E+\lambda\cdot\theta^{\dagger}$
up to the (not unitary) gauge transformation.

\begin{lem}
If $\eta>0$ is sufficiently small,
then there exists a positive constant $C$
such that the following holds for any
$\lambda,\lambda'\in\Delta(\lambda_0,\eta)$:
\begin{equation} \label{eq;8.18.15}
 |d''(\lambda)-d''(\lambda')|_{h,\poin}
\leq |\lambda-\lambda_0|\cdot C.
\end{equation}
Note that $d''(\lambda)-d''(\lambda')$ are $(0,1)$-forms.
\end{lem}
\pf
We have the following:
\[
 d''(\lambda)-d''(\lambda')
=(\lambda-\lambda')\cdot(\phi_2+\phi_3)
+\lambda\cdot(\psi(\lambda)-\psi(\lambda')).
\]
Thus we obtain the result from Lemma \ref{lem;9.10.40}
and the estimates for $\phi_i$ $(i=2,3)$.
\hfill\qed

\vspace{.1in}
Let $p_{\lambda}$ denote the projection
$\Delta(\lambda_0,\eta)\times (X-D)\lrarr X-D$,
and then we have the $C^{\infty}$-bundle
$p_{\lambda}^{-1}E$ on $\Delta(\lambda_0,\eta)\times (X-D)$.
We have the naturally defined operator:
\[
 \delbar_{\lambda}+d''(\lambda):
 C^{\infty}(p_{\lambda}^{-1}(E))
\lrarr
 C^{\infty}\bigl(p_{\lambda}^{-1}(E)\otimes\Omega^{0,1}_{ 
 \nbigx(\lambda_0,\eta)-\nbigd(\lambda_0,\eta)  }\bigr).
\]

\begin{lem} \label{lem;9.10.86}
The operator  $\delbar_{\lambda}+d''(\lambda)$
gives a holomorphic structure,
i.e.,
$(\delbar_{\lambda}+d''(\lambda))^2=0$.
\end{lem}
\pf
Note that $g(\lambda-\lambda_0)$ is holomorphic
with respect to $\lambda$.
Then the claim can be checked by a direct calculation.
\hfill\qed

\vspace{.1in}

We use the notation in the subsubsection \ref{subsubsection;10.12.1}.
We put as follows:
\[
 \big\langle
 f_1,f_2
 \big\rangle_{\lambda,\veca,N}
:=
 \int\big(f_1,f_2\big)_{\veca,N}\dvol
+\int\big(d''(\lambda)f_1,\,\,d''(\lambda)f_2\bigr)_{\veca,N}\dvol,
\quad\quad
 ||f||^2_{\lambda,\veca,N}
:=\big\langle
   f,f
  \big\rangle_{\lambda,\veca,N}.
\]
\begin{lem}
For any $\lambda\in\Delta(\lambda_0,\eta)$,
the norms
$||\cdot||_{\lambda,\veca,N}$
and $||\cdot||_{\lambda_0,\veca,N}$
are equivalent.
\end{lem}
\pf
It follows from the inequality (\ref{eq;8.18.15})
\hfill\qed

\vspace{.1in}

Hence the completions with respect to the norms
$||\cdot||_{\lambda,\veca,N}$ are independent of a choice of
$\lambda\in \Delta(\lambda_0,\eta)$
for some sufficiently small positive number $\eta$.
Let $A^{0,q}_{\veca,N}(\nbige^{\lambda_0})$ be
the completion of the space $A_c^{0,q}(\nbige^{\lambda_0})$
with respect to the norm
$||\cdot||_{\veca,N}$.
Then we have the family of complexes
$\bigl(A^{0,\cdot}_{\veca,N}(\nbige^{\lambda_0}),d''(\lambda)\bigr)$
$(\lambda\in\Delta(\lambda_0,\eta))$.

%% file: 15.3.tex

\subsubsection{Extension of holomorphic sections}

\begin{lem} \label{lem;9.10.80}
There exists a positive number $\eta>0$
and the family of linear morphisms
$G(\lambda):
  \ker(d''(\lambda_0))\lrarr A^{0,0}_{\veca,N}(\nbige^{\lambda_0})$
depending $\lambda\in\Delta(\lambda_0,\eta)$
satisfying the following:
\begin{itemize}
\item
The vanishing
$H^i(A^{0,\cdot}_{\veca,N}(\nbige^{\lambda_0}),d''(\lambda))=0$
holds for any $i>0$
and for any $\lambda\in \Delta(\lambda_0,\eta)$.
\item
The morphism
$G(\lambda)$
satisfies the conditions 
{\rm\ref{8.18.21}}, {\rm\ref{8.18.22}} and {\rm\ref{8.18.23}}
in Lemma {\rm\ref{lem;9.10.50}}.
It gives the trivialization of the family
$\bigl\{
 \ker d''(\lambda)\,|\,\lambda\in\Delta(\lambda_0,\eta)
 \bigr\}$,
namely
 $G(\lambda)$ gives the homeomorphism
of $\ker d''(\lambda_0)$ and $\ker d''(\lambda)$
for any point $\lambda\in\Delta(\lambda_0,\eta)$.
\end{itemize}
\end{lem}
\pf
Note that $d''(\lambda_0)=\delbar_{\nbigelambdazero}$,
and hence the conditions in Lemma \ref{lem;9.10.75}
is satisfied due to Lemma \ref{lem;9.10.72}.
Then we obtain the result due to Lemma \ref{lem;9.10.75}.
\hfill\qed

\vspace{.1in}

Recall that we have the $C^{\infty}$-bundle
$p_{\lambda}^{-1}(E)$ with the hermitian metric
$h_{\veca,N}=
 h\cdot \prod_{i=1}^l |z_i|^{a_i}\cdot\bigl(-\log|z_i|\bigr)^N$.
We have the holomorphic structure
$\delbar_{\lambda}+d''(\lambda)$
(Lemma \ref{lem;9.10.86}).

\begin{cor}
For any section $f$ of $\nbige^{\lambda_0}$ over $\nbigx^{\lambda_0}$,
we have a holomorphic section
$\tilde{f}$ of the holomorphic bundle
$\bigl(p_{\lambda}^{-1}E,\delbar_{\lambda}+d''(\lambda)\bigr)$ 
over $\Delta(\lambda_0,\eta)\times (X-D)$
such that
$\tilde{f}_{|\{\lambda\}\times (X-D)}
\in A^{0,0}_{\veca,N}(\nbigelambdazero)$.
\end{cor}
\pf
Let $G(\lambda)$ be the family of the morphism
given in Lemma \ref{lem;9.10.80}.
We put $\tilde{f}_{|\{\lambda\}\times (X-D)}:=G(\lambda)(f)$.
We have the absolute convergent series
in $A^{0,0}_{\veca,N}(\nbigelambdazero)$:
\begin{equation} \label{eq;9.10.87}
 \tilde{f}=\sum (\lambda-\lambda_0)^i\cdot f_i,
\quad\quad
 f_i\in A^{0,0}_{\veca,N}(\nbigelambdazero).
\end{equation}
By our construction,
it is clear that 
the restrictions
$\tilde{f}_{\{\lambda\}\times (X-D)}$
are contained in
$A^{0,0}_{\veca,N}(\nbigelambdazero)$.
Since (\ref{eq;9.10.87}) is absolute convergent,
we also have the following finiteness:
\begin{equation}\label{eq;9.10.82}
 \int_{\Delta(\lambda_0,\epsilon)\times (X-D)}
 |\tilde{f}|_h^2
 \cdot\prod_{i=1}^l |z_i|^{a_i}(-\log|z_i|)^N\dvol
< \sum \frac{\pi}{h+1}\epsilon^{2(h+2)}||f_i||<\infty.
\end{equation}
The finiteness (\ref{eq;9.10.82}) implies
that $\tilde{f}$ can be regarded as
$L^2$-section on
$\Delta(\lambda_0,\eta)\times (X-D)$
with respect to the metric $h_{\veca,N}$.

We have clearly $\delbar_{\lambda}\tilde{f}=0$.
We have $d''(\lambda)(\tilde{f}_{|\lambda})=0$
for any point $\lambda\in\Delta(\lambda_0,\eta)$,
by our construction.
Then we obtain $d''(\lambda)(\tilde{f})=0$
in the distribution sense, due to Fubini's theorem.
Hence we obtain 
$\bigl(\delbar_{\lambda}+d''(\lambda)\bigr)\tilde{f}=0$
in the distribution sense.
Thus we can conclude that $\tilde{f}$ is holomorphic section
with respect to $\delbar_{\lambda}+d''(\lambda)$.
\hfill\qed

\vspace{.1in}

We put $F(\lambda):=g(\lambda-\lambda_0)^{-1}\cdot\tilde{f}$.
Then $F$ gives a section of the $C^{\infty}$-bundle
$p_{\lambda}^{-1}(E)$ over $\Delta(\lambda_0,\eta)\times X$.

\begin{lem} \label{lem;9.10.90}
 For any positive number $\epsilon$,
 there exists a positive number $\eta$ satisfying the following:
\begin{itemize}
\item
 $F$ is holomorphic with respect to the holomorphic structure
 $\delbar_{\nbige}=\delbar_{\lambda}+\delbar_E+\lambda\cdot\theta^{\dagger}$,
 i.e.,
 $F$ gives a holomorphic section of $\nbige$
 over $\Delta(\lambda_0,\eta)\times (X-D)$.
\item
 For any $\lambda\in\Delta(\lambda_0,\eta)$
 and for any $\epsilon'>0$
 there exists a positive constant $C$ satisfying the following
 inequality:
 \[
  \bigl|F_{\{\lambda\}\times (X-D)}\bigr|_h=
 C\cdot\Bigl(\prod |z_i|^{-a_i-\epsilon-\epsilon'}\Bigr).
 \]
\end{itemize}
\end{lem}
\pf
Since $g$ is holomorphic with respect to $\lambda$,
we have the following relation from (\ref{eq;9.10.88}):
\[
 \delbar_{\lambda}+d''(\lambda)
=g(\lambda-\lambda_0)\cdot
 \Bigl(
 \delbar_{\lambda}+\delbar_E+\lambda\cdot\theta^{\dagger}
 \Bigr)\cdot
g(\lambda-\lambda)^{-1}.
\]
Then the holomorphic property of $F$ with respect to 
$\delbar_{\nbige}$ follows from
the holomorphic property of $\tilde{f}$
with respect to $\delbar_{\lambda}+d''(\lambda)$.

For any positive number $\epsilon$,
there exist positive constants $\eta>0$ and $C_1$
such that the following holds
for any $\lambda\in\Delta(\lambda_0,\eta)$:
\begin{equation} \label{eq;9.10.89}
 \bigl|g(\lambda-\lambda_0)\bigr|_h\leq
 C_1\cdot \Bigl(\prod |z_i|^{-\epsilon}\Bigr).
\end{equation}
Since $\tilde{f}_{\{\lambda\}\times (X-D)}$
is an element of $A^{0,0}_{\veca,N}(\nbigelambdazero)$,
we obtain the following finiteness
for any $\lambda\in\Delta(\lambda_0,\eta)$
from (\ref{eq;9.10.89}):
\[
 \int_{\{\lambda\}\times (X-D)}
 |F(\lambda)|_h^2
\cdot\prod |z_i|^{a_i+\epsilon}\dvol<\infty
\]
It implies the second claim.
\hfill\qed

\vspace{.1in}
We formulate the result in this subsubsection.

\begin{prop} \label{prop;9.10.60}
Let $\vecb=(b_1,\ldots,b_l)$ be an element of $\real^l$.
Assume that $b_i\not \in \Par(\nbigelambdazero,i)$.
Then there exists a positive constant $\eta$
satisfying the following:
\begin{itemize}
\item
For any holomorphic section $f$ of
$\prolongg{\vecb}{\nbigelambdazero}$,
there exists a holomorphic section $F$
of $\prolongg{\vecb}{\nbige}$ over $\Delta(\lambda_0,\eta)$
such that $F_{|\{\lambda\}\times X}=f$.
\end{itemize}
\end{prop}
\pf
It immediately follows from Lemma \ref{lem;9.10.90}.
\hfill\qed

%% file: 15.4.tex

\subsubsection{Prolongation of $\nbige$ around $\lambda_0$}

\begin{condition}\label{condition;a11.16.1}
Let $\vecb=(b_1,\ldots,b_l)$ be an element of $\real^l$
such that 
$b_i\not\in\Par(\nbige^{\lambda_0},i)$ for any $i$.
\hfill\qed
\end{condition}

\begin{thm}
Let $\vecb$ be an element of $\real^l$ as in
Condition {\rm\ref{condition;a11.16.1}}.
Then there exists a positive number $\eta$ such that
$\prolongg{\vecb}{\nbige}$ is locally free
over $\Delta(\lambda_0,\eta)\times X$.
\end{thm}
\pf
By considering the tensor product
of $\harmonicbundle$ and the model bundle $L(\vecu)$,
we may assume that the residue of $\tr(\theta)$ is trivial.
We can also assume that $\vecb=0$.
Note we have $0\not\in\Par(\nbige^{\lambda_0},i)$
due to our assumption.

Let $\vecv$ be a frame of $\prolong{\nbigelambdazero}$,
compatible with $\lefttop{i}F$ $(i=1,\ldots,l)$.
We put $\lefttop{i}a_j:=\lefttop{i}\deg(v_j)$.
We have $-1<\lefttop{i}a_j<0$.
We put as follows:
\[
 \veca_j:=
 \bigl(
 \lefttop{1}a_j,
 \lefttop{2}a_j,\ldots,
 \lefttop{l}a_j
 \bigr),
\quad\quad
 \vecc:=\sum_j \veca_j,
\quad\quad
 \lefttop{i}c=\sum_j^{\rank E} \lefttop{i}a_j.
\]
Then $v_j$ is a holomorphic section of
$\prolongg{\veca_j}{\nbigelambdazero}$,
and we have
$\det(\prolong{\nbigelambdazero})
=\prolongg{\vecc}{\det(\nbigelambdazero)}$.

By using Proposition \ref{prop;9.10.60},
we can take $\epsilon>0$ and $\eta>0$,
and holomorphic sections $\tilde{v}_j$
of $\prolongg{\veca_j+\epsilon\cdot\vecdelta}{\nbigelambdazero}$
over $\Delta(\lambda_0,\eta)\times X$
such that $\tilde{v}_{j\,|\,\{\lambda_0\}\times X}=v_j$.
In particular, $\tilde{v}_j$ give holomorphic sections
of $\prolong{\nbige}$ over $\Delta(\lambda_0,\eta)\times X$.
Hence we obtain the tuple of sections $\tilde{\vecv}:=(\tilde{v}_i)$
of $\prolong{\nbige}$ over $\Delta(\lambda_0,\eta)\times X$.

Then we have the following inequality
for any sufficiently small $\epsilon>0$
and for any $\lambda\in \Delta(\lambda_0,\eta)$:
\[
 \lefttop{i}\deg(\Omega(\tilde{\vecv})_{|\lambda})
< \sum_j \lefttop{i}a_j+r\cdot \epsilon
< \lefttop{i}c+1.
\]
Hence $\Omega(\tilde{\vecv})_{|\{\lambda\}\times X}$
is a holomorphic section of
$\prolongg{\vecc}{\det(\nbigelambda)}
=\prolongg{\vecc}{\det(\nbige)}_{|\{\lambda\}\times X}$.
Since we have assumed that the residue of $\tr(\theta)$ is trivial,
we have $\prolongg{\vecc}{\det(\nbige)}_{|\lambda}=
 \det(\prolong{\nbigelambda})$.
Since we have $\Omega(\tilde{\vecv})_{|(\lambda_0,O)}\neq 0$,
the section
$\Omega(\tilde{\vecv})$ gives a frame
around $(\lambda_0,O)$.

Then we can use the standard argument as follows
(the last argument in the proof of Lemma \ref{lem;d11.14.51}).
Let $f$ be a holomorphic section of $\prolong{\nbige}$
on a neighbourhood $U$ of $(\lambda_0,O)$.
On $\bigl( U\cap (X-D)\bigr) \times\Delta(\lambda_0,\epsilon)$,
we have the following description:
\[
 f=\sum f_i\cdot \tilde{v}_i.
\]
As usual, we consider the section
$f\wedge \tilde{v}_2\wedge\cdots\wedge\tilde{v}_r
=f_1\cdot\Omega(\tilde{\vecv})$
of $\det(\prolong{\nbige})=\prolongg{\vecc}{\det(\nbige)}$
over $U$,
and we can derive that $f_1$ is holomorphic.
Similarly, we can show that $f_i$ are holomorphic for any $i$.
It implies that $\tilde{\vecv}$ gives a frame of
$\prolong{\nbige}$ on a neighbourhood of $(\lambda_0,O)$.
\hfill\qed

%% file: a55.tex

\subsubsection{The parabolic structures
$\lefttop{i}\Fzero(\prolongg{\vecb}{\nbige})$}

\label{subsubsection;10.18.20}

Let $\vecb$ be an element of $\real^l$ as
in Condition \ref{condition;a11.16.1}.
Let $\vecv=(v_j)$ be a frame of $\prolongg{\vecb}{\nbigelambdazero}$
compatible with the parabolic structure.
For each $v_j$, we obtain the element
$\veca(v_j)=(a_1(v_j),\ldots,a_l(v_j))\in\real^l$,
where we put $a_i(v_j):=\lefttop{i}\deg^F(v_j)$.

Let us pick a positive number $\epsilon$ as follows:
\begin{equation}\label{eq;a11.16.2}
  0<\epsilon<
 \frac{1}{3}
 \min \Bigl(
 \bigcup_{i=1}^l
 \bigl\{
 |a_1-a_2|\,\big|\,
 a_1,a_2\in\Par(\nbigelambdazero,i),\,\,
 a_1\neq a_2
 \bigr\}
 \Bigr).
\end{equation}

We can pick a positive number $\epsilon_0$
such that there exist sections $\tilde{v}_j$ of
$\prolongg{\veca(v_j)+\epsilon\cdot\vecdelta}{\nbige}$
on $\nbigx(\lambda_0,\epsilon_0)$,
such that $\tilde{v}_{j|\nbigx^{\lambda_0}}=v_j$.
If $\epsilon_0$ is sufficiently small,
we may assume the following:
\begin{condition}\label{condition;a11.16.3}
 Let $u_1,u_2$ be elements of $\KMS(\nbige^0,i)$ such that
 $\paramap(\lambda_0,u_1)<\paramap(\lambda_0,u_2)$.
 Note we have
 $\paramap(\lambda_0,u_1)+\epsilon<\paramap(\lambda_0,u_2)$,
 due to {\rm(\ref{eq;a11.16.2})}.
 Then the inequality
 $\paramap(\lambda,u_1)+\epsilon<\paramap(\lambda',u_2)$
 holds for any $\lambda,\lambda'\in\Delta(\lambda_0,\epsilon_0)$.
\end{condition}

On $\nbigd_i(\lambda_0,\epsilon_0)$,
the filtration $\lefttop{i}F^{(\lambda_0)}$
is obtained as follows:
\[
 \lefttop{i}F^{(\lambda_0)}_a
 \bigl(\prolongg{\vecb}{\nbige}_{|\nbigd_i(\lambda_0,\epsilon_0)}\bigr)
:=
 \big\langle
  \tilde{v}_{j\,|\,\nbigd_i(\lambda_0,\epsilon_0)}\,
 \big|\,\lefttop{i}\deg(v_j)\leq a'
 \big\rangle,
\quad\quad
 a':=\max\bigl\{
 x\in\Par(\prolongg{\vecb}{\nbige},i)\,\big|\,
 x\leq a
 \bigr\}.
\]
For any $\lambda\in\Delta(\lambda_0,\eta)$ and
$c\in \Par(\prolongg{\vecb}{\nbigelambdazero},i)$,
we put as follows:
\[
 d(\lambda,\lambda_0,c)
:=\max\bigl\{
 \paramap(\lambda,u)\,\big|\,
 u\in \KMS(\nbige^0,i),\,\,
 \paramap(\lambda_0,u)=c
 \bigr\}.
\]
Recall that we have the parabolic filtration
$\lefttop{i}F$ of $\prolongg{\vecb}{\nbigelambda}$.
So we compare the two vector subbundle
$\lefttop{i}F^{(\lambda_0)}_{c\,|\,\{\lambda\}\times D_i}$
and $\lefttop{i}F_{d(\lambda,\lambda_0,c)}$
of $\prolongg{\vecb}{\nbigelambda}_{|\nbigd^{\lambda}_i}$.

\begin{lem} \label{lem;9.10.91}
For any $\lambda\in\Delta(\lambda_0,\epsilon_0)$,
we have the equality
$\lefttop{i}F^{(\lambda_0)}_{c\,|\,\{\lambda\}\times D_i}
=\lefttop{i}F_{d(\lambda,\lambda_0,c)}$.
In particular, the filtration
$\lefttop{i}F^{(\lambda_0)}$ is independent of choices of
a compatible frame $\vecv$ and an extension $\tilde{\vecv}$.
\end{lem}
\pf
Since both of 
$\lefttop{i}F^{(\lambda_0)}_{c\,|\,\{\lambda\}\times D_i}$
and $\lefttop{i}F_{d(\lambda,\lambda_0,c)}$ are vector subbundles,
we have only to show
$\lefttop{i}F^{(\lambda_0)}_{c\,|\,(\lambda,P)}
=\lefttop{i}F_{d(\lambda,\lambda_0,c)\,|\,(\lambda,P)}$ 
for any $P\in D_i^{\circ}$.
For that purpose, we have only to consider
the restriction $\nbige$ to the curve $\pi_i^{-1}(P)$.
Thus we can restrict our attention to the case
$\dim(X)=1$, which is assumed in the following of the proof.

In the case $\deg^F(v_i)=c$,
we have the following inequality 
for any $\lambda\in\Delta(\lambda_0,\epsilon_0)$,
due to Condition \ref{condition;a11.16.3}:
\[
 \deg^F\bigl(\tilde{v}_{i\,|\,\nbigxlambda}\bigr)\leq
 c+\epsilon
 <
 \min\bigl\{
 d\in\Par(\nbigelambda,i)\,\big|\,
 d>d(\lambda,\lambda_0,c)
 \bigr\}.
\]
It implies
$\tilde{v}_{i\,|\,(\lambda,O)}$
is contained in the space
$F_{d(\lambda,\lambda_0,c)}(\prolongg{b}{\nbigelambda}_{|O})$.
Hence we obtain 
$\Fzero_{c\,|\,(\lambda,O)}\subset
  F_{d(\lambda,\lambda_0,c)}(\prolongg{b}{\nbigelambda}_{|O})$.

Due to our construction of the filtration $\Fzero$,
we have the following equality:
\[
 \rank\Fzero_c
=\sum_{b-1<c'\leq c}
 \multiplicity(\lambda_0,c')
=\sum_{\substack{ u\in \KMS(\nbige^0),\\ b-1<\paramap(\lambda_0,u)\leq c}}
 \multiplicity(0,u).
\]
On the other hand,
we have the following equality:
\[
 \rank F_{d(\lambda,\lambda_0,c)}(\prolongg{b}{\nbigelambda}_{|O})
=\sum_{b-1<c'\leq d(\lambda,\lambda_0,c)}
 \multiplicity(\lambda,c')
=\sum_{\substack{ u\in \KMS(\nbige^0),\\ 
 b-1<\paramap(\lambda,u)\leq d(\lambda,\lambda_0,c)}}
 \multiplicity(0,u)
=\sum_{\substack{ u\in \KMS(\nbige^0),\\
 b-1<\paramap(\lambda_0,u)\leq c}}
 \multiplicity(0,u).
\]
Note $\paramap(\lambda_0,u)>b-1$ if and only if
$\paramap(\lambda,u)>b-1$ due to our assumption
$b\not\in \Par(\nbigelambdazero,i)$.
Thus we obtain
$\Fzero_{c}(\prolongg{b}{\nbige})_{|(\lambda,O)}=
  F_{d(\lambda,\lambda_0,c)}(\prolongg{b}{\nbigelambda}_{|O})$.
\hfill\qed

\vspace{.1in}
On $\nbigd_i(\lambda_0,\epsilon_0)$, we obtain the following vector bundle:
\[
 \lefttop{i}\Gr^{F^{(\lambda_0)}}_c
 \bigl(\prolongg{\vecb}{\nbige}_{|\nbigd_i(\lambda_0,\epsilon_0)}\bigr)
:=
 \frac{\lefttop{i}F^{(\lambda_0)}_c
 \bigl(\prolongg{\vecb}{\nbige}_{|\nbigd_i(\lambda_0,\epsilon_0)}\bigr)  }
 {\lefttop{i}F^{(\lambda_0)}_{<c}
 \bigl(\prolongg{\vecb}{\nbige}_{|\nbigd_i(\lambda_0,\epsilon_0)}\bigr) }.
\]
\begin{cor}
We have the following on $\nbigd_i^{\lambda}$
for any $\lambda\in\Delta(\lambda_0,\eta)$:
\begin{equation} \label{eq;9.11.3}
 \lefttop{i}\Gr^{F^{(\lambda_0)}}_{c}
 \bigl(\prolongg{\vecb}{\nbige}_{|\nbigd_i(\lambda_0,\epsilon_0)}\bigr)
   _{|\nbigd_i^{\lambda}}
\simeq
 \frac{\lefttop{i}F_{d(\lambda,\lambda_0,c)}
  \bigl(\prolongg{\vecb}{\nbigelambda}\bigr)}
   {\sum_{b<c}\lefttop{i}F_{d(\lambda,\lambda_0,b)}
 \bigl(\prolongg{\vecb}{\nbigelambda}\bigr) }.
\end{equation}
\end{cor}
\pf
It immediately follows from Lemma \ref{lem;9.10.91}.
\hfill\qed

%% file: 15.8.tex

\subsubsection{Regularity of $\DD$ and $\DD^{\lambda}$}

\begin{lem} \label{lem;9.10.95}
$\DD$ is the regular $\lambda$-connection,
namely,
if $f$ is a section of $\prolongg{\vecb}{\nbige}$,
then $\DD f$ is a section of
$\prolongg{\vecb}{\nbige}\otimes\Omega^{1,0}_X(\log D)$.
\end{lem}
\pf
Let us consider the case $\lambda_0\neq 0$.
We may assume that any $\lambda\in\Delta(\lambda_0,\epsilon_0)$
is generic.
In this case,
the prolongment
$\prolongg{\vecb}{\nbige}_{|X^{\ast}(\lambda_0,\epsilon_0)}$
can be essentially regarded as a quasi canonical prolongment.
Thus $\DD f_{|X^{\ast}(\lambda_0,\epsilon_0)}$
gives a section of
$\prolongg{\vecb}{\nbige}\otimes{\Omega^{1,0}(\log \nbigd)
 _{|\nbigx^{\ast}(\lambda_0,\epsilon_0)}}$.
Hence $\DD f$ gives a section of
$\prolongg{\vecb}{\nbige}\otimes{\Omega^{1,0}(\log \nbigd)}$
over $\nbigx(\lambda_0,\epsilon_0)- \nbigd^{\lambda_0}$.
Note that the codimension of $\nbigd^{\lambda_0}$
in $\nbigx(\lambda_0,\epsilon_0)$ is two.
Thus $\DD f$ gives a section of
$\prolongg{\vecb}{\nbige}\otimes{\Omega^{1,0}(\log \nbigd)}$
over $\nbigx(\lambda_0,\epsilon_0)$.

Then we can show the claim in the case $\lambda=0$ by the same argument.
\hfill\qed

\begin{cor}
$\DD^{\lambda}$ is the regular $\lambda$-connection.
\end{cor}
\pf
Let $\vecb=(b_1,\ldots,b_l)$ be an element of $\real^l$
and let $f$ be a section of $\prolongg{\vecb}{\nbigelambda}$.
We may assume that $b_i\not\in\Par(\nbigelambda,i)$
for any $i$.
We can take a holomorphic section $F$
of $\prolongg{\vecb}{\nbige}$ over $\Delta(\lambda,\eta)\times X$
such that $F_{|\nbigxlambda}=f$.
Due to Lemma \ref{lem;9.10.95},
$\DD F$ is a section of
$\prolongg{\vecb}{\nbige}\otimes
 \Omega_X(\log D)$.
Since we have $\DD F_{|\nbigxlambda}=\DD^{\lambda}f$,
we obtain the result.
\hfill\qed

\vspace{.1in}

\subsubsection{The residue and the $\lambda$-connection on the divisors}

On $\nbigd_i(\lambda_0,\epsilon_0)$,
we have $\prolongg{\vecb}{\nbige}_{|\nbigd_i(\lambda_0,\epsilon_0) }$.
Then we have the endomorphisms $\Res_i(\DD)$,
which preserve the parabolic filtrations
due to Lemma \ref{lem;9.10.95}.
We also have the induced $\lambda$-connection
$\lefttop{i}\DD$ 
of $\prolongg{\vecb}{\nbige}_{|\nbigd_i(\lambda_0,\epsilon_0) }$,
which is defined as follows: 
For any $f\in\prolongg{\vecb}{\nbige}_{|\nbigd_i(\lambda_0,\epsilon_0)}$,
pick $F\in\prolongg{\vecb}{\nbige}$
such that $F_{|\nbigd_i(\lambda_0,\epsilon_0) }=f$.
Then we put
$\lefttop{i}\DD(f):=\DD F_{|\nbigd_i(\lambda_0,\epsilon_0) }$.
\begin{lem}
It is well defined.
\end{lem}
\pf
Assume $F_{|\nbigd_i(\lambda_0,\epsilon_0) }=0$.
Then we have the description $F=z_i\cdot G$
for some $G\in \prolongg{\vecb}{\nbige}$.
We have the following:
\[
 \DD(z_i\cdot G)=
 \lambda\cdot dz_i\cdot G
+z_i\cdot \DD G.
\]
Thus $\pi(\DD(z_i\cdot G)_{|\nbigd_i(\lambda_0,\epsilon_0)  })=0$,
where $\pi$ denotes the projection
$\Omega_X(\log D)_{|\nbigd_i(\lambda_0,\epsilon_0)  }
\lrarr
 \Omega_{D_i}(\log D\cap D_i)$.
Hence we are done.
\hfill\qed

\vspace{.1in}

Let $\vecv$ be a frame of $\prolongg{\vecb}{\nbige}$,
which is compatible with the parabolic
filtrations $\bigl(\lefttop{i}F\,|\,i=1,\ldots,l\bigr)$.
Then we obtain the $\lambda$-connection form $\nbiga$
determined by $\DD\vecv=\vecv\cdot\nbiga$.
We develop $\nbiga$ as $\nbiga=\sum \nbiga^k\frac{dz_k}{z_k}$.
\begin{lem}
Then we have the following formula:
\[
 \Res_i(\DD)\vecv_{|\nbigd_i(\lambda_0,\epsilon_0)}
=\vecv_{|\nbigd_i(\lambda_0,\epsilon_0)}
 \cdot \nbiga^i_{|\nbigd_i(\lambda_0,\epsilon_0) },
\quad\quad
 \lefttop{i}\DD\vecv_{|\nbigd_i(\lambda_0,\epsilon_0)}
=\vecv_{|\nbigd_i(\lambda_0,\epsilon_0)}
 \cdot\sum_{k\neq i}
 \nbiga^k_{|\nbigd_i(\lambda_0,\epsilon_0) }
 \frac{dz_k}{z_k}
\]
\end{lem}
\pf
It immediately follows from the definitions.
\hfill\qed

\begin{lem} \label{lem;9.10.96}
$\lefttop{i}\DD$ and $\Res_i(\DD)$ preserve
the filtration $\lefttop{i}F^{(\lambda_0)}$.
\end{lem}
\pf
Let $f$ be a section of
$\lefttop{i}F^{(\lambda_0)}_{c}\bigl(
 \prolongg{\vecb}{\nbige}_{|\nbigd_i(\lambda_0,\epsilon_0)}
 \bigr)$.
Let $\vecb'$ be the element of $\real^l$
as follows:
\[
 q_j(\vecb')=
\left\{
 \begin{array}{ll}
 b_j & (j\neq i)\\
 c   & (j=i).
 \end{array}
\right.
\]
Then by definition of the filtration $\Fzero$,
we can take a holomorphic section $F\in\prolongg{\vecb'}{\nbige}$
such that $\pi(F_{|\nbigd_i(\lambda_0,\epsilon_0)})=f$,
where $\pi$ denotes the natural morphism
$\prolongg{\vecb'}{\nbige}_{|\nbigd_i(\lambda_0,\epsilon_0) }
\lrarr 
 \prolongg{\vecb}{\nbige}_{|\nbigd_i(\lambda_0,\epsilon_0) }$.
By definition,
$\lefttop{i}\DD(f)$ is the image of
$\DD F_{|\nbigd_i(\lambda_0,\epsilon_0) }$
via the morphism $\pi'$,
where $\pi'$ denotes the tensor product of
$\pi$ and the identity morphism of
$\Omega^{1,0}_{D_i} \bigl(\log (D\cap D_i)\bigr)$.
Thus $\lefttop{i}\DD$ preserves the filtration $\lefttop{i}F$.
By a similar argument,
we can show that $\Res_i(\DD)$ preserves the filtration
$\lefttop{i}F$.
\hfill\qed

\begin{cor}
 Assume
 $\lefttop{k}\deg^{F^{(\lambda_0)}}(v_i)
 <
  \lefttop{k}\deg^{F^{(\lambda_0)}}(v_j)$.
Then we have
$\nbiga^h_{i\,j\,|\,\nbigd_{k}}=0$ for any $h$.
\hfill\qed
\end{cor}

\begin{lem}
The endomorphism
$\Res_i(\DD)$ is flat with respect to
$\lefttop{i}\DD$.
\end{lem}
\pf
It immediately follows from the flatness of $\DD$.
\hfill\qed

%% file: 15.5.tex

\subsubsection{The $\EE$-decomposition}

\label{subsubsection;10.18.21}

For $\lambda\in\cnum_{\lambda}$,
we put
$\nbigt(\lambda,c,i):=\bigl\{u\in \KMS(\nbige^0,i)\,\big|\,
 \paramap(\lambda,u)=c\,\bigr\}$.
Since $\Res_i(\DD^{\lambda})$
preserves the parabolic filtration $F$
of $\prolongg{\vecb}{\nbigelambda}_{|D_i}$,
due to Lemma \ref{lem;9.10.96},
we have the induced action of $\Res_i(\DD^{\lambda})$
on $\lefttop{i}\Gr^{F}_c(\prolongg{\vecb}{\nbigelambda})$.

\begin{lem} \label{lem;9.11.2}
We have the generalized eigen decomposition
of $\lefttop{i}\Gr^F_c$ with respect to $\Res(\DD^{\lambda})$
as follows:
\[
 \lefttop{i}\Gr^F_c(\prolongg{\vecb}{\nbigelambda})
=\bigoplus_{u\in \nbigt(\lambda,c,i)}
 \EE\bigl(\lefttop{i}\Gr^F_c(\prolongg{\vecb}{\nbigelambda}),
 \eigenmap(\lambda,u)
 \bigr).
\]
The rank of $\EE\bigl(\lefttop{i}\Gr^F_c(\prolongg{\vecb}{\nbigelambda}),
 \eigenmap(\lambda,u)\bigr)$
is $\multiplicity(0,u)$.
\end{lem}
\pf
We have only to show 
$\rank \EE\bigl(\Gr^F_c(\prolongg{\vecb}{\nbigelambda})_{|P},
 \eigenmap(\lambda,u)\bigl)
=\multiplicity(0,u)$
for any $P\in D_i^{\circ}$.
It follows from the result in the case of curves
(Corollary \ref{cor;9.11.1}).
\hfill\qed

\vspace{.1in}
Let $c$ and $d$ be real numbers such that $d<c$.
The residue $\Res(\DD)$ induces the endomorphism
of
$\lefttop{i}F_c\big/\lefttop{i}F_d(\prolongg{\vecb}{\nbigelambda})$.
\begin{cor} \label{cor;9.11.4}
The set of the eigenvalues of $\Res(\DD)$
on
 $\lefttop{i}F_c\big/\lefttop{i}F_d(\prolongg{\vecb}{\nbigelambda})$
is as follows:
\[
 S=
 \Bigl\{\eigenmap(\lambda,u)\,
 \Big|\,u\in \bigcup_{d<a\leq c}\nbigt(\lambda,a,i)
 \Bigr\}.
\]
The rank of the generalized eigenspace corresponding to
$\beta\in S$ is given as follows:
\[
 \sum_{d<a\leq c}
 \sum_{\substack{u\in \nbigt(\lambda,a,i)\\
 \eigenmap(\lambda,u)=\beta
 }} \multiplicity(0,u).
\]
\end{cor}
\pf
It immediately follows from Lemma \ref{lem;9.11.2}.
\hfill\qed

\vspace{.1in}

Let us consider the generalized eigen decomposition
of $\prolongg{\vecb}{\nbigelambda}_{|D_i}$
with respect to $\Res_i(\DD^{\lambda})$.
We put as follows for any point $\lambda\in\cnum_{\lambda}$
and for any element $\vecb\in\real^l$:
\[
 \nbigk(\lambda,\vecb,i):=
 \bigl\{
 u\in \KMS(\nbige^0,i)\,\big|\,
 b_i-1<\paramap(\lambda,u)\leq b_i
 \bigr\}.
\]
Here $b_i$ denotes the $i$-th component of $\vecb$.

\begin{cor} \label{cor;9.11.5}
The set of the eigenvalues of the endomorphism
$\Res(\DD^{\lambda})$ on $\prolongg{\vecb}{\nbigelambda}_{|D_i}$
is given as follows:
\[
 \Sp(\nbigelambda,i)
=\bigl\{
 \eigenmap(\lambda,u)
 \,\big|\,
 u\in \nbigk(\lambda,\vecb,i)
 \bigr\}
\]
The rank 
of the generalized eigenspace corresponding to
$\beta\in\Sp(\nbigelambda,i)$ is given as follows:
\[
 \multiplicity(\lambda,\beta,i)
:=\sum_{\substack{u\in \nbigk(\lambda,\vecb,i)\\
 \eigenmap(\lambda,u)=\beta}}
 \multiplicity(0,u).
\]
\end{cor}
\pf
This is the special case of Corollary \ref{cor;9.11.4}.
\hfill\qed

\vspace{.1in}
\noindent
{\bf Notation}
We put
$\lefttop{i}\EE(\prolongg{\vecb}{\nbigelambda},\beta)=
 \EE\bigl(\prolongg{\vecb}{\nbigelambda}_{|D_i},\beta\bigr)$,
for simplicity.
\hfill\qed

\vspace{.1in}

Let us pick any point $\lambda_0\in\cnum_{\lambda}$.
Then there exist small positive numbers $\eta_2$ and $\epsilon_2$
such that
we have the following decomposition
into vector bundles on $\nbigd_i(\lambda_0,\eta_2)$,
due to Corollary \ref{cor;9.11.5}:
\[
 \prolongg{\vecb}{\nbige}_{|\nbigd_i(\lambda_0,\eta_2)}
=\bigoplus_{\beta\in\Sp(\nbigelambdazero,i)}
 \EE_{\epsilon_2}\bigl(
 \Res(\DD),\beta \bigr).
\]
\noindent
{\bf Notation}
We put as follows:
$ \lefttop{i}\EEzero\bigl(
 \prolongg{\vecb}{\nbige},\beta
 \bigr)
=\EE_{\epsilon_2}\bigl(
 \Res(\DD),\beta \bigr)$,
for simplicity.
\hfill\qed

\begin{lem} \label{lem;9.11.7}
Assume $\eta_2$ and $\epsilon_2$ are sufficiently small.
Let $\lambda$ be any element of $\Delta(\lambda_0,\eta_2)$.
Then we have the following decomposition:
\[
 \lefttop{i}\EEzero\bigl(
 \prolongg{\vecb}{\nbige},\beta
 \bigr)_{|\nbigd_i^{\lambda}}
=\bigoplus_{\substack{u\in \nbigk(\lambda,\vecb,i),\\
 \eigenmap(\lambda_0,u)=\beta
 }} \lefttop{i}\EE\bigl(
 \prolongg{\vecb}{\nbigelambda}_{|\nbigd_i(\lambda_0,\eta_2)},
 \eigenmap(\lambda,u) \bigr).
\]
\end{lem}
\pf
Since we have assumed $b_i\not\in\Par(\nbigelambda,i)$,
we have $\nbigk(\lambda,\vecb,i)=\nbigk(\lambda_0,\vecb,i)$.
If $\epsilon_2$ is sufficiently small,
$|\eigenmap(\lambda,u)-\beta|<\epsilon_2$
implies $\eigenmap(\lambda_0,u)=\beta$.
If $\eta_2$ is sufficiently small,
we have $\eigenmap(\lambda,u)\neq \eigenmap(\lambda,u')$
for $u,u'\in \nbigk(\lambda,\vecb,i)$ such that $u\neq u'$.
Thus we are done.
\hfill\qed

\begin{lem} \label{lem;9.11.6}
The filtration $\lefttop{i}\Fzero$ and
the decomposition $\lefttop{i}\EEzero$
are compatible.
\end{lem}
\pf
Since $\Res_i(\DD^{\lambda})$ preserves the parabolic filtration
$\lefttop{i}F$,
the parabolic filtration $\lefttop{i}F$
and the decomposition $\lefttop{i}E$ are compatible.
Then Lemma \ref{lem;9.11.6} follows from Lemma \ref{lem;9.11.7}
and Lemma \ref{lem;9.10.91}.
\hfill\qed

\begin{lem}
$\lefttop{i}\DD$ and $\Res_i(\DD)$ preserve
$\lefttop{i}\EE^{(\lambda_0)}$.
\end{lem}
\pf
As for $\Res_i(\DD)$, it is clear.
Since $\Res_i(\DD)$ is flat with respect to
$\lefttop{i}\DD$,
the generalized eigen decomposition is preserved
by $\lefttop{i}\DD$.
\hfill\qed

\vspace{.1in}

Let pick any point $\lambda_0\in\cnum_{\lambda}$.
Since $\Res_i(\DD)$ preserves the filtration $\lefttop{i}\Fzero$
due to Lemma \ref{lem;9.10.96},
we have the induced action of $\Res_i(\DD)$
on $\lefttop{i}\Gr^{F^{(\lambda_0)}}_c(\prolongg{\vecb}{\nbige})$.
Then there exist small positive numbers $\eta_2$ and $\epsilon_2$
such that
we have the following decomposition
into vector bundles on $\Delta(\lambda_0,\eta_2)\times D_i$,
due to Lemma \ref{lem;9.11.2}:
\[
 \lefttop{i}\Gr^{F^{(\lambda_0)}}_c(\prolongg{\vecb}{\nbige})
=\bigoplus_{u\in\nbigt(\lambda_0,c,i)}
 \EE_{\epsilon_2}
 \bigl(\lefttop{i}\Gr^{F^{(\lambda_0)}}_c(\prolongg{\vecb}{\nbige}),\,\,
 \eigenmap(\lambda_0,u)
 \bigr).
\]

\begin{lem} \label{lem;8.19.10}
Assume $\epsilon_2$ and $\eta_2$ are sufficiently small.
Let $\lambda$ be any point of $\Delta^{\ast}(\lambda_0,\eta_2)$
and $P$ be any point of $D_i$.
The subspace
$\EE_{\epsilon_2}
 \bigl(\lefttop{i}\Gr^{\Fzero}_c(\prolongg{\vecb}{\nbige}),
 \eigenmap(\lambda_0,u)
 \bigr)_{|(\lambda,P)}$
of $\lefttop{i}\Gr^{\Fzero}_c(\prolongg{\vecb}{\nbige})_{|(\lambda,P)}$
is the generalized eigenspace of the induced endomorphism
$\Res(\DD^{\lambda})$ corresponding to
the eigenvalue $\eigenmap(\lambda,u)$.
\end{lem}
\pf
Since we have the isomorphism (\ref{eq;9.11.3}),
the set of 
the eigenvalue of the induced morphism $\Res(\DD^{\lambda})$ 
on $\lefttop{i}\Gr^{\Fzero}_c(\prolongg{\vecb}{\nbige})_{|(\lambda,P)}$
is as follows, due to Corollary \ref{cor;9.11.4}:
\[
 \bigl\{
 \eigenmap(\lambda,u)\,\big|\,
 u\in \nbigt(\lambda_0,c,i)
 \bigr\}.
\]
Let $u$ and $u'$ be elements of $\nbigt(\lambda_0,c,i)$.
If $\eta_2$ is sufficiently small,
$|\eigenmap(\lambda,u)-\eigenmap(\lambda,u')|<\epsilon_2$
if and only if
$\eigenmap(\lambda,u)=\eigenmap(\lambda,u')$.
Since we have $\paramap(\lambda_0,u)=\paramap(\lambda_0,u')=c$,
the condition $\eigenmap(\lambda,u)=\eigenmap(\lambda,u')$
implies $u=u'$.
Hence the condition
$|\eigenmap(\lambda,u)-\eigenmap(\lambda,u')|<\epsilon_2$
implies $\eigenmap(\lambda,u)=\eigenmap(\lambda,u')$.
we obtain the result.
\hfill\qed

\vspace{.1in}
\noindent
{\bf Notation}
Let
$\lefttop{i}\Gr^{F^{(\lambda_0)},\EE^{(\lambda_0)}}_{\kmsmap(\lambda_0,u)}$
denote the space
$\lefttop{i}\EEzero
 \bigl(\lefttop{i}\Gr^{\Fzero}_{\paramap(\lambda_0,u)},
 \eigenmap(\lambda_0,u)\bigr)$
for $u\in \KMS(\nbige^0,i)$.
\hfill\qed


\subsubsection{Weak norm estimate}
\label{subsubsection;a12.1.5}

Pick a frame $\vecv$ of $\prolongg{\vecb}{\nbige}$
over $\nbigx(\lambda_0,\epsilon_0)$,
which is compatible with $F^{(\lambda_0)}$ and $\EE^{(\lambda_0)}$.
For each $v_j$ and for each $i$,
we have 
the unique element $u_i(v_j)\in\KMS(\nbige^0,i)$
satisfying the following:
\[
 \lefttop{i}\deg^{F^{(\lambda_0)},\EE^{(\lambda_0)}}
 (v_{j\,|\,\nbigx^{\lambda_0}})
=\kmsmap(\lambda_0,u_i(v_j)).
\]

\begin{lem}
Then $\lefttop{i}\deg^{F^{(\lambda)},\EE^{(\lambda)}}
 (v_{j\,|\,\nbigx^{\lambda}})
=\kmsmap(\lambda,u_i(v_j)) $ also holds for any $\lambda$.
\end{lem}
\pf
We have the following for some sufficiently small positive number
$\epsilon_3$:
\[
\max\Bigl\{
 \bigl|
 \lefttop{i}\deg^{F}(v_{j\,|\,\nbigx^{\lambda}})
-\lefttop{i}\deg^{F}(v_{j\,|\,\nbigx^{\lambda_0}})
 \bigr|,\,\,
  \bigl|
 \lefttop{i}\deg^{\EE}(v_{j\,|\,\nbigx^{\lambda}})
-\lefttop{i}\deg^{\EE}(v_{j\,|\,\nbigx^{\lambda_0}})
 \bigr|
 \Bigr\}<\epsilon.
\]
Since the unique element of $\KMS(\nbigelambda,i)$ satisfying
such condition is $\kmsmap(\lambda,u_i(v_j))$,
we are done.
\hfill\qed

\vspace{.1in}

\label{subsubsection;9.11.11}

Let us consider the $C^{\infty}$-frame $\vecv'$
of $\nbige$ over $\Delta(\lambda_0,\eta_2)\times (X-D)$,
given as follows:
\[
 v_j':=
 v_j\cdot
 \prod_{i=1}^l |z_i|^{\paramap(\lambda,u_i(v_j))},
\quad\quad
 \vecv':=(v_j').
\]
By a standard argument,
we obtain the following.
\begin{prop} \label{prop;10.16.2}
The frame $\vecv'$ is adapted up to log order.
\end{prop}
\pf
It is easy to see that $H(h,\vecv')\leq
C\cdot\bigl(-\sum_{i=1}^l\log|z_i|\bigr)^M$
for some positive constants $M$ and $C$,
which follows from Lemma \ref{lem;8.19.10}.
By considering the dual frame $\vecv^{\lor\,\prime}$,
we obtain $H(h,\vecv')\geq
C'\cdot\bigl(-\sum_{i=1}^l\log|z_i|\bigr)^{-M'}$, as usual
(See the proof of Proposition \ref{prop;9.11.10}, for example).
Thus we are done.
\hfill\qed

%% file: 15.6.tex

\subsubsection{Some functoriality}

Let consider the functoriality for tensor product.
Let $(E^{(a)},\delbar_{E^{(a)}},h^{(a)},\theta^{(a)})$ $(a=1,2)$
be tame harmonic bundles. We obtain the deformed holomorphic bundles
$\nbige^{(a)}$.
Let pick $\lambda_0\in\cnum$.
For simplicity, we consider the following situation:
\begin{itemize}
\item
 $0\not\in\Par(\prolong{\nbige^{(a)\,\lambda}},i)$
for $a=1,2$ and for $i=1,\ldots,l$.
\item
For any $\lambda\in\Delta(\lambda_0,\epsilon_0)$,
the sets $\Par(\prolong{\nbige^{(a)\,\lambda}},i)$ $(a=1,2)$ are
$\eta_{a\,i}$-small such that $\eta_{1\,i}+\eta_{2\,i}<1$
for $i=1,\ldots,l$.
\end{itemize}

Take $\vecv^{(a)}$ be frames of $\prolong{\nbige^{(a)}}$
over $\Delta(\lambda_0,\epsilon_0)$,
which are compatible with the filtration
$\lefttop{i}F^{(\lambda_0)}$
and the decomposition $\lefttop{i}\EE^{(\lambda_0)}$
$(i=1,\ldots,l)$.

\begin{lem}
$\vecv^{(1)}\otimes\vecv^{(2)}$ is a frame of
$\prolong{\bigl( \nbige^{(1)}\otimes\nbige^{(2)}\bigr)}$,
which is compatible with the filtration
$F^{(\lambda_0)}$ and the decomposition $\EE^{(\lambda_0)}$
\end{lem}
\pf
We obtain the frame $\vecv^{(a)\,\prime}$
from $\vecv^{(a)}$,
which is adapted up to log order.
Then $\vecv^{(1)\,\prime}\otimes\vecv^{(2)\,\prime}$
is adapted up to log order,
as in the subsubsection \ref{subsubsection;9.11.11}.
Due to our assumption,
we obtain that $\vecv^{(1)}\otimes\vecv^{(2)}$
gives a frame of $\prolong{\big(\nbige^{(1)}\otimes\nbige^{(2)}\big)}$.
Then the compatibilities with
the filtration $F$ and the decomposition $\EE$ are clear.
\hfill\qed

\begin{cor}
 Let $\eta$ be a positive number and $R$ be positive integer
 such that $R\cdot \eta<1$.
 Assume the following for simplicity:
 For any $\lambda\in\Delta(\lambda_0,\epsilon_0)$,
 the sets $\Par(\prolong{\nbigelambda},i)$ are $\eta$-small.

Then we have the following canonical isomorphism:
\[
 \prolong{\bigwedge^R\nbige}\simeq
 \bigwedge^R\prolong{\nbige},
\quad\quad
 \prolong{\Sym^R\nbige}\simeq
 \Sym^R(\prolong{\nbige}).
\]
\hfill\qed
\end{cor}

%% file: 15.7.tex

\subsubsection{KMS-spectrum}

Let $I$ be a subset of $I$
and $\vecc$ be an element of $\real^I$.
Due to the compatibility of the filtrations
$\bigl(\lefttop{i}\Fzero\,|\,i\in I\bigr)$,
we have the following vector bundles
on $\nbigd_I(\lambda_0,\epsilon_0)$:
\[
 \lefttop{I}F^{(\lambda_0)}_{\vecc}
 =\bigcap_{i\in I}\lefttop{i}F^{(\lambda_0)}
 _{c_i\,|\,\nbigd_I(\lambda_0,\epsilon_0)},
\quad\quad
 \lefttop{I}\Gr^{F^{(\lambda_0)}}_{\vecc}=
 \frac{\lefttop{I}F^{(\lambda_0)}_{\vecc}}
  {\sum_{\vecc'\lneq \vecc}
     \lefttop{I}F^{(\lambda_0)}_{\vecc'}}.
\]
On the vector bundle $\lefttop{I}\Gr^{F^{(\lambda_0)}}_{\vecc}$,
we have the induced endomorphisms
$\lefttop{I}\Gr^{F^{(\lambda_0)}}_{\vecc}\Res_i(\DD)$ $(i\in I)$.
Hence we obtain the tuple of endomorphisms:
\[
 \lefttop{I}\Gr^{F^{(\lambda_0)}}_{\vecc}(\Res_I(\DD))
:=
 \bigl(
  \lefttop{I}\Gr^{F^{(\lambda_0)}}_{\vecc}(\Res_i(\DD))\,
 \big|\,i\in I
 \bigr).
\]
We have the subset
$\Sp\bigl(\lefttop{I}\Gr^{\Fzero}_{\vecc}(\Res_I(\DD))\bigr)
 \subset\cnum^I$,
and we have the following decomposition:
\[
 \lefttop{I}\Gr^{\Fzero}_{\vecc}
=\bigoplus_{\vecgamma\in
   \Sp\bigl(\lefttop{I}\Gr^{\Fzero}_{\vecc}(\Res_I(\DD))\bigr)}
 \lefttop{I}\EEzero\bigl(
  \lefttop{I}\Gr^{F^{(\lambda_0)}}_{\vecc},\vecgamma
 \bigr),
\quad\quad
 \lefttop{I}\EEzero\bigl(
  \lefttop{I}\Gr^{F^{(\lambda_0)}}_{\vecc},\vecgamma
 \bigr)
=\bigcap_{ i\in I}
 \lefttop{i}
 \EEzero\bigl(
 \lefttop{I}\Gr^{F^{(\lambda_0)}}_{\vecc},q_i(\vecgamma)
 \bigr).
\]
For a pair $\vecu=(\vecc,\vecgamma)$
such that $\vecc\in \prod_{i\in I}\Par(\nbigelambdazero,i)$
and $\vecgamma\in \prod_{i\in I}\Sp(\nbigelambdazero,i)$,
we put as follows:
\[
 \lefttop{I}\Gr^{\Fzero,\EEzero}_{\vecu}(\nbige)
=
\lefttop{I}\EE^{\Fzero}(
  \lefttop{I}\Gr^{\Fzero}_{\vecc}\Res_I(\DD),\vecgamma
 ).
\]
We obtain the following subset:
\[
 \KMS(\nbigelambda,I)
:=\big\{\vecu\in\real^I\times\cnum^I
 \,\big|\,
\lefttop{I}\Gr^{\Fzero,\EEzero}_{\vecu}(\nbige)\neq 0
 \bigr\}
\subset \prod_{i\in I}\KMS(\nbigelambda,i).
\]
For any element $\vecu\in\KMS(\nbigelambda,I)$,
we put as follows:
\[
 \multiplicity(\lambda,\vecu):=
 \dim \lefttop{I}\Gr^{\Fzero,\EEzero}_{\vecu}(\nbige).
\]
The number $\multiplicity(\lambda,\vecu)$ is called
the multiplicity of $\vecu$.
We have the natural $\seisuu^I$-action 
on $\KMS(\nbigelambda,I)$ which preserves the multiplicities.
We put as follows:
\[
 \overline{\KMS}(\nbigelambda,I)
=\KMS(\nbigelambda,I)\big/\seisuu^I,
\quad\quad
 \KMS(\prolongg{\vecb}{\nbigelambda},I)
:=\bigl\{\vecu\in\KMS(\nbigelambda,I)\,\big|\,
 b_i-1<c_i\leq b_i \bigr\}.
\]

Recall that we have the morphism
$\kmsmap(\lambda):\KMS(\nbige^0,i)\lrarr \KMS(\nbigelambda,i)$.
It induces the morphism
$\kmsmap(\lambda):
 \prod_{i\in I}\KMS(\nbige^0,i)\lrarr
 \prod_{i\in I}\KMS(\nbigelambda,i)$.

\begin{prop} \label{prop;9.11.12}
The morphism above induces
$\kmsmap_I(\lambda):\KMS(\nbige^0,I)\lrarr\KMS(\nbigelambda,I)$,
and we have equality
$\multiplicity(0,\vecu)
=\multiplicity(\lambda,\kmsmap_I(\lambda,\vecu))$.
\end{prop}
\pf
Let $\vecu$ be an element of $\prod_{i\in I}\KMS(\nbige^0,i)$,
then we obtain
$\kmsmap(\lambda_0,\vecu)\in \prod_{i\in I}\KMS(\nbigelambdazero,i)$.
Assume $\kmsmap(\lambda_0,\vecu)\in \KMS(\nbigelambdazero,i)$.
Then we have the vector bundle over $\nbigd_I(\lambda_0,\epsilon_0)$:
\[
 \Gr^{F^{(\lambda_0)},\EE^{(\lambda_0)}}_{\kmsmap(\lambda_0,\vecu)}(\nbige).
\]
We use the following lemma.
\begin{lem}\label{lem;9.11.13}
The following holds for any $\lambda\in\Delta(\lambda_0,\epsilon_0)$:
\[
 \Gr^{F^{(\lambda_0)},\EE^{(\lambda_0)}}_{\kmsmap(\lambda_0,\vecu_0)}
  (\nbige)_{|\nbigd_I^{\lambda}}
=\Gr^{F,\EE}_{\kmsmap(\lambda,\vecu)}(\nbigelambda).
\]
\end{lem}
\pf
We have only to check the case $\dim(X)=1$.
It follows from Lemma \ref{lem;8.19.10}.
\hfill\qed

\vspace{.1in}
To show Proposition \ref{prop;9.11.12},
let us consider the following condition for $\lambda$ and
$\vecu\in \prod_{i\in I}\KMS(\nbige^0,i)$:
\begin{quote}
 $P(\lambda,\vecu)$:
 $\kmsmap(\lambda,\vecu)\in \KMS(\nbigelambda,I)$.
\end{quote}

\begin{lem}
Let us pick $\vecu\in\prod_{i\in I}\KMS(\nbige^0,i)$.
For any $\lambda$, there exists a positive number $\eta(\lambda)$
such that the following two conditions are equivalent:
\begin{itemize}
\item
$P(\lambda,\vecu)$ holds
\item
$P(\lambda',\vecu)$ holds
for some $\lambda'\in\Delta(\lambda,\eta(\lambda))$.
\end{itemize}
\end{lem}
\pf
It follows from Lemma \ref{lem;9.11.13}.
\hfill\qed

\vspace{.1in}
Let us return to the proof of Proposition \ref{prop;9.11.12}.
Let us take any point $\lambda\in\cnum_{\lambda}$.
We have the line $\ell=\{t\cdot\lambda\,|\,0\leq t\leq 1\}$
in the plane $\cnum_{\lambda}$.
We can pick a finite subset $S=\{t_1,\ldots,t_h\}\subset [0,1]$
such that
$\ell\subset
 \bigcup_{t_i\in S}\Delta\bigl(t_i\cdot\lambda,\eta(t_i\cdot\lambda)\bigr)$.
Then it follows the equivalence of
$P(0)$ and $P(\lambda)$.
Thus the proof of Proposition \ref{prop;9.11.12} is accomplished.
\hfill\qed

%% file: a57.tex

\subsubsection{The induced vector bundles $\lefttop{\lbar}\nbigg_{\vecu}$
 on $\nbigd_{\lbar}$}
\label{subsubsection;10.16.10}

Let $\vecu$ be an element of $\KMS(\nbige^0,\lbar)$.
Let us pick $\lambda_0\in\cnum_{\lambda}$.
We put
$(\vecb,\vecbeta):=\kmsmap(\lambda_0,\vecu)\in\real^l\times\cnum^l$.
Let us pick any sufficiently small $\epsilon_1>0$ such that
$b_i+\epsilon_1\not\in\Par(\nbigelambdazero,i)$
for any $i$.
By considering the prolongment
$\prolongg{\vecb+\epsilon_1\cdot\vecdelta}{\nbige}$,
we obtain the vector bundle over
$D_{\lbar}\times\Delta(\lambda_0,\epsilon_0)$:
\[
 \lefttop{\lbar}\nbigg_{\vecu}^{(\lambda_0)}:=
 \lefttop{\lbar}\Gr^{\Fzero,\EEzero}_{
 \kmsmap(\lambda_0,\vecu) }(\nbige).
\]
Clearly,
it does not depend on a choice of $\epsilon_1$
on a neighbourhood of $\lambda_0$.

Let us pick $\lambda\in\Delta(\lambda_0,\epsilon_0)$.
We put
$\bigl(\vecb(\lambda),\vecbeta(\lambda)\bigr):=\kmsmap(\lambda,\vecu)$.
Let us pick $\epsilon_0'$ as
$\Delta^{\ast}(\lambda_0,\epsilon_0)\supset
 \Delta(\lambda,\epsilon_0')$.
We have the following vector bundles on
$\nbigd^{\circ}_{\lbar}(\lambda,\epsilon_0')$:
\[
 \lefttop{\lbar}\nbigg^{(\lambda)}_{\vecu},
\quad\quad
 \lefttop{\lbar}
 \nbigg^{(\lambda_0)}_{\vecu\,|\,\Delta(\lambda,\epsilon_0')}.
\]
\begin{lem}
The vector bundles $\lefttop{\lbar}\nbigg^{(\lambda)}_{\vecu}$
and $\lefttop{\lbar}
 \nbigg^{(\lambda_0)}_{\vecu\,|\,\Delta(\lambda,\epsilon_0')}$
are naturally isomorphic.
\end{lem}
\pf
Note that we have the following decomposition:
\[
 \lefttop{\lbar}\EE^{(\lambda_0)}\bigl(
 \prolongg{\vecb+\epsilon\cdot\vecdelta}{\nbige}_{
 |\nbigd_{\lbar}(\lambda_0,\epsilon_0)
 },\,\vecbeta
 \bigr)_{|\nbigd_{\lbar}^{\circ\,\ast}(\lambda_0,\epsilon_0)}
=\bigoplus_{
\substack{
 \vecu'\in\KMS(\nbige^0,\lbar),\\
 \kmsmap(\lambda_0,\vecu')=\vecbeta}}
 \lefttop{\lbar}\EE\bigl(
 \prolongg{\vecb+\epsilon\cdot\vecdelta}{\nbige}_{
 \nbigd_{\lbar}^{\circ\,\ast}(\lambda_0,\epsilon_0)},\,
 \eigenmap(\lambda,\vecu')
 \bigr).
\]
We have the following:
\[
 \lefttop{\lbar}\nbigg^{(\lambda_0)}_{
  |\nbigd^{\circ}_{\lbar}(\lambda,\epsilon_0')}=
 \lefttop{\lbar}\Gr^{F^{(\lambda_0)}}_{\vecb}
 \bigl(\EE^{(\lambda_0)}(\prolongg{\vecb+\epsilon\cdot\vecdelta}
 {\nbige}_{\nbigd^{\circ}_{\lbar}(\lambda_0,\epsilon_0) },\,\vecbeta
 )
 \bigr)_{|\nbigd^{\circ}_{\lbar}(\lambda,\epsilon_0')}
\simeq
 \EE^{(\lambda)}\bigl(\prolongg{\vecb+\epsilon\cdot\vecdelta}{
 \nbige }_{|\nbigd^{\circ}_{\lbar}(\lambda,\epsilon_0')},\,
 \eigenmap(\lambda,\vecu) \bigr).
\]
On the other hand,
we have $\lefttop{\lbar}\nbigg^{(\lambda)}_{\vecu}= 
 \EE^{(\lambda)}\bigl(\prolongg{\vecb+\epsilon\cdot\vecdelta}{
 \nbige }_{|\nbigd^{\circ}_{\lbar}(\lambda,\epsilon_0')},\,
 \eigenmap(\lambda,\vecu) \bigr)$.
Thus we are done.
\hfill\qed

\vspace{.1in}

Thus we obtain the vector bundle
$\lefttop{\lbar}\nbigg_{\vecu}$ over $\nbigd_{\lbar}$.
When we distinguish the dependence of $\lefttop{\lbar}\nbigg_{\vecu}$
on the harmonic bundle $\harmonicbundle$,
we use the notation
$\lefttop{\lbar}\nbigg_{\vecu}\harmonicbundle$,
or simply $\lefttop{\lbar}\nbigg_{\vecu}(E)$.

\subsubsection{The induced frame}

Let $\vecv=(v_i)$ be a compatible frame of
$\prolongg{\vecb}{\nbige}$ over $\nbigx(\lambda_0,\epsilon_0)$.
For each section $v_i$,
we have the element $\vecu(v_i)\in \KMS(\nbige^0,\lbar)$
such that
$\deg^{\EEzero,\Fzero}(v_i)=\kmsmap(\lambda_0,\vecu(v_i))$.
For any element $\vecu\in\KMS\bigl(\nbige^0,\lbar\bigr)$,
we put as follows:
\[
 \vecv_{\vecu}
:=\bigl(z^{-\vecn}\cdot v_i\,\big|\,
 \vecu(v_i)+\vecn=\vecu
 \bigr).
\]
Note that $\vecn$ is determined by
$\vecu$ and $\vecb$, not by $v_i$.

\begin{lem} \label{lem;10.24.1}
The tuple of sections $\vecv_{\vecu}$
induces the frame of $\lefttop{\lbar}\nbigg^{(\lambda_0)}_{\vecu}$.
\end{lem}
\pf
First we remark the following:
Let $\vecn$ be an element of $\seisuu^l$.
We put $\tilde{\vecv}:=\bigl(z^{-\vecn}v_i\bigr)$.
Then $\tilde{\vecv}$ gives the frame of
$\prolongg{\vecb+\vecn}{\nbige}$
over $\nbigx(\lambda_0,\epsilon_0)$.

It is easy to see that $\vecv_{\vecu}$ induces
the tuple of sections of
$\lefttop{\lbar}\nbigg^{(\lambda_0)}_{\vecu}$.
By using the remark above,
we can show that $\vecv_{\vecu}$ is a frame.
\hfill\qed

\subsubsection{The nilpotent map $\nbign_{i,\vecu}$ and the pairing with
   the dual}
\label{subsubsection;10.24.2}

We have the endomorphism 
$\Res_i(\DDlambda)$ on $\nbigg_{\vecu}$.
The unique eigenvalue of $\Res_i(\DDlambda)$
on $\lefttop{\lbar}\nbigg_{\vecu\,|\,\lambda}$ is
$\eigenmap\bigl(\lambda,q_i(\vecu)\bigr)$.
Hence the nilpotent part
$\nbign_{i\,\vecu}:=
 \Res_i(\DDlambda)-\eigenmap\bigl(\lambda,q_i(\vecu)\bigr)$ 
gives the holomorphic endomorphism.

Let $\lambda_0$ be an element of $\cnum_{\lambda}$
and $\epsilon_0>0$.
We have the naturally defined morphism
$\prolongg{\vecb}{\nbige}\otimes
 \big(\prolongg{-\vecb+(1-\epsilon)\vecdelta}\nbige^{\lor}\big)
\lrarr\nbigo$
over $\Delta(\lambda_0,\epsilon_0)$,
which is the morphism of $\lambda$-connections.
It is easy to see that
$\lefttop{i}\EE^{(\lambda_0)}$ and
$\lefttop{i}F^{(\lambda_0)}$ are preserved
over $\nbigd_i(\lambda_0,\epsilon_0)$.
Thus we obtain the following morphism
over $\nbigd_{\lbar}(\lambda_0,\epsilon_0)$:
\[
 \lefttop{\lbar}\Gr^{F^{(\lambda_0)},\EE^{(\lambda_0)}}
 _{\kmsmap(\lambda_0,\vecu)}(\nbige)
\otimes
 \lefttop{\lbar}\Gr^{F^{(\lambda_0)},\EE^{(\lambda_0)}}
 _{\kmsmap(\lambda_0,-\vecu)}(\nbige^{\lor})
\lrarr\nbigo_{\nbigd_{\lbar}(\lambda_0,\epsilon_0)}.
\]
Then we obtain the morphism
$S:
 \lefttop{\lbar}\nbigg_{\vecu}(E)
  \otimes\lefttop{\lbar}\nbigg_{-\vecu}(E^{\lor})
 \lrarr\nbigo_{\nbigd_{\lbar}}$.

\begin{lem}
We have the equality
$S\bigl(\nbign_{i\,\vecu}\otimes\id\bigr)
+S\bigl(\id\otimes\nbign_{i,-\vecu}\bigr)=0$.
\end{lem}
\pf
It follows from
$S\bigl(\DD\otimes\id\bigr)+S\bigl(\id\otimes\DD\bigr)
=\DD\circ S$.
\hfill\qed

\subsubsection{Functoriality for dual}
\label{subsubsection;10.26.20}

Due to the results in the subsubsection \ref{subsubsection;10.24.2},
we obtain the naturally defined morphism
$\lefttop{\lbar}\nbigg_{-\vecu}(E^{\lor})
\lrarr
 \lefttop{\lbar}\nbigg_{\vecu}(E)^{\lor}$.
\begin{lem}
The naturally defined morphism
$\lefttop{\lbar}\nbigg_{-\vecu}(E^{\lor})
\lrarr
 \lefttop{\lbar}\nbigg_{\vecu}(E)^{\lor}$
is isomorphic.
\end{lem}
\pf
Let $\vecv$ be a frame of $\prolongg{\vecb}{\nbige}$,
which is compatible with $\EEzero$ and $\Fzero$.
The dual frame $\vecv^{\lor}$ gives the frame
of $\prolongg{-\vecb+(1-\epsilon)\vecdelta}{\nbige}$
for some positive constant $\epsilon$.
It is also compatible with $\EEzero$ and $\Fzero$.
By using Lemma \ref{lem;10.24.1},
$\vecv_{\vecu}$ and $\bigl(\vecv^{\lor}\bigr)_{-\vecu}$ 
induce the frames of $\lefttop{\lbar}\nbigg_{\vecu}(E)$
and $\lefttop{\lbar}\nbigg_{-\vecu}(E^{\lor})$ respectively.
Thus we are done.
\hfill\qed

\subsubsection{Functoriality for tensor products}
\label{subsubsection;10.26.21}

Let $(E_i,\delbar_{E_i},\theta_i,h_i)$ $(i=1,2)$ be
tame harmonic bundles over $X-D$.
Let $\vecu_i$ $(i=1,2)$ be elements of
$KMS(\nbige^0_i,\lbar)$.
We have the induced morphisms on $\nbigx(\lambda_0,\epsilon_0)$:
\[
 \lefttop{\lbar}F_{\paramap(\lambda_0,\vecu_1)}(\nbige_1)
\otimes
 \lefttop{\lbar}F_{\paramap(\lambda_0,\vecu_2)}(\nbige_2)
\lrarr
 \lefttop{\lbar}F_{\paramap(\lambda_0,\vecu_1+\vecu_2)}
 (\nbige_1\otimes\nbige_2),
\]
\[
 \lefttop{\lbar}\Gr^F_{\paramap(\lambda_0,\vecu_1)}(\nbige_1)
\otimes
 \lefttop{\lbar}\Gr^F_{\paramap(\lambda_0,\vecu_2)}(\nbige_2)
\lrarr
 \lefttop{\lbar}\Gr^F_{\paramap(\lambda_0,\vecu_1+\vecu_2)}
 (\nbige_1\otimes\nbige_2).
\]
Since the morphism is compatible with the residues
of the $\lambda$-connections,
we obtain the following induced morphisms:
\[
 F_{\vecu_1,\vecu_2}:
 \lefttop{\lbar}\nbigg_{\vecu_1}^{(\lambda_0)}
 \otimes\lefttop{\lbar}\nbigg_{\vecu_2}^{(\lambda_0)}
\lrarr
 \lefttop{\lbar}\nbigg_{\vecu_1+\vecu_2}^{(\lambda_0)}.
\]

\begin{lem}
The morphism $F_{\vecu_1,\vecu_2}$ is compatible with
the nilpotent morphisms
of the induced vector bundles.
\end{lem}
\pf
It is clear from our construction.
\hfill\qed

\vspace{.1in}

Let $\vecb_i$ be elements of $\real^l$.
\begin{lem}
Let $\vecu$ be an element of $\KMS(\nbige^0_1\otimes\nbige^0_2,\lbar)$.
We have the isomorphism:
\[
 \bigoplus_{\substack{
 \vecu_i\in\KMS(\prolongg{\vecb_i}{\nbige^0},\lbar),\\
 \vecu_1+\vecu_2=\vecu
 }
 }
 \lefttop{\lbar}\nbigg_{\vecu_1}(\nbige_1)
\otimes
 \lefttop{\lbar}\nbigg_{\vecu_2}(\nbige_2)
\simeq
 \lefttop{\lbar}\nbigg_{\vecu}(\nbige_1\otimes\nbige_2).
\]
\end{lem}
\pf
We can show the lemma by using Lemma \ref{lem;10.24.1}.
\hfill\qed

\subsubsection{Functoriality for pull backs}
\label{subsubsection;10.26.6}

We put $X^{(1)}:=\Delta_z^n$ and $X^{(2)}:=\Delta_{\zeta}^{m}$.
We put $D^{(2)}_i:=\{\zeta_i=0\}$ 
and $D^{(2)}:=\bigcup_{i=1}^l D_i^{(2)}$.
Let $\harmonicbundle$ be a tame harmonic bundle
over $X^{(2)}-D^{(2)}$.

Let $\vecc=\big(c_{j\,i}\,\big|\,1\leq j\leq n,\,\,1\leq i\leq m\big)$
be an element of $\seisuu_{>0}^{n\cdot m}$.
Let us consider the morphism
$\psi:X_1\lrarr X_2$ given as follows:
\[
 \psi^{\ast}(\zeta_i):=\prod_{j=1}^n z_j^{c_{j\,i}}.
\]
We put $D^{(1)}:=\psi^{-1}D^{(2)}$.
We obtain the harmonic bundle
$\psi^{-1}\harmonicbundle$ over $X^{(1)}-D^{(1)}$.

Let us pick a point $\lambda_0\in\cnum_{\lambda}$
and a sufficiently small positive number $\epsilon_0$
such that $\prolongg{\vecb}{\nbige}$ is locally free
on $\nbigx^{(2)}(\lambda_0,\epsilon_0)$
for some $\vecb\in\real^{l}$.
Let $\vecv=(v_i)$ be a frame of $\prolongg{\vecb^{(2)}}{\nbige}$
which is compatible with $\Fzero$ and $\EEzero$.
We have the elements 
$\vecu(v_i)\in\KMS(\nbige^0,\lbar)$
for each $v_i$ satisfying the following:
\[
 \deg^{\EEzero,\Fzero}(v_i)=\kmsmap(\lambda_0,\vecu(v_i)).
\]
We put as follows:
\[
 v_i':=v_i\cdot \prod_{k=1}^l|\zeta_k|^{\paramap(\lambda,u_k(v_i))},
\quad
 \vecv'=(v_i').
\]
We have already seen that $C^{\infty}$-frame $\vecv'$ on
$\nbigx^{(2)}(\lambda_0,\epsilon_0)-\nbigd^{(2)}(\lambda_0,\epsilon_0)$
is adapted up to log order.

We obtain the holomorphic frame
$\psi^{-1}\vecv$ 
and the $C^{\infty}$-frame
$\psi^{-1}\vecv'$
of $\psi^{-1}\nbige$ 
over
$\nbigx^{(1)}(\lambda_0,\epsilon_0)
-\nbigd^{(1)}(\lambda_0,\epsilon_0)$.
Note that $\psi^{-1}\vecv'$ is adapted up to log order.

We put as follows:
\[
 \vecc\cdot\vecu(v_i)
:=\Bigl(
 \sum_{k}c_{1\,k}\cdot u_k(v_i),
  \ldots,
 \sum_kc_{n\,k}\cdot u_k(v_i)
 \Bigr)\in (\real\times\cnum)^n.
\]
The elements $\vecn(v_i)\in\seisuu^{n}$ are determined by
the following conditions:
\[
 \vecb^{(1)}-\vecdelta<
 \paramap(\lambda_0,\vecc\cdot\vecu(v_i))+\vecn_i
 \leq \vecb^{(1)}.
\]
We put as follows:
\[
 w_i:=\psi^{\ast}v_i\cdot \prod z_j^{-n_j},
\quad
 \vecw=(w_i).
\]
Then $\vecw$ is a tuple of sections
of $\prolongg{\vecb^{(1)}+\eta\cdot\vecdelta}{\psi^{\ast}\nbige}$
over $\nbigx^{(1)}(\lambda_0,\epsilon_1)$.
Here $\eta$ denotes any small positive number,
and a small positive number $\epsilon_1$ depends $\eta$.

We put
$\vecd(w_i)
:=\paramap\bigl(\lambda_0,\vecc\cdot\vecu(v_i)\bigr)+\vecn_i$.
We put as follows:
\[
 w'_i:=w_i\cdot \prod_{j=1}^n|z_j|^{d_j(w_i)},
\quad
 \vecw':=(w_i').
\]
Thus we obtain $C^{\infty}$-frame $\vecw'$
of $\nbige$ over
$\nbigx^{(1)}(\lambda_0,\epsilon_1)-\nbigd^{(1)}(\lambda_0,\epsilon_1)$.
\begin{lem}
The tuple $\vecw'$ is adapted up to log order.
\end{lem}
\pf
It is easy to see that we have some $C^{\infty}$-functions
$f_i$ such that
$w'_i=f_i\cdot \psi^{-1}v'_i$ and $|f_i|=1$ hold.
Since $\psi^{-1}\vecv$ is adapted up to log order,
we obtain the lemma.
\hfill\qed

\begin{lem}\label{lem;10.26.5}
The tuple $\vecw$ gives the frame of
$\prolongg{\vecb^{(1)}+\eta\cdot\vecdelta}{\psi^{\ast}\nbige}$.
It is compatible with $\Fzero$ and $\EEzero$.
\end{lem}
\pf
The first claim and the compatibility with $\FFzero$
follow from Lemma \ref{lem;9.5.10}
and Lemma \ref{lem;10.10.31}.

Let $A=\sum A_k\cdot d\zeta_k/\zeta_k$ denote the $\lambda$-connection
form of $\DD$ with respect to $\vecv$,
i.e., $\DD\vecv=\vecv\cdot A$.
Then we have
$\psi^{\ast}\DD\psi^{\ast}\vecv=\psi^{\ast}\vecv\cdot
 \psi^{\ast}A$.
We have the following:
\[
 \psi^{\ast}A=
 \sum_i\sum_j \frac{dz_j}{z_j}\cdot c_{j\,i}\cdot \psi^{\ast}A_i.
\]
Thus we obtain the following:
\[
 \psi^{\ast}\DD\cdot\vecw
=\vecw\cdot
  \sum_j
\Bigl(
\sum_k c_{j\,k}\cdot \psi^{\ast}A_k
+N_j
\Bigr)\cdot \frac{dz_j}{z_j}.
\]
Here $N$ denotes the diagonal matrix whose $i$-th component
is $-n_j(v_i)$.
Then it is easy to see that $\vecw$ is compatible with 
the decompositions $\lefttop{i}\EEzero$.
\hfill\qed

\vspace{.1in}

We have the naturally defined morphism
$\psi^{\ast}\lefttop{\lbar}\nbigg_{\vecu}(E)\lrarr
 \lefttop{\lbar}\nbigg_{\vecc\cdot\vecu}(\psi^{\ast}E)$.
\begin{lem} \label{lem;10.26.12}
The following naturally defined morphism is isomorphic:
\[
 \bigoplus_{\vecc\cdot\vecu=\vecu_1}
 \psi^{\ast}\lefttop{\lbar}\nbigg_{\vecu}(E)
\lrarr
 \lefttop{\lbar}\nbigg_{\vecu_1}(\psi^{\ast}(E)).
\]
\end{lem}
\pf
It follows from Lemma \ref{lem;10.24.1}
and Lemma \ref{lem;10.26.5}.
\hfill\qed

\begin{cor} \label{cor;10.26.35}
The correspondence $\vecu\longmapsto\vecc\cdot\vecu$
induces the surjective morphism
$\KMS(\nbige^0,\lbar)\lrarr \KMS(\psi^{\ast}\nbige^0,\nbar)$.
We have 
$ \sum_{\vecc\cdot\vecu=\vecu_1}
 \multiplicity(0,\vecu)
=\multiplicity(0,\vecu_1)$.
\hfill\qed
\end{cor}

%% file: 16.tex

\subsubsection{Preliminary}

\label{subsubsection;10.18.16}

We denote the space of the multi-valued flat sections
by $H(\nbigelambda)$.
Let $\lefttop{i}M^{\lambda}$ denote the monodromy
of $\nbigelambda$ with respect to the loop around
the divisor $\nbigd_i^{\lambda}$.
We obtain the tuple of endomorphisms
$\vecM=(\lefttop{1}M,\ldots,\lefttop{l}M)$.
Then we obtain the generalized eigen decomposition:
\[
 H(\nbigelambda)=
\bigoplus_{\vecomega\in\Sp(\vecM)}
 \EE\bigl(H(\nbigelambda),\vecomega\bigr).
\]
We denote the restriction of $\vecM$ to
$\EE (H(\nbigelambda),\vecomega)$
by $\vecM_{\vecomega}$.
We obtain the tuple of endomorphisms
$\vecN_{\vecomega}
=(\lefttop{1}N_{\vecomega},\ldots,\lefttop{l}N_{\vecomega})$
given as follows:
\[
 \lefttop{i}N_{\vecomega}
:=\frac{-1}{2\pi\sqrt{-1}}
 \log \lefttop{i}M^u_{\vecomega}.
\]
Here $\lefttop{i}M^u_{\vecomega}$ denotes the unipotent part
of $\lefttop{i}M_{\vecomega}$.

Let $\vecb=(b_1,\ldots,b_l)$ be an element of $\real^l$
and $\vecomega=(\omega_1,\ldots,\omega_l)$ be an element of
$\cnum^l$.
We put as follows:
\[
 \vecalpha(\vecb,\vecomega):=
 \bigl(
 \alpha(b_1,\omega_1),\ldots,\alpha(b_l,\omega_l)
 \bigr).
\]
Here $\alpha(b,\omega)$ for $(b,\omega)\in\real\times\cnum$
is given in the page \pageref{page;9.11.15}.

\subsubsection{The increasing order and the filtration $\lefttop{i}\nbigf$}

\label{subsubsection;10.18.17}

Let $s$ be an element of $\EE(H(\nbigelambda),\vecomega)$.
Let $P$ be an element of $D_i^{\circ}$.
Then we put $\lefttop{i}\ord(s):=\ord(s_{|\pi_i^{-1}(P)})$.

\begin{lem} \label{lem;9.11.28}
The number $\lefttop{i}\ord(s)$ is independent of a choice of $P$.
\end{lem}
\pf
Let $P$ and $P'$ be two points of $D_i^{\circ}$.
Let $\gamma$ be a path in $D_i^{\circ}$ connecting
$P$ and $P'$.

Since $s$ is holomorphic with respect to $\delbar_{\nbigelambda}$,
we have the equality
$ d(s,s)_h=
(\del_{\nbige^{\lambda}}s,s)_h+(s,\del_{\nbige^{\lambda}}s)
=2\Re\bigl(
 (\bar{\lambda}+\lambda^{-1})\cdot(\theta s,s)
 \bigr)$.
Hence we obtain the following equality:
\[
 d \log |s|_h^2
=-2\Re\Bigl(
 (\bar{\lambda}+\lambda^{-1})
\cdot
 \frac{(\theta s,s)_h}{|s|^2_h}.
 \Bigr)
\]
Let $q_i$ be the projection of $X=\Delta^n$ onto the $i$-th component.
Let $Q$ be any point of $\Delta-\{O\}$,
and then $q_i^{-1}(Q)$ is a hyperplane of $X$.
Due to the estimate of $\theta$ (Lemma \ref{lem;9.7.105}, for example),
there exists a positive constant $C$ which is independent of $Q$,
satisfying the following inequality on $q_i^{-1}(Q)$:
\[
 \bigl|d \log |s|^2_{h\,|\,q_i^{-1}(Q)} \bigr|\leq
 C\cdot \Bigl(\sum_{j\neq i} \frac{1}{|z_j|}\Bigr).
\]
Thus we obtain the following inequality
on the path $\gamma_Q:=q_i^{-1}(Q)\cap \pi_i^{-1}(\gamma)$
for some constant $C$ which is independent of $Q$:
\[
 \big| d\log |s_{|\gamma_Q}|^2_h\big|\leq C.
\]
Hence we obtain the following estimate,
which is independent of $Q$:
\[
 \Bigl|
 \log\bigl|s_{|\pi_i^{-1}(P)}(Q)\bigr|_h-
 \log\bigl|s_{|\pi_i^{-1}(P')}(Q)\bigr|_h
\Bigr|
 <C'
\]
It implies our claim.
\hfill\qed

\vspace{.1in}

Let $s$ be an element of $\EE(H(\nbigelambda),\vecomega)$
and $\vecb$ be an element of $\real^l$:
We put as follows:
\[
 F(s,\vecb):=
 \exp\Bigl(
  \sum_{i=1}^l \log z_i\cdot
  \bigl(
   \alpha(b_i,\omega_i)+\lefttop{i}N_{\vecomega}
  \bigr)
 \Bigr)\cdot s.
\]
Then $F(s,\vecb)$ is a holomorphic section of $\nbigelambda$
over $X-D$.
We put as follows:
\[
 c_i:=-\lefttop{i}\ord(s)-\Re(\alpha(b_i,\omega_i)).
\]
We put $\vecc:=(c_i)$.

\begin{lem}
$F(s,\vecb)$ is a holomorphic section of $\prolongg{\vecc}{\nbigelambda}$.
\end{lem}
\pf
Due to the result in the case of curves (Lemma \ref{lem;9.6.3}),
we obtain the following for any point $P\in D_i^{\circ}$:
\[
 -\ord\bigl(F(s,\vecb)_{|\pi_i^{-1}(P)}\bigr)\leq c_i
\]
Then we obtain the result due to Corollary \ref{cor;11.28.15}.
\hfill\qed

\begin{cor}
 $F\bigl(s,-\vecord(s)\bigr)$
 is a holomorphic section of $\prolong{\nbigelambda}$.
\hfill\qed
\end{cor}

We put as follows:
\[
 \lefttop{i}\nbigf_{a}
 \EE(H(\nbigelambda),\vecomega)
=\bigl\{s\in \EE(H(\nbigelambda),\vecomega)\,
 \big|\,-\lefttop{i}\ord(s)\leq a\bigr\}.
\]

\begin{lem}
The monodromies preserves the filtration $\lefttop{i}\nbigf_a$.
\end{lem}
\pf
It immediately follows from the definition of the filtration
$\lefttop{i}\nbigf$ and Lemma \ref{lem;9.11.28}.
\hfill\qed

\subsubsection{Functoriality of $\lefttop{i}\nbigf$}

For a positive integer $c$,
we have the morphism
$\psi_{c\cdot\vecdelta_i}:X\lrarr X$
given by $(z_1,\ldots,z_n)\longmapsto
 (z_1,\ldots,z_{i-1},z_i^c,z_{i+1},\ldots,z_n)$.
We have the natural isomorphism
$\psi_{c\cdot\vecdelta}^{\ast}:
 H(\nbigelambda)\simeq H(\psi_{c\cdot\vecdelta_i}^{-1}\nbigelambda)$
via the pull back.

\begin{lem}\label{lem;9.11.35}
Under the isomorphism, we have the following:
\[
 \psi_{c\cdot\vecdelta}^{\ast}
\bigl(
 \lefttop{i}\nbigf_a
 H(\nbigelambda)
=\lefttop{i}\nbigf_{c\cdot a}
 H\bigl(
 \psi_{c\cdot\vecdelta}^{-1}\nbigelambda
 \bigr)
\]
\end{lem}
\pf
It can be reduced to the case of curves
(Lemma \ref{lem;a11.17.1}).
\hfill\qed

\vspace{.1in}

We have the natural isomorphism
$H(\nbige^{(1)\,\lambda}\otimes\nbige^{(2)\,\lambda})
\simeq
 H(\nbige^{(1)\,\lambda})\otimes 
 H(\nbige^{(2)\,\lambda})$.
We have the following isomorphism:
\[
 \EE\bigl(H(\nbige^{(1)\,\lambda} \otimes\nbige^{(2)\,\lambda}\bigr),
 \vecomega)
\simeq
 \bigoplus_{\vecomega_1\cdot\vecomega_2=\vecomega}
 \EE\bigl(H(\nbige^{(1)\lambda}),\vecomega_1\bigr)
 \otimes
 \EE\bigl(H(\nbige^{(2)\,,\lambda}),\vecomega_2\bigr).
\]
Then we have the two filtrations on
$\EE\bigl(
 H(\nbige^{(1)\,\lambda}\otimes\nbige^{(2)\,\lambda}),
 \vecomega\bigr)$.
One is $\lefttop{i}\nbigf$ for
$\EE\bigl(
 H(\nbige^{(1)\,\lambda}\otimes\nbige^{(2)\,\lambda}),\vecomega\bigr)$.
The other is induced by $\lefttop{i}\nbigf$
for $\EE\bigl(H(\nbige^{(b)\lambda}),\vecomega\bigr)$ for $b=1,2$.

\begin{lem} \label{lem;10.12.6}
They are same. Namely the following holds:
\[
 \lefttop{i}
 \nbigf_a\Bigl(
 \EE\bigl(
 H(\nbige^{(1)\,\lambda} \otimes\nbige^{(2)\,\lambda},\vecomega\bigr)
 \Bigr)
\simeq
 \bigoplus_{\vecomega_1\cdot\vecomega_2=\vecomega}
 \sum_{a_1+a_2\leq a}
 \lefttop{i}\nbigf_{a_1}\Bigl(
 \EE\bigl(H(\nbige^{(b)\lambda}),\vecomega_1\bigr)
 \Bigr)
 \otimes
 \lefttop{i}\nbigf_{a_2}\Bigl(
 \EE\bigl(H(\nbige^{(b)\lambda}),\vecomega_2\bigr)
 \Bigr).
\]
\end{lem}
\pf
Due to the result in the case of curves,
we obtain the following:
\begin{multline}
  \lefttop{i}\nbigf_a
 \Bigl(
 \bigoplus_{q_i(\vecomega)=\omega}
 \EE\bigl(
 H(\nbige^{(1)\,\lambda} \otimes\nbige^{(2)\,\lambda}),\vecomega\bigr)
 \Bigr) \\
\simeq
 \bigoplus_{\omega_1\cdot\omega_2=\vecomega}
 \sum_{a_1+a_2\leq a}
 \lefttop{i}\nbigf
 \Bigl(
 \bigoplus_{q_i(\vecomega_1)=\omega_1}
 \EE\bigl(H(\nbige^{(1)\,\lambda}),\vecomega_1\bigr)
 \Bigr)
\otimes
 \lefttop{i}\nbigf
 \Bigl(
 \bigoplus_{q_i(\vecomega_2)=\omega_2}
 \EE\bigl(H(\nbige^{(2)\,\lambda}),\vecomega_2\bigr)
 \Bigr).
\end{multline}
Since the monodromy endomorphisms $M^{\lambda}_j$ $(j\neq i)$
preserve the filtration $\lefttop{i}\nbigf$,
we obtain the results.
\hfill\qed

\vspace{.1in}

We have the following isomorphism:
\begin{equation} \label{eq;9.11.20}
 H\Bigl(\bigotimes^R\nbigelambda\Bigr)
\simeq
 \bigotimes^RH(\nbigelambda).
\end{equation}
We have the following isomorphism:
\begin{equation} \label{eq;9.11.21}
 \EE\Bigl(H\Bigl(\bigotimes^R\nbigelambda\Bigr),\vecomega\Bigr)
\simeq
 \bigoplus_{f\in\nbigs(\vecomega,R)}
 \bigotimes_{\vecomega'\in\Sp(\vecM^{\lambda})}
 \bigotimes^{f(\vecomega')}
 \EE\bigl(H(\nbigelambda),\vecomega'\bigr).
\end{equation}
Here we put as follows:
\[
 \nbigs(\vecomega,R):=
 \Bigl\{
 f:\Sp(\vecM^{\lambda})\lrarr\seisuu_{\geq 0}\,
 \Big|\,
 \sum f(\vecomega')=R,\quad
 \prod \vecomega^{\prime\,f(\vecomega')}=\vecomega
 \Bigr\}.
\]
We naturally have the filtrations on the both sides
of (\ref{eq;9.11.20}) and (\ref{eq;9.11.21}).

\begin{cor}
The filtrations of the both sides
{\rm(\ref{eq;9.11.20})} and {\rm(\ref{eq;9.11.21})}
are preserved by the isomorphisms.
\hfill\qed
\end{cor}

We have the isomorphism
$H(\nbige^{\lor\,\lambda})\simeq H(\nbigelambda)^{\lor}$.
preserving the $\EE$-decomposition.
\begin{equation} \label{eq;10.12.5}
 \EE(H(\nbige^{\lor\,\lambda}),\vecomega)
\simeq
 \EE(H(\nbigelambda)^{\lor},\vecomega).
\end{equation}
On the left hand side of (\ref{eq;10.12.5}),
we have the filtration $\lefttop{i}\nbigf$.
On the right hand side of (\ref{eq;10.12.5}),
we have the induced filtration $\lefttop{i}{\nbigf}^{\lor}$
by $\lefttop{i}\nbigf$,
given as follows:
\[
 \lefttop{i}{\nbigf}^{\lor}_aH(\nbigelambda)^{\lor}:=
 \bigl\{
 f\in H(\nbigelambda)^{\lor}\,\big|\,
 f \bigl(
 \lefttop{i}\nbigf_bH(\nbigelambda)
 \bigr)
 \subset
 \lefttop{i}\nbigf_{a+b}H(\nbigelambda),
\,\,(\forall b\in\real)
 \bigr\}.
\]
 
\begin{lem}
The isomorphism preserves the filtrations.
\end{lem}
\pf
It can be reduced to the case of curves
(Lemma \ref{lem;9.11.23})
by an argument similar to the proof of Lemma
\ref{lem;10.12.6}.
\hfill\qed

%% file: 16.1.tex

\subsubsection{The case $\lambda$ is generic}

Assume that $\lambda$ is generic.
Recall that
$\eigenmap^f(\lambda):\KMSoverline(\nbige^0,i)
\lrarr
 \Sp^f(\nbigelambda,i)$ is isomorphic for any $i$.

\begin{lem} \label{lem;9.11.26}
 The filtration $\lefttop{i}\nbigf$ on
 $\EE(H(\nbigelambda),\vecomega)$ is trivial
 in the following sense:
\begin{quote}
For each $i$, we have the unique element
$u_i\in\KMS(\prolong{\nbige^0},i)$
such that $\omega_i=\eigenmap^f(\lambda,u_i)$.
Then $\lefttop{i}\Gr^{\nbigf}_b\EE(H(\nbigelambda),\vecomega)\neq 0$
if and only if
$b=\paramap^f(\lambda,u_i)$.
\end{quote}
\end{lem}
\pf
It follows from the result in the case of curves
(the subsubsection \ref{subsubsection;9.11.25}).
\hfill\qed

\begin{cor}
If $\lambda$ is generic,
then the filtrations
$(\lefttop{i}\nbigf\,|\,i=1,\ldots,l)$ are compatible.
\end{cor}
\pf
It immediately follows Lemma \ref{lem;9.11.26} .
\hfill\qed

\vspace{.1in}

Let $\lambda\in\cnum_{\lambda}^{\ast}$ be generic,
and let $\vecomega=(\omega_1,\ldots,\omega_l)$
be an element of $\Sp(\vecM^{\lambda})$.
We have the unique element $u_i\in\nbigk(\nbige,\lambda,0,i)$
such that $\eigenmap^f(\lambda,u_i)=\omega_i$.
\begin{lem}
Let $\vecbeta$ be an element of $\cnum^l$
whose $i$-th component is
$\beta_i=\eigenmap(\lambda,u_i)$.
Recall that we have the space
$\EE(\prolong{\nbigelambda}_{|O},\vecbeta)$,
which is a generalized eigenspace
of the tuple $\Res_{\,\lbar}(\DD^{\lambda})$.
Then we have the isomorphism:
\[
 \EE(\prolong{\nbigelambda}_{|O},\vecbeta)
\simeq
 \EE(H(\nbigelambda),\vecomega).
\]

\end{lem}
\pf
Let $\vecs$ be base of $H(\nbigelambda)$,
which is compatible with $\EE$.
We put $\vecomega(s_j):=\vecdeg^{\EE}(s_j)$
and $\vecb(s_j):=\vecdeg^{\nbigf}(s_j)$.

We put $v_j:=F(s_j,\vecb(s_j))$.
Then we obtain the tuple $\vecv=(v_i)$ of sections
of $\prolong{\nbigelambda}$.
Due to the result in the case of curves,
$\vecv_{|\pi_i^{-1}(P)}$ is a frame of
$\prolong{\nbigelambda_{|\pi_i^{-1}(P)} }$.
Thus $\vecv$ is a frame of $\prolong{\nbigelambda}$.
The frames $\vecv$ and $\vecs$ induce
the isomorphism desired.
\hfill\qed

\begin{cor}
Assume that $\lambda$ is generic.
For any $\vecomega\in\Sp(\vecM^{\lambda})$,
we have the unique element $\vecu\in\KMSoverline(\nbige^0,\lbar)$
such that $\eigenmap^f(\lambda,\vecu)=\vecomega$.
\hfill\qed
\end{cor}

\begin{cor}\label{cor;a11.17.2}
Assume $\vecomega=\eigenmap^f(\lambda,\vecu)$
and $\vecb=\paramap^f(\lambda,\vecu)$.
We have the following equalities:
\[
 \dim \lefttop{\lbar}\nbigf_{\vecb}\Bigl(
 \EE(H(\nbigelambda),\vecomega)\Bigr)
=\left\{
 \begin{array}{ll}
 \multiplicity(0,\vecu) &(\vecb\geq \paramap^f(\lambda,\vecu))\\
 \mbox{{}}\\
 0  & (\mbox{\rm otherwise}).
 \end{array}
 \right.
\]
\[
 \dim \lefttop{\lbar}\Gr^{\nbigf}_{\vecb}\Bigl(
 \EE(H(\nbigelambda),\vecomega)\Bigr)
=\left\{
 \begin{array}{ll}
 \multiplicity(0,\vecu) &(\vecb= \paramap^f(\lambda,\vecu))\\
 \mbox{{}}\\
 0  & (\mbox{\rm otherwise}).
 \end{array}
 \right.
\]
\hfill\qed
\end{cor}

%% file: 16.2.tex

\subsubsection{The family of the space of the multi-valued flat sections}
   
Let $\nbigh=\nbigh(E)$ be the holomorphic vector bundle
over $\cnum^{\ast}_{\lambda}$,
obtained by $\nbigh_{|\lambda}=H(\nbigelambda)$.
We have the monodromy endomorphisms
$\vecM=(M_i\,|\,i=1,\ldots,l)$,
where $M_i$ denotes the monodromy with respect to $D_i$.

Let $\lambda_0$ be an element of $\cnum_{\lambda}^{\ast}$.
We put $S_0:=\Sp(\vecM^{\lambda_0})$.
Let $\epsilon_0$ and $\epsilon_1$ be sufficiently small numbers.
Then we obtain the following decompositions, (see (\ref{eq;a12.9.1})
 for the notation $\EE_{\epsilon_1}$):
\[
 \nbigh_{|\Delta(\lambda_0,\epsilon_0)}
=
 \bigoplus_{\vecomega\in S_0}
 \nbigh^{(\lambda_0)}_{\vecomega},
\quad\quad
 \nbigh^{(\lambda_0)}_{\vecomega}
 :=\lefttop{\lbar}\EE_{\epsilon_1}\bigl(H(\nbigelambda),\vecomega\bigr).
\]
We put 
$ \nbigs(\vecomega):=
\bigl\{\vecu\in\KMSoverline({\nbige^0},\lbar)\,\big|\,
 \eigenmap^f(\lambda_0,\vecu)=\vecomega
 \bigr\}$.
We may assume that
any $\lambda\in\Delta^{\ast}(\lambda_0,\epsilon_0)$ is generic.
Then we have the following decomposition
on $\Delta^{\ast}(\lambda_0,\epsilon_0)$:
\[
 \nbigh^{(\lambda_0)}_{\vecomega\,|\,\Delta^{\ast}(\lambda_0,\epsilon_0)}
=\bigoplus_{\vecu\in \nbigs(\vecomega)}
 \nbigh_{\eigenmap^f(\lambda,\vecu)},
\quad\quad
 \nbigh_{\eigenmap^f(\lambda,\vecu)\,|\,\lambda}
=\EE(H(\nbigelambda),\eigenmap^f(\lambda,\vecu)).
\]

As in the case of the curves,
we consider the filtration $\lefttop{i}\nbigf^{(\lambda_0)}$.
We put as follows:
\begin{equation}\label{eq;a12.9.6}
 \lefttop{i}\nbigf^{(\lambda_0)}_d
 \nbigh^{(\lambda_0)}_{\vecomega\,|\,\Delta^{\ast}(\lambda_0,\epsilon_0)}
=\bigoplus_{
 \substack{
 \vecu\in\nbigs(\vecomega)\\
 \paramap^f(\lambda,q_i(\vecu))\leq d
 }
 }\nbigh_{\eigenmap^f(\lambda,\vecu)}.
\end{equation}
Since it is given as the sum of the generalized eigenspaces,
the filtration $\lefttop{i}\nbigf^{(\lambda_0)}$
can be prolonged to the filtration of
$\nbigh^{(\lambda_0)}_{\omega}$,
which we denote by $\lefttop{i}\nbigf^{(\lambda_0)}$.

Note the following isomorphism for any $P\in D_i^{\circ}$:
\[
 H(\nbigelambda)\simeq
 H(\nbigelambda_{|\pi_i^{-1}(P)}).
\]
We have the filtration $\nbigf^{(\lambda_0)}$
on the right hand side.

\begin{lem}\mbox{{}}
\begin{itemize}
\item
 Under the isomorphism above, we have
 $\nbigh^{(\lambda_0)}_{\omega}=
 \bigoplus_{q_i(\vecomega)=\omega}\nbigh_{\vecomega}^{(\lambda_0)}$.
\item
 $\nbigf^{(\lambda_0)}_d\nbigh_{\omega}=
 \bigoplus_{q_i(\vecomega)=\omega}
 \lefttop{i}\nbigf^{(\lambda_0)}\nbigh^{(\lambda_0)}_{\vecomega}$.
\item
 In particular,
 we have 
 $\bigl(
 \lefttop{i}\nbigf_d^{(\lambda_0)}\nbigh_{\vecomega}^{(\lambda_0)}
 \bigr)_{|\lambda_0}
 =\lefttop{i}\nbigf\bigl(\EE\bigl(H(\nbigelambdazero),\vecomega\bigr)\bigr)$.
\hfill\qed
\end{itemize}
\end{lem}

%% file: a55.1.tex

\subsubsection{The dimension and the virtual dimension}

For any $\lambda$,
pick $\vecomega\in\Sp(\vecM^{\lambda})$
and $\veca\in\real^l$.
We put as follows:
\[
 d(\lambda,\vecomega,\veca):=
 \dim\lefttop{\lbar}\nbigf_{\veca}
 \EE(H(\nbigelambda),\vecomega),
\quad\quad
 v.d(\lambda,\vecomega,\veca):=
 \sum_{\substack{
 \vecu\in \nbigs(\vecomega)\\
 \paramap^f(\lambda,\vecu)\leq\veca
 } }
 \multiplicity(0,\vecu).
\]

\begin{lem}\label{lem;10.12.10}
Let $\lambda_0$ be any element of $\cnum_{\lambda}$.
We have the inequality
$d(\lambda_0,\vecomega,\veca)\geq v.d(\lambda_0,\vecomega,\veca)$.
\end{lem}
\pf
Let us pick a sufficiently small positive number $\epsilon$.
Then we have the following inequality,
for any $\lambda\in\Delta(\lambda_0,\epsilon_0)$.
\begin{equation}\label{eq;9.11.30}
 d(\lambda_0,\vecomega,\veca)=
\dim\Bigl(
 \bigcap_{i=1}^l
 \lefttop{i}\nbigf^{(\lambda_0)}_{a_i}
 \nbigh^{(\lambda_0)}_{\vecomega\,|\,\lambda_0}
 \Bigr)
\geq
 \dim\Bigl(
\bigcap_{i=1}^l
 \lefttop{i}\nbigf^{(\lambda_0)}_{a_i}
 \nbigh^{(\lambda_0)}_{\vecomega\,|\,\lambda}
 \Bigr).
\end{equation}
For any generic $\lambda\in\Delta^{\ast}(\lambda_0,\epsilon_0)$,
we have the following equality:
\[
 \lefttop{i}\nbigf^{(\lambda_0)}_{a_i}
 \nbigh^{(\lambda_0)}_{\vecomega\,|\,\lambda}
=\bigoplus_{
 \substack{\vecu\in\nbigs(\vecomega)\\
 \paramap^f(\lambda_0,u_i)\leq a_i
 }
 }\nbigh^{(\lambda_0)}_{\eigenmap^f(\lambda,\vecu)\,|\,\lambda}.
\]
Thus we obtain the following:
\[
 \bigcap_{i=1}^l
 \lefttop{i}\nbigf^{(\lambda_0)}_{a_i}
 \nbigh^{(\lambda_0)}_{\vecomega\,|\,\lambda}
=\bigoplus_{
 \substack{\vecu\in\nbigs(\vecomega)\\
 \paramap^f(\lambda_0,\vecu)\leq \veca
 }
 }\nbigh^{(\lambda_0)}_{\eigenmap^f(\lambda,\vecu)\,|\,\lambda}.
\]
Therefore we obtain the following equality,
due to Corollary \ref{cor;a11.17.2}.
\begin{equation} \label{eq;9.11.31}
 \dim\bigcap_{i=1}^l
 \lefttop{i}\nbigf^{(\lambda_0)}_{a_i}
 \nbigh^{(\lambda_0)}_{\vecomega\,|\,\lambda}
 =\sum_{
 \substack{
 \vecu\in\nbigs(\vecomega)\\
 \paramap^f(\lambda_0,\vecu)\leq \veca
 }}
 \dim\nbigh_{\eigenmap^f(\lambda,\vecu)\,|\,\lambda}
=\sum_{
 \substack{
 \vecu\in\nbigs(\vecomega)\\
 \paramap^f(\lambda_0,\vecu)\leq \veca
 }}
 \multiplicity(0,\vecu)
=v.d(\lambda_0,\vecomega,\veca).
\end{equation}
We obtain the result from (\ref{eq;9.11.30}) and (\ref{eq;9.11.31}).
\hfill\qed

%% file: 16.3.tex

\subsubsection{Preliminary proposition}

Let us consider the compatibility of the filtrations
$\lefttop{i}\nbigf$ $(i=1,\ldots,l)$.
We may assume that $l=n$.
We put $X^{(1)}:=\{(z_1,\ldots,z_n)\,|\,z_{n-1}=z_n\}
\subset X$.
We put $D_i^{(1)}:=D_i\cap X^{(1)}$ for $i\leq n-1$.
We have the natural isomorphism
$X^{(1)}\simeq \Delta^{n-1}=\{(z_1,\ldots,z_{n-1})\in\Delta^{n-1}\}$.
We have the natural isomorphism
$H(\nbigelambda)\lrarr H(\nbigelambda_{|X^{(1)}})$.

We have the tuple of monodromies
$\vecM:=(M_1,\ldots,M_n)$
and $\vecM^{(1)}:=(M_1^{(1)},\ldots,M_{n-1}^{(1)})$.
Here we have the following:
\[
 M_i^{(1)}=
 \left\{
 \begin{array}{ll}
 M_i & (i\leq n-2)\\
 M_{n-1}\circ M_n &(i=n-1).
 \end{array}
 \right.
\]
For any $\vecomega\in\Sp(\vecM)$,
$\phi(\vecomega)\in \Sp(\vecM^{(1)})$ is given as follows:
\[
 q_i(\phi(\vecomega))=
 \left\{
 \begin{array}{ll}
 \omega_i &(i\leq n-2)\\
 \omega_{n-1}\cdot\omega_n &(i=n-1).
 \end{array}
 \right.
\]

The map $\phi:\real^n\lrarr\real^{n-1}$
is defined as follows:
\[
 q_i(\phi(\veca))=
\left\{
 \begin{array}{ll}
 a_i &(i\leq n-2)\\
 a_{n-1}+a_n & (i=n-1).
 \end{array}
\right.
\]

We have the filtrations $\lefttop{i}\nbigf$
 $(i=1,\ldots,n)$
on $\EE(H(\nbigelambda),\vecomega)$.
For any $\veca\in\real^{n}$,
we have the subspace
$\lefttop{\nbar}\nbigf_{\veca}
:=\bigcap_{i=1}^n
 \lefttop{i}\nbigf_{a_i}$
of $H(\nbigelambda)$.
We put as follows:
\[
 \lefttop{\nbar}\nbigf'_{\veca}
=\sum_{\substack{\vecb\lneq \veca}}
 \lefttop{\nbar}\nbigf_{\vecb}.
\]
Here $\vecb\lneq\veca$ means $\vecb\leq \veca$ and $\vecb\neq\veca$.

We also have the filtrations 
on $\EE\bigl(H(\nbigelambda_{|X^{(1)}}),\vecomega^{(1)}\bigr)$,
which we denote by $\lefttop{i}\nbigf^{(1)}$
$(i=1,\ldots,n-1)$.
Similarly, we have
$\lefttop{\nminusitibar}\nbigf^{(1)}_{\veca^{(1)}}$
and 
$\lefttop{\nminusitibar}\nbigf^{\prime\,(1)}_{\veca^{(1)}}$
for any $\veca^{(1)}\in\real^{n-1}$.

\begin{lem}
We have the following implication:
\[
 \lefttop{\nbar}\nbigf_{\veca}\cap
 \EE(H(\nbigelambda),\vecomega)
\subset
 \lefttop{\nminusitibar}\nbigf^{(1)}_{\phi(\veca)}
 \cap
 \EE\bigl(H(\nbigelambda_{|X^{(1)}}),\phi(\vecomega)\bigr).
\]
\end{lem}
\pf
Let $s$ be an element of the left hand side.
Let $P$ be any point of $D_i^{\circ}$.
Then we have the following for any $i=1,\ldots,n$:
\[
  -\ord(F(s,\vecb)_{|\pi_i^{-1}(P)})
\leq
 a_i-\Re\bigl(\alpha(b_i,\omega_i)\bigr).
\]

Let $C$ be the subset of $X^{(1)}\simeq \Delta^{n-1}$
such that
$C=\pi_{n-1}^{-1}(P)$
for a point
$P\in D^{(1)}_{n-1}
 -\bigcup_{i<n-1}D^{(1)}_{n-1}\cap D^{(1)}_i$.
The inclusion $C\subset X$ is obtained by
the diagonal embedding
$\{P\}\times \Delta_{n-1}\lrarr \{P\}\times\Delta_{n-1}\times\Delta_n$.
We have the following,
due to Corollary \ref{cor;11.28.15}:
\[
 -\ord(F(s,\vecb)_{|C})
\leq
 a_{n-1}+a_n-\Re(\alpha(b_{n}+b_{n-1},\omega_n\cdot\omega_{n-1})).
\]
It implies
$-\ord(s_{|C})\leq a_{n-1}+a_n$
due to Lemma \ref{lem;9.6.3}.
Hence we have the inequality:
\[
 \lefttop{n-1}\deg^{\nbigf^{(1)}}(s)
\leq a_{n-1}+a_n.
\]
Thus we are done.
\hfill\qed

\begin{prop} \label{prop;9.11.40}
The following holds.
\begin{description}
\item[(A)]
 $d(\lambda,\vecomega,\veca)=v.d(\lambda,\vecomega,\veca)$.
\item[(B)]
 The following morphism is isomorphic:
\[
 \sum_{\substack{
 \phi(\veca)\leq \veca^{(1)}\\
 \phi(\vecomega)=\vecomega^{(1)}
 } }
 \Bigl(
 \lefttop{\nbar}\nbigf_{\veca}
 \cap
 \EE\bigl(H(\nbigelambda),\vecomega\bigr)
\Bigr)
\lrarr
 \lefttop{\nminusitibar}
 \nbigf^{(1)}_{\veca^{(1)}}
 \cap
 \EE\bigl(H(\nbigelambda_{|X^{(1)}}),\vecomega^{(1)}\bigr).
\]
\item[(C)]
The following morphism is injective.
\[
 \frac{\lefttop{\nbar}\nbigf_{\veca}
   \cap\EE\bigl(H(\nbigelambda),\vecomega\bigr)}
 {\lefttop{\nbar}\nbigf'_{\veca}
  \cap \EE\bigl(H(\nbigelambda),\vecomega\bigr)}
\lrarr
 \frac{
 \lefttop{\nminusitibar}\nbigf^{(1)}_{\phi(\veca)}
 \cap
 \EE\bigl(H(\nbigelambda_{|X^{(1)}}),\phi(\vecomega)\bigr)}
 {\lefttop{\nminusitibar}\nbigf^{\prime\,(1)}_{\phi(\veca)}
 \cap
 \EE\bigl(H(\nbigelambda_{|X^{(1)}}),\phi(\vecomega)\bigr)}.
\]
\end{description}
\end{prop}
The proposition will be proved in the subsubsections
\ref{subsubsection;10.12.50}--\ref{subsubsection;10.12.51}.

\vspace{.1in}

Before entering the proof,
we simplify the problem.
Considering the morphism
$\psi_c^{-1}(z_1,\ldots,z_{n})=(z_1,\ldots,z_{n-1},z_n^c)$,
we may assume the following for some $\epsilon>0$.
(Note Lemma \ref{lem;9.11.35}):\\
\begin{description}
\item[{Assumption D}]
 For any $a,b\in\Par^f(\nbigelambda,n)$ such that $a\neq b$,
and for any $\lambda'\in\Delta(\lambda,\epsilon)$,
\[
 |a-b|>\sum_{c\in\Par^f(\nbigelambda,n-1)}|c|.
\]
\end{description}

\begin{lem}\mbox{{}}
Under the assumption D, the following holds:
\begin{enumerate}
\item
 The morphism
$\Par^f(\nbigelambda,n)\times\Par^f(\nbigelambda,n-1)
\lrarr \real$ given by $(a,b)\longmapsto a+b$
is injective.
In particular,
it induces the total order $\leq_1$ on
the set
$\Par^f(\nbigelambda,n)\times\Par^f(\nbigelambda,n-1)$.
\item
We have the natural orders on $\Par^f(\nbigelambda,i)$
$(i=n,n-1)$,
and we obtain the lexicographic order $\leq_2$
on $\Par^f(\nbigelambda,n)\times\Par^f(\nbigelambda,n-1)$.
We have the equality $\leq_1=\leq_2$.
\end{enumerate}
\end{lem}
\pf
It immediately follows from our assumption above.
\hfill\qed

\subsubsection{A proof of the claims $(A)$ and $(B)$ of
Proposition \ref{prop;9.11.40}}
\label{subsubsection;10.12.50}

Let $\veca^{(1)}$ be an element of $\real^{n-1}$,
whose $i$-th components are $a_{i}^{(1)}$.
Then the element
$(a_{n}^{\circ},a_{n-1}^{\circ})\in 
 \Par^f(\nbigelambda,n)\times\Par^f(\nbigelambda,n-1)$
is determined as follows:
\[
 (a_{n}^{\circ},a_{n-1}^{\circ}):=
 \max\bigl\{(b_n,b_{n-1})\in
 \Par^f(\nbigelambda,n)\times\Par^f(\nbigelambda,n-1)
 \,\big|\,b_n+b_{n-1}\leq a_{n-1}^{(1)}\bigr\}.
\]

\begin{lem}
We have the following equality:
\begin{equation}\label{eq;a11.17.5}
 \sum_{\phi(\veca)\leq\veca^{(1)}}
 \lefttop{\nbar}\nbigf_{\veca}
=
\Bigl(
\lefttop{\nminusnibar}\nbigf_{\veca'}
\cap\lefttop{n-1}\nbigf_{a^{\circ}_{n-1}}
\cap\lefttop{n}\nbigf_{a^{\circ}_n}
\Bigr)
+
\Bigl(
 \lefttop{\nminusnibar}\nbigf_{\veca'}
\cap
 \lefttop{n}\nbigf_{<a^{\circ}_n}
\Bigr).
\end{equation}
Here we put $\veca':=(a_1^{(1)},\ldots,a^{(1)}_{n-2})\in\real^{n-2}$.
\end{lem}
\pf
It is clear that
$\lefttop{\nminusnibar}\nbigf_{\veca'}
\cap\lefttop{n-1}\nbigf_{a^{\circ}_{n-1}}
\cap\lefttop{n}\nbigf_{a^{\circ}_n}$
is contained in the left hand side of (\ref{eq;a11.17.5}).
Under the assumption $(D)$,
we have
$\lefttop{n}\nbigf_{<a_n^{\circ}}=
 \lefttop{n-1}\nbigf_{b}\cap\lefttop{n}\nbigf_{a_n^{\circ}-\eta}$
for some real numbers $b$ and $\eta$
such that $b+a_n^{\circ}-\eta<a_{n-1}^{(1)}$.
Thus $\lefttop{\nminusnibar}\nbigf_{\veca'} \cap
 \lefttop{n}\nbigf_{<a^{\circ}_n}$
is contained in the left hand side of (\ref{eq;a11.17.5}).
Thus we obtain the implication $\supset$.

Next, we would like to show the implication $\subset$.
We have only to show 
that $\lefttop{\nbar}\nbigf_{\veca}\EE\bigl(H(\nbigelambda),\vecomega\bigr)$
is contained in the right hand side
when we have $\phi(\veca)\leq\veca^{(1)}$.

Under the assumption $D$,
the condition $\phi(\veca)\leq \veca^{(1)}$
implies the following:
\[
 a_i\leq a^{(1)}_i,\,\,(i\leq n-2),
\quad\mbox{ and }\quad
 \left\{
 \begin{array}{l}
 a_{n}< a_n^{\circ},\\
 \,\,\,\,\mbox{\rm or},\\
 a_n=a_n^{\circ},\, a_{n-1}\leq a_{n-1}^{\circ}.
 \end{array}
 \right.
\]
Then the implication $\subset$ immediately follows.
\hfill\qed

\vspace{.1in}

We will show the claims $(A)$ and $(B)$
in Proposition \ref{prop;9.11.40} by an induction on $n$.
We assume that the claims $(A)$ and $(B)$ for $n-1$,
and we will show that the claims $(A)$ and $(B)$ hold for $n$.

Let us consider the following claims:
\begin{description}
\item[$(A,a,\leq)$]
 $(A)$ holds for any $\veca$ such that $a_{n}+a_{n-1}\leq a$.
\item[$(A,a,<)$]
 $(A)$ holds for any $\veca$ such that $a_n+a_{n-1}<a$.
\item[$(B,a,\leq)$]
 $(B)$ holds for any $\veca^{(1)}$ such that $a^{(1)}_{n-1}\leq a$.
\item[$(B,a,<)$]
 $(B)$ holds for any $\veca^{(1)}$ such that $a^{(1)}_{n-1}<a$.
\end{description}
If $a$ is sufficiently negative,
then $(A,a,\leq)$ and $(B,a,\leq)$ are true trivially.

\begin{lem} \label{lem;9.11.50}
$(A,a,<)$ implies $(B,a,\leq)$ and $(A,a,\leq)$.
\end{lem}
\pf
We have the following implication:
\begin{equation}\label{eq;a11.17.6}
 \Bigl(
 \lefttop{\nminusnibar}{\nbigf}_{\veca'}
 \cap
 \lefttop{n-1}\nbigf_{a^{\circ}_{n-1}}
 \cap
 \lefttop{n}\nbigf_{a^{\circ}_n}
 \Bigr)
+\Bigl(
 \lefttop{\nminusnibar}\nbigf_{\veca'}
 \cap
 \lefttop{n}\nbigf_{<a^{\circ}_n}
 \Bigr)
\subset
 \lefttop{\nminusitibar}\nbigf^{(1)}_{\veca^{(1)}}.
\end{equation}
We also have the following:
\[
 \Bigl(
 \lefttop{\nminusnibar}\nbigf_{\veca'}
\cap
 \lefttop{n-1}\nbigf_{a^{\circ}_{n-1}}
\cap
 \lefttop{n}\nbigf_{a^{\circ}_n}
 \Bigr)
\cap
 \Bigl(
 \lefttop{\nminusnibar}\nbigf_{\veca'}
 \cap
 \lefttop{n}\nbigf_{<a^{\circ}_n}
 \Bigr)
=
 \lefttop{\nminusnibar}\nbigf_{\veca'}
 \cap\lefttop{n-1}\nbigf_{a^{\circ}_{n-1}}
 \cap\lefttop{n}\nbigf_{<a^{\circ}_n}.
\]
Thus we have the following equality:
\begin{multline}\label{eq;9.11.45}
\dim\Bigl(
 \lefttop{\nminusnibar}\nbigf_{\veca'}
\cap\lefttop{n-1}\nbigf_{a^{\circ}_{n-1}}
\cap\lefttop{n}\nbigf_{a^{\circ}_n}
\cap \EE\bigl(H(\nbigelambda),\vecomega\bigr)
 \Bigr)
\leq 
\dim \bigl(
 \lefttop{\nminusitibar}\nbigf^{(1)}_{\veca^{(1)}}
 \cap \EE\bigl(H(\nbigelambda),\vecomega\bigr)
 \bigr)\\
+\dim\bigl(
 \lefttop{\nminusnibar}\nbigf_{\veca'}
 \cap \lefttop{n-1}\nbigf_{a^{\circ}_{n-1}}
 \cap \lefttop{n}\nbigf_{<a_n^{\circ}}
 \cap \EE\bigl(H(\nbigelambda),\vecomega\bigr)
 \bigr)
-\dim\Bigl(
 \lefttop{\nminusnibar}\nbigf_{\veca''}\cap
 \lefttop{n}\nbigf_{<a_n^{\circ}}
 \cap\EE\bigl(H(\nbigelambda),\vecomega\bigr)
  \Bigr).
\end{multline}

By using the assumption of the induction on $n$,
or by using $(A,a,<)$,
we obtain the following equality:
\begin{equation} \label{eq;9.11.46}
\begin{array}{ll}
{\displaystyle
 \dim\Bigl(
 \lefttop{\nminusnibar}\nbigf_{\veca'}
 \cap
 \lefttop{n}\nbigf_{<a^{\circ}_n}
\cap \EE\bigl(H(\nbigelambda),\vecomega \bigr)
 \Bigr)
=\sum_{\vecu\in S_1(\vecomega)}\multiplicity(0,\vecu)},\\
 \mbox{{}}\\
 S_1(\vecomega)=\bigl\{
\vecu\in\nbigs(\vecomega)\,\big|\,
 \paramap^f(\lambda,u_i)\leq a_i^{(1)}\,(i\leq n-2),\,\,
 \paramap^f(\lambda,u_{n})<a_n^{\circ}
 \bigr\}
\end{array}
\end{equation}
By using $(A,a,<)$, we obtain the following:
\begin{equation}\label{eq;9.11.47}
\begin{array}{l}
{\displaystyle
 \dim\Bigl(
 \lefttop{\nminusnibar}\nbigf_{\veca'}
\cap\lefttop{n-1}\nbigf_{a^{\circ}_{n-1}}
\cap\lefttop{n}\nbigf_{<a^{\circ}_n}
\cap \EE\bigl(H(\nbigelambda),\vecomega\bigr)
 \Bigr)
=\sum_{ \vecu\in S_2(\vecomega)}\multiplicity(0,\vecu),}\\
\mbox{{}}\\
 {\displaystyle
 S_2(\vecomega)=\bigl\{\vecu\in\nbigs(\vecomega)\,\big|\,
 \paramap^f(\lambda,u_i)\leq a_i^{(1)},\,\,
 \paramap^f(\lambda,u_{n-1})\leq a_{n-1}^{\circ},\,\,
 \paramap^f(\lambda,u_{n})<a_n^{\circ}
 \bigr\}.}
\end{array}
\end{equation}
Due to the assumption of the induction on $n$,
we have the following:
\begin{equation} \label{eq;9.11.48}
\begin{array}{l}
{\displaystyle
  \dim\bigl(
 \lefttop{\nminusitibar}\nbigf^{(1)}_{\veca^{(1)}}
 \cap
 \EE\bigl(
 H(\nbigelambda),\vecomega^{(1)}
 \bigr)
 \bigr)
=\sum_{\substack{
 \vecu\in\nbigs(\vecomega^{(1)})\\
 \paramap^f(\lambda,\vecu)\leq \veca^{(1)} }}
 \multiplicity(0,\vecu)
=\sum_{\phi(\vecomega)=\omega^{(1)}}
 \sum_{ \vecu\in S_3(\vecomega)}
 \multiplicity(0,\vecu),}  \\
\mbox{{}}\\
{\displaystyle
 S_3(\vecomega):=\bigl\{\vecu\in\nbigs(\vecomega)\,\big|\,
 \paramap^f(\lambda,u_i)\leq a_i^{(1)}\,\,(i\leq n-2),\,\,
 \paramap^f(\lambda,u_{n-1})+\paramap^f(\lambda,u_n)
 \leq a_{n-1}^{(1)}
 \bigr\}. }
\end{array}
\end{equation}
We put 
$ S_4(\vecomega):=
 S_3(\vecomega)-\bigl(S_1(\vecomega)-S_2(\vecomega)\bigr)$.
It is easy to check the following:
\[
 S_4=\bigl\{
 \vecu\in\nbigs(\vecomega)\,\big|\,
 \paramap^f(\lambda,u_i)\leq a_i\,\,(i\leq n-2),\,\,
 \paramap^f(\lambda,u_i)\leq a^{\circ}_{i}\,\, (i=n-1,n)
 \bigr\}.
\]
Then we obtain the following inequality by a direct calculation
from (\ref{eq;9.11.45}), (\ref{eq;9.11.46}), (\ref{eq;9.11.47})
and (\ref{eq;9.11.48}):
\begin{multline} \label{eq;a11.17.9}
 \sum_{\phi(\vecomega)=\vecomega'}
 \dim
 \Bigl(
 \lefttop{\nminusnibar}\nbigf_{\veca'}
 \cap
 \lefttop{n-1}\nbigf_{a^{\circ}_{n-1}}
 \cap
 \lefttop{n}\nbigf_{a^{\circ}_n}
 \cap \EE\bigl(H(\nbigelambda),\vecomega\bigr)
 \Bigr)
\leq
 \sum_{\phi(\vecomega)=\omega^{(1)}}
 \sum_{\vecu\in S_4(\vecomega) }\multiplicity(0,\vecu) \\
=\sum_{\phi(\vecomega)=\vecomega^{(1)}}
 v.d(\veca,\vecomega).
\end{multline}
Here $\veca$ is determined by
$q_{i}(\veca)=a_{i}^{(1)}$ for $i\leq n-2$
and $q_i(\veca)=a_i^{\circ}$ for $i=n-1,n$.
On the other hand,
we have already known the inequality (Lemma \ref{lem;10.12.10}):
\begin{equation}\label{eq;a11.17.10}
 \dim
 \Bigl(
 \lefttop{\nminusnibar}\nbigf_{\veca'}
 \cap
 \lefttop{n-1}\nbigf_{a^{\circ}_{n-1}}
 \cap
 \lefttop{n}\nbigf_{a^{\circ}_n}
 \cap \EE\bigl(H(\nbigelambda),\vecomega\bigr)
 \Bigr)
\geq
 v.d(\lambda,\veca,\vecomega).
\end{equation}
From (\ref{eq;a11.17.9}) and (\ref{eq;a11.17.10}),
we can conclude that the equality in (\ref{eq;a11.17.10}) holds,
which implies $(A,a,\leq)$.
We can also conclude that 
the equality in (\ref{eq;9.11.45}) holds,
which implies that
the equality in (\ref{eq;a11.17.6}) holds.
Thus we obtain $(B,a,\leq)$.
Thus the proof of Lemma \ref{lem;9.11.50}
is accomplished.
\hfill\qed

\vspace{.1in}

For any $a$ and some $\epsilon>0$,
the following implications are clear:
\[
\begin{array}{l}
 (B,a,\leq)\Longrightarrow(B,a+\epsilon,<)\\
 (A,a,\leq)\Longrightarrow(A,a+\epsilon,<).
\end{array}
\]
Hence we obtain $(A,a)$ and $(B,a)$ for any $a$.
Thus the induction on $n$ can proceed.
Namely the proof of the claims $(A)$ and $(B)$ 
of Proposition \ref{prop;9.11.40} is accomplished.

\subsubsection{A proof of the claim $(C)$ of Proposition \ref{prop;9.11.40}}
\label{subsubsection;10.12.51}

Now we shall prove the claim $(C)$.
The following inclusion is surjective,
due to $(B)$:
\[
 \sum_{\substack{
 \phi(\vecb)\lneq \phi(\veca)
 }}
 \lefttop{\nbar}\nbigf_{\vecb}
\cap
 \EE\bigl(H(\nbigelambda),\vecomega\bigr)
\lrarr
 \lefttop{\nminusitibar}\nbigf^{\prime\,(1)}_{\phi(\veca)}
\cap
 \EE\bigl(H(\nbigelambda),\vecomega\bigr).
\]
Thus it is isomorphic.

\begin{lem}
Let $\vecb$ be an element of $\real^n$.
Assume $\phi(\vecb)\lneq \phi(\veca)$.
We put $\veca'=(a_1,\ldots,a_{n-2})\in\real^{n-2}$.
Then we have one of the following:
\begin{itemize}
\item
 $\lefttop{\nbar}\nbigf_{\vecb}\subset
 \lefttop{\nminusnibar}\nbigf_{\veca'}
\cap\lefttop{n}\nbigf_{<a_n}
 $.
\item
 $\lefttop{\nbar}\nbigf_{\vecb}\subset 
 \lefttop{\nminusnibar}\nbigf_{\veca'}
\cap\lefttop{n-1}\nbigf_{<a_{n-1}}
\cap\lefttop{n}\nbigf_{a_n}
 $.
\item
 $\lefttop{\nbar}\nbigf_{\vecb}\subset 
 \lefttop{\nminusnibar}\nbigf'_{\veca'}
\cap\lefttop{n-1}\nbigf_{a_{n-1}}
\cap\lefttop{n}\nbigf_{a_n}
 $.
\end{itemize}
\end{lem}
\pf
The condition implies
$b_i\leq a_i$ $(i\leq n-2)$ and
$b_{n-1}+b_n\leq a_{n-1}+a_n$,
and at least one of the inequalities is not equality.
Then we have at least one of the following:
\begin{itemize}
\item
 $b_n\leq a_n$.
\item
 $b_n=a_n$ and $b_{n-1}\leq a_{n-1}$.
\item
 $b_n=a_n$, $b_{n-1}=a_{n-1}$,
 $(b_1,\ldots,b_{n-2})\leq \veca'$
 and $(b_1,\ldots,b_{n-2})\neq \veca'$.
\end{itemize}
Then the claim follows.
\hfill\qed

\begin{lem} \label{lem;9.11.57}
We have the following:
\begin{equation}\label{eq;9.11.56}
 \lefttop{\nbar}\nbigf_{\veca}\cap
 \Bigl(
 \sum_{\substack{
 \phi(\vecb)\lneq\phi(\veca)
 }
 }\lefttop{\nbar}\nbigf_{\vecb}
 \Bigr)
=\sum_{\substack{
 \vecb\lneq\veca
 }}
 \lefttop{\nbar}\nbigf_{\vecb}
=\lefttop{\nbar}\nbigf'_{\veca}.
\end{equation}
\end{lem}
\pf
We have the following:
\begin{equation}\label{eq;9.11.55}
 \sum_{\substack{
 \phi(\vecb)\lneq\phi(\veca)
 }
 }\lefttop{\nbar}\nbigf_{\vecb}
=
\Bigl(
 \lefttop{\nminusnibar}\nbigf_{\veca'}
\cap\lefttop{n}\nbigf_{<a_n}
\Bigr)
+
\Bigl(
\lefttop{\nminusnibar}\nbigf_{\veca'}
\cap \lefttop{n-1}\nbigf_{<a_{n-1}}
\cap \lefttop{n}\nbigf_{a_n}
+
\lefttop{\nminusnibar}\nbigf'_{\veca'}
\cap \lefttop{n-1}\nbigf_{a_{n-1}}
\cap \lefttop{n}\nbigf_{a_n}
\Bigr).
\end{equation}
Note that the second term in the right hand side
of (\ref{eq;9.11.55})
is contained in $\lefttop{\nbar}\nbigf_{\veca}$.
Let us pick elements:
$ x\in
 \lefttop{\nminusnibar}\nbigf_{\veca'}
  \cap
 \lefttop{n-1}\nbigf_{<a_{n-1}}\cap\lefttop{n}\nbigf_{a_n}
 +
 \lefttop{\nminusnibar}\nbigf'_{\veca'}
  \cap
 \lefttop{n-1}\nbigf_{a_{n-1}}\cap\lefttop{n}\nbigf_{a_n}$
and 
$ y\in \lefttop{\nminusnibar}\nbigf_{\veca'}
 \cap(\lefttop{n}\nbigf_{<a_n})$.
Assume $x+y\in\lefttop{\nbar}\nbigf_{\veca}$.
Then we obtain
$y\in
 \lefttop{\nbar}\nbigf_{\veca}\cap
 \lefttop{\nminusnibar}\nbigf_{\veca'}\cap
 \lefttop{n}\nbigf_{<a_n}$.
Hence we have
$y\in \lefttop{\nminusnibar}\nbigf_{\veca'}
 \cap \lefttop{n-1}\nbigf_{a_{n-1}}
 \cap \lefttop{n}\nbigf_{<a_n}$.
Thus the left hand side of (\ref{eq;9.11.56}) is as follows:
\[
  \lefttop{\nminusnibar}\nbigf_{\veca'}
 \cap \lefttop{n-1}\nbigf_{a_{n-1}}
 \cap \lefttop{n}\nbigf_{<a_n}
+
\Bigl(
\lefttop{\nminusnibar}\nbigf_{\veca'}
\cap \lefttop{n-1}\nbigf_{<a_{n-1}}
\cap \lefttop{n}\nbigf_{a_n}
+
\lefttop{\nminusnibar}\nbigf'_{\veca'}
\cap \lefttop{n-1}\nbigf_{a_{n-1}}
\cap \lefttop{n}\nbigf_{a_n}
\Bigr).
\]
It is same as the right hand side of (\ref{eq;9.11.56}).
Thus we are done.
\hfill\qed

\vspace{.1in}

The claim $(C)$ immediately follows Lemma \ref{lem;9.11.57}.
Thus the proof of Proposition \ref{prop;9.11.40} is accomplished.
\hfill\qed

\subsubsection{A consequence}

Let $C_0$ be the diagonal curve:
$C_0:=\{(z,\ldots,z)\in\Delta^n\}$.
The restriction gives the natural isomorphism
$ H(\nbigelambda)\lrarr
 H(\nbigelambda_{|C_0})$.
We have the filtration on the right hand side,
which we denote by $\lefttop{C_0}\nbigf$.

\begin{cor}\label{cor;10.12.30}
The following morphism is injective:
\[
 \frac{\lefttop{\nbar}\nbigf_{\veca}\cap
 \EE\bigl(H(\nbigelambda),\vecomega\bigr)}
{\lefttop{\nbar}\nbigf'_{\veca}\cap
 \EE\bigl(H(\nbigelambda),\vecomega\bigr) }
\lrarr
 \frac{\lefttop{C_0}\nbigf_{|\veca|}\cap
 \EE\bigl(H(\nbigelambda),\vecomega\bigr)}
{\lefttop{C_0}\nbigf'_{|\veca|}\cap
 \EE\bigl(H(\nbigelambda),\vecomega\bigr)
}.
\]
Here $|\veca|$ denotes $\sum_{i=1}^n a_i$
for $\veca=(a_1,\ldots,a_n)$.
\end{cor}
\pf
We have only to use the claims $(C)$ in Proposition \ref{prop;9.11.40}
inductively.
\hfill\qed

%% file: 16.4.tex

\subsubsection{The compatibility of the filtrations
$\big(\lefttop{i}\nbigf\,\big|\,i=1,\ldots,n\big)$}

Let $s$ be an element of
$\lefttop{\nbar}\nbigf_{\veca}\cap\EE\bigl(H(\nbigelambda),\vecomega\bigr)$.
Then we obtain the section
$F(s,\veca)$ of $\prolong{\nbigelambda}$.
For any $i$,
we have the element 
$\vecu\in\nbigk(\nbige,\lambda,0,\lbar)$,
satisfying $\kmsmap^f(\lambda,\vecu)=(\veca,\vecomega)$.
Due to the result in the case of curves (Lemma \ref{lem;9.6.3}),
we have
$\lefttop{i}\deg^F F(s,\veca)\leq
 \paramap(\lambda,u_i)$.
We put $c_i:=\lefttop{i}\deg^F F(s,\veca)$
and $\vecc:=(c_1,\ldots,c_n)$.
We obtain the following map:
\[
 \Phi:
 \lefttop{\nbar}\Gr^{\nbigf}_{\veca}
 \EE(H(\nbigelambda),\vecomega)
\lrarr
 \lefttop{\nbar}\Gr^F_{\vecc}
 (\prolong{\nbigelambda}_{|O}),
\quad
 s\longmapsto F(s,\vecu)(O).
\]

\begin{lem}
The map $\Phi$ is injective.
\end{lem}
\pf
Let $\vecv$ be a frame of $\prolong{\nbigelambda}$,
which is compatible with $(\lefttop{i}F\,|\,i=1,\ldots,n)$.
We describe as follows:
\[
 F(s,\veca):=\sum f_j\cdot v_j.
\]
Assume that $\Phi(s)=0$.
We have $f_j(O)=0$ unless $\vecdeg(v_j)\lneq\vecc$.
In this case,
we have $-\ord(s_{|C_0})<|\veca|$.
It implies
$s_{|C_0}\in \lefttop{C_0}\nbigf'_{|\veca|} H(\nbigelambda_{|C_0})$.
Then we obtain $s=0$ in
$\lefttop{\nbar}\Gr_{\veca}\EE(H(\nbigelambda),\vecomega)$,
due to Corollary \ref{cor;10.12.30}.
\hfill\qed

\vspace{.1in}

We put 
$\gamma_i:=
 \eigenmap(\lambda,u_i)$
and 
$\vecgamma:=(\gamma_1,\ldots,\gamma_n)$.

\begin{lem}
We have $\Image(\Phi)\subset
 \lefttop{\nbar}\EE\bigl(
 \lefttop{\nbar}\Gr_{\vecc}^{F}(\Res_{\nbar} \DD^{\lambda}),
 \vecgamma \bigr)$.
\end{lem}
\pf
We have only to check that
$F(s,\veca)_{|(\lambda,P)}$
is contained in
$\lefttop{i}\EE\bigl(
 \lefttop{i}\Gr^{F}_{c_i}
 (\prolong{\nbigelambda}_{D_i}), 
 \gamma_i
 \bigr)$ for any $P\in D_i^{\circ}$.
It follows from the result in the case of curves
(Corollary \ref{cor;9.11.60}).
\hfill\qed

\vspace{.1in}

Thus we obtain the morphisms
$\Phi_{(\veca,\vecomega)}:
 \lefttop{\nbar}\Gr^{\nbigf}_{\veca}\bigl(
 \lefttop{\nbar}\EE\bigl(H(\nbigelambda),\vecomega\bigr)
 \bigr)
\lrarr
 \lefttop{\nbar}\EE(
 \lefttop{\nbar}\Gr^{F}_{\vecc}(
 \prolong{\nbigelambda}_{|O}
 ),\vecgamma)$.
Then we obtain the following injection:
\[
\bigoplus_{(\veca,\vecomega)}
\Phi_{(\veca,\vecomega)}:
 \bigoplus_{(\veca,\vecomega)}
 \lefttop{\nbar}\Gr^{\nbigf}_{\veca}
 \lefttop{\nbar}\EE\bigl(H(\nbigelambda),\vecomega\bigr)
\lrarr
 \bigoplus_{(\vecc,\vecgamma)}
 \lefttop{\nbar}\EE\bigl(
 \lefttop{\nbar}\Gr^F_{\vecc}(\Res_{\nbar}\DD^{\lambda}),\vecgamma\bigr).
\]

\begin{prop}
The morphisms $\Phi_{(\veca,\vecomega)}$ are isomorphic.
\end{prop}
\pf
Since $\Phi_{(\veca,\vecomega)}$ is injective,
we obtain the following inequalities:
\begin{equation} \label{eq;9.11.65}
 \rank \nbigelambda
\leq
\sum_{(\veca,\vecomega)}
\dim\lefttop{\nbar}\Gr^{\nbigf}_{\veca}
\EE(H(\nbigelambda),\vecomega)
\leq
 \sum_{(\vecc,\vecgamma)}
 \dim\EE(\lefttop{\nbar}\Gr^F_{\vecc}(\Res_{\nbar}\DD^{\lambda}),\vecgamma)
=\rank \nbigelambda.
\end{equation}
Then the proposition immediately follows.
\hfill\qed

\begin{thm}
 The tuple of the filtrations
 $\bigl(\lefttop{i}\nbigf\,\big|\,i=1,\ldots,n\bigr)$ is compatible
 in the sense of Definition {\rm\ref{df;b11.12.1}}.
\end{thm}
\pf
Due to (\ref{eq;9.11.65}),
we have the following equality:
\[
 \sum _{\veca}
 \rank \lefttop{\nbar}\Gr^{\nbigf}_{\veca}
 H(\nbigelambda)
=\rank H(\nbigelambda).
\]
It implies the compatibility of the filtrations,
due to Lemma \ref{lem;10.12.31}.
\hfill\qed

\begin{cor}\mbox{{}}
\begin{itemize}
\item
 The tuple of the filtrations
 $\bigl(\lefttop{i}\nbigf^{(\lambda_0)}\,\big|\,i=1,\ldots,l\bigr)$
 are compatible. 
\item
 We have the following decomposition on
 $\Delta^{\ast}(\lambda_0,\epsilon_0)$:
 \[
  \lefttop{I}\nbigf^{(\lambda_0)}_{\vecb}
 \nbigh^{(\lambda_0)}_{\vecomega\,|\,
 \Delta^{\ast}(\lambda_0,\epsilon_0)
 }=
\bigoplus_{\substack{
\vecu\in S(\vecomega,\vecb)
 }}
 \nbigh_{\eigenmap^f(\lambda,\vecu)}.
 \]
Here we put
$S(\vecomega,\vecb)
:=\bigl\{
 \vecu\in \KMSoverline(\nbige^0,I)\,\big|\,
 \eigenmap^f(\lambda_0,\vecu)=\vecomega,\,\,
 \paramap^f(\lambda_0,\vecu)\leq\vecb\bigr\}$.
\hfill\qed
\end{itemize}
\end{cor}

%% file: a57.1.tex

\subsubsection{Weak norm estimate}

Let $\vecs=(s_i)$ be a frame of $\nbigh_{|\Delta(\lambda_0,\epsilon_0)}$,
which is compatible with $\EEzero$ and $\nbigfzero$.
For each $s_i$, we have the element $\vecu(s_i)\in\KMS(\nbige^0,\lbar)$
such that
$\deg^{\EEzero,\nbigfzero}(s_i)
=\kmsmap^f(\lambda_0,\vecu(s_i))$.
We put as follows:
\begin{equation} \label{eq;10.24.5}
 \vecs':=\bigl(s'_i\bigr),
\quad
 s'_i:=
 s_i\cdot\prod_{j=1}^l |z_j|^{\paramap^f(\lambda,u_j(s_i))}.
\end{equation}

\begin{lem} \label{lem;10.24.6}
$\vecs'$ gives a $C^{\infty}$-frame of
 $\nbige_{|\nbigx(\lambda_0,\epsilon_0)-\nbigd(\lambda_0,\epsilon_0)}$,
which is adapted up to log order.
\end{lem}
\pf
It is easy to see the following inequality for some positive
constants $C_1$ and $M_1$:
\[
 H(h,\vecs')\leq
 C_1\cdot \Bigl(
 -\sum_{j=1}^l\log|z_j|
 \Bigr)^{M_1}.
\]
Let $\vecs^{\lor}$ be the dual frame of $\vecs$.
Due to the result in the case of curves,
we have $\vecu(s^{\lor}_i)=-\vecu(s_i)$.

Let $\vecs^{\lor\,\prime}$ be the modification
as in (\ref{eq;10.24.5}).
We obtain the inequality
$H(h,\vecs^{\lor\,\prime})
 \leq C_2\cdot\bigl(-\sum_{j=1}^l\log|z_j|\bigr)^{M_2}$
for some positive constants $C_2$ and $M_2$.
Since we have $H(h,\vecs')\cdot H(h,\vecs^{\lor\,\prime})=1$,
we obtain the result.
\hfill\qed

%% file: 16.6.tex

\subsubsection{The induced vector bundle $\lefttop{I}\nbigg_{\vecu}(\nbigh)$}

\label{subsubsection;10.18.9}

Let $\vecu$ be an element of $\KMSoverline(\nbige^0,I)$.
We put as follows:
\[
 \lefttop{I}\nbigg^{(\lambda_0)}_{\vecu}\nbigh
:=\lefttop{I}
 \Gr^{\nbigf^{(\lambda_0)}}_{\paramap^f(\lambda_0,\vecu)}
 \bigl(
\nbigh^{(\lambda_0)}_{\eigenmap^f(\lambda_0,\vecu)}\bigr).
\]
Let $\lambda\in\Delta^{\ast}(\lambda_0,\epsilon_0)$ be generic.
Let us pick $\epsilon'_0>0$ such that
$\Delta(\lambda,\epsilon_0')\subset
 \Delta^{\ast}(\lambda_0,\epsilon_0)$.

\begin{lem} \label{lem;10.12.50}
We have the following isomorphism:
\[
 \lefttop{I}\nbigg^{(\lambda_0)}_{\vecu}
 \nbigh_{|\Delta(\lambda,\epsilon_0')}
\simeq
 \lefttop{I}\nbigg^{(\lambda)}_{\vecu}\nbigh.
\]
\end{lem}
\pf
For $\vecomega=\eigenmap^f(\lambda_0,\vecu)$,
we have the following decomposition on
$\Delta^{\ast}(\lambda_0,\epsilon_0)$:
\[
 \nbigh^{(\lambda_0)}_{
  \vecomega\,|\,\Delta^{\ast}(\lambda_0,\epsilon_0)}
=\bigoplus_{
 \eigenmap^f(\lambda_0,\vecu)=\vecomega
 }
 \nbigh_{\eigenmap^f(\lambda,\vecu)}.
\]
Then we obtain the natural isomorphisms:
\[
 \lefttop{I}\nbigg_{\vecu}^{(\lambda_0)}
 \nbigh_{|\Delta(\lambda,\epsilon_0')}
\simeq
 \lefttop{I}\Gr^{\nbigf^{(\lambda_0)}}_{\vecomega}
 \nbigh_{\eigenmap^f(\lambda_0,\vecu)\,|\,
 \Delta(\lambda,\epsilon_0')}
\simeq
 \nbigh_{\eigenmap^f(\lambda,\vecu)\,|\,\Delta(\lambda,\epsilon_0')}
\simeq
 \lefttop{I}\nbigg^{(\lambda)}_{\vecu}\nbigh.
\]
Thus we are done.
\hfill\qed

\vspace{.1in}

Due to Lemma \ref{lem;10.12.50}, we obtain the vector bundle
$\lefttop{I}\nbigg_{\vecu}\nbigh$ over $\cnum_{\lambda}^{\ast}$.

\subsubsection{The induced pairing and nilpotent maps}

\label{subsubsection;10.24.6}

We have the natural pairing:
$ \nbigh(E)\otimes\nbigh(E^{\lor})
\lrarr\nbigo_{\cnum_{\lambda}^{\ast}}$.
Pick a point $\lambda_0\in\cnum_{\lambda}^{\ast}$
and a sufficiently small positive number $\epsilon_0>0$.
On $\Delta(\lambda_0,\epsilon_0)$,
we have the filtrations 
$\lefttop{i}\nbigf^{(\lambda_0)}$
and the decompositions $\lefttop{i}\EE^{(\lambda_0)}$,
which was preserved by the pairing.
Hence we obtain the following induced pairing:
\[
S:\lefttop{I}\nbigg_{\vecu}\nbigh(E)\otimes
 \lefttop{I}\nbigg_{-\vecu}\nbigh(E^{\lor})
 \lrarr\nbigo_{\cnum_{\lambda}^{\ast}}.
\]

The monodromies $M_i$ induces the endomorphisms
$M_{i\,\vecu}$
on $\lefttop{I}\nbigg_{\vecu}\nbigh$.
We denote the unipotent part 
by $M_{i\,\vecu}^u$.
Then we obtain the following nilpotent maps:
\[
 \nbign_{\vecu,i}:=
\frac{-1}{2\pi\sqrt{-1}}
\log M^u_{i\,\vecu}.
\]

\begin{lem}
We have the relation
$S\bigl(\nbign_{\vecu\,i}\otimes id\bigr)
+S\bigl(id\otimes\nbign_{-\vecu\,i}\bigr)=0$.
\end{lem}
\pf
Since the monodromy endomorphisms preserve the pairings,
the claim is obtained.
\hfill\qed

%% file: a57.2.tex

\subsubsection{Functoriality for the dual}
\label{subsubsection;10.26.22}

Due to the subsubsection \ref{subsubsection;10.24.6},
we obtain the naturally defined morphism
$\lefttop{\lbar}\nbigg_{-\vecu}(\nbigh(E^{\lor}))
\lrarr
\lefttop{\lbar}\nbigg_{\vecu}(\nbigh)^{\lor}$.

\begin{lem} \label{lem;10.24.10}
The morphism
$\lefttop{\lbar}\nbigg_{-\vecu}(\nbigh)(E^{\lor})
\lrarr
\lefttop{\lbar}\nbigg_{\vecu}(\nbigh(E))^{\lor}$
is isomorphic.
\end{lem}
\pf
It follows from Lemma \ref{lem;10.24.6}.
\hfill\qed

\subsubsection{Functoriality for tensor products}
\label{subsubsection;10.26.23}

The morphism $\nbigh(E_1)\otimes\nbigh(E_2)\lrarr\nbigh(E_1\otimes E_2)$
preserves the $\nbigfzero$ and $\EEzero$.
We have the naturally defined morphism:
\begin{equation}\label{eq;10.24.7}
\bigoplus_{\substack{
 \vecu_i\in\KMSoverline{(\nbige^0,\lbar)},\\
 \vecu_1+\vecu_2=\vecu}
 }
 \lefttop{\lbar}\nbigg_{\vecu_1}\nbigh(E_1)
\otimes
 \lefttop{\lbar}\nbigg_{\vecu_2}\nbigh(E_2)
\lrarr
 \lefttop{\lbar}\nbigg_{\vecu_1+\vecu_2}\nbigh\bigl(E_1\otimes E_2\bigr).
\end{equation}

\begin{lem} \label{lem;10.24.11}
The morphism {\rm (\ref{eq;10.24.7})} is isomorphic.
\end{lem}
\pf
Under the isomorphism
$\nbigh(E_1)\otimes\nbigh(E_2)\simeq\nbigh(E_1\otimes E_2)$,
the induced filtration on the left hand side
and the filtration on the right hand side are same,
due to the result in the case of curves.
Then we obtain the lemma.
\hfill\qed

\subsubsection{Functoriality for pull backs}
\label{subsubsection;10.26.24}

We use the setting in the subsubsection \ref{subsubsection;10.26.6}.
We have the naturally defined isomorphism 
$\nbigh(E)\simeq \nbigh(\psi^{\ast}(E))$.

\begin{lem}\label{lem;10.26.10}
We have the naturally defined isomorphism:
\[
\bigoplus_{\psi_{\vecc}^{\ast}(\vecu)=\vecu_1}
 \lefttop{\lbar}\nbigg_{\vecu}\bigl(\nbigh(E)\bigr)
 \simeq
 \lefttop{\lbar}\nbigg_{\vecu_1}\bigl(\nbigh(\psi^{\ast}E)\bigr).
\]
Here $\psi_{\vecc}^{\ast}(\vecu)$ denotes the element
whose $i$-th component is $c_i\cdot u_i$ for $\vecu=(u_1,\ldots,u_l)$.
\end{lem}
\pf
It is easy to see that the isomorphism
is strictly compatible with the filtration $\nbigfzero$.
It is also easy to check that the isomorphism is compatible
with the decompositions $\EEzero$.
Thus we obtain the isomorphism desired.
\hfill\qed

%% file: 16.5.tex

\subsubsection{The compatibility of the decompositions}

\label{subsubsection;10.15.3}

Let $\vecc$ be an element of $\real^l$.
Let $\vecs$ be a base of $H(\nbigelambda)$,
compatible with $\EE$ and $\nbigf$.
We put $\veca_i:=\vecdeg^{\nbigf}(s_i)$.
We put $v_i:=F(s_i,\veca_i-\vecc)$,
and then $\vecv=(v_i)$ is a tuple of sections
of $\prolongg{\vecc}{\nbigelambda}$.

\begin{lem}
The tuple $\vecv$ gives a frame of $\prolongg{\vecc}{\nbigelambda}$.
The restriction of $\vecv$ to $\nbigd_i(\lambda_0,\epsilon_0)$
is compatible with the decomposition
$\lefttop{i}\EE$ and the filtration $\lefttop{i}F$
of the vector bundle
$\prolongg{\vecc}{\nbigelambda}_{|\nbigd_i^{\lambda}}$.
\end{lem}
\pf
We have only to show that
$\vecv_{|\pi_j^{-1}(P)}$ gives a frame of
$\prolongg{c_j}{\nbige^{\lambda}_{|\pi_j^{-1}(P)}}$
for any point $P\in D_j^{\circ}$.
Then it follows from the result in the case of curves
(the subsubsection \ref{subsubsection;a11.15.3}).
\hfill\qed

\vspace{.1in}

By the action of monodromies,
we have the following decomposition:
\begin{equation} \label{eq;10.15.1}
 \nbigelambda
=\bigoplus_{\vecomega\in\Sp(\vecM^{\lambda})}
 \nbigelambda_{\vecomega},
\quad\quad
 \nbigelambda_{\vecomega\,|\,P}
=\EE\bigl(\nbigelambda_{|P},\,\,\vecomega\bigr).
\end{equation}

\begin{lem}
The decomposition {\rm(\ref{eq;10.15.1})} is prolonged
to the decomposition of the vector bundle
$\prolongg{\vecc}{\nbigelambda}$.
Namely, we have the following decomposition of
$\prolongg{\vecc}{\nbigelambda}$:
\[
 \prolongg{\vecc}{\nbigelambda}=
\bigoplus_{\vecomega}
 \prolongg{\vecc}{\nbigelambda_{\vecomega}}.
\]
\end{lem}
\pf
We have only to use the frame $\vecv$ above.
\hfill\qed

\vspace{.1in}

Then we obtain the decomposition
of the restriction of $\prolongg{\vecc}{\nbige}$ to $D_i$:
\[
 \prolongg{\vecc}{\nbigelambda}_{|D_i}
=\bigoplus_{\vecomega\in \Sp(\vecM^{\lambda})}
 \prolongg{\vecc}{\nbigelambda}_{\vecomega\,|\,D_i}.
\]
On the other hand,
we have the decomposition of $\prolongg{\vecc}{\nbige}_{|D_i}$
induced by the action of the residue $\Res_i(\DDlambda)$:
\[
 \prolongg{\vecc}{\nbigelambda}_{|D_i}
 =\bigoplus_{\beta\in\Sp(\prolongg{\vecc}{\nbigelambda},i)}
 \EE\bigl(\prolongg{\vecc}{\nbigelambda}_{|D_i},\beta\bigr).
\]
The relation of the two decompositions are given
in the following lemma.
\begin{lem} \label{lem;10.15.5}
We have the following equality:
\begin{equation} \label{eq;10.15.2}
 \bigoplus_{q_i(\vecomega)=\omega}
 \prolongg{\vecc}{\nbigelambda}_{\vecomega\,|\,D_i}
=\bigoplus_{\beta\in L(\lambda,\omega)}
 \EE(\prolongg{\vecc}{\nbigelambda}_{|D_i},\beta),
\end{equation}
Here we put
$L(\lambda,\omega):=
 \bigl\{\beta\in\Sp(\prolongg{\vecc}{\nbigelambda},i)\,\big|\,
 \exp(-2\pi\sqrt{-1}\lambda^{-1}\cdot\beta)=\omega
 \bigr\}$.
\end{lem}
\pf
Since the both sides of (\ref{eq;10.15.2}) are 
vector subbundles of $\prolongg{\vecc}{\nbigelambda}_{|D_i}$,
we have only to show the equality 
for the fibers over the points $P\in D_i^{\circ}$.
Thus we have only to consider the case $\dim(X)=1$.
Then it follows from Corollary \ref{cor;a11.17.20}.
\hfill\qed

\begin{cor}
Let $I$ be a subset of $\lbar$.
Let $\vecomega^{\circ}=(\omega^{\circ}_i\,|\,i\in I)$
be an element of $\Sp^f(\nbigelambda,I)$.
Then we have the following:
\[
 \bigoplus_{q_i(\vecomega)=\omega_i^{\circ}}
 \prolongg{\vecc}{\nbigelambda}_{\vecomega\,|\,D_I}
=\bigoplus_{\vecbeta\in L(\lambda,\vecomega^{\circ})
 }
 \EE\bigl(\prolongg{\vecc}{\nbigelambda}_{|D_I},\vecbeta\bigr).
\]
Here we put
$L(\lambda,\vecomega^{\circ})
 :=
\bigl\{
 \vecbeta\in\Sp(\prolongg{\vecc}{\nbigelambda},I)\,\big|\,
 \exp\bigl(-2\pi\sqrt{-1}\lambda^{-1}\cdot\beta_i\bigr)=\omega_i^{\circ}
 \,\,(i\in I)
\bigr\}$.
\hfill\qed
\end{cor}

\subsubsection{The compatibility of the filtrations}

We have the filtration
$\lefttop{i}\nbigf\bigl(
 \nbigelambda_{\vecomega}\bigr)$
of $\nbigelambda_{\vecomega}$
induced by
the filtration
$\lefttop{i}\nbigf\bigl(
 \EE\bigl(H(\nbigelambda),\vecomega\bigr)\bigr)$
of $\EE\bigl(H(\nbigelambda),\vecomega\bigr)$.

\begin{prop}\mbox{{}} \label{prop;10.15.10}
\begin{itemize}
\item
 The vector subbundle
 $\lefttop{i}\nbigf_b\bigl(
 \nbigelambda_{\vecomega}\bigr)$ of $\nbigelambda$
 can be prolonged to the subbundle
 $\prolongg{\vecc}{
 \lefttop{i}\nbigf_b\bigl(
 \nbigelambda_{\vecomega}\bigr)
 }$
 of $\prolongg{\vecc}{\nbigelambda_{\vecomega}}$.
\item
 The family
 $\prolongg{\vecc}{
 \lefttop{i}\nbigf\bigl(
 \nbigelambda_{\vecomega}\bigr)}
 =\bigl\{
 \prolongg{\vecc}{
 \lefttop{i}\nbigf_a\bigl(
 \nbigelambda_{\vecomega}\bigr)}
 \,\big|\,a\in\real\bigr\}$
 gives the filtration of the vector bundle
 $\prolongg{\vecc}{\nbigelambda_{\vecomega}}$
in the category of vector bundles.
\item
 The tuple of the filtrations 
$\bigl(
 \prolongg{\vecc}{\lefttop{i}\nbigf\bigl(
 \nbigelambda_{\vecomega}\bigr)}
 \,\big|\,
 i=1,\ldots,l \bigr)$ 
 of the vector bundle $\prolongg{\vecc}{\nbigelambda}$
 is compatible.
\end{itemize}
\end{prop}
\pf
It is easy to check the claims
by using the frame $\vecv$ given in the first part
of the subsubsection \ref{subsubsection;10.15.3}.
\hfill\qed

\vspace{.1in}

Then we have the filtrations
$\prolongg{\vecc}{
 \lefttop{i}\nbigf\bigl(\nbigelambda_{\vecomega}\bigr)}_{|D_i}$
of the vector bundle
$\prolongg{\vecc}{\nbigelambda}_{|D_i}$.
On the other hand,
we have the parabolic filtration $\lefttop{i}F$
of $\prolongg{\vecc}{\nbigelambda}_{|D_i}$.
The relation of two filtrations are given
in the proposition \ref{prop;10.15.4}

We recall that the number
$\alpha(a,\omega)\in \cnum$ for $(a,\omega)\in\real\times\cnum$
is determined by the following conditions:
\[
 \exp\bigl(-2\pi\sqrt{-1}\cdot\alpha(a,\omega)\bigr)=\omega,
\quad\quad
 a\leq \Re\bigl(\alpha(a,\omega)\bigr)<a+1.
\]
We put $d(a,\omega):=a-\Re(\alpha(a,\omega))\in\real$.

\begin{prop}\label{prop;10.15.4}
We have the following equality:
\[
 \bigoplus_{q_i(\vecomega)=\omega}
 \prolongg{\vecc}{\lefttop{i}\nbigf_{b}\bigl(
 \nbigelambda_{\vecomega}\bigr)_{|D_i}}
=\lefttop{i}
 F_{d(b-c_i,\omega)}\bigl(
 \EE\bigl(\prolongg{\vecc}{\nbigelambda}_{|D_i},
 \lambda\cdot\alpha(b-c_i,\omega)\bigr)\bigr)
\oplus
\bigoplus_{\beta\in K(\lambda,\omega,b)
 }
 \EE\bigl(\prolongg{\vecc}\nbigelambda_{|D_i},\beta\bigr).
\]
Here we put
$K(\lambda,\omega,b):=
 \bigl\{
 \beta\in\Sp(\nbigelambda,i)\,\big|\,
 \exp\bigl(-2\pi\sqrt{-1}\lambda^{-1}\cdot\beta\bigr)=\omega,\,\,
 \Re\bigl(\lambda^{-1}\cdot\beta\bigr)
 <\Re\bigl(\alpha(b-c_i,\omega)\bigr)
 \bigr\}$.
\end{prop}
\pf
The claim can be easily reduced to the case $\dim(X)=1$,
as in the proof of Lemma \ref{lem;10.15.5}.
\hfill\qed

\subsubsection{The induced vector bundle
 $\lefttop{\lbar}\nbigg_{\vecu}(\nbigelambda)$}

Let $\vecu$ be an element of $\KMS(\nbige^0,\lbar)$.
We put as follows:
\[
 \vecb:=\paramap^f(\lambda,\vecu),\quad
 \vecomega:=\eigenmap^f(\lambda,\vecu),\quad
 \vecc:=\paramap(\lambda,\vecu).
\]
We obtain the holomorphic vector bundle
$\lefttop{\lbar}\nbigg_{\vecu}(\nbigelambda)$
over $X$,
given as follows:
\[
 \lefttop{\lbar}\nbigg_{\vecu}\bigl(
 \nbigelambda\bigr)
:=\frac{\prolongg{\vecc}{
 \lefttop{\lbar}\nbigf_{\vecb}\bigl(\nbigelambda_{\vecomega}\bigr)}}
 {\sum_{\vecb'\lneq\vecb}
 \prolongg{\vecc}{
 \lefttop{\lbar}\nbigf_{\vecb'}\bigl(
 \nbigelambda_{\vecomega}
 \bigr)}}.
\]
\begin{lem}
We have the natural isomorphism:
\[
 \lefttop{\lbar}\nbigg_{\vecu}\bigl(
 \nbigelambda\bigr)_{|D_{\lbar}}
 \simeq
 \lefttop{\lbar}\nbigg_{\vecu\,|\,\{\lambda\}\times D_{\lbar}}.
\]
\end{lem}
\pf
It follows from Lemma \ref{lem;10.15.5},
Proposition \ref{prop;10.15.4} and the definitions.
\hfill\qed

\vspace{.1in}

We have the induced $\lambda$-connection $\DDlambda$
of $\lefttop{\lbar}\nbigg_{\vecu}(\nbigelambda)$,
which is flat and regular.
Then we obtain the residues $\Res_i(\DDlambda)$.

\begin{cor}
The endomorphism $\Res_i(\DDlambda)$ has the unique
eigenvalue $\eigenmap\bigl(\lambda,q_i(\vecu)\bigr)$.
\hfill\qed
\end{cor}

Let
$H\bigl(\lefttop{\lbar}\nbigg_{\vecu}(\nbigelambda)\bigr)$
be the space of the multi-valued flat sections of
$\lefttop{\lbar}\nbigg_{\vecu}(\nbigelambda)$.
Naturally we have the following isomorphism:
\[
 H\bigl(\lefttop{\lbar}\nbigg_{\vecu}(\nbigelambda)\bigr)
\simeq
 \lefttop{\lbar}\Gr^{\nbigf}_{\vecb}
 \EE\bigl(H(\nbigelambda),\vecomega\bigr)
=\nbigg_{\vecu}\bigl(\nbigh\bigr)_{|\lambda}.
\]

%% file: 17.tex

\subsubsection{Prolongation of the decompositions and the filtrations}

\label{subsubsection;10.18.15}

The decomposition
$\nbigh=\bigoplus_{\vecomega\in\Sp(\vecM^{\lambda_0})}
 \nbigh^{(\lambda_0)}_{\vecomega}$
and the filtration $\lefttop{i}\nbigf^{(\lambda_0)}$ $(i=1,\ldots,l)$
on $\nbigh_{|\Delta(\lambda_0,\epsilon_0)}$
induce those on $\nbige_{|
 \nbigx(\lambda_0,\epsilon_0)-\nbigd(\lambda_0,\epsilon_0)}$.
Namely, we have the decomposition
$\nbige=
 \bigoplus_{\vecomega\in\Sp(\vecM^{\lambda_0})}
 \nbige^{(\lambda_0)}_{\vecomega}$
and the filtration
$\lefttop{i}\nbigf^{(\lambda_0)}(\nbige_{\vecomega})$
over $\nbigx(\lambda_0,\epsilon_0)-\nbigd(\lambda_0,\epsilon_0)$.

Let $\vecb$ be an element of $\real^l$
such that $b_i\not\in\Par(\nbige^{\lambda_0},i)$.
If $\epsilon_0$ is sufficiently small,
we have the locally free sheaf $\prolongg{\vecb}{\nbige}$
over $\nbigx(\lambda_0,\epsilon_0)$.

\begin{prop}\mbox{{}} \label{prop;10.15.16}
\begin{enumerate}
\item \label{8.20.1}
 The vector subbundle
 $\lefttop{i}\nbigf^{(\lambda_0)}_c\nbige^{(\lambda_0)}_{\vecomega}$
 of ${\nbige}$
 is prolonged to the subbundle
$\prolongg{\vecb}{\lefttop{i}\nbigf_c^{(\lambda_0)}
 \nbige^{(\lambda_0)}_{\vecomega} }$
of $\prolongg{\vecb}{\nbige}$ on $\nbigx(\lambda_0,\epsilon_0)$.
\item \label{8.20.2}
 We have
 $\prolongg{\vecb}{\nbige^{(\lambda_0)}_{\vecomega\,|\,\nbigx^{\lambda_0}}}=
 \prolongg{\vecb}{\nbige^{\lambda_0}_{\vecomega}}$.
On the other hand,
 we have the following decomposition,
 for any $\lambda\in\Delta^{\ast}(\lambda_0,\epsilon_0)$:
 \[
  \prolongg{\vecb}{\nbige^{(\lambda_0)}
   _{\vecomega\,|\,\nbigx^{\lambda}}}=
 \bigoplus_{\vecu\in\nbigs(\vecomega)}
 \prolongg{\vecb}{\nbigelambda_{\eigenmap^f(\lambda_0,\vecu)}}.
 \]
Here we put $\nbigs(\vecomega):=\bigl\{
\vecu\in \KMS(\nbige^0,\lbar)\,\big|\,
 \eigenmap^f(\lambda_0,\vecu)=\vecomega \bigr\}$.
\item \label{8.20.3}
 We have
 $\prolongg{\vecb}{\lefttop{i}\nbigf^{(\lambda_0)}\bigl(
  \nbige^{(\lambda_0)}_{\vecomega}\bigr)_{|\nbigx^{\lambda_0}}}
=\prolongg{\vecb}{\lefttop{i}\nbigf\bigl(
 \nbige^{\lambda_0}_{\vecomega}\bigr)}$.
 On the other hand,
  we have the decomposition,
 for any $\lambda\in\Delta^{\ast}(\lambda_0,\epsilon_0)$:
\[
 \prolongg{\vecb}{\lefttop{i}\nbigf^{(\lambda_0)}_c
  \nbige^{(\lambda_0)}_{\vecomega\,|\,\lambda}}
=\bigoplus_{\vecu\in \nbigs(\vecomega,b,i)}
 \prolongg{\vecb}{\lefttop{i}\nbigf_{\paramap^f(\lambda,q_i(\vecu))}
 \bigl(
 \nbigelambda_{\eigenmap^f(\lambda,\vecu)}}\bigr)
=\bigoplus_{\vecu\in \nbigs(\vecomega,b,i)}
 \prolongg{\vecb}{\nbigelambda_{\eigenmap^f(\lambda,\vecu)}}.
\]
Here we put
$\nbigs(\vecomega,b,i):=\bigl\{\vecu\in\KMS(\nbige^0,\lbar)\,\big|\,
 \eigenmap^f(\lambda_0,\vecu)=\vecomega,\,\,
 \paramap^f\bigl(\lambda_0,q_i(\vecu)\bigr)\leq b
 \bigr\}$.
\item \label{8.20.4}
 The tuple of the filtrations
 $\bigl(
 \prolongg{\vecb}{\lefttop{i}\nbigf^{(\lambda_0)}}\,\big|\,
 i=1,\ldots,l \bigr)$ is compatible.
\end{enumerate}
\end{prop}
\pf
Once we know the claim \ref{8.20.1},
then the rests follow from the result at the specializations
(Lemma \ref{lem;10.15.5},
Proposition \ref{prop;10.15.10}
and Proposition \ref{prop;10.15.4}).
Hence we have only to check the claim \ref{8.20.1}.
We may assume that $\vecb=0$.

\begin{lem} \label{lem;10.15.11}
Let $\epsilon_i>0$ be positive numbers such that
$\rank(E)\cdot\epsilon_i<1$.
Assume that $\Par(\prolong{\nbigelambda},i)$ are $\epsilon_i$-small
for any $\lambda\in\Delta(\lambda_0,\epsilon_0)$.
Then the filtration
$\lefttop{i}\nbigf^{(\lambda_0)}_a\nbige_{\vecomega}$
can be prolonged to the subbundle of $\prolong{\nbige}$
over $\nbigx(\lambda_0,\epsilon_i)$.
\end{lem}
\pf
Let $\vecs$ be a frame of $\nbigh$ over $\Delta(\lambda_0,\epsilon_0)$
compatible with $\EE^{(\lambda_0)}$ and $\nbigf^{(\lambda_0)}$.
We put
$R:=
 \rank\bigl(\lefttop{i}\nbigf^{(\lambda_0)}_a\bigl(
 \nbige^{(\lambda_0)}_{\vecomega}\bigr)\bigr)$.
Due to the assumption of Lemma \ref{lem;10.15.11},
we have
$\bigwedge^R\bigl(
 \prolong{\nbige}\bigr)
=\Bigl(
 \prolong{
 \bigwedge^R\nbige}\Bigr)$
over $\nbigx(\lambda_0,\epsilon_0)$.

From the frame $\vecs$ of $\nbigh$ over $\Delta(\lambda_0,\epsilon_0)$,
we obtain the naturally induced frame
$\tilde{\vecs}=(\tilde{s}_j)$ of $\bigwedge^R\nbigh$
over $\Delta(\lambda_0,\epsilon_0)$,
which is compatible with $\EE^{(\lambda_0)}$ and $\nbigf^{(\lambda_0)}$.
There exists $j_0$ such that 
$\tilde{s}_{j_0}$ gives a frame of
the line bundle
$\bigwedge^R\bigl(
 \lefttop{i}\nbigf^{(\lambda_0)}_a\nbigh^{(\lambda_0)}_{\vecomega}
 \bigr)$.

There exists $\vecu_0\in \KMS\bigl(\prolong{\nbige^0},\lbar\bigr)$
such that
$\vecdeg^{\nbigf^{(\lambda_0)},\EE^{(\lambda_0)}}
 (\tilde{s}_{j_0})=\kmsmap^f(\lambda_0,\vecu_0)$.
We have the following:
\[
 \vecu_0=
 \sum_{\vecu\in \nbigs(\vecomega,a,i)}
 \vecu,
\quad
\nbigs(\vecomega,a,i):=
\bigl\{
 \vecu\in \KMS(\prolong{\nbige^0},\lbar)\,\big|\,
 \eigenmap^f(\lambda_0,\vecu)=\vecomega,\,\,
 \paramap^f(\lambda_0,q_i(\vecu))\leq a
\bigr\}.
\]
We put as follows:
$ v=\exp\Bigl(
 \log z\cdot
 \bigl(
 \lambda^{-1}\cdot\eigenmap(\lambda,\vecu_0)
+\nu\bigl(\paramap(\lambda,\vecu_0)\bigr)
 \bigr)
 \Bigr)\cdot \tilde{s}_{j_0}$.
Then $v$ is a section of $\prolong{\Bigl(\bigwedge^R\nbige\Bigr)}$
over $\nbigx(\lambda_0,\epsilon_0)$.
It is easy to see that Lemma \ref{lem;10.15.11} can be reduced 
to the following lemma.
\begin{lem}\label{lem;10.15.15}
 There exists $\eta>0$ and a neighbourhood $U$ of $O$ in $X$,
such that
 $v_{|(\lambda,P)}\neq 0$ for
 any $(\lambda,P)\in \Delta(\lambda_0,\eta)\times U$.
\end{lem}
\pf
The section $v_{|\nbigx^{\lambda_0}}$ gives an element of
the frame of $\prolong{\nbigelambdazero}$ induced by $\tilde{s}$,
as in the subsubsection \ref{subsubsection;10.15.3}.
Then  we obtain Lemma \ref{lem;10.15.15} and Lemma \ref{lem;10.15.11}.
\hfill\qed

\vspace{.1in}

Let us return to the proof of Proposition \ref{prop;10.15.16}.
Note the following:
We can pick $\vecc\in\seisuu_{>0}^l$
and $\vecu\in \real^l$
such that $\prolong{\psi_{\vecc}^{-1}\nbige\otimes L(\vecu)}$
satisfies the condition 
of Lemma \ref{lem;10.15.11}
on $\Delta(\lambda_0,\epsilon_0)$.
Hence
$\psi_{\vecc}^{-1}(\lefttop{i}\nbigf_a^{(\lambda_0)}\nbige_{\vecomega})$
can be prolonged to the $\mu_{\vecc}$-equivariant subbundle
of $\prolong{\psi_{\vecc}^{-1}\nbige}$.
Then we pick the equivariant frame,
and take the descent of the frame.
(See the argument in the subsubsection \ref{subsubsection;c11.16.2}).
Then it follows that
$\lefttop{i}\nbigf^{(\lambda_0)}\nbige^{(\lambda_0)}_{\vecomega}$
can be prolonged to the subbundle
of $\prolong{\nbige}_{\vecomega}^{(\lambda_0)}$.
Thus we obtain the claim \ref{8.20.1} of Proposition \ref{prop;10.15.16}.
\hfill\qed

%% file: 17.1.tex

\subsubsection{The induced bundle $\lefttop{\lbar}\nbigg_{\vecu}(\nbige)$}

\label{subsubsection;10.18.10}

Let $\vecu$ be an element of $\KMS(\nbige^0,\lbar)$.
Let $\lambda_0$ be an element of $\cnum_{\lambda}^{\ast}$.
We put $(\vecb,\vecomega)=\kmsmap^f(\lambda_0,\vecu)$.
Pick a sufficiently small $0<\epsilon<1$
such that
$c_i=\paramap(\lambda_0,u_i)+\epsilon'\not\in\Par(\nbigelambdazero,i)$
for any $i$
and for any $0<\epsilon'\leq \epsilon$.
We put $\vecc'=\paramap(\lambda_0,\vecu)+\epsilon\cdot\vecdelta$.
Pick sufficiently small $\epsilon_0$.
Then we have the following vector bundle
over $\nbigx(\lambda_0,\epsilon_0)$:
\[
 \lefttop{\lbar}\nbigg^{(\lambda_0)}_{\vecu}(\nbige)
:=
 \lefttop{\lbar}
  \Gr^{\nbigf^{(\lambda_0)}}_{\vecb}\bigl(
 \prolongg{\vecc'}{\nbige^{(\lambda_0)}_{\vecomega}}
 \bigr).
\]

The following lemma is clear from our construction.
\begin{lem}
It is independent of a choice of $\epsilon$
on a neighbourhood of $\lambda_0$.
\hfill\qed
\end{lem}

We may assume that
any point $\lambda\in\Delta^{\ast}(\lambda_0,\epsilon_i)$
is generic.
Pick a positive number $\epsilon_0'$ such that
$\Delta(\lambda,\epsilon_0')\subset
 \Delta^{\ast}(\lambda_0,\epsilon_0)$.
We may assume that we have the vector bundle
$\nbigg^{(\lambda)}_{\vecu}(\nbige)$
on $\nbigx(\lambda,\epsilon_0')$.

\begin{lem}
We have the following natural isomorphism:
\[
 \nbigg^{(\lambda_0)}_{\vecu}
 (\nbige)_{|\nbigx(\lambda,\epsilon_0')}
\simeq
 \nbigg^{(\lambda)}_{\vecu}(\nbige).
\]
\end{lem}
\pf
We have the following decomposition:
\[
 \prolongg{\vecc'}{\nbige^{(\lambda_0)}}
  _{|\nbigx(\lambda,\epsilon_0')}
=\bigoplus_{\vecu'\in\nbigs(\vecomega)}
 \prolongg{\vecc'}{
 \nbige^{(\lambda)}_{\eigenmap^f(\lambda,\vecu')}},
\quad
\nbigs(\vecomega):=\bigl\{\vecu'\in\KMSoverline(\nbige^0,\lbar)\,\big|\,
 \eigenmap^f(\lambda_0,\vecu')=\vecomega
 \bigr\}.
\]
We have the following isomorphism:
\[
 \nbigg^{(\lambda_0)}_{\vecu}(\nbige)_{|\nbigx(\lambda,\epsilon_0')}
\simeq
 \lefttop{\lbar}\prolongg{\vecc'}{
 \Gr^{\nbigf^{(\lambda_0)}}_{\vecb}\bigl(
 \nbige^{(\lambda_0)}_{\vecomega}\bigr)
 }_{|\nbigx(\lambda,\epsilon_0')}
\simeq
 \prolongg{\vecc'}{
 \nbige^{(\lambda)}_{\eigenmap^f(\lambda,\vecu)}}
\simeq
 \nbigg^{(\lambda)}_{\vecu}(\nbige).
\]
Hence we are done.
\hfill\qed

\vspace{.1in}

Hence we obtain the vector bundle
$\lefttop{\lbar}\nbigg_{\vecu}(\nbige)$
over $\cnum_{\lambda}^{\ast}\times X$,
and the induced regular $\lambda$-connection $\DD$
on $\lefttop{\lbar}\nbigg_{\vecu}(\nbige)$.

\begin{lem}
Let $\vecu$ be an element of $\KMS(\nbige^0,\lbar)$.
\begin{itemize}
\item
We have the following isomorphism:
\[
 \bigl(\lefttop{\lbar}\nbigg_{\vecu}(\nbige)_{|\nbigd_{\lbar}},\,
 \Res_i(\DD)\bigr)
\simeq
 \bigl(\lefttop{\lbar}\nbigg_{\vecu},\,\,
 \nbign_i+\eigenmap(\lambda,u_i)\bigr).
\]
\item
Taking the multi-valued flat sections of
$\lefttop{\lbar}\nbigg_{\vecu}(\nbige)$,
we obtain the vector bundle $\lefttop{\lbar}\nbigg_{\vecu}(\nbigh)$
over $\cnum^{\ast}_{\lambda}$.
\end{itemize}
\end{lem}
\pf
It immediately follows from our construction.
\hfill\qed

\subsubsection{Pairing}

From the natural pairing
$\nbige\otimes\nbige^{\lor}\lrarr\nbigo_{\nbigx-\nbigd}$,
we obtain the following morphism
on $\nbigx(\lambda_0,\epsilon_0)$ for
$\lambda_0\in\cnum_{\lambda}^{\ast}$ and for any sufficiently small
$\epsilon_0>0$:
\[
 \prolongg{\vecb}{\bigl(\nbige\bigr)}
\otimes
 \prolongg{-\vecb+(1-\epsilon)\vecdelta}{\bigl(\nbigelor\bigr)}
\lrarr
\nbigo_{\nbigx}.
\]
Since it preserves the filtrations $\nbigf^{(\lambda_0)}$
and the decompositions $\EE^{(\lambda_0)}$,
we obtain the following morphism of regular $\lambda$-connections:
\[
 \nbigg_{\vecu}\nbige\otimes
 \nbigg_{-\vecu}\nbige^{\lor}
\lrarr
 \nbigo_{\nbigx}.
\]

\subsubsection{Functoriality}
\label{subsubsection;10.26.25}

We have the naturally defined morphism
$\lefttop{\lbar}\nbigg_{-\vecu}\bigl(\nbige^{\lor}\bigr)
 \lrarr
 \lefttop{\lbar}\nbigg_{\vecu}\bigl(\nbige\bigr)^{\lor}$.
We also have the morphism
$\lefttop{\lbar}\nbigg_{\vecu_1}(\nbige_1)\otimes
 \lefttop{\lbar}\nbigg_{\vecu_2}(\nbige_2)
\lrarr
 \lefttop{\lbar}\nbigg_{\vecu_1+\vecu_2}(\nbige_1\otimes\nbige_2)$.
\begin{lem}
The morphism
$\lefttop{\lbar}\nbigg_{-\vecu}\bigl(\nbige^{\lor}\bigr)
 \lrarr
 \lefttop{\lbar}\nbigg_{\vecu}\bigl(\nbige\bigr)^{\lor}$
is isomorphic.
The following morphism is isomorphic:
\[
\bigoplus_{\substack{
 \vecu_i\in\KMS(\prolongg{\vecb_1}{\nbige_i}),\\
 \vecu_1+\vecu_2=\vecu
 }}
 \lefttop{\lbar}\nbigg_{\vecu_1}(\nbige_1)\otimes
 \lefttop{\lbar}\nbigg_{\vecu_2}(\nbige_2)
\lrarr
 \lefttop{\lbar}\nbigg_{\vecu_1+\vecu_2}(\nbige_1\otimes\nbige_2).
\]
\end{lem}
\pf
It follows from Lemma \ref{lem;10.24.10}
and Lemma \ref{lem;10.24.11}.
\hfill\qed

\vspace{.1in}
We also have the functoriality for the pull backs.
We use the setting in the subsubsection \ref{subsubsection;10.26.6}.
\begin{lem}\label{lem;10.26.11}
We have the naturally defined isomorphism as follows:
\[
 \bigoplus_{\vecc\cdot\vecu=\vecu_1}
 \lefttop{\lbar}\nbigg_{\vecu}\bigl(\nbige\bigr)
 \simeq
 \lefttop{\lbar}\nbigg_{\vecu_1}\bigl(\psi^{\ast}\nbige\bigr).
\]
\end{lem}
\pf
It follows from Lemma \ref{lem;10.26.10}.
\hfill\qed

%% file: 17.2.tex

\subsubsection{The induced morphism $\Phi_{\vecu}^{\can}$}

\label{subsubsection;10.18.36}

From the regular $\lambda$-connection
$\bigl(\lefttop{\lbar}\nbigg_{\vecu}(\nbige),\DD\bigr)$,
we obtain the isomorphism
$\Phi_{\vecu}^{\can}:
 \lefttop{\lbar}\nbigg_{\vecu}\nbigh\lrarr
 \lefttop{\lbar}\nbigg_{\vecu\,|\,\cnum_{\lambda}^{\ast}}$,
which we will explain.
For any holomorphic section $s$ of
$\lefttop{\lbar}\nbigg_{\vecu}(\nbigh)$
over an open subset $U\subset\cnum^{\ast}$,
we put as follows:
\[
 v=\exp\Bigl(\sum_j
 \log z_j\cdot
 \bigl(
  \lambda^{-1}\cdot\eigenmap(\lambda,u_j)
 +\nbign_{j\,\vecu} \bigr)
 \Bigr)\cdot s.
\]
Then it gives the holomorphic section of
the vector bundle
$\lefttop{\lbar}\nbigg_{\vecu}\bigl(\nbige\bigr)$
over $U\times X$.
Then the restriction $v_{|U\times\{O\}}$ is a section
of $\lefttop{\lbar}\nbigg_{\vecu}$ over $U$.
We put
$\Phi_{\vecu}^{\can}(s):=v_{|U\times\{O\}}$,
and then
we obtain the isomorphism $\Phi_{\vecu}^{\can}$ desired.

The morphism $\Phi_{\vecu}^{\can}$ be also seen as follows:
Let $s$ be a section of $\lefttop{\lbar}\nbigg_{\vecu}\nbigh$.
We have the expression as follows:
\[
 s=\exp\Bigl(
 -\sum_j
\log z_j\cdot
 \bigl(
  \lambda^{-1}\cdot\eigenmap(\lambda,u_j)
 +\nbign_{j,\vecu}
 \bigr)
 \Bigr)\cdot v
=
 \sum_J(\log z)^J\cdot v_J.
\]
Here $J=(j_1,\ldots,j_l)$ denotes multi-indices
and $\bigl(\log z\bigr)^J=\prod_{h=1}^l \bigl(\log z_h\bigr)^{j_h}$.
\begin{lem}
We have $\Phi^{\can}_{\vecu}(s)=v_{0}(O)$.
\end{lem}
\pf
It immediately follows from our construction.
\hfill\qed

\begin{lem} \label{lem;10.26.13}
The isomorphism
$\Phi_{\vecu}^{\can}$ preserves the morphisms $\nbign_{\vecu\,i}$
and the pairing.
It is compatible with tensor products and duals.
It is also compatible with pull backs
as in the subsubsection {\rm\ref{subsubsection;10.26.6}}.
\end{lem}
\pf
It is clear from our construction.
\hfill\qed

\subsubsection{The induced morphisms $\Phi_{\vecu,P,O}$}

\label{subsubsection;10.18.35}

For any point $P\in X$,
we have the isomorphisms of
$\lefttop{\lbar}\nbigg_{\vecu}$
and $\lefttop{\lbar}\nbigg_{\vecu}(\nbige)_{|\cnum^{\ast}\times\{P\}}$.
For simplicity of notation, we denote
$\lefttop{\lbar}\nbigg_{\vecu}(\nbige)_{|\cnum^{\ast}\times\{P\}}$
by $\lefttop{\lbar}\nbigg_{\vecu}(\nbige)_{|P}$.

Take a normalizing frame $\vecv$ of
$\lefttop{\lbar}\nbigg_{\vecu}(\nbige)$,
namely we take a holomorphic frame $\vecv$ of
$\lefttop{\lbar}\nbigg_{\vecu}(\nbige)$
such that $\DD\vecv=\vecv\cdot\sum A_i\frac{dz_i}{z_i}$
holds for some constant matrices $A_i$.
For any point $P,Q\in X$,
the trivialization $\vecv$
gives the isomorphism
$\Phi_{\vecu,P,Q}:
 \lefttop{\lbar}\nbigg_{\vecu}(\nbige)_{|P}
\lrarr
 \lefttop{\lbar}\nbigg_{\vecu}(\nbige)_{|Q}$,
by the correspondence
$v_{i\,|\,P}\longmapsto v_{i\,|\,Q}$.
If we fix the coordinate,
the $\Phi_{\vecu,P,Q}$ does not depend on a choice of
normalizing frame.
Note that we have the isomorphism:
\[
 \lefttop{\lbar}\nbigg_{\vecu}(\nbige)_{|O}
\simeq
 \lefttop{\lbar}\nbigg_{\vecu\,|\,\cnum^{\ast}}.
\]
Thus we obtain the isomorphism $\Phi_{\vecu,P,O}$
of 
$\lefttop{\lbar}\nbigg_{\vecu}(\nbige)_{|P}$
and $\lefttop{\lbar}\nbigg_{\vecu}$.

\begin{lem} \label{lem;10.26.14}
The isomorphism
$\Phi_{\vecu,P,O}$ preserves the morphisms $\nbign_{\vecu\,i}$
and the pairing.
It is compatible with tensor product and dual.
It is also compatible with pull backs
as in the subsubsection {\rm\ref{subsubsection;10.26.6}}.
\hfill\qed
\end{lem}

\begin{rem}
In our previous paper {\rm\cite{mochi}},
we used only the morphism $\Phi_{\vecu,P,O}$
and did not use the morphism $\Phi^{\can}_{\vecu}$.
\hfill\qed
\end{rem}

%% file: 21.tex

\subsubsection{Conjugate}

Let $X$ be a complex manifold.
We denote the conjugate of $X$ by $X^{\dagger}$.
We put $\nbigx^{\dagger}=\cnum_{\mu}\times X^{\dagger}$.
Let $\harmonicbundle$ be a tame harmonic bundle over $X$.
Then we obtain the tame harmonic bundle $(E,\del_E,h,\theta^{\dagger})$
over $X^{\dagger}$,
and thus 
the deformed holomorphic bundle
$\nbige^{\dagger}$ over $\nbigx^{\dagger}-\nbigd^{\dagger}$
and the $\mu$-connection $\DD^{\dagger}$
on $\nbigx^{\dagger}$.
We also have the associated flat connections $\DD^{\dagger\,f}$.

Let $\sigma:\cnum_{\mu}\lrarr\cnum_{\lambda}$ be the morphism
given by $\mu\longmapsto -\overline{\mu}$.
It induces the anti-holomorphic map
$\nbigx^{\dagger}\lrarr\nbigx$.

Let $U$ be an open subset of $\nbigx^{\dagger}$.
We have the isomorphism $\sigma:U\lrarr\sigma(U)$.
Let $\vecv$ be a frame of $\prolongg{\vecb}{\nbige}$
over $\sigma(U)$.
Then we put as follows:
\begin{equation} \label{eq;10.16.1}
\vecv^{\dagger}:=
 \sigma^{\ast}\Bigl(
 \vecv\cdot \overline{H(h,\vecv)}^{-1}
 \Bigr).
\end{equation}
Then $\vecv^{\dagger}$ is a tuple of $C^{\infty}$-sections
of $\nbige^{\dagger}$ on $U$.
\begin{lem}\mbox{{}} \label{lem;10.16.7}
\begin{enumerate}
\item
The tuple $\vecv^{\dagger}$ is a holomorphic frame
of $\nbige^{\dagger}$ on $U$.
\item \label{number;10.16.8}
Let $\nbiga$ be the $\lambda$-connection one form of $\DD$
with respect to the frame $\vecv$.
Then the $\mu$-connection one form of $\DD^{\dagger}$
with respect to the frame $\vecv^{\dagger}$
is given by $\sigma^{\ast}\Bigl({}^t\overline{\nbiga}\Bigr)$.
\end{enumerate}
\end{lem}
\pf
It can be checked by direct calculations.
(See the subsubsection 3.1.6 in our previous paper \cite{mochi},
for example).
\hfill\qed

\subsubsection{The comparison of the flat connections}

\label{subsubsection;a12.9.5}

We identify $\cnum^{\ast}_{\mu}$ and $\cnum^{\ast}_{\lambda}$
by the relation $\mu=\lambda^{-1}$.
It induces the identification of the $C^{\infty}$-manifolds
$\nbigx^{\dagger\,\shikaku}=\nbigx^{\shikaku}$.

We have the holomorphic family of the flat connections
$(\nbige^{\shikaku},\DD^{f})$
on $\nbigx^{\shikaku}-\nbigd^{\shikaku}$.
We also have the holomorphic family of the flat connections
$(\nbige^{\dagger\,\shikaku},\DD^{\dagger\,f})$
on $\nbigx^{\shikaku\,\dagger}-\nbigd^{\shikaku\,\dagger}$.
\begin{lem} \label{lem;10.18.1}
Under
the identification $\nbigx^{\dagger\,\shikaku}=\nbigx^{\shikaku}$
given above,
we have
$\bigl(\nbige^{\dagger\,\shikaku},\DD^{\dagger\,f}\bigr)
=\bigl(\nbige^{\shikaku},\DD^{f}\bigr)$
over $\nbigx^{\shikaku}$.
\end{lem}
\pf
By definition of $\DD$ and $\DD^{\dagger}$,
we obtain the following:
\[
 \DD^f=
 \delbar_E+\lambda\theta^{\dagger}+
 \lambda^{-1}\cdot\Bigl(
 \lambda\del_E+\theta
 \Bigr)
=\del_E+\mu\theta+
\mu^{-1}\Bigl(
 \mu\cdot\delbar_E+\theta^{\dagger}
 \Bigr)
=\DD^{\dagger\,f}.
\]
Thus we are done.
\hfill\qed

\vspace{.1in}

Let $M_i$ denote the monodromy endomorphism of $\nbige$
with respect to the loop $\gamma_i$
around the divisor $D_i$ with the anti-clockwise direction:
\begin{equation}
 \gamma_i:[0,1]\lrarr
 \bigl(z_1,\ldots,z_{i-1},e^{2\pi\sqrt{-1}t}\cdot z_i,z_{i+1},
 \ldots,z_n
 \bigr).
\end{equation}
Let $M_i^{\dagger}$ denote the monodromy endomorphism of $\nbige^{\dagger}$
of $\gamma_i^{-1}$.
The following lemma immediately follows from Lemma \ref{lem;10.18.1}.
\begin{lem} \label{lem;10.18.28}
We have $M_i^{-1}=M_i^{\dagger}$.
\hfill\qed
\end{lem}

\subsubsection{The variation of pure twistor structures and the conjugate}

Due to Lemma \ref{lem;10.18.1},
we obtain the patched object
as in the subsubsection \ref{subsubsection;10.18.2}.
Thus we obtain the variation of pure twistor structures.
It can be simply described as follows:
Let $p:X\times\proj^1\lrarr X$ denote the projection.
We put $\nbige^{\sankaku}:=p^{-1}(E)$.
The differential operator
$\DD^{\sankaku}:
 C^{\infty}\bigl(X\times\proj^1,\nbige^{\sankaku}\bigr)
\lrarr
 C^{\infty}\bigl(X\times\proj^1,\nbige^{\sankaku}\otimes\xi\Omega_X^1\bigr)$
is given as follows:
\[
 \DD^{\sankaku}:=
 \bigl(\delbar_E+\theta\bigr)\otimes\sqrt{-1}\cdot f_0^{(1)}
+\bigl(\del_E+\theta^{\dagger}\bigr)\otimes f_{\infty}^{(1)}.
\]
\begin{lem}
$\bigl(\nbige^{\sankaku},\DD^{\sankaku}\bigr)$
is a variation of pure twistor structures.
\end{lem}
\pf
It can be checked by a direct calculation
(Lemma \ref{lem;a12.1.1}).
\hfill\qed

\vspace{.1in}

We obtain the conjugate
$\bigl(\sigma^{\ast}\nbige^{\sankaku},
 \DD^{\sankaku}_{\sigma^{\ast}\nbige^{\sankaku}}\bigr)$
(see the subsubsection \ref{subsubsection;10.18.3}).
In this case, we have $\sigma^{\ast}\nbige^{\sankaku}=p^{-1}(E)$.
\begin{lem}
The $\proj^1$-holomorphic structure
$d''_{\lambda}$ is given as follows:
\[
 d''_{\lambda}\sigma^{\ast}g
=\sigma^{\ast}\Bigl(
 \frac{\del g}{\del \bar{\lambda}}
 \Bigr)\cdot(-d\bar{\lambda}).
\]
\end{lem}
\pf
By definition,
we have the following:
\[
 d''_{\lambda}\sigma^{\ast}g
=\varphi_0\sigma^{\ast}\bigl(\delbar_{\lambda}g\bigr)
=\sigma^{\ast}\Bigl(
 \frac{\del g}{\del \bar{\lambda}}
 \Bigr)\cdot(-d\bar{\lambda}).
\]
Thus we are done.
\hfill\qed

\begin{lem}\label{lem;10.18.8}
Let $g$ be a section of $\nbige$.
The $C^{\infty}$-sections $A_i$, $B_i$, $C_i$ and $D_i$
are determined as follows:
\[
 \del_E g=\sum A_i\cdot dz_i,
\quad
 \theta^{\dagger}\cdot g=\sum B_i\cdot d\overline{z}_i,
\quad
 \delbar_E g=\sum C_i\cdot d\overline{z}_i,
\quad
 \theta \cdot g=\sum D_i\cdot dz_i.
\]
Then we have the following formula:
\begin{equation} \label{eq;a12.1.2}
 \DD(\sigma^{\ast}g)=
 \sum_i \Bigl(
 \sigma^{\ast}A_i\cdot d\bar{z}_i\otimes\sqrt{-1}f_0^{(1)}
-\sigma^{\ast}B_i\cdot dz_i\otimes \sqrt{-1}f_0^{(1)}
+\sigma^{\ast}C_i\cdot dz_i\otimes f_{\infty}^{(1)}
-\sigma^{\ast}D_i\cdot d\bar{z}_i\otimes f_{\infty}^{(1)}
 \Bigr).
\end{equation}
\end{lem}
\pf
We have
$\DD(\sigma^{\ast}g)=\varphi_0\sigma^{\ast}(\DD g)$
by definition.
We can check the formula (\ref{eq;a12.1.2})
by using Lemma \ref{lem;10.18.6}.
\hfill\qed

\subsubsection{Polarization}

For any sections $f$ and $\sigma^{\ast}(g)$
of $\nbige^{\sankaku}$ and $\sigma^{\ast}\nbige^{\sankaku}$,
we have the $C^{\infty}$-function
$S\bigl(f,\sigma^{\ast}(g)\bigr):=
 h\bigl(f(\lambda,x),g(-\bar{\lambda},x)\bigr)$
Thus we obtain the pairing
$S:\nbige^{\sankaku}\otimes\sigma^{\ast}\nbige^{\sankaku}
\lrarr \Tate(0)$.

\begin{lem}
The pairing $S$ is a morphism of $\proj^1$-holomorphic bundles.
\end{lem}
\pf
We have the following equality:
\[
  \delbar_{\lambda} S(f,\sigma^{\ast}g)
=h\bigl(
 \delbar_{\lambda}f(\lambda,x),g(-\bar{\lambda},x)
 \bigr)
+h\bigl(f(\lambda,x),
 \del_{\lambda}\bigl(g(-\bar{\lambda},x)\bigr)
 \bigr).
\]
The first term in the right hand side can be rewritten
as $S\bigl(\delbar_{\lambda}f,\sigma^{\ast}g\bigr)$.
The second term in the right hand side can be rewritten as follows:
\[
 h\Bigl(f(\lambda,x), 
 \frac{\del g}{\del \bar{\lambda}}(-\bar{\lambda},x)
\Bigr)\cdot (-d\bar{\lambda})
=S\Bigl(
 f,\sigma^{\ast}\frac{\del g}{\del \bar{\lambda}}
 \Bigr)(-d\bar{\lambda})
=S\Bigl(
 f,d''\sigma^{\ast}g
 \Bigr).
\]
Thus we are done.
\hfill\qed

\begin{lem} \label{lem;10.18.7}
The pairing $S$ is a morphism of variation of pure twistors.
\end{lem}
\pf
We have the following equalities:
\begin{multline}
 \delbar_X S\bigl(f,\sigma^{\ast}g\bigr)\otimes\sqrt{-1}\cdot f_0^{(1)}
=h\bigl(\delbar_Ef(\lambda,x),
 g(-\bar{\lambda},x)
 \bigr)\otimes\sqrt{-1}\cdot f_0^{(1)}
+h\bigl(f(\lambda,x),
 \del_E g(-\bar{\lambda},x)
 \bigr)\otimes\sqrt{-1}\cdot f_0^{(1)} \\
=h\bigl((\delbar_E+\theta) f(\lambda,x),
 g(-\bar{\lambda},x)
 \bigr)
 \otimes\sqrt{-1}\cdot f_0^{(1)}
+h\bigl(f(\lambda,x),
 (\del_E-\theta^{\dagger})g(-\bar{\lambda},x)
 \bigr)
 \otimes\sqrt{-1}\cdot f_0^{(1)}\\
=S\bigl((\delbar_E+\theta)f,\sigma^{\ast}g\bigr)
+
 \sum_i
 S(f,\sigma^{\ast}A_i)\cdot d\bar{z}_i\otimes\sqrt{-1}f_0^{(1)}
-\sum_i S(f,\sigma^{\ast}B_i)\cdot dz_i\otimes\sqrt{-1}f_0^{(1)}.
\end{multline}
Here we have used Lemma \ref{lem;10.18.8}.
On the other hand,
we also have the following:
\begin{multline}
 \del_X S\bigl(f,\sigma^{\ast}g\bigr)\otimes f_{\infty}^{(1)}
=h\bigl(\del_Ef,g(-\bar{\lambda},x)\bigr)\otimes f_{\infty}^{(1)}
+h\bigl(f,\delbar_E g(-\bar{\lambda},x)\bigr)\otimes f_{\infty}^{(1)}\\
=h\bigl((\del_E+\theta^{\dagger})f,
 g(-\bar{\lambda},x)
 \bigr)
 \otimes f_{\infty}^{(1)}
+h\bigl(f,(\delbar_E-\theta)g(-\bar{\lambda},x)
 \bigr)\otimes f_{\infty}^{(1)}\\
=S\bigl((\del_E+\theta^{\dagger})f,\sigma^{\ast}g\bigr)
+\sum_i S\bigl(f,\sigma^{\ast}C_i\bigr)
 \cdot d\bar{z}_i\otimes f_{\infty}^{(1)}
-\sum_i S\bigl(f,\sigma^{\ast}D_j\bigr)\cdot
 dz_j\otimes f_{\infty}^{(1)}.
\end{multline}
Then Lemma \ref{lem;10.18.7} immediately follows.
\hfill\qed

\begin{cor}[Simpson]
The tuple $(\nbige^{\sankaku},\DD^{\sankaku},h)$ is
a variation of polarized pure twistor structures.
\hfill\qed
\end{cor}

\begin{lem}
We obtain the isomorphism
$\dualhenomap:\sigma^{\ast}\nbige^{\sankaku}\simeq\nbige^{\lor\,\sankaku}$
of the variation of pure twistors.
In particular,
we obtain the isomorphisms
$\sigma^{\ast}\nbige^{\dagger}\simeq \nbige^{\lor}$
and $\sigma^{\ast}\nbige\simeq \nbige^{\lor\,\dagger}$.
We also denote them by $\dualhenomap$.
\end{lem}
\pf
Since the pairing $S$ is perfect,
it induces the isomorphism $\dualhenomap$.
\hfill\qed

%% file: a53.tex

\subsubsection{Compatible frame and the KMS-structure of the conjugate}

We put $X=\Delta^n$, $D_i:=\{z_i=0\}$ and $D=\bigcup_{i=1}^l D_i$.
Let us pick a point $\lambda_0\in\cnum_{\lambda}$.
Let us pick a sufficiently small positive number $\epsilon_0$
such that the sheaf $\prolongg{\vecb}{\nbige}$
on $\nbigx(\lambda_0,\epsilon_0)$
is locally free.
Let $\vecv=(v_i)$ be a frame of $\prolongg{\vecb}{\nbige}$,
which is compatible with $\EEzero$ and $\Fzero$.
For each $v_i$, 
we have the element $\vecu(v_i)\in\KMS(\nbige^0,\lbar)$
such that the following holds:
\[
 \kmsmap(\lambda_0,\vecu(v_i))=\deg^{\EEzero,\Fzero}(v_i)
\in \KMS(\nbige^{\lambda_0},\lbar).
\]
Let $u_j(v_i)\in \KMS(\nbige^0,j)$ denote the $j$-th component
of $\vecu(v_i)$.

We denote the restriction
$\vecv_{|\nbigx(\lambda_0,\epsilon_0)-\nbigd(\lambda_0,\epsilon_0)}$
by $\vecv$, for simplicity.
Then we obtain the holomorphic frame $\vecv^{\dagger}$
of $\nbige^{\dagger}$ over
$\nbigx^{\dagger}(-\overline{\lambda}_0,\epsilon_0)$,
which is given by (\ref{eq;10.16.1}).
(Note
$\sigma(\nbigx^{\dagger\,}(-\bar{\lambda}_0,\epsilon_0))
 =\nbigx(\lambda_0,\epsilon_0)$.)
We put as follows:
\begin{equation} \label{eq;10.25.25}
\begin{array}{ll}
 v_i':=v_i\cdot\prod_{j=1}^l|z_j|^{\paramap(\lambda,u_j(v_i))},
 &
 \vecv':=\bigl(v_i'\bigr),\\
 \mbox{{}}\\
 v_i^{\dagger\,\prime}
:=v_i^{\dagger}\cdot\prod_{j=1}^l|z_j|^{-\paramap(\lambda,u_j(v_i))},
 &
 \vecv^{\dagger\,\prime}=\bigl(v_i^{\dagger\,\prime}\bigr).
\end{array}
\end{equation}
Then $\vecv'$ is a $C^{\infty}$-frame of
$\nbige$ over
$\nbigx(\lambda_0,\epsilon_0)-\nbigd(\lambda_0,\epsilon_0)$,
and $\vecv^{\dagger\,\prime}$ is a $C^{\infty}$-frame of
$\nbige^{\dagger}$ over
$\nbigx^{\dagger}(-\overline{\lambda}_0,\epsilon_0)
 -\nbigd^{\dagger}(-\overline{\lambda}_0,\epsilon_0)$.

\begin{lem}
The frames $\vecv'$ and
$\vecv^{\dagger\,\prime}$ are adapted up to log order.
\end{lem}
\pf
The adaptedness of $\vecv'$ up to log order
has already been shown in Proposition \ref{prop;10.16.2}
(the subsubsection \ref{subsubsection;a12.1.5}).
Let $L$ be the diagonal matrix
such that
$L_{i\,i}:=\prod_j |z_j|^{\paramap(\lambda,u_j(v_i)}$.
Then we have the relations
$\vecv^{\dagger\,\prime}=\vecv^{\dagger}\cdot L^{-1}$
and
$\vecv'=\vecv\cdot L$.
Then we obtain the following:
\[
 \vecv^{\dagger\,\prime}
=\vecv^{\dagger}\cdot L^{-1}
=\vecv\cdot\overline{H(h,\vecv)}^{-1}\cdot L^{-1}
=\vecv'\cdot L^{-1}\cdot \overline{H(h,\vecv)}^{-1}\cdot L^{-1}
=\vecv'\cdot \overline{(L\cdot H(h,\vecv)\cdot L)}^{-1}
=\vecv'\cdot \overline{H(h,\vecv')}^{-1}.
\]
Since $\vecv'$ is adapted up to log order,
and since $H(h,\vecv')$ and $H(h,\vecv')^{-1}$
is bounded up to log order,
$\vecv^{\dagger\,\prime}$ is adapted up to log order.
\hfill\qed

\vspace{.1in}
Recall that
we put $\vecu^{\dagger}:=(\overline{\vecalpha},-\vecb)$,
for any element $\vecu=(\vecalpha,\vecb)\in\cnum^l\times\real^l$
(the subsubsection \ref{subsubsection;10.16.5}).

\begin{cor} \mbox{{}} \label{cor;10.16.9}
\begin{itemize}
\item
There exists a positive number $\epsilon>0$
such that $\vecv^{\dagger}$ is a frame of
$\prolongg{-\vecb+(1-\epsilon)\cdot\vecdelta}
   {\nbige^{\dagger}}$
over $\nbigx^{\dagger}(-\overline{\lambda}_0,\epsilon_0)$.
\item
The frame $\vecv^{\dagger}$ is compatible with the parabolic filtration
$F^{(-\overline{\lambda}_0)}$
and the decomposition $\EE^{(-\overline{\lambda}_0)}$.
\item
 We have the following:
\[
 \deg^{\EE^{(-\overline{\lambda}_0)},F^{(-\overline{\lambda}_0)}}
 (v_i^{\dagger})
=\kmsmap(-\overline{\lambda_0},\vecu(v_i)^{\dagger}).
\]
\end{itemize}
\end{cor}
\pf
By using Lemma \ref{lem;9.5.10}, we obtain the first claim.
By using Lemma \ref{lem;10.10.31},
we obtain that the frame $\vecv^{\dagger}$ is compatible with
the filtration $F^{(-\overline{\lambda}_0)}$,
and we have the following:
\[
 \deg^{F^{(-\overline{\lambda}_0)}}(v_i^{\dagger})
=-\paramap(\lambda_0,\vecu(v_i))
=\paramap\bigl(-\overline{\lambda}_0,\vecu(v_i)^{\dagger}\bigr).
\]
Here we have used Lemma \ref{lem;10.16.6}.

Due to the claim \ref{number;10.16.8} in Lemma \ref{lem;10.16.7},
we obtain that the frame $\vecv^{\dagger}$
is compatible with $\EE^{(-\overline{\lambda}_0)}$.
Moreover we have the following:
\[
 \deg^{\EE^{(-\overline{\lambda}_0)}}(v_i)
=\overline{\eigenmap\bigl(\lambda_0,\vecu(v_i)\bigr)}
=\eigenmap\bigl(-\overline{\lambda}_0,\vecu(v_i)^{\dagger}\bigr).
\]
Thus we obtain the second and the third claims in
Corollary \ref{cor;10.16.9}.
\hfill\qed

\begin{cor}
By the correspondence $\vecu\longmapsto \vecu^{\dagger}$,
we have the isomorphism preserving the multiplicity:
\[
 \KMS(\nbigelambda,\lbar)
\lrarr \KMS(\nbige^{\dagger\,-\bar{\lambda}},\lbar).
\]
In particular,
we have the isomorphism
$\KMS(\nbige^0,\lbar)
\lrarr \KMS(\nbige^{\dagger\,0},\lbar)$.
\hfill\qed
\end{cor}

\begin{lem} \label{lem;10.18.11}
Via the isomorphism
$\dualhenomap:\sigma^{\ast}\nbige^{\dagger}\simeq\nbige^{\lor}$,
we have $\dualhenomap(\vecv^{\dagger})=\vecv^{\lor}$.
\end{lem}
\pf
It can be shown by an elementary linear algebraic argument.
(See the subsubsection 3.1.6 in the previous paper, for example.)
\hfill\qed

\subsubsection{The conjugate and the dual}

Let $\lambda_0$ be a point of $\cnum_{\lambda}$.
Let $\vecb=(b_1,\ldots,b_l)$ be an element of $\real^l$ such that
$b_i\not\in \KMS(\nbige^{\lor\,\lambda_0},i)$ for any $i$.
We take a sufficiently small positive number $\epsilon_0$,
then we have the locally free sheaf
$\prolongg{\vecb}{\nbige^{\lor}}$ on $\nbigx(\lambda_0,\epsilon_0)$.

\begin{cor} \label{cor;10.18.12}
The sheaf $\prolongg{\vecb}{\nbige^{\dagger}}$
 on $\nbigx^{\dagger}(-\overline{\lambda}_0,\epsilon_0)$
is locally free.
We have the isomorphism
$\dualhenomap:
 \sigma^{\ast}\bigl(\prolongg{\vecb}{\nbige^{\dagger}}\bigr)
\simeq
 \prolongg{\vecb}{\nbige^{\lor}}$.
\end{cor}
\pf
It follows from Lemma \ref{lem;10.18.11}.
\hfill\qed

\vspace{.1in}
The morphism $\dualhenomap$ induces the
isomorphism $\dualhenomap_{|\nbigd_i(\lambda_0,\epsilon_0)}$:
\[
 \dualhenomap_{|D_i}:
 \sigma^{\ast}
  \prolongg{\vecb}{\nbige^{\dagger}}_{|\nbigd_i(-\bar{\lambda}_0,\epsilon_0)}
\lrarr
  \prolongg{\vecb}{\nbige^{\lor}}_{|\nbigd_i(\lambda_0,\epsilon_0)}.
\]
Recall that we have the decomposition $\lefttop{i}\EEzero$
and the filtration $\lefttop{i}\Fzero$
of $\prolongg{\vecb}{\nbige^{\lor}}_{|\nbigd_i(-\bar{\lambda}_0,\epsilon_0)}$
given in the subsubsection \ref{subsubsection;10.18.20}
and the subsubsection \ref{subsubsection;10.18.21}.
Similarly,
we have the decomposition and the filtration of
$\sigma^{\ast}
  \prolongg{\vecb}{\nbige^{\dagger}}_{|\nbigd_i(\lambda_0,\epsilon_0)}$.

\begin{lem}\label{lem;10.18.13}
The morphism $\dualhenomap_{|D_i}$ preserves the filtrations
and the decompositions.
In particular,
the morphism $\dualhenomap_{|D_i}$ induces
the isomorphism
$\KMS(\prolongg{\vecb}{\nbige^{\dagger,-\overline{\lambda}}},i)
\lrarr
\KMS(\prolongg{\vecb}{\nbige^{\lor,\lambda}},i)$
given by the correspondence
$(\alpha,a)\longmapsto (-\overline{\alpha},a)$.
\end{lem}
\pf
It can be shown by using the comparison
of $\deg^{\EE,F}(v_i^{\dagger})$ and $\deg^{\EE,F}(v_i^{\lor})$.
\hfill\qed

\begin{rem}\label{rem;a12.1.6}
{\rm
When we have $\deg^{\EE^{(-\lambdabar_0)}}(v^{\dagger}_i)=\alpha$,
we obtain
$\deg^{\EE^{(\lambda_0)}}(\sigma^{\ast}(v_i^{\dagger}))
  =-\overline{\alpha}$,
due to the relation $\dualhenomap(v_i^{\dagger})=v_i^{\lor}$.
It can be directly seen, which we explain in the following.
For simplicity, we assume
that $\dim(X)=1$ and that we have the equality
$\DD^{\dagger}v_i^{\dagger}=
 \alpha\cdot v_i^{\dagger}\cdot d\zbar/\zbar$.
In that case, we have the following:
\[
 \DD^{\sankaku}(v_i^{\dagger})=
 \alpha\cdot v_i^{\dagger}\cdot
 d\zbar/\zbar\otimes f_{\infty}^{(1)}.
\]
Then we obtain the following:
\begin{multline}
 \DD^{\sankaku}\sigma^{\ast}(v_i^{\dagger})
=\varphi_0\sigma^{\ast}\bigl(
 \DD^{\sankaku}v_i^{\dagger}
 \bigr)
=\varphi_0\sigma^{\ast}\bigl(
 \alpha\cdot v_i^{\dagger}\cdot
 d\zbar/\zbar\otimes f_{\infty}^{(1)}
 \bigr)
=\overline{\alpha}
\cdot \sigma^{\ast}(v_i^{\dagger})
\cdot dz/z\otimes (-\sqrt{-1})\cdot f_{0}^{(1)} \\
=-\overline{\alpha}
\cdot \sigma^{\ast}(v_i^{\dagger})
\cdot dz/z\otimes\sqrt {-1}\cdot f_0^{(1)}.
\end{multline}
Since we have
$\DD^{\sankaku}(\sigma^{\ast}v_i^{\dagger})
=\DD(\sigma^{\ast}v_i^{\dagger})\otimes \sqrt{-1}f_0^{(1)}$,
we obtain 
$\DD(\sigma^{\ast}v_i^{\dagger})=
-\overline{\alpha}\cdot \sigma^{\ast}(v_i)\cdot dz/z$.
\hfill\qed
}
\end{rem}

Let us pick a point $\lambda_0\in\cnum^{\ast}_{\lambda}$.
Then we have the filtration and the decomposition
of $\prolongg{\vecb}{\nbige^{\lor}}$ on $\nbigx(\lambda_0,\epsilon_0)$,
given in the subsubsection \ref{subsubsection;10.18.15}.
Similarly, we have the filtration and the decomposition
of $\prolongg{\vecb}{\nbige^{\dagger}}$
on $\nbigx^{\dagger}(-\overline{\lambda}_0,\epsilon_0)$.

\begin{lem} \label{lem;10.18.18}
The filtrations and the decompositions are preserved
by the morphism $\dualhenomap$.
\end{lem}
\pf
Recall that the decompositions are induced by the monodromy
endomorphisms.
Since the $\dualhenomap$ preserves the flat connection,
the decompositions are preserved.

Recall that the restriction of the filtrations
to $\nbigx^{\ast}(\lambda_0,\epsilon_0)$
have the canonical splittings, given by the generalized eigenspaces
of the monodromy actions.
(Here we may assume that 
any point $\lambda$ of $\Delta^{\ast}(\lambda_0,\epsilon_0)$ is generic.)
Thus the restriction of the filtrations
to $\nbigx^{\ast}(\lambda_0,\epsilon_0)$ are preserved.
Then it follows that the filtrations are preserved
on whole $\nbigx(\lambda_0,\epsilon_0)$.
\hfill\qed

\vspace{.1in}
The isomorphism $\dualhenomap$
induces the isomorphism
$\sigma^{\ast}\nbigh^{\dagger}\simeq \nbigh^{\lor}$,
which we denote also by $\dualhenomap$.
Let us pick a point $\lambda_0\in\cnum_{\lambda}$
and a sufficiently small positive number $\epsilon_0$.
Then we have the filtrations and the decompositions
of $\nbigh^{\lor}$ on $\Delta_{\lambda}(\lambda_0,\epsilon_0)$,
given in the subsubsections
\ref{subsubsection;10.18.16}--\ref{subsubsection;10.18.17}.
Similarly, we have the filtrations and the decompositions
of $\nbigh^{\dagger}$ on $\Delta_{\mu}(-\overline{\lambda}_0,\epsilon_0)$.

\begin{lem} \label{lem;10.18.22}
The morphism $\dualhenomap$ preserves
the filtrations and the decompositions.
\end{lem}
\pf
It can be shown by an argument 
similar to Lemma \ref{lem;10.18.18}.
\hfill\qed

\subsubsection{The induced objects and the pairing}

Let 
$\vecu^{\dagger}$ be an element of $\KMS(\nbige^{\dagger\,0},\lbar)$.

\begin{itemize}
\item
By applying the constructions
in the subsubsection \ref{subsubsection;10.16.10},
we obtain the vector bundle
$\lefttop{\lbar}\nbigg^{\dagger}_{\vecu^{\dagger}}$
on $\nbigd_{\lbar}^{\dagger}$.
We also have the endomorphisms $\Res_i(\DD^{\dagger})$ 
and the nilpotent parts $\nbign^{\dagger}_i$ $(i\in \lbar)$.
\item
By applying the construction in the subsubsection 
\ref{subsubsection;10.18.9},
we obtain the holomorphic vector bundle
$\lefttop{\lbar}\nbigg^{\dagger}_{\vecu^{\dagger}}\nbigh^{\dagger}$
on $\cnum^{\ast}_{\mu}$.
We also have the monodromy endomorphisms,
and the nilpotent parts $\nbign^{\dagger}_i$.
\item
By applying the construction in the subsubsection
\ref{subsubsection;10.18.10},
we obtain the holomorphic vector bundle
$\lefttop{\lbar}\nbigg^{\dagger}_{\vecu^{\dagger}}(\nbige^{\dagger})$
on $\nbigx^{\dagger\,\shikaku}$.
We have the holomorphic family of the flat connections
$\DD^{\dagger\,f}$.
\end{itemize}

\begin{cor} \label{cor;10.18.23}
The morphism $\dualhenomap$ induces the following isomorphisms:
\[
 \sigma^{\ast}\lefttop{\lbar}\nbigg^{\dagger}_{\vecu^{\dagger}}(E)
\simeq
 \lefttop{\lbar}\nbigg_{-\vecu}(E^{\lor}),
\quad
 \sigma^{\ast}
 \lefttop{\lbar}\nbigg^{\dagger}_{\vecu^{\dagger}}\nbigh^{\dagger}(E)
\simeq
 \lefttop{\lbar}\nbigg_{-\vecu}\nbigh(E^{\lor}),
\quad
 \sigma^{\ast}
 \lefttop{\lbar}\nbigg^{\dagger}_{\vecu^{\dagger}}(\nbige^{\dagger})
\simeq
 \lefttop{\lbar}\nbigg_{-\vecu}(\nbige^{\lor}).
\]
In the first and the second isomorphisms,
the isomorphisms reverse the signature of the nilpotent maps.
In the third isomorphism,
the isomorphism preserves the family of the flat connections.
\end{cor}
\pf
It immediately follows from Lemma \ref{lem;10.18.13},
Lemma \ref{lem;10.18.18} and Lemma \ref{lem;10.18.22}.
(See also Remark \ref{rem;a12.1.6}.)
\hfill\qed

\begin{cor}
We have the naturally defined pairings:
\[
 \begin{array}{l}
 \lefttop{\lbar}\nbigg_{\vecu}
 \otimes
 \sigma^{\ast}\lefttop{\lbar}\nbigg^{\dagger}_{\vecu^{\dagger}}
\lrarr\nbigo_{\nbigd_{\lbar}},\\
 \mbox{{}}\\
 \lefttop{\lbar}\nbigg_{\vecu}\nbigh
 \otimes
 \sigma^{\ast}\lefttop{\lbar}\nbigg^{\dagger}_{\vecu^{\dagger}}
 \nbigh^{\dagger}
\lrarr\nbigo_{\cnum_{\lambda}^{\ast}},\\
 \mbox{{}}\\
 \lefttop{\lbar}\nbigg_{\vecu}(\nbige)
 \otimes
 \sigma^{\ast}\lefttop{\lbar}\nbigg^{\dagger}_{\vecu^{\dagger}}
 (\nbige^{\dagger})
\lrarr\nbigo_{\nbigx^{\shikaku}}.
 \end{array}
\]
\end{cor}
\pf
It immediately follows from Corollary \ref{cor;10.18.23}.
\hfill\qed

%% file: a53.1.tex

\subsubsection{The identification of the flat bundles with
 filtrations and the decompositions}

Let $M_i$ denote the monodromy endomorphism of $\nbige$
with respect to the loop $\gamma_i$
around the divisor $D_i$ with the anti-clockwise direction:
\begin{equation} \label{eq;10.18.27}
 \gamma_i:[0,1]\lrarr
 \bigl(z_1,\ldots,z_{i-1},e^{2\pi\sqrt{-1}t}\cdot z_i,z_{i+1},
 \ldots,z_n
 \bigr).
\end{equation}
Let $M_i^{\dagger}$ denote the monodromy endomorphism of $\nbige^{\dagger}$
of $\gamma_i^{-1}$.
Recall Lemma \ref{lem;10.18.28}.

Let us pick a point $\lambda_0\in\cnum^{\ast}_{\lambda}$
and a small neighbourhood $U\subset\cnum^{\ast}_{\lambda}$ of $\lambda_0$.
Then we have the filtration and the decomposition
of $\nbige$ on $U\times X$,
given in the subsubsection \ref{subsubsection;10.18.15}.
We have the point
$\lambda_0^{-1}\in \cnum_{\mu}^{\ast}$
and the neighbourhood $U'$ in $\cnum_{\mu}^{\ast}$,
which is same as $U$ by the identification
$\lambda=\mu^{-1}$.
Then we have the filtration and the decomposition
of $\nbige^{\dagger}$
on $X^{\dagger}\times \sigma(U)$ similarly.
As is noted in the subsubsection \ref{subsubsection;a12.9.5},
we have
$\bigl(\nbige^{\shikaku},\DD^f\bigr)
=\bigl(\nbige^{\dagger\,\shikaku},\DD^{\dagger\,f}\bigr)$
as flat bundles.

\begin{cor} \label{cor;10.18.29}
The identification $\nbige^{\dagger\,\shikaku}=\nbige^{\shikaku}$
on $X\times U$
preserves the filtrations and the decompositions.
\end{cor}
\pf
The decomposition is obtained from the generalized eigen decomposition
of the monodromy endomorphisms.
Thus the decompositions are preserved,
due to Lemma \ref{lem;10.18.28}.
We put $U^{\ast}=U-\{\lambda_0\}$.
We recall that the restriction of the filtrations
to $(X-D)\times U^{\ast}$
have the splittings given by the generalized eigen decompositions
as in (\ref{eq;a12.9.6}).
We also recall the relation
$\paramap^f(\lambda,u)=\paramap^f(\lambda^{-1},u^{\dagger})$
(Lemma \ref{lem;a12.1.10}).
Then we obtain that the restriction of the filtrations
to $(X-D)\times U$
are preserved due to Lemma \ref{lem;10.18.28}.
Then it follows that the filtrations are preserved.
\hfill\qed

\vspace{.1in}

By considering the spaces of the multi-valued flat sections,
we obtain the vector bundle 
$\nbigh(E,\delbar_E,\theta^{\dagger},h)$
on $\cnum^{\ast}_{\mu}$.
We denote it by $\nbigh^{\dagger}(E)$ for simplicity.
Namely, $\nbigh^{\dagger}_{|\mu}$ denotes
the space of the multi-valued flat sections of $\nbige^{\dagger\,\mu}$.
We have the monodromy endomorphisms $M^{\dagger}_i$
of $\nbigh^{\dagger}$
with respect to the loop $\gamma_i^{-1}$.
We also have the monodromy endomorphisms $M_i$
of $\nbigh$ with respect to the loop $\gamma_i$.

\begin{lem} \label{lem;10.18.26}
Under the identification $\cnum^{\ast}_{\mu}=\cnum^{\ast}_{\lambda}$ above,
we have $\nbigh^{\dagger}(E)=\nbigh(E)$.
We have $M_i^{-1}=M_i^{\dagger}$.
\end{lem}
\pf
It follows from the coincidence of the flat connections
$\DD^{\mu,f}=\DD^{\lambda,f}$.
\hfill\qed

\vspace{.1in}

Let us pick a point $\lambda_0\in\cnum_{\lambda}$
and an appropriate neighbourhood $U$ of $\lambda_0$.
Then we have the filtrations and the decompositions
of $\nbigh$ on $U$,
given in the subsubsections
\ref{subsubsection;10.18.16}--\ref{subsubsection;10.18.17}.
Similarly, we have the filtrations and the decompositions
of $\nbigh^{\dagger}$ on $U$.

\begin{cor} \label{cor;10.18.30}
The identification $\nbigh=\nbigh^{\dagger}$
preserves the filtrations and the decompositions.
\end{cor}
\pf
It can be shown by an argument similar to
the proof of Corollary \ref{cor;10.18.29}.
\hfill\qed

\subsubsection{The identification of the induced objects}

We have the holomorphic bundles
$\lefttop{\lbar}\nbigg_{\vecu}\nbigh$ on $\cnum_{\lambda}^{\ast}$.
We have the tuple of the monodromy endomorphisms
$\vecM=(M_1,\ldots,M_l)$.
Here $M_i$ denote the monodromy along the loop $\gamma_i$
given in (\ref{eq;10.18.27}).
We also have the holomorphic bundles
$\lefttop{\lbar}\nbigg^{\dagger}_{\vecu^{\dagger}}$
on $\cnum_{\mu}^{\ast}$.
We have the tuple of the monodromy endomorphisms
$\vecM^{\dagger}=(M_1^{\dagger},\ldots,M_l^{\dagger})$.
Here $M_i^{\dagger}$ denotes the monodromy along the loop $\gamma_i^{-1}$.

We identify $\cnum_{\lambda}^{\ast}$ and $\cnum_{\mu}^{\ast}$
by the relation $\lambda=\mu^{-1}$.

\begin{lem} \label{lem;10.18.31}
 We have the natural identification
$\lefttop{\lbar}\nbigg_{\vecu}\nbigh=
 \lefttop{\lbar}\nbigg^{\dagger}_{\vecu^{\dagger}}\nbigh$.
We also have $M^{-1}=M^{\dagger}$.
\end{lem}
\pf
Due to Lemma \ref{lem;10.18.26},
we have the natural identification
$\nbigh=\nbigh^{\dagger}$ over $\cnum_{\lambda}^{\ast}$,
on which we have $M^{-1}=M^{\dagger}$.
Due to Corollary \ref{cor;10.18.30} and
Lemma \ref{lem;a12.1.10},
we obtain the result.
\hfill\qed

\vspace{.1in}

We have the holomorphic vector bundle
$\lefttop{\lbar}\nbigg_{\vecu}(\nbige)$
on $\nbigx^{\shikaku}$.
We denote the restriction of $\lefttop{\lbar}\nbigg_{\vecu}(\nbige)$
to $\nbigx^{\shikaku}-\nbigd^{\shikaku}$
by the same notation.
Then we have the holomorphic family of the regular connections
$\DD^{f}$.

Similarly, we also have the $C^{\infty}$-bundle
$\lefttop{\lbar}\nbigg^{\dagger}_{\vecu^{\dagger}}(\nbige^{\dagger})$
with the holomorphic family of the flat connections $\DD^{\dagger\,f}$
on $\nbigx^{\dagger\,\shikaku}-\nbigd^{\dagger\,\shikaku}$.
By the relation $\lambda=\mu^{-1}$,
we have
$\nbigx^{\shikaku}-\nbigd^{\shikaku}
=\nbigx^{\dagger\,\shikaku}-\nbigd^{\dagger\,\shikaku}
=(X-D)\times \cnum_{\lambda}^{\ast}$.

\begin{lem}\label{lem;10.18.32}
We have the natural identification
$\lefttop{\lbar}\nbigg^{\dagger}_{\vecu^{\dagger}}\bigl(
 \nbige^{\dagger}\bigr)
=\lefttop{\lbar}\nbigg_{\vecu}\bigl(\nbige\bigr)$
and $\DD^{f}=\DD^{\dagger\,f}$
on $\nbigx^{\shikaku}-\nbigd^{\shikaku}$.
\end{lem}
\pf
Recall Lemma \ref{lem;10.18.1}
and Corollary \ref{cor;10.18.29}.
Then Lemma  \ref{lem;10.18.32}
can be shown by an argument similar to the proof of
Lemma \ref{lem;10.18.31}.
\hfill\qed

%% file: 20.5.tex

\subsubsection{The vector bundle $S_{\vecu}(E,P)$}

Let $\vecu\in\KMS(\nbige^0,\lbar)$.
Then we obtain the holomorphic vector bundle
$\lefttop{\lbar}\nbigg_{\vecu}$ on $\nbigd_{\lbar}$.
For simplicity,
we denote the restriction of
$\lefttop{\lbar}\nbigg_{\vecu}$
to $\cnum_{\lambda}\times\{O\}$
by the same notation.
Namely, we have the holomorphic bundle
$\lefttop{\lbar}\nbigg_{\vecu}$
on $\cnum_{\lambda}\simeq\cnum_{\lambda}\times\{O\}$.
Similarly,
we have the holomorphic bundle
 $\lefttop{\lbar}\nbigg^{\dagger}_{\vecu^{\dagger}}$ on $\cnum_{\mu}$.

Let $P$ be a point of $X-D$.
Then we have the isomorphism
$ \Phi_{P,O}:
 \lefttop{\lbar}\nbigg(\nbige)_{|\{P\}\times\cnum_{\lambda}^{\ast}}
 \simeq
 \lefttop{\lbar}\nbigg_{\vecu|\cnum^{\ast}_{\lambda}}$
over $\cnum^{\ast}_{\lambda}$,
given in the subsubsection \ref{subsubsection;10.18.35}.
Similarly, we have the isomorphism 
$\Phi^{\dagger}_{P,O}:
 \lefttop{\lbar}\nbigg_{\vecu^{\dagger}}(\nbige^{\dagger})
 _{|\{P\}\times\cnum_{\mu}^{\ast}}
\simeq
 \lefttop{\lbar}\nbigg_{\vecu^{\dagger}|\cnum_{\mu}^{\ast}}$
over $\cnum^{\ast}_{\mu}$.
From the morphisms $\Phi_{P,O}$ and $\Phi_{P,O}^{\dagger}$,
we obtain the isomorphism:
\[
 \Phi_{P,O}^{\dagger\,-1}\circ\Phi_{P,O}:
 \lefttop{\lbar}\nbigg_{\vecu\,|\,\cnum_{\lambda}^{\ast}}
\lrarr
 \lefttop{\lbar}\nbigg^{\dagger}_{\vecu^{\dagger}\,|\,\cnum_{\mu}^{\ast}}.
\]
Then we obtain the vector bundle,
which we denote by $S_{\vecu}(E,P)$
or simply by $S(P)$.

\subsubsection{The vector bundle $S^{\can}_{\vecu}(E)$}
\label{subsubsection;a12.5.5}

Similarly we have the following isomorphisms,
given in the subsubsection \ref{subsubsection;10.18.36}:
\[
 \Phi^{\can}:\lefttop{\lbar}\nbigg_{\vecu}\nbigh
  \lrarr\lefttop{\lbar}\nbigg_{\vecu\,|\,\cnum^{\ast}_{\lambda}},
\quad\quad
\Phi^{\dagger\,\can}:
 \lefttop{\lbar}\nbigg^{\dagger}_{\vecu^{\dagger}}(\nbigh^{\dagger})
  \lrarr\lefttop{\lbar}
 \nbigg^{\dagger}_{\vecu^{\dagger}\,|\,\cnum^{\ast}_{\mu}}.\\
\]
Since we have the canonical identification
$\lefttop{\lbar}\nbigg_{\vecu}\nbigh
=\lefttop{\lbar}\nbigg^{\dagger}_{\vecu^{\dagger}}\nbigh^{\dagger}$,
we obtain the isomorphism
$\Phi^{\dagger\,\can}\circ
 (\Phi^{\can})^{-1}:
 \lefttop{\lbar}\nbigg_{\vecu\,|\,\cnum^{\ast}_{\lambda}}
\lrarr
 \lefttop{\lbar}\nbigg^{\dagger}_{\vecu^{\dagger}\,|\,\cnum^{\ast}_{\mu}}$.
Thus we obtain the vector bundle,
which we denote by $S^{\can}_{\vecu}(E)$.

\subsubsection{Pairing}

\begin{lem}
We have the natural isomorphisms
$\sigma^{\ast} S^{\can}_{\vecu}(E)\simeq S^{\can}_{\vecu}(E)^{\lor}$
and
$\sigma^{\ast} S_{\vecu}(E,P)\simeq S_{\vecu}(E,P)^{\lor}$.
\end{lem}
\pf
It follows from Corollary \ref{cor;10.18.23}
and our construction.
\hfill\qed

\begin{cor}
Let $\vecu$ be an element of $\KMS(\nbige^0,\lbar)$.
We have the naturally induced pairings:
\[
 \begin{array}{l}
 S^{\can}_{\vecu}(E)\otimes\sigma^{\ast} S^{\can}_{\vecu}(E)
\lrarr \Tate(0),\\
\mbox{{}}\\
 S_{\vecu}(E,P)\otimes\sigma^{\ast}S_{\vecu}(E,P)
\lrarr \Tate(0).
 \end{array}
\]
They are perfect.
\hfill\qed
\end{cor}

\subsubsection{Nilpotent maps}

We have the nilpotent part of the residues
$\nbign_{\vecu\,i\,|\,\cnum_{\lambda}^{\ast}}$
on $\lefttop{\lbar}\nbigg_{\vecu\,|\,\cnum_{\lambda}^{\ast}}$.
We also have
$\nbign^{\dagger}_{\vecu^{\dagger}\,i\,|\,\cnum_{\lambda}^{\ast}}$
on $\lefttop{\lbar}\nbigg^{\dagger}_{\vecu^{\dagger}\,|\,\cnum_{\mu}^{\ast}}$.

\begin{lem}
Due to the isomorphisms
$\Phi^{\dagger\,\can\,-1}\circ\Phi^{\can}$
or
$\Phi^{\dagger\,-1}\circ\Phi_{P,O}$,
we have the following:
\[
 \lambda^{-1}\cdot\nbign_{\vecu\,i\,|\,\cnum_{\lambda}^{\ast}}
=-\mu^{-1}\cdot\nbign^{\dagger}_{\vecu^{\dagger}\,i\,|\,\cnum_{\mu}^{\ast}}.
\]
\end{lem}
\pf
We have the relation:
\[
 \exp\Bigl(
 2\pi\sqrt{-1}\lambda^{-1}\cdot\nbign_{\vecu\,i\,|\,\cnum_{\lambda}^{\ast}}
 \Bigr)
=\exp\Bigl(-2\pi\sqrt{-1}
 \mu^{-1}\cdot\nbign^{\dagger}_{\vecu^{\dagger}\,i\,|\,\cnum_{\mu}^{\ast}}
 \Bigr).
\]
Thus we are done.
\hfill\qed

\vspace{.1in}

Thus we obtain the following morphisms:
\[
 \nbign_i^{\sankaku}:
 S_{\vecu}(E,P)\lrarr S_{\vecu}(E,P)\otimes \Tate(-1),
\quad
 \nbign_i^{\sankaku}:
 S^{\can}_{\vecu}(E)\lrarr S^{\can}_{\vecu}(E)\otimes \Tate(-1).
\]
Here we put 
$\nbign^{\sankaku}_{i\,|\,\cnum_{\lambda}}
 :=\nbign_{\vecu\,i}\otimes t_0^{(-1)}$
and
$\nbign^{\sankaku}_{i\,|\,\cnum_{\mu}}
 :=\nbign^{\dagger}_{\vecu^{\dagger}\,i}\otimes t_{\infty}^{(-1)}$.
Note that we have the relation
$t_0^{(-1)}=-\lambda^2\cdot t_{\infty}^{(-1)}$.
Thus $\nbign^{\sankaku}$ is well defined.

\begin{lem}
The isomorphisms
$\sigma^{\ast}S^{\can}(E)\simeq S^{\can}_{-\vecu}(E^{\lor})$
and $\sigma^{\ast}S_{\vecu}(E,P)\simeq S_{-\vecu}(E^{\lor},P)$
preserves the nilpotent morphisms.
\end{lem}
\pf
It follows from Corollary \ref{cor;10.18.23}
and the relation $\varphi_0(\sigma^{\ast}(t_{\infty}^{(-1)}))=-t_0^{(-1)}$.
\hfill\qed

\begin{cor}
We have the relation
$S\bigl(\nbign^{\sankaku}_i\otimes \id\bigr)
+S\bigl(\id\otimes\sigma^{\ast}\nbign^{\sankaku}_i\bigr)=0$.
\hfill\qed
\end{cor}

%% file: a57.3.tex

\subsubsection{Functoriality}

The functoriality of the induced vector bundle
can be easily obtained from the functorialities
of $\lefttop{\lbar}\nbigg_{\vecu}(E)$,
$\lefttop{\lbar}\nbigg_{\vecu}(\nbigh(E))$,
$\lefttop{\lbar}\nbigg_{\vecu}(\nbige)$,
and the morphisms $\Phi^{\can}$ and $\Phi_{P,O}$.
The lemmas in this subsubsection follows from
the results in the subsubsections 
\ref{subsubsection;10.26.20}--\ref{subsubsection;10.26.6},
\ref{subsubsection;10.26.22}--\ref{subsubsection;10.26.24},
\ref{subsubsection;10.26.25}
and
\ref{subsubsection;10.18.36}--\ref{subsubsection;10.18.35}.

We have the naturally defined morphisms:
\begin{equation}\label{eq;10.24.15}
S^{\can}_{-\vecu}(E^{\lor})\lrarr S^{\can}_{\vecu}(E)^{\lor},
\quad
S_{-\vecu}(E^{\lor},P)\lrarr S_{\vecu}(E,P)^{\lor}.
\end{equation}
We have the naturally defined nilpotent maps
$\nbign_i^{\lor}$ on $S^{\can}_{\vecu}(E)^{\lor}$
and $S_{\vecu}(E,P)^{\lor}$.
\begin{lem}
The morphisms {\rm(\ref{eq;10.24.15})} are isomorphic.
They are compatible with the pairing.
The signature of the nilpotent map is reversed.
\hfill\qed
\end{lem}

Let $\vecb_i$ $(i=1,2)$ be elements of $\real^l$.
We have the naturally defined morphism:
\begin{equation} \label{eq;10.24.16}
\begin{array}{l}
{\displaystyle
\bigoplus_{\substack{
 \vecu_i\in\KMS(\prolongg{\vecb_i}{\nbige^0},\lbar),\\
 \vecu_1+\vecu_2=\vecu
 }}
 S^{\can}_{\vecu_1}(E_1)\otimes S^{\can}_{\vecu_2}(E_2)
\lrarr S^{\can}_{\vecu_1+\vecu_2}(E_1\otimes E_2), }
 \\ \mbox{{}}\\
 {\displaystyle
 \bigoplus_{\substack{
 \vecu_i\in\KMS(\prolongg{\vecb_i}{\nbige^0},\lbar),\\
 \vecu_1+\vecu_2=\vecu
 }}
 S_{\vecu_1}(E_1,P)\otimes S_{\vecu_2}(E_2,P)
\lrarr
 S_{\vecu_1+\vecu_2}(E_1\otimes E_2,P).
 }
\end{array}
\end{equation}

\begin{lem}
The morphisms {\rm(\ref{eq;10.24.16})} are isomorphic.
They are compatible with the pairings and the nilpotent maps.
\hfill\qed
\end{lem}

We also have the functoriality for the pull backs.
We use the setting in the subsubsection \ref{subsubsection;10.26.6}.
\begin{lem}\label{lem;10.26.50}
We have the naturally defined isomorphism:
\[
 \bigoplus_{\vecc\cdot\vecu=\vecu_1}
 S_{\vecu}^{\can}(E)
 \simeq
 S_{\vecu_1}^{\can}(\psi^{\ast}E),
\quad
  \bigoplus_{\vecc\cdot\vecu=\vecu_1}
 S_{\vecu}(E,\psi(P))
 \simeq
 S_{\vecu_1}(\psi^{\ast}E,P).
\]
They are compatible with the nilpotent maps and the pairings.
\hfill\qed
\end{lem}

%% file: a58.2.tex

\subsubsection{The construction}

Let us consider the case $X=\Delta$ and $D=\{O\}$.
In this case, we have one nilpotent map
$\nbign^{\sankaku}$ on $S^{\can}_u(E)$ and $S_u(E,P)$.
Due to Simpson,
we know that the conjugacy classes of $\nbign^{\sankaku}_{|\lambda}$
are independent of a choice of $\lambda\in\proj^1$
(Corollary \ref{cor;9.11.1}).
Thus the weight filtrations are the filtration
in the category of vector bundles.
Thus we obtain the associated graded bundle
$\Gr^W_hS^{\can}_u(E)$ and $\Gr^W_h S_u(E,P)$.

We have another construction of
$\Gr^W_hS^{\can}_u(E)$ and $\Gr^W_hS_u(E,P)$.
\begin{itemize}
\item
We have the vector bundles
$\Gr^W_h \nbigg_u(E)$ on $\cnum_{\lambda}$
and $\Gr^W_h\nbigg_{u^{\dagger}}^{\dagger}(E)$
on $\cnum_{\mu}$.
\item
We have the vector bundles
$\Gr^W_h\nbigg_u(\nbigh)$ on $\cnum_{\lambda}$
and $\Gr^W_h\nbigg^{\dagger}_{u^{\dagger}}(\nbigh^{\dagger})$
on $\cnum_{\mu}$.
\item
We have the vector bundles
$\Gr^W_h\nbigg_u(\nbige)$ on $\nbigx^{\shikaku}$
and $\Gr^W_h\nbigg^{\dagger}_{u^{\dagger}}(\nbige^{\dagger})$
on $\nbigx^{\shikaku\,\dagger}$.
We have the canonical isomorphisms:
\[
 \Gr^W_h\nbigg_u(\nbige)_{|\cnum_{\lambda}^{\ast}\times\{O\}}
\simeq
 \Gr^W_h\nbigg_u(E)_{|\cnum_{\lambda}^{\ast}},
\quad
 \Gr^W_h\nbigg^{\dagger}_{u^{\dagger}}(\nbige^{\dagger})
 _{|\cnum_{\mu}^{\ast}\times\{O\}}
\simeq
 \Gr^W_h\nbigg^{\dagger}_{u^{\dagger}}(E)_{|\cnum_{\mu}^{\ast}}.
\]
\end{itemize}

We have the family of the induced flat connections
$\DD^f$ and $\DD^{\dagger\,f}$
on $\Gr^W_h\nbigg_u(\nbige)$ 
and $\Gr^W_h\nbigg^{\dagger}_{u^{\dagger}}(\nbige^{\dagger})$
respectively.

\begin{lem}
The monodromy endomorphisms of 
$\Gr^W_h\nbigg_u(\nbige)$ 
and $\Gr^W_h\nbigg^{\dagger}_{u^{\dagger}}(\nbige^{\dagger})$
with respect to the flat connections
are of the form: $F(\lambda)\times$ identity.
\hfill\qed
\end{lem}

As in the cases of $S^{\can}_u(E)$,
we obtain the isomorphisms:
\[
 \Phi^{\can}:\Gr^W_h\nbigg_u(\nbigh)\lrarr
 \Gr^W_h\nbigg_u(E)_{|\cnum_{\lambda}^{\ast}},
\quad
 \Phi^{\can\,\dagger}:
 \Gr^W_h\nbigg^{\dagger}_{u^{\dagger}}(\nbigh^{\dagger})
\lrarr
 \Gr^W_h\nbigg^{\dagger}_{u^{\dagger}}(E)_{|\cnum_{\mu}^{\ast}}.
\]
Thus we obtain the gluing of
$\Gr^W_h\nbigg_u(E)$ and
$\Gr^W_h\nbigg^{\dagger}_{u^{\dagger}}(E)$
via $\Phi^{\can\,\dagger}\circ \Phi^{\can\,-1}$.
Thus we obtain the vector bundle,
which is naturally isomorphic to 
$\Gr^WS^{\can}_u(E)$.

Similarly we obtain the gluing
$\Gr^W_h\nbigg_u(E)$ and
$\Gr^W_h\nbigg^{\dagger}_{u^{\dagger}}(E)$
via $\Phi^{\dagger}_{O,P}\circ \Phi^{-1}_{O,P}$.
The resulted vector bundle is naturally isomorphic to
$\Gr_h^WS_u(E,P)$.

%% file: a57.4.tex

\subsubsection{The gluing matrices}
\label{subsubsection;10.25.30}

For simplicity we put as follows:
\[
  \Gr^W\nbigg(E):=
 \bigoplus_{k\in\seisuu}\bigoplus_{u\in \KMS(\prolong{\nbige^0})}
 \Gr^W_k\nbigg_u(E),
\quad
 \Gr^W\nbigg^{\dagger}(E):=
 \bigoplus_{k\in\seisuu}\bigoplus_{u\in \KMS(\prolong{\nbige^0})}
 \Gr^W_k\nbigg^{\dagger}_{u^{\dagger}}(E).
\]

Let $\vecw$ be a frame of
$\Gr^W\nbigg(E)$ compatible with the grading.
We denote the degree of $w_i$ by
$u(w_i)\in \KMS(\prolong{\nbige^0})$.
Let $\vecw^{\dagger}$ be a frame of
$\Gr^W\nbigg^{\dagger}(E)$.
We denote the degree of $w_i^{\dagger}$ by 
$u(w_i^{\dagger})\in\KMS(\prolong{\nbige^0})$.
\begin{rem}
Note $u(w_i^{\dagger})$ denotes an element of
$\KMS(\prolong{\nbige^0})$ not $\KMS(\nbige^{\dagger\,0})$.
\hfill\qed
\end{rem}

Let $\vecw$ be a frame of $\Gr^W\nbigg(E)$.
We have the dual frame $\vecw^{\dagger}$
of $\bigl(\Gr^W\nbigg(E)\bigr)^{\lor}$.
Via the isomorphism
$\sigma^{\ast}\Gr^W\nbigg^{\dagger}(E)
\simeq\bigl(
 \Gr^W\nbigg(E)\bigr)^{\lor}$,
we obtain $\vecw^{\dagger}$.
\begin{lem}
We have $u(w_i)=u(w_i^{\dagger})$.
\hfill\qed
\end{lem}

By gluing
$(\Phi^{\can\,\dagger})
 \circ\Phi^{\can\,-1}$,
we obtain the relation
$\vecw=\vecw^{\dagger}\cdot A^{\can}$
for some holomorphic function
$A^{\can}:\cnum_{\lambda}^{\ast}\lrarr GL(r)$.
By gluing $\Phi_{P,O}^{\dagger}\circ\Phi_{P,O}^{-1}$,
we obtain the relation
$\vecw=\vecw^{\dagger}\cdot A_{P,O}$.
We would like to give a method to calculate $A^{\can}$
and $A_{P,O}$.

We can take the normalizing frame
$\vecv$ of $\Gr^W\nbigg(\nbige)$
which is a lift of $\vecw$.
We can also take the normalizing frame $\vecv^{\dagger}$ of
$\Gr^W\nbigg^{\dagger}(\nbige^{\dagger})$,
which is a lift of $\vecw^{\dagger}$.
Since we have the identification
$\Gr^W\nbigg(\nbige)
=\Gr^W\nbigg^{\dagger}(\nbige^{\dagger})$
over $\nbigx^{\shikaku}-\nbigd^{\shikaku}$
as $C^{\infty}$-bundles,
we have the relation:
\[
 v_i^{\dagger}=\sum_j v_j\cdot J_{j\,i},
\quad
 \vecv=\vecv^{\dagger}\cdot J.
\]

\begin{lem}
We have $J_{j\,i}=0$ unless
$u(w_i^{\dagger})=u(w_j)$.
\end{lem}
\pf
It follows from the compatibility of
the $\vecw$ and $\vecw^{\dagger}$ with the grading.
\hfill\qed

\begin{lem} \label{lem;10.24.20}
In the case $u(w_i^{\dagger})=u(w_j)=u$,
there exist holomorphic functions
$K_{j\,i}$ on $\cnum_{\lambda}^{\ast}$ such that
the following holds:
\[
 J_{j\,i}=
 \exp\Bigl(
 -\eigenmap(\lambda,u)\cdot\log|z|^2
 \Bigr)\cdot K_{j\,i}.
\]
\end{lem}
\pf
It follows from a direct calculation.
\hfill\qed

\begin{cor} \label{cor;10.25.51}
There exists the $GL(r)$-valued holomorphic function
$K$ on $\cnum^{\ast}_{\lambda}$ such that 
$J(\lambda,z)= C(\lambda,z)\cdot K(\lambda)$.
Here $C(\lambda,z)$ is given as follows:
 \begin{equation}\label{eq;10.25.50}
 C(\lambda,z)=
 \bigoplus_{\substack{u\in\KMS(\nbige^0),\\ k\in\seisuu}}
  \exp\Bigl(
 -\eigenmap(\lambda,u)\cdot\log|z|^2
 \Bigr)
\cdot\id_{\Gr^{W}_k\nbigg_{u}}.
\end{equation}
Here $C(\lambda,z)$ is regarded as the endomorphism
of $\bigoplus \Gr^{W}_k\nbigg_{u}$
via the frame $\vecw$.
\end{cor}
\pf
It immediately follows from \ref{lem;10.24.20}.
\hfill\qed

\begin{lem}\label{lem;a12.9.10}\mbox{{}}
\begin{itemize}
\item
Via the gluing
$\Phi^{\can\,\dagger}
\circ
 \Phi^{\can\,-1}$,
we have the relation
$\vecw^{\dagger}=\vecw\cdot K$.
\item
Via the gluing
$\Phi_{O,P}^{\dagger}\circ
 \Phi_{O,P}^{-1}$,
we have the relation
$\vecw^{\dagger}
=\vecw\cdot C(\lambda,P)\cdot K(\lambda)$.
Here $C(\lambda,P)$ denote the matrix given as follows:
\begin{equation}\label{eq;10.25.45}
 C(\lambda,P)=\bigoplus_{\substack{u\in\KMS(\nbige^0),\\k\in\seisuu}}
  \exp\Bigl(
 -\eigenmap(\lambda,u)\cdot\log|z(P)|^2
 \Bigr)
\cdot\id_{\Gr^W_k\nbigg_{u}(E)}.
\end{equation}
Here $C(\lambda,P)$ is regarded as the endomorphism
of $\bigoplus \Gr^{W}_k\nbigg_u$ via the frame $\vecw$.
\end{itemize}
\end{lem}
\pf
It follows from the definitions and Lemma \ref{lem;10.24.20}.
\hfill\qed

\begin{cor}
The vector bundles
$\Gr^W_hS_u^{\can}(E)$ and $\Gr^W_hS_{u}(E,P)$
are isomorphic
for any $u\in\KMS(\nbige^0)$ and $h\in\seisuu$.
\hfill\qed
\end{cor}

\subsubsection{Local lifting and the gluing matrices}

\label{subsubsection;10.25.2}

Let $\vecw$ be a frame of $\Gr^W\nbigg(E)$ over
$\cnum_{\lambda}$,
and let $\vecv$ be a normalizing frame of
the bundle $\Gr^W\nbigg(\nbige)$
on $\nbigx^{\shikaku}$ as in the previous subsubsection
\ref{subsubsection;10.25.30}.
Let us pick a point $\lambda_0\in\cnum_{\lambda}^{\ast}$,
and let $U(\lambda_0)$ be an appropriate neighbourhood of $\lambda_0$
in $\cnum_{\lambda^{\ast}}$.

Let take a non-negative number $\epsilon=\epsilon(\lambda_0)$
satisfying the following:
\begin{itemize}
\item
 In the case $0\not\in \KMS(\nbige^{\lambda_0})$,
 we put $\epsilon=0$.
\item
 In the case $0\in\KMS(\nbige^{\lambda_0})$,
 $\epsilon$ is taken any positive number
 such that $\{r\,|\,0<r\leq \epsilon\}\cap\KMS(\nbige^{\lambda_0})=\emptyset$.
\end{itemize}

For each $w_i$,
the integer $\nu(w_i)=\nu(w_i,\lambda_0)$ is determined 
by the condition
$-1+\epsilon(\lambda_0)
 <\paramap(\lambda_0,u(w_i))+\nu(w_i)\leq \epsilon(\lambda_0)$.

Let $\tilde{\vecv}$ be a frame of $\prolongg{\epsilon}{\nbige}$
on $X\times U(\lambda_0)$
satisfying the following:
\begin{itemize}
\item
 The frame 
 $\tilde{\vecv}$ is compatible with
 the filtrations $\nbigfzero$, the decompositions $\EEzero$
 and the weight filtration $W$.
\item
 We put $\bar{v}_i:=\tilde{v}_i\cdot z^{\nu(w_i)}$.
 Then the tuple $\bar{\vecv}=(\bar{v}_i)$
 gives a frame of $\nbige$ on
 $(X-D)\times U$
 which is compatible with $\nbigfzero$, $\EEzero$ and $W$.
 Then $\bar{\vecv}$ induces the frame
 of $\Gr^W\nbigg(\nbige)$.
\item
 The induced frame of $\bar{\vecv}$ is same as
 $\vecv$.
\end{itemize}
Such $\tilde{\vecv}$ is called a local lift
of $\vecw$ around $\lambda_0$.

\begin{rem} \label{rem;b12.9.30}
In the case $\lambda_0=0$,
a local lift $\tilde{\vecv}$ of $\vecw$
is a frame of $\prolongg{\epsilon(0)}{\nbige}$
on $U(0)\times X$ satisfying the following:
\begin{itemize}
\item
$\tilde{\vecv}$ is compatible with 
the decomposition $\EE^{(0)}$, 
the parabolic filtration $F^{(0)}$,
and the weight filtration $W$.
\item
Then $\tilde{\vecv}$ induces the frame of
$\Gr^{W}\nbigg(E)$ on a neighbourhood of $U(0)$.
The induced frame is same as $\vecw$.
\hfill\qed
\end{itemize}
\end{rem}

\vspace{.1in}

Let $\vecw$ be a frame of $\Gr^{W}\nbigg^{\dagger}$
on $\cnum_{\mu}$,
and let $\vecv^{\dagger}$ be a normalizing frame of
the bundle $\Gr^W\nbigg^{\dagger}(\nbige^{\dagger})$ over
$X^{\dagger}\times\cnum_{\mu}^{\ast}$
as in the subsubsection \ref{subsubsection;10.25.30}.
Let $U(\lambda_0)$ be as above.
Note that $\lambda_0$ gives the point
$\mu_0=\lambda_0^{-1}\in\cnum_{\mu}^{\ast}$.
Let $U'$ denote the neighbourhood of $\lambda_0^{-1}$
corresponding to $U(\lambda_0)$ via the isomorphism
$\cnum_{\lambda}^{\ast}\simeq\cnum_{\mu}^{\ast}$.

We have the non-negative number
$\epsilon'=\epsilon\bigl(-\overline{\lambda}_0^{-1}\bigr)$.
For each $w_i^{\dagger}$,
the integer $\nu(w_i^{\dagger})$ is determined by the condition
$-\epsilon'<
 \paramap(\lambda_0^{-1},u(w_i^{\dagger})^{\dagger})
 +\nu(w_i^{\dagger})
 \leq 1-\epsilon'$.
We can take a frame $\tilde{\vecv}^{\dagger}$
of $\prolongg{1-\epsilon'}{\nbige^{\dagger}}$
satisfying the following:
\begin{itemize}
\item
 The frame $\tilde{\vecv}^{\dagger}$ is compatible with
 $\nbigf^{(\mu_0)}$, $\EE^{(\mu_0)}$ and $W$.
\item
 We put $\bar{v}_i^{\dagger}=
 \bar{z}^{\nu(w_i^{\dagger})}\cdot
 \tilde{v}_i^{\dagger}$.
 Then the frame $\bar{\vecv}^{\dagger}$ is compatible with
 $\nbigf^{(\lambda_0)}$, $\EEzero$ and $W$.
 Then $\bar{\vecv}^{\dagger}$ induces the frame of
 $\Gr^W\nbigg^{\dagger}(E)$.
\item
 The induced frame of $\bar{\vecv}^{\dagger}$ is same as $\vecv^{\dagger}$.
\end{itemize}
Such $\tilde{\vecv}^{\dagger}$ is called a local lift of
$\vecv^{\dagger}$ around $\mu_0$.

\begin{rem}
As in Remark {\rm\ref{rem;b12.9.30}},
a local lift at $\mu_0=0$ is also defined,
similarly.
\hfill\qed
\end{rem}

On $(X-D)\times U=(X^{\dagger}-D^{\dagger})\times U'$,
we have $\nbige=\nbige^{\dagger}$ as
$U$-holomorphic bundles.
Hence we obtain the relation:
\[
 \tilde{v}^{\dagger}_i
=\sum \tilde{J}_{j\,i}\cdot \tilde{v}_j.
\]
Here $\tilde{J}_{j\,i}$ are $U$-holomorphic.

\begin{lem}\label{lem;10.25.5}
In the case $\tilde{J}_{j\,i}\neq 0$,
we have the following:
\[
 \eigenmap^f\bigl(\lambda_0^{-1},u(w^{\dagger}_i)^{\dagger}\bigr)
=\eigenmap^f\bigl(\lambda_0,u(w_j)\bigr)^{-1},
\quad
 \paramap^f\bigl(\lambda_0^{-1},u(w^{\dagger}_i)^{\dagger}\bigr)
\geq
 \paramap^f\bigl(\lambda_0,u(w_j)\bigr),
\quad
 \deg^W(w_i^{\dagger})
\geq
 \deg^W(w_j).
\]
Here $\vecw$ and $\vecw^{\dagger}$ as in the subsubsection
{\rm\ref{subsubsection;10.25.30}}.
\end{lem}
\pf
It follows from the compatibility 
of the $\tilde{\vecv}$, $\tilde{\vecv}^{\dagger}$
and the filtrations and the decompositions.
\hfill\qed

\begin{lem}\label{lem;10.25.6}
In the case
$u(w^{\dagger}_i)
=u(w^{\dagger}_j)$
and $\deg^W(w^{\dagger}_i)=\deg^W(w_j)$,
we have the relation:
\[
 \tilde{J}_{j\,i}=
 J_{j\,i}\cdot z^{\nu(w_i)}\cdot \bar{z}^{-\nu(w_i^{\dagger})}.
\]
\end{lem}
\pf
It follows from the condition
that $\tilde{\vecv}$ and $\tilde{\vecv}^{\dagger}$
are lifts of $\vecv$ and $\vecv^{\dagger}$
respectively.
\hfill\qed

%% file: a88.1.tex

\subsubsection{Preliminary}
\label{subsubsection;a12.10.20}

We put $X=\Delta$ and $D=\{O\}$.
Let $\harmonicbundle$ be a tame harmonic bundle over $X-D$.
We can take a model bundle
 $(E_0,\delbar_{E_0},h_0,\theta_0)$.
We denote the deformed holomorphic bundle of $E_0$
by $\nbige_0$.

Let $\vece$ and $\vece_0$ be holomorphic frames
of $\prolong{E}$ and $\prolong{E_0}$ respectively,
which are compatible with $F$, $\EE$ and $W$.
We take the isomorphism
$\Phi:\prolong{E}_0\lrarr \prolong{E}$
via the condition $\Phi(\vece_0)=\vece$.
It satisfies the following:
\begin{itemize}
\item
The morphism $\Phi$ is compatible with $\EE$, $F$ and $W$
at the origin $O$.
\item
We have the induced isomorphism
$\Gr^F\Phi:\Gr^F(\prolong{E}_0)\lrarr \Gr^F(\prolong{E})$
and the endomorphisms
$\Gr^F\bigl(\Res(\theta)\bigr)$ $\Gr^F\bigl(\Res(\theta_0))\bigr)$.
Under the isomorphism $\Gr^F$,
we have 
$\Gr^F\bigl(\Res(\theta)\bigr)=\Gr^F\bigl(\Res(\theta_0)\bigr)$.
\end{itemize}
Recall that
the morphism $\Phi$ and $\Phi^{-1}$ give
the isomorphism of $E$ and $E_0$
which are bounded with respect to the metrics $h$ and $h_0$,
due to Simpson.
(See the subsubsection 4.3.3 in the previous paper, for example.)

We use the setting of the subsubsection \ref{subsubsection;10.25.2}.
Let $\vecw$ be a frame of $\Gr^W\nbigg(E)$ over $\cnum_{\lambda}$
and $\vecv$ be the normalizing frame of $\Gr^W\nbigg(\nbige)$
over $\nbigx$ as in the subsubsection \ref{subsubsection;10.25.2}.
Let $\lambda_0$ be a point of $\cnum_{\lambda}$.
Let $\epsilon=\epsilon(\lambda_0)$ be the non-negative number
and $U(\lambda_0)$ be an appropriate neighbourhood of $\lambda_0$
as in the subsubsection \ref{subsubsection;10.25.2}.
Let $\vecv^{(\lambda_0)}$ be a local lift of
$\vecw$ on $X\times U(\lambda_0)$.
(Note that we use the  notation $\vecv$
 instead of $\tilde{\vecv}$, for simplicity).
In the case $\lambda_0=0$,
we may assume
$\vecv^{(0)}_{|X\times\{0\}}=\vece$.

On the other hand,
let $\vecw_0$ be a frame of
$\Gr^{W}\nbigg(E_0)$ over $\cnum_{\lambda}$.
We may assume that
we can take the canonical frame $\vecv_0^{(\lambda_0)}$
of $\prolongg{\epsilon(\lambda_0)}{\nbige_0}$
on $X\times U(\lambda_0)$,
which induce $\vecw_0$
(See the subsubsection \ref{subsubsection;b12.9.50}
for a canonical frame).
It is compatible with $\EE^{(\lambda_0)}$, $\nbigf^{(\lambda_0)}$
and $W$ of $\nbige_0$,
in the case $\lambda_0\neq 0$.
In the case $\lambda_0=0$,
it is compatible with $\EE^{(0)}$, $F^{(0)}$ and $W$.

We put as follows:
\begin{equation} \label{eq;10.25.10}
\begin{array}{l}
 \vecv^{(\lambda_0)\,\prime}=
 \bigl(v_i^{(\lambda_0)\,\prime}\bigr),
\quad
 v_i^{(\lambda_0)\prime}
:=v_i^{(\lambda_0)}\cdot |z|^{b(v_i)}
\cdot\bigl(-\log|z|\bigr)^{-\frac{1}{2}k(v_i)},\\
\mbox{{}}\\
 \vecv^{(\lambda_0)\,\prime}_0=
 \bigl(v_{0\,i}^{(\lambda_0)\,\prime}\bigr),\quad
 v_{0\,i}^{(\lambda_0)\prime}
:=v_{0\,i}^{(\lambda_0)}\cdot |z|^{b(v_{0\,i})}
 \cdot\bigl(-\log|z|\bigr)^{-\frac{1}{2}k(v_{0\,i})}.
\end{array}
\end{equation}
Here we put
$b(v_i)(\lambda):=\deg^F(v_{i\,|\,\lambda})$
and $k(v_i):=\deg^{W}(v_i)$,
and similar to $b(v_{0\,i})$ and $k(v_{0\,i})$.

We put $r:=\rank(E)$.
We have the $GL(r)$-valued $C^{\infty}$-functions $B^{(\lambda_0)}$
and $B_0^{(\lambda_0)}$
determined by the following conditions:
\begin{equation}\label{eq;d12.9.1}
 \vecv'=\vece'\cdot B^{(\lambda_0)},
\quad\quad
 \vecv_0'=\vece_0'\cdot B^{(\lambda_0)}_0.
\end{equation}

The  functions $I^{(\lambda_0)\prime}_{i\,j}$ are determined 
as follows:
\[
 \Phi(v^{(\lambda_0)\prime}_{0\,i})
=\sum I_{j\,i}^{(\lambda_0)\prime}\cdot
 v_j^{(\lambda_0)\,\prime}.
\]
Then we obtain the $M(r)$-valued function
$I^{(\lambda_0)\,\prime}:=\bigl(I_{j\,i}^{(\lambda_0)\,\prime}\bigr):
 U(\lambda_0)\times(X-D)\lrarr M(r)$.
The following lemma is easy to see.

\begin{lem}
In the case
 $\deg^{\EEzero,\Fzero}(v_{0\,i}^{(\lambda_0)})
 =\deg^{\EEzero,\Fzero}(v_j^{(\lambda_0)})$
and $k(v_{0\,i}^{(\lambda_0)})=k(v_j^{(\lambda_0)})$,
the function $I_{i\,j}^{(\lambda_0)\,\prime}$ is holomorphic
with respect to the variable $\lambda$.
\end{lem}
\pf
If we denote
$\Phi(v_{0\,i}^{(\lambda_0)})
 =\sum I_{j\,i}^{(\lambda_0)}
 \cdot v_j^{(\lambda_0)}$,
then $I_{j\,i}^{(\lambda_0)}$ are holomorphic
with respect to $\lambda$.
Then the lemma immediately follows.
\hfill\qed

\begin{lem} \label{lem;9.15.1} \mbox{{}}
\begin{itemize}
\item
 $I^{(\lambda_0)\,\prime}$ and $I^{(\lambda_0)\,\prime\,-1}$
 are bounded.
\item
 We have the estimate
 $|I^{(\lambda_0)\,\prime}_{j\,i}|\leq C\cdot (-\log|z|)^{-1}$
 unless
 $\deg^{\EE^{(\lambda_0)},F^{(\lambda_0)}}(v_{0\,i}^{(\lambda_0)})
 =\deg^{\EE^{(\lambda_0)},F^{(\lambda_0)}}(v_j^{(\lambda_0)})$.
\item
 In the case $\deg^W\bigl(v^{(\lambda_0)}_{0\,i}\bigr)
 \neq \deg^W\bigl(v_j^{(\lambda_0)}\bigr)$,
 we have the following finiteness:
\[
 \int_{U(\lambda_0)}
 || I^{(\lambda_0)\,\prime}_{j\,i}||_W<\infty.
\]
\end{itemize}
\end{lem}
\pf
It can be shown by an argument
similar to the proof of Proposition \ref{prop;9.5.15}.
\hfill\qed

%% file: a88.2.tex

\subsubsection{The construction of the isomorphism $\Psi$}
\label{subsubsection;10.19.15}

Let $\Gr^W\nbigg(E)$ and $\Gr^W\nbigg(E_0)$
be as in the subsubsection \ref{subsubsection;10.25.30}.
In this subsubsection, we would like to construct
the holomorphic isomorphism
$\Psi:
 \Gr^W\nbigg(E)\lrarr
 \Gr^W\nbigg(E_0)$.

For any point $\lambda_0\in\cnum_{\lambda}$,
we take neighbourhood $U(\lambda_0)$
as in the subsubsection \ref{subsubsection;10.25.2}.
Then we obtain the covering
$\bigl\{U(\lambda_0)\,\big|\,
 \lambda_0\in\cnum_{\lambda} \bigr\}$ 
of the complex plane $\cnum_{\lambda}$.
Then we can take a discrete subset $S$ of $\cnum_{\lambda}$
such that
$\bigl\{U(\lambda_0)\,\big|\,
 \lambda_0\in S\bigr\}$ is a covering of $\cnum_{\lambda}$.
We may assume that $0\in S$.

On each $X\times U(\lambda_0)$ $(\lambda_0\in S)$,
we have the frames $\vecv^{(\lambda_0)}$ 
and $\vecv^{(\lambda_0)}_0$ of $\prolongg{\epsilon}{\nbige}$
and $\prolongg{\epsilon}{\nbige_0}$ respectively,
as in the subsubsection \ref{subsubsection;a12.10.20}.
We recall that we assume that 
$\vecv^{(0)}_{|\nbigx^{0}}=\vece$,
$\vecv^{(0)}_{0|\nbigx^0}=\vece_0$
and $\Phi(\vece_0)=\vece$.
In the following, we identify $E$ and $E_0$
via the isomorphism $\Phi$
as $C^{\infty}$-vector bundles on $X-D$.

\begin{lem} \label{lem;9.15.4}
We can take
a sequence of subsets
$\bigl\{U_N\subset X-D\,\big|\,N=1,2,\ldots\bigr\}$
satisfying the following conditions:
\begin{enumerate}
\item
 The volume of $U_N$ with respect to the measure
${\displaystyle
 \frac{|dz|\cdot |d\bar{z}|}{|z|^2\cdot(-\log|z|)\cdot\log(-\log|z|)}}$
is infinite.
\item
 $U_N\supset U_{N+1}$.
\item \label{number;9.15.3}
 Let $P$ be any point of $U_N$
 and $\lambda_0$ be any point of
 $\bigl\{\lambda\in S\,\big|\,|\lambda|<N\bigr\}$.
 Then the inequality
 $\big|I^{(\lambda_0)\,\prime}_{j\,i}(P)\big| \leq N^{-1}$
 holds for any $i$ and $j$ such that
 $\deg^{\Fzero,\EEzero,W}\bigl(v_{0\,i}^{(\lambda_0)}\bigr)
 \neq
 \deg^{\Fzero,\EEzero,W}\bigl(v_j^{(\lambda_0)}\bigr)$.
\item \label{number;9.15.10}
 Let $\vece$ be a frame of $\prolong{E}$ as above.
For any point $P\in U_N$,
the inequality $|h(e_i',e_j')(P)|\leq N^{-1}$
holds,
in the case $\deg^{\EE,F,W}(e_i)\neq \deg^{\EE,F,W}(e_j)$.
\item \label{number;9.15.11}
 For any point $P\in U_N$, the following inequalities hold:
\[
 |z|\cdot(-\log|z|)\cdot
  \bigl|
(\theta_0-\theta)
 \bigr|_h(P)\leq N^{-1},
\quad\quad
  |z|\cdot(-\log|z|)\cdot
 \bigl|
 (\theta_0^{\dagger}-\theta^{\dagger})
 \bigr|_h(P)\leq N^{-1}.
\]
\end{enumerate}
\end{lem}
\pf
It follows from Lemma \ref{lem;9.15.1},
the asymptotic orthogonality
(the subsubsections \ref{subsubsection;a12.1.10},
 \ref{subsubsection;a12.1.11}
 and \ref{subsubsection;a12.1.12}),
Lemma \ref{lem;9.15.2},
and our choice $\Gr^F\bigl(\Res(\theta-\theta_0)\bigr)=0$.
\hfill\qed

\vspace{.1in}

Let us pick points $P_N\in U_N$ such that $\lim_{N\to \infty}P_N=O$.

\begin{lem} \label{lem;10.19.3}
Note that we have the following
for any $\lambda_0\in S$:
\begin{enumerate}
\item \label{number;10.19.4}
 For any point $\lambda\in U(\lambda_0)$,
 we have
 $\lim_{N\to\infty} I_{j\,i}^{(\lambda_0)\,\prime}(\lambda,P_N)=0$
 in the case 
 $\deg^{\EE,F,W}(v_{0\,i}^{(\lambda_0)})
 \neq \deg^{\EE,F,W}(v_j^{(\lambda_0)})$.
 It follows from the claim {\rm\ref{number;9.15.3}}
 in Lemma {\rm\ref{lem;9.15.4}}.
\item \label{number;10.19.5}
 The $M(r)$-valued functions
 $I^{(\lambda_0)\,\prime}_{|U(\lambda_0)\times\{P_N\}}$
 and
 $I^{(\lambda_0)\,\prime\,-1}_{|U(\lambda_0)\times\{P_N\}}$
 are bounded.
 In the case
 $\deg^{\EEzero,\Fzero\,W}(v_{0\,i})=\deg^{\EEzero,\Fzero\,W}(v_j)$,
 the function
 $I_{\infty\,i\,j\,|\,
 U(\lambda_0)\times\{P_N\}}^{(\lambda_0)\,\prime}$
 is holomorphic with respect to the variable $\lambda$.
\item \label{number;10.19.6}
 The sequences of hermitian matrices
 $\bigl\{H(h,\vece')_{|P_N}\bigr\}$ and
 $\bigl\{H(h,\vece')^{-1}_{|P_N}\bigr\}$
 are bounded.
\hfill\qed
\end{enumerate}
\end{lem}

\begin{lem}
By taking a subsequence $\{N_i\}$ of $\{N\}$,
we may assume that
$\bigl\{I^{(\lambda_0)\,\prime}
  _{|U(\lambda_0)\times\{P_{N_i}\}}
 \bigr\}$,
$\bigl\{I^{(\lambda_0)\,\prime\,-1}
 _{|U(\lambda_0)\times\{P_{N_i}\}}
 \bigr\}$,
$\bigl\{H(h,\vece')_{|P_{N_i}}\bigr\}$
and $\bigl\{H(h,\vece')^{-1}_{|P_{N_i}}\bigr\}$
are convergent for any $\lambda_0\in S$.
\end{lem}
\pf
Since 
$\bigl\{H(h,\vece')_{|P_N}\bigr\}$
and $\bigl\{H(h,\vece')^{-1}_{|P_N}\bigr\}$
are bounded sequences of hermitian matrices,
we can take a convergent subsequence,
when we fix $\lambda_0$.
Due to the claims \ref{number;10.19.4} and \ref{number;10.19.5}
in Lemma \ref{lem;10.19.3},
we can take a convergent subsequences of
$\bigl\{I^{(\lambda_0)\,\prime\,-1}
  _{|U(\lambda_0)\times\{P_N\}}\bigr\}$
when we fix $\lambda_0$.
Then we obtain the lemma by using the standard diagonal argument.
\hfill\qed

\vspace{.1in}

We denote the limit by $I^{(\lambda_0)\prime}_{\infty}$
and $H^{(\lambda_0)}_{\infty}$.
Then $I^{(\lambda_0)\prime}_{\infty}$
is holomorphic $M(r)$-valued function
on $U(\lambda_0)$
(note the claim \ref{number;10.19.5} in Lemma \ref{lem;10.19.3}).
Due to our construction and the claim \ref{number;10.19.4}
in Lemma \ref{lem;10.19.3},
the $M(r)$-valued function
$I^{(\lambda_0)\,\prime}_{\infty}$ can be regarded
as the direct sum:
\[
 \bigoplus_{(u,k)\in
 \nbigk(\nbige,\lambda_0)\times\seisuu}
 I_{u,k}^{(\lambda_0)},\quad
 I_{u,k}^{(\lambda_0)}:=
 \Bigl(
 I^{(\lambda_0)\,\prime}_{\infty\,\,i,j}\,\,
 \Big|\,\,
 \deg^{\Fzero,\EEzero,W}\bigl(v_{0\,i}^{(\lambda_0)}\bigr)
=\deg^{\Fzero,\EEzero,W}\bigl(v_j^{(\lambda_0)}\bigr)
=\bigl(\kmsmap(\lambda_0,u),k\bigr)
 \Bigr).
\]

Then the frames $\vecw$, $\vecw_0$
and the function $I_{\infty}^{(\lambda_0)\,\prime}$
induce the isomorphism defined on
$U(\lambda_0)$:
\[
 \Psi^{(\lambda_0)}:
 \Gr^{W}\nbigg(E_0)_{|U(\lambda_0)}
\lrarr
 \Gr^W\nbigg(E)_{|U(\lambda_0)}.
\]
Note that we use the following easy lemma implicitly.
\begin{lem}
In the case
 $\deg^{\Fzero,\EEzero}(v_{0\,i})=\deg^{\Fzero,\EEzero}(v_j)$,
we have 
$\nu(w_{0\,i},\lambda_0)=\nu(w_j,\lambda_0)$.
(See the subsubsection {\rm\ref{subsubsection;10.25.2}}
for $\nu(w_j,\lambda_0)$).
\hfill\qed
\end{lem}

It is easy to check the following lemma.
\begin{lem}\mbox{{}} \label{lem;c12.9.10}
\begin{itemize}
\item
Fix a subsequence $\{P_l\}$.
Then the morphism $\Psi^{(\lambda_0)}$ is independent of
choices of $\vecv^{(\lambda_0)}$, $\vecv_0^{(\lambda_0)}$,
 $\vecw$ and $\vecw_0$.
\item
In the case 
$A:=U(\lambda_0)\cap U(\lambda_1)\neq \emptyset$,
we have
$\Psi^{(\lambda_0)}_{|A}
=\Psi^{(\lambda_1)}_{|A}$.
\item
In particular,
we obtain the global isomorphism on $\cnum_{\lambda}$:
\[
  \Psi:
\Gr^W\nbigg(E_0)\lrarr\Gr^W\nbigg(E).
\]
Once we fix a subsequence $\{P_l\}$,
then the morphism $\Psi$ is independent of choices of
$\vecw$, $\vecw_0$, $\vecv^{(\lambda_0)}$, $\vecv^{(\lambda_0)}_0$
and $S$.
\hfill\qed
\end{itemize}
\end{lem}

In the following,
we may assume that $\Psi(\vecw_0)=\vecw$.
We can also take $S$ as $\cnum_{\lambda}$
not a discrete subset.

On the other hand $H^{(\lambda_0)}_{\infty}$
is a positive definite hermitian matrix,
and it can be regarded as a direct sum
of the hermitian matrices $H_{u}^{(\lambda_0)}$
$(u\in \KMS(\prolongg{\epsilon}{\nbige^0}))$,
due to the claim \ref{number;10.19.6} in Lemma \ref{lem;10.19.3}.
We have the induced hermitian metric $h_{u,k}$
on the vector space $\Gr^W_k\nbigg_{u}(E)_{|0}$
induced by $H_{u,k}$.

\subsubsection{Modification of the model bundle}
\label{subsubsection;10.25.35}

We put $\theta=f_0\cdot dz/z$ and
$\theta^{\dagger}=f_0^{\dagger}\cdot d\bar{z}/\bar{z}$.
We pick $\rho$ as in (\ref{eq;9.15.5}) in the page \pageref{eq;9.15.5}.
Then the sequence of the endomorphisms
$\bigl\{(-\log|z|)\cdot\bigl(f_0(P_N)-\rho(P_N)\bigr)\bigr\}$
and
$\bigl\{(-\log|z|)\cdot
\bigl(f_0^{\dagger}(P_N)-\rho^{\dagger}(P_N)\bigr)\bigr\}$
converges to the morphisms,
due to the condition \ref{number;9.15.11} in Lemma \ref{lem;9.15.4}
and Proposition \ref{prop;9.6.77}:
\[
 \begin{array}{l}
 f_{\infty}:\Gr^W_k\nbigg_{u}(E)_{|0}
 \lrarr
 \Gr^W_{k-2}\nbigg_{u}(E)_{|0},\\
 \mbox{{}}\\
 f^{\dagger}_{\infty}:
 \Gr^W_k\nbigg_{u}(E)_{|0}
\lrarr
 \Gr^W_{k+2}\nbigg_{u}(E)_{|0}.
 \end{array}
\]
By our construction
$f_{\infty}$ and $f^{\dagger}_{\infty}$ are mutually adjoint
with respect to the metric
$h_{\infty}=\bigoplus h_{u,k}$.

On the other hand,
we have the metric $h_{0\,u\,k}$
on $\Gr^W_k(\nbigg_{u}(E_0))$.
We also have the morphisms
$f_{0\,\infty}$ and $f_{0\,\infty}^{\dagger}$,
which coincides $f_{\infty}$ and $f^{\dagger}_{\infty}$
under the isomorphism $\Psi$,
due to the conditions \ref{number;9.15.10} and \ref{number;9.15.11}
in Lemma \ref{lem;9.15.4}.

In the following,
we identify $\Gr^W_k(\nbigg_{u}(E_0))$
and $\Gr^W_k(\nbigg_{u}(E))$
via the isomorphism $\Psi$.
We also identify $(f_{\infty},f_{\infty}^{\dagger})$
and $(f_{0\,\infty},f_{0\,\infty}^{\dagger})$.

We have the primitive decomposition
$\Gr^W=\bigoplus P_hGr_k^W$
with respect to the morphisms $f_{\infty}=f_{0\,\infty}$.

From the construction of the model bundle,
the primitive decomposition is orthogonal
with respect to the hermitian metric $h_{0\,\infty}$.
In particular, 
the morphism
$f_{\infty}^{\dagger}=f_{0\,\infty}^{\dagger}$
preserves the primitive decomposition.

\begin{lem}
The primitive decomposition
is orthogonal with respect to the hermitian metric $h_{\infty}$.
\end{lem}
\pf
Since we have $f_{\infty}^{\dagger}=f_{0\,\infty}^{\dagger}$,
the morphism $f_{\infty}^{\dagger}$
preserves the primitive decomposition.
It implies the orthogonality of the primitive decomposition
with respect to the metric $h_{\infty}$.
\hfill\qed

\vspace{.1in}

It is easy to see that
the model bundle $(E_0,\delbar_{E_0},h_0,\theta_0)$
is isomorphic to the following:
\[
 \bigoplus_{k,u}
 (P_k\Gr^W_k,H_{0,u,k})\otimes Mod(k)\otimes L(u).
\]
We have the other model bundle $(E_1,\delbar_{E_1},h_1,\theta_1)$
given as follows:
\[
 \bigoplus_{k,u}
 (P_k\Gr^W_k,H_{u,k})\otimes
 Mod(k)\otimes L(u).
\]
Note that we have the natural isomorphism
of the deformed holomorphic bundles
$\nbige_0$ and $\nbige_1$ compatible with $\lambda$-connections,
which is mutually bounded,
and thus we have $\nbigg_{u}(E_1)\simeq \nbigg_u(E_0)$,
and then $H_{1\,u\,h}=H_{u\,h}$.
Thus we may assume $H_{0\,u\,h}=H_{u\,h}$
from the beginning.

%% file: a88.3.tex

\subsubsection{The morphism $\Psi^{\dagger}$}
\label{subsubsection;10.25.1}

Since we have the isomorphism
$\sigma^{\ast}\Gr^W\nbigg
\simeq
 \bigl(\Gr^W\nbigg\bigr)^{\lor}$,
the isomorphism $\Psi$ induces
the isomorphism
$\Gr^W\nbigg^{\dagger}(E_0)\lrarr\Gr^W\nbigg^{\dagger}(E)$.
It can be regarded as follows (Lemma \ref{lem;c12.9.50}):

Let $\vecw^{\dagger}$ be the frame of
$\Gr^{W}\nbigg^{\dagger}(E)$,
which is induced by the dual frame $\vecw^{\lor}$
of $\bigl(\Gr^{W}\nbigg(E)\bigr)^{\lor}$
via the isomorphism
$\sigma^{\ast}\Gr^W\nbigg
\simeq
 \bigl(\Gr^W\nbigg\bigr)^{\lor}$.
Similarly we have the frame $\vecw_0^{\dagger}$
of $\Gr^W\nbigg^{\dagger}(E_0)$.

On the other hand,
we have the frame $\vecv^{(\lambda_0)\,\dagger}$
of $\prolongg{1-\epsilon(\lambda_0)}{\nbige^{\dagger}}$
over $\sigma(U(\lambda_0))\times X$
given as in (\ref{eq;10.16.1}):
\begin{equation}\label{eq;c12.9.1}
 \vecv^{(\lambda_0)\,\dagger}
=\sigma^{\ast}\Bigl(
 \vecv^{(\lambda_0)}\cdot
 \overline{H(h,\vecv^{(\lambda_0)})}^{-1}
 \Bigr).
\end{equation}
It is easy to see that $\vecv^{(\lambda_0)\,\dagger}$
is a local lift of $\vecw^{\dagger}$
around $-\overline{\lambda}_0$.
Similarly we obtain the local lift 
$\vecv^{(\lambda_0)\,\dagger}_0$
of $\vecw_0^{\dagger}$ around $-\overline{\lambda}_0$.

From (\ref{eq;c12.9.1}),
we obtain the following:
\begin{multline}\label{eq;c12.9.61}
 \vecv^{(\lambda_0)\,\dagger\,\prime}
=\sigma^{\ast}\Bigl(
 \vecv^{(\lambda_0)\,\prime}
 \cdot \overline{H(h,\vecv^{(\lambda_0)\,\prime})}^{-1}
 \Bigr)
=\sigma^{\ast}\Bigl( 
 \vece'\cdot B^{(\lambda_0)}
 \cdot \overline{H(h,\vecv^{(\lambda_0)\,\prime})}^{-1}
\Bigr) \\
=\vece'\cdot\sigma^{\ast}\Bigl(
 \overline{H(h,\vece')}^{-1}
\cdot
 \lefttop{t}\overline{B}^{(\lambda_0)\,-1}
 \Bigr).
\end{multline}
Similarly we have the following:
\begin{equation}\label{eq;c12.9.60}
 \vecv_0^{(\lambda_0)\,\dagger\,\prime}
=\vece_0'\cdot
 \sigma^{\ast}\Bigl(
 \overline{H(h,\vece_0')}^{-1}
\cdot
 \lefttop{t}\overline{B}_0^{(\lambda_0)\,-1}
 \Bigr).
\end{equation}
We determine the $GL(r)$-valued function $C^{(\lambda_0)}$
on $\sigma(U(\lambda_0))\times (X-D)$
by the following condition:
\[
 \Phi\bigl(\vecv_0^{(\lambda_0)\,\dagger\,\prime}\bigr)
=\vecv^{(\lambda_0)\,\dagger\,\prime}\cdot
 C^{(\lambda_0)}.
\]

\begin{lem}\label{lem;c12.9.50}
Let $\{P_l\}$ be a sequence as in Lemma {\rm\ref{lem;c12.9.10}}.
The sequence $\big\{C^{(\lambda_0)}(\lambda,P_l)\big\}$
converges to the identity matrix.
\end{lem}
\pf
Recall that we have the following,
due to $\Psi(\vecw_0)=\vecw$:
\[
\lim_{N\to\infty}
 \bigl(
 B^{(\lambda_0)\,-1}\cdot B_0^{(\lambda_0)}
 \bigr)(\lambda,P_N)
=\lim _{N\to \infty}
 I^{(\lambda_0)\,\prime}(\lambda,P_N)=
 \mbox{\rm the identity matrix}
\]
We also have the following,
due to the modification
in the subsubsection \ref{subsubsection;10.25.35}:
\[
 \lim_{N\to\infty}
 H(h,\vece_0')^{-1}_{|P_N}\cdot H(h,\vece')_{|P_N}=
\mbox{the identity matrix}.
\]
We also have the boundedness 
of the sequences
$\bigl\{H(h,\vece')_{|P_N}\bigr\}$,
$\bigl\{H(h_0,\vece_0')_{|P_N}\bigr\}$,
$\bigl\{B^{(\lambda_0)\,-1}(\lambda,P_N)\bigr\}$
and $\bigl\{B^{(\lambda_0)}_0(\lambda,P_N)\bigr\}$.
Then the lemma immediately follows
from (\ref{eq;c12.9.61}) and (\ref{eq;c12.9.60}).
\hfill\qed

%% file: a58.tex

\subsubsection{The isomorphism of induced vector bundles}

In the subsubsection \ref{subsubsection;10.19.15},
we constructed the isomorphism
$\Psi:
 \Gr^W\nbigg(E_0)\lrarr
 \Gr^W\nbigg(E)$.
In the subsubsection \ref{subsubsection;10.25.1},
we obtained the isomorphism
$\Psi^{\dagger}:
 \Gr^W\nbigg^{\dagger}(E_0)
\lrarr
 \Gr^W\nbigg^{\dagger}(E)$.
They induce the isomorphisms
for any $k\in\seisuu$ and $u\in\KMS(\nbige^0)$:
\[
 \Psi:\Gr^W_k\nbigg_u(E_0)\lrarr\Gr^W_k\nbigg_u(E),
\quad
 \Psi^{\dagger}:
 \Gr^W_k\nbigg^{\dagger}_{u^{\dagger}}(E_0)
\lrarr
 \Gr^W_k\nbigg^{\dagger}_{u^{\dagger}}(E).
\]

\begin{prop} \label{prop;10.25.20}
Let $k$ be an integer and $u$ be an element of
$\KMS(\nbige^0)$.
The isomorphisms $\Psi$ and $\Psi^{\dagger}$ induce
the isomorphisms
$\Psi^{\sankaku}:\Gr^W_k S^{\can}_u(E_0)\lrarr \Gr^W S^{\can}_u(E)$
and
$\Psi^{\sankaku}:\Gr^W_k S_u(E_0,P)\lrarr \Gr^W_k S_u(E,P)$.
\end{prop}
\pf
We have only to show that the morphisms
are compatible with the gluings.
For simplicity, we put as follows:
\[
 \Gr^W S^{\can}(E)
:=\bigoplus_{k\in\seisuu}\bigoplus_{u\in\KMS(\prolong{\nbige^0})}
 \Gr^W_kS^{\can}_u(E).
\]
Similarly,
we have $\Gr^W S^{\can}(E_0)$.
We have only to show that
$\Psi$ and $\Psi^{\dagger}$
induce the isomorphism
of $\Gr^WS^{\can}(E_0)\lrarr\Gr^W S^{\can}(E)$.

Let $\vecw$ and $\vecw_0$
be frames of $\Gr^{W}\nbigg(E)$ and $\Gr^W\nbigg(E_0)$ respectively
such that
$\Psi(\vecw_0)=\vecw$.
We have the induced frames $\vecw^{\dagger}$ and $\vecw^{\dagger}_0$
of $\Gr^W\nbigg^{\dagger}(E)$ and $\Gr^W\nbigg^{\dagger}(E_0)$
respectively.
We have $\Psi^{\dagger}(\vecw_0^{\dagger})=\vecw^{\dagger}$.
We use $u(w_i)=u(w_{0\,i})=u(w_i^{\dagger})=u(w_{0\,i}^{\dagger})$
without mention.
We also use
$\deg^W(w_i)=\deg^W(w_{0\,i})$
and $\deg^{W}(w_i^{\dagger})=\deg^{W}(w_{0\,i}^{\dagger})$.

\begin{rem}
We also have
$\deg^{W}(w_i)=-\deg^{W}(w_i^{\dagger})$.
\hfill\qed
\end{rem}

We have the normalizing frames $\vecv$ (resp. $\vecv_0$)
of $\Gr^W\nbigg(\nbige)$ (resp. $\Gr^W\nbigg(\nbige_0)$)
over $\nbigx^{\shikaku}$,
which is a lift of $\vecw$ (resp. $\vecw_0$).
We also have the normalizing frame
$\vecv^{\dagger}$ (resp. $\vecv^{\dagger}_0$)
of $\Gr^W\nbigg^{\dagger}(\nbige^{\dagger})$
(resp. $\Gr^W\nbigg^{\dagger}(\nbige^{\dagger}_0)$),
which is a lift of $\vecw^{\dagger}$ (resp. $\vecw^{\dagger}_0$).

\begin{rem}
Note that we use $\vecv$ as a normalizing frame of
$\Gr^W\nbigg(\nbige)$, not a local lift.
\hfill\qed
\end{rem}

We have the relations
$\vecv=\vecv^{\dagger}\cdot J$
and
$\vecv_0=\vecv_0^{\dagger}\cdot J_0$.
It is clear that Proposition \ref{prop;10.25.20}
can be reduced to the following lemma,
due to Lemma \ref{lem;a12.9.10}.

\begin{lem}\label{lem;10.25.21}
We have $J=J_0$.
\end{lem}
\pf
Let $\lambda_0$ be any point of $\cnum_{\lambda}$.
We can take local lifts $\tilde{\vecv}^{(\lambda_0)}$
and $\tilde{\vecv}_0^{(\lambda_0)}$
of $\vecw$ and $\vecw_0$ on $U(\lambda_0)$,
which induce the normalizing frames
$\vecv$ and $\vecv_0$.
We take $C^{\infty}$-frames
$\tilde{\vecv}^{(\lambda_0)\,\prime}$
and $\tilde{\vecv}_0^{(\lambda_0)\,\prime}$
as in (\ref{eq;10.25.25}).
We have the functions $B^{(\lambda_0)}$ and $B_0^{(\lambda_0)}$
from $(X-D)\times U(\lambda_0)$
to $GL(r)$
determined by the conditions
$\tilde{\vecv}^{(\lambda_0)\,\prime}=\vece'\cdot B^{(\lambda_0)}$
and
$\tilde{\vecv}_0^{(\lambda_0)\,\prime}=\vece_0'\cdot B_0^{(\lambda_0)}$,
as in (\ref{eq;d12.9.1}).

On the other hand, we put
$ \tilde{\vecv}^{\dagger\,(\lambda_0^{-1})}
:=\bigl(
 \tilde{\vecv}^{(-\overline{\lambda}^{-1}_0)}
 \bigr)^{\dagger}$
as in (\ref{eq;c12.9.1}).
Similarly
we obtain $\tilde{\vecv}_0^{\dagger\,(\lambda_0^{-1})}$.
\begin{lem}
The frames
$\tilde{\vecv}^{\dagger\,(\lambda_0^{-1})}$
and $\tilde{\vecv}_0^{\dagger\,(\lambda_0^{-1})}$
are local lifts of $\vecw^{\dagger}$ and $\vecw^{\dagger}_0$
respectively.
They induce the normalizing frames
$\vecv^{\dagger\,(\lambda_0^{-1})}$
and $\vecv_0^{\dagger\,(\lambda_0^{-1})}$
respectively.
\end{lem}
\pf
It is easy to check the claim from the definitions.
\hfill\qed

\vspace{.1in}

We obtain $\tilde{\vecv}^{\dagger\,(\lambda_0^{-1})\,\prime}$
and $\tilde{\vecv}_0^{\dagger\,(\lambda_0^{-1})\,\prime}$
as in (\ref{eq;10.25.25}).
Then we have the following relations:
\[
 \tilde{\vecv}^{\dagger\,(\lambda_0^{-1})\,\prime}
=\tilde{\vecv}^{(\lambda_0)\,\prime}\cdot \tilde{J}^{(\lambda_0)},
\quad
 \tilde{\vecv}_0^{\dagger\,(\lambda_0^{-1})\,\prime}
=\tilde{\vecv}_0^{(\lambda_0)\,\prime}\cdot \tilde{J}^{(\lambda_0)}_0.
\]
On the other hand, we have the following equalities:
\begin{multline}
 \tilde{\vecv}^{\dagger\,(\lambda_0^{-1})\,\prime}
=\sigma^{\ast}\bigl(
 \tilde{\vecv}^{(-\overline{\lambda}^{-1}_0)\,\prime}
\cdot \overline{H(h,\tilde{\vecv}^{(-\overline{\lambda}^{-1}_0)\,\prime})^{-1}}
 \bigr)
=\vece'\cdot\overline{H(h,\vece')}^{-1}
 \cdot
 \sigma^{\ast}
 \bigl(
  \overline{ \lefttop{t}B^{(-\overline{\lambda}^{-1}_0)}}
 \bigr)^{-1} \\
=\tilde{\vecv}^{(\lambda_0)\,\prime}\cdot
 B^{(\lambda_0)\,-1}\cdot
 \overline{H(h,\vece')}^{-1}
 \cdot
 \sigma^{\ast}
 \bigl(
  \overline{ \lefttop{t}B^{(-\overline{\lambda}^{-1}_0)}}
 \bigr)^{-1}.
\end{multline}
Hence we obtain the relation:
\[
 \tilde{J}^{(\lambda_0)}
=B^{(\lambda_0)\,-1}\cdot
 \overline{H(h,\vece')}^{-1}
 \cdot
 \sigma^{\ast}
 \bigl(
  \overline{ \lefttop{t}B^{(-\overline{\lambda}^{-1}_0)}}
 \bigr)^{-1}.
\]
Similarly we have the relation:
\[
\tilde{J}_0^{(\lambda_0)}
=B_0^{(\lambda_0)\,-1}\cdot
 \overline{H(h_0,\vece_0')}^{-1}
 \cdot
 \sigma^{\ast}
 \bigl(
  \overline{ \lefttop{t}B_0^{(-\overline{\lambda}^{-1}_0)}}
 \bigr)^{-1}.
\]

\begin{lem} \label{lem;10.25.40}
Let $\{P_n\}$ be as in the subsubsection 
{\rm \ref{subsubsection;10.19.15}}.
We have the following:
\[
 \lim_{N\to\infty} \Bigl(
 \bigl(
B_0^{(\lambda_0)}\bigr) ^{-1}
 \cdot B^{(\lambda_0)}
 \Bigr)_{|U(\lambda_0)\times\{P_N\}}=
 \mbox{{\rm the identity matrix}},
\]
\[
 \lim_{N\to\infty} \Bigl(
 \bigl(
 B_0^{(-\overline{\lambda}^{-1}_0)}
 \bigr)^{-1}
 \cdot B^{(-\overline{\lambda}^{-1}_0)}
 \Bigr)_{|\sigma(U(\lambda))\times\{P_N\}}=
 \mbox{{\rm the identity matrix}.}
\]
\end{lem}
\pf
We have the relation:
$ \tilde{\vecv}^{(\lambda_0)\,\prime}
=\tilde{\vecv}^{(\lambda_0)\,\prime}_0\cdot
 \bigl(
 B_0^{(\lambda_0)}\bigr)^{-1}
 \cdot
 B^{(\lambda_0)}$.
Thus the first claim follows from our construction
of the morphism $\Psi$
and $\Psi(\vecw_0)=\vecw$.
The second claim can be shown similarly.
\hfill\qed

\begin{lem} \label{lem;10.25.41}
Let $\{P_N\}$ be as in the subsubsection
{\rm\ref{subsubsection;10.19.15}}.
We have the following:
\[
 \lim_{N\to\infty} H(h,\vece')_{|P_N}\cdot H(h_0,\vece_0')^{-1}_{|P_N}
=\mbox{{\rm the identity matrix}}.
\]
\end{lem}
\pf
It follows from our construction
(See the subsubsection \ref{subsubsection;10.25.35}).
\hfill\qed

\begin{lem}
We have the following:
\[
\lim_{N\to\infty} 
 \Bigl(
 \bigl(
 \tilde{J}^{(\lambda_0)}_0
 \bigr)^{-1}\cdot \tilde{J}^{(\lambda_0)}
 \Bigr)_{|U(\lambda_0)\times\{P_N\}}
=\mbox{{\rm the identity matrix}}
\]
\end{lem}
\pf
It follows from Lemma \ref{lem;10.25.40},
Lemma \ref{lem;10.25.41}
and the boundedness of
the sequences $\bigl\{H(h,\vece')_{|P_N}\bigr\}$ etc.
\hfill\qed

\vspace{.1in}

Let us return to the proof of Lemma \ref{lem;10.25.21}.
In the case 
$u(w_i^{\dagger})=u(w_j)=u$
and $\deg^{W}(w_i^{\dagger})=\deg^{W}(w_j)=h$,
we can develop
the $(i,j)$-component of
$\bigl(
 \tilde{J}^{(\lambda_0)}\bigr)^{-1}\cdot \tilde{J}^{(\lambda_0)}$
can be developed as follows,
by using Lemma \ref{lem;10.25.5}:
\[
 \Bigl(
 \bigl(
 \tilde{J}_0^{(\lambda_0)}\bigr)^{-1}
 \cdot\tilde{J}^{(\lambda_0)}
 \Bigr)_{i\,j}
=\sum_{\substack{u(w_k^{\dagger})=u,\\ \deg^{W}(w_k^{\dagger})=h}}
 \bigl(
 \tilde{J}^{(\lambda_0)}_0\bigr)^{-1}_{i,k}
 \cdot\tilde{J}^{(\lambda_0)}_{k,j}.
\]
Due to Lemma \ref{lem;10.25.5}
and Lemma \ref{lem;10.25.6},
we have the following, 
in the case $u(w_i^{\dagger})=u(w_j)=u$,
and $\deg^W(w_i^{\dagger})=\deg^W(w_j)$:
\[
 \tilde{J}^{(\lambda_0)}_{i\,j}=J_{i\,j}
 \cdot z^{\nu(w_i)}\cdot \bar{z}^{\nu(w_j^{\dagger})},
\quad
 \bigl(
 \tilde{J}^{(\lambda_0)}_0
 \bigr)^{-1}_{i\,j}
=\bigl(J_0^{-1}\bigr)_{i\,j}
\cdot z^{-\nu(w_i)}\cdot \bar{z}^{-\nu(w_j^{\dagger})}.
\]
Then we obtain the following:
\[
 \lim_{N\to\infty}
 \Bigl(
 J_0^{-1}\cdot J
 \Bigr)_{|U(\lambda_0)\times\{P_N\}}
=\mbox{{\rm the identity matrix}}.
\]
Due to Corollary \ref{cor;10.25.51},
we have the $GL(r)$-valued holomorphic functions
$K$ and $K_0$ on $\cnum^{\ast}$
such that $J(\lambda,z)=K(\lambda)\cdot C(\lambda,z)$
and $J_0(\lambda,z)=K_0(\lambda)\cdot C(\lambda,z)$,
where $C(\lambda,z)$ is given by the formula (\ref{eq;10.25.50}).
Thus we obtain the following:
\[
\mbox{{\rm the identity matrix}}
=\lim_{N\to\infty}
 \bigl(
 K_0(\lambda)\cdot C(\lambda,P_N)
 \bigr)^{-1}
\cdot
 K(\lambda)\cdot C(\lambda,P_N)
=\lim_{N\to\infty}K_0(\lambda)^{-1}\cdot K(\lambda)
=K_0(\lambda)^{-1}\cdot K(\lambda).
\]
Note $K(\lambda)$, $K_0(\lambda)$ and $C(\lambda,z)$ are commutative.
Thus we obtain $K(\lambda)=K_0(\lambda)$.
Then the equality $J(\lambda,z)=J_0(\lambda,z)$ holds.
Therefore we obtain Lemma \ref{lem;10.25.21},
and thus Proposition \ref{prop;10.25.20}.
\hfill\qed

%% file: a58.1.tex

\subsubsection{Theorem (the one dimensional case)}

\begin{prop}
The tuples $(\Gr^WS^{\can}_u(E),W,N^{\sankaku},S)$
and $(\Gr^WS^{\can}_u(E_0),W,N^{\sankaku},S)$
are isomorphic.
The tuples $(\Gr^WS_u(E_0,P),W,N_0^{\sankaku},S_0)$
and $(\Gr^{W}S_u(E,P),W,N^{\sankaku},S)$ are isomorphic.
\end{prop}
\pf
We have already constructed the isomorphism
$\Psi^{\sankaku}:S^{\can}_u(E)\lrarr S^{\can}$.
Under the isomorphism,
the restrictions of
the morphisms
$N_0^{\sankaku}:
 \Gr^W_k S^{\can}_u(E_0)\lrarr 
 \Gr^W_{k-2} S^{\can}_u(E_0)\otimes \Tate(-1)$
and
$N^{\sankaku}:
 \Gr^W_k S^{\can}_u(E)\lrarr 
 \Gr^W_{k-2} S^{\can}_u(E)\otimes \Tate(-1)$
to the fibers on $0\in\proj^1$ are same.
Since $\Gr^W_kS^{\can}_u(E_0)$
and $\Gr^W_{k-2}S^{\can}_u(E_0)\otimes\Tate(-1)$
are pure twistors of weight $k$,
we obtain the coincidence
$N_0^{\sankaku}=N^{\sankaku}$ over $\proj^1$.
Similarly we obtain the coincidence
$S_0=S$ over $\proj^1$.
\hfill\qed

\vspace{.1in}

Then we obtain the following immediately.
\begin{thm}\label{thm;9.15.16}
The tuples
 $(S^{\can}_{u}(E),W,\vecN^{\sankaku},S)$ and
 $(S_u(E,P),W,\vecN^{\sankaku},S)$ $(P\in X-D)$
are polarized mixed twistor of $(0,1)$-type.
\hfill\qed
\end{thm}

%% file: 21.4.tex

\subsubsection{Preliminary}

Let $S_i$ $(i=1,\ldots,l)$ be finite subsets of $\real\times\cnum$.
For any element $\vecc\in\seisuu_{>0}^l$,
we have the map
$\varphi_{\vecc}:\prod_{i=1}^l S_i\lrarr\real\times\cnum$
as follows:
\[
 \varphi_{\vecc}(\vecu)=\sum c_i\cdot u_i.
\]
Let $T$ denotes the set of the elements
$\vecc\in\seisuu_{>0}^l$
such that $\varphi_{\vecc}$ are injective.

\begin{lem}\label{lem;10.1.21}
The set $T$ is Zariski dense in $\cnum^l$.
Namely, there does not exists a closed subset $Z$ of $\cnum^l$
defined as a zero set of some polynomials
such that $T\subset Z$.
\end{lem}
\pf
We put as follows:
\[
\begin{array}{l}
 \eta_i:=\min\bigl\{|d-d'|\,\big|\,d,d'\in S_i,\,\,\,d\neq d'\bigr\}
 \\
 \mbox{{}}\\
 \xi_i:=\max\bigl\{|d-d'|\,\big|\,d,d'\in S_i\bigr\}.
\end{array}
\]
Then the set $T$ contains any elements $\vecc=(c_i)$
such that
$c_i\cdot\xi_i<3^{-1}\cdot c_{i+1}\cdot\eta_{i+1}$.
Then it is easy to see that $T$ is Zariski dense in $\cnum^n$.
\hfill\qed

\vspace{.1in}

We put
$S(\vecc):=\bigl\{\varphi_{\vecc}(\vecu)\,\big|\,
 \vecu\in \prod_{i=1}^l S_i\bigr\}$.
A complex number $\lambda$ is called generic
with respect to $S_i$ $(i=1,\ldots,l)$
and $S(\vecc)$ $(\vecc\in T)$,
if the maps
$\eigenmap(\lambda):S_i\lrarr \cnum$
and $\eigenmap(\lambda):S(\vecc)\lrarr\cnum$
are injective.
\begin{lem}
The set of the complex numbers,
which are not generic with respect to
$S_i$ $(i=1,\ldots,l)$
and $S(\vecc)$ $(\vecc\in T)$,
is discrete in $\cnum^{\ast}$.
\end{lem}
\pf
It can be shown by an argument similar to
the proof of Lemma \ref{lem;10.11.21}.
\hfill\qed

\subsubsection{Weak result}

We put $X=\Delta^n$ and $D=\bigcup_{i=1}^n D_i$.
Let $\harmonicbundle$ be a tame harmonic bundle over $X-D$.
We put $S_i:=\KMS(\prolong{\nbige^0},i)$
and we use Lemma \ref{lem;10.1.21}.
Let us take a point $\lambda_0\in\cnum_{\lambda}^{\ast}$
which is generic
with respect to $S_i$ $(i=1,\ldots,n)$
and $S(\vecc)$ $(\vecc\in T)$.

Let us consider the morphism $\phi_{\vecc}:\Delta\lrarr X$
given by $z\longmapsto \bigl(z_1^{c_1},\ldots,z_n^{c_n}\bigr)$.
The image is denoted by $C(\vecc)$.
Via the pull back,
we obtain the harmonic bundle
$\phi_{\vecc}^{-1}\harmonicbundle$
over the punctured disc $\Delta-\{O\}$.

\begin{lem}\mbox{{}}\label{lem;a12.1.20}
\begin{itemize}
\item
The morphism
$\KMS(\nbige^{0},\lbar)
 \lrarr\KMS(\phi^{\ast}_{\vecc}\nbige^{0})$
is given by the correspondence
$\vecu\longmapsto \vecc\cdot\vecu=\sum c_i\cdot u_i$.
\item
The morphism is injective.
\item
$\lambda_0$ is generic with respect to
the sets $\KMS(\phi^{\ast}_{\vecc}\nbige^0)$
for any $\vecc\in T$.
\end{itemize}
\end{lem}
\pf
The first claim follows from Corollary \ref{cor;10.26.35}.
The following map is injective,
due to our choices of $\vecc$ and $\lambda_0$:
\[
 \prod_{i=1}^l\KMS(\nbige^0,i)
\lrarr \cnum,
\quad
 \vecu\longmapsto
 \lambda_0^{-1}\cdot\eigenmap(\lambda_0,\vecc\cdot\vecu).
\]
Then the second and the third claims are obtained.
\hfill\qed

\vspace{.1in}
Let $\vecu$ be an element of $\KMS(\nbige^0,\lbar)$.
We have the nilpotent maps
$\nbign^{\lambda_0}_i$ $(i=1,\ldots,n)$ on
$\lefttop{\lbar}\nbigg_{\vecu\,|\,(\lambda_0,O)}$
induced by the residues.
For any element $\vecd\in\cnum^n$,
we put $\nbign^{\lambda_0}(\vecd):=\sum d_i\cdot\nbign^{\lambda_0}_i$.
Then $\nbign^{\lambda_0}(\vecd)$ is the endomorphism
of $\lefttop{\lbar}\nbigg_{\vecu\,|\,(\lambda_0,O)}$.
We also put $\nbign^{\lambda_0}(\vect):=\sum t_i\cdot \nbign^{\lambda_0}_i$
for variables $t_i$.
Then $\nbign^{\lambda_0}(\vect)$ is the endomorphism
of 
$\lefttop{\lbar}\nbigg_{\vecu\,|\,(\lambda_0,O)}
\otimes_{\cnum}\cnum(t_1,\ldots,t_n)$.
Recall that if the conjugacy classes
of $\nbign^{\lambda_0}(\vecd)$ and $\nbign^{\lambda_0}(\vect)$
are same, then $\vecd$ is called generic
with respect to the tuple
$\bigl(\nbign_1^{\lambda_0},\ldots,\nbign_l^{\lambda_0}\bigr)$.
The following lemma is easy to see.
\begin{lem}\label{lem;10.26.55}
Let $T_1$ denote the subset of $T$,
which consists of the elements $\vecc$
which are generic with respect to
$(\nbign_1^{\lambda_0},\ldots,\nbign^{\lambda_0}_l)$.
Then $T_1$ is Zariski dense in $\cnum^l$.
\hfill\qed
\end{lem}

\begin{lem} \label{lem;9.15.15}
We have the isomorphisms:
\[
 S^{\can}_{\vecu}(E)\simeq
 S^{\can}_{\vecc\cdot\vecu}\bigl(\phi_{\vecc}^{-1}E\bigr),
\quad
 S_{\vecu}\bigl(E,\phi_{\vecc}(P)\bigr)\simeq
 S_{\vecc\cdot\vecu}\bigl(\phi_{\vecc}^{-1}E,P\bigr).
\]
The isomorphisms are compatible with the pairings.
Under the isomorphisms,
the nilpotent induced by the residue on 
$S_{\vecc\cdot\vecu}^{\can}\bigl(\phi_{\vecc}^{-1}E\bigr)$
and 
$S_{\vecc\cdot\vecu}\bigl(\phi_{\vecc}^{-1}E,P\bigr)$
are given by
$\nbign^{\sankaku}(\vecc)=\sum c_i\cdot\nbign^{\sankaku}_i$.
\end{lem}
\pf
It follows from Lemma \ref{lem;10.26.50}
and the injectivity in Lemma \ref{lem;a12.1.20}.
\hfill\qed

\vspace{.1in}
Let $W(\vecc)$ denote the weight filtration on
$S_{\vecu}^{\can}(E)$ or $S_{\vecu}(E,P)$
induced by the nilpotent map $\nbign^{\sankaku}(\vecc)$.

\begin{cor}
Let $\vecc$ be an element of $T$.
The filtered vector bundles
$(S^{\can}_{\vecu}(E),W(\vecc))$ and $(S_{\vecu}(E,P),W(\vecc))$
are the mixed twistor structure.
\end{cor}
\pf
It follows from Lemma \ref{lem;9.15.15} and Theorem \ref{thm;9.15.16}.
\hfill\qed

\begin{lem}\label{lem;10.26.60}
Let $\vecc$ be an element $T_1$.
Then $\vecc$ is generic
with respect to the tuple
 $(\nbign_1^{\lambda},\ldots,\nbign_n^{\lambda})$
for any $\lambda\in\proj^1$.
\end{lem}
\pf
For any elements $\vecc_i\in T_1$, $(i=1,2)$
the conjugacy classes of $\nbign^{\lambda_0}(\vecc_1)$
and $\nbign^{\lambda_0}(\vecc_2)$ are same.
For any $\vecc\in S$ and for any $\lambda\in\proj^1$,
the conjugacy classes of
$\nbign^{\lambda_0}(\vecc)$ and $\nbign^{\lambda}(\vecc)$
are same.
Thus the conjugacy classes of
$\nbign^{\lambda}(\vecc_i)$ $(i=1,2)$ are same.
Since $T_1$ is Zariski dense in $\cnum^l$,
we can conclude that $\vecc$ $(\vecc\in T_1)$
are generic with respect to
the tuple $(\nbign_1^{\lambda},\ldots,\nbign_l^{\lambda})$.
\hfill\qed

\begin{lem} \label{lem;10.26.61}
For any $i$,
we have
$\nbign_i^{\sankaku}\cdot W_h(\vecc)\subset
 W_{h-1}(\vecc)\otimes\Tate(-1)$.
\end{lem}
\pf
Due to Lemma \ref{lem;10.26.60},
we may apply a lemma of Cattani-Kaplan
(Proposition 1.9 in \cite{ck}. 
It is not difficult to prove directly.)
\hfill\qed

\begin{lem}
For any $i$,
we have
$\nbign_i^{\sankaku}\cdot W_h(\vecc)\subset
 W_{h-2}(\vecc)\otimes\Tate(-1)$.
\end{lem}
\pf
We put $\theta=\sum f_i\cdot dz_i\cdot z_i^{-1}$,
and we put $\rho_i$ as in the page \pageref{eq;9.15.17}.
We have the following inequality
for some positive constant $C$, due to Simpson's Main estimate:
\begin{equation}\label{eq;10.26.51}
\bigl|(f_i-\rho_i)_{|C(\vecc)}\bigr|_h
\leq C\cdot |z|^{-1}\cdot\bigl(-\log|z|\bigr)^{-1}.
\end{equation}
Let us consider the restriction of $\nbign^{\sankaku}_i$
to the fibers over $0\in\proj^1$.
We denote the restriction by
$\nbign^{\sankaku}_{i\,|\,0}$.
Since $\nbign^{\sankaku}_{i\,|\,0}$ is given
by $(f_i-\rho_i)_{(0,O)}$,
we obtain
$\nbign^{\sankaku}_{i\,|\,0}\cdot
 W_{h}(\vecc)_{|0}\subset 
 W_{h-2}(\vecc)_{|0}$,
due to (\ref{eq;10.26.51}) and the norm estimate
in one dimensional case.
Similarly, we obtain
$\nbign^{\sankaku}_{i\,|\,\infty}\cdot W_{h}(\vecc)_{|\infty}\subset 
 W_{h-2}(\vecc)_{|\infty}$.

Due to Lemma \ref{lem;10.26.61},
we have
$\nbign^{\sankaku}\cdot W_h(\vecc)\subset W_{h-1}(\vecc)\otimes
\Tate(-1)$.
Let us consider the induced morphism
$\nbign^{\sankaku}:
 \Gr^{W(\vecc)}_h\lrarr \Gr^{W(\vecc)}_{h-1}\otimes\Tate(-1)$.
We have already shown the vanishing of the induced morphism
at $(x=0,\infty)$.
Note that $\Gr^{W(\vecc)}_h$ and 
$\Gr^{W(\vecc)}_{h-1}\otimes\Tate(-1)$
are pure twistors of weight $h$ and $h+1$ respectively.
Thus the vanishing at $x=0,\infty$
implies vanishing over $\proj^1$.
Thus we obtain the implication
$\nbign_{i}\cdot W_{h}(\vecc)\subset 
 W_{h-2}(\vecc)$.
\hfill\qed

\begin{cor}
For any element $\vecd\in\cnum^n$
which is generic with respect to
the tuple $(\nbign^{\sankaku}_1,\ldots,\nbign^{\sankaku}_n)$,
we have $W(\vecd)=W(\vecc)$,
where $\vecc$ is taken as in Lemma {\rm\ref{lem;10.26.55}}.
\hfill\qed
\end{cor}

%% file: a59.2.tex
\subsubsection{Preliminary norm estimate}

\label{subsubsection;10.26.85}

Let $\harmonicbundle$ be a tame harmonic bundle
over $X-D$.
Let $\lambda$ be generic,
and we take a normalizing frame $\vecv$ of
$\prolong{\nbigelambda}$,
which is compatible with $\EE$ and $F$.

We have the matrices $A^k\in M_r(\cnum)$ $(k=1,\ldots,n)$
determined as follows:
\[
 \DD^{\lambda}\vecv
=\vecv\cdot \sum_{k=1}^n A^k\cdot\frac{dz_k}{z_k}.
\]
Let $f_{A_k}$ denote the endomorphism induced by $A_k$
for the frame $\vecv$.
We denote the nilpotent part of $A_k$ by $N_k$.
We impose the following assumption in this subsubsection.

\begin{assumption}\label{assumption;10.26.65}
The conjugacy classes of $N_k$ are independent of
$k=1,\ldots,n$.
\hfill\qed
\end{assumption}

\vspace{.1in}

Under the assumption \ref{assumption;10.26.65},
we obtain the weight filtration $W^{(0)}$ of the vector bundle
$\prolong{\nbigelambda}$
induced by $N_k$, which is independent of a choice of $k$.
We may assume that $\vecv$ is compatible with the filtration $W^{(0)}$.
We put $b_i(v_j):=\lefttop{i}\deg^{F}(v_j)$
and $\vecb(v_j)=\bigl(b_i(v_j)\,\big|\,i=1,\ldots,n\bigr)$.
We put $\beta_i(v_j)=\lefttop{i}\deg^{\EE}(v_j)$
and $\vecbeta(v_j)=\bigl(\beta_i(v_j)\,\big|\,i=1,\ldots,n\bigr)$.
We put $k(v_j):=\deg^{W}(v_j)$.

Then we obtain the $C^{\infty}$-frame
$\vecv'=(v_i')$ of $\nbige^{\lambda}$,
given as follows:
\[
 v_i':=
 v_i\cdot \prod_{h=1}^{n}|z_h|^{b_h(v_i)}
 \cdot\Bigl(
 -\sum_{h=1}^n\log|z_h|^2
 \Bigr)^{-k(v_i)/2}.
\]
\begin{lem}\label{lem;10.26.75}
The frame $\vecv'$ is adapted over $X-D$.
\end{lem}
\pf
We put as follows:
\[
 Y_m:=
\Bigl\{
 (z_1,\ldots,z_n)\in X-D\,\Big|\,
 |z_h|=1,\,\,(h\leq m)
\Bigr\}.
\]
We consider the following claim:
\begin{quote}
$(P_m)$:
 The restriction of the frame
 $\vecv'_{|Y_m}$ is adapted over $Y_m$.
\end{quote}

We show the claim $(P_m)$ by a descending induction on $m$.
Since $Y_n$ is compact, the claim $(P_n)$ holds.
We assume that $(P_{m+1})$ holds,
and we will derive $(P_m)$.

Let us pick the elements
$\vecu_1,\ldots,\vecu_a\in \KMS(\nbige^0,\nbar)$
such that
$\bigl\{\kmsmap(\lambda,\vecu_i)\,\big|\,i=1,\ldots,a\bigr\}
=\KMS(\prolong{\nbigelambda},\nbar)$.
Then we have the generalized eigen decomposition
$\nbigelambda=\bigoplus_i \nbigelambda_{\vecu_i}$
of $f_{A_k}$ $(k=1,\ldots,n)$.
Here $\nbigelambda_{\vecu_i}$ denote the subbundle
corresponding to the eigenvalues
$\eigenmap(\lambda,\vecu_i)
=\bigl(\eigenmap(\lambda,q_1(\vecu_i)),
 \ldots,
 \eigenmap(\lambda,q_n(\vecu_i))\bigr)$.
Correspondingly we have the decomposition
of the endomorphisms $f_{A_m}$:
\[
 f_{A_m}=\bigoplus_{\vecu_i}
\Bigl(
 \eigenmap(\lambda,q_m(\vecu_i))+N_{m\,\vecu_i}
\Bigr).
\]

We take the model bundle $Mod(N_{m\,\vecu_i})$
corresponding to the nilpotent map $N_{m\,\vecu_i}$.
We obtain the deformed holomorphic bundle
$\nbigelambda_{\vecu_i}$,
the $\lambda$-connection $\DD^{\lambda}_{\vecu_i}$,
the metric $h_{1\,\vecu_i}$,
and the canonical frame $\vecv_{1\,\vecu_i}$.

Let $\phi$ be the holomorphic map
$X\lrarr \Delta$ given by
$(z_1,\ldots,z_n)\longmapsto \prod_{i=1}^n z_i$.
Then we obtain the harmonic bundle
$\phi^{-1}Mod(N_{\vecu_i})$,
the deformed holomorphic bundle
$\phi^{-1}\nbigelambda_{\vecu_i}$,
the $\lambda$-connection $\DD^{\lambda}_{1\,\vecu_i}$,
the metric $h_{1\,\vecu_i}$,
and the canonical frame $\phi^{-1}\vecu_i$,
over $X-D$.

On the other hand,
we have the model bundle $L(\vecu_i)$ over $X-D$
of rank $1$.
We obtain the deformed holomorphic bundle
$\nbigl^{\lambda}(\vecu_i)$,
the $\lambda$-connection $\DD^{\lambda}_{2\,\vecu_i}$,
the metric $h_{2\,\vecu_i}$,
and the canonical frame $e_{\vecu_i}$.

Then we obtain the following:
\[
 \nbigelambda_{0\,\vecu_i}:=
 \nbigl^{\lambda}(\vecu_i)\otimes
 \phi^{-1}\nbigelambda_{\vecu_i},
\,\,
 \DD_{0\,\vecu_i}^{\lambda}:=\DD_2^{\lambda}\otimes\phi^{-1}\DD_1^{\lambda},
\,\,
 h_{0\,\vecu_i}:=h_{2\,\vecu_i}\otimes\phi^{-1}h_{1\,\vecu_i},
\,\,
 \vecv_{0\,\vecu_i}
 :=e_{\vecu_0}\otimes\phi^{-1}\vecv_{1\,\vecu_i}.
\]
By taking a direct sum,
we obtain $\nbigelambda_0$, $\DD_0^{\lambda}$,
$h_{0}$ and $\vecv_0$.

Moreover, by taking the $dz_m$-component,
we obtain the $\lambda$-connections $\gminiq_m(\DD_0^{\lambda})$
and $\gminiq_m(\DD^{\lambda})$
along the $z_m$-direction.
Note the following relations due to our construction:
\begin{equation} \label{eq;9.15.30}
 \gminiq_m(\DD^{\lambda})\vecv=
 \vecv\cdot A_m\cdot\frac{dz_m}{z_m},
\quad\quad
 \gminiq_m(\DD_0^{\lambda})\vecv_0=
 \vecv_0\cdot A_m\cdot\frac{dz_m}{z_m}.
\end{equation}

Let consider the morphism
$\Phi:\nbigelambda_0\lrarr\nbigelambda$
given by the frames $\vecv_0$ and $\vecv$.
The equalities (\ref{eq;9.15.30}) implies that
$\Phi$ is flat with respect to the $\lambda$-connections
along the $z_m$-direction.
Moreover, 
the restriction $\Phi_{|Y_{m+1}}$ and the inverse 
$\Phi^{-1}_{|Y_{m+1}}$ are bounded,
due to our assumption of the induction.
Thus we obtain the boundedness of $\Phi_{|Y_m}$ and the inverse.
Then it is easy to derive the adaptedness
of $\vecv'$ on $Y_{m}$ 
and the induction can proceed.
See the subsection 6.1 in \cite{mochi}
 for more detail of the argument.
\hfill\qed

%% file: a59.1.tex

\subsubsection{Preliminary constantness of the filtrations}
\label{subsubsection;11.9.1}

We use the setting in the subsubsection \ref{subsubsection;10.26.85}.
Let us pick the elements
$\vecu_1,\ldots,\vecu_a\in \KMS(\nbige^0,\nbar)$
such that
$\bigl\{\eigenmap(\lambda,\vecu_i)\,\big|\,i=1,\ldots,a\bigr\}
=\KMS(\prolong{\nbigelambda},\nbar)$.
Then we have the generalized eigen decomposition
$\prolong{\nbigelambda}=\bigoplus_i \prolong{\nbigelambda_{\vecu_i}}$
of $f_{A_k}$ $(k=1,\ldots,n)$.
Here $\prolong{\nbigelambda_{\vecu_i}}$ denote the subbundle
corresponding to the eigenvalues
$\eigenmap(\lambda,\vecu_i)
=\bigl(\eigenmap(\lambda,q_1(\vecu_i)),
 \ldots,
 \eigenmap(\lambda,q_n(\vecu_i))\bigr)$.
Correspondingly we have the decomposition
of $A_k$:
\[
 A_k=\bigoplus_{\vecu_i}
\Bigl(
 \eigenmap\bigl(\lambda,q_k(\vecu_i)\bigr)+N_{k\,\vecu_i}
\Bigr).
\]
The frame $\vecv$ induces the frame $\vecv_{\vecu_i}$
of $\prolong{\nbigelambda_{\vecu_i}}$.

Let $\veca$ be any element of $\real_{>0}^n$.
Let $\hyperh$ denote the upper half plane
$\{\zeta\,|\,\Image(\zeta)>0\}$.
We put $y:=-\Image(\zeta)$.
Let us consider the following morphism
$\psi_{\veca}:\hyperh\lrarr X-D$,
given by $z_i=\exp\bigl(\sqrt{-1}\cdot a_i\cdot\zeta\bigr)$
for $i=1,\ldots,n$.
We obtain the harmonic bundle
$\psi_{\veca}^{-1}\harmonicbundle$.
We have the deformed holomorphic bundle
$\psi_{\veca}^{-1}\nbigelambda
=\bigoplus_{\vecu_i}\psi_{\veca}^{-1}\nbigelambda_{\vecu_i}$.
We have the holomorphic frame
$\vecw_{\vecu}:=\psi_{\veca}^{-1}(\vecv_{\vecu})$.
We have the following relation:
\begin{equation} \label{eq;10.26.70}
 \DD\vecw_{\vecu}=\vecw_{\vecu}\cdot
\Bigl(
 \sum_k \sqrt{-1}a_k\cdot A^k_{\vecu}
\Bigr)\cdot d\zeta
=\vecw_{\vecu}\cdot
 \sum_k
 a_k\cdot
 \Bigl(
 \eigenmap\bigl(\lambda,q_k(\vecu)\bigr)
+N_{\vecu,k}
\Bigr)\cdot
\sqrt{-1}\cdot d\zeta.
\end{equation}

We put $N_{\vecu}(\veca)=\sum_i a_i N_{\vecu,i}$,
and we take a model bundle
$Mod(N_{\vecu}(\veca))$
on $\Delta^{\ast}_z$ for $N_{\vecu}(\veca)$.
We put $\vecu(\veca):=\sum_k a_k\cdot q_k(\vecu)$,
and we take a model bundle $L(\vecu(\veca))$ on $\Delta^{\ast}_z$.
We put $E_{1,\vecu}:=L(\vecu(\veca))\otimes Mod(N_{\vecu}(\veca))$,
and we denote the deformed holomorphic bundle
by $\nbige^{\lambda}_{1,\vecu}$.
We have the metric $h_{1\,\vecu}$ 
and the $\lambda$-connection $\DD_{1\,\vecu}^{\lambda}$
on $\nbige^{\lambda}_{1,\vecu}$.
We have the normalizing frame $\vecv_{1\,\vecu}$,
such that
the following relation holds on $\Delta^{\ast}_z$:
\[
 \DD^{\lambda}_{\vecu}\vecv_{1\,\vecu}
=\vecv_{1\,\vecu}\cdot 
\bigl(
 \eigenmap\bigl(\lambda,\vecu(\veca)\bigr)
+N_{\vecu}(\veca)
\bigr)\cdot\frac{dz}{z}.
\]

Let $\phi$ denote the holomorphic map
$\hyperh\lrarr \Delta^{\ast}$ given by
$z=\exp(\sqrt{-1}\zeta)$.
We put $\nbige^{\lambda}_{0\,\vecu}:=\phi^{\ast}\nbige_{1\,\vecu}$
and the pull backs of
$h_{1\,\vecu}$, $\DD^{\lambda}_{1\,\vecu}$, and $\vecv_{1\,\vecu}$
are denoted by $h_{0,\vecu}$, $\DD^{\lambda}_{0\,\vecu}$
and $\vecv_{0\,\vecu}$ respectively.
On the upper half plane $\hyperh$,
we have the following:
\begin{equation}\label{eq;10.26.71}
\DD^{\lambda}_{0\,\vecu}\phi^{\ast}(\vecv_{0\,\vecu})
=\phi^{\ast}(\vecv_{0\,\vecu})\cdot
 \bigl(
 \eigenmap\bigl(\lambda,\vecu(\veca)\bigr)
+N_{\vecu}(\veca)
 \bigr)\cdot\sqrt{-1}d\zeta.
\end{equation}

Then we have the isomorphism
$\Phi_{\vecu}:\nbige^{\lambda}_{0\,\vecu}
\lrarr \psi_{\veca}^{-1}\nbige^{\lambda}_{\vecu}$
given by the frames $\vecv_{0\,\vecu}$ and $\vecv_{\vecu}$.
Then $\Phi_{\vecu}$ and $\Phi^{-1}_{\vecu}$ 
are compatible with $\lambda$-connections,
due to (\ref{eq;10.26.70}) and (\ref{eq;10.26.71}).

We put $\nbigelambda_0:=\bigoplus_{\vecu_i}\nbigelambda_{0\,\vecu}$.
We have the induced metric $h_0$.
We obtain the isomorphism
$\Phi:=\bigoplus_{\vecu}\Phi_{\vecu}$
from $\nbige^{\lambda}_{0}$ to $\psi_{\veca}^{-1}\nbigelambda$.
Then $\Phi$ and $\Phi^{-1}$ are compatible with
the $\lambda$-connections.
We regard them as the flat sections of
$Hom\bigl(\nbige^{\lambda}_{0},
 \psi_{\veca}^{-1}\nbigelambda\bigr)$
and
$Hom\bigl(\psi_{\veca}^{-1}\nbigelambda,
 \nbigelambda_{0}\bigr)$.
The metrics $h_{0}$ and $\psi_{\veca}^{-1}h$
induce the metrics
$h_{2}$ and $h_{3}$ of
$Hom\bigl(\nbige^{\lambda}_{0},
 \psi_{\veca}^{-1}\nbigelambda\bigr)$
and
$Hom\bigl(\psi_{\veca}^{-1}\nbigelambda,\nbigelambda_{0}
 \bigr)$.

\begin{lem}
We have the following estimate
for some positive constants $C_1$ and $C_2$:
\begin{equation} \label{eq;10.26.90}
 \max\Bigl\{
 \log |\Phi|_{h_2},\,\,
  \log |\Phi^{-1}|_{h_3}
\Bigr\}
\leq C_1+C_2\cdot \log y.
\end{equation}
\end{lem}
\pf
For each element $v_i\in \vecv_{\vecu}$,
we put as follows:
\[
 v_i':=v_i\cdot\prod_{j=1}^{l}|z_j|^{\paramap(\lambda_0,q_j(\vecu))}.
\]
Then $\vecv'=(v_i')$ is $C^{\infty}$-frame
of $\nbigelambda$ over $X-D$,
which is adapted up to log order.
In particular,
we have the following inequality on $\hyperh$,
for some positive constants $C_i$ $(i=1,2)$ and $M$:
\[
 C_1\cdot y^{-M}
\leq
 H(h,\psi_{\veca}^{-1}\vecv')
\leq
 C_2\cdot y^M
\]
We also have the following equality for $v_i\in\vecv_{\vecu}$:
\[
 \psi_{\veca}^{-1}v_i'
=w_i\cdot
 \exp\Bigl(
 -y\cdot \paramap(\lambda_0,\vecu(\veca))
 \Bigr).
\]

On the other hand,
we put as follows, for $v_{1\,j}\in\vecv_{1\,\vecu}$:
\[
 v'_{1\,j}:=v_{1\,j}\cdot |z|^{\paramap(\lambda_0,\vecu(\veca))}.
\]
Then the $C^{\infty}$-frame 
$\vecv'_1=(v'_{1\,j})$ of $\bigoplus_{\vecu}\nbigelambda_{1\,\vecu}$
is adapted up to log order.
We put $v'_{0\,j}:=\phi^{-1}v'_{1\,j}$
and $\vecv'_0:=(v'_{0\,j})$.
Then we obtain the following inequalities
for some positive constants $C_i'$ $(i=1,2)$ and $M'$:
\[
 C_1'\cdot y^{-M'}
\leq H(h_0,\vecv'_0)
\leq C_2'\cdot y^{M'}.
\]
We also have the following equality,
for $v_{0\,i}\in \vecv_{0\,\vecu}$:
\[
 v_{0\,i}'=v_{0\,i}\cdot
 \exp\Bigl(
 -y\cdot \paramap(\lambda_0,\vecu(\veca))
 \Bigr).
\]
Since $\Phi$ and $\Phi^{-1}$ are given
by the frames $\psi_{\veca}^{-1}\vecv_{\vecu}$
and $\vecv_{0\,\vecu}$,
we obtain the estimate
$\max\bigl\{
 |\Phi|_{h_2},\,|\Phi^{-1}|_{h_3}
 \bigr\}\leq C''\cdot y^{M''}$ for some positive constants
$C''$ and $M''$.
Thus we are done.
\hfill\qed

\begin{lem}
The functions $\log|\Phi|_{h_2}$
and $\log|\Phi^{-1}|_{h_3}$ are subharmonic.
Namely we have the following inequalities:
\begin{equation}\label{eq;10.26.81}
\Delta \log|\Phi|_{h_2}^2\leq 0,
\quad\quad
 \Delta \log|\Phi^{-1}|_{h_3}^2\leq 0.
\end{equation}
\end{lem}
\pf
It can be shown by an argument similar to the proof of
Lemma 4.1 in \cite{s2}.
\hfill\qed

\vspace{.1in}
The functions $\log |\Phi|_{h_2}^2$ and $\log |\Phi^{-1}|_{h_3}^2$
can be regarded as follows:
We have the holomorphic bundle
$\nbigelambda_1:=\bigoplus_{\vecu_i}\nbigelambda_{1\,\vecu_i}$
on $\Delta^{\ast}$.
We also have the frame $\vecv_1$ induced by
the frames $\vecv_{1\,\vecu}$.
We denote the projection of $\Delta^{\ast}\times (X-D)$
onto the $i$-th component by $q_i$.
We have the holomorphic bundles
$q_1^{-1}\nbige^{\lambda}_{1}$
and $q_2^{-1}\nbigelambda$.
The frames $\vecv_1$ and $\vecv$ give
the isomorphisms $\Phi'$ and $\Phi^{\prime\,-1}$
of $q_1^{-1}\nbige^{\lambda}_{1}$
and $q_2^{-1}\nbigelambda$.
The metrics $h_1$ and $h$ induce
the metrics $h_2'$ and $h_3'$ of 
the bundles
$Hom(q_1^{\ast}\nbige^{\lambda}_1,q_2^{\ast}\nbigelambda)$
and
$Hom(q_2^{\ast}\nbigelambda,q_1^{\ast}\nbigelambda_1)$.
The morphisms $\phi$ and $\psi_{\veca}$ induce
the morphism $F:\hyperh\lrarr \Delta^{\ast}\times(X-D)$.
Then we have the following equalities:
\begin{equation} \label{eq;10.26.80}
 \log |\Phi|_{h_2}^2=F^{-1}\Bigl(
 \log |\Phi'|_{h_2'}^2
 \Bigr),
\quad\quad
 \log |\Phi^{-1}|_{h_3}^2=F^{-1}\Bigl(
 \log |\Phi^{\prime\,-1}|_{h_3'}^2
 \Bigr).
\end{equation}

We obtain the functions $G_1:=\Xi\bigl(\log|\Phi|_{h_2}\bigr)$
and $G_2:=\Xi\bigl(\log|\Phi^{-1}|_{h_3}\bigr)$
on $\real_{>0}$ (the subsubsection \ref{subsubsection;10.26.95}).
\begin{lem} \label{lem;10.26.93}
There exists a positive constants $C$
such that the following inequality holds on
the upper half plane $\hyperh$:
\[
\max\bigl\{
 G_1,G_2
 \bigr\}
\leq C.
\]
\end{lem}
\pf
It follows from
(\ref{eq;10.26.81}), (\ref{eq;10.26.80}),
Lemma \ref{lem;10.26.91}
that the functions $G_i$ $(i=1,2)$ are convex below.
Then
Lemma \ref{lem;10.26.93} follows from
(\ref{eq;10.26.90}) and Lemma \ref{lem;10.26.92}.
\hfill\qed

\vspace{.1in}

We take the metric $h_{4\,\vecu}$ of
$\psi_{\veca}^{-1}\nbigelambda_{\vecu}$.
For any $w_i,w_j\in\vecw_{\vecu}$,
we put as follows:
\[
 h_{4\,\vecu}(w_i,w_j)=\delta_{i\,j}\cdot 
 \exp\Bigl(
 2\cdot y\cdot \paramap\bigl(\lambda,\vecu(\veca)\bigr)
 \Bigr)\cdot y^{k(v_i)}.
\]
Here $\delta_{i\,j}$ denotes $1$ in the case $(i=j)$ and $0$
in the case $i\neq j$.
Recall that we put $k(v_i):=\deg^W(v_i)$.
Taking the direct sum of $h_{4\,\vecu}$,
we have the induced metric $h_4$ of 
$\psi_{\veca}^{-1}\nbigelambda$.

\begin{lem}
The metrics $\psi_{\veca}^{-1}h$ and $h_4$ are mutually bounded
\end{lem}
\pf
It immediately follows from Lemma \ref{lem;10.26.75}.
\hfill\qed

\vspace{.1in}

We have the weight filtration $W^{(1)}$ on
$\nbigelambda_{1\,\vecu}$
induced by the logarithms of the monodromy.
We take the normalizing frame $\tilde{\vecv}_{0\,\vecu}$
of $\nbigelambda_{1\,\vecu}$,
which is compatible with $W^{(1)}$.
Let $\tilde{\vecv}_{0\,\vecu}$ denote the pull back 
of $\tilde{\vecv}_{1\,\vecu}$ via the morphism $\phi$.
We take the metric $h_{5\,\vecu}$ of $\nbigelambda_{0\,\vecu}$.
For any $\tilde{v}_{0\,i},\tilde{v}_{0\,j}\in\tilde{\vecv}_{0\,\vecu}$,
we put as follows:
\[
 h_{5\,\vecu}\bigl(\tilde{v}_{0\,i},\tilde{v}_{0\,j}\bigr)=
 \delta_{i\,j}\cdot
 \exp\Bigl(
 2\cdot y\cdot\paramap\bigl(\lambda,\vecu(\veca)\bigr)
 \Bigr)\cdot y^{k(\tilde{v}_{0\,i})}.
\]
Then we have the induced metric
$h_5$ of $\nbigelambda_0$,
and the metrics $h_5$ and $h_0$ are mutually bounded.

The metrics $h_4$ and $h_5$ induce the metrics
$h_6$ and $h_7$
of 
$Hom(\nbige^{\lambda}_0,\psi_{\veca}^{-1}\nbigelambda)$
and
$Hom(\psi_{\veca}^{-1}\nbigelambda,\nbigelambda_0)$.
Then $h_6$ and $h_2$ are mutually bounded,
and $h_7$ and $h_3$ are mutually bounded.
Hence we obtain the following inequality on $\hyperh$,
for some positive constants $C_4$:
\[
 \max\Bigl\{
 \Xi(
 \log
 |\Phi|_{h_6}^2,
 ),
 \Xi(
 \log|\Phi^{-1}|_{h_7}^2
 )
 \Bigr\}
\leq C_4.
\]
Since the functions
$\log |\Phi|_{h_6}$ and $\log|\Phi^{-1}|_{h_7}$ are independent of
the real part of $\zeta$,
we have the equalities
$ \Xi(\log |\Phi|_{h_6})=\log |\Phi|_{h_6}$
and
$ \Xi(\log|\Phi^{-1}|_{h_7})=\log|\Phi^{-1}|_{h_7}$.
Thus we obtain the following inequalities
on $\hyperh$, for some positive constant $C_5$:
\begin{equation} \label{eq;10.26.100}
 \max\Bigl\{
 |\Phi|_{h_6},\,\,
 |\Phi^{-1}|_{h_7}
 \Bigr\}\leq C_5.
\end{equation}

\begin{lem} \label{lem;9.16.1}
Under Assumption {\rm\ref{assumption;10.26.65}},
the weight filtration
of $\nbign(\veca)=\sum a_i\cdot \nbign^{\lambda}_i$
are independent of a choice of $\veca\in\real_{>0}^n$
for any $\lambda\in\proj^1$.
\end{lem}
\pf
Let $f$ be a holomorphic section of
$\prolongg{\paramap(\lambda_0,\vecu(\veca))}{
 \nbigelambda_{1\,\vecu}}$.
It gives the section of
$W^{(1)}_k\nbigelambda_{1\,\vecu}$
if and only if
we have the following estimate:
\begin{equation}\label{eq;a12.1.30}
 \bigl|
 \phi^{\ast}f
 \bigr|_{h_{5\,\vecu}}
=O\Bigl(
 \exp\bigl(-y\cdot\paramap(\lambda,\vecu(\veca))\bigr)
 \cdot y^{k/2}
 \Bigr).
\end{equation}
Since $h_{4\,\vecu}$ and $h_{5\,\vecu}$ are mutually bounded,
(\ref{eq;a12.1.30}) is equivalent to the following:
\[
  \bigl|
 \phi^{\ast}f
 \bigr|_{h_{5\,\vecu}}
=O\Bigl(
 \exp\bigl(-y\cdot\paramap(\lambda,\vecu(\veca))\bigr)
 \cdot y^{k/2}
 \Bigr).
\]
It is equivalent to
$\phi^{\ast}f$ gives the section of
$W^{(0)}_k\bigl(\nbigelambda_{\vecu}\bigr)$.
Since $W^{(1)}_k$ is the weight filtration
of $\nbign(\veca)$,
we obtain the result.
\hfill\qed

%% file: a59.3.tex

\subsubsection{Constantness of the filtration at generic $\lambda$}

Let $\harmonicbundle$ be a tame harmonic bundle.
Let $\lambda$ be generic,
and $\vecv$ be a normalizing frame of $\prolong{\nbigelambda}$,
compatible with $F$ and $\EE$.
We put $b_i(v_j)=\lefttop{i}\deg^F(v_j)$
and $\vecb(v_j)=(b_i(v_j)\,|\,i=1,\ldots,n)$.
We put $\beta_i(v_j)=\lefttop{i}\deg^{\EE}(v_j)$
and $\vecbeta(v_i)=(\beta_i(v_j)\,|\,i=1,\ldots,n)$.

Then we obtain the matrices $A^k\in M_r(\cnum)$
determined as follows:
\[
 \DD^{\lambda}\vecv
=\vecv\cdot \sum_{k=1}^n A^k\frac{d\zeta_k}{\zeta_k}.
\]
We denote the nilpotent part of $A_k$ by $N_k$.

Let us pick $\vecc_1,\ldots,\vecc_n\in\seisuu_{>0}^n$.
We denote the tuple $(\vecc_1,\ldots,\vecc_n)$
by $\vecC$.
Let $\phi_{\vecC}:X\lrarr X$ be the morphism
as follows:
\[
 \phi_{\vecC}^{-1}(\zeta_k)
=\prod_{j=1}^n z_h^{c_{h\,k}}.
\]
Then we have the following:
\[
 \phi_{\vecC}^{-1}\Bigl(
 \sum A^k\frac{d\zeta_k}{\zeta_k}
 \Bigr)
=\sum_{h=1}^n
 \Bigl(
 \sum_{k=1}^n
 c_{h\,k}\cdot A^k
 \Bigr)
 \frac{dz_h}{z_h}.
\]

We put $\vecc_h\cdot \vecb(v_j):=\sum c_{h\,i}\cdot b_i(v_j)$.
We decompose as follows:
\[
 \vecc_h\cdot \vecb(v_j)
=n_{h\,j}+\kappa_{h\,j},\quad\quad
 n_{h\,j}\in\seisuu,\,\,\,
 -1<\kappa_{h\,j}\leq 0.
\]
We put $\tilde{v}_j:=v_j\cdot \prod_{h=1}^nz^{n_{h\,j}}$,
and $\tilde{\vecv}=(\tilde{v}_j)$.
Then $\tilde{\vecv}$ is a normalizing frame of
$\prolong{\phi_{\vecc}^{-1}\nbigelambda}$.
We have
$ \tilde{\DD}\tilde{\vecv}
=\tilde{\vecv}\cdot \sum_{h}\tilde{A}^h\cdot dz_h/z_h$
for $\tilde{A}^h\in M_r(\cnum)$.
The components $\tilde{A}^{h}_{i\,j}$ is given as follows:
\[
 \tilde{A}^h_{i\,j}
=\sum _k c_{h\,k}\cdot A^k_{i\,j}
+\delta_{i\,j}\cdot n_{h\,j}.
\]
Hence the nilpotent part $\tilde{N}_h$
of $\tilde{A}^h$ is given by $N(\vecc_h):=\sum c_{h\,k}\cdot N_k$.
If $\vecc_h$ $(h=1,\ldots,n)$ are generic
with respect to the tuple
$(\nbign_1,\ldots,\nbign_n)$,
then the conjugacy classes of $\tilde{N}_h$
are independent of $h=1,\ldots,n$.

\begin{lem} \label{lem;9.16.2}
The weight filtration of $\nbign(\veca)$
are independent of
a choice of $\veca\in\real_{>0}^n$.
\end{lem}
\pf
We can pick generic elements
$\vecc_1,\ldots,\vecc_n\in\seisuu^{n}_{>0}$,
such that 
$\veca=\sum a_{h}'\vecc_h$ for some $\veca'=(a_h')\in\real_{>0}^n$.
We have the following relation:
\[
 N(\veca)
=\sum_k a_k N_k
=\sum_{k,h}c_{h\,k}\cdot a'_h\cdot N_k
=\sum _h a'_h\cdot N(\vecc_h).
\]
Thus Lemma \ref{lem;9.16.2} can be reduced to
Lemma \ref{lem;9.16.1}.
\hfill\qed

%% file: 21.5.tex

\subsubsection{Theorem (the higher dimensional case)}

Let $V$ be a pure twistor of weight $n$.
Let $S:V\otimes\sigma(V)\lrarr \Tate(-n)$ denote
the $(-1)^n$-symmetric pairing.
We say that $S$ is a semi-polarization,
if the induced hermitian pairing
on $H^0\bigl(\proj^1,V\otimes\nbigo(-n)\bigr)$ is 
positive semi-definite.

We have the induced pairing
$S:P_h\Gr^{W}_{-h}S^{\can}_{\vecu}\otimes
   P_{h}\Gr^W_{h}S^{\can}_{\vecu}\lrarr \Tate(0)$.
\begin{lem} \label{lem;9.15.21}
For any $\vecd\in\real^n_{>0}$,
$S\bigl(N(\vecd)^h\otimes\id\bigr)$ gives a semi-polarization
on $P\Gr^{W}_h$.
Here $P\Gr^W_h$ denotes the primitive part
for $N(\vecd)$.
\end{lem}
\pf
Let $q_i$ denote the projection of $\proj^1\times\real_{>0}^n$
onto the $i$-th component.
We have the $C^{\infty}$-vector bundle
$q_1^{\ast}\Gr^W_h$ on $\proj^1\times\real_{>0}^n$.
Since the conjugacy classes of $N(\veca)$ $(\veca\in\real^n_{>0})$
are independent of a choice of $\veca$,
we obtain the vector bundle
$Pq_1^{\ast}\Gr^W_h$,
by taking the primitive parts for $N(\veca)$.
For generic $\veca\in\rnum^{n}$,
the pairing $S\bigl(N(\vecd)^h\otimes\id\bigr)$
gives the polarization,
due to Theorem \ref{thm;9.15.16}.
Then it follows that the pairing
$S\bigl(N(\vecd)^h\otimes\id\bigr)$
is a semi-polarization for any elements $\vecd\in\real^n_{>0}$.
\hfill\qed

\begin{lem} \label{lem;9.15.20}
The conjugacy classes of $N(\vecd)$ are independent 
of a choice of $\vecd\in\real_{>0}^n$.
In particular,
$W(\veca)=W(\veca')$ for any $\veca,\veca'\in\real_{>0}^n$.
\end{lem}
\pf
It follows from Lemma \ref{lem;9.16.2}.
\hfill\qed

\begin{cor} \label{cor;9.15.25}
For any element $\vecd\in\real^n_{>0}$,
$S\bigl(N(\vecd)^h\otimes\id\bigr)$ gives a polarization
on the primitive part $P\Gr^{W}_h$
for $N(\vecd)$.
\end{cor}
\pf
Due to Lemma \ref{lem;9.15.20},
the pairing
$S\bigl(N(\veca)^h\otimes\id\bigr)$ is perfect.
Then it gives the polarization
due to Lemma \ref{lem;9.15.21}.
\hfill\qed

\begin{cor} \label{cor;9.15.26}
The tuple
 $(\Gr^WS^{\can}_{u}(E),W^{(0)},\vecN^{\sankaku\,(0)},S^{(0)})$
is a split nilpotent orbit.
Here $\vecN^{\sankaku\,(0)}$ denotes the tuple
$(N^{\sankaku}_1,\ldots,N^{\sankaku}_l)$.
\end{cor}
\pf
It follows from Corollary \ref{cor;9.15.25}.
\hfill\qed

\begin{thm} \label{thm;10.26.110}
The tuples
 $(S^{\can}_{u}(E),W,\vecN^{\sankaku},S)$ and
 $(S_u(E,P),W,\vecN^{\sankaku},S)$ $(P\in X-D)$
are polarized mixed twistor structures of $(0,l)$-type.
Here $\vecN^{\sankaku}$ denotes the tuple
of nilpotent maps $\nbign_i^{\sankaku}$ $(i\in\lbar)$.
\end{thm}
\pf
It follows from Corollary \ref{cor;9.15.26}
and Lemma \ref{lem;a12.1.35}.
\hfill\qed

%% file: 22.tex

\subsubsection{Strong sequential compatibility on
$S^{\can}_{\vecu}(E)$ and $S_{\vecu}(E,P)$}

From a harmonic bundle $\harmonicbundle$
over $X-D$,
we obtain the polarized limiting mixed twistor structure
\[
 (S^{\can}_{\vecu}(E),\nbign^{\sankaku}_i,S),
\quad
 (S_{\vecu}(E,P),\nbign^{\sankaku}_i,S).
\]

\begin{cor} \label{cor;10.26.112}
The tuple of nilpotent maps
$(\nbign^{\sankaku}_1,\ldots,\nbign^{\sankaku}_n)$
on $S^{\can}_{\vecu}(E)$
or $S_{\vecu}(E,P)$
is strongly sequentially compatible.
\end{cor}
\pf
It follows from Theorem \ref{thm;10.26.110}
and Lemma \ref{lem;10.26.111}.
\hfill\qed

\begin{rem}
From the associated graded mixed twistor structure
$(\Gr^WS^{\can}_{\vecu}(E),W,\vecN,S)$,
we obtain the patched objects,
which gives the harmonic bundle.
In the case where $\harmonicbundle$ is nilpotent
and with trivial parabolic structure,
it is same as the limiting CVHS in our previous paper.
\hfill\qed
\end{rem}

\begin{rem}\label{rem;b12.5.1}
From the limiting mixed twistor structure
$(S^{\can}_{\vecu}(E),\nbign_i,S)$,
we obtain the variation of $\proj^1$-holomorphic bundle $\nbigv$.
We also obtain the pairing $\nbigv\otimes\sigma(\nbigv)\lrarr \Tate(0)$.
Does it give a variation of polarized pure twistor structures?
If it is true, it gives a partial generalization
of Schmid's nilpotent orbit theorem.

The replacement of a variation of Hodge structures
to the nilpotent orbit seems fundamental
in the study of Cattani-Kaplan-Schmid and Kashiwara-Kawai.
In this study, we do not use such replacement.
\hfill\qed
\end{rem}

\begin{cor} \label{cor;c11.18.1}
Let $\veca$ be an element of $\seisuu_{\geq \,0}^l$
On $\lefttop{\lbar}\nbigg_{\vecu\,|\,(\lambda,P)}$,
the conjugacy classes of $\prod_i\nbign^{a_i}_{i\,|\,(\lambda,P)}$
is independent of a choice of
$(\lambda,P)\in \cnum\times D_{\lbar}$.
\end{cor}
\pf
When we fix a point $P\in D_{\lbar}$,
the independence follows from the fact that
$\prod_i\nbign_i^{a_i}$ induces the morphisms of mixed twistor structures.
When we fix $\lambda\in\cnum^{\ast}$,
we can derive the independence by using the normalizing frame.
\hfill\qed

\vspace{.1in}
Similarly we can show the following:
\begin{lem}
Let $\veca$ and $\vecb$ be elements of $\seisuu_{\geq 0}^l$.
Then $\Image\prod_{i=1}\nbign_i^{a_i}\cap\Ker\prod_{i=1}\nbign_i^{b_i}$
gives the vector subbundle of $\lefttop{\lbar}\nbigg_{\vecu}$
on $\nbigd_{\lbar}$.
\hfill\qed
\end{lem}

\begin{cor}\label{cor;10.26.115}
On $\lefttop{\lbar}\nbigg_{\vecu}$,
the tuple of nilpotent maps $\nbign_1,\ldots,\nbign_l$
are strongly sequentially compatible.
\end{cor}
\pf
This is a direct corollary of
Corollary \ref{cor;10.26.112}.
\hfill\qed

%% file: 22.1.tex

\subsubsection{Sequential compatibility of the tuple
$\bigl(\nbign_i,
 \lefttop{i}\Fzero,
 \lefttop{i}\EEzero\,\big|\,
 i\in\lbar
 \bigr)$}

Let $\lambda_0$ be an element of $\cnum_{\lambda}$.
Let $\epsilon_0$ be a sufficiently small positive number.
Then we may assume that
$\prolongg{\vecb}{\nbige}$
over $\nbigx(\lambda_0,\epsilon_0)$ is locally free,
and that  we have $\lefttop{i}\Fzero$,
$\lefttop{i}\EEzero$ and $\Res_i(\DD)$
on $\prolong{\nbige}_{|\nbigd_i}$.
We have already known that
the tuple
$\bigl(
 \lefttop{i}\Fzero,\lefttop{i}\EEzero\,\big|\,
 i=1,\ldots,l
 \bigr)$ is compatible
in the sense of Definition \ref{df;10.11.90}.

\begin{lem}
Let $P$ be a point of 
$D_{\mbar}^{\circ}=
   D_{\mbar}-\bigcup_{i>m}(D_{i}\cap D_{\mbar})$.
\begin{itemize}
\item
For any $m_1\leq m$,
the tuple
$(\nbign_1,\ldots,\nbign_{m_1},\lefttop{m_1+1}F^{(\lambda_0)},
\ldots, \lefttop{m}F^{(\lambda_0)})$
is sequentially compatible
on the bundle $\lefttop{\mbariti}\Gr^{F(\lambda_0)}
 \lefttop{\mbariti}\EE^{(\lambda_0)}(
  \prolong{\nbige}_{|P\times \Delta(\lambda_0,\epsilon_0)},\vecbeta)$.
\item
Let $\veca$ be an element of $\seisuu_{\geq\,0}^{m_1}$.
The conjugacy classes of $\prod_i\nbign^{a_i}_{i\,|\,(\lambda,P)}$
on $\lefttop{\mbariti}\Gr^{\Fzero}
 \lefttop{\mbar}\EEzero\bigl(
 \prolong{\nbige}_{|(\lambda,P)}\bigr)$
are independent of a choice of $\lambda\in\cnum_{\lambda}$.
\end{itemize}
\end{lem}
\pf
We know that
$\lefttop{\mbar}\nbign_1$, $\ldots$,
$\lefttop{\mbar}\nbign_{m}$
are sequentially compatible
on
$\lefttop{\mbar}\Gr^{F^{(\lambda_0)}}
 \lefttop{\mbar}\EE^{(\lambda_0)}(
 \prolong{\nbige}_{|P\times\Delta(\lambda_0,\epsilon_0)},\vecbeta
 )$
and their conjugacy classes are independent of a choice of $\lambda$
(Corollary \ref{cor;c11.18.1} and Corollary  \ref{cor;10.26.115}).

On $\Delta^{\ast}(\lambda_0,\epsilon_0)$,
the vector bundle
$\lefttop{\mbariti}\Gr^{F^{(\lambda_0)}}
 \lefttop{\mbar}\EE^{(\lambda_0)}(
 \prolong{\nbige}_{|P\times\Delta(\lambda_0,\epsilon_0)},\vecbeta
 )$
is decomposed into the generalized eigen bundles
of the endomorphisms $\Res_i(\DD)$ $(i=m_1+1,\ldots,m)$.
It satisfies the following:
\begin{itemize}
\item
The decomposition gives the splitting
of the filtrations $\lefttop{i}F$ $(i=m_1+1,\ldots,m)$.
\item
The decomposition is compatible with
$\nbign_i$ $(i=1,\ldots,m_1)$.
\end{itemize}

For any $\lambda\in\Delta(\lambda_0,\epsilon_0)$,
let $\nbigr$ denote the local ring
of $\nbigo_{\Delta(\lambda_0,\epsilon_0)}$ at $\lambda$.
We put as follows:
\[
 V=\lefttop{\mbariti}\Gr^{F^{(\lambda_0)}}
   \lefttop{\mbar}\EE(
 \prolong{\nbige}_{|P\times\Delta(\lambda_0,\epsilon_0)}
 )
\otimes_{\nbigo_{\Delta(\lambda_0,\epsilon_0)}}\nbigr.
\]
Then we obtain the naturally induced filtrations
$\lefttop{i}F$ $(i=m_1+1,\ldots,m)$ and
the nilpotent maps
$\nbign_i$ $(i=1,\ldots,m_1)$.
Then we have only to apply 
Lemma \ref{lem;c11.12.20} and
Proposition \ref{prop;10.26.120}.
\hfill\qed

\vspace{.1in}
On $D_{\mbar}(\lambda_0,\epsilon_0)$,
we have the vector bundle
$\lefttop{\mbar}\Gr^{F^{(\lambda_0)}}_{\veca}
 \EE^{(\lambda_0)}(\prolong{\nbige}_{|\nbigd_{\mbar}},\vecbeta)$
and
$\nbign_1,\ldots,\nbign_m$.
\begin{lem}\label{lem;10.26.121}\mbox{{}}
\begin{itemize}
\item
Let $\veca$ be an element of $\seisuu_{\geq\,0}^l$.
The conjugacy classes of $\prod_i\nbign^{a_i}_{i\,|\,(\lambda,P)}$ are 
independent of a choice of $(\lambda,P)\in\nbigd_{\mbar}$.
\item
The tuple of the nilpotent maps
$\nbign_1,\ldots,\nbign_m$ is 
sequentially compatible.
\end{itemize}
\end{lem}
\pf
Let $P$ be a point of $D_{\mbar}$.
We have already seen the result
on $\Delta(\lambda_0,\epsilon_0)\times \{P\}$.

We put
$I_{\vech\,|\,(\lambda,P)}
:=\bigcap_{j=1}^m W(\jbar)_{h_j\,|\,(\lambda,P)}$,
and then we have only to show that
$\bigl\{I_{\vech\,|\,(\lambda,P)}\,
 \big|\,(\lambda,P)\in \nbigd_{\mbar}
 \bigr\}$
forms a vector bundle.
For that purpose,
we have only to show
that the ranks of $I_{\vech\,|\,(\lambda,P)} $
are independent of $(\lambda,P)$.

If we fix $P$,
then they are independent of a choice of $\lambda$
due to the previous lemma.

Let pick a generic $\lambda$.
Then we have a normalizing frame for $\prolong{\nbige^{\lambda}}$
Due to the normalizing frame,
we obtain the isomorphism for any $P_1,P_2\in D_{\mbar}$:
\[
 \bigl(\prolong{\nbigelambda}_{|P_1},
 \Res_i(\DD^{\lambda})\,
 \big|\,i\in\mbar
 \bigr)
 \simeq
 \bigl(\prolong{\nbigelambda}_{|P_2},
 \Res_i(\DD^{\lambda})\,\big|\,i\in\mbar
 \bigr).
\]
Thus we are done.
\hfill\qed

\vspace{.1in}
We obtain the tuple 
$\bigl(\nbign_i,
 \lefttop{i}\Fzero,
 \lefttop{i}\EEzero\,\big|\,
 i\in\lbar
 \bigr)$
as in the subsubsection \ref{subsubsection;10.26.123}.

\begin{thm}
The tuple
$\bigl(\nbign_i,
 \lefttop{i}\Fzero,
 \lefttop{i}\EEzero\,\big|\,
 i\in\lbar
 \bigr)$
is sequentially compatible.
\end{thm}
\pf
It follows from Lemma \ref{lem;10.26.121}.
\hfill\qed

%% file: a33.2.tex

\subsubsection{Decomposition}

First we have a remark.
\begin{rem}\label{rem;9.23.12}
Let $P$ be a point of $D_I^{\circ}$.
Let $\vecu$ be an element of $\KMS(\nbige^0,I)$.
Let $q_I:X\lrarr D_I$ be the projection.
By considering the restriction of
$\harmonicbundle$ to $q_I^{-1}(P)$,
we obtain the Pol-MTS $(S^{\can}_{\vecu},W,\vecN,S)$
of type $(0,|I|)$.
Note that we can apply Proposition {\rm\ref{prop;9.23.10}}
to the tuple $(S^{\can}_{\vecu},W,\vecN,S)$.
\end{rem}

Let us consider the vector bundle
$\lefttop{I}\nbigg_{\vecu}(\nbige)$
and the nilpotent maps
$\nbign_i$ $(i\in I)$ on $\nbigd_I$.
Let $J$ and $K$ be subsets of $I$
such that $J\cap K=\emptyset$.
We put $\nbign_K:=\prod_{k\in K}\nbign_k$.
Due to the limiting mixed twistor theorem,
the conjugacy classes of $\nbign_{K\,|\,(\lambda,P)}$
are independent of $(\lambda,P)\in \nbigd_I$.
Thus we obtain the vector bundle
$\Image(\nbign_K)$.
For any element $i\in J$, we put $J':=J-\{i\}$.
We have the morphisms
$\var_i:\nbigv_J(\Image\nbign_K)\lrarr\nbigv_{J'}(\Image\nbign_K)$
and $\can_i:\nbigv_{J'}(\Image\nbign_K)\lrarr\nbigv_{J}(\Image\nbign_K)$
For any element $i\in K$, we put $K':=K-\{i\}$.
We have the morphisms
$\var_i:\nbigv_J(\Image\nbign_K)\lrarr\nbigv_{J}(\Image\nbign_{K'})$
and $\can_i:\nbigv_{J}(\Image\nbign_{K'})\lrarr\nbigv_{J}(\Image\nbign_K)$
(See the subsubsection \ref{subsubsection;9.23.11}).

\begin{lem}
Let $P$ be a point of $D_i^{\circ}$.
Then we have the following decomposition:
\[
 \Gr^{W(N)}(\nbigv_J(\Image(\nbign_K)))_{|\cnum_{\lambda}\times\{P\}}
=\Image \can_{i\,|\,\cnum_{\lambda}\times\{P\}}
\oplus
 \Ker \var_{i\,|\,\cnum_{\lambda}\times\{P\}}.
\]
\end{lem}
\pf
We have only to apply the limiting mixed twistor theorem
and Proposition \ref{prop;9.23.10}.
(See Remark (\ref{rem;9.23.12})).
\hfill\qed

\begin{lem} \label{lem;9.23.16}
Let $P$ be a point of $D_I$.
\begin{itemize}
\item
 The numbers $\dim\Image(\can_{i\,|\,(\lambda,P)})$
and $\dim\Ker(\var_{i\,|\,(\lambda,P)})$
are independent of $\lambda$.
\item
We have the decomposition
$\Gr^{W(N)}(\nbigv_J(\Image\nbign_K))_{|(P,\lambda)}
=\Image \can_{i\,|\,\cnum_{\lambda}\times\{P\}}\oplus 
 \ker \var_{i\,|\,\cnum_{\lambda}\times\{P\}}$.
\end{itemize}
\end{lem}
\pf
We have the subset $I_0\subset \lbar$ such that
$P\in D_K^{\circ}$.
Let $\vecu_1$ be any element of $\KMS(\nbige^0,I_0)$.
The ranks of $\can_i$ and $\var_i$
on $\lefttop{I_0}\nbigg_{\vecu_1}$ are independent of $\lambda$
due to the limiting mixed twistor theorem,
and we have the decomposition
of $\lefttop{I_0}\nbigg_{\vecu_1}$.
Then we have only to apply Lemma \ref{lem;9.23.15}.
\hfill\qed

\begin{prop} \label{prop;9.23.21}
We have the decomposition of the vector bundle:
\[
 \Gr^{W(N)}(\nbigv_J(\Image\nbign_K))
=\Image (\can_i)\oplus\Ker (\var_i).
\]
\end{prop}
\pf
We have only to show that
$\Image (\can_i)$ and $\Ker (\var_i)$ are vector bundles.
Namely we have only to show
the ranks of
the morphisms
$\can_{i\,|\,(\lambda,P)}$
and $\var_{i\,|\,(\lambda,P)}$
are independent of $(\lambda,P)\in \nbigd_i$.
When we fix a point $P$, it follows from Lemma \ref{lem;9.23.16}.
When we pick a generic $\lambda$,
we can show the numbers are independent of $P$
by using the normalizing frame.
(See the last part of the proof of Lemma \ref{lem;10.26.121}.)
\hfill\qed

%% file: a84.1.tex

\subsubsection{The induced polarized pure twistor structure }
\label{subsubsectio;b12.3.100}

Let us consider the case $l=n=|I|$.
Due to the limiting mixed twistor theorem
and Kashiwara's lemma (Corollary \ref{cor;9.23.7}),
we obtain the polarized mixed twistor structure
of type $(n-1,n+1)$:
\[
 \bigl(
 \nbigv_{\nbar}(S^{\can}_{\vecu}),
 W,\nbigvecn,\nbigv_{\nbar}(S)
 \bigr).
\]

\begin{lem}\label{lem;b12.3.170}
We put as follows:
\[
 \nbigc_{\nbar}:=
 \Ker\Bigl(
 P_h\Gr^{W(N)}_h\nbigv_{\nbar}(S^{\can}_{\vecu}(E))
\lrarr
 \bigoplus_{|I|=n-1}
 P_h\Gr^{W(N)}_h\nbigv_I(S^{\can}_{\vecu}(E))
 \Bigr).
\]
Then $\nbigc_{\nbar}$ is pure twistor of weight $h+n-1$.
The pairing $\nbigv_{\nbar}(S)(N^h\otimes\id)$
gives the polarization of $\nbigc_{\nbar}$.
\end{lem}
\pf
It follows from Saito's lemma (Lemma \ref{lem;9.26.20}).
\hfill\qed

%% file: 25.tex

\subsubsection{The $\lambda$-connection form for pull back}
\label{subsubsection;10.27.9}

We put $X=\Delta^n_{\zeta}$,
and $D=\bigcup_{i=1}^l D_i$.
Let $\harmonicbundle$ be a tame harmonic bundle over $X-D$.
Let $\lambda_0$ be a point
and $\epsilon_0$ be a sufficiently small number.
We assume that
we have $(\prolong{\nbige},\DD)$
over $\nbigx(\lambda_0,\epsilon_0)$.

Let $\vecv=(v_i)$ be a holomorphic frame of $\prolong{\nbige}$
compatible with $\Fzero$ and $\EEzero$.
For each $v_i$,
we have the element $\vecu(v_i)\in\KMS(\nbige^0,\lbar)$
such that
$\deg^{\EEzero,\Fzero}(v_i)=\kmsmap(\lambda_0,\vecu(v_i))$.
We put as follows:
\[
 v_i':=v_i\cdot
 \prod_{j=1}^l|\zeta_j|^{\paramap(\lambda_0,\vecu(v_i))},
\quad\quad
 \vecv'=(v_i').
\]

We pick any element
$\vecc:=\bigl(c_{j\,i}\,|\,j\in\mbar,i\in\nbar\bigr)
  \in \seisuu_{\geq 0}^{m\cdot n}$.
We put $\tilde{X}=\Delta^m$.
Then we have the morphism
$\phi_{\vecc}:\tilde{X}\lrarr X$
determined as follows:
\[
 \phi_{\vecc}^{\ast}\zeta_i
=\prod_{j=1}^m z_j^{c_{j\,i}}.
\]
Then we obtain the $C^{\infty}$-frame
$\phi_{\vecc}^{\ast}\vecv'$ of
$\phi_{\vecc}^{\ast}\nbige$ over $\tilde{X}-\tilde{D}$,
which is adapted up to log order.
We have the following equalities:
\[
 \phi_{\vecc}^{\ast}(v_i')
=\phi_{\vecc}^{\ast}(v_i)
 \cdot
 \prod_{j=1}^l\Bigl(
 \prod_{h=1}^m |z_h|^{c_{h\,j}}
 \Bigr)^{\paramap(\lambda,u_j(v_i)}
=\phi_{\vecc}^{\ast}(v_i)
 \cdot
 \prod_{h=1}^m|z_h|^{\paramap\bigl(\lambda,\sum c_{h\,j}\cdot u_j(v_i)\bigr)}
=\phi_{\vecc}^{\ast}(v_i)\times
 \prod_{h=1}^m
 |z_h|^{\paramap\bigl(\lambda,\vecc_h\cdot \vecu(v_i)\bigr)}.
\]
Here we put $\vecc_h\cdot \vecu(v_i)=\sum c_{h\,j}\cdot u_j(v_i)$.

For some small positive number $\epsilon_0'$
such that $\epsilon_0'<\epsilon_0$,
we may assume that
$\nu\bigl(\paramap\bigl(\lambda,\vecc_h\cdot\vecu(v_i)\bigr)\bigr)$
does not independent of $\lambda\in\Delta(\lambda_0,\epsilon'_0)$
for any $v_i$ and for any $h$.
We put
$\nu(v_i,h)
:=\nu\bigl(\paramap\bigl(\lambda,\vecc_h\cdot\vecu(v_i)\bigr)\bigr)
\in\seisuu$.
(See the subsubsection \ref{subsubsection;b11.11.10}).
We put $\kappa(v_i,h):=
 \kappa\bigl(\paramap\bigl(\lambda,\vecc_h\cdot\vecu(v_i)\bigr)\bigr)$.
Then we put as follows:
\[
 u_i:=\phi_{\vecc}^{\ast}v_i
\cdot \prod_{h=1}^m z_h^{\nu(v_i,h)}.
\]
Then $u_i$ is a holomorphic section 
of $\prolong{\phi_{\vecc}^{\ast}\nbige}$.
We put as follows:
\[
 u_i':=u_i\cdot
 \prod_{h=1}^m |z_h|^{\kappa(v_i,h)},
\quad\quad
 \vecu':=(u_i').
\]
Then $\vecu'$ is adapted up to log order
over $\nbigx(\lambda_0,\epsilon_0')-\nbigd(\lambda_0,\epsilon_0')$,
by our construction.
Hence $\vecu=(u_i)$ is a holomorphic frame 
of $\prolong{\phi_{\vecc}^{\ast}\nbige}$,
which is compatible with the filtrations $F^{(\lambda_0)}$,
due to Lemma \ref{lem;9.5.10} and Lemma \ref{lem;10.10.31}.
We also have the following:
\begin{equation}\label{eq;10.27.10}
\lefttop{h}\deg^{F^{(\lambda_0)}}(u_i)
=\kappa\bigl(v_i,h\bigr).
\end{equation}

Let $A=\sum_k A^k\cdot\zeta^{-1}_k\cdot d\zeta_k$
denote the $\lambda$-connection form
of $\DD$ with respect to the frame $\vecv$,
i.e.,
$\DD\vecv=\vecv\cdot A$ holds.
We put $\tilde{\DD}:=\phi_{\vecc}^{\ast}\DD^{\ast}$.
Then we have the following equalities:
\[
 \tilde{\DD}\phi_{\vecc}^{\ast}v_i
=\sum_j \phi_{\vecc}^{\ast}v_j\cdot
 \Bigl(
  \sum_k \phi_{\vecc}^{\ast}A_{j\,i}^k
 \cdot\phi_{\vecc}^{\ast}\frac{d\zeta_k}{\zeta_k}
 \Bigr)
=\sum_j\phi_{\vecc}^{\ast}v_j
\cdot
 \Bigl(
 \sum_{h}\frac{dz_h}{z_h}
 \cdot \sum_k c_{h\,k}\cdot\phi_{\vecc}^{\ast}A^k_{j\,i}
 \Bigr).
\]
Thus we obtain the following:
\begin{multline} \label{eq;10.26.130}
 \tilde{\DD}u_i
=\tilde{\DD}\Bigl(
 \phi_{\vecc}^{\ast}v_i
 \cdot
 \prod_{h=1}^m z_h^{
 \nu(v_i,h)}
 \Bigr)
=\tilde{\DD}\bigl(
 \phi_{\vecc}^{\ast}v_i\bigr)
 \cdot \prod_{h=1}^m z_h^{\nu(v_i,h)}
+\phi_{\vecc}^{\ast}v_i\cdot
 \prod_{h=1}^m z_h^{\nu(v_i,h)}
\cdot
 \Bigl(
 \sum_{h=1}^m \nu(v_i,h)
 \cdot\frac{dz_h}{z_h}
 \Bigr) \\
=\sum_{j}u_j\cdot
 \Bigl(
 \sum_{h=1}^m \frac{dz_h}{z_h}
  \sum_{k=1}^n c_{h\,k}\cdot\phi^{\ast}A^k_{j\,i}
 \Bigr)\cdot\prod_{h=1}^m
 z^{\nu(v_i,h)-\nu(v_j,h)}
+u_i\cdot \Bigl(
 \sum_{h=1}^m \nu(v_i,h)\cdot
 \frac{dz_h}{z_h}
 \Bigr).
\end{multline}

\subsubsection{Diagonal case}
\label{subsubsection;10.27.65}

Let us consider the case that
$\vecc$ is a diagonal matrix whose $i$-th diagonal component is $c_i$,
that is,
$\varphi_{\vecc}^{\ast}\zeta_i=z_i^{c_i}$.
Then we obtain the following formula from (\ref{eq;10.26.130}):
\begin{multline}
 \tilde{\DD}u_i
=\sum_j u_j\cdot
 \Bigl(
 \sum_{h=1}^m\frac{dz_h}{z_h}\cdot c_h
 \cdot\phi_{\vecc}^{\ast}A^h_{j\,i}
 \Bigr)\cdot
 \prod_{p=1}^m
 z^{\nu(c_p\cdot b_p(v_i))-\nu(c_p\cdot b_p(v_j))}
+u_i\cdot\sum_{h=1}^m \nu(c_h\cdot b_h(v_i))\cdot\frac{dz_h}{z_h}
 \\
=:\sum u_j\cdot\tilde{A}^h_{j\,i}\frac{dz_h}{z_h}.
\end{multline}
Here we put $b_p(v_i):=\paramap(\lambda_0,u_p(v_i))$.

\begin{lem}\mbox{{}}
We have the following vanishings:
\begin{enumerate}
\item
 In the case $\lefttop{k}\deg^{F^{(\lambda_0)}}(v_i)<
 \lefttop{k}\deg^{F^{(\lambda_0)}}(v_j)$,
 we have $\tilde{A}^h_{j\,i\,|\,\nbigd_k}=0$,
 for any $h$.
\item
 If $c_k$ is sufficiently large
 and if
 the inequality $\lefttop{k}\deg^{F^(\lambda_0)}(v_i)>
 \lefttop{k}\deg^{F^{(\lambda_0)}}(v_j) $ holds,
then we have $\tilde{A}^h_{j\,i\,|\,\nbigd_k}=0$.
\end{enumerate}
\end{lem}
\pf
In the case $\lefttop{k}\deg^{\Fzero}(v_i)<\lefttop{k}\deg^{\Fzero}(v_j)$,
we have $A^h_{j\,i\,|\,\nbigd_k}=0$.
Hence $A_{j\,i}^h\cdot \zeta_k^{-1}$ is holomorphic,
and thus
$\phi_{\vecc}^{\ast}A_{j\,i}^h\cdot z_k^{-c_k}$
is holomorphic.
Since we have the following inequality:
$ -c_k<\nu\bigl(c_k\cdot b_k(v_i)\bigr)
 -\nu\bigl(c_k\cdot b_k(v_j)\bigr)$,
we obtain the first claim.

Let us show the second claim.
If $c_k$ is sufficiently large
and the inequality $\lefttop{k}\deg^{F^(\lambda_0)}(v_i)>
 \lefttop{k}\deg^{F^{(\lambda_0)}}(v_j) $ holds,
then we obtain the inequality
$\nu(c_k\cdot b_k(v_i))>\nu(c_k\cdot b_k(v_j))$.
Since $A_{j\,i}^k$ is holomorphic,
we obtain the result.
\hfill\qed

\vspace{.1in}

Assume that $c_k$ $(k=1,\ldots,l)$ are sufficiently large.
Then we have the following formula:
\begin{multline}\label{eq;10.27.1}
 \Res_k(\tilde{\DD})u_i
=\sum u_j\cdot\tilde{A}^k_{j\,i\,|\,\nbigd_k}
=u_i\cdot \nu(c_k\cdot b_k(v_i))+ \\
\sum_{\lefttop{k}\deg^{\Fzero}(v_i)=\lefttop{k}\deg^{\Fzero}(v_j)}
 \!\!\!\!
 u_j\cdot c_k\cdot
 \phi_{\vecc,k}^{\ast}(
 A_{j\,i\,|\,\nbigd_k}^k
 )\cdot
 \prod_{a\neq k}
 z_a^{\nu(c_a\cdot b_a(v_i))-\nu(c_a\cdot b_a(v_j))}.
\end{multline}
Here
$\phi_{\vecc,k}:\tilde{\nbigd}_k\lrarr\nbigd_k$
denotes the restriction of $\phi$ to $\tilde{\nbigd}_k$,
given by
$\phi_{\vecc,k}^{\ast}\zeta_i=z_i^{c_i}$
for $i\neq k$.

We obtain the $\lambda$-connection
$\lefttop{k}\tilde{\DD}$
of
$\prolong{\psi_{\vecc}^{-1}\nbige}_{|\nbigd_k(\lambda_0,\epsilon_0')}$,
given as follows:
\begin{multline}\label{eq;10.27.2}
 \lefttop{k}\tilde{\DD}u_i
=\sum_j u_j\cdot
 \Bigl(
 \sum_{h\neq k}\tilde{A}^h_{j\,i\,|\,\nbigd_k}
 \cdot\frac{dz_h}{z_h}
 \Bigr) \\
=\sum_{\lefttop{k}\deg(v_i)=\lefttop{k}\deg(v_j)}
 u_j\cdot
 \Bigl(
 \sum_{h\neq k}\frac{dz_h}{z_h}
 \cdot c_h\cdot \phi_{\vecc,k}^{\ast}
 A^h_{j\,i\,|\,\nbigd_k}
 \Bigr)\cdot
 \prod_{p=1}^l z_p^{\nu(c_p\cdot b_p(v_i))-\nu(c_p\cdot b_p(v_j)) }
+u_i\cdot \sum_{h\neq k}
  \nu(c_h\cdot b_h(v_i))\frac{dz_h}{z_h}.
\end{multline}

Let $\lefttop{k}\nbigk_{u}$
denote the vector subbundle of 
$\prolong{\phi_{\vecc}^{\ast}\nbige}_{|\nbigd_k}$,
generated by
$\vecu_u:=
 \bigl\{u_i\,\big|\,\deg^{\Fzero,\EEzero}(v_i)=\kmsmap(\lambda_0,u)\bigr\}$.
Then we have the following decomposition $\lefttop{k}\nbigk$
for $k=1,\ldots,l$:
\[
\prolong{\phi_{\vecc}^{\ast}\nbige}_{|\nbigd_k}
=\bigoplus_u \lefttop{k}\nbigk_{u}.
\]

\begin{lem}\mbox{{}} \label{lem;10.27.5}
\begin{itemize}
\item
 The vector subbundle $\lefttop{k}\nbigk_u$ is preserved
 by the residue $\Res_k(\tilde{\DD})$
 and the induced $\lambda$-connection $\lefttop{k}\DD$.
\item
 We have the following formula:
 \[
  \Res_k(\tilde{\DD})_{|\lefttop{k}\nbigk_u}(\vecu_u)
=\vecu_u\cdot
 \Bigl(
 \nu\bigl(c_k\cdot \paramap(\lambda_0,u)\bigr)
+\phi_{\vecc,k}^{\ast}\nbigr_u\cdot
 \prod_{p\neq k}z_p^{\nu(c_p\cdot b_p(v_i))-\nu(c_p\cdot b_p(v_j))}
 \Bigr).
 \]
Here $\nbigr_u$ denotes the representation matrix
of $\Res_k(\DD)$ on
the vector subbundle
$\langle v_{i}\,|\,\lefttop{k}\deg^{\FEzero}(v_i)=\kmsmap(\lambda_0,u)
 \rangle \subset \prolong{\nbige}_{|\nbigd_k(\lambda_0,\epsilon_0')}$
with respect to the frame
$\vecv_u=
 \{ v_{i}\,|\,\lefttop{k}\deg^{\FEzero}(v_i)=\kmsmap(\lambda_0,u) \}$.
\item
In particular,
the eigenvalue of $\Res_k(\widetilde{\DD})_{|\lefttop{k}\nbigk_u}$
is $\nu(c_k\cdot \paramap(\lambda_0,u))+\eigenmap(\lambda,u)$.
\end{itemize}
\end{lem}
\pf
The claims immediately follow from
(\ref{eq;10.27.1}) and (\ref{eq;10.27.2}).
\hfill\qed

\begin{lem}
The decompositions $\lefttop{k}\nbigk$ $(k=1,\ldots,l)$ are
compatible
in the sense of Definition {\rm\ref{df;10.11.90}}.
\end{lem}
\pf
It immediately follows from our construction.
\hfill\qed

\begin{lem}
We have the following decomposition
on $\nbigd_k(\lambda_0,\epsilon_0)$:
\begin{equation} \label{eq;10.27.3}
  \lefttop{k}\Fzero_b
 \prolong{\phi_{\vecc}^{\ast}\nbige}_{
  |\nbigd_k(\lambda_0,\epsilon_0)}
=\bigoplus_{\kappa(c_k\cdot \paramap(\lambda_0,u))\leq b}
 \lefttop{k}\nbigk_{u}.
\end{equation}
Namely the decompositions $(\lefttop{k}\nbigk\,|\,k=1,\ldots,l)$
gives the splitting of the filtrations
$\bigl(\lefttop{k}\Fzero\,\big|\,k=1,\ldots l\bigr)$.
\hfill\qed
\end{lem}

\begin{lem}\label{lem;10.27.60}
When $c_k$ $(k=1,\ldots,l)$ are sufficiently large,
$\phi^{\ast}\nbige$ is CA in the sense of 
Definition {\rm\ref{df;10.27.4}} below.
\end{lem}
\pf
It immediately follows from the second claim
of Lemma \ref{lem;10.27.5}.
\hfill\qed

\subsubsection{Convenient}

Let $\lambda_0$ be a point of $\cnum_{\lambda}$,
and $\epsilon_0$ be a sufficiently small positive number.
Recall that
we always have the following decomposition on
$\nbigd_i^{\ast}(\lambda_0,\epsilon_0)$:
\[
 \prolong{\nbige}_{|\nbigd_i^{\ast}(\lambda_0,\epsilon_0)}
=
 \bigoplus
 \EE\bigl(\prolong{\nbige}_{|\nbigd_i^{\ast}(\lambda_0,\epsilon_0)},
 \eigenmap(\lambda,u)\bigr)
 ).
\]

\begin{df} \label{df;10.27.4}
$\nbige$ is called convenient at $\lambda_0$,
if the following holds:
\begin{description}
\item[$(A)$]
The decomposition above is prolonged to
the decomposition
$\prolong{\nbige}_{|\nbigd_k(\lambda_0,\epsilon_0)}
=\bigoplus \lefttop{k}\nbigk_u$.
Moreover the tuple of the decompositions
$\bigl(\lefttop{k}\nbigk\,\big|\,k=1,\ldots,l\bigr)$ 
is compatible.

\item[$(B)$]
There exists a sequence of positive numbers $\eta_1>\eta_2\cdots>\eta_l>0$:
\begin{itemize}
\item
 $\sum \eta_i<1/2$
\item
 $\Par(\prolong{\nbige},i)$ is $\eta_i$-small.
\item
$\min\bigl\{|a-b|\,\big|\,a\neq b\in\Par(\prolong{\nbige^0},i)\cup\{0\}
  \bigr\}
 >2\cdot \sum_{j> i} \eta_j$.
\end{itemize}
 
\end{description}

$\nbige$ is called $CA$ (resp. $CB$)
if the condition $(A)$ (resp. $(B)$) holds.
\hfill\qed
\end{df}

Note that the decompositions $\lefttop{k}\nbigk$
is uniquely determined if $\nbige$ is convenient at $\lambda_0$.

\begin{lem} \label{lem;10.27.8}
Assume that $\nbige$ is CA at $\lambda_0$.
Then the induced connection
$\lefttop{i}\DD$ and the residue $\Res_i(\DD)$
preserve the decomposition $\lefttop{i}\nbigk$.
\end{lem}
\pf
Since the restrictions of 
$\lefttop{i}\DD$ and $\Res_i(\DD)$
to $\nbigd_i^{\ast}(\lambda_0,\epsilon_0)$
preserve
$\lefttop{i}\nbigk_{u\,|\,\nbigd_i^{\ast}(\lambda_0,\epsilon_0)}$,
they preserve $\lefttop{i}\nbigk_u$
on $\nbigd(\lambda_0,\epsilon_0)$.
\hfill\qed

\begin{lem} \label{lem;10.27.7}
Assume that $\nbige$ is CA at $\lambda_0$.
We have the following decomposition:
\begin{equation} \label{eq;10.27.6}
 \lefttop{i}\Fzero_b\EE^{(\lambda_0)}
 \bigl(\prolong{\nbige}_{|\nbigd_i(\lambda_0,\epsilon_0)},\beta
 \bigr)
=\bigoplus_{\substack{
 \eigenmap(\lambda_0,u)=\beta \\
 -1<\paramap(\lambda_0,u)\leq b
 }  }
 \lefttop{i}\nbigk_u.
\end{equation}
\end{lem}
\pf
By definition, the restrictions of the both sides of
(\ref{eq;10.27.6}) to $\nbigd^{\ast}_i(\lambda_0,\epsilon_0)$ are same.
Then Lemma \ref{lem;10.27.7} immediately follows.
\hfill\qed

\vspace{.1in}

Hence, if $\nbige$ is CA at $\lambda_0$,
we obtain the naturally defined isomorphism:
\[
 \lefttop{i}\nbigk_u\simeq
  \lefttop{i}\Gr^{\Fzero}_{\paramap(\lambda_0,u)}
 \lefttop{i}
 \EEzero\bigl(\prolong{\nbige}_{|\nbigd_i(\lambda_0,\epsilon_0)},
 \eigenmap(\lambda_0,u)\bigr).
\]
Thus the nilpotent parts $\nbign_i$ of the residues $\Res_i(\DD)$
are well defined as elements of
$\End\bigl(\prolong{\nbige}_{|\nbigd_i(\lambda_0,\epsilon_0)}\bigr)$
for any $i=1,\ldots,l$.

\begin{lem}
Assume that $\nbige$ is CA at $\lambda_0$.
Then the tuple
$\bigl(\lefttop{k}\nbigk_u,\nbign_i\,\big|\,k\in\lbar,i\in\lbar\bigr)$
is strongly sequentially compatible
in the sense of Definition {\rm\ref{df;b11.12.25}}.
\hfill\qed
\end{lem}

Let $\vecv$ be a frame of $\prolong{\nbige}$,
which is compatible with
the decompositions $(\lefttop{i}\nbigk\,|\,i=1,\ldots,l)$.
Let $\sum A^h\cdot z_h^{-1}dz_h$ denote
the $\lambda$-connection form
of $\DD$ with respect to the frame $\vecv$:
\[
 \DD\vecv=
 \vecv\cdot
 \Bigl(
 \sum_h A^h\cdot\frac{dz_h}{z_h}
 \Bigr).
\]
\begin{lem}
Assume that $\nbige$ is CA at $\lambda_0$.
In the case
$\lefttop{k}\deg^{\nbigk}(v_i)\neq \lefttop{k}\deg^{\nbigk}(v_j)$,
we have $A^h_{i\,j\,|\,\nbigd_k}=0$.
\end{lem}
\pf
It follows from Lemma \ref{lem;10.27.8}.
\hfill\qed

\vspace{.1in}

Assume that $\nbige$ is CA at $\lambda_0$,
and $\lefttop{k}\nbigk$ $(k=1,\ldots,l)$ be the decompositions
as in Definition \ref{df;10.27.4}.
Then we put as follows, for any subset $I\subset\lbar$:
\[
 \lefttop{I}\nbigk_{\vecu}:=
 \bigcap_{i\in I}\lefttop{i}\nbigk_{q_i(\vecu)}.
\]
Then we obtain the decomposition of
$\prolong{\nbige}_{|\nbigd_I(\lambda_0,\epsilon_0)}
=\bigoplus\lefttop{I}\nbigk_{\vecu}$.

%% file: 25.1.tex

\subsubsection{Functoriality for some pull back $(I)$}
\label{subsubsection;10.27.15}

Let us return to the functoriality.
Let $\eta$ be a positive number such that
$(1+\eta)^{l-1}\cdot 2/3<1$.
We put
$\tilde{X}:=\Delta(1+\eta)^{l-1}\times\Delta(2/3)\times\Delta^{n-l}$.
Let us consider the morphism $\phi:\tilde{X}\lrarr X$
given as follows:
\[
 \phi^{\ast}(\zeta_i)=
\left\{
\begin{array}{ll}
 \prod_{j=i}^lz_j & (i\leq l)\\
 \mbox{{}}\\
 z_i & (i>l).
\end{array}
\right.
\]
We put $\tilde{D}:=\phi^{-1}(D)$.

\begin{lem}
Assume that $\nbige$ is convenient at $\lambda_0$.
Then we have the natural isomorphism
$\prolong{\bigl(
 \phi^{\ast}\nbige\bigr)}
 \simeq
 \phi^{\ast}\bigl(
 \prolong{\nbige}\bigr)$.

Let $\vecv$ be a frame of $\prolong{\nbige}$
which is compatible with the decompositions
$\bigl(\lefttop{i}\nbigk\,\big|\,i\in\lbar\bigr)$.
In this case,
$u_i=\phi^{\ast}v_i$ holds
(see the subsubsection {\rm\ref{subsubsection;10.27.9}}),
and the degrees $\lefttop{h}\deg^{\Fzero}(u_i)$ are given as follows:
\[
 \lefttop{h}\deg^{\Fzero}(u_i)=
 \sum_{k\leq h}\lefttop{k}\deg^{\Fzero}(v_i).
\]
\end{lem}
\pf
If follows from the formula (\ref{eq;10.27.10})
for the degree of $u_i$.
\hfill\qed

\vspace{.1in}

Let $\sum_k A^k\cdot d\zeta_k/\zeta_k$
be the $\lambda$-connection form of $\DD$
with respect to the frame $\vecv$,
i.e.,
$\DD\vecv=\vecv\cdot(\sum_k A^k\cdot d\zeta_k/\zeta_k)$
holds.
Then we have the following:
\[
 \tilde{\DD}\vecu
=\vecu\cdot\sum_h\Bigl(
 \sum_{k\leq h}
 \phi^{\ast}A^k
 \Bigr)\cdot\frac{dz_k}{z_k}
=:\vecu\cdot\sum_h \tilde{A}^h\cdot\frac{dz_k}{z_k}.
\]
Here we have the following:
\[
 \tilde{A}^h=
\left\{
\begin{array}{ll}
 \sum_{k=1}^h\phi^{\ast}A^k & (h\leq l)\\
 \mbox{{}}\\
  A^h & (h>l).
\end{array}
\right.
\]
In particular, we obtain the following:
\[
 \tilde{A}^h_{|\tilde{\nbigd}_i(\lambda_0,\epsilon_0)}
=\sum_{k\leq h}
 \phi_i^{\ast}A^k_{|\nbigd_{\ibar}(\lambda_0,\epsilon_0)}.
\]
Here $\phi_{i}$ denotes the morphism
$\tilde{\nbigd}_i\lrarr \nbigd_{\ibar}$
is given as follows:
\[
 (\lambda,z_1,\ldots,z_{i-1},0,z_{i+1},\ldots,z_n)
\longmapsto
 \bigl(\lambda,\mbox{$\prod_{j=i+1}^lz_{j}$},
 \,\ldots,\,z_{l-1}\cdot z_l,\,z_l,\,\ldots,z_n\bigr).
\]
In particular,
we have the following formula:
\[
 \Res_i(\tilde{\DD})
=\sum_{k\leq i}
 \phi_i^{\ast}\Res_k(\DD)_{|\nbigd_{\ibar}}.
\]
We obtain the decomposition
of the vector bundle
$ \prolong{\phi^{\ast}\nbige}_{|\tilde{\nbigd}_i(\lambda_0,\epsilon_0)}
=\bigoplus
 \phi^{\ast}\lefttop{\ibar}\nbigk_{\vecu}$.
We put as follows:
\[
 \lefttop{i}\tilde{\nbigk}_{u}
=\bigoplus_{\phi_i^{\ast}(\vecu)=u}
 \phi^{\ast}\lefttop{\ibar}\nbigk_{\vecu}.
\]
Here we put $\phi^{\ast}_i(\vecu):=\sum_{j\leq i} q_j(\vecu)$:
Then we obtain the decomposition:
$ \prolong{\phi^{\ast}\nbige}_{|\tilde{\nbigd}_i(\lambda_0,\epsilon_0)}
=\bigoplus
 \lefttop{\ibar}\tilde{\nbigk}_{u}$.
Clearly the tuple of the decompositions
$\bigl(\lefttop{i}\tilde{\nbigk}\,\big|\,i=1,\ldots,l\bigr)$
is compatible.

\begin{lem}
The eigenvalue of $\Res_i(\tilde{\DD})$ on
$\lefttop{i}\tilde{\nbigk}_{u}$ is 
$\eigenmap(\lambda,u)$.
In particular,
$\phi^{\ast}\nbige$ is CA.
\end{lem}
\pf
It immediately follows from our construction.
\hfill\qed

\vspace{.1in}
The following lemma is also seen easily from our construction.
\begin{lem}
The nilpotent part 
$\tilde{\nbign}_i$ of the residue $\Res_i(\tilde{\DD})$
is given by
$\phi_i^{\ast}\nbign(\ibar)
=\sum_{j\leq i}\phi_i^{\ast}\nbign_j$.
\hfill\qed
\end{lem}

%% file: 25.2.tex

\subsubsection{Functoriality for some pull back $(II)$}

\label{subsubsection;10.27.36}

We put $\tilde{X}=\Delta_{z}^{n-1}$.
Let us consider the morphism
$\phi:\tilde{X}\lrarr X$,
given as follows:
\[
 \phi^{\ast}(\zeta_i)=
 \left\{
 \begin{array}{ll}
 z_1 & (i=1,2)\\
 z_{i-1}      & (i\geq 3).
 \end{array}
 \right.
\]

\begin{lem}
Assume $\nbige$ is convenient at $\lambda_0$.
Then we have the natural isomorphism
$\prolong{\phi^{\ast}\nbige}\simeq
 \phi^{\ast}\prolong{\nbige}$.

Let $\vecv$ be a frame of frame of $\prolong{\nbige}$
compatible with the decompositions $\lefttop{k}\nbigk$ $(k=1,\ldots,l)$.
In this case,
we have $u_i=\phi^{\ast}v_i$
(see the subsubsection {\rm\ref{subsubsection;10.27.9}}),
and the degrees $\lefttop{h}\deg^{\Fzero}(u_i)$ are given as follows:
\[
 \lefttop{h}\deg^{\Fzero}(u_i)
=\left\{
 \begin{array}{ll}
 \lefttop{1}\deg^{\Fzero}(v_i)
+\lefttop{2}\deg^{\Fzero}(v_i), & (h=1),\\
 \mbox{{}}\\
 \lefttop{h+1}\deg^{\Fzero}(v_i), & (h\geq 2).
 \end{array}
 \right.
\]
\end{lem}
\pf
If follows from the formula (\ref{eq;10.27.10})
for the degree of $u_i$.
\hfill\qed

\vspace{.1in}

Assume that $\nbige$ is convenient at $\lambda_0$.
Let $\sum_k \nbiga^k d\zeta_k/\zeta_k $ denote
the $\lambda$-connection of $\DD$ with respect to
the frame $\vecv$,
i.e.,
$\DD\vecv=\vecv\cdot(\sum_k \nbiga^k \cdot d\zeta_k/\zeta_k)$
holds.
Then we have the following formula:
\[
 \tilde{\DD}\vecu
=\vecu\cdot
 \Bigl(
 (\phi^{\ast}\nbiga^1+\phi^{\ast}\nbiga^2)\cdot\frac{dz_1}{z_1}
+\sum_{h=1}^{n-1}\phi^{\ast}\nbiga^{h+1}\cdot\frac{dz_h}{z_h}
 \Bigr).
\]
Hence we have the following formula:
\[
 \Res_h(\tilde{\DD})=
 \left\{
\begin{array}{ll}
 \phi^{\ast}\Res_1(\DD)_{|\nbigd_{\nibar}(\lambda_0,\epsilon_0)}
+\phi^{\ast}\Res_2(\DD)_{|\nbigd_{\nibar}(\lambda_0,\epsilon_0)}
 & (h=1)\\
 \mbox{{}}\\
 \phi^{\ast}\Res_{h+1}(\DD) & (h\geq 2).
\end{array}
 \right.
\]
Thus the nilpotent parts are as follows:
\[
 \tilde{\nbign}_i
=\left\{
 \begin{array}{ll}
 \phi^{\ast}(\nbign_1+\nbign_2) & (i=1)\\
 \mbox{{}}\\
 \phi^{\ast}\nbign_{i+1} & (i\geq 2).
 \end{array}
 \right.
\]
Hence $\tilde{\nbign}(\ibar)=\phi^{\ast}\tilde{\nbign}(\iitibar)$.
Let $\tilde{W}(\ibar)$ denote the weight filtration
of $\tilde{\nbign}(\ibar)$,
and then we have
$\tilde{W}(\ibar)=\phi^{\ast}W(\iitibar)$.

On the divisor $\tilde{\nbigd}_1(\lambda_0,\epsilon_0)$,
we put as follows:
\[
 \lefttop{1}\tilde{\nbigk}_u=
 \bigoplus_{u_1+u_2=u}\phi^{\ast}\lefttop{\nibar}\nbigk_{(u_1,u_2)}.
\]
On the divisors $\tilde{\nbigd}_i(\lambda_0,\epsilon_0)$ $(i\geq 2)$,
we put $\lefttop{i}\tilde{\nbigk}_u=\phi^{\ast}(\lefttop{i+1}\nbigk_u)$.
Then the decompositions $\lefttop{i}\tilde{\nbigk}$
satisfies the condition $(A)$ in Definition \ref{df;10.27.4}.
It is also easy to check the condition $(B)$
in Definition \ref{df;10.27.4}.
Hence we have the following lemma.
\begin{lem}
$\phi^{\ast}\nbige$ is convenient at $\lambda_0$.
\hfill\qed
\end{lem}

%% file: 26.tex

\subsubsection{Statement}

Let $\harmonicbundle$ be a tame harmonic bundle
over $X-D$.
Let $\lambda_0$ be a point of $\cnum_{\lambda}$
and $\epsilon_0$ be a small positive number.
Assume that $\prolong{\nbige}$ is convenient at $\lambda_0$
and locally free over $\nbigx(\lambda_0,\epsilon_0)$.
Let $\vecv$ be a frame of $\prolong{\nbige}$
compatible with
the tuple
$\bigl(\lefttop{i}\nbigk,W(\mbar)\,\big|\,
 i\in\lbar,\,m\in\lbar\bigr)$.

Since $\vecv$ is compatible with
$\EEzero$ and $\Fzero$,
we have the element $\vecu(v_i)$ for each $v_i$.
We put as follows:
\[
 b_m(v_i):=
 \paramap\bigl(\lambda,q_m(\vecu(v_i))\bigr)
\quad
 h_m(v_i):=
 \frac{1}{2}\cdot
 \bigl(\deg^{W(\mbar)}(v_i)-\deg^{W(\underline{m-1})}(v_i)\bigr).
\]
We put as follows:
\[
 v_i':=
 v_i\cdot\prod_{m=1}^l
 |\zeta_m|^{b_m(v_i)}\cdot\bigl(-\log|\zeta_m|\bigr)^{-h_m(v_i)}.
\]
Then we obtain the $C^{\infty}$-frame $\vecv'=(v_i')$
of $\nbige$ over
$\nbigx(\lambda_0,\epsilon_0)-\nbigd(\lambda_0,\epsilon_0)$.

We consider the following subsets
$Z$ and $\del Z$:
\begin{equation} \label{eq;10.27.70}
 \begin{array}{l}
 Z:=\bigl\{
 (\zeta_1,\ldots,\zeta_n)\in X-D\,\big|\,
 |\zeta_{j-1}|\leq |\zeta_{j}|\leq 2^{-1},\,\,j\in\lbar
 \bigr\},\\
 \mbox{{}}\\
 \del Z:=\bigl\{
 (\zeta_1,\ldots,\zeta_n)\in X-D\,\big|\,
 |\zeta_j|=2^{-1},\,\,j\in\lbar
 \bigr\}.
 \end{array}
\end{equation}
The purpose of this subsection is to show the following proposition.
\begin{prop} \label{prop;9.16.71}
 On the region $Z$,
 the $C^{\infty}$-frame  $\vecv'$ is adapted.
 Moreover there exists positive constants $C_i$ $(i=1,2)$,
 depending on $H(h,\vecv')_{|\del Z}$,
 such that $C_1\leq H(h,\vecv')\leq C_2$
 holds on the region $Z$.
\end{prop}

Let $\phi:\tilde{X}\lrarr X$ be a morphism
given in the subsubsection \ref{subsubsection;10.27.15}.
We obtain the holomorphic frame $\vecu=\phi^{\ast}\vecv$
of $\prolong{\phi^{\ast}\nbige}$.
We also obtain the $C^{\infty}$-frame
$\vecu'=(u_i')=\phi^{\ast}\vecv'$ over $\tilde{X}-\tilde{D}$.
It is easy to see the following:
\[
 u_i'
= u_i\cdot \prod_{m=1}^l
 |z_m|^{b_m(u_i)}\cdot
 \Bigl(
-\sum_{t\geq m}\log|z_t|
 \Bigr)^{-h_m(u_i)}
=u_i\cdot\prod_{m=1}^l
 |z_m|^{\sum_{j\leq m}b_j(v_i)}\cdot
 \Bigl(
 -\sum_{t\geq m}\log |z_t|
 \Bigr)^{-h_m(v_i)}.
\]
Here we have $b_m(u_i)=\sum_{j\leq m}b_j(v_i)$
and $h_m(u_i)=h_m(v_i)$.

Let us consider the following subsets
$\tilde{Z}$ and $\del\tilde{Z}$:
\[
\begin{array} {l}
 \tilde{Z}:=\bigl\{
 (z_1,\ldots,z_n)\in\tilde{X}-\tilde{D}\,\big|\,
 |z_i|\leq 1,(i\leq l-1),\,\,
 |z_{l}|\leq 2^{-1}
 \bigr\},\\
\mbox{{}}\\
 \del\tilde{Z}:=\bigl\{
 (z_1,\ldots,z_n)\in\tilde{X}-\tilde{D}\,\big|\,
 |z_i|=1,\,(i\leq l-1),\,\,
 |z_{l}|= 2^{-1}
 \bigr\}.
\end{array}
\]
Note the $\phi(\tilde{Z})=Z$ and $\phi(\del\tilde{Z})=\del Z$.

It is easy to see that
Proposition \ref{prop;9.16.71}
is equivalent to
Proposition \ref{prop;9.16.70}.
\begin{prop}\label{prop;9.16.70}
The $C^{\infty}$-frame $\vecu'$ is adapted on
$\tilde{X}-\tilde{D}$.
Moreover there exist positive constant $C_i$ $(i=1,2)$,
depending only on $H(h,\vecu')_{|\del\tilde{Z}}$,
such that $C_1\leq H(h,\vecu')\leq C_2$
over $\tilde{Z}$.
\end{prop}

We will show Proposition \ref{prop;9.16.70},
or equivalently Proposition \ref{prop;9.16.71}
in the following subsubsections.
We use an induction on the dimension of $X$.
We assume that the propositions hold in the case $\dim X\leq n-1$,
and we will prove the propositions hold in the case $\dim X=n$.
The hypothesis of the induction will be used in
Lemma \ref{lem;10.27.25}.

\subsubsection{Step 1. Independence of a choice of compatible frames}

\begin{lem}
Let $\vecv$ be a frame of $\prolong{\nbige}$
over $\nbigx(\lambda_0,\epsilon_0)$,
which is compatible with the tuple
$\bigl(\lefttop{k}\nbigk,W(\mbar)\,\big|\,
 k\in\lbar,\,\,m\in\lbar\bigr)$.

Assume that the claim in Proposition {\rm\ref{prop;9.16.71}}
holds for $\vecv$.
Then the same claim holds for any other frame
of \,\,$\prolong{\nbige}$
over $\nbigx(\lambda_0,\epsilon_0)$,
which is compatible with the tuple
$\bigl(\lefttop{k}\nbigk,W(\mbar)\,\big|\,
 k\in\lbar,\,\,m\in\lbar\bigr)$.
\end{lem}
\pf
Let $\vecv^{(1)}$ be other frame of $\prolong{\nbige}$
over $\nbigx(\lambda_0,\epsilon_0)$,
which is compatible with the tuple
$\bigl(\lefttop{k}\nbigk,W(\mbar)\,\big|\,
 k\in\lbar,\,\,m\in\lbar\bigr)$.
We have the relation of the form:
\[
 v_i^{(1)}
=\sum B_{j\,i}\cdot v_j.
\]
Here $B_{j\,i}$ are holomorphic on $X$
and $B_{j\,i\,|\nbigd_{\kbar}}=0$
unless the following holds:
\begin{equation} \label{eq;8.30.1}
 \lefttop{\kbar}\deg^{\nbigk}(v_i^{(1)})
=\lefttop{\kbar}\deg^{\nbigk}(v_j^{(1)}),
\quad\quad
 \deg^{W(\mbar)}(v_i^{(1)})
\leq
 \deg^{W(\mbar)}(v_j^{(1)})\,\,\,\,
 (\forall m\leq k)
\end{equation}
We have the induced relation
$u_i^{(1)}=\sum \phi^{\ast}B_{j\,i}\cdot u_j=
 \sum \tilde{B}_{j\,i}\cdot u_j$.
Then we have $\tilde{B}_{j\,i\,|\,\tilde{\nbigd}_k}=0$
unless (\ref{eq;8.30.1}) holds for $i$ and $j$.
We also have the induced relation
$u_i^{(1)\,\prime}=
 \sum \tilde{B}'_{j\,i}\cdot u_j'$,
and then we have the following:
\[
 \tilde{B}_{j\,i}'
=\tilde{B}_{j\,i}\cdot
 \prod_{p}|z_p|^{b_p(u_i^{(1)})-b_p(u_j)}
\cdot
 \prod_p\Bigl(
 -\sum_{t\geq p}\log|z_t|
 \Bigr)^{-h_p(u_i^{(1)})+h_p(u_j)}.
\]

Once we obtain the boundedness of $\tilde{B}'$,
then we obtain the boundedness of $\tilde{B}^{\prime\,-1}$
by symmetry.
It implies the equivalence of the adaptedness
of $\vecu^{(1)\,\prime}$ and $\vecu'$.
So we have only to prove the boundedness of $\tilde{B}'$.

\vspace{.1in}
\noindent
(i)
Note that 
we have $\tilde{B}_{j\,i\,|\tilde{\nbigd}_p}=0$
in the case $b_p(u_i^{(1)})-b_p(u_j)<0$.
We also have
$-1<b_p(u_i^{(1)})-b_p(u_j)$,
due to convenience of $\nbige$ at $\lambda_0$.
Thus $\tilde{B}'_{j\,i\,|\tilde{\nbigd_p}}=0$.

\vspace{.1in}
\noindent
(ii)
We have the following equality:
\[
 \prod_p\Bigl(
 -\sum_{t\geq p}\log |z_t|
 \Bigr)^{-h_p(u_i^{(1)})+h_p(u_j)}
=\prod_p\Bigl(
 \frac{-\sum_{t\geq p}\log |z_t|}{-\sum_{t\geq p+1}\log|z_t|}
 \Bigr)^{a_p}.
\]
Here we put as follows:
\[
 a_p=-\frac{1}{2}
 \Bigl(\deg^{W(\pbar)}(u_i^{(1)})-\deg^{W(\pbar)}(u_j)\Bigr).
\]
In the case $a_p\leq 0$,
it is easy to see that
$(-\sum_{t\geq p}\log|z_t|)^{a_p}\cdot
 (-\sum_{t\geq p+1}\log|z_t|)^{-a_p}$
is bounded.
In the case $a_p>0$,
we have $\tilde{B}_{j\,i\,|\,\tilde{\nbigd}_p}=0$,
and we have the following inequality on $\tilde{Z}$,
for some positive constant $C$:
\[
 \left(
 \frac{-\sum_{t\geq p+1}\log|z_t|-\log|z_p|}
 { -\sum_{t\geq p+1}\log|z_t|  }
 \right)^{a_p}
\leq
 \left(
 1+\frac{-\log|z_p|}{C}
 \right)^{a_p}
\]
Here we have used
$-\sum_{t\geq p+1}\log|z_t|\geq -\log|z_l|\geq C$
for some positive constant $C$.

From (i) and (ii),
the boundedness of $\tilde{B}'_{i\,j}$ follows immediately.
\hfill\qed

\subsubsection{Step 2. Strongly compatible frame $\vecv$}

Thus we pick a frame $\vecv$ which is strongly compatible
with
 $\bigl(
 \lefttop{k}\nbigk,W(\mbar)\,\big|\,
 k\in\lbar,\,\,m\in\lbar
 \bigr)$
in the sense of Corollary \ref{cor;10.27.26}.
Namely we take a frame compatible frame $\vecv$
as the following condition is satisfied:

\begin{condition}\label{condition;10.27.28}
A compatible frame $\vecv$ consists of
sections $v_{h,\vecu,\vech,i}$.
The following holds:
\[
 N(\itibar)v_{h,\vecu,\vech,i}
=\left\{
 \begin{array}{ll}
 v_{h,\vecu,\vech-2\vecdelta_1,i} & (-h+2\leq q_1(\vech)\leq h,\,\,
 h-q_1(\vech) \mbox{ is even})\\
 \mbox{{}}\\
 0 & (\mbox{otherwise}).
 \end{array}
 \right.
\]
We have $\lefttop{j}\deg^{\nbigk}(v_{h,\vecu,\vech,i})=q_j(\vecu)$
and $\deg^{W(\mbar)}(v_{h,\vecu,\vech,i})=q_m(\vech)$.
\hfill\qed
\end{condition}
Then we obtain the frame $\vecu=\bigl(u_{h,\vecu,\vech,i}\bigr)$
of $\prolong{\phi^{\ast}\nbige}$.

Let $q_1$ denote the projection of $\tilde{X}$
onto the first component $\Delta(1+\eta)$.
We have the naturally defined projection
$\Omega^{1,0}_{\tilde{X}}\lrarr q_1^{\ast}\Omega^{1,0}_{\Delta(1+\eta)}$.

Let $\gminiq_1(\tilde{\DD})$ denote the composite of the following:
\[
\begin{CD}
 \phi^{\ast}\nbige
 @>{\tilde{\DD}}>>
 \phi^{\ast}\nbige\otimes\Omega^{1,0}_{\tilde{X}}
 @>>>
 \nbige\otimes q_1^{\ast}\Omega_{\Delta(1+\eta)}^{1,0}.
\end{CD}
\]

\begin{lem}
In general,
let $\vecv^{(1)}$ be a frame of $\prolong{\nbige}$,
which is not necessarily compatible.
For the frame $\vecu^{(1)}=\phi^{\ast}\vecv^{(1)}$
of $\prolong{\phi^{\ast}\nbige}$,
we have the following implication:
\[
 \DD v^{(1)}_i=\sum u^{(1)}_j\cdot\nbiga^k_{j\,i}\frac{d\zeta_k}{\zeta_k}
\Longrightarrow
 \gminiq_1(\tilde{\DD})u^{(1)}_i
=\sum u^{(1)}_j\cdot\tilde{\nbiga}_{j\,i}^1\cdot\frac{dz_1}{z_1}
=\sum u^{(1)}_j\cdot \phi^{\ast}\nbiga^1_{j\,i}\cdot\frac{dz_1}{z_1}.
\]
\end{lem}
\pf
It can be shown by a direct calculation.
\hfill\qed

\vspace{.1in}

In particular,
if $\vecv$ is as in Condition \ref{condition;10.27.28}.
we have the following
on $\tilde{\nbigd}_k(\lambda_0,\epsilon_0)$ $(2\leq k\leq l)$.
(Here we use the notation $\tilde{\nbigd}_k$ instead of
$\tilde{\nbigd}_k(\lambda_0,\epsilon_0)$
for simplicity of the notation.):
\begin{equation} \label{eq;10.27.30}
 \gminiq_1(\tilde{\DD})
 u_{h,\vecu,\vech,i\,|\,
  \tilde{\nbigd}_k}
=\left\{
 \begin{array}{ll}
 \Bigl(
 \eigenmap(\lambda,u_1)\cdot u_{h,\vecu,\vech,i}
+u_{h,\vecu,\vech-2\vecdelta_1,i}\Bigr)\cdot
 {\displaystyle\left(\frac{dz_1}{z_1}\right)
 _{|\tilde{\nbigd}_k}},
 &(-h+2\leq \gminiq_1(\vech)\leq h,\,\,
 h-q_1(\vech) \mbox{ even }),\\
 \mbox{{}}\\
 \eigenmap(\lambda,u_1)\cdot u_{h,\vecu,\vech,i}
 {\displaystyle\left(
      \frac{dz_1}{z_1}\right)_{|\tilde{\nbigd}_k}}
 & (q_1(\vech)=-h)\\
\mbox{{}}\\
 0 & (\mbox{otherwise}).
 \end{array}
 \right.
\end{equation}

On the divisor $\tilde{\nbigd}_1$,
we have the following formula:
\begin{equation} \label{eq;10.27.31}
 \Res(\gminiq_1\tilde{\DD})
 u_{h,\vecu,\vech,i\,|\,
  \tilde{\nbigd}_1}
=\left\{
 \begin{array}{ll}
 \Bigl(
 \eigenmap(\lambda,u_1)\cdot u_{h,\vecu,\vech,i}
+u_{h,\vecu,\vech-2\vecdelta_1,i}\Bigr)
 _{|\tilde{\nbigd}_1},
 &(-h+2\leq \gminiq_1(\vech)\leq h,\,\,
 h-q_1(\vech) \mbox{ even }),\\
 \mbox{{}}\\
 \eigenmap(\lambda,u_1)\cdot u_{h,\vecu,\vech,i
  \,|\,\tilde{\nbigd}_1}
 & (q_1(\vech)=-h)\\
\mbox{{}}\\
 0 & (\mbox{otherwise}).
 \end{array}
 \right.
\end{equation}

\subsubsection{Model bundle and the comparing morphism}

For $h=q_1(\vech)\geq 0$, $\vecu$ and $i$,
we take the vector subspace of
$\prolong{\nbige}_{|(\lambda,O)}$ as follows:
\[
 V_{h,\vecu,\vech,i}:=
 \bigoplus _{a=0}^h
 \cnum\cdot v_{h,\vecu,\vech-2a\vecdelta_1,i}
\subset\prolong{\nbige}_{|(\lambda_0,O)}.
\]
Then we have the decomposition:
$\prolong{\nbige}_{|(\lambda_0,O)}
=\bigoplus V_{h,\vecu,\vech,i}$.
The decomposition is compatible with $\Res_1(\DD)$ and $\nbign_1$.
Let $\nbign_{h,\vecu,\vech,i}$
denote the restriction of $\nbign_1$
to $V_{h,\vecu,\vech,i}$.

Then we pick the model bundle
$ E(V_{h,\vecu,\vech,i},\nbign_{h,\vecu,\vech,i})
\otimes L(u_1)$
over $\Delta^{\ast}$.
We obtain the deformed holomorphic bundle
$\nbige_{h,\vecu,\vech,i}$ over $\Delta^{\ast}\times\cnum_{\lambda}$.
The direct sum of $\nbige_{h,\vecu,\vech,i}$
is denoted by $\nbige_0$.
We have the natural metric $h_0$ on $\nbige_0$,
which is a direct sum of the metrics
$h_{h,\vecu,\vech,i}$.

On $\Delta(\lambda_0,\epsilon_0)\times\Delta$,
we have the prolongation
$\prolong{\nbige_0}$
and the canonical frame $\vecu^0=(u^0_{h,\vecu,\vech,h,i})$.
Then we have the following:
\begin{equation} \label{eq;10.27.32}
 \DD_0 u^0_{h,\vecu,\vech,i}:=
 \left\{
 \begin{array}{ll}
 \bigl(
 \eigenmap(\lambda,u_1)\cdot u_{h,\vecu,\vech,i}
+u_{h,\vecu,\vech-2\vecdelta_1,i}
 \bigr)\cdot{\displaystyle\frac{dz_1}{z_1}}
 & (-h+2\leq q_1(\vech)\leq h,\,\,h-q_1(\vech)\mbox{ even })\\
 \mbox{{}}\\
 \bigl(
 \eigenmap(\lambda,u_1)\cdot u_{h,\vecu,\vech,i}\bigr)
\cdot{\displaystyle\frac{dz_1}{z_1}},
 & (q_1(\vech)=-h),\\
\mbox{{}}\\
 0 & (\mbox{otherwise}).
 \end{array}
 \right.
\end{equation}
We also have the following formula:
\begin{equation}\label{eq;10.27.33}
 \Res(\DD_0) u^0_{h,\vecu,\vech,i}:=
 \left\{
 \begin{array}{ll}
 \eigenmap(\lambda,u_1)\cdot u_{h,\vecu,\vech,i}
+u_{h,\vecu,\vech-2\vecdelta_1,i}
 & (-h+2\leq q_1(\vech)\leq h,\,\,h-q_1(\vech)\mbox{ even })\\
 \mbox{{}}\\
 \eigenmap(\lambda,u_1)\cdot u_{h,\vecu,\vech,i}
 & (q_1(\vech)=-h),\\
\mbox{{}}\\
 0 & (\mbox{otherwise}).
 \end{array}
 \right.
\end{equation}

We have the holomorphic vector bundle
$\phi^{\ast}q_1^{\ast}\prolong{\nbige_0}$
over $\Delta(\lambda_0,\epsilon_0)\times\tilde{X}$.
We have the $\lambda$-connection
$\tilde{\DD}_0:=\phi^{\ast}q_1^{\ast}\DD_0$.
The frames
$\phi^{\ast}q_1^{\ast}\vecu^0$ and $\vecu$
give the isomorphism
$\Phi:\phi^{\ast}q_1^{\ast}\prolong{\nbige}_0\lrarr
 \prolong{\nbige}$
over $\tilde{\nbigx}(\lambda_0,\epsilon_0)$.

\begin{lem}
We have the following equality on 
the divisor $\bigcup_{k=2}^l\tilde{\nbigd}_k$:
\[
 \Phi\circ
 \gminiq_1\bigl(\tilde{\DD}_0\bigr)
-\gminiq_1(\tilde{\DD})\circ\Phi=0
\]
We have the following equality on the divisor $\tilde{\nbigd_1}$:
\[
 \Phi\circ\Res\bigl(\gminiq_1(\tilde{\DD}_0)\bigr)
-\Res\bigl(\gminiq_1(\tilde{\DD})\bigr)\circ\Phi=0.
\]
\end{lem}
\pf
It immediately follows from 
(\ref{eq;10.27.30}),
(\ref{eq;10.27.31}),
(\ref{eq;10.27.32})
and (\ref{eq;10.27.33}).
\hfill\qed

\subsubsection{Step 3. The metric on $\phi^{\ast}q_1^{\ast}\nbige_0$}

We have the metric $\phi^{\ast}q_1^{\ast}h_{h,\vecu,\vech,i}$
on $\phi^{\ast}q_1^{\ast}\nbige_{h,\vecu,\vech,i}$.
We put as follows:
\[
 \tilde{h}_{h,\vecu,\vech,i}:=
 \phi^{\ast}q_1^{\ast}\bigl(h_{h,\vecu,\vech,i}\bigr)
\cdot
 \prod_{k=2}^l|z_k|^{-\sum_{2\leq t\leq k}\paramap(\lambda,u_t)}
 \times
 \prod_{k=2}^l
 \Bigl(
 -\sum_{m\geq k}\log|z_m|^2
 \Bigr)^{q_k(\vech)}.
\]
The metrics $\tilde{h}_{h,\vecu,\vech,i}$ induce the metric
$\tilde{h}_0$ on $\phi^{\ast}q_1^{\ast}\nbige_0$.

We put as follows:
\[
 \tilde{u}^{0\,\prime}_{h,\vecu,\vech,i}
:=
 \tilde{u}^0_{h,\vecu,\vech,i}\cdot
 \prod_{k=1}^l|z_k|^{\sum_{1\leq t\leq k}\paramap(\lambda,u_t)}
\times
 \prod_{k=1}^l\Bigl(
 -\sum_{m\geq k}\log|z_t|^2
 \Bigr)^{-q_k(\vech)}.
\]
\begin{lem}
Then the $C^{\infty}$-frame $\tilde{\vecu}^{0\,\prime}$ 
over
$\tilde{\nbigx}(\lambda_0,\epsilon_0)
-\tilde{\nbigd}(\lambda_0,\epsilon_0)$
is adapted with respect to the metric $\tilde{h}_0$.
\end{lem}
\pf
We put as follows:
\[
 u^{0\,\prime}_{h,\vecu,\vech,i}
:=
 u^0_{h,\vecu,\vech,i}
\cdot
 |z|^{\paramap(\lambda,u_1)}
\cdot
 (-\log|z_1|)^{-q_1(\vech)}.
\]
Then the $C^{\infty}$-metric
$\vecu^{0\,\prime}=(u^{0\,\prime}_{h,\vecu,\vech,i})$
over $\Delta^{\ast}$
is adapted with respect to $h_0$.
Then the adaptedness of $\vecu^{0\,\prime}$ immediately follows.
\hfill\qed

\vspace{.1in}

We put as follows:
\[
 \lefttop{\heartsuit}\tilde{u}^{0\,\prime}_{h,\vecu,\vech,i}
:=
 \tilde{u}^0_{h,\vecu,\vech,i}
\times
 \prod_{k=2}^l|z_k|^{\sum _{1\leq t\leq k}\paramap(\lambda,u_t)}
\times
 \prod_{k=3}^l
 \Bigl(
 -\sum_{m\geq k}\log|z_m|^2
 \Bigr)^{-q_k(\vech)}
\times
 \Bigl(
 -\sum_{m\geq 2}\log|z_m|^2
 \Bigr)^{-q_1(\vech)-q_2(\vech)}.
\]
We obtain the frame 
$C^{\infty}$-frame
$\lefttop{\heartsuit}\tilde{\vecu}^{0\,\prime}$
of $\nbige_0$ over
$\tilde{\nbigx}(\lambda_0,\epsilon_0)
-\tilde{\nbigd}(\lambda_0,\epsilon_0)$.
We put as follows:
\begin{equation} \label{eq;10.27.35}
Y:=
 \bigl\{(z_1,\ldots,z_n)\in\tilde{Z}\,\big|\,
 |z_1|=1 
 \bigr\}.
\end{equation}
\begin{cor}\label{cor;10.27.40}
The restriction
$\lefttop{\heartsuit}\tilde{\vecu}^{0\,\prime}_{|Y}$
is adapted with respect to the metric $\tilde{h_0}_{|Y}$.
\hfill\qed
\end{cor}

\subsubsection{Step 4. The end of the proof}
We put as follows:
\[
 \lefttop{\heartsuit}u'_{h,\vecu,\vech,i}
:=u_{h,\vecu,\vech,i}
 \cdot
 \prod_{k=2}^l |z_k|^{\sum_{1\leq t\leq k}\paramap(\lambda,u_t)}
\times
 \prod_{k=3}^l
 \Bigl(
 -\sum_{m\geq k}\log |z_m|^2
 \Bigr)^{-q_k(\vech)}
\times
 \Bigl(-\sum_{m\geq 2}\log|z_m|^2
 \Bigr)^{-q_1(\vech)-q_2(\vech)}.
\]

\begin{lem}\label{lem;10.27.25}
On the set $Y$ given in {\rm(\ref{eq;10.27.35})},
the $C^{\infty}$-frame
$\lefttop{\heartsuit}\vecu'$ is adapted with respect to
the metric $\phi^{\ast}h$.
\end{lem}
\pf
We put $\tilde{X}_a:=\{(a,z_2,\ldots,z_n)\in\tilde{X}\}$
and $X_a:=\{(\zeta_1,\ldots,\zeta_n)\in X\,|\,\zeta_1=a\zeta_2\}$.
Due to the result in the subsubsection \ref{subsubsection;10.27.36},
the restriction $\nbige_{|X_a}$ is convenient at $\lambda_0$,
and the frame $\vecv_{|X_a}$ is compatible.
Hence we obtain the result due to the assumption of the induction.
\hfill\qed

\begin{cor}
The restriction $\Phi_{|Y}$ is bounded
over the set $Y$.
\end{cor}
\pf
It immediately follows from
Corollary \ref{cor;10.27.40}
and Lemma \ref{lem;10.27.25}.
\hfill\qed

\vspace{.1in}
Then, by using the method explained in the subsection 6.1
of our previous paper,
we obtain the boundedness of $\Phi$
on the region $\tilde{Z}$.
Thus the induction can proceed,
and therefore we obtain
Proposition
\ref{prop;9.16.71} and Proposition \ref{prop;9.16.70}.
\hfill\qed

%% file: 26.1.tex

We put $X:=\Delta^n$, $D_i:=\{z_i=0\}$
and $D=\bigcup_{i=1}^l D_i$.
Let $\harmonicbundle$ be a tame harmonic bundle over $X-D$.
Let us pick any point $\lambda_0\in\cnum$.
We pick $\vecb\in\real^l$ such that
$q_i(\vecb)\not\in\KMS(\nbige^{\lambda_0},i)$
for $i=1,\ldots,l$.
Let pick a sufficiently small positive number $\epsilon_0$
such that $\prolongg{\vecb}{\nbige}$ is locally free
over the closure of $\Delta(\lambda_0,\epsilon_0)$.

Let $\vecv=(v_i)$ be a frame of $\prolongg{\vecb}{\nbige}$,
which is compatible with 
$\EE^{(\lambda_0)}$, $\Fzero$ and $W$
(Corollary\ref{cor;10.9.12}).
For each $v_i$,
we have the element $\vecu(v_i)\in\KMS(\nbige^0,\lbar)$
such that
$\deg^{\EEzero,\Fzero}(v_i)=
 \kmsmap\bigl(\lambda_0,\vecu(v_i)\bigr)$.
We put as follows, for each $v_i$:
\[
 b_j(v_i)=\paramap\bigl(\lambda,q_j(\vecu(v_i))\bigr),
\quad
 h_j(v_i)=\frac{1}{2}
 \bigl(
 \deg^{W(\jbar)}(v_i)
-\deg^{W(\underline{j-1})}(v_i)
 \bigr).
\]
Then we put as follows:
\[
 v_i':=v_i\cdot
 \prod_j \Bigl(
 |z_j|^{b_j(v_i)}
 \cdot
 (-\log|z_j|)^{-h_j(v_i)}
 \Bigr).
\]
Then we obtain the $C^{\infty}$-frame $\vecv'=(v_i')$
on $\nbigx(\lambda_0,\epsilon_0)-\nbigd(\lambda_0,\epsilon_0)$.

For any positive number $C$,
we put as follows:
\[
 Z(C):=
 \bigl\{
 (z_1,\ldots,z_n)\in X-D\,\big|\,
 |z_{i-1}|^C\leq |z_{i}|,\,\,
 (i\in\lbar)
 \bigr\}.
\]

\begin{thm} \label{thm;10.27.50}
Let $C$ be any positive number.
Then the $C^{\infty}$-frame $\vecv'$ is adapted
over the region $Z(C)\times\Delta(\lambda_0,\epsilon_0)$.
\end{thm}
\pf
First we give some easy reductions.

\begin{lem}
For the proof of Theorem {\rm\ref{thm;10.27.50}}
We may assume $\vecb=0$.
\end{lem}
\pf
We take the model bundle $L(\vecb)$ over $X-D$,
and prolongment 
$\prolongg{-\vecb}{\nbigl(\vecb)}$
of the deformed holomorphic bundle
over $\nbigx$.
We have the canonical frame $e$ of $\prolongg{-\vecb}{\nbigl(\vecb)}$
such that $|e|_h=\prod_{j=1}^l |z_j|^{q_j(\vecb)}$.

We have the naturally defined isomorphism
$\prolongg{\vecb}{\nbige}\otimes \prolongg{-\vecb}{\nbigl(\vecb)}
\simeq
 \prolong{\bigl(\nbige\otimes\nbigl(\vecb)\bigr)}$.
Once we show the claim of Theorem \ref{thm;10.27.50}
for $\prolong{\bigl(\nbige\otimes\nbigl(\vecb)\bigr)}$,
then we obtain the claim for $\prolongg{\vecb}{\nbige}$, too.
\hfill\qed

\begin{lem} \label{lem;10.27.51}
Assume that we have already shown the following claim:
\begin{quote}
(P):
The $C^{\infty}$-frame
$\vecv'$ is adapted
over $Z(C)\times\Delta(\lambda_0,\epsilon_0')$
for some positive number $\epsilon_0'$.
(The number $\epsilon_0'$ can be smaller than $\epsilon_0$).
\end{quote}
Then Theorem {\rm\ref{thm;10.27.50}} is obtained.
\end{lem}
\pf
Let $\lambda_1$ be any point of 
the closure $\overline{\Delta(\lambda_0,\epsilon_0)}$
of $\Delta(\lambda_0,\epsilon_0)$.
Due to the assumption of Lemma \ref{lem;10.27.51},
we may assume that we have some positive number $\epsilon_1'$
such that $\vecv'$ is adapted over
$Z(C)\times\Delta(\lambda_1,\epsilon_1)$.
We may assume that we can take a finite subset
$S\subset\Delta(\lambda_0,\epsilon_0)$
such that
$\Delta(\lambda_0,\epsilon_0)\subset
 \bigcup_{\lambda_1\in S}\Delta(\lambda_1,\epsilon_1)$.
Then the adaptedness of $\vecv'$
over $\Delta(\lambda_1,\epsilon_1)$ for $\lambda_1\in S$ implies
the adaptedness of $\vecv'$
over $\Delta(\lambda_0,\epsilon_0)$.
\hfill\qed

\vspace{.1in}
Let us return to the proof of Theorem \ref{thm;10.27.50}.
Note 
we may freely replace a positive number $\epsilon_0$
with a smaller one,
due to Lemma \ref{lem;10.27.51}.
Let $\eta_1$ be a positive number such that
$\eta_1\cdot\rank\nbige<1/3$.

\begin{lem}\label{lem;a12.1.50}
We can pick elements $(a_i,c_i)\in\seisuu_{>0}\times\{r\,|\,-1<r<0\}$
and numbers $\eta_i$ $(i\in\lbar)$
as follows, inductively:
\begin{itemize}
\item
 First we take $(a_1,c_1)$ satisfying the following:
\begin{itemize}
\item
We put
 $S_1:=\bigl\{
 a_1+\kappa(c_1\cdot b)\,\big|\,
 b\in\Par\bigl(\prolong{\nbigelambdazero},1\bigr)
 \bigr\}$.
Then $S_1$ is $\eta_1$-small,
and we have $0\not\in S_1$.
\end{itemize}
Then we put as follows:
\[
  \eta_2:=\frac{1}{3}
 \min\bigl\{|a-b|\,\big|\,a,b\in S_1,\,\,a\neq b\bigr\}.
\]
\item
 Suppose that we have already pick $(a_j,c_j)$ $(j< i)$
and $\eta_j$ $(j\leq i)$.
Then we take $(a_i,c_i)$ as follows:
 \begin{itemize}
  \item The inequality $c_i>C\cdot c_{i-1}$ holds.
  \item We put
   $S_i:=\bigl\{
 a_i+\kappa(c_i\cdot b)\,\big|\,
 b\in\Par(\prolong{\nbigelambdazero},i)
 \bigr\}$.
Then the set $S_i$ is $\eta_i$-small and we have $0\not\in S_i$.
 \end{itemize}
Then we put as follows:
 \[
  \eta_{i+1}:=\frac{1}{3}
 \min\bigl\{|a-b|\,\big|\,a,b\in S_i,\,\,a\neq b\bigr\}.
 \]
\end{itemize}
Moreover, we may assume that $c_i$ is sufficiently large
with respect to $\KMS(\prolong{\nbigelambdazero},i)$
for each $i$,
in the sense of Definition {\rm\ref{df;9.5.11}}.
\hfill\qed
\end{lem}

Let $\veca=(a_1,\ldots,a_l)$ be such an element of $\real^l$
as in Lemma \ref{lem;a12.1.50}.
We can take a small positive number $\epsilon_0'$
such that
$\prolong{\bigl(
 \phi_{\vecc}^{\ast}\nbige\otimes L(\veca)\bigr)}$
is locally free on 
$\Delta(\lambda_0,\epsilon_0')\times X$.
Due to our choice of $(c_1,\ldots,c_l)$,
$\prolong{\bigl(
\phi_{\vecc}^{\ast}\nbige\otimes L(\veca)\bigr)}$
is convenient.
(See Lemma \ref{lem;10.27.60} and Definition \ref{df;10.27.4}).

Let $e$ be the canonical base of 
the deformed holomorphic bundle $\nbigl(\veca)$
of the model bundle $L(\veca)$.
We have $|e|=\prod_{j=1}^l|z_j|^{a_j}$.
We put $e':=e\cdot \prod_{j=1}^l|z_j|^{-a_j}$,
and then we have $|e'|=1$.
On the other hand,
we have the frame $\vecu=(u_i)$
of $\prolong{\nbige}$ over
$\nbigx(\lambda_0,\epsilon_0')$
as in the subsubsection \ref{subsubsection;10.27.65}:
\[
 u_i:=\phi_{\vecc}^{\ast}v_i\cdot
 \prod_{j}z_j^{\nu(c_j\cdot b_j(v_i))}.
\]
Then the tensor product
$\vecu\otimes e=(u_i\otimes e)$ is a frame of
the vector bundle
$\prolong{\phi_{\vecc}^{\ast}\nbige\otimes L(\veca)}$,
which is compatible with
$\bigl(
 \lefttop{k}\nbigk,W(\mbar)\,\big|\,
 k\in\lbar,\,\,m\in\lbar
 \bigr)$.

We take the $C^{\infty}$-frame $\vecu'=(u_i')$
of $\phi_{\vecc}^{\ast}\nbige$
over $\nbigx(\lambda_0,\epsilon_0')-\nbigd(\lambda_0,\epsilon_0')$:
\[
 u_i':=u_i\cdot \prod_{j=1}^l
 |z_j|^{\kappa(c_j\cdot b_j(v_i))}.
\]
Due to Proposition \ref{prop;9.16.71},
 we obtain the adaptedness of
the $C^{\infty}$-frame
$\vecu'\otimes e'=(u_i'\otimes e')$
on the region $Z$, given in (\ref{eq;10.27.70}).
It is easy to see that
we have the relation
$u_i'=v_i'\cdot \omega_i$ for some
$C^{\infty}$-function $\omega_i$ on $Z$
such that $|\omega_i|=1$.
Thus the frame $\phi_{\vecc}^{\ast}\vecv'$ is adapted
on the region $Z$.
It implies the adaptedness of the frame $\vecv'$
on the following region: 
\[
  Z(c_1,\ldots,c_l):=
 \bigl\{
(z_1,\ldots,z_n)\in X-D\,\big|\,
 |z_{i-1}|^{c_{i-1}}<|z_{i}|^{c_{i}},\,\,
 i\in\lbar
 \bigr\}.
\]
Since we have $C\cdot c_i<c_{i+1}$ due to our choice
of $c_1,\ldots,c_l$,
we have the implication
$Z(C)\subset Z(c_1,\ldots,c_l)$.
Thus we are done.
\hfill\qed

%% file: a60.tex

\subsubsection{Preliminary}

Let $\hyperh$ denote the upper half plane.
We use the complex coordinate $(\zeta_1,\ldots,\zeta_n)$
of $\hyperh^n$.
We also use the real coordinate $\zeta_i=x_i+\sqrt{-1}y_i$.

\begin{lem}\label{lem;10.28.5}
Let $(a_k,n_k)$ $(k=1,\ldots,l)$ be elements of $\real\times\seisuu$.
We have the following equality:
\begin{equation}\label{eq;10.28.1}
 \prod_{k=1}^l y_k^{-a_k+a_{k-1}+n_k}
=\prod_{k=1}^{l-1}
 \left(
 \frac{y_k}{y_{k+1}}
 \right)^{-a_k+\sum_{i\leq k}n_i}
 \times y^{-a_l+\sum_{i\leq l}n_i}.
\end{equation}
\end{lem}
\pf
We use an induction on $l$.
We assume that the following equality holds:
\begin{equation}\label{eq;10.28.2}
 \prod_{k=1}^{l-1}
 y_k^{-a_k+a_{k-1}+n_k}
=\prod_{k=1}^{l-2}
  \left(
 \frac{y_k}{y_{k+1}}
 \right)^{-a_k+\sum_{i\leq k}n_i}
 \times y^{-a_{l-1}+\sum_{i\leq l-1}n_i}.
\end{equation}
By a direct calculation, we have the following:
\begin{equation}\label{eq;10.28.3}
 y_{l-1}^{-a_{l-1}+\sum_{i\leq l-1}n_i}
\times y_l^{-a_l+a_{l-1}+n_l}
=\left(
 \frac{y_{l-1}}{y_l}
 \right)^{-a_{l-1}+\sum_{i\leq l-1}n_i}
\times y_l^{-a_l+\sum_{i\leq l}n_i}.
\end{equation}
From the equalities (\ref{eq;10.28.2}) and (\ref{eq;10.28.3}),
we obtain the equality (\ref{eq;10.28.1}).
Thus the induction can proceed.
\hfill\qed

\vspace{.1in}

Let $C_i$ $(i=1,2,3)$ be positive numbers.
We put as follows:
\[
 \tilde{Z}(C_1,C_2,C_3)
:=\Bigl\{
 (\zeta_1,\ldots,\zeta_l)\in\hyperh^l
 \,\Big|\,
 |x_i|\leq C_1,\,\,
 y_{i+1}\leq C_2\cdot y_i,\,\,(i\in\lbar),\,\,
 y_n\geq C_3
 \Bigr\}.
\]
Let $(a_k,n_k)$ $(k=1,\ldots,l)$ be elements of $\real\times\seisuu$.
Let us consider the following function:
\[
 F_l:=\prod_{k=1}^ly_k^{-a_k+a_{k-1}}\cdot
 \bigl|x_k+\sqrt{-1}y_k\bigr|^{n_k}.
\]

\begin{lem}\label{lem;10.28.4}
Let $C_b$ $(b=1,2,3)$ be any positive numbers.
Assume $a_k\geq \sum_{i\leq k}n_i$.
The function $F_l$ is bounded
on the region $\tilde{Z}(C_1,C_2,C_3)$.
\end{lem}
\pf
We have only to show the boundedness
of $\prod_{k=1}^ly_k^{-a_k+a_{k-1}+n_k}$
over $\tilde{Z}(C_1,C_2,C_3)$.
Then Lemma \ref{lem;10.28.4} follows from
Lemma \ref{lem;10.28.5}.
\hfill\qed

\subsubsection{Norm estimate in the case where $\lambda$ is fixed}

Let us consider the norm estimate for flat sections.
For simplicity,
we consider the case $X=\Delta^n$ and $D=\bigcup_{i=1}^n D_i$,
where we put $D_i:=\{z_i=0\}$.
Let $\harmonicbundle$ be a tame harmonic bundle
over $X-D$.
We have the universal covering $\pi:\hyperh^n\lrarr X-D$,
given by $\zeta_i\longmapsto\exp\bigl(\sqrt{-1}\zeta_i\bigr)$.
Let $\lambda$ be a point of $\cnum^{\ast}_{\lambda}$.
Let us consider the norm estimate of flat sections
of $\pi^{-1}\nbigelambda$.

Let $\vecs=(s_i)$ be a frame of
$H(\nbigelambda)$, which is compatible with
$\bigl(\lefttop{i}\EE,\lefttop{j}\nbigf,W(\mbar)\,\big|\,
 i\in\nbar,\,j\in\nbar,\,m\in\nbar
 \bigr)$.
We have the elements
$\vecu(s_i)\in\KMSoverline(\nbige^0,\nbar)$,
such that
$\deg^{\EE,\nbigf}(s_i)=\kmsmap^f(\lambda,\vecu(s_i))$.
Let $M_k^u$ denote the unipotent part of the monodromy
of $\DD^{\lambda,f}$,
and we put as follows:
\[
 N_k:=\frac{-1}{2\pi\sqrt{-1}}
 \log M_k^u.
\]
For any $\vecn\in\seisuu_{\geq 0}^l$,
we put as follows:
\[
 N^{\vecn}:=\prod_{k=1}^l N_k^{n_k}.
\]
The matrix $b(\vecn):=\bigl(b(\vecn)_{j\,i}\bigr)$
is determined by the relation
$N^{\vecn}\vecs=\vecs\cdot b(\vecn)$,
i.e.,
$N^{\vecn}s_i=\sum b(\vecn)_{j\,i}s_j$.

\begin{lem}\label{lem;10.28.6}
Assume that $b(\vecn)_{j\,i}\neq 0$.
Then we have the following:
\begin{itemize}
\item
 $\deg^{\EE}(s_i)=\deg^{\EE}(s_j)$.
\item
 $\lefttop{k}\deg^{\nbigf}(s_i)\geq \lefttop{k}\deg^{\nbigf}(s_j)$
 for any $k=1,\ldots,n$.
\item
Let $l$ be any integer such that $1\leq l\leq n$.
 In the case 
 $\lefttop{k}\deg^{\EE,\nbigf}(s_i)= \lefttop{k}\deg^{\EE\nbigf}(s_j)$
 for $k\leq l$,
 we also have the following, for any $k\leq l$:
\[
 \deg^{W(\kbar)}(s_i)\geq
 \deg^{W(\kbar)}(s_j)+2\sum_{t\leq k}n_t.
\]
\end{itemize}
\end{lem}
\pf
It immediately follows from our choice of $\vecs$.
\hfill\qed

\vspace{.1in}
We put $v_i:=F\bigl(s_i,\paramap^f(\lambda,\vecu(s_i))\bigr)$.
Then $\vecv=(v_i)$ is a frame of $\prolong{\nbigelambda}$
which is compatible with
$\bigl(\lefttop{i}\EE,\lefttop{j}F,W(\mbar)\,\big|\,
 i\in\nbar,\,j\in\nbar,\,m\in\nbar
 \bigr)$.
We have the elements
$\vecu(v_i)\in \KMS(\nbige^0,\nbar)$
such that $\deg^{\EE,F}(v_i)=\kmsmap(\lambda,\vecu(v_i))$.

We put $\alpha_k(s_j):=\lambda^{-1}\cdot\eigenmap(\lambda,u_k(v_j))$.
Let $\vecn!$ denote the number $\prod_{i=1}^nn_i!$
for $\vecn=(n_1,\ldots,n_n)$.
Then we have the following:
\begin{equation}
 v_i=\prod_{k=1}^n z_k^{\alpha_k(s_i)}
 \cdot\sum_j\frac{b(\vecn)_{j\,i}}{\vecn!}
 \cdot\prod_{k=1}^n
 \bigl(\log z_k\bigr)^{n_k}\cdot s_j.
\end{equation}
Note we also have the following relation:
\[
 \paramap(\lambda,u_k(v_j))
=\paramap^f(\lambda,u_k(s_j))-\Re\bigl(
 \alpha_k(s_j)
 \bigr).
\]
We put $h_k(v_j):=2^{-1}\bigl(
 \deg^{W(\kbar)}(v_j)-\deg^{W(\underline{k-1})}(v_j)
\bigr)$
and
$h_k(s_j):=2^{-1}\bigl(
\deg^{W(\kbar)}(s_j)-\deg^{W(\underline{k-1})}(s_j)
 \bigr)$.
Note we have $h_k(s_j)=h_k(v_j)$.

We put as follows:
\[
 \begin{array}{l}
 v_j':=
 v_j\cdot \prod_{k=1}^l |z_k|^{\paramap(\lambda,u_k(v_j))}
 \cdot\bigl(
 -\log|z_k|
 \bigr)^{-h_k(v_j)},\\
 \mbox{{}}\\
 s_j':=
 s_j\cdot \prod_{k=1}^l |z_k|^{\paramap^f(\lambda,u_k(s_j))}
 \cdot\bigl(
 -\log|z_k|
 \bigr)^{-h_k(s_j)}.
 \end{array}
\]
Then we obtain the frames $\vecv'=(v_i')$
and $\vecs'=(s_i')$.

We put as follows:
\[
 B_{j\,i}:=
 \sum_{\vecn}\frac{b(\vecn)_{j\,i}}{\vecn!}
 \prod_{k}
 z_k^{\paramap^f(\lambda,u_k(s_i))
 -\paramap^f(\lambda,u_k(s_j))}
\times\bigl(-\log|z_k|\bigr)^{-h_k(s_i)+h_k(s_j)}
 \times\bigl(-\log z_k\bigr)^{n_k}.
\]
Thus we obtain the matrix valued function $B$.
\begin{lem}
We have the relation
$\vecv'=\vecs'\cdot B$.
\end{lem}
\pf
We have the following:
\[
  v_i'=
\sum_{j,\vecn}
\frac{b(\vecn)_{j\,i}}{\vecn!}
 \prod_kz_k^{
 \alpha_k(s_i)+
 \paramap(\lambda,u_k(v_i))-\paramap^f(\lambda,u_k(s_j))}
 \bigl(-\log|z_k|\bigr)^{-h_k(v_i)+h_k(s_j)}
 \times (\log z)^{n_k}\cdot s_j'.
\]
We have 
$ \alpha_k(s_i)+\paramap(\lambda,u_k(v_i))
=\paramap^f(\lambda,u_k(s_i))$
and $h_k(v_i)=h_k(s_i)$.
Thus we obtain the result.
\hfill\qed

\vspace{.1in}

\begin{lem}\mbox{{}}
\begin{itemize}
\item We have $B_{i\,i}=1$.
\item
Assume $B_{j\,i}\neq 0$,
then we have
$\lefttop{k}\deg^{\nbigf}(s_i)\geq \lefttop{k}\deg^{\nbigf}(s_j)$.
\end{itemize}
Namely the matrix $B$ is triangular,
and the diagonal components are $1$.
\end{lem}
\pf
It immediately follows from Lemma \ref{lem;10.28.6}.
\hfill\qed

\begin{lem}\label{lem;10.28.7}
The matrix valued functions $B$ and $B^{-1}$ are bounded
over $\tilde{Z}(C_1,C_2,C_3)$.
\end{lem}
\pf
Assume $b(\vecn)_{j\,i}\neq 0$.
Then we have
$\paramap^f(\lambda,u_k(s_i))-\paramap^f(\lambda,u_k(s_j))\geq 0$
for any $k$.
Moreover let $h$ be the number
such that $\paramap^f(\lambda,u_k(s_i))-\paramap^f(\lambda,u_k(s_j))=0$
for any $k< h$
and that $\paramap^f(\lambda,u_h(s_i))-\paramap^f(\lambda,u_h(s_j))\neq 0$.

\begin{lem} \label{lem;10.28.8}
Let $h$ be as above.
On the region $\tilde{Z}(C_1,C_2,C_3)$,
we have the boundedness of the following:
\[
 \prod_{k\geq h}z_k^{\paramap^f\bigl(\lambda,u_k(s_i)-u_k(s_j)\bigr)}
\times\bigl(-\log|z_k|\bigr)^{-h_k(s_i)+h_k(s_j)}
\times (\log z_k)^{n_k}.
\]
\end{lem}
\pf
We have only to compare the order of
$|z_{h}|^{\paramap^f\bigl(\lambda,u_h(s_i)-u_h(s_j)\bigr)}$
and $\prod_{k\geq h}(-\log|z_k|)^M$.
Note there exists a positive constant
we have $|z_h|^{C}\leq |z_k|$
for any $k\geq h$ over $Z(C_1,C_2,C_3)$.
\hfill\qed

\begin{lem} \label{lem;10.28.9}
Let $h$ be as above.
We have the boundedness of the following function
over $\tilde{Z}(C_1,C_2,C_3)$:
\[
 \prod_{k=1}^{h-1}
 (-\log|z_k|)^{-h_k(s_i)+h_k(s_j)}
\times(\log z_k)^{n_k}.
\]
\end{lem}
\pf
We put
$a_k:=2^{-1}\bigl(\deg^{W(\kbar)}(s_i)-\deg^{W(\kbar)}(s_j)\bigr)$.
Then we have
$-h_k(s_i)+h_k(s_j)=-a_k+a_{k-1}$
and $-a_k+\sum_{i\leq k}n_i\leq 0$ for any $k\leq h$.
Then we obtain the desired boundedness 
from Lemma \ref{lem;10.28.4}.
\hfill\qed

\vspace{.1in}
Let us return to the proof of Lemma \ref{lem;10.28.7}.
The boundedness of $B$ immediately follows from Lemma \ref{lem;10.28.8}
and Lemma \ref{lem;10.28.9}.
Since $B$ is triangular such that the diagonal components are 1,
the boundedness of $B^{-1}$
follows from the boundedness of $B$.
\hfill\qed

\begin{thm}
The frame $\vecs'$ is adapted on the region
$\tilde{Z}(C_1,C_2,C_3)$.
\end{thm}
\pf
It follows from the adaptedness of $\vecv'$
and the boundedness of $B$ and $B^{-1}$
on the region $\tilde{Z}(C_1,C_2,C_3)$.
\hfill\qed

%% file: a74.tex

We recall the nearby cycle functor
for $\nbigr$-modules introduced by Sabbah,
with some minor generalization.
We recommend the reader to read the readable paper \cite{sabbah}.
In particular, see the chapters 0 and 1 for the basic
of $\nbigr$-modules,
and see the chapter 3 in \cite{sabbah}
for $V$-filtration and the nearby cycle functor.
For most of the definitions and the lemmas
contained in this section,
the reader can find the counterpart in \cite{sabbah}.
We also use some results in \cite{sabbah}
without mention.

We consider the right $\nbigr$-modules
in this subsection.

\begin{rem}
Sabbah kindly informed the author on the revision of his paper
{\rm\cite{sabbah2}},
in which the generalization is discussed.
We keep this section for our reference
in the later discussion.
\hfill\qed
\end{rem}

%% file: 30.tex

\subsubsection{$V$-filtration}

Let $\nbigm$ be an $\nbigr_{X\times\cnum}$-module.
Let $t$ be the coordinate of $\cnum$.

\begin{df}
A $V$-filtration at $\lambda_0$ is defined to be
a filtration $\Ulambdazero$ of $\nbigm$ indexed by $\seisuu$,
defined on $\nbigx(\lambda_0,\epsilon_0)$ for some $\epsilon_0>0$
indexed by $\seisuu$
satisfying the following:
\[
 \Ulambdazero_a(\nbigm)\cdot V_m(\nbigr)
\subset
 \Ulambdazero_{a+m}(\nbigm).
\]
(See the section {\rm 3.1.a.} in {\rm\cite{sabbah}} for $V_m(\nbigr)$.)
A $V$-filtration $\Ulambdazero$ at $\lambda_0$
is called monodromic,
if there exists a monic $b(s)\in\cnum[\lambda][s]$
such that
\begin{itemize}
\item
 $b\bigl(t\del_t-k\cdot\lambda\bigr)$ acts trivially on
 $\Ulambdazero_k(\nbigm)\big/\Uzero_{k-1}(\nbigm)$
 for any $k\in\seisuu$.
\item
 $\gcd\Bigl(b\bigl(s-k\lambda\bigr),\,\,
 b\bigl(s-l\lambda\bigr)\Bigr)\in\cnum[\lambda]-\{0\}$
 if $k\neq l$.\hfill\qed
\end{itemize}
\end{df}

Later we will consider a refinement of $V$-filtration,
which is the filtration indexed by $\real$.
It will be also called by $V$-filtration.
Recall the definition of good $V$-filtration.
\begin{df}
Let $\bigl(\nbigm,\Ulambdazero\nbigm\bigr)$ be
a $V$-filtered $\nbigr_{\nbigx}$-module
defined over $\nbigx(\lambda_0,\epsilon_0)$.
It is called good, if the following holds:
\begin{itemize}
\item
 For any compact subset $K\subset \nbigx(\lambda_0,\epsilon_0)$,
 there exists $k_0\geq 0$ such that the following holds:
 \[
  \Ulambdazero_{-k}\nbigm
=\Ulambdazero_{-k_0}\nbigm\cdot
 t^{k-k_0},
\quad\quad
 \Ulambdazero_k\nbigm
=\sum_{0\leq j\leq k-k_0}
 \Ulambdazero_{k_0}\nbigm\cdot \deldel^j_t.
 \]
\item
 $\Ulambdazero_l\nbigm$ is $V_0\nbigr$-coherent.
\hfill\qed
\end{itemize}
\end{df}

Let $A$ be a finite subset of $(\real/\seisuu)\times\cnum$.
Let $\pi:\real\times\cnum\lrarr(\real/\seisuu)\times\cnum$
be the projection.
We put $\tilde{A}:=\pi^{-1}(A)$.
For any real number $c\in\real$,
we put as follows:
\[
 A_c^{(\lambda_0)}:=
\bigl\{u\in\tilde{A}\,\big|\,
 c-1<\paramap(\lambda_0,u)\leq c
 \bigr\}.
\]

\begin{df} \label{df;9.17.5}
A coherent $\nbigr$-module
$\nbigm$ is called specializable along $X_0$ at $\lambda_0$ for $A$,
if the following holds:
\begin{itemize}
\item
 There exists a good $V$-filtration $\Ulambdazero(\nbigm)$
 at $\lambda_0$.
\item
 We have a map $f:A\lrarr\seisuu_{\geq 0}$,
 which can be regarded as the function
 on $\tilde{A}$ or $A_c^{(\lambda_0)}$ $(c\in\real)$,
and we put as follows:
\[
 b_U(s):=\prod_{u\in A_0^{(\lambda_0)}}
  \bigl(s+\eigenmap(\lambda,u)\bigr)^{f(u)}.
\]
Then 
$ b_U(t\!\cdot\!\deldel_t-k\lambda)$
acts trivially on $\Ulambdazero_k\big/\Ulambdazero_{k-1}$.
\end{itemize}
In such case,
we say $\nbigm$ is specializable for $(X_0,A,\lambda_0)$.
If we would like to distinguish a $V$-filtration $\Uzero$
and a function $f$ as above,
we say that $\nbigm$ is specializable for
$(X_0,A,\lambda_0,\Uzero,f)$.
\hfill\qed
\end{df}

Note the following relation:
\[
 b_U(s-k\lambda)=\prod_{u\in A_k^{(\lambda_0)}}
 \bigl(s+\eigenmap(\lambda,u)\bigr)^{f(u)}.
\]

The following lemma is clear.
\begin{lem}
If $\nbigm$ is specializable for $(X_0,\lambda_0,A)$,
then the $\nbigr$-submodules and the $\nbigr$-quotient modules are also
specializable for  $(X_0,\lambda_0,A)$.
\hfill\qed
\end{lem}

Recall the following lemma.
\begin{lem}
Assume that $\nbigm$ is specializable for
$(X_0,\lambda_0,A,\Uzero,f)$.
Let $m$ be a local section of $\nbigm$ around $\lambda_0$.
Then there exists a finite subset $S\subset \tilde{A}$
such that 
we have $m\cdot b_S(t\deldel_t)\in m\cdot V_{-1}(\nbigr)$.
Here we put as follows:
\[
 b_S(s):=\prod_{u\in S}\bigl(s+\eigenmap(\lambda,u)\bigr)^{f(u)}.
\]
\end{lem}
\pf
$m\cdot \nbigr$ is a submodule of $\nbigm$.
Thus $\Ulambdazero$ induces a good $V$-filtration
of $m\cdot\nbigr$.
On the other hand, the good $V$-filtration of $\nbigr$
induces the good $V$-filtration of $m\cdot\nbigr$.
Since the two good $V$-filtrations are equivalent,
the result holds.
\hfill\qed

\subsubsection{Refinement and the decomposition}

On $\Ulambdazero_k\big/\Ulambdazero_{k-1}$,
we have the filtration $\bigl\{F_c\,\big|\,k-1\leq c\leq k\bigr\}$
defined as follows:
\[
 F_c:=
 \ker\Bigl(
 \prod_{\substack{
 u\in A_k^{(\lambda_0)},\\
 \paramap(\lambda_0,u)\leq c}}
 \bigl(t\deldel_t+\eigenmap(\lambda,u)\bigr)^{f(u)}
 \Bigr)
\subset \Uzero_k\big/\Uzero_{k-1}.
\]
Then we obtain the refinement of the filtration $\Ulambdazero$
as follows:
Let $c$ be any real number.
We take the integer $k$ satisfying $k-1< c\leq k$.
Let $\pi_k$ denote the naturally defined projection
$\Ulambdazero_k\lrarr \Ulambdazero_k\big/\Ulambdazero_{k-1}$.
Then we put $\Ulambdazero_c:=\pi_k^{-1}(F_c)$.
In the following, we put as follows
for any real number $c\in\real$:
\[
 \Gr^{\Ulambdazero}_c(\nbigm)
:=\Ulambdazero_{c}/\Ulambdazero_{<c}.
\]

For any real number $c$, we put as follows:
\[
 \nbigk(A,c,\lamdazero)
:=\bigl\{u\in\tilde{A}\,\big|\,\paramap(\lambda_0,u)=c\bigr\},
\quad\quad
 b_c(s):=
 \prod_{ u\in \nbigk(A,c,\lambdazero) }
\bigl(s+\eigenmap(\lambda,u)\bigr)^{f(u)}.
\]
Then $b_c(t\deldel_t)$ acts trivially on
$\Gr^{\Ulambdazero}_c(\nbigm)$.

\begin{lem}
We have the decomposition
on a neighbourhood of $\lambda_0$:
\begin{equation} \label{eq;9.17.1}
 \Gr^{\Ulambdazero}_c(\nbigm)
=\bigoplus_{u\in\nbigk(A,c,\lambda_0)}
 \EE\bigl(-\eigenmap(\lambda,u)\bigr).
\end{equation}
Here $\EE\bigl(-\eigenmap(\lambda,u)\bigr)$ denotes
the kernel of $\bigl(t\deldel_t+\eigenmap(\lambda,u)\bigr)^N$
for any sufficiently large integer $N$.
\end{lem}
\pf
Let $u_i$ $(i=1,2)$ be elements of
$\nbigk(A,c,\lambda_0)$.
Note we have $\eigenmap(\lamdazero,u_1)\neq \eigenmap(\lamda_0,u_2)$
if $u_1\neq u_2$.
Then the decomposition (\ref{eq;9.17.1}) immediately follows.
\hfill\qed

\vspace{.1in}
We put as follows:
\[
 \psi^{U^{(\lambda_0)}}_{t,u}(\nbigm):=\EE(-\eigenmap(\lambda,u)).
\]

For $\lambda\in\cnum_{\lamda}$,
the function
$\phi_{\lambda}:\nbigk(A,c,\lamda_0)\lrarr \real$
is defined by
$u\longmapsto \paramap(\lambda,u)$.

\begin{lem}
Assume $|\lambda-\lambda_0|$ is sufficiently small.
For any $c,c'\in \real$ such that $c\neq c'$,
we have the following:
\[
 \phi_{\lamda}\bigl(\nbigk(A,c,\lamda_0)\bigr)
\cap
 \phi_{\lamda}\bigl(\nbigk(A,c',\lamda_0)\bigr)=\emptyset.
\]
\end{lem}
\pf
Since $A$ is finite,
the set
$\bigl\{c\in\real\,
  \big|\,\nbigk(A,c,\lamda_0)\neq \emptyset\bigr\}$
is a discrete and periodic subset of $\real$.
Then the lemma immediately follows.
\hfill\qed

\begin{lem} \label{lem;9.17.10}
Assume that $\nbigm$ is specializable for $(X_0,A,\lamda_0)$.
If $|\lambda_1-\lamda_0|$ is sufficiently small,
then $\nbigm$ is specializable 
for $(X_0,A,\lambda_1)$.
If there exists the $V$-filtration $\Uzero$ at $\lambda_0$
such that $\Gr^{U^{(\lambda_0)}}$ is strict,
then there exists the $V$-filtration $U^{(\lambda_1)}$
at $\lambda_1$ such that $\Gr^{U^{(\lambda_1)}}$ is strict.
\end{lem}
\pf
For any real number $d\in\real$,
we put as follows:
\[
 S(d):=\bigl\{c\in\real\,\big|\,\exists u\in \nbigk(A,c,\lambdazero),\,\,
 \paramap(\lambda_1,u)=d \bigr\}.
\]
Since $|\lambda_1-\lambda_0|$ is small,
we have $\big|S(d)\big|\leq 1$.
First, let us consider the case
$S(d)=\{c\}$.
Let $\pi_c:\Ulambdazero_c\lrarr \Gr^{\Ulambdazero}_c$
be the projection,
and we put as follows:
\[
 U^{(\lamda_1)}_d:=
 \pi_c^{-1}\Bigl(
 \bigoplus_{\substack{
 u\in\nbigk(A,c,\lamda_0),\\
 \paramap(\lambda_1,u)\leq d
 }}
 \psi^{U^{(\lambda_0)}}_{t,u}
 \Bigr).
\]
Let us consider the case $S(d)=\emptyset$.
In that case,
we put $d_0:=\max\{d'\leq d\,|\,S(d')\neq \emptyset\}$,
and we put
$U_{d}^{(\lambda_1)}:=U_{d_0}^{(\lambda_1)}$.
Then it is easy to check that
$U^{(\lambda_1)}$ is the filtration we desired.
\hfill\qed

\begin{rem}
In the proof, the construction of $U^{(\lambda_1)}$ from
 $U^{(\lambda_0)}$
is also given.
\hfill\qed
\end{rem}

\subsubsection{Lemmas for uniqueness}

\begin{lem} \label{lem;9.17.19}
Assume that $\nbigm$ is specializable
for $(X_0,A,\lambda_0,\Uzero,f)$ and
$(X_0,A',\lambda_0,U^{\prime\,(\lambda_0)},f'  )$.
Assume the strictness of
$\Gr^{U^{\prime\,(\lambda_0)}}(\nbigm)$.
Then $\Uzero_a\subset U_a^{\prime\,(\lambda_0)}$
for any $a\in\real$.
\end{lem}
\pf
In this proof, we omit to denote $(\lambda_0)$.
We straightforwardly follow the argument given in
the proof of (2) Lemma 3.3.4. in  \cite{sabbah}.

Since $U$ and $U'$ are good,
there exists $l\geq 0$
such that $U'_{c-l}\subset U_c\subset U'_{c+l}$
for any $c$.
Pick $m\in U_c\nbigm$,
and we take $d$ as
$m\in U_d'-U'_{<d}$.
Note that $d\leq c+l$.

There exists a finite subset $S\subset \tilde{A}$
satisfying the following:
\begin{itemize}
\item
 For any element $u\in S$, the inequality $\paramap(\lambda_0,u)\leq c$
 holds.
\item
 We put $B_U(s):=\prod_{u\in S}\bigl(s+\eigenmap(\lambda,u)\bigr)^{f(u)}$.
Then $m\cdot B_U(t\deldel_t)\in U_{<d-l}\subset U'_{<d}$.
\end{itemize}

On the other hand,
we put as follows:
\[
 B_{U}'(s):=
\prod_{u\in\nbigk(A',d,\lambda_0)}
\bigl(s+\eigenmap(\lambda,u)\bigr)^{f'(u)}.
\]
Then we have $m\cdot B_{U}'(t\deldel_t)\in U'_{<d}$.

Assume $c<d$,
and we will derive a contradiction.
Since we have $\paramap(\lambda_0,u)=d$ for any
element $u\in\nbigk(A',d,\lambda_0)$,
we have
$\bigl\{\eigenmap(\lambda,u)\,\big|\,u\in S\bigr\}
\cap
 \bigl\{\eigenmap(\lambda,u)\,\big|\,u\in\nbigk(A',\lambda_0,d)\bigr\}
=\emptyset$,
where we regard $\eigenmap(\lambda,u)$ as
functions of $\lambda$.
Thus we have $\gcd(B_U,B_{U'})\in \cnum[\lambda]-\{0\}$.
Hence there exists an element $g(\lambda)\in\cnum[\lambda]-\{0\}$ 
such that
$m\cdot g(\lambda)\in U'_{<d}$.
Due to the strictness of $\Gr^{U'}_{d}(\nbigm)$,
we obtain $m\in U'_{<d}$,
which contradicts our choice of $d$.
Hence $c\geq d$.
It implies $U_c\subset U'_c$.
\hfill\qed

\begin{lem} \label{lem;9.17.7}
Assume that $\nbigm$ is specializable
for $(X_0,A,\lambda_0,\Uzero,f)$ and
$(X_0,A',\lambda_0,U^{\prime\,(\lambda_0)},f'  )$.
Assume the strictness of
$\Gr^{\Uzero}(\nbigm)$ and
$\Gr^{U^{\prime\,(\lambda_0)}}(\nbigm)$.
Then we have $\Uzero=U^{\prime\,(\lambda_0)}$.
If we have $f(u)\neq 0$ for any $u\in A$,
then $A$ is contained in $A'$.
\end{lem}
\pf
By using Lemma \ref{lem;9.17.19},
we obtain $\Uzero=U^{\prime\,(\lambda_0)}$.
Moreover,
we obtain the two decomposition:
\[
 \Gr^{U}_c\nbigm
=\bigoplus_{
  u\in \nbigk(A,c,\lambda_0)}
 \EE\bigl(-\eigenmap(\lambda,u)\bigr)
=\bigoplus_{
 u\in \nbigk(A',c,\lambda_0)
 }
 \EE\bigl(-\eigenmap(\lambda,u)\bigr).
\]
Thus we obtain the second claim.
\hfill\qed

\vspace{.1in}

The second claim in Lemma \ref{lem;9.17.7}
implies that we can take the unique minimal $A$
if we impose the strictness to $\Gr^{\Uzero}(\nbigm)$.

\begin{lem} \label{lem;9.17.25}
Assume the following:
\begin{itemize}
\item
$\nbigm$ is specializable for
$(X_0,A,\lambda_0,\Uzero,f)$ and
$(X_0,A',\lambda_0,U^{\prime\,(\lambda_0)},f')$.
\item
The map $\eigenmap(\lambda_0):
 \tilde{A}\cup \tilde{A}'\lrarr \cnum$ is injective.
\end{itemize}
Then we have $\Uzero=U^{\prime\,(\lambda_0)}$.

Moreover, if we have
$f(u)\neq 0$ 
for any $u\in A$,
then we have $A\subset A'$.
\end{lem}
\pf
We omit to denote $(\lambda_0)$.
We put $A''=A'\cup A$.
The argument is essentially same as the proof of
Lemma \ref{lem;9.17.19}.
Let us pick a section $g\in U_c$.
We have the real number $d$ determined by the condition
$g\in U'_d-U'_{<d}$.
Then we have the following:
\begin{equation} \label{eq;9.17.20}
 g\cdot b_1\bigl(t\cdot\deldel_t\bigr) \in U'_{<d},
\quad\quad
b_1(s):=
 \prod_{u\in\nbigk(A',\lambda_0,d)}
 \bigl(s+\eigenmap(\lambda,u)
 \bigr).
\end{equation}
Since there exists a positive integer $l$ such that
$U_{c-l}\subset U'_{<d}$,
there exists a finite subset
$S\subset\bigl\{u\in \tilde{A}\,\big|\,\paramap(\lambda_0,u)\leq
c\bigr\}$,
such that the following holds:
\begin{equation} \label{eq;9.17.21}
 g\cdot b_2\bigl(t\cdot\deldel_t\bigr) \in U'_{<d},
\quad\quad
b_2(s):=
\prod_{u\in S}\bigl(s+\eigenmap(\lambda,u)
 \bigr).
\end{equation}
Assume $c<d$, and we will derive a contradiction.
We have $S\cap \nbigk(A',\lambda_0,d)=\emptyset$.
Then we obtain the following,
due to the injectivity assumption:
\begin{equation} \label{eq;9.19.1}
 \bigl\{
 -\eigenmap(\lambda_0,u)\,\big|\,
 u\in S
 \bigr\}
\cap
  \bigl\{
 -\eigenmap(\lambda_0,u)\,\big|\,
 u\in \nbigk(A',\lambda_0,d)
 \bigr\}
=\emptyset.
\end{equation}
It implies $\gcd\bigl(b_1(s),b_2(s)\bigr)=1$.
Then we obtain $f\in U'_{<d}$
from (\ref{eq;9.17.20}) and (\ref{eq;9.17.21}),
but it contradicts our choice of $d$.
Thus we obtain $c\geq d$.
By symmetry we obtain $c=d$,
and thus $U=U'$.

By the same argument as the proof of Lemma \ref{lem;9.17.7},
we obtain the second claim.
\hfill\qed

\subsubsection{A lemma for strict compatibility of the morphism
  and $\Uzero$}

\begin{prop} \label{prop;9.17.30}
Let $\nbigm$ and $\nbign$ be specializable
for $(X_0,A,\lambda_0,\Uzero(\nbigm),f_{\nbigm})$
and $(X_0,A',\lambda_0,\Uzero(\nbign),f'_{\nbign})$
respectively.
Let $\phi:\nbigm\lrarr\nbign$ be a morphism
of $\nbigr$-modules.
Assume $\Gr^{\Uzero}(\nbign)$ is strict.
\begin{enumerate}
\item
 $\phi\bigl( \Uzero_c(\nbigm)\bigr)$
 is contained in $\Uzero_c(\nbign)$.
\item
 Assume that $\Gr^{\Uzero}(\nbigm)$ is strict,
 and that the induced morphism
 $\Gr^{\Uzero}_c(\phi):
 \Gr^{\Uzero}_c(\nbigm)\lrarr \Gr^{\Uzero}_c(\nbign)$
 is strict, i.e.,
 $\Cok(\Gr^{\Uzero}_c(\phi))$ is strict.
Then $\phi$ is strict with respect to the filtrations
$\Uzero(\nbigm)$ and $\Uzero(\nbign)$,
namely,
we have
$\Uzero_c(\nbigm)\cap \Image(\phi)=
\phi\bigl(\Uzero_c(\nbign)\bigr)$.
\end{enumerate}
\end{prop}
\pf
Let us consider the image $\Image(\phi)$.
The good $V$-filtration $\Uzero(\nbigm)$ induces
the good $V$-filtration $U^{(1)}$ on $\Image(\phi)$
via the surjection
$\nbigm\lrarr\Image(\phi)$.
It satisfies the conditions
in Definition \ref{df;9.17.5} for $(X_0,A,\lambda_0)$.

The good $V$-filtration $\Uzero(\nbign)$ induces
the good $V$-filtration $U^{(2)}$ on $\Image(\phi)$
via the inclusion $\Image(\phi)\subset\nbign$.
It satisfies the conditions in Definition \ref{df;9.17.5}
for $(X_0,A',\lambda_0)$.
Moreover $\Gr^{U^{(2)}}(\Image(\phi))$ is strict.
Hence
we obtain $U^{(1)}_c\subset U^{(2)}_c$,
due to Lemma \ref{lem;9.17.19}.
It implies that the morphism $\phi$ preserves the filtration.
Thus we obtain the first claim of Proposition \ref{prop;9.17.30}.

To show the second claim of Proposition \ref{prop;9.17.30},
let us consider the induced morphism
$\phi':
 \Uzero_d/\Uzero_{<c}(\nbigm)\lrarr\Uzero_d/\Uzero_{<c}(\nbign)$.
\begin{lem} \label{lem;9.17.32}
The morphism $\phi'$ is strict and we have the following:
\begin{equation}\label{eq;9.17.31}
  \Image(\phi')\cap 
\Bigl(
\Uzero_{d'}/\Uzero_{<c}(\nbign)\Bigr)
=\phi'\Bigl(\Uzero_{d'}/\Uzero_{<c}(\nbigm)
 \Bigr).
\end{equation}
\end{lem}
\pf
We note that there exists
an open dense subset $Y\subset\Delta(\lambda_0,\epsilon_0)$ 
such that
the restrictions of $U^{(2)}$ and $U^{(1)}$ to
$Y\times X$
are same, which is due to Lemma \ref{lem;9.17.25}.
Thus the equality (\ref{eq;9.17.31}) holds
on $Y\times X_0$.
Because $\Gr^{\Uzero}(\nbigm)$ and $\Gr^{\Uzero}(\nbign)$ are strict,
and because the morphisms $\Gr^{\Uzero}(\phi)$ is strict,
we can derive that the equality (\ref{eq;9.17.31}) holds
on $\nbigx_0(\lambda_0,\epsilon_0)$.
\hfill\qed

\begin{cor} \label{cor;9.17.26}
Let $c$ and $d$ be real numbers such that $c<d$.
Let $h$ be an element of $\Uzero_d(\nbigm)$
such that $\phi(h)\in\Uzero_c(\nbign)$.
Then there exists an element $h_1\in\Uzero_c(\nbigm)$
such that
$\phi(h-h_1)\in \Uzero_{<c}(\nbign)$.
\end{cor}
\pf
It immediately follows from Lemma \ref{lem;9.17.32}.
\hfill\qed

\vspace{.1in}

Since the $V$-filtrations $U^{(i)}$ $(i=1,2)$ on $\Image(\phi)$
are good,
there exists a positive number $l_0$
such that
$U^{(2)}_{c}\subset U^{(1)}_{c+l_0}$ for any real number $c$.
It implies the following:
For any element $g\in U^{(2)}_c$,
there exists an element
$\bar{g}\in \Uzero_{c+l_0}(\nbigm)$
such that $g=\phi(\bar{g})$.

\begin{lem}
There exists an element $\bar{g}_1\in \Uzero_c(\nbigm)$
such that $\phi(\bar{g}-\bar{g}_1)$
is contained in $\Uzero_{c-l_0-1}(\nbign)$.
\end{lem}
\pf
Due to Corollary \ref{cor;9.17.26},
we have an element $\bar{g}_2\in \Uzero_{c}(\nbigm)$
such that
$\phi(\bar{g}-\bar{g}_2)\in \Uzero_{<c}(\nbign)$.
By using Corollary \ref{cor;9.17.26} inductively,
we obtain the element desired.
Note that the set
$\bigl\{d\in\real\,\big|\,\Gr^{\Uzero}_d(\nbigm)\neq 0\bigr\}
 \cup\bigl\{d\in \real\,\big|\,\Gr^{\Uzero}_{d}(\nbign)\neq 0\bigr\}$
is discrete in $\real$.
\hfill\qed

Let us return to the proof of Proposition \ref{prop;9.17.30}.
Since $\phi(\bar{g}-\bar{g}_1)$ is contained in
$U^{(2)}_{c-l_0-1}$, we can pick an element
$\bar{g}_3\in \Uzero_{c-1}(\nbigm)$ such that
$\phi(\bar{g}_3)=\phi(\bar{g}-\bar{g}_1)$.
Then the element
$\bar{g}_3+\bar{g}_1$
is contained in $\Uzero_c(\nbigm)$,
and it satisfies the following:
\[
 \phi(\bar{g}_3+\bar{g}_1)=\phi(\bar{g})=g.
\]
Thus $g$ is contained in $U^{(1)}_c$,
namely we obtain $U^{(2)}=U^{(1)}$.
It means the strictness of $\phi$
with respect to the filtrations
$\Uzero(\nbigm)$ and $\Uzero(\nbign)$.
Therefore the proof of Proposition \ref{prop;9.17.30}
is accomplished.
\hfill\qed

\subsubsection{$\psi_{t,u}^{(\lambda_0)}$ and $\psi_{t,u}$}

\begin{df}
When $\nbigm$ is specializable for $(X_0,A,\lambda_0,\Uzero,f)$
such that $\Gr^{\Uzero}$ is strict,
we put as follows:
\[
 \psi_{t,u}^{(\lambda_0)}(\nbigm)
:=\psi_{t,u}^{U^{(\lambda_0)}}(\nbigm)
\subset
\Gr^{\Ulambdazero}_{\paramap(\lamda_0,u)}.
\]
It is well defined due to Lemma {\rm\ref{lem;9.17.19}}.
Clearly, it is strict.
\hfill\qed
\end{df}

We have the nilpotent map
$t\deldel_t+\eigenmap(\lambda,u)$
on $\psi_{t,u}^{(\lambda_0)}(\nbigm)$.

\begin{lem}
We have the decomposition:
\begin{equation} \label{eq;9.17.38}
 \Gr^{\Ulambdazero}_c(\nbigm)
=\bigoplus_{u\in\nbigk(A,\lambda_0,c)}
 \psi^{(\lambda_0)}_{t,u}(\nbigm).
\end{equation}
\end{lem}
\pf
It is just a reformulation of (\ref{eq;9.17.1}).
\hfill\qed

\begin{lem} \label{lem;9.17.12}
Assume that
$\Delta(\lambda_1,\epsilon_1)\subset\Delta(\lambda_0,\epsilon_0)$.
We have the following:
\begin{equation} \label{eq;9.17.11}
 \psi^{(\lambda_0)}_{t,u}\nbigm_{|\Delta(\lambda_1,\epsilon_1)}
=\psi^{(\lambda_1)}_{t,u}\nbigm.
\end{equation}
\end{lem}
\pf
The construction of $U^{(\lambda_1)}$ from $U^{(\lambda_0)}$
is given in the proof of Lemma \ref{lem;9.17.10}.
Then (\ref{eq;9.17.11}) can be checked easily
by using the uniqueness
of the $V$-filtration $U^{(\lambda_1)}$ such that
$\Gr^{U^{(\lambda_1)}}$ is strict.
\hfill\qed

\vspace{.1in}

\begin{lem}\label{lem;9.17.35}
Assume the following:
\begin{itemize}
\item
$\nbigm$ is specializable for 
$\bigl(X_0,A(\lambda_0),\lambda_0,\Uzero,f^{(\lambda_0)}\bigr)$
for any $\lambda_0\in\cnum_{\lambda}$. 
\item
 $\Gr^{\Uzero}$ is strict for any $\lambda_0\in\cnum_{\lambda}$.
\item
 $f^{(\lambda_0)}(u)\neq 0$ for any $u\in A(\lambda_0)$. 
\end{itemize}
Then the following holds:
\begin{enumerate}
\item
If $\Uzero$ is defined over $\Delta(\lambda_0,\epsilon_0)$
and $U^{(\lambda_1)}$ is defined over $\Delta(\lambda_1,\epsilon_1)$.
Then we have $U^{(\lambda_0)}=U^{(\lambda_1)}$
on $\Delta(\lambda_0,\epsilon_0)\cap\Delta(\lambda_1,\epsilon_1)$.
\item
$A(\lambda_0)$ does not depend on $\lambda_0\in\cnum_{\lambda}$.
\end{enumerate}
\end{lem}
\pf
It immediately follows from Lemma \ref{lem;9.17.7}.
\hfill\qed

\begin{df} \label{df;9.17.36}
If the assumption of Lemma {\rm\ref{lem;9.17.35}} is satisfied,
the $\nbigr$-module $\nbigm$ is called specializable along $X_0$.
In that case, the set $A$ is denoted by
$\KMSoverline(\nbigm,X_0)$ or $\KMSoverline(\nbigm,t)$.
The set $\tilde{A}$ is denoted by $\KMS(\nbigm,X_0)$ or $\KMS(\nbigm,t)$.
\hfill\qed
\end{df}

\begin{rem}
Note that the strictness of $\Gr^{\Uzero}$ $(\lambda_0\in\cnum)$
is contained in Definition {\rm\ref{df;9.17.36}}.
\hfill\qed
\end{rem}

\begin{df}
Assume $\nbigm$ is specializable along $X_0$.
Then $\bigl\{\psi^{(\lambda_0)}_{u,t}\nbigm\,\big|\,\lambda_0\in\cnum\bigr\}$
determines the globally defined $\nbigr_{\nbigx_0}$-module,
due to Lemma {\rm\ref{lem;9.17.12}}.
We denote it by $\psi_{t,u}\nbigm$.
\hfill\qed
\end{df}

Let $\nbigm$ and $\nbign$ be strictly specializable along $X_0$.
Let $\phi:\nbigm\lrarr\nbign$ be a morphism of $\nbigr$-modules.
\begin{lem}
The morphism $\phi$ preserves the $V$-filtrations $\Uzero$ 
at $\lambda_0$ for any $\lambda_0\in\cnum_{\lambda}$.
\end{lem}
\pf
It immediately follows from the first claim in
Proposition \ref{prop;9.17.30}.
\hfill\qed

\vspace{.1in}

Then we obtain the induced morphism 
$\Gr^{\Uzero}(\phi):\Gr^{\Uzero}(\nbigm)\lrarr\Gr^{\Uzero}(\nbign)$.
\begin{lem}
It induces the morphisms
$\psi_{t,u}^{(\lambda_0)}(\phi):
 \psi^{(\lambda_0)}_{t,u}(\nbigm)
 \lrarr\psi^{(\lambda_0)}_{t,u}(\nbign)$.
They can be glued, and we obtain the morphism
$\psi_{t,u}(\phi):\psi_{t,u}(\nbigm)\lrarr\psi_{t,u}(\nbign)$.
\end{lem}
\pf
Since the decomposition (\ref{eq;9.17.38})
is obtained as a generalized eigen decomposition,
the first claim is clear.
The second claim is also clear from the construction
of $U^{(\lambda_1)}$ from $U^{(\lambda_0)}$
given in Lemma \ref{lem;9.17.10}.
\hfill\qed

\vspace{.1in}

Let $\vecdelta_0$ denote the element $(1,0)\in\real\times\cnum$.
Then we have the naturally induced morphisms:
\[
 t:\psi_{t,u}^{(\lambda_0)}\nbigm
   \lrarr
   \psi_{t,u-\vecdelta_0}^{(\lambda_0)}\nbigm.
\]
\[
 \deldel_t:
 \psi_{t,u}^{(\lambda_0)}\nbigm
 \lrarr
 \psi_{t,u+\vecdelta_0}^{(\lambda_0)}\nbigm.
\]
In particular,
we put as follows:
\[
\begin{array}{l}
 \can=\deldel_t:\psi^{(\lambda_0)}_{t,-\vecdelta_0}\nbigm
 \lrarr\psi^{(\lambda_0)}_{t,0}\nbigm,\\
 \mbox{{}}\\
 \var=t:\psi^{(\lambda_0)}_{t,0}\nbigm\lrarr
 \psi^{(\lambda_0)}_{t,-\vecdelta_0}\nbigm.
\end{array}
\]
If $\nbigm$ is specializable along $X_0$,
then we have
$t:\psi_{t,u}\nbigm\lrarr\psi_{t,u-\vecdelta_0}\nbigm$
and
$\deldel_t:\psi_{t,u}\nbigm\lrarr\psi_{t,u+\vecdelta_0}\nbigm$.
In particular,
we have
$\can:\psi_{t,-\vecdelta_0}\nbigm
 \lrarr\psi_{t,0}\nbigm$
and
$\var:\psi_{t,0}\nbigm\lrarr
 \psi_{t,-\vecdelta_0}\nbigm$.

\subsubsection{Strictly specializable}

\begin{df}
An $\nbigr$-module
$\nbigm$ is called strictly specializable along $X_0$,
if the following holds:
\begin{enumerate}
\item
 It is specializable along $X_0$.
\item
 For any $\lambda_0\in\cnum_{\lambda}$ 
 and for any $c<0$,
 the morphism
 $t:U_c^{(\lambda_0)}\nbigm\lrarr U^{(\lambda_0)}_{c-1}\nbigm$
 is isomorphic.
\item
 For any $\lambda_0\in\cnum_{\lambda}$
 and for any $c>-1$,
 the morphism
 $\deldel_t:\Gr^{U^{(\lambda_0)}}_c\nbigm\lrarr
 \Gr^{U^{(\lambda_0)}}_{c+1}\nbigm$ is isomorphic.
\hfill\qed
\end{enumerate}
\end{df}

\begin{df}\mbox{{}}
\begin{itemize}
\item
Let $\nbigm$ and $\nbign$ be strictly specializable
along $X_0$.
A morphism $\phi:\nbigm\lrarr\nbign$ is called
strictly specializable
if $\psi_{t,u}(\phi):\psi_{t,u}\nbigm\lrarr\psi_{t,u}\nbign$ is strict,
i.e.,
the cokernel $\Cok\bigl(\psi_{t,u}(\phi)\bigr)$ is strict.
\item
We have the category with
strictly specializable $\nbigr$-modules along $X_0$
and strictly specializable morphisms.
We denote it by $\nbigs^2(X,t)$.
Note that $X_0=\{t=0\}$.
\item
Let $f$ be a holomorphic function on $X$.
An $\nbigr$-module $\nbigm$ is called strictly specializable along $f$
if $i_{f\,\ast}\nbigm$ is strictly specializable.
\hfill\qed
\end{itemize}
\end{df}
Here $i_f$ denotes the naturally defined inclusion
$X\lrarr X\times\cnum$.

\begin{prop} \label{prop;a11.23.1}
Assume $\nbigm$ is strictly specializable along $X_0$.
\begin{enumerate}
\item It we have a direct sum decomposition $\nbigm=\nbigm_1\oplus\nbigm_2$,
 then $\nbigm_i$ $(i=1,2)$ are strictly specializable.
\item
 Assume that $\nbigm$ is supported in $\nbigx_0$.
 Then we have $U^{(\lambda_0)}_{<0}\nbigm=0$
 for any $\lambda_0\in\cnum_{\lambda}$.
 We also have $\psi_{t,u}\nbigm=0$ if
 $u$ does not contained in $\seisuu_{\geq\,0}\times\{0\}$.
\item \label{number;10.2.1}
 The following conditions are equivalent.
   \begin{itemize}
  \item
    $\var:\psi_{t,0}\nbigm\lrarr\psi_{t,-\vecdelta_0}\nbigm$
    is injective.
  \item
    Let $\nbigm'$ be a submodule of $\nbigm$ such that
    the support of $\nbigm'$ is contained in $\nbigx_0$.
    Then $\nbigm=\nbigm'$ or $0$.
  \item
    Let $\nbigm'$ be a submodule of $\nbigm$ such that
    the support of $\nbigm'$ is contained in $\nbigx_0$.
    Assume that $\nbigm'\in \nbigs^2(X,t)$.
    Then $\nbigm'=\nbigm$ or $0$.
   \end{itemize}
 \item \label{number;10.2.2}
  Assume $\can:\psi_{t,-\vecdelta_0}\nbigm\lrarr\psi_{t,0}\nbigm$ is
      surjective.
 Let $\nbigm''\in \nbigs^2(X,t)$ be a quotient of $\nbigm$ such that
 the support of $\nbigm''$ is contained in $\nbigx_0$.
 Then $\nbigm''=0$.
 \item
  The following conditions are equivalent.
 \begin{itemize}
  \item $\psi_{t,0}\nbigm=\Image\can\oplus\Ker\var$.
  \item We have the decomposition $\nbigm=\nbigm'\oplus\nbigm''$,
 where the support of $\nbigm''$ is contained in $X_0$
 and $\nbigm'$ has neither submodules or quotients
 contained in $\nbigs^2(X,t)$,
 whose support is contained in $X_0$.
 \end{itemize}
 \end{enumerate}
\end{prop}
\pf
The proof of the claims 1. 3. 4, 5.
are same as those the proof of Proposition 3.3.9. in \cite{sabbah}.
Let us see the claim 2.

(i) Since the multiplication
 $t\cdot $ is injective on $U^{(\lambda_0)}_{<0}$,
we obtain $U^{(\lambda_0)}_{<0}=0$.
In particular, $\psi^{(\lambda_0)}_{t,u}\nbigm=0$
if $\paramap(\lambda_0,u)<0$.

(ii)
Assume that $\paramap(\lambda_0,u)\geq 0$ is not integer.
We can take $l\in\seisuu_{>0}$ such that
$-1<\paramap(\lambda_0,u)-l<0$.
Then we obtain the surjection:
\[
 \deldel_t^l:
 \psi^{(\lambda_0)}_{t,u-l\vecdelta_0}\nbigm
\lrarr
 \psi^{(\lambda_0)}_{t,u}\nbigm.
\]
Thus $\psi^{(\lambda_0)}_{t,u}\nbigm=0$ in this case.

(iii)
Assume $u\not\in \seisuu_{\geq 0}\times\{0\}$,
and $\paramap(\lambda_0,u)=0$.
Note that $\eigenmap(\lambda_0,u)\neq 0$.
Then the composite of the morphisms
\[
\begin{CD}
 \psi^{(\lambda_0)}_{t,u}\nbigm
 @>{t}>>
 \psi^{(\lambda_0)}_{t,u-\vecdelta_0}\nbigm
 @>{\deldel_t}>> \psi^{(\lambda_0)}_{t,u}\nbigm
\end{CD}
\]
is isomorphic.
Then we obtain $\psi^{(\lambda_0)}_{t,u}\nbigm $
in this case.

(iv)
If $u\not\in\seisuu_{\geq\,0}\times\{0\}$
and $\paramap(\lambda_0,u)\geq 0$,
then it can be reduced to the case (iii)
by the argument in (ii).

Then we obtain
$V_0\nbigm
=\psi^{(\lambda_0)}_{t,0}\nbigm
=\ker(t:\nbigm\lrarr\nbigm)$.
Since
$\deldel_t^k:\psi^{(\lambda_0)}_{t,0}
\lrarr \psi^{(\lambda_0)}_{t,k\vecdelta_0}$
is isomorphic,
we obtain $\nbigm=i_{+}\psi^{(\lambda_0)}_{t,0}\nbigm$.
\hfill\qed

\vspace{.1in}

Let $\nbigs^2_{X_0}(X,t)$ denote the subcategory
of $\nbigs^2(X,t)$,
whose objects have the supports contained in $\nbigx_0$.
\begin{cor}
We have the equivalence
$ \nbigs^{2}_{\nbigx_0}(X,t)
\simeq
 \bigl(\mbox{\rm strict }\nbigr_{\nbigx_0}\mbox{\rm -modules}\bigr)$.
\hfill\qed
\end{cor}

\begin{df}\mbox{{}}
\begin{itemize}
\item
$\nbigm$ is strictly $S$-decomposable along $X_0$,
if it is strictly specializable along $X_0$
and $\psi_{t,0}=\Image(\can)\oplus\Ker(\var)$.

\item
$\nbigm$ is $S$-decomposable at $P$,
if for any holomorphic function $f$,
and $i_{f\,+}\nbigm$ is strictly $S$-decomposable
at $(x,0)$.

\item
$\nbigm$ is strictly $S$-decomposable 
if $\nbigm$ is strictly $S$-decomposable at any $x\in X$.
\hfill\qed
\end{itemize}
\end{df}

\begin{lem}\label{lem;a11.23.2}
Let $\nbigm$ and $\nbigm'$ be  $\nbigr$-modules,
which are strictly $S$-decomposable along $X_0$.
Let $f:\nbigm'\lrarr\nbigm$ be a morphism.
Assume the following:
\begin{itemize}
\item
For any $\lambda_0\in\cnum_{\lambda}$,
there exists a number $h(\lambda_0)<0$
such that
the induced morphism
$\Vzero_{h(\lambda_0)}(\nbigm')\lrarr \Vzero_{h(\lambda_0)}(\nbigm)$
is isomorphic.
\item
We have $\psi_{t,0}(\nbigm)=\Image(\can)$
and $\psi_{t,0}(\nbigm')=\Image(\can)$.
\end{itemize}
Then $f$ is isomorphic around a neighbourhood of $X_0$.
\end{lem}
\pf
We have only to show that
$\Vzero_h(\nbigm')\lrarr\Vzero_h(\nbigm)$ is isomorphic
for any $h$ and for any $\lambda_0$.

In the case $h<h(\lambda_0)$,
the coincidence $\Vzero_h(\nbigm')\lrarr\Vzero_h(\nbigm)$
follows from
$\Vzero_{h(\lambda_0)}(\nbigm')
\simeq
 \Vzero_{h(\lambda_0)}(\nbigm)$
and 
the uniqueness of the $V$-filtration
whose associated graded module is strict.

In the case $h(\lambda_0)\leq h<0$,
we have a large integer $N$ such that
$h-N<h(\lambda_0)$.
We have the following commutative diagramm:
\[
 \begin{CD}
 \Vzero_h(\nbigm') @>>> \Vzero_h(\nbigm)\\
 @V{t^N}VV @V{t^N}VV \\
 \Vzero_{h-N}(\nbigm') @>>> \Vzero_{h-N}(\nbigm).
 \end{CD}
\]
Since the both of the vertical arrows and the lower horizontal arrow
are isomorphic,
the upper vertical arrow is also isomorphic.

In particular, we know 
$\Vzero_{<0}(\nbigm')\simeq\Vzero_{<0}(\nbigm)$.

Since $\nbigm$ and $\nbigm'$ are strictly $S$-decomposable,
and since
we have
$\psi_{t,0}=\Image(\can)$ for both of $\nbigm$ and $\nbigm'$,
they are generated by $\Vzero_{<0}$.
Thus we obtain the surjectivity of the morphism $f$.
Since we have $\Vzero_{<0}(\ker f)=0$,
the support of $\ker(f)$ is contained in $X_0$.
Then we obtain $\ker(f)=0$ due to Proposition \ref{prop;a11.23.1}.
\hfill\qed

\begin{cor} \label{cor;a11.23.3}
Let $\nbigm$ and $\nbigm'$ be  $\nbigr$-modules,
which are strictly $S$-decomposable along $X_0$,
and $\psi_{t,0}(\nbigm)=\Image(\can)$
and $\psi_{t,0}(\nbigm')=\Image(\can)$ hold.
Let $f:\nbigm'\lrarr\nbigm$ be a morphism.
Assume that the support of the cokernel of $f$
is contained in $X_0$.
Then $f$ is isomorphic.
\end{cor}
\pf
For any $\lambda_0$,
there exists a sufficiently negative number $h(\lambda_0)$
such that $\Vzero_{h(\lambda_0)}\cok(f)=0$.
It means
$\Vzero_{h(\lambda_0)}(\nbigm')\simeq \Vzero_{h(\lambda_0)}(\nbigm)$.
Thus we can apply Lemma \ref{lem;a11.23.2}.
\hfill\qed

%% file: a34.1.tex

\subsubsection{$\tildepsizero(\nbigm)$ and $\tildepsi(\nbigm)$}

Let $\nbigm$ be an $\nbigr_{\nbigx}$-module,
which is strictly specializable along $X_0$.

\begin{df}
Let $u=(a,\alpha)$ be an element of $\KMS(\nbigm,t)$
such that $u\not\in\seisuu_{\geq\,0}\times \{0\}$.
The $\nbigr_{\nbigx_0}$-module $\tildepsizero_{t,u}(\nbigm)$
is defined as follows:
Let us pick an integer $b$ such that
$\paramap(\lambda_0,u-b\cdot\vecdelta_0)
=\paramap(\lambda_0,u)-b<0$.
Then we put 
$\tildepsizero_{t,u}(\nbigm):=
 \psizero_{t,u-b\cdot\vecdelta_0}(\nbigm)$.

It is well defined in the following sense:
For any non-negative integer $N$,
we have the canonical isomorphism:
\[
 t^N:
 \psizero_{t,u-b\cdot\vecdelta_0}
\lrarr
 \psizero_{t,u-(b+N)\cdot\vecdelta_0}.
\]
\hfill\qed
\end{df}

Let $u$ be an element of $\KMS(\nbigm,t)$.
Let us consider the following set:
\[
 S(u):=
 \bigl\{
 \lambda\in\cnum^{\ast}\,\big|\,
 \exists b\in \seisuu,\,\,\mbox{\rm s.t.}\,\,
 \eigenmap\bigl(\lambda,u-b\cdot\vecdelta_0\bigr)=0,\,\,
 \paramap\bigl(\lambda_0,u-b\cdot\vecdelta_0\bigr)\geq 0
 \bigr\}.
\]
\begin{lem}
The set $S(u)$ is discrete in $\cnum_{\lambda}$.
\end{lem}
\pf
The set $S(u)$ is contained in the following set:
\[
 \bigcup_{b\in \seisuu}
 \bigl\{
 \lambda\,\big|\,
 \eigenmap(\lambda,u-b\vecdelta_0)=0
 \bigr\}.
\]
Then the discreteness of $S(u)\cap \cnum^{\ast}$
in $\cnum^{\ast}$ can be shown
by an argument similar to the proof of
Lemma \ref{lem;10.11.21}.

In a neighbourhood $U$ of $\lambda=0$,
$S(u)\cap U$ is contained in the following finite set:
\[
 \bigcup_{0\leq b\leq \paramap(0,u)+1}
 \bigl\{
 \lambda\,\big|\,
 \eigenmap(\lambda,u-b\vecdelta_0)=0
 \bigr\}.
\]
Then we obtain the discreteness of $S(u)$
in $\cnum_{\lambda}$.
\hfill\qed

\begin{lem} \label{lem;9.25.10}
Let $\lambda_0$ be an element of $\cnum^{\ast}-S$.
Then we have the canonical isomorphism
$\psizero_{t,u}\simeq \tildepsizero_{t,u}$.
\end{lem}
\pf
Since $\lambda_0$ is an element of $\cnum^{\ast}$,
the eigenvalues of the endomorphism $s=t\cdot\deldel_t$
cannot be $0$.
It implies that the morphisms
$t:\psizero_{t,u}\lrarr\psizero_{t,u-\vecdelta_0}$
are isomorphic for any $u$.
Thus we have the isomorphism
$t^N:\psizero_{t,u}\lrarr \psizero_{t,u-N\cdot\vecdelta_0}$.
It gives the isomorphism desired.
\hfill\qed

\vspace{.1in}

We have the following straightforward corollary of
Lemma \ref{lem;9.25.10}.

\begin{cor} \label{cor;9.25.11}
Let $u$ be an element of $\KMS(\nbigm,t)$.
Let $\lambda_0$ and $\lambda_1$ be a point.
Assume that
$\Delta(\lambda_0,\epsilon_0)\cap\Delta(\lambda_1,\epsilon_1)
\subset \cnum-S(u)$.
Then we have the canonical isomorphism
$\tildepsi^{(\lambda_0)}_{t,u}(\nbigm)
\simeq
 \tildepsi^{(\lambda_1)}_{t,u}(\nbigm)$.
\hfill\qed
\end{cor}

\begin{df}
The $\nbigr_{\nbigx_0}$-module
$\tildepsi_{t,u}(\nbigm)$
is defined by
$\tildepsi_{t,u}(\nbigm)_{|\Delta(\lambda_0,\epsilon_0)}:=
 \tildepsizero_{t,u}(\nbigm)$.

It is well defined due to Corollary {\rm\ref{cor;9.25.11}}.
\hfill\qed
\end{df}

\begin{lem}
We have the canonical inclusion
$\psizero_{t,u}\lrarr \tildepsizero_{t,u}$.
\end{lem}
\pf
Recall that we have assumed 
$u\not\in\seisuu_{\geq\,0}\times\{0\}$.
Then the induced morphism
$t:\psizero_{t,u}\lrarr\psizero_{t,u-\vecdelta_0}$
is injective.
Thus we have the injection
$t^N:\psizero_{t,u}\lrarr\psizero_{t,u-N\cdot\vecdelta_0}$.
It gives the desired inclusion.
\hfill\qed

\begin{cor}
We have the canonical inclusion
$\psi_{t,u}(\nbigm)\lrarr \tildepsi_{t,u}(\nbigm)$.
\hfill\qed
\end{cor}

%% file: a74.1.tex

We recall the specialization of sesqui-linear pairing
introduced by Sabbah, with minor generalization.
We recommend the reader to read
the sections 1.5--1.7 and the section 3.5--3.7
of \cite{sabbah}.
We consider the left $\nbigr$-modules in this subsection.

%% file: 30.1.tex

\subsubsection{$\overline{\nbigr}$ and $\overline{\nbigr}$-module}

We put $\deldelbar_i:=-\lambda^{-1}\cdot\delbar_i$,
and
$\overline{\nbigr}_{\nbigx}:=
 \nbigo_{\nbigx^{\dagger\,\shikaku}}[\deldelbar_i]$,
where $\nbigx^{\dagger}=X^{\dagger}\times\cnum^{\ast}$.
We use the map $\sigma:\cnum^{\ast}\lrarr\cnum^{\ast}$
given by $\sigma(\lambda)=-\overline{\lambda}^{-1}$.

Let $U$ be a subset of $\cnum_{\lambda}^{\ast}$.
Let $\nbigm$ be a left $\nbigr$-module over $X\times U$.
Then the left $\overline{\nbigr}$-module structure
on $\sigma^{\ast}\nbigm$ on $X^{\dagger}\times\sigma(U)$ is given.
Let $f$ be a section of $\nbigo_{X^{\dagger}\times \sigma(U)}$,
and $v$ be a section of $\sigma^{\ast}(\nbigm)$.
\[
 f\bullet v:=\overline{\sigma^{\ast}(f)}\cdot v,
\quad\quad
 \deldelbar_i\bullet v:=
 \deldel_i\cdot v.
\]
Note the following relation:
\begin{multline}
 \deldelbar_i\bullet  (f\bullet v)
=\deldelbar_i\bullet
 (\overline{\sigma^{\ast}(f)}\cdot v)
=\deldel_i\cdot(\overline{\sigma^{\ast}(f)}\cdot v)
=\overline{\sigma^{\ast}(f)}\cdot \deldel_i\cdot v
+(\lambda\del_i\overline{\sigma^{\ast}(f)})\cdot v\\
=f\bullet (\deldelbar_i\bullet v)
+ \overline{\sigma^{\ast}
 (-\lambda^{-1}\cdot \delbar_i f)} 
 \cdot v
=f\bullet (\deldelbar_i\bullet v)
+\overline{\sigma^{\ast}
 \bigl(\deldelbar_i(f)
 \bigr)}
 \cdot v
=f\bullet (\deldelbar_i\bullet v)
+\bigl(\deldelbar_i(f)\bigr)\bullet v.
\end{multline}
Thus we obtain the well defined left $\overline{\nbigr}$-module
structure on $\sigma^{\ast}\nbigm$.
We often use the notation $\overline{\nbigm}$
instead of $\sigma^{\ast}\nbigm$.

\vspace{.1in}
\noindent
{\bf Notation}
We use the following notation:
Let $\pi_{U}:X\times U\lrarr X$ denote the projection.
Then we put $\nbigm_U:=\pi_{U\,\ast}\nbigm$,
and $\overline{\nbigm}_U:=\pi_{U\,\ast}\sigma^{\ast}\nbigm$.

\subsubsection{Preliminary I}

Let $\lambda_0$ be a point of $\cnum_{\lambda}$.
In this subsubsection, $U$ denotes an open subset
$\Delta(\lambda_0,\epsilon_0)$
for some small positive number $\epsilon_0$.
Let us consider the case $X=X_0\times\cnum$.
Let $t$ be the coordinate of $\cnum$.
Let $\nbigm'$ and $\nbigm''$ be objects of $\nbigs^2(X,t)$.
Let $C_U:\nbigm'_U\lrarr \overline{\nbigm''}_U\lrarr \distribution^U_{X}$
be a sesqui-linear pairing.
We recall the construction of Sabbah to obtain
the specialization along $t$:
\[
 \psi_{t,u}^{(\lambda_0)}C_U:
 \psi_{t,u}^{(\lambda_0)}\nbigm'_{U}
 \otimes
 \overline{\psi_{t,u}^{(\sigma(\lambda_0))}\nbigm''}_U
\lrarr
 \distribution^{U}_{X_0}.
\]

Let $W_0$ be an open subset of $X_0$.
Let $m$ be a section of $\nbigm'_{U}$ 
and $\mu$ be a section of $\nbigm''_U$
on $W=W_0\times \Delta_t$.
Let us pick a $C^{\infty}$ $(n-1,n-1)$-form $\phi$ on $W_0$
whose support is compact,
and a $C^{\infty}$-function $\chi$ on $\Delta_t$
such that $\chi=1$ around the origin $O\in\Delta$
and that the support of $\chi$ is compact.

For any integer $k\in\seisuu$,
we put as follows:
\[
 \nbigi^{(k)}_{C(m,\bar{\mu}),\phi}(s)
:=\left\{
 \begin{array}{ll}
 \bigl\langle C(m,\bar{\mu}),\,\,|t|^{2s}t^k\cdot\chi(t)\cdot\phi\wedge
 \frac{i}{2\pi}dt\wedge d\bar{t}\bigr\rangle,
 & (k\geq 0),\\
\mbox{{}}\\
 \bigl\langle
 C(m,\bar{\mu}),\,\,|t|^{2s}\bar{t}^{|k|}\cdot\chi(t)\cdot\phi\wedge
 \frac{i}{2\pi}dt\wedge d\bar{t}\bigr\rangle,
 & (k<0).
 \end{array}
 \right.
\]

Then $\nbigi^{(k)}_{C(m,\bar{\mu}),\phi}(s)$
is a $H(\AAA\cap U)$-valued holomorphic function
defined on the half plane
$\bigl\{s\in\cnum\,\big|\,\Re(s)>\sigma_0-2^{-1}|k|\bigr\}$,
where $\sigma_0$ denotes some real number.

\begin{lem}\label{lem;9.17.50}
In the case $k\geq 0$,
we have the following:
\begin{equation}\label{eq;9.17.42}
 \lambda\cdot
 \bigl(
 s+k+\lambda^{-1}\eigenmap(\lambda,u)
 \bigr)\cdot
 \nbigi^{(k)}_{C(m,\bar{\mu}),\phi}(s)
=\nbigi^{(k)}_{C(m',\bar{\mu}),\phi}
+F.
\end{equation}
Here we put $m':=(-\deldel_tt+\eigenmap(\lambda,u))\cdot m$,
and $F$ denotes an entire function of the variable $s$.
\end{lem}
\pf
Let us consider the following:
\begin{equation} \label{eq;9.17.40}
 \Big\langle
 \bigl(-\deldel_tt+\eigenmap(\lambda,u)\bigr)\cdot
 C(m,\bar{\mu}),\,
 |t|^{2s}t^k\phi\!\cdot\!\chi(t)\!\cdot\!\frac{i}{2\pi}dt\wedge d\bar{t}
 \Big\rangle
=\Big\langle
 C(m,\bar{\mu}),\,\,\,
  \bigl(t\deldel_t+\eigenmap(\lambda,u)\bigr)\cdot
 |t|^{2s}t^k\phi\!\cdot\!\chi(t)\!\cdot\!
 \frac{i}{2\pi}dt\wedge d\bar{t}
 \Big\rangle.
\end{equation}
The left hand side of (\ref{eq;9.17.40}) can be rewritten as follows:
\[
 {\rm L.H.S.}=
 \Big\langle
 C\bigl((-\deldel_tt+\eigenmap(\lambda,u))\cdot m,
 \bar{\mu}\bigr),\,\,
 |t|^{2s}t^k\phi\cdot\chi(t)\frac{i}{2\pi}dt\wedge d\bar{t}
 \Big\rangle
=\Big\langle
 C(m',\overline{\mu}),\,\,
 |t|^{2s}t^k\cdot\phi\cdot\chi(t)\frac{i}{2\pi}dt\wedge d\bar{t}
 \Big\rangle.
\]
The right hand side of (\ref{eq;9.17.40}) can be rewritten as follows:
\begin{multline} \label{eq;9.17.41}
 {\rm R.H.S.}=
 \Big\langle
 C(m,\bar{\mu}),\,\,
 \bigl(
 (t\deldel_t+\eigenmap(\lambda,u))\cdot
 |t|^{2s}t^k
 \bigr)
 \cdot\phi\cdot
 \chi(t)\cdot\frac{i}{2\pi}dt\wedge d\bar{t}
 \Big\rangle \\
+\Big\langle
 C(m,\bar{\mu}),\,\,
 |t|^{2s}t^{k+1}\cdot\phi\cdot\lambda\cdot \del_t\chi(t)\cdot
 \frac{i}{2\pi}dt\wedge d\bar{t}
 \Big\rangle.
\end{multline}

Since we have $\del_t\chi(t)=0$ around $t=0$,
the second term in (\ref{eq;9.17.41}) is entire.
The first term in (\ref{eq;9.17.41}) is as follows:
\[
 \bigl(
 (s+k)\cdot \lambda+\eigenmap(\lambda,u)
 \bigr)\cdot
 \Big\langle
 C(m,\bar{\mu}),\,\,
 |t|^{2s}\cdot t^k\cdot \phi\cdot\chi(t)
 \frac{i}{2\pi}dt\wedge d\bar{t}
 \Big\rangle.
\]
Then (\ref{eq;9.17.42}) follows immediately.
\hfill\qed

\begin{lem}
In the case $k\geq 0$,
we have the following equality:
\begin{equation} \label{eq;9.17.43}
 -\lambda^{-1}\cdot
\Bigl(s+k+\frac{\eigenmap(\lambda,u)}{\lambda}\Bigr)\cdot
\nbigi^{(-k)}_{C(m,\bar{\mu}),\phi}(s)
=\nbigi^{(-k)}_{C(m,\bar{\mu}'),\phi}
+F.
\end{equation}
Here we put $\mu'=\bigl(-\deldel_tt+\eigenmap(\lambda,u)\bigr)\cdot\mu$,
and $F$ denotes an entire function of $s$.
\end{lem}
\pf
Similarly we consider the following:
\begin{multline}\label{eq;9.17.45}
 \Big\langle
 \bigl(-\deldelbar_t\bar{t}
+\sigma^{\ast}\overline{\bigl(\eigenmap(\lambda,u)\bigr)}\bigr)\cdot
 C(m,\bar{\mu}) ,\,\,
 |t|^{2s}\bar{t}^k\!\cdot\! \phi\!\cdot\!\chi(t)\!\cdot\!
 \frac{i}{2\pi}dt\wedge d\bar{t}
 \Big\rangle \\
=\Big\langle
 C(m,\bar{\mu}),\,\,\,
 \bigl(
 \bar{t}\deldelbar_t+\sigma^{\ast}\overline{\bigl(
 \eigenmap(\lambda,u)\bigr)}
 \bigr)\cdot
 |t|^{2s}\bar{t}^k\!\cdot\!\phi\!\cdot\!\chi(t)
  \frac{i}{2\pi}dt\wedge d\bar{t}
 \Big\rangle.
\end{multline}

The left hand side of (\ref{eq;9.17.45}) is as follows:
\[
 {\rm L.H.S.}=
 \Big\langle
 C\bigl(m,\,\,
 \overline{(-\deldel_tt+\eigenmap(\lambda,u))\!\cdot\! \mu}\bigr),\,\,
 |t|^{2s}\bar{t}^k\!\cdot\!\phi\!\cdot\!\chi(t)
 \frac{i}{2\pi}dt\wedge d\bar{t}
 \Big\rangle.
\]
The right hand side of (\ref{eq;9.17.45}) is as follows:
\begin{multline} \label{eq;9.17.46}
 {\rm R.H.S.}=
 \Big\langle
 C(m,\bar{\mu}),\,\,\,
 \bigl(
\bigl( \bar{t}\deldelbar_t+
 \sigma^{\ast}\bigl(\overline{\eigenmap(\lambda,u)}\bigr)
 \bigr)
\cdot
 |t|^{2s}\bar{t}^k\bigr)
\cdot\phi\cdot\chi(t) \cdot
 \frac{i}{2\pi}dt\wedge d\bar{t}
 \Big\rangle \\
-\lambda^{-1}
 \Big\langle
 C(m,\bar{\mu}),\,\,\,
 |t|^{2s}\bar{t}^{k+1}\cdot\phi\cdot
 \delbar_t \chi(t)\cdot
 \frac{i}{2\pi}dt\wedge d\bar{t}
 \Big\rangle.
\end{multline}
The second term of (\ref{eq;9.17.46}) is entire.
To see the first term of (\ref{eq;9.17.46}),
note the following:
\[
 \bigl(
 \bar{t}\deldelbar_t
 +\sigma^{\ast}\bigl(\overline{\eigenmap(\lambda,u)}\bigr)
 \bigr)\cdot
 |t|^{2s}\bar{t}^k
= \bigl(
 -\lambda^{-1}\cdot(s+k)
+\sigma^{\ast}\bigl(\overline{\eigenmap(\lambda,u)}\bigr)
 \bigr)\cdot
 |t|^{2s}\bar{t}^k
\]
Here we have $\sigma^{\ast}(\bar{\lambda})=-\lambda^{-1}$.
Thus it is same as the following:
\[
 -\lambda^{-1}\cdot |t|^{2s}\bar{t}^k
\Big(s+k
+\sigma^{\ast}\overline{\left(\frac{\eigenmap(\lambda,u)}{\lambda}\right)}
\Big)
=-\lambda^{-1}\cdot|t|^{2s}\bar{t}^k\cdot
\Big(s+k+\frac{\eigenmap(\lambda,u)}{\lambda}\Big)
\]
Here we have used the following equality:
\[
 \sigma^{\ast}
 \overline{\left(\frac{\eigenmap(\lambda,u)}{\lambda}\right)}
=\bar{\alpha}\cdot(-\lambda)-a-\alpha\cdot(-\lambda^{-1})
=-\bar{\alpha}\cdot\lambda-a+\alpha\cdot\lambda^{-1}
=\frac{\eigenmap(\lambda,u)}{\lambda}.
\]
Then (\ref{eq;9.17.43}) follows immediately.
\hfill\qed

\subsubsection{Preliminary II}
\label{subsubsection;a11.18.2}

Let $m$ be an element of $U^{(\lambda_0)}_c\nbigm$
such that $0\neq \pi_c(m)\in \psi^{(\lambda_0)}_{t,u}$
via the projection $\pi_c:\Uzero_c\nbigm\lrarr \Gr^{\Uzero}_c\nbigm$.
Let $b_m(s)$ be the Bernstein polynomial of $m$ at $\lambda_0$,
i.e. $b_m(-\deldel_t\cdot t)m\in V_{-1}\nbigr\cdot m$.
Then $b_m(s)$ is of the following form:
\[
 b_m(s)=
 \bigl(s+\eigenmap(\lambda,u_0)\bigr)^{\nu(u_0)}
\cdot
 \prod_{u\in S_0}
 \bigl(s+\eigenmap(\lambda,u)\bigr)^{\nu(u)}.
\]
Here  $S_0$ denotes a finite subset of $\real\times\cnum$
such that $\paramap(\lambda_0,u)<c$ for any $u\in S_0$.
Then we put as follows for any positive integer $\sigma$:
\[
 B_m^{(\sigma)}(s):=
 \prod_{\nu=0}^{\sigma}
 b_m(s+\nu\lambda).
\]

\begin{lem} \label{lem;9.19.4}
There exists a finite subset $S_1(\sigma)\subset \real\times\cnum$
such that the following holds:
\begin{itemize}
\item
 $B_m^{(\sigma)}(s)=\prod_{u\in S_1(\sigma)}
  \bigl(s+\eigenmap(\lambda,u)\bigr)^{\nu'(u)}$.
\item
 For any $u\in S_1(\sigma)$, we have $\paramap(\lambda_0,u)\leq c$.
 If $\paramap(\lambda_0,u)=c$, then $u=u_0$.
\end{itemize}
Moreover,
there exists a positive number $C$
such that $\nu'(u)\leq C$ for any $u\in\bigcup_{\sigma}S_1(\sigma)$.
\end{lem}
\pf
It is clear from our construction.
\hfill\qed

\begin{lem} \label{lem;9.17.51}
We have the following equality:
\begin{equation} \label{eq;9.19.3}
 \Bigl(
 \prod_{u\in S_1(\sigma)}
 \lambda\cdot\bigl(s+k+\lambda^{-1}\cdot\eigenmap(\lambda,u)\bigr)^{\nu'(u)}
 \Bigr)
 \cdot \nbigi^{(k)}_{C(m,\bar{\mu}),\phi}(s)
=\nbigi^{(k)}_{C(m',\bar{\mu}),\phi}
+\mbox{\rm an entire function}.
\end{equation}
Here $m'= B_m^{(\sigma)}(-\deldel_tt)\cdot m$.
The first term in the right hand side is
holomorphic
on the half plane
$\bigl\{s\in\cnum\,\big|\,\Re(s)>\sigma_0-\sigma-2^{-1}|k|\bigr\}$.
\end{lem}
\pf
The equality (\ref{eq;9.19.3}) follows from Lemma \ref{lem;9.17.50}.
By the construction of $B_m^{(\sigma)}$,
we have the following for some $P\in V_0(\nbigr_{\nbigx})$:
\[
 B_m^{(\sigma)}(-\deldel_tt)\cdot m
=t^{\sigma}P\cdot m.
\]
Hence $\nbigi^{(k)}_{C(m',\bar{\mu}),\phi}$ 
in Lemma \ref{lem;9.17.51} is holomorphic 
on the half plane
$\bigl\{s\in\cnum\,\big|\,\Re(s)>\sigma_0-\sigma-2^{-1}|k|\bigr\}$.
\hfill\qed

\vspace{.1in}

Let $Z(f)$ denote the zero set of a holomorphic function $f$.
\begin{lem}\mbox{{}}\label{lem;a11.18.1}
We regard $\nbigi^{(k)}_{C(m,\bar{\mu}),\phi}(s)$
as a `function' of $(s,\lambda)\in \cnum_s\times U$.
\begin{itemize}
\item
$\nbigi^{(k)}_{C(m,\bar{\mu}),\phi}$ is meromorphic 
on $\cnum_s\times U$.
\item
There exists a discrete subset $S_2$
of $\real\times\cnum$ such that
the pole of $\nbigi^{(k)}_{C(m,\bar{\mu}),\phi}$
is contained in the following:
\[
 \bigcup_{u\in S_2}
 Z\bigl(s+k+\lambda^{-1}\eigenmap(\lambda,u)\bigr).
\]
The order of the poles are bounded.
For any element $u\in S_2$,
we have $\paramap(\lambda_0,u)\leq c$.
If $\paramap(\lambda_0,u)=c$,
then $u=u_0$.
\end{itemize}
\end{lem}
\pf
It immediately follows from Lemma \ref{lem;9.19.4}
and Lemma \ref{lem;9.17.51}.
\hfill\qed

\subsubsection{Preliminary III}
\label{subsubsection;a11.18.3}

Let $\mu$ be an element of $U^{(\sigma(\lambda_0))}_d\nbigm$
such that
$0\neq \pi_d(\mu)\in \psi_{t,u_1}^{(\sigma(\lambda_0))}\nbigm$
via the projection
$U^{(\sigma(\lambda_0))}_d\nbigm
\lrarr
\Gr^{U^{(\sigma(\lambda_0))}}_d\nbigm$.
Let $b_{\mu}$ be a Bernstein polynomial of $\mu$
at $\sigma(\lambda_0)$.
Then it is of the following form:
\[
 b_{\mu}(s)=
 \bigl(s+\eigenmap(\lambda,u_1)\bigr)^{\nu(u_1)}\cdot
 \prod_{u\in S_3}\bigl(s+\eigenmap(\lambda,u)\bigr)^{\nu(u)}.
\]
Here $S_3$ is a subset of $\real\times\cnum$.
For any element $u\in S_3$, we have $\paramap(-\sigma({\lambda}_0),u)<d$.
Then we put as follows for any positive integer $\sigma$:
\[
 B_{\mu}^{(\sigma)}(s)
:=\prod_{\nu=0}^{\sigma}
 b_{\mu}(s+\nu\lambda).
\]
The following lemma is clear.
\begin{lem}
There exists a finite subset $S_4(\sigma)\subset\real\times\cnum$
satisfying the following:
\begin{itemize}
\item
$ B^{(\sigma)}_{\mu}(s)
=\prod_{u\in S_4(\sigma)}\bigl(s+\eigenmap(\lambda,u)\bigr)^{\nu'(u)}$.
\item
For any element $u\in S_4(\sigma)$,
we have $\paramap(\sigma({\lambda}_0),u)\leq d$.
If $\paramap(\sigma({\lambda}_0),u)=d$,
then $u=u_1$.
\end{itemize}
Moreover,
there exists a positive number $C$
such that $\nu'(u)\leq C$
for any elements $u\in \bigcup_{\sigma}S_4(\sigma)$.
\hfill\qed
\end{lem}

\begin{lem}
$\prod_{u\in S_1}
 \bigl(
 -\lambda^{-1}\cdot
 (s+k+\lambda^{-1}\eigenmap(\lambda,u))^{\nu'(u)}
 \bigr)\cdot
 \nbigi^{(-k)}_{C(m,\bar{\mu}),\phi}(s)$
is holomorphic on the half plane
$\bigl\{s\in\cnum\,\big|\,\Re(s)>\sigma_0-\sigma-2^{-1}|k|\bigr\}$.
\end{lem}
\pf
It can be shown by an argument similar to the proof
of Lemma \ref{lem;9.17.51}.
\hfill\qed

\begin{cor}\label{cor;04.2.17.2}
$\nbigi^{(-k)}_{C(m,\bar{\mu}),\phi}$ is meromorphic
on $\cnum\times U$.
There exists a discrete subset $S_5\subset\real\times\cnum$
such that
the pole of $\nbigi^{(-k)}_{C(m,\bar{\mu}),\phi}$
is contained in
$\bigcup_{u\in S_5}Z\bigl(s+k+\lambda^{-1}\eigenmap(\lambda,u)\bigr)$.
The orders of the poles are bounded.
For any $u\in S_5$, we have
$\paramap(\sigma(\lambda_0),u)\leq d$.
If $\paramap(\sigma({\lambda}_0),u)=d$,
then $u=u_0$.
\end{cor}
\pf
Similar to Lemma \ref{lem;a11.18.1}.
\hfill\qed

\subsubsection{The construction of the specialization $\psizero_{t,u}C$}

Let $[m]$ be a section of $\psizero_{t,u_0}\nbigm'$ on $W_0\times U$,
and $m$ be a section of $\nbigm'$ on $W_0\times\Delta_t\times U$
such that $\pi_c(m)=[m]$.
Here we put $c=\paramap(\lambda_0,u_0)$
and $\pi_c$ denotes the projection 
$\Uzero_c\nbigm'\lrarr \Gr^{\Uzero}_c\nbigm'$,
as in the subsubsection \ref{subsubsection;a11.18.2}.
Let $[\mu]$ be a section of $\psi^{(\sigma(\lambda_0))}_{t,u_0}\nbigm''$
on $W_0\times\sigma(U)$, and $\mu$ be a section of $\nbigm''$
on $W_0\times\Delta_t\times\sigma(U)$
such that $\pi_d(\mu)=[\mu]$.
Here we put $d=\paramap\bigl(\sigma(\lambda_0),u_0\bigr)$
and $\pi_d$ denotes the projection
$U^{(\sigma(\lambda_0))}_d\nbigm''
 \lrarr \Gr^{U^{(\sigma(\lambda_0))}}_d\nbigm''$
as in the subsubsection \ref{subsubsection;a11.18.3}.

Then we put as follows:
\begin{equation}\label{eq;9.19.5}
\big\langle
 \psi^{(\lambda_0)}_{t,u_0}C\bigl([m],[\bar{\mu}]\bigr),\,
 \phi
\big\rangle
:=
 \Res_{Z(s+\eigenmap(\lambda,u_0))}
 \bigl(
 \nbigi^{(0)}_{C(m,\bar{\mu}),\phi}(s)
 \bigr).
\end{equation}
Here the residue at $Z(s+\eigenmap(\lambda,u_0))$
means the coefficient of $(s+\eigenmap(\lambda,u_0))^{-1}$
for the development $\sum a_i\cdot (s+\eigenmap(\lambda,u_0))^i$.

Recall that we have a discrete subset $S$ of $\real\times\cnum$
such that
the poles of $\nbigi^{(0)}_{C(m,\bar{\mu}),\phi}(s)$
is contained in
$\bigcup_{u\in S}Z(s+\lambda^{-1}\cdot e(\lambda,u))$.
\begin{lem} \label{lem;04.2.17.5}
We may assume
$\eigenmap(\lambda_0,u)\neq \eigenmap(\lambda_0,u_0)$
for any $u\in S-\{u_0\}$.
\end{lem}
\pf
We have the following equality in general:
\begin{multline}\label{eq;04.2.17.3}
 \paramap(\lambda_0,u)
+\paramap(\sigma(\lambda_0),u)
=a+2\Re(\lambda_0\cdot\bar{\alpha})
+a+2\Re(-\overline{\lambda}_0^{-1}\cdot\bar{\alpha})
=2\Re\Bigl(
 a+\lambda_0\cdot\bar{\alpha}
-\lambda_0^{-1}\cdot \alpha
 \Bigr) \\
=2\Re\bigl(\lambda_0^{-1}\cdot\eigenmap(\lambda_0,u)\bigr).
\end{multline}
Assume that $\eigenmap(\lambda_0,u)=\eigenmap(\lamda_0,u_0)$.
Then we obtain the equality from (\ref{eq;04.2.17.3}):
\begin{equation}\label{eq;04.2.17.1}
 \paramap(\lambda_0,u)+\paramap(\sigma(\lambda_0),u)
=\paramap(\lambda_0,u_0)+\paramap(\sigma(\lambda_0),u_0).
\end{equation}

In the case $\paramap(\lambda_0,u)>\paramap(\lambda_0,u_0)$,
we can exclude $u$ from $S$
because of Lemma \ref{lem;a11.18.1}.
In the case $\paramap(\lambda_0,u)<\paramap(\lambda_0,u_0)$,
we have $\paramap(\sigma(\lambda_0),u)>\paramap(\sigma(\lambda_0),u_0)$
due to the equality (\ref{eq;04.2.17.1}).
Then we can exclude $u$ due to 
Corollary \ref{cor;04.2.17.2}.
\hfill\qed

\begin{lem}
The right hand side of {\rm(\ref{eq;9.19.5})} gives
a holomorphic function on a neighbourhood of $\lambda_0$ in $U$.
\end{lem}
\pf
It immediately follows from Lemma \ref{lem;04.2.17.5}.
\hfill\qed

\begin{lem} \label{lem;a11.18.5}
{\rm(\ref{eq;9.19.5})} is well defined.
\end{lem}
\pf
Let $m_1$ denote another lift of $[m]$.
Then the non-trivial pole of $\nbigi^{(0)}_{C(m-m_1,\bar{\mu}),\phi}$
is contained in the following:
\[
 \bigcup_{\paramap(\lambda_0,u)<c}Z(s+\eigenmap(\lambda,u)).
\]
Thus the residue at $Z(s+\eigenmap(\lambda,u_0))$ is $0$.
Hence (\ref{eq;9.19.5}) is independent of a choice of the lift of $[m]$.
Similarly, it can be shown that 
(\ref{eq;9.19.5}) is independent of a choice of the lift of $[\mu]$.
\hfill\qed

\vspace{.1in}
Thus we obtain the specialization morphism
\[
 \psizero_{t,u}C:
 \psizero_{t,u}\nbigm'_U
\otimes
 \overline{\psi_{t,u}^{(\sigma(\lambda_0))}\nbigm''}_U
\lrarr\distribution^U_{X_0}.
\]

We put $N=-\deldel_tt+\eigenmap(\lambda,u_0)$,
which induces the  nilpotent map on
$\psi^{(\lambda_0)}_{t,u_0}\nbigm$
and $\psi_{t,u_0}^{(\sigma(\lambda_0))}\nbigm$.

\begin{lem} \label{lem;9.17.60}
We have the following equality:
\[
 \psi^{(\lambda_0)}_{t,u_0}C\bigl(N[m],[\bar{\mu}]\bigr)
=(i\lambda)^2\cdot
  \psi^{(\lambda_0)}_{t,u_0}C\bigl([m],\bar{N}[\bar{\mu}]\bigr).
\]
\end{lem}
\pf
By definition, we have the following:
\[
 \big\langle
 \psi^{(\lambda_0)}_{t,u_0} C(N[m],[\bar{\mu}]),\,
 \phi
 \big\rangle
=\Res_{Z(s+\eigenmap(\lambda,u_0))}
 \Big\langle
 C(Nm,\bar{\mu}),
 |t|^{2s}\cdot\phi\cdot\chi\frac{i}{2\pi}dt\wedge d\bar{t}
 \Big\rangle
\]
We have the following equality:
\begin{multline}
 \Big\langle
 C(Nm,\bar{\mu}),
 |t|^{2s}\cdot\phi\cdot\chi\frac{i}{2\pi}dt\wedge d\bar{t}
 \Big\rangle 
=\Big\langle
 C(m,\bar{\mu}),\,\,\,
 \bigl(t\deldel_t+\eigenmap(\lambda,u)\bigr)
 \cdot
 (|t|^{2s}\cdot\phi\cdot\chi)
 \frac{i}{2\pi}dt\wedge d\bar{t}
 \Big\rangle \\
=\Big\langle
 C(m,\bar{\mu}),\,\,
\bigl(
 (t\deldel_t+\eigenmap(\lambda,u))
 |t|^{2s}\bigr)\cdot
 \phi\cdot\chi
 \frac{i}{2\pi}dt\wedge d\bar{t}
 \Big\rangle
+\mbox{\rm an entire function}.
\end{multline}
On the other hand, we have the following:
\[
 \bigl\langle
 \psi^{(\lambda_0)}_{t,u_0}C\bigl([m],\bar{N}[\bar{\mu}]\bigr),\,
 \phi
 \bigr\rangle
=\Res_{Z(s+\eigenmap(\lambda,u_0))}
 \Big\langle
 C(m,\bar{N}\bar{\mu}),\,\,
 |t|^{2s}\phi\!\cdot\!\chi\!\cdot\!\frac{i}{2\pi}dt\wedge d\bar{t}
 \Big\rangle
\]
We have the following equalities:
\begin{multline}
 \Big\langle
 C(m,\bar{N}\bar{\mu}),\,\,
 |t|^{2s}\phi\!\cdot\!\chi\frac{i}{2\pi}dt\wedge d\bar{t}
 \Big\rangle 
=\Big\langle
 C(m,\overline{(-\deldel_tt+\eigenmap(\lambda,u_0))\cdot\mu}),\,\,
 |t|^{2s}\phi\chi\frac{i}{2\pi}dt\wedge d\bar{t}
 \Big\rangle \\
=\Big\langle
 C(m,\bar{\mu}),\,\,
 \bigl(-\deldelbar_t\bar{t}
 +\sigma^{\ast}\bigl(\overline{\eigenmap(\lambda,u_0)}\bigr)
 \bigr)
 \cdot
 |t|^{2s}\phi\!\cdot\!
 \chi\!\cdot\!
 \frac{i}{2\pi}dt\wedge d\bar{t}
 \Big\rangle \\
=\Big\langle
 C(m,\bar{\mu}),\,\,
 \bigl(
 \bigl(-\deldelbar_t\bar{t}
 +\sigma^{\ast}\bigl(\overline{\eigenmap(\lambda,u_0)}\bigr)
 \bigr)\cdot
 |t|^{2s}
 \bigr)
 \!\cdot\!\phi\!\cdot\!\chi
 \frac{i}{2\pi}dt\wedge d\bar{t}
 \Big\rangle
+\mbox{\rm an entire function}.
\end{multline}
We have the following relation:
\[
  \bigl(t\deldel_t+\eigenmap(\lambda,u)\bigr)
 \cdot
 |t|^{2s}
=\bigl(s\lambda+\eigenmap(\lambda,u)\bigr)\cdot|t|^{2s}
=\lambda\cdot\bigl(s+\lambda^{-1}\cdot\eigenmap(\lambda,u)\bigr)
 \cdot |t|^{2s}.
\]
We also have the following:
\begin{multline}
 \bigl(\bar{t}\deldelbar_t
 +\varphi\bigl({\eigenmap(\lambda,u)}\bigr)
 \bigr)\cdot
 |t|^{2s}
=\bigl(s\cdot\sigma^{\ast}{\bar{\lambda}}+
 \sigma^{\ast}\bigl(
 \overline{\eigenmap(\lambda,u)}
 \bigr)\bigr)
 \cdot |t|^{2s}
=-\lambda^{-1}\cdot|t|^{2s}\cdot
 \bigl(s+\sigma^{\ast}\bigl(
 \overline{\lambda^{-1}\cdot\eigenmap(\lambda,u)}
 \bigr)
 \bigr) \\
=-\lambda^{-1}\cdot|t|^{2s}\cdot
 \bigl(s+\lambda^{-1}\cdot\eigenmap(\lambda,u)
 \bigr)
\end{multline}
Thus we obtain the relation desired.
\hfill\qed

%% file: a34.2.tex

\subsubsection{The induced pairing on $\tildepsizero_{t,u}(\nbigm)$}

The pairing on $\psi_{t,u}(\nbigm)$ is not so good,
if $u$ is not contained in $\real\times\{0\}\subset\real\times\cnum$.
We modify it to construct the induced pairing:
\[
 \tildepsizero_{t,u}C:
 \tildepsizero_{t,u}\nbigm\otimes
 \overline{
 \tildepsi^{(\sigma(\lambda_0))}_{t,u}\nbigm}
\lrarr\distribution^{U}_{X_0}.
\]

\begin{lem}\label{lem;9.25.13}
Let $m$ be a section of $\psizero_{t,u}(\nbigm)$
and $\mu$ be a section of $\psi^{(\sigma(\lambda_0))}_{t,u}(\nbigm)$.
Then we have the equality:
\begin{equation}\label{eq;9.25.12}
 \psizero_{t,u-\vecdelta_0}C\bigl(
 [t\cdot m],\overline{[t\cdot\mu]}
 \bigr)
=\psizero_{t,u}C\bigl(
 [m],\overline{[\mu]}
 \bigr).
\end{equation}
\end{lem}
\pf
We have the following equality:
\begin{multline}
 \nbigi^{(0)}_{C(t\cdot m,\overline{t\cdot\mu}),\varphi}(s)
=\big\langle
 C(t\cdot m,\overline{t\cdot\mu}),
 \varphi\wedge |t|^{2s}\cdot \chi(t)\cdot \frac{i}{2\pi}dt\wedge d\bar{t}
 \big\rangle
=\big\langle
 C(m,\overline{\mu}),
 \varphi\wedge |t|^{2(s+1)}\cdot \chi(t)\cdot \frac{i}{2\pi}dt\wedge d\bar{t}
 \big\rangle \\
=\nbigi^{(0)}_{C(m,\bar{\mu}),\varphi}(s+1).
\end{multline}
Then the formula (\ref{eq;9.25.12}) immediately follows.
\hfill\qed

\vspace{.1in}

Let $m$ be a section of $\tildepsizero_{t,u}(\nbigm)$,
and $\mu$ be a section of $\tildepsi^{(\sigma(\lambda_0))}_{t,u}(\nbigm)$.
We pick a sufficiently large integer $N$
and $m_1\in\psizero_{t,u-N\vecdelta_0}(\nbigm)$
and $\mu_1\in \psi^{(\sigma(\lambda_0))}_{t,u-N\vecdelta_0}(\nbigm)$
corresponding to $m$ and $\mu$ respectively.
Then the pairing 
$\tildepsizero C(m,\bar{\mu})$
is defined to be
$\psizero C(m_1,\overline{\mu_1})$.

\begin{cor}\label{cor;a11.18.6}
It is well defined.
\end{cor}
\pf
It immediately follows from the definition
of $\tildepsizero(\nbigm)$
and Lemma \ref{lem;9.25.13}.
\hfill\qed

\subsubsection{The induced pairings $\psi_{t,u}C$ and $\tildepsi_{t,u}C$}

\begin{cor}
We obtain the induced pairings:
\[
\begin{array}{l}
 \tildepsi_{t,u} C:
 \tildepsi_{t,u}(\nbigm')_{\AAA}
\otimes
 \overline{\tildepsi_{t,u}(\nbigm'')}_{\AAA}
\lrarr
 \distribution_{X_0}^{\AAA},\\
\mbox{{}}\\
 \psi_{t,u}C:\psi_{t,u}(\nbigm')_{\AAA}
\otimes
 \overline{\psi_{t,u}(\nbigm'')}_{\AAA}
\lrarr
 \distribution_{X_0}^{\AAA}.
\end{array}
\]
\end{cor}
\pf
It follows from Lemma \ref{lem;a11.18.5}
and Corollary \ref{cor;a11.18.6}.
\hfill\qed

\begin{cor}\label{cor;a11.18.7}
We have the following relations:
\[
 \tildepsi_{t,u_0}C\bigl(N[m],[\bar{\mu}]\bigr)
=(i\lambda)^2\cdot
  \tildepsi_{t,u_0}C\bigl([m],\bar{N}[\bar{\mu}]\bigr), 
\quad
\psi_{t,u_0}C\bigl(N[m],[\bar{\mu}]\bigr)
=(i\lambda)^2\cdot
  \psi_{t,u_0}C\bigl([m],\bar{N}[\bar{\mu}]\bigr).
\]
Here we put $N:=-\deldel_tt+\eigenmap(\lambda,u)$.
\end{cor}
\pf
It immediately follows from Lemma \ref{lem;9.17.60}.
\hfill\qed

\begin{cor}
We have the induced pairing:
\[
 \tildepsi_{t,u}C:
 P_h\Gr^W_h\tildepsi_{t,u}\nbigm'_{\AAA}
\otimes
 \overline{
 P_h\Gr^W_{-h}\tildepsi_{t,u}\nbigm''}_{\AAA}
\lrarr\distribution^{\AAA}_{X_0}.
\]
Here $W$ denotes the weight filtration induced by $N$
in Corollary {\rm\ref{cor;a11.18.7}},
and $P_k\Gr^W_h$ denote the primitive part of the associated
graded modules.
\end{cor}
\pf
It immediately follows from Corollary \ref{cor;a11.18.7}.
\hfill\qed

\begin{cor}
We have the specialization of a $\nbigr$-triple
$(\nbigm',\nbigm'',C)$,
where $\nbigm'$ and $\nbigm''$ are strictly specializable.
\hfill\qed
\end{cor}

%% file: a35.tex

\subsubsection{Uniqueness of the pairing}

Let $Z$ be a closed irreducible subset of $X$.
Let $\nbigm'$ and $\nbigm''$ be strictly specializable
$\nbigr$-modules, whose supports are $Z$.
We assume that
the morphisms $\can$ for $\nbigm'$ and $\nbigm''$
are surjective.
Let $\lambda_0$ be a point of $\cnum^{\ast}$,
and let $U$ be an open subset $\Delta(\lambda_0,\epsilon_0)$
for some sufficiently small $\epsilon_0>0$.
Let $C_a:\nbigm'_{U}\otimes\overline{\nbigm''}_{\sigma(U)}
 \lrarr\distribution^{U}_{X}$
$(a=1,2)$ be sesqui-linear pairings.

Let $Z'\subset Z$ be a Zariski open subset.

\begin{lem}\label{lem;9.29.1}
Assume that $\lambda_0$ is generic with respect to
$\KMS(\nbigm',t)\cup\KMS(\nbigm'',t)\cup\{(0,0)\}$.
Assume that we have $C_1=C_2$ on $Z'$.
Then we have $C_1=C_2$ on $Z$.
\end{lem}
\pf
We follow the argument of Sabbah (Proposition 3.6 in \cite{sabbah}).
Since it is a local property,
we may assume that there exits a holomorphic function $f$
such that $f^{-1}(0)\cap Z\supset Z-Z'$.
We have only to show the coincidence of the pairings:
\[
 i_{f\,+}C_a:
 i_{f\,+}\nbigm'_U\otimes
 \overline{i_{f\,+}\nbigm''}_U\lrarr\distribution^{U}_X,
\quad
 (a=1,2).
\]
Thus we may assume that $X=X_0\times \cnum$
and $f=t$ is the coordinate of $\cnum$.

Let $m'$ and $m''$ be sections of
$\nbigm'_{U}[\deldel_t]$ and $\nbigm''_{\sigma(U)}[\deldel_t]$
respectively.
We put
$A(m',\overline{m''}):=C_1(m',\overline{m''})-C_2(m',\overline{m''})$.
Since the support of $A(m',\overline{m''})$ is contained in $X_0$,
we have the following development:
\[
 A(m',m'')=
 \sum_{a+b\leq p}
 \eta_{a,b}\cdot\deldel_t^a\cdot\deldelbar^b_t\cdot\delta_{X_0}.
\]
Here $\eta_{a,b}$ denotes the $H(U\cap\AAA)$-valued distributions 
on $X_0$.
We have only to show $\eta_{a,b}=0$ for any $a$ and $b$.

Let us consider the case $m'\in \Vzero_{<0}(\nbigm')$.
We have a finite subset
$S\subset \KMS(\nbigm,0,t)$ such that the following holds:
\begin{itemize}
\item
$\paramap(\lambda_0,u)<0$ for any element $u\in S$.
\item
We put $B(x):=\prod_{u\in S}\bigl(x+\eigenmap(\lambda,u)\bigr)$.
Then we have $B(-\deldel_tt)\cdot m'=P\cdot t^{p+1}\cdot m'$.
Here $P$ denotes an element of $V_0\nbigr_{X\times\cnum}$.
\end{itemize}
Then we have the following vanishing:
\[
 B(-\deldel_tt)\cdot
 A(m',m'')
=P\cdot t^{p+1}\cdot A(m',m'')=0.
\]

Note that we have the following:
\[
 \bigl(-\deldel_tt+\eigenmap(\lambda,u)\bigr)
 \cdot \deldel_t^a\cdot
 \delta_{X_0}
=\bigl(
 a\lambda+\eigenmap(\lambda,u)
 \bigr)\cdot
 \deldel_t^a\cdot
 \delta_{X_0}
=\eigenmap(\lambda,u+a\cdot\vecdelta_{0})
 \cdot\deldel_t^a\cdot\delta_{X_0}.
\]
Note we have $\paramap(\lambda_0,u-a\vecdelta_0)<0$.
Then we have $a\lambda_0+\eigenmap(\lambda_0,u)\neq 0$
due to the genericness of $\lambda_0$.
Thus we obtain $\eta_{a,b}=0$ in the case $m'\in \Vzero_{<0}(\nbigm')$.
Since $\nbigm'$ is generated by $\Vzero_{<0}\nbigm'$ around $\lambda_0$,
the general case can be reduced to the case $m'\in \Vzero_{<0}(\nbigm')$.
\hfill\qed

\begin{prop}\label{prop;a11.23.10}
Assume that we have $C_1=C_2$ on $Z'$.
Then we have $C_1=C_2$ on $Z$.
\end{prop}
\pf
Let $\phi$ be a test function.
Since $C_i(m',\overline{m''})(\phi)$ are holomorphic functions
on $\AAA$,
we have only to show that the coincidence on a neighbourhood
of a generic $\lambda$.
Thus the proposition follows from Lemma \ref{lem;9.29.1}.
\hfill\qed

%% file: a74.2.tex

We put $X=\Delta^n$ and $D=\bigcup_{i=1}^l D_i$.
Let $\harmonicbundle$ be a tame harmonic bundle over $X-D$.
The sheaf $\nbige$ with $\lambda$-connection $\DD$
can be naturally regarded as
the left $\nbigr_{\nbigx-\nbigd}$-module.
By tensoring $\bigwedge^n \Omega_{X-D}$,
we obtain the right $\nbigr_{\nbigx-\nbigd}$-module.
In this section, we use the right $\nbigr$-module structure,
if we do not mention.
For simplicity,
we omit to denote `$\otimes \bigwedge^n\Omega_{X-D}$'.
The author hopes that there are no confusion.

%% file: a75.tex

\subsubsection{Definition}
\label{subsubsection;a12.2.7}

Let $\lambda_0$ be a point of $\cnum_{\lambda}$.
Pick $\vecb\in\real^l$
such that $b_i\not\in\Par(\nbigelambdazero,i)$ for each $i$.

\begin{description}
\item [$(C1)$]
Let $\epsilon_0$ be a sufficiently small positive number
such that we have $\prolongg{\vecb}{\nbige}$
on $\nbigx(\lambda_0,\epsilon_0)$.
\end{description}

We have the natural inclusion $j:X-D\lrarr X$,
which induces $j:\nbigx-\nbigd\lrarr\nbigx$.
We have the subsheaf $\naiveprolong{\nbige}^{(\lambda_0)}$
of $j_{\ast}\nbige_{|\nbigx(\lambda_0,\epsilon_0)}$
given as follows:
\[
 \naiveprolong{\nbige}^{(\lambda_0)}
:=\bigcup_{a=1}^{\infty}
 \Bigl(\prod_{i=1}^lz_i\Bigr)^{-a}\cdot
 \prolongg{\vecb}{\nbige}.
\]
The following lemma is easy to see.
\begin{lem}\mbox{{}} \label{lem;10.2.3}
\begin{itemize}
\item
$\naiveprolong{\nbige}^{(\lambda_0)}$
is characterized as follows:
\[
 \Gamma\bigl(U,\naiveprolong{\nbige}^{(\lambda_0)}\bigr)
=\Bigl\{
 f\in \Gamma(U,j_{\ast}\nbige),\,\Big|\,
 |f|\leq C_1\cdot \prod_{i=1}^l|z_i|^{-C_2},\,\,\,\,\,
 \bigl(\exists C_1, C_2>0\bigr)
 \Bigr\},
\quad
\bigl(U\subset\nbigx(\lambda_0,\epsilon_0)\bigr).
\]
\item
Let $\vecb'$ and $\epsilon_0'$ be other choices.
We pick $0<\epsilon_0''<\min(\epsilon_0',\epsilon_0'')$.
Then we obtain two sheaves
$\naiveprolong{\nbige}^{(\lambda_0)}(\vecb,\epsilon_0)$
and $\naiveprolong{\nbige}^{(\lambda_0)}(\vecb',\epsilon_0')$.
Then we have the following:
\[
 \naiveprolong{\nbige}^{(\lambda_0)}(\vecb,\epsilon_0)_{|
 \Delta(\lambda_0,\epsilon_0'')}
=\naiveprolong{\nbige}^{(\lambda_0)}(\vecb',\epsilon_0')_{|
 \Delta(\lambda_0,\epsilon_0'')}.
\]
\hfill\qed
\end{itemize}
\end{lem}

In the following, we omit to denote $\vecb$
and $\epsilon_0$.

\begin{cor} \label{cor;a10.2.4}
Pick $\lambda_1\in\Delta(\lambda_0,\epsilon_0)$,
and $\epsilon_1$ such that
$\Delta(\lambda_1,\epsilon_1)\subset \Delta(\lambda_0,\epsilon_0)$.
Then we have the following:
\[
 \naiveprolong{\nbige}^{(\lambda_0)}_{|\Delta(\lambda_1,\epsilon_1)}
=\naiveprolong{\nbige}^{(\lambda_1)}.
\]
\end{cor}
\pf It follows from Lemma \ref{lem;10.2.3}.
\hfill\qed

\begin{df}
We define the sheaf $\naiveprolong{\nbige}$
by the following condition:
\[
 \naiveprolong{\nbige}_{|\Delta(\lambda_0,\epsilon_0)}
=\naiveprolong{\nbige}^{(\lambda_0)}.
\]
It is well defined,
due to Corollary {\rm\ref{cor;a10.2.4}}.
\hfill\qed
\end{df}

If $\epsilon_0$ is sufficiently small,
we may assume that the following condition is satisfied:
\begin{description}
\item[(C2)]
 We put
$0\neq\eta:=
 \min
 \bigcup_i\big\{|a-b|\,\big|\,a,b\in\Par(\nbigelambdazero,i),a\neq b\big\}$.
We pick $\eta_1<3^{-1}\eta$.
Then for any $u\in \KMS(\nbige^0,i)$ and
for any $\lambda\in\Delta(\lambda_0,\epsilon_0)$,
the inequality
$\big|\paramap(\lambda,u)-\paramap(\lambda_0,u)\big|<\eta_1$
holds.
\end{description}

Note the following elementary lemma.
\begin{lem} \label{lem;9.19.10}
Let $\epsilon_0$ be a positive number satisfying the condition $(C2)$.
For any element $u\in \KMS(\nbige^0,i)$ such that
$\kmsmap(\lambda_0,u)\in \KMS(\prolongg{\vecb}{\nbige^{\lambda_0}},i)$,
and for any $\lambda\in\Delta(\lambda_0,\epsilon_0)$,
we have the following:
\[
 \begin{array}{l}
 \paramap(\lambda_0,u)<0 \Longrightarrow
 \paramap(\lambda,u)<0,\\
\mbox{{}}\\
 \paramap(\lambda_0,u)>0 \Longrightarrow
 \paramap(\lambda,u)>0.
 \end{array}
\]
\end{lem}
\pf
Obvious.
\hfill\qed

\vspace{.1in}

Let $\epsilon_0$ be a positive number satisfying the condition $(C2)$.
Then we have the parabolic structure $\lefttop{i}F=
 \bigl(\lefttop{i}\Fzero_c\,\big|\,
 b_i-1\leq c\leq b_i \bigr)$
of $\prolongg{\vecb}{\nbige}$.
We shall define the subsheaf
$\naiveprolongg{\vecc}{\nbige}^{(\lambda_0)}$
as follows:

\vspace{.1in}
\noindent
1.
For $\vecc=(c_i)\in\real^l$ such that
$c_i\in\Par\bigl(\prolongg{\vecb}{\nbige}^{\lambda_0},i\bigr)$,
we put as follows:
\[
 \naiveprolongg{\vecc}{\nbige}^{(\lambda_0)}
:=\bigl\{f\in\prolongg{\vecb}{\nbige}\,\big|\,
 f_{|\nbigd_i}\in\lefttop{i}\Fzero_{c_i}
 \bigr\}.
\]

\begin{lem}
We have
$\naiveprolongg{\vecc}{\nbige}^{(\lambda_0)}
=\prolongg{\vecc+\eta_1\vecdelta}{\nbige}$.
In particular,
it is locally free on $\Delta(\lambda_0,\epsilon_0)$.
\end{lem}
\pf
It is clear from our choice of $\epsilon_0$ and $\eta_1$.
\hfill\qed

\vspace{.1in}

In particular,
we obtain the parabolic filtration $\lefttop{i}\Fzero$
and the decomposition $\lefttop{i}\EE^{(\lambda_0)}$
of $\naiveprolongg{\vecc}{\nbige^{(\lambda_0)}}$.

\vspace{.1in}
\noindent
2. For $\vecc\in\prod_{i=1}^l[b_i-1,b_i]$,
 we put $\vecc'=(c_i')$,
 $c_i':=
 \max\bigl\{
 x\in \Par(\prolongg{\vecb}{\nbige}^{(\lambda_0)},i)
  \,\big|\,x\leq c_i
 \bigr\}$,
and we put as follows:
\[
 \naiveprolongg{\vecc}{\nbige^{(\lambdazero)}}:=
 \naiveprolongg{\vecc'}{\nbige^{(\lambdazero)}}.
\]

\noindent
3. For general $\vecc\in\real^l$
we take $\vecn\in\seisuu^l$
such that $\vecc-\vecn=\vecc'\in\prod_{i=1}^l[b_i-1,b_i]$,
and then we put as follows:
\[
  \naiveprolongg{\vecc}{\nbige^{(\lambdazero)}}:=
 \Bigl(\prod_{i=1}^lx_i^{-n_i}\Bigr)\cdot
 \naiveprolongg{\vecc'}{\nbige^{(\lambdazero)}}.
\]

\begin{rem}
As before, the subsheaf
$\naiveprolongg{\vecc}{\nbige}^{(\lambda_0)}$
is given independently of choices of $\vecb$ and $\epsilon_0$.
\hfill\qed
\end{rem}

\begin{lem}
We have the following:
\[
 \naiveprolongg{\vecc}{\nbige}^{(\lambda_0)}
 _{|\Delta(\lambda_1,\epsilon_1)}
=\naiveprolongg{\vecd}{\nbige}^{(\lambda_1)}
\]
Here we put $\vecd=(d_i)$ and
$d_i=\max\bigl\{
 \paramap(\lambda_1,u)\,\big|\,
 u\in \KMS(\nbige^0,i),\,\,\paramap(\lambda_0,u)\leq c_i
 \bigr\}$.
\end{lem}
\pf
It immediately follows from our construction.
\hfill\qed

\vspace{.1in}

Since $\nbige$ is a right $\nbigr_{\nbigx-\nbigd}$-module
(see the remark in the first part of this section),
$j_{\ast}\nbige$ is a $\nbigr_{\nbigx}$-module.
The following lemma is clear.
\begin{lem} \label{lem;a11.19.1}
$\naiveprolong{\nbige}$ is the $\nbigr_{\nbigx}$-submodule
of $j_{\ast}\nbige$.
We have the following implication:
\[
 \naiveprolongg{\vecc}{\nbige^{(\lambdazero)}}\cdot \deldel_i
 \subset
 \naiveprolongg{\vecc+\vecdelta_i}{\nbige^{(\lambdazero)}},
\quad\quad
 \naiveprolongg{\vecc}{\nbige^{(\lambdazero)}}\cdot z_i
 \subset
 \naiveprolongg{\vecc-\vecdelta_i}{\nbige^{(\lambdazero)}}.
\]
\hfill\qed
\end{lem}

\subsubsection{The filtrations $\lefttop{i}\Vzero$ of
 $\naiveprolong{\nbige^{(\lambdazero)}}$}

We put $\vecdelta:=\overbrace{(1,\ldots,1)}^l$.
For $\vecb\in\real^l$,
we put as follows on $\Delta(\lambda_0,\epsilon_0)$:
\[
 \lefttop{\lbar}V^{(\lambda_0)}_{\vecb}\bigl(
 \naiveprolong{\nbige^{(\lambdazero)}}\bigr)
:=
 \naiveprolongg{\vecb+\vecdelta}{\nbige^{(\lambda_0)}}.
\]

For any $I\subset \lbar$
and for any $\vecc\in\real^I$,
we put as follows:
\[
 \lefttop{I}V^{(\lambda_0)}_{\vecc}\bigl(
 \naiveprolong{\nbige^{(\lambdazero)}}\bigr)
 =
 \bigcup_{q_I(\vecb)=\vecc}
 \lefttop{\lbar}V^{(\lambda_0)}_{\vecb}\bigl(
 \naiveprolong{\nbige^{(\lambdazero)}}\bigr).
\]
The following equality is clear:
\[
 \bigcap_{i\in I}\lefttop{i}V^{(\lambda_0)}_{c_i}\bigl(
 \naiveprolong{\nbige^{(\lambdazero)}}\bigr)
=\lefttop{I}V^{(\lambda_0)}_{\vecc}\bigl(
 \naiveprolong{\nbige^{(\lambdazero)}}\bigr).
\]

Let $I\sqcup J=\lbar$ be a decomposition.
Let $\vecc$ and $\vecd$ be elements of
$\real^I$ and $\real^J$ respectively.
On the vector bundle
$\lefttop{I}\Gr^{\Vzero}_{\vecc}
 \lefttop{J}\Vzero_{\vecd}\bigl(\naiveprolong{\nbige}\bigr)$,
we have the induced filtrations
$\lefttop{j}\Vzero$ $(j\in J)$.
The following lemma is easy to see.
\begin{lem}\label{lem;a12.2.10}
Let $\vecb'$ be an element of $\real^J$
such that $\vecb'\leq \vecb$.
We have the natural isomorphism:
\[
 \lefttop{J}\Gr^{\Vzero}_{\vecb'}
 \Bigl(
 \lefttop{I}\Gr^{\Vzero}
 \lefttop{J}\Vzero_{\vecb}
 \bigl(
 \naiveprolong{\nbige}^{(\lambda_0)}
 \bigr)
 \Bigr)
\simeq
 \lefttop{\lbar}\Gr^{\Vzero}_{\vecc+\vecb'}
 \bigl(
 \naiveprolong{\nbige}^{(\lambda_0)}
 \bigr).
\]
\end{lem}
\pf
It can be reduced to the compatibility of the parabolic filtrations.
\hfill\qed

\vspace{.1in}
Let $i$ be an element of $I$.
The actions of $z_i$ and $\deldel_i$ induce
the morphisms:
\begin{equation}\label{eq;a12.2.1}
 \begin{array}{l}
 z_i:
 \lefttop{I}\Gr^{\Vzero}_{\vecc}\lefttop{J}\Vzero_{\vecd}
 \bigl(\naiveprolong{\nbige}\bigr)
\lrarr
 \lefttop{I}\Gr^{\Vzero}_{\vecc-\vecdelta_i}
 \lefttop{J}\Vzero_{\vecd}
 \bigl(\naiveprolong{\nbige}\bigr),\\
 \mbox{{}}\\
 \deldel_i:
  \lefttop{I}\Gr^{\Vzero}_{\vecc}\lefttop{J}\Vzero_{\vecd}
 \bigl(\naiveprolong{\nbige}\bigr)
\lrarr
 \lefttop{I}\Gr^{\Vzero}_{\vecc+\vecdelta_i}
 \lefttop{J}\Vzero_{\vecd}
 \bigl(\naiveprolong{\nbige}\bigr),\\
 \end{array}
\end{equation}
\begin{lem} \mbox{{}}\label{lem;a12.2.11}
\begin{itemize}
\item
The induces morphism $z_i$ in {\rm(\ref{eq;a12.2.1})} is isomorphic.
It is strict with respect to the filtrations
$\lefttop{j}\Vzero_j$  $(j\in J)$.
\item
The induced morphism $\deldel_i$ in {\rm(\ref{eq;a12.2.1})}
is strict with respect to the filtrations
$\lefttop{j}\Vzero_j$ $(j\in J)$.
\end{itemize}
\end{lem}
\pf
The first claim is clear.
The second claim can be reduced to the strictness
of the residue $\Res_i(\DD)$ with respect to
the parabolic filtrations $\lefttop{j}\Fzero$ $(j\in J)$.
The generalized eigen decompositions
of the residues $\Res_j(\DD)$ $(j\in J)$
gives the splitting of $\lefttop{j}\Fzero$ $(j\in J)$
on the open dense subset of generic points in
$\Delta(\lambda_0,\epsilon_0)$.
Hence we can show the strictness of
$\Res_i(\DD)$ with respect to the filtrations
$\lefttop{j}F$ $(j\in J)$
by using Proposition \ref{prop;c11.12.60}.
\hfill\qed

%% file: a75.2.tex

\subsubsection{The right action of $z_i\deldel_i$}

Let $I$ be a subset of $\lbar$,
and $J=\lbar-I$.
Let $\vecc$ and $\vecd$ be elements of
$\real^I$ and $\real^{J}$ respectively.
On $\nbigd_I(\lambda_0,\epsilon_0)$,
we have the sheaf
$\lefttop{I}\Gr^{V^{(\lambda_0)}}_{\vecc}
 \lefttop{J}\Vzero_{\vecd}\bigl(\naiveprolong{\nbige}\bigr)$.

\begin{lem} \label{lem;a11.19.2}
 We have the natural isomorphism
$ \lefttop{I}\Gr^{\Vzero}_{\vecc}
 \lefttop{J}\Vzero_{\vecd}\bigl(\naiveprolong{\nbige}\bigr)
\simeq\lefttop{I}\Gr^{\Fzero}_{\vecc+\vecdelta_I}
  \bigl(\naiveprolongg{\vecc+\vecd+\vecdelta}{\nbige}_{|\nbigd_I}\bigr)$.
In particular,
$\lefttop{I}\Gr^{\Vzero}_{\vecc}\lefttop{J}\Vzero_{\vecd}
\bigl(\naiveprolong{\nbige}\bigr)$ is a vector bundle
over $\nbigd_I(\lambda_0,\epsilon_0)$.
\end{lem}
\pf
It is clear from our construction of the filtration
$\lefttop{I}\Vzero$.
\hfill\qed

\vspace{.1in}
The right actions of $z_i\deldel_i$ $(i\in I)$
on $\naiveprolong{\nbige}$ induces the
endomorphisms of
$\lefttop{I}\Gr^{\Vzero}_{\vecc}\lefttop{J}\Vzero_{\vecd}
\bigl(\naiveprolong{\nbige}\bigr)$.
(See Lemma \ref{lem;a11.19.1}).

\begin{lem}
The endomorphism induced by $z_i\cdot\deldel_i$
is same as $-\Res_i(\DD)-\lambda$,
under the isomorphism in Lemma {\rm\ref{lem;a11.19.2}}.
\end{lem}
\pf
The right action of $z_i\cdot\deldel_i$
corresponds to the left action of
$-\deldel_i\cdot z_i=-z_i\cdot\deldel_i-\lambda$.
Then the claim immediately follows.
\hfill\qed

\begin{lem}
The eigen functions of the right action of $z_i\deldel_i$ $(i\in I)$
on $\lefttop{I}\Gr^{\Vzero}_{\vecc}\lefttop{J}\Vzero_{\vecd}
 \bigl(\naiveprolong{\nbige}\bigr)$
are as follows:
\[
\bigl\{
 -\eigenmap(\lambda,u)\,\big|\,
 u\in\KMS(\nbige^0,i),\,\,\paramap(\lambda_0,u)=c_i
\bigr\}.
\]
\end{lem}
\pf
 Hence the eigenvalue of the right action of $z_i\deldel_i$
 is described as $-\eigenmap(\lambda,u)-\lambda$
 for $u\in\KMS(\nbige^0,i)$ such that
 $\paramap(\lambda,u)=c_i+1$.
 Thus it is described as
 $-\eigenmap(\lambda,u)$ for some $u\in \KMS(\nbige^0,i)$
 such that $\paramap(\lambda,u)=c_i$.
\hfill\qed

\vspace{.1in}

In the following,
we put as follows:
\[
\begin{array}{l}
 \nbigk(\nbige,\lambda_0,i,c)
:=
 \bigl\{u\in\KMS(\nbige^0,i)\,\big|\,\paramap(\lambda_0,u)=c\bigr\},\\
\mbox{{}}\\
 \nbigk(\nbige,\lambda_0,I,c)
:=
 \bigl\{\vecu\in\KMS(\nbige^0,I)\,\big|\,
 \paramap(\lambda_0,\vecu)=\vecc\bigr\}.
\end{array}
\]

\subsubsection{$\lefttop{I}\tilde{T}(\vecb,\vecd)$ and some properties}
We put as follows:
\[
 \lefttop{I}\psizero_{\vecu}
 \lefttop{J}\Vzero_{\vecd}
 (\naiveprolong{\nbige})
:=\lefttop{I}\EE\bigl(-\eigenmap(\lambda,\vecu)\bigr).
\]
Then we have the decomposition:
\[
 \lefttop{I}\Gr^{\Vzero}_{\vecc}
 \lefttop{J}\Vzero_{\vecd}
 \bigl(
 \naiveprolong{\nbige}
 \bigr)
=\bigoplus_{\vecu\in\nbigk(\nbige,\lambda_0,I,\vecc)}
 \lefttop{I}\EE\bigl(-\eigenmap(\lambda,\vecu)\bigr)
=\bigoplus_{\vecu\in\nbigk(\nbige,\lambda_0,I,\vecc)}
  \lefttop{I}\psizero_{\vecu}
 \lefttop{J}\Vzero_{\vecd}
 (\naiveprolong{\nbige}).
\]

\begin{lem}
We have the following induced morphisms for $i\in I$:
\[
 z_i:
 \lefttop{I}\psizero_{\vecu}\lefttop{J}\Vzero_{\vecd}
 (\naiveprolong{\nbige})
\lrarr
 \lefttop{I}\psizero_{\vecu-\vecdelta_{0\,i}}
 \lefttop{J}\Vzero_{\vecd}(\naiveprolong{\nbige}).
\]
\[
 \deldel_i:
 \lefttop{I}\psizero_{\vecu}
 \lefttop{J}\Vzero_{\vecd}(\naiveprolong{\nbige})
\lrarr
 \lefttop{I}\psizero_{\vecu+\vecdelta_{0,i}}
 \lefttop{J}\Vzero_{\vecd}(\naiveprolong{\nbige}).
\]
Here $\vecdelta_{0,i}$ denotes the element of
$\bigl(\real\times\cnum\bigr)^I$
determined by
$q_i(\vecdelta_{0,i})=(1,0)$ and
$q_j(\vecdelta_{0,i})=(0,0)$ $(j\neq i)$.
\end{lem}
\pf
It follows from the relations
$(z_i\cdot\deldel_i)\cdot z_i
=z_i\cdot\bigl(z_i\cdot\deldel_i\bigr)+\lambda\cdot z_i$
and
$(z_i\cdot\deldel_i)\cdot\deldel_i
=\deldel_i\cdot(z_i\cdot\deldel_i)-\lambda\cdot\deldel_i$.
\hfill\qed

\vspace{.1in}

Let $\vecb=(b_i\,|\,i\in I)$ be an element of $\real^I$.
The elements $\veca=(a_i\,|\,i\in I)\in\real^I_{<0}$
and $\vecn=(n_i\,|\,i\in I)\in\seisuu^I_{\geq 0}$
are taken as follows:
If $b_i\geq 0$,
the numbers $a_i$ and $n_i$ are determined by
the conditions $-1\leq a_i<0$ and $n_i:=b_i-a_i$.
If $b_i< 0$, we put $a_i:=b_i$ and $n_i:=0$.
If $\paramap(\lambda_0,\vecu)=\vecb$,
we obtain the morphism:
\[
 \prod_i\deldel_i^{n_i}:
 \lefttop{I}\psizero_{\vecu'}
 \lefttop{J}\Vzero_{\vecd}\bigl(\naiveprolong{\nbige}\bigr)
\lrarr
 \lefttop{I}\psizero_{\vecu}
 \lefttop{J}\Vzero_{\vecd}\bigl(\naiveprolong{\nbige}\bigr).
\]
Here we put $\vecu'=\vecu-\sum_i n_i\cdot\vecdelta_{0,i}$.
The image of the morphism is denoted by
$\lefttop{I}\tilde{T}^{(\lambda_0)}(\vecu,\vecd)$.
By our construction,
we have the surjection,
in the case $\paramap(\lambda_0,u_i)\geq -1$:
\[
 \deldel_i:
 \lefttop{I}\tilde{T}^{(\lambda_0)}(\vecu,\vecd)
\lrarr
 \lefttop{I}\tilde{T}^{(\lambda_0)}(\vecu+\vecdelta_{0,i},\vecd).
\]
The multiplication of $z_i$ induces
the morphism
$\lefttop{I}\tilde{T}^{(\lambda_0)}(\vecu,\vecd)
\lrarr
\lefttop{I}\Gr^{\Vzero}_{\vecb-\vecdelta_i}\lefttop{J}\Vzero_{\vecd}
 \bigl(\naiveprolong{\nbige}\bigr)$.

\begin{lem}
The image is contained in
$ \lefttop{I}\tilde{T}^{(\lambda_0)}(\vecu-\vecdelta_{0,i},\vecd)$.
\end{lem}
\pf
We decompose $\vecu=\vecu'+\sum n_i\vecdelta_{0,i}$ as above.
In the case $n_i=0$, the claim is clear.
In the case $n_i>0$,
the claim follows from the relation
$ \deldel_i^{n_i}z_i
=\bigl(z_i\deldel_i+\lambda\cdot n_i\bigr)\deldel_i^{n_i-1}$.
\hfill\qed

\vspace{.1in}
The following lemma is clear.
\begin{lem}
The morphism
$z_i:
 \lefttop{I}\tilde{T}(\vecu,\vecd)
 \lrarr
 \lefttop{I}\tilde{T}(\vecu-\vecdelta_{0,i},\vecd)$
is injective.
It is isomorphic in the following cases:
\begin{itemize}
\item
 $\paramap(\lambda_0,u_i)<0$.
\item
 $\paramap(\lambda_0,u_i)=0$ and $u_i\neq (0,0)$.
\item
 $\paramap(\lambda_0,u_i)>0$ and $\eigenmap(\lambda_0,u_i)\neq 0$.
\hfill\qed
\end{itemize}
\end{lem}

\vspace{.1in}

The unique eigenvalue of the right action of $z_i\deldel_i$ on
$\lefttop{I}\tilde{T}^{(\lambda_0)}(\vecu,\vecd)$ is
$-\eigenmap(\lambda,u_i)$,
by our construction.

\begin{lem} \label{lem;a11.19.5}
The morphism
$\deldel_i:\lefttop{I}\tilde{T}^{(\lambda_0)}(\vecu,\vecd)
\lrarr \lefttop{I}\tilde{T}^{(\lambda_0)}(\vecu+\vecdelta_{0,i},\vecd)$
is injective unless $u_i=(-1,0)$.
\end{lem}
\pf
Let us consider the following morphisms:
\[
 \begin{CD}
 \lefttop{I}\tilde{T}^{(\lambda_0)}(\vecu,\vecd)
 @>{\deldel_i}>>
 \lefttop{I}\tilde{T}^{(\lambda_0)}(\vecu+\vecdelta_{0,i},\vecd)
 @>{z_i}>>
 \lefttop{I}\tilde{T}^{(\lambda_0)}(\vecu,\vecd).
 \end{CD}
\]
The composite is the right action of
$\deldel_iz_i=z_i\deldel_i+\lambda$.
The eigenfunction of $\deldel_iz_i$
on $\lefttop{I}\tilde{T}^{(\lambda_0)}(\vecu,\vecd)$
is given by
$\eigenmap(\lambda,u_i)+\lambda$.
In the case $u_i\neq (-1,0)$,
we have $-\eigenmap(\lambda,u_i)+\lambda
=-\alpha+(a+1)\lambda+\bar{\alpha}\lambda^2
\neq 0$.
Hence the composite $\deldel_iz_i$ is injective,
and thus $\deldel_i$ is injective.
\hfill\qed

\vspace{.1in}
The following lemma can be shown similarly.
\begin{lem}
The induced morphism $\deldel_i$
is isomorphic in the following cases:
\begin{itemize}
\item
$\paramap(\lambda_0,u_i)>-1$.
\item
$\paramap(\lambda_0,u_i)=-1$ and $u_i\neq (-1,0)$.
\item
 $\paramap(\lambda_0,u_i)<-1$ and
 $-\eigenmap(\lambda_0,u)+\lambda\neq 0$.
\hfill\qed
\end{itemize}
\end{lem}

\vspace{.1in}

\begin{lem}
The sheaf $\lefttop{I}\tilde{T}^{(\lambda_0)}(\vecu,\vecd)$
is a locally free $\nbigo_{\nbigd_I(\lambda_0,\epsilon_0)}$-module.
\end{lem}
\pf
We may assume $I=\kbar$ for some $k\leq l$.
We put $J:=\lbar-\kbar$.
Due to Lemma \ref{lem;a11.19.5},
we may assume that $u_i=(0,0)$ $(i\leq m)$
and $\paramap(\lambda_0,u_i)<0$ $(m<i\leq k)$
for some $m\leq k$.
We put $\vecu':=\vecu-\sum_{i=1}^m\vecdelta_{0,i}$.
Then $\lefttop{I}\tilde{T}^{(\lambda_0)}(\vecu,\vecd)$
is the image of the morphism:
\begin{equation}\label{eq;a11.19.10}
\prod_{i=1}^m\deldel_i:
 \lefttop{\kbar}\psizero_{\vecu'}\lefttop{J}\Vzero
\bigl(\naiveprolong{\nbige}\bigr)
\lrarr
 \lefttop{\kbar}\psizero_{\vecu}\lefttop{J}\Vzero
\bigl(\naiveprolong{\nbige}\bigr).
\end{equation}
On $\lefttop{\kbar}\psizero_{\vecu}\lefttop{J}\Vzero
 \bigl(\naiveprolong{\nbige}\bigr)$,
we have the nilpotent endomorphisms
$\nbign_i:=\Res_i(\DD)+\lambda$ $(i\in \mbar)$.
It is clear that
the image of the morphism (\ref{eq;a11.19.10})
is same as the image of
$\prod_{i=1}^m \nbign_i$.
Then it is locally free due to the limiting mixed twistor theorem.
(See Lemma \ref{lem;10.26.121}).
\hfill\qed

\vspace{.1in}

For $\vecb\in\real^I$ and $\vecd\in\real^{\lbar-I}$,
we put as follows:
\[
 \lefttop{I}\tilde{T}^{(\lambda_0)}(\vecb,\vecd)
:=\bigoplus_{\vecu\in\nbigk(\nbige,\lambda_0,I,\vecb)}
 \lefttop{I}\tilde{T}^{(\lambda_0)}(\vecu,\vecd).
\]

\begin{cor}\label{cor;c11.19.10}
The sheaf \,$\lefttop{I}\tilde{T}^{(\lambda_0)}(\vecb,\vecd)$
is a locally free $\nbigo_{\nbigd_I(\lambda_0,\epsilon_0)}$-module.
The surjective morphism
$\deldel_i:\lefttop{I}\tilde{T}^{(\lambda_0)}(\vecb,\vecd)
\lrarr \lefttop{I}\tilde{T}^{(\lambda_0)}(\vecb+\vecdelta_i,\vecd)$
is isomorphic in the case $b_i>-1$.

In the case $\vecb\in\real^I_{<0}$,
we have
$\lefttop{I}\tilde{T}^{(\lambda_0)}(\vecb,\vecd)
\simeq
 \lefttop{I}\Gr^{\Vzero}_{\vecb}
 \lefttop{J}\Vzero_{\vecd}\bigl(
 \naiveprolong{\nbige}\bigr)$.
\hfill\qed
\end{cor}

%% file: a75.5.tex

\subsubsection{A lemma}
\label{subsubsection;c11.24.1}

Let $k$ and $m$ be integers such that $1\leq m\leq k\leq l$.
Let $\vecc$ be an element of $\real^{k}$
such that $c_i=0$ $(1\leq i\leq m-1)$, $c_m=-1$
and $c_i<0$ $(m<i\leq k)$.
We put $J:=\lbar-\kbar$,
and let $\vecd$ be an element of $\real^J$.
For any subset $K\subset\kbar$, we put $K':=\lbar-K$.
We have the projection:
\[
 \pi_K:\lefttop{\lbar}\Vzero_{\vecc+\vecd}
 \bigl(\naiveprolong{\nbige}\bigr)
\lrarr
 \lefttop{K}\Gr^{\Vzero}_{q_K(\vecc)}
 \lefttop{J}\Vzero_{q_{K'}(\vecc+\vecd)}
 \bigl(\naiveprolong{\nbige}\bigr).
\]

The following lemma will be used in the proof
of Lemma \ref{lem;d11.19.2}.
\begin{lem} \label{lem;9.19.41}
Let $s$ be a section of the following locally free
sheaf on $\nbigd_{\lbar}$:
\begin{equation}\label{eq;a11.19.20}
 \Ker\Bigl(
 \lefttop{\kbar}\tilde{T}^{(\lambda_0)}(\vecc,\vecd)
 \stackrel{\deldel_m}{\lrarr}
 \lefttop{\kbar}\tilde{T}^{(\lambda_0)}(\vecc+\vecdelta_m,\vecd)
 \Bigr).
\end{equation}
Then we can take a section
$g\in
 \lefttop{\lbar}
  \Vzero_{\vecc+\vecd}(\naiveprolong{\nbige})$
satisfying the following:
\begin{itemize}
\item
$ \pi_{\kbar}(g)=s
 \in\lefttop{\kbar}\tilde{T}^{(\lambda_0)}
 (\vecc,\vecd)
\subset
 \lefttop{\kbar}\Gr^{\Vzero}_{\vecc}
 \lefttop{J}\Vzero_{\vecd}\bigl(\naiveprolong{\nbige}\bigr)$.
\item
 For any subset $K\subset\kbar$,
 we have the following:
\begin{equation}\label{eq;a11.19.21}
 \pi_K(g_{|\nbigd_K})\in
\left\{
 \begin{array}{ll}
 \ker\bigl(\Res_m(\DD)\bigr)\cap
 \Image\Bigl(\prod_{i\in K,i\leq m-1}(\Res_i(\DD)+\lambda)\Bigr),
 & (m\in K),\\
 \mbox{{}}\\
 \Image\Bigl(
 \prod_{i\in K,i\leq m-1}
 (\Res_i(\DD)+\lambda)  \Bigr), &(m\not\in K).
\end{array}
\right.
\end{equation}
Here the right hand side of {\rm(\ref{eq;a11.19.21})}
denote the subbundles of
$\lefttop{K}\Gr^{\Vzero}_{q_K(\vecc)}
 \lefttop{K'}\Vzero_{q_{K'}(\vecc+\vecd)}\bigl(
 \naiveprolong{\nbige}
 \bigr)$.
\end{itemize}
\end{lem}
\pf
Before entering the proof,
we give a few remarks.

\begin{enumerate}
\item
The subbundle (\ref{eq;a11.19.20})
of $\lefttop{\kbar}\Gr^{\Vzero}_{\vecc}
 \lefttop{J}\Gr^{\Vzero}_{\vecd}\bigl(\naiveprolong{\nbige}\bigr)$
is same as
$\Ker\bigl(\Res_m(\DD)\bigr)
 \cap\Image\bigl(\prod_{i=1}^{m-1}(\Res_i(\DD)+\lambda)\bigr)$.
\item
By tensoring the model bundle of rank 1,
we can reduce the problem to the case
$c_i=0$ $(i<m)$, $=-1$ $(m\leq i\leq k)$.
Hence we consider only the case in the following.
\end{enumerate}

For the proof of Lemma \ref{lem;9.19.41},
we need some preparation.
Let us take a frame $\vecv$ of
$\lefttop{\lbar}\Vzero_{\vecc+\vecd}\bigl(
 \naiveprolong{\nbige}
 \bigr)$, which is compatible with $\Fzero$ and $\EEzero$.
Since we have
$\lefttop{\lbar}\Vzero_{\vecc+\vecd}\bigl(
 \naiveprolong{\nbige}
 \bigr)=
 \naiveprolongg{\vecc+\vecd+\vecdelta_{\lbar}}{\nbige^{(\lambda_0)}}$
by definition,
we have the following:
\[
 \begin{array}{ll}
 0<\lefttop{j}\deg^{\Fzero}(v_i)\leq 1, & (1\leq j<m),\\
\mbox{{}}\\
 -1<\lefttop{j}\deg^{\Fzero}(v_i)\leq 0, & (m\leq j\leq k).
 \end{array}
\]
Let us take positive integer $b_j$ $(j\in\kbar)$.
Then we obtain the integers $n_j(v_i)$
determined as follows:
\[
 \begin{array}{ll}
 0<-n_j(v_i)+b_j\cdot \lefttop{j}\deg^{\Fzero}(v_i)\leq 1,
 & (1\leq j<m)\\
 \mbox{{}}\\
 -1<-n_j(v_i)+b_j\cdot \lefttop{j}\deg^{\Fzero}(v_i)\leq 0,
 & (m\leq j\leq k).\\
 \end{array}
\]
Note the following:
\[
 \begin{array}{ll}
 0\leq n_j(v_i)\leq b_j-1, & (1\leq j<m),\\
 \mbox{{}}\\
 -b_j+1\leq n_j(v_i)\leq 0, & (m\leq j\leq k).
 \end{array}
\]
If $b_j$ is sufficiently large,
we have
$n_j(v_i)\neq n_j(v_p)$
in the case
$\lefttop{j}\deg^{\Fzero}(v_i)\neq 
 \lefttop{j}\deg^{\Fzero}(v_p)$.

Let us consider the morphism $\psi_{\vecb}:X\lrarr X$
given by the correspondence:
\[
(z_1,\ldots,z_n)\longmapsto
 (z_1^{b_1},\ldots,z_k^{b_k},z_{k+1},\ldots,z_n).
\]
We put
$\tilde{v}_i:=\psi_{\vecb}^{-1}(v_i)\cdot\prod_{j}z_j^{n_j(v_i)}$,
and $\tilde{\vecv}:=(\tilde{v}_i)$.
Then $\tilde{\vecv}$ gives the frame
of $\lefttop{\lbar}\Vzero_{\vecc+\vecd}\bigl(
 \naiveprolong{\psi_{\vecb}^{-1}\nbige}
 \bigr)$,
which is compatible with $\EEzero$ and $\Fzero$.

Recall that we have the natural $\prod_{i=1}^k\mu_{b_i}$-action
on $X$
(see the subsection \ref{subsection;a11.9.53} for the notation).
Let $\omega_i$  be the generator of $\mu_{b_i}$,
and then the action is given by
$\omega_i\cdot\bigl(z_1,\ldots,z_n\bigr)
=\bigl(z_1,\ldots,z_{i-1},\omega_i\cdot z_i,z_{i+1},\ldots,z_n\bigr)$.
The action is naturally lifted to the action
on $\lefttop{\lbar}\Vzero_{\vecc+\vecd}\bigl(
 \naiveprolong{\psi_{\vecb}^{-1}\nbige}
 \bigr)$.

On the divisor $\nbigd_i(\lambda_0,\epsilon_0)$,
we have the decomposition:
\[
 \lefttop{\lbar}\Vzero_{\vecc+\vecd}\bigl(
 \naiveprolong{\psi_{\vecb}^{-1}\nbige}
 \bigr)_{|\nbigd_i(\lambda_0,\epsilon_0)}
=\bigoplus_{a\in S(i)}\lefttop{i}U_a.
\]
Here $\omega_i$ acts as $\omega_i^a$ on $\lefttop{i}U_a$,
and we put as follows:
\[
 S(i):=
\left\{
 \begin{array}{ll}
 \bigl\{ a\in\seisuu\,\big|\, 0\leq a\leq b_i-1 \bigr\},
 & (1\leq i<m),\\
 \mbox{{}}\\
 \bigl\{ a\in\seisuu\,\big|\, -b_i+1\leq a\leq 0 \bigr\},
 & (m\leq i\leq k).
 \end{array}
\right.
\]
The tuple of the decompositions
$\bigl(\lefttop{i}U\,\big|\,i\in\kbar\bigr)$
is compatible.

We put as follows:
\[
 n_i:=\left\{
 \begin{array}{ll}
 b_i-1, & (1\leq i<m),\\
 \mbox{{}}\\
 0 & (m\leq i\leq k).
 \end{array}
 \right.
\]
We obtain the element
$\vecn_K:=\bigl(n_i\,\big|\,i\in K\bigr)\in\seisuu^K$.
We have the vector bundle
$\lefttop{K}U_{\vecn_K}
:=\bigcap_{k\in K}\lefttop{k}U_{n_k\,|\,\nbigd_K}$.
Then the frames $\vecv$ and $\tilde{\vecv}$
induces the isomorphism:
\begin{equation} \label{eq;a11.19.25}
 \lefttop{K}U_{\vecn_K}
\simeq
 \lefttop{K}\Gr^{\Vzero}_{q_K(\vecc)}
 \lefttop{K'}\Vzero_{q_{K'}(\vecc+\vecd)}\bigl(
 \naiveprolong{\nbige}
 \bigr).
\end{equation}
On the right hand side of (\ref{eq;a11.19.25}),
we have the nilpotent endomorphisms
$\Res_i(\DD)+\lambda$ $(i\in K,i<m)$.
On the left hand side,
the corresponding morphism $F_i$
is given by
$b_i^{-1}\cdot\bigl(\Res_i(\psi_{\vecb}^{\ast}\DD)+\lambda\bigr)$.
In the case $m\in K$,
we also have the nilpotent endomorphism
$\Res_m(\DD)$ on the right hand side of (\ref{eq;a11.19.25}).
The corresponding morphism $F_m$
is given by $b_m^{-1}\cdot\Res_m(\psi_{\vecb}^{\ast}\DD)$.

In particular,
we have the isomorphism
$\lefttop{\kbar}U_{\vecn_{\kbar}}
\simeq
 \lefttop{\kbar}\Gr^{\Vzero}_{\vecc}
 \lefttop{J}\Vzero_{\vecd}\bigl(\naiveprolong{\nbige}\bigr)$.
Let $\tilde{s}$ denote the section
of $\lefttop{\kbar}U_{\vecn_{\kbar}}$
corresponding to $s$.

Then $\tilde{s}$ is contained in the following:
\[
 \ker\bigl(F_{m\,|\,\nbigd_{\kbar}}\bigr)\cap
 \Image\Bigl(\prod_{i\leq m-1}F_{i\,|\,\nbigd_{\kbar}}\Bigr).
\]

\begin{lem}
We can take an equivariant section
$\tilde{g}\in 
 \lefttop{\lbar}_{\vecc+\vecd}\bigl(
 \naiveprolong{\psi_{\vecb}^{-1}\nbige}
 \bigr)$
such that
 $\tilde{g}_{|\nbigd_{\kbar}}=\tilde{s}$
and that the following holds:
 \[
  \tilde{g}_{|\nbigd_K}
 \in
 \left\{
 \begin{array}{ll}
 \ker(F_m)\cap\Image\Bigl(\prod_{i\leq m-1,i\in K}F_i\Bigr)
 \cap \lefttop{K}U_{\vecn_K}
 & (m\in K)\\
 \mbox{{}}\\
 \Image\Bigl(\prod_{i\leq m-1,i\in J}F_i\Bigr) 
 \cap \lefttop{K}U_{\vecn_K}
 & (m\not\in K).
 \end{array}
 \right.
 \]
Here we put $\vecn_K:=\bigl(n_k\,\big|\,k\in K\bigr)$.
\end{lem}
\pf
It follows from Lemma \ref{lem;10.10.2}.
\hfill\qed

\vspace{.1in}

We put $g_1:=\tilde{g}\cdot \prod_{i=1}^{n_i}z_i^{-n_i}$,
which may be assumed invariant under the action of $\mu_{\vecb}$.
Then we have the section $g$ of
$\lefttop{\lbar}V_{\vecc+\vecd}(\naiveprolong{\nbige})$,
such that $g_1=\psi_{\vecb}^{-1}g$.
Due to the isomorphism (\ref{eq;a11.19.25}),
it can be checked that the section $g$ has the desired properties.
Thus the proof of Lemma \ref{lem;9.19.41} is accomplished.
\hfill\qed


%% file: a75.1.tex

\subsubsection{Preliminary for the construction of $\rmoduleprolong{E}$}

Let $\vecb$ be an element of $\real^l$.
We have the vector bundle
$\naiveprolongg{\vecb}{\nbige^{(\lambdazero)}}
   _{|\nbigd_i(\lambda_0,\epsilon_0)}$
over $\nbigd_i(\lambda_0,\epsilon_0)$ $(i=1,\ldots,l)$.
We often denote $\nbigd_i(\lambda_0,\epsilon_0)$
by $\nbigd_i$ for simplicity of notation.
We have the action of the residue $\Res_i(\DD)$
on $\naiveprolongg{\vecb}{\nbige^{(\lambdazero)}}
    _{|\nbigd_i(\lambda_0,\epsilon_0)}$.
Note that the right action of $-z_m\cdot\deldel_m$
induces the endomorphism $\Res_i(\DD)+\lambda$
on $\naiveprolongg{\vecb}{\nbigelambdazero}_{|\nbigd_i(\lambda_0,\epsilon_0)}$.

We put as follows
(see the subsubsection \ref{subsubsection;a12.2.5}
for $\EE_{\eta}$):
\[
 \lefttop{i}\EEzero\bigl(\naiveprolongg{\vecb}{\nbige^{(\lambdazero)}}
 _{|\nbigd_i(\lambda_0,\epsilon_0)},\beta\bigr)
:=\EE_{\eta}(\Res_i(\DD),\beta).
\]
Then we have the decomposition:
\[
 \naiveprolongg{\vecb}{\nbige^{(\lambda_0)}}_{|\nbigd_i}
=\bigoplus_{\beta\in \Sp(\prolongg{\vecb}{\nbigelambdazero},i)}
 \lefttop{i}\EEzero\bigl(
 \naiveprolongg{\vecb}{\nbige}^{(\lambda_0)}_{|\nbigd_i},\beta
\bigr).
\]

Let us consider the case $\vecb=\vecdelta=(\overbrace{1,\ldots,1}^l)$.
We have the decomposition
$\naiveprolongg{\vecdelta}{\nbige}^{(\lambda_0)}_{|\nbigd_i}
=\lefttop{i}\EE_0^{(\lambda_0)}
\oplus
 \lefttop{i}\EE_1^{(\lambda_0)}$.
Here we put as follows:
\[
 \lefttop{i}\EE^{(\lambda_0)}_0
=\lefttop{i}\EE^{(\lambda_0)}\bigl(
 \naiveprolongg{\vecdelta}{\nbige}^{(\lambda_0)}
 _{|\nbigd_i(\lambda_0,\epsilon_0)},
 -\lambda_0 \bigr),
\quad\quad
 \lefttop{i}\EE^{(\lambda_0)}_1
=\bigoplus_{\beta\neq -\lambda_0}
 \lefttop{i}\EE^{(\lambda_0)}\bigl(
 \naiveprolongg{\vecdelta}{\nbige}^{(\lambda_0)}
 _{|\nbigd_i(\lambda_0,\epsilon_0)},
 \beta\bigr).
\]

On $\lefttop{i}\EE^{(\lambda_0)}_1$,
the morphism $\Res_i(\DD)+\lambda$ is invertible.
We denote the inverse
by $(\Res_i(\DD)+\lambda)^{-1}$.
We have the filtration
$\lefttop{i}\Fzero$ on $\lefttop{i}\EE^{(\lambda_0)}_{j}$ $(j=0,1)$.
We have the following naturally defined projections:
\[
\begin{array}{l}
 \pi_i:
 \naiveprolongg{\vecdelta}{\nbige^{(\lambda_0)}}_{|\nbigd_i}
\lrarr
 \lefttop{i}\Gr^{\Fzero}_1\lefttop{i}\EE^{(\lambda_0)}_0,
\quad\quad
 \pi_i':
 \naiveprolongg{\vecdelta}{\nbige^{(\lambda_0)}}_{|\nbigd_i}
\lrarr
 \lefttop{i}\EE^{(\lambda_0)}_1.
\end{array}
\]

We put as follows:
\[
 \nbigk_0:=
 \bigl\{f\in\naiveprolongg{\vecdelta}{\nbige^{(\lambda_0)}}\,\big|\,
 \pi_i(f_{|\nbigd_i})=0\,\,\,(i\leq l)
 \bigr\}.
\]
We put as follows, for any integer $m\geq 1$:
\[
 \nbigk_m:=
 \bigl\{f\in\nbigk_0
 \,\big|\,
 \pi'_{i}(f_{|\nbigd_i})=0\,\,(i\leq m)
 \bigr\}.
\]

Let $\vecv=(v_j)$ be a frame of
$\naiveprolongg{\vecdelta}{\nbige^{(\lambda_0)}}$
compatible with $\EE^{(\lambda_0)}$ and $\Fzero$.
For each $v_j$,
we put as follows:
\[
\begin{array}{l}
 S_0(v_j):=
 \bigl\{i\,\big|\,
 \lefttop{i}\deg^{\Fzero,\EEzero}(v_j)
 =(1,-\lambda_0)
 \bigr\},
 \\
 \mbox{{}}\\
 S_m(v_j):=
 S_0(v_j)\sqcup
 \bigl\{i\,\big|\,i\leq m,\,\,
 \lefttop{i}\deg^{\EEzero}(v_j)\neq -\lambda_0\bigr\}.
\end{array}
\]
We put as follows:
\[
 \tilde{v}_{m,j}:=v_j\cdot
 \prod_{i\in S_m(v_j)}z_i,
\quad\quad
 \tilde{\vecv}_m=\bigl(\tilde{v}_{m,j}\bigr).
\]

\begin{lem}
$\tilde{\vecv}_m$ is a frame of $\nbigk_m$.
In particular,
$\nbigk_m$ is locally free.
\end{lem}
\pf
It immediately follows from the definition
of $\nbigk_m$.
\hfill\qed

\vspace{.1in}

We have
$ \nbigk_0\supset\nbigk_1\supset\cdots\supset
 \nbigk_l$.
We also have
$\nbigk_l\subset
 \naiveprolongg{(1-\epsilon)\vecdelta}{\nbige}^{(\lambda_0)}$.

\begin{lem} \label{lem;9.18.1}
For any section $f\in\nbigk_{m-1}$,
there exists $g\in\nbigk_m$ such that
$f-g\deldel_m\in\nbigk_m$.
\end{lem}
\pf
We can regard $f$ as a section of $\naiveprolongg{\vecdelta}{\nbige}$.
The restriction
$f_{|\nbigd_m}\in \naiveprolongg{\vecdelta}{\nbige}_{|\nbigd_m}$
is decomposed as
$f_{|\nbigd_m}=f_0+f_1$,
where $f_j\in \lefttop{m}\EEzero_j$
$(j=0,1)$
such that $\pi_m(f_0)=0$ in
$\Gr^{\Fzero}_1\bigl(\lefttop{m}\EE_0\bigr)$.

For any $i<m$,
the restriction
$f_{1\,|\,\nbigd_i\cap\nbigd_m}$
is contained in
$\lefttop{i}\EEzero_{0\,|\,\nbigd_i\cap\nbigd_m}
\cap
 \lefttop{m}\EEzero_{1\,|\,\nbigd_i\cap\nbigd_m}$,
and we have $\pi_i(f_{1\,|\,\nbigd_i\cap\nbigd_m})=0$
in
$\lefttop{i}\Gr^{\Fzero}_1(\lefttop{i}\EE_0)_{|\nbigd_i\cap\nbigd_m}$.

We put $g_1:=\bigl(\Res_m(\DD)+\lambda\bigr)^{-1}f_1\in
 \lefttop{m}\EE_1$.
Then for any $i<m$,
we have
$g_{1\,|\,\nbigd_i\cap\nbigd_m}\in
 \lefttop{i}\EE_{0\,|\,\nbigd_i\cap\nbigd_m}
 \cap
 \lefttop{m}\EE_{1|\nbigd_i\cap\nbigd_m}$
and
$\pi_i(g_{1\,|\,\nbigd_i\cap\nbigd_m})=0$
in $\lefttop{i}\Gr^{\Fzero}(\lefttop{i}\EE_0)_{|\nbigd_i\cap\nbigd_m}$.

We can pick $\tilde{g}\in\nbigk_{m-1}$
such that $\tilde{g}_{|\nbigd_m}=g_1$.
We put $g:=-\tilde{g}\cdot z_m\in\nbigk_{m-1}(-\nbigd_m)
\subset\nbigk_m$.
Then we have
$f-g\cdot \deldel_m\in \nbigk_{m-1}$,
and we have the following:
\[
 (f-g\cdot\deldel_m)_{|\nbigd_m}
=f_{|\nbigd_m}
-\bigl(
 \tilde{g}\cdot (-z_m\cdot\deldel_m)
 \bigr)_{|\nbigd_m}
=f_0+f_1-(\Res_m(\DD)+\lambda)\cdot g_1=f_0.
\]
Hence $f-\deldel_m g\in\nbigk_m$.
\hfill\qed

\begin{cor} \label{cor;10.2.4}
For any section $f\in\nbigk_0$,
there exist sections $g_1,\ldots,g_l\in\nbigk_l$
and $f_0\in\nbigk_l$
such that the following holds:
\[
 f=\sum_{m=1}^lg_m\cdot\prod_{j=1}^m \deldel_j
+f_0.
\]
\end{cor}
\pf
We have only to use Lemma \ref{lem;9.18.1}
inductively.
\hfill\qed

\subsubsection{The definition of the prolongment $\rmoduleprolong{E}$}

Let $\epsilon_0>0$ be as in $(C2)$.
The subsheaf $\gbige^{(\lambda_0)}$ of
$\naiveprolong{\nbige}^{(\lambda_0)}$
is obtained, which we will explain in the following.
We put as follows, for any element $\vecb\in\real_{<0}^l$:
\[
 \lefttop{\lbar}V^{(\lambda_0)}_{\vecb}
 \gbige^{(\lambda_0)}
:=\lefttop{\lbar}V^{(\lambda_0)}_{\vecb}
 \naiveprolong{\nbige}.
\]
Then we put as follows:
\[
 \lefttop{\lbar}\Vzero_{<0}\rmoduleprolong{E}^{(\lambda_0)}
:=\bigcup_{\vecb\in\real_{<0}^l}
 \lefttop{\lbar}\Vzero_{\vecb}\rmoduleprolong{E}^{(\lambda_0)}.
\]
The sheaf $\gbige^{(\lambda_0)}$ is the $\nbigr$-submodule
of $\naiveprolong{\nbige}^{(\lambda_0)}$
generated by $\lefttop{\lbar}V^{(\lambda_0)}_{<0}\gbige^{(\lambda_0)}$.

\begin{lem}
For $\lambda_1\in\Delta(\lambda_0,\epsilon_0)$,
we have the following:
\[
 \gbige^{(\lambda_0)}_{|\Delta(\lambda_1,\epsilon_1)}
=\gbige^{(\lambda_1)}.
\]
\end{lem}
\pf
Due to the condition $(C2)$ on $\epsilon_0$
and Lemma \ref{lem;9.19.10},
we obtain 
$\lefttop{\lbar}V^{(\lambda_0)}_{<0}
  \gbige^{(\lambda_0)}_{|\Delta(\lambda_1,\epsilon_1)}
\subset
 \lefttop{\lbar}V^{(\lambda_1)}_{<0}
  \gbige^{(\lambda_1)}$.
It implies the implication
$\gbige^{(\lambda_0)}_{|\Delta(\lambda_0,\epsilon_0)}
\subset\gbige^{(\lambda_1)}$.

Let us see
$\lefttop{\lbar}V^{(\lambda_1)}_{<0}\gbige^{(\lambda_1)}
\subset \gbige^{(\lambda_0)}_{|\Delta(\lambda_1,\epsilon_1)}$.
Due to the condition $(C2)$ on $\epsilon_0$,
we have
$\lefttop{\lbar}V^{(\lambda_1)}_{<0}\gbige^{(\lambda_1)}
\subset \nbigk_{0\,|\,\Delta(\lambda_1,\epsilon_1)}$.
Hence any sections
$f\in \lefttop{\lbar}V^{(\lambda_1)}_{<0}\gbige^{(\lambda_1)}$
is described as follows, due to Corollary \ref{cor;10.2.4}:
\[
 f=f_0+
 \sum_{m=1}^lg_m\cdot
 \prod_{j=m}^l\deldel_j,\quad
\quad
\bigl(
 f_0,g_m\in\lefttop{\lbar}V^{(\lambda_0)}_{<0}\gbige^{(\lambda_0)}
 _{|\Delta(\lambda_1,\epsilon_1)}
\bigr).
\]
Thus we obtain $f\in \gbige^{(\lambda_0)}_{|\Delta(\lambda_1,\epsilon_1)}$.
\hfill\qed

\begin{cor}
We obtain the global $\nbigr_{\nbigx}$-module $\gbige$.
It is strict.
\hfill\qed
\end{cor}

\subsubsection{The filtrations $\lefttop{i}\Vzero$ of $\gbige$}

For $\vecb\in\real^l$,
we put as follows:
\[
 \lefttop{\lbar}V^{(\lambda_0)}_{\vecb}\gbige
:=\sum_{\substack{
 \vecc<0,\vecn\in\seisuu_{\geq 0}\\
 \vecc+\vecn\leq \vecb}}
 \lefttop{\lbar}V^{(\lambda_0)}_{\vecc}\gbige\cdot
 \deldel^{\vecn}.
\]
For any subset $I\subset \lbar$,
we put as follows:
\[
 \lefttop{I}V^{(\lambda_0)}_{\veca}
 \gbige
:=\bigcup_{q_I(\vecb)=\veca}
 \lefttop{\lbar}V^{(\lambda_0)}_{\vecb}\gbige.
\]
In particular,
$\lefttop{\{i\}}V^{(\lambda_0)}_{\veca}\gbige$
is denoted by $\lefttop{i}V^{(\lambda_0)}_{\veca}\gbige$.

\subsubsection{Remark}
\label{subsubsection;b12.3.1}

Let $(E,\delbar_E,h,\theta)$ be a tame harmonic bundle
over $X-\bigcup_{i=1}^l D_i$.
Let $(E',\delbar_{E'},h',\theta')$ be the restriction
$(E,\delbar_E,h,\theta)$ to $X-\bigcup_{i=1}^n D_i$.
Then we obtain the two $\nbigr$-modules $\gbige(E)$
and $\gbige(E')$
on $\nbigx$.
They can be regarded as the subsheaf
of $i_{\ast}\nbige'$,
where $i$ denote the inclusion $X-\bigcup_{i=1}^lD_i\subset X$.
It is easy to see that we have
$\lefttop{\lbar}V^{(\lambda_0)}_{<0}\gbige(E)=
 \lefttop{\nbar}V^{(\lambda_0)}_{<0}\gbige(E')$
for any $\lambda_0$,
which implies $\gbige(E)=\gbige(E')$.
Hence we may assume $l=n$, if we need it.

%% file: a75.6.tex

\subsubsection{$\lefttop{I}T\kakkolambdazero(\vecc,\vecd)$
 and the easy properties}

\label{subsubsection;a11.23.20}

Let us pick a decomposition $I\sqcup J=\lbar$.
For any elements $\vecc\in\real^{I}$
and $\vecd\in\real^{J}_{<0}$,
we put as follows:
\[
 \lefttop{I}T^{(\lambda_0)}(\vecc,\vecd)
:=\frac{\lefttop{\lbar}\Vzero_{\vecc+\vecd}\gbige}
  {\sum_{\vecb\in S}\lefttop{\lbar}\Vzero_{\vecb}\gbige},
\]
Here we put
$S:=\{\vecb\in\real^l\,|\,q_I(\vecb)\lneq\vecc, q_J(\vecb)\leq \vecd\}$.

\begin{lem}
Let $\vecb=(b_1,\ldots,b_l)$ be an element of $\real^l$.
If we have $b_i+1\geq 0$ for some $i\in\lbar$,
then we have 
$\lefttop{\lbar}\Vzero_{\vecb+\vecdelta_{i}}(\gbige)
=\lefttop{\lbar}\Vzero_{\vecb}(\gbige)\cdot\deldel_i
+\lefttop{\lbar}\Vzero_{\vecb'}(\gbige)$.
Here we put $\vecb':=(b_1,\ldots,b_{i-1},-\epsilon,b_{i+1},\ldots,b_l)$
for some positive number $\epsilon$.
\end{lem}
\pf
The implication $\supset$ is clear.
We show the implication $\subset$.
Let $f$ be an element of
$\lefttop{\lbar}V_{\vecc}^{(\lambda_0)}\gbige$
such that $\vecc+\vecn\leq \vecb+\vecdelta_{i}$
and $\vecc\in\real_{<0}^I$.
In the case $q_i(\vecn)\geq 1$,
we have
$f\cdot \deldel^{\vecn-\vecdelta_{i}}\in\lefttop{\lbar}\Vzero_{\vecb}$,
which implies
$f\cdot\deldel^{\vecn}$ is contained in
$\lefttop{\lbar}\Vzero_{\vecb}\cdot\deldel_i$.
In the case  $q_i(\vecn)=0$,
we have the following:
\[
 q_j(\vecc+\vecn)\leq q_j(\vecb),\,
(j\neq i),
\quad
 q_i(\vecc)<0.
\]
Thus $f\cdot\deldel^{\vecn}$ is contained in
$\lefttop{\lbar}\Vzero_{\vecb'}$.
Therefore we are done.
\hfill\qed

\begin{cor}\label{cor;b11.19.1}
Let $\vecc$ be an element of $\real^I$ such that
$q_i(\vecc)+1\geq 0$ for some $i\in I$.
Let $\vecd\in\real^{J}_{<0}$.
Then the induced morphism
$\deldel_i:\lefttop{I}T^{(\lambda_0)}(\vecc,\vecd)
\lrarr
 \lefttop{I}T^{(\lambda_0)}(\vecc+\vecdelta_{i},\vecd)$
is surjective.
\hfill\qed
\end{cor}

Since we have the inclusion
$\lefttop{\lbar}\Vzero_{\vecb}(\gbige)
\lrarr\lefttop{\lbar}\Vzero_{\vecb}\bigl(\naiveprolong{\nbige}\bigr)$
for any $\vecb\in\real^l$,
we have the naturally defined morphism
$f_{\vecc,\vecd}:
 \lefttop{I}T^{(\lambda_0)}(\vecc,\vecd)
\lrarr
 \lefttop{I}\Gr^{\Vzero}_{\vecc}
 \lefttop{J}\Vzero_{\vecd}\bigl(\naiveprolong{\nbige}\bigr)$.

\begin{lem} \label{lem;9.18.10}
$\Image(f_{\vecc,\vecd})
=\lefttop{I}\tilde{T}^{(\lambda_0)}(\vecc,\vecd)$.
\end{lem}
\pf
In the case $\vecc<0$, the claim immediately follows from the definition
of $\lefttop{I}T^{(\lambda_0)}$
and $\lefttop{I}\tilde{T}^{(\lambda_0)}$,
namely, both of them are same as
$\lefttop{I}\Gr^{\Vzero}_{\vecc}\lefttop{J}\Gr_{\vecd}^{\Vzero}
 \bigl(\naiveprolong{\nbige}\bigr)$.

We have the following commutative diagramm:
\begin{equation} \label{eq;9.19.15}
 \begin{CD}
 \lefttop{I}T^{(\lambda_0)}(\vecc,\vecd)
 @>{f_{\vecc,\vecd}}>>
 \lefttop{I}\Gr^{\Vzero}_{\vecc}
 \lefttop{J}\Gr^{\Vzero}_{\vecd}(\naiveprolong{\nbige})\\
 @V{\deldel_i}VV @V{\deldel_i}VV \\
 \lefttop{I}T^{(\lambda_0)}(\vecc+\vecdelta_{i},\vecd)
 @>{f_{\vecc+\vecdelta_{i},\vecd}}>> 
 \lefttop{I}\Gr^{\Vzero}_{\vecc+\vecdelta_{i}}
 \lefttop{J}\Vzero_{\vecd}(\naiveprolong{\nbige}).
\end{CD}
\end{equation}
If $q_i(\vecc)+1\geq 0$,
the left $\deldel_i$ in the diagramm (\ref{eq;9.19.15})
is surjective, due to Corollary \ref{cor;b11.19.1}.
As for the right $\deldel_i$,
we have 
$\deldel_i(\lefttop{I}\tilde{T}^{(\lambda_0)}(\vecc,\vecd))
=\lefttop{I}\tilde{T}^{(\lambda_0)}(\vecc+\vecdelta_i,\vecd)$.
Then
it is easy to show the following implication:
\[
 \Image f_{(\vecc,\vecd)}=
 \lefttop{I}\tilde{T}^{(\lambda_0)}(\vecc,\vecd)\,\,
\Longrightarrow\,\,
 \Image(f_{\vecc+\vecdelta_i,\vecd})
 =\lefttop{I}\tilde{T}^{(\lambda_0)}(\vecc+\vecdelta_i,\vecd).
\]
Hence the general case can be reduced to the case $\vecc<0$.
\hfill\qed

\begin{rem}
Later we see that $f_{(\vecc,\vecd)}$ is isomorphic
(Lemma {\rm\ref{lem;9.18.11}}).
\hfill\qed
\end{rem}

%% file: 31.1.tex

\subsubsection{Preliminary}
\label{subsubsection;d11.19.5}

For any subset $I\subset\lbar$,
let $\vecdelta_I\in\real^l$ denote
the element determined by
the conditions
$q_i(\vecdelta_I)=1$ $(i\in I)$ and $q_i(\vecdelta_I)=0$ $(i\not\in I)$.
In particular, we put
$\vecdelta_{\mbar}:=(\overbrace{1,\ldots,1}^m,0,\ldots,0)$
and
$\vecdelta_m:=(\overbrace{0,\ldots,0}^{m-1},1,0,\ldots,0)$.

Let $m$ be an integer such that $1\leq m\leq l$.
Let $\vecb$ be an element of $\real_{\leq 0}^{l}$
such that $b_i=0$ for $i\leq m$.
For any subset $J\subset\lbar$,
we denote $q_J(\vecb)$ by $\vecb_J$.

Let $I$ be a subset of $\lbar$.
We put $J:=\lbar-I$.
We have the endomorphisms
$\Res_i(\DD)_{|\nbigd_I}+\lambda$
$(i\in I)$
of $\lefttop{I}\Gr^{\Vzero}_{0}
 \lefttop{J}\Vzero_{\vecb_J}\bigl(\naiveprolong{\nbige}\bigr)$.
We put as follows,
for $I\subset\lbar$ and $m\leq l$:
\[
 \lefttop{I}Q_m\bigl(
 \lefttop{\lbar}\Vzero_{\vecb}\bigl(
 \naiveprolong{\nbige}
 \bigr)
 \bigr):=
 \Image\Bigl(
 \prod_{\substack{i\in I,\\i\leq m}}
 \bigl(\Res_i(\DD)+\lambda\bigr)
 \Bigr)
\subset
 \lefttop{I}\Gr^{\Vzero}_{0}
 \lefttop{J}\Vzero_{\vecb_J}
 \bigl(\naiveprolong{\nbige}\bigr).
\]
We often omit to use the notation $\lefttop{I}Q_m$
instead of 
$\lefttop{I}Q_m\bigl(
 \lefttop{\lbar}\Vzero_{\vecb}\bigl(
 \naiveprolong{\nbige}\bigr)
 \bigr)$,
if there are no confusion.
We have the generalized eigen decomposition
for the tuple of the endomorphisms
$\Res_I(\DD)=\bigl(\Res_i(\DD)\,\big|\,i\in  I\bigr)$:
\begin{equation} \label{eq;9.19.11}
 \lefttop{I}\Gr^{\Vzero}_{0}
 \lefttop{J}\Vzero_{\vecb_J}
 \bigl(\naiveprolong{\nbige}\bigr)
=\bigoplus_{\vecbeta\in\cnum^I}
 \lefttop{I}\EE^{(\lambda_0)}
 \bigl(\Res_I(\DD),\vecbeta
 \bigr).
\end{equation}
Let $\nbign_{i\,\vecbeta}$
denote the nilpotent part of the restriction $\Res_i(\DD)$
to $ \lefttop{I}\EE
 \bigl(\Res_I(\DD),\vecbeta
 \bigr)$.
The decomposition (\ref{eq;9.19.11}) induces the decomposition:
\[
 \lefttop{I}Q_m=
\bigoplus_{\vecbeta\in\cnum^I}\lefttop{I}Q_{m,\vecbeta},
\quad\quad
 \lefttop{I}Q_{m,\vecbeta}
=\Image\Bigl(
 \prod_{\substack{i\in I,\,i\leq m \\ \beta_i=-\lambda_0}}
 \nbign_{i,\vecbeta}
 \Bigr).
\]

\begin{lem}
$\lefttop{I}Q_{m}\bigl(
 \lefttop{\lbar}\Vzero_{\vecb}\bigl(
 \naiveprolong{\nbige}
 \bigr)
 \bigr)$
is a vector subbundle
of $\lefttop{I}\Gr^{\Vzero}_{0}
 \lefttop{J}\Vzero\bigl(\naiveprolong{\nbige}\bigr)$.
\end{lem}
\pf
The conjugacy classes of the nilpotent maps
$\prod_{i\in I,i\leq m\,\beta_i=-\lambda_0}
  \nbign_{i,\vecbeta\,|\,(\lambda,P)}$
is independent of $(\lambda,P)$,
which follows from a limiting mixed twistor theorem.
(See Lemma \ref{lem;10.26.121}).
Then the lemma follows.
\hfill\qed

\vspace{.1in}

Let $m$ be an integer such that $0\leq m\leq l$.
Let $\vecb$ be an element of $\real_{\leq 0}^{l}$
such that $b_i=0$ for $i\leq m$.
For any subset $I\subset\mbar$,
we put $J:=\lbar-I$.
We have the projection
$\tilde{\pi}_I:
 \lefttop{\lbar}\Vzero_{\vecb}\bigl(
 \naiveprolong{\nbige}
 \bigr)
\lrarr
 \lefttop{I}\Gr^{\Vzero}_{0}
 \lefttop{J}\Vzero_{\vecb_J}\bigl(
 \naiveprolong{\nbige}
 \bigr)$.

Then we put as follows,
for $\vecd=\sum_{i=m+1}^l d_i\cdot\vecdelta_i$ $(-1\leq d_i<0)$:
\[
 \nbigl_{m,\vecd}
:=
 \bigl\{f\in \lefttop{\lbar}\Vzero_{\vecd}\bigl(
 \naiveprolong{\nbige}\bigr)
 \,\,\big|\,\,
 \tilde{\pi}_I(f)\in
 \lefttop{I}Q_m,\,\,\forall I\subset\mbar
 \bigr\}.
\]

\begin{lem} \label{lem;9.18.3}
For $m\geq 1$,
we have the following equality for some positive number $\epsilon$:
\begin{equation}\label{eq;a12.2.6}
 \nbigl_{m,\vecd}
=\nbigl_{m-1,\vecd-\vecdelta_m}\cdot\deldel_m
+\nbigl_{m-1,\vecd-\epsilon\vecdelta_m}
\end{equation}
\end{lem}
\pf
Let us see the implication $\supset$.
Clearly we have
$\nbigl_{m-1,\vecd-\epsilon\vecdelta_m}
\subset
 \nbigl_{m,\vecd}$.
Let $f$ be a section of $\nbigl_{m-1,\vecd-\vecdelta_m}$.
Since $\nbigl_{m-1,\vecd-\vecdelta_m}$ is contained in
$\lefttop{\lbar}\Vzero_{\vecd-\vecdelta_m}(\naiveprolong{\nbige})$,
we have
$f\cdot \deldel_m \in
 \lefttop{\lbar}\Vzero_{\vecd}(\naiveprolong{\nbige})$.
For any subset $I\subset \mbar$ such that $m\not\in I$,
we have the following:
\[
 \tilde{\pi}_I(f)\in
 \lefttop{I}Q_{m-1}\bigl(
 \lefttop{\lbar}\Vzero_{\vecd-\vecdelta_m}\bigl(
 \naiveprolong{\nbige}
 \bigr)
 \bigr),
\quad\quad
 \tilde{\pi}_I\Bigl((f\cdot\deldel_m)\Bigr)
\in\lefttop{I}Q_{m-1}\bigl(
 \lefttop{\lbar}\Vzero_{\vecd}
 \bigl(\naiveprolong{\nbige}\bigr)
 \bigr)
=\lefttop{I}Q_m\bigl(
 \lefttop{\lbar}\Vzero_{\vecd}
 \bigl(\naiveprolong{\nbige}\bigr)
 \bigr).
\]
We put $\tilde{f}:=f\cdot z_m^{-1}$,
and then we have 
$ (f\cdot \deldel_m)_{|\nbigd_m}
=-\bigl(\Res_m(\DD)+\lambda\bigr)\tilde{f}_{|\nbigd_m}$.
For a subset $I\subset\lbar$ such that $m\not\in I$,
we put $I'=I\sqcup\{m\}$.
Note we have the following:
\[
\tilde{\pi}_{I'}(\tilde{f}_{|\nbigd_{I'}})
\in
 \lefttop{I'}Q_{m-1}\bigl(
 \lefttop{\lbar}\Vzero_{\vecd}\bigl(
 \naiveprolong{\nbige}
 \bigr)
 \bigr).
\]
Thus we obtain
$\tilde{\pi}_{I'}\bigl((f\cdot\deldel_m)_{|\nbigd_{I'}}\bigr)
\in \lefttop{I'}Q_m\bigl(
 \lefttop{\lbar}\Vzero_{\vecd}\bigl(
 \naiveprolong{\nbige}
 \bigr)
 \bigr)$.
Thus we obtain the implication $\supset$
in (\ref{eq;a12.2.6}).

\vspace{.1in}

Let us show the implication $\subset$ in (\ref{eq;a12.2.6}).
It will be accomplished after Lemma \ref{lem;9.18.4}.
For $\vecc=(c_1,\ldots,c_m)\in\seisuu_{>0}^m$,
we have the morphism $\phi_{\vecc}:X\lrarr X$
given by $(z_1,\ldots,z_n)\longmapsto
 (z_1^{c_1},\ldots,z_m^{c_m},z_{m+1},\ldots,z_n)$.
Let $\vecv$ be a holomorphic frame
of $\lefttop{\lbar}\Vzero_{\vecd}\bigl(
 \naiveprolong{\nbige}\bigr)$,
which is compatible with $\Fzero$ and $\EEzero$.

Since we have $\lefttop{\lbar}\Vzero_{\vecd}\bigl(
 \naiveprolong{\nbige}
 \bigr)=\naiveprolongg{\vecd+\vecdelta_{\lbar}}{\nbige}^{(\lambda_0)}$,
we have the positivity
$0<b_j(v_i):=\lefttop{j}\deg^{\Fzero}(v_i)$
for $j\in\mbar$.
We have the integers $N_j(v_i)\in\seisuu_{\geq\,0}$ 
for any $j\in\mbar$,
determined as follows:
\[
0<c_j\cdot b_j(v_i)-N_j(v_i)\leq 1.
\]
If $c_j$ is sufficiently large,
we have $N_j(v_i)\neq N_j(v_k)$
in the case $b_j(v_i)\neq b_j(v_k)$.

\begin{lem}
If $\lefttop{j}\deg^{\Fzero}(v_i)=1$,
we have $N_j(v_i)=c_j-1$.
\end{lem}
\pf
Note $0<c_j\cdot 1-(c_j-1)=1\leq 1$.
\hfill\qed

\vspace{.1in}

Let us pick $\epsilon_0>0$ satisfying $(C2)$
for $\phi_{\vecc}^{-1}\nbige$
(the subsubsection \ref{subsubsection;a12.2.7}).
We put as follows:
\[
 \tilde{v}_i:=\phi_{\vecc}^{-1}(v_i)\cdot\prod_{j=1}^m z_j^{N_j(v_i)},
\quad\quad
 \tilde{\vecv}=(\tilde{v}_i).
\]
Then $\tilde{\vecv}$ is a holomorphic frame of
$\lefttop{\lbar}\Vzero_{\vecd}\bigl(
 \naiveprolong{
 \phi_{\vecc}^{-1}\nbige}
 \bigr)$,
which is compatible with $\EEzero$ and $\Fzero$.

We have the $\mu_{\vecc}=\prod_{i=1}^m\mu_{c_i}$-action 
on $X$ as usual.
The action is lifted to the action on 
$\lefttop{\lbar}\Vzero_{\vecd}\bigl(
 \naiveprolong{
 \phi_{\vecc}^{-1}\nbige}
 \bigr)$.
For $j\in\mbar$,
we have the weight decomposition:
\[
 \lefttop{\lbar}\Vzero_{\vecd}\bigl(
 \naiveprolong{
 \phi_{\vecc}^{-1}\nbige}
 \bigr) _{|\nbigd_j}
=\bigoplus_{a=0}^{c_j-1}\lefttop{j}U_a.
\]
Here the action of $\mu_{c_j}$ on $\lefttop{j}U_a$ is of weight $a$.
If $c_j$ is sufficiently large,
we have the map
$\varphi_j:\{a\,|\,\lefttop{j}U_a\neq 0\}\lrarr 
\Par\bigl(
 \prolongg{\vecdelta_{\lbar}+\vecd}{
 \phi_{\vecc}^{\ast}\nbige^{(\lambdazero)}},
 j\bigr)$,
and the decomposition gives the splitting of the filtration
$\lefttop{j}\Fzero$:
\[
 \lefttop{j}\Fzero_b
=\bigoplus_{\varphi_j(a)\leq b}\lefttop{j}U_a.
\]
(See the subsubsection \ref{subsubsection;c11.16.1}
for such a splitting.)

\vspace{.1in}
We put $N_i=c_i-1$ $(i\in I)$,
and we put $\vecN_I:=(N_i\,|\,i\in I)\in\seisuu^I$.
Then the frames $\vecv$ and $\tilde{\vecv}$ induce the isomorphism,
for $I\subset\mbar$:
\begin{equation} \label{eq;9.19.21}
 \lefttop{I}U_{\vecN_I}
=\bigcap_{i\in I}\lefttop{i}U_{N_i}
\simeq
 \lefttop{I}\Gr^{\Vzero}_{0}
 \lefttop{J}\Vzero_{\vecd}
\bigl(\naiveprolong{\nbige}\bigr).
\end{equation}
We have the endomorphisms
$ \Res_i\bigl(\DD\bigr)+\lambda$ on the right hand side
of (\ref{eq;9.19.21}).
On the other hand,
we have the endomorphism
$c_i^{-1}\cdot\bigl(\Res(\phi_{\vecc}^{\ast}\DD+\lambda)\bigr)$
on the left hand side.
\begin{lem}
Under the isomorphism {\rm(\ref{eq;9.19.21})},
the endomorphism
$c_i^{-1}\cdot\bigl(\Res(\phi_{\vecc}^{\ast}\DD+\lambda)\bigr)$
corresponds to
$ \Res_i\bigl(\DD\bigr)+\lambda$.
\end{lem}
\pf
If $\DD\vecv=\vecv\cdot\sum_i\nbiga_i\cdot dz_i/z_i$,
then we have the following:
\[
 \bigl(
\phi_{\vecc}^{-1}\DD\bigr)\tilde{\vecv}
=\tilde{\vecv}\cdot
 \Bigl(
 \sum_i
 \phi_{\vecc}^{-1}\nbiga_i \cdot c_i\cdot\frac{dz_i}{z_i}
 +\sum \nbign_i\cdot\lambda\frac{dz_i}{z_i}
 \Bigr).
\]
Here $\nbign_i$ denote the diagonal matrices
such that $(\nbign_i)_{j\,j}=N_i(v_j)$,
and we put $c_i=1$ for $i\geq m+1$,
for simplicity.

We put
$\tilde{\vecv}_I
 :=(\tilde{v}_i\,|\,\lefttop{I}\deg^{\Fzero}(v_i)=\vecdelta_I)$.
Then we have the following:
\[
 \Res_i(\phi_{\vecc}^{-1}\DD)\tilde{\vecv}_I
=\tilde{\vecv}_I\cdot
 \Bigl(
 c_i\cdot\phi_{\vecc}^{-1}\nbiga_{i|\lefttop{i}U_{N_i}}
+(c_i-1)\cdot\lambda
 \Bigr).
\]
Hence we obtain the following:
\[
 \Bigl(
 \Res_i(\phi_{\vecc}^{-1}\DD)+\lambda
 \Bigr)\tilde{\vecv}_I
=\tilde{\vecv}_I\cdot
 c_i\cdot
 \Bigl(
 \phi_{\vecc}^{-1}\nbiga_{|\lefttop{i}U_{N_i}}+\lambda
 \Bigr).
\]
Thus we are done.
\hfill\qed

\vspace{.1in}

Let us consider a section $f$ of $\nbigl_{m,\vecd}$
and consider the pull back $\phi_{\vecc}^{-1}f$.
Recall we put $N_j=c_j-1$
for $j=1,\ldots,m$.
We put as follows:
\[
 \tilde{f}:=
 \phi_{\vecc}^{-1}f\cdot
 \prod_{j=1}^m z_j^{N_j}.
\]
Then we obtain the section
$\tilde{f}\in
 \lefttop{\lbar}\Vzero_{\vecd}
 \bigl(
 \naiveprolong{
 \phi_{\vecc}^{\ast}\nbige}\bigr)$.

\begin{lem}
\mbox{{}}
\begin{itemize}
\item
We have $\tilde{f}_{|\nbigd_i}\in \lefttop{i}\Fzero_{d_i}$
for any $i>m$.
\item
We have $\tilde{f}_{|\nbigd_i}\in \lefttop{i}U_{N_i}$
for any $i\leq m$.
\end{itemize}
\end{lem}
\pf
The first claim is clear.
Note that our $\vecc$ is an element of $\seisuu_{>0}^m$.

We have the description $f=\sum f_j\cdot v_j$,
and then we have the following:
\begin{equation} \label{eq;9.19.20}
 \tilde{f}=
 \sum \phi_{\vecc}^{-1}(f_j)\cdot
 \prod_{i=1}^m z_i^{N_i-N_i(v_j)}\cdot \tilde{v}_j.
\end{equation}
The second claim immediately follows from (\ref{eq;9.19.20}).
\hfill\qed

\vspace{.1in}

We take the subbundles
$\lefttop{I}\tilde{Q}_m$ of
$\lefttop{I}U_{\vecN_I}$ for $I\subset\mbar$:
\[
 \lefttop{I}\tilde{Q}_{m}
:=\Image\Bigl(
 \prod_{i\in I}\bigl(
 \Res_i(\phi_{\vecc}^{-1}\DD)+\lambda\bigr)
 \Bigr)\cap \lefttop{I}U_{\vecN_I}.
\]
We also take the vector subbundles
$\lefttop{I}\tilde{R}_m$
of $\lefttop{I}U_{\vecN_I}$
for any subset $I\subset\mbar$ such that $m\in I$:
\[
 \lefttop{I}\tilde{R}_m:=
\Image\Bigl(
 \prod_{i\in I,i<m}
 \bigl( \Res_i(\phi_{\vecc}^{-1}\DD)+\lambda
 \bigr)
 \Bigr)\cap\lefttop{I}U_{\vecN_I}.
\]
In the case $m\in I$,
we have the equivariant inclusion
$\lefttop{I}\tilde{Q}_m\subset\lefttop{I}\tilde{R}_m$
and the naturally defined equivariant surjection:
\begin{equation}\label{eq;c11.18.1}
 \Res_m(\phi_{\vecc}^{-1}\DD)+\lambda:
\lefttop{I}\tilde{R}_m\lrarr
\lefttop{I}\tilde{Q}_m.
\end{equation}

We have
$\tilde{f}_{|\nbigd_I}\in\lefttop{I}\tilde{Q}_m$
 for any subset $I\subset\mbar$,
for we have $\lefttop{I}\tilde{Q}_m\simeq \lefttop{I}Q_m$
under the isomorphism 
$\lefttop{I}U_{\vecN_I}\simeq
\lefttop{I}\Gr^{\Vzero}_{0}
\lefttop{J}\Vzero_{\vecd}
 \bigl(
 \naiveprolong{\nbige}
 \bigr)$
given in (\ref{eq;9.19.21}).

\begin{lem} \label{lem;9.19.40}
For any subset $I$ such that $m\in I$,
we can pick the equivariant section $g_I\in\lefttop{I}\tilde{R}_m$,
satisfying the following conditions:
\[
 c_m^{-1}\cdot
 \bigl(
 \Res_m(\phi_{\vecc}^{-1}\DD)+\lambda\bigr)\cdot g_I
=\tilde{f}_{|\nbigd_I},
\quad\quad
 g_{I\,|\,\nbigd_{I'}}=g_{I'},\,\,(I\subset I').
\]
\end{lem}
\pf
In the case $I=\mbar$,
we have the surjection
$\lefttop{\mbar}\tilde{R}_m\lrarr\lefttop{\mbar}\tilde{Q}_m$,
thus we have only to take an appropriate lift.

Assume that we have already picked $g_{I'}$
for any subsets $I\subsetneq I'$,
and we will construct $g_I$.
We would like to take a lift of $\tilde{f}_{|\nbigd_I}$
for the morphism (\ref{eq;c11.18.1}).
The data 
$(g_{I'}\,|\,I'\supsetneq I)$ gives
an equivariant lift of $\tilde{f}_{|\del\nbigd_I}$,
where we put
$\del\nbigd_I:=
 \bigl(\bigcup_{I\subsetneq I'}\nbigd_{I'}\bigr)$.
We have only to extend it equivariantly.
\hfill\qed

\begin{lem}
We may take an equivariant section
$\tilde{g}\in
 \lefttop{\lbar}\Vzero_{\vecd}\bigl(
 \naiveprolong{\phi_{\vecc}^{\ast}\nbige}
 \bigr)$,
satisfying the following:
\begin{itemize}
\item
 $\tilde{g}_{|\nbigd_I}\in\lefttop{I}U_{\vecN_I}$
 for $I\subset \mbar$.
\item
 $\tilde{g}_{|\nbigd_I}=g_I$
 for $I\subset\mbar$ such that $m\in I$.
\item
 $\tilde{g}$ satisfies the same transformation rule
 as $\tilde{f}$.
\end{itemize}
\end{lem}
\pf
Note the following:
If $m\not\in I\subset\mbar$,
then we have
$g_{I\sqcup\{m\}}\in
 \lefttop{I}U_{\vecN_I\,|\,\nbigd_{I\sqcup\{m\}}}$.
Then the lemma can be shown by an argument
similar to the proof of Proposition \ref{prop;10.9.10}
\hfill\qed

\vspace{.1in}

Note that $\tilde{g}\cdot\prod_{j=1}^m z^{-N_j}$
is equivariant and satisfying the same transformation
as $\phi_{\vecc}^{-1}f$.
Thus we have some section $g'\in\nbige$
such that
$\phi_{\vecc}^{-1}g'=\tilde{g}\cdot \prod_{j=1}^mz_j^{-N_j}$.
We put $g=-z_m\cdot g'$.

\begin{lem} \label{lem;9.18.5}
$g\in\nbigl_{m-1,\vecd-\vecdelta_m}$.
\end{lem}
\pf
For $i>m$,
we have $\lefttop{i}\deg(g)\leq d_i$,
which is clear from our construction.

Since we have $\lefttop{m}\deg(g')\leq 1$,
we have $\lefttop{m}\deg(g)\leq 0$.
For $I\subset\mminusitibar$,
we have the following,
under the isomorphism given in (\ref{eq;9.19.21}):
\[
 \lefttop{I}Q_{m-1}
 \bigl(\lefttop{\lbar}\Vzero_{\vecd}
 \bigl(\naiveprolong{\nbige}\bigr)
 \bigr)
=\lefttop{I}Q_m
 \bigl(\lefttop{\lbar}\Vzero_{\vecd}
 \bigl(\naiveprolong{\nbige}\bigr)
 \bigr)
\simeq \lefttop{I}\tilde{Q}_m.
\]
Hence we have
$\tilde{\pi}_I(g)=-\tilde{\pi}_I(g'\cdot z_m)\in
 \lefttop{I}Q_{m-1}\bigl(
 \lefttop{\lbar}\Vzero\bigl(
 \naiveprolong{\nbige}
 \bigr)
  \bigr)$.
Thus we obtain the result.
\hfill\qed

\begin{lem} \label{lem;9.18.4}
We have $f-g\cdot\deldel_m \in \nbigl_{m-1,\vecd-\epsilon\cdot\vecdelta_m}$
for some $\epsilon>0$.
\end{lem}
\pf
We have the following:
\[
 (g\cdot\deldel_m)_{|\nbigd_m}=
\bigl(g'\cdot(-z_m\cdot\deldel_m)\bigr)_{|\nbigd_m}
=(\Res_m(\DD)+\lambda)g'_{|\nbigd_m}.
\]
Thus we have
$\tilde{\pi}_{I}\bigl((g\cdot\deldel_m -f)_{|\nbigd_I} \bigr)=0$
for $I\subset\mbar$ such that $m\in I$,
by our construction.
In the case $m\not\in I$,
since $\deldel_m$ preserves $\lefttop{I}Q_{m-1}$,
it is easy to see that
$\tilde{\pi}_I(f-g\cdot\deldel_m)\in Q_{m-1}$.
Thus we are done.
\hfill\qed

\vspace{.1in}
The implication $\subset$ in the claim of
Lemma \ref{lem;9.18.3} immediately follows from
Lemma \ref{lem;9.18.5} and Lemma \ref{lem;9.18.4}.
Thus the proof of Lemma \ref{lem;9.18.3} is accomplished.
\hfill\qed

\begin{cor}\label{cor;d11.19.20}
We have the following, for some positive number $\epsilon$:
\[
 \nbigl_{m,\vecd}
=\sum_{I\subset\mbar}
  \nbigl_{0,\vecd-\vecdelta_I-\epsilon\cdot\vecdelta_{\mbar-I}}\cdot
\Bigl(
\prod_{i\in I}\deldel_i
 \Bigr).
\]
\end{cor}
\pf
We have only to use Lemma \ref{lem;9.18.3} inductively.
\hfill\qed

%% file: a75.3.tex

\subsubsection{The injectivity of $f_{\vecc,\vecd}$}

Recall that $f_{\vecc,\vecd}$ is surjective
(Lemma \ref{lem;9.18.10}).

\begin{lem} \label{lem;9.18.11}
The morphism $f_{\vecc,\vecd}$ is injective.
As a result it gives the isomorphism:
\[
  \lefttop{I}T^{(\lambda_0)}(\vecc,\vecd)
\stackrel{\sim}{\lrarr}
 \lefttop{I}\tilde{T}^{(\lambda_0)}(\vecc,\vecd).
\]
\end{lem}
\pf
First we note that
$f_{\vecc,\vecd}$ 
for $\vecc<0$ is isomorphic by definition.

We have the following commutative diagramm due to Lemma \ref{lem;9.18.10}:
\begin{equation} \label{eq;9.19.16}
 \begin{CD}
 \lefttop{I}T^{(\lambda_0)}(\vecc,\vecd)
 @>{f_{\vecc,\vecd}}>>
 \lefttop{I}\tilde{T}^{(\lambda_0)}(\vecc,\vecd)\\
 @V{\deldel_i}VV @V{\deldel_i}VV \\
 \lefttop{I}T^{(\lambda_0)}(\vecc+\vecdelta_i,\vecd)
 @>{f_{\vecc+\vecdelta_i,\vecd}}>>
 \lefttop{I}\tilde{T}^{(\lambda_0)}(\vecc+\vecdelta_i,\vecd).\\
 \end{CD}
\end{equation}
Let us assume $q_i(\vecc)+1> 0$.
Then the right $\deldel_i$ in the diagramm (\ref{eq;9.19.16})
is isomorphic,
(Corollary \ref{cor;c11.19.10}),
and the left $\deldel_i$ is surjective
(Corollary \ref{cor;b11.19.1}).
Moreover,
we have already known $f_{\vecc,\vecd}$ and $f_{\vecc+\vecdelta_i,\vecd}$
are surjective, due to Lemma \ref{lem;9.18.10}.
Hence,
under the assumption $q_i(\vecc)+1> 0$,
if $f_{\vecc,\vecd}$ is isomorphic,
then $f_{\vecc+\vecdelta_i,\vecd}$ and the left $\deldel_i$ 
in the diagramm (\ref{eq;9.19.16}) is isomorphic.

Thus we have only to show the injectivity of $f_{\vecc,\vecd}$
in the case $c_j\leq 0$ for any $j$.
For such $\vecc$, we put $ m(\vecc):=\#\{i\,|\,c_i=0\}$,
and we use an induction on $m(\vecc)$.

If $m(\vecc)=0$, then we obtain $\vecc<0$,
and thus $f_{\vecc,\vecd}$ is isomorphic.
Assume that we have already known the
injectivity of $f_{\vecc,\vecd}$ for $m(\vecc)\leq m-1$,
and then we may derive the injectivity of $f_{\vecc,\vecd}$
for $m(\vecc)=m$.

We may assume that $I=\kbar$.
We put $J:=\lbar-I$.
Let $\vecc$ be an element of $\real_{\leq 0}^{k}$
such that $c_i=0$ $(i\leq m-1)$, $c_m=-1$,
and $c_i<0$ $(m<i\leq k)$.
Let $\vecd$ be an element of $\real^{J}$.
Due to the hypothesis of the induction,
$f_{\vecc,\vecd}$ is isomorphic.
Once we prove that
$f_{\vecc+\vecdelta_m,\vecd}$ is isomorphic,
then the induction can proceed.

Hence Lemma \ref{lem;9.18.11} is reduced
to the following lemma.

\begin{lem} \label{lem;d11.19.2}
Under the isomorphism $f_{\vecc,\vecd}$,
we have the following:
\begin{equation}\label{eq;d11.19.1}
  \ker\Bigl(
 -\deldel_m:\lefttop{\kbar}
 \tilde{T}^{(\lambda_0)}(\vecc,\vecd)
\lrarr
 \lefttop{\kbar}\tilde{T}^{(\lambda_0)}(\vecc+\vecdelta_m,\vecd)
 \Bigr)
=
 \ker\Bigl(
 -\deldel_m:\lefttop{\kbar}
 T^{(\lambda_0)}(\vecc,\vecd)
\lrarr
 \lefttop{\kbar}T^{(\lambda_0)}(\vecc+\vecdelta_m,\vecd)
 \Bigr).
\end{equation}
\end{lem}
\pf
The implication $\supset$ is clear from the commutative diagramm
(\ref{eq;9.19.16}).
Let us show the implication $\subset$.
Let $s$ be a section of the left hand side of (\ref{eq;d11.19.1}).
Let take a section $g$ of 
$\lefttop{\lbar}
  \Vzero_{\vecc+\vecd}(\naiveprolong{\nbige})$ as in
Lemma \ref{lem;9.19.41}.
The section $g$ induces
the section of $\lefttop{I}T^{(\lambda_0)}(\vecc,\vecd)$,
which is the inverse image $f_{\vecc,\vecd}^{-1}(s)$.
Hence we have only to show that
the element 
$g\cdot\deldel_m$
induces the trivial element of
$\lefttop{I}T^{(\lambda_0)}(\vecc+\vecdelta_m,\vecd)$.

We put $\vecc':=\vecc+\vecdelta_m$.
Then $g\cdot\deldel_m$ is a section of
$\lefttop{\lbar}\Vzero_{\vecc'+\vecd}\bigl(\naiveprolong{\nbige}\bigr)$.
For any subset $I\subset\mbar$,
we have the projection
$\tilde{\pi}_I:\lefttop{\lbar}\Vzero_{\vecc'+\vecd}\bigl(
 \naiveprolong{\nbige}\bigr)
\lrarr
 \lefttop{I}\Gr^{\Vzero}_{0}
 \lefttop{J}\Vzero_{q_J(\vecc'+\vecd)}\bigl(
 \naiveprolong{\nbige}\bigr)$,
as in the subsubsection \ref{subsubsection;d11.19.5}.
By our choice of $g$,
we have the following, in the case $I\subset\mbar$:
\[
 \tilde{\pi}_I\bigl(-\big(g\cdot\deldel_m\big)_{|\nbigd_I}\bigr)
\in\Image\Bigl(
 \prod_{\substack{i\in I,\\ i\leq m-1}}(\Res_i(\DD)+\lambda)
 \Bigr).
\]
Moreover, we have 
$ \tilde{\pi}_I\bigl(-\big(g\cdot\deldel_m\big)_{|\nbigd_I}\bigr)=0$
in the case $m\in I\subset\mbar$.

We put
$\vecd':=(c_{m+1},\ldots,c_{k},d_{k+1},\ldots,d_l)
 \in \real^{\lbar-\mbar}$.
Then we obtain
$g\cdot\deldel_m\in\nbigl_{m-1,\vecd'-\epsilon\cdot\vecdelta_m}$.
Due to Corollary \ref{cor;d11.19.20},
we obtain the following:
\[
 g\cdot\deldel_m\in 
 \sum_{I\subset\underline{m-1}}
 \nbigl_{0,\vecd'-\vecdelta_I-\epsilon\cdot\vecdelta_{I^c}}
 \cdot\prod_{i\in I}\deldel_i.
\]
Here we put $I^c:=\mbar-I$.
Note
$\nbigl_{0,\vecd'-\epsilon\cdot\vecdelta_{I^c}-\vecdelta_I}
\cdot\deldel_I
\subset
 \lefttop{\lbar}\Vzero_{\vecc'+\vecd}(\gbige)$.
Hence we obtain $g\cdot\deldel_m \in
\lefttop{\lbar}\Vzero_{\vecc+\vecdelta_m+\vecd-\epsilon\vecdelta_m}
 \bigl(\gbige\bigr)$,
which means the vanishing
$g\cdot\deldel_m =0$ in
$\lefttop{\kbar}T^{(\lambda_0)}(\vecc,\vecd)$.
Therefore we obtain the implication
$\subset$ in (\ref{eq;d11.19.1}).
Thus we obtain Lemma \ref{lem;d11.19.2}
and Lemma \ref{lem;9.18.11}.
\hfill\qed

%% file: a77.2.tex

\subsubsection{$\lefttop{\lbar}\psi_{\vecu}(\gbige)$
 and $\lefttop{\lbar}\tildepsi_{\vecu}(\gbige)$}
\label{subsubsection;b11.23.20}

In particular,
we have the isomorphism
$\lefttop{\lbar}T^{(\lambda_0)}(\vecc)\simeq
 \lefttop{\lbar}\tilde{T}^{(\lambda_0)}(\vecc)$
for any $\vecc\in\real^l$.
For any $\vecu\in\prod_{i=1}^l\nbigk(\nbige,\lambda_0,i,c_i)$,
we put as follows:
\[
 \lefttop{\lbar}\psizero_{\vecu}(\gbige):=
 \lefttop{\lbar}\tilde{T}^{(\lambda_0)}(\vecu)
\subset\lefttop{\lbar}\tilde{T}^{(\lambda_0)}(\vecc).
\]
We have another characterization.
The action of
$\prod_{u\in\nbigk(\nbige,\lambda,i,c_i)}
 \bigl(z_i\cdot\deldel_i+\eigenmap(\lambda,u)\bigr)^N$
on $\lefttop{\lbar}\Gr^{\Vzero}_{\vecc}(\gbige)$
is $0$ if $N$ is sufficiently large.
Hence we obtain the following:
\[
 \lefttop{\lbar}\psizero_{\vecu}(\gbige)
=\bigcap_{i=1}^l
 \Ker\bigl(z_i\deldel_i+\eigenmap(\lambda,u_i)\bigr)^N.
\]
We have the decomposition:
\[
 \lefttop{\lbar}\Gr^{\Vzero}_{\vecc}=
\bigoplus_{\vecu\in\prod_i\nbigk(\nbige,\lambda,i,c_i)}
 \lefttop{\lbar}\psizero_{\vecu}(\gbige).
\]

For any $\vecu\in\prod_i\nbigk(\nbige,\lambda,i,c_i)$,
we take a sufficiently large integer $N$
such that $\paramap(\lambda_0,u_i)-N\vecdelta_{0}<0$
for any $i$.
We put
$\lefttop{\lbar}\tildepsi_{\vecu}(\gbige)
:=\lefttop{\lbar}\psi_{\vecu-N\cdot\vecdelta_{0,\lbar}}(\gbige)$.
Here $\vecdelta_{0,\lbar}$ denote the element
$(\overbrace{\vecdelta_0,\ldots,\vecdelta_0}^l)\in
 (\cnum\times\real)^l$.
It is well defined in the sense
that we have the canonical isomorphism
$\prod_{i}z_i^{N'-N}:
 \lefttop{\lbar}\tildepsi_{\vecu-N\cdot\vecdelta_{0,\lbar}}(\gbige)
:=\lefttop{\lbar}\psi_{\vecu-N'\cdot\vecdelta_{0,\lbar}}(\gbige)$
for two choices $N$ and $N'$.

Recall that we have the coherent sheaves
$\lefttop{\lbar}\nbigg_{\vecu}(E)$ $(\vecu\in \KMS(\nbige^0,\lbar))$
on $\cnum_{\lambda}$
(the subsubsection \ref{subsubsection;10.16.10}).
We have
$\lefttop{\lbar}\nbigg^{(\lambda_0)}_{\vecu+\vecdelta_{0,\lbar}}(E)
\simeq
 \lefttop{\lbar}\psizero_{\vecu}(\gbige)$
by definition,
when $\paramap(\lambda_0,\vecu)<0$.
Since we have the canonical isomorphism
$\prod_{i} z_i^N:
 \lefttop{\lbar}\nbigg_{\vecu}^{(\lambda_0)}(E)
\lrarr
 \lefttop{\lbar}\nbigg_{\vecu-N\cdot\vecdelta_{0,\lbar}}(E)$,
we obtain the isomorphism
$\lefttop{\lbar}\tildepsizero_{\vecu}(\gbige)
\simeq
 \lefttop{\lbar}\nbigg_{\vecu}^{(\lambda_0)}(E)$
for any $\vecu\in\prod_i\nbigk(\nbige,\lambda,i,c_i)$.

We have the nilpotent part $\nbign_i$ of the residue
$\Res_i(\DD)$ on $\lefttop{\lbar}\nbigg^{\lambda_i}_{\vecu}(E)$.
On the other hand,
we have the nilpotent part $\nbign_i$ of
$-z_i\cdot\deldel_i$ on $\lefttop{\lbar}\tildepsizero_{\vecu}(\gbige)$.
\begin{lem}
The isomorphism preserves the nilpotent morphisms $\nbign_i$.
\end{lem}
\pf
The $\Res_i(\DD)$ corresponds to the left action
of $z_i\cdot\deldel_i=\deldel_i\cdot z_i-\lambda$,
which corresponds to the right action of
$-z_i\cdot\deldel_i-\lambda$.
\hfill\qed

In particular, we obtain the following.
\begin{lem}\mbox{{}}
\begin{itemize}
\item
When $\Delta(\lambda_1,\epsilon_1)\subset\Delta(\lambda_0,\epsilon_0)$,
we have the canonical isomorphism
$\lefttop{\lbar}\tildepsizero_{\vecu}(\gbige)_{|\nbigx(\lambda_1,\epsilon_1)}
\simeq
 \lefttop{\lbar}\tildepsi^{(\lambda_1)}_{\vecu}(\gbige)_{|
 \nbigx(\lambda_1,\epsilon_1)}$.
\item
Hence we obtain the globally defined $\nbigr$-module
$\lefttop{\lbar}\tildepsi_{\vecu}(\gbige)$
on $\nbigd_{\lbar}$.
\item
We have the isomorphism
$\lefttop{\lbar}\nbigg_{\vecu+\vecdelta_{0,\lbar}}(E)
\simeq
 \lefttop{\lbar}\tildepsi_{\vecu}(\gbige)$.
It preserves the nilpotent parts of
$\Res(\DD)$ and $-z_i\cdot\deldel_i$.
\hfill\qed
\end{itemize}
\end{lem}

%% file: a76.1.tex

\subsubsection{Preliminary}

\label{subsubsection;a12.2.30}

For $\vecb\in\real^I$,
we put as follows:
\[
\begin{array}{l}
 \nbigs(\vecb):=
 \bigl\{\vecc\in\real^I
 \,\big|\, \vecc\leq \vecb\bigr\},\\
 \mbox{{}}\\
 \nbigs^0(\vecb):=
 \bigl\{\vecc\in\nbigs(\vecb)\,\big|\,\vecc\neq \vecb\bigr\}.
\end{array}
\]

\begin{df}
Let $S\subset\real^I\cap \KMS(\nbigelambdazero,I)$ be a finite subset.
It is called primitive,
if the following holds:
\begin{itemize}
\item
 For any $\vecb,\vecb'\in\nbigs$ such that $\vecb\neq \vecb'$,
 we have  $\vecb\not\in\nbigs(\vecb')$.
\hfill\qed
\end{itemize}
\end{df}

For a finite subset $S\subset\real^I\cap\KMS(\nbigelambdazero,I)$,
we put as follows:
\[
 \nbigs(S):=\bigcup_{\vecb\in S}\nbigs(\vecb),
\quad
 \nbigs^0(S):=\bigcup_{\vecb\in S}\nbigs^0(\vecb).
\]

The following lemma can be checked elementarily.
\begin{lem}\label{lem;a11.20.10}
Let $S\subset \seisuu^I$ and $\vecn\in\seisuu^I$
such that $\vecn\not\in\nbigs(S)$.
Then there does not exist the following decomposition:
\[
 \prod_{i\in I}z_i^{-n_i}
=\sum_{\vecb\in S}
 f_{\vecb}(z)\cdot \prod_{i\in I}z_i^{-b_i}.
\]
Here $f_{\vecb}$ $(\vecb\in S)$ denote holomorphic functions
on $\Delta^I$.
\hfill\qed
\end{lem}

\subsubsection{The filtrations $\lefttop{I}\nbigv_{S,\vecd}(\gbige)$}

For $\vecd\in\real^{\lbar-I}_{<0}$,
we put as follows:
\[
 \lefttop{I}\nbigv_{S,\vecd}(\gbige)
:=\sum_{\vecb\in S}\lefttop{\lbar}\Vzero_{\vecb+\vecd}(\gbige)
=\sum_{\vecb\in\nbigs(S)}\lefttop{\lbar}\Vzero_{\vecb+\vecd}(\gbige),
\quad\quad
 \lefttop{I}\nbigv'_{S,\vecd}(\gbige)
:=\sum_{\vecb\in\nbigs^0(S)}\lefttop{\lbar}\Vzero_{\vecb+\vecd}(\gbige).
\]
Similarly we put as follows, for any $\vecd\in\real^{\lbar-I}$:
\[
 \lefttop{I}\nbigv_{S,\vecd}\bigl(\naiveprolong{\nbige}\bigr)
:=\sum_{\vecb\in S}\lefttop{I}\Vzero_{\vecb+\vecd}
  \bigl(\naiveprolong{\nbige}\bigr),
\quad
 \lefttop{\lbar}\nbigv'_{S,\vecd}\bigl(\naiveprolong{\nbige}\bigr)
:=\sum_{\vecb\in\nbigs'(S)}
 \lefttop{\lbar}\Vzero_{\vecb+\vecd}\bigl(\naiveprolong{\nbige}\bigr).
\]

\begin{lem}
There exists a finite subset $S'$,
such that
$\lefttop{I}\nbigv'_{S,\vecd}\bigl(\gbige\bigr)=
 \lefttop{I}\nbigv_{S',\vecd}\bigl(\gbige\bigr)$
and
$\lefttop{I}\nbigv'_{S,\vecd}\bigl(\naiveprolong{\nbige}\bigr)
=\lefttop{I}\nbigv_{S',\vecd}\bigl(\naiveprolong{\nbige}\bigr)$.
\end{lem}
\pf
It follows
from the discreteness of the set $\Par(\nbige^{\lambda_0},I)$.
\hfill\qed

\vspace{.1in}

We have the naturally defined morphism:
\[
 \lefttop{I}\nbigv_{S,\vecd}\bigl(\gbige\bigr)
\lrarr
 \lefttop{I}\nbigv_{S,\vecd}\bigl(\naiveprolong{\nbige}\bigr).
\]
If $S$ consists of the unique element $\vecc$,
we have the following by definition:
\[
 \frac{\lefttop{I}\nbigv_{\vecc,\vecd}\bigl(\gbige\bigr)}
  {\lefttop{I}\nbigv'_{\vecc,\vecd}\bigl(\gbige\bigr)}
\simeq
 \lefttop{I}T(\vecc,\vecd),\quad\quad
 \frac{\lefttop{I}\nbigv_{\vecc,\vecd}
 \bigl(\naiveprolong{\nbige}\bigr)}
 {\lefttop{I}\nbigv'_{\vecc,\vecd}
 \bigl(\naiveprolong{\nbige}\bigr)}
\simeq
 \lefttop{I}\Gr^{\Vzero}_{\vecc}\lefttop{J}\Vzero_{\vecd}
 \bigl(
 \naiveprolong{\nbige}\bigr).
\]
Here we put $J:=\lbar-I$.

Let $I$ be a subset of $\lbar$.
We put $J:=\lbar-I$.
Let $\vecd$ be an element of $\real^J$.
Let $\vecv$ be a frame of
$\lefttop{\lbar}\Vzero_{\vecd}\bigl(\naiveprolong{\nbige}\bigr)$,
which is compatible with $\Fzero$ and $\EEzero$.
We put $\lefttop{j}\deg^{\Vzero}(v_i):=\lefttop{j}\deg^{\Fzero}(v_i)-1$,
and
$\vecd(v_i):=
 \bigl(\lefttop{j}\deg^{\Vzero}(v_i)\,\big|\,j\in I
 \bigr)$.
For any element $\vecc\in\real^I$,
we have the element $\vecn(v_i,\vecc)\in\seisuu^I$ determined
by the condition
$\vecc-\vecdelta<\vecn(v_i,\vecc)+\vecd(v_i)\leq \vecc$.
When we put $\tilde{v}_i:=v_i\cdot z^{-\vecn(v_i,\vecc)}$,
the tuple $\tilde{\vecv}:=\bigl(\tilde{v}_i\bigr)$
is a frame of $\lefttop{\lbar}\Vzero_{\vecc+\vecd}\bigl(
 \naiveprolong{\nbige} \bigr)$,
which is compatible with $\Fzero$ and $\EEzero$.
Here $z^{-\vecn(v_i,\vecc)}$ denotes $\prod_{j}z_j^{-n_j(v_i,\vecc)}$.

\begin{lem} \label{lem;9.18.15}
Let $\vecc$ be an element of $\real^I$ and
$S$ be a subset of $\real^I$
such that $\vecc\not\in\nbigs(S)$.
For any section $f\in
\lefttop{I}\nbigv_{\vecc,\vecd}(\naiveprolong{\nbige})
\cap\lefttop{I}\nbigv_{S,\vecd}(\naiveprolong{\nbige})$,
we have $f=0$ in
$\lefttop{I}\Gr^{\Vzero}_{\vecc}\lefttop{J}\Vzero_{\vecd}
 (\naiveprolong{\nbige})$.
\end{lem}
\pf
We put $\tilde{v}_i:=v_i\cdot z^{-\vecn(v_i,\vecc)}$.
Then we have the description
$f=\sum \tilde{f}_i\cdot\tilde{v}_i=
   \sum \tilde{f}_i\cdot z^{-\vecn(v_i,\vecc)}\cdot v_i$,
for holomorphic functions $\tilde{f}_i$ on $\nbigx(\lambda_0,\epsilon_0)$.
Assume that $f\neq 0$ in
$\lefttop{I}\Gr^{\Vzero}_{\vecc}\lefttop{J}V_{\vecd}\bigl(
 \naiveprolong{\nbige}\bigr)$,
and we will derive a contradiction.
Under the assumption, we have the following:
\[
 \Bigl(
 \sum_{\vecd(v_a)+\vecn(v_a,\vecc)=\vecc}
 \tilde{f}_a\cdot\tilde{v}_a
 \Bigr)_{|\,\nbigd_I(\lambda_0,\epsilon_0)}\neq 0.
\]
Then there exists $i_0$ such that
$\vecd(v_{i_0})+\vecn(v_{i_0},\vecc)=\vecc$
and
$\tilde{f}_{i_0\,|\,\nbigd_I(\lambda_0,\epsilon_0)}\neq 0$.

On the other hand,
since we have $f\in\lefttop{I}\nbigv_{S,\vecd}(\naiveprolong{\nbige})$
by our assumption,
we have the description $f=\sum_{\vecb\in S} f_{\vecb}$
for some sections
$f_{\vecb}\in \lefttop{\lbar}V_{\vecb+\vecd}
 \bigl(\naiveprolong{\nbige}\bigr)$.
For each $f_{\vecb}$,
we have the description
$f_{\vecb}=\sum f_{\vecb,i}\cdot z^{-\vecn(v_i,\vecb)}\cdot v_i$
for holomorphic functions $f_{\vecb,i}$
on $\nbigx(\lambda_0,\epsilon_0)$.

As a result, we have the equality:
\begin{equation}\label{eq;a11.20.2}
 \tilde{f}_{i_0}\cdot z^{-\vecn(v_{i_0},\vecc)}
=\sum_{\vecb\in S} f_{\vecb,i_0}\cdot z^{-\vecn(v_{i_0},\vecb)}.
\end{equation}

\begin{lem}\label{lem;a11.20.1}
\[
 \vecn(v_{i_0},\vecc)\not\in
 \bigcup_{\vecc\in S}\nbigs\bigl(\vecn(v_{i_0},\vecb)\bigr).
\]
\end{lem}
\pf
Assume that $\vecn(v_{i_0},\vecc)\in\nbigs\bigl(\vecn(v_{i_0},\vecb)\bigr)$
for some $\vecb\in S$.
Then we have the following:
\[
 \vecc-\vecd(v_{i_0})
=\vecn(v_{i_0},\vecc)
\leq
 \vecn(v_{i_0},\vecb)\leq \vecb-\vecd(v_{i_0}).
\]
Thus we obtain $\vecc\leq \vecb$,
which contradicts our assumption $\vecc\not\in\nbigs(S)$
of Lemma \ref{lem;9.18.15}.
Hence we obtain Lemma \ref{lem;a11.20.1}.
\hfill\qed

\vspace{.1in}

Thus the equality (\ref{eq;a11.20.2}) and Lemma \ref{lem;a11.20.10}
contradicts,
and the proof of Lemma \ref{lem;9.18.15}
is accomplished.
\hfill\qed

\vspace{.1in}

For a primitive subset $S\subset \seisuu^I_{\geq\,0}$,
we have the natural surjection:
\begin{equation}\label{eq;a11.20.3}
 \bigoplus_{\vecb\in S}
 \lefttop{I}\Gr^{\Vzero}_{\vecb}
 \lefttop{J}\Vzero_{\vecb}\bigl(
 \naiveprolong{\nbige}\bigr)
\lrarr
 \frac{\lefttop{I}\nbigv_{S,\vecd}(\naiveprolong{\nbige})}
 {\lefttop{I}\nbigv'_{S,\vecd}(\naiveprolong{\nbige})}.
\end{equation}

\begin{cor}
The morphism {\rm(\ref{eq;a11.20.3})} is isomorphic.
\end{cor}
\pf
We have only to show the injectivity of (\ref{eq;a11.20.3}).
Let us consider sections
$f_{\vecb}\in
  \lefttop{I}\Vzero_{\vecb}
  \lefttop{J}\Vzero_{\vecd}\bigl(\naiveprolong{\nbige}\bigr)$
$(\vecb\in S)$
such that the summation $g=\sum_{\vecb\in S} f_{\vecb}$ is contained in
$\lefttop{I}\nbigv'_{\vecb}\bigl(\naiveprolong{\nbige}\bigr)$.
For any element $\vecb_1\in S$,
we have $f_{\vecb_1}=g-\sum_{\vecb\in S-\{\vecb_1\}}f_{\vecb}$.
Hence
there exists a finite subset $S'$ such that
$ \vecb_1\not\in\nbigs(S')$
and 
$f_{\vecb_1}\in
 \lefttop{I}\nbigv_{\vecb_1,\vecd}\bigl(\naiveprolong{\nbige}\bigr)
 \cap
 \lefttop{I}\nbigv_{S',\vecd}\bigl(\naiveprolong{\nbige}\bigr)$.
Due to Lemma \ref{lem;9.18.15},
we obtain $f_{\vecb_1}=0$ in
$\lefttop{I}\Gr^{\Vzero}_{\vecb_1}
 \lefttop{J}\Vzero_{\vecd}(\naiveprolong{\nbige})$.
It implies the injectivity of the morphism (\ref{eq;a11.20.3}).
\hfill\qed

%% file: a76.tex

\subsubsection{The comparison of
 $\lefttop{\lbar}\nbigv_{S}(\gbige)$
and $\lefttop{\lbar}\nbigv_S\bigl(\naiveprolong{\nbige}\bigr)$}

For a primitive subset $S$,
we have the following naturally defined commutative diagramm:
\begin{equation} \label{eq;9.19.50}
 \begin{CD}
 \bigoplus_{\vecb\in S}
 \lefttop{I}T^{(\lambda_0)}(\vecb,\vecd)
 @>{g'}>{{\rm injective}}>
 \bigoplus_{\vecb\in S}
 \lefttop{I}\Gr^{\Vzero}_{\vecb}
 \lefttop{J}\Vzero_{\vecd}\bigl(\naiveprolong{\nbige}\bigr)\\
 @V{\rm sur}V{f}V @V{\simeq}V{f'}V \\
 {\displaystyle
 \frac{\lefttop{I}\nbigv_{S,\vecd}(\gbige)}
 {\lefttop{I}\nbigv'_{S,\vecd}(\gbige)}  }
 @>{g}>>
 {\displaystyle
 \frac{\lefttop{I}\nbigv_{S,\vecd}\bigl( \naiveprolong{\nbige}\bigr)}
 {\lefttop{I}\nbigv'_{S,\vecd}\bigl( \naiveprolong{\nbige}\bigr)}
 }.
 \end{CD}
\end{equation}

\begin{lem} \label{lem;a11.21.1}
The morphism $g$ in {\rm(\ref{eq;9.19.50})} is injective,
and the morphism $f$ in {\rm (\ref{eq;9.19.50})} is isomorphic.
We have
$\Image(g)=f'\bigl(\bigoplus_{\vecb\in S}\lefttop{I}T_{\vecb}\bigr)$.
\end{lem}
\pf
It immediately follows the diagramm (\ref{eq;9.19.50}).
Note we have
$\Image(g')=\bigoplus_{\vecb\in S}\lefttop{I}T_{\vecb}$.
\hfill\qed

Let $S$ be a primitive subset of $\real^l$.
\begin{prop} \label{prop;b11.21.5}
We have
$\lefttop{\lbar}\nbigv_S\bigl(\gbige\bigr)=
 \lefttop{\lbar}\nbigv_S\bigl(\naiveprolong{\nbige}\bigr)
\cap \gbige$.
\end{prop}
\pf
The implication $\subset$ is clear.
Thus we have only to show the implication $\supset$.
Let $N$ be a large number such that
$-N\cdot\vecdelta<\vecb$ for any element $\vecb\in S$.

Let $f$ be a section of
$\gbige\cap\lefttop{\lbar}\nbigv_S\bigl(\naiveprolong{\nbige}\bigr)$.
We have the following description:
\[
 f=\sum_{\vecb\in T}f_{\vecb}.
\]
Here $T$ denotes some primitive subset of $\real^l$,
and $f_{\vecb}$ are sections of
$\lefttop{\lbar}\Vzero_{\vecb}\bigl(\gbige\bigr)$.
If $T$ is contained in $\nbigs(S)$,
then $f$ is contained in $\lefttop{\lbar}\nbigv_S(\gbige)$,
and thus there are nothing to show.
We will give an algorithm to replace the primitive subset $T$
when $T$ is not contained in $\nbigs(S)$.

In general, we put 
$P(\vecb):=\bigl\{i\in \lbar\,\big|\,b_i>-N\bigr\}$
for any element $\vecb\in\real^l$.
Then we put
$T_j:=\bigl\{\vecb\in T\,\big|\,|P(\vecb)|=j\bigr\}$.
We divide $T_j$ into $T_j^{\star}\sqcup T_j^{\times}$,
where we put as follows:
\[
 T_j^{\star}:=\bigl\{\vecb\in T_j\,\big|\,\vecb\in\nbigs(S)\bigr\},
\quad
 T_j^{\times}:=T_j-T_j^{\star}.
\]
\begin{lem}
We have $T_0=T_0^{\star}$.
\end{lem}
\pf
It follows from our choice of $N$.
\hfill\qed

\vspace{.1in}

We divide $T_j^{\times}$ as follows:
\[
 T_j^{\times}=\coprod_{\substack{I\subset\lbar,\\ |I|=j}}
 T_I^{\times},
\quad\quad
 T_I^{\times}:=\bigl\{\vecb\in T_j^{\times}\,\big|\,
 P(\vecb)=I \bigr\}.
\]
We have the naturally defined morphism
$\pi_I:T_I^{\times}\lrarr \real_{>-N}^I$.
We put $S_I^{\times}:=\pi_I\bigl(T_I^{\times}\bigr)$.
For any element $\vecc\in S_I^{\times}$,
we put $T_I^{\times}(\vecc):=\pi_I^{-1}(\vecc)$.
Then we have the decomposition:
\[
 f=\sum_{j=1}^l\Bigl(
 \sum_{\vecb\in T_j^{\star}}f_{\vecb}
+\sum_{|I|=j}
 \sum_{\vecc\in S_I^{\times}}
 \sum_{\vecb\in T_I^{\times}(\vecc)}
 f_{\vecb}
 \Bigr).
\]

Note that $T\subset\nbigs(S)$ is equivalent to
$\bigcup_m T_m^{\times}=\emptyset$.
Hence we put as follows,
when $T\not\subset\nbigs(S)$:
\[
 m_0(T):=\max\bigl\{m\,\big|\,T_m^{\times}\neq\emptyset\bigr\}.
\]
Let $I$ be a subset of $\lbar$ such that
$|I|=m_0(T)$ and $S_I^{\times}\neq\emptyset$.
For any element $\vecc\in S_I^{\times}$,
we put $|\vecc|:=\sum_{i\in I} (c_i+N)>0$.
We put as follows:
\[
 q_0(T):=\max\Bigl(
 \bigcup_{|I|=m_0(T)}\bigl\{|\vecc|\,\,\big|\,\,\vecc\in S_I^{\times}\bigr\}
 \Bigr).
\]
We also put as follows:
\[
 r_0(T):=\Bigl|
 \bigcup_{|I|=m_0(T)}\bigl\{
 \vecc\in S_I^{\times}\,\big|\,
 |\vecc|=q_0
 \bigr\}
 \Bigr|.
\]

Let take $\vecc\in S_I^{\times}$ such that $|\vecc|=q_0(T)$.
We have the section:
\[
 F_{\vecc}:=
 \sum_{\vecb\in T_I^{\times}(\vecc)} f_{\vecb}
 \in
 \lefttop{I}\Vzero_{\vecc}
 \lefttop{J}\Vzero_{-N\cdot\vecdelta_J}\bigl(\gbige\bigr)
\subset
 \lefttop{I}\Vzero_{\vecc}
 \lefttop{J}\Vzero_{-N\cdot\vecdelta_J}\bigl(
 \naiveprolong{\nbige}
 \bigr).
\]
Here we put $J:=\lbar-I$.

We have the induced section
$[F_{\vecc}]$ of
$\lefttop{I}T(\vecc,-N\cdot\vecdelta_J)
=\lefttop{I}\tilde{T}(\vecc,-N\cdot\vecdelta_J)
 \subset
\lefttop{I}\Gr^{\Vzero}_{\vecc}
 \lefttop{J}\Vzero_{-N\cdot\vecdelta_J}\bigl(\naiveprolong{\nbige}\bigr)$.

\begin{lem}\label{lem;b11.21.2}
Assume $[F_{\vecc}]\neq 0$.
Then we have some primitive subset $U_0\subset\real^J_{\leq -N}$
and the decomposition
$[F_{\vecc}]=\sum_{\vecd\in U_0}G_{\vecc,\vecd}$,
where $G_{\vecc,\vecd}$ are sections of
$\lefttop{I}T(\vecc,\vecd)\subset
 \lefttop{I}\Gr^{\Vzero}_{\vecc}\lefttop{J}\Vzero_{\vecd}
 \bigl(\naiveprolong{\nbige}\bigr)$,
such that
the induced section $[G_{\vecc,\vecd}]$ is not $0$
in $\lefttop{I}\Gr^{\Vzero}_{\vecc}
 \lefttop{I}\Gr^{\Vzero}_{\vecd}
 \bigl(\naiveprolong{\nbige}\bigr)
=\lefttop{\lbar}\Gr^{\Vzero}_{\vecc+\vecd}\bigl(
 \naiveprolong{\nbige}
 \bigr)$.
(Note Lemma {\rm\ref{lem;a12.2.10}}.)
\end{lem}
\pf
The filtrations
$\lefttop{j}\Vzero\bigl(\naiveprolong{\nbige}\bigr)$
$(j\in J)$
induces the filtrations $\lefttop{j}\Vzero$ of
$\lefttop{I}\Gr^{\Vzero}_{\vecc}\Vzero_{-N\cdot\vecdelta}\bigl(
 \naiveprolong{\nbige}\bigr)$.
Since the induced morphisms $\deldel_i$ $(i\in I)$
are strict with respect to the filtrations
$\lefttop{j}\Vzero$ $(j\in J)$
(see Lemma \ref{lem;a12.2.11}),
we have the following equalities
for any $\vecd\leq -N\cdot\vecdelta$:
\[
 \lefttop{J}\Vzero_{\vecd}\cap
 \lefttop{I}\tilde{T}(\vecc,-N\cdot\vecdelta_J)
=\lefttop{I}\tilde{T}(\vecc,\vecd).
\]
Then the claim follows immediately.
\hfill\qed

\begin{cor}
Assume $[F_{\vecc}]\neq 0$ in $\lefttop{I}T(\vecc,-N\cdot\vecdelta)$.
Then we have some primitive subset $U_0\subset\real^J_{\leq -N}$
and the decomposition $F_{\vecc}=\sum_{\vecd\in U_0} F_{\vecc,\vecd}$,
where $F_{\vecc,\vecd}$ is a section of
$\lefttop{I}\Vzero_{\vecc}\lefttop{J}\Vzero_{\vecd}
\bigl(\gbige\bigr)$ such that
the induced section $[F_{\vecc,\vecd}]$ is not $0$
in $\lefttop{\lbar}\Gr^{\Vzero}_{\vecc+\vecd}(\gbige)$.
\end{cor}
\pf
Due to Lemma \ref{lem;b11.21.2},
we have the decomposition
$F_{\vecc}=\sum_{\vecd\in U_0}F_{\vecc,\vecd}$
such that $[F_{\vecc,\vecd}]\neq 0$
in $\lefttop{\lbar}\Gr^{\Vzero}_{\vecc+\vecd}
   \bigl(\naiveprolong{\nbige}\bigr)$.
Then we obtain the claim due to Lemma \ref{lem;9.18.11}.
\hfill\qed

\vspace{.1in}

Then we obtain the decomposition:
\[
 f=\sum_{\vecb\in T-T_I^{\times}(\vecc)} f_{\vecb}
 +\sum_{\vecd\in U_0}F_{\vecc,\vecd}.
\]
Recall we have $|I|=m_0(T)$,
which we denote by $m_0$ for simplicity.
\begin{lem}
For any $\vecd\in U_0$,
there exists an element $\vecb\in \bigcup_{m_0\leq m}T_m^{\star}$,
such that $\vecc+\vecd\leq \vecb$.
\end{lem}
\pf
Let us consider the following two cases:
\begin{description}
\item[Case 1.]
 There exists an element $\vecb\in T-T_I^{\times}(\vecc)$
 such that $\vecc+\vecd\leq\vecb$.
\item[Case 2.]
 There does not exist such element.
\end{description}
In the case 1,
we have $\vecc\leq q_I(\vecb)$.
Hence we have $\vecb\in \bigcup_{m_0<m}T_m$
or $\vecb\in T_{m_0}$.
Due to our choice of $m_0$,
we have $\bigcup_{m_0<m}T_m=\bigcup_{m_0<m}T_m^{\star}$.
Due to our choice of $m_0$ and $\vecc$,
we have the following:
\[
 \bigl\{
 \vecb\in T_{m_0}-T_{I}^{\times}(\vecc)\,\big|\,
 q_I(\vecb)\geq \vecc
 \bigr\}
\subset T_{m_0}^{\star}.
\]
Thus we obtain $\vecb\in \bigcup_{m_0\leq m}T_m^{\star}$
in this case.

Let us consider the case 2.
Then we obtain $F_{\vecc,\vecd}=0$ in
$\lefttop{\lbar}\Gr^{\Vzero}_{\vecc+\vecd}
 \bigl(\naiveprolong{\nbige}\bigr)$
due to Lemma \ref{lem;9.18.15}.
It contradicts our choice of $F_{\vecc,\vecd}$.
Hence we can conclude that the case 2 does not happen.
\hfill\qed

\vspace{.1in}

\begin{cor}
In the case $[F_{\vecc}]$ is not $0$
in $\lefttop{I}T(\vecc,\vecd)$,
we have some sections
$H_{\vecb}$ $(\vecb\in \bigcup T^{\star}_m)$
satisfying the following:
\begin{itemize}
\item
 $H_{\vecb}\in
 \lefttop{I}\Vzero_{\vecc}\lefttop{J}\Vzero_{-N\cdot\vecdelta_J}
 \bigl(\gbige\bigr)
 \cap \lefttop{\lbar}\Vzero_{\vecb}(\gbige)$.
\item
 $[F_{\vecc}-\sum_{\vecb} H_{\vecb}]=0$
 in $\lefttop{I}T(\vecc,-N\cdot\vecdelta_J)$.
\hfill\qed
\end{itemize}
\end{cor}

We put $\tilde{F}_{\vecc}:=F_{\vecc}-\sum_{\vecb} H_{\vecb}$.
For any element $\vecb\in T^{\star}:=\bigcup T^{\star}_m$,
we put $\tilde{f}_{\vecb}:=f+H_{\vecb}$.
For any element
$\vecb\in T-\bigl(T^{\star}\cup T_I^{\times}(\vecc)\bigr)$,
we put $\tilde{f}_{\vecb}:=f_{\vecb}$.
Then we have the following decomposition:
\[
 f=\sum_{T-T_I^{\times}(\vecc)}\tilde{f}_{\vecb}
 +\tilde{F}_{\vecc}.
\]
Since the induced section $[\tilde{F}_{\vecc}]$
is $0$ in $\lefttop{I}T(\vecc,-N\cdot\vecdelta_J)$,
there is some primitive set $U_1\subset\real^l$
and the decomposition:
$\tilde{F}_{\vecc}
=\sum_{\vecb'\in U_1} \overline{F}_{\vecb'}$
such that the following conditions hold
for any element $\vecb'\in U_1$:
\[
 q_J(\vecb')\leq -N\cdot\vecdelta_J,
\quad
 \sum_{i\in I}(b_i+N)<r_0.
\]
We obtain the following decomposition:
\[
 f=\sum_{T-T_I^{\times}(\vecc)}\tilde{f}_{\vecb}
 +\sum_{\vecb'\in U_1}\overline{F}_{\vecb'}.
\]
Then it is easy to obtain 
a decomposition $\sum_{\vecb\in T'}f'_{\vecb}$
for a primitive subset $T'\subset\real^l$
satisfying the following in the lexicographic order:
\[
 \bigl(m_0(T'),r_0(T'),q_0(T')\bigr)
 <\bigl(m_0(T),r_0(T),q_0(T)\bigr).
\]
It is easy to see that the reduction procedure above
can stop in finite times.
Thus the proof of Proposition \ref{prop;b11.21.5}
is accomplished.
\hfill\qed

\begin{cor}
We have
$\lefttop{\lbar}V_{\vecb}\gbige
=\gbige\cap\lefttop{\lbar}V_{\vecb}\bigl(\naiveprolong{\nbige}\bigr)$.
\hfill\qed
\end{cor}

%% file: a76.2.tex

\subsubsection{The distributivity of the filtrations
 $\bigl(
 \lefttop{i}\Vzero\bigl(\gbige\bigr)\,\big|\,i\in\lbar
 \bigr)$}

\begin{cor}
We have
$\lefttop{I}\Vzero_{\vecb}(\gbige)
=\lefttop{I}\Vzero_{\vecb}\bigl(\naiveprolong{\nbige}\bigr)\cap\gbige$.
In particular,
we have
$\lefttop{i}\Vzero_{b}(\gbige)
=\lefttop{i}\Vzero_b\bigl(
 \naiveprolong{\nbige}\bigr)\cap\gbige$.
\end{cor}
\pf
For any section
$f\in\lefttop{I}\Vzero_{\vecb}
  \bigl(\naiveprolong{\nbige}\bigr)\cap\gbige$,
there exists an element $\vecc$ such that
$q_I(\vecc)=\vecb$ and
$f\in\lefttop{\lbar}\Vzero_{\vecc}\bigl(\naiveprolong{\nbige}\bigr)
 \cap\gbige$.
Hence
$f\in \lefttop{\lbar}\Vzero_{\vecc}(\gbige)$,
and thus
$f\in \lefttop{I}\Vzero_{\vecb}(\gbige)$.
Thus we obtain the implication $\supset$.
The inverse implication $\subset$ can be shown similarly.
\hfill\qed

\begin{cor} \label{cor;9.18.31}
We have
$ \lefttop{I}\Vzero_{\vecb}(\gbige)
=\bigcap_{i\in I}\lefttop{i}\Vzero_{\vecb}(\gbige)$.
\end{cor}
\pf
If $f\in\bigcap_{i\in I}\lefttop{i}\Vzero_{b_i}(\gbige)$,
then we obtain
$f\in \bigcap_{i\in I}\lefttop{i}\Vzero_{b_i}\bigl(\naiveprolong{\nbige}\bigr)
 \cap\gbige$.
It implies
$f\in \lefttop{I}\Vzero_{\vecb}\bigl(
 \naiveprolong{\nbige}\bigr)
 \cap\gbige$.
Thus we obtain $f\in\lefttop{I}\Vzero_{\vecb}(\gbige)$.
It implies the implication $\supset$.
The reverse implication is clear.
\hfill\qed

\vspace{.1in}
The following lemma can be shown similarly.
\begin{lem}
Let $I$ be a subset of $\lbar$ such that $i\not\in I$.
Let $S$ be a primitive subset of $\seisuu^I_{\geq 0}$.
Then we have the following:
\[
 \lefttop{i}\Vzero_{b}(\gbige)
 \cap
 \Bigl(\sum_{\vecc\in S}\lefttop{I}\Vzero_{\vecc}(\gbige)\Bigr)
=\sum_{\vecc\in S}\lefttop{I\sqcup\{i\}}\Vzero_{(\vecc,b)}(\gbige).
\]
\end{lem}
\pf
We have only to use Proposition \ref{prop;b11.21.5}.
\hfill\qed

\begin{cor}
We have the natural isomorphism:
\begin{equation} \label{eq;9.18.20}
 \lefttop{i}\Gr^{\Vzero}_b
 \lefttop{I}\Gr^{\Vzero}_{\vecc}(\gbige)
\simeq
 \lefttop{I\sqcup\{i\}}
 \Gr^{\Vzero}_{(\vecc,b)}(\gbige).
\end{equation}
\end{cor}
\pf
The left hand side of (\ref{eq;9.18.20}) is 
isomorphic to the following, by definition:
\[
 {\rm L.H.S.}
=\frac{\lefttop{i}\Vzero_b\lefttop{I}\Vzero_{\vecc}(\gbige)}
 {\lefttop{i}\Vzero_{<b}\lefttop{I}V_{\vecc}(\gbige)
 +
\lefttop{i}\Vzero_b\cap
 \sum_{\vecc'\lneq \vecc}
 \lefttop{I}\Vzero_{\vecc'}(\gbige)
 }
\simeq
 \frac{
 \lefttop{i}\Vzero_{b}\lefttop{I}\Vzero_{\vecc}(\gbige) }
 {\lefttop{i}\Vzero_{<b}\lefttop{I}\Vzero_{\vecc}(\gbige)
+\sum_{\vecc'\lneq\vecc}
 \lefttop{I}\Vzero_{\vecc'}\lefttop{i}\Vzero_b(\gbige )}.
\]
The last term is the right hand side of (\ref{eq;9.18.20}).
Thus we are done.
\hfill\qed

\begin{cor}\label{cor;04.2.24.1}
The $\nbigr$-module $\gbige$ is coherent and holonomic.
\end{cor}
\pf
Let us consider the filtrations $\lefttop{i}\tilde{V}^{(\lambda_0)}$ 
of $\gbige$ indexed by $\real_{\geq\,-1}$,
given as follows:
\[
 \lefttop{i}\tilde{V}^{(\lambda_0)}_a
:=
\left\{
 \begin{array}{ll}
 \lefttop{i}\Vzero_a & (a\geq -1)\\
 \mbox{{}}\\
 0 & (a<-1).
 \end{array}
\right.
\]
We also consider the filtrations $\lefttop{i}F$
of $\nbigr$ indexed by $\seisuu_{\geq\,0}$,
where $\lefttop{i}F_b$ is generated by
the sections of $\nbigr$ whose degree with respect to $\deldel_i$
are less than $b$.
Clearly we have the relation:
\[
 \lefttop{i}F_b(\nbigr)\cdot
 \lefttop{i}\tilde{V}^{(\lambda_0)}_a
\subset
 \lefttop{i}\tilde{V}^{(\lambda_0)}_{a+b}.
\]
We obtain the $\lefttop{\lbar}\Gr^{F}(\nbigr)$-module
$\lefttop{\lbar}\Gr^{\tilde{V}^{(\lambda_0)}}(\gbige)$.
We have the following:
\[
 \lefttop{\lbar}\Gr^{\tilde{V}^{(\lambda_0)}}(\gbige)
=\bigoplus_{I\sqcup J=\lbar}
 \bigoplus_{\veca\in\real_{\geq\,-1}^{I}}
 \lefttop{I}\Gr^{\Vzero}_{\veca}
 \lefttop{J}\Vzero_{-\vecdelta_J}(\gbige).
\]
Hence it is coherent as a $\lefttop{\lbar}\Gr^{F}(\nbigr)$-module.
By using a standard theory of filtered ring sheaves 
(see the appendices in \cite{bjork} or \cite{kashiwara_text}),
we can show that $\gbige$ is a coherent $\nbigr$-module.
Similarly, we can show that
the characteristic variety of $\gbige$
is contained in
$(N_X^{\ast}X\bigcup_I N^{\ast}_{D_I}X)\times\cnum$,
where $N_Y^{\ast}X$ denotes the cotangent bundle of 
$Y$ in $X$.
Then we obtain that $\gbige$ is holonomic.
\hfill\qed

\vspace{.1in}

Let $\nbigr_i$ denote the subsheaves of $\nbigr$
generated by $\deldel_j$ $(j\neq i)$ as algebra.
We have the natural action of
$\nbigr_i$ on $\lefttop{i}\Vzero_a\gbige$.

\begin{cor}
$\gbige$ is regular along $z_i=0$,
namely,
$\lefttop{i}\Vzero_a(\gbige)$ is
coherent as an $\nbigr_i$-module.
\end{cor}
\pf
It can be shown similarly to
Corollary \ref{cor;04.2.24.1}.
\hfill\qed

\begin{lem} \label{lem;a12.2.15}
Let $I\sqcup J=\lbar$ be a decomposition.
For any $\vecc\in\real^I$ and $\vecd\in\real_{<0}^J$,
the $\nbigr$-module
$\lefttop{I}\Gr^{\Vzero}_{\vecc}\lefttop{J}\Vzero_{\vecd}(\gbige)$
is strict.
\end{lem}
\pf
We have the injection
$\lefttop{I}\Gr^{\Vzero}_{\vecc}
 \lefttop{J}\Vzero_{\vecd}(\gbige)
\lrarr
 \lefttop{I}\Gr^{\Vzero}_{\vecc}
 \lefttop{J}\Vzero_{\vecd}\bigl(
 \naiveprolong{\nbige}\bigr)$.
Then the claim immediately follows.
\hfill\qed

\begin{lem} \label{lem;9.18.30}
Let $I\sqcup J=\lbar$ be a decomposition.
For any $\vecc\in\real^I$ and $\vecd\in\real^J$,
the $\nbigr$-module
$\lefttop{I}\Gr^{\Vzero}_{\vecc}\lefttop{J}\Vzero_{\vecd}(\gbige)$
is strict.
\end{lem}
\pf
In the case $I=\lbar$,
the claim follows from Lemma \ref{lem;a12.2.15}.
We use a descending induction on $|I|$.
We assume that the claim holds for any $I$ such that $|I|=m$,
and then we will show the claim for $I$ such that $|I|=m-1$.
We may assume that $I=\underline{m-1}$.

We put $M_-(\vecd):=\{i\in J\,\,|\,d_i<0\}$.
We use a descending induction on  $|M_-(\vecd)|$.
In the case $|M_-(\vecd)|=|J|$,
$\lefttop{\mminusitibar}\Gr^{\Vzero}_{\vecc}
 \lefttop{\mminusitibar^c}\Vzero_{\vecd}(\gbige)$
is strict due to Lemma \ref{lem;a12.2.15}.
Assume we have proved the strictness holds
in the case $|M_-(\vecd)|>m$,
and we will prove that the strictness holds
in the case $|M_-(\vecd)|=m$.
Let pick an element $i\in J$ such that $d_i\geq 0$.
We put $J':=J-\{i\}$ and $I':=I\sqcup\{i\}$.
Let $\pi:\real^J\lrarr\real^{J'}$ denote the projection.
By the hypothesis of the induction on $|I|$,
$\lefttop{I'}\Gr^{\Vzero}_{(\vecc,d)}
 \lefttop{J'}\Vzero_{\pi(\vecd)}(\gbige)$
is strict for any $d\in\real$.
By the hypothesis of the induction on $|M_-(\vecd)|$,
$\lefttop{I}\Gr^{\Vzero}_{(\vecc,b)}
 \lefttop{J}\Vzero_{\vecd-N\vecdelta_i}(\gbige)$ is strict
if $N$ is sufficiently large.
Thus we are done.
\hfill\qed

\begin{cor}
$\lefttop{I}\Gr^{\Vzero}_{\vecc}(\gbige)$ is strict.
\end{cor}
\pf
Since
$\lefttop{I}\Gr^{\Vzero}_{\vecc}(\gbige)$
is an inductive limit
of $\lefttop{I}\Gr^{\Vzero}_{\vecc}
 \lefttop{J}V_{\vecd}(\gbige)$ ($\vecd\in\real^{\lbar-I}$),
the corollary immediately follows from Lemma \ref{lem;9.18.30}.
\hfill\qed

\begin{cor}\label{cor;b12.6.1}
$\lefttop{I}\Gr^{\Vzero}_{\vecc}\gbige$ is
strictly specializable
along $\nbigd_i$ at $\lambda_0$.
\hfill\qed
\end{cor}

\begin{cor}
The $\nbigr$-module $\gbige$ is strictly $S$-decomposable
along $\nbigd_i$,
and we have $\psi_{z_i,0}(\gbige)=\Image(\can)$.
\end{cor}
\pf
We have already obtained strict specializability along $\lambda_0$
(Corollary \ref{cor;b12.6.1}).
We have the following commutative diagramm:
\[
\begin{CD}
 \psi_{z_i,0}(\gbige)
 @>>>
 \lefttop{i}\Gr^{\Vzero}_0(\gbige)
 @>>>
 \lefttop{i}\Gr^{\Vzero}_0\bigl(\naiveprolong{\nbige}\bigr)\\
 @V{\var}VV @VVV @VVV \\
 \psi_{z_i,-\vecdelta_0}(\gbige)
 @>>>
 \lefttop{i}\Gr^{\Vzero}_{-1}(\gbige)
 @>>>
 \lefttop{i}\Gr^{\Vzero}_{-1}\bigl(\naiveprolong{\nbige}\bigr).
\end{CD}
\]
The horizontal arrows are injective,
and the right vertical arrows are isomorphic.
Hence we have $\Ker(\var)=0$.

Let $I$ be a subset of $\lbar$ such that $i\in I$.
We put $J:=\lbar-I$.
Let $\vecc$ be any element of $\real^I$
such that $q_i(\vecc)=0$,
and $\vecd$ be any element of $\real_{<0}^J$.
Then the natural morphism
$\lefttop{I}T(\vecc-\vecdelta_i,\vecd)\lrarr
 \lefttop{I}T(\vecc,\vecd)$ is surjective.
Then we can show
that
$\lefttop{i}\Gr^{\Vzero}_{-1}(\gbige)
 \lrarr
 \lefttop{i}\Gr^{\Vzero}_0(\gbige)$ is surjective,
by an argument similar to the proof of
Lemma \ref{lem;9.18.30}.
Thus we obtain
$\psi_{z_i,0}(\gbige)=\Image(\can)$.
\hfill\qed

\vspace{.1in}

Recall that we have $\lefttop{\nbar}\tildepsi_{\vecu}(\gbige)$
(the subsubsection \ref{subsubsection;b11.23.20}).
\begin{lem}\label{lem;b11.23.10}
We have the isomorphism:
\[
 \lefttop{\lbar}\tildepsi_{\vecu}(\gbige)
\simeq
 \lefttop{1}\tildepsi_{u_1}\Bigl(
 \lefttop{2}\tildepsi_{u_2}\Bigl(
 \cdots
 \lefttop{l}\tildepsi_{u_l}\bigl(
 \gbige
 \bigr)
 \Bigr)
 \Bigr).
\]
\end{lem}
\pf
We have only to show that there exists
the canonical isomorphism as follows,
in the case $\paramap(\lambda_0,\vecu)<0$:
\begin{equation} \label{eq;b11.23.100}
 \lefttop{\lbar}\psizero_{\vecu}(\gbige)
\simeq
 \lefttop{1}\psizero_{u_1}\Bigl(
 \lefttop{2}\psizero_{u_2}\Bigl(
 \cdots
 \lefttop{l}\psizero_{u_l}\bigl(
 \gbige
 \bigr)
 \Bigr)
 \Bigr).
\end{equation}
Both of (\ref{eq;b11.23.100}) are naturally isomorphic
to $\lefttop{\lbar}\nbigg_{\vecu+\vecdelta_{0,\lbar}}(E)$.
\hfill\qed

%% file: 31.5.tex

\subsubsection{Primitive Decomposition of sections of
 $\naiveprolong{\nbige}$ and $\gbige$}
\label{subsubsection;9.21.5}

The following lemma can be shown elementarily.
\begin{lem}\label{lem;a11.22.1}
Let $f$ be a holomorphic function on $\Delta_z^n\times\Delta_w^m$.
There exists the unique primitive subset $S\subset\seisuu^l_{\geq\,0}$
such that we have a holomorphic decomposition
$f=\sum_{\vecp\in S}z^{\vecp}\cdot a_{\vecp}(z,w)$
such that $a_{\vecp}(0,w)\neq 0$.
\hfill\qed
\end{lem}

First let us consider the primitive sections of
$\naiveprolong{\nbige}$.

\begin{df}
Let $I$ be a subset of $\lbar$.
A section $f\in\naiveprolong{\nbige}$ is called $I$-primitive,
if $f\in\lefttop{I}\Vzero_{\vecb}\bigl(\naiveprolong{\nbige}\bigr)$
and $f\neq 0$ in
$\lefttop{I}\Gr^{\Vzero}_{\vecb}\bigl(\naiveprolong{\nbige}\bigr)$.
\mbox{{}}\hfill\qed
\end{df}

The following lemmas are easy to see.
\begin{lem}
Let $\vecv=(v_i)$ be a frame of $\naiveprolongg{\vecb}{\nbige}$,
which is compatible with $\EEzero$ and $\Fzero$.
For any $v_i$, for any element $\vecc\in\seisuu^l$
and for any subset $I\subset\lbar$,
the section $v_i\cdot \prod_{j=1}^l z_j^{-c_j}\cdot a$, $(a(0)\neq 0)$,
is $I$-primitive.
\hfill\qed
\end{lem}

\begin{lem} \label{lem;9.18.40}
Assume that $f$ is $I$-primitive
and $\lefttop{I}\deg^{V}(f)=\vecc$.
For any section
$g\in\lefttop{I}\nbigv'_{\vecc}\bigl(\naiveprolong{\nbige}\bigr)$,
then $f+g$ is $I$-primitive
such that $\lefttop{I}\deg^{V}(f+g)=\vecc$.
\hfill\qed
\end{lem}

\begin{lem}
Let $f$ be a section.
There exists a finite subset $S\subset\real^I$
and a decomposition:
$f=\sum_{\vecb\in S}f_{\vecb}$.
Here $f_{\vecb}$ is $I$-primitive
such that $\lefttop{I}\deg(f_{\vecb})=\vecb$.

Let $\max(S)$ denote the set of the maximal elements of $S$.
Then there exist the $I$-primitive
sections $g_{\vecb}$ for $\vecb\in S$
such that $\lefttop{I}\deg(g_{\vecb})=\vecb$
and $f=\sum_{\vecb\in \max(S)}g_{\vecb}$.
\hfill\qed
\end{lem}

\begin{df}
An $I$-primitive decomposition of $f$ is a decomposition
$f=\sum_{\vecb\in S}f_{\vecb}$
such that the following holds:
\begin{itemize}
\item
 The subset $S\subset\real^I$ is primitive.
\item
 The sections $f_{\vecb}$ are $I$-primitive
 such that $\lefttop{I}\deg(f_{\vecb})=\vecb$.
\hfill\qed
\end{itemize}
\end{df}

\begin{lem}
For any section $f$ of $\naiveprolong{\nbige}$,
there exists an $I$-primitive decomposition.
\end{lem}
\pf
We have a development $f=\sum f_i\cdot v_i$,
and we have a primitive decomposition of $f_i$,
as in Lemma \ref{lem;a11.22.1}.
Then we obtain the decomposition
$f=\sum_{\vecb\in S'}f_{\vecb}$
such that $f_{\vecb}$ is $I$-primitive
with $\lefttop{I}\deg(f_{\vecb})=\vecb$.
Then we have only to apply Lemma \ref{lem;9.18.40}.
\hfill\qed

\begin{lem}\label{lem;9.21.16}
Let $f=\sum_{\vecb\in S}f_{\vecb}$ be an
$I$-primitive decomposition.
Then the set $S$ is canonically determined.
For any $\vecb\in S$,
the section
$[f_{\vecb}]$ of
$\lefttop{I}\Gr^V_{\vecb}\bigl(\naiveprolong{\nbige}\bigr)$
is canonically determined.
\end{lem}
\pf
Assume we have two primitive decomposition of $f$:
\[
 f=\sum_{\vecb\in S}f_{\vecb}=\sum_{\vecb'\in S'}f'_{\vecb'}.
\]
Assume $\vecb\in S$ and $\vecb\not\in S'$,
and we will derive a contradiction.

Let us consider $A=\bigl\{\vecb'\in S'\,\big|\,\vecb\leq \vecb'\bigr\}$.
Assume $A\neq\emptyset$.
Note that $\vecb'\in A$ is not contained in $S$.
From the equality above,
there exists $S''$ such that
$\vecb\not\in \nbigs(S'')$ and
$f_{\vecb}\in \lefttop{I}\nbigv_{S''}$.
It implies
$[f_{\vecb'}]=0$
in $\lefttop{I}\Gr^V_{\vecb'}\bigl(\naiveprolong{\nbige}\bigr)$,
and thus we have arrived at the contradiction.

If $A=\emptyset$,
by a similar argument,
we obtain $[f_{\vecb}]=0$
in $\lefttop{I}\Gr^V_{\vecb}\bigl(\naiveprolong{\nbige}\bigr)$,
thus we have arrived at the contradiction.
Thus we obtain $\vecb\in S'$.
By symmetry, we obtain $S=S'$.

We have the following:
\[
 f_{\vecb}-f_{\vecb}'
=\sum_{\vecb'\in S-\{\vecb\}}
 (f_{\vecb'}-f'_{\vecb'}).
\]
Thus we obtain $[f_{\vecb}-f_{\vecb}']=0$
in $\lefttop{I}\Gr^{V}_{\vecb}(\naiveprolong{\nbige})$.
\hfill\qed

\begin{df}
The set $S$ above is denoted by $\lefttop{I}\Prim(f)$.
For any element $\vecb\in\lefttop{I}\Prim(f)$,
we put as follows:
\[
 \lefttop{I}P_{\vecb}(f)
:=[f_{\vecb}]
\in\lefttop{I}\Gr^{V}_{\vecb}\bigl(\naiveprolong{\nbige}\bigr).
\]
\hfill\qed
\end{df}

The following lemma is easy to see.
\begin{lem} \label{lem;9.18.32}
$f\in \lefttop{I}\nbigv_{S}\bigl(\naiveprolong{\nbige}\bigr)$
if and only if
$\lefttop{I}\Prim(f)\subset \nbigs(S)$.
\hfill\qed
\end{lem}

We have the natural inclusion
$\iota:\gbige\lrarr\naiveprolong{\nbige}$.
For any section $f\in \gbige$,
we put $\lefttop{I}\Prim(f):=\lefttop{I}\Prim(\iota(f))$,
and we put as follows:
\[
 \lefttop{I}\Prim(f)
:=\lefttop{I}\Prim\bigl(\iota(f)\bigr)
 \in \lefttop{I}\Gr^V_{\vecb}(\gbige)
 \subset \lefttop{I}\Gr^V_{\vecb}(\naiveprolong{\nbige}).
\]

\begin{lem}
$f\in \lefttop{I}\nbigv_S(\gbige)$
if and only if
$\lefttop{I}\Prim(f)\subset\nbigs(S)$.
\end{lem}
\pf
It follows from Corollary \ref{cor;9.18.31} and
Lemma \ref{lem;9.18.32}.
\hfill\qed

%% file: a79.tex


\begin{prop}\label{prop;a11.26.2}
Let $\tilde{\gbige}$ be a coherent $\nbigr$-module
on $\nbigx$ satisfying the following conditions:
\begin{enumerate}
\item
 The restriction $\tilde{\gbige}_{|\nbigx-\nbigd}$ is isomorphic to 
 $\nbige$.
\item\label{number;d11.26.2}
 We have the injection $\tilde{\gbige}\lrarr \iota_{\ast}\nbige$,
 where $\iota$ denotes the open immersion $X-D\lrarr X$.
\item\label{number;d11.26.3}
 The $\nbigr$-module $\tilde{\gbige}$ is holonomic.
 It is also regular and strictly $S$-decomposable
 along $z_i=0$ $(i=1,\ldots,n)$.
\item\label{number;d11.26.4}
 Let $\lambda\in\cnum^{\ast}$ be generic
 with respect to $\bigcup_i\KMS(\tilde{\gbige},i)$.
 Then the specialization $\tilde{\gbige}_{|\nbigxlambda}$
 is strictly $S$-decomposable along $z_i=0$ $(i=1,\ldots,n)$.
 We also have
 $\psi_{z_i,0}\tilde{\gbige}_{|\nbigxlambda}=\Image (\can)$.
\item \label{number;d11.26.5}
 Let $\lefttop{i}\Vzero$ be the $V$-filtration along $z_i=0$ at
     $\lambda$ such that
 $\lefttop{i}\Gr^{\Vzero}\bigl(\tilde{\gbige}\bigr)$ is strict.
 We put $\lefttop{\nbar}\Vzero_{<0}
 :=\bigcap_{i=1}^n\lefttop{i}\Vzero_{<0}$.
Then $\lefttop{\nbar}\Vzero_{<0}$ is a coherent locally free
 $\nbigo_{|\nbigx}$-module,
and $\tilde{\gbige}$ is generated by $\lefttop{\nbar}\Vzero_{<0}$.
\end{enumerate}
Then we have the natural isomorphism
$\tilde{\gbige}\simeq \gbige$.
\end{prop}
\pf
First we remark that $\gbige$ satisfies the conditions above.
We have the isomorphism
$\tilde{\gbige}_{|\nbigx-\nbigd}\simeq \nbige\simeq
 \gbige_{|\nbigx-\nbigd}$,
and we can regard $\tilde{\gbige}$ and $\gbige$
as the $\nbigr$-submodules of $\iota_{\ast}\nbige$,
where $\iota$ denotes the open immersion $X-D\lrarr X$.
We will show that they are same as the $\nbigr$-submodules.

Let $\lambda\in\cnum^{\ast}$ be generic with respect to
$\KMS(\tilde{\gbige},z_i)\cup\KMS(\gbige,z_i)$.
Let $\pi_i:X\lrarr D_i$ denote the naturally defined projection.
We put 
$X_i:=\pi_i^{-1}\left(D_i-\bigcup_{j\neq i}D_i\cap D_j\right)$
and $\nbigx_i^{\lambda}:=X_i\times\{\lambda\}$.
Let us consider the restriction
$\tilde{\gbige}_{|\nbigx_i^{\lambda}}$
and
$\gbige_{|\nbigx_i^{\lambda}}$,
which we denote by
$\tilde{\gbige}^{\lambda}_i$
and $\gbige^{\lambda}_i$ respectively.

\begin{lem}
The $D$-modules
$\tilde{\gbige}^{\lambda}_{i}$
and $\gbige^{\lambda}_i$ are regular holonomic.
They are prolongation of $\nbigelambda_{|X_i\setminus D_i}$.
\end{lem}
\pf
They are regular along $z_i=0$.
The restrictions to $X_i\setminus D_i$ are isomorphic to
$\nbigelambda_{|X_i\setminus D_i}$.
Then the claims immediately follow.
\hfill\qed

\vspace{.1in}

Let $L$ be the local system corresponding for
$\nbigelambda_{|X_i\setminus D_i}$.
Let $\tilde{\nbigf}$ and $\nbigf$ be the perverse sheaves
corresponding to 
$\tilde{\gbige}^{\lambda}_i$ and $\gbige^{\lambda}_i$.
Since $\gbige^{(\lambda)}_i$ and $\tilde{\gbige}^{\lambda}_i$
are strictly $S$-decomposable,
and since $\psi_{z_i,0}=\Image\can$
for both of them,
$\tilde{\nbigf}$ and $\nbigf$ are the intermediate extensions
of $L$.
Hence we have the isomorphism
$\tilde{\nbigf}\simeq\nbigf$,
extending the canonical isomorphism
$\tilde{\nbigf}_{|X_i\setminus D_i}
 \simeq L
 \simeq \nbigf_{|X_i\setminus D_i}$.
Due to Riemann-Hilbert correspondence,
we obtain the isomorphism
$\tilde{\gbige}^{\lambda}_{i}\simeq
 \gbige^{\lambda}_i$,
extending
$\tilde{\gbige}^{\lambda}_{|X_i\setminus D_i}\simeq 
 \nbige\simeq \gbige^{\lambda}_{|X_i\setminus D_i}$.

\begin{lem}\label{lem;c12.4.1}
The KMS-spectrum $\KMS(\gbige,z_i)$ and
$\KMS(\tilde{\gbige},z_i)$ are same.
\end{lem}
\pf
Let $\lambda$ be any generic point.
The $V$-filtration $\lefttop{i}V^{(\lambda)}$ of
$\tilde{\gbige}$ and $\gbige$ along $z_i=0$
induce the Kashiwara-Malgrange filtrations $V$
of $\tilde{\gbige}^{\lambda}_i$ and $\gbige^{\lambda}_i$
along $z_i=0$.
Let us consider the following sets:
\[
 \Sp\bigl(\tilde{\gbige}^{\lambda}_i,z_i\bigr)
:=\bigl\{
 \alpha\in\cnum\,\big|\,
 \Gr^V\bigl(\tilde{\gbige}^{\lambda}_i\bigr)\neq 0
  \bigr\},
\quad\quad
 \Sp\bigl(\gbige^{\lambda}_i\bigr)
:=\bigl\{
 \alpha\in\cnum\,\big|\,
 \Gr^{V}\bigl(\gbige^{\lambda}_i\bigr)\neq 0
 \bigr\}.
\]
Then we have
$\Sp\bigl(\tilde{\gbige}^{\lambda}_i,z_i\bigr)
=\Sp\bigl(\gbige^{\lambda}_i,z_i\bigr)$
for any generic $\lambda$,
due to the uniqueness property of
the graduation of Kashiwara-Malgrange filtration.
Since the set of generic $\lambda$ is uncountable,
we obtain the coincidence $\KMS(\tilde{\gbige},z_i)=\KMS(\gbige,z_i)$.
\hfill\qed

\begin{lem} \label{lem;a11.26.1}
Let $\lambda\in\cnum^{\ast}$ be generic.
The $V$-filtrations of $\tilde{\gbige}^{\lambda}_i$
and $\gbige^{\lambda}_i$, induced by $\lefttop{i}V^{(\lambda)}$, are same.
\end{lem}
\pf
Let $\lefttop{i}V$ denote the induced $V$-filtrations
$\lefttop{i}V^{(\lambda)}$
of $\tilde{\gbige}^{\lambda}_i$ and $\gbige^{\lambda}_i$.
We put $A:=\KMS(\gbige,z_i)=\KMS(\tilde{\gbige},z_i)$.
We put
$A(c):=
 \bigl\{u\in\KMS(\gbige,z_i)\,\big|\,
    \paramap(\lambda,u)=c\bigr\}$.
Note we have the following:
\[
 \lefttop{i}\Gr^{V}_c(\tilde{\gbige}^{\lambda})
=\bigoplus_{u\in A(c)}
 \EE(z_i\deldel_i,-\eigenmap(\lambda,u)),
\quad\quad
 \lefttop{i}\Gr^{V}_c(\gbige^{\lambda})
=\bigoplus_{u\in A(c)}
 \EE(z_i\deldel_i,-\eigenmap(\lambda,u)).
\]
Then the coincidence of the $V$-filtration
follows from the uniqueness of the Kashiwara-Malgrange filtration.
\hfill\qed

\vspace{.1in}

Let $\lambda\in\cnum^{\ast}$ be generic.
Let $V$ be the induced $V$-filtration
of $\gbige^{\lambda}_i$ or $\tilde{\gbige}^{\lambda}_i$
by $\lefttop{i}V^{(\lambda)}$.

\begin{lem}
Let $f$ be a section of
$V_{\vecb}(\tilde{\gbige}^{\lambda}_i)$.
Then $f$ naturally gives the section of
$V_{\vecb}(\gbige^{\lambda}_i)$.
\end{lem}
\pf
It follows from Lemma \ref{lem;a11.26.1}.
\hfill\qed

\vspace{.1in}

Let us pick a point $\lambda_0\in\cnum_{\lambda}$
and sufficiently small positive number $\epsilon_0$.
We will restrict our attention to the restrictions
of $\gbige$ and $\tilde{\gbige}$ to $\nbigx(\lambda_0,\epsilon_0)$.
\begin{lem}
Let $f$ be a section of $\lefttop{\nbar}\Vzero_{<0}(\tilde{\gbige})$.
Then $f$ naturally gives the section of
$\lefttop{\nbar}\Vzero_{<0}(\gbige)$.
\end{lem}
\pf
Let $\lambda\in\Delta(\lambda_0,\epsilon_0)$ be generic.
Then the restriction $f_{|\nbigx_i^{\lambda}}$
gives the section of
$V_{<0}(\tilde{\gbige}^{\lambda}_i)
=V_{<0}(\gbige^{\lambda}_i)$.
Due to Corollary \ref{cor;11.28.15},
$f$ gives the section of
$\lefttop{\nbar}\Vzero_{<0}\bigl(\gbige\bigr)$.
\hfill\qed

\vspace{.1in}

Therefore we obtain the inclusion
$j:\lefttop{\nbar}\Vzero_{<0}(\tilde{\gbige})
\subset
 \lefttop{\nbar}\Vzero_{<0}(\gbige)$.
Let $\lambda\in\Delta(\lambda_0,\epsilon_0)$ be generic.
\begin{lem}
The restriction of $j$ to generic $\lambda$
is isomorphic.
\end{lem}
\pf
The restriction of $j$ to $\bigcup_{i=1}^n \nbigx_i^{\lambda}$
is isomorphic due to Lemma \ref{lem;a11.26.1}.
Since $\nbigx-\bigcup_{i=1}^n \nbigx_i^{\lambda}$
is of codimension 2 in $\nbigxlambda$
and since both of $\lefttop{n}V(\tilde{\gbige}^{\lambda})$
and $\lefttop{n}V(\gbige^{\lambda})$ are locally free,
the restriction of  $j$ to generic $\lambda$ is isomorphic.
\hfill\qed

\vspace{.1in}

Hence the inclusion $j$ is isomorphic
outside 
the closed subset whose codimension is larger than $2$.
Since both of 
$\lefttop{\nbar}\Vzero_{<0}(\tilde{\gbige})$
and $\lefttop{\nbar}\Vzero_{<0}(\gbige)$
are locally free,
we can conclude that the inclusion $j$ is isomorphic.

Since $\tilde{\gbige}_{|\nbigx(\lambda_0,\epsilon_0)}$
and $\gbige_{|\nbigx(\lambda_0,\epsilon_0)}$
are generated by $\lefttop{\nbar}\Vzero_{<0}$,
they are same as the submodule of
$\iota_{\ast}\nbige_{|\nbigx(\lambda_0,\epsilon_0)}$.
Thus the proof of Proposition \ref{prop;a11.26.2}
is accomplished.
\hfill\qed

%% file: a78.tex

We use the right $\nbigr$-module structures
in this section.
We put $X=\Delta^n$,
$D_i:=\{z_i=0\}$ and $D=\bigcup_{i=1}^n D_i$.
Let $\harmonicbundle$ be a tame harmonic bundle over $X-D$.
Note the remark in the subsubsection \ref{subsubsection;b12.3.1}.

%% file: 32.tex

\subsubsection{Preliminary}

Let $l$ be a positive integer such that $l\leq n$.
Let us pick an element $\vecm=(m_1,\ldots,m_l)\in\seisuu_{>0}^l$,
and we put $g=\prod_{i=1}^{l}z_i^{m_i}$.
\begin{rem}
The meaning of $l$ is different from 
that in the previous section.
\hfill\qed
\end{rem}

We obtain $\gbige$ and $\naiveprolong{\nbige}$,
and thus $i_{g\,\ast}\gbige=\gbige[\deldel_t]$
and $i_{g\,\ast}\bigl(\naiveprolong{\nbige}\bigr)=
 \naiveprolong{\nbige}[\deldel_t]$.
We regard them as a $\nbigr_{\nbigx\times\cnum}$-module.
We have the following formulas (cf. \cite{saito2}):
\begin{equation} \label{eq;9.21.1}
 \begin{array}{ll}
 (u\otimes\deldel_t^j)\cdot a
=u\cdot a\otimes\deldel_t^j,
 &
 (u\otimes\deldel_t^j)\cdot\deldel_i
=u\deldel_i\otimes\deldel_t^j
-u(\del_ig)\otimes\deldel_t^{j+1},\\
\mbox{{}}\\
 (u\otimes \deldel_t^j)\cdot t
=u\cdot g\otimes \deldel_t^j
+j\cdot\lambda\cdot u\otimes\deldel_t^{j-1},
 &
 (u\otimes\deldel_t^j)\cdot\deldel_t
=u\otimes\deldel_t^{j+1}.
 \end{array}
\end{equation}

We put $s_i:=z_i\deldel_i$ and $s:=t\deldel_t$.
Then we have the following:
\[
 (u\otimes\deldel^j)\cdot(s_i+m_i s)
=us_i\otimes \deldel_t^j
+m_i\cdot\lambda\cdot j\cdot u\otimes\deldel_t^j.
\]
Here we put $m_i=0$ for $i>l$.
In the following,
$(u\otimes 1)s^j$ is denoted by $u\otimes s^j$.

\begin{lem}
We have
$ (u\otimes s^j)\cdot s_i
=(u\otimes 1)s_is^j
=us_i\otimes s^j-m_iu\otimes s^{j+1}$.
\end{lem}
\pf
It can be checked by a direct calculation
from (\ref{eq;9.21.1}).
\hfill\qed

\vspace{.1in}

It is easy to check
$\naiveprolong{\nbige}[\deldel_t]=\naiveprolong{\nbige}[s]$.
\begin{lem}
For any section $s$ of $i_{g\,\ast}\naiveprolong{\nbige}$,
we have the following:
\[
 (u\otimes s^j)\cdot\deldel_i
=(u\cdot z_i^{-1}\otimes s^j)\cdot s_i
=u\cdot\deldel_i\otimes s^j
-m_i\cdot u\cdot z_i^{-1}\otimes s^{j+1}.
\]
\hfill\qed
\end{lem}

\subsubsection{The filtration $\Uzero$ and the endomorphisms}

Following Saito \cite{saito2},
we introduce the filtration $\Uzero$.
For any negative real number $b$, we put as follows:
\[
 U_b^{(\lambda_0)}(i_{g\,\ast}\gbige)
:=\bigl(
 \lefttop{\lbar}\Vzero_{b\cdot\vecm}\gbige\otimes 1
 \bigr)\cdot\nbigr_{\nbigx}.
\]
For any real number $b$,
we put as follows:
\[
 U_b^{(\lambda_0)}(i_{g\,\ast}\gbige)
 :=\sum_{\substack{
 \vecc<0,\\
 \vecc+j\cdot\vecm\leq b\cdot\vecm}}
 \bigl(
 \lefttop{\lbar}\Vzero_{\vecc}\gbige
 \otimes \deldel_t^j
 \bigr)\cdot \nbigr_{\nbigx}
=\sum_{\substack{b'+j\leq b,\\ b'<0}}
 \Uzero_{b'}(\gbige[\deldel_t])\cdot\deldel_t^j.
\]
For simplicity, we use the following notation:
\[
 \Psi^{(\lambda_0)}_b:=
 \Gr^{U^{(\lambda_0)}}_b\bigl(
 \gbige[\deldel_t]
 \bigr).
\]
It is the $\nbigr_{\nbigx}$-module.

\begin{lem}
The following immediately follows from the construction.
\begin{itemize}
\item
 $U^{(\lambda_0)}_a\bigl(
 i_{g\,\ast}\gbige\bigr)
 \cdot V_m\nbigr
 \subset U^{(\lambda_0)}_{a+m}\bigl(
 i_{g\,\ast}\gbige\bigr)$.
\item
 The induced morphism $t:U_b^{(\lambda_0)}(i_{g\,\ast}\gbige)
 \lrarr U_{b-1}^{(\lambda_0)}(i_{g\,\ast}\gbige)$
 is onto if $b<0$.
\item
 The induced morphism
 $\deldel_t:\Psi^{(\lambda_0)}_b\lrarr \Psi^{(\lambda_0)}_{b+1}$
is onto if $b\geq -1$.
\hfill\qed
\end{itemize}
\end{lem}

\begin{lem}
The following immediately follows from the definition.
\begin{itemize}
\item $U_b^{(\lambda_0)}$ is a coherent $V_0\nbigr_{\nbigx}$-module.
\item
 The filtration $\Uzero$ is a good $V$-filtration.
\hfill\qed
\end{itemize}
\end{lem}

For $b=0$, we also put as follows:
\[
 U^{\prime\,(\lambda_0)}_0(\gbige[\deldel_t])
:=\bigl(
 \lefttop{\lbar}\Vzero_{0}\gbige\otimes 1
 \bigr)\cdot\nbigr_{\nbigx},
\quad
 \Psi^{\prime\,(\lambda_0)}
:=\frac{U^{\prime\,(\lambda_0)}_0(\gbige[\deldel_t]) }
{ \Uzero_{<0}(\gbige[\deldel_t])}.
\]
\begin{lem}
We have 
$\Uzero_0(\gbige[\deldel_t])
\subset U^{\prime\,(\lambda_0)}_0(\gbige[\deldel_t])$.
In particular, we have the inclusion
$\Psizero_0\subset \Psi^{\prime\,(\lambda_0)}_0$.
\end{lem}
\pf
It immediately follows from the definition.
\hfill\qed

\begin{lem}
The induced morphism
$t:\Psi^{\prime\,(\lambda_0)}\lrarr\Psizero_{-1}$ is isomorphic.
\end{lem}
\pf
It can be checked directly from the definition.
\hfill\qed

\begin{cor}\label{cor;9.22.15}
$\Psizero_0$ is isomorphic to the image of the induced morphism
$ \deldel_t: \Psizero_{-1}\lrarr  \Psi^{\prime\,(\lambda_0)}_0$.
The induced morphism
$t:\Psizero_0\lrarr\Psizero_{-1}$ is injective.
\hfill\qed
\end{cor}

\subsubsection{The endomorphisms}

Let us pick a large integer $N$ such that $N\geq \rank\nbige$.

\begin{lem}
For any section $f$
of $\lefttop{\lbar}\Vzero_{\vecb}\bigl(\naiveprolong{\nbige}\bigr)$,
we have the following:
\[
 f\cdot
 \prod_{i=1}^l\prod_{u\in\nbigk(\nbige,\lambda_0,i,b_i)}
 \bigl(s_i+\eigenmap(\lambda,u)\bigr)^N
\in \lefttop{\lbar}\Vzero_{\vecb-\epsilon\vecdelta}
(\naiveprolong{\nbige}).
\]
Here $\epsilon$ denotes some positive number.
\hfill\qed
\end{lem}

For any section
$f=f\otimes 1\in
 \lefttop{\lbar}\Vzero_{b\cdot\vecm}(\naiveprolong{\nbige})\otimes 1$,
we put as follows:
\[
 \phi_{1,h}(f\otimes 1):=(f\otimes 1)\cdot s_h,
\quad\quad
 \phi_{2\,h}(f\otimes 1):=
 (f\cdot s_h)\otimes 1.
\]

\begin{lem}
The following claims can be checked by a direct calculation.
\begin{itemize}
\item
We have the following equality in $\naiveprolong{\nbige}[s]$
for any $h\leq l$:
\begin{equation}\label{eq;a11.24.1}
 (f\otimes 1)\cdot s=
\frac{1}{m_h}
 \bigl(
 \phi_{2,h}(f\otimes 1)
-\phi_{1,h}(f\otimes 1)
 \bigr).
\end{equation}
\item We have the following:
\[
\begin{array}{l}
 \phi_{1,h}(f\otimes 1)\in
 \lefttop{\lbar}V_{\vecb-\vecdelta_h}
 (\naiveprolong{\nbige})
 \cdot\deldel_h,\\
\mbox{{}}\\
 \prod_{u\in\nbigk(\nbige,\lambda_0,h,b\cdot m_h)}
 \bigl(\phi_{2,h}+\eigenmap(\lambda,u)\bigr)^N
 (f\otimes 1)\in\lefttop{\lbar}V_{\vecb-\epsilon\vecdelta_h}
 (\naiveprolong{\nbige}).
\end{array}
\]
\item
$\phi_{1,h}$ and $\phi_{2,h}$ are commutative.
\hfill\qed
\end{itemize}
\end{lem}

\begin{lem} \label{lem;9.21.2}
We have the following:
\[
 (f\otimes 1)\cdot
 \prod_{u\in\nbigk(\nbige,\lambda_0,h,b\cdot m_h)}
 \Bigl(
 s+\frac{\eigenmap(\lambda,u)}{m_h}
 \Bigr)^N
\in\sum_{j\geq 0}
 \lefttop{\lbar}V_{\vecb-\epsilon\cdot\vecdelta_h}(
 \naiveprolong{\nbige})\cdot\deldel_h^j.
\]
\end{lem}
\pf
We have only to apply the lemmas above to the following:
\[
 \prod_{u\in \nbigk(\nbige,\lambda_0,h,b\cdot m_h)}
 \Bigl(
 \frac{\phi_{2,h}+\eigenmap(\lambda,u)+\phi_{1,h}}{m_h}
 \Bigr)^N
(f\otimes 1).
\]
\hfill\qed

\begin{cor} \label{cor;9.18.50}
For any section $f$ of $\lefttop{\lbar}\Vzero_{b\cdot \vecm}\gbige$,
we have the following:
\[
 (f\otimes 1)\cdot
 \prod_{h=1}^{l}\prod_{u\in \nbigk(\nbige,\lambda_0,h,bm_h)}
 \Bigl(
 s+\frac{\eigenmap(\lambda,u)}{m_h}
 \Bigr)^N
\in \Uzero_{<b}(i_{g\,\ast}\gbige).
\]
\end{cor}
\pf
It immediately follows from Lemma \ref{lem;9.21.2}.
\hfill\qed

\vspace{.1in}
We put as follows:
\begin{equation}\label{eq;04.2.4.5}
\begin{array}{l}
 \KMS(i_{g\,\ast}\gbige^0):=
\bigcup_{h=1}^l
 \bigl\{u\in\real\times\cnum\,\big|\,
 m_h\cdot u\in \KMS(\nbige^0,h)
 \bigr\},
 \\
\mbox{{}}\\
 \nbigk(i_{g\,\ast}\gbige,\lambda_0,b)
:=\bigl\{u\in \KMS(i_{g\,\ast}\gbige^0)
  \,\big|\, \paramap(\lambda_0,u)=b
 \bigr\}.
\end{array}
\end{equation}

\begin{cor} \label{cor;9.18.51}
If $N>l\cdot\rank(\nbige)$,
we have the following:
\[
 (f\otimes 1)\cdot
 \prod_{u\in\nbigk(\gbige[\deldel_t],\lambda_0,b)}
 \bigl(s+\eigenmap(\lambda,u)\bigr)^N
\in \Uzero_{<b}(i_{g\,\ast}\gbige).
\]
\end{cor}
\pf
It immediately follows from Corollary \ref{cor;9.18.50}.
\hfill\qed

\vspace{.1in}

We put as follows:
\begin{equation} \label{eq;9.22.17}
 F_b(x):=\prod_{u\in\nbigk(\gbige[\deldel_t],\lambda_0,b)}
 \bigl(x+\eigenmap(\lambda,u)\bigr)^N.
\end{equation}

\begin{lem} \label{lem;9.22.16}
The action of $F_b(s)$ vanishes on $\Psizero_b$ for any $b$.
\end{lem}
\pf
In the case $b<0$,
it immediately follows from Corollary \ref{cor;9.18.51}.
By using the surjectivity of the induced morphisms
$\deldel_t:\Psizero_{c-1}\lrarr \Psizero_c$ $(c\geq 0)$
and the obvious relation
$\deldel_t\circ F_{c-1}=F_c\circ\deldel_t$ for any real number $c$,
we can reduce the general case to the case $b<0$.
\hfill\qed

\subsubsection{The filtration $\lefttop{\nbar}\Vzero$ on $\Uzero_b$}
\label{subsubsection;a11.24.15}
We put as follows:
\[
 \lefttop{\nbar}V_0\nbigr_{\nbigx}:=
 \nbigo_{\nbigx}[s_1,\ldots,s_n].
\]

For any real number $b<0$ and for any element
$\vecc\in\real^{l}_{\leq 0}\times\real^{n-l}$,
we put as follows:
\[
 \lefttop{\nbar}\Vzero_{\vecc}
 \Uzero_b\bigl(\gbige[\deldel_t]\bigr)
:=\bigl(
 \lefttop{\nbar}\Vzero_{\vecc+b\cdot\vecm}
 \bigl(\gbige\bigr)
 \otimes 1
 \bigr)\cdot \lefttop{\nbar}V_0\nbigr_{\nbigx}.
\]

\begin{lem}
For any negative number $b$ and any element
$\vecc\in\real_{\leq 0}^l\times\real^{n-l}$,
we have the following:
\[
 \lefttop{\nbar}\Vzero_{\vecc}
 \Uzero_b\gbige[\deldel_t]
=\lefttop{\nbar}V_{\vecc+b\cdot\vecm}\bigl(\gbige\bigr)
\otimes \cnum[s].
\]
Both of them are generated by
$\lefttop{\nbar}\Vzero_{\vecc+b\cdot\vecm}(\gbige)$
over $\lefttop{\nbar}V_0\nbigr_{\nbigx}$.
\hfill\qed
\end{lem}

For any element $\vecc\in\real^n$, we put as follows:
\[ 
 \lefttop{\nbar}
 \Vzero_{\vecc}\bigl(\Uzero_b(\gbige[\deldel_t])\bigr)
:=\sum_{(\vecn,\veca)\in S}
 \lefttop{\nbar}\Vzero_{\veca}
 \bigl(\Uzero_b(\gbige[\deldel_t])\bigr)\cdot\deldel^{\vecn},\\
\]
Here we put
$S:=\bigl\{
 (\vecn,\veca)\in\seisuu_{\geq\,0}^l
  \times(\real_{\leq 0}^l\times\real^{n-l})
 \,\big|\,
  \vecn+\veca\leq \vecc
 \bigr\}$.
The following lemma is easy to see.
\begin{lem}
We have
$ \lefttop{\nbar}\Vzero_{\vecc}
 \Uzero_b\gbige[\deldel_t]
\subset
 \lefttop{\nbar}\Vzero_{\vecc+b\cdot\vecm}
 (\naiveprolong{\nbige})\otimes\cnum[s]$.
\hfill\qed
\end{lem}

The submodule
$\lefttop{\nbar}\Vzero_{\vecc}\Uzero_b\gbige[\deldel_t]$
of $\Uzero_b\gbige[\deldel_t]$
induces the submodule of
$\Gr^{\Uzero}_b\gbige[\deldel_t]$.
It is denoted by
$\lefttop{\nbar}V_{\vecc}
 (\Psizero_b)$.

%% file: 32.1.tex

\subsubsection{Reductions I, II and III}

\label{subsubsection;9.18.61}

We put
$\vecd_0:=-\epsilon\cdot\vecdelta_{\nbar-\lbar}
=(\overbrace{0,\ldots,0}^l,-\epsilon,\ldots,-\epsilon)$
for some sufficiently small positive number $\epsilon$.
In the following,
let $\supp(\vecn)$ denote the set $\bigl\{i\,\big|\,n_i\neq 0\bigr\}$
for an element $\vecn\in\real^l$.
We use the coordinate
$(z_1,\ldots,z_{l},x_1,\ldots,x_{n-l})$.
Note we use the notation in the subsubsection
\ref{subsubsection;a12.2.30}.

\vspace{.1in}

\noindent
{\bf Reduction I.}\\
Let $f$ be a section of $\Uzero_b\bigl(\gbige[\deldel_t]\bigr)$
for $b<0$.
We have a development:
\[
 f=\sum_{h,\vecn}f_{h,\vecn}\otimes s^h\deldel_z^{\vecn},
\quad\quad
 \bigl(f_{h,\vecn}\in\lefttop{\lbar}V_{b\cdot\vecm}(\gbige)
 \bigr).
\]
We take a primitive decomposition of 
each section $f_{h,\vecn}$ of $\lefttop{\lbar}V_{b\cdot\vecm}(\gbige)$
given in the subsubsection \ref{subsubsection;9.21.5}.
Namely,
we take a decomposition of $f_{h,\vecn}$ as follows:
\[
 f_{h,\vecn}=\sum_{\vecb\in S(h,\vecn)}
 f_{h,\vecn,\vecb}\otimes s^h\cdot\deldel_z^{\vecn}.
\]
Here $S(h,\vecn)$ is a finite subset of $\real^n$
such that $q_{\lbar}\bigl(S(h,\vecn)\bigr)\subset\nbigs(b\cdot\vecm)$,
and $f_{h,\vecn,\vecb}$ are $\nbar$-primitive sections
of $\gbige$ such that
$\lefttop{\nbar}\deg^{\Vzero}(f_{h,\vecn,\vecb})=\vecb$.

Then we obtain the following decomposition:
\begin{equation} \label{eq;9.18.60}
 f=
 \sum_{h,\vecn}
 \sum_{\vecb\in S(h,\vecn)}
 f_{h,\vecn,\vecb}\otimes s^h\cdot\deldel_z^{\vecn}.
\end{equation}
Here $f_{h,\vecn,\vecb}$ are $\nbar$-primitive
sections of $\gbige$ such that
$\lefttop{\nbar}\deg^{\Vzero}(f_{h,\vecn,\vecb})=\vecb$.

The procedure to take a development
as in (\ref{eq;9.18.60})
is called Reduction I.
Note the following:
if $f_{h,\vecn}\in\lefttop{\nbar}\Vzero_{\vecc}$,
then $S(h,\vecn)\subset \nbigs(\vecc)$.

\vspace{.1in}
\noindent
{\bf Reduction II.}\\
Let $\vecb$ be an element of $\real^n$
such that $q_{\lbar}(\vecb)\leq b\cdot\vecm$.
Let us consider a section 
of the form $f_{\vecb}\otimes s^h\deldel_z^{\vecn}$
for $\nbar$-primitive $f_{\vecb}$
such that $\lefttop{\nbar}\deg^{\Vzero}(f_{\vecb})=\vecb$.
We consider the following condition.
\begin{description}
\item[(A2)]
 There exists $j_0\in\gminis(\vecn)$
 such that $q_{j_0}(\vecb)\leq b\cdot m_{j_0}-1$.
\end{description}
If $(A2)$ is satisfied,
we have
$q_{\lbar}(\vecb+\vecdelta_{j_0})\leq b\cdot \vecm$,
and thus
$f_{\vecb}\cdot z_{j_0}^{-1}$
is contained in
$\lefttop{\lbar}V_{b\cdot\vecm}\gbige$.
Thus we have the following equalities
in $\Uzero_b\gbige[\deldel_t]$:
\[
\begin{array}{ll}
 f_{\vecb}\otimes s^h\deldel_z^{\vecn} &
=\bigl(f_{\vecb}\cdot z_{j_0}^{-1}\bigr)\otimes s^h s_{j_0}
 \cdot\bigl(\deldel_z^{\vecn-\vecdelta_{j_0}}\bigr)\\
 &
=\bigl(f_{\vecb}\cdot z_{j_0}^{-1}\bigr)\cdot s_{j_0}
\otimes s^h\cdot
 \deldel_z^{\vecn-\vecdelta_{j_0}}
-m_{j_0}\cdot
 \bigl(f_{\vecb}\cdot z_{j_0}^{-1}\bigr)\otimes
 s^{h+1}\cdot\deldel_z^{\vecn-\vecdelta_{j_0}}.
\end{array}
\]
We note
$(f_{\vecb}\cdot z_{j_0}^{-1})\cdot s_{j_0}$
and $f_{\vecb}\cdot z_{j_0}^{-1}$
are sections of
$\lefttop{\nbar}V_{\vecb+\vecdelta_{j_0}}\gbige$.
We also note 
$\vecn-\vecdelta_{j,0}\lneq \vecn$.

\vspace{.1in}

\noindent
{\bf Reduction III.}\\
Let $\vecb$ be an element of $\real^n$ such that
$q_{\lbar}(\vecb)\leq b\cdot\vecm$.
Let us consider a section of the form
$f_{\vecb}\otimes s^h\deldel_z^{\vecn}$
for $\nbar$-primitive $f_{\vecb}$
such that $\lefttop{\nbar}\deg^{V}(f_{\vecb})=\vecb$.
We put
$ K(\vecb):=\bigl\{i\,\big|\,i\leq l,\,\,q_i(\vecb)=b\cdot m_i\bigr\}$,
and we consider the following condition:
\begin{description}
\item[(A3)]
 $h\geq |K(\vecb)|$.
\end{description}

Assume that $(A3)$ is satisfied.
Since we have
$f_{\vecb}\cdot \prod_{k\in K}z_k
 \in\lefttop{\nbar}\Vzero_{\vecb-\vecdelta_K}
\subset \lefttop{\lbar}\Vzero_{<b\cdot \vecm}$,
we have the following:
\begin{equation} \label{eq;9.18.55}
 f_{\vecb}\otimes s^{h-|K(\vecb)|}\!\cdot\!\prod_{k\in K(\vecb)}s_k
 =\Bigl(
 f_{\vecb}\cdot\prod_{k\in K(\vecb)}z_i\Bigr)
 \otimes s^{h-|K(\vecb)|}
 \!\cdot\!\prod_{k\in K(\vecb)}\deldel_k
 \in
 \Uzero_{<b}\gbige[\deldel_t].
\end{equation}
In general, we have the following equality for any $u$,
by using the formula (\ref{eq;a11.24.1}):
\begin{equation}\label{eq;9.18.56}
u\otimes s^{h-|K|}\cdot\prod_{k\in K}s_k
=\sum_{K'\subset K}
 u\cdot \prod_{k\in K'}s_k\otimes s^{h-|K'|}\cdot
 \prod_{k\in K-K'}(-m_k).
\end{equation}
Thus we obtain the following:
\begin{lem}\label{lem;a12.2.50}
If $(A3)$ is satisfied,
we have the following decomposition modulo $\Uzero_{<b}$
\begin{equation} \label{eq;9.18.57}
 f_{\vecb}\otimes s^h\deldel_z^{\vecn}
\equiv
 -\sum_{\emptyset\neq K'\subset K(\vecb)}
 (-1)^{|K'|}\cdot
 \Bigl(
 f_{\vecb}\cdot\prod_{k\in K'}\widetilde{s}_k
 \Bigr)\otimes s^{h-|K'|}
\cdot\deldel_z^{\vecn},
\quad\quad
(\modulo \Uzero_{<b}).
\end{equation}
Here we put $\widetilde{s}_i:=m_i^{-1}\cdot s_i$.
In the decomposition {\rm(\ref{eq;9.18.57})}, the following holds:
\begin{itemize}
\item
 $f_{\vecb}\cdot\prod_{k\in K'}\widetilde{s}_k$
 is contained in
 $\lefttop{\nbar}\Vzero_{\vecb}\gbige$.
\item
 We have $h-|K'|<h$.
\end{itemize}
\end{lem}
\pf
It immediately follows
from (\ref{eq;9.18.55}) and (\ref{eq;9.18.56}).
\hfill\qed

\subsubsection{$(B1)$-decomposition}

\begin{df}
A section 
$f$ of $\Uzero_b\gbige[\deldel_t]$ is called $(B1)$,
if it has the following decomposition mod $\Uzero_{<b}$:
\begin{equation}\label{eq;9.2.1}
 f\equiv
 \sum_{h,\vecn}
 \sum_{\vecb\in S(h,\vecn)}
 f_{h,\vecn,\vecb} \otimes s^h\deldel_z^{\vecn}.
\end{equation}
\begin{description}
\item[(B1.1)]
 $f_{h,\vecn,\vecb}$ is $\nbar$-primitive
 such that $\lefttop{\nbar}\deg^{V}(f_{h,\vecn,\vecb})=\vecb$.
\item[(B1.2)]
 Neither $(A2)$ nor $(A3)$ hold.
 i.e., the following holds for each $(h,\vecn,\vecb)$
such that $f_{h,\vecn,\vecb}\neq 0$:
\begin{itemize}
\item
 For any $j\in\supp(\vecn)$,
 we have $b\cdot m_j-1<q_j(\vecb)\leq b\cdot m_j$.
\item
 $h<|K(\vecb)|$.
\end{itemize}
\hfill\qed
\end{description}
\end{df}

\begin{lem} \label{lem;9.18.61}
Any section $f$ of $\Uzero_b\gbige[\deldel_t]$ $(b<0)$
is $(B1)$.
Moreover the condition $(B1.3)$ holds:
\begin{description}
\item[(B1.3)]
 If $f\in
 \sum_{\vecc\in S}\lefttop{\nbar}\Vzero_{\vecc}\Uzero_b\gbige[\deldel_t]$,
 we can take
 $S(h,\vecn)\subset\bigcup_{\vecc\in S}\nbigs(\vecc-\vecn)$
 for any $h$ and any $\vecn$.
\end{description}
\end{lem}
\pf
By using the reduction $I$,
we can take the development (\ref{eq;9.18.60})
satisfying $(B1.1)$ and $(B1.3)$.
Then we have only to show that
each $f_{h,\vecn,\vecb}\otimes s^h\deldel_z^{\vecn}$
satisfies the condition $(B1)$.

We consider the following partial order on the set
$\seisuu_{\geq\,0}^{l}\times\seisuu_{\geq\,0}
=\bigl\{(\vecn,h)\bigr\}$:
\[
 (\vecn',h')\leq (\vecn,h)
\Longleftrightarrow
 \vecn'\leq\vecn\,\mbox{ and } \,h'\leq h.
\]

We put as follows:
\[
\begin{array}{l}
 \nbigs(\vecn,h)
:=\bigl\{(\vecn',h')\in\seisuu^l_{\geq\,0}\times\seisuu_{\geq\,0}\,
  \big|\,(\vecn',h')\leq (\vecn,h)\bigr\},\\
 \mbox{{}}\\
 \nbigs'(\vecn,h)
:=\bigl\{(\vecn',h')\in\nbigs(\vecn,h)\,\bigr|\,(\vecn',h')\neq (\vecn,h)
 \bigr\}.
\end{array}
\]
Then we consider the following claims:
\begin{description}
\item[$P(\vecn,h)$:]
 For any $(\vecn',h')\in \nbigs(\vecn,h)$,
 any section
 $f_{h',\vecn',\vecb'}\otimes s^h\deldel_z^{\vecn}$
with $(B1.1)$ satisfies $(B1.2)$ and $(B1.3)$.
\item[$Q(\vecn,h)$:]
 For any $(\vecn',h')\in \nbigs'(\vecn,h)$,
 any section
 $f_{h',\vecn',\vecb'}\otimes s^h\deldel_z^{\vecn}$
with $(B1.1)$ satisfies $(B1.2)$ and $(B1.3)$.
\end{description}
We have only to show the claims $P(\vecn,h)$ hold
for any $(\vecn,h)\in\seisuu^{l}_{\geq 0}\times \seisuu_{\geq\,0}$.

In the case $(\vecn,h)=(0,0)$,
the condition $(B1.2)$ is trivial.
Thus the claim $P(0,0)$ holds.
The claim $Q(\vecn,h)$ is equivalent to
$\bigcup_{(\vecn',h')\in\nbigs'(\vecn,h)}P(\vecn',h')$.
Hence we have only to show $Q(\vecn,h)\Longrightarrow P(\vecn,h)$.

We make the following procedure to
$f_{h,\vecn,\vecb}\otimes s^h\deldel_z^{\vecn}$:
\begin{itemize}
\item
 If $(A2)$ for $f_{h,\vecn,\vecb}\otimes s^h\deldel_z^{\vecn}$
 holds,
 then we make $\Red II+\Red I+\Red III +\Red I$.
\item
 If $(A3)$ for $f_{h,\vecn,\vecb}\otimes s^h\deldel_z^{\vecn}$ holds,
 then we make $\Red III + \Red I$.
\item
 If neither $(A2)$ nor $(A3)$ hold,
 then $(B1.2)$ for $f_{h,\vecn,\vecb}$ is already satisfied.
 Thus we do nothing.
 Note $(B1.3)$ also holds in this case.
\end{itemize}
By this procedure,
we obtain the development of
$f_{h,\vecn,\vecb}\otimes s^h\deldel_z^{\vecn}$
modulo $\Uzero_{<b}\gbige[\deldel_t]$:
\begin{equation} \label{eq;a12.6.2}
 f_{h,\vecn,\vecb}\otimes s^h\deldel_z^{\vecn}
\equiv\sum_{\substack{
 (\vecn',h')\in\nbigs'(\vecn,h),\\
 \vecb'\in S(\vecn',h')}}
 f'_{h',\vecn',\vecb'}\otimes s^{h'}\deldel_z^{\vecn'}.
\end{equation}
We may apply $Q(\vecn,h)$ 
to each term in the right hand side of (\ref{eq;a12.6.2}).
Thus the induction can proceed,
and we obtain Lemma \ref{lem;9.18.61}.
\hfill\qed

\begin{df}
A development satisfying $(B1.1)$, $(B1.2)$ and $(B1.3)$
is called a $(B1)$-development.
\hfill\qed
\end{df}

%% file: 32.2.tex

\subsubsection{Reductions I' and II'}

\noindent
{\bf Reduction I'}\\
For any section $f\in \Uzero_{<b}\bigl(\gbige[\deldel_t]\bigr)$,
we have a development:
\[
 f=\sum_{h\in\seisuu_{\geq \,0}}
 \sum_{\vecn\in\seisuu_{\geq\,0}^l}
  f_{h,\vecn}\otimes s^h\deldel_z^{\vecn},
\quad\quad
\bigl(f_{h,\vecn}\in \lefttop{\lbar}\Vzero_{<b\cdot\vecm}(\gbige)\bigr).
\]
By taking an $\nbar$-primitive decomposition of each section
$f_{h,\vecn}$ of $\lefttop{\lbar}\Vzero_{<b\cdot\vecm}\gbige$,
we obtain a decomposition of the following form,
as in Reduction I in
the subsubsection \ref{subsubsection;9.18.61}:
\begin{equation}
 f=\sum_{h,\vecn}\sum_{\vecb\in S(h,\vecn)}
 f_{h,\vecn,\vecb}\otimes s^h\deldel_z^{\vecn}.
\end{equation}
Here $S(h,\vecn)$ denote subsets of $\real^n$,
and $f_{h,\vecn,\vecb}$ are $\nbar$-primitive sections
of $\gbige$
such that $\lefttop{\nbar}\deg^{\Vzero}(f_{h,\vecn,\vecb})=\vecb$.

\vspace{.1in}
\noindent
{\bf Reduction II'}\\
Let us consider a section
of $\Uzero_{<b}\gbige[\deldel_t]$
of the form
$f_{\vecb}\otimes s^h\cdot\deldel_z^{\vecn}$
for  $\nbar$-primitive $f_{\vecb}$ such that
$\lefttop{\nbar}(f_{\vecb})=\vecb$.
Let us consider the following condition:
\begin{description}
\item[$(A2')$]
 There exists $j_0\in\gminis(\vecn)$
 such that $q_{j_0}(\vecb)<b\cdot m_{j_0}-1$.
\end{description}

\begin{lem}
If $(A2')$ is satisfied,
we have the following:
\[
 f_{\vecb}\otimes s^h\deldel_z^{\vecn}
=(f_{\vecb}\cdot z_{j_0}^{-1})\cdot s_{j_0}
 \otimes s^h\cdot\deldel_z^{\vecn-\vecdelta_{j_0}}
-m_{j_0}\cdot(f_{\vecb}\cdot z_{j_0}^{-1})\otimes
 s^{h+1}\cdot\deldel_z^{\vecn-\vecdelta_{j_0}}.
\]
\end{lem}
\pf
Similar to Reduction II in the subsubsection
\ref{subsubsection;9.18.61}.
\hfill\qed

\subsubsection{Reduction (B1')}

\begin{df}
A section $f$ of $\Uzero_{<b}\gbige[\deldel_t]$ is called $(B1')$,
if it has a decomposition as follows:
\begin{equation}\label{eq;9.2.2}
 f=
 \sum_{h,\vecn} \sum_{\vecb\in S(h,\vecn)}
 f_{h,\vecn,\vecb}\otimes s^h\deldel_z^{\vecn}.
\end{equation}
\begin{description}
\item[(B1'.1)]
 $f_{h,\vecn,\vecb}$ is $\nbar$-primitive
 and $\lefttop{\nbar}\deg^V(f_{h,\vech,\vecb})=\vecb$.
\item[(B1'.2)]
 $(A2')$ is not satisfied,
i.e.,
 For any $j\in \gminis(\vecn)$,
 the inequalities $b\cdot m_j-1\leq q_j(\vecb)<b\cdot m_j$ hold.
\hfill\qed
\end{description}
\end{df}

\begin{lem}
Any section $f\in \Uzero_{<b}\gbige[\deldel_t]$ is $(B1')$.
\end{lem}
\pf
The argument is essentially same as the proof of
Lemma \ref{lem;9.18.61}.
By using $\Red I'$,
there exists a decomposition (\ref{eq;9.2.2})
of $f$ satisfying $(B1'.1)$.
We have only to show that
$f_{h,\vecn,\vecb}\otimes s^h\deldel_z^{\vecn}$
with $(B1'.1)$ satisfies $(B1'.2)$.

Let us consider the following claims:
\begin{description}
\item[$P(\vecn)$:]
 For any $\vecn'\in\nbigs(\vecn)$,
 any section of the form $f_{h,\vecn',\vecb}\otimes s^h\deldel_z^{\vecn}$ 
 satisfies $(B1')$.
\item[$Q(\vecn)$:]
 For any $\vecn'\in\nbigs'(\vecn)$,
 any section of the form $f_{h,\vecn',\vecb}\otimes s^h\deldel_z^{\vecn}$ 
 satisfies $(B1')$.
\end{description}
We have only to show 
that the claims $P(\vecn)$ holds for any $\vecn\in\seisuu_{\geq \,0}$.
In the case $\vecn=0$,
the claim $P(0)$ is trivial.
The claim $Q(\vecn)$ is equivalent to
$\bigcup_{\vecn'\in\nbigs'(\vecn)}P(\vecn')$.
Thus we have only to show
$Q(\vecn)\Longrightarrow P(\vecn)$.

We make the following procedure to 
a section of the form
$f_{h,\vecn,\vecb}\otimes s^h\deldel_z^{\vecn}$:
\begin{itemize}
\item
 If $(A2')$ for $f_{h,\vecn,\vecb}\otimes s^h\deldel_z^{\vecn}$ holds,
 then we make $\Red II' +\Red I'$.
\item
 If $(A2')$ does not hold,
 $(B1'.2)$ for $f_{h,\vecn,\vecb}\otimes s^h\deldel_z^{\vecn}$
 is already satisfied,
 and thus we do nothing.
\end{itemize}
As in the proof of Lemma \ref{lem;9.18.61},
we obtain a development of
$f_{h,\vecn,\vecb}\otimes s^h\deldel_z^{\vecn}$,
such that $Q(\vecn)$ can be applied to each term.
Thus we are done.
\hfill\qed

%% file: 32.3.tex

\subsubsection{Preliminary}

For any section $f$ of $\Uzero_b\gbige[\deldel_t]$,
we take a $(B1)$-development modulo $\Uzero_{<b}\gbige[\deldel_t]$:
\begin{equation}\label{eq;9.3.1}
f\equiv
 \sum_{h,\vecn}
 \sum_{\vecb\in S(h,\vecn)}
 f_{h,\vecn,\vecb}\otimes s^h\deldel_z^{\vecn}.
\end{equation}

We put as follows:
\[
 \nbigt(f):=
 \Bigl\{\vecn+\vecb\in\real^{n}\,\Big|\,
 \sum_h f_{h,\vecn,\vecb}\otimes s^h\neq 0
 \Bigr\}.
\]
Then we obtain the set of the maximal elements of $\nbigt(f)$,
which is denoted by $\max\nbigt(f)$.

\begin{lem} \label{lem;9.21.10}
Let $\vecc$ be an element of $\nbigt(f)$.
Then we have a decomposition $\vecc=\vecn+\vecb$
such that $f_{h,\vecn,\vecb}\neq 0$ in {\rm (\ref{eq;9.3.1})}.
Such $\vecn$ and $\vecb$ is determined uniquely for $\vecc$.
\end{lem}
\pf
Recall that $\vecn\in\seisuu^{l}_{\geq \,0}$ and
$\vecb\in \real^{n}$.
They satisfy the following conditions:
\begin{itemize}
\item $\vecb+\vecn=\vecc$ holds.
\item We have $q_i(\vecb)=q_i(\vecc)$ in the case
 $q_i(\vecc)\leq b\cdot m_i$, or $i>l$.
\item We have $b\cdot m_i-1<q_i(\vecb)\leq b\cdot m_i$.
\end{itemize}
It is easy to see that these conditions determine
uniquely $\vecb$ and $\vecn$.
\hfill\qed

\vspace{.1in}
Due to Lemma \ref{lem;9.21.10},
we can use the notation `$(\vecn,\vecb)$' to describe
an element of $\nbigt(f)$.
For maximal element $(\vecn,\vecb)$ of $\nbigt(f)$,
we put as follows:
\[
 P_{(\vecn,\vecb)}(f):=
 \sum_{h<|K(\vecb)|}[f_{h,\vecn,\vecb}]\otimes s^h
 \in \bigoplus_{h<|K(\vecb)|}
 \bigl(\lefttop{\nbar}\Gr^{\Vzero}_{\vecb}\gbige\bigr)
 \otimes s^h.
\]

\begin{prop} \label{prop;9.18.70}
The set $\max\nbigt(f)$ is canonically determined for a section $f$
of $\Uzero_b\bigl(\gbige[\deldel_t]\bigr)$.
For each element $(\vecn,\vecb)\in\max\nbigt(f)$,
the section $P_{(\vecn,\vecb)}(f)$ is canonically determined.
\end{prop}
\pf
Let us take other $(B1)$-development of $f$:
\[
 f\equiv
 \sum_{h,\vecn}
 \sum_{\vecb\in S'(h,\vecn)} f'_{h,\vecn,\vecb}\otimes s^h\deldel_z^{\vecn}.
\]
We obtain the sets $\nbigt'(f)$ and $\max\nbigt'(f)$.
We put as follows:
\[
 F:=
\sum_{h,\vecn}
 \sum_{\vecb\in S'(h,\vecn)}f'_{h,\vecn,\vecb}\otimes s^h\deldel_z^{\vecn}
-\sum_{h,\vecn}
 \sum_{\vecb\in S(h,\vecn)}f_{h,\vecn,\vecb}\otimes s^h\deldel_z^{\vecn}
 \in \Uzero_{<b}\bigl(\gbige[\deldel_t]\bigr).
\]
Let $(\vecn_0,\vecb_0)$ be an element of
$\max\bigl(\nbigt(f)\cup\nbigt'(f)\bigr)$.
We have only to show the following:
\begin{equation} \label{eq;9.21.15}
 \sum_h
 \bigl(f_{h,\vecn_0,\vecb_0}-f'_{h,\vecn_0,\vecb_0}\bigr)\otimes s^h
\equiv 0,\quad
\mbox{ in } \sum_{h<|K(\vecb)|}
 \lefttop{\nbar}\Gr_{\vecb}^{\Vzero}(\gbige)\otimes s^h.
\end{equation}
Here we put $f_{h,\vecn_0,\vecb_0}=0$
(resp. $f'_{h,\vecn_0,\vecb_0}=0$)
if $(\vecn_0,\vecb_0)\not\in\nbigt(f)$,
(resp. if $(\vecn_0,\vecb_0)\not\in\nbigt'(f)$).
We assume that the equation (\ref{eq;9.21.15}) does not hold,
and we will derive a contradiction.

We have the natural inclusion
$\gbige[\deldel_t]\subset\naiveprolong{\nbige}[\deldel_t]
 =\naiveprolong{\nbige}[s]$.
We regard $F$ as a section of $\naiveprolong{\nbige}[s]$.
Note the development of a section $u\in\naiveprolong{\nbige}[\deldel_t]$,
of the following form:
\begin{equation} \label{eq;a11.24.5}
 u\otimes s^h\deldel^{\vecn}
=uz^{-\vecn}\otimes s^{h+|\vech|}
\cdot
 \prod(-m_i)^{n_i}
+\sum_{h'<h+|\vecn|}a_{h'}\otimes s^{h'}.
\end{equation}
Here $a_{h'}$ denote sections of $\naiveprolong{\nbige}$.
If $u$ is a section of
$\lefttop{\nbar}\Vzero_{\vecb}\bigl(
\naiveprolong{\nbige}\bigr)$,
then $u\cdot z^{-\vecn}$ and
$a_{h'}$ are sections of
$\lefttop{\nbar}\Vzero_{\vecb+\vecn}\bigl(
\naiveprolong{\nbige}\bigr)$.

Then we have the following:
\[
 (f_{h,\vecn,\vecb}-f'_{h,\vecn,\vecb})\otimes s^h\deldel^{\vecn}
\in
 \bigoplus_{h\leq h'\leq h+|\vecn|}
 \lefttop{\nbar}\Vzero_{\vecb+\vecn}
 \bigl(\naiveprolong{\nbige}\bigr)\otimes s^{h'}.
\]
We put $S''(h,\vecn):=S'(h,\vecn)\cup S(h,\vecn)$.
Then we have the following:
\[
 F\in
 \sum_{\vecn,h}\sum_{\vecb\in S''(h,\vecn)}
 \Bigl[
 \bigoplus_{h'<|K(\vecb)|+|\vecn|}
 \lefttop{\nbar}\Vzero_{\vecb+\vecn}
 \bigl(\naiveprolong{\nbige}\bigr)\otimes s^{h'}
\Bigr]=:\nbigl.
\]
Since $\vecn_0+\vecb_0$ is assumed to be maximal,
we have the following naturally defined projection:
\[
\begin{CD}
\nbigl @>{\pi_{\vecn_0+\vecb_0}}>>
 {\displaystyle \bigoplus_{h'}}
 \lefttop{\nbar}\Gr_{\vecb_0+\vecn_0}^{\Vzero}
  \bigl(\naiveprolong{\nbige}\bigr)
 \otimes s^{h'}.
\end{CD}
\]
We have the following:
\[
 \pi_{\vecn_0+\vecb_0}(F)
=\pi_{\vecn_0+\vecb_0}
\Bigl(
 \sum_h
 (f_{h,\vecn_0,\vecb_0}-f'_{h,\vecn_0,\vecb_0})
  \otimes
 s^h\deldel^{\vecn_0}
\Bigr).
\]
We put as follows:
\[
 h_0:=\max\bigl\{h\,\,\big|\,\,
 \bigl[f_{h,\vecn_0,\vecb_0}-f'_{h,\vecn_0,\vecb_0}\bigr]\neq 0
 \mbox{ in }\lefttop{\nbar}\Gr^{\Vzero}_{\vecb}(\gbige)
 \bigr\}.
\]
By our choice of $h_0$,
the coefficient of $s^{h'}$ in $\pi_{\vecn_0+\vecb_0}(F)$
is as follows,
due to (\ref{eq;a11.24.5}):
\begin{equation} \label{eq;9.21.22}
  \left\{
 \begin{array}{ll}
 0, & (h'>h_0+|\vecn_0|), \\
 \mbox{{}}\\
 \big[z^{-\vecn_0}\cdot
  (f_{h_0,\vecn_0,\vecb_0}-f'_{h_0,\vecn_0,\vecb_0})\big],
  & (h'=h_0+|\vecn_0|).
 \end{array}
 \right.
\end{equation}

On the other hand,
we have a $(B1')$-development of $F\in \Uzero_{<b}$:
\[
 F=\sum_{h,\vecn}\sum_{\vecb\in S_0'(h,\vecn)}
 F_{h,\vecn,\vecb}\otimes s^h\deldel_z^{\vecn}
\]
Here $F_{h,\vecn,\vecb}$ is $\nbar$-primitive
section of $\gbige$,
such that $\lefttop{\nbar}\deg^{\Vzero}(F_{h,\vecn,\vecb})=\vecb$.
As before we have the following, under the inclusion
$\gbige[\deldel_t]\subset\naiveprolong{\nbige}[s]$:
\[
 F\in
 \sum_{h,\vecn}\sum_{\vecb\in S_0'(h,\vecn)}
 \lefttop{\nbar}V_{\vecb+\vecn}
 \bigl(\naiveprolong{\nbige}\bigr)\otimes s^h
=:\nbigl_0'
\]
We put
 $\nbigt_0(F):=\bigl\{\vecn+\vecb\in\real^n\,\big|\,
 \sum_h F_{h,\vecn,\vecb}\otimes s^h\neq 0
 \bigr\}$.

\begin{lem} \label{lem;9.21.20}
For any element $\vecc\in\nbigt_0(F)$,
the decomposition $\vecc=\vecn_1+\vecb_1$ is uniquely determined
by the following conditions:
\begin{itemize}
\item
 $\vecn_1\in\seisuu_{\geq \,0}$
 and $\vecb\in\real^n$.
\item
 $q_i(\vecb)=q_i(\vecc)$ if $q_i(\vecc)<b\cdot m_i$,
 or if $i>l$.
\item
 $b\cdot m_i-1\leq q_i(\vecb)<b\cdot m_i$.
\end{itemize}
\end{lem}
\pf
Similar to Lemma \ref{lem;9.21.10}.
\hfill\qed

\vspace{.1in}

Let $\vecn_1+\vecb_1\in \nbigt_0(F)$ be a maximal element.
We have the naturally defined morphism:
\[
\begin{CD}
 \nbigl_0'@>{\pi_{\vecn_1+\vecb_1}}>>
 \bigoplus_{h'}\lefttop{\nbar}\Gr_{\vecb_1+\vecn_1}
 (\naiveprolong{\nbige})\otimes s^{h'}.
\end{CD}
\]
We have
$\pi_{\vecn_1+\vecb_1}(F)
=\pi_{\vecn_1+\vecb_1}\bigl(
 \sum_h F_{h,\vecn_1,\vecb_1}\otimes s^h\deldel_z^{\vecn_1}
 \bigr)$.

\begin{lem}
$\vecn_0+\vecb_0$ is an element of $\max\nbigt_0(F)$.
\end{lem}
\pf
It follows from Lemma \ref{lem;9.21.16}.
\hfill\qed

\vspace{.1in}

Let us pick the decomposition $\vecn_1+\vecb_1=\vecn_0+\vecb_0$,
uniquely given in Lemma \ref{lem;9.21.20}

\begin{lem} \label{lem;9.21.21}
We have the following relation:
\[
 \vecn_1=\vecn_0+\vecdelta_{K(\vecb_0)},
\quad\quad
 \vecb_1=\vecb_0-\vecdelta_{K(\vecb_0)}.
\]
\end{lem}
\pf
First we remark $q_i(\vecb_0)\leq b\cdot m_i$ for any $i=1,\ldots,l$.
In the case $q_i(\vecb_0)<b\cdot m_i$,
it is easy to check that $q_i(\vecb_0)=q_i(\vecb_1)$.
In the case $q_i(\vecb_0)=b\cdot m_i$,
we have $q_i(\vecb_1)=b\cdot m_i-1$ and $q_i(\vecn_1)=1$,
because of the condition $q_i(\vecb_1)<b\cdot m_i$.
Thus we are done.
\hfill\qed

\vspace{.1in}

We put $h_1:=\max\{h\,|\,F_{h,\vecn_1,\vecb_1}\neq 0\}$.
The coefficient of $s^{h'}$ in $\pi_{\vecn_0+\vecb_0}(F)$
is as follows:
\begin{equation}\label{eq;a11.24.10}
 \left\{
\begin{array}{ll}
 [z^{-\vecn_1}\cdot F_{h_1,\vecn_1,\vecb_1}]\neq 0
 & (h'=h_1+|\vecn_1|)\\
 \mbox{{}}\\
 0 & (h'>h_1+|\vecn_1|).
\end{array}
 \right.
\end{equation}
Note the following inequality due to Lemma \ref{lem;9.21.21}
and $h_0<|K(\vecb_0)|$:
\[
 h_1+|\vecn_1|
=h_1+|\vecn_0|+|K(\vecb_0)|
\geq |\vecn_0|+|K(\vecb_0)|
 >|\vecn_0|+h_0
\]
Then (\ref{eq;9.21.22}) and (\ref{eq;a11.24.10}) contradict.
Therefore the equation (\ref{eq;9.21.15}) holds,
and thus the proof of Proposition \ref{prop;9.18.70} is accomplished.
\hfill\qed

\begin{cor} \label{cor;9.22.1}
For any section 
$f\in \Psizero_b=\Gr^{\Uzero}_b\bigl(\gbige[\deldel_t]\bigr)$,
the set $\max\nbigt(f)$ is canonically determined.
For any elements $(\vecn,\vecb)\in\max\nbigt(f)$,
the sections
$P_{(\vecn,\vecb)}(f)\in \lefttop{\nbar}\Gr^{\Vzero}_{\vecc}\Psizero_b$
are canonically determined.
\hfill\qed
\end{cor}

\noindent
{\bf Notation}
For any element $\vecc\in\real^n$,
the decomposition $\varpi(\vecc)+\vartheta(\vecc)=\vecc$
is determined by the following conditions:
\begin{itemize}
\item $\varpi(\vecc)\in\real^n$ and
 $\vartheta(\vecc)\in\real^l$.
\item
  We have $q_i(\varpi(\vecc))=q_i(\vecc)$
 in the case $q_i(\vecc)\leq 0$, or in the case $i>l$.
\item
 We have $-1<q_i(\varpi(\vecc))\leq 0$
 in the case $q_i(\vecc)> 0$ and $i\leq l$.
\end{itemize}
We also use the notation
$M(\vecc):=\bigl\{i\,\big|\,i\leq l,\,q_i(\vecc_i)=0\bigr\}$
for any $\vecc\in\real^n$.

Recall that we put
$\Psizero_b:=\Gr^{\Uzero}_b\bigl(
 \gbige[\deldel_t]\bigr)$.
We have the naturally defined morphism
(see the subsubsection \ref{subsubsection;a11.24.15}
 for the definition of $\lefttop{\nbar}\Vzero$ on $\Uzero_b$
 and $\Psi^{(\lambda_0)}_b$):
\[
 \lefttop{\nbar}\Gr^{\Vzero}_{\varpi(\vecc)+b\cdot\vecm}
\bigl(\gbige
 \bigr)[s]\cdot\deldel_z^{\vartheta(\vecc)}
\lrarr
\lefttop{\nbar}\Gr^{\Vzero}_{\vecc}\bigl(
 \Psizero_b
\bigr).
\]
It induces the following morphism:
\begin{equation} \label{eq;9.21.25}
 \frac{
 \lefttop{\nbar}\Gr^{\Vzero}_{\varpi(\vecc)+b\cdot\vecm}
 \bigl(\gbige
 \bigr)[s]}
{\prod_{i\in M(\varpi(\vecc))}(s-\widetilde{s}_i^{\circ})}
\cdot\deldel_z^{\vartheta(\vecc)}
\lrarr
 \lefttop{\nbar}\Gr^{\Vzero}_{\vecc}
 \Psizero_b.
\end{equation}
Here the action of $\widetilde{s}_i^{\circ}$ on
$\lefttop{\nbar}\Gr^{\Vzero}_{\varpi(\vecc)+b\cdot\vecm}(\gbige)[s]$
is induced by the action of $\widetilde{s}_i$ on
$\lefttop{\nbar}\Gr^{\Vzero}_{\varpi(\vecc)+b\cdot\vecm}(\gbige)$.

\begin{cor} \label{cor;9.22.5}
The morphism {\rm (\ref{eq;9.21.25})} is isomorphic.
In particular,
$\lefttop{\nbar}\Gr^{\Vzero}_{\vecc}
 \Psizero_b$ is a locally free $\nbigo_{\nbigd_{\nbar}}$-module.
\end{cor}
\pf
By using Proposition \ref{prop;9.18.70},
it is easy to see that the morphism {\rm (\ref{eq;9.21.25})} is isomorphic.
\hfill\qed

\subsubsection{$\nbar$-primitive decomposition}

\begin{df}
Let $f$ be an element of $\Psizero_b$ for $b<0$.
We put as follows:
\[
 \lefttop{\nbar}\Prim(f):=
 \bigl\{\vecc\in\real^n\,\big|\,
 \vecc+b\cdot\vecm\in\max\bigl(\nbigt(f)\bigr)\bigr\}.
\]
For any element $\vecc\in\lefttop{\nbar}\Prim(f)$,
we put as follows:
\[
 P_{\vecc}(f):=
 \Bigl[
 \sum_h
 f_{h,\vartheta(\vecc),\varpi(\vecc)+b\cdot\vecm}\otimes
 s^h\cdot\deldel_z^{\vartheta(\vecc)}
 \Bigr]
 \in \lefttop{\nbar}\Gr^{\Vzero}_{\vecc}\Psizero_b.
\]
It is well defined due to Corollary {\rm\ref{cor;9.22.1}}.
\hfill\qed
\end{df}

\begin{df}
Let $\vecb$ be an element of $\real^n$.
A section $f$ of $\Psizero_b$ is called $\nbar$-primitive
such that $\lefttop{\nbar}\deg^{\Vzero}(f)=\vecb$,
if the following holds:
\begin{itemize}
\item $f\in \lefttop{\nbar}\Vzero_{\vecb}\Psizero_b$.
\item $[f]\neq 0$ in $\lefttop{\nbar}\Gr^{\Vzero}_{\vecb}\Psizero_b$.
\hfill\qed
\end{itemize}
\end{df}

\begin{df}
For any section $f$ of $\Psizero_b$,
an $\nbar$-primitive development is defined to be
a development $f=\sum_{\vecb\in S} f_{\vecb}$,
where $S$ denotes a finite subset of $\real^n$
and $f_{\vecb}$ $(\vecb\in S)$ are $\nbar$-primitive
such that $\lefttop{\nbar}\deg^{\Vzero}(f_{\vecb})=\vecb$.
\hfill\qed
\end{df}

\begin{lem}
For a section $f\in \Gr^U_b\gbige[\deldel_t]$,
we have a $\nbar$-primitive development:
\[
 f=\sum_{\vecc\in\lefttop{\nbar}\Prim(f)}f_{\vecc},
\quad\quad
 f_{\vecc}\in
 \lefttop{\nbar}\Vzero_{\vecc}\Gr^{\Uzero}_b\gbige[\deldel_t].
\]
Here $f_{\vecc}$ is primitive
such that $\lefttop{\nbar}\deg(f_{\vecc})=\vecc$.
We have
$[f_{\vecc}]=P_{\vecc}(f)$
in $\lefttop{\nbar}\Gr^{\Vzero}_{\vecc}\Psizero_b$.
\hfill\qed
\end{lem}

\begin{df}
Let $I$ be a subset of $\nbar$.
For any element $\vecb\in\real^I$,
we put as follows:
\[
 \lefttop{I}\Vzero_{\vecb}(\Psizero_b):=
 \bigcup_{\substack{\vecc\in\real^n,\\ q_I(\vecc)=\vecb}}
 \lefttop{\nbar}\Vzero_{\vecc}(\Psizero_b).
\]
For any subset $S\subset\real^I$,
we put $\lefttop{I}\Vzero_{S}=\sum_{\vecb\in S}\lefttop{I}V_{\vecb}$.

We use the notation $\lefttop{i}\Vzero_{b}$
instead of $\lefttop{\{i\}}\Vzero_b$.
\hfill\qed
\end{df}

\begin{lem}
Let $b$ be a negative real number.
Let $S$ be a finite subset of $\real^n$.
We have the following equivalence:
\[
 f\in \lefttop{\nbar}\Vzero_{S}\bigl(\Psizero_b\bigr)
\Longleftrightarrow
 \lefttop{\nbar}\Prim(f)\subset \nbigs(S).
\]
\end{lem}
\pf
$\Longleftarrow$ is clear.
We show $\Longrightarrow$.
We have a development
$f=\sum_{\vecb\in S} f_{\vecb}$ such that
$f_{\vecb}\in \lefttop{\nbar}V_{\vecb}\Psizero_b$.
Due to $(B1.3)$,
we obtain $\Prim(f_{\vecb})\subset\nbigs(\vecb)$.
Thus we obtain the result.
\hfill\qed

\begin{cor} \label{cor;9.18.80}
Let $I$ be a subset of $\nbar$.
Let $S$ be a finite subset of $\real^I$.
We have the following equivalence:
\[
 f\in\sum_{\vecb\in S}
 \lefttop{I}V_{\vecb}\Psizero_b
\Longleftrightarrow
 q_I\bigl(
 \lefttop{\nbar}\Prim(f)\bigr)
 \subset\nbigs(S).
\]
\end{cor}
\pf
$\Longleftarrow$ is clear.
We show $\Longrightarrow$.
For any element $\vecb\in S$,
we pick the element $\tilde{\vecb}\in \real^n$
appropriately satisfying $q_I(\tilde{\vecb})=\vecb$
and the following:
\[
 f\in \lefttop{\nbar}\Vzero_{\tilde{S}}\Psizero_b,
\quad\quad
 \tilde{S}:=\bigl\{
 \tilde{\vecb}\,\big|\,\vecb\in S
 \bigr\}.
\]
Then we have $\lefttop{\nbar}\Prim(f)\subset \nbigs(\tilde{S})$.
It implies
$q_I\bigl(\lefttop{\nbar}\Prim(f)\bigr)\subset\nbigs(S)$.
\hfill\qed

\begin{cor}
Let $I$ be a subset of $\nbar$,
and $\vecb=\bigl(b_i\,|\,i\in I\bigr)$ be an element of $\real^I$.
We have the following:
\[
 \bigcap_{i\in I}
 \lefttop{i}V_{b_i}\Psizero_b
=\lefttop{I}V_{\vecb}\Psizero_b
\]

Let $I$ be a subset of $\nbar$
and $S$ be a finite subset of $\real^I$.
Let $i$ be an element of $\nbar-I$, and $c$ be a real number.
We put $I_1=I\sqcup\{i\}$.
We have the naturally defined subset
$S_1=\bigl\{(\vecb,c)\,\big|\,\vecb\in S\bigr\}
\subset \real^{I_1}$.
Then we have the following:
\[
  \lefttop{i}V_c\Psizero_b
\cap
 \Bigl(
 \lefttop{I}V_{S}\Psizero_b
 \Bigr)
=\sum_j \lefttop{I_1}V_{S_1}\Psizero_b.
\]
\end{cor}
\pf
It is easy to check by using Corollary \ref{cor;9.18.80}.
\hfill\qed

\begin{cor}
Let $I$ be a subset of $\nbar$,
and $i$ be an element of $\nbar-I$.
Let $\vecb$ be an element of $\real^I$ and $c$ be a real number.
We put $I_0=I\sqcup\{i\}$,
and we have the naturally defined element
$\vecb_0=(\vecb,c)\in \real^{I_0}$.
Then we have
$\lefttop{i}\Gr^V_c\lefttop{I}\Gr^{V}_{\vecb}\Psizero_b
\simeq
 \lefttop{I_0}\Gr^{V}_{\vecb_0}\Psizero_b$.
\hfill\qed
\end{cor}


%% file: b32.1.tex

\subsubsection{
The local freeness of the sheaf
$\lefttop{I}\Gr^{\Vzero}_{\vecc}\lefttop{J}\Vzero_{\vecd}\Psizero_b$
 $(\vecd<0)$}

\label{subsubection;9.22.30}

Let $I$ be a subset of $\nbar$,
and we put $J:=\nbar-I$.
We put $I_0:=I\cap\lbar$ and $I_1:=\nbar-I_0$.
We also put $J_0:=J\cap\lbar$ and $J_1:=\nbar-J_0$.

For any element $\vecc\in\real^I$,
we take $\varpi(\vecc)\in\real_{\leq 0}^{I_0}\times\real^{I_1}$
and $\vartheta(\vecc)\in\seisuu_{\geq 0}^{I_0}$
determined by the following conditions:
\begin{itemize}
\item
 $\varpi(\vecc)+\vartheta(\vecc)=\vecc$.
\item
 $q_i(\varpi(\vecc))=q_i(\vecc)$ in the cases
 $q_i(\vecc)\leq 0$ or $i\in I_1$.
\item
 $-1<q_i(\varpi(\vecc))\leq 0$
 in the case $q_i(\vecc)\geq 0$ 
 and $i\in I_0$.
\end{itemize}
We put
$M(\varpi(\vecc)):=\bigl\{i\in I_0\,\big|\,q_i(\varpi(\vecc))=0\bigr\}$.

For any elements $\vecc\in\real^I$ and $\vecd\in\real_{<0}^J$,
we have the following morphisms:
\[
 \lefttop{\nbar}\Vzero_{\varpi(\vecc)+\vecd+b\cdot\vecm}
 \bigl(\gbige
 \bigr)[s]\cdot\deldel_z^{\vartheta(\vecc)}
\lrarr
 \lefttop{\nbar}\Vzero_{\vecc+\vecd}
 \Uzero\bigl(\gbige[\deldel_t]\bigr)
\lrarr
 \lefttop{\nbar}\Vzero_{\vecc+\vecd}\Psizero_b
\lrarr
 \lefttop{I}\Gr^{\Vzero}_{\vecc}
 \lefttop{I}\Vzero_{\vecd}\Psizero_b.
\]
The composite of the morphisms is denoted by $\tilde{\Phi}$.
Note the following:
\begin{itemize}
\item
 $\lefttop{\nbar}\Vzero_{\varpi(\vecc)+\vecd+b\vecm}
 \bigl(\gbige\bigr)[s]\cdot
 \prod_{i\in M(\varpi(\vecc))}(s-\tilde{s}^{\circ}_i)
 \cdot\deldel_z^{\vecn}
 \subset \Uzero_{<b}\gbige[\deldel_t]$
 (Lemma \ref{lem;a12.2.50}).
 Here the action of $\tilde{s}_i^{\circ}$ 
 on $\lefttop{\nbar}\Vzero_{\varpi(\vecc)+\vecd+b\vecm}(\gbige)[s]$
 is given by
 the action of $\tilde{s}_i$
 on $\lefttop{\nbar}\Vzero_{\varpi(\vecc)+\vecd+b\vecm}(\gbige)$.
\item
 The image of
$\lefttop{I}\Vzero_{\lneq\varpi(\vecc)+q_I(b\vecm)}
 \lefttop{J}\Vzero_{\vecd+q_J(b\vecm)}
 \bigl(\gbige\bigr)[s]$
via $\tilde{\Phi}$ is contained
in $\lefttop{I}V_{\lneq\vecc}\lefttop{J}V_{\vecd}\Psizero_b$.
\end{itemize}
Thus the following morphism $\Phi$ is induced:
\begin{equation} \label{eq;9.22.3}
 \Phi:
 \frac{
 \lefttop{I}T^{(\lambda_0)}
 \bigl(\varpi(\vecc)+b\vecm_I,\, \vecd+b\vecm_J\bigr)[s]}
{\prod_{i\in M(\varpi(\vecc))}(s-\tilde{s}^{\circ}_i)}
\cdot\deldel_z^{\vartheta(\vecc)}
\lrarr
 \lefttop{I}\Gr^{\Vzero}_{\vecc}\lefttop{J}\Vzero_{\vecd}
 \Psizero_b.
\end{equation}
See the subsubsection \ref{subsubsection;a11.23.20}
for $\lefttop{I}T^{(\lambda_0)}$.
The action of $\tilde{s}_i^{\circ}$
on
$\lefttop{I}T^{(\lambda_0)}
 \bigl(\varpi(\vecc)+b\vecm_I,\, \vecd+b\vecm_J\bigr)[s]$
is induced by the action of $\tilde{s}_i$ on
$\lefttop{I}T^{(\lambda_0)}
 \bigl(\varpi(\vecc)+b\vecm_I,\, \vecd+b\vecm_J\bigr)$.

\begin{prop} \label{prop;9.22.6}
The morphism $\Phi$ is isomorphic.
\end{prop}
\pf
First we show the surjectivity.
Let $f$ be an element of
$\lefttop{\nbar}V_{\vecc+\vecd}\Uzero_b(\gbige[\deldel_t])$.
Let us take an $\nbar$-primitive decomposition of $f$:
\[
 f\equiv\sum f_{h,\vecn,\vecb}\otimes s^h\cdot\deldel_z^{\vecn}.
\]
In the development,
we have $\vecb\leq \varpi(\vecc)+\vecd+b\cdot\vecm$
and $K(\vecb)\subset M(\varpi(\vecc))$.
Hence we obtain the surjectivity.

Let us show the injectivity.
We have the natural projection:
\begin{equation} \label{eq;9.22.2}
 \pi:
 \lefttop{I}\Gr^{\Vzero}_{\varpi(\vecc)+b\vecm_I}
 \lefttop{J}\Vzero_{\vecd+b\vecm_J}(\gbige)[s]
 \cdot\deldel_z^{\vartheta(\vecc)}
\lrarr
\frac{
 \lefttop{I}\Gr^{\Vzero}_{\varpi(\vecc)+b\vecm_I}
 \lefttop{J}\Vzero_{\vecd+b\vecm_J}(\gbige)[s]
 \cdot\deldel_z^{\vartheta(\vecc)}
}{\prod_{i\in M(\varpi(\vecc))}(s-\tilde{s}^{\circ}_i)}.
\end{equation}
Note that the right hand side of (\ref{eq;9.22.2})
is same as the left hand side of (\ref{eq;9.22.3}).
Let us pick a non trivial section $f$
of the right hand side of (\ref{eq;9.22.2}).
We can pick a section $\tilde{f}$
of the left hand side of (\ref{eq;9.22.2})
satisfying the following:
\begin{itemize}
\item $\pi(\tilde{f})=f$.
\item
 We have a decomposition
 $\tilde{f}
 =\sum_{\vecb\in S_0}
 \tilde{f}_{\vecb,h}\otimes s^h\cdot\deldel_z^{\vartheta(\vecc)}$.
Here $S_0$ denotes a primitive subset of $\real_{<0}^J$,
and $\tilde{f}_{\vecb,h}$ $(\vecb\in S_0)$ are $\nbar$-primitive
such that
$\lefttop{\nbar}\deg^{\Vzero}(\tilde{f}_{\vecb,h})
=\varpi(\vecc)+b\cdot\vecm+\vecb$.
\end{itemize}
To show the injectivity of $\Phi$,
we have only to show $\tilde{\Phi}(\tilde{f})\neq 0$.

We have
 $\tilde{\Phi}(\tilde{f})
 \in
 \lefttop{I}\Gr^{\Vzero}_{\vecc}\lefttop{J}\Vzero_{S_0}\Psizero_b$.
For any element $\vecb\in S_0$,
we have the following morphisms:
\[
 \lefttop{I}\Gr^{\Vzero}_{\vecc}\lefttop{J}\Vzero_{S_0}
 \Psizero_b
\lrarr
 \lefttop{\nbar}\Gr^{\Vzero}_{\vecc+\vecb}
 \Psizero_b
\simeq
 \frac{
 \lefttop{\nbar}\Gr^{\Vzero}_{\varpi(\vecc)+\vecd+b\vecm}(\gbige)[s]
 }{
 \prod_{i\in M(\varpi(\vecc))}(s-\tilde{s}_i^{\circ})
 }\cdot\deldel_z^{\vartheta(\vecc)}.
\]
Here the left morphism is the naturally defined projection,
and the right morphism is the isomorphism given in 
Corollary \ref{cor;9.22.5}.
The composite is denoted by $\pi_{\vecb}$.
Then we have the following,
which can be checked directly from the definition:
\[
 \pi_{\vecb}\bigl(
 \tilde{\Phi}(\tilde{f})
 \bigr)
=\sum_h [\tilde{f}_{\vecb,h}]\otimes s^h\cdot\deldel_z^{\vartheta(\vecc)}.
\]
Then we obtain $\pi_{\vecb}\bigl(\tilde{\Phi}(\tilde{f})\bigr)\neq 0$.
Thus we obtain the injectivity of $\tilde{\Phi}$,
and thus the proof of 
Proposition \ref{prop;9.22.6} is accomplished.
\hfill\qed

\begin{cor}\label{cor;9.22.7}
Let $I$ and $J$ be subsets of $\nbar$ such that $\nbar=I\sqcup J$.
Let $b$ be a negative number,
$\vecc$ be an element of $\real^I$,
$\vecd$ be an element of $\real_{<0}^J$.
Then the sheaf
$\lefttop{I}\Gr^{\Vzero}_{\vecc}\lefttop{J}\Gr_{\vecd}\Psizero_{b}$
is an $\nbigo_{\nbigd_{I}}$-locally free sheaf.
In particular, it is strict.
\hfill\qed
\end{cor}

%% file: b32.2.tex

\subsubsection{Strictness of $\Psizero_b$
 and the strict $S$-decomposability $\gbige[\deldel_t]$ along $t=0$}

\begin{lem} \label{lem;9.22.8}
Let $b$ be a negative number.
Let $I$ and $J$ be subsets of $\nbar$
such that $I\sqcup J=\nbar$.
Let $\vecc$ and $\vecd$ be elements of $\real^I$ and $\real^J$
respectively.
Then the sheaf
$\lefttop{I}\Gr^{\Vzero}_{\vecc}\lefttop{J}\Vzero_{\vecd}
  \bigl(\Psizero_b\bigr)$
is strict.
\end{lem}
\pf
For any element $\vecd\in\real^J$, we put as follows:
$M_+(\vecd):=\bigl\{i\,\big|\,q_i(\vecd)\geq 0\bigr\}$.
Note that $|M_+(\vecd)|\leq |J|$.
Let us consider the following claims:
\begin{description}
\item[$P(m,i)$:]
 The strictness holds for $(J,\vecd)$
 such that $\bigl(|J|,|M_+(\vecd)|\bigr)\leq (m,i)$.
\item[$Q(m,i)$:]
 The strictness holds for $(J,\vecd)$
 such that $\bigl(|J|,|M_+(\vecd)|\bigr)\lneq (m,i)$.
\end{description}
We have already known that $P(m,0)$ holds
for any $m$ due to Corollary \ref{cor;9.22.7}.
We have only to show
the implication $Q(m,i)\Longrightarrow P(m,i)$
in the case $i>0$.

Let $(J,\vecd)$ be the tuple such that
$\bigl(|J|,|M_+(\vecd)|\bigr)=(m,i)$ such that $i>0$.
Let us pick the element $j_0\in J$ such that
$q_{j_0}(\vecd)>0$.
 We put $I_1:=I\sqcup\{j_0\}$ and $J_1:=J-\{j_0\}$.
 For an element $\vecc\in\real^{I}$ and a real number $a$,
 we have the naturally defined element
 $(\vecc,a)\in\real^{I_1}$.
 We have the naturally defined projection
 $\pi:\real^{J}\lrarr\real^{J_1}$.

We have the induced filtration $\lefttop{j_0}\Vzero\Psizero_b$
on $\lefttop{I}\Gr^{\Vzero}_{\vecc}\lefttop{J}\Vzero_{\vecd}\Psizero_b$.
Note the following:
\begin{itemize}
\item
 We have
 $\lefttop{j_0}\Gr^{\Vzero}_a\lefttop{I}\Gr^{\Vzero}_{\vecc}
 \lefttop{J}\Vzero_{\pi(\vecd)}\bigl(\Psizero_b\bigr)
 \simeq
 \lefttop{I_1}\Gr^{\Vzero}_{(\vecc,a)}
  \lefttop{J_1}\Vzero_{\pi(\vecd)}\bigl(\Psizero_b\bigr)$.
 It is strict due to the hypothesis of the induction
 $Q(m,i)$.
\item
 Let $N$ be the real number such that $N>q_{j_0}(\vecd)$.
 Then
 $\lefttop{I}\Gr^{\Vzero}_{\vecc}
 \lefttop{J}\Vzero_{\vecd-N\vecdelta_{j_0}}\bigl(
 \Psizero_b\bigr)$ is strict,
 due to the hypothesis of the induction
 $Q(m,i)$.
\end{itemize}
Then the strictness of $(J,\vecd)$ follows easily.
Thus we obtain Lemma \ref{lem;9.22.8}.
\hfill\qed

\begin{cor}\label{cor;a11.24.30}
Let $b$ be a negative number.
Let $I$ be a subset of $\nbar$,
and $\vecc$ be an element of $\real^I$.
Then the sheaf $\lefttop{I}\Gr^{\Vzero}_{\vecc}\bigl(\Psizero_b\bigr)$
is strict.
\end{cor}
\pf
We put $J:=\nbar-I$.
Since $\lefttop{I}\Gr^{\Vzero}_{\vecc}\bigl(\Psizero_b\bigr)$
is the inductive limit of
$\lefttop{I}\Gr^{\Vzero}_{\vecc}\lefttop{J}\Vzero_{\vecd}
   \bigl(\Psizero_b\bigr)$  $(\vecd\in\real^J)$.
Thus the corollary follows from Lemma \ref{lem;9.22.8}.
\hfill\qed

\begin{cor}\label{cor;9.22.10}
Let $b$ be a negative number.
The sheaf $\Psizero_b$ is strict.
\end{cor}
\pf
It immediately follows from Corollary \ref{cor;a11.24.30}.
\hfill\qed

\begin{prop} \label{prop;9.22.9}
Let $b$ be any real number.
The sheaf $\Psizero_b$ is strict.
\end{prop}
\pf
In the case $b<0$, we have already shown the claim in Corollary
\ref{cor;9.22.10}.
Let us consider the case $b=0$.
Then we have the injection
$t:\Psizero_0\lrarr \Psizero_{-1}$
(Corollary \ref{cor;9.22.15}).
Thus the strictness of $\Psizero_0$ follows.

Let us consider the following claims:
\begin{description}
\item[$P(a)$:] The strictness of $\Psizero_b$ holds
for any $b\leq a$.
\item[$Q(a)$:] The strictness of $\Psizero_b$ holds
for any $b< a$.
\end{description}
We have already shown that $P(0)$.
Since the set $\{a\,|\,\Psizero_a\neq 0\}$ is discrete,
we have only to show the implication
$Q(a)\Longrightarrow P(a)$ in the case $a>0$,
which we will show in the following.

Let us consider the following morphisms:
\begin{equation} \label{eq;9.22.18}
\begin{CD}
\Psizero_{a-1}@>{\deldel_t}>>\Psizero_a @>{t}>>\Psizero_{a-1}.
\end{CD}
\end{equation}
By our construction of the filtration $\Uzero$,
the morphism $\deldel_t$ is surjective.
Let us consider the composite,
which is the endomorphism of $\Psizero_{a-1}$
induced by the multiplication of
$\deldel_t\cdot t=t\cdot \deldel_t+\lambda$ from the right.
Due to Lemma \ref{lem;9.22.16},
we have $F_{a-1}(t\cdot\deldel_t)=F_{a}(\deldel_t\cdot t)=0$
on $\Psizero_{a-1}$.
Here $F_a$ is given in (\ref{eq;9.22.17}).
Note that 
$\eigenmap(\lambda,u)$ does not vanish identically,
if we have $\paramap(\lambda_0,u)=a>0$.
Since $\Psizero_{a-1}$ is strict due to our assumption
$Q(a)$, the composite of the morphisms
in (\ref{eq;9.22.18}) is injective. 
Thus the morphism $\deldel_t$ is isomorphic.
It means the strictness of $\Psizero_a$.
Namely the claim $P(a)$ holds,
and thus the proof of Proposition \ref{prop;9.22.9}
is accomplished.
\hfill\qed

\vspace{.1in}

Recall that
the the following endomorphism
identically vanishes
on $\Psizero_b$ for any $b$:
\[
 \prod_{u\in \nbigk(\gbige[\deldel_t],\lambda_0,b)}
 \bigl(s+\eigenmap(\lambda,u)\bigr)^N.
\]
Here $N$ denotes a sufficiently large integer.
Recall that we put as follows for any element
$u\in \nbigk(\gbige[\deldel_t],\lambda_0,b)$:
\[
 \psizero_{t,u}(\gbige[\deldel_t])
 =\Ker\bigl(s+\eigenmap(\lambda,u)\bigr)^N.
\]

\begin{cor}
We have the following:
\begin{enumerate}
\item
The sheaf $\psizero_{t,u}(\gbige[\deldel_t])$ is strict
for any $u$.
\item
The morphism
$t:\psizero_{t,u}(\gbige[\deldel_t])\lrarr
 \psizero_{t,u-\vecdelta_0}(\gbige[\deldel_t])$
is injective, and it is
isomorphic in the cases
\begin{enumerate}
\item
$\paramap(u)<0$.
\item
$\paramap(u)=0$ and $u\neq (0,0)$.
\end{enumerate}
\item
The morphism
$\deldel_t:\psizero_{t,u}(\gbige[\deldel_t])
\lrarr\psizero_{t,u+\vecdelta_0}(\gbige[\deldel_t])$
is isomorphic in the cases
\begin{enumerate}
\item $\paramap(u)>-1$.
\item $\paramap(u)=-1$ and $u\neq (-1,0)$.
\end{enumerate}
If we have $\paramap(u)<-1$,
the morphism $\deldel_t$ is injective.
If we have $u=(-1,0)$,
then the morphism $\deldel_t$ is surjective.
\hfill\qed
\end{enumerate}
\end{cor}

In particular, we obtain the following proposition.
\begin{prop}
The sheaf $\gbige[\deldel_t]$ is strictly $S$-decomposable
along $t=0$.
\hfill\qed
\end{prop}

%% file: 32.4.tex

\subsubsection{The endomorphisms of
$\lefttop{I}\Gr^{\Vzero}_{\vecc}\psi_{t,u}\gbige[\deldel_t]$
and $\lefttop{I}\psi_{\vecu}\psi_{t,u}\gbige[\deldel_t]$}

\label{subsubsection;a11.26.15}

Let pick an element $i\in\nbar$.
In the following, we put $m_i=0$ if $i>l$.
We see the action of $s_i$ on 
$\lefttop{i}\Gr_c^{\Vzero}\psi_{t,u}\gbige[\deldel_t]$,
which is induced by the multiplication from the right.
Assume $b:=\paramap(\lambda_0,u)<0$ and $c\leq 0$.
A section of
$\lefttop{i}\Gr_c^{\Vzero}\psi_{t,u}\gbige[\deldel_t]$
can be described as a sum of
sections of the form
$(f\otimes s^h)\deldel_z^{\vecn}$ such that
$q_i(\vecn)=0$,
where $f\in \lefttop{i}\Vzero_{c+b\cdot m_i}\gbige$.
We have the following formula:
\[
 (f\otimes s^h)\deldel_z^{\vecn}\cdot s_i
=(fs_i\otimes s^h)\deldel_z^{\vecn}
-m_i\cdot (f\otimes s^{h+1})\deldel_z^{\vecn}.
\]

We put as follows:
\[
\begin{array}{l}
 g_{1\,i}\bigl((f\otimes s^h)\deldel_z^{\vecn}\bigr)
:=fs_i\otimes s^h\deldel_z^{\vecn},\\
 \mbox{{}}\\
 g_{2\,i}\bigl((f\otimes s^h)\deldel_z^{\vecn}\bigr)
:=(f\otimes s^h\deldel_z^{\vecn})s_i\\
\mbox{{}}\\
g_s\bigl((f\otimes s^h)\deldel_z^{\vecn}\bigr)
:=f\otimes s^{h+1}\deldel_z^{\vecn}.
\end{array}
\]
In the case $i>l$, we have $g_{1\,i}=g_{2\,i}$.

\begin{lem}
We have
$g_{2\,i}=g_{1\,i}-m_i\cdot g_s$.
The morphisms $g_{2\,i}$, $g_{1\,i}$ and $g_s$ commute.
\end{lem}
\pf
The first claim is clear from the definition.
The commutativity of $g_{2\,i}$ and $g_s$
are easy to see.
Since $g_{1\,i}$ can be described as a linear combination
of $g_{2\,i}$ and $g_s$,
we obtain the second claim.
\hfill\qed

\begin{lem} \label{lem;9.22.21}
On $\lefttop{i}\Gr^{\Vzero}_c\psi_{t,u}(\gbige[\deldel_t])$,
we have the following vanishing
for any sufficiently large integer $N$:
\[
 \prod_{u_1\in\nbigk(\nbige,i,\lambda_0,c+m_ib)}
 \bigl(g_{1\,i}+\eigenmap(\lambda,u_1)\bigr)^N=0,
\quad\quad
 \bigl(g_s+\eigenmap(\lambda,u)\bigr)^N=0.
\]
\end{lem}
\pf
For any section $f$ of $\lefttop{i}\Vzero_{c+m_ib}(\gbige)$,
we have the following:
\[
 f\cdot \prod_{u_1\in\nbigk(\nbige,i,\lambda_0,c+m_ib)}
 (t\cdot \deldel_t+\eigenmap(\lambda,u_1))^N
\in \lefttop{i}\Vzero_{<c+m_ib}(\gbige).
\]
It means the first vanishing.
The second vanishing follows from the definition
of $\psi_{t,u}(\gbige[\deldel_t])$.
\hfill\qed

\vspace{.1in}
We put as follows:
\[
 G_{i,c,N}(x):=
 \prod_{u_1\in\nbigk(\nbige,i,\lambda_0,c+m_ib)}
 \bigl(x+\eigenmap(\lambda,u_1-m_iu)\bigr)^N.
\]
Here $x$ denotes a variable.

\begin{lem} \label{lem;9.22.22}
For a sufficiently large $N$,
we have $G_{i,c,N}(g_{2\,i})=0$
on $\lefttop{i}\Gr_c^{\Vzero}\psi_{t,u}(\gbige[\deldel_t])$:
\end{lem}
\pf
We have the equality
$ g_{2\,i}+\eigenmap(\lambda,u_1-m_iu)
=g_{1\,i}+\eigenmap(\lambda,u_1)
-m_i\cdot \bigl(g_s+\eigenmap(\lambda,u)\bigr)$.
In the case $c\leq 0$,
we obtain $G_{i,c,N}(g_{2\,i})$ by the equality above
and Lemma \ref{lem;9.22.21}.
Since we have the surjection
$\deldel_t^M:\lefttop{i}\Gr^{\Vzero}_{c-M}\lrarr
\lefttop{i}\Gr^{\Vzero}_c$
and the relation
$G_{i,c,N}(g_{2\,i})\circ \deldel_i^M
=\deldel_i^M\circ G_{i,c-M,N}(g_{2\,i})$,
the general case can be reduce to the case $c<0$.
\hfill\qed

\begin{lem}
On $\lefttop{i}\Gr^{\Vzero}_c\lefttop{I}\Gr^{\Vzero}\psi_{t,u}
(\gbige[\deldel_t])$,
we have 
$G_{i,c,N}(g_{2\,i})=0$
for any sufficiently large integer $N$:
\end{lem}
\pf
It follows from Lemma \ref{lem;9.22.22}
and the isomorphism
$ \lefttop{i}\Gr^{\Vzero}_c
 \lefttop{I}\Gr^{\Vzero}\psi_{t,u}
=\lefttop{I}\Gr^{\Vzero}
 \lefttop{i}\Gr_c^{\Vzero}\psi_{t,u}$.
\hfill\qed

\begin{cor}
Let $I$ be a subset of $\nbar$.
Let $\vecc=(c_i\,|\,i\in I)$ be an element of 
$\real^I$.
On $\lefttop{I}\Gr^{\Vzero}_{\vecc}\psi_{t,u}(\gbige[\deldel_t])$,
we have $G_{i,c_i,N}(g_{2\,i})=0$
for any sufficiently large integer $N$.
\hfill\qed
\end{cor}

\begin{lem}
For any element $u_1\in\nbigk(\nbige,i,\lambda_0,c+m_ib)$,
we have $\paramap(\lambda_0,u_1-m_iu)=c$.
\end{lem}
\pf
It can be shown by a direct calculation.
Note the equality
$b=\paramap(\lambda_0,u)$.
\hfill\qed

\begin{cor}
The $\nbigr_{\nbigd_I}$-module
$\lefttop{I}\Gr^{\Vzero}_{\vecc}\psi_{t,u}(\gbige[\deldel_t])$
is strictly specializable along $z_i=0$ $(i\in\nbar-I)$ at $\lambda_0$.
The $V$-filtrations at $\lambda_0$
are given by $\lefttop{i}\Vzero$.
\hfill\qed
\end{cor}

We put 
$\nbigk(\psi_{t,u}\gbige,i,\lambda_0,c)
:=\bigl\{
  u_1-m_i\cdot u\,\big|\,u_1\in\nbigk(\nbige,i,\lambda_0,c+m_ib)
  \bigr\}$,
where we put $b:=\paramap(\lambda_0,u)$.
For any element $\vecc=(c_i\,|\,i\in I)$,
and for any element
$\vecu\in \prod_{i\in I}\nbigk(\psi_{t,u}\gbige,i,\lambda_0,c_i)$,
the submodule
$\lefttop{I}\psizero_{\vecu}\bigl(\psi_{t,u}(\gbige[\deldel_t])\bigr)$ of
$\lefttop{I}\Gr^{\Vzero}_{\vecc}
 \psi_{t,u}(\gbige[\deldel_t])$
is defined as follows,
for any sufficiently large integer $N$:
\[
 \lefttop{I}\psizero_{\vecu}
\bigl(
 \psi_{t,u}(\gbige[\deldel_t])
\bigr)
:=
 \bigcap_{i\in I}\Ker\bigl(g_{2\,i}+\eigenmap(\lambda,q_i(\vecu))\bigr)^N.
\]
We put as follows:
\[
 \KMS(\psi_{t,u}\gbige[\deldel_t],i)
:=\bigl\{
 u_1-m_i\cdot u\,\big|\,
 u_1\in \KMS(\nbige^0,i)
 \bigr\}.
\]

\begin{lem}\mbox{{}}
\begin{enumerate}
\item
Let $I$ be a subset of $\nbar$.
Let $\vecu$ be an element of
$\prod_{i\in I}\KMS(\psi_{t,u}\gbige[\deldel_t],i)$.
Then
$\bigl\{\lefttop{I}\psizero_{\vecu}\psi_{t,u}\gbige[\deldel_t]\,\big|\,
 \lambda_0\in\cnum_{\lambda}
 \bigr\}$
give the $\nbigr_{\nbigd_I}$-module
$\lefttop{I}\psi_{\vecu}\psi_{t,u}\gbige[\deldel_t]$.
\item
Let $i$ be an element of $\nbar-I$.
Then $\lefttop{I}\psi_{\vecu}\psi_{t,u}\gbige[\deldel_t]$
is strictly specializable along $z_i=0$.
\item
Let $i$ be an element of $\nbar-I$,
and $v$ be an element of
$\KMS(\psi_{t,u}\gbige[\deldel_t],i)$.
We put $I_1:=I\sqcup\{i\}$.
We have the naturally defined element
$(\vecu,v)\in \prod_{j\in I_1}\KMS(\psi_{t,u}\gbige[\deldel_t],j)$.
Then we have the following isomorphism:
\[
 \psi_{z_i,v}\lefttop{I}\psi_{\vecu}\psi_{t,u}\gbige[\deldel_t]
\simeq
 \lefttop{I_1}\psi_{(\vecu,v)}\psi_{t,u}\gbige[\deldel_t].
\]
\end{enumerate}
\end{lem}
\pf
It can be shown by an inductive argument.
\hfill\qed

%% file: 32.5.tex

\subsubsection{The decompositions of
$\lefttop{I}\Gr^V_{\vecb_I}\lefttop{J}V_{\vecb_J}(\Psizero_b)$}
\label{subsubsection;a11.26.6}

Let $I\sqcup J=\nbar$ be a decomposition of $\nbar$.
We put $I_0:=I\cap \lbar$ and $I_1:=I-I_0$.
For any element $\vecb\in\real^n$,
we put $\vecb_Y:=q_Y(\vecb)$ for $Y=I,J$.
For any element $\vecb\in\real^n$ such that $\vecb_J<0$,
we put
$ \lefttop{I}\nbigl_{\vecb,b}:=
 \lefttop{I}\Gr^{\Vzero}_{\vecb_I}
 \lefttop{J}\Vzero_{\vecb_J}\bigl(
 \Psizero_b\bigr)$.
We have the right action of $s_i$ $(i\in I)$.
We put $\tilde{s}_i:=m_i^{-1}\cdot s_i$
for $i\in I_0$.
We would like to see the generalized eigen decomposition
of $ \lefttop{I}\nbigl_{\vecb,b}:=
 \lefttop{I}\Gr^{\Vzero}_{\vecb_I}
 \lefttop{J}\Vzero_{\vecb_J}\bigl(
 \Psizero_b\bigr)$
with respect to the actions of $\tilde{s}_i$ $(i\in M(\varpi(\vecb_I)))$.

We have the following isomorphism
(Proposition \ref{prop;9.22.6}):
\[
\lefttop{I}\nbigl_{\vecb,b}
\simeq
 \frac{
 \lefttop{I}T^{(\lambda_0)}
 \bigl(\varpi(\vecb_I)+b\cdot\vecm_I,\, \vecb_J+b\cdot\vecm_J\bigr)[s]}
{\prod_{i\in M(\varpi(\vecb_I))}(s-\tilde{s}^{\circ}_i)}
\cdot\deldel_z^{\vartheta(\vecb_I)}.
\]
Here the action of $\tilde{s}^{\circ}_i$
on $\lefttop{I}T^{(\lambda_0)}
 \bigl(\varpi(\vecb_I)+b\vecm_I,\, \vecb_J+b\vecm_J\bigr)[s]$
is induced by the action of $\tilde{s}_i$
on $\lefttop{I}T^{(\lambda_0)}
 \bigl(\varpi(\vecb_I)+b\vecm_I,\, \vecb_J+b\vecm_J\bigr)$.
We use a similar convention in the following.
\begin{rem}
The actions of $\tilde{s}_i$ and $\tilde{s}_i^{\circ}$
are different.
The relation is $\tilde{s}_i=\tilde{s}_i^{\circ}-s$.
\hfill\qed
\end{rem}

On
$\lefttop{I}T^{(\lambda_0)}
 \bigl(\varpi(\vecb_I)+b\cdot\vecm_I,\vecb_J+b\vecm_J\bigr)$,
we have the action of the tuple of endomorphisms
\[
 \widetilde{\vecs}:=
 \bigl(\widetilde{s}_i\,\big|\,i\in I_0\bigr).
\]
We have the generalized eigen decomposition:
\[
 \lefttop{I}T^{(\lambda_0)}
 \bigl(\varpi(\vecb_I)+b\vecm_I,\vecb_J+b\vecm_J\bigr)
:=\bigoplus_{\vecu\in S(\varpi(\vecb_I),b,\vecm)}
 \EE\bigl(\widetilde{\vecs},\,\,
 -\eigenmap(\lambda,\vecu)
 \bigr).
\]
Here we put
$ S(\varpi(\vecb_I),b,\vecm):=
\bigl\{\bigl(u_i\,\big|\,i\in I_0\bigr)
 \,\big|\, m_i\cdot u_i\in
 \nbigk\bigl(
 \nbige,\lambda_0,i,b\cdot m_i+q_i(\varpi(\vecb_I))
 \bigr)\bigr\}$.

We put
$\nbigq(I,\vecb,\vecu,\vecm):=
 \EE\bigl(\widetilde{\vecs},\,\,
 -\eigenmap(\lambda,\vecu)  \bigr)$,
and then we have the following decomposition:
\[
 \lefttop{I}T^{(\lambda_0)}
  \bigl(\varpi(\vecb_I)+b\vecm_I,\vecb_J+b\vecm_J\bigr)
=\bigoplus_{\vecu\in S(\varpi(\vecb_I),b,\vecm)}
 \nbigq(I,\vecb,\vecu,\vecm).
\]
Then we obtain the following decomposition:
\begin{equation} \label{eq;9.18.90}
 \lefttop{I}\nbigl_{\vecb,b}
\simeq
 \bigoplus_{\vecu\in S(\varpi(\vecb_I),b,\vecm)}
\frac{
 \nbigq(I,\vecb,\vecu,\vecm)[s]}
 {\prod_{i\in M(\varpi(\vecb_I))}(s-\tilde{s}^{\circ}_i)}
\cdot\deldel_z^{\vartheta(\vecb_I)}.
\end{equation}

Let $\vecu=(u_i\,|\,i\in I_0)$ be an element of
$S(\varpi(\vecb_I),b,\vecm)$.
Since we have $q_i(\varpi(\vecb_I))=0$ for any
elements $i\in M(\varpi(\vecb_I))$,
we have $\paramap(\lambda_0,u_i)=b$ for $i\in M(\varpi(\vecb_I))$.
Thus,
if we have $\eigenmap(\lambda_0,u_i)=\eigenmap(\lambda_0,u_j)$
$\bigl(i,j\in M(\varpi(\vecb_I))\bigr)$,
then we obtain $u_i=u_j$.
We also remark that
we have $\paramap(\lambda_0,u_i)\neq b$
for any $i\in I_0-M(\varpi(\vecb_I))$.
We obtain the decomposition of $M(\varpi(\vecb_I))$
as follows:
\begin{equation} \label{eq;9.18.91}
 M(\varpi(\vecb_I))
 =\coprod_{u\in\nbigk(\gbige[\deldel_t],\lambda_0,b)}
M(\vecu,u),
\quad\quad
 M(\vecu,u):=
 \bigl\{
 i\in M(\varpi(\vecb_I))\,\big|\,
 u_i=u\bigr\}
=\bigl\{
 i \in I_0\,\big|\,u_i=u
 \bigr\}.
\end{equation}
Corresponding to the decomposition (\ref{eq;9.18.91}),
we have the following:
\[
 \prod_{i\in M(\varpi(\vecb_I))}
 (s-\tilde{s}^{\circ}_i)
=\prod_{u\in\nbigk(\gbige[\deldel_t],\lambda_0,b)}
 \prod_{i\in M(\vecu,u)}
 (s-\tilde{s}^{\circ}_i).
\]
Then we have the following isomorphism:
\[
\frac{
 \nbigq(I,\vecb,\vecu,\vecm)[s]}
 {\prod_{i\in M(\varpi(\vecb_I))}(s-\tilde{s}^{\circ}_i)}
\simeq
 \bigoplus_{u\in \nbigk(\gbige[\deldel_t],\lambda_0,b)}
 \frac{\nbigq(I,\vecb,\vecu,\vecm)[s]}
 {\prod_{i\in M(\vecu,u)} (s-\tilde{s}^{\circ}_i)}.
\]
We put as follows:
\[
\nbigq(I,\vecb,\vecu,\vecm,u):=
\frac{\nbigq(I,\vecb,\vecu,\vecm)[s]}
 {\prod_{i\in M(\vecu,u)}(s-\tilde{s}^{\circ}_i)}
\]
Then the eigenvalue of the action $s$ 
on $\nbigq(I,\vecb,\vecu,\vecm,u)$ is $-\eigenmap(\lambda,u)$,
and the eigenvalue of $\tilde{s}^{\circ}_j=m_j^{-1}\cdot g_{1\,j}$ is
$-\eigenmap(\lambda,u_j)$.
Note that the eigen functions of
$m_i^{-1}\cdot g_{1\,i}$ is $-\eigenmap(\lambda,u)$
for any $i\in M(\vecu,u)$.
We put $N:=s+\eigenmap(\lambda,u)$
and $\tilde{N}_i:=\tilde{s}_i+\eigenmap(\lambda,u)$
for $i\in M(\vecu,u)$.
Then we have the following:
\begin{equation}\label{eq;b12.6.50}
 \nbigq(I,\vecb,\vecu,\vecm,u)\simeq
\frac{\nbigq(I,\vecb,\vecu,\vecm)[N]}
 {\prod_{i\in M(\vecu,u)}(N-\tilde{N}_i)}.
\end{equation}

We obtain the decomposition, as follows:
\begin{equation}\label{eq;b12.3.10}
 \lefttop{I}\nbigl_{\vecb,b}
\simeq
 \bigoplus_{(u,\vecu)}
 \nbigq(I,\vecb,\vecu,\vecm,u)\cdot\deldel_z^{\vartheta(\vecb_I)}.
\end{equation}
Here
$(u,\vecu)$ runs through the set
$\nbigk(\gbige[\deldel_t],\lambda_0,b)\times
 S(\varpi(\vecb_I),b,\vecm)$.
We also have the decomposition:
\begin{equation} \label{eq;b12.3.11}
 \lefttop{I}\Gr^V_{\vecb_1}
 \lefttop{J}V_{\vecb_2}\psizero_u(\gbige[\deldel_t])
=\bigoplus_{\vecu\in S(\varpi(\vecb_I),b,\vecm)}
 \nbigq(I,\vecb,\vecu,u,\vecm)\cdot\deldel_z^{\vartheta(\vecb_I)}.
\end{equation}

Let us consider the action of $\tilde{s}_i$ $(i\in I_0)$
on $\nbigq(I,\vecb,\vecu,\vecm,u)$.
Recall that we have the relation
$\tilde{s}_i=\tilde{s}_i^{\circ}-s$.

\begin{lem}\mbox{{}}
\begin{itemize}
\item
Let $i$ be an element of $I_0$.
The eigenvalue of $\tilde{s}_i$ on
$\nbigq(I,\vecb,\vecu,\vecm,u)$ is
$-\eigenmap(\lambda,u_i-u)$.
\item
The decompositions {\rm (\ref{eq;b12.3.10})} and
{\rm (\ref{eq;b12.3.11})} are generalized eigen decompositions
with respect to the endomorphism
$\tilde{s}_i$ $(i\in I_0)$ and $s$.
\end{itemize}
\end{lem}
\pf
The second claim immediately follows from the first claim.
The eigenvalue of $s$ is $-\eigenmap(\lambda,u)$
and the eigenvalue of $\tilde{s}_i^{\circ}$
is $-\eigenmap(\lambda,u_i)$.
Thus we obtain the first claim.
\hfill\qed

\begin{cor}\label{cor;b12.3.20}
Let $i$ be an element of $M(\varpi(\vecb))$.
Then $\tilde{s}_i$ is nilpotent on 
$\nbigq(I,\vecb,\vecu,\vecm,u)$ if and only if
$i$ is contained in $M(\vecu,u)$.
\hfill\qed
\end{cor}

%% file: a33.3.tex

\subsubsection{The decomposition of
$\Gr^{W(N)}\nbigq(I,\vecb,\vecu,\vecm,u)$}

\label{subsubsection;a11.26.40}

Recall we put $I_0=I\cap \lbar$ and $I_1=I-I_0$.
We assume that $\vecb_I\leq 0$.
From our construction,
$\nbigq(I,\vecb,\vecu,\vecm)$
is a subbundle of
$\lefttop{I}\Gr^{\Vzero}_{\vecb_I+b\cdot\vecm_I}
 \lefttop{J}\Vzero_{<\vecb_J+b\vecm_J}
 (\naiveprolong{\nbige})$ over $\nbigd_{I}(\lambda_0,\epsilon_0)$.
Moreover, they are contained in
the $-\eigenmap(\lambda,\vecu)$-part
$\EE\bigl(\tilde{\vecs},-\eigenmap(\lambda,\vecu)\bigr)$,
in the generalized eigen decomposition with respect to
the tuple of the morphisms
$\tilde{\vecs}=\bigl(\tilde{s}_i\,\big|\,i\in M(\vecb_I)\bigr)$.
Let $\tilde{\nbign}_i$ denote the nilpotent part of $\tilde{s}_i$
on $\EE(\tilde{\vecs},-\eigenmap(\lambda,\vecu))$.

We also have the action of $s_i$ $(i\in I_1)$
on $\EE\bigl(\tilde{\vecs},-\eigenmap(\lambda,\vecu)\bigr)$.
We denote the nilpotent part of $s_i$ by $\nbign_i$.
We put $R(\vecu):=\{i\in I_1\,|\,u_i=(0,0)\}\subset I_1$.
Then we have $s_i=\nbign_i$ for $i\in R(\vecu)$.
We put $\nbign_{R(\vecu)}:=\prod_{i\in R(\vecu)}\nbign_i$.

\begin{lem}
We have $\nbigq(I,\vecb,\vecu,\vecm)=\Image\nbign_{R(\vecu)}$.
\end{lem}
\pf
It can be directly checked by our construction,
by using Lemma \ref{lem;9.18.11}.
\hfill\qed

\vspace{.1in}

Let $j$ be an element of $R(\vecu)$.
Note 
$\nbigq(I,\vecb-\vecdelta_j,\vecu-\vecdelta_{0,j},\vecm)
\simeq \Image\nbign_{R(\vecu)-\{j\}}$.
Then the naturally defined morphism
$z_j:\nbigq(I,\vecb,\vecu,\vecm)\lrarr
 \nbigq(I,\vecb-\vecdelta_j,\vecu-\vecdelta_{0,j},\vecm)$
is isomorphic to the natural inclusion.
The naturally defined morphism
$\deldel_j:\nbigq(I,\vecb-\vecdelta_j,\vecu-\vecdelta_{0,j},\vecm)
\lrarr \nbigq(I,\vecb,\vecu,\vecm)$
is isomorphic to the morphism $\tilde{\nbign}_j$.
Then we obtain the naturally induced morphisms,
for any element $j\in R(\vecu)$:
\[
 z_j:
 \nbigq(I,\vecb,\vecu,\vecm,u)\lrarr
 \nbigq(I,\vecb-\vecdelta_j,\vecu-\vecdelta_{0,j},\vecm,u),
\quad
 \deldel_j:\nbigq(I,\vecb-\vecdelta_j,\vecu-\vecdelta_{0,j},\vecm,u)
\lrarr \nbigq(I,\vecb,\vecu,\vecm,u).
\]

Let $j$ be an element of $M(\varpi(\vecb_I))\subset I_0$.
Note that the multiplication $\cdot z_j$ from the right
induces the isomorphism
$\phi:\nbigq(I,\vecb,\vecu,\vecm)\lrarr
 \nbigq(I,\vecb-\vecdelta_j,\vecu-\vecdelta_{0,j},\vecm)$.
Then we obtain the following isomorphism:
\begin{equation} \label{eq;9.23.20}
 \nbigq(I,\vecb-\vecdelta_j,\vecu-\vecdelta_{0,j},\vecm,u)
\simeq
 \frac{\nbigq(I,\vecb,\vecu,\vecm)[N]}{
 \prod_{i\in M'(\vecu,u)}(N-N_i).
 }
\end{equation}
Here we put $M'(\vecu,u):=M(\vecu,u)-\{j\}$.
Under the isomorphism {\rm(\ref{eq;9.23.20})},
the morphism $z_j:
 \nbigq(I,\vecb,\vecu,\vecm,u)\lrarr
 \nbigq(I,\vecb-\vecdelta_j,\vecu-\vecdelta_{0,j},\vecm,u)$
is induced by the identity of
$\nbigq(I,\vecb,\vecu,\vecm,u)[N]$.

\begin{lem}
The morphism $\deldel_j:
 \nbigq(I,\vecb-\vecdelta_j,\vecu-\vecdelta_{0,j},\vecm,u)
\lrarr \nbigq(I,\vecb,\vecu,\vecm,u)$
is induced by the 
morphism $m_j\cdot(\tilde{N}_j-N)$ on
$\nbigq(I,\vecb,\vecu,\vecm,u)[N]$.
\end{lem}
\pf
Let $\sum_{j=1}^h\phi(f_i)\cdot N^i$
be a section of the right hand side of
the isomorphism (\ref{eq;9.23.20}).
We have the following:
\begin{multline}
 \Bigl(
 \sum_{j=1}^h\phi(f_i)\cdot N^i
 \Bigr)\cdot\deldel_j
=\Bigl(
 \sum_{j=1}^hf_i\cdot N^i
 \Bigr)\cdot s_j
=\sum_{j=1}^h (f_i\cdot s_j)\cdot N^i
-m_j\cdot \sum_{j=1}^h f_i\cdot N^i\cdot s \\
=\sum_{j=1}^h f_i\cdot (s_j-m_j\cdot s)\cdot N^i
=m_j\sum_{j=1}^h f_i\cdot (\tilde{N}_j-N)\cdot N^i.
\end{multline}
Thus we are done.
\hfill\qed

\begin{lem}\label{lem;a11.26.45}
Assume that $\vecb_I\leq 0$ and $\vecb_J<0$.
Let $i$ be an element of
$M(\varpi(\vecb_I))\cup R(\vecu)$.
Then we have the following decomposition:
\[
  \Gr^{W(N)}\bigl(
 \nbigq(I,\vecb,\vecu,\vecm,u)
 \bigr)
=\Image(\can_i)\oplus \Ker(\var_i).
\]
\end{lem}
\pf
It follows from Proposition \ref{prop;9.23.21}.
\hfill\qed

\vspace{.1in}

Let $j$ be an element of $I_1$.
Then we have the induced morphism:
\begin{equation}\label{eq;9.23.22}
\deldel_j:
 \nbigq(I,\vecb-\vecdelta_j,\vecu-\vecdelta_{0,j},\vecm)
\lrarr
 \nbigq(I,\vecb,\vecu,\vecm).
\end{equation}
\begin{lem} \label{lem;9.23.23}
In the case $q_j(\vecb)>0$,
the morphism {\rm (\ref{eq;9.23.22})} is isomorphic.
Thus the induced morphism
$\deldel_j:
 \nbigq(I,\vecb-\vecdelta_j,\vecu-\vecdelta_{0,j},\vecm,u)
\lrarr
 \nbigq(I,\vecb,\vecu,\vecm,u)$
is isomorphic.
\end{lem}
\pf
It is clear from our construction.
\hfill\qed

\begin{lem}\label{lem;9.23.24}
Let $j$ be an element of $I_0$.
Assume that $q_j(\vecb)>0$.
Then we have
$\nbigq(I,\vecb,\vecu,\vecm)
=\nbigq(I,\vecb-\vecdelta_j,\vecu-\vecdelta_{0,j},\vecm)$
and 
$\nbigq(I,\vecb,\vecu,\vecm,u)
=\nbigq(I,\vecb-\vecdelta_j,\vecu-\vecdelta_{0,j},\vecm,u)$.
\end{lem}
\pf
In the case $q_j(\vecb)>0$,
we have $\varpi(\vecb_I)=\varpi(\vecb_I-\vecdelta_{j})$.
Then the lemma is clear from our construction.
\hfill\qed

\begin{prop} \label{prop;9.23.25}
Assume $\vecb_J<0$.
Let $i$ be an element of
$M(\varpi(\vecb_I))\cup R(\vecu)$.
Then we have the following decomposition:
\[
  \Gr^{W(N)}\bigl(
 \nbigq(I,\vecb,\vecu,\vecm,u)
 \bigr)
=\Image(\can_i)\oplus \Ker(\var_i).
\]
\end{prop}
\pf
We have already checked such decomposition
in the case $\vecb_I\leq 0$.
By using the isomorphisms given in Lemma \ref{lem;9.23.23}
and Lemma \ref{lem;9.23.24},
the general case can be reduced to the case $\vecb_I\leq 0$.
\hfill\qed

%% file: a84.tex

\subsubsection{The special case}
\label{subsubsection;b12.3.150}

Let us consider the case $l=n$,
$I=\nbar$ and $\vecb=0$.
In the case we have $\varpi(\vecb_I)=0=\vartheta(\vecb_I)$.
We also have the following:
\[
 S(0,b,\vecm)=
\bigl\{
 \bigl(
 u_i\,\big|\,i\in\vecn
 \bigr)\,\big|\,
 m_i\cdot u_i\in \nbigk(\nbige,\lambda_0,i,b\cdot m_i)
\bigr\}.
\]
Then we have the following decomposition:
\[
 \lefttop{\nbar}\Gr^{\Vzero}_0
 \bigl(\psi_{t,u}\bigr)
=\bigoplus_{\vecu\in S(0,b,\vecm)}
 \nbigq\bigl(
 \nbar,0,\vecu,\vecm,m
 \bigr),
\quad
 \nbigq(\nbar,0,\vecu,\vecm,u)
=\frac{\nbigq(\nbar,0,\vecu,\vecm)[N]
 }{\prod_{i\in M(\vecu,u)}(N-\tilde{N}_i)}.
\]
\begin{lem}
The endomorphism $\tilde{s}_i$ on
$\nbigq(\nbar,0,\vecu,\vecm,u)$ is nilpotent for any $i\in\nbar$,
if and only if
$u_i=u$ holds for any $i\in\nbar$.
\end{lem}
\pf
It follows from Corollary \ref{cor;b12.3.20}.
\hfill\qed

\vspace{.1in}

Let us consider the case
$\vecu=\bigl(\overbrace{u,\ldots,u}^n\bigr)$.
We put
$\hat{\vecu}=\bigl(m_i\cdot u\,\big|\,i\in\nbar\bigr)+\vecdelta_{0,\nbar}
 \in (\real\times\cnum)^n$.
Note that $\hat{\vecu}\in \prod_{i=1}^n\KMS(\nbige^0,i)$.
\begin{lem}
We have the natural isomorphism:
$\nbigq(\nbar,0,\vecu,\vecm)\simeq
 \lefttop{\nbar}\nbigg_{\hat{\vecu}}(E)$.
(See the subsubsection {\rm\ref{subsubsection;10.16.10}}
for $\lefttop{\nbar}\nbigg_{\hat{\vecu}}(E)$).
\end{lem}
\pf
It can be checked directly from the definitions.
\hfill\qed

\vspace{.1in}
Hence we obtain the isomorphism:
\[
 \nbigq(\nbar,0,\vecu,\vecm,u)
\simeq
 \nbigv_{\nbar}\bigl(
 S_{\hat{\vecu}}^{\can}(E)
 \bigr)_{|\cnum_{\lambda}}.
\]
See the subsubsection \ref{subsubsection;9.13.50}
for $\nbigv_{\nbar}\bigl(S_{\hat{\vecu}}^{\can}(E)\bigr)$.
In particular,
we obtain the subbundle:
\[
 \nbigc_{\nbar\,|\,\cnum_{\lambda}}
\subset P_h\Gr^{W(N)}_h\nbigq(\nbar,0,\vecu,\vecm,u).
\]

%% file: a83.1.tex

In this section,
we use the right $\nbigr$-module structure.

%% file: a80.1.tex

\subsubsection{$\lefttop{I}\nbigl$ and
 the filtrations $\lefttop{\nbar}\Vzero$}

In this subsection,
$\tildepsi_{t,u}(\gbige[\deldel_t])$ is often denoted by
$\tildepsi_{t,u}$, for simplicity.
We put as follows:
\[
\begin{array}{l}
 \Par_{b\cdot m_i}\bigl(
 \nbigelambdazero,i
 \bigr):=
 \bigl\{
 \alpha-b\cdot m_i\,\big|\,
 \alpha\in \Par(\nbigelambdazero,i)
 \bigr\},\\
 \mbox{{}}\\
 \Par_{b\cdot m_i}^-\bigl(
 \nbigelambdazero,i
 \bigr):=
 \bigl\{
 r\in  \Par_{b\cdot m_i}\bigl(
 \nbigelambdazero,i
 \bigr)\,\big|\,r<0
 \bigr\},\\
 \mbox{{}}\\
\Par_{b\cdot m_i}^{\geq 0}\bigl(
 \nbigelambdazero,i
 \bigr):=
 \bigl\{
 r\in  \Par_{b\cdot m_i}\bigl(
 \nbigelambdazero,i
 \bigr)\,\big|\,r\geq 0
 \bigr\}.\\
\end{array}
\]

We put as follows, for any subset $I\subset\nbar$:
\[
 \lefttop{I}\nbigl:=
 \lefttop{I}\Vzero_{<0}(\tildepsi_{t,u})\Big/
 \sum_{I'\supsetneq I}\lefttop{I'}\Vzero_{<0}(\tildepsi_{t,u}).
\]
Then $\lefttop{I}\nbigl$ is 
an $\nbigo_{\nbigd_J(\lambda_0,\epsilon_0)}$-sheaf,
where we put $J:=\nbar-I$.
For any element $\veca\in \real^n$,
we put as follows:
\[
 \lefttop{\nbar}\Vzero_{\veca}\bigl(\lefttop{I}\nbigl\bigr)
:=\Image\Bigl(
 \lefttop{\nbar}\Vzero_{\veca}(\tildepsi_{t,u})\cap
 \lefttop{I}\Vzero_{<0}(\tildepsi_{t,u})
\lrarr
 \lefttop{I}\nbigl
 \Bigr).
\]
Then 
$\lefttop{\nbar}\Vzero_{\veca}\bigl(\lefttop{I}\nbigl\bigr)$
is a coherent $\nbigo_{\nbigd_J(\lambda_0,\epsilon_0)}$-sheaf.

\begin{lem}\mbox{{}}
\begin{itemize}
\item
 Let $j$ be an element of $\nbar-I$.
 In the case $a_j<0$,
 we have $\lefttop{\nbar}\Vzero_{\veca}(\lefttop{I}\nbigl)=0$.
\item
 Let $j$ be an element of $I$.
 In the case $a_j\geq 0$,
 we put as follows:
\[
 a_i':=\left\{
 \begin{array}{ll}
 a_i & (i\neq j)\\ \mbox{{}}\\
 \max 
 \Par^-_{b\cdot m_j}(\nbigelambdazero,j),
 & (i=j).
 \end{array}
 \right.
\]
Then we have
$\lefttop{\nbar}\Vzero_{\veca'}(\lefttop{I}\nbigl)
=\lefttop{\nbar}\Vzero_{\veca}(\lefttop{I}\nbigl)$.
\end{itemize}
\end{lem}
\pf
The claims are clear from our construction
of $\lefttop{I}\Vzero_{<0}(\tildepsi_{t,u})$.
\hfill\qed

\vspace{.1in}

Let $\pi$ denote the projection
$\lefttop{I}\Vzero_{<0}(\tildepsi_{t,u})\lrarr\lefttop{I}\nbigl$.
For any section $f$ of $\lefttop{I}\nbigl$,
we take a section $F$ of $\lefttop{I}\Vzero_{<0}(\tildepsi_{t,u})$
such that $\pi(F)=f$.

\begin{lem}\mbox{{}}
\begin{itemize}
\item
The set 
$\lefttop{\nbar}\Prim(f):=
 \lefttop{\nbar}\Prim(F)\cap\bigl(
 \real^I_{<0}\times\real^J_{\geq 0}\bigr)$
is canonically determined for $f$.
\item
For any element $\veca\in\lefttop{\nbar}\Prim(f)$,
the section $P_{\veca}(f):=P_{\veca}(F)$
of $\lefttop{\nbar}\Gr^{\Vzero}_{\veca}(\tildepsi_{t,u})$
is canonically determined for $f$.
\end{itemize}
\end{lem}
\pf
Since we have
$F-F'\in \sum_{I'\supsetneq I}\lefttop{I'}\Vzero_{<0}(\tildepsi_{t,u})$
for other choice of $F'$,
the claims are clear.
\hfill\qed

\vspace{.1in}
Let $S$ be a subset of $\real^n$.
We put as follows:
\[
 \lefttop{\nbar}V^{(\lambda_0)}_{S}(\lefttop{I}\nbigl)
:=\sum_{\veca\in S}
 \lefttop{\nbar}\Vzero_{\veca}\bigl(\lefttop{I}\nbigl\bigr).
\]

\begin{lem} \label{lem;a11.26.5}
The projection
$\pi:\lefttop{I}\Vzero_{<0}(\tildepsi_{t,u})
\lrarr\lefttop{I}\nbigl$ induces the surjection:
\[
 \lefttop{\nbar}V^{(\lambda_0)}_S(\tildepsi_{t,u})
 \cap \lefttop{I}\Vzero_{<0}(\tildepsi_{t,u})
\lrarr
 \lefttop{\nbar}V^{(\lambda_0)}_S\bigl(\lefttop{I}\nbigl\bigr).
\]
\end{lem}
\pf
It follows from
$\lefttop{\nbar}V^{(\lambda_0)}_S(\tildepsi_{t,u})\cap
 \lefttop{I}V_{<0}(\tildepsi_{t,u})
=\sum_{\veca\in S}\lefttop{\nbar}
 \Vzero_{\veca}(\tildepsi_{t,u})
  \cap \lefttop{I}\Vzero_{<0}(\tildepsi_{t,u})$.
\hfill\qed

\vspace{.1in}

Let $I\sqcup J=\nbar$ be a decomposition.

\begin{cor}
Let $f$ be a section of $\lefttop{I}\nbigl$.
Then $f$ is contained in
$\lefttop{\nbar}V^{(\lambda_0)}_S\bigl(\lefttop{I}\nbigl\bigr)$
if and only if
$\lefttop{\nbar}\Prim(f)$ is contained in $\nbigs(S)$.
\hfill\qed
\end{cor}

\subsubsection{Local freeness}

Let $I\sqcup J=\nbar$ be a decomposition.

\begin{cor} \label{cor;a11.26.16}
Let $\vecb$ be an element of $\real^I_{<0}$
and $\vecc$ be an element of $\real^J_{\geq \,0}$.
Then we have the naturally defined isomorphism:
\[
 \lefttop{J}\Gr^{\Vzero}_{\vecc}
 \lefttop{I}\Vzero_{\vecb}(\tildepsi_{t,u})
\simeq
 \lefttop{J}\Gr^{\Vzero}_{\vecc}
 \lefttop{I}\Vzero_{\vecb}(\lefttop{I}\nbigl).
\]
In particular,
$\lefttop{J}\Gr^{\Vzero}_{\vecc}
 \lefttop{I}\Vzero_{\vecb}(\lefttop{I}\nbigl)$
is coherent and locally free
$\nbigo_{\nbigd_J}$-module.
\end{cor}
\pf
We have the naturally defined morphism
from the left hand side to the right hand side.
By using Lemma \ref{lem;a11.26.5},
we can check that the morphism is isomorphic.
The local freeness follows from 
the result in the subsubsection \ref{subsubsection;a11.26.6}.
\hfill\qed

\begin{lem}\label{lem;a11.26.10}
Let $S$ be a primitive subset of $\real^J_{\geq 0}$.
The $\nbigo_{\nbigd_J}$-module
$\lefttop{J}\Vzero_{S}
 \lefttop{I}\Vzero_{\vecb}(\lefttop{I}\nbigl)$
is locally free and coherent.
\end{lem}
\pf
We have only to check the claims for
primitive subsets $S$,
which are contained in
$\prod_{j\in J}\Par^{\geq\,0}_{b\cdot m_j}(\nbigelambdazero,j)$.
For such a primitive subset $S$,
we put as follows:
\[
 r(S):=\max\bigl\{|\vecc|\,\big|\,\vecc\in S\bigr\}.
\]
Here we put $|\veca|=\sum a_i$.
We use an induction
on the number $r(S)$.

We have the following exact sequence for some
$S'\subset
 \prod_{j\in J}\Par^{\geq\,0}_{b\cdot m_j}(\nbigelambdazero,j)$:
\begin{equation}\label{eq;a12.2.100}
 0\lrarr
 \lefttop{J}\Vzero_{S'}\lefttop{I}\Vzero_{\vecb}
 \bigl(\lefttop{I}\nbigl\bigr)
 \lrarr
 \lefttop{J}\Vzero_{S}\lefttop{I}\Vzero_{\vecb}
 \bigl(\lefttop{I}\nbigl\bigr)
 \lrarr
 \bigoplus_{\veca\in S}
 \lefttop{I}\Gr^{\Vzero}_{\veca}
 \lefttop{J}\Vzero_{\vecb}\bigl(\lefttop{I}\nbigl\bigr)
 \lrarr 0.
\end{equation}
Then we have $r(S')<r(S)$.
Due to the hypothesis of the induction,
we may assume that $\lefttop{J}\Vzero_{S'}\lefttop{I}\Vzero_{\vecb}
 \bigl(\lefttop{I}\nbigl\bigr)$
is locally free and coherent.
The third term in the sequence (\ref{eq;a12.2.100})
is locally free due to Corollary \ref{cor;a11.26.16}.
Then we obtain the local freeness of
$\lefttop{J}\Vzero_{S}\lefttop{I}\Vzero_{\vecb}
 \bigl(\lefttop{I}\nbigl\bigr)$.
\hfill\qed

\begin{cor}
The $\nbigo_{\nbigd_J}$-module
 $\lefttop{I}\Vzero_{\vecb}(\lefttop{I}\nbigl)$
 is locally free of infinite rank.
\hfill\qed
\end{cor}

\subsubsection{Filtrations $\lefttop{i}\Vzero$ and the compatibility}

Let $K$ be a subset of $\nbar$
and $\vecb$ be an element of $\real^K$.
We put as follows:
\[
 \lefttop{K}\Vzero_{\vecb}\bigl(\lefttop{I}\nbigl\bigr)
:=\bigcup_{q_K(\veca)=\vecb}
 \lefttop{\nbar}\Vzero_{\veca}\bigl(\lefttop{I}\nbigl\bigr).
\]
In particular,
we put $\lefttop{i}\Vzero_b:=\lefttop{\{i\}}\Vzero_b$.

Let $S$ be a finite subset of $\real^K$.
We put as follows:
\[
 \lefttop{K}\Vzero_S\bigl(\lefttop{I}\nbigl\bigr)
:=\sum_{\veca\in S}
 \lefttop{K}\Vzero_{\veca}\bigl(\lefttop{I}\nbigl\bigr).
\]

\begin{lem}\mbox{{}} \label{lem;a11.26.11}
\begin{itemize}
\item
The projection $\pi:\lefttop{I}\Vzero_{<0}\lrarr\lefttop{I}\nbigl$
induces the surjection:
\[
 \lefttop{K}\Vzero_S\cap\lefttop{I}\Vzero_{<0}(\tildepsi_{t,u})
\lrarr
 \lefttop{K}\Vzero_S\bigl(\lefttop{I}\nbigl\bigr).
\]
\item
Let $f$ be a section of $\lefttop{I}\nbigl$.
Then $f$ is contained in
$\lefttop{K}\Vzero_S\bigl(\lefttop{I}\nbigl\bigr)$
if and only if
$q_K\bigl(\lefttop{\nbar}\Prim(f)\bigr)$
is contained in $\nbigs(S)$.
\hfill\qed
\end{itemize}
\end{lem}

We have the filtrations
$\bigl(\lefttop{j}\Vzero\,\big|\,j\in J\bigr)$
on $\lefttop{J}\Vzero_s\lefttop{I}\Vzero_{\vecb}\bigl(
 \lefttop{I}\nbigl\bigr)$
in the category of the vector bundles.

\begin{lem}
The tuple of the filtrations
$\bigl(\lefttop{j}\Vzero\,\big|\,j\in J\bigr)$
is compatible
in the sense of Definition {\rm \ref{df;b11.12.3}}.
\end{lem}
\pf
From the proof of Lemma \ref{lem;a11.26.10},
we have the equality:
\[
 \sum \rank \lefttop{J}\Gr^{\Vzero}_{\veca}
 \lefttop{I}\Vzero_{\vecb}\bigl(\lefttop{I}\nbigl\bigr)
=\rank \lefttop{J}\Vzero_S\lefttop{I}\Vzero_{\vecb}
 \bigl(\lefttop{I}\nbigl\bigr).
\]
From Lemma \ref{lem;a11.26.11},
we have the following:
\[
 \bigcap_{j\in J}\lefttop{j}\Vzero_{c_j}
 \lefttop{I}\Vzero_{\vecb}\bigl(\lefttop{I}\nbigl\bigr)
=\lefttop{J}\Vzero_{\vecc}\lefttop{I}\Vzero_{\vecb}
 \bigl(\lefttop{I}\nbigl\bigr).
\]
Then we obtain the compatibility
by using Lemma \ref{lem;10.12.31}.
\hfill\qed

\subsubsection{The actions of $g_{1,i}$ and the generic splitting}

On 
$\lefttop{J}\Vzero_S\lefttop{I}\Vzero_{\vecb}
 \bigl(\lefttop{I}\nbigl\bigr)$
and
$\lefttop{J}\Gr^{\Vzero}_{\vecc}\lefttop{I}\Vzero_{\vecb}
 \bigl(\lefttop{I}\nbigl\bigr)$,
we have the action of the tuple
$\bigl(g_{1,i}\,\big|\,i\in J\bigr)$,
given in the subsubsection \ref{subsubsection;a11.26.15}.

\begin{lem} \label{lem;a11.26.17}
The following endomorphism vanishes
on $\lefttop{J}\Gr^{\Vzero}_{\vecc}\lefttop{I}\Vzero_{\vecb}
 \bigl(\lefttop{I}\nbigl\bigr)$:
\[
 \prod_{u_1\in\nbigk(\nbige,i,\lambda_0,c_i+m_ib)}
 \bigl(g_{1,i}+\eigenmap(\lambda,u_1)\bigr)^N.
\]
\end{lem}
\pf
It follows from Lemma \ref{lem;9.22.21}
and Corollary \ref{cor;a11.26.16}.
\hfill\qed

\begin{lem} \label{lem;a11.26.20}
Let $\lambda_1\in\Delta(\lambda_0,\epsilon_0)$  be generic,
and $\epsilon_1$ be a positive number such that
$\Delta(\lambda_1,\epsilon_1)\subset \Delta(\lambda_0,\epsilon_0)$
and that any $\lambda\in\Delta(\lambda_1,\epsilon_1)$ is generic.
Let us consider the endomorphisms
of $\bigl(g_{1,i}\,\big|\,i\in J\bigr)$
on
$\lefttop{J}\Vzero_S\lefttop{I}\Vzero_{\vecb}
    \bigl(\lefttop{I}\nbigl\bigr)
 _{|\Delta(\lambda_1,\epsilon_1)}$.
Then the generalized eigen decomposition
for $\bigl(g_{1,i}\,\big|\,i\in J\bigr)$
gives the splitting of the filtrations
$\bigl(\lefttop{i}\Vzero\,\big|\,i\in J\bigr)$.
\end{lem}
\pf
It follows from Lemma \ref{lem;a11.26.17}.
\hfill\qed

\subsubsection{The weight filtration of $N=s+\eigenmap(\lambda,u)$}

We put $s:=t\cdot\deldel_t$
and $N:=s+\eigenmap(\lambda,u)$.
Since the eigenfunction of $s$ on
$\psi_{t,u}$ is $-\eigenmap(\lambda,u)$,
the induced actions of $N$
are nilpotent
on the vector bundles $\lefttop{I}\nbigl$,
$\lefttop{J}\Vzero_{S}\lefttop{I}\Vzero_{\vecb}
 \bigl(\lefttop{I}\nbigl\bigr)$
and
$\lefttop{J}\Gr^{\Vzero}_{\vecc}
 \lefttop{I}\Vzero_{\vecb}\bigl(\lefttop{I}\nbigl\bigr)$.

\begin{lem} \label{lem;a11.26.25}
The conjugacy classes of $N$
on 
$\lefttop{J}\Gr^{\Vzero}_{\vecc}
 \lefttop{I}\Vzero_{\vecb}\bigl(\lefttop{I}\nbigl\bigr)_{|(\lambda,P)}$
are independent of a choice of
$(\lambda,P)\in \nbigd_J(\lambda_0,\epsilon_0)$.
\end{lem}
\pf
Recall the descriptions in the subsubsection
\ref{subsubsection;a11.26.6}.
Then the independence of $\lambda$ follows from
the limiting mixed twistor theorem.
When we fix a generic $\lambda$,
we can show the independence of $P$
by using the normalizing frame.
Thus we obtain the independence of $(\lambda,P)$.
\hfill\qed

\begin{cor}
The weight filtration
$W\bigl(N,\lefttop{J}\Gr^{\Vzero}_{\vecc}\lefttop{I}\Vzero_{\vecb}
 \bigl(\lefttop{I}\nbigl\bigr)
 \bigr)$
is a filtration in the category of vector bundles.
\hfill\qed
\end{cor}

\begin{lem} \label{lem;a11.26.30}
Let us consider the action of $N$
on $\lefttop{J}\Vzero_{S}\lefttop{I}\Vzero_{\vecb}
 \bigl(\lefttop{I}\nbigl\bigr)$.
\begin{itemize}
\item
 The conjugacy classes of $N$
 on $\lefttop{J}\Vzero_{S}\lefttop{I}\Vzero_{\vecb}
 \bigl(\lefttop{I}\nbigl\bigr)_{|(\lambda,P)}$
 are independent of a choice of
 $(\lambda,P)\in\nbigd_J(\lambda_0,\epsilon_0)$.
\item
 The filtration
 $W\bigl(N,\lefttop{J}\Vzero_{S}\lefttop{I}\Vzero_{\vecb}
 \bigl(\lefttop{I}\nbigl\bigr)\bigr)$
 is a filtration in the category of vector bundles.
\item
 The projection
 $\lefttop{J}\Vzero_{S}\lefttop{I}\Vzero_{\vecb}\bigl(
 \lefttop{I}\nbigl\bigr)\lrarr
 \lefttop{J}\Gr^{\Vzero}_S\lefttop{I}\Vzero_{\vecb}\bigl(
 \lefttop{I}\nbigl
 \bigr)$ induces the surjective morphism:
\[
 W\bigl(N,
 \lefttop{J}\Vzero_{S}\lefttop{I}\Vzero_{\vecb}\bigl(
 \lefttop{I}\nbigl\bigr)
 \bigr)
 \lrarr
 W\bigl(N,
 \lefttop{J}\Gr^{\Vzero}_S\lefttop{I}\Vzero_{\vecb}\bigl(
 \lefttop{I}\nbigl
 \bigr)
 \bigr).
\]
\item
In the case $S\subset\nbigs(S')$,
we have the following equality:
\[
 W\bigl(
 N,\lefttop{J}\Vzero_S\lefttop{I}\Vzero_{\vecb}
 \bigl(\lefttop{I}\nbigl\bigr)
 \bigr)
= W\bigl(
 N,\lefttop{J}\Vzero_{S'}\lefttop{I}\Vzero_{\vecb}
 \bigl(\lefttop{I}\nbigl\bigr)
 \bigr)
\cap
 \lefttop{J}\Vzero_S\lefttop{I}\Vzero_{\vecb}
 \bigl(\lefttop{I}\nbigl\bigr).
\]
\end{itemize}
\end{lem}
\pf
Let us see the first claim.
The independence of $\lambda$ for a fixed $P$
can be shown by using the generic splitting
(Lemma \ref{lem;a11.26.20}),
Lemma \ref{lem;c11.12.20} and Lemma \ref{lem;a11.26.25}.
The independence of $P$ for a fixed generic $\lambda$
can be shown by using the normalizing frame.
Thus we obtain the first claim.

The other claims can be shown similarly.
\hfill\qed

\begin{lem} \label{lem;a12.2.110}
We have the following:
\begin{equation} \label{eq;a11.26.31}
 \sum_{\veca\in S}
 W_h\bigl(N,
 \lefttop{J}\Vzero_{\veca}\lefttop{I}\Vzero_{\vecb}(\lefttop{I}\nbigl)
 \bigr)
=W_h\bigl(N,
 \lefttop{J}\Vzero_S\lefttop{I}\Vzero_{\vecb}(\lefttop{I}\nbigl)
 \bigr).
\end{equation}
\end{lem}
\pf
We have only to check the claims for
primitive subsets $S$,
which are contained in
$\prod_{j\in J}\Par^{\geq\,0}_{b\cdot m_j}(\nbigelambdazero,j)$.
For such a primitive subset $S$,
we put as follows:
\[
 r(S):=\max\bigl\{|\vecc|\,\big|\,\vecc\in S\bigr\}.
\]
Here we put $|\veca|=\sum a_i$.
We use an induction
on the number $r(S)$.
The left hand side and the right hand side
of (\ref{eq;a11.26.31})
are denoted by $\nbiga$ and $\nbigb$ respectively.

Let us consider the exact sequence:
\[
\begin{CD}
 0@>>>
 \lefttop{J}\Vzero_{S'}\lefttop{I}\Vzero_{\vecb}\bigl(
 \lefttop{I}\nbigl
 \bigr)
@>>>
 \lefttop{J}\Vzero_S\lefttop{I}\Vzero_{\vecb}\bigl(
 \lefttop{I}\nbigl
 \bigr)
@>{\pi_S}>>
 \bigoplus_{\veca\in S}
 \lefttop{J}\Gr^{\Vzero}_{\veca}\lefttop{I}\Vzero_{\vecb}
 \bigl(\lefttop{I}\nbigl\bigr)
@>>> 0
\end{CD}
\]
Here we have $r(S')<r(S)$.

Due to Lemma \ref{lem;a11.26.30},
we have the following:
\[
 \pi_S(\nbiga)=\pi_S(\nbigb)
=\bigoplus_{\veca\in S}
 W_h\bigl(
 N,\lefttop{J}\Gr^{\Vzero}_{\veca}
 \lefttop{I}\Vzero_{\vecb}(\lefttop{I}\nbigl)
 \bigr).
\]
We also obtain the following, due to Lemma \ref{lem;a11.26.30}:
\begin{equation}\label{eq;a11.26.38}
 \ker\pi_S\cap\nbigb
=W_h\bigl(N,
 \lefttop{J}\Vzero_{S'}\lefttop{I}\Vzero_{\vecb}
 \bigl(\lefttop{I}\nbigl\bigr)
 \bigr).
\end{equation}
We have the following for $\ker\pi_S\cap\nbiga$:
\begin{multline}\label{eq;a11.26.35}
 \ker \pi_S\cap\nbiga
=\sum_{\veca\in S}\Bigl(
 W_h\bigl(N, \lefttop{J}\Vzero_{\veca}\lefttop{I}\Vzero_{\vecb}
 \bigl(\lefttop{I}\nbigl\bigr)
 \bigr)
 \cap
 \lefttop{J}\Vzero_{S'}\lefttop{I}\Vzero_{\vecb}
 \bigl(\lefttop{I}\nbigl\bigr)
 \Bigr) \\
=\sum_{\veca\in S}\Bigl(
 W_h\bigl(N, \lefttop{I}\Vzero_{\vecb}(\lefttop{I}\nbigl)\bigr)
\cap
\lefttop{J}\Vzero_{\veca}\lefttop{I}\Vzero_{\vecb}
 \bigl(\lefttop{I}\nbigl\bigr)
 \cap
 \lefttop{J}\Vzero_{S'}\lefttop{I}\Vzero_{\vecb}
 \bigl(\lefttop{I}\nbigl\bigr)
 \Bigr).
\end{multline}
The following lemma can be easily shown.
\begin{lem}
We have the subset
$S(\veca) \subset
 \prod_j\Par_{b\cdot m_j}^{\geq\,0}(\nbigelambdazero,j)$
satisfying the following:
\[
 \nbigs(S(\veca))
\cap\prod_{j\in J}\Par^{\geq\,0}_{b\cdot m_j}\bigl(
 \nbigelambdazero,j
 \bigr)
 =\nbigs(\veca)\cap \nbigs(S')
 \cap\prod_{j\in J}\Par^{\geq\,0}_{b\cdot m_j}\bigl(
 \nbigelambdazero,j
 \bigr).
\]
We also have the following:
\[
 \nbigs\Bigl(
 \bigcup_{\veca\in S}S(\veca')
 \Bigr)\cap\prod_{j\in J}\Par^{\geq\,0}_{b\cdot m_j}\bigl(
 \nbigelambdazero,j
 \bigr)
=\nbigs(S')\cap\prod_{j\in J}\Par^{\geq\,0}_{b\cdot m_j}\bigl(
 \nbigelambdazero,j
 \bigr).
\]
\hfill\qed
\end{lem}

Then the right hand side of (\ref{eq;a11.26.35})
can be rewritten as follows:
\begin{equation}\label{eq;a11.26.36}
 \sum_{\veca\in S}
 \Bigl(
 W_h(N,\lefttop{I}\Vzero_{\vecb}(\lefttop{I}\nbigl))
\cap
 \lefttop{J}\Vzero_{S(\veca)}\lefttop{I}\Vzero_{\vecb}
 \bigl(\lefttop{I}\nbigl\bigr)
 \Bigr).
\end{equation}
Since we have $r(S(\veca))<r(S)$,
(\ref{eq;a11.26.36}) can be rewritten as follows:
\begin{equation}\label{eq;a11.26.37}
\sum_{\veca\in S}
\sum_{\vecc\in S(\veca)}
 W_h\bigl(N,\lefttop{J}\Vzero_{\vecc}\lefttop{I}\Vzero_{\vecb}
 \bigl(\lefttop{I}\nbigl\bigr)\bigr)
=\sum_{\vecc\in S'}
 W_h\bigl(N,\lefttop{J}\Vzero_{\vecc}\lefttop{I}\Vzero_{\vecb}
 \bigl(\lefttop{I}\nbigl\bigr)\bigr).
\end{equation}
Since $r(S')<r(S)$,
the right hand sides of (\ref{eq;a11.26.37})
and (\ref{eq;a11.26.38}) are same.
Then we obtain
$\ker \pi_S\cap\nbiga=\ker\pi_S\cap\nbigb$.
Thus we obtain $\nbiga=\nbigb$.
\hfill\qed

%% file: a80.2.tex
\subsubsection{The decomposition}

We reformulate the result
in the subsubsection \ref{subsubsection;a11.26.40}.

Let $S$ be a primitive subset of $\real^n$.
Let $i$ be an element of $\nbar$.
Let $I$ be a subset of $\nbar$ such that $i\in I$.
We put $J':=(\nbar-I)\cup\{i\}$.
We have the $\nbigo_{|\nbigd_{J'}}$-locally free sheaf
$\lefttop{i}\Gr^{\Vzero}_{-1}\lefttop{\nbar}\Vzero_{S}
 \bigl(\lefttop{I}\nbigl\bigr)$.
We also have the following direct summand:
\[
 \lefttop{i}\psi_{-\vecdelta_0}
 \lefttop{\nbar}\Vzero_S\bigl(\lefttop{I}\nbigl\bigr)
:=\EE\bigl(z_i\deldel_i,-1\bigr).
\]
We put $I':=I-\{i\}$. Note we have $J'\sqcup I'=\nbar$.
We put $S':=\bigl\{\veca+\vecdelta_i\,\big|\,\veca\in S\bigr\}$.
Then we have
$\nbigo_{\nbigd_J'}$-locally free sheaf
$\lefttop{i}\Gr_0\lefttop{\nbar}\Vzero_{S'}(\lefttop{I'}\nbigl)$.
We also have the following direct summand:
\[
 \lefttop{i}\psi_0\lefttop{\nbar}\Vzero_{S'}
\bigl(\lefttop{I'}\nbigl\bigr):=
\EE\bigl(z_i\deldel_i,0\bigr).
\]

The multiplication $z_i$ induces
the morphism:
\[
 \var_i:
 P_h\Gr^{W(N)}_h\lefttop{i}\psi_{0}\lefttop{\nbar}\Vzero_{S'}
 \bigl(\lefttop{I'}\nbigl\bigr)
\lrarr
 P_h\Gr^{W(N)}_h\lefttop{i}\psi_{-1}\lefttop{\nbar}\Vzero_S
 \bigl(\lefttop{I}\nbigl\bigr).
\]
The multiplication of $\deldel_i$ induces the morphism:
\[
 \can_i:
 P_h\Gr^{W(N)}_h \lefttop{i}\psi_{-1}\lefttop{\nbar}\Vzero_S
 \bigl(\lefttop{I}\nbigl\bigr)
\lrarr
 P_h\Gr^{W(N)}_h\lefttop{i}\psi_{0}\lefttop{\nbar}\Vzero_{S'}
 \bigl(\lefttop{I'}\nbigl\bigr).
\]

\begin{lem}
We have the decomposition:
\[
 P_h\Gr^{W(N)}_h\lefttop{i}\psi_{0}\lefttop{\nbar}\Vzero_{S'}
 \bigl(\lefttop{I'}\nbigl\bigr)
=
 \Image(\can_i)\oplus\Ker(\var_i).
\]
\end{lem}
\pf
When we take the $\lefttop{J'}\Gr^{\Vzero}$,
the result follows from Lemma \ref{lem;a11.26.45}.
Then we obtain the result by using
the generic splitting (Lemma \ref{lem;a11.26.20})
and Lemma \ref{lem;9.23.15}.
\hfill\qed

\begin{cor}\label{cor;c11.26.20}
We have the decomposition:
\[
 P_h\Gr^{W(N)}_h\lefttop{i}\psi_0(\lefttop{I}\nbigl)
=\Image(\can_i)\oplus \Ker\var_i.
\]
In particular,
$\Image(\can_i)$ and $\Ker(\var_i)$
are subbundles of $P_h\Gr^{W(N)}_h\lefttop{i}\psi_0(\lefttop{I}\nbigl)$
on $\nbigd_{J'}$.
\hfill\qed
\end{cor}

\begin{cor}\label{cor;b12.3.50}
$\Image(\var_i)$ is a subbundle
of $P_h\Gr^{W(N)}_h\lefttop{i}\psi_{-\vecdelta_0}(\lefttop{I}\nbigl)$
on $\nbigd_{J'}$.
\hfill\qed
\end{cor}

%% file: a83.2.tex

\subsubsection{Preliminary}

Let $\nbiga$ be an abelian category.
Let $C$ be an object of $\nbiga$,
and let $f:C\lrarr C$ be a nilpotent endomorphism of $C$.
Recall that we obtain the weight filtration
$W(f)$, characterized by the following conditions
(see (1.6) in \cite{d3}):
\begin{itemize}
\item
 $f\bigl(W_h(f)\bigr)\subset W_{h-2}(f)$.
\item
 The induced morphism
 $f^h:\Gr^{W(f)}_h\lrarr \Gr^{W(f)}_{-h}$ is isomorphic 
for any $h\geq 0$.
\end{itemize}

The weight filtration has a functoriality.
\begin{lem}
Let $C_i$ $(i=1,2)$ be objects of $\nbiga$.
Let $f_i$ be nilpotent endomorphisms of $C_i$.
Let $\phi:C_1\lrarr C_2$ be a morphism
such that the following diagramm is commutative:
\[
 \begin{CD}
 C_1 @>{\phi}>> C_2 \\
 @V{f_1}VV @V{f_2}VV\\
 C_1 @>{\phi}>> C_2.
 \end{CD}
\]
Then the morphism $\phi$  preserves the weight filtrations,
i.e.,
$\phi\bigl(W_h(f_1)\bigr)\subset \phi\bigl(W_h(f_2)\bigr)$
\end{lem}
\pf
Assume that $f_1^{d+1}=f_2^{d+1}=0$.
Recall that we have
$W_d(f_i)=C_i$, $W_{d-1}(f_i)=\ker(f_i^{d})$,
$W_{-d}=\Image(f_i^d)$ and $W_{-d-1}=0$.
Then it is easy to see that
$W_d$ and $W_{-d}$ are preserved by $\phi$.
Then it is easy to derive that
$\phi$ preserves $W_h$  for any $h$.
(See the way of the recursive construction of the weight filtration
in (1.6) in \cite{d3}.)
\hfill\qed

\vspace{.1in}
Let $0\lrarr C_1\stackrel{a}{\lrarr} C_2\stackrel{b}{\lrarr} C_3\lrarr 0$
be an exact sequence in the abelian category $\nbiga$.
Let $f_i$ be nilpotent endomorphisms of $C_i$.
Let $W(C_i)$ be the filtrations of $C_i$
which are preserved by the morphisms $a$ and $b$.
Moreover we assume that
the induced sequence
$0\lrarr\Gr^{W}(C_1)\lrarr \Gr^{W}(C_2)\lrarr \Gr^{W}(C_3)\lrarr 0$
is exact,
i.e.,
the morphisms $a$ and $b$ are strict
with respect to the filtrations.

\begin{lem}
Let $i$ be either $1$ or $3$.
We put $S:=\{1,2,3\}-\{i\}$.
Assume that $W(C_j)$ are the weight filtrations of $f_j$
for $j\in S$.
Then $W(C_i)$ is also the weight filtration of $f_i$
 for $i$.
\end{lem}
\pf
We have only to show that $W(C_i)$ satisfies the axioms of
the weight filtration.
Since $f_2\bigl( W_h(C_2)\bigr)\subset W_{h-2}(C_2)$,
we obtain $f_i\bigl(W_h(C_i)\bigr)\subset W_{h-2}(C_i)$.
Due to the exact sequence
$0\lrarr\Gr^{W}(C_1)\lrarr \Gr^{W}(C_2)\lrarr \Gr^{W}(C_3)\lrarr 0$,
we obtain that the morphisms
$f_i^h:\Gr^W_h(C_i)\lrarr \Gr^{W}_{-h}(C_i)$
are isomorphic for any $h$.
\hfill\qed

\vspace{.1in}
Let $0\lrarr C_1\stackrel{a}{\lrarr} C_2\stackrel{b}{\lrarr} C_3\lrarr 0$
be an exact sequence in the abelian category $\nbiga$.
Let $f_i$ be nilpotent endomorphisms of $C_i$.

\begin{cor}\mbox{{}}\label{cor;a12.2.120}
\begin{itemize}
\item
If the morphism $a$ is strict
with respect to the filtrations
$W(f_1)$ and $W(f_2)$,
then $b$ is also strict with respect to
the filtrations $W(f_2)$ and $W(f_3)$.
\item
If the morphism $b$ is strict
with respect to the filtrations
$W(f_2)$ and $W(f_3)$,
then $a$ is also strict with respect to
the filtrations $W(f_1)$ and $W(f_2)$.
\hfill\qed
\end{itemize}
\end{cor}

%% file: a80.3.tex

\subsubsection{Definition of $\FFzero$ and
the decomposition of $\Gr^{\FFzero}$}

We introduce the filtration $\FFzero$ of $\psizero_{t,u}$
as follows:
\[
 \FFzero_m:=
 \sum_{|I|\geq l-m}\lefttop{I}\Vzero_{<0}(\psi_{t,u}),
\quad\quad
0=\FFzero_0\subset\FFzero_1\subset \cdots\subset \FFzero_l=\psizero_{t,u}.
\]
We put $\GGzero_m:=\psi_{t,u}/\FFzero_m$

\begin{lem} \label{lem;9.20.16}
We have $\Gr^{\FFzero}_m=\bigoplus_{|I|=l-m}\lefttop{I}\nbigl$.
\end{lem}
\pf
We have the naturally defined surjective morphism
$\bigoplus_{|I|=l-m}\lefttop{I}\nbigl\lrarr
 \Gr^{\FFzero}_m$.
Let $\bigl(\lefttop{I}f\,\big|\,|I|=l-m\bigr)$ be an element of
$ \bigoplus_{|I|=l-m}\lefttop{I}\nbigl$
such that
$\sum \lefttop{I}f\in\FFzero_{m-1}$.
Then we have the following:
\[
 \lefttop{I}f\in
 \Bigl(
 \sum_{\substack{
 |I'|\geq l-m,\\ I'\neq I
 }}\lefttop{I'}\Vzero_{<0}(\psi_{t,u})
 \Bigr)
\cap \lefttop{I}\Vzero_{<0}(\psi_{t,u})
=\sum_{\substack{
 |I'|\geq l-m,\\
 I'\neq I
 }}\lefttop{I\cup I'}\Vzero_{<0}(\psi_{t,u}).
\]
If $|I'|\geq l-m$ and $I'\neq I$,
then $|I\cup I'|\geq l-m+1$ and $I\cup I'\neq I$.
Hence we obtain $\lefttop{I}f=0$.
Thus the lemma is proved.
\hfill\qed

\subsubsection{Compatibility of the weight filtration $(I)$}

We put as follows:
\[
 \lefttop{\nbar}\Vzero_S\bigl(\Gr^{\FFzero}_m\bigr)
:=\Image\bigl(
 \lefttop{\nbar}\Vzero_S\cap \FFzero_m\lrarr\Gr^{\FFzero}_m
 \bigr),
\quad\quad
 \lefttop{\nbar}\Vzero_S\bigl(\GGzero_m\bigr)
:=\Image\bigl(
 \lefttop{\nbar}\Vzero_S\lrarr \GGzero_m
 \bigr).
\]
Then we obtain the following exact sequence:
\begin{equation}
 \begin{CD}
 0 @>>>
 \lefttop{\nbar}\Vzero_S
 \bigl(\Gr^{\FFzero}_m\bigr)
 @>>>
 \lefttop{\nbar}\Vzero_S\bigl(
 \GGzero_{m-1}
 \bigr)
 @>>>
 \lefttop{\nbar}\Vzero_S\bigl(
 \GGzero_{m}
 \bigr)
 @>>> 0.
 \end{CD}
\end{equation}

\begin{lem}\label{lem;a12.2.130}
We have the following:
\[
 W_h\bigl(N,\lefttop{\nbar}\Vzero_S(\GGzero_{m-1})
 \bigr)
\cap
 \lefttop{\nbar}\Vzero_{S}\bigl(
 \Gr^{\FFzero}_{m}\bigr)
=W_h\bigl(N,\lefttop{\nbar}\Vzero_S\bigl(\Gr^{\FFzero}_m\bigr)\bigr).
\]
\end{lem}
\pf
Due to Lemma \ref{lem;9.20.16},
we have only to show the following,
for $|I|=n-m$:
\begin{equation}\label{eq;a11.26.50}
 W_h\bigl(N,\lefttop{\nbar}\Vzero_S(\GGzero_{m-1})\bigr)
\cap
\lefttop{\nbar}\Vzero_S(\lefttop{I}\nbigl)
=
 W_h\bigl(
 N,\lefttop{\nbar}\Vzero_S(\lefttop{I}\nbigl)
 \bigr).
\end{equation}
The implication $\supset$ is clear from the functoriality
of the weight filtration.
Let us show the implication $\subset$.
When we restrict the both sides of (\ref{eq;a11.26.50})
to $\nbigd_I-\bigcup_{I'\supsetneq I}\nbigd_{I'}$,
then the equality in (\ref{eq;a11.26.50}) holds.
Since the right hand side is a subbundle
of $\lefttop{\nbar}\Vzero_S(\lefttop{I}\nbigl)$
(Lemma \ref{lem;a11.26.30}),
we obtain the implication $\subset$.
\hfill\qed

\vspace{.1in}

We put as follows:
\[
 W_h'\bigl(N,\lefttop{\nbar}\Vzero_S(\GGzero_m)\bigr)
:=\Image\Bigl(
 W_h\bigl(
 \lefttop{\nbar}\Vzero_S(\GGzero_{m-1})\bigr)
\lrarr
 \lefttop{\nbar}\Vzero_S(\GGzero_m)
 \Bigr).
\]
\begin{lem}\label{lem;a11.26.80}
We have
$W_h(N,\lefttop{\nbar}\Vzero_S(\GGzero_{m}))
=W'_h(N,\lefttop{\nbar}\Vzero_S(\GGzero_{m}))$
\end{lem}
\pf
It follows from Lemma \ref{lem;a12.2.130}
and Corollary \ref{cor;a12.2.120}.
\hfill\qed

\begin{cor}\label{cor;a11.26.81}
The projection $\psi_{t,u}\lrarr \GGzero_m$
induces the surjection:
\[
 W_h\bigl(N,\lefttop{\nbar}\Vzero_S(\psi_{t,u})\bigr)
\lrarr
 W_h\bigl(N,\lefttop{\nbar}\Vzero_S(\GGzero_m)\bigr).
\]
\hfill\qed
\end{cor}

\begin{cor}
We have the following:
\[
 W_h\bigl(N,\lefttop{\nbar}\Vzero_S(\psizero_{t,u})\bigr)
 \cap\FFzero_m
=W_h\bigl(N,\Vzero_S(\FFzero_m)\bigr).
\]
\end{cor}
\pf
It follows from Lemma \ref{lem;a11.26.80}
and Corollary \ref{cor;a12.2.120}.
\hfill\qed

\subsubsection{The compatibility of the weight filtration $(II)$}

Let $S$ and $S'$ be primitive subsets
such that $S\subset\nbigs(S')$.
Let us consider the following diagramm:
\begin{equation}
\begin{CD}
 0 @>>>
 \lefttop{\nbar}\Vzero_{S'}\bigl(\Gr^{\FFzero}_m\bigr)
 @>>>
 \lefttop{\nbar}\Vzero_{S'}\bigl(\GGzero_{m-1}\bigr)
 @>>>
 \lefttop{\nbar}\Vzero_{S'}\bigl(\GGzero_{m}\bigr)
 @>>> 0\\
 @. @AAA @AAA @AAA @.\\
 0 @>>>
 \lefttop{\nbar}\Vzero_{S}\bigl(\Gr^{\FFzero}_m\bigr)
 @>>>
 \lefttop{\nbar}\Vzero_{S}\bigl(\GGzero_{m-1}\bigr)
 @>{\pi_1}>>
 \lefttop{\nbar}\Vzero_{S}\bigl(\GGzero_{m}\bigr)
 @>>> 0\\
\end{CD}
\end{equation}
The vertical arrows are injective.

\begin{lem}\label{lem;b11.26.21}
We have the following:
\[
 W_h\bigl(N,\lefttop{\nbar}\Vzero_{S'}(\GGzero_m)\bigr)\cap
 \Vzero_{S}(\GGzero_m)
=W_h\bigl(N,\lefttop{\nbar}\Vzero_{S}(\GGzero_m)
 \bigr).
\]
\end{lem}
\pf
We use a descending induction on $m$.
The claim in the case $m=n$ is trivial.
We assume that the claim holds for $m$,
and we will show the claim for $m-1$.

The implication $\supset$ follows from the functoriality
of the weight filtration.
Let us show the implication $\subset$.
We have only to show the following two equalities:
\begin{equation}\label{eq;b11.26.1}
 \Bigl[
W_h\bigl(N,\lefttop{\nbar}\Vzero_{S'}(\GGzero_{m-1})\bigr)
 \cap \lefttop{\nbar}\Vzero_{S}(\GGzero_{m-1})
 \Bigr]
 \cap \lefttop{\nbar}\Vzero_{S}(\Gr^{\FFzero}_m)
=W_h(N,\lefttop{\nbar}\Vzero_{S'}(\GGzero_{m}))
 \cap \lefttop{\nbar}\Vzero_{S}(\Gr^{\FFzero}_m).
\end{equation}
\begin{equation}\label{eq;b11.26.2}
 \pi_1\Bigl(
 W_h(N,\lefttop{\nbar}\Vzero_{S'}(\GGzero_{m-1}))
 \cap \lefttop{\nbar}\Vzero_{S}(\GGzero_{m-1})
 \Bigr)
\subset\pi_1\Bigl(
 W_h(N,\lefttop{\nbar}\Vzero_{S'}(\GGzero_{m}))
 \Bigr).
\end{equation}
The equality (\ref{eq;b11.26.1}) can be shown as follows:
\begin{multline}
 \Bigl[
W_h\bigl(N,\lefttop{\nbar}\Vzero_{S'}(\GGzero_{m-1})\bigr)
 \cap \lefttop{\nbar}\Vzero_{S}(\GGzero_{m-1})
 \Bigr]
 \cap \lefttop{\nbar}\Vzero_{S}(\Gr^{\FFzero}_m)\\
=
 W_h\bigl(N,\lefttop{\nbar}\Vzero_{S'}(\GGzero_{m-1})\bigr)
 \cap \lefttop{\nbar}\Vzero_{S'}(\Gr^{\FFzero}_{m})
 \cap \lefttop{\nbar}\Vzero_{S}(\Gr^{\FFzero}_m) \\
=W_h(N,\lefttop{\nbar}\Vzero_{S'}(\Gr^{\FFzero}_m))
 \cap \lefttop{\nbar}\Vzero_{S}(\Gr^{\FFzero}_m)
=W_h(N,\lefttop{\nbar}\Vzero_{S}(\Gr^{\FFzero}_m)) \\
=W_h(N,\lefttop{\nbar}\Vzero_{S'}(\GGzero_{m}))
 \cap \lefttop{\nbar}\Vzero_{S}(\Gr^{\FFzero}_m).
\end{multline}

The implication (\ref{eq;b11.26.2}) can be shown as follows:
\begin{multline}
 \pi_1\Bigl(
 W_h\bigl(N,\lefttop{\nbar}\Vzero_{S'}(\GGzero_{m-1})\bigr)
 \cap \lefttop{\nbar}\Vzero_{S}(\GGzero_{m-1})
 \Bigr)
\subset
 \pi_1\Bigl(
 W_h\bigl(N,\lefttop{\nbar}\Vzero_{S'}(\GGzero_{m-1})\bigr)
 \Bigr)
 \cap 
 \pi_1\Bigl(
 \lefttop{\nbar}\Vzero_{S}(\GGzero_{m-1})
 \Bigr) \\
=W_h\bigl(N,\lefttop{\nbar}\Vzero_{S'}(\GGzero_{m})\bigr)
\cap
 \lefttop{\nbar}\Vzero_{S}(\GGzero_{m})
=W_h(N,\lefttop{\nbar}\Vzero_{S}(\GGzero_{m})).
\end{multline}
Thus we are done.
\hfill\qed

\begin{cor}\label{cor;a12.3.10}
We have the following:
\[
 W_h\bigl(N,\lefttop{\nbar}\Vzero_{S'}(\psi_{t,u})\bigr)\cap
 \lefttop{\nbar}\Vzero_{S}
=W_h\bigl(N,\lefttop{\nbar}\Vzero_{S}(\psi_{t,u})\bigr).
\]
In particular, we have the following:
\[
 W_h(N,\psi_{t,u})\cap
 \lefttop{\nbar}\Vzero_{S}
=W_h\bigl(N,\lefttop{\nbar}\Vzero_{S}(\psi_{t,u})\bigr).
\]
\hfill\qed
\end{cor}

\begin{lem} \label{lem;b11.26.20}
We have the following:
\[
 \sum_{\veca\in S}W_h\bigl(N,\lefttop{\nbar}\Vzero_{\veca}(\GGzero_m)\bigr)
=W_h\bigl(
 N,\lefttop{\nbar}\Vzero_{S}(\GGzero_m)
 \bigr).
\]
\end{lem}
\pf
We use a descending induction on $m$.
In the case $m=n$, the claim is trivial.
We assume that the claim for $m$ holds,
and we will prove the claim for $m-1$.
We use the following exact sequence:
\[
 \begin{CD}
 0 @>>> \Gr^{\FFzero}_m
 @>>> \GGzero_{m-1} @>{\pi_2}>> \GGzero_m @>>> 0.
 \end{CD}
\]

We have only to show the following two claims:
\begin{equation}\label{eq;b11.26.5}
 \sum_{\veca\in S}W_h\bigl(N,\lefttop{\nbar}\Vzero_{\veca}(\GGzero_m)\bigr)
\cap \Gr^{\FFzero}_m
=W_h\bigl(
 N,\lefttop{\nbar}\Vzero_{S}(\GGzero_m)
 \bigr)\cap \Gr^{\FFzero}_m.
\end{equation}
\begin{equation}\label{eq;b11.26.6}
\pi_2\Bigl(
 \sum_{\veca\in S}
 W_h\bigl(N,\lefttop{\nbar}\Vzero_{\veca}(\GGzero_m)\bigr)
 \Bigr)
=\pi_2\Bigl(
 W_h\bigl(
 N,\lefttop{\nbar}\Vzero_{S}(\GGzero_m)
 \bigr)
 \Bigr).
\end{equation}

Let us show (\ref{eq;b11.26.5}).
We have the following:
\begin{multline} \label{eq;b11.26.7}
 \sum_{\veca\in S} W_h(N,\Vzero_{\veca}\Gr^{\FFzero}_m)
\subset
 \sum_{\veca\in S}W_h\bigl(N,\lefttop{\nbar}\Vzero_{\veca}(\GGzero_m)\bigr)
\cap \Gr^{\FFzero}_m \\
\subset
 W_h\bigl(
 N,\lefttop{\nbar}\Vzero_{S}(\GGzero_m)
 \bigr)\cap \Gr^{\FFzero}_m
=W_h\bigl(
 N,\lefttop{\nbar}\Vzero_S(\Gr^{\FFzero}_m)
 \bigr).
\end{multline}
Due to Lemma \ref{lem;a12.2.110}, we have already known
$\sum_{\veca\in S} W_h(N,\Vzero_{\veca}\Gr^{\FFzero}_m)
=W_h\bigl(
 N,\lefttop{\nbar}\Vzero_S(\Gr^{\FFzero}_m)
 \bigr)$,
we obtain (\ref{eq;b11.26.5}).

The equality (\ref{eq;b11.26.6}) can be shown as follows:
\begin{multline}
 \pi_2\Bigl(
 \sum_{\veca\in S}W_h(N,\lefttop{\nbar}\Vzero_{\veca}(\GGzero_m))
 \Bigr)
=\sum_{\veca\in S}W_h(N,\lefttop{\nbar}\Vzero_{\veca}(\GGzero_{m-1}))
=W_h\bigl(
 N,\lefttop{\nbar}\Vzero_{S}(\GGzero_{m-1})
 \bigr) \\
=\pi_2\Bigl(
 W_h\bigl(
 N,\lefttop{\nbar}\Vzero_{S}(\GGzero_m)
 \bigr)
 \Bigr).
\end{multline}
Hence we are done.
\hfill\qed

\begin{cor}
We have the following:
\[
 W_h\bigl(N,\lefttop{\nbar}\Vzero_S(\psi_{t,u})\bigr)
=\sum_{\veca\in S}
 W_h\bigl(N,\lefttop{\nbar}\Vzero_{\veca}(\psi_{t,u})\bigr).
\]
\hfill\qed
\end{cor}

%% file: a80.4.tex

\begin{lem} \label{lem;9.20.23}
We have the following:
\[
 W_h\bigl(N,\lefttop{I}\Vzero_S(\psi_{t,u})\bigr)
=W_h(N,\psi_{t,u})\cap \lefttop{I}\Vzero_S(\psi_{t,u})
=\sum_{\veca\in S}\lefttop{I}\Vzero_{\veca}(\psi_{t,u})
 \cap W_h(N,\psi_{t,u}).
\]
\end{lem}
\pf
It easily follows from Lemma \ref{lem;b11.26.21}
and Lemma \ref{lem;b11.26.20}.
\hfill\qed

\begin{lem}
We have the following:
\[
 W_h(N,\psi_{t,u})\cap
\Bigl(
 \lefttop{I}\Vzero_S(\psi_{t,u})
+\lefttop{i}\Vzero_{\vecc}(\psi_{t,u})
\Bigr)
=\lefttop{I}\Vzero_{S}(\psi_{t,u})\cap W_h(N)
+\lefttop{i}\Vzero_{\vecc}\cap W_h(N).
\]
\end{lem}
\pf
Similar to Lemma \ref{lem;9.20.23}.
\hfill\qed

\vspace{.1in}

\begin{lem}\label{lem;a12.3.2}
Let $S$ and $S'$ be primitive subsets of $\real^n$
such that $S\subset\nbigs(S')$.
The following sequence is exact:
\[
 \begin{CD}
 0 @>>>
 W_{h}\bigl(N,\lefttop{\nbar}\Vzero_{S}(\psi_{t,u})\bigr)
 @>>>
 W_{h}\bigl(N,\lefttop{\nbar}\Vzero_{S'}(\psi_{t,u})\bigr)
 @>>>
 W_{h}\bigl(
 N,\lefttop{\nbar}\Vzero_{S'}\big/\lefttop{\nbar}\Vzero_S
 \bigr)
 @>>>0.
 \end{CD}
\]
\end{lem}
\pf
We have the exact sequence:
\begin{equation}\label{eq;a12.3.1}
 \begin{CD}
 0 @>>>
 \lefttop{\nbar}\Vzero_S(\psi_{t,u})
 @>{a}>>
 \lefttop{\nbar}\Vzero_{S'}(\psi_{t,u})
 @>>>
 \lefttop{\nbar}\Vzero_{S'}\big/
 \lefttop{\nbar}\Vzero_{S}
 @>>> 0.
 \end{CD}
\end{equation}
Due to Lemma \ref{lem;9.20.23},
the morphism $a$ in (\ref{eq;a12.3.1})
is strict with respect to the weight filtrations.
Hence the sequence (\ref{eq;a12.3.1}) is strict
with respect to the weight filtrations
due to Corollary \ref{cor;a12.2.120}.
Thus we are done.
\hfill\qed

%% file: a80.5.tex

\subsubsection{The filtrations on $\Gr^{W(N)}(\psi_{t,u})$}

We put as follows, for any $\vecc\in\real^I$:
\[
 \lefttop{I}V_{\vecc}^{(\lambda_0)}(\Gr^{W(N)}_h(\psi_{t,u}))=
\Image\Bigl(
 \lefttop{I}\Vzero_{\vecc}\cap W_h(N,\psi_{t,u})
\lrarr \Gr^{W(N)}_{h}(\psi_{t,u})
 \Bigr).
\]
We have the following isomorphism:
\[
 \lefttop{I}\Vzero_{\vecc}(\Gr^{W(N)}_h(\psi_{t,u}))
\simeq
 \frac{\lefttop{I}\Vzero_{\vecc}\cap W_h(N,\psi_{t,u})}
 {\lefttop{I}\Vzero_{\vecc}\cap W_{h-1}(N,\psi_{t,u})}
\simeq
 \Gr^{W(N)}_{h}(\lefttop{I}\Vzero_{\vecc}).
\]
We put
$\lefttop{I}\Vzero_S\Gr^{W(N)}_h(\psi_{t,u}):=
\sum_{\veca\in S}\lefttop{I}\Vzero_{\veca}\Gr^{W(N)}_h(\psi_{t,u})$.

\begin{lem}
The projection $W_h(N,\psi_{t,u})\lrarr \Gr^{W(N)}_{h}(\psi_{t,u})$
induces the surjection:
\[
\lefttop{I}V^{(\lambda_0)}_{S}(\psi_{t,u})
\lrarr
\lefttop{I}V^{(\lambda_0)}_{S}(\Gr^{W(N)}(\psi_{t,u})).
\]
\end{lem}
\pf
It follows from Lemma \ref{lem;9.20.23}.
\hfill\qed

\begin{lem}\label{lem;b11.26.50}
Let $I$ be a subset of $\nbar$ and $i$ be an element of
$\nbar-I$.
We have the following:
\[
 \Bigl[
 \lefttop{I}\Vzero_S\cap\lefttop{i}\Vzero_c
 \Bigr]
 \bigl(
 \Gr^{W(N)}(\psi_{t,u})\bigr)
=\sum_{\veca\in S}\lefttop{I\sqcup\{i\}}\Vzero_{(\veca,c)}
 \bigl(
 \Gr^{W(N)}(\psi_{t,u})\bigr).
\]
\end{lem}
\pf
The implication $\supset$ is clear.
We show the implication $\subset$.
For any section
$f\in\lefttop{I}\Vzero_S\Gr^{W(N)}(\psi_{t,u})
    \cap \lefttop{i}\Vzero_c\Gr^{W(N)}(\psi_{t,u})$,
we pick sections
$ f_1\in \lefttop{I}\Vzero_S(\psi_{t,u})\cap W_h(N,\psi_{t,u})$
and $ f_2\in\lefttop{i}\Vzero_c(\psi_{t,u})\cap W_h(N,\psi_{t,u})$
such that $\rho_h(f_1)=\rho_h(f_2)=f$.
Here $\rho_h$ denotes the projection
$W_h(N)\lrarr \Gr^{W(N)}_h$.

Then we have the following:
\begin{multline}
 f_1-f_2\in W_{h-1}(N,\psi_{t,u})\cap
 \bigl(
 \lefttop{I}\Vzero_S(\psi_{t,u})+\lefttop{i}\Vzero_c(\psi_{t,u})
 \bigr) \\
=\bigl(W_{h-1}(N,\psi_{t,u})\cap \lefttop{I}\Vzero_S(\psi_{t,u})\bigr)
+\bigl(W_{h-1}(N,\psi_{t,u})\cap\lefttop{i}V_{c}(\psi_{t,u})\bigr).
\end{multline}
Hence we obtain $f_1-f_2=g_1-g_2$
for some sections
$g_1\in\lefttop{I}\Vzero_S(\psi_{t,u})\cap W_{h-1}(N,\psi_{t,u})$
and $g_2\in\lefttop{i}\Vzero_c(\psi_{t,u})\cap W_{h-1}(N,\psi_{t,u})$.
Then we obtain the following:
\begin{multline}
 f_1-g_1=f_2-g_2\in
\lefttop{I}\Vzero_S(\psi_{t,u})\cap W_h(N,\psi_{t,u})
 \cap \lefttop{i}\Vzero_c(\psi_{t,u})\cap W_h(N,\psi_{t,u})\\
=\lefttop{I}\Vzero_S(\psi_{t,u})\cap
  \lefttop{i}\Vzero_c(\psi_{t,u})\cap W_h(N,\psi_{t,u}).
\end{multline}
Then the implication $\subset$ follows easily.
\hfill\qed

\begin{cor}
Let $I$ be a subset of $\nbar$
and $\veca$ be an element of $\real^I$.
Let $i$ be an element of $\nbar-I$
and $c$ be a real number.
We have the isomorphism:
$\lefttop{i}\Gr^{\Vzero}_{c}
 \lefttop{I}\Gr^{\Vzero}_{\veca}
 \Gr^{W(N)}_h(\psi_{t,u})
\simeq
 \lefttop{I\sqcup\{i\}}\Gr^{\Vzero}_{(\veca,c)}
 \Gr^{W(N)}_h(\psi_{t,u})$.
\end{cor}
\pf
It follows from Lemma \ref{lem;b11.26.50}.
\hfill\qed

\begin{lem} \label{lem;9.22.40}
Let $S$ and $S'$ be primitive subsets of $\real^n$
such that $S\subset\nbigs(S')$.
We have the following isomorphism:
\begin{equation} \label{eq;9.20.24}
 \frac{\lefttop{\nbar}\Vzero_{S'}\bigl(\Gr^{W(N)}_h(\psi_{t,u})\bigr)}
 {\lefttop{\nbar}\Vzero_S
  \bigl(\Gr^{W(N)}_h(\psi_{t,u})\bigr) }
\simeq
 \Gr^{W(N)}_h
 \Bigl(
 \lefttop{\nbar}\Vzero_{S'}\bigl(\psi_{t,u}\bigr)
\Big/
 \lefttop{\nbar}\Vzero_S
  \bigl(\psi_{t,u}\bigr)
 \Bigr).
\end{equation}
\end{lem}
\pf
The left hand side of (\ref{eq;9.20.24})
can be rewritten as follows:
\begin{equation}\label{eq;a12.3.8}
\frac{\lefttop{\nbar}\Vzero_{S'}\bigl(\Gr^{W(N)}_h(\psi_{t,u})\bigr)}
{\lefttop{\nbar}\Vzero_S\bigl(\Gr^{W(N)}_h(\psi_{t,u})\bigr)}
\simeq
 \frac{\lefttop{\nbar}\Vzero_{S'}(\psi_{t,u})\cap W_h(N,\psi_{t,u})}
 {\lefttop{\nbar}\Vzero_S(\psi_{t,u})\cap W_h(N,\psi_{t,u})
 +W_{h-1}(N,\psi_{t,u})\cap\lefttop{\nbar}\Vzero_{S'}(\psi_{t,u})}.
\end{equation}
On the other hand,
the right hand side of (\ref{eq;9.20.24})
can be rewritten as follows:
\begin{equation}\label{eq;a12.3.7}
 {\rm R.H.S.}
\simeq
 \Cok\Bigl(
 \Gr^{W(N)}_h(\lefttop{\nbar}\Vzero_{S}(\psi_{t,u}))
\lrarr
 \Gr^{W(N)}_h(\lefttop{\nbar}\Vzero_{S'}(\psi_{t,u}))
 \Bigr).
\end{equation}
We have the following commutative diagramm:
\begin{equation}\label{eq;a12.3.5}
 \begin{CD}
 @. 0 @. 0 @. 0 @. \\
 @. @VVV @VVV @VVV \\
 0 @>>> \lefttop{\nbar}\Vzero_{S}\cap W_{h-1}(N)
   @>>> \lefttop{\nbar}\Vzero_{S'}\cap W_{h-1}(N)
   @>>> W_{h-1}\bigl(N,
 \lefttop{\nbar}\Vzero_{S'}\big/\lefttop{\nbar}\Vzero_{S}\bigr)
   @>>> 0 \\
 @. @VVV @VVV @VVV @.\\
 0 @>>> \lefttop{\nbar}\Vzero_{S}\cap W_{h}(N)
   @>>> \lefttop{\nbar}\Vzero_{S'}\cap W_{h}(N)
   @>>> W_{h}\bigl(N,
 \lefttop{\nbar}\Vzero_{S'}\big/\lefttop{\nbar}\Vzero_{S}\bigr)
   @>>> 0 \\
  @. @VVV @VVV @VVV @.\\
 0 @>>> \lefttop{\nbar}\Vzero_S
 \bigl(\Gr^{W(N)}_h\bigr)
   @>>> \lefttop{\nbar}\Vzero_{S'}
 \bigl(\Gr^{W(N)}_h\bigr)
   @>>> \lefttop{\nbar}\Vzero_{S'}\bigl(\Gr^{W(N)}_h\bigr)
 \big/\lefttop{\nbar}\Vzero_S\bigl(\Gr^{W(N)}_h\bigr)
   @>>> 0 \\
 @. @VVV @VVV @VVV \\
 @. 0 @. 0 @. 0 \\
 \end{CD}
\end{equation}
Here we omit to denote ``$\psi_{t,u}$''.
The first and the second rows 
are exact due to Lemma \ref{lem;a12.3.2}.
The third row is exact by definition.
The first and the second columns are also exact by definition.
Hence we obtain the exactness of the third columns.
Therefore we obtain the isomorphism:
\begin{equation}\label{eq;a12.3.6}
 \Cok\Bigl(
 W_{h-1}\bigl(
 N,\lefttop{\nbar}\Vzero_{S'}\big/\lefttop{\nbar}\Vzero_S
 \bigr)
\lrarr
  W_{h}\bigl(
 N,\lefttop{\nbar}\Vzero_{S'}\big/\lefttop{\nbar}\Vzero_S
 \bigr)
 \Bigr)
\simeq
 \Cok\bigl(
 \lefttop{\nbar}\Vzero_S\bigl(\Gr^{W(N)}_h\bigr)
 \lrarr
 \lefttop{\nbar}\Vzero_{S'}\bigl(\Gr^{W(N)}_h\bigr)
 \bigr).
\end{equation}
The right hand side of (\ref{eq;a12.3.6})
is isomorphic to the right hand side of
(\ref{eq;9.20.24}) due to (\ref{eq;a12.3.7}).

On the other hand,
since the first and the second rows in (\ref{eq;a12.3.5}) are exact,
the left hand side of (\ref{eq;a12.3.6})
is isomorphic to the left hand side of (\ref{eq;9.20.24}),
due to (\ref{eq;a12.3.8}).
Thus we are done.
\hfill\qed

\begin{lem} \label{lem;b11.26.100}
Let $I\sqcup J=\nbar$ be a decomposition.
Let $\vecc$ be an element of $\real^J_{<0}$,
and $\vecb$ be an element of $\real^I$.
Then 
$\lefttop{I}\Gr^{\Vzero}_{\vecb}
 \lefttop{J}\Gr^{\Vzero}_{\vecc}
 \Gr^{W(N)}_h(\psi_{t,u})$ is strict.
\end{lem}
\pf
Due to Lemma \ref{lem;9.22.40},
$\lefttop{I}\Gr^{\Vzero}_{\vecb}
 \lefttop{J}\Gr^{\Vzero}_{\vecc}
 \Gr^{W(N)}_h(\psi_{t,u})$
is isomorphic to
$\Gr^{W(N)}_h
 \lefttop{I}\Gr^{\Vzero}_{\vecb}
 \lefttop{J}\Vzero_{\vecc}(\psi_{t,u})$.
It is strict,
for $W(N)$ on 
$\lefttop{I}\Gr^{\Vzero}_{\vecb}\lefttop{J}\Vzero_{\vecc}$
is a filtration in the category of vector bundles.
\hfill\qed

\begin{lem}
Let $I\sqcup J=\nbar$ be a decomposition.
Let $\vecb$ be an element of $\real^I$,
and $\vecc$ be an element of $\real^J$.
Then
$\lefttop{I}\Gr^{\Vzero}_{\vecb}
 \lefttop{J}\Vzero_{\vecc}\Gr^{W(N)}_h(\psi_{t,u})$
is strict.
\end{lem}
\pf
We put $M_-(\vecc):=\bigl\{i\,\big|\,c_i<0\bigr\}$.
Note we have $|M_-(\vecc)|\leq n-|I|$.
We use a descending induction
on $(|I|,|M_-(\vecc)|)$.

In the case $(|I|,|M_-(\vecc)|)=(n,0)$
or $(|I|,|M_-(\vecc)|)=(m,n-m)$,
the claim holds due to Lemma \ref{lem;b11.26.100}.
Let us consider the case
$|I|=m$ and $|M_-(\vecc)|<n-m$.
We may assume that $J=\{1,\ldots,n-m\}$
and $I=\nbar-J$.
We may assume that $c_1\geq 0$.

We have the filtration $\lefttop{1}\Vzero$
on
 $\lefttop{I}\Gr^{\Vzero}_{\vecb}
 \lefttop{J}\Vzero_{\vecc}
 \Gr^{W(N)}_h(\psi_{t,u})$.
For any $d$,
we may assume that
the following is strict:
\[
 \lefttop{1}\Gr^{\Vzero}
 \bigl(
 \lefttop{I}\Gr^{\Vzero}_{\vecb}
 \lefttop{J}\Vzero_{\vecc}
 \Gr^{W(N)}_h(\psi_{t,u})
 \bigr).
\]
For $d<0$,
we may assume that
$\lefttop{1}\Vzero_{d}
 \bigl(
 \lefttop{I}\Gr^{\Vzero}_{\vecb}
 \lefttop{J}\Vzero_{\vecc}
 \Gr^{W(N)}_h(\psi_{t,u})
 \bigr)$ is strict.
Then we obtain the strictness
of $\bigl(
 \lefttop{I}\Gr^{\Vzero}_{\vecb}
 \lefttop{J}\Vzero_{\vecc}
 \Gr^{W(N)}_h(\psi_{t,u})
 \bigr)$.
\hfill\qed

\begin{cor}
$\Gr^{W(N)}_h\psi_{t,u}$ 
and $\lefttop{i}\Gr^{\Vzero}\Gr^{W(N)}_h\psi_{t,u}$
are strict.
\hfill\qed
\end{cor}

\begin{cor}
$\Gr^{W(N)}_h\psi_{t,u}$ 
is strictly specializable along $z_i=0$
for any $i=1,\ldots,n$.
\hfill\qed
\end{cor}

%% file: a80.6.tex

\subsubsection{Preliminary}

Let $m$ be an integer such that $0\leq m\leq n$.
Let $c$ be a real number.
Let $h$ be an integer.
Let $S$ be a primitive subset of $\real^n$.
We put as follows:
\[
\begin{array}{l}
 \nbiga(m,c,S,h):=
 \Bigl(
 \FFzero_{<m}\cap\lefttop{i}\Vzero_c
 +\FF_m\cap \lefttop{i}\Vzero_{<c}
 \Bigr)
 \cap \lefttop{\nbar}\Vzero_S
 \cap W_h(N,\psi_{t,u}),\\
\mbox{{}}\\
 \tilde{\nbiga}(m,c,S,h):=
 \Bigl(
 \FFzero_{<m}\cap \lefttop{i}\Vzero_{c}
 +\FF_{m}\cap \lefttop{i}\Vzero_{<c}
 \Bigr)\cap
 \lefttop{\nbar}\Vzero_S\cap
 \Gr^{W(N)}_h(\psi_{t,u}).
\end{array}
\]
\begin{lem}\label{lem;a12.3.17}
We have the following:
\begin{equation}\label{eq;a12.3.16}
 \nbiga(m,c,s,h)
=W\Bigl(
 N,
 \bigl(
 \FFzero_{<m}\cap \lefttop{i}\Vzero_c
+\FFzero_m\cap \lefttop{i}\Vzero_{<c}
 \bigr)
\cap \lefttop{\nbar}\Vzero_S(\psi_{t,u})
 \Bigr).
\end{equation}
\end{lem}
\pf
The proof is essentially same as the arguments for
Corollary \ref{cor;a12.3.10}.
We only indicate an outline.

We put as follows:
\[
 \nbigy:=
 \frac{\lefttop{i}\Vzero_c\cap \lefttop{\nbar}\Vzero_S}
 {\bigl(\FFzero_{m-1}\cap \lefttop{i}\Vzero_c
 +\FFzero_m\cap \lefttop{i}\Vzero_{<c}\bigr)\cap \Vzero_S}.
\]
Then we have the following exact sequence:
\[
 \begin{CD}
 0 @>>>
 \lefttop{i}\Vzero_{<c}\cap \lefttop{\nbar}\Vzero_S
 \bigl(\Gr^{\FFzero}_m\bigr)
 @>>>
 \lefttop{i}\Vzero_c\cap\lefttop{\nbar}\Vzero_S
 \bigl(\GGzero_{m-1}\bigr)
 @>>>
 \nbigy @>>>0.
 \end{CD}
\]
By an argument similar to the proof of Lemma \ref{lem;a12.2.130},
we obtain the following:
\[
 W_h\Bigl(
 N,\lefttop{i}\Vzero_{<c}\cap \lefttop{\nbar}\Vzero_S
 \bigl(\Gr^{\FFzero}_m\bigr)
 \Bigr)
=W_h\Bigl(N,
 \lefttop{i}\Vzero_c\cap\lefttop{\nbar}\Vzero_S
 \bigl(\GGzero_{m-1}\bigr)
 \Bigr)
\cap
 \Bigl(
 \lefttop{i}\Vzero_{<c}\cap \lefttop{\nbar}\Vzero_S
 \bigl(\Gr^{\FFzero}_m\bigr)\Bigr).
\]
Then by an argument similar to the proof of
Lemma \ref{lem;a11.26.80},
we can show that
the projection 
$\lefttop{i}\Vzero_c\cap\lefttop{\nbar}\Vzero_S\bigl(\GGzero_{m-1}\bigr)
\lrarr
 \nbigy$ induces the surjection:
\[
 W_h\bigl(N,
 \lefttop{i}\Vzero_c\cap\lefttop{\nbar}\Vzero_S\bigl(\GGzero_{m-1}\bigr)
 \bigr)
\lrarr
 W_h\bigl(N,\nbigy\bigr).
\]
Then it immediately follows that
the projection
$\lefttop{i}\Vzero_{c}\cap\lefttop{\nbar}\Vzero_S(\psi_{t,u})
\lrarr \nbigy$
induces the surjection:
\begin{equation}\label{eq;a12.3.11}
 W_h\bigl(N,
 \lefttop{i}\Vzero_{c}\cap\lefttop{\nbar}\Vzero_S(\psi_{t,u})
 \bigr)
\lrarr
 W_h\bigl(N,\nbigy \bigr).
\end{equation}
Let us consider the following exact sequence:
\begin{equation}\label{eq;a12.3.15}
 \begin{CD}
 0
@>>>
 \bigl(
 \FFzero_{m-1}\cap\lefttop{i}\Vzero_c
 \cap\FFzero_m\cap\lefttop{i}\Vzero_{<c}\bigr)
 \cap\lefttop{\nbar}\Vzero_S(\psi_{t,u})
@>{a}>>
 \lefttop{i}\Vzero_c\lefttop{\nbar}\Vzero_S(\psi_{t,u})
@>{b}>>
 \nbigy
@>>> 0.
 \end{CD}
\end{equation}
The surjectivity of (\ref{eq;a12.3.11}) implies
that the morphism $b$ in (\ref{eq;a12.3.15})
is strict with respect to the weight filtrations.
Due to Corollary \ref{cor;a12.2.120},
the sequence (\ref{eq;a12.3.15}) is strict
with respect to the weight filtration.
In particular,
the morphism $a$ is strict with respect to
the weight filtration.
Thus we obtain the equality (\ref{eq;a12.3.16}).
\hfill\qed

\vspace{.1in}
We put as follows:
\[
\begin{array}{ll}
 \nbigb(m,c,S,h):=
 \Bigl(\FF_m\cap \lefttop{i}\Vzero_c\Bigr)
 \cap \lefttop{\nbar}\Vzero_S
 \cap W_h(N,\psi_{t,u}),
 &
 \tilde{\nbigb}(m,c,S,h):
 \Bigl(
 \FF_m\cap\lefttop{i}\Vzero_c\Bigr)\cap
 \lefttop{\nbar}\Vzero_S\cap
 \Gr^{W(N)}_h(\psi_{t,u}), \\
\mbox{{}} \\
 \nbigc(m,c,S,h):=
 W_h\bigl(N,
 \lefttop{\nbar}\Vzero_S\lefttop{i}\Gr^{\Vzero}_c\Gr^{\FFzero}_m\bigr),
 &
 \tilde{\nbigc}(m,c,S,h):=
 \Gr^{\FFzero}_m
 \lefttop{i}\Gr^{\Vzero}_c
 \lefttop{\nbar}\Vzero_S\bigl(
 \Gr^{W(N)}_h(\psi_{t,u})
 \bigr).
\end{array}
\]

\begin{lem}\label{lem;a12.3.21}
The following sequence is exact:
\[
 \begin{CD}
 0 @>>>
 \nbiga(m,c,S,h) @>>>
 \nbigb(m,c,S,h) @>>>
 \nbigc(m,c,S,h) @>>> 0.
 \end{CD}
\]
\end{lem}
\pf
We have the following exact sequence:
\begin{multline}\label{eq;a12.3.18}
 0 \lrarr
 \bigl(
 \FFzero_{<m}\cap\lefttop{i}\Vzero_c
 +\FFzero_m\cap\lefttop{i}\Vzero_{<c}
 \bigr)\cap\lefttop{\nbar}\Vzero_S
\lrarr
 \FFzero_m\cap\lefttop{i}\Vzero_c\cap\lefttop{\nbar}\Vzero_S \\
\lrarr
 \lefttop{\nbar}\Vzero_S\lefttop{i}\Gr^{\Vzero}_c
 \Gr^{\FFzero}_m
\lrarr 0.
\end{multline}
Due to Lemma \ref{lem;a12.3.17} and Corollary \ref{cor;a12.2.120}
the sequence (\ref{eq;a12.3.18}) is strict
with respect to the weight filtration.
Due to Corollary \ref{cor;a12.3.10},
we have the following:
\[
 \nbigb(m,c,S,h)=
W_h\bigl(
 N,\FFzero_m\cap \lefttop{i}\Vzero_c\cap\lefttop{\nbar}\Vzero_S
 \bigr).
\]
Thus we are done.
\hfill\qed

\vspace{.1in}

Then we have the following diagramm:
\begin{equation}\label{eq;a12.3.20}
 \begin{CD}
 @. 0 @. 0 @. 0 @. \\
 @. @VVV @VVV @VVV @. \\
0 @>>> \nbiga(m,c,S,h-1)
  @>>> \nbigb(m,c,S,h-1)
  @>>> \nbigc(m,c,S,h-1)
  @>>> 0 \\
 @. @VVV @VVV @VVV @. \\
 0 @>>> \nbiga(m,c,S,h)
  @>>> \nbigb(m,c,S,h)
  @>>> \nbigc(m,c,S,h)
  @>>> 0 \\
 @. @VVV @VVV @VVV @. \\
 0 @>>> \tilde{\nbiga}(m,c,S,h)
  @>>> \tilde{\nbigb}(m,c,S,h)
  @>>> \tilde{\nbigc}(m,c,S,h)
  @>>> 0 \\
 @. @VVV @VVV @VVV @. \\
 @. 0 @. 0 @. 0 @.\\
 \end{CD}
\end{equation}
The first and the second rows are exact due to
Lemma \ref{lem;a12.3.21}.
The third row is exact by definition.
The first and the second columns are also exact by definition.
Therefore we obtain the exactness
of the third columns.

\begin{cor}
We have the following isomorphism:
\begin{equation}\label{eq;a12.3.25}
\tilde{\nbigc}(m,c,S,h)
\simeq
\Cok
\Bigl(
 \nbigc(m,c,S,h-1)
\lrarr
 \nbigc(m,c,S,h)
\Bigr).
\end{equation}
\end{cor}
\pf
It follows from the exactness of 
the third column in the diagramm (\ref{eq;a12.3.20}).
\hfill\qed

\begin{cor}\label{cor;c11.26.1}
We have the natural isomorphism:
\begin{equation}\label{eq;a12.3.26}
 \Gr^{\FFzero}_m
 \lefttop{i}\Gr^{\Vzero}_c
 \lefttop{\nbar}\Vzero_{S}
 \Gr^{W(N)}_h(\psi_{t,u})
\simeq
 \Gr^{W(N)}_h
 \Gr^{\FFzero}_m
 \lefttop{i}\Gr^{\Vzero}_c
 \lefttop{\nbar}\Vzero_S(\psi_{t,u}).
\end{equation}
\end{cor}
\pf
(\ref{eq;a12.3.26}) is just a reformulation of
(\ref{eq;a12.3.25}).
\hfill\qed

\subsubsection{The filtrations $\FFzero$ 
 on $\lefttop{i}\psi_{0}P_h\Gr^{W(N)}_h(\psi_{t,u})$
 and $\lefttop{i}\psi_{-\vecdelta_0}P_h\Gr^{W(N)}_h(\psi_{t,u})$}

We have the induced filtrations $\FFzero$
on $\lefttop{i}\psi_{0}P_h\Gr^{W(N)}_h(\psi_{t,u})$
and $\lefttop{i}\psi_{-\vecdelta_0}P\Gr^{W(N)}_h(\psi_{t,u})$.

\begin{lem}
We have the following isomorphism:
\[
 \Gr^{\FFzero}_m\lefttop{i}\psi_0P_h\Gr^{W(N)}_h(\psi_{t,u})
\simeq
 P_h\Gr^{W(N)}_h\lefttop{i}\psi_{0}\Gr^{\FFzero}_m(\psi_{t,u})
=\bigoplus_{
 \substack{|I|=n-m,\\ i\not\in I}}
 P_h\Gr^{W(N)}_h\lefttop{i}\psi_0\lefttop{I}\nbigl.
\]
We have the following isomorphism:
\[
 \Gr^{\FFzero}_m\bigl(
 \lefttop{i}\psi_{-\vecdelta_0}P_h\Gr^{W(N)}_h(\psi_{t,u})
 \bigr)
\simeq
 P_h\Gr^{W(N)}_h
 \lefttop{i}\psi_{-\vecdelta_0}
 \Gr^{\FFzero}_m(\psi_{t,u})
=\bigoplus_{\substack{|I|=n-m+1,\\ i\in I}}
 P_h\Gr^{W(N)}_h
 \lefttop{i}\psi_{-\vecdelta_0}(\lefttop{I}\nbigl).
\]
\end{lem}
\pf
It follows from Corollary \ref{cor;c11.26.1}.
\hfill\qed

\subsubsection{The morphisms $\can_i$ and $\var_i$}

The morphisms $\can_i$ and $\var_i$ induce the following:
\[
 \can_i:
 \FFzero_m\lefttop{i}\psi_{-\vecdelta_0}
 P_h\Gr^{W(N)}_h(\psi_{t,u})
\lrarr
 \FFzero_{m+1}\lefttop{i}\psi_0
P_h\Gr^{W(N)}_h(\psi_{t,u}).
\]
\[
 \var_i:
 \FFzero_{m+1}\lefttop{i}\psi_0
P_h\Gr^{W(N)}_h(\psi_{t,u})
\lrarr 
 \FFzero_m\lefttop{i}\psi_{-\vecdelta_0}
 P_h\Gr^{W(N)}_h(\psi_{t,u}).
\]
Thus we obtain the induced morphisms:
\[
 \can_i(m):
 \frac{\FFzero_{n-1}}{\FFzero_{m-1}}
 \lefttop{i}\psi_{-\vecdelta_0}
 P_h\Gr^{W(N)}_h(\psi_{t,u})
\lrarr
 \frac{\FFzero_{n}}{\FFzero_{m}}
 \lefttop{i}\psi_0
P_h\Gr^{W(N)}_h(\psi_{t,u}).
\]
\[
 \Gr_m(\can_i):
 \Gr^{\FFzero}_m\lefttop{i}\psi_{-\vecdelta_0}
 P_h\Gr^{W(N)}_h(\psi_{t,u})
\lrarr
 \Gr^{\FFzero}_{m+1}\lefttop{i}\psi_{0}
 P_h\Gr^{W(N)}_h(\psi_{t,u}).
\]
\[
 \var_i(m):
 \frac{\FFzero_{n}}{\FFzero_{m}}
 \lefttop{i}\psi_0
P_h\Gr^{W(N)}_h(\psi_{t,u})
\lrarr
 \frac{\FFzero_{n-1}}{\FFzero_{m-1}}
 \lefttop{i}\psi_{-\vecdelta_0}
 P_h\Gr^{W(N)}_h(\psi_{t,u}).
\]
\[
  \Gr_m(\var_i):
 \Gr^{\FFzero}_{m+1}\lefttop{i}\psi_{0}
 P_h\Gr^{W(N)}_h(\psi_{t,u})
\lrarr
 \Gr^{\FFzero}_m\lefttop{i}\psi_{-\vecdelta_0}
 P_h\Gr^{W(N)}_h(\psi_{t,u}).
\]

In the following diagramm,
we omit to denote $\lefttop{i}\psi_{0}P_h\Gr^{W(N)}_h(\psi_{t,u})$:
\begin{equation} \label{eq;c11.26.5}
 \begin{array}{ccccccc}
 0 \lrarr
 & \Gr^{\FFzero}_{m+1}
 & \lrarr
 & \FFzero_n/\FFzero_m
 & \lrarr
 & \FFzero_n/\FFzero_{m+1}
 & \lrarr 0 \\
 & \bigcup & & \bigcup & & \bigcup \\
 0\lrarr
 & \Image\Gr_m(\can_i)
 & \stackrel{a}{\lrarr}
 & \Image(\can_i(m))
 & \stackrel{b}{\lrarr}
 & \Image(\can_i(m+1))
 & \lrarr 0.
 \end{array}
\end{equation}
The upper sequence is exact by definition.
\begin{lem} \label{lem;c11.26.15}
The lower sequence is exact.
\end{lem}
\pf
The surjectivity of $b$ and the injectivity of $a$
are easy.
Let us show $\Ker(b)=\Image(a)$.
We have only to show $\ker(b)\subset\Image(a)$.
The argument is similar to the proof of
Lemma \ref{lem;a12.2.130}.

Recall that we have the decomposition:
\[
 \Image(\Gr_m(\can_i))=
 \bigoplus_{\substack{|I|=n-m-1\\ i\not\in I}}
 \nbiga_I.
\]
Here $\nbiga_I$ denotes a subbundle of
$P_h\Gr^{W(N)}_h\lefttop{i}\psi_0\lefttop{I}\nbigl$
over $\nbigd_{\nbar-I}$
(Corollary \ref{cor;c11.26.20}).
It is easy to see that the
restrictions of
$\Image(\Gr_m(\can_i))$
and $\Image(\can_i(m))$ to $\nbigx-\nbigd^{[m]}$ are same.
Here we put
$\nbigd^{[m]}:=\bigcup_{|I'|=n-m}\nbigd_{\nbar-I'}$.
Thus we obtain
$\ker(b)\subset \Image(a)$.
\hfill\qed

\vspace{.1in}
In the following diagram,
we omit to denote
$\lefttop{i}\psi_{-\vecdelta_0}P_h\Gr^{W(N)}_h(\psi_{t,u})$:
\[
 \begin{array}{ccccccc}
 0 \lrarr
 & \Gr^{\FFzero}_{m}
 & \lrarr
 & \FFzero_{n-1}/\FFzero_{m-1}
 & \lrarr
 & \FFzero_{n-1}/\FFzero_{m}
 & \lrarr 0 \\
 & \bigcup & & \bigcup & & \bigcup \\
 0\lrarr
 & \Image\Gr_m(\var_i)
 & \stackrel{a}{\lrarr}
 & \Image(\var_i(m))
 & \stackrel{b}{\lrarr}
 & \Image(\var_i(m+1))
 & \lrarr 0.
 \end{array}
\]
The upper sequence is exact by definition.

\begin{lem}\label{lem;c11.26.16}
The lower sequence is also exact.
\end{lem}
\pf
It can be shown by an argument similar to the proof
of Lemma \ref{lem;c11.26.15}.
We use Corollary \ref{cor;b12.3.50}
instead of Corollary \ref{cor;c11.26.20}.
\hfill\qed

\vspace{.1in}
In the following diagramm,
we omit to denote $\lefttop{i}\psi_0P_h\Gr^{W(N)}_h(\psi_{t,u})$:
\[
  \begin{array}{ccccccc}
 0 \lrarr
 & \Gr^{\FFzero}_{m+1}
 & \lrarr
 & \FFzero_n/\FFzero_m
 & \lrarr
 & \FFzero_n/\FFzero_{m+1}
 & \lrarr 0 \\
 & \bigcup & & \bigcup & & \bigcup \\
 0\lrarr
 & \Ker\Gr_m(\var_i)
 & \stackrel{a}{\lrarr}
 & \Ker(\var_i(m))
 & \stackrel{b}{\lrarr}
 & \Ker(\var_i(m+1))
 & \lrarr 0.
 \end{array}
\]
The upper sequence is exact by definition.

\begin{lem}\label{lem;b12.3.51}
The lower sequence is also exact.
\end{lem}
\pf
It immediately follows from Lemma \ref{lem;c11.26.16}.
\hfill\qed

\begin{lem}
We have the following decomposition:
\[
 \frac{\FFzero_n}{\FFzero_m}
 \bigl(\lefttop{i}\psi_0P_h\Gr^{W(N)}_h(\psi_{t,u})\bigr)
=\Image(\can_i(m))\oplus
 \Ker(\var_i(m)).
\]
\end{lem}
\pf
Recall that we have the decomposition
of $\Gr^{\FFzero}_m(\lefttop{i}\psi_0P_h\Gr^{W(N)}_h(\psi_{t,u}))$.
(Corollary \ref{cor;c11.26.20}).
We also have Lemma \ref{lem;c11.26.15} and
Lemma \ref{lem;b12.3.51}.
Then the lemma can be shown by an easy descending induction
on $m$.
\hfill\qed

\begin{cor} \label{cor;d11.26.1}
We have the decomposition:
\[
 \lefttop{i}\psi_0P_h\Gr^{W(N)}_h(\psi_{t,u})
=\Image(\can_i)\oplus
 \Ker(\var_i).
\]
As a result,
$P_h\Gr^{W(N)}_h(\psi_{t,u})$ is 
strictly $S$-decomposable along $z_i=0$.
\hfill\qed
\end{cor}

%% file: a80.7.tex

\subsubsection{Some properties of the components}

\begin{prop}\label{prop;b12.6.150}
We have the decomposition:
\[
 P_h\Gr^{W(N)}\psi_{t,u}=
 \bigoplus_{I\subset\nbar}
 \nbigm_I.
\]
Here $\nbigm_I$ denotes the $\nbigr$-module
satisfying the following:
\begin{itemize}
\item
 The support of $\nbigm_I$ is $\nbigd_I$,
 we have the injection
 $\nbigm_I\lrarr \iota_{I\,\ast}(\nbigm_{I|\nbigd_I^{\circ}})$.
 Here we put
 $\nbigd_I^{\circ}:=\nbigd_I-\bigcup_{I'\supsetneq I}\nbigd_{I'}$,
and $\iota_{I}$ denotes the open immersion
 $\nbigd_I^{\circ}\lrarr\nbigd_I$.
\item
 $\nbigm_I$ are holonomic.
\item
 $\nbigm_I$ are regular and strictly $S$-decomposable along $z_i=0$
 for any $i$.
\end{itemize}
\end{prop}
\pf
It immediately follows from Corollary \ref{cor;d11.26.1}
and Proposition \ref{prop;a11.23.1}.
\hfill\qed

\begin{lem} \label{lem;a11.30.60}
For any subset $I\subset\nbar$,
the $\nbigr$-module $\nbigm_I$ satisfies
the conditions {\rm\ref{number;d11.26.2}--\ref{number;d11.26.5}}
in Proposition {\rm\ref{prop;a11.26.2}}.
\end{lem}
\pf
It is easy to check
the conditions \ref{number;d11.26.2} and \ref{number;d11.26.3}.
(The condition \ref{number;d11.26.2} is replaced
by the injectivity
$\nbigm_I\lrarr
 \iota_{I\,\ast}\bigl(\nbigm_{I\,|\,\nbigd_I^{\circ}}\bigr)$.)

Let $\lambda\in\cnum^{\ast}_{\lambda}$ be generic.
By the same argument as those in
the section \ref{section;b11.24.1}--\ref{section;b12.2.0},
we can show that
$\Gr^{W(N)}_h(\psi_{t,u})_{|\nbigxlambda}$ 
is also strictly $S$-decomposable along $z_i=0$
for $i=1,\ldots,n$.
In that case,
it is easy to see the specialization 
$\nbigm_{I\,|\,\nbigxlambda}$ gives the $I$-component
of $\Gr^{W(N)}_h(\psi_{t,u})_{|\nbigxlambda}$.
Hence $\nbigm_{I\,|\,\nbigxlambda}$ 
are also strictly $S$-decomposable along $z_i=0$,
i.e.,
the condition \ref{number;d11.26.4} is satisfied.

Let us check the condition \ref{number;d11.26.5}.
We put $J:=\nbar-I$.
Let us consider $\lefttop{I}\psi_{0}\nbigm_I$,
which is the $\nbigr$-module on $\nbigd_I$.
Let us pick any point $\lambda_0\in\cnum_{\lambda}$.
We denote $\nbigd_I(\lambda_0,\epsilon_0)$
by $\nbigd_I$ for simplicity.
Around $\lambda_0$,
the sheaf $\lefttop{I}\psi_0\nbigm_I$
is a direct summand of
$\Gr^{W(N)}\lefttop{I}\Gr^{\Vzero}_0(\psi_{t,u})$.
We have the induced filtrations
$\lefttop{i}\Vzero$ $(i\in J)$
on $\lefttop{I}\psi_0\nbigm_I$.
Note that $\lefttop{i}\Vzero$ is the $V$-filtration
along $z_i=0$.
Since 
$\lefttop{J}\Vzero_{<0}\lefttop{I}\Gr^{\Vzero}_{0}
 \Gr^{W(N)}_h(\psi_{t,u})$ is 
coherent and locally free $\nbigo_{\nbigd_I}$-module,
and since
$\lefttop{J}\Vzero_{<0}\nbigm_I$ is a direct summand
of $\lefttop{J}\Vzero_{<0}\lefttop{I}\Gr^{\Vzero}_{0}
 \Gr^{W(N)}_h(\psi_{t,u})$,
the $\nbigo_{\nbigd_I}$-module
$\lefttop{J}\Vzero_{<0}\nbigm_I$ is also
coherent and locally free.

By our construction,
It is easy to see that
$\lefttop{J}\Vzero_{<0}
 \lefttop{I}\Gr^{\Vzero}_0\Gr^{W(N)}(\psi_{t,u})$
generates $\lefttop{I}\Gr^{\Vzero}_0\Gr^{W(N)}(\psi_{t,u})$.
Hence we can conclude that
$\lefttop{J}\Vzero_{<0}\lefttop{I}\psi_0\nbigm_I$
generate $\lefttop{I}\psi_0\nbigm_I$.
Thus we are done.
\hfill\qed

%% file: a81.2.tex

\subsubsection{The isomorphism}

\label{subsubsection;a11.30.20}

Let us see the component $\nbigm_{\nbar}$.
Since the support of $\nbigm_{\nbar}$ is
$\cnum_{\lambda}\times \{O\}$,
and since $\nbigm_{\nbar}$ is specially $S$-decomposable
along $z_i$ $(i\in \nbar)$,
$\nbigm_{\nbar}$ is isomorphic to the push forward
of $\lefttop{\nbar}\psi_0\nbigm_{\nbar}$.
We also have the following:
\begin{equation}\label{eq;a11.30.10}
\lefttop{\nbar}\psi_0\nbigm_{\nbar}
\simeq
 \Ker\Bigl(
 \lefttop{\nbar}\tildepsi_{0}
 \bigl(
 P_h\Gr^{W(N)}_h(\tildepsi_{t,u}\gbige[\deldel_t])
 \bigr)
\lrarr
 \bigoplus_{i}
 \lefttop{\nbar}\tildepsi_{-\vecdelta_{0,i}}
 \bigl(
 P_h\Gr^{W(N)}_h(\tildepsi_{t,u}\gbige[\deldel_t])
 \bigr)
\Bigr)
\end{equation}

On the other hand,
$\lefttop{\nbar}\psi_0\bigl(
 P_h\Gr^{W(N)}_h
 \bigl(
 \tildepsi_{t,u}
 \bigr)
 \bigr)$
is naturally isomorphic to
$P_h\Gr^{W(N)}_h\nbigq(\nbar,0,\vecu,\vecm,u)$
(see the subsubsection \ref{subsubsection;b12.3.150}).
Here we put $\vecu=(\overbrace{u,\ldots,u}^n)$.
We put $\nbigc_{\nbar\,0}:=\nbigc_{\nbar\,|\,\cnum_{\lambda}}$.

\begin{lem}\label{lem;b12.3.165}
We have the natural isomorphism:
\[
 \nbigc_{\nbar\,0}\simeq
 \lefttop{\nbar}\psi_0\bigl(\nbigm_{\nbar}
 \bigr).
\]
\end{lem}
\pf
We have only to compare
the right hand side of
(\ref{eq;a11.30.10})
and $\nbigc_{\nbar\,0}$,
which can be checked directly from the definitions.
\hfill\qed

%% file: a85.1.tex

In this section,
we use the left $\nbigr$-module structure on $\nbige$
and $\gbige$.

%% file: a34.tex

\subsubsection{The pairing on $\nbige$ and the prolongation}

\label{subsubsection;9.25.2}

We put $X:=\Delta^n$, $D_i:=\{z_i=0\}$ and $D:=\bigcup_{i=1}^n D_i$.
Let $\harmonicbundle$ be a tame harmonic bundle over $X-D$.
Let us recall the sesqui-linear pairing
$C_0:\nbige_{\AAA}\otimes\overline{\nbige_{\AAA}}
 \lrarr \distribution_X^{\AAA}$
given by Sabbah.
Let $f$ be a section of $\nbige_{|\nbigx(\lambda_0,\epsilon_0)}$
and $g$ be a section of
$\nbige_{|\sigma(\nbigx(\lambda_0,\epsilon_0))}$.
Then the pairing is defined as follows:
\begin{equation}\label{eq;a11.22.3}
 C_0(f,\bar{g}):=h\bigl(f,\sigma^{\ast}(g)\bigr).
\end{equation}
We prolong it to the sesqui-linear pairing $\gbigc$
for the $\nbigr$-module $\gbige$.
We have only to consider the pairings of the sections
of $\gbige_{|\Delta(\lambda_0,\epsilon_0)}$
and $\overline{\gbige_{|\sigma(\Delta(\lambda_0,\epsilon_0))}}$.

Let $f$ be a section of
$V^{(\lambda_0)}_{<0}\bigl(\gbige_{|\nbigx(\lambda_0,\epsilon_0)}\bigr)$
and $g$ be a section of
$V^{(\sigma(\lambda_0))}_{<0}\bigl(
 \gbige_{|\sigma(\nbigx(\lambda_0,\epsilon_0))}\bigr)$.
Then we have the following estimate
of $C^{\infty}$-functions on $(X-D)\times\Delta(\lambda_0,\epsilon_0)$,
for some positive constants $C$ and $\epsilon$:
\[
 \bigl|h\bigl(f,\sigma^{\ast}(g)\bigr)\bigr|\leq
 C\cdot \prod_{i=1}^n |z_i|^{-2+\epsilon}.
\]
Thus the functions $h\bigl(f,\sigma^{\ast}(g)\bigr)$
is an $L^1$-function with respect to the standard metric
$\sum_{i=1}^n |dz_i\cdot d\bar{z}_i|$.
Hence $h\bigl(f,\sigma^{\ast}(g)\bigr)$ naturally gives
a distribution.
We denote it by $\gbigc(f,\bar{g})$.

For an element $\vecp=(p_1,\ldots,p_n)\in\seisuu_{\geq\,0}^n$,
we denote $\prod_{i=1}^n\deldel_i^{p_i}$ by
$\deldel_z^{\vecp}$.
Let us consider sections of the form
$\deldel_z^{\vecn_1}\cdot f$
and $\deldel_z^{\vecn_2}\cdot g$,
where $f$ and $g$ are sections
of $V^{(\lambda_0)}_{<0}\bigl(\gbige_{|\nbigx(\lambda_0,\epsilon_0)}\bigr)$
and
$V^{(\sigma(\lambda_0))}_{<0}\bigl(
 \gbige_{|\sigma(\nbigx(\lambda_0,\epsilon_0))}\bigr)$
respectively,
and $\vecn_1$ and $\vecn_2$ be elements of $\seisuu_{\geq\,0}^{n}$.
(Note we use the left $\nbigr$-modules in this subsection.)
Since $\gbigc(f,g)$ gives a distribution
on $X$,
we obtain the following distribution:
\begin{equation} \label{eq;9.25.1}
 \gbigc\bigl(\deldel_z^{\vecn_1}\!\cdot\! f,\,\,
\overline{\deldel_z^{\vecn_2}\!\cdot\! g}
 \bigr):=
 \deldel_z^{\vecn_1}\cdot\deldelbar_z^{\vecn_2}\cdot
 \gbigc(f,\bar{g}).
\end{equation}

Since $\gbige_{|\nbigx(\lambda_0,\epsilon_0)}$
and $\gbige_{|\sigma(\nbigx(\lambda_0,\epsilon_0))}$
are generated by
$V^{(\lambda_0)}_{<0}\bigl(\gbige_{|\nbigx(\lambda_0,\epsilon_0)}\bigr)$
and
$V^{(\sigma(\lambda_0))}_{<0}\bigl(
 \gbige_{|\sigma(\nbigx(\lambda_0,\epsilon_0))}\bigr)$,
the formula (\ref{eq;9.25.1})
gives the pairing $\gbigc$ of $\gbige$.

\subsubsection{The uniqueness}

Let $C:\gbige_{\AAA}\otimes\overline{\gbige_{\AAA}}
\lrarr \distribution_X$ be a sesqui-linear pairing.

\begin{lem} \label{lem;9.25.3}
Assume that $C$ vanishes on $X-D$.
Then it vanishes on $X$.
\end{lem}
\pf
Similar to the proof of Proposition \ref{prop;a11.23.10}.
\hfill\qed

\begin{cor} \label{cor;a11.23.6}
Let $C'$ be a pairing of $\gbige$.
Assume that the restriction of $C'$
to $X-D$ coincides with $C_0$ given by Sabbah
{\rm(\ref{eq;a11.22.3})}.
Then $C'$ is same as
$\gbigc$ given in the subsubsection
{\rm \ref{subsubsection;9.25.2}}.
\end{cor}
\pf
It immediately follows from Lemma \ref{lem;9.25.3}.
\hfill\qed

%% file: a85.tex

\subsubsection{The definition of
 $\lefttop{\nbar}\overline{\nbigg}_{\vecu}(\nbige)$}

\label{subsubsection;a12.4.1}

Let $\vecu$ be an element of
$\KMSoverline(\nbige^0,\nbar)$.
We put as follows:
\[
 \lefttop{\nbar}\overline{\nbigg}_{\vecu}(\nbige):=
\lefttop{\nbar}\Gr^{\nbigfzero}_{\paramap^f(\lambda_0,\vecu)}
\lefttop{\nbar}\EEzero(\eigenmap^f(\lambda_0,\vecu)).
\]
It is easy to see that 
we have the globally defined
flat bundle
$\lefttop{\nbar}\overline{\nbigg}_{\vecu}(\nbige)$
on $(X-D)\times\cnum^{\ast}_{\lambda}$
satisfying the following:
\[
 \lefttop{\nbar}\overline{\nbigg}_{\vecu}(\nbige)_{|
 \Delta(\lambda_0,\epsilon_0)\times(X-D)}
=\lefttop{\nbar}\overline{\nbigg}_{\vecu}^{(\lambda_0)}(\nbige)
 _{|\Delta(\lambda_0,\epsilon_0)\times(X-D)}.
\]

\begin{rem}\label{rem;c12.3.20}
Note that $\lefttop{\nbar}\overline{\nbigg}_{\vecu}$
is not same as $\lefttop{\nbar}\nbigg_{\vecu_1}$
which is defined for $\vecu_1\in\KMS(\nbige^0,\nbar)$
(the subsubsection {\rm\ref{subsubsection;10.18.10}}).
Let $\pi:\KMS(\nbige^0,\nbar)\lrarr\KMSoverline(\nbige^0,\nbar)$
be the projection.
Then we have the following relation:
\begin{equation}\label{eq;a12.4.2}
  \lefttop{\nbar}\nbigg_{\vecu\,|\,(X-D)\times\cnum^{\ast}_{\lamda}}
=\lefttop{\nbar}\overline{\nbigg}_{\pi(\vecu)}.
\end{equation}
\hfill\qed
\end{rem}

%% file: a83.tex

\subsubsection{Some vanishing}

Let $s_i$ and $\bar{s}_i$ 
denote the left action of
$-\deldel_i\cdot z_i$ and $-\deldelbar_i\cdot \bar{z}_i$
($i=1,\ldots,n$) respectively.
Let us pick sections $f_1$ of $\gbige_{|\nbigx(\lambda_0,\epsilon_0)}$
and $f_2$ of $\gbige_{|\sigma(\nbigx(\lambda_0,\epsilon_0))}$
satisfying
$\bigl(s_i+\eigenmap(\lambda,u_a)\bigr)^N\cdot f_a=0$
for some $u_a\in\KMS(\nbige^0,i)$ $(a=1,2)$
and for some sufficiently large integer $N$.
Note that we also have
$(\bar{s}_2+\overline{\sigma^{\ast}\eigenmap(\lambda,u_2)})^N
\cdot \bar{f}_2=0$.

We put $F:=\gbigc(f_1,\overline{f_2})$.
We have the following:
\begin{equation}\label{eq;9.25.15}
 \Bigl(
 -\frac{\del}{\del z_i}z_i
 +\frac{\eigenmap(\lambda,u_1)}{\lambda}
 \Bigr)^N F=0.
\end{equation}
\begin{equation}\label{eq;9.25.16}
\Bigl(
-\frac{\del}{\del \bar{z}_i} \bar{z}_i
+\sigma^{\ast}
 \Bigl(
 \frac{\eigenmap(\lambda,u_2)}{\lambda}
 \Bigr)
 \Bigr)^N F=0.
\end{equation}
From the equality (\ref{eq;9.25.15}),
$F$ is of the following form:
\[
 F=z_i^{\lambda^{-1}\cdot\eigenmap(\lambda,u_1)-1}
 \cdot
 \sum_k a_k\cdot(\log z_i)^k.
\]
Here we have $\del a_k/\del z_i=0$.

From the formula (\ref{eq;9.25.16}),
we obtain the following:
\begin{equation}\label{eq;9.25.17}
 \Bigl(
-\frac{\del}{\del \bar{z}_i}
 \bar{z}_i
+\frac{\eigenmap(\lambda,u_2)}{\lambda}
 \Bigr)^N F=0.
\end{equation}
Here we have used the equalities:
\[
 \sigma^{\ast}(\lambda)=-\lambda^{-1},
\quad
 \sigma^{\ast} \overline{\Bigl(
 \frac{\eigenmap(\lambda,u)}{\lambda}
\Bigr)}
=\frac{\eigenmap(\lambda,u)}{\lambda}.
\]
The equation (\ref{eq;9.25.17}) implies
that $F$ is of the following form:
\[
 F=\bar{z}_i^{\lambda^{-1}\cdot\eigenmap(\lambda,u_2)-1}
 \cdot\sum_k b_k\cdot(\log \bar{z}_i)^k.
\]
Here we have $\del b_k/\del \bar{z}_i=0$.

\begin{lem} \label{lem;9.25.18}
Assume that $F$ is not $0$.
Then we have $u_1=u_2+b\cdot(1,0)$ for some integer $b$,
and $F$ is of the following form on $\nbigx(\lambda_0,\epsilon_0)$:
\[
 F=|z_i|^{2\lambda^{-1}\eigenmap(\lambda,u)-2}
 \cdot \bar{z}_i^{-b}
 \cdot
 \Bigl(
 \sum_k a_k\cdot \bigl(\log|z_i|^2\bigr)^k
 \Bigr).
\]
Here $a_k$ are independent of the variables $z_i$ and $\bar{z}_i$,
and we put $u=u_1$.
\end{lem}
\pf
Note that $F$ is a $C^{\infty}$-function on $X-D$.
In particular, it is not multi-valued.
It implies that
$\lambda^{-1}\eigenmap(\lambda,u_1)-\lambda^{-1}\eigenmap(\lambda,u_2)$
is a integer.
It implies $u_1-u_2\in\seisuu\times\{0\}$.
Thus we obtain the first claim.
The second claim also follows from the uni-valuedness
of $F$.
\hfill\qed

%% file: a84.2.tex

\subsubsection{The filtration and the decomposition}

Let $\lambda_0\in\cnum_{\lambda}$ be generic.
Recall that we have the generalized eigen decomposition
with respect to the monodromy actions,
for some small positive number $\epsilon_0$:
\[
 \nbige_{|(X-D)\times\Delta(\lambda_0,\epsilon_0)}
=\bigoplus_{\vecu\in\KMSoverline(\nbige^0,\nbar)}
 \nbige^{(\lambda_0)}_{\vecu}.
\]
Here $\nbige^{(\lambda_0)}_{\vecu}$
denotes the generalized eigen space 
corresponding to $\eigenmap^f(\lambda,\vecu)$.

\begin{lem}\label{lem;c12.3.1}
Let $U$ denote the disc $\Delta(\lambda_0,\epsilon_0)$.
The restriction of the sesqui-linear pairing $C_0$
to
$\nbige^{(\lambda_0)}_{\vecu_1\,U}
 \otimes
 \overline{\nbige^{(\sigma(\lambda_0))}_{\vecu_2\,\sigma(U)}}$
is trivial unless $\vecu_1=\vecu_2$.
\end{lem}
\pf
It immediately follows from Lemma \ref{lem;9.25.18}.
\hfill\qed

\vspace{.1in}

Let $\lambda_0$ be a point of $\cnum_{\lambda}$,
which is not necessarily generic.
For a sufficiently small positive number $\epsilon_0$,
we have the filtrations $\lefttop{i}\nbigfzero$ $(i\in\nbar)$
and the decomposition $\lefttop{i}\EEzero$ $(i\in \nbar)$
of $\nbige_{|\Delta(\lambda_0,\epsilon_0)\times(X-D)}$.
We may assume that any point
$\lambda\in\Delta^{\ast}(\lambda_0,\epsilon_0)$
is generic.

\begin{lem} \label{lem;c12.3.10}
Let $U$ denote the disc $\Delta(\lambda_0,\epsilon_0)$.
Let $\vecalpha=(\alpha_i)$ and $\vecbeta=(\beta_i)$ be elements of
$\cnum^{\ast\,n}$.
The restriction of the sesqui-linear pairing $C_0$
to $\EEzero(\vecalpha)\otimes
 \overline{\EEzero(\vecbeta)}$
is trivial unless
$\alpha_i=\bar{\beta}_i^{-1}$ for any $i$.
\end{lem}
\pf
Let $\lambda_1$ be any point of $\Delta^{\ast}(\lambda_0,\epsilon_0)$,
and $\epsilon_1$ be a positive number 
such that
$U_1:=\Delta(\lambda_1,\epsilon_1)
  \subset\Delta^{\ast}(\lambda_0,\epsilon_0)$.
Then we have the following decomposition:
\[
 \EEzero(\vecalpha)_{|U_1\times(X-D)}
=\bigoplus_{\substack{\vecu\in\KMSoverline(\nbige^0,\nbar)\\
 \eigenmap^f(\lambda_0,\vecu)=\vecalpha
 } }
 \nbige^{(\lambda_1)}_{\vecu},
\quad
 \EE^{(\sigma(\lambda_0))}(\vecbeta)_{|\sigma(U_1)\times(X-D)}
=\bigoplus_{\substack{\vecu\in\KMSoverline(\nbige^0,\nbar)\\
 \eigenmap^f(\sigma(\lambda_0),\vecu)=\vecbeta}
 }
 \nbige^{(\sigma(\lambda_1))}_{\vecu}.
\]
In the case
$\eigenmap^f(\lambda_0,\vecu_1)\neq
 \overline{\eigenmap^f(\sigma(\lambda_0),\vecu_2)}^{-1}$,
we have $\vecu_1\neq\vecu_2$.
Hence we obtain the vanishing on $U_1$
due to Lemma \ref{lem;c12.3.1}.
Then we obtain the vanishing on $U$
by using the argument given in the last part of
Lemma \ref{lem;9.25.3}.
\hfill\qed

\vspace{.1in}

Let us consider the set
$S(\alpha):=
\bigl\{u\in \KMSoverline(\nbige^0,i)\,\big|\,
  \eigenmap^f(\lambda_0,u)=\alpha\bigr\}$.
We have the two maps
\[
 \paramap^f(\lambda_0):S(\alpha)\lrarr\real,
\quad
 \paramap^f\bigl(\sigma(\lambda_0)\bigr):S(\alpha)\lrarr\real.
\]

\begin{lem} \label{lem;a11.22.10}
Let $u_1$ and $u_2$ be elements of $S(\alpha)$.
Then we have $\paramap^f(\lambda_0,u_1)>\paramap^f(\lambda_0,u_2)$
if and only if
we have
 $\paramap^f\bigl(\sigma(\lambda_0),u_1\bigr)<
  \paramap^f\bigl(\sigma(\lambda_0),u_2\bigr)$.
\end{lem}
\pf
It can be checked by a direct calculation.
\hfill\qed

\vspace{.1in}

\begin{lem} \label{lem;a12.6.3}
Let $u$ be an element of $S(\alpha)$ above.
We put $x:=\paramap^f(\lambda_0,u)$.
Let us consider the restriction of $C_0$
to the following:
\begin{equation}\label{eq;c12.3.2}
 \lefttop{i}\nbigfzero_{x}\EEzero(\vecalpha)
\otimes
 \overline{\lefttop{i}\nbigfzero_y\EEzero(\vecbeta)}.
\end{equation}
Here we assume $\vecalpha=\bar{\vecbeta}^{-1}$.

In the case $y<\paramap^f(\sigma(\lambda_0),u)$,
then the restriction of $C_0$ to {\rm(\ref{eq;c12.3.2})}
is $0$.
\end{lem}
\pf
Let us take $U_1=\Delta(\lambda_0,\epsilon_1)$
as in the proof of Lemma \ref{lem;c12.3.10}.
We have the following decompositions:
\[
 \lefttop{i}\nbigfzero_{x}
 \EEzero(\vecalpha)_{|U_1\times(X-D)}
=\bigoplus_{\paramap^f(\lambda_0,q_i(\vecu'))\leq x}
 \nbige^{(\lambda_1)}_{\vecu'},
\quad\quad
 \lefttop{i}\nbigfzero_y
 \EE^{(\sigma(\lambda_0))}_{|\sigma(U_1)\times(X-D)}
=\bigoplus_{\paramap^f(\sigma(\lambda_0),q_i(\vecu'))\leq y}
 \nbige^{(\sigma(\lambda_1))}_{\vecu'}.
\]
In the case 
$\paramap^f(\lambda_0,u_1)\leq \paramap^f(\lambda_0,u)$
and
$\paramap^f(\sigma(\lambda_0),u_2)\leq y<
 \paramap^f(\sigma(\lambda_0),u)$,
we have $u_1\neq u_2$
due to Lemma \ref{lem;a11.22.10}.
Then we obtain the vanishing 
as in the proof of Lemma \ref{lem;c12.3.10}.
\hfill\qed

\vspace{.1in}

\begin{lem}
The pairing
$C_0:\nbige_{\AAA}\otimes\overline{\nbige}_{\AAA}
\lrarr\distribution^{\AAA}_{X-D}$
induces the pairing
\begin{equation}\label{lem;c12.3.15}
 C_0:\lefttop{\nbar}\overline{\nbigg}_{\vecu\,\AAA}
 \otimes
 \overline{\lefttop{\nbar}\overline{\nbigg}_{\vecu\,\AAA}}
\lrarr
\distribution_{X-D}^{\AAA}.
\end{equation}
\end{lem}
\pf
It follows from Lemma \ref{lem;c12.3.10}.
(See the subsubsection \ref{subsubsection;a12.4.1}
 for the definition of
$\lefttop{\nbar}\overline{\nbigg}_{\vecu}(\nbige)$.)
\hfill\qed

\subsubsection{In the case $\paramap(\lambda_0,\vecu)<\vecdelta_{\nbar}$
 and $\paramap(\sigma(\lambda_0),\vecu)<\vecdelta_{\nbar}$}

Let $\vecu$ be an element of $\KMS(\nbige^0,\nbar)$.
Let $\lambda_0$ be a point of $\cnum_{\lambda}^{\ast}$.
We denote a small disc $\Delta(\lambda_0,\epsilon_0)$
by $U$.

Let $F$ and $G$ be sections of
$\lefttop{\nbar}\nbigg_{\vecu}(\nbige)$
over $X\times U$ and $X\times\sigma(U)$ respectively.
Recall we have the relation (\ref{eq;a12.4.2}).
Thus we can naturally regard
$F$ and $G$
as sections of
$\lefttop{\nbar}\overline{\nbigg}_{\pi(\vecu)}(\nbige)$
over $(X-D)\times U$ and $(X-D)\times \sigma(U)$
respectively.
Here $\pi$ denotes the natural projection
$\KMS(\nbige^0,\nbar)\lrarr\KMSoverline(\nbige^0,\nbar)$.
Hence we have the pairing
$C_0(F,G)$,
which gives an element of
$C^{\infty}(X-D,\nbigo(U))$.

On the other hand,
we can take the section $\tilde{F}$ of
$\nbige_{|(X-D)\times U}$
such that $\pi_1(\tilde{F})=F$.
Here $\pi_1$ denotes the following projection:
\[
 \pi_1:
 \lefttop{\nbar}\nbigfzero_{\paramap^f(\lambda_0,\vecu)}
 \EEzero(\eigenmap^f(\lambda_0,\vecu))
\lrarr
 \lefttop{\nbar}\Gr^{\nbigfzero}_{\paramap^f(\lambda_0,\vecu)}
 \EEzero(\eigenmap^f(\lambda_0,\vecu)).
\]
We can pick $\tilde{F}$
such as it is the section of
$\naiveprolongg{\paramap(\lambda_0,\vecu)}{\nbige}_{|X\times U}$.
(See the subsubsection \ref{subsubsection;10.18.10}).
In the case $\paramap(\lambda_0,\vecu)<\vecdelta$,
$\tilde{F}$ naturally gives the section of
$\Vzero_{\paramap(\lambda_0,\vecu)-\vecdelta_{\nbar}}(\gbige)
   _{|X\times U}$.

Similarly we can pick a section $\tilde{G}$
of $\Vzero_{\paramap(\sigma(\lambda_0),\vecu)-\vecdelta_{\nbar}}(\gbige)
   _{|X\times \sigma(U)}$
for $G$.

\begin{lem}\label{lem;a12.4.21}
$C_0(F,G)$ naturally gives the
$\nbigo(\AAA)$-valued distribution 
$\Phi$ on $X$,
and we have
$\Phi=\gbigc(\tilde{F},\tilde{G})$.
\end{lem}
\pf
It is clear from the definitions
of $\gbigc(\tilde{F},\tilde{G})$.
\hfill\qed

%% file: a85.2.tex

\subsubsection{The induced sesqui-linear pairing
 on $\lefttop{\nbar}\nbigg_{\hat{\vecu}}(E)$}
\label{subsubsection;a12.4.20}

Let $\vecu$ be an element of $\KMS(\nbige^0,\nbar)$.
We put $\hat{\vecu}=\vecu+\vecdelta_{0,\nbar}$.
Let $f$ be sections of
$\lefttop{\nbar}\nbigg_{\hat{\vecu}}(E)(\AAA)$.
We have the sections 
$F$ of $\lefttop{\nbar}\nbigg_{\hat{\vecu}}(\nbige)_{|X\times\AAA}$
satisfying the following:
\begin{itemize}
\item
 We have
 $F_{|\{O\}\times \AAA}=f$
 under the isomorphism
 $\lefttop{\nbar}\nbigg_{\hat{\vecu}}(\nbige)
  _{\{O\}\times\cnum_{\lambda^{\ast}}}
 \simeq \lefttop{\nbar}\nbigg_{\hat{\vecu}}(E)_{|\cnum_{\lambda}^{\ast}}$.
\item
 We have the following vanishing, for any $i$
and for any sufficiently large integer $N$:
\[
 \Bigl(
 -\deldel_iz_i+\eigenmap(\lambda,u_i)
 \Bigr)^NF=0.
\]
\end{itemize}
The section $F$ is called the lift of $f$.

Let $f$ and $g$ be sections of
$\lefttop{\nbar}\nbigg_{\hat{\vecu}}(E)$
over $\AAA$.
Let $F$ and $G$ be the lifts of $f$ and $G$
respectively.
We can naturally regard $F$ and $G$
as the sections of
$\lefttop{\nbar}\overline{\nbigg}_{\pi(\hat{\vecu})}(\nbige)$
over $(X-D)\times \AAA$.
Thus we have the pairing
$C_0(F,G)$.

\begin{lem}\label{lem;a12.4.5}
The function $C_0(F,G)$ is of the following form:
\begin{equation}\label{eq;a12.4.6}
 C_0(F,G)
=\prod_{i=1}^n |z_i|^{2\lambda^{-1}\eigenmap(\lambda,u_i)-2}
 \cdot 
 \Bigl(
 \sum_{\vecn\in\seisuu_{\geq 0}^n}
 a_{\vecn}(\lambda)\cdot \prod_i\bigl(\log|z_i|^2\bigr)^{n_i}
 \Bigr).
\end{equation}
Here $a_{\vecn}(\lambda)$ denote holomorphic functions
on $\AAA$.
\end{lem}
\pf
It immediately follows from Lemma \ref{lem;9.25.18}.
\hfill\qed

\vspace{.1in}
\begin{df}\label{df;a12.4.22}
We put as follows:
\[
 \lefttop{\nbar}\Psi_{\hat{\vecu}}(C_0)(f,g):=a_0.
\]
Here $a_0$ is given in the development {\rm(\ref{eq;a12.4.6})}.
Thus we obtain the sesqui-linear pairing:
\[
 \lefttop{\nbar}\Psi_{\hat{\vecu}}(C_0):
 \lefttop{\nbar}\nbigg_{\hat{\vecu}}(E)_{\AAA}
\otimes
 \overline{\lefttop{\nbar}\nbigg_{\hat{\vecu}}(E)_{\AAA}}
\lrarr\nbigo(\AAA).
\]
\hfill\qed
\end{df}

\vspace{.1in}
Let $N$ be any integer.
The multiplication of $(\prod_{i=1}^n z_i)^N$
induces the isomorphism
$\lefttop{\nbar}\nbigg_{\hat{\vecu}}(E)
\lrarr
 \lefttop{\nbar}\nbigg_{\hat{\vecu}-N\cdot\vecdelta_{0,\nbar}}(E)$.

\begin{lem}\label{lem;a12.4.10}
Under the isomorphism
$\lefttop{\nbar}\nbigg_{\hat{\vecu}}(E)
\lrarr
 \lefttop{\nbar}\nbigg_{\hat{\vecu}-N\cdot\vecdelta_{0,\nbar}}(E)$
above,
we have
$\lefttop{\nbar}\Psi_{\hat{\vecu}}(C_0)
=\lefttop{\nbar}\Psi_{\hat{\vecu}-N\vecdelta_{0,\nbar}}(C_0)$.
\end{lem}
\pf
It is clear from the definition.
\hfill\qed

\subsubsection{Comparison of the induced sesqui-linear pairings}

Let $\vecu$ be an element of $\KMS(\nbige^0,\nbar)$
such that any $i$-th components are not contained in
$\seisuu\times\{0\}$.
We put $\hat{\vecu}=\vecu+\vecdelta_{0,\nbar}$.
We have the isomorphism
$\lefttop{\nbar}\tildepsi_{\vecu}(\gbige)
\simeq
 \lefttop{\nbar}\nbigg_{\hat{\vecu}}(E)$
by definition
(see the subsubsection \ref{subsubsection;b11.23.20}).
Thus we obtain the following sesqui-linear pairing:
\[
 \lefttop{\nbar}\Psi_{\hat{\vecu}}(C_0):
 \lefttop{\nbar}\tildepsi_{\vecu}(\gbige)_{\AAA}
\otimes
 \overline{\lefttop{\nbar}\tildepsi_{\vecu}(\gbige)_{\AAA}}
\lrarr \nbigo(\AAA).
\]

On the other hand,
we have the following pairing,
due to  Lemma \ref{lem;b11.23.10}:
\[
 \lefttop{\nbar}\tildepsi_{\vecu}(\gbigc):
 \lefttop{\nbar}\tildepsi_{\vecu}(\gbige)_{\AAA}
\otimes
 \overline{\lefttop{\nbar}\tildepsi_{\vecu}(\gbige)_{\AAA}}
\lrarr \nbigo(\AAA).
\]

\begin{prop}\label{prop;a12.4.25}
We have
$\lefttop{\nbar}\tildepsi_{\vecu}(\gbigc)
=
\lefttop{\nbar}\Psi_{\hat{\vecu}}(C_0)$.
\end{prop}
\pf
We have only to compare them
on a neighbourhood of any point $\lambda_0\in \AAA$.
Due to Lemma \ref{lem;a12.4.10}
and the definition of $\lefttop{\nbar}\tildepsi_{\vecu}\gbigc$,
we may assume $\paramap(\lambda_0,\vecu)<0$ and
$\paramap(\sigma(\lambda_0),\vecu)<0$.

Let $U$ denote a small disc
$\Delta(\lambda_0,\epsilon_0)$.
Let $f$ and $g$ be sections of
$\lefttop{\nbar}\psi_{\vecu}(\gbige)$
on $U$ and $\sigma(U)$ respectively.
Let $F$ and $G$ be the lifts of $f$ and $g$
in the sense of the subsubsection \ref{subsubsection;a12.4.20}.
Due to Lemma \ref{lem;a12.4.21},
Lemma \ref{lem;a12.4.5} and Definition \ref{df;a12.4.22},
the proposition can be reduced
to the following lemma.

\begin{lem}\label{lem;b11.23.50}
We have $\lefttop{\nbar}\psizero_{\vecu}\gbigc(f,\bar{g})=a_0$.
Here $a_0$ is given in the development
{\rm(\ref{eq;a12.4.6})}.
\end{lem}
\pf
By the definition of $\lefttop{\nbar}\psizero_{\vecu}\gbigc$,
we have the following:
\begin{multline}
 \lefttop{\nbar}\psizero_{\vecu}\gbigc(f,\bar{g})
=\Res_{s_1+\alpha(\lambda,u_1)}
 \Res_{s_2+\alpha(\lambda,u_2)}
 \cdots
 \Res_{s_n+\alpha(\lambda,u_n)}
 \Bigl\langle
 \gbigc(F,\bar{G}),
 \chi\cdot
 \prod_{i=1}^n|z_i|^{2s_i}
 \frac{\sqrt{-1}}{2\pi}dz_i\wedge d\bar{z_i}
 \Bigr\rangle
\end{multline}
Here we put
$\alpha(\lambda,u_i)
:=\lambda^{-1}\cdot\eigenmap(\lambda,u_i)$,
and $\chi$ denotes a $C^{\infty}$-function on $\cnum^n$
as follows:
\begin{condition} \mbox{{}}\label{condition;b11.23.61}
\begin{itemize}
\item
 The support of $\chi$ is compact.
\item
 $\chi$ is identically $1$ around the origin.
\item
 $\chi$ depends only on $|z_i|$ $(i=1,\ldots,n)$.
\end{itemize}
\end{condition}
Then Lemma \ref{lem;b11.23.50} can be reduced to the following lemma.

\begin{lem}\label{lem;b11.23.66}
We have the following formula:
\begin{equation}\label{eq;b11.23.62}
 \Res_{s+\alpha(\lambda,u)}
 \int\chi\cdot |z|^{s+2\alpha(\lambda,u)-2}
 \cdot(\log|z|^2)^m
 \cdot\frac{\sqrt{-1}}{2\pi}   dz\wedge d\bar{z}
=\left\{
 \begin{array}{ll}
 1 & (m=0)\\
 \mbox{{}}\\
 0 & (m\neq 0).
 \end{array}
 \right.
\end{equation}
Here we put $\alpha(\lambda,u):=\eigenmap(\lambda,u)$
for $u\in\real\times\cnum$,
and $\chi$ denotes a $C^{\infty}$-function on $\cnum^1$
as in Condition {\rm\ref{condition;b11.23.61}}.
\end{lem}
\pf
We put $T:=s+\alpha(\lambda,u)$,
the left hand side of (\ref{eq;b11.23.62}) can be rewritten as follows:
\[
\Res_{T=0}
 \int
 \chi\cdot |z|^{2T-2}\cdot(\log |z|^2)^m
 \cdot \frac{\sqrt{-1}}{2\pi}dz\wedge d\bar{z}.
\]
Then we obtain the formula (\ref{eq;b11.23.62})
from Corollary \ref{cor;b11.23.65}.
Thus we obtain Lemma \ref{lem;b11.23.66},
Lemma \ref{lem;b11.23.50}
and
thus Proposition \ref{prop;a12.4.25}.
\hfill\qed

%% file: a85.3.tex

\subsubsection{Preliminary}\label{subsubsection;b11.23.70}

Let $V$ be a vector bundle over $\proj^1$,
and let $S:V\otimes\sigma(V)\lrarr \Tate(0)$ is a perfect symmetric pairing.
Since it is symmetric, it satisfies
$\sigma^{\ast}S(f,\sigma^{\ast}(g))=
 S(g,\sigma^{\ast}(f))$.
We have the perfect strict $\nbigr$-triple
$\Theta(V^{\lor})=\bigl(V_0,\sigma^{\ast}(V^{\lor}_{\infty}),C_V\bigr)$
(the subsubsection \ref{subsubsection;b11.23.1}).
The perfect pairing
$S:V_0\otimes \sigma^{\ast}(V_{\infty})\lrarr\nbigo$,
induces the isomorphism
$\rho_2:\sigma^{\ast}(V_{\infty}^{\lor})\lrarr V_0$.

\begin{lem}\label{lem;b11.23.2}
For any $F\in V_0(\AAA)$ and
for any $G\in \sigma(V_{\infty}^{\lor})(\AAA)$,
we have the following:
\[
 C_V(F,\bar{G})=S(F,\sigma^{\ast}\rho_2(G)).
\]
\end{lem}
\pf
Let $\langle\cdot,\cdot\rangle$ denote
the pairing of the vector bundle and its dual.
We have the following equalities:
\[
 C_V(F,\bar{G})=\langle F,\sigma^{\ast}G\rangle
=\sigma^{\ast}\langle \sigma^{\ast}F,G\rangle
=\sigma^{\ast}S\bigl(\rho_2(G),\sigma^{\ast}F\bigr)
=S(F,\sigma^{\ast}\rho_2(G)).
\]
Thus we are done.
\hfill\qed

\vspace{.1in}
For sections $F,G\in V_0(\AAA)$,
we put $C_S(F,\bar{G}):=S(F,\sigma^{\ast}G)$,
which gives the pairing
$V_0(\AAA)\otimes \overline{V_0(\AAA)}\lrarr \nbigo(\AAA)$.
Thus we obtain the perfect $\nbigr$-triple
$(V_0,V_0,C_S)$.
The isomorphism
$\rho_2:\sigma^{\ast}(V_{\infty}^{\lor})
\simeq V_0$
induces the isomorphism
$\Theta(V^{\lor})\simeq(V_0,V_0,C_S)$,
due to Lemma \ref{lem;b11.23.2}.

We have the isomorphism
$f_S:V^{\lor}\lrarr \sigma^{\ast}\bigl(V^{\lor}\bigr)$
induced by $S$.
It induces the isomorphism
$\Theta(f_S):\Theta(V^{\lor})
 \simeq \Theta\bigl(\sigma^{\ast}\bigl(V^{\lor}\bigr)\bigr)
\simeq \Theta(V^{\lor})^{\ast}$.
We also have the natural morphism
$g_S:(V_0,V_0,C_S)\lrarr (V_0,V_0,C_S)^{\ast}=(V_0,V_0,C_S)$.
The following lemma is clear from our construction.
\begin{lem}
Under the isomorphism
$\Theta\bigl(V^{\lor}\bigr)\simeq (V_0,V_0,C_S)$,
we have $\Theta(f_S)=g_S$.
\hfill\qed
\end{lem}

%% file: a85.4.tex

\subsubsection{The $\nbigr$-triple 
   $\Theta\bigl(S^{\can}_{-\hat{\vecu}}(E)\bigr)$}

Let $\vecu$ be an element of $\prod_{i}\KMS(\nbige^0,i)$
such that
any $i$-th component is not contained in
$\seisuu_{\geq\,0}\times\{0\}$.
For simplicity we put $\hat{\vecu}:=\vecu+\vecdelta_{0,\nbar}$.
Applying the construction in the subsubsection
\ref{subsubsection;b11.23.70}
for $V=S^{\can}_{\hat{\vecu}}(E)$
and the pairing
$S:S^{\can}_{\hat{\vecu}}(E)\otimes
   \sigma^{\ast}S^{\can}_{\hat{\vecu}}(E)  \lrarr \Tate(0)$,
we obtain the isomorphism:
\[
 \Theta\bigl(
 S^{\can}_{-\hat{\vecu}}(E^{\lor})
 \bigr)
\simeq
 \bigl(
 \lefttop{\nbar}\nbigg_{\hat{\vecu}}(E),
 \lefttop{\nbar}\nbigg_{\hat{\vecu}}(E),
 C_S
 \bigr).
\]
Here $C_S$ can be calculated as follows:
Let $f$ and $g$ be sections
of $\lefttop{\nbar}\nbigg^{(\lambda_0)}_{\hat{\vecu}}(E)$
and $\lefttop{\nbar}\nbigg^{(\sigma(\lambda_0))}_{\hat{\vecu}}(E)$.
We take the corresponding multi-valued holomorphic sections
$F_0$ and $G_0$ of
$\lefttop{\nbar}\nbigg^{(\lambda_0)}_{\hat{\vecu}}(\nbige)$
and $\lefttop{\nbar}\nbigg^{(\sigma(\lambda_0))}_{\hat{\vecu}}(\nbige)$,
namely,
we have
$\Phi^{\can}_{\hat{\vecu}}(F_0)=f$ and
$\Phi^{\can}_{\hat{\vecu}}(G_0)=g$.
\begin{lem}
We have
$C_S(f,\bar{g})=h(F_0,\sigma^{\ast}G_0)$.
Here $h$ denotes the hermitian metric 
of the given tame harmonic bundle
$\harmonicbundle$.
\end{lem}
\pf
It is clear from our construction.
See Lemma \ref{lem;b11.23.2}.
\hfill\qed

\subsubsection{The isomorphism}

Since we have the isomorphism
$\lefttop{\nbar}\tildepsizero_{\vecu}(\gbige)
\simeq
 \lefttop{\nbar}\nbigg_{\hat{\vecu}}(E)$
(the subsubsection \ref{subsubsection;b11.23.20}),
we would like to compare the sesqui-linear pairings
and the polarizations under the isomorphism.

\begin{lem}
We have $C_S(f,\bar{g})=\lefttop{\nbar}\tildepsizero_{\vecu}\gbigc(f,\bar{g})$.
\end{lem}
\pf
We have the formula of the following form:
\[
 F=\prod_{i=1}^n z_i^{-\lambda^{-1}\eigenmap(\lambda,u_i)-1}
 \Bigl(
 F_0+\sum_{\vecn\neq 0}F_{\vecn}\prod_{i=1}^n
 \bigl(\log z_i\bigr)^{n_i}
 \Bigr)
\]
\[
G=\prod_{i=1}^n \bar{z}_i^{-\lambda^{-1}\eigenmap(\lambda,u_i)-1}
 \Bigl(
 G_0+\sum_{\vecn\neq 0}G_{\vecn}\prod_{i=1}^n
 \bigl(\log \bar{z}_i\bigr)^{n_i}
 \Bigr).
\]
Here $F_{\vecn}$ and $G_{\vecn}$ denote multi-valued flat sections.
Then the term $a_0$ in the development (\ref{eq;a12.4.6})
is $h(F_0,\sigma^{\ast}G_0)$.
Thus we are done.
\hfill\qed

\begin{cor}
We have the isomorphism
$\Theta\bigl(S_{-\hat{\vecu}}(E^{\lor})\bigr)
\simeq
 \lefttop{\nbar}\tildepsi_{\vecu}
 \bigl(
 \gbige,\gbige,\gbigc
 \bigr)$,
which induces the isomorphism
of the polarized mixed twistor structure.
\hfill\qed
\end{cor}

\begin{cor}\label{cor;b12.3.250}
We have
$S(f,\sigma^{\ast}g)=\lefttop{\nbar}\tildepsi_{\vecu}
 \gbigc(f,\bar{g})$.
\hfill\qed
\end{cor}

%% file: a34.6.tex

This subsection is a continuation
of the subsubsection \ref{subsubsection;a11.30.20}.

\subsubsection{Preliminary}

\label{subsubsection;b11.9.1}

Let us calculate the induced sesqui-linear pairing
$\tildepsi_{t,u}\bigl((\sqrt{-1}N)^h\otimes\id\bigr)$
on the component $\nbigm_{\nbar}$ of
$P_h\Gr^{W(N)}_h\tildepsi_{t,u}\gbige[\deldel_t]$
(see the subsubsection \ref{subsubsection;a11.30.20}).
We would like to show that
the $\nbigr$-triple
$\bigl(
 \nbigm_{\nbar},\nbigm_{\nbar},
 \tildepsi_{t,u}\gbigc\bigl((\sqrt{-1}N)^h\otimes\id\bigr)
 \bigr)$
is a polarized pure twistor structure of weight $h$
in the sense of Sabbah.
(Note the support is $\{O\}$ in this case.)
Since $\nbigm_{\nbar}$ is isomorphic to
the push-forward of $\nbigc_{\nbar\,0}:=\nbigc_{\nbar\,|\,\cnum_{\lambda}}$
(Lemma \ref{lem;b12.3.165}),
we  calculate the sesqui-linear pairings of
the sections of $\nbigc_{\nbar\,0}$.
We would like to use
Lemma \ref{lem;b12.3.170}
and Lemma \ref{lem;a11.30.41}.

\begin{rem} \label{rem;a12.3.280}
{\rm
In Lemma {\rm \ref{lem;b12.3.170}},
the result is stated for our polarized mixed twistor structure.
In particular, the nilpotent maps 
are obtained from the residues $\Res_i(\DD)$ there.
In this section,
we use the nilpotent parts of
the left action of $-\deldel_iz_i$.
Hence the signatures are reversed.
\hfill\qed
}
\end{rem}

\vspace{.1in}

We put as follows:
\[
 L_k(z):=\frac{\bigl(\log|z|^2\bigr)^{k}}{k!}.
\]

Let $f$ be a section of
$\nbigc_{\nbar}(\AAA)$
and $g$ be a section of
$\nbigc_{\nbar}(\AAA)$.
We denote the nilpotent part of $-\deldel_tt$ by $N$.
We denote the nilpotent part of $-\deldel_iz_i$ by $N_i$.
We also denote the nilpotent part of $m_i^{-1}(-\deldel_iz_i)$
by $\tilde{N}_i$.

We have the description
$f=\sum_{a=0}^{n-1} (\sqrt{-1} N)^af_a$
and $g=\sum_{b=0}^{n-1}(-\sqrt{-1}N)^bg_b$,
where $f_a$ and $g_b$ are sections of
$\lefttop{\vecn}\nbigg_{\hat{\vecu}}(E)$,
where we put
$\hat{\vecu}=\bigl(m_i\cdot u\,\big|\,i\in\nbar\bigr)
 +\vecdelta_{0,\nbar}$.

Let us consider the pairing
$\tildepsi_{t,u}\gbigc\bigl((\sqrt{-1}N)^k\cdot f,\bar{g}\bigr)$.
First we calculate the function on a neighbourhood
of a point $\lambda_0\in\AAA$.
We may assume that $\paramap(\lambda_0,u)<0$.
We put $u_i=m_i\cdot u$.

\vspace{.1in}

Let us pick a point $\lambda_0$ of $\AAA$.
We can take the  holomorphic sections
$F_a$ of
$\lefttop{\nbar}\nbigg_{\hat{\vecu}} (\nbige)$
$(a=0,\ldots,n-1)$,
such that the following holds:
\begin{itemize}
\item
$F_{a\,|\,\AAA\times \{O\}}=f_a$
\item
For any $a$ and $i$ and for sufficiently large $M$,
$(-\deldel_iz_i+\eigenmap(\lambda,u_i))^MF_a=0$ holds.
\end{itemize}
(See the subsubsection \ref{subsubsection;10.18.10}
for $\lefttop{\nbar}\nbigg_{\hat{\vecu}}(\nbige)$).

Similarly we can pick the sections $G_b$ of
$\lefttop{\nbar}\nbigg_{\hat{\vecu}}
   ^{(\sigma(\lambda_0))}(\nbige)$
for $g_b$ $(b=0,\ldots,n-1)$.

Note the following equality:
\[
 \tildepsi_{t,u}\bigl(
 (\sqrt{-1}N)^{i+k} f_i,\,\,
 \overline{(-\sqrt{-1}N)^j\cdot g_j}\bigr)
=(-\lambda^2)^{-j}
 \tildepsi_{t,u}\bigl(
 (\sqrt{-1}N)^{i+j+k}\cdot f_i,
 \,\,\overline{g_j}\bigr).
\]

We have the following equality:
\begin{multline}
\Big\langle
 i_+\gbigc\Bigl(
  \bigl(s+\eigenmap(\lambda,u)\bigr)^{i+j+k}\cdot F_i,
 \overline{G_j}
 \Bigr),\,\,
 |t|^{2T}\chi(t)\cdot\varphi\wedge\frac{i}{2\pi}dt\wedge d\bar{t}
\Big\rangle \\
=\Big\langle
 i_+\gbigc\Bigl(
 F_i,
 \overline{G_j}
 \Bigr),\,\,
 \bigl(
 \bigl(t\deldel_t+\eigenmap(\lambda,u)\bigr)^{i+j+k}
 |t|^{2T}
 \bigr)\cdot
 \chi(t)\cdot\varphi\wedge\frac{i}{2\pi}dt\wedge d\bar{t}
\Big\rangle
+\Big\langle
 i_+\gbigc\Bigl(
 F_i,
 \overline{G_j}
 \Bigr),\,\,
 |t|^{2T}\cdot
 \chi_1(t)\cdot\varphi\wedge\frac{i}{2\pi}dt\wedge d\bar{t}
\Big\rangle.
\end{multline}
Here $\chi_1(t)$ denotes a $C^{\infty}$-function
which vanishes identically around $t=0$.
The last term is entire function of $T$,
thus we can ignore it.
We have the following (see the section 3.7 in \cite{sabbah}):
\begin{multline}
\Big\langle
 i_+\gbigc\Bigl(
 F_i,
 \overline{G_j}
 \Bigr),\,\,
 \bigl(
 \bigl(t\deldel_t+\eigenmap(\lambda,u)\bigr)^{i+j+k}
 |t|^{2T}\bigr)
\cdot
 \chi(t)\cdot\varphi\wedge\frac{i}{2\pi}dt\wedge d\bar{t}
\Big\rangle \\
=\lambda^{i+j+k}\cdot 
 \Bigl(
 T+\frac{\eigenmap(\lambda,u)}{\lambda}
 \Bigr)^{i+j+k}\cdot
 \Big\langle
 \gbigc(F_i,\overline{G_j}),\,\,\,\,
 \prod_{i=1}^n|z_i|^{2m_iT}\cdot
 \chi\Bigl(
 \prod_{i=1}^n|z_i|^{2m_i}
 \Bigr)
 \wedge \varphi
 \Big\rangle.
\end{multline}
We may assume that
$\chi\bigl(
 \prod_{i=1}^n|z_i|^{2m_i}
 \bigr)=1 $
on a neighbourhood of the support of $\varphi$.
Thus we may ignore $\chi\bigl(
 \prod_{i=1}^n|z_i|^{2m_i} \bigr)$ in the following.

Due to our choice of $F_i$ and $G_j$,
there exist holomorphic functions $a_{\vecn,k}(\lambda)$
$\bigl(k\in\seisuu,\vecn\in\seisuu^n_{\geq\,0}\bigr)$
such that the following holds
(see Lemma \ref{lem;9.25.18}):
\[
 \sum_{a+b=k-h}(-1)^b\cdot\lambda^{k-2b}
 \cdot\gbigc(F_a,G_b)
=\prod_{i=1}^n
 |z_i|^{2\lambda^{-1}\cdot\eigenmap(\lambda,u_i)-2}
\cdot
 \sum_{\vecn}a_{\vecn,k}
 \prod_{i=1}^n
 L_{n_i}(z_i).
\]
\begin{lem}
We have the following,
for any $\vecl=(l_1,\ldots,l_n)\in\seisuu^n_{\geq\,0}$:
\[
 (-\lambda)^{|\vecn|}
 a_{\vecn,k}=
 \sum
 (-1)^{b}\cdot\lambda^{k-2b}
\sum_{a+b=k}
 \lefttop{\nbar}\tildepsi_{\vecu}
 \gbigc\Bigl(f_a\cdot\!\! \prod_{i=1}^nN_i^{n_i},\,\,
 g_b\Bigr)
\]
\end{lem}
\pf
We have the following:
\[
 \sum_{a+b=k}(-1)^b\lambda^{k-2b}\cdot
 \gbigc\Bigl(
 \prod_{i=1}^n(-\deldel_iz_i+\eigenmap(\lambda,u_i))^{l_i}
 \cdot F_a,G_b
 \Bigr)
=\prod_{i=1}^n
 |z_i|^{2\lambda^{-1}\cdot\eigenmap(\lambda,u_i)-2}
\cdot
 \sum_{\vecn\geq \vecl}a_{\vecn,k}
 \prod_{i=1}^n
 L_{n_i-l_i}(z_i)
\times(-\lambda)^{|\vecl|}
\]
Then we obtain the result
by using  Corollary \ref{cor;b11.23.65} inductively.
\hfill\qed

\vspace{.1in}

Let us consider the following function:
\begin{multline} \label{eq;9.26.5}
 \sum_{a,b}(-1)^b\cdot\lambda^{a+h-b}\cdot
 (\sqrt{-1})^{a+b+h}
 \Bigl(
 T+\frac{\eigenmap(\lambda,u)}{\lambda}
 \Bigr)^{a+b+h}\cdot
 \gbigc(F_a,\overline{G_b})
\cdot \prod_{i=1}^n|z_i|^{2m_iT} \\
=\sum_k
 (\sqrt{-1})^{k}
 \Bigl(
 T+\frac{\eigenmap(\lambda,u)}{\lambda}
 \Bigr)^k\cdot
\sum_{a+b=k-h}(-1)^b\cdot\lambda^{k-2b}\cdot
 \gbigc\bigl(F_a,\overline{G_b}\bigr)
 \cdot\prod_{i=1}^n|z_i|^{2m_iT} \\
=\sum_k(\sqrt{-1})^{k}
 \Bigl(T+\frac{\eigenmap(\lambda,u)}{\lambda}
 \Bigr)^k\cdot
 \prod_{i=1}^n|z_i|^{
 2m_i(T+\lambda^{-1}\eigenmap(\lambda,u))-2}\cdot
 \sum_{\vecn}a_{\vecn,k}
 \prod_{i=1}^n
 L_{n_k}(z_i).
\end{multline}
We put $\tau=T+\lamda^{-1}\cdot\eigenmap(\lambda,u)$,
and then the right hand side of (\ref{eq;9.26.5})
can be rewritten as follows:
\begin{equation}
 \sum_k (\sqrt{-1})^k\cdot\tau^k\cdot \prod_{i=1}^n|z_i|^{2m_i\tau-2}
 \cdot
 \sum_{\vecn}a_{\vecn,k}\cdot
 \prod_{i=1}^n
 L_{n_i}(z_i).
\end{equation}
Then we have the following:
\begin{multline}\label{eq;b12.3.200}
 \Big\langle
 \tildepsi_{t,u}\bigl(
 (\sqrt{-1}N)^{i+k}f_i,\,\,
 \overline{
 (-\sqrt{-1}N)^{j}g_j
 }
 \bigr),
 \varphi
 \Big\rangle\\
=\Res_{\tau=0}
 \Big\langle
 \sum_k (\sqrt{-1})^k\cdot\tau^k\cdot
 \prod_{i=1}^n|z_i|^{2m_i\tau-2}\cdot L_{n_i}(z_i)
\cdot\sum_{\vecn}a_{\vecn,k},\,\,
\varphi
 \Big\rangle.
\end{multline}
We put as follows:
\[
 \nbigu:=\Bigl\{(\vecn,k)\in\seisuu_{\geq \,0}^n\times\seisuu_{\geq\,0}
 \,\Big|\,
 |\vecn|=k-n+1\Bigr\}.
\]
Then the distribution in the right hand side of (\ref{eq;b12.3.200})
 can be regarded as the product
of the delta function at $\{O\}$
and the following holomorphic function of $\lambda$
(Lemma \ref{lem;a11.30.20}):
\begin{multline}
 \sum_{(\vecn,k)\in \nbigu}
 (\sqrt{-1})^k
 a_{\vecn,k}\cdot (-1)^{|\vecn|}\cdot\prod_{i=1}^{n}m_i^{-n_i-1} \\
=\prod_{i=1}^n m_i^{-1}\cdot
 (\sqrt{-1})^{n-1}
\cdot
 \sum_{(\vecn,k)\in\nbigu}
  \lambda^{k-|\vecn|}
\sum_{a+b=k}
 (-\lambda^{-2})^b\cdot
 \lefttop{\nbar}\tildepsi_{\vecu}
 \gbigc\Bigl(\prod_{i=1}^n(\sqrt{-1}\tilde{N}_i)^{n_i}\cdot f_a,
 \,\,
 g_b\Bigr) \\
=\prod_{i=1}^nm_i^{-1}\cdot
 \sum_{(\vecn,k)\in\nbigu}
 (\sqrt{-1}\lambda)^{n-1-2b}\cdot
  \lefttop{\nbar}\tildepsi_{\vecu}
 \gbigc\Bigl(\prod_{i=1}^n(\sqrt{-1}\tilde{N}_i)^{n_i}\cdot f_a,
 \,\,
 g_b\Bigr).
\end{multline}
Due to Corollary \ref{cor;b12.3.250}, we obtain the following:
\begin{lem} \label{lem;a11.30.40}
Let $f=\sum f_i\cdot (\sqrt{-1}N)^i$
and $g=\sum g_i\cdot (-\sqrt{-1}N)^j$ be sections of
$\nbigc_{\nbar}(\AAA)$.
We have the following formula:
\[
 \tildepsi_{t,u}\gbigc\bigl( (\sqrt{-1}N)^h\cdot f,\overline{g}
 \bigr)\\
=\prod_{a=1}^n m_a^{-1}\cdot
 \sum_{-|\vecn|-n+i+j+h=-1}
 (\sqrt{-1}\lambda)^{n-1-2j}
 \cdot
 S\Bigl(
 \prod_{a=1}^n(\sqrt{-1}\tilde{N}_a)^{n_a}
 f_i,\sigma^{\ast}{g}_j
 \Bigr)
\]
Here $S$ denotes the pairing of
$\lefttop{\nbar}\nbigg_{\vecu+\vecdelta_{0,\nbar}}$.
\hfill\qed
\end{lem}

\begin{cor}\label{cor;a11.30.50}
The $\nbigr$-triple $(\nbigm_{\nbar},\nbigm_{\nbar},C')$ is
the polarized pure twistor module of weight $h$.
\end{cor}
\pf
It follows from Lemma \ref{lem;a11.30.40},
Lemma \ref{lem;a11.30.41},
Lemma \ref{lem;b12.3.170},
and Remark \ref{rem;a12.3.280}.
\hfill\qed

\subsubsection{Consequences}
\label{subsubsection;04.2.4.10}

We have the pairing:
\[
 \Gr^{W(N)}\tildepsi_{t,u}\gbigc\bigl( 
 (\sqrt{-1}N)^h\cdot,\cdot 
 \bigr):
 P_h\Gr^{W(N)}_h\tildepsi_{t,u}\gbige_{\AAA}
\otimes
\overline{
 P_h\Gr^{W(N)}_{h}
 \tildepsi_{t,u}\gbige_{\AAA}}
\lrarr
 \distribution^{\AAA}_X.
\]

Let $I$ be a subset of $\nbar$.
Let $\nbigm_{I}$ denote the component of
$\Gr^{W(N)}_h\tildepsi_{t,u}\gbige$,
whose strict support is $\nbigd_I$.
Let $\gbigc_{I}$ denote the restriction
of $\Gr^{W(N)}\tildepsi_{t,u}\gbigc\bigl((\sqrt{-1}N)^h\cdot,\cdot\bigr)$
to $\nbigm_{I\,\AAA}\otimes\overline{\nbigm_{I\,\AAA}}$.
Then we obtain the $\nbigr$-triple
$\bigl(\nbigm_{I},\nbigm_{I},\gbigc_{I}\bigr)$.

Since the support of $\nbigm_I$ is $\nbigd_I$,
and since $\nbigm_I$ is strictly $S$-decomposable
along $z_i=0$ for any $i\in\nbar$,
we have the $\nbigr$-triple 
$(\nbigm_I',\nbigm_I',\gbigc_I)$ on $\nbigd_I$
such that whose push out with respect to
the inclusion $\nbigd_I\lrarr \nbigx$
is $(\nbigm_I,\nbigm_I,\gbigc_I)$.
It is easy to see that $\nbigm_I$ is smooth
on $\nbigd_I^{\circ}=\nbigd_I-\bigcup_{I'\supsetneq I}\nbigd_{I'}$.

\begin{lem}
The restriction of $(\nbigm_I,\nbigm_I,\gbigc_I)$
to $\nbigd_I^{\circ}$ is a variation of pure twistors
of weight $h$.
\end{lem}
\pf
Let $q_I$ denote the projection of $X$ onto $D_I$,
and the induced projection of $\nbigx$ onto $\nbigd_I$.
Let $Q$ be a point of $D_I^{\circ}$.
By considering the restriction of
$(\nbigm_I,\nbigm_I,\gbigc_{I})$
to the hyperplane $q_I^{-1}(Q)$,
and by applying Corollary \ref{cor;a11.30.50},
we obtain the result.
\hfill\qed

\vspace{.1in}
Then we obtain the harmonic bundle
$(E_{I},\delbar_{E_I},\theta_I,h_I)$
on $D_I^{\circ}$,
such that the following holds:
\[
 \bigl(
 \nbige(E_I),
 \nbige(E_I),
 C_0
 \bigr)
\otimes\nbigo(h)
 \simeq
 \bigl(
 \nbigm_I,\nbigm_I,\gbigc_I
 \bigr)_{|\,D^{\circ}_I}.
\]
Then we obtain the $\nbigr$-triple 
$(\gbige(E_I),\gbige(E_I),\gbigc)$
on $D_I$.

\begin{prop}\label{prop;b12.4.1}
The $\nbigr$-triple
$(\gbige(E_I),\gbige(E_I),\gbigc)\otimes\nbigo(h)$
is naturally isomorphic to
$(\nbigm_I,\nbigm_I,\gbigc_I)$.
\end{prop}
\pf
Due to Lemma \ref{lem;a11.30.60} and
Proposition \ref{prop;a11.26.2},
$\gbige(E_I)$ and $\nbigm_I$ are naturally isomorphic.
Due to Corollary \ref{cor;a11.23.6},
we obtain the coincidence of the sesqui-linear pairings
under the isomorphism.
Thus we are done.
\hfill\qed

%% file: a77.tex

\subsubsection{The definitions}
\label{subsubsection;04.2.16.200}

Let $X$ be a complex manifold and $Z$ be an irreducible closed subset of
$X$. 

\begin{df} \label{df;a11.23.60}
A generically defined variation of polarized pure twistor structures 
of weight $w$ on $Z$ is defined to be the data as follows:
\begin{itemize}
\item
 A Zariski open smooth subset $U$ of $Z$.
\item
 A smooth polarized pure twistor structure $\nbigt$ on $U$.
 (See the section {\rm 2.2} in {\rm \cite{sabbah}} for
 a smooth polarized pure twistor structure.
 It is equivalent to a polarized variation of pure twistor structures,
 given by Simpson.)
\item
 There exists a resolution $\varphi:\tilde{Z}\lrarr Z$
 satisfying the following:
 \begin{itemize}
  \item
  $\tilde{D}:=\varphi^{-1}(Z-U)$ is a normal crossing divisor of
  $\tilde{Z}$ 
  \item
  The corresponding harmonic bundle to
  $\varphi^{-1}\nbigt$ is tame
  with respect to the divisor $\tilde{D}$.
 \end{itemize}
\end{itemize}

If the corresponding harmonic bundle to $\varphi^{-1}\nbigt$
is tame and pure imaginary (see {\rm\cite{mochi3}}),
then 
$\nbigt$ is called pure imaginary.
 \hfill\qed
\end{df}

\begin{df}
Let $(U,\nbigt)$ and $(U',\nbigt')$  be generically defined
tame smooth pure twistor structures on $Z$.
They are called equivalent, if 
there exists a generically defined tame smooth pure twistor structures
$(U'',\nbigt'')$ on $Z$ satisfying 
$U''\subset U\cap U'$,
$\nbigt_{|U''}\simeq\nbigt''$ and $\nbigt'_{|U''}\simeq \nbigt''$.
\hfill\qed
\end{df}

The special case $w=0$ is as follows.
\begin{df}
A generically defined tame harmonic bundle on $Z$
is defined to be a data as follows:
\begin{itemize}
\item
 A Zariski open subset $U$ of $Z$.
 It is smooth.
\item 
 A harmonic bundle $(E,\delbar_E,\theta,h)$
 defined on $U$.
\item
 There exists a resolution $\varphi:\tilde{Z}\lrarr Z$
 satisfying the following:
 \begin{itemize}
  \item
  $\tilde{D}:=\varphi^{-1}(Z-U)$ is a normal crossing divisor of
  $\tilde{Z}$ 
  \item
  $\varphi^{-1}\harmonicbundle$ is a tame harmonic bundle
  on $\tilde{Z}-\tilde{D}$.
 \end{itemize}
\end{itemize}

If $\varphi^{-1}\harmonicbundle$ is tame and pure imaginary,
then $\harmonicbundle$ is called pure imaginary.
\hfill\qed
\end{df}

\begin{df}
 Let $\bigl(U,\harmonicbundle\bigr)$
 and $\bigl(U',(E',\delbar_{E'},h',\theta')\bigr)$ be
 generically defined tame harmonic bundles on $Z$.
 They are called equivalent if the following holds:
\begin{itemize}
\item
 There exists a generically defined tame harmonic bundle
 $\bigl(U'',(E'',\delbar_{E''},h'',\theta'')\bigr)$ on $Z$
 satisfying the following:
\[
  U''\subset U\cap U',
\quad
 \harmonicbundle_{|U''}=(E'',\delbar_{E''},h'',\theta''),
\quad
 (E',\delbar_{E'},h',\theta')_{|U''}=
 (E'',\delbar_{E''},h'',\theta'').
\]
\hfill\qed
\end{itemize}
\end{df}

\subsubsection{Statement}
\label{subsubsection;b12.5.21}

Let $\VPTgen(Z,w)$ denote the set of the equivalence classes
of generically defined tame variation of polarized pure twistor
structures of weight $w$ on $Z$.
Let $\MPT(Z,w)$ denote the set of the equivalence classes
of polarized regular pure twistor $D$-module of weight $w$
whose strict support is $Z$.

The following theorem is one of the main result
in this paper.
\begin{thm}\label{thm;b12.5.20}
We have the bijective correspondence
$\VPTgen(Z,w)\simeq \MPT(Z,w)$.
\end{thm}

Now it is a rather formal consequence of our study.
The map from $\MPT(Z,w)\lrarr \VPTgen(Z,w)$ is given
in the subsection \ref{subsection;04.1.27.35}.
The map from $\VPTgen(Z,w)\lrarr \MPT(Z,w)$ is given
in the subsection \ref{subsection;a11.23.51}
and the subsection \ref{subsection;a11.23.52}.
It is clear from our construction that the composite
$\VPTgen(Z,w)\lrarr \MPT(Z,w)\lrarr \VPTgen(Z,w)$ is identity.
It is shown in
the subsection \ref{subsection;a11.23.52}
that the composite $\MPT(Z,w)\lrarr \VPTgen(Z,w)\lrarr \MPT(Z,w)$
is identity.

%% file: d5.tex

\subsubsection{The statement}
\label{subsubsection;04.1.27.70}

Let $X$ be a complex manifold,
and $Z$ be a closed irreducible subvariety of $X$.
Let $\nbigt$ be a regular pure twistor $D$-module of weight $n$
whose strict support is $Z$.
Due to the result of Sabbah,
we have a smooth Zariski open subset $U$ of $Z$
and a harmonic bundle $\harmonicbundle$ defined on $U$
the restriction of $\nbigt$
to the open subset $X-(Z-U)$ is
push-forward of $\nbigm(E)\otimes\nbigo(n)$,
via the inclusion $U\lrarr X-(Z-U)$.

\begin{prop}\label{prop;04.1.27.52}
The harmonic bundle $\harmonicbundle$ is tame.
\end{prop}
The proof of Proposition \ref{prop;04.1.27.52} is given
in the subsubsections
\ref{subsubsection;04.1.27.50}--\ref{subsubsection;04.1.27.55}.

\subsubsection{One dimensional case}
\label{subsubsection;04.1.27.50}

Let us consider the case $\dim Z=1$.
Let $Q$ be a point of $Z-U$.
Since the property is local,
we may assume the following:
\begin{condition}\mbox{{}}\label{condition;04.2.4.50}
\begin{itemize}
\item
 $X=\Delta^n$ and $Z$ is irreducible.
\item
 Let $q_1$ denote the projection $\Delta^n\lrarr \Delta$
 onto the first component.
 We have $q_1(Q)=O$.
 The restriction of $q_1$ to $Z$ is proper and surjective.
 The restriction of $q_1^{-1}(\Delta^{\ast})\cap Z\lrarr \Delta^{\ast}$
 is a covering map.
\end{itemize}
\end{condition}
Since the restriction $q_{1\,|\,Z}$ is proper,
we have the push-forward $q_{1\,+}\bigl(\nbigt\bigr)$,
which is a regular pure twistor $D$-module.
Let $\nbigt_1=(\nbigm_1,\nbigm_1,C)$ denote the component of
$q_{1,\,+}\bigl(\nbigt\bigr)$ whose strict support is
$\Delta$.
Then there exists the harmonic bundle
$(E_1,\delbar_{E_1},\theta_1,h_1)$ on $\Delta^{\ast}$
such that
the restriction
$\nbigt_{1\,|\,\Delta^{\ast}}\simeq
 \nbigm(E_1)\otimes\nbigo(n)$.
Let us take a $V$-filtration of $\nbigm_1$ along $z_1=0$
at $\lambda=0$ on $\Delta\times\Delta(\epsilon)$ for some positive
number $\epsilon>0$.
Due to the regularity,
$V_{a}$ is a coherent $\nbigo_{\Delta\times\Delta(\epsilon)}$-module
for $a<0$.
We have
$\deldel_{z_1}\cdot V_{a}\subset V_{a+1}$.
By considering the specialization at $\lambda=0$,
we obtain that $(E_1,\delbar_{E_1},\theta_1,h_1)$ is tame.

Let us consider the normalization $\pi:\tilde{Z}\lrarr Z$.
Then we obtain the morphism
$F:=q_1\circ\pi:\tilde{Z}\lrarr\Delta$.
Then $F^{-1}\bigl(E_1,\delbar_{E_1},h_1,\theta_1\bigr)$ is tame.
Since $\pi^{-1}\bigl(E,\delbar_E,h,\theta\bigr)$ is 
a direct summand of $F^{-1}\bigl(E_1,\delbar_{E_1},h_1,\theta_1\bigr)$,
we can conclude that 
$\pi^{-1}\bigl(E,\delbar_E,\theta,h\bigr)$ is tame on $\tilde{Z}$,
namely,
$\harmonicbundle$ is tame.

\subsubsection{The general case}
\label{subsubsection;04.1.27.55}

Let $\pi:\tilde{Z}\lrarr Z$ be a birational morphism
such that $\tilde{Z}$ is smooth
and that $\tilde{D}:=\pi^{-1}(Z-Z')$
is the normal crossing divisor of $\tilde{Z}$.
Let $C$ be a smooth curve in $\tilde{Z}$,
which meets the smooth part of $\tilde{D}$ transversally.
Then we obtain the curve $\pi(C)$ in $Z$.
Let us specialize $\nbigt$ to $\pi(C)$,
which gives the pure twistor $D$-module 
$\nbigt'$ whose strict support is $\pi(C)$.
The restriction of $\nbigt'$ to $\pi(C)\cap Z'$
corresponds to the restriction of the harmonic bundle
$\harmonicbundle_{|\pi(C)\cap Z'}$.
The pull back of $\harmonicbundle_{|\pi(C)\cap Z'}$
is isomorphic to the restriction 
$\pi^{-1}\harmonicbundle_{C\setminus\tilde{D}}$.
Then we can obtain the tameness of
$\pi^{-1}\harmonicbundle_{C\setminus\tilde{D}}$,
due to the result in the one dimensional case
(the subsubsection \ref{subsubsection;04.1.27.50}).
Then we obtain the tameness 
of $\pi^{-1}\harmonicbundle$,
due to Corollary \ref{cor;04.1.27.51}.
Thus the proof of Proposition \ref{prop;04.1.27.52} is accomplished.
\hfill\qed

%% file: a35.1.tex

\subsubsection{A prolongment of a tame variation
of polarized pure twistor structure
with normal crossing divisor}

Let $X$ be a complex manifold
and $D$ be a normal crossing divisor.
Let $\harmonicbundle$ be a tame harmonic bundle over $X-D$.
Then we have the $\nbigr$-triple
$(\nbige,\nbige,C_0)$ on $X-D$,
as in the section 2.2 of \cite{sabbah}.
Then we obtain the $\nbigr$-triple
$(\gbige,\gbige,\gbigc)$ on $X$.
In the case $X=\Delta^n$, $D_i=\{z_i=0\}$ and $D=\bigcup_{i=1}^l D_i$,
the construction of $(\gbige,\gbige,\gbigc)$ is given
in the subsection \ref{subsection;a11.23.40}
and the subsubsection \ref{subsubsection;9.25.2}.
The local construction can be globalized.

Let $\nbigt$ be a tame variation of polarized pure twistor structures
of weight $w$ on $X-D$.
Let $\harmonicbundle$ denote the corresponding tame harmonic bundle
such that
$\nbigt\simeq(\nbige,\nbige,C_0)\otimes\nbigo_{\proj^1}(w)$.
Then we have the natural prolongment
$\gbigt:=(\gbige,\gbige,\gbigc)\otimes\nbigo_{\proj^1}(w)$.

\begin{thm}\label{thm;11.5.2}
The $\nbigr$-triple
$\gbigt=(\gbige,\gbige,\gbigc)\otimes\nbigo_{\proj^1}(w)$
is a regular pure twistor $D$-module of weight $w$.
\end{thm}
\pf
We have only to consider the case $w=0$.
We use the right $\nbigr$-structures.
To distinguish it, we use the notation
$\gbige\otimes\omega_X$
instead of $\gbige$,
where we put $\omega_X:=\wedge ^n\Omega^{1,0}_X$.
So we show the following:
Let $\harmonicbundle$ be a tame harmonic bundle
over $X-D$.
Then the $\nbigr$-triple
$\nbigm(E):=
 \bigl(\gbige\otimes\omega_X,\gbige\otimes\omega_X,\gbigc\bigr)$
is a regular pure twistor $D$-module.

We use an induction on the dimension of $X$.
We assume that the claim has already been shown in the case $\dim (X)<n$,
and we will show the claim in the case $\dim(X)=n$.

Since the property is local,
we may assume $X=\Delta^n$ and $D=\bigcup_{i=1}^nD_i$,
where $D_i=\{z_i=0\}$.
Let $g$ be a holomorphic function on $X$.
Then we have only to show that
$\Gr^{W(N)}_h\tildepsi_{g,u}\nbigm(E)$
satisfies the conditions
Definition \ref{df;11.5.10} and Definition \ref{df;11.5.11},
which are Definition 4.1.2 and Definition 4.2.2
in \cite{sabbah}.

Let $\pi:\widetilde{X}\lrarr X$ be a birational morphism
satisfying the following:
\begin{itemize}
\item
 We put  $\tilde{g}:=\pi\circ g$,
and $\tilde{D}:=\pi^{-1}(D)\cup\tilde{g}^{-1}(0)$.
Then $\tilde{D}$ is normal crossing.
\item
 The restriction $\pi_{|\tilde{X}-\tilde{D}}$
 gives isomorphism $\tilde{X}-\tilde{D}\lrarr X-\bigl(D\cup g^{-1}(0)\bigr)$.
\end{itemize}

Let $(\tilde{E},\delbar_{\tilde{E}},\tilde{h},\tilde{\theta})$
denote the pull back of $\harmonicbundle_{|X-\big(D\cup g^{-1}(0)\bigr)}$.
Then we obtain the $\nbigr$-triple
$\nbigm(\tilde{E})
=(\tildegbige\otimes\omega_{\tilde{X}},
 \tildegbige\otimes\omega_{\tilde{X}},\tilde{\gbigc})$
on $\tilde{X}$.
Since $P_h\Gr^{W(N)}_h\tildepsi_{\tilde{g},u}\nbigm(\tilde{E})$
are isomorphic to the tensor product
of the Tate objects and the $\nbigr$-triples induced by 
a tame harmonic bundles
(Proposition \ref{prop;b12.4.1}),
they are regular pure twistor $D$-modules of weight $h$
due to the hypothesis of the induction.

We obtain the $\nbigr$-triple
$\nbigh^i\pi_{\ast}\nbigm(\tilde{E})=(\nbigm_1,\nbigm_1,C)$.
Due to the argument of Sabbah-Saito
we know the following:
\begin{lem}\mbox{{}}
\begin{enumerate}
\item
 The $\nbigr$-modules
 $\nbigm_i$ are strictly $S$-decomposable along $g=0$,
 $g\cdot \prod_{i=1}^n z_i=0$,
 and $\prod_{i=1}^n z_i=0$.
\item
 The $\nbigr$-triples
 $P_h\Gr^{W(N)}_h\tildepsi_{g,u}\nbigh^i\pi_{\ast}\nbigm(\tilde{E})$
 are regular pure twistor $D$-modules with the polarization.
\end{enumerate}
\end{lem}
\pf
See the section 6.2.6 in \cite{sabbah}.
\hfill\qed

\vspace{.1in}

We take the component $\nbigt=(\nbigm_2,\nbigm_2,C)$
of $\nbigh^0\pi_{\ast}\nbigm(\tilde{E})$
such that $\nbigm_2$ does not have non-trivial
submodules whose supports are contained in
$g^{-1}(0)$.
Note that the regularity of $\nbigm_2$
along $g=0$ follows from the regularity of
$\widetilde{\gbige}$ along $\tilde{g}=0$.

\begin{lem} \label{lem;a11.23.7}
To check $P_h\Gr^{W(N)}_h\tildepsi_{g,u}\nbigm(E)$
satisfies the conditions in the definitions {\rm\ref{df;11.5.10}}
and {\rm\ref{df;11.5.11}},
we have only to show $\nbigt\simeq \nbigm(E)$.
\end{lem}
\pf
It follows from the fact
that $P_h\Gr^{W(N)}_h\tildepsi_{g,u}\nbigt$
is a direct summand of
$P_h\Gr^W{(N)}_h\tildepsi_{t,u}\nbigh^0\pi_{\ast}\nbigm(\tilde{E})$.
\hfill\qed

\vspace{.1in}

Let us show that $\nbigt$ is isomorphic to $\nbigm(E)$.
We put $X_1:=X-\bigl(D\cup g^{-1}(0)\bigr)$.
Let $\iota$ denote the open immersion $X_1\lrarr X$.
Then $\gbige\otimes\omega_X$ and
$\nbigm_2$ 
are contained in
$\iota_{\ast}\nbige\otimes\omega_X$.

We use the following remark.
We take the metric $g_1$ of $\Omega^{1,0}_X\bigl(\log D\bigr)$
such that $z_i^{-1}\cdot dz_i$ $(i\in\nbar)$ are
the orthonormal frame.
Let $Q$ be any point of $\tilde{X}$.
We take a neighbourhood of $Q$
and a coordinate $(\zeta_1,\ldots,\zeta_n)$ on $U$
such that $U\cap \tilde{D}=\bigcup_{i=1}^l\bigl\{\zeta_i=0\bigr\}$.
We also take a metric $g_2$ of
$\Omega^{1,0}_{\widetilde{X}}\bigl(\log\tilde{D}\bigr)$ on $U$,
such that $\zeta_i^{-1}\cdot d\zeta_i$ are the orthonormal frame.
Then the naturally defined morphism
$\pi^{\ast}\Omega^{1,0}_X(\log D)\lrarr 
\Omega^{1,0}_{\widetilde{X}}(\log \tilde{D})$
is bounded on $U$.
The metric $h$ on $E$ and the metrics $g_a$ $(a=1,2)$
induce the metrics
$h_1$ on $E\otimes\omega_X$
and $h_2$ on $\pi^{\ast}E\otimes\omega_{\widetilde{X}}$ restricted to $U$.
The naturally induced morphism
$\pi^{\ast}(E\otimes\omega_X)\lrarr
\pi^{\ast}E\otimes\omega_{\widetilde{X}}$
is bounded on $U$.

Let $f\cdot dz_1\wedge\cdots\wedge dz_n$
be a section of $\gbige\otimes\omega_X$.
The fact $f\cdot dz_1\wedge\cdots\wedge dz_n\in \Vzero_{<0}\otimes\omega_X$
is equivalent to the following estimate 
for some positive numbers $\epsilon$ and $C$:
\[
 \bigl|f\cdot dz_1\wedge\cdots\wedge dz_n\bigr|_{h_1}
\leq C\cdot \prod_{i=1}^n|z_i|^{\epsilon}.
\]
Then we obtain the following estimate on $U$
for some positive numbers $C'$ and $\epsilon'$:
\[
 \bigl|\pi^{\ast}(f)\cdot\pi^{\ast}(dz_1\wedge\cdots \wedge dz_n)\bigr|_{h_2}
\leq C'\cdot \prod_{i=1}^n|\zeta_i|^{\epsilon'}.
\]
It means that
the restriction
$\pi^{\ast}(f)\cdot \pi^{\ast}(dz_1\wedge\cdots\wedge dz_n)$ to $U$
naturally induces the section of
$(\lefttop{\nbar}\Vzero_{<0}\tilde{\gbige}\otimes\omega_{\tilde{X}})_{|U}$.
Hence $f\cdot dz_1\wedge\cdots\wedge dz_n$
naturally induces the section of
$\nbigm_2$.
Since $\gbige\otimes\omega_X$
is generated by $\lefttop{\nbar}\Vzero_{<0}\gbige\otimes\omega_X$,
we obtain the inclusion
$\psi:\gbige\otimes\omega_X\subset\nbigm_2$.

\begin{lem} \label{lem;a11.23.5}
The restriction of the inclusion to $X-D$ is isomorphic,
namely $(\gbige\otimes\omega_X)_{|X-D}\simeq \nbigm_{2\,|\,X-D}$.
\end{lem}
\pf
By our construction,
we have
$(\gbige\otimes\omega_X)_{|X_1}\simeq
 \nbigm_{2\,|\,X_1}$.
Hence the support of the cokernel of the morphism
$(\gbige\otimes\omega_X)_{|X-D}\lrarr \nbigm_{2\,|\,X-D}$
is contained in $g^{-1}(0)\cap(X-D)$.
Since both of $(\gbige\otimes\omega_X)_{|X-D}$ and $\nbigm_{2\,|\,X-D}$
are strictly $S$-decomposable along $g^{-1}(0)\cap(X-D)$,
we obtain $(\gbige\otimes\omega_X)_{|X-D}\simeq \nbigm_{2\,|\,X-D}$
due to Corollary \ref{cor;a11.23.3}.
\hfill\qed

\begin{lem}
The inclusion $\gbige\otimes\omega_X\lrarr \nbigm_2$ is isomorphic.
\end{lem}
\pf
Both of $\gbige\otimes\omega_X[\deldel_t]$ and $\nbigm_2[\deldel_t]$ are 
strictly $S$-decomposable along $\prod_{i=1}^n z_i$.
Then we obtain the result
due to Lemma \ref{lem;a11.23.5} and Corollary \ref{cor;a11.23.3}.
\hfill\qed

\vspace{.1in}

The coincidence of the sesqui-linear pairings
are obtained by the uniqueness
(Proposition 3.6.6. in \cite{sabbah} and 
 Corollary \ref{cor;a11.23.6} in this paper).
Hence $\nbigt$ and $\nbigm(E)$ are isomorphic.
Due to Lemma \ref{lem;a11.23.7},
the proof of Theorem \ref{thm;11.5.2} is accomplished.
\hfill\qed

%% file: a77.1.tex

\subsubsection{The existence of prolongment}
\label{subsubsection;04.2.4.25}

Let $X$ and $Z$ be as in the subsection \ref{subsection;a11.23.50}.
Let $(U, \nbigt)$ be a generically defined
variation of polarized pure twistor structures of weight $w$.
We take a resolution
$\varphi:\tilde{U}\lrarr U$ as in Definition \ref{df;a11.23.60}.
Due to Theorem \ref{thm;11.5.2},
we obtain the polarized regular pure twistor $D$-module $\gbigt'$,
which is a prolongment of $\varphi^{-1}\nbigt$.
Let us take the component $\gbigt$ of
$\nbigh^0(\varphi_{\ast}\gbigt')$
whose strict support is $Z$,
which is a polarized regular pure twistor $D$-module.
By our construction,
we have the naturally defined isomorphism
$\gbigt_{|U}\simeq \nbigt$,
namely $\gbigt$ gives a prolongment of $\nbigt$.

%% file: a35.2.tex

\subsubsection{The statement of the uniqueness}

\label{subsubsection;9.30.11}

Let $X$ be a complex manifold and $Z$ be an irreducible
closed subset of $X$.
Let $\nbigt=(\nbigm',\nbigm'',C)$ be a regular pure twistor $D$-module
whose strict support is $Z$.
Due to Proposition \ref{prop;04.1.27.52},
we have a smooth Zariski open subset $U$ of $Z$
and a tame harmonic bundle $\harmonicbundle$
defined over $U$,
such that
the restriction of $\nbigt$
to the open subset $X-(Z-U)$ is
push-forward of $\nbigm(E)$ via the inclusion $U\lrarr X-(Z-U)$.
We put $Z-U=Y$.

Let take a birational morphism $\varphi:\tilde{Z}\lrarr Z$
such that $\tilde{Z}$ is smooth
and $\tilde{D}:=\varphi^{-1}(Y)$ is a normal crossing divisor.
We have the harmonic bundle
$\varphi^{-1}(E,\delbar_E,h,\theta)
=(\tilde{E},\delbar_{\tilde{E}},\tilde{h},\tilde{\theta})$.
Then we obtain the regular pure twistor $D$-module
$\nbigm(\tilde{E})=(\tilde{\gbige},\tilde{\gbige},\tilde{\gbigc})$.
We have the push-forward
$\nbigh^0\bigl(
 \varphi_{\ast}\nbigm(\tilde{E})\bigr)$.
Let $\tildenbigt=(\tildenbigm',\tildenbigm'',\tilde{C})$
denote the component of $\nbigh^0\bigl(
 \varphi_{\ast}\nbigm(\tilde{E})\bigr)$
whose strict support is $Z$.
We have the natural isomorphism
$\nbigt_{X-Y}\simeq \tildenbigt_{X-Y}$,
because both of them are push-forward
of $\nbigm(E)$.
We would like to show that the isomorphism
can be prolonged to the isomorphism
of $\nbigt$ and $\tildenbigt$.

\begin{thm}\label{thm;11.5.6}
We have the natural isomorphism
$\nbigt\simeq\tildenbigt$.
\end{thm}

We will show the theorem in the following subsubsections.
Before entering the proof,
we give some remarks.

Since the claim is local on $X$,
we may assume that $X=\Delta^n$.
We may also assume that we have a holomorphic function
$g$ such that $Y\subset g^{-1}(0)$.

We put $\omega_U:=\bigwedge^{\dim U}\Omega_U^{1,0}$.
We denote the immersion $U\lrarr X$ by $\iota_U$.
Then $\nbigm'$, $\nbigm''$, $\tildenbigm'$ and $\tildenbigm''$
are contained in $\iota_{U\,\ast}\nbige\otimes\omega_U$,
where $\nbige\otimes\omega_U$ is regarded as the right $\nbigr$-module.

\begin{lem}
We have only to show that
$\nbigm'$ (resp. $\nbigm''$)
and $\tildenbigm'$ (resp. $\tildenbigm''$) 
are same as the $\nbigr$-submodules of
$\iota_{U\,\ast}\nbige\otimes\omega_U$.
\end{lem}
\pf
It follows from the uniqueness of the pairing
(Proposition \ref{prop;a11.23.10})
\hfill\qed

\begin{lem}\label{lem;04.2.16.100}
Let $\lambda\in\cnum_{\lambda}$ be generic.
We have $\nbigm'_{|\nbigx^{\lambda}}=\widetilde{\nbigm}'_{|\nbigx^{\lambda}}$
in $\bigl(\iota_{U\,\ast}\nbige\otimes\omega_U\bigr)_{|\nbigx^{\lambda}}$.
\end{lem}
\pf
Both of 
$\nbigm'_{|\nbigx^{\lambda}}$ and $\widetilde{\nbigm}'_{|\nbigx^{\lambda}}$
can be regarded as the regular holonomic $D$-modules.
The restrictions to $X-Y$ is same
as the push-forward of
$\bigl(\nbigelambda,\DD^{\lambda,f}\bigr)$.
By the Riemann-Hilbert correspondence,
the regular holonomic $D$-modules
$\nbigm'_{|\nbigx^{\lambda}}$ and $\widetilde{\nbigm}'_{|\nbigx^{\lambda}}$
correspond to the perverse sheaves $\nbigf_1$ and $\nbigf_2$,
and the flat bundle $(\nbigelambda,\DD^{\lambda,f})$
corresponds to the local system $L$ on $U$.
Then both of $\nbigf_i$ $(i=1,2)$ are the intermediate extensions
of $L$,
due to the strictly $S$-decomposability of 
$\nbigm'_{|\nbigx^{\lambda}}$ and $\widetilde{\nbigm}'_{|\nbigx^{\lambda}}$.
Thus $\nbigf_1$ and $\nbigf_2$ are isomorphic.
Therefore we have the isomorphism
$\nbigm'_{|\nbigxlambda}\simeq\widetilde{\nbigm}'_{|\nbigxlambda}$,
whose restriction to to $X-Y$ can be regarded
as the identity of $\nbige^{\lambda}\otimes\omega_U$.
Hence we obtain the result.
\hfill\qed

\vspace{.1in}

We will use the following elementary lemma.
\begin{lem} \label{lem;9.30.3}
Let $B$ be a connected complex manifold.
Let $f$ be a holomorphic function
on $B\times\Delta^{\ast}$.
Assume the following:
\begin{itemize}
\item
For any point $b\in B$,
the restriction $f_{|\{b\}\times \Delta^{\ast}}$
is meromorphic at the origin $\{O\}$.
In particular, we obtain the order of the pole $p(b)$.
\end{itemize}
Then the numbers $p(b)$ are bounded of $b\in B$.
\end{lem}
\pf
We have the development
\[
 f_{|\{b\}\times\Delta^{\ast}}
=\sum f_n(b)\cdot z^{n}.
\]
The coefficients $f_n(b)$ are calculated as follows:
\[
 f_n(b)=
 \int_{|z|=1/2} z^{-n-1}\cdot f_{|\{b\}\times\{|z|=1/2\}} dz.
\]
Hence the dependence of $f_n(b)$ $(n\in\seisuu)$ on $b\in B$
are holomorphic.
Assume that the set $S:=\{n<0\,|\,f_n\not\equiv 0\}$ is infinite.
We have $B\neq \bigcup_{n\in S}f_n^{-1}(0)$.
For any point $b\in B-\bigcup_{n\in S}f_n^{-1}(0)$,
the restriction of $f$ to $\{b\}\times\Delta^{\ast}$
is transcendental,
which contradicts to the assumption of the lemma.
Thus the set $S$ is finite.
\hfill\qed

\subsubsection{Step 1. In the case $\dim Z=1$}

We may assume that $\tilde{Z}=\Delta=\bigl\{w\,\big|\,|w|<1\bigr\}$,
and that $\varphi=(\varphi_1,\ldots,\varphi_n):\Delta\lrarr X=\Delta^n$
gives the homeomorphism
$\tilde{Z}\lrarr Z$.
We may also assume that $\varphi_n(w)=w^l$ for some integer $l$
and that $U\simeq\tilde{Z}-\{0\}=\Delta^{\ast}$.

Let $\lambda_0$ be a point of $\cnum_{\lambda}$.
We pick a frame $\vecv=(v_i)$ of $V_{<0}\tildegbige$
on $\Delta(\lambda_0,\epsilon_0)\times\Delta$.
Each $v_i$ naturally gives the section of
$\tildenbigm'$.
We denote it by $\varphi_{\ast}(v_i)$.
Then we obtain the submodule of $\tildenbigm'$
generated by $\varphi_{\ast}\vecv$.
We denote the submodule by $\varphi_{\ast}\vecv\cdot\nbigr$.

\begin{lem}
We have $\varphi_{\ast}\vecv\cdot\nbigr=\tildenbigm'$.
\end{lem}
\pf
By definition, we have
$\tildenbigm'=
R\varphi_{\ast}(
 \tildegbige\otimes^L_{\nbigr_{\Delta}}\varphi^{-1}\nbigr_{\nbigx}
 )$.
Since $\varphi$ is homeomorphic, we have
$R\varphi_{\ast}=\varphi_{\ast}$.
Thus $\tildenbigm'$ is same as the component
of
$\varphi_{\ast}(\tildegbige\otimes_{\nbigr_{\Delta}}\varphi^{-1}\nbigr_X)$.
Thus it is generated by $\varphi_{\ast}(\vecv)$.
\hfill\qed

\begin{lem}
Let $f$ be a section of $\iota_{U\,\ast}\nbige\otimes\omega_U$.
Then we have the unique description as follows:
\begin{equation} \label{eq;9.30.1}
 f=\sum_{\vecn\in\seisuu_{\geq\,0}^{n-1}}
 v_i\cdot f_{\vecn,i}\cdot \prod_{i=1}^{n-1}\deldel_i^{n_i}.
\end{equation}
Here $f_{\vecn,i}$ denote holomorphic functions on
$U\simeq\Delta^{\ast}$.
\end{lem}
\pf
Note that $\del/\del z_1,\ldots,\del/\del z_{n-1}$ gives the frame of
the normal bundle of $U$ in $X-Y$.
Then the claim is clear.
\hfill\qed

\begin{cor} \label{cor;9.30.6}
The morphism
$z_n^N:\iota_{U\,\ast}\nbige\otimes\omega_U
 \lrarr\iota_{U\,\ast}\nbige\otimes\omega_U$
is injective for any positive integer $N$.
\end{cor}
\pf
It is clear from the description (\ref{eq;9.30.1}).
\hfill\qed

\vspace{.1in}
For a section $f$ of $\iota_{U\,\ast}\nbige\otimes\omega_U$
or $\iota_{U\,\ast}(\nbige\otimes\omega_U)_{|\nbigx^{\lambda}}$,
we often use the notation  $u_{\vecn}(f)$
to denote $\sum_{i} v_i\cdot f_{\vecn,i}$.
Here $f_{\vecn,i}$ are given in the development (\ref{eq;9.30.1}).
Then $u_{\vecn}(f)$ can be regarded as the holomorphic sections of
$\tildenbige\otimes\Omega_{\Delta^{\ast}}$ 
over $\Delta(\lambda,\epsilon_0)\times\Delta^{\ast}$
or $\Delta^{\ast}$.
We also note that 
the set $\{\vecn\,|\,u_{\vecn}(f)\neq 0\}$ is finite.

\begin{lem} \label{lem;9.30.2}
Let $f$ be a section of $\tildenbigm'$, or 
$\tildenbigm'_{|\nbigx^{\lambda}}$.
The holomorphic functions $f_{\vecn,i}$ on $\Delta^{\ast}$ in 
the development {\rm (\ref{eq;9.30.1})} are meromorphic
at $\Delta(\lambda_0,\epsilon)\times\{O\}$ or $\{O\}$,
i.e.,
they are not transcendental.

In particular, $u_{\vecn}(f)$ are sections
of $\naiveprolong{\tildenbige}$ on
$\Delta(\lambda_0,\epsilon_0)\times\Delta$ or $\Delta$.
\end{lem}
\pf
We have
$\varphi_{\ast}(\del/\del w)=\sum a_j\cdot\varphi^{-1}(\del/\del z_j)$
for holomorphic functions $a_j$ on $\Delta$.
We have $a_n=l\cdot w^{l-1}$.
Then $\varphi^{-1}(\del/\del z_n)$ can be described
as the linear combination of
$\del/\del w$ and $\varphi^{-1}(\del/\del z_j)$ $(j=1,\ldots,n-1)$
with the meromorphic coefficients.
Since we have $\tildenbigm'=\varphi_{\ast}\vecv\cdot\nbigr$,
we obtain the result.
\hfill\qed

\vspace{.1in}

Let us consider the $V$-filtrations at $\lambda_0$
of $\nbigm'$ and $\tildenbigm'$
along $z_n=0$.
Let us pick a section of $f$ of $\Vzero_{<0}\nbigm'$.

\begin{lem}
For any sufficiently large integer $N$,
the section $f\cdot z_n^N$ is contained in $\Vzero_{<0}\tildenbigm'$.
\end{lem}
\pf
Since $f$ can be regarded as the section of
$\iota_{U\,\ast}\nbige\otimes\omega_U$,
we have the development
$f=\sum v_i\cdot f_{i\,\vecn}\cdot\prod_{i=1}^{n-1}\deldel_i^{n_i}$.
Let $\lambda\in\Delta(\lambda_0,\epsilon_0)$ be generic.
Then the restriction
$f_{|\nbigx^{\lambda}}$
is contained in
$\Vzero_{<0}\nbigm'_{|\nbigx^{\lambda}}=
 \Vzero_{<0}\tildenbigm'_{|\nbigx^{\lambda}}$
(Lemma \ref{lem;04.2.16.100}).
Due to Lemma \ref{lem;9.30.2},
the restrictions
$f_{i\,\vecn\,|\,\{\lambda\}\times\Delta^{\ast}}$
is meromorphic.
Due to Lemma \ref{lem;9.30.3},
the order of the poles are bounded independently of a choice of
generic $\lambda\in\Delta(\lambda_0,\epsilon_0)$.
Hence we obtain some sufficiently large integer $N$
such that $f_{i\,\vecn}\cdot \varphi^{\ast}z_n^N$ 
are holomorphic on $\Delta(\lambda_0,\epsilon_0)\times \Delta$
for any $i$ and $\vecn$,
due to Hartogs' Theorem.
Then we have $\sum u_{\vecn}(f)\cdot \varphi^{\ast}z_n^N$
are sections of $\Vzero_{<0}\tildegbige$.
In particular, we have
$f\cdot \varphi(z_n^N)\in\tildegbige\otimes\varphi^{-1}\nbigr_{\nbigx}$.
Hence
we obtain $f\cdot z_n^N\in\tildenbigm'$.
If we take larger $N$,
then we may assume that $f\cdot z_n^N\in \Vzero_{<0}\tildenbigm'$.
\hfill\qed

\vspace{.1in}

Since $f\cdot z_n^N$ is contained in
$\Vzero_{<0}\nbigm'\cap \Vzero_{<0}\tildenbigm'$,
we have
$(f\cdot z_n^N)\cdot\nbigr\subset\nbigm'\cap\tildenbigm'$.
The $V$-filtrations
$\Vzero(\nbigm')$ and $\Vzero(\tildenbigm')$
induce the filtrations $U^{(1)}$ and $U^{(2)}$
on $(f\cdot z_n^N)\cdot\nbigr$.

\begin{lem} \label{lem;9.30.5}
We have $U^{(1)}=U^{(2)}$
\end{lem}
\pf
Both of $U^{(1)}$ and $U^{(2)}$
are good and monodromic.
Moreover $\Gr^{U^{(a)}}$ $(a=1,2)$ are strict.
Thus we have $U^{(1)}=U^{(2)}$,
due to Lemma \ref{lem;9.17.19}.
\hfill\qed

\begin{lem} \label{lem;9.30.7}
We have $\nbigm'\subset\tildenbigm'$.
\end{lem}
\pf
If $f$ is contained in $\Vzero_{b}\nbigm'$ for a negative number $b<0$,
we have $f\cdot z_n^N\in \Vzero_{b-N}(\nbigm')$.
Due to Lemma \ref{lem;9.30.5},
we have $f\cdot z_n^N\in\Vzero_{b-N}(\tildenbigm')$.
Since $z_n^N:\Vzero_b(\tildenbigm')\lrarr\Vzero_{b-N}(\tildenbigm')$
are surjective for $b<0$,
we have a section $\tilde{f}$
of $\Vzero_{b}(\tildenbigm')$ such that
$\tilde{f}\cdot z_n^N=f\cdot z_n^N$.
Due to the injectivity of the morphism $z_n^N$
(Corollary \ref{cor;9.30.6}),
we have $\tilde{f}=f$.
It means that $f$ is contained in $\Vzero_b(\tildenbigm')$.
Thus we obtain
$\Vzero_{<0}(\nbigm')\subset\Vzero_{<0}(\tildenbigm')$.
Since $\nbigm'$ is generated by $\Vzero_{<0}(\nbigm')$,
we obtain Lemma \ref{lem;9.30.7}.
\hfill\qed

\vspace{.1in}

\begin{lem}\label{lem;a11.23.20}
We have $\nbigm'=\tildenbigm'$.
\end{lem}
\pf
Let $i$ denote the inclusion $\nbigm'\subset\tildenbigm$.
Since $i$ is isomorphic on $X-Y$,
the support of the cokernel of $i$ is contained in $Y=Z-U$.
Thus we have $\Vzero_{b}(\nbigm')=\Vzero_{b}(\tildenbigm')$
for any sufficiently negative number $b$.
Since both of the filtrations $\Vzero(\nbigm')$
and $\Vzero(\tildenbigm')$ are strictly specializable,
it is easy to derive that $\nbigm'=\tildenbigm'$
(Corollary \ref{cor;a11.23.3}).
\hfill\qed

\subsubsection{Step 2. In the case we have a good normal frame}
\label{subsubsection;a11.23.35}

We put $r:=\dim X-\dim Z$.
Then we can pick holomorphic functions $f_1,\ldots,f_r$ 
such that $Z$ is one of the irreducible components
of $\bigcap_{i=1}^r f_i^{-1}(0)$ with the reduced structure.
We put as follows:
\[
 w_j:=
\sum_{i=1}^n
 \frac{\del f_j}{\del z_i}
\frac{\del}{\del z_i}
\]

The following lemma is clear.
\begin{lem}\label{lem;a11.23.30}
There exists a Zariski open subset $U'$ of $U$
such that the restrictions of $w_1,\ldots,w_r$ to $U'$
give the frame of $N_{U'} X$.
\hfill\qed
\end{lem}

In this subsubsection, we impose the following assumption:
\begin{condition}\mbox{{}}\label{condition;b12.4.20}
$U'=U$.
\end{condition}

Under the assumption, we have the following lemma.
\begin{lem}
Let $f$ be a section of $\iota_{U\,\ast}(\nbige\otimes\omega_U)$.
We have the following unique description:
\begin{equation}\label{eq;a11.23.15}
 f=\sum_{\vecn\in\seisuu_{\geq\,0}^{r}}
 u_{\vecn}\cdot \prod_{j=1}^r w_j^{n_j}.
\end{equation}
Here $u_{\vecn}$ $(\vecn\in\seisuu_{\geq\,0}^r)$
denote holomorphic sections of $\nbige\otimes\omega_U$
on $U\simeq \tilde{Z}-\tilde{D}$.
The set $\{\vecn\,|\,u_{\vecn}\neq 0\}$ is finite.
\hfill\qed
\end{lem}

Let $f$ be a section of
$\nbigm'\subset\iota_{U\,\ast}(\nbige\otimes\omega_U)$.
We have the development of $f$
as in (\ref{eq;a11.23.15}).
We may assume that there is a holomorphic function $g$ on $X$
such that $Z-U\subset g^{-1}(0)$,
by shrinking $X$.

\begin{lem}
\label{lem;9.30.10}
Let $\tilde{C}\subset\tilde{Z}$ be a curve
transversal with the smooth part of $\tilde{D}$.
Then the restrictions $u_{\vecn\,|\,\tilde{C}}$
give meromorphic sections of
$\bigl(
 \naiveprolong{\tilde{\nbige}}\otimes\omega_{\widetilde{Z}}
\bigr)_{|\tilde{C}}$.
\end{lem}
\pf
We use an induction on the dimension of $Z$.
The case $\dim\tilde{Z}=1$ has been done in the Step 1
(the previous subsubsection).
We assume that the claim in the case $\dim Z\leq n-1$ holds,
and we will show the claim in the case $\dim Z=n$.

Let us consider the image $C=\varphi(\tilde{C})\subset Z$.
We may assume that
$\tilde{C}\cap \tilde{D}$ consists of  one point $\tilde{P}$.
We put $P=\varphi(\tilde{P})$.
In the following argument, we can shrink $X$
if we need to do so.

We take a holomorphic function $F$ on $X$ 
such that $C\cap U$ is contained in the smooth part of $F^{-1}(0)$
and $Z\not\subset F^{-1}(0)$.
Let us consider the inclusion $\iota_F:X\lrarr X\times\cnum$
of the graph,
and the $V$-filtration $\Vzero$ of $\iota_F\nbigm'$
along $t=0$, where $t$ denotes the coordinate of $\cnum$.

\begin{lem}
For any sufficiently large natural number $N$,
the section $f\cdot g^N\otimes 1$ is contained in
$\Vzero_{-1}(\iota_{F}\nbigm')$.
\end{lem}
\pf
It is easy to see that the restriction $f_{U}\otimes 1$ is contained
in $\Vzero_{-1}$.
Assume that $f\otimes 1$ is contained
in $\Vzero_{a}(\iota_F\nbigm')$ $(a>-1)$.
Then it induces the sections of
$\Gr^W\Gr^{\Vzero}_a(\iota_F\nbigm')$,
whose support is contained in $U-Z\subset g^{-1}(0)$.
Hence we can kill it by multiplying the power $g$.

By using an inductive argument,
we can take a sufficiently large number $N$
such that $f\cdot g^N\otimes 1$ is contained in $\Vzero_{-1}$.
\hfill\qed

\vspace{.1in}

Let $Z_1$ denote the irreducible component of
$F^{-1}(0)$ containing $C$.
Let $\nbigt_1=(\nbigm'_1,\nbigm''_1,C_1)$
denote the component of $\psi_{F,-1}\nbigt$
whose strict support is $Z_1$.

Then $f\cdot g^N$ induces the section $f_1$ of $\nbigm_1'$,
and we have
$f_{1\,|\,U\cap Z_1}=(f\cdot g^N)_{|U\cap Z_1}/dF$.
We have a development $f_1=\sum u_{\vecn}'\cdot \prod w_j^{n_j}$,
Then it is easy to see that
$u_{\vecn\,|\,C}$ can be holomorphically reconstructed 
from $u_{\vecn}'$.

Let $\varphi_1$ denote a resolution
$\widetilde{Z}_1\lrarr Z_1$
such that $\widetilde{D}_1=\widetilde{Z}_1-\varphi_1^{-1}(Z_1\cap U)$
is normal crossing.
We have the curve $\tilde{C}_1\subset \widetilde{Z}_1$
corresponding to $C$.
We may assume that $\tilde{C}_1$ transversal with $\tilde{D}_1$.
Due to the hypothesis of the induction,
the sections
$(u'_{\vecn\,|\,U\cap Z_1}/dF)_{|C}$
naturally give the meromorphic sections of
$\bigl(
\naiveprolong{\varphi_1^{-1}(\nbige)}\otimes\omega_{\widetilde{Z}_1}
 \bigr)_{|\widetilde{C}_1}$.
Then it is easy to derive that
the sections $u_{\vecn\,|\,\widetilde{C}}$
give the meromorphic sections of
$\bigl(
 \naiveprolong{\tilde{\nbige}}\otimes\omega_{\widetilde{Z}}
\bigr)_{|\tilde{C}}$.
\hfill\qed

\vspace{.1in}
Let us pick a point $P$ of $\tilde{Z}$.
We can take a neighbourhood $\nbigu$ of $\tilde{Z}$
on which we have a coordinate $(\zeta_1,\ldots,\zeta_n)$
such that
$\nbigu\cap \tilde{D}\simeq
\bigcup_{i=1}^l \{\zeta_i=0\}$
and $\zeta_i(P)=0$.
Such an open subset is called an admissible open set.
We remark that we can consider
the filtrations $\lefttop{\lbar}\Vzero$
on each admissible open subset.
Let $\tilde{D}_{\nbigu,i}$ denote the component of
$\nbigu\cap\tilde{D}$ corresponding to $\{\zeta_i=0\}$.
Let $\pi_i$ denote the projection
of $\nbigu$ onto $\{\zeta_i=0\}$.

We pick a frame $\vecv=(v_i)$ of $\Vzero_{<0}\tildegbige$
on $\nbigu$.
Let us consider the restriction of $u_{\vecn}$
to $\nbigu\cap U$,
where $u_{\vecn}$ are obtained
in the development (\ref{eq;a11.23.15})
for the section $f\in\nbigm'$.
Then we have the description
$u_{\vecn}=\sum \alpha_{i\,\vecn}v_{i}$,
where $\alpha_{i\,\vecn}$ are holomorphic 
on $\nbigu-(\nbigu\cap \tilde{D})$.
Due to Lemma \ref{lem;9.30.10},
the restrictions of $\alpha_{i\,\vecn}$
to the curves $\{\lambda\}\times\pi_p^{-1}(Q)$
are meromorphic for any point
$(\lambda,Q)
 \in\Delta(\lambda_0,\epsilon_0)\times \tilde{D}^{\circ}_{p,\nbigu}$.
Here we put
$\tilde{D}_{\nbigu,p}^{\circ}:=
\tilde{D}_{\nbigu,p}
 -\bigcup_{r\neq p}(\tilde{D}_{\nbigu,p}\cap \tilde{D}_{\nbigu,r})$.
Due to Lemma \ref{lem;9.30.3},
the degree $\lefttop{i}\deg^{\Vzero}(\alpha_{j,\vecn})$
are bounded, independently of $(\lambda,Q)$.
It implies the following lemma.
\begin{lem}
The restrictions of $u_{\vecn}$ to $\nbigu\cap U$
give sections of $\naiveprolong{\tildenbige}\otimes\omega_{\widetilde{Z}}$.
\hfill\qed
\end{lem}

We may take a finite covering of $\widetilde{Z}$
by admissible open subsets.
Hence we obtain the following:
\begin{lem} \label{lem;b12.4.5}
$u_{\vecn}$ give sections of
$\naiveprolong{\tildenbige}\otimes\omega_{\widetilde{Z}}$.
\hfill\qed
\end{lem}

\begin{lem} \label{lem;b12.4.10}
For any sufficiently large number $N$,
the sections
$u_{\vecn}\cdot g^N$ $(\vecn\in\seisuu^n)$
are contained in $\tilde{\gbige}$.
Moreover
they are contained in
$\lefttop{\lbar}\Vzero_{<0}(\tilde{\gbige})$
on each admissible open subset.
\end{lem}
\pf
It immediately follows from Lemma \ref{lem;b12.4.5}.
\hfill\qed

\begin{cor}\label{cor;b12.4.15}
Let $f$ be a section of $\nbigm'$.
There exists a positive integer
such that
$f\cdot g^N$ naturally gives the section of
$\widetilde{\nbigm'}$
\end{cor}
\pf
It immediately follows from Lemma \ref{lem;b12.4.10}.
\hfill\qed

\vspace{.1in}

Let us show the equality
$\nbigm'=\widetilde{\nbigm}'$.
To show it,
we consider the inclusion 
$X\lrarr X\times\cnum$ for the graph of $g$.
Let $t$ denote the coordinate of $\cnum$.
We have only to show
$i_{g\,\ast}\nbigm'=i_{g\,\ast}\widetilde{\nbigm}'$.

We denote the $V$-filtration along $t=0$
of $i_{g\,\ast}\nbigm'$ and $i_{g\,\ast}\widetilde{\nbigm}'$
by $\Vzero$.

Let $f$ be a section of
$\Vzero_{<0}\bigl(i_{g\,\ast}\nbigm'\bigr)$.
We have the development:
\[
 f=\sum f_i\cdot \deldel_t^i.
\]
Due to Corollary \ref{cor;b12.4.15},
we have an integer $N$ such that
$f\cdot t^N\in \Vzero_{<0}\bigl(i_{g\,\ast}\widetilde{\nbigm}'\bigr)$.
Note $t=g$ on the graph of $g$.

The rests are completely same as the Step 1.
Since $f\cdot t^N$ is contained in
$\Vzero_{<0}i_{g\,\ast}\nbigm'\cap\Vzero_{<0}i_{g\,\ast}\widetilde{\nbigm}'$,
we have
$(f\cdot t^N)\cdot\nbigr\subset
 \nbigm'\cap\widetilde{\nbigm}'$.
The $V$-filtrations
$\Vzero(\nbigm')$ and $\Vzero(\widetilde{\nbigm}')$
induce the same $V$-filtration
on $(f\cdot t^N)\cdot\nbigr$
(Lemma \ref{lem;9.30.5}).

Since the multiplication of $g$
on $\nbigm'$ and $\widetilde{\nbigm}'$ are injective,
it can be checked easily that
the multiplication of $t$ on
$i_{g\,\ast}\nbigm'$ and $i_{g\,\ast}\widetilde{\nbigm}'$
are injective.
Then we obtain the implication
$\nbigm'\subset\widetilde{\nbigm}'$
by the same argument as the proof of Lemma \ref{lem;9.30.7},
and we obtain the equality
$\nbigm'=\widetilde{\nbigm}'$
by the same argument as the proof of Lemma \ref{lem;a11.23.20}.
Thus we are done
when Condition \ref{condition;b12.4.20} is satisfied.

\subsubsection{The end of the proof of Theorem \ref{thm;11.5.6}}

Let us consider the general case.
Let $U'$ be a Zariski open subset of $U$
as in Lemma \ref{lem;a11.23.30}.
We put $Y':=Z-U'$.
Let us take a birational morphism
$\varphi:\tilde{Z}'\lrarr Z$
such that $\tilde{Z}'$ is smooth and
$\tilde{D}':=\varphi^{\prime\,-1}(Y')$ is a normal crossing divisor.
We put
$(\tilde{E}',\delbar_{\tilde{E}'},\tilde{h}',\tilde{\theta}')
:=\varphi^{\prime\,-1}(\harmonicbundle)$.
Then we obtain the regular pure twistor $D$-module
$\nbigm(\tilde{E}'):=(\tilde{\gbige}',\tilde{\gbige}',\tilde{\gbigc}')$.
Taking the the push-forward via $\varphi'$,
and taking the component whose strict support is $X$,
we obtain $\tilde{\nbigt}$ as is explained in the subsubsection
\ref{subsubsection;9.30.11}.
Due to the result in the subsubsection \ref{subsubsection;a11.23.35},
we have the naturally defined isomorphism
$\tilde{\nbigt}'\simeq\nbigt$.
By applying the same consideration to $\tilde{\nbigt}$,
we have the isomorphism $\tilde{\nbigt}'\simeq\tilde{\nbigt}$.
Thus we obtain the isomorphism
$\tilde{\nbigt}\simeq\nbigt$.
\hfill\qed

%% file: d7.tex

\subsubsection{The correspondence in the pure imaginary case}

Let $X$ be a complex manifold,
and $Z$ be an irreducible closed subset of $X$.
Let $\VPTgen^{pi}(Z,w)$ denote the set of the equivalence classes
of generically defined tame variation of polarized pure twistor
structures of weight $w$ on $Z$.
Let $\MPT^{pi}(Z,w)$ denote the set of the equivalence classes
of polarized regular pure twistor module of weight $w$
whose strict support is $Z$.
(See the subsubsection \ref{subsubsection;04.2.16.200}
 for the definitions.)
Restricting the correspondence of Theorem \ref{thm;b12.5.20}
to the pure imaginary cases,
we have the following theorem.

\begin{thm}\label{thm;04.2.4.1}
We have the bijective correspondence
$\VPTgen^{pi}(Z,w)\simeq \MPT^{pi}(Z,w)$.
\end{thm}

We have only to show the following:
\begin{itemize}
\item
The image of  $\VPTgen^{pi}(Z,w)$
via the map $\VPTgen(Z,w)\lrarr\MPT(Z,w)$
(the subsection \ref{subsection;a11.23.51})
is contained in $\MPT^{pi}(Z,w)$.
The sketch of the proof is given in the subsubsection 
\ref{subsubsection;04.2.4.30}.
\item
The image of $\MPT^{pi}(Z,w)$
via the map $\MPT(Z,w)\lrarr\VPTgen(Z,w)$
(the subsection \ref{subsection;04.1.27.35})
is contained in $\VPTgen^{pi}(Z,w)$.
\end{itemize}

Both of them are the minor modifications
of the previous argument.
Hence we give only a sketch of the proof.

\subsubsection{The prolongation preserves the pure imaginary property}

\label{subsubsection;04.2.4.30}

Let $U$ be a smooth Zariski open subset of $Z$.
Let $\harmonicbundle$ be a tame pure imaginary harmonic bundle on $U$.
We have the prolongment
$\nbigt=\bigl(\gbige,\gbige,\gbigc\bigr)$ as pure twistor $D$-module on $Z$
of $\harmonicbundle$.

\begin{prop}\label{prop;04.2.4.20}
The pure twistor $D$-module
$\nbigt\otimes\nbigo(w)$ is pure imaginary.
\end{prop}
\pf
We use an induction on the dimension of $X$.
Since the property is local,
we can assume that
$X=\Delta^n=\bigl\{(z_1,\ldots,z_n)\,\big|\, |z_i|<1\bigr\}$,
$D_i=\bigl\{z_i=0\bigr\}$
and $D=\bigcup_{i=1}^l D_i$.

\begin{lem}\label{lem;04.2.4.21}
In the case $g=\prod_{i=1}^n z_i^{m_i}$,
we have the following:
\begin{itemize}
\item
 $\tilde{\psi}_{g,u}\nbigt=0$ unless $u\in\real\times(\sqrt{-1}\real)$.
\item
 The twistor $D$-module $P_h\Gr^{W(N)}_h\tildepsi_{t,u}\nbigt$ is
 pure imaginary.
\end{itemize}
\end{lem}
\pf
From (\ref{eq;04.2.4.5}), Lemma \ref{lem;9.22.16}
and the strictness of $\Psi_b$ (Proposition \ref{prop;9.22.9}),
we obtain
$\psi_{t,u}\bigl(i_{g\,\ast}\gbige\bigr)=0$
unless $u\in\real\times\bigl(\sqrt{-1}\real\bigr)$.
Then the first claim immediately follows.
From Lemma \ref{lem;9.22.22} and the strictness of
$\psi_{t,u}\bigl(i_{g\,\ast}\gbige\bigr)$,
we obtain the vanishing
$\lefttop{i}\psi_{v}\psi_{t,u}\bigl(i_{g\,\ast}\gbige\bigr)=0$
unless $v\in\real\times\bigl(\sqrt{-1}\real\bigr)$.
It follows that
$\lefttop{i}\psi_v\Gr^{W(N)}\psi_{t,u}\bigl(i_{g\,\ast}\gbige\bigr)=0$
unless $v\in\real\times\bigl(\sqrt{-1}\real\bigr)$.
Recall that we have the decomposition
$P_h\Gr^{W(N)}_h\psi_{t,u}\bigl(i_{g\,\ast}\gbige\bigr)
=\bigoplus_{I}\nbigm_I$ as in Proposition \ref{prop;b12.6.150}.
Then we obtain
$\psi_{z_i,v}\nbigm=0$ unless
$v\in\real\times\bigl(\sqrt{-1}\real\bigr)$.
It means that
the tame harmonic bundle corresponding to
the pure twistor $D$-modules
$\bigl(\nbigm_I,\nbigm_I,\gbigc_I\bigr)$
given as in the subsubsection \ref{subsubsection;04.2.4.10}
are pure imaginary.
Due to the hypothesis of our induction for the proof of
Proposition \ref{prop;04.2.4.20},
we obtain that the twistor $D$-modules
$\bigl(\nbigm_I,\nbigm_I,\gbigc_I\bigr)$
are pure imaginary.
Thus we can conclude that
$P_h\Gr^{W(N)}_h\tildepsi_{t,u}\nbigt$ are pure imaginary.
\hfill\qed

\vspace{.1in}

Let us return to the proof of Proposition \ref{prop;04.2.4.20}.
Let $g$ be a general holomorphic function on $X$.
Let us take the birational map $\pi:\widetilde{X}\lrarr X$
such that $\bigl(g\circ\pi\bigr)^{-1}(0)$ is normal crossing
and that $\widetilde{X}$ is smooth.
Note that $\pi^{-1}\harmonicbundle$ is pure imaginary
(Lemma 7.1 in \cite{mochi3}).
Let $\tilde{\nbigt}$ be the pure twistor $D$-module
corresponding to $\pi^{-1}\harmonicbundle$.
We put $\tilde{g}:=g\circ\pi$.
From Lemma \ref{lem;04.2.4.21},
we know that
$\tildepsi_{\tilde{g},u}\tilde{\nbigt}=0$
unless $u\in\real\times\bigl(\sqrt{-1}\real\bigr)$
and that
$P_h\Gr^{W(N)}_h\tildepsi_{\tilde{g},u}\tilde{\nbigt}$ is pure imaginary.
Note that
$\nbigt$ is the direct summand of
$\pi_{+}\tilde{\nbigt}$ whose strict support is $X$.
Then we obtain
$\tildepsi_{g,u}\nbigt=0$
unless $u\in\real\times\bigl(\sqrt{-1}\real\bigr)$
from Theorem 3.3.15 in \cite{sabbah2},
and we obtain that
$P_h\Gr^{W(N)}_h\tildepsi_{g,u}\nbigt$ is pure imaginary
from the argument in the section 6.2.b in \cite{sabbah}.
Thus the proof of Proposition \ref{prop;04.2.4.20} is accomplished.
\hfill\qed

\vspace{.1in}
We can show that
the image of $\VPTgen^{pi}(Z,w)$ via the map
$\VPTgen^{pi}(Z,w)\lrarr \MPT(Z,w)$ is contained in
$\MPT^{pi}(Z,w)$,
by using Proposition \ref{prop;04.2.4.20}
and the argument in the subsubsection \ref{subsubsection;04.2.4.25}.

\subsubsection{The generic part is also pure imaginary}

Let $\nbigt=\bigl(\gbige,\gbige,\gbigc\bigr)$ is a 
pure imaginary pure twistor $D$-module whose strict support is $Z$.
We obtain the harmonic bundle $\harmonicbundle$
on a Zariski dense subset $U$ of $Z$,
as in the subsubsection \ref{subsubsection;04.1.27.70}.
Let us show that $\harmonicbundle$ is pure imaginary.

First we consider the case $\dim Z=1$.
Let $Q$ be a point of $Z-U$.
Since the property is local,
we may assume the condition \ref{condition;04.2.4.50}.
Let $(E_1,\delbar_{E_1},\theta_1,h_1)$
be as in the subsubsection
\ref{subsubsection;04.1.27.70}.
Due to the result of the section 5 in \cite{sabbah2},
$(E_1,\delbar_{E_1},\theta_1,h_1)$ is pure imaginary.
Or, we can directly check it
by seeing the monodromy of
$\bigl(\nbige_1^{\lambda},\DD^{\lambda,f}\bigr)$
for generic $\lambda$.
Let $\pi:\tilde{Z}\lrarr Z$ be a resolution,
and then $\pi^{-1}\harmonicbundle$ is 
a direct summand of 
$\bigl(q_1\circ\pi\bigr)^{-1}\bigl(E_1,\delbar_{E_1},\theta_1,h_1\bigr)$.
Since
$\bigl(q_1\circ\pi\bigr)^{-1}\bigl(E_1,\delbar_{E_1},\theta_1,h_1\bigr)$
is pure imaginary,
we obtain that $\harmonicbundle$ is pure imaginary.

Let us consider the general case.
Let $\pi:\tilde{Z}\lrarr Z$ be a resolution
such that $\tilde{D}:=\pi^{-1}\bigl(Z-U\bigr)$ is a normal crossing 
divisor of $\tilde{Z}$.
Let $C$ be a curve contained in $\tilde{Z}$,
which intersects with the smooth part of $\tilde{D}$ transversally.
Let us consider the specialization of
$\nbigt$ to $\pi(C)$,
which is also pure imaginary.
Then we obtain that 
$\pi^{-1}\harmonicbundle_{|C\setminus \tilde{D}}$
is pure imaginary from the result in the one dimensional case,
which is already discussed above.
Therefore we obtain that $\harmonicbundle$ is pure imaginary.

\vspace{.1in}

Then it is easy to derive that
the image of $\MPT^{pi}(Z,w)$ via the map
$\MPT(Z,w)\lrarr \VPTgen(Z,w)$ is contained in
$\VPTgen^{pi}(Z,w)$.
In all, the proof of Theorem \ref{thm;04.2.4.1} is accomplished.
\hfill\qed

%% file: d8.tex

Let $X$ be a complex manifold,
and $Z$ be an irreducible closed subset of $X$.
Let $\nbigt=(\nbigm',\nbigm'',\gbigc)$
be a pure imaginary pure twistor $D$-module.
$\Xi_{Dol}(\nbigt)$ is defined to be
the $D$-module $\nbigm'_{|\lambda=1}$.
We refer the 4.1.i in \cite{sabbah2}
for the fundamental property of the functor
$\Xi_{Dol}$.
In particular, we recall the following:
\begin{prop}
[Sabbah \cite{sabbah2}] \mbox{{}}
 $\Xi_{Dol}(\nbigt)$ is a semisimple regular holonomic $D$-module.
\hfill\qed
\end{prop}
We have a Zariski open subset $U$ of $Z$
such that the restriction $\nbigt_{|U}$ is a variation 
of polarized pure twistor structures.
By considering the restriction to $\lambda=1$,
we obtain the flat bundle, or equivalently the local system $L$.
We remark that it is semisimple.
Let $\nbigf$ be the regular holonomic $D$-module on $Z$,
which corresponds to the intermediate extension of $L$
via the Riemann-Hilbert correspondence.
Since $\Xi_{Dol}(\nbigt)$ is semisimple,
we have the isomorphism
$\nbigf\simeq \Xi_{Dol}(\nbigt)$.

Let $\RHD^{ss}(Z)$ denote the set of the equivalence classes of
the semisimple regular holonomic $D$-modules
whose strict support is $Z$.
As is explained above, we have the map
$\Xi_{Dol}:\MPT^{pi}(Z,0)\lrarr \RHD^{ss}(Z)$.

\begin{thm}[The conjecture of Sabbah] \label{thm;04.2.5.20}
The map $\Xi_{Dol}:\MPT^{pi}(Z,0)\lrarr \RHD^{ss}(Z)$
is surjective.
\end{thm}
\pf
Let $\nbigf$ be the element of $\RHD^{ss}(Z)$.
We have a smooth Zariski open subset $U$ of $Z$
such that $\nbigf_{|U}$ corresponds to a flat bundle.
We can take a tame pure imaginary pluri-harmonic metric $h$
of $\nbigf_{|U}$ (Theorem 6.2 in \cite{mochi3}).
Hence we obtain the tame pure imaginary harmonic bundle
$\harmonicbundle$ on $U$.
Let $\nbigt=(\gbige,\gbige,\gbigc)$
denote the pure imaginary pure twistor $D$-module
corresponding to $\harmonicbundle$
given in Theorem \ref{thm;04.2.4.1}.
Since the restrictions of
the semisimple regular holonomic $D$-modules
$\nbigf$ and $\Xi_{Dol}(\nbigt)$ to $U$
are isomorphic to $\nbigf_{|U}$,
they are isomorphic.
Thus we obtain the surjectivity desired.
\hfill\qed

\begin{rem}
We formulate the theorem as the surjectivity
between the sets of the equivalence classes,
$\Xi_{Dol}$ is given as the functor
which is compatible with the vanishing cycle functors
and the push-forward (see {\rm\cite{sabbah2}}).
As a result, we obtain the regular holonomic version
of Kashiwara's conjecture,
combining the results of Sabbah and us.
\hfill\qed
\end{rem}

%% file: a62.tex

We recall the definitions of pure twistor $D$-modules
and its polarization, due to Sabbah, given in Chapter 4
in his paper \cite{sabbah}.
The definitions
are based on the work of Saito on his celebrated pure Hodge modules,
as Sabbah himself noted in \cite{sabbah}.
We shall consider only the regular holonomic case.

\subsubsection{Pure twistor $D$-modules}

We only consider the regular case.
Let $X$ be an $n$-dimensional complex manifold.
Let $w$ be an integer.
The definition of pure twistor $D$-modules of weight $w$
is given as follows, inductively.
\begin{df}[Definition {\rm 4.1.2. \cite{sabbah}}]\label{df;11.5.10}
The category $MT^{(r)}_{\leq d}(X,w)$ is defined to be
the full subcategory of $\rtriplecat(X)$,
whose object $(\nbigm',\nbigm'',C)$ satisfies the following:
\begin{description}
\item[(HSD)]
 The $\nbigr$-modules $\nbigm'$ and $\nbigm''$
 are holonomic and strictly $S$-decomposable.
 The dimension of their strict support is less than $d$.
\item[(REG)]
 Let $U$ be any open subset of $X$,
 and $f$ be any holomorphic function on $U$.
 Then $\nbigm'_{|U}$ and $\nbigm''_{|U}$ are regular
 along $\{f=0\}$.
\item[($MT_{0}$)]
 Let us consider the case $d=0$.
 Let $\{x_i\}$ be the union of the strict supports of
 $\nbigm'$ and $\nbigm''$.
 Then we have
 $\bigl(\nbigm'_{|\{x_i\}},\nbigm''_{|\{x_i\}},C\bigr)
 =i_{\{x_i\}\,+}\bigl(\nbigh',\nbigh'',C_0\bigr)$,
where $\bigl(\nbigh',\nbigh'',C_0\bigr)$
 is a pure twistor structure of dimension $0$ and weight $w$.
\item[($MT_{>0}$)]
 Let us consider the case $d>0$.
 Let $U$ be any open subset of $X$,
 and $f$ be any holomorphic function on $U$.
 Let $u$ be any element of
 $\real\times\cnum-\bigl(\seisuu_{\geq 0}\times\{0\}\bigr)$,
 and $l$ be any integer.
 Then the induced $\nbigr$-triple
 $\Gr^W_l\psitilde_{f,u}(\nbigm',\nbigm'',C):=
 \bigl(
 \Gr^W_{-l}\psitilde_{f,u}(\nbigm'),
 \Gr^W_l\psitilde_{f,u}(\nbigm''),
 \Gr^W_l\psitilde_{f,u}C
 \bigr)$
 is the object of $MU^{(r)}_{\leq d-1}(U,w+l)$.

\mbox{{}}
\hfill\qed
\end{description}
\end{df}

\begin{rem}\label{rem;11.5.1}
{\rm
Since our index set of the KMS-structure is 
not $\real\times\{0\}$
but $\real\times\cnum$,
our pure twistor $D$-module is wider than 
those given by Sabbah's definition.
However it is just a minor modification,
and the results and the proofs in the section 4.1. in \cite{sabbah}
are valid, when we change as follows:
\begin{itemize}
\item
``$\alpha\in \closedopen{-1}{0}$'' $\Longrightarrow$ 
``$u\in\real\times\cnum-\bigl(\seisuu_{\geq 0}\times\{0\}\bigr)$''.
\item
``$\alpha\in [-1,0]$'' $\Longrightarrow$ 
``$u\in\real\times\cnum$''.
\item
``$\psi_{t,\alpha}$'' $\Longrightarrow$ ``$\psitilde_{t,u}$''.
\item
``$z_o\in \cnum^{\ast}$'' $\Longrightarrow$
``any generic $\lambda_0$''.
This case appears only in Proposition 4.1.21 (1) in \cite{sabbah}.
\hfill\qed
\end{itemize}
}
\end{rem}

For example,
Proposition 4.1.3 is changed as follows:
\begin{lem}
Let $(\nbigm',\nbigm'',C)$ is an object of
$MT^{(r)}_{\leq d}(X,w)$.
Then $\nbigm'$ and $\nbigm''$ are strict.

Let $U$ be any open subset of $X$,
and $f$ be any holomorphic function on $U$.
Let $u$ be any element of $\real\times\cnum$.
Then $\tilde{\psi}_{t,u}(\nbigm')$ and $\tilde{\psi}_{t,u}(\nbigm'')$
are strict.
\hfill\qed
\end{lem}

Since the paper \cite{sabbah} is very well written,
and since the modification is quite minor,
we do not reproduce the propositions and the proof, here.
We recommend the reader to see \cite{sabbah}.

\subsubsection{Polarization}

We recall only the definition of polarization of
pure twistor $D$-modules,
which is given in the subsection 4.2 in \cite{sabbah}
\begin{df}[Definition {\rm 4.2.2} in {\rm\cite{sabbah}}]\label{df;11.5.11}
Let $\nbigt$ be an object of $MT^{(r)}_{\leq d}(X,w)$.
A polarization of $\nbigt$ is 
a sesqui-linear Hermitian duality
$\nbigs:\nbigt\lrarr\nbigt^{\ast}(-w)$ of weight $w$
satisfying the following:
\begin{description}
\item[($MTP_0$)]
Let us consider the case $d=0$.
Let $\{x_i\}$ be the strict support of $\nbigt$.
We may assume that
$\nbigt=i_{\{x_i\}\,+}\bigl(\nbigh',\nbigh'',C_0\bigr)$.
Then we have $\nbigs=i_{\{x_i\}\,+}\nbigs_0$
for some polarization $\nbigs_0$
of $\bigl(\nbigh',\nbigh'',C_0\bigr)$.

\item[($MTP_{>0}$)]
Let us consider the case $d>0$.
Let $U$ be any open subset of $X$,
and $f$ be any holomorphic function on $U$.
Let $u$ be an element of
$\real\times\cnum-\{\seisuu_{\geq \,0}\times\{0\}\}$,
and $l$ be any non-negative integer.
Then $P\Gr^W_l\tildepsi_{t,u}S$ gives a polarization
of $P_l\Gr^W_l\tildepsi_{t,u}\nbigt$.
\hfill\qed
\end{description}
\end{df}

\begin{rem}{\rm
Again the propositions and the proofs
in the subsection 4.2 in \cite{sabbah}
are valid, when we change them as in Remark \ref{rem;11.5.1}.
}
\end{rem}

%% file: d6.tex

\subsubsection{The definition}

We have the inclusion of the set of the pure imaginary numbers
$\sqrt{-1}\real\subset\cnum$.
\begin{df}[Pure imaginary pure twistor $D$-module]
The category $MT^{sabbah}_{\leq d}(X,w)$ is defined to be
the full subcategory of $\nbigr$-Triples(X),
whose objects $(\nbigm',\nbigm'',C)$ satisfy the following conditions,
in addition to the conditions in Definition {\rm\ref{df;11.5.10}}:
\begin{description}
\item[(PI)]
 We have 
 $\widetilde{\psi}_{t,u}(\nbigm')=\widetilde{\psi}_{t,u}(\nbigm'')=0$ 
 unless $u\in\real\times (\sqrt{-1}\real)$.
\end{description}
Such objects are called pure imagianry pure twistor $D$-module.

A polarization of a pure imaginary pure twistor $D$-module
is defined by Definition {\rm\ref{df;11.5.11}}.
\hfill\qed
\end{df}

\begin{rem}{\rm
Sabbah give the definition of pure twistor $D$-modules in \cite{sabbah2},
which is our pure imaginary pure twistor $D$-module.
Since his purpose is to attack Kashiwara's conjecture,
it is natural to restrict the attention to pure imaginary case.
On the other hand, 
our definition may be natural from the view point of
Simpson's Meta Theorem.
\hfill\qed
}
\end{rem}

%% file: a62.1.tex

Let $C$ be a holomorphic curve.

Generalizing the result in Chapter 5 of \cite{sabbah},
we can show the following:
\begin{thm} \label{thm;11.5.5}
Let $w$ be any integer.
The variation of polarized pure twistor structures of weight $w$
which are generically defined over $C$,
are bijectively corresponds to
the regular pure twistor $D$-modules of weight $w$
whose strict support is $C$.
\end{thm}

\begin{rem}{\rm
Theorem {\rm\ref{thm;11.5.5}} is a special case
of Theorem {\rm\ref{thm;11.5.2}}
and Theorem \ref{thm;11.5.6}.
Although we use the decomposition theorem
in the proof Theorem \ref{thm;11.5.6},
we do not need the decomposition theorem in the proof,
if the dimension of the base manifold is $1$,
for we do not have to consider the blow up
in that case.
However, we give
an explanation of the correspondence
and an outline of the proof.
\hfill\qed
}
\end{rem}

\noindent
{\bf Outline of the proof of Theorem \ref{thm;11.5.5}}
We only consider the case $w=0$,
for we have only to consider the tensor product with $\nbigo(w)$
to obtain the other cases.

Let $\harmonicbundle$ be a tame harmonic bundle
defined over a Zariski open subset $C'\subset C$.
Then we obtain the pure twistor $D$-module 
$\nbigt(E):=(\gbige,\gbige,\gbigc)$,
as proved in Theorem \ref{thm;11.5.2}.
It is easy to see that
the restriction $\nbigt(E)_{|C'}$ gives
the harmonic bundle $\harmonicbundle$.

On the other hand,
Let $\nbigt=(\nbigm,\nbigm,C,Id)$ be a regular pure twistor $D$-module
on $C$.
Then we have a Zariski open subset $C'$,
such that $\nbigt_{|C'}$ gives a harmonic bundle $\harmonicbundle$.
Due to the regularity, it is tame.
Then we obtain the pure twistor $D$-modules
$\nbigt(E)$.

\begin{lem}
We have the naturally defined isomorphism
$\nbigt\lrarr\nbigt(E)$.
\end{lem}
\pf
On the Zariski open subset $C'$,
we have the isomorphism
$\nbigt_{|C'}\lrarr \nbigt(E)_{|C'}$,
due to our construction of $\nbigt(E)_{|C'}$.
We shall prolong the morphism defined on $C'$
to the morphism defined on $C$.
For that purpose,
we have only to discuss
on a neighbourhood of a point $P\in C-C'$.
So we can assume that $C=\Delta$ and $C'=\Delta^{\ast}$.
Let $z$ be the coordinate of $\Delta$.

Let us pick a point $\lambda_0\in\cnum_{\lambda}$ which is generic
with respect to the sets $KMS(\nbigm,z)$ and $KMS(\gbige,z)$.
Then the restrictions
$\nbigm_{|\{\lambda_0\}\times C}$ and $\gbige_{|\{\lambda_0\}\times C}$
can be regarded as regular holonomic $D$-modules,
which are strictly $S$-decomposable.
Let $L$ be a local system
corresponds to the flat connection
$V=\nbigm_{|\{\lambda_0\}\times \Delta^{\ast}}
=\gbige_{|\{\lambda_0\}\times\Delta^{\ast}}$.
By Riemann-Hilbert correspondence,
both of regular holonomic $D$-modules
$\nbigm_{|\lambda_0}$ and $\gbige_{|\lambda_0}$
correspond to the perverse sheaf which is
the intermediate extension of $L$.
Thus we have the isomorphism
$\nbigm_{|\lambda_0}\simeq\gbige_{|\lambda_0}$.
If we regard them as the submodules
of $\iota_{\ast}V$,
then we have $\nbigm_{|\lambda_0}=\gbige_{|\lambda_0}$.
The $V$-filtrations are also same.
(See the proof of Lemma \ref{lem;c12.4.1} and
 Lemma \ref{lem;a11.26.1}.)

Let us pick a point $\lambda_0\in\cnum_{\lambda}$,
which is not necessarily generic.
Let $f$ be a section of
$V^{(\lambda_0)}_{<0}\nbigm$ on
$\Delta(\lambda_0,\epsilon_0)\times \Delta$.
Let us pick any point $\lambda\in\Delta(\lambda_0,\epsilon_0)$,
which is generic.
Then the restriction $f_{|\lambda}$ gives
a section of $\Vzero_{<0}\gbige_{|\lambda}$.

The section $f$ also gives the holomorphic section
of $\nbige$ over $\Delta(\lambda_0,\epsilon_0)\times\Delta^{\ast}$,
which we denote by $f$.
Since $f_{|\lambda}$ gives a section of
$\Vzero_{<0}\gbige_{|\lambda}$,
we have $-\ord(f_{|\lambda})<1$.
Due to Corollary \ref{cor;11.28.15},
we can conclude that $f$ gives the section
of $\Vzero_{<0}\gbige$.
Thus we obtain the morphism
$\Vzero_{<0}\nbigm\lrarr \Vzero_{<0}\gbige$
on $\Delta(\lambda_0,\epsilon_0)\times \Delta$.
Since $\nbigm$ is generated by
$\Vzero_{<0}\nbigm$,
we obtain the morphism $\nbigm\lrarr\gbige$
defined over $\Delta(\lambda_0,\epsilon_0)\times \Delta$,
and thus the morphism defined over $\cnum_{\lambda}\times\Delta$.

Once we obtain the morphism $\nbigm\lrarr\gbige$
whose restriction to $C'$ is isomorphic,
the morphism is isomorphic on $C$,
due to the strictly $S$-decomposability of the both sides.
It is easy to see that the sesqui-linear pairings also coincide
(Corollary \ref{cor;a11.23.6}).

\hfill\qed

%% file: a62.2.tex

Let $X$ and $Y$ be complex manifolds,
and $f$ be a projective morphism from $X$ to $Y$.
Let $c$ be the first Chern class of a relatively ample line bundle
on $X$ with respect to $f$.
The decomposition theorem and the hard Lefschetz theorem
for pure twistor $D$-modules is as follows.

\begin{thm}[Theorem 6.1.1 in \cite{sabbah}]
Let $(\nbigt,\nbigs)$ be a polarized pure twistor $D$-module
of weight $w$ on $X$.
Then
$\bigl(\bigoplus_i f_{+}^i\nbigt,\nbigl_c,\bigoplus_i f^i_{+}\bigr)$
is a polarized graded Lefschetz twistor $D$-modules.
\end{thm}

The decomposition theorem for polarized pure twistor $D$-module
was shown by Sabbah in \cite{sabbah}
(based on Saito's argument in \cite{saito1}).
Although our pure twistor $D$-modules
(Definition \ref{df;11.5.10} and Definition \ref{df;11.5.11})
are wider than those given in \cite{sabbah},
the argument of Sabbah-Saito essentially works.
Roughly speaking,
their argument is divided into the following two steps.
\begin{description}
\item[Step 1]
By using the induction on the dimensions
of $\Supp(\nbigt)$ and $f\bigl(\Supp(\nbigt)\bigr)$,
the problem is reduced to the decomposition theorem
for a pure twistor $D$-modules on a smooth projective curve.
(See the sections 6.2.b and 6.2.c in \cite{sabbah}.
See also the section 5.3 in \cite{saito1}.)
\item[Step 2]
We prove the decomposition theorem
for the pure twistor $D$-module on a smooth projective curve.
(See the sections 6.2.a in \cite{sabbah}.
See also the sections 6, 7 and 11 in \cite{z}.)
\end{description}

Since we do not have to change the argument for the step 1,
we only give an argument for the step 2.
Although it is just a minor modification 
of the argument given in \cite{z} and \cite{sabbah},
we give some detail.

%% file: a61.tex

\subsubsection{Preliminary calculation of cohomology}

Following Zucker,
we put $\gbigm_1:=\gbiga/\gbigb$,
where $\gbiga$ and $\gbigb$ are given as follows:
\[
\begin{array}{l}
 \gbiga:=\Bigl\{
 \mbox{\rm measurable function $f$}\,\Big|\,
 \int_0^A |f(r)|^2\cdot |\log r|\cdot r\cdot dr<\infty
 \,\,\mbox{\rm for some $0<A<1$}
 \Bigr\},\\
\mbox{{}}\\
\gbigb:=
\Bigl\{
 f\in\gbiga\,\Big|\,
 f=u'\,\,\mbox{\rm weakly}\,\,
 \int_0^A |u|^2\cdot |\log r|^{-1}r^{-1}dr
 \,\,\mbox{\rm for some $0<A<1$}
 \Bigr\}.
\end{array}
\]

In the following,
we use the Poincar\`{e} metric and the induced volume form
on $\Delta^{\ast}$.

Let $V$ be a $C^{\infty}$-vector bundle of rank 1
on $\Delta^{\ast}$
with a trivialization $e$.
Let $\nabla$ be the flat connection of $V$
such that $\nabla(e)=\alpha\cdot e\cdot dz/z$
for some complex number $\alpha$.
Let $h$ be a metric of $V$
such that
$h(e,e)= r^{-2a}\cdot|\log r|^{k}$
for some real number $a$ and an integer $k$.

Let $\nbigl^p(V)_{(2)}$ be the sheaf
of germs of locally $L^2$-section of $V\otimes\Omega^p$
for which $\nabla(\phi)$ is $L^2$-section of $V\otimes\Omega^{p+1}$.
Then we obtain the complex of sheaves $\nbigl^{\cdot}(V)_{(2)}$:
\[
 \nbigl^0(V)_{(2)}
\lrarr
 \nbigl^1(V)_{(2)}
\lrarr
 \nbigl^2(V)_{(2)}.
\]

Let $\nbigh^i(\nbigl^{\cdot}(V)_{(2)})$
denote the stalk of the $i$-th cohomology sheaves
of $\nbigl^{\cdot}(V)_{(2)}$
at $O$.

\begin{lem}\label{lem;11.4.3}
In the case $\alpha\in\cnum-\seisuu$,
we have $\nbigh^i(\nbigl^{\cdot}(V)_{(2)})=0$
for $i=0,1,2$.
\end{lem}
\pf
The claim $i=0$ can be shown by a direct calculation.
The claims for $i=1,2$ are shown by Zucker
in Prop (11.3) in his paper.
Note that the condition $-1<\alpha=-a<0$ is imposed
in the section 11 in his paper.
However the only assumption $\alpha\not\in\seisuu$
is used in the proof of of Proposition (11.2).
\hfill\qed

\vspace{.1in}

Let $H(V)$ denote the space of the flat sections of $V$.

\begin{lem} \label{lem;11.4.2}
In the case $\alpha=0$,
we have the following:
\begin{equation}\label{eq;11.3.10}
 \nbigh^0\bigl(\nbigl^{\cdot}(V)_{(2)}\bigr)
=\left\{
 \begin{array}{ll}
 H(V) & (a<0,\,\,\mbox{\rm or }\,\, a=0,k\leq 0)\\
 \mbox{{}}\\
 0    & (\mbox{\rm otherwise}).
 \end{array}
 \right.
\end{equation}
\begin{equation}
\nbigh^1\bigl(\nbigl^{\cdot}(V)_{(2)}\bigr)
=\left\{
 \begin{array}{ll}
 {\displaystyle\frac{dt}{t}\otimes H(V)},
 & (a<0,\,\,\mbox{\rm or } a=0,k\leq -2)\\
 \mbox{{}} \\
 \gbigm_1\cdot dr\otimes H(V),
 & (a=0,k=1)\\
 \mbox{{}}\\
 0 & (\mbox{\rm otherwise})
 \end{array}
 \right.
\end{equation}
\begin{equation}
\nbigh^2\bigl(\nbigl^{\cdot}(V)_{(2)}\bigr)
=\left\{ 
 \begin{array}{ll}
 {\displaystyle \gbigm_1\cdot dr\wedge \frac{dt}{t}\otimes H(V)},
 & (a=0,k=-1)\\
 \mbox{{}}\\
 0 & (\mbox{\rm otherwise}).
 \end{array}
 \right.
\end{equation}
\end{lem}

The proof of Lemma \ref{lem;11.4.2},
due to Zucker,
will be given in the subsubsections
\ref{subsubsection;a11.9.16}--\ref{subsubsection;a11.9.15}.
Before entering the proof,
we give the formulas of the norms
of the sections of $V\otimes\Omega^{i}$.
\begin{equation} \label{eq;11.3.1}
 \int
 |f\cdot e|^2
 \dvol
=\int|f|^2\cdot r^{-2a}\cdot |\log r|^k\cdot
 \frac{dr\cdot d\theta}{r\cdot |\log r|^2}
=\int |f|^2\cdot r^{-2a-1}\cdot
 |\log r|^{k-2}\cdot dr\cdot d\theta.
\end{equation}
\begin{equation}\label{eq;11.3.2}
\int
 |f\cdot dr\cdot e|^2
 \dvol
=\int|f|^2\cdot r^2\cdot |\log r|^2\cdot
 r^{-2a}|\log r|^k\cdot
 \frac{dr\cdot d\theta}{r\cdot |\log r|^2}
=\int |f|^2\cdot r^{-2a+1}\cdot
 |\log r|^{k}\cdot dr\cdot d\theta.
\end{equation}
\begin{equation} \label{eq;11.3.3}
\int |f\cdot d\theta\cdot e|^2\cdot \dvol
=\int |f|^2\cdot |\log r|^2\cdot r^{-2a}\cdot |\log r|^k
 \cdot \frac{dr\cdot d\theta}{r|\log r|^2}
=\int |f|^2\cdot r^{-2a-1}|\log r|^k
 \cdot dr\cdot d\theta.
\end{equation}
\begin{equation} \label{eq;11.3.4}
\int |f\cdot dr\cdot d\theta\cdot e|^2\dvol
=\int |f|^2r^{-2a-1}\cdot |\log r|^{k+2}\cdot dr\cdot d\theta.
\end{equation}
Here $f$ denotes $C^{\infty}$-functions,
and the metrics of $V\otimes\Omega^i$ are induced
by the metric $h$ of $V$ and the Poincar\`{e} metric.

We denote the $L^2$-norm by $||\cdot||_{(2)}$.

\subsubsection{The calculation of $\nbigh^0$}
\label{subsubsection;a11.9.16}

From (\ref{eq;11.3.1}),
we have $\int |e|^2\dvol<\infty$
if and only if
we have $a<0$ or $a=0,k\leq 0$.
Thus (\ref{eq;11.3.10}) immediately follows.

\subsubsection{The calculation of $\nbigh^1$}

We need a preparation.
Let $\omega=f\cdot dr+g\cdot d\theta$ be a $C^{\infty}$-section
of $E\otimes\Omega^1$,
whose support is compact.
We have the Fourier developments:
\[
 f=\sum f_n\cdot e^{\sqrt{-1}n\theta},
\quad
 g=\sum g_n\cdot e^{\sqrt{-1}n\theta}.
\]
We take $\epsilon_n$ as follows:
\[
 d\omega=
 \sum (g_n'-\sqrt{-1}n\cdot f_n)\cdot e^{\sqrt{-1}n\theta}
 dr\wedge d\theta
=\sum \epsilon_n\cdot e^{\sqrt{-1}n\theta} dr\wedge d\theta.
\]
We put $u_n:=(n\sqrt{-1})^{-1}g_n$ for $n\neq 0$,
and we put $u:=\sum_{n\neq 0}u_n\cdot e^{\sqrt{-1}n\theta}$.
Then the $n$-th Fourier coefficient $(du)_n$
of $du$ is as follows:
\[
 (du)_n:=
 u_n'+\sqrt{-1}n\cdot u_n\cdot d\theta
=(n\sqrt{-1})^{-1}g_n'\cdot dr+g_n\cdot d\theta
=f_n\cdot dr+g_n\cdot d\theta
-\frac{\sqrt{-1}}{n}\epsilon_n\cdot dr.
\]
Thus we obtain the following:
\[
 U:=du-\omega-(f_0\cdot dr+g_0\cdot d\theta)
=-\sqrt{-1}\sum_{n\neq 0}n^{-1}\cdot
 \epsilon_n\cdot e^{\sqrt{-1}n\theta}\cdot dr.
\]
\begin{lem}
We have $||U\cdot e||_{(2)}\leq ||d\omega\cdot e||_{(2)}$.
\end{lem}
\pf
From (\ref{eq;11.3.2}) and (\ref{eq;11.3.4}),
we obtain the following:
\[
 ||U\cdot e||_{(2)}^2
=\int |U\cdot e|^2\dvol
=\sum_{n\neq 0}
 n^{-2}\int |\epsilon_n|^2\cdot r^{-2a+1}
 |\log r|^{k}\cdot dr.
\]
\[
 ||d\omega\cdot e||^2_{(2)}
=\int |d\omega\cdot e|^2\cdot \dvol
=\sum_{n\neq 0}
 \int |\epsilon_n|^2\cdot r^{-2a+1}
 \cdot |\log r|^k\cdot dr.
\]
Thus we are done.
\hfill\qed

\vspace{.1in}

Let $\omega=f\cdot dr+g\cdot d\theta$ be a section of $\nbigl^1$
such that $d\omega=0$.
We have the Fourier development
$f=\sum f_n\cdot e^{\sqrt{-1}n\theta}$
and
$g=\sum g_n\cdot e^{\sqrt{-1}n\theta}$.
We can take a sequence
$\{\omega^{(k)}\,|\,k=1,2,\ldots\}$
satisfying the following:
\begin{itemize}
\item
 $\omega^{(k)}$ is $C^{\infty}$-section of
 $V\otimes\Omega^1$, whose support is compact.
\item
 $\{\omega^{(k)}\}$ converges to $\omega$ in $L^2$.
\item
 $\{d\omega^{(k)}\}$ converges to $0$ in $L^2$.
\end{itemize}

For each $\omega^{(k)}$,
we take $u^{(k)}\in \nbigl^0(V)_{(2)}$ as above.
Then $\{u^{(k)}\}$ converges to
$\sum_{n\neq 0}n^{-1}\cdot g_n\cdot e^{\sqrt{-1}n\theta}$,
and we have the following:
\[
 ||du^{(k)}-\omega^{(k)}-
 (f_0^{(k)}dr+g_0^{(k)}d\theta)||_{(2)}
\leq
 ||d\omega^{(k)}||_{(2)}
\to 0,\quad (k\to\infty).
\]
Hence $\omega-(f_0\cdot dr+g_0\cdot d\theta)$
is coboundary.

We also have $(d\omega^{(k)})_0=g_0^{(k)}\cdot dr\wedge d\theta$.
Thus we have
$\lim_{k\to\infty} g_0^{(k)}=g_0$ and
$\lim_{k\to\infty} g_0^{(k)\,\prime}=0$.
Thus $g_0$ is constant.
Moreover
we have
$\int |g_0\cdot d\theta\cdot e|^2\dvol<\infty$
if and only if
we have $a<0$ or $a=0,k<-1$.
In the case $a<0$ or $a=0,k<-1$,
we have $(\log r)\cdot e\in L^2$
and $d\log r\cdot e\in L^2$.
Thus $g_0\frac{dt}{t}\cdot e$
and $\sqrt{-1}g_0\cdot d\theta\cdot e$
are same modulo $d\cdot\nbigl^0$.

Let us consider the equation $u_0'=f_0$.
From (\ref{eq;11.3.2}),
we have 
$\int |f_0|\cdot r^{-2a+1}\cdot |\log r|^k dr<\infty$.

\begin{lem} \label{lem;11.3.20}
If $a<0$ or $a=0,k<1$,
then $f_0\cdot dr$ is a coboundary.
\end{lem}
\pf
We put $\int_A^r f_0\cdot d\rho$.
We have the following:
\begin{multline}
 |u_0|^2
\leq \int_r^A|f_0|^2\rho\cdot |\log \rho|^{1-\epsilon}
 \cdot
 \int_{r}^A\rho^{-1}|\log \rho|^{-1+\epsilon}d\rho
=\epsilon^{-1}
 \int_r^A|f_0|^2\cdot\rho\cdot|\log \rho|^{1-\epsilon} d\rho
 \cdot
 (\log A^{\epsilon}-\log r) \\
\leq 2\epsilon^{-1}\cdot
 \int_r^A|f_0|^2\cdot\rho\cdot|\log \rho|^{1-\epsilon} d\rho
 \cdot
 |\log r|.
\end{multline}
Here $\epsilon$ is a positive number.
In the case $a=0$,
we impose the condition $0<\epsilon<1-k$.
Then we obtain the following:
\begin{multline}
||u_0\cdot e||^2_{(2)}
=
\int_0^A|u_0|^2 r^{-2a-1}|\log r|^{k-2}dr
\leq
 2\epsilon^{-1}\cdot
\int_0^A\int_{r}^A|f_0|^2\cdot\rho|\log\rho|^{1-\epsilon} d\rho
\cdot
 |\log r|^{\epsilon}r^{-2a-1}|\log r|^{k-2}dr
 \\
=2\epsilon^{-1}\cdot
 \int_0^A |f_0|^2\cdot \rho|\log\rho|^{1-\epsilon}d\rho
\cdot
 \int_0^{\rho}r^{-2a-1}|\log r|^{k-2+\epsilon}dr
\leq
 C\cdot \int_0^A|f_0|^2\rho^{-2a+1}|\log \rho|^{k-1}\cdot d\rho \\
\leq 
 C\cdot \int_0^A|f_0|^2\rho^{-2a+1}|\log \rho|^{k}\cdot d\rho
=C\cdot ||f_0\cdot dr\cdot e||_{(2)}^2.
\end{multline}
Here $C$ is positive constant depending
only on $\epsilon$, $a$ and $k$.
Thus Lemma \ref{lem;11.3.20} follows.
\hfill\qed

\begin{lem}
If $a>0$ or $a=0,k>1$,
then $f_0\cdot dr$ is a coboundary.
\end{lem}
\pf
We put $u_0:=\int_0^rf_0\cdot d\rho$.
Then we have the following:
\[
 |u_0|^2\leq\int_0^r|f_0|^2\cdot\rho\cdot|\log\rho|^{1+\epsilon}\cdot d\rho
 \cdot\int_0^r\rho^{-1}\cdot|\log \rho|^{-1-\epsilon}d\rho
=\epsilon^{-1}\int_0^r|f_0|^2\cdot\rho\cdot|\log\rho|^{1+\epsilon} d\rho
 \cdot |\log r|^{-\epsilon}.
\]
Here $\epsilon$ denotes a positive number.
In the case $a=0$, we impose the condition
$0<\epsilon<k-1$.
We obtain the following:
\begin{multline}
||u_0\cdot e||_{(2)}^2
=\int_0^A |u_0|^{2}\cdot r^{-2a-1}\cdot |\log r|^{k-2}dr
\leq
 \epsilon^{-1}\cdot\int_0^A\int_0^r
 |f_0|^2\cdot|\log \rho|^{1+\epsilon}
 \cdot d\rho\cdot
 |\log r|^{k-2-\epsilon}\cdot r^{-2a-1}dr \\
=\epsilon^{-1}\cdot
 \int_0^A|f_0|^2\cdot\rho\cdot|\log\rho|^{1+\epsilon}d\rho
\cdot
 \int_{\rho}^A r^{-2a-1}|\log r|^{k-2-\epsilon} dr 
\leq 
 C\int_0^A|f_0|^2\cdot\rho^{-2a+1}|\log \rho|^{k-1}d\rho \\
\leq
 C\cdot \int_0^A|f_0|^2\cdot\rho^{-2a+1}|\log\rho|^kd\rho
=||f_0\cdot dr\cdot e||_{(2)}^2.
\end{multline}
Thus we are done.
\hfill\qed

\vspace{.1in}
In the case $a=0,k=1$,
the group $\gbigm_1\otimes dr\otimes H(V)$ remains.

\subsubsection{The calculation of $\nbigh^2$}
\label{subsubsection;a11.9.15}

We put $\omega=f\cdot dr\wedge d\theta$.
We have the Fourier development
$f=\sum_n f_n\cdot e^{\sqrt{-1}n\theta}$.
For $n\neq 0$,
we put $g_n:=(n\sqrt{-1})^{-1}f_n$
and $g=\sum_{n\neq 0}g_n\cdot e^{\sqrt{-1}n\theta}$.
Then we have
$d(g\cdot dr)=\omega-f_0\cdot dr\wedge d\theta$.
We also have the following:
\begin{multline}
||g\cdot dr||_{(2)}^2
=2\pi\sum_{n\neq 0}
\int |g_n|^2\cdot r^{-2a}\cdot |\log r|^k\cdot r\cdot dr
=2\pi\sum_{n\neq 0} n^{-2}
\int |f_n|^2\cdot r^{-2a+1}\cdot|\log r|^{k}dr \\
\leq
 C\sum_{n\neq 0}n^{-2}\int |f_n|^2\cdot r^{-2a+1}|\log r|^{k+2}dr
\leq 
C\cdot ||\omega||_{(2)}^2.
\end{multline}
Hence $\omega-f_0\cdot dr\wedge d\theta$ is a coboundary.

Let us consider the equation
$d(\eta_0)=f_0\cdot dr\wedge d\theta\cdot e$ 
for $\eta_0=h_0\cdot d\theta\cdot e$.
In the case $a<0$ or $a=0,k<-1$,
we put $h_0=\int_A^r f_0\cdot d\rho$.
In the case $a>0$ or $a_0,k>-1$,
we put $h_0=\int_0^r f_0\cdot d\rho$.
Then, as before $h_0$ is the  $L^2$-section.
Thus $f_0\cdot dr\wedge d\theta$ is a coboundary.

In the case $a=0,k=-1$,
the group $\gbigm_1\cdot dr\wedge d\theta$ remains.
Note we have $\sqrt{-1}dr\wedge d\theta=dr\wedge dt/t$.

Thus the proof of Lemma \ref{lem;11.4.2} is finished.
\hfill\qed

%% file: a61.1.tex

\subsubsection{The complexes
 $\nbigl^{\cdot}(E)$ and
 $\bigl(\nbigelambda\otimes\Omega^{\cdot,0}\bigr)_{(2)}$}

Let $E$ be a holomorphic vector bundle over $\Delta^{\ast}$.
Let $\DDlambda$ be the $\lambda$-connection on $E$,
and $h$ be a hermitian metric of $E$.
Let $\nbigl^p(E)$ denote the sheaf
on $\Delta$ of germs of locally $L^2$-sections $\phi$
of $E\otimes\Omega^p$
for which $\DD^{\lambda}\phi$ are $L^2$-sections of
$E\otimes\Omega^{p+1}$.
Then we obtain the complex of sheaves
$\nbigl^{\cdot}(E)$
induced by $\DDlambda$.

Let us consider a tame harmonic bundle
$\harmonicbundle$ over $\Delta^{\ast}$.
We obtain $(\nbigelambda,\DDlambda,h)$
over $\Delta^{\ast}$.
We have the projection:
\[
 \prolong{\nbigelambda}_{|O}
\stackrel{\pi}{\lrarr}
 \Gr^F_0(\nbigelambda_{|O})
=\bigoplus_{\beta\in\cnum} 
 \EE\bigl(
 \Gr^F_0\bigl(\nbigelambda_{|O}\bigr),\beta
 \bigr).
\]
We have the weight filtration $W$
on $\EE\bigl(\Gr^F_0\bigl(\nbigelambda_{|O}\bigr),\beta\bigr)$,
induced by the nilpotent part of the residue.
We put as follows:
\[
 I:=W_1\EE\bigl(\Gr^F_0\bigl(\nbigelambda_{|O}\bigr),0\bigr)
\oplus
 \bigoplus_{\beta\neq 0}
 W_{-1}\EE\bigl(
 \Gr^F_0\bigl(\nbigelambda_{|O}\bigr),\beta
 \bigr).
\]
We put as follows:
\[
 \nbigelambda_{(2)}:=
 \bigl\{
 f\in\prolong{\nbigelambda}\,\big|\,
 \pi(f_{|O})\in I
 \bigr\}.
\]
We also put as follows:
\[
 \bigl(
 \nbigelambda\otimes\Omega^{1,0}
 \bigr)_{(2)}
:=\bigl\{
 f\in\bigl(
 \prolong{\nbigelambda}\otimes\Omega^{1,0}
 \bigr)_{(2)}\,\big|\,
 \pi(f)\in W_{-1}\Gr^F_0\bigl(\nbigelambda_{|O}\bigr)
 \bigr\}.
\]
The $\lambda$-connection $\DDlambda$
induces the morphism
$\nbigelambda_{(2)}\lrarr
\bigl(\nbigelambda\otimes\Omega^{1,0}\bigr)_{(2)}$.
We denote the complex by
$\bigl(
 \nbigelambda\otimes\Omega^{\cdot,0}
 \bigr)_{(2)}$.

\begin{lem}
We have the following:
\[
 \nbigelambda_{(2)}=\bigl\{
 f\in\nbigelambda\,\big|\,
 ||f||_{(2)}<\infty,\,\,
 || \DDlambda f||_{(2)}<\infty
 \bigr\},
\]
\[
 \nbigelambda\otimes\Omega^{1,0}_{(2)}
=\bigl\{
 f\in\nbigelambda\otimes\Omega^{1,0}\,\big|\,
 ||f||_{(2)}<\infty
 \bigr\}.
\]
\end{lem}
\pf
It is clear from our construction
and the norm estimate.
\hfill\qed

\vspace{.1in}
We have the naturally defined morphism
$\Psi:\bigl(\nbigelambda\otimes\Omega^{\cdot,0}\bigr)_{(2)}
\lrarr
 \nbigl^{\cdot}(\nbigelambda)$.

\begin{prop}\label{prop;11.4.5}
The morphism $\Psi$ is quasi isomorphic.
\end{prop}
The proposition \ref{prop;11.4.5} is also essentially due to Zucker
\cite{z}.
We divide the proof into the two cases.
The case $\lambda=0$ is discussed
in the subsubsections
\ref{subsubsection;d12.4.1}--\ref{subsubsection;a11.9.22}.
The case $\lambda\neq 0$ is discussed
in the subsubsections
\ref{subsubsection;d12.4.2}--\ref{subsubsection;a11.9.31}.

\subsubsection{The case $\lambda=0$, Preliminary}
\label{subsubsection;d12.4.1}

Let us consider the case $\lambda=0$.
We have $\prolong{E}$ over $\Delta$,
the Higgs field $\theta=f\cdot dz/z$
and the metric $h$.
We have the decomposition (the subsubsection \ref{subsubsection;a11.9.17}):
\[
 E=\bigoplus_{\alpha\in\Sp(f_{|O})}E_{\alpha},
\quad
 f=\bigoplus_{\alpha\in \Sp(f_{|O})}f_{\alpha},\quad
 f_{\alpha}\in \End(E_{\alpha}).
\]
We may assume that there exists a 
sufficiently small positive number $\epsilon$
such that 
$|\alpha-\beta|<\epsilon$
for any eigenvalue of $f_{\alpha\,|\,P}$ $(P\in \Delta^{\ast})$.

Let $\vecv=(v_i)$ be a frame of $\prolong{E}$
compatible with the following:
\begin{itemize}
\item
 The decomposition $E=\bigoplus E_{\alpha}$.
\item
 The parabolic filtration $F$ of $E_{|O}$.
\item
 The weight filtration on $\Gr^F(E_{|O})$.
\end{itemize}
For each $v_i$,
we put as follows:
\[
 b(v_i):=\deg^F(v_i),
\quad
 \alpha(v_i):=\deg^{\EE}(v_i),
\quad
 h(v_i):=\frac{\deg^W(v_i)}{2}.
\]
If we put
$v_i':=v_i\cdot (-\log |z_i|)^{-h(v_i)}\cdot 
 |z|^{b(v_i)}$
and $\vecv':=(v_i')$.
Recall that $\vecv'$ is adapted over $\Delta^{\ast}$
due to Simpson
(see the subsubsection 4.3.3 in \cite{mochi}).
Let us consider the metric $\tilde{h}$
determined as follows:
\[
 \tilde{h}(v_i,v_j):=
 \delta_{i\,j}\cdot |z|^{-2b(v_i)}
 \cdot\bigl(-\log |z|^2\bigr)^{-2h(v_i)}.
\]
Then the metrics $\tilde{h}$ and $h$ are mutually bounded.
Since the complexes
$\nbigl^{\cdot}(E,h)$ and $\nbigl^{\cdot}(E,\tilde{h})$
coincide,
we may use the metric $\tilde{h}$
in the following argument.

We have the decompositions:
\[
 \nbigl^{\cdot}(E)_{(2)}
=\bigoplus_{\alpha\in\cnum}\nbigl^{\cdot}(E_{\alpha})_{(2)},
\quad
 \bigl(
 E\otimes\Omega^{\cdot,0}
 \bigr)_{(2)}
=\bigoplus_{\alpha\in\cnum}
 \bigl(
 E_{\alpha}\otimes\Omega^{\cdot,0}
 \bigr)_{(2)}.
\]
The decompositions are compatible with the morphism
$\Psi$.
The claim of Proposition \ref{prop;11.4.5} for $\lambda=0$
immediately follows from the following lemma.
\begin{lem}\label{lem;a11.9.20}
The morphism 
$\Psi_{\alpha}:
 \bigl(
 E_{\alpha}\otimes\Omega^{\cdot,0}
 \bigr)_{(2)}
\lrarr
 \nbigl^{\cdot}(E_{\alpha})_{(2)}
 $
is quasi isomorphic.
\end{lem}
Lemma \ref{lem;a11.9.20} is proved in the subsubsections
\ref{subsubsection;a11.9.21}--\ref{subsubsection;a11.9.22}.

\subsubsection{The proof of Lemma \ref{lem;a11.9.20} ($\alpha=0$)}
\label{subsubsection;a11.9.21}

In the case $\alpha=0$,
the morphism $\theta_0:E_0\lrarr E_0\otimes\Omega^{1,0}$ is bounded.
Let us consider the sheaf
$\nbigl^{p,q}(E_0)_{(2)}$ on $\Delta$
of germs of locally $L^2$-section of $E_0\otimes \Omega^{p,q}$,
for which $\delbar \phi$ is $L^2$.
Then we obtain the complex of sheaves
$\nbigl^{p,\cdot}(E_0)_{(2)}$:
\[
 \nbigl^{p,0}(E_0)_{(2)}
\stackrel{\delbar}{\lrarr}
 \nbigl^{(p,1)}(E_0)_{(2)}.
\]

We denote the cohomology sheaves of
$\nbigl^{p,\cdot}(E_0)_{(2)}$ by $\nbigh^{p,\cdot}$.

\begin{lem} \label{lem;11.4.1}
We have the following isomorphisms:
\[
 \begin{array}{l}
 \nbigh^{0,0}\simeq E_{0\,(2)},
\quad
 \nbigh^{1,0}\simeq \bigl(E_0\otimes\Omega^{1,0}\bigr)_{(2)},\\
 \mbox{{}}\\
 \nbigh^{0,1}\simeq\gbigm_1\otimes d\bar{t}\otimes
 \Gr^W_1\Gr^{F}_0\bigl( E_{0\,|\,O} \bigr), \\
 \mbox{{}}\\
 \nbigh^{1,1}\simeq \gbigm_1\otimes d\bar{t}\otimes
 {\displaystyle \frac{dt}{t}}
 \otimes \Gr^W_{-1}\Gr^F_0(E_{0\,|\,O}).
 \end{array}
\]
\end{lem}
\pf To show Lemma \ref{lem;11.4.1},
we need some preparation.
We consider the lexicographic order on
$\real\times \seisuu$.
Let us consider the filtration $\tilde{F}$
of $E_0$ in the category of vector bundles
indexed by $\real\times\seisuu$,
given as follows:
\[
 \tilde{F}_{(b,k)}(E_0)
:=\big\langle
 v_i\,\big|\,
 \deg^{\EE}(v_i)=0,\,\,
 \bigl(\deg^F(v_i),\deg^W(v_i)\bigr)
 \leq (b,k)
  \big\rangle.
\]
Then we obtain the vector bundle
$V_{b,k}=\Gr^{\tilde{F}}_{(b,k)}(E_0)$
on $\Delta^{\ast}$.
The tuple
$\vecv_{0\,b,k}
 =\{v_i\,|\,\alpha(v_i)=0,b(v_i)=b,2h(v_i)=k\}$
naturally induces the frame of $V_{b,k}$.
We denote it by
$\tilde{\vecv}_{0\,b,k}=\bigl(\tilde{v}_i\,|\,
 \alpha(v_i)=0,b(v_i)=b,2h(v_i)=k \bigr)$.

We have the naturally induced metric $\tilde{h}_{b,k}$
on $V_{b,k}$,
for which we have
$\tilde{h}_{b,k}\bigl(\bar{v}_i,\bar{v}_j\bigr)
=\delta_{i\,j}\cdot |z|^{-2b}\cdot \bigl|\log |z|\bigr|^{k}$.
For the complex $\bigl(\nbigl^{p,\cdot}(V_{b,h}),\delbar\bigr)$,
we recall Proposition 6.4 and Proposition 11.5 in \cite{z}.
\begin{lem}
We have the following isomorphisms:
\[
 \nbigh^0\bigl(
 \nbigl^{p,\cdot}(V_{b,k})
 \bigr)
\simeq
 \bigl(V_{b,k}\otimes\Omega^{p,0}\bigr)_{(2)},
\]
\[
 \nbigh^1\bigl(
 \nbigl^{0,\cdot}(V_{b,k})
 \bigr)
=\left\{
 \begin{array}{ll}
 \gbigm_1\otimes V_{b,k\,|\,O}\cdot d\bar{t}
 & (b_0,k=1), \\
 \mbox{{}}\\
 0 & (\mbox{\rm otherwise}),
 \end{array}
 \right.
\]
\[
 \nbigh^2\bigl(
 \nbigl^{1,\cdot}(V_{b,k})
 \bigr)
=\left\{
 \begin{array}{ll}
 \gbigm_1\otimes V_{b,k\,|\,O}\cdot d\bar{t}\frac{dt}{t}
 & (b=0,k=-1),\\
 \mbox{{}}\\
 0 & (\mbox{\rm otherwise}).
 \end{array}
 \right.
\]
\end{lem}
\pf
The case $b=0$ is shown in Proposition 6.4 in \cite{z}.
and the case $-1<b<0$ is shown in Proposition 11.5 in \cite{z}.
The general case can be reduced to the two cases above.
\hfill\qed

\vspace{.1in}
We have the naturally defined inclusion
$\bigl(
\Gr^{\tilde{F}}_{(b,k)}(E_0)\otimes\Omega^{i,0}
\bigr)_{(2)}\lrarr
\nbigl^{i,\cdot}\bigl(\Gr^{\tilde{F}}_{(b,k)}(E_0)\bigr)$,
and we have the isomorphism:
\[
 \bigl(
\Gr^{\tilde{F}}_{(b,k)}(E_0)\otimes\Omega^{i,0}
\bigr)_{(2)}
\simeq
 \nbigh^0\bigl(
 \nbigl^{i,\cdot}\bigl(\Gr^{\tilde{F}}_{(b,k)}(E_0)\bigr)
 \bigr).
\]
We also have the isomorphism
$(E_0\otimes\Omega^{i,0})_{(2)}
 \simeq\nbigh^0\bigl(\nbigl^{i,\cdot}(E_0)\bigr)_{(2)}$.
Then Lemma \ref{lem;11.4.1} is obtained
by an easy spectral sequence argument.
\hfill\qed

\vspace{.1in}
We have the morphism of the complexes:
$\theta_0:\nbigl^{0,\cdot}(E_0)_{(2)}
  \lrarr\nbigl^{1,\cdot}(E_0)_{(2)}$.
It is easy to see that
it induces the isomorphism
$\nbigh^1\bigl(
 \nbigl^{0,\cdot}(E_0)_{(2)}
 \bigr)\lrarr
 \nbigh^1\bigl(
 \nbigl^{1,\cdot}(E_0)_{(2)}
 \bigr)$.
Then it follows from an easy spectral sequence argument,
that the inclusion
$\bigl(E_0\otimes\Omega^{\cdot\,0}\bigr)_{(2)}
\lrarr \nbigl^{\cdot}(E_0)_{(2)}$
is quasi isomorphic.

\subsubsection{The proof of Lemma \ref{lem;a11.9.20} $(\alpha\neq 0)$}
\label{subsubsection;a11.9.22}

This is the easier part in the proof of Lemma \ref{lem;a11.9.20}.
In the case $\alpha\neq 0$,
the morphism
$\theta_{\alpha}:E_{\alpha}\lrarr E_{\alpha}\otimes\Omega^{1,0}$
is invertible.
Let $F:E_{\alpha}\otimes\Omega^{1,0}\lrarr E_{\alpha}$
denote the inverse.
Then there exists a positive constant $C$
such that
$|F|_h\leq C\cdot (-\log |z|)^{-1}$.

The morphism $F$ induces the morphism
$E_{\alpha}\otimes\Omega^{1,1}\lrarr E_{\alpha}\otimes\Omega^{0,1}$,
which we denote also by $F$.

\begin{lem}
We have $\nbigh^2\bigl(\nbigl^{\cdot}(E_{\alpha})\bigr)=0$.
\end{lem}
\pf
Let $g$ be an $L^2$-section of $E_{\alpha}\otimes\Omega^{1,1}$.
Then $F(g)$ is the $L^2$-section of $E_{\alpha}\otimes\Omega^{0,1}$
such that
$\theta_{\alpha}(F(g))=g$.
It implies $\nbigh^2(\nbigl^{\cdot}(E_{\alpha}))=0$.
\hfill\qed

\vspace{.1in}

\begin{lem}
We have $\nbigh^1\bigl(\nbigl^{\cdot}(E_{\alpha})_{(2)}\bigr)=0$
and $\nbigh^1\bigl(
 \bigl(E_{\alpha}\otimes\Omega^{\cdot,0}\bigr)_{(2)}\bigr)=0$.
\end{lem}
\pf
The vanishing
$\nbigh^1\bigl(
 \bigl(E_{\alpha}\otimes\Omega^{\cdot,0}\bigr)_{(2)}\bigr)=0$
can be checked directly.
Let us show the vanishing
$\nbigh^1\bigl(\nbigl^{\cdot}(E_{\alpha})_{(2)}\bigr)=0$.
Let $g_1\cdot d\bar{z}+g_2\cdot dz$
be an $L^2$-section of $E_{\alpha}\otimes\Omega^1$
such that the following holds:
\[
 (\delbar+\theta)(g_1\cdot d\bar{z}+g_2\cdot dz)
=\theta(g_1)\cdot d\bar{z}+\delbar g_2\cdot dz=0.
\]
We put $g_3:=F(g_2\cdot dz)$,
which is a $L^2$-section of $E_{\alpha}$.
Then we have $\theta(g_3)=g_2\cdot dz$ by our construction.
We also have the following:
\[
 \theta_{\alpha}(\delbar g_3)
=-\delbar\theta_{\alpha}(g_3)
=-\delbar(g_2)\cdot dz=\theta_{\alpha}(g_1)\cdot d\bar{z}.
\]
Since $\theta_{\alpha}$ is invertible,
we obtain $\delbar g_3=g_1$.
Thus we obtain 
the vanishing $\nbigh^1\bigl(\nbigl^{\cdot}(E_{\alpha})_{(2)}\bigr)$.
\hfill\qed

\begin{lem}
We have the vanishings
$\nbigh^0\bigl(\nbigl^{\cdot}(E_{\alpha})_{(2)}\bigr)=0$
and
$\nbigh^0\bigl(
 \bigl(E_{\alpha}\otimes\Omega^{\cdot,0}\bigr)_{(2)}
 \bigr)=0$.
\end{lem}
\pf
It immediately follows from the invertibility
of $\theta_{\alpha}$.
\hfill\qed

\vspace{.1in}
Therefore the proof of Lemma \ref{lem;a11.9.20} is finished.
\hfill\qed

\subsubsection{The case $\lambda\neq 0$, Preliminary}
\label{subsubsection;d12.4.2}

Let us continue to prove the case $\lambda\neq 0$
in Proposition \ref{prop;11.4.5}.
We consider the complex induced by the flat connection
$\DD^{\lambda,f}$ instead of $\DD^{\lambda}$.
Clearly, the complexes are naturally quasi isomorphic.

Let us consider the space $H\bigl(\nbigelambda\bigr)$
of the multi-valued flat section of $\nbigelambda$.
We have the filtration $\nbigf H\bigl(\nbigelambda\bigr)$
by an increasing order,
and we have the generalized eigen decomposition
$H\bigl(\nbigelambda\bigr)
=\bigoplus_{\omega\in \Sp^f(\nbigelambda)}
 \EE\bigl(H\bigl(\nbigelambda\bigr),\omega\bigr)$.
We also have the weight filtration $W$
on $\Gr^{\nbigf}\bigl(H\bigl(\nbigelambda\bigr)\bigr)$.

Let $\vecs$ be a frame of $H\bigl(\nbigelambda\bigr)$,
which is compatible with the decomposition $\EE$,
the filtration $\nbigf$,
and the filtration $W$.
We put as follows:
\[
 \omega(s_i):=\deg^{\EE}(s_i),
\quad
 b(s_i):=\deg^{\nbigf}(s_i),
\quad
 h(s_i):=\frac{1}{2}\deg^W(s_i).
\]
We put $v_i:=F(s_i,0)$,
and $\vecv=(v_i)$.
Then we have $\DD \vecv=\vecv\cdot A\cdot dz/z$
for some constant matrix $A\in M(r)$,
whose eigenvalues $\alpha$ satisfy
the condition $0\leq \Re(\alpha)<1$.

Let $\alpha(v_i)$ denote the complex number
satisfying
$\exp\bigl(-2\pi\sqrt{-1}\cdot \alpha(s_i)\bigr)=\omega(s_i)$
and
$0\leq \Re\alpha(s_i)<1$.
We also put as follows:
\[
 b(v_i):=b(s_i)-\Re\alpha(v_i),
\quad
 h(v_i):=h(s_i).
\]

We put $v_i':=v_i\cdot |z|^{b(v_i)}\cdot (-\log |z|)^{-h(v_i)}$
and $\vecv':=(v_i')$.
\begin{lem}
The $C^{\infty}$-frame $\vecv'$ over $\Delta^{\ast}$
is adapted.
\end{lem}
\pf
It is easy to reduce the claim
to the norm estimate in one dimensional case.
\hfill\qed

\vspace{.1in}

Let us consider the metric $\tilde{h}$
defined as follows:
\[
 \tilde{h}(v_i,v_j)=\delta_{i\,j}\cdot
 |z|^{-2b(v_i)}\bigl(-\log|z|\bigr)^{2h(v_i)}.
\]
Then the metrics $\tilde{h}$ and $h$ are mutually bounded.
Thus we use $\tilde{h}$ in the following.

We have the decomposition
$\nbigelambda=\bigoplus_{\alpha\in\cnum}\nbige^{\lambda}_{\alpha}$,
where $\nbige^{\lambda}_{\alpha}$ is given as follows:
\[
 \nbige^{\lambda}_{\alpha}
:=\big\langle
 v_i\,\big|\,\alpha(v_i)=\alpha
 \big\rangle.
\]
The morphism $\Psi$ is compatible with the decomposition,
i.e.,
$\Psi$ is a direct sum of the morphisms
$\Psi_{\alpha}:
 \bigl(\nbigelambda_{\alpha}\otimes\Omega^{\cdot,0}\bigr)_{(2)}
\lrarr
 \nbigl^{\cdot}(\nbigelambda_{\alpha})$.
The claim of Proposition \ref{prop;11.4.5} for $\lambda\neq 0$
immediately follows from the following lemma.
\begin{lem}\label{lem;a11.9.30}
The $\Psi_{\alpha}$ is quasi isomorphic.
\end{lem}

Lemma \ref{lem;a11.9.30}, which is essentially due to Zucker,
is proved in the subsubsections
\ref{subsubsection;d12.4.3}--\ref{subsubsection;a11.9.31}.

Before entering the proof of Lemma \ref{lem;a11.9.30},
we need some preparation.
We have the filtration $\tilde{\nbigf}$ of $\nbige^{\lambda}_{\alpha}$
given as follows:
\[
 \tilde{\nbigf}_{b}\nbige^{\lambda}_{\alpha}
=\big\langle
 v_i\,\big|\,\alpha(v_i)=\alpha,
 b(v_i)\leq b
 \big\rangle.
\]
We have the naturally induced frame
$\vecv_{b,\alpha}
 :=\bigl(v_{i}\,\big|\,\alpha(v_i)=\alpha,b(v_i)=b\bigr)$.
We also have the weight filtration $W$
on $\Gr^{\tilde{\nbigf}}_b\bigl(\nbige^{\lambda}_{\alpha}\bigr)$:
\[
 W_h\Gr^{\tilde{\nbigf}}_b\bigl(\nbige^{\lambda}_{\alpha}\bigr)
=\bigl\langle
 v_i\,\big|\,
 \alpha(v_i)=\alpha,\,b(v_i)=b,\,2h(v_i)\leq h
 \bigr\rangle.
\]
Then we obtain the vector bundle
$V_{b,h}:=\Gr^W_h\Gr^{\tilde{\nbigf}}_b
 \bigl(\nbige^{\lambda}_{\alpha}\bigr)$.
We have the naturally induced frame
$\tilde{\vecv}_{b,h}:=
 \bigl(\tilde{v}_i\,\big|\,
 \alpha(v_i)=\alpha,\,b(v_i)=b,\,2h(v_i)=h
 \bigr)$.
Then we have the induced metric $\tilde{h}_{b,h}$
on $V_{b,h}$,
for which we have
$\tilde{h}_{b,h}\bigl(\tilde{v}_i,\tilde{v}_j\bigr)
=\delta_{i\,j}\cdot |z|^{-2b}\cdot (-\log|z|)^{h}$.
We also have the naturally induced connection
$\nabla$, for which we have
$\nabla(\tilde{v}_i)=\tilde{v}_i\cdot \alpha\cdot dz/z$.

\subsubsection{Proof of Lemma \ref{lem;a11.9.30} (The case $\alpha=0$)}
\label{subsubsection;d12.4.3}

Due to Lemma \ref{lem;11.4.2},
we have the following:
\[
 \nbigh^0\bigl(
 \nbigl^{\cdot}(V_{b,h})_{(2)}
 \bigr)
=\left\{
 \begin{array}{ll}
  H(V_{b,h}), &(b<0,\,\,\mbox{\rm or } b=0,h<0).
 \mbox{{}}\\
 0 & (\mbox{\rm otherwise}).
 \end{array}
 \right.
\]
\[
 \nbigh^1\bigl(
 \nbigl^{\cdot}(V_{b,h})_{(2)}
 \bigr)
=\left\{
 \begin{array}{ll}
 {\displaystyle \frac{dt}{t}\otimes V_{b,h}}
 & (b<0,b=0,h\leq -2), \\
 \mbox{{}}\\
 \gbigm_1\cdot dr\otimes e, & (b_0,h=1)\\
 \mbox{{}}\\
 0 & (\mbox{\rm otherwise}).
 \end{array}
 \right.
\]
\[
 \nbigh^2\bigl(
 \nbigl^{\cdot}(V_{b,h})_{(2)}
 \bigr)
=\left\{
 \begin{array}{ll}
 {\displaystyle\gbigm_1\cdot dr\wedge \frac{dt}{t}\otimes e}
 & (b=0,h=-1),\\
 \mbox{{}}\\
 0 & (\mbox{\rm otherwise}).
 \end{array}
 \right.
\]
On the other hand,
the following can be checked directly:
\[
 \nbigh^0\bigl(
 \nbigl(V_{b,h}\otimes\Omega^{\cdot,0})_{(2)}
 \bigr)
=\left\{
 \begin{array}{ll}
 H\bigl(V_{b,h}\bigr),
 & (b<0,\,\mbox{\rm or }b=0,h\leq 0),\\
 \mbox{{}}\\
 0 & (\mbox{\rm otherwise}).
 \end{array}
 \right.
\]
\[
 \nbigh^1\bigl(
 \bigl(V_{b,h}\otimes\Omega^{\cdot,0}
 \bigr)_{(2)}
 \bigr)
=\left\{
 \begin{array}{ll}
 {\displaystyle
 \frac{dt}{t}\otimes H\bigl(V_{b,h}\bigr),
 } & (b<0,\,\,\mbox{\rm or } b=0,h\leq -2),\\
 \mbox{{}}\\
 0 & (\mbox{\rm otherwise}).
 \end{array}
 \right.
\]
We have the spectral sequence:
\[
 E_1^{p,q}=\nbigh^{p+q}\bigl(
 \Gr^W_{-p}\nbigl^{\cdot}\bigl(
 \Gr^{\nbigf}_b(\nbige^{\lambda}_0)
 \bigr)_{(2)}
 \bigr)
\Longrightarrow
 \nbigh^{p+q}\bigl(
 \nbigl^{\cdot}\bigl(
 \Gr^{\nbigf}_{b}(\nbige^{\lambda}_0)
 \bigr)_{(2)}
 \bigr)
\]
\[
 E_1^{p,q}=\nbigh^{p+q}\bigl(
 \Gr^W_{-p}
 \bigl(
 \Gr^{\nbigf}_b(\nbige^{\lambda}_0)
 \otimes\Omega^{\cdot,0}
 \bigr)_{(2)}
 \bigr)
\Longrightarrow
 \nbigh^{p+q}\bigl(
 \bigl(
 \Gr^{\nbigf}_b\bigl(\nbige^{\lambda}_0\bigr)\otimes\Omega^{\cdot,0}
 \bigr)_{(2)}
 \bigr).
\]

We have $E_1^{p,q}=E_2^{p,q}$ for both of them.
In the calculation of $E_3^{p,q}$,
the $\gbigm_1\cdot dr\otimes e$
in $\nbigh^1(\nbigl^{\cdot}(V_{b,h})_{(2)})$
$(b=0,h=1)$
and
$\gbigm_1\cdot dr\wedge dt/t\otimes e$
in $\nbigh^2(\nbigl^{\cdot}(V_{b,h})_{(2)})$
$(b=0,h=-1)$ are canceled.

The morphism
$\bigl(
\Gr^{\nbigf}_b(\nbige^{\lambda}_0)\otimes\Omega^{\cdot,0}
\bigr)_{(2)}
\lrarr \nbigl^{\cdot}\bigl(
 \Gr^{\nbigf}_b(\nbige^{\lambda}_0)
 \bigr)_{(2)}$ induces
the morphisms of the spectral sequences.
Then the induced morphisms
are isomorphic at the $E_3$-level.
As a consequence,
the morphism
$\bigl(
\Gr^{\nbigf}_b(\nbige^{\lambda}_0)\otimes\Omega^{\cdot,0}
\bigr)_{(2)}
\lrarr \nbigl^{\cdot}\bigl(
 \Gr^{\nbigf}_b(\nbige^{\lambda}_0)
 \bigr)_{(2)}$ is quasi isomorphic.

As a result,
we obtain that
$\bigl(
 \nbige^{\lambda}_0\otimes\Omega^{\cdot,0}
 \bigr)_{(2)}
\lrarr
 \nbigl^{\cdot}(\nbige^{\lambda}_0)_{(2)}$ is
quasi isomorphic.

\subsubsection{The proof of Lemma \ref{lem;a11.9.30}
 (The case $\alpha\neq 0$)}
\label{subsubsection;a11.9.31}

Due to Lemma \ref{lem;11.4.3},
we have the vanishing for $i=0,1,2$
\[
 \nbigh^{i}\bigl(
 \nbigl^{\cdot}
 \bigl(\Gr^{W}_h\Gr^{\nbigf}_b(\nbige_{\alpha}^{\lambda})
 \bigr)_{(2)}
 \bigr)=0.
\]
The following vanishings can be checked by a direct calculation:
\[
 \nbigh^i\bigl(
 \bigl(
 \Gr^W_h\Gr^{\nbigf}_b\bigl(\nbige_{\alpha}^{\lambda} \bigr)
 \bigr)_{(2)}
 \bigr).
\]
Thus the morphism 
$\bigl(\nbigelambda\otimes\Omega^{\cdot,0}\bigr)_{(2)}
\lrarr
 \nbigl^{\cdot}(\nbigelambda)_{(2)}$
is quasi isomorphic.
Therefore 
we obtain Lemma \ref{lem;a11.9.30},
and thus Proposition \ref{prop;11.4.5}.
\hfill\qed

%% file: a61.2.tex

\subsubsection{First replacement}

Let $C$ be a quasi projective curve over $\cnum$,
and $\overline{C}$ be the smooth completion.
We take a Kahler metric of $C$
which is equivalent to the Poincar\`{e} metric
around the points $\overline{C}-C$.
Then $C$ is a complete Kahler manifold.

Let $\harmonicbundle$ be a tame harmonic bundle over $C$.
Then we have the complex of sheaves
$\nbigl^{\cdot}(\nbigelambda)_{(2)}$
and
$\bigl(\nbigelambda\otimes\Omega^{\cdot,0}\bigr)_{(2)}$
over $\overline{C}$.

\begin{prop}
The naturally defined morphisms
$\bigl(\nbigelambda\otimes\Omega^{\cdot,0}\bigr)_{(2)}
\lrarr
 \nbigl^{\cdot}(\nbigelambda)_{(2)}$
are quasi isomorphic.
\end{prop}
\pf
It follows from Proposition \ref{prop;11.4.5}.
\hfill\qed

\subsubsection{Second replacement}

It is not so clear how we consider a family version
of $\nbigl^{\cdot}(\nbigelambda)_{(2)}$ $(\lambda\in\cnum)$.
So we replace them, as is explained in the following.
Let $\DD^{\lambda\,\ast}$ denote the formal adjoint
of the differential operator $\DD^{\lambda}$.
Recall the formula of the Laplacians:
\begin{equation}\label{eq;a12.9.15}
 \square_{\lambda}=
 \DD^{\lambda}\circ\DD^{\lambda\,\ast}+
 \DD^{\lambda\,\ast}\circ\DD^{\lambda}
=\bigl(1+|\lambda|^2\bigr)\square_{0}=
 (\delbar+\theta)\circ(\delbar+\theta)^{\ast}
+(\delbar+\theta)^{\ast}\circ(\delbar+\theta).
\end{equation}

Let $\tilde{\nbigl}^{0}(\nbigelambda)_{(2)}$
denote the sheaf of germs of locally $L^2$-sections $\phi$
of $\nbigelambda$ for which $\square_{\lambda} \phi$ is $L^2$.
Let $\tilde{\nbigl}^1(\nbigelambda)_{(2)}$
denote the sheaf of germs of locally $L^2$-sections $\phi$
of $\nbigelambda\otimes\Omega^1$ for which
$\DD^{\lambda}\phi$ and $\DD^{\lambda\,\ast}\phi$ are $L^2$.
Let $\tilde{\nbigl}^2(\nbigelambda)_{(2)}$ denote 
the sheaf of germs of locally $L^2$-sections $\phi$
of $\nbigelambda\otimes\Omega^2$.
Then we obtain the complex of sheaves
$\tilde{\nbigl}^{\cdot}(\nbigelambda)_{(2)}$.
We have the naturally defined morphism
$\tilde{\nbigl}^{\cdot}(\nbigelambda)_{(2)}
\lrarr
 \nbigl^{\cdot}(\nbigelambda)_{(2)}$.

The following lemma is easy to see.
\begin{lem}
The sheaves $\nbigl^{i}(\nbigelambda)_{(2)}$
and $\tilde{\nbigl}^i(\nbigelambda)_{(2)}$
are soft.
\hfill\qed
\end{lem}

\begin{lem}
The morphism
$\Gamma\bigl(
 \overline{C},\tilde{\nbigl}^{\cdot}(\nbigelambda)_{(2)}
 \bigr)
\lrarr
 \Gamma\bigl(
 \overline{C},\nbigl^{\cdot}(\nbigelambda)_{(2)}
 \bigr)$ is quasi isomorphic.
\end{lem}
\pf
We put
$V:=\bigoplus
 L^2(\nbigelambda\otimes\Omega^{i}\bigr)$.
Since
the operator $\DD^{\lambda}+\DD^{\lambda\,\ast}$
on $V$ is self adjoint.
Thus we have the orthogonal decomposition
$L^2\bigl(
 \nbigelambda\otimes\Omega^i
 \bigr)=\gbigh^i\oplus
 \Image\bigl(\DD^{\lambda}+\DD^{\lambda\,\ast}\bigr)$.
Then we obtain the isomorphisms:
\[
 H^i\bigl(
 \Gamma\bigl(\tilde{\nbigl}(\nbigelambda)_{(2)}\bigr)
 \bigr)
\simeq
 H^i\bigl(
 \Gamma\bigl(
 \nbigl(\nbigelambda)_{(2)}
 \bigr)
 \bigr)
\simeq
 \gbigh^i.
\]
Thus we are done.
\hfill\qed

\vspace{.1in}
The morphism
$\bigl(\nbigelambda\otimes\Omega^{\cdot,0}\bigr)_{(2)}
\lrarr
 \nbigl^{\cdot}(\nbigelambda)_{(2)}$
is decomposed into
the morphisms
$\bigl(\nbigelambda\otimes\Omega^{\cdot,0}\bigr)_{(2)}
\lrarr
 \tilde{\nbigl}^{\cdot}(\nbigelambda)_{(2)}
\lrarr
 \nbigl^{\cdot}(\nbigelambda)_{(2)}$.
Thus we obtain the following.
\begin{prop}\label{prop;e12.4.10}
We have the natural isomorphism
of $\hyperh^i\bigl(
 \bigl(
\nbigelambda\otimes\Omega^{\cdot,0}
 \bigr)_{(2)}
 \bigr)\simeq
 H^i\bigl(
\Gamma\bigl(
 \tilde{\nbigl}^{\cdot}(\nbigelambda)_{(2)}
 \bigr)
 \bigr)$.
\hfill\qed
\end{prop}

%% file: a86.tex

\subsubsection{Family of $L^2$-complexes} 

Let $C$ be a quasi projective curve
and $\overline{C}$ be the smooth completion.
We take a complete Kahler metric of $C$
which is equivalent to Poincar\'{e} metric around
$\overline{C}-C$.
Then we have the naturally defined measure
on $\cnum_{\lambda}\times C$.

Let $\harmonicbundle$ be a tame harmonic bundle over $C$.
Let $p_{\lambda}:\cnum_{\lambda}\times C\lrarr C$
denote the projection.
Then we have the $C^{\infty}$-bundle $\nbige=p_{\lambda}^{-1}E$.
We have the naturally defined metric
on $\nbige\otimes\Omega_C^{p}$.

Let $\nbigl^2\bigl((\lambda,P),\nbige\otimes\Omega_C^p\bigr)$
denote the space of germs of $L^2$-sections
of $\nbige\otimes\Omega^p_C$ at the point $(\lambda,P)$.
Let $f$ be an element of
$\nbigl^2\bigl((\lambda,P),\nbige\otimes\Omega_C^p\bigr)$.
Let us consider the following condition.
\begin{condition}\label{condition;11.4.10}
There exist open subsets $\lambda\in U_1\subset \cnum_{\lambda}$
and $P\in U_2\subset \overline{C}$,
and an element $L^2(U_1\times U_2)$
satisfying the following:
\begin{enumerate}
\item \label{number;11.4.11} The germ of $F$ at $(\lambda,P)$ is $f$.
\item \label{number;11.4.12} Due to Fubini's theorem,
 $F$ induces the measurable function
 $\Phi_F^{U_1\times U_2}:U_1\lrarr L^2(U_2,E\otimes\Omega^p_C)$.
 Then the function $\Phi_F^{U_1\times U_2}$ is holomorphic.
\end{enumerate}
\end{condition}

Note the following easy lemma.
\begin{lem}
Let $f$ be an element of
$\nbigl^2\bigl((\lambda,P),\nbige\otimes\Omega^p_C\bigr)$.
Assume that there exist $U_i$ $(i=1,2)$ and $F$
as in Condition {\rm\ref{condition;11.4.10}}.
Let $\lambda\in U_1'$ and $P\in U_2'$ be open subsets,
and $F'$ be an element of $L^2(U_1'\times U_2')$
whose germ is $f$.
Then there exist 
open subsets $\lambda\in U_1''$ and $P\in U_2''$ be open subsets
such that
and $F'_{|B}$ satisfies the conditions
{\rm\ref{number;11.4.11}} and {\rm\ref{number;11.4.12}}
in Condition {\rm\ref{condition;11.4.10}}.
\hfill\qed
\end{lem}

\begin{df}\mbox{{}}
\begin{itemize}
\item
An element
$f\in\nbigl^2\bigl((\lambda,P),\nbige\otimes\Omega^{p}_C\bigr)$
is called $\lambda$-holomorphic,
if Condition {\rm\ref{condition;11.4.10}} holds
for $f$.
\item
Let $\nbigc\bigl((\lambda,P),\nbige\otimes\Omega^p_C\bigr)$
denote the subspace of
$\nbigl^2\bigl((\lambda,P),\nbige\otimes\Omega^p_C\bigr)$,
which consists of the $\lambda$-holomorphic germs.
\hfill\qed
\end{itemize}
\end{df}

\begin{lem} \label{lem;11.4.17}
The following sequence is exact.
\[
\begin{CD}
 0@>>>
 \nbigc\bigl((\lambda_0,P),\nbige\otimes\Omega^p_{C}\bigr)
 @>{\lambda-\lambda_0}>>
 \nbigc\bigl((\lambda_0,P),\nbige\otimes\Omega^p_C\bigr)
 @>>>\nbigl^2(P,\nbigelambda\otimes\Omega^{p}_C)
 @>>> 0.
 \end{CD}
\]
\hfill\qed
\end{lem}

Let $f$ be an element of
$\nbigc\bigl((\lambda_0,P),\nbige\otimes\Omega^p_C\bigr)$.
Let $U_i$ $(i=1,2)$ and $F$ be as in Condition \ref{condition;11.4.10},
and then we obtain the holomorphic function
$\Phi^{U_1\times U_2}_{F}:U_1\lrarr
L^2\bigl(U_2,E\otimes\Omega^p_C\bigr)$.
Let $U_i'$ $(i=1,2)$ and $F'$ be also
as in Condition \ref{condition;11.4.10},
and then we obtain the holomorphic function
$\Phi^{U_1'\times U_2'}_{F'}:U_1'\lrarr
L^2\bigl(U_2',E\otimes\Omega^p_C\bigr)$.
Let us take a sufficiently small open subset $U_1''\subset U_1\cap U_1'$,
and let us put $U_2''=U_2\cap U_2'$.
Then both of $\Phi^{U_1\times U_2}_F$ and $\Phi^{U_1'\times U_2'}_{F'}$
induce the holomorphic morphisms
$U_1''\lrarr L^2(U_2'',E\otimes\Omega^p_C)$,
and they coincide.
In this sense,
the germ $f\in\nbigc\bigl((\lambda_0,P),E\otimes\Omega^{p}_C\bigr)$
determines the germ $\Phi_f$ of holomorphic function
to the space $\nbigl^2(P,E\otimes\Omega^p_C)$
at $\lambda_0$.

We put as follows:
\[
\begin{array}{l}
 \tilde{\nbigl}^2(P,E\otimes\Omega^0):=
 \bigl\{
 f\in\nbigl^2(P,E\otimes\Omega^0)\,\big|\,
 \square f\in\nbigl^2\bigl(P,E\otimes\Omega^0\bigr)
 \bigr\},\\
\mbox{{}}\\
 \tilde{\nbigl}^2(P,E\otimes\Omega^2):=
 \nbigl^2(P,E\otimes\Omega^2).
\end{array}
\]
We also put as follows:
\[
 \tilde{\nbigl}^2(P,E\otimes\Omega^1):=
 \left\{
 f\in\nbigl^2(P,E\otimes\Omega^1)\,\left|\,
 \begin{array}{ll}
 (\delbar+\theta)f\in \nbigl^2(P,E\otimes\Omega^2),\\
 (\delbar+\theta)^{\ast}f\in \nbigl^2(P,E\otimes\Omega^0)
 \end{array}
 \right.
 \right\}.
\]
\begin{lem}\label{lem;11.4.15}
Let $\lambda$ be any element of $\cnum_{\lambda}$.
We have the following:
\[
 \tilde{\nbigl}^2(P,E\otimes\Omega^1):=
 \left\{
 f\in\nbigl^2(P,E\otimes\Omega^1)\,\left|\,
 \begin{array}{ll}
 \DD^{\lambda}f\in \nbigl^2(P,E\otimes\Omega^2),\\
 \bigl(\DD^{\lambda}\bigr)^{\ast}f\in \nbigl^2(P,E\otimes\Omega^0)
 \end{array}
 \right.
 \right\}.
\]
\end{lem}
\pf
It follows from (\ref{eq;a12.9.15}).
\hfill\qed

\vspace{.1in}

Let us consider the following condition 
for an element
$f\in \nbigc\bigl((\lambda_0,P),\nbige\otimes\Omega^p_C\bigr)$:
\begin{description}
\item[($\star$)]
 $\Phi_f$ is the germ of a holomorphic function
 to $\tilde{\nbigl}^2\bigl(P,E\otimes\Omega^p\bigr)$
 at $P$.
\end{description}

Then we put as follows:
\[
 \tilde{\nbigc}\bigl(
 (\lambda_0,P),\nbige\otimes\Omega^i_C
 \bigr):=
\bigl\{
 f\in  \nbigc\bigl((\lambda_0,P),\nbige\otimes\Omega^p_C\bigr)
 \,\big|\,
 \mbox{\rm $f$ satisfies $(\star)$}
 \bigr\}.
\]
\begin{lem}\label{lem;11.4.16}
The $\lambda$-connection induces the complex:
\[
 \tilde{\nbigc}\bigl((\lambda,P),E\otimes\Omega^0\bigr)
\lrarr
 \tilde{\nbigc}\bigl((\lambda,P),E\otimes\Omega^1\bigr)
\lrarr
 \tilde{\nbigc}\bigl((\lambda,P),E\otimes\Omega^2\bigr).
\]
\end{lem}
\pf
It follows from Lemma \ref{lem;11.4.15}.
\hfill\qed

\vspace{.1in}

Let $\nbigs(E\otimes\Omega^2)$ denote the sheaf
of germs $\tilde{\nbigc}$ on $\cnum_{\lambda}\times \overline{C}$.
Due to Lemma \ref{lem;11.4.16},
the $\lambda$-connection $\DD$
induces the complex
$\nbigs(E\otimes\Omega^{\cdot})$.

\begin{lem}
We have the exact sequence of the complexes:
\[
\begin{CD}
 0 @>>> \nbigs(E\otimes\Omega^{\cdot})
 @>{\lambda-\lambda_0}>>
 \nbigs(E\otimes\Omega^{\cdot})
 @>>>
 \tilde{\nbigl}(\nbigelambda\otimes\Omega^{\cdot})_{(2)}
 @>>> 0.
\end{CD}
\]
\end{lem}
\pf
It follows from Lemma \ref{lem;11.4.17}.
\hfill\qed

\vspace{.1in}

The following lemma can be checked directly.
\begin{lem}
The sheaves $\nbigs(E\otimes\Omega^p)$ $(p=0,1,2)$
are $f$-soft for the projection
$\cnum_{\lambda}\times \overline{C}\lrarr\cnum_{\lambda}$.
(See {\rm\cite{ks}} for the definition of $f$-soft.)

In particular, we have
$R^if_{\ast}\nbigs(E\otimes\Omega^p)=0$
for $i\neq 0$.
\hfill\qed
\end{lem}

Hence $Rf_{\ast}\nbigs(E\otimes\Omega^p_C)$ is canonically
quasi isomorphic to
$f_{\ast}\nbigs(E\otimes\Omega^p_C)$.
The sheaf $f_{\ast}\nbigs(E\otimes\Omega^p)$
is naturally isomorphic to
the sheaf of $\tilde{L}^2(E\otimes\Omega^p_C)$-valued holomorphic functions
over $\cnum_{\lambda}$,
i.e.,
it corresponds to the trivial vector bundle
$\cnum_{\lambda}\times \tilde{L}^2(E\otimes\Omega^p_C)$.

\begin{lem}\label{lem;e12.4.20}
The complex $f_{\ast}\nbigs(E\otimes\Omega^p_C)$
is quasi isomorphic to
$\bigoplus_{i}\gbigh^i[-i]$.
\end{lem}
\pf
We regard
$\bigoplus_{i}\gbigh^i[-i]$
as the complex with the trivial differential.
We have the naturally defined morphism
of $\bigoplus_{i}\gbigh^i[-i]$
to $f_{\ast}\nbigs(E\otimes\Omega^p_C)$.
By using the argument in the subsubsection \ref{subsubsection;9.10.51},
we can show that the morphism is quasi isomorphic.
\hfill\qed

%% file: a86.1.tex

\subsubsection{The sub-complex $\nbigq^{(\lambda_0)}$
 of $\gbige\otimes\Omega^{\cdot,0}$}

Let $\lambda_0$ be a point of $\cnum_{\lambda}$.
Let $P$ be a point of $\bar{C}-C$.
We pick an appropriate coordinate around $P$,
then we can pick an embedding $\Delta\lrarr \overline{C}$
such that the image of $O$ is $P$.
Let us consider the restriction
of $\harmonicbundle$ to $\Delta^{\ast}$.

We take a sufficiently small positive number $\epsilon_0$,
then we have the filtration
$\Vzero\bigl(\naiveprolong{\nbige}\bigr)$.
We have the projections:
\[
 \begin{array}{l}
 \pi_1:
 \Vzero_{\leq -1}\bigl(\naiveprolong{\nbige}\bigr)
 \lrarr
 \Gr^{\Vzero}_{-1}\bigl(\naiveprolong{\nbige}\bigr),\\
 \mbox{{}}\\
 \pi_2:
 \Vzero_{\leq 0}\bigl(
 \naiveprolong{\nbige}\bigr)\otimes\Omega^{1,0}_C
\lrarr
 \Gr^{\Vzero}_0\bigl(\naiveprolong{\nbige}\bigr)
 \otimes\Omega^{1,0}_C.
 \end{array}
\]
We have the generalized decomposition
with respect to the action of $-\deldel_z z$:
\[
\begin{array}{l}
 \Gr^{\Vzero}_{-1}\bigl(\naiveprolong{\nbige}\bigr)
=\bigoplus_{\paramap(\lambda_0,u)=-1}
 \EE\bigl(\Gr^{\Vzero}_{-1}, -\eigenmap(\lambda,u)
 \bigr),\\
 \mbox{{}}\\
\Gr^{\Vzero}_{-0}\bigl(\naiveprolong{\nbige}\bigr)
 \otimes\Omega^{1,0}_C
=\bigoplus_{\paramap(\lambda_0,u)=0}
 \EE\bigl(\Gr^{\Vzero}_0,-\eigenmap(\lambda,u)\bigr).
\end{array}
\]
We put as follows:
\[
 \nbigq^{(\lambda_0)\,0}:=
 \pi_1^{-1}\Bigl(
 W_{1}\EE\bigl(\Gr^{\Vzero}_{-1}\bigl(\naiveprolong{\nbige}\bigr),
 -1\bigr)
 \Bigr),
\quad\quad
 \nbigq^{(\lambda_0),1}:=
 \pi_2^{-1}\Bigl(
 W_{-1}\EE\bigl(\Gr^{\Vzero}_{0}\bigl(\naiveprolong{\nbige}\bigr),
 0\bigr)
\otimes \Omega^{1,0}_C
 \Bigr).
\]

We consider the above procedure for any point of $\bar{C}-C$.
Then we obtain the complex
$\nbigq^{(\lambda_0)\,\bullet}$
on $\bar{C}$.

\begin{lem} \label{lem;f12.4.1}
We have the inclusion
$\nbigq^{(\lambda_0)\,\bullet}\lrarr
 \gbige\otimes\Omega^{\bullet,0}$,
which is quasi isomorphic.
\end{lem}
\pf
The following morphism is isomorphic
due to the strictly specializability of $\gbige$:
\[
 \frac{\gbige}{\Vzero_{-1}(\gbige)}
\lrarr
 \frac{\gbige\otimes\Omega^{1,0}}
 {\Vzero_{0}\gbige\otimes\Omega^{1,0}}.
\]
The following morphism is isomorphic,
which we can show by using Proposition \ref{prop;b11.21.5}:
\[
 \frac{\Gr^{\Vzero}_{-1}\gbige}{
 W_1\EE\bigl(\Gr^{\Vzero}_{-1},-1\bigr)}
\lrarr
  \frac{\Gr^{\Vzero}_{0}\gbige\otimes\Omega^{1,0}}{
 W_{-1}\EE\bigl(\Gr^{\Vzero}_{0},0\bigr)\otimes\Omega^{1,0}}.
\]
Then we obtain the result.
\hfill\qed

\vspace{.1in}
Let us consider the case
$\Delta(\lambda_1,\epsilon_1)\subset
 \Delta(\lambda_0,\epsilon_0)$.
\begin{lem}\label{lem;f12.4.2}
We have the following commutative diagramm:
\[
 \begin{CD}
 \nbigq^{(\lambda_0)}_{|\nbigx(\lambda_1,\epsilon_1)}
 @>>>
 \bigl(
 \gbige\otimes\Omega^{\cdot,0}\bigr)_{|\nbigx(\lambda_1,\epsilon_1)}\\
 @VVV @VVV \\
 \nbigq^{(\lambda_1)}
 @>{=}>> 
\bigl(
 \gbige\otimes\Omega^{\cdot,0}\bigr)_{|\nbigx(\lambda_1,\epsilon_1)}.
 \end{CD}
\]
\end{lem}
\pf
It can checked easily from the definition.
\hfill\qed

\subsubsection{The comparison of
$\nbigq^{(\lambda_0)}$ and $\nbigs(\nbige\otimes\Omega^{\cdot})$}

We have the natural inclusion:
\begin{equation}
 \nbigq^{(\lambda_0)}
\lrarr
 \nbigs(\nbige\otimes\Omega^{\cdot})_{|
 \nbigx(\lambda_0,\epsilon_0)}.
\end{equation}
Hence we obtain the morphism:
\begin{equation}\label{eq;e12.4.1}
Rf_{\ast}\bigl(\nbigq^{(\lambda_0)}
 \bigr)
\lrarr
Rf_{\ast}\bigl(
 \nbigs(\nbige\otimes\Omega^{\cdot})_{|
 \nbigx(\lambda_0,\epsilon_0)}
 \bigr).
\end{equation}
First let us consider the specialization
of (\ref{eq;e12.4.1})
at $\lambda\in\Delta(\lambda_0,\epsilon_0)$.

\begin{lem}\label{lem;e12.4.25}
The following morphism is quasi isomorphic:
\[
 Rf_{\ast}\nbigq^{(\lambda_0)}_{|\nbigxlambda}
\lrarr
 Rf_{\ast}\nbigs(E\otimes\Omega^{\cdot})_{|\nbigxlambda}
=Rf_{\ast}\tilde{\nbigl}^2\bigl(E\otimes\Omega^{\cdot}\bigr).
\]
\end{lem}
\pf
We have the natural inclusion:
\begin{equation}\label{eq;e12.4.5}
 \nbigq^{(\lambda_0)}_{|\nbigxlambda}
\lrarr
 \bigl(
 \nbigelambda\otimes\Omega^{\cdot,0}\bigr)_{(2)}.
\end{equation}
It is easy to check the morphism (\ref{eq;e12.4.5})
is quasi isomorphic.
Then the lemma immediately follows from 
Proposition \ref{prop;e12.4.10}.
\hfill\qed

\vspace{.1in}

We have the following commutative diagram:
\[
 \begin{CD}
 0 @>>>
 \nbigq^{(\lambda_0)\,\cdot}
 @>{\lambda-\lambda_0}>>
 \nbigq^{(\lambda_0)\,\cdot}
 @>>> 
 \nbigq^{(\lambda_0)}_{|\nbigxlambdazero}
 @>>>0 \\
 @VVV @VVV @VVV @VVV @VVV \\
 0 @>>>
 \nbigs(\nbige\otimes\Omega^{\cdot})
 @>{\lambda-\lambda_0}>>
 \nbigs(\nbige\otimes\Omega^{\cdot})
 @>>>
 \tilde{\nbigl}(\nbigelambdazero\otimes\Omega^{\cdot})_{(2)}
 @>>> 0.
 \end{CD}
\]

Let us consider the following induced commutative diagramm:
\begin{equation}\label{eq;a11.9.40}
 \begin{CD}
 \nbigh^i\bigl(
 Rf_{\ast}\bigl(
 \nbigq^{(\lambda_0)\,\cdot}
 \bigr)
 @>>>
 \nbigh^i\bigl(
 Rf_{\ast}\bigl(
 \nbigq^{(\lambda_0)}_{|\nbigxlambdazero}
 \bigr) \\
 @V{\varphi}VV @VVV \\
 \nbigh^i(f_{\ast}\nbigs(\nbige\otimes\Omega^{\cdot}))
 @>>>
 \nbigh^i\bigl(
 \tilde{\nbigl}^2(\nbigelambdazero\otimes\Omega^{\cdot})_{(2)}
 \bigr)
 \end{CD}
\end{equation}

\begin{lem}
$\nbigh^i\bigl(
 Rf_{\ast}\bigl(
 \nbigq^{(\lambda_0)}
 \bigr)$
and $\nbigh^i\bigl( f_{\ast}\nbigs(\nbige\otimes\Omega^{\cdot})\bigr)$
are isomorphic.
\end{lem}
\pf
We will use a descending induction on $i$.
We assume that the isomorphism for $i+1$ is shown,
and we will show that the morphism for $i$ is isomorphic.

Due to Lemma \ref{lem;e12.4.20},
the multiplication of $(\lambda-\lambda_0)$
on
$\nbigh^{i+1}\bigl(Rf_{\ast}\nbigs(\nbige\otimes\Omega^{\cdot})\bigr)$
is injective.
Hence we obtain that the upper vertical arrow
in (\ref{eq;a11.9.40}) is surjective.
Due to Lemma \ref{lem;e12.4.25},
the right vertical arrow is isomorphic.
Then we obtain the following:
\[
 \nbigh^i(f_{\ast}\nbigs(\nbige\otimes\Omega^{\cdot}))
=(\lambda-\lambda_0)\cdot
 \nbigh^i(f_{\ast}\nbigs(\nbige\otimes\Omega^{\cdot}))
+\Image\varphi.
\]
Then we obtain the surjectivity of $\varphi$
due to Nakayama's Lemma.

Let us show the injectivity of $\varphi$.
Let $f$ be a section of $\nbigh^i\bigl(\nbigq^{(\lambda_0)}\bigr)$
such that $\varphi(f)=0$.
Since the right vertical arrow is isomorphic,
$f$ is of the form $(\lambda-\lambda_0)\cdot g$
for some section $g$ of $\nbigh^i\bigl(\nbigq^{(\lambda_0)}\bigr)$.
Since the multiplication of $(\lambda-\lambda_0)$ 
on $\nbigh^i\bigl(Rf_{\ast}\nbigs(E\otimes\Omega^{\cdot})\bigr)$
is injective,
we obtain $\varphi(g)=0$.
Thus we obtain
$\Ker(\varphi)=(\lambda-\lambda_0)\cdot\Ker(\varphi)$.
Then we obtain $\Ker(\varphi)=0$
due to Nakayama's Lemma.
\hfill\qed

\begin{lem}\label{lem;f12.4.3}
Let us consider the case
$\Delta(\lambda_1,\epsilon_1)\subset
 \Delta(\lambda_0,\epsilon_0)$.
We have the following commutative diagramm:
\[
 \begin{CD}
 \nbigh^i\bigl(
 Rf_{\ast}
 \nbigq^{(\lambda_0)}
 \bigr)_{|\Delta(\lambda_1,\epsilon_1)}
 @>>>
 \nbigh^i\bigl(
 Rf_{\ast}
 \nbigs(E\otimes\Omega^{\cdot})
 \bigr)_{|\Delta(\lambda_1,\epsilon_1)}\\
 @VVV @VVV \\
 \nbigh^i\bigl(
 Rf_{\ast}
 \nbigq^{(\lambda_1)}
 \bigr)
 @>>>
 \nbigh^i\bigl(
 Rf_{\ast}
 \nbigs(E\otimes\Omega^{\cdot})
 \bigr)_{|\Delta(\lambda_1,\epsilon_1)}.
 \end{CD}
\]
The horizontal arrows are isomorphic.
\end{lem}
\pf
It is clear from the definition.
\hfill\qed

\begin{lem}
We have the canonical isomorphisms:
\[
 \nbigh^i\bigl(
 Rf_{\ast}\bigl(\gbige\otimes\Omega^{\cdot,0}\bigr)
 \bigr)
\simeq
 \nbigh^i\bigl(
 f_{\ast}\nbigs(E\otimes\Omega^{\cdot})
 \bigr)
\simeq
 \gbigh^i\otimes\nbigo_{\cnum_{\lambda}}.
\]
\end{lem}
\pf
It follows from 
Lemma \ref{lem;f12.4.1},
Lemma \ref{lem;f12.4.2},
Lemma \ref{lem;e12.4.25}
and Lemma \ref{lem;f12.4.3}.
\hfill\qed

\begin{cor}
The induced $\nbigr$-triple
$\bigl(
 R^if_{\ast}\bigl(\gbige\otimes\Omega^{\cdot}\bigr),
 R^if_{\ast}\bigl(\gbige\otimes\Omega^{\cdot}\bigr),
 C
 \bigr)$
is pure twistor of weight $i$.
The polarization is naturally given.
The hard Lefschetz theorem holds.
\end{cor}
\pf
The twistor property and the positivity
can be shown by an argument similar 
to the proof of Theorem 2.24 in the paper of Sabbah \cite{sabbah}.
Since we can take the harmonic repesentatives of the cohomology classes,
the hard Lefschetz theorem clearly holds.
\hfill\qed